\documentclass[english,twopage]{book}
\usepackage{mathptmx}

\usepackage[T1]{fontenc}
\usepackage[latin9]{luainputenc}
\usepackage{fancyhdr}
\usepackage{mathrsfs}
\usepackage{amsmath}
\usepackage{amsthm}
\usepackage{amssymb}
\usepackage{makeidx}
\makeindex

\makeatletter

\newcommand{\noun}[1]{\textsc{#1}}

\numberwithin{equation}{section}
\numberwithin{figure}{section}
\theoremstyle{plain}
\newtheorem{thm}{\protect\theoremname}[section]
  \theoremstyle{definition}
  \newtheorem{defn}[thm]{\protect\definitionname}
  \theoremstyle{plain}
  \newtheorem{cor}[thm]{\protect\corollaryname}
  \theoremstyle{plain}
  \newtheorem{lem}[thm]{\protect\lemmaname}
  \theoremstyle{plain}
  \newtheorem{prop}[thm]{\protect\propositionname}
  \theoremstyle{definition}
  \newtheorem{xca}[thm]{\protect\exercisename}
  \theoremstyle{remark}
  \newtheorem{rem}[thm]{\protect\remarkname}

\let\OrgIndex\index

\renewcommand*{\index}[1]{\OrgIndex{#1}}

\makeatother

\usepackage{babel}
  \providecommand{\corollaryname}{Corollary}
  \providecommand{\definitionname}{Definition}
  \providecommand{\exercisename}{Exercise}
  \providecommand{\lemmaname}{Lemma}
  \providecommand{\propositionname}{Proposition}
  \providecommand{\remarkname}{Remark}
\providecommand{\theoremname}{Theorem}

\begin{document}

\title{\textbf{\huge{}Foundations of Constructive Probability Theory}}

\author{Yuen-Kwok Chan, \thanks{Mortgage Analytics, Citigroup, (Retired); all opinions expressed by
the author are his own.} \thanks{The author is grateful to the late Prof E.Bishop for teaching him
Constructive Mathematics, to the late Profs R. Getoor and R. Blumenthal
for teaching him Probability and for mentoring, to the late Profs
R. Pyke and W. Birnbaum and the other statisticians in the Mathematics
Department of Univ of Washington, circa 1970's, for their moral support.
The author is also thankful to the constructivists in the Mathematics
Department of New Mexico State University, circa 1975, for hosting
a sabbatical visit and for valuable discussions, especially to Profs
F. Richman, D. Bridges, M. Mandelkern, W. Julian, and the late Prof.
R.Mines.} \thanks{Contact: chan314@gmail.com}}

\date{June 2019}

\maketitle
\tableofcontents{}

\part{Introduction and Preliminaries}

\chapter{Introduction}

The present work on probability theory is an outgrowth of the constructive
analysis in \cite{Bishop67} and \cite{BishopBridges85}. 

Perhaps the simplest explanation of constructive mathematics is by
way of focusing on the following two commonly used theorems. The first,
the \emph{\index{principle of finite search}principle of finite search},
states that, given a finite sequence of 0-or-1 integers, either all
members of the sequence are equal to $0$, or there exists a member
which is equal to $1$. We use this theorem without hesitation because,
given the finite sequence, a finite search would determine the result.

The second theorem, which we may call the \emph{\index{principle of infinite search}principle
of infinite search}, states that, given an infinite sequence of 0-or-1
integers, either all members of the sequence are equal to $0$, or
there exists a member which is equal to $1$. The name ``infinite
search'' is perhaps unfair, but it brings into sharp focus that the
computational meaning of this theorem is not clear. The theorem is
tantamount to an infinite loop in computer programming.

Most mathematicians acknowledge the important distinction between
the two theorems, but regard the principle of infinite search as an
expedient tool to prove theorems, with the belief that theorems so
proved can then be specialized to constructive theorems, when necessary. 

Contrary to this belief, many classical theorems proved directly or
indirectly via the principle of infinite search are actually equivalent
to the latter, and, as such, can never have a constructive proof.
Oftentimes, not even the numerical meaning of the theorems in question
is clear. 

We believe that, for the constructive formulations and proofs of even
the most abstract theorems, the easiest way is a disciplined and systematic
approach, by using only finite searches and by quantifying mathematical
objects and theorems at each and every step, with natural numbers
as a starting point. The above cited references show that this approach
is not only possible, but fruitful.

It should be emphasized that we do not claim that theorems whose proofs
require the principle of infinite search are untrue or incorrect.
They are certainly correct and consistent derivations from commonly
accepted axioms. There is indeed no reason why we cannot discuss such
classical theorems alongside their constructive counterparts. The
term ``non-constructive mathematics'' is not meant to be pejorative.
We will use, in its place, the more positive term ``classical mathematics''. 

Moreover, it is a myth that constructivists use a different system
of logic. The only logic we use is everyday logic; no formal language
is needed. The present author considers himself a mathematician who
is neither interested in, nor equipped to comment on, the formalization
of mathematics, classical or constructive. 

Since a constructively valid argument is also correct from the classical
view point, a reader of the classical persuasion should have no difficulties
understanding our proofs. Proofs using only finite searches are surely
agreeable to any reader who is accustomed to infinite searches. 

Indeed, the author would consider the present book a success if the
reader, but for this introduction and occasional remarks in the text,
finishes reading without realizing that this is a constructive treatment.
At the same time, we hope that a reader of the classical persuasion
might consider the more disciplined approach of constructive mathematics
for his or her own research an invitation to a challenge. 

We hasten to add that we do not think that finite computations in
constructive mathematics are the end. We would prefer a finite computation
with $n$ steps to one with $n!$ steps. We would be happy to see
a systematic and general development of mathematics which is not only
constructive, but also computationally efficient. That admirable goal
will however be left to abler hands. 

Probability theory, rooted in applications, can naturally be expected
to be constructive. Indeed, the crowning achievements of probability
theory \textemdash{} the laws of large numbers, the central limit
theorems, the analysis of Brownian Motion processes and their stochastic
integrals, and that of Levy processes, to name just a few \textemdash{}
are exemplars of constructive mathematics. Kolmogorov, the grandfather
of modern probability theory, actually took an interest in the formalization
of general constructive mathematics.

On the other hand, many a theorem in modern probability actually implies
the principle of infinite search. The present work attempts a systematic
constructive development. Each existence theorem will be a construction.
The input data, the construction procedure, and the output objects
are the essence and integral parts of the theorem. Incidentally, by
inspecting each step in the procedure, we can routinely observe how
the output varies with the input. Thus a continuity theorem in epsilon-delta
terms routinely follows an existence theorem. For example, we will
construct a Markov process from a given semigroup, and prove that
the resulting Markov process varies continuously with the semigroup,
in epsilon-delta terms often derived from the Borell-Cantelli lemma. 

The reader with the probability literature will notice that our constructions
resemble Kolmogorov's construction of the Brownian motion process,
which is replete with Borel-Cantelli estimates and rates of convergence.
This is in contrast to popular proofs of existence via Prokhorov's
Theorem. The reader can regard Part III of this book, Chapters 6-11,
the part on stochastic processes, as an extension of Kolmogorov's
constructive methods to stochastic processes: Danielle-Kolmogorov-Skorokhod
construction of random fields, measurable random fields, a.u. continuous
processes, a.u. càdlàg processes, martingales, strong Markov processes,
and Feller processes, all with locally compact state spaces. 

Such a systematic, constructive, and general treatment of stochastic
processes, we believe, has not previously been attempted.

The purpose of this book is twofold. A student with a general mathematics
background can use it at the first-year graduate-school level can
use it as an introduction to probability or to constructive mathematics,
and an expert in probability can use it as a reference for further
constructive development in his or her own research specialties. 

Part II of this book, Chapters 3-5, is a re-packaging and expansion
 of the measure theory in \cite{BishopBridges85}. This is so we can
have a self-contained probability theory in terms familiar to probabilists. 

For expositions of constructive mathematics, see the first chapters
of the last cited reference. See also \cite{Richman82} and \cite{Stolzenberg70}.
We give a synopsis in the next chapter, along with basic notations
and terminologies.

\chapter{Preliminaries}

\section*{Natural numbers}

We start with the natural numbers as known in elementary schools.
All mathematical objects are constructed from natural numbers, every
theorem ultimately a calculation on the natural numbers. From natural
numbers are constructed the integers and the rational numbers, along
with the arithmetical operations, in the manner taught in elementary
schools. 

We claim to have a natural number only when we have provided a finite
method to calculate it, i.e. to find its decimal representation. This
is the fundamental difference from classical mathematics, which requires
no such finite method; an infinite procedure in a proof is considered
just as good in classical mathematics.

The notion of a finite natural number is so simple and so immediate
that no attempt is needed to define them in even simpler terms. A
few examples would suffice as clarification: $1,2,3$ are natural
numbers. So are $9^{9}$ and $9^{9^{9}}$; the multiplication method
will give, at least in principle, their decimal expansion in a finite
number of steps. On the other hand, the ``truth value'' of a particular
mathematical statement is a natural number only if a finite method
has been supplied which, when carried out, would conclusively prove
or disprove the statement.

\section*{Calculations and theorems}

An algorithm or a calculation means any finite, step-by-step procedure.
A mathematical object is defined when we specify the calculations
that need to be done to produce this object. We say that we have proved
a theorem if we have provided a step-by-step method that translates
the calculations doable in the hypothesis to a calculation in the
conclusion of the theorem. The statement of the theorem is merely
a summary of the algorithm contained in the proof.

Although we do not, for good reasons, write mathematical proofs in
a computer language, the reader would do well to compare constructive
mathematics to the development of a large computer software library,
successive objects and library functions being built from previous
ones, each with a guarantee to finish in a finite number of steps.

\section*{Proofs by contradiction}

There is a trivial form of proofs by contradiction which is valid
and useful in constructive mathematics. Suppose we have already proved
that one of two given alternatives, $A$ and $B$, must hold, meaning
that we have given a finite method, which, when unfolded, gives either
a proof for $A$ or a proof for $B$. Suppose subsequently we also
prove that $A$ is impossible. Then we can conclude that we have a
proof of $B$; we need only exercise said finite method, and see that
the resulting proof is for $B$. 

\section*{Recognizing non-constructive theorems}

Consider the simple theorem ``if $a$ is a real number, then $a\leq0$
or $0<a$'', which may be called the principle of excluded middle
for real numbers. We can see that this theorem implies the principle
of infinite search by the following argument. Let $(x)_{i=1,2,...}$
be any given sequence of 0-or-1 integers. Define the real number $a=\sum_{i=1}^{\infty}x_{i}2^{-i}$.
If $a\leq0$, then all members of the given sequence are equal to
0; if $0<a,$ then some member is equal to 1. Thus the theorem implies
the principle of infinite search, and therefore cannot be constructively
valid. 

Any theorem that implies this limited principle of excluded middle
cannot have a constructive proof. This observation provides a quick
way for the constructive analyst to recognize certain theorems as
nonconstructive. Then we can proceed to find constructive substitutes.

For the aforementioned principle of excluded middle of real numbers
itself, a useful constructive substitute is the theorem ``if $a$
is a real number, then, for arbitrarily small $\varepsilon>0$, we
have $a<\varepsilon$ or $0<a$''. Heuristically, this is a recognition
that a general real number $a$ can be computed with arbitrarily small,
but nonzero, error.

\section*{Prior knowledge}

We assume that the reader of this book has familiarity of calculus
and metric spaces, and has had an introductory course in probability
theory at the level of \cite[Feller]{FellerI71} or \cite[Ross]{Ross03}.
We recommend prior reading of the first four chapters of \cite{BishopBridges85},
which contain the basic treatment of the real numbers, set theory,
and metric spaces. We will also require some rudimentary knowledge
of complex numbers and complex analysis.

The reader should have no difficulty in switching back and forth between
constructive mathematics and classical mathematics, any more than
in switching back and forth between classical mathematics and computer
programming. Indeed, the reader is urged to read, concurrently with
this book if not before, the many classical texts in probability.

\section*{Notations and conventions}

If $x,y$ are mathematical objects, we write $x\equiv y$ to mean
``$x$ is defined as $y$'', ``$x$, which is defined as $y$'',
``$x,$ which has been defined earlier as $y$'', or any other grammatical
variation depending on the context. 

\subsection*{Numbers}

Unless otherwise indicated, $N,Q$, and $R$ will denote the set of
integers, the set of rational numbers in the decimal or binary system,
and the set of real numbers respectively. We will also write $\{1,2,\cdots\}$
for the set of positive integers. The set $R$ is equipped with the
Euclidean metric. Suppose $a,b,a_{i}\in R$ for $i=m,m+1\cdots$ for
some $m\in N$. We will write $\lim_{i\rightarrow\infty}a_{i}$ for
the limit of the sequence $a_{m},a_{m+1},\cdots$ if it exists, without
explicitly referring to $m$. We will write $a\vee b,a\wedge b,a_{+},a_{-}$
for $\max(a,b),\min(a,b),a\vee0,a\wedge0$ respectively. The sum $\sum_{i=m}^{n}a_{i}\equiv a_{m}+\cdots+a_{n}$
is understood to be $0$ if $n<m$. The product $\prod_{i=m}^{n}a_{i}\equiv a_{m}\cdots a_{n}$
is understood to be $1$ if $n<m$. Suppose $a_{i}\geq0$ for $i=m,m+1\cdots$
. We write $\sum_{i=m}^{\infty}a_{i}<\infty$ if and only if $\sum_{i=m}^{\infty}|a_{i}|<\infty$,
in which case $\sum_{i=m}^{\infty}a_{i}$ is taken to be $\lim_{n\rightarrow\infty}\sum_{i=m}^{n}a_{i}$
. In other words, unless otherwise specified, convergence of a series
of real numbers means absolute convergence.

\subsection*{Sets and functions}

In general, a set is a collection of objects equipped with an equality
relation. To define a set is to specify how to construct an element
of the set, and how to prove that two elements are equal. A set is
also called a family.

The usual set-theoretic notations are used. Let two subsets $A$ and
$B$ of a set $\Omega$ be given. We will write $A\cup B$ for the
union, and $A\cap B$ or $AB$ for the intersection. We write $A\subset B$
if each member $\omega$ of $A$ is a member of $B$. We write $A\supset B$
for $B\subset A$, The \index{set-theoretic complement}\emph{set-theoretic
complement} of a subset $A$ of the set $\Omega$ is defined as the
set $\{\omega\in\Omega:\omega\in A\mbox{ implies a contradiction}\}$.
We write $\omega\notin A$ if $\omega\in A$ implies a contradiction.
The set $\Omega$ is said to be \emph{discrete}\index{discrete set}
if, for each $\omega,\omega'\in\Omega$ either (i) $\omega=\omega'$
or (ii) the assumption that $\omega=\omega'$ leads to a contradiction.
Given a set $\Omega$, the subset $\phi$ defined as the set-theoretic
complement of $\Omega$ and is called the empty set. The notion of
a set-theoretic complement is otherwise rarely used in the present
book. Instead, we will use heavily the notions of a metric complement
or a measure-theoretic complement, to be defined later. 

Suppose $A,B$ are sets. A finite, step by step, method $X$ which
produces an element $X(x)\in B$ given any $x\in A$ is called an
\index{operation}\emph{operation} from $A$ to $B$. The element
$X(x)$ need not be unique. Two different applications of the operation
$X$ with the same input element $x$ can produce different outputs.
An example of an operation is $[\cdot]_{1}$, which assigns to each
$a\in R$ an integer $[a]_{1}\in(a,a+2)$. This operation is a substitute
of the classical operation $[\cdot],$ and will be used frequently
in the present work.

Suppose $\Omega,\Omega'$ are sets. Suppose $X$ is an operation which,
for each $\omega$ in some non-empty subset $A$ of $\Omega$, constructs
a unique member $X(\omega)$ in $\Omega'$. Then the operation $X$
is called a \emph{function}\index{function} from $\Omega$ to $\Omega'$,
or simply a function on $\Omega$. The subset $A$ is called the \emph{domain}\index{domain of a function}
of $X$. A function is also called a \emph{mapping}\index{}. We then
write $X:\Omega\rightarrow\Omega'$, and write $domain(X)$ for the
set $A$. Thus a function $X$ is an operation which has the additional
property that if $\omega_{1}=\omega_{2}$ in $domain(X)$, then $X(\omega_{1})=X(\omega_{2})$
in $\Omega'$. The non-empty $domain(X)$ is not required to be the
whole set $\Omega$. This  is convenient when we work with functions
defined only almost everywhere, in a sense to be made precise later.
Separately, we sometimes use the expression $\omega\rightarrow X(\omega)$
for the a function $X$ whose domain is understood. For example, the
expression $\omega\rightarrow\omega^{2}$ stands for the function
$X:R\rightarrow R$ defined by $X(\omega)\equiv\omega^{2}$ for each
$\omega\in R$. 

Let $X:\Omega\rightarrow\Omega'$ be a function, and let $A$ be a
subset of $\Omega$ such that $A\cap domain(X)$ is non-empty. Then
the \index{restriction of a function}\emph{restriction} $X|A$ of
$X$ to $A$ is defined as the function from $A$ to $\Omega'$ with
$domain(X|A)\equiv A\cap domain(X)$ and $(X|A)(\omega)$ for each
$\omega\in domain(X|A)$. The set
\[
B\equiv\{\omega'\in\Omega':\omega'=X(\omega)\:for\:some\:\omega\in domain(X)\}
\]
is called the \index{range of a function}\emph{range} of the function
$X$, and is denoted by $range(X)$. 

A function $X:A\rightarrow B$ is called a \emph{\index{surjection}surjection
}if $range(X)=B$; in that case, there exists an operation $Y:B\rightarrow A$,
not necessarily a function, such that $X(Y(b))=b$ for each $b\in B$.
The function $X$  is called an \index{injection}\emph{injection}
if for each $a,a'\in domain(X)$ with $X(a)=X(a')$ we have $a=a'$.
It is called a bijection \index{bijection} if $domain(X)=A$ and
if $X$ is both a surjection and an injection. 

Let $X:B\rightarrow A$ be a surjection with $domain(X)=B$. Then
the triple $(A,B,X)$ is called an \emph{indexed set}\index{indexed set}\index{indexed family}.
In that case, we write $X_{b}\equiv X(b)$ for each $b\in B$. We
will, by abuse of notations, call $A$ or $\{X_{b}:b\in B\}$ an \emph{indexed
set}, and write $A\equiv\{X_{b}:b\in B\}$. We will call $B$ the
index set, and say that $A$ is indexed by the members $b$ of $B$.

A set $A$ is said to be \emph{finite}\index{finite set} if there
exists a bijection $v:\{1,\cdots,n\}\rightarrow A$, for some $n\geq1$,
in which case we write $|A|\equiv n$ and call it the \index{size of a finite set}
\emph{size }of $A$. We will then call $v$ an \index{enumeration}\emph{enumeration}
of the set $A$, and call the pair $(A,v)$ an \emph{\index{enumerated set}enumerated
set}. When the enumeration $v$ is understood from context, we will
abuse notations and simply call the set $A\equiv\{v_{1},\cdots,v_{n}\}$
an enumerated set. 

A set $A$ is said to be \emph{countable }\index{countable set} if
there exists a surjection $v:\{1,2,\cdots\}\rightarrow A$. A set
$A$ is said to be \emph{countably infinite}\index{countably infinite set}
if there exists a bijection $v:\{1,2,\cdots\}\rightarrow A$. We will
then call $v$ an \emph{enumeration}\index{enumeration} of the set
$A$, and call the pair $(A,v)$ an \emph{\index{enumerated set}}enumerated
set. When the enumeration $v$ is understood from context, we will
abuse notations and simply call the set $A\equiv\{v_{1},v_{2},\cdots\}$
an enumerated set. 

Suppose $X:\Omega\rightarrow\Omega'$ and $X':\Omega'\rightarrow\Omega''$
are such that the set $A$ defined by $A=\{\omega\in domain(X):X(\omega)\in domain(X')\}$
is non-empty. Then the \index{composite function}\emph{composite
function} $X'\circ X:\Omega\rightarrow\Omega''$ is defined to have
$domain(X'\circ X)=A$ and $(X'\circ X)(\omega)=X'(X(\omega))$ for
$\omega'\in A$. The alternative notations $X'(X)$ will also be used
for $X'\circ X$.

Henceforth, unless otherwise indicated, we write $X(\omega)$ only
with the implicit condition that $\omega\in domain(X)$. 

Two functions $X,Y$ are considered equal, $X=Y$ in symbols, if 
\[
domain(X)=domain(Y)
\]
and $X(\omega)=Y(\omega)$ for each $\omega\in domain(X)$. When emphasis
is needed, this equality will be referred to as the \emph{\index{set-theoretic equality of functions}set-theoretic
equality}, in contradistinction to almost everywhere equality, to
be defined later. 

Let $\Omega$ be a set and let $n\geq1$ be arbitrary integer. A function
$\omega:\{1,\cdots,n\}\rightarrow\Omega$ which assigns to each $i\in\{1,\cdots,n\}$
an element $\omega(i)\equiv\omega_{i}\in\Omega$ is called \emph{a
finite sequence}\index{sequence} of elements in $\Omega$. A function
$\omega:\{1,2,\cdots,\}\rightarrow\Omega$ which assigns to each $i\in\{1,2,\cdots\}$
an element $\omega(i)\equiv\omega_{i}\in\Omega$ is called \emph{an
infinite sequence}\index{sequence} of elements in $\Omega$. We will
then write $\omega\equiv(\omega_{1},\cdots,\omega_{n})\equiv$ or
$(\omega_{i})_{i=1,\cdots,n}$, in the first case, and write $(\omega_{1},\omega_{2},\cdots)$
or $(\omega_{i})_{i=1,2,\cdots,}$, in the second case, for the sequence
$\omega$. If, in addition, $j$ is a sequence of integers in $domain(\omega)$,
with that $j_{k}<j_{h}$ for each $k<h$ in $domain(j)$, then the
sequence $\omega\circ j:domain(j)\rightarrow\Omega$ is called a \emph{subsequence}\index{subsequence}
of $\omega$. Throughout this book, we will write a subscripted symbol
$a_{b}$ interchangeably with $a(b)$ to lessen the burden on subscripts.
Thus, $a_{b(c)}$ stands for of $a_{b_{c}}$. Similarly, $\omega_{j_{k}}\equiv\omega_{j(k)}\equiv\omega(j(k))$
for each $k\in domain(j)$, and we write $(\omega_{j(1)},\omega_{j(2)},\cdots)$
or $(\omega_{j(k)})_{k=1,2,\cdots}$, or simply $(\omega_{j(k)}),$
for the subsequence when the domain of $j$ is clear. If $(\omega_{1},\cdots,\omega_{n})$
is a sequence, we will write $\{\omega_{1},\cdots,\omega_{n}\}$ for
the range of $\omega$. Thus an element $\omega_{0}\in\Omega$ is
in $\{\omega_{1},\cdots,\omega_{n}\}$ if and only if there exists
$i=1,\cdots,n$ such that $\omega_{0}=\omega_{i}$.

Suppose $(\omega_{i})_{i=1,2,\cdots,}$ and $(\omega'_{i})_{i=1,2,\cdots,}$
are two infinite sequences. We will write $(\omega_{i},\omega'_{i})_{i=1,2,\cdots}$
for the merged sequence $(\omega_{1},\omega'_{1},\omega_{2},\omega'_{2},\cdots)$.
Similar notations for several sequences. 

Let $(\Omega_{n})_{n=0,1,\cdots}$ be a sequence of non-empty sets.
Consider any $0\leq n\leq\infty$, i.e $n$ is a non-negative integer
or the symbol $\infty$. We will let $\Omega^{(n)}$ denote the Cartesian
product $\prod_{j=0}^{n}\Omega_{j}$. Consider $0\leq k<\infty$ with
$k\leq n$. The \emph{\index{coordinate function}coordinate function}
$\pi_{k}$ is the function with $domain(\pi_{k})=\Omega^{(n)}$ and
$\pi_{k}(\omega_{0},\omega_{1},\cdots)=\omega_{k}$. If $\Omega_{n}=\Omega$
for each $n\geq0$, then we will write $\Omega^{n}$ for $\Omega^{(n)}$
for each $n\geq0$. Let $X$ be a function on $\Omega_{k}$ and let
$Y$ be a function on $\Omega^{(k)}$. When confusion is unlikely,
we will use the same symbol $X$ also for the function $X\circ\pi_{k}$
on $\Omega^{(n)}$, which depends only on the $k$-th coordinate.
Likewise we will use $Y$ also for the function $Y\circ(\pi_{0},\cdots,\pi_{k})$
on $\Omega^{(n)}$, which depends only on the first $k+1$ coordinates.
Thus every function on $\Omega_{k}$ or $\Omega^{(k)}$ is identified
with a function on $\Omega^{(\infty)}$. Accordingly, sets of functions
on $\Omega_{k},\Omega^{(k)}$ are regarded also as sets of functions
on $\Omega^{(n)}$. 

Let $M$ the family of all real-valued functions on $\Omega$, equipped
with the set-theoretic equality for functions. Suppose $X,Y\epsilon M$
and suppose $f$ is a function on $R\times R$ such that the set 
\[
D\equiv\{\omega\in domain(X)\cap domain(Y):(X(\omega)Y(\omega))\in domain(f)\}
\]
is non empty. Then $f(X,Y)$ is defined as the function with $domain(f(X,Y))\equiv D$
and $f(X,Y)(\omega)\equiv f(X(\omega),Y(\omega))$ for each $\omega\in D$.
The definition extends to a sequence of functions in the obvious manner.
Suppose $(X_{i})_{i=m,m+1,\cdots}$ is a sequence in $F$. Suppose
the set 
\[
D\equiv\{\omega\in\cap_{i=m}^{\infty}domain(X_{i}):\sum_{i=m}^{\infty}|X_{i}(\omega)|<\infty\}
\]
is non-empty, then $\sum_{i=m}^{\infty}X_{i}$ is defined as the function
with $domain(\sum_{i=m}^{\infty}X_{i})\equiv D$ and with value $\sum_{i=m}^{\infty}X_{i}(\omega)$
for each $\omega\in D$. Unless otherwise specified, convergence for
series means absolute convergence. 

Suppose $X,Y\epsilon M$ and $A$ is a subset of $\Omega$, and suppose
$a\in R$. We say $X\leq Y$ on $A$ if (i) $A\cap domain(X)=A\cap domain(Y)$
and (ii) $X(\omega)\leq Y(\omega)$ for each $\omega\in A\cap domain(X)$.
If $X\leq Y$ on $\Omega$ we will simply write $X\leq Y$. Thus $X\leq Y$
implies $domain(X)=domain(Y)$. We write $X\leq a$ if $X(\omega)\leq a$
for each $\omega\in domain(X)$. We will write 
\[
(X\leq a)\equiv\{\omega\in domain(X):X(\omega)\leq a\}.
\]
We make similar definitions when the relation $\leq$ is replaced
by $<,\geq,>,$ or $=$. We say $X$ is non-negative if $X\geq0$. 

Suppose $a\in R$. We will abuse notations and write $a$ also for
the constant function $X$ with $domain(X)=\Omega$ and with $X(\omega)=a$
for each $\omega\in domain(X)$. 

Let $X$ be a function on the product set $\Omega'\times\Omega''$.
Let $\omega'\in\Omega'$ be such that $(\omega',\omega'')\in domain(X)$
for some $\omega''\in\Omega''$. Define the function $X(\omega',\cdot)$
on $\Omega''$ by
\[
domain(X(\omega',\cdot))\equiv\{\omega''\in\Omega'':(\omega',\omega'')\in domain(X)\}
\]
Form$X(\omega',\cdot)(\omega'')\equiv X(\omega',\omega'')$. Similarly,
let $\omega''\in\Omega''$ be  such that $(\omega',\omega'')\in domain(X)$
for some $\omega'\in\Omega'$. Define the function $X(\cdot,\omega'')$
on $\Omega'$ by 
\[
domain(X(\cdot,\omega''))\equiv\{\omega'\in\Omega':(\omega',\omega'')\in domain(X)\}.
\]
and $X(\cdot,\omega'')(\omega')\equiv X(\omega',\omega'')$. Given
a function $X$ on the Cartesian product $\Omega'\times\Omega''\times\cdots\times\Omega^{(n)}$,
for each $(\omega',\omega'',\cdots,\omega^{(n)})\in domain(X)$, we
define similarly the functions $X(\cdot,\omega'',\omega''',\cdots,\omega^{(n)})$,
$X(\omega',\cdot,\omega''',\cdots,\omega^{(n)}),\cdots$,$X(\omega',\omega'',\cdots,\omega^{(n-1)},\cdot)$
on the sets $\Omega',\Omega'',\cdots,\Omega^{(n)}$ respectively.

Let $M',M''$ denote the families of all real-valued functions on
two sets $\Omega',\Omega''$ respectively, and let $L''$ be a subset
of $M'$. Suppose
\begin{equation}
T:\Omega'\times L''\rightarrow R\label{eq:T:Omega' X L'' into R}
\end{equation}
is a real-valued function. We can define a function 
\[
T^{*}:L''\rightarrow M'
\]
with 
\[
domain(T^{*})\equiv\{X''\in L'':domain(T(\cdot,X''))\mbox{ is non-empty}\}
\]
and by $T^{*}(X'')\equiv T(\cdot,X'')$. When there is no risk of
confusion, we write $T$ also for the function $T^{*}$, $TX''$ for
$T(\cdot,X'')$, and write
\[
T:L''\rightarrow M'
\]
interchangeably with the expression (\ref{eq:T:Omega' X L'' into R}).
Thus the duality
\begin{equation}
T(\cdot,X'')(\omega')\equiv T(\omega',X'')\equiv T(\omega',\cdot)(X'').\label{eq:duality of T}
\end{equation}

\subsection*{Metric spaces}

We recommend prior reading of the first four chapters of \cite{BishopBridges85},
which contain the basic treatment of the real numbers, set theory,
and metric spaces. We will use without comment theorems about metric
spaces and continuous functions from these chapters. The definitions
and notations, with few exceptions, are familiar to readers of classical
texts. A summary of these definitions follows.

Let $(S,d)$ be a metric space. If $J$ is a subset of $S$, its \index{metric complement}\emph{metric
complement} is the set $\{x\in S:d(x,y)>0\mbox{ for all }y\in J\}$,
Unless otherwise specified, $J_{c}$ will denote the metric complement
of $J$. A condition is said to hold for \emph{all but countably many}\index{all but countably many}
members of $S$ if it holds for each member in the metric complement
$J_{c}$ of some countable subset $J$ of $S.$ We will say that two
elements $x,y\in S$ are \emph{unequal}\index{unequal elements of metric space},
and write $x\neq y$, if $d(x,y)>0$.

We will call a subset $A$ of $S$ \index{metrically discrete subset}\emph{metrically
discrete} if, for each $x,y\in A$ we have $x=y$ or $d(x,y)>0$.
Classically each subset $A$ of $S$ is metrically discrete. 

Let $(f_{n})_{n=1,2,\cdots}$ be a sequence of functions from a set
$\Omega$ to $S$ such that the set
\[
D\equiv\{\omega\in\cup_{n=1}^{\infty}\cap_{i=n}^{\infty}domain(f_{i}):\lim_{i\rightarrow\infty}f_{i}(\omega)\:exists\;in\;S\}
\]
is non-empty, then $\lim_{i\rightarrow\infty}f_{i}$ is defined as
the function with $domain(\lim_{i\rightarrow\infty}$$f_{i})\equiv D$
and with value $\lim_{i\rightarrow\infty}f_{i}(\omega)$ for each
$\omega\in D$. We emphasize that $\lim_{i\rightarrow\infty}f_{i}$
is well defined only if it can be shown that $D$ is non-empty. Note
that for each $\omega\in D$, the value $f_{i}(\omega)$ is defined
in $S$ for each $i\geq n$ for some $n\geq1$, but not necessary
for any $i<n$. 

A function $f:S\rightarrow S'$ is said to be uniformly continuous
on a subset $A\subset domain(f)$, relative to the metrics $d,d'$
on $S,S'$ respectively, if there exists an operation $\delta:(0,\infty)\rightarrow(0,\infty)$
such that $d'(f(x),f(y))<\varepsilon$ for each $x,y\in A$ with $d(x,y)<\delta(\varepsilon)$.
When there is need to be precise as to the metrics $d,d'$, we will
say that $f:(S,d)\rightarrow(S',d')$ is uniformly continuous on $A$.
The operation $\delta$ is called a \emph{modulus of continuity}\index{modulus of continuity}
of $f$ on $A$. If there exists a coefficient $c\geq0$ such that
$d'(f(x),f(y))\leq cd(x,y)$ for all $x,y\in A$, then the function
$f$ is said to be \emph{Lipschitz continuous} \index{Lipschitz continuity}
on $A$, and the constant $c$ is then called a Lipschitz constant
of $f$ on . In that case, we will say simply that  $f$ has Lipschitz
constant $c$.

A metric space $(S,d)$ is said to be \emph{totally bounded} if, for
each $\varepsilon>0$, there exists a finite subset $A\subset S$
such that for each $x\in S$ there exists $y\in A$ with $d(x,y)<\varepsilon$.
The subset $A$ is then called an \index{epsilon-approximation of a totally bounded metric space}
\emph{$\varepsilon$-approximation} of $S$. A \emph{compact metric
space}\index{compact metric space} $K$ is defined as a complete
and totally bounded metric space. 

A subset $A\subset S$ is said to be \index{bounded set}\emph{bounded}
if there exists $x\in S$ and $a>0$ such that $A\subset(d(\cdot,x)\leq a)$.
A subset $S'\subset S$ is said to be \emph{locally compact\index{locally compact metric space}}
if every bounded subset of $S'$ is contained in some compact subset.
The metric space $(S,d)$ is said to be locally compact if the subset
$S$ is locally compact. A function $f:(S,d)\rightarrow(S',d')$ is
said to be \index{continuous function on locally compact space}continuous
if $domain(f)=S$ and if it is uniformly continuous on each compact
subset $K$ of $S$. 

Suppose $(S_{n},d_{n})_{n=1,2,\cdots}$ is a sequence of metric spaces.
For each integer $n\geq1$, define 
\[
d^{(n)}(x,y)\equiv(\bigotimes_{i=1}^{n}d_{i})(x,y)\equiv(d_{1}\otimes\cdots\otimes d_{n})(x,y)\equiv\bigvee_{i=1}^{n}d_{i}(x_{i},y_{i})
\]
for each $x,y\in\prod_{i=1}^{n}S_{i}$. Then $(S^{(n)},d^{(n)})\equiv\bigotimes_{i=1}^{n}(S_{i},d_{i})\equiv(\prod_{i=1}^{n}S_{i},\bigotimes_{i=1}^{n}d_{i})$
is a metric space called the \index{product metric space}\emph{product
metric space} of $S_{1},\cdots,S_{n}$. Define the infinite product
metric $\bigotimes_{i=1}^{\infty}d_{i}$ on $\prod_{i=1}^{\infty}S_{i}$
by 
\[
d^{(\infty)}(x,y)\equiv(\bigotimes_{i=1}^{\infty}d_{i})(x,y)\equiv\sum_{i=1}^{\infty}2^{-i}(1\wedge d_{i}(x_{i},y_{i}))
\]
for each $x,y\in\prod_{i=1}^{\infty}S_{i}$. Define the infinite product
metric space 
\[
(S^{(\infty)},d^{(\infty)})\equiv\bigotimes_{i=1}^{\infty}(S_{i},d_{i})\equiv(\prod_{i=1}^{\infty}S_{i},\bigotimes_{i=1}^{\infty}d_{i}).
\]
Suppose, in addition, $(S_{n},d_{n})$ is a copy of the same metric
space $(S,d)$ for each $n\geq1$. Then we simply write $(S^{n},d^{n})\equiv(S^{(n)},d^{(n)})$
and $(S^{\infty},d^{\infty})\equiv(S^{(\infty)},d^{(\infty)})$. Thus,
in this case, 
\[
d(x,y)\equiv\bigvee_{i=1}^{n}d(x_{i},y_{i})
\]
for each $x=(x_{1},\cdots,x_{n}),y=(y_{1},\cdots,y_{n})\in S^{n}$,
and
\[
d^{\infty}(x,y)\equiv\sum_{i=1}^{\infty}2^{-i}(1\wedge d_{i}(x_{i},y_{i})).
\]
for each $x=(x_{1},x_{2},\cdots),y=(y_{1},y_{2},\cdots)\in S^{\infty}$. 

If, in addition, $(S_{n},d_{n})$ is locally compact for each $n\geq1$,
then the finite product space $(S^{(n)},d^{(n)})$ is locally compact
for each $n\geq1$, while the infinite product space $(S^{(\infty)},d^{(\infty)})$
is complete but not necessarily locally compact. If $(S_{n},d_{n})$
is a compact for each $n\geq1$, then both the finite and infinite
product spaces are compact.

Suppose $(S,d)$ is a  metric space. We will write $C_{u}(S,d)$,
or simply $C_{u}(S)$, for the space of real-valued functions functions
on $(S,d)$ with $domain(f)=S$ which are uniformly continuous on
each bounded subset of $S$. We will write $C_{ub}(S,d)$, or simply
$C_{ub}(S)$, for the subspace of $C_{u}(S)$ whose members are bounded.
Let $x_{\circ}$ be an arbitrary, but fixed, reference point in $(S,d)$.
A continuous function $f$ on $(S,d)$ is then said to vanish at infinity
if, for each $\varepsilon>0$, there exists $a>0$ such that $|f|\leq\varepsilon$
for each $x\in S$ with $d(x,x_{\circ})>a$. Write $C_{0}(S,d)$,
or simply $C_{0}(S)$, for the space of continuous functions on $(S,d)$
which vanish at infinity. A real-valued function $f$ on $S$ is said
to have a subset $A\subset S$ as\emph{ support}\index{support} if
$x\in domain(f)$ and $|f(x)|>0$ together imply $x\in A$. Then we
also say that $f$ is supported by $A$, or that $A$ supports $f$.
 We will write $C(S,d)$, or simply $C(S)$, for the subspace of $C_{u}(S,d)$
whose members have bounded supports. In the case where $(S,d)$ is
locally compact, $C(S)$ consists of continuous functions on $(S,d)$
with compact supports. Summing up,
\[
C(S)\subset C_{0}(S)\subset C_{ub}(S)\subset C_{u}(S).
\]

Suppose a subset $A$ of $R$ is nonempty. A number $b\in R$ is called
a lower bound of $A$, and $A$ said to bounded from below, if $b\leq a$
for each $a\in A$. A lower bound $b$ of $A$ is called the \emph{greatest
lower bound}\index{greatest lower bound}, or \emph{infimum}\index{infimum},
of $A$ if $b\geq b'$ for each lower bound $b'$ of $A$. In that
case, we write $\inf A\equiv b$. 

Similarly, a number $b\in R$ is called an upper bound of $A$, and
$A$ said to be bounded from above, if $b\geq a$ for each $a\in A$.
An upper bound $b$ of $A$ is called the \emph{least upper bound}\index{least upper bound},
or \emph{supremum}\index{supremum}, of $A$ if $b\leq b'$ for each
upper bound $b'$ of $A$. In that case, we write $\sup A\equiv b$. 

There is no constructive general proof for the existence of an infimum
for an subset of $R$ that is bounded from below. Existence needs
to be proved before each usage for each special case, much as in the
case of limits. In that regard, \cite{BishopBridges85} proves that,
if a non-empty subset $A$ of $R$ is totally bounded, then both $\inf A$
and $\sup A$ exist.

Suppose $f$ is a continuous function on a compact metric space $(K,d)$.
Then the last cited text proves that $\inf_{K}f\equiv\inf\{f(x):x\in K\}$
and $\sup_{K}f\equiv\sup\{f(x):x\in K\}$ exist.

\subsection*{Miscellaneous}

The symbols $\Rightarrow,\Leftarrow,$ and $\Leftrightarrow$ will
in general stand for ``only if'', ``if'', and ``if and only if''
respectively. An exception will be made where the symbol $\Rightarrow$
is used for weak convergence, defined later. The intended meaning
will be clear from context.

We will often write ``$x,y,\cdots,z\in A$'' as an abbreviation
for ``$\{x,y,\cdots,z\}\subset A$''. 

Unless it is otherwise indicated by context, the symbols $i,j,k,m,n,p$
will denote integers, the symbols $a,b$ will denote real numbers,
and the symbols $\varepsilon,\delta$ positive real numbers. For example
the statement ``for each $i\geq1"$ will mean ``for each integer
$i\geq1$''.

Suppose $(a_{n})_{n=1,2,\cdots}$ is a sequence of real numbers. Then
$a_{n}\rightarrow a$ stands for $\lim_{n\rightarrow\infty}a_{n}=a$.
We write  $a_{n}\uparrow a$ if $(a_{n})$ is a nondecreasing sequence
and $a_{n}\rightarrow a$. Similarly, we write $a_{n}\downarrow a$
if $(a_{n})$ is a nonincreasing sequence and $a_{n}\rightarrow a$.
More generally, suppose $f$ is a function on some subset $A\subset R$.
Then $f(x)\rightarrow a$ stands for $\lim_{x\rightarrow x_{0}}f(x)=a$
where $x_{0}$ can stand for a real number or for one of the symbols
$\infty$ or $-\infty$.

We use the common ``big O'' and ``small o'' notation. Suppose
$f$ and $g$ are functions on some subset $A\subset R$. Let $x_{0}$
stand for a real number or for one of the symbols $\infty$ or $-\infty$.
If for some $c>0$, we have $|f(x)|\leq c|g(x)|$ for all $x\in A$
in some neighborhood $B$ of $x_{0}$, then we write $f(x)=O(g(x))$.
If for each $c>0$, we have $|f(x)|\leq c|g(x)|$ for each $x\in A$
in some neighborhood $B$ of $x_{0}$, then we write $f(x)=o(g(x))$.
A subset $B\subset R$ is a neighborhood of $x_{0}$ if there exists
an open interval $(a,b)$ such that either (i) $x_{0}\in(a,b)$, (ii)
$b=x_{0}=\infty$, or (iii) $a=x_{0}=-\infty$. 

Finally, we use the symbol $\square$ to mark the end of a proof or
a definition.

\part{Probability Theory}

\chapter{Partitions of Unity}

In the Introduction, we summarized the basic concepts and theorems
about metric spaces from \cite{BishopBridges85}. Locally compact
metric spaces were introduced. They can be regarded as a simple, but
wide ranging, generalization of the real line. Most, if not all, metric
spaces in the present book are locally compact.

In the present chapter, we will define binary approximations and partitions
of unity for a locally compact metric space $(S,d)$. Roughly speaking,
a binary approximation is a digitization of $(S,d)$, a generalization
of the binary numbers which digitize the space $R$ of real numbers.
A partition of unity is then a sequence in $C(S,d)$ which serves
as a basis for $C(S,d)$ in the sense that each $f\in C(S,d)$ can
be approximated by linear combinations of members in the partition
of unity. 

We first cite a theorem from \cite{BishopBridges85} which guarantees
an abundance of compact subsets. 
\begin{thm}
\label{Thm.  Abundance of Compact Sets}\textbf{\emph{ (Abundance
of compact sets). }}Let $f:K\rightarrow R$ be a continuous function
on a compact metric space $(K,d)$ with $domain(f)=K$. Then, for
all but countably many real numbers $\alpha>\inf_{K}f$, the set $(f\leq\alpha)\equiv\{x\in K:f(x)\leq\alpha\}$
is compact.
\end{thm}
\begin{proof}
See Theorem (4.9) in Chapter 4 of \cite{BishopBridges85}.
\end{proof}
Classically, the set $(f\leq\alpha)$ is compact for each $\alpha\geq\inf_{K}f$,
without exception. Such a general statement would however imply the
principle of infinite search, and is therefore nonconstructive. Theorem
\ref{Thm.  Abundance of Compact Sets} above is  sufficient for all
our purposes. 
\begin{defn}
\label{Def. Convention (f<=00003Da) for compactness} \textbf{(Convention
for compact sets $(f\leq a)$).} We hereby adopt the convention that,
if the compactness of the set $(f\leq\alpha)$ is required in a discussion,
compactness has been explicitly or implicitly verified, usually by
proper prior selection of the constant $\alpha$, enabled by Theorem
\ref{Thm.  Abundance of Compact Sets}. $\square$
\end{defn}
The following corollary guarantees an abundance of compact neighborhoods
of a compact set.
\begin{cor}
\label{Cor. Abundance of Compact Nbds, balls} \textbf{\emph{(Abundance
of compact neighborhoods).}} Let $(S,d)$ be a locally compact metric
space, and let $K$ be a compact subset of $S$. Then the subset 
\[
K_{r}\equiv(d(\cdot,K)\leq r)\equiv\{x\in S:d(x,K)\leq r\}
\]
is compact for all but countably many $r>0$. 
\end{cor}
\begin{proof}
Let $n\geq1$ be arbitrary. Then $K_{n}\equiv(d(\cdot,K)\leq n)$
is a bounded set. Since $(S,d)$ is locally compact, there exists
a compact set $S_{n}$ such that $K_{n}\subset S_{n}\subset S$. The
continuous function $f$ on $(S_{n},d)$ defined by $f\equiv d(\cdot,K)$
has infimum $0$. Hence, by Theorem \ref{Thm.  Abundance of Compact Sets},
the set $\{x\in S_{n}:d(x,K)\leq r\}$ is compact for all but countably
many $r>0$. On the other hand, for all $r<n$ we have 
\[
K_{r}=K_{r}K_{n}\subset K_{r}S_{n}=\{x\in S_{n}:d(x,K)\leq r\}.
\]
Thus $K_{r}$ is compact for all $r\in(0,n)\cap A_{n}$, where $A_{n}$
contains all but countably many $r>0$. Define $A\equiv\bigcap_{n=1}^{\infty}A_{n}$.
Then $A$ contains all but countably many $r>0$. Now let $r\in(0,\infty)\cap A$
be arbitrary. Then $r\in(0,n)\cap A_{n}$ for some $n\geq1$, whence
$K_{r}$ is compact. 
\end{proof}
Separately, the next elementary metric space lemma will be convenient.
\begin{lem}
\label{Lem. If S compact, then subspace of fcts depending on finitely many coordinates is  dense in C(S^inf)}
\textbf{\emph{(If $(S,d)$ is compact, then the subspace of $C(S^{\infty},d^{\infty})$
whose members depend on finitely many coordinates is dense).}} Suppose
$(S,d)$ is a compact metric space. 

Let $n\geq1$ be arbitrary. Define the truncation function $j_{n}^{*}:S^{\infty}\rightarrow S^{\infty}$
by
\[
j_{n}^{*}(x_{1},x_{2},\cdots)\equiv(x_{1},x_{2},\cdots,x_{n},x_{\circ},x_{\circ},\cdots)
\]
for each $(x_{1},x_{2},\cdots)\in S^{\infty}$. Then $j_{n}^{*}\circ j_{m}^{*}=j_{n}^{*}$
for each $m\geq n$. Let 
\begin{equation}
L_{0,n}\equiv\{f\in C(S^{\infty},d^{\infty}):f=f\circ j_{n}^{*}\}.\label{eq:temp-205}
\end{equation}
Let $L_{0,\infty}\equiv\bigcup_{n=1}^{\infty}L_{0,n}$. Then $L_{0,n}\subset L_{0,n+1}$.
Moreover, the following holds.

1. $L_{0,n}$ and $L_{0,\infty}$ are linear subspaces of $C(S^{\infty},d^{\infty})$,
and consist of functions which depend, respectively, on the first
$n$ and on finitely many coordinates.

2. The subspace $L_{0,\infty}$ is dense in $C(S^{\infty},d^{\infty})$
relative to the supremum norm $\left\Vert \cdot\right\Vert $. Specifically,
let $f\in C(S^{\infty},d^{\infty})$ be arbitrary, with a modulus
of continuity $\delta_{f}$. Then $f\circ j_{n}^{*}\in L_{0,n}$.
Moreover, for each $\varepsilon>0$ we have $\bigl\Vert f-f\circ j_{n}^{*}\bigr\Vert\leq\varepsilon$
if $n>-\log_{2}(\delta_{f}(\varepsilon))$. In particular, if $f$
has Lipschitz constant $c>0$, then $\bigl\Vert f-f\circ j_{n}^{*}\bigr\Vert\leq\varepsilon$
if $n>\log_{2}(c\varepsilon^{-1})$.
\end{lem}
\begin{proof}
Let $m\geq n\geq1$ and $w\in S^{\infty}$ be arbitrary. Then, for
each $(x_{1},x_{2},\cdots)\in S^{\infty}$, we have
\[
j_{n}^{*}(j_{m}^{*}(x_{1},x_{2},\cdots))=j_{n}^{*}(x_{1},x_{2},\cdots,x_{m},x_{\circ},x_{\circ},\cdots)
\]
\[
=(x_{1},x_{2},\cdots,x_{n},x_{\circ},x_{\circ},\cdots)=j_{n}^{*}(x_{1},x_{2},\cdots).
\]
Hence $j_{n}^{*}\circ j_{m}^{*}=j_{n}^{*}$.

1. It is clear from the defining equality \ref{eq:temp-205} that
$L_{0,n}$ is a linear subspace of $C(S^{\infty},d^{\infty})$. Let
$f\in L_{0,n}$ be arbitrary. Then $f=f\circ j_{n}^{*}=f\circ j_{n}^{*}\circ j_{m}^{*}=f\circ j_{m}^{*}$.
Hence $f\in L_{0,m}$. Thus $L_{0,n}\subset L_{0,m}$. Consequently,
$L_{0,\infty}\equiv\bigcup_{p=1}^{\infty}L_{0,p}$ is a union of a
nondecreasing sequence of linear subspaces of $C(S^{\infty},d^{\infty})$,
and is therefore also a linear subspace of $C(S^{\infty},d^{\infty})$.

2. Let $f\in C(S^{\infty},d^{\infty})$ be arbitrary, with a modulus
of continuity $\delta_{f}$. Let $\varepsilon>0$ be arbitrary. Suppose
$n>-\log_{2}(\delta_{f}(\varepsilon))$. Then $2^{-n}<\delta_{f}(\varepsilon)$.
Let $(x_{1},x_{2},\cdots)\in S^{\infty}$ be arbitrary. Then
\[
d^{\infty}((x_{1},x_{2},\cdots),j_{n}^{*}(x_{1},x_{2},\cdots))
\]
\[
=d^{\infty}((x_{1},x_{2},\cdots),(x_{1},x_{2},\cdots,x_{n},x_{\circ},x_{\circ},\cdots))
\]
\[
\equiv\sum_{k=1}^{n}2^{-k}\widehat{d}(x_{k},x_{k})+\sum_{k=n+1}^{\infty}2^{-k}\widehat{d}(x_{k},x_{\circ})\leq0+2^{-n}<\delta_{f}(\varepsilon),
\]
where $\widehat{d}\equiv1\wedge d$. Hence 
\[
|f(x_{1},x_{2},\cdots)-f\circ j_{n}^{*}(x_{1},x_{2},\cdots)|<\varepsilon,
\]
where $(x_{1},x_{2},\cdots)\in S^{\infty}$ is arbitrary. We conclude
that $\left\Vert f-f\circ j_{n}^{*}\right\Vert \leq\varepsilon$,
as alleged. 
\end{proof}

\section{Binary Approximations}

Let $(S,d)$ be an arbitrary locally compact metric space. Then $S$
has a countable dense subset. A binary approximation, defined below
in this section, is a structured and well-quantified countable dense
subset.

Recall that (i) $|A|$ denotes the number of elements in an arbitrary
finite set $A$, (ii) a subset $A$ of $S$ is said to be metrically
discrete if, for each $y,z\in A$, either $y=z$ or $d(y,z)>0$, and
(iii) a finite subset $A$ of $K\subset S$ is called an $\varepsilon-$approximation
of $K$ if for each $x\in K$ there exists $y\in A$ with that $d(x,y)<\varepsilon$.
Classically, each subset of $(S,d)$ is metrically discrete. 

$\square$
\begin{defn}
\label{Def. Binary approximationt and Modulus of local compactness}
\textbf{(Binary approximation and modulus of local compactness).}
Let $(S,d)$ be a locally compact metric space, with an arbitrary,
but fixed, reference point $x_{\circ}$. Let $A_{0}\equiv\{x_{\circ}\}\subset A_{1}\subset A_{2}\subset\cdots$
be a sequence of metrically discrete and finite subsets of $S$. For
each $n\geq1$, let $\kappa_{n}\equiv|A_{n}|$. Suppose
\begin{equation}
(d(\cdot,x_{\circ})\leq2^{n})\subset\bigcup_{x\in A(n)}(d(\cdot,x)\leq2^{-n})\label{eq:temp-505}
\end{equation}
and
\begin{equation}
\bigcup_{x\in A(n)}(d(\cdot,x)\leq2^{-n+1})\subset(d(\cdot,x_{\circ})\leq2^{n+1})\label{eq:temp-502}
\end{equation}
for each $n\geq1$. Then the sequence $\xi\equiv(A_{n})_{n=1,2,\cdots}$
of subsets is called a \emph{binary approximation}\index{binary approximation}
for $(S,d)$ relative to $x_{\circ}$, and the sequence of integers
\[
\left\Vert \xi\right\Vert \equiv(\kappa_{n})_{n=1,2,\cdots}\equiv(|A_{n}|)_{n=1,2,\cdots}
\]
is called the \emph{modulus of local compactness}\index{modulus of local compactness}
of $(S,d)$ corresponding to $\xi$. 

Thus a binary approximation is an expanding sequence of $2^{-n}$-approximation
for $(d(\cdot,x_{\circ})\leq2^{n})$ as $n\rightarrow\infty$. The
next proposition shows that the definition is not vacuous. $\square$

First note that $\bigcup_{n=1}^{\infty}A_{n}$ is dense in $(S,d)$
in view of relation \ref{eq:temp-505}. In the case where $(S,d)$
is compact, for $n\geq1$ so large that $S=(d(\cdot,x_{\circ})\leq2^{n})$,
relation \ref{eq:temp-505} says that we need at most $\kappa_{n}$
points to make a $2^{-n}$-approximation of $S$. The number $\log\kappa_{n}$
is thus a bound for Kolmogorov's $2^{-n}$-\emph{entropy\index{Kolmogorov's varepsilon-entropy@Kolmogorov's $\varepsilon$-entropy}}
of the compact metric space $(S,d)$, which represents the informational
content in a $2^{-n}-$approximation of $S$. (See \cite{Lorentz66}
for a definition of $\varepsilon$-entropy\index{epsilon-entropy}). 
\end{defn}
\begin{lem}
\label{Lem. Existence of metrically discrete epsilon approximations}
\textbf{\emph{(Existence of metrically discrete $\varepsilon$-approximations).}}
Let $K$ be a compact subset of the locally compact metric space $(S,d)$.
Let $A_{0}$ be a metrically discrete finite subset of $K$ . Let
$\varepsilon>0$ be arbitrary. Then the following holds.

1. There exists a metrically discrete finite subset $A_{1}$ of $K$
such that (i) $A_{0}\subset A_{1}$ and (ii) $A_{1}$ is an $\varepsilon$-approximation
of $K$.

2. In particular, there exists a metrically discrete finite set $A_{1}$
which is an $\varepsilon$-approximation of $K$.
\end{lem}
\begin{proof}
1. Let $A\equiv\{y_{1},\cdots,y_{m}\}$ be an $\varepsilon$-approximation
of $K$. We apply a weeding procedure on $A$. Successively examine
each $y_{i}\in A$. Either (i) $d(y_{i},x)>0$ for each $x\in A_{0}$,
and $d(y_{i},y_{j})>0$ for each $j$ with $1\leq j\leq i-1$, or
(ii) $d(y_{i},x)<\varepsilon$ for some $x\in A_{0}$ or $d(y_{i},y_{j})<\varepsilon$
for some $j$ with $1\leq j\leq i-1$. In case (ii), discard the point
$y_{i}$ from $A$, decrease $m$ by 1, relabel the thus-far surviving
points as $y_{1},\cdots,y_{m}$, redefine $A\equiv\{y_{1},\cdots,y_{m}\}$,
and repeat the procedure. Observe that, after at most $m$ steps,
each surviving member in the set $A$, if any, is of positive distance
to each member of $A_{0}$, and is of positive distance to any other
surviving member of $A$. Now define $A_{1}\equiv A_{0}\cup A$. In
view of the above observation , the set $A_{0}\cup A$ is metrically
discrete. In other words, $A_{1}$ is metrically discrete. 

Let $B$ be the set of discarded points. Let $y\in K$ be arbitrary.
Since $A\cup B$ is equal to the original $\varepsilon$-approximation
of $K$, we see that $y$ is of distance less than $\varepsilon$
to some member of $A$, or it is of distance less than $\varepsilon$
to some member of $B$. At the same time, each member of $B$ is of
distance less than $\varepsilon$ to some members of $A_{0}\cup A$,
according to condition (ii). Hence, in either case, $y$ is of distance
less than $2\varepsilon$ to some members of $A_{0}\cup A$, thanks
to the triangle inequality. We conclude that $A_{1}\equiv A_{0}\cup A$
is a $2\varepsilon$-approximation of $K$. Assertion 1 has been proved.

2. Assertion 2 is a special case of Assertion 1 by taking $A_{0}\equiv\phi$.
\end{proof}
\begin{prop}
\textbf{\emph{\label{Prop. Existence of separant} (Existence of binary
approximations).}} Each locally compact metric space $(S,d)$ has
a binary approximation.
\end{prop}
\begin{proof}
Let $x_{\circ}\in S$ be an arbitrary, but fixed, reference point.
Proceed inductively on $n\geq1$ to construct a metrically discrete
and finite subset $A_{n}$ of $S$ to satisfy relations \ref{eq:temp-505}
and \ref{eq:temp-502}. 

Let $n\geq0$ be arbitrary. Let $A_{0}\equiv\{x_{\circ}\}$. Suppose
the set $A_{n}$ has been constructed for some $n\geq0$, such that,
if $n\geq1$, then (i) $A_{n}$ is metrically discrete and finite,
and (ii) relations \ref{eq:temp-505} and \ref{eq:temp-502} are satisfied.
Proceed to construct $A_{n+1}$.

To that end, write $\varepsilon\equiv2^{-n-2}$, and take any $r\in[2^{n+1},2^{n+1}+\varepsilon)$
such that 
\[
K\equiv(d(\cdot,x_{\circ})\leq r)
\]
is compact. This is possible in view of Corollary \ref{Cor. Abundance of Compact Nbds, balls}.
If $n=0$, then $A_{n}\equiv\phi\subset K$ trivially. If $n\geq1$,
then, according to the induction hypothesis, $A_{n}$ is metrically
discrete, and by relation \ref{eq:temp-502}, we have
\[
A_{n}\subset\bigcup_{x\in A(n)}(d(\cdot,x)\leq2^{-n+1})\subset(d(\cdot,x_{\circ})\leq2^{n+1})\subset K.
\]
Hence we can apply Lemma \ref{Lem. Existence of metrically discrete epsilon approximations}
to construct a $2^{-n-1}$ approximation $A_{n+1}$ of $K$ which
is metrically discrete and finite. We conclude that
\[
(d(\cdot,x_{\circ})\leq2^{n+1})\subset K\subset\bigcup_{x\in A(n+1)}(d(\cdot,x)\leq2^{-n-1})
\]
proving relation \ref{eq:temp-505} for $n+1$.

Now let 
\[
y\in\bigcup_{x\in A(n+1)}(d(\cdot,x)\leq2^{-n})
\]
be arbitrary. Then $d(y,x)\leq2^{-n}$ for some $x\in A_{n+1}\equiv A_{n}\cup A\subset K$.
Therefore 
\[
d(x,x_{\circ})\leq r<2^{n+1}+\varepsilon.
\]
Consequently
\[
d(y,x_{\circ})\leq2^{-n}+2^{n+1}+\varepsilon\equiv2^{-n}+2^{n+1}+2^{-n-2}\leq2^{n+2}.
\]
Thus 
\[
\bigcup_{x\in A(n+1)}(d(\cdot,x)\leq2^{-n})\subset(d(\cdot,x_{\circ})\leq2^{n+2}),
\]
proving relation \ref{eq:temp-502} for $n+1$. Induction is completed.
Thus the sequence $\xi\equiv(A_{n})_{n=1,2,\cdots}$ satisfies all
the conditions in Definition \ref{Def. Binary approximationt and Modulus of local compactness}
to be a binary approximation of $(S,d)$.
\end{proof}
\begin{defn}
\label{Def. Finite product and power of binary approx of locally compact S}
\textbf{(Finite product and power of binary approximations).} Let
$n\geq1$ be arbitrary. For each $i=1,\cdots,n$, let $(S_{i},d_{i})$
be a locally compact metric space, with a reference point $x_{i,\circ}\in S_{i}$
and with a binary approximation $\xi_{i}\equiv(A_{i,p})_{p=1,2,\cdots}$
relative to $x_{i,\circ}$. Let $(S^{(n)},d^{(n)})\equiv(\prod_{i=1}^{n}S_{i},\bigotimes_{i=1}^{n}d_{i})$
be the product metric space, with $x_{\circ}^{(n)}\equiv(x_{1,\circ},\cdots,x_{n,\circ})$
designated as the reference point in $(S^{(n)},d^{(n)})$.

For each $p\geq1$, let $A_{p}^{(n)}\equiv A_{1,p}\times\cdots\times A_{n,p}$.
The next lemma proves that $(A_{p}^{(n)})_{p=1,2,\cdots}$ is a binary
approximation of $(S^{(n)},d^{(n)})$ relative to $x_{\circ}^{(n)}$.
We will call $\xi^{(n)}\equiv(A_{p}^{(n)})_{p=1,2,\cdots}$ the \index{product binary approximation, finite}\emph{product
binary approximation} of  $\xi_{1},\cdots,\xi_{n}$, and write $\xi^{(p)}\equiv\xi_{1}\otimes\cdots\otimes\xi_{p}$.
If $(S_{i},d_{i})=(S,d)$ for some locally compact metric space, with
$x_{i,\circ}=x_{\circ}$ and $\xi_{i}=\xi$ for each $i=1,\cdots,n$,
we will call $\xi^{(n)}$ the $n$-th \index{power of a binary approximation, finite}\emph{power}
of $\xi$, and write $\xi^{n}\equiv\xi^{(n)}$. $\square$
\end{defn}
\begin{lem}
\label{Lem. Product separant is a separant} \textbf{\emph{(Finite
product binary approximation is indeed a binary approximation).}}
Use the assumptions and notations in Definition \ref{Def. Finite product and power of binary approx of locally compact S}.
Then $\xi^{(n)}$ is indeed a binary approximation of $(S^{(n)},d^{(n)})$
relative to $x_{\circ}^{(n)}$. Let $\left\Vert \xi_{i}\right\Vert \equiv(\kappa_{i,p})_{p=1,2,\cdots}\equiv(|A_{i,p}|)_{p=1,2,\cdots}$
be the modulus of local compactness of $(S_{i},d_{i})$ corresponding
to $\xi_{i}$, for each $i=1,\cdots,n$. Let $\left\Vert \xi^{(n)}\right\Vert $
be the modulus of local compactness of $(S^{(n)},d^{(n)})$ corresponding
to $\xi^{(n)}$. Then $\left\Vert \xi^{(n)}\right\Vert =(\prod_{i=1}^{n}\kappa_{i,p})_{p=1,2,\cdots}$. 

In particular, if $\xi_{i}\equiv\xi$ for each $i=1,\cdots,n$, for
some binary approximation $\xi$ of some locally compact metric space
$(S,d)$, then $\left\Vert \xi^{n}\right\Vert =(\kappa_{p}^{n})_{p=1,2,\cdots}$.
\end{lem}
\begin{proof}
Recall that $A_{p}^{(n)}\equiv A_{1,p}\times\cdots\times A_{n,p}$
for each $p\geq1$. Hence $A_{1}^{(n)}\subset A_{2}^{(n)}\subset\cdots$. 

1. Let $p\geq1$ be arbitrary. Let
\[
x\equiv(x_{1},\cdots,x_{n}),y\equiv(y_{1},\cdots,y_{n})\in A_{p}^{(n)}\equiv A_{1,p}\times\cdots\times A_{n,p}
\]
be arbitrary. For each $i=1,\cdots,n$, because $(A_{i,q})_{q=1,2,\cdots}$
is a binary approximation, the set $A_{i,p}$ is metrically discrete.
Hence either (i) $x_{i}=y_{i}$ for each $i=1,\cdots,n$, or (ii)
$d_{i}(x_{i},y_{i})>0$ for some $i=1,\cdots,n$. In case (i) we have
$x=y$. In case (ii) we have 
\[
d^{(n)}(x,y)\equiv\bigvee_{j=1}^{n}d_{j}(x_{j},y_{j})\geq d_{i}(x_{i},y_{i})>0.
\]
Thus $A_{p}^{(n)}$ is metrically discrete. 

2. Next note that
\[
(d^{(n)}(\cdot,x_{\circ}^{(n)})\leq2^{p})\equiv\{(y_{1},\cdots,y_{n})\in S^{(n)}:\bigvee_{i=1}^{n}d_{i}(y_{i},x_{i,\circ})\leq2^{p}\}
\]
\[
=\bigcap_{i=1}^{n}\{(y_{1},\cdots,y_{n})\in S^{(n)}:d_{i}(y_{i},x_{i,\circ})\leq2^{p}\}
\]
\begin{equation}
\subset C\equiv\bigcap_{i=1}^{n}\bigcup_{z(i)\in A(i,p)}\{(y_{1},\cdots,y_{n})\in S^{(n)}:d_{i}(y_{i},z_{i})\leq2^{-p}\},\label{eq:temp-255}
\end{equation}
where the last inclusion is due to relation \ref{eq:temp-505} applied
to the binary approximation $(A_{i,q})_{q=1,2,\cdots}$. Basic Boolean
operations yield 
\[
C=\bigcup_{(z(1),\cdots,z(n))\in A(1,p)\times\cdots\times A(n,p)}\bigcap_{i=1}^{n}\{(y_{1},\cdots,y_{n})\in S^{(n)}:d_{i}(y_{i},z_{i})\leq2^{-p}\}
\]
\[
=\bigcup_{(z(1),\cdots,z(n))\in A(1,p)\times\cdots\times A(n,p)}\{(y_{1},\cdots,y_{n})\in S^{(n)}:\bigvee_{i=1}^{n}d_{i}(y_{i},z_{i})\leq2^{-p}\}
\]
\begin{equation}
\bigcup_{x\in A_{p}^{(n)}}(d^{(n)}(\cdot,x)\leq2^{-p}).\label{eq:temp-503}
\end{equation}
Thus relation \ref{eq:temp-505} has been verified for the sequence
$\xi^{(n)}\equiv(A_{q}^{(n)})_{q=1,2,\cdots}$. 

Reversing direction, we have, similarly, 
\[
\bigcup_{x\in A_{p}^{(n)}}(d^{(n)}(\cdot,x)\leq2^{-p+1})
\]
\[
=\bigcap_{i=1}^{n}\bigcup_{z(i)\in A(i,p)}\{(y_{1},\cdots,y_{n})\in S^{(n)}:d_{i}(y_{i},z_{i})\leq2^{-p+1}\}
\]
\[
\subset\bigcap_{i=1}^{n}\{(y_{1},\cdots,y_{n})\in S^{(n)}:d_{i}(y_{i},x_{i,\circ})\leq2^{p+1}\}
\]
\[
=\{(y_{1},\cdots,y_{n})\in S^{(n)}:\bigvee_{i=1}^{n}d_{i}(y_{i},x_{i,\circ})\leq2^{p+1}\}
\]
\[
=(d^{(n)}(\cdot,x_{\circ}^{(n)})\leq2^{p+1}),
\]
which verifies relation \ref{eq:temp-502} for the sequence $\xi^{(n)}\equiv(A_{q}^{(n)})_{q=1,2,\cdots}.$.
Thus all the conditions in Definition \ref{Def. Binary approximationt and Modulus of local compactness}
have been proved for the sequence $\xi^{(n)}$ to be a binary approximation
of $(S^{(n)},d^{(n)})$ relative to $x_{\circ}^{(n)}$. Moreover
\[
\left\Vert \xi^{(n)}\right\Vert \equiv(|A_{q}^{(n)}|)_{q=1,2,\cdots}=(\prod_{i=1}^{n}|A{}_{i,q}|)_{q=1,2,\cdots}\equiv(\prod_{i=1}^{n}\kappa_{i,q})_{=1,2,\cdots}.
\]
\end{proof}
We now extend the construction of product binary approximations to
the infinite product space $(S^{\infty},d^{\infty})$ in the case
where $(S,d)$ is compact. As usual, $\widehat{d}\equiv1\wedge d$.
\begin{defn}
\label{Def. Countable product of binary approxximations; compact}
\textbf{(Countable product of binary approximation for compact space).}
Suppose $(S,d)$ is a compact metric space, with a reference point
$x_{\circ}\in S$, and with a binary approximation $\xi\equiv(A_{n})_{n=1,2,\cdots}$
relative to $x_{\circ}$. Let $(S^{\infty},d^{\infty})$ be the countable
power of metric space $(S,d)$, with $x_{\circ}^{\infty}\equiv(x_{\circ},x_{\circ},\cdots)$
designated as the reference point in $(S^{\infty},d^{\infty})$.

For each $n\geq1$, define the subset 
\[
B_{n}\equiv A_{n+1}^{n+1}\times\{x_{\circ}^{\infty}\}=\{(x_{1},\cdots,x_{n+1},x_{\circ},x_{\circ}\cdots):x_{i}\in A_{n+1}\;\mathrm{for\;each}\;i=1,\cdots,n+1\}.
\]
The next lemma proves that $\xi^{\infty}\equiv(B_{n})_{n=1,2,}$ is
a binary approximation of $(S^{\infty},d^{\infty})$ relative to $x_{\circ}^{\infty}$.
We will call $\xi^{\infty}$ the countable \index{countable power binary approximation}\emph{
power of the binary approximation}  $\xi$. $\square$
\end{defn}
\begin{lem}
\label{Lem. Countable Product of binary approximation is a binary approxt; compacts}
\textbf{\emph{(Countable product binary approximation for infinite
product of compact metric spaces is indeed a binary approximation).}}
Suppose $(S,d)$ is a compact metric space, with a reference point
$x_{\circ}\in S$, and with a binary approximation $\xi\equiv(A_{n})_{n=1,2,\cdots}$
relative to $x_{\circ}$. Without loss of generality, assume that
$d\leq1$. Then the sequence $\xi^{\infty}\equiv(B_{n})_{n=1,2,}$
in Definition \ref{Lem. Countable Product of binary approximation is a binary approxt; compacts}
is indeed a binary approximation of $(S^{\infty},d^{\infty})$ relative
to $x_{\circ}^{\infty}$. 

Let $\left\Vert \xi\right\Vert \equiv(\kappa_{n})_{n=1,2,\cdots}\equiv(|A_{n}|)_{n=1,2,\cdots}$
denote the modulus of local compactness of $(S,d)$ corresponding
to $\xi$. Then the modulus of local compactness of $(S^{\infty},d^{\infty})$
corresponding to $\xi^{\infty}$ is given by 
\[
\left\Vert \xi^{\infty}\right\Vert =(\kappa_{n+1}^{n+1})_{n=1,2,\cdots}.
\]
\end{lem}
\begin{proof}
Let $n\geq1$ be arbitrary. 

1. Let
\[
x\equiv(x_{1},\cdots,x_{n+1},x_{\circ},x_{\circ}\cdots),y\equiv(y_{1},\cdots,y_{n+1},x_{\circ},x_{\circ}\cdots)\in B_{n}
\]
be arbitrary. Since $A_{n+1}$ is metrically discrete we have either
(i) $x_{i}=y_{i}$ for each $i=1,\cdots,n$+1, or (ii) $\widehat{d}(x_{i},y_{i})>0$
for some $i=1,\cdots,n+1$. In case (i) we have $x=y$. In case (ii)
we have 
\[
d^{\infty}(x,y)\equiv\sum_{j=1}^{\infty}2^{-j}\widehat{d}(x_{j},y_{j})\geq2^{-i}\widehat{d}(x_{i},y_{i})>0.
\]
Thus we see that $B_{n}$ is metrically discrete. 

2. Next, let $y\equiv(y_{1},y_{2},\cdots)\in S^{\infty}$ be arbitrary.
Let $j=1,\cdots,n+1$ be arbitrary. Then 
\[
y_{j}\in(d(\cdot,x_{\circ})\leq2^{n+1})\subset\bigcup_{z\in A(n+1)}(d(\cdot,z)\leq2^{-n-1}),
\]
where the first containment relation is a trivial consequence of the
hypothesis that $d\leq1$, and the second is an application of relation
\ref{eq:temp-505}. Hence there exists some $u_{j}\in A_{n+1}$ with
$d(y_{j},u_{j})\leq2^{-n-1}$. It follows that 
\[
u\equiv(u_{1},\cdots,u_{n+1},x_{\circ},x_{\circ},\cdots)\in B_{n},
\]
and
\[
d^{\infty}(y,u)\leq\sum_{j=1}^{n+1}2^{-j}\widehat{d}(y_{j},u_{j})+\sum_{j=n+2}^{\infty}2^{-j}
\]
\[
\leq\sum_{j=1}^{n+1}2^{-j}2^{-n-1}+2^{-n-1}<2^{-n-1}+2^{-n-1}=2^{-n}.
\]
We conclude that
\[
(d^{\infty}(\cdot,x_{\circ}^{\infty})\leq2^{n})=S^{\infty}\subset\bigcup_{u\in B(n)}(d^{\infty}(\cdot,u)\leq2^{-n}).
\]
where the equality is trivial because $d^{\infty}\leq1$. Thus relation
\ref{eq:temp-505} is verified for the sequence $(B_{n})_{n=1,2,\cdots}$.
At the same time, we have trivially 
\[
\bigcup_{u\in B(n)}(d^{\infty}(\cdot,u)\leq2^{-n+1})\subset S^{\infty}=(d^{\infty}(\cdot,x_{\circ}^{\infty})\leq2^{n+1}).
\]
Thus all the conditions in Definition \ref{Def. Binary approximationt and Modulus of local compactness}
have been verified for the sequence $\xi^{\infty}\equiv(B_{n})_{n=1,2,\cdots}.$
to be a binary approximation of $(S^{\infty},d^{\infty})$ relative
to $x_{\circ}^{\infty}$. Moreover,
\[
\left\Vert \xi^{\infty}\right\Vert \equiv(|B_{n}|)_{n=1,2,\cdots}=(|A_{n+1}^{n+1}|)_{n=1,2,\cdots}\equiv(\kappa_{n+1}^{n+1})_{n=1,2,\cdots}.
\]
\end{proof}

\section{Partitions of Unity}

In this section, we define and construct a partition of unity determined
by a binary approximation of a locally compact metric space. Versions
of partitions of unity abound in the literature, providing approximate
linear bases in the analysis of linear spaces of functions. The present
version, roughly speaking, furnishes an approximate linear basis for
$C(S)$, the space of continuous functions with compact supports on
a locally compact metric space. 

First we list, without proof, an elementary lemma for ease of later
reference.
\begin{lem}
\label{Ex.  Lipshitz Constants} \textbf{\emph{(Elementary lemma for
Lipschitz continuous functions).}}\emph{ }Let $(S,d)$ be an arbitrary
metric space. A real-valued function $f$ on S is said to be \emph{Lipschitz
continuous\index{Lipschitz continuous function}}, with Lipschitz
constant $c\geq0$ if $|f(x)-f(y)|\leq cd(x,y)$ for each $x,y\in S$.
We will then also say simply that the function has Lipschitz constant
$c$.

Let $x_{\circ}\in S$ be an arbitrary, but fixed, reference point.
Let $f,g$ be real-valued functions with Lipschitz constants $a,b$
respectively on $S$. Then the following holds. 

1. $d(\cdot,x_{\circ})$ has Lipschitz constant $1$.

2. $\alpha f+\beta g$ has Lipschitz constant $|\alpha|a+|\beta|b$
for each $\alpha,\beta\in R$. 

3. $f\vee g$ and $f\wedge g$ have Lipschitz constant $a\vee b$. 

4. $1\wedge(1-cd(\cdot,x_{\circ}))_{+}$ has Lipschitz constant $c$
for each $c>0$, 

5. If $\left\Vert f\right\Vert \vee\left\Vert g\right\Vert \leq1$
then $fg$ has Lipschitz constant $a+b$, 

6. Suppose $(S',d')$ is a locally compact metric space. Suppose $f'$
is a real-valued functions on $S'$, with Lipschitz constant $a'>0$.
Suppose $\left\Vert f\right\Vert \vee\left\Vert f'\right\Vert \leq1$.
Then $f\otimes f':S\times S'\rightarrow R$ has Lipschitz constant
$a+a'$ where $S\times S'$ is equipped with the product metric $\overline{d}\equiv d\otimes d'$,
and where $f\otimes f'(x,x')\equiv f(x)f'(x')$ for each $(x,x')\in S\times S'$. 

7. Assertion 6 above can be generalized to a $p$-fold product $f\otimes f'\otimes\cdots\otimes f^{(p)}$. 
\end{lem}
The next definition and proposition are essentially  Proposition 6.15
in $\qquad\qquad$ \cite{BishopBridges85}. 
\begin{defn}
\label{Def.. Epsilon-parittion of unity} \textbf{($\varepsilon$-partition
of unity).} Let $A$ be an arbitrary metrically discrete and finite
subset of a locally compact metric space $(S,d)$. Because the set
$A$ is finite, we can write $A=\{x_{1},\cdots,x_{\kappa}\}$ for
some sequence $x\equiv(x_{1},\cdots,x_{\kappa})$ where $x:\{1,\cdots,\kappa\}\rightarrow A$
is an enumeration of the finite set $A$. Thus $|A|\equiv\kappa$.
Let $\varepsilon>0$ be arbitrary. Define, for each $k=1,\cdots,\kappa$,
\begin{equation}
\eta_{k}\equiv1\wedge(2-\varepsilon^{-1}d(\cdot,v_{k}))_{+}\in C(S),\label{eq:temp-233-1-2}
\end{equation}
and
\begin{equation}
g_{k}^{+}\equiv\eta_{1}\vee\cdots\vee\eta_{k}\in C(S).\label{eq:temp-130-1-2}
\end{equation}
In addition, define $g_{0}^{+}\equiv0$ and, for each $k=1,\cdots,\kappa$,
define
\begin{equation}
g_{v(k)}\equiv g_{k}^{+}-g_{k-1}^{+}.\label{eq:temp-506-2}
\end{equation}
Then the subset $\{g_{x}:x\in A\}$ of $C(S)$ is called the $\varepsilon$\emph{-partition
of unity}\index{varepsilon-partition of unity@$\varepsilon$-partition of unity}
of $(S,d)$, determined by the enumerated set $A$. The members of
$\{g_{x}:x\in A\}$ are called the \emph{basis functions}\index{basis functions of an varepsilon-partition of unity@basis functions of an $\varepsilon$-partition of unity}
of the $\varepsilon$\emph{-}partition of unity\emph{.} $\square$ 
\end{defn}
\begin{prop}
\label{Prop. Properties of  epsilon partition of unity } \textbf{\emph{(Properties
of $\varepsilon$-partition of unity).}} Let $A=\{x_{1},\cdots,x_{\kappa}\}$
be an arbitrary metrically discrete and enumerated finite subset of
a locally compact metric space $(S,d)$. Let $\varepsilon>0$ be arbitrary.
Let $\{g_{x}:x\in A\}$ be the $\varepsilon$\emph{-}partition of
unity determined by the enumerated set $A$. Then the following holds.

1. $g_{x}$ has values in $[0,1]$ and has $(d(\cdot,x)<2\varepsilon)$
as support, for each $x\in A$.

2. $\sum_{x\in A}g_{x}\leq1$ on $S$.

3. $\sum_{x\in A}g_{x}=1$ on $\bigcup_{x\in A}(d(\cdot,x)\leq\varepsilon)$.

4. For each $x\in A$, the functions $g_{x}$, $\sum_{y\in A;y<x}g_{y}$,
and $\sum_{y\in A}g_{y}$ have Lipschitz constant $2\varepsilon^{-1}$.
Here $y<x$ means $y=x_{i}$ and $x=x_{j}$ for some $i,j\in\{1,\cdots,\kappa\}$
with $i<j$.
\end{prop}
\begin{proof}
1. Use the notations in Definition \ref{Def.. Epsilon-parittion of unity}.
Let $k=1,\cdots,\kappa$ be arbitrary. Suppose $y\in S$ is such that
$g_{v(k)}(y)>0$. By the defining equality \ref{eq:temp-506-2}, it
follows that $g_{k}^{+}(y)>g_{k-1}^{+}(y)$. Hence $\eta_{k}(y)>0$
by equality \ref{eq:temp-130-1-2}. Equality \ref{eq:temp-233-1-2}
then implies that $d(y,v_{k})<2\varepsilon$. In short $g_{v(k)}$
has $(d(\cdot,v_{k})<2\varepsilon)$ as support. In general $g_{v(k)}\geq0$
in view of equalities \ref{eq:temp-233-1-2}, \ref{eq:temp-130-1-2},
and \ref{eq:temp-506-2}.

2. $\sum_{x\in A}g_{x}=g_{\kappa}^{+}\equiv\eta_{1}\vee\cdots\vee\eta_{\kappa}\leq1$.
Condition 2 is verified. Consequently $g_{x}\leq1$ for each $x\in A$.

3. Suppose $y\in S$ is such that $d(y,v_{k})\leq\varepsilon$ for
some $k=1,\cdots,\kappa$. Then $\eta_{k}(y)=1$ according to equality
\ref{eq:temp-233-1-2}. Hence $\sum_{x\in A}g_{x}(y)\equiv g_{k}^{+}(y)=1$
by equality \ref{eq:temp-130-1-2}.

4. Now let $k=1,\cdots,\kappa$ be arbitrary. Refer to Exercise \ref{Ex.  Lipshitz Constants}
for basic operations of Lipschitz constants. Then, in view of the
defining equality \ref{eq:temp-233-1-2}, the function $\eta_{k}$
has Lipschitz constant $\varepsilon^{-1}$. Hence $g_{k}^{+}\equiv\eta_{0}\vee\cdots\vee\eta_{k}$
has Lipschitz constant $\varepsilon^{-1}$. In particular, $\sum_{y\in A}g_{y}\equiv g_{\kappa}^{+}$
has Lipschitz constant $\varepsilon^{-1}$. Moreover, for each $k=1,\cdots,\kappa$,
the function 
\[
\sum_{y\in A;y<v(k)}g_{y}\equiv\sum_{i=1}^{k-1}g_{v(i)}=g_{k}^{+}
\]
has Lipschitz constant $\varepsilon^{-1}$ whence $g_{v(k)}\equiv g_{k}^{+}-g_{k-1}^{+}$
has Lipschitz constant $2\varepsilon^{-1}$. In other words, for each
$x\in A$, the functions $\sum_{y\in A}g_{y}$, $\sum_{y\in A;y<x}g_{y}$,
and $g_{x}$ have Lipschitz constant $c\equiv2\varepsilon^{-1}$. 
\end{proof}
Recall that if $f\in C(S)$ then $\sup_{x\in S}|f(x)|$ exists and
is denoted by $\left\Vert f\right\Vert $.
\begin{defn}
\label{Def. Partition of unity for locally compact (S,d)} \textbf{(Partition
of unity of locally compact metric space).} Let $(S,d)$ be a locally
compact metric space, with a reference point $x_{\circ}\in S$. Let
the nondecreasing sequence $\xi\equiv(A_{n})_{n=1,2,\cdots}$ of enumerated
finite subsets of $(S,d)$ be a binary approximation of $(S,d)$ relative
to $x_{\circ}$.

For each $n\geq1$, let $\{g_{n,x}:x\in A_{n}\}$ be the $2^{-n}$\emph{-}partition
of unity of $(S,d)$ determined by $A_{n}$. Then the sequence 
\[
\pi\equiv(\{g_{n,x}:x\in A_{n}\})_{n=1,2,\cdots}
\]
is called a \emph{partition of unity}\index{partition of unity of locally compact metric space}
of $(S,d)$ determined by the binary approximation $\xi$. $\square$
\end{defn}
\begin{prop}
\label{Prop. Properties of parittion of unity-1} \textbf{\emph{(Properties
of partition of unity).}} Let $\xi\equiv(A_{n})_{n=1,2,\cdots}$ be
a binary approximation of the locally compact metric space $(S,d)$
relative to a reference point $x_{\circ}$. Let $\pi\equiv(\{g_{n,x}:x\in A_{n}\})_{n=1,2,\cdots}$
be the partition of unity determined by $\xi$. Let $n\geq1$ be arbitrary.
Then the following holds.

1. $g_{n,x}\in C(S)$ has values in $[0,1]$ and has support $(d(\cdot,x)\leq2^{-n+1})$,
for each $x\in A_{n}$.

2. $\sum_{x\in A(n)}g_{n,x}\leq1$ on $S$.

3. $\sum_{x\in A(n)}g_{n,x}=1$ on $\bigcup_{x\in A(n)}(d(\cdot,x)\leq2^{-n})$

4. For each $x\in A_{n}$, the functions $g_{n,x}$, $\sum_{y\in A(n);y<x}g_{n,y}$,
and $\sum_{y\in A(n)}g_{n,y}$ have Lipschitz constant $2^{n+1}$.

5. For each $x\in A_{n}$, 
\begin{equation}
g_{n,x}=\sum_{y\in A(n+1)}g_{n,x}g_{n+1,y}\label{eq:temp-504-1}
\end{equation}
on $S$.
\end{prop}
\begin{proof}
Assertions 1-4 are restatements of their counterparts in Proposition
\ref{Prop. Properties of  epsilon partition of unity } for the case
$\varepsilon\equiv2^{-n}$.

5. Now let $x\in A_{n}$ be arbitrary. By Assertion 1,
\[
(g_{n,x}>0)\subset(d(\cdot,x)\leq2^{-n+1}).
\]
At the same time
\[
(d(\cdot,x)\leq2^{-n+1})\subset(d(\cdot,x_{\circ})\leq2^{n+1})
\]
\[
\subset\bigcup_{y\in A(n+1)}(d(\cdot,y)\leq2^{-n-1})\subset(\sum_{y\in A(n+1)}g_{n+1,y}=1)
\]
where the first inclusion is by relation \ref{eq:temp-502}, the second
by relation \ref{eq:temp-505} applied to $n+1$, and the third by
Assertion 3 applied to $n+1$. Combining, 
\[
(g_{n,x}>0)\subset(\sum_{y\in A(n+1)}g_{n+1,y}=1).
\]
The desired equality \ref{eq:temp-504-1} in Assert 5 follows. 
\end{proof}
\begin{prop}
\label{Prop. Approx  by  interpolation}\emph{ }\textbf{\emph{(Approximation
by interpolation). }}Let $A$ be an arbitrary metrically discrete
enumerated finite subset of a locally compact metric space $(S,d)$.
Let $\varepsilon>0$ be arbitrary. Let $\{g_{x}:x\in A\}$ be an $\varepsilon$-partition
of unity corresponding to $A$.

Let $f\in C(S)$ be arbitrary, with a modulus of continuity $\delta_{f}$
and with $\bigcup_{x\in A}(d(\cdot,x)\leq\varepsilon)$ as support.
Let $\alpha>0$ be arbitrary. Suppose $\varepsilon<\frac{1}{2}\delta_{f}(\frac{1}{3}\alpha)$.
Then

\emph{
\begin{equation}
\left\Vert f-\sum_{x\in A}f(x)g_{x}\right\Vert \leq\alpha\label{eq:temp-151}
\end{equation}
on $S$.}
\end{prop}
\begin{proof}
For abbreviation, write $h\equiv\sum_{x\in A}f(x)g_{x}$. Let $y\in S$
be arbitrary. 

1. Suppose $g_{x}(y)>0$ for some $x\in A$. Since $g_{x}$ has $(d(\cdot,x)<2\varepsilon)$
as support, it follows that $d(y,x)<2\varepsilon<\delta_{f}(\frac{1}{3}\alpha)$.
Hence 
\begin{equation}
|f(y)-f(x)|g_{x}(y)<\frac{1}{3}\alpha g_{x}(y).\label{eq:temp-156}
\end{equation}

2. Suppose $|f(y)-f(x)|g_{x}(y)>\frac{1}{3}\alpha g_{x}(y)$ for some
$x\in A$. Then $g_{x}(y)>0$, leading to inequality \ref{eq:temp-156}
by Step 1, a contradiction. Hence 
\begin{equation}
|f(y)-f(x)|g_{x}(y)\leq\frac{1}{3}\alpha g_{x}(y)\label{eq:temp-63}
\end{equation}
for each $x\in A$. 

3. Either $|f(y)|>0$ or $|f(y)|<\frac{1}{3}\alpha$. First suppose
$|f(y)|>0$. Then $y\in\bigcup_{x\in A}(d(\cdot,x)\leq\varepsilon)$
since the latter set supports $f$, by hypothesis. Hence $\sum_{x\in A}g_{x}(y)=1$
by Condition 3 of Definition \ref{Def.. Epsilon-parittion of unity}.
Therefore
\[
|f(y)-h(y)|=|\sum_{x\in A}f(y)g_{x}(y)-\sum_{x\in A}f(x)g_{x}(y)|
\]
\[
\leq\sum_{x\in A}|f(y)-f(x)|g_{x}(y)<\sum_{x\in A}\frac{1}{3}\alpha g_{x}(y)\leq\frac{1}{3}\alpha
\]
where the second inequality follows from inequality \ref{eq:temp-63}.

4. Now suppose $|f(y)|<\frac{1}{3}\alpha$. Then
\[
|f(y)-h(y)|<\frac{1}{3}\alpha+\sum_{x\in A}|f(x)|g_{x}(y).
\]
Suppose the summand corresponding to some $x\in A$ is greater than
$0$. Then $g_{x}(y)>0$. Hence inequality \ref{eq:temp-156} in Step
1 holds. Consequently
\begin{equation}
|f(y)-f(x)|g_{x}(y)<\frac{1}{3}\alpha g_{x}(y).\label{eq:temp-156-1}
\end{equation}
\[
|f(y)-h(y)|<\frac{1}{3}\alpha+\sum_{x\in A}|f(x)|g_{x}(y)
\]
\[
\leq\frac{1}{3}\alpha+\sum_{x\in A}(|f(y)|+\frac{1}{3}\alpha)g_{x}(y)\leq\frac{1}{3}\alpha+\frac{2}{3}\alpha\sum_{x\in A}g_{x}(y)\leq\alpha.
\]

Combining, we see that $|f(y)-h(y)|\leq\alpha$ for arbitrary $x\in S$.
\end{proof}
\begin{prop}
\label{Prop. Approx  by  Lipschitz continuous functios}\emph{ }\textbf{\emph{(Approximation
by Lipschitz continuous function). }}Let $\xi\equiv(A_{n})_{n=1,2,\cdots}$
be a binary approximation of the locally compact metric space $(S,d)$
relative to a reference point $x_{\circ}$. Let $\pi\equiv(\{g_{n,x}:x\in A_{n}\})_{n=1,2,\cdots}$
be the partition of unity determined by $\xi$. Let $f\in C(S)$ be
a arbitrary, with a modulus of continuity $\delta_{f}$, and with
$\left\Vert f\right\Vert \leq1$.

Let $\alpha>0$ be arbitrary. Let $n\geq1$ be so large that\emph{
(i)} $f$ has the set $(d(\cdot,x_{\circ})\leq2^{n})$ as support,
and \emph{(ii)} $2^{-n}<\frac{1}{2}\delta_{f}(\frac{1}{3}\alpha)$.
Then there exists $g\in C(S)$ with Lipschitz constant $2^{n+1}|A_{n}|$,
such that $\left\Vert f-g\right\Vert \leq\alpha$. Specifically, we
can take 
\[
g\equiv\sum_{x\in A(n)}f(x)g_{n,x}.
\]
\end{prop}
\begin{proof}
By the definition of a partition of unity, the set $A_{n}$ is a $2^{-n}$\emph{-}partition
of unity of $(S,d)$. By hypothesis, the function $f\in C(S)$ has
support 
\[
(d(\cdot,x_{\circ})\leq2^{n})\subset\bigcup_{x\in A(n)}(d(\cdot,x)\leq2^{-n}),
\]
where the displayed relation is according to Proposition \ref{Prop. Properties of  epsilon partition of unity }.
At the same time, $2^{-n}<\frac{1}{2}\delta_{f}(\frac{1}{3}\alpha)$
by hypothesis. Hence Proposition \ref{Prop. Approx  by  interpolation}
implies that $\left\Vert f-g\right\Vert \leq\alpha,$ where
\[
g\equiv\sum_{x\in A(n)}f(x)g_{x}\in C(S)
\]
. Again, according to Proposition \ref{Prop. Properties of  epsilon partition of unity },
each of the functions $g_{x}$ in the last sum has Lipschitz constant
$2^{n+1}$, while $f(x)$ is bounded by $1$ by hypothesis. Hence,
using basic properties of Lipschitz constants in Exercise \ref{Ex.  Lipshitz Constants},
we conclude that the function $g$ has Lipschitz constant $\sum_{x\in A(n)}|f(x)|2^{n+1}\leq2^{n+1}|A_{n}|$,
as desired. 
\end{proof}

\section{One-point Compactification }

The infinite product of a locally compact metric space is not necessarily
locally compact, while the infinite product of a compact metric space
remains compact. For that reason, we will find it sometimes useful
to embed a locally compact metric space into a compact metric space
such that, while the metric is not preserved, the continuous functions
are. This is made precise in the present section as a first application
of partitions of unity.

The next definition is essentially from \cite{BishopBridges85}. 
\begin{defn}
\label{Def. One point compactification} \textbf{(One-point compactification).}
A \emph{one-point compactification}\index{one-point compactification}
of a locally compact metric space $(S,d)$ is a metric space $(\overline{S},\overline{d})$
with an element $\triangle$, called the \emph{\index{point at infinity}point
at infinity},  such that the following five conditions hold.

1. $\widetilde{S}\equiv S\cup\{\Delta\}$ is dense in $(\overline{S},\overline{d})$.
Moreover, $\overline{d}\leq1$.

2. For each compact subset $K$ of $(S,d)$, there exists $c>0$ such
that $\overline{d}(x,\Delta)\geq c$ for each $x\in K$. 

3. Let $K$ be an arbitrary compact subset  of $(S,d)$. Let $\varepsilon>0$
be arbitrary. Then there exists $\delta_{K}(\varepsilon)>0$ such
that for each $y\in K$ and $z\in S$ with $\overline{d}(y,z)<\delta_{K}(\varepsilon)$,
we have $d(y,z)<\varepsilon$. In particular,  the identity mapping
$\bar{\iota}:(S,\overline{d})\rightarrow(S,d)$ is uniformly continuous
on each compact subset of $S$.

4. The identity mapping $\iota:(S,d)\rightarrow(S,\overline{d})$,
defined by $\iota(x)\equiv x$ for each $x\in S$, is uniformly continuous
on $(S,d)$. In other words, for each $\varepsilon>0$, there exists
$\delta_{\overline{d}}(\varepsilon)>0$ such that $\overline{d}(x,y)<\varepsilon$
for each $x,y\in S$ with $d(x,y)<\delta_{\overline{d}}(\varepsilon)$.

5. For each $n\geq1$, we have 
\[
(d(\cdot,x_{\circ})>2^{n+1})\subset(\overline{d}(\cdot,\Delta)\leq2^{-n}).
\]
Thus, as a point $x\in S$ moves away from $x_{\circ}$ relative to
$d$, it converges to the point $\Delta$ at infinity relative to
$\overline{d}$. $\square$
\end{defn}
The next proposition clarifies the relation between continuous functions
on $(S,d)$ and continuous functions on $(\overline{S},\overline{d})$.
First some notations.
\begin{defn}
\label{Def. Restriction of a family of functions)} \textbf{(Restriction
of a family of functions). }Let $A,A'$ be arbitrary sets and let
$B$ be an arbitrary subset of $A$. Recall that the restriction of
a function $f:A\rightarrow A'$ to a subset $B\subset A$ is denoted
by $f|B$. Suppose $F$ is a family of functions from $A$ to $A'$
and suppose $B\subset A$. Then we call the family 
\[
F|B\equiv\{f|B:f\in F\}
\]
the \index{restriction of a family of functions} \emph{restriction
of $F$ to} $B$. 
\end{defn}
$\square$

Recall that $C_{ub}(S,d)$ denotes the space of bounded and uniformly
continuous functions on a locally compact metric space $(S,d)$. 
\begin{prop}
\label{Prop.Relationof Conintuous functions on (S,d)  to those on    (Sbar,dbar)}
\textbf{\emph{(Continuous functions on $(S,d)$ and continuous functions
on $(\overline{S},\overline{d})$).}} Let $(S,d)$ be a locally compact
metric space, with a fixed reference point $x_{\circ}\in S$. Let
$(\overline{S},\overline{d})$ be a one-point compactification of
$(S,d)$. Then the following holds.

1. Each compact subset $K$ of $(S,d)$ is also a compact subset of
$(\overline{S},\overline{d})$. 

2. $C(S,d)\subset C(\overline{S},\overline{d})|S\subset C_{ub}(S,d)$.
Moreover, if $f\in C(\overline{S},\overline{d})$ has a modulus of
continuity $\overline{\delta}$, then $\bar{f}|S\in C_{ub}(S,d)$
has the same modulus of continuity $\overline{\delta}$. 
\end{prop}
\begin{proof}
1. Suppose $K$ is a compact subset of $(S,d)$. By Conditions 3 and
4 of Definition \ref{Def. One point compactification}, the identity
mapping $\iota:(K,d)\rightarrow(K,\overline{d})$ and its inverse
$\bar{\iota}:(K,\overline{d})\rightarrow(K,d)$ are uniformly continuous.
Hence, since by assumption $(K,d)$ is compact, so is $(K,\overline{d})$.

2. First consider each $f\in C(S,d)$. Let the compact subset $K$
of $(S,d)$ be a support of $f$. Extend $f$ to a function $\tilde{f}$
on $\widetilde{S}\equiv S\cup\{\Delta\}$ by defining $\tilde{f}(\triangle)\equiv0$
and $\tilde{f}(x)\equiv x$ for each $x\in S$. We will show that
$\tilde{f}$ is uniformly continuous on $(\widetilde{S},\overline{d})$.
To that end, let $\varepsilon>0$ be arbitrary. Let $\delta>0$ be
so small that $|\tilde{f}(x)-\tilde{f}(y)|<\varepsilon$ for each
$x,y\in S$ with $d(x,y)<\delta$. Then, by Condition 2 in Definition
\ref{Def. One point compactification}, we have $\overline{\delta}\equiv\delta\wedge\overline{d}(K,\Delta)>0$.
Now consider each $x,y\in\widetilde{S}$ with $\overline{d}(x,y)<\overline{\delta}$.
Suppose, for the sake of a contradiction, that $|\tilde{f}(x)-\tilde{f}(y)|>\varepsilon$.
Either (i) $x=\Delta$ or (ii) $x\in S$. Consider case (i). Then
$\tilde{f}(x)=0$. Hence $|\tilde{f}(y)|>0$. Therefore $y\in S$
and $|f(y)|\equiv|\tilde{f}(y)|>0$. Since $K$ is a support of $f$
we see that $y\in K$. Combining, $\overline{d}(y,x)\geq\overline{d}(K,\Delta)\geq\overline{\delta}$,
a contradiction. Thus $x\in S$. Similarly $y\in S$. Therefore, by
the definition of $\delta$, we have $|f(x)-f(y)|<\varepsilon$, again
a contradiction. Summing up, we see that $|f(x)-f(y)|\leq\varepsilon$.
Since $\varepsilon>0$ and $x,y\in\widetilde{S}$ with $\overline{d}(x,y)<\overline{\delta}$
are arbitrary, $\tilde{f}$ is a uniformly continuous function on
$(\widetilde{S},\overline{d})$. As such $\tilde{f}$ can be extended
by continuity to a function $\bar{f}\in C(\overline{S},\overline{d})$,
thanks to the denseness of $\widetilde{S}$ in $(\overline{S},\overline{d})$.
Since $f\in C(S,d)$ is arbitrary and since $f=\bar{f}|S$, we have
proved that $C(S,d)\subset C(\overline{S},\overline{d})|S$.

Now consider each $\bar{f}\in C(\overline{S},\overline{d})$, with
a modulus of continuity $\overline{\delta}$. Then $\bar{f}$ is bounded
since $C(\overline{S},\overline{d})$ is compact. Let $\varepsilon>0$
and $x,y\in S$ be arbitrary with $d(x,y)<\overline{\delta}(\varepsilon)$.
Then, by condition 4 in Definition \ref{Def. One point compactification},
we have $\overline{d}(x,y)\leq1\wedge d(x,y)<\overline{\delta}(\varepsilon)$.
Hence $|\bar{f}(x)-\bar{f}(y)|<\varepsilon$. Since $\varepsilon>0$
is arbitrary, we conclude that $\bar{f}|S\in C_{ub}(S,d)$, with modulus
of continuity also given by $\overline{\delta}$. Summing up, we have
proved that $C(\overline{S},\overline{d})|S\subset C_{ub}(S,d)$.
\end{proof}
The next theorem constructs a one-point compactification. The proof
follows the lines of Theorem 6.8 in Chapter 4 of \cite{BishopBridges85}.

\begin{thm}
\label{Thmf. Construction of one-point compactificaton from separant}
\textbf{\emph{(Construction of a one-point compactification from a
binary approximation).}} Let $(S,d)$ be a locally compact metric
space. Let the sequence $\xi\equiv(A_{n})_{n=1,2,\cdots}$ of subsets
be a binary approximation of $(S,d)$ relative to $x_{\circ}$. Then
there exists a one-point compactification $(\overline{S},\overline{d})$
of $(S,d)$, such that the following conditions hold.

(i). For each $p\geq1$ and for each $y,z\in S$ with 
\[
d(y,z)<p^{-1}2^{-p-1},
\]
we have 
\[
\overline{d}(y,z)<2^{-p+1}.
\]

(ii). For each $n\geq1$ and for each $y\in(d(\cdot,x_{\circ})\leq2^{n})$
and for each $z\in S$ with 
\[
\overline{d}(y,z)<2^{-n-1}|A_{n}|^{-2},
\]
we have
\[
d(y,z)<2^{-n+2}.
\]

The one-point compactification $(\overline{S},\overline{d})$ constructed
in the proof is said to be \emph{de}term\emph{ined by the binary approximation}
$\xi$.
\end{thm}
\begin{proof}
Let $\pi\equiv(\{g_{n,x}:x\in A_{n}\})_{n=1,2,\cdots}$ be the partition
of unity of $(S,d)$ determined by $\xi$. Let $n\geq1$ be arbitrary.
Then $\{g_{n,x}:x\in A_{n}\}$ is a $2^{-n}$-partition of unity corresponding
to the metrically discrete and enumerated finite set $A_{n}$. Moreover,
by Proposition \ref{Prop. Properties of parittion of unity-1}, $g_{n,x}$
has Lipschitz constant $2^{n+1}$ for each $x\in A_{n}$.

1. Define 
\[
\widetilde{S}\equiv\{(x,i)\in S\times\{0,1\}:i=0\;or\;(x,i)=(x_{\circ},1)\}.
\]
and define $\Delta\equiv(x_{\circ},1)$. Identify each $x\in S$ with
$\bar{x}\equiv(x,0)\in\widetilde{S}$. Thus $\widetilde{S}=S\cup\{\Delta\}$.
Extend each function $f\in C(S)$ to a function on $\widetilde{S}$
by defining $f(\Delta)\equiv0$. In particular $g_{n,x}(\Delta)\equiv0$
for each $x\in A_{n}$. Define 
\begin{equation}
\overline{d}(y,z)\equiv\sum_{n=1}^{\infty}2^{-n}|A_{n}|^{-1}\sum_{x\in A(n)}|g_{n,x}(y)-g_{n,x}(z)|\label{eq:temp-216}
\end{equation}
for each $y,z\in\widetilde{S}$. Then $\overline{d}(y,y)=0$ for each
$y\in\widetilde{S}$. Symmetry and triangle inequality of the function
$\overline{d}$ are immediate consequences of equality \ref{eq:temp-216}.
Moreover, $\overline{d}\leq1$ since the functions $g_{n,x}$ have
values in $[0,1]$.

2. Let $y\in S$ be arbitrary, and let $K$ be an arbitrary compact
subset of $(S,d)$. Suppose $y\in K$. Let $n\geq1$ be so large that
\[
y\in K\subset(d(\cdot,x_{\circ})\leq2^{n}).
\]
Then 
\[
y\in\bigcup_{x\in A(n)}(d(\cdot,x)\leq2^{-n})\subset(\sum_{x\in A(n)}g_{n,x}=1).
\]
where the membership relation of on the left-hand side is by expression
\ref{eq:temp-505} in Definition \ref{Def. Binary approximationt and Modulus of local compactness},
and where the inclusion on the right-hand side is according to Assertion
3 of Proposition \ref{Def.. Epsilon-parittion of unity}. Hence the
defining equality \ref{eq:temp-216} yields 
\begin{equation}
\overline{d}(y,\Delta)\geq2^{-n}|A_{n}|^{-1}\sum_{x\in A(n)}g_{n,x}(y)=2^{-n}|A_{n}|^{-1},\label{eq:temp-177}
\end{equation}
establishing Condition 2 in Definition \ref{Def. One point compactification}.

3. Let $n\geq1$ be arbitrary. Let $y\in(d(\cdot,x_{\circ})\leq2^{n})$
and $z\in S$ be arbitrary such that
\[
\overline{d}(y,z)<\delta_{\xi,n}\equiv2^{-n-1}|A_{n}|^{-2}.
\]
 As seen in Step 2, 
\[
\sum_{x\in A(n)}g_{n,x}(y)=1.
\]
Hence there exists $x\in A_{n}$ such that
\begin{equation}
g_{n,x}(y)>\frac{1}{2}|A_{n}|^{-1}>0.\label{eq:temp-217}
\end{equation}
At the same time,
\[
|g_{n,x}(y)-g_{n,x}(z)|\leq\sum_{u\in A(n)}|g_{n,u}(y)-g_{n,u}(z)|
\]
\[
\leq2^{n}|A_{n}|\overline{d}(y,z)<2^{n}|A_{n}|\delta_{\xi,n}\equiv\frac{1}{2}|A_{n}|^{-1}.
\]
Hence inequality \ref{eq:temp-217} implies that $g_{n,x}(z)>0$.
Consequently, $y,z\in(d(\cdot,x)<2^{-n+1})$. Thus $d(y,z)<2^{-n+2}$.
This establishes Assertion (ii) of the theorem.

Now let $K$ be an arbitrary compact subset of $(S,d)$ and let $\varepsilon>0$
be arbitrary. Let $n\geq1$ be so large that $K\subset(d(\cdot,x_{\circ})\leq2^{n})$
and that $2^{-n+2}<\varepsilon$. Let $\delta_{K}(\varepsilon)\equiv\delta_{\xi,n}$.
Then, by the preceding paragraph, for each $y\in K$ and $z\in S$
with $\overline{d}(y,z)<\delta_{K}(\varepsilon)\equiv\delta_{\xi,n}$,
we have $d(y,z)<\varepsilon$. Condition 3 in Definition \ref{Def. One point compactification}
has been verified. 

In particular, suppose $y,z\in\widetilde{S}$ are such that $\overline{d}(y,z)=0$.
Then either $y=z=\Delta$ or $y,z\in S$, in view of inequality \ref{eq:temp-177}.
Suppose $y,z\in S$. Then the preceding paragraph applied to the compact
set $K\equiv\{y,z\}$, implies that $d(y,z)=0$. Since $(S,d)$ is
a metric space, we conclude that $y=z$. In view of the last paragraph
of Step 1 above, $(\widetilde{S},\overline{d})$ is a metric space.

4. Recall that $g_{n,x}$ has values in $[0,1]$, and, as remarked
above, has Lipschitz constant $2^{n+1}$, for each $x\in A_{n}$,
for each $n\geq1$. Let $p\geq2$ be arbitrary. Let $y,z\in S$ be
such that $d(y,z)<p^{-1}2^{-p-1}$. Then 
\[
\overline{d}(y,z)\equiv\sum_{n=1}^{\infty}2^{-n}|A_{n}|^{-1}\sum_{x\in A(n)}|g_{n,x}(y)-g_{n,x}(z)|
\]
\[
\leq\sum_{n=1}^{p}2^{-n}2^{n+1}d(y,z)+2^{-p}
\]
\begin{equation}
<p2p^{-1}2^{-p-1}+2^{-p}=2^{-p}+2^{-p}=2^{-p+1}.\label{eq:temp-270}
\end{equation}
Since $2^{-p+1}$ is arbitrarily small, we see that the identity mapping
$\iota:(S,d)\rightarrow(S,\overline{d})$ is uniformly continuous.
This establishes Condition 4 in Definition \ref{Def. One point compactification}. 

5. Let $n\geq1$ be arbitrary. Consider each $y\in(d(\cdot,x_{\circ})>2^{n+1}).$
Let $m\geq n$ be arbitrary. Then 
\[
y\in(d(\cdot,x_{\circ})>2^{m+1})\subset\bigcap_{x\in A(m)}(d(\cdot,x)\geq2^{-m+1})
\]
by relation \ref{eq:temp-502} in Definition \ref{Def. Binary approximationt and Modulus of local compactness}
of a binary approximation. For each $x\in A_{m}$, since $g_{m,x}$
has support $(d(\cdot,x)\geq2^{-m+1})$, we infer $g_{m,x}(y)=0$.
Hence the defining equality \ref{eq:temp-216} reduces to
\[
\overline{d}(y,\triangle)\equiv\sum_{m=1}^{\infty}2^{-m}|A_{m}|^{-1}\sum_{x\in A(m)}g_{m,x}(y)
\]
\[
=\sum_{m=n+1}^{\infty}2^{-m}|A_{m}|^{-1}\sum_{x\in A(m)}g_{m,x}(y)
\]
\begin{equation}
\leq\sum_{m=n+1}^{\infty}2^{-m}=2^{-n}.\label{eq:temp-216-1}
\end{equation}
Since $y\in(d(\cdot,x_{\circ})>2^{n+1})$ is arbitrary, we conclude
that 
\begin{equation}
(d(\cdot,x_{\circ})>2^{n+1})\subset(\overline{d}(\cdot,\Delta)\leq2^{-n}).\label{eq:temp-235}
\end{equation}
This proves Condition 5 in Definition \ref{Def. One point compactification}. 

6. We will prove next that $(\widetilde{S},\overline{d})$ is totally
bounded. To that end, let $p\geq1$ be arbitrary. Let $m\equiv m_{p}\equiv[(p+2)+\log_{2}p]_{1}$.
Recall here that $[\cdot]_{1}$ is the operation which assigns to
each $a\in[0,\infty)$ and integer $[a]_{1}$ in $(a,a+2)$. Then
\[
2^{-m}<\overline{\delta}_{p}\equiv p^{-1}2^{-p-2}.
\]
Note that
\[
\widetilde{S}\equiv S\cup\{\Delta\}\subset(d(\cdot,x_{\circ})<2^{m})\cup(d(\cdot,x_{\circ})>2^{m-1})\cup\{\Delta\}
\]
\[
\subset\bigcup_{x\in A(m)}(d(\cdot,x)\leq2^{-m})\cup(\overline{d}(\cdot,\Delta)\leq2^{-m+2})\cup\{\Delta\}
\]
where the second inclusion is due to relation \ref{eq:temp-505},
and to relation \ref{eq:temp-235} applied to $m-2$. Continuing,
\[
\widetilde{S}\subset\bigcup_{x\in A(m)}(d(\cdot,x)\leq\overline{\delta}_{p})\cup(\overline{d}(\cdot,\Delta)<p^{-1}2^{-p})\cup\{\Delta\}.
\]
\[
\subset\bigcup_{x\in A(m)}(\overline{d}(\cdot,x)<2^{-p})\cup(\overline{d}(\cdot,\Delta)<2^{-p})\cup\{\Delta\},
\]
thanks to the inequality \ref{eq:temp-270} in Step 4. Consequently,
the set 
\[
\overline{A}_{p}\equiv A_{m(p)}\cup\{\Delta\}
\]
is a metrically discrete $2^{-p}$-approximation of $(\widetilde{S},\overline{d})$.
Since $2^{-p}$ is arbitrarily small, the metric space $(\widetilde{S},\overline{d})$
is totally bounded. Hence its completion $(\overline{S},\overline{d})$
is compact, and $\widetilde{S}$ is dense in $(\overline{S},\overline{d})$,
proving Condition 1 in Definition \ref{Def. One point compactification}.
Note that, since $\widetilde{S}\equiv S\cup\{\Delta\}$ is a dense
subset of $(\overline{S},\overline{d})$, the sequence $\overline{A}_{p}$
is a $2^{-p}$-approximation of $(\overline{S},\overline{d})$.

Summing up, $(\overline{S},\overline{d})$ satisfies all the conditions
in Definition \ref{Def. One point compactification} to be a one-point
compactification of $(S,d)$.
\end{proof}
Proposition \ref{Prop.Relationof Conintuous functions on (S,d)  to those on    (Sbar,dbar)}
established the relation of continuity on $(S^{n},d^{n})$ to continuity
on $C(\overline{S}^{n},\overline{d}^{n})$ in the case $n=1$. The
next lemma generalizes to the case where $n\geq1$. 
\begin{cor}
\emph{\label{Cor. Extension of f in C(S**n) to point of infinity}}
\textbf{\emph{(}}\textbf{\textup{Extension}}\textbf{\emph{ of each
$f\in C(S^{n},d^{n})$ to $(\overline{S}^{n},\overline{d}^{n})$).}}
Let $n\geq1$ be arbitrary. Then
\[
C(S^{n},d^{n})\subset C(\overline{S}^{n},\overline{d}^{n})|S^{n}\subset C_{ub}(S^{n},d^{n}).
\]
\end{cor}
\begin{proof}
1. Let $h\in C(S^{n},d^{n})$ be arbitrary with a modulus of continuity
$\delta_{h}$.\emph{ }Then there exists\emph{ }$r>0$ with such that
$K_{r}\equiv(d(x_{\circ},\cdot)\leq r)$ is compact in $(S,d)$, and
such that $K_{r}^{n}$ is a support of $h$. Let $s>1$ be such that
$K\equiv(d(K_{r},\cdot)\leq s)$\emph{ }is compact in $(S,d)$. Then
$K_{r},K$ are compact subsets of $C(\overline{S},\overline{d})$,
according to Proposition \ref{Prop.Relationof Conintuous functions on (S,d)  to those on    (Sbar,dbar)}.
By Definition \ref{Def. One point compactification}, for each $\varepsilon>0$,
there exists $\delta_{K}(\varepsilon)\in(0,1)$ such that, for each
$x,y\in K$ with $\overline{d}(x,y)<\delta_{K}(\varepsilon)$, we
have 
\begin{equation}
d(x,y)<\varepsilon.\label{eq:temp-174}
\end{equation}

Now let $\varepsilon'\in(0,1)$ be arbitrary. Write $\varepsilon\equiv1\wedge\delta_{h}(\varepsilon')$
and define $\overline{\delta}_{K}(\varepsilon')\equiv\delta_{K}(\varepsilon)$
Let $u\equiv(x_{1},\cdots,x_{n}),v\equiv(y_{1},\cdots,y_{n})\in\overline{S}^{n}$
be arbitrary such that 
\begin{equation}
\overline{d}^{n}(u,v)\equiv\bigvee_{i=1}^{n}\overline{d}(x_{i},y_{i})<\overline{\delta}_{K}(\varepsilon')\equiv\delta_{K}(\varepsilon)\equiv\delta_{K}(1\wedge\delta_{h}(\varepsilon')).\label{eq:temp-184-1}
\end{equation}
We will prove that 
\[
|h(u)-h(v)|\leq\varepsilon'.
\]
First note that, by inequality \ref{eq:temp-184-1}, we have 
\begin{equation}
d(x_{i},y_{i})<\varepsilon\equiv1\wedge\delta_{h}(\varepsilon')\label{eq:temp-158}
\end{equation}
for each $i=1,\cdots,n$. Suppose, for the sake of a contradiction,
that 
\begin{equation}
|h(u)-h(v)|>\varepsilon'.\label{eq:temp-181}
\end{equation}
Then $h(u)>0$ or $h(v)>0$. Suppose $h(u)>0$. Then $(x_{1},\cdots,x_{n})\equiv u\in K_{r}^{n}$
since $K_{r}^{n}$ contains a support of $h$. Let $i=1,\cdots,n$
be arbitrary. Then $x_{i}\in K_{r}$, whence, by inequality \ref{eq:temp-158},
\[
y_{i}\in(d(\cdot,K_{r})\leq1)\subset(d(\cdot,K_{r})\leq s)\equiv K.
\]
Thus $x_{i},y_{i}\in K$. At the same time, $\overline{d}(x_{i},y_{i})<\delta_{K}(\varepsilon)$
by inequality \ref{eq:temp-184-1}. Consequently, inequality \ref{eq:temp-174}
holds for $x_{i},y_{i}$. Combining,
\[
d^{n}(u,v)\equiv\bigvee_{i=1}^{n}d(x_{i},y_{i})<\varepsilon\leq\delta_{h}(\varepsilon').
\]
Since $\delta_{h}$ is a modulus of continuity of $h\in C(S^{n},d^{n})$,
it follows that 
\begin{equation}
|h(u)-h(v)|<\varepsilon',\label{eq:temp-265}
\end{equation}
a contradiction to inequality \ref{eq:temp-181}. Similarly, the assumption
$h(v)>0$ also leads to a contradiction. Summing up, the assumption
of inequality \ref{eq:temp-181} leads to a contradiction. Hence 
\[
|h(u)-h(v)|\leq\varepsilon',
\]
where $u,v\in\overline{S}^{n}$ are arbitrary with $\overline{d}^{n}(u,v)<\overline{\delta}_{K}(\varepsilon')$.
In other words, $h$ is uniformly continuous on $(\overline{S}^{n},\overline{d}^{n})$,
with modulus of continuity $\overline{\delta}_{K}$. 

2. Conversely, let $\bar{h}\in C(\overline{S}^{n},\overline{d}^{n})$
be arbitrary. By Definition \ref{Def. One point compactification}
of the compactification, the identity mapping $\iota:(S,d)\rightarrow(S,\overline{d})$
is uniformly continuous. Hence so is the identity mapping $\iota^{n}:(S^{n},d^{n})\rightarrow(S^{n},\overline{d}^{n})$.
Therefore $\bar{h}|S^{n}=\bar{h}\circ\iota^{n}$ is bounded and uniformly
continuous on $(S^{n},d^{n})$.
\end{proof}
\begin{cor}
\label{Cor. Compactification of binary approximation} \textbf{\emph{(Compactification
of binary approximation). }}Use the same notations and assumptions
as in Theorem \ref{Thmf. Construction of one-point compactificaton from separant}.
In particular, let $\xi\equiv(A_{n})_{n=1,2,\cdots}$ be a binary
approximation of $(S,d)$ relative to the reference point $x_{\circ}$.
For each $n\geq1$, let $A_{n}\equiv\{x_{n,1},\cdots,x_{n,\kappa(n)}\}$.
Thus $\left\Vert \xi\right\Vert \equiv(|A_{n}|)_{n=1.2.\cdots}=(|\kappa_{n}|)_{n=1.2.\cdots}$.

Let $p\geq1$ be arbitrary. Write $m_{p}\equiv[(p+2)+\log_{2}p]_{1}$.
Define 
\[
\overline{A}_{p}\equiv A_{m(p)}\cup\{\Delta\}\equiv\{x_{m(p),1},\cdots,x_{m(p),\kappa(m(p))},\Delta\}.
\]
Then $\overline{\xi}\equiv(\overline{A}_{p})_{p=1,2,\cdots}$ is a
binary approximation of $(\overline{S},\overline{d})$ relative to
$x_{\circ}$, called the \emph{\index{compactification of a binary approximation}
compactification} of $\xi$. Thus the corresponding modulus of local
compactness of $(\overline{S},\overline{d})$ is given by
\[
\left\Vert \overline{\xi}\right\Vert \equiv(|\overline{A}_{p}|)_{p=1,2,\ldots}=(\kappa_{m(p)}+1)_{p=1,2,\ldots}=(|A_{m(p)}|+1)_{p=1,2,\ldots}
\]
and is therefore determined by $\left\Vert \xi\right\Vert $.
\end{cor}
\begin{proof}
Let $p\geq1$ be arbitrary. According to Step 6 of the proof of Theorem
\ref{Thmf. Construction of one-point compactificaton from separant},
the finite set $\overline{A}_{p}$ is a metrically discrete $2^{-p}$-approximation
of $(\overline{S},\overline{d})$. Hence
\[
(\overline{d}(\cdot,x_{\circ})\leq2^{p})\subset\overline{S}\subset\sum_{x\in\overline{A}(p)}(\overline{d}(\cdot,x)\leq2^{-p}).
\]
At the same time, Condition 1 of Definition \ref{Def. One point compactification}
says that $\overline{d}\leq1$. Hence
\[
\bigcup_{x\in\overline{A}(p)}(\overline{d}(\cdot,x)\leq2^{-p+1})\subset\overline{S}\subset(\overline{d}(\cdot,x_{\circ})\leq1)\subset(\overline{d}(\cdot,x_{\circ})\leq2^{p+1}).
\]
Thus all the conditions in Definition \ref{Def. Binary approximationt and Modulus of local compactness}
have been verified for $\overline{\xi}\equiv(\overline{A}_{p})_{p=1,2,\cdots}$
to be a binary approximation of $(\overline{S},\overline{d})$ relative
to $x_{\circ}$. 
\end{proof}

\chapter{Integration and Measure}

We introduce next the Riemann-Stieljes integral on $R$. Then we give
a general treatment of integration- and measure theory in terms of
Daniell integrals, adapted from \cite{BishopBridges85}. The standard
graduate course in measure theory usually starts with a chapter of
measurable sets, before defining a measure. In contrast, the Daniell
integration theory starts with the integration and the integrable
functions. Thus we discuss the computation of the integration early
on. We remark that it is possible to adhere to the traditional approach
of starting with measurable sets. (See \cite{BishopCheng72}). However,
Daniell integrals are more natural, and cleaner, in the present context.

\section{The Riemann-Stieljes Integral}
\begin{defn}
\label{Def. Distribution Func} \textbf{(Distribution function).}
A \emph{distribution function}\index{distribution function} is a
nondecreasing real-valued function $F$ whose $domain(F$) is dense
in $R$. $\square$
\end{defn}
Let $F$ be a distribution function, and let $X\in C(R)$. 

By a \emph{partition}\index{partition} of $R$ we mean a finite and
increasing sequence ($x_{0}$,$\cdots,x_{n})$ in $domain(F)$. One
partition is said to be a \emph{refinement}\index{refinement of a partition}
of another if the former contains the latter as a subsequence. For
any partition ($x_{1}$,$\cdots,x_{n}$), define its \emph{mesh}\index{mesh}
as $\bigvee_{i=1}^{n}(x_{i}-x_{i-1})$ and define the \emph{Riemann-Stieljes
sum}\index{Riemann-Stieljes sum} as

\[
S(x_{0},\cdots,x_{n})\equiv\sum_{i=1}^{n}X(x_{i})(F(x_{i})-F(x_{i-1}))
\]

\begin{thm}
\label{Thm.  Existence of RS integral} \textbf{\emph{(Exisence of
Riemann-Stieljes integral). }}For any $X\in C(R)$, the Riemann-Stieljes
sum converges as the mesh of the partition ($x_{0}$,$\cdots,x_{n}$)
approaches 0 with $x_{0}\rightarrow-\infty$ and $x_{n}\rightarrow+\infty$
. The limit will be called the \emph{Riemann-Stieljes integral}\index{Riemann-Stieljes integral}
of $X$ with respect to the function $F$, and will be denoted by
$\int_{-\infty}^{+\infty}X(x)dF(x)$, or more simply by $\int X(x)dF(x)$.
\end{thm}
\begin{proof}
Suppose $X$ vanishes outside the compact interval $[a,b]$ where
$a,b\in domain(F)$. Let $\varepsilon>0$. Consider a partition ($x_{0}$,$\cdots,x_{n}$)
with (i) $x_{0}<a-2<b+2<x_{n}$ and (ii) it has mesh less than $1\wedge\delta_{X}(\varepsilon)$
where $\delta_{X}$ is a modulus of continuity for $X$. 

Let $i$ be any index with $0<i\leq n$. Suppose we insert $m$ points
between $(x_{i-1},x_{i})$ and make a refinement $(\cdots,x_{i-1},y_{1},\cdots,y_{m-1},x_{i},\cdots)$.
Let $y_{0}$ and $y_{m}$ denote $x_{i-1}$ and $x_{i}$ respectively.
Then the difference in Riemann-Stieljes sums for the new and old partitions
is bounded by 
\[
|X(x_{i})(F(x_{i})-F(x_{i-1}))-\sum_{j=1}^{m}X(y_{j})(F(y_{j})-F(y_{j-1})|
\]
\[
=|\sum_{j=1}^{m}(X(x_{i})-X(y_{j}))(F(y_{j})-F(y_{j-1})|
\]
\[
\leq|\sum_{j=1}^{m}\varepsilon(F(y_{j})-F(y_{j-1})|=\varepsilon(F(x_{i})-F(x_{i-1})
\]
Moreover, the difference is $0$ if $x_{i}<a-2$ or $x_{i-1}>b$+2.
Since $x_{i}-x_{i-1}<1$, the difference is 0 if $x_{i-1}<a-1$ or
$x_{i}>b+1$.

Since any refinement of ($x_{0}$,$\cdots,x_{n}$) can be obtained
by inserting points between the pairs $(x_{i-1},x_{i})$, we see that
the Riemann-Stieljes sum of any refinement differs from that for ($x_{0}$,$\cdots,x_{n}$)
by at most $\sum\varepsilon(F(x_{i})-F(x_{i-1}))$ where the sum is
over all $i$ for which $a<x_{i-1}$ and $x_{i}<b$. The difference
is therefore at most $\varepsilon(F(b)-F(a))$.

Consider a second partition ($u_{0}$,$\cdots,u_{p}$) satisfying
the conditions (i) and (ii). Because the domain of $F$ is dense,
we can find a third partition ($v_{0}$,$\cdots,v_{q}$) satisfying
the same conditions and the additional condition that $|v_{k}-x_{i}|>0$
and $|v_{k}-u_{j}|>0$ for all $i,j,k$. Then ($v_{0}$,$\cdots,v_{q}$)
and ($x_{0}$,$\cdots,x_{n}$) have a common refinement, namely the
merged sequence rearranged in increasing order. So their Riemann-Stieljes
sums differ from each other by at most $2\varepsilon(F(b)-F(a))$
by the first part of this proof. Similarly, the Riemann-Stieljes sum
for ($u_{0}$,$\cdots,u_{p}$) differs from that of ($v_{0}$,$\cdots,v_{q}$)
by at most $2\varepsilon(F(b)-F(a))$. Hence the Riemann-Stieljes
sums for ($u_{0}$,$\cdots,u_{p}$) and ($x_{0}$,$\cdots,x_{n}$)
differ by at most $4\varepsilon(F(b)-F(a)).$ 

Since $\varepsilon$ is arbitrary, the asserted convergence is proved. 
\end{proof}
\begin{thm}
\label{Thm. Basics of he-Riemann-Stieljes integral} \textbf{\emph{(Basic
properties of the Riemann-Stieljes integral). }}The Riemann-Stieljes
integral is linear on $C(R)$. It is also positive: if $\int X(x)dF(x)>0$,
then there exists $x\in R$ such that $X(x)>0.$
\end{thm}
\begin{proof}
Linearity follows trivially from the defining formulas. Suppose $a,b\in domain(F)$
are such that $X$ vanishes outside $[a,b]$. If the integral is greater
than some positive number $c$, then the Riemann-Stieljes sum $S(x_{0},\cdots,x_{n})$
for some partition with $x_{1}=a$ and $x_{n}=b$ is greater than
$c$. If follows that $X(x_{i})(F(x_{i})-F(x_{i-1}))$ is greater
than or equal to $c/n$ for some index $i$. Hence $X(x_{i})(F(b)-F(a))\geq X(x_{i})(F(x_{i})-F(x_{i-1}))\geq c/n$.
This implies $X(x_{i})>c/(n(F(b)-F(a))>0$.
\end{proof}
In the special case where $domain(F)=R$ and $F(x)=x$ for each $x\in R$,
the Riemann-Stieljes sums and Riemann-Stieljes integral are called
the Riemann sums and the Riemann integral respectively. 

\section{Integration on Locally Compact Metric Spaces}

In this section, the Riemann-Stieljes integration is generalized to
a locally compact metric space\emph{ }$(S,d)$. 

Classically, integration is usually defined in terms of a measure,
a function on a family of subsets which is closed relative to the
operations of countable unions, countable intersections, and relative
complements. In the case of a metric space, one such family can be
generated via these three operations from the family of all open subsets.
Members of the family thus generated are called Borel sets. In the
special case of $R$, the open sets can in turn be generated from
a countable subfamily of intervals in successive partitions of $R$,
wherein ever smaller intervals cover any compact interval in $R$.
The intervals in the countable family can thus serve as building blocks
in the analysis of measures on $R$.

The Daniell integration theory is a more natural choice for the constructive
development. Integrals of functions, rather than measures of sets,
are the starting point. In the special case of a locally compact metric
space $(S,d)$, the family $C(S)$ supplies the basic integrable functions.
The family $C(S)$ can be generated, via linear operations and uniform
convergence, from a countable subfamily obtained from successive partitions
of the unit function $1$ by non-negative members of $C(S)$, wherein
members with ever smaller compact supports sum to 1 on any given compact
subset in $S$. The functions in this countable subfamily can then
serve as building blocks in the analysis of integrations on $S$. 
\begin{defn}
\label{Def. integration on loc compact space} \textbf{(Integration
on a locally compact metric space).} An \emph{integration} on a locally
compact metric space\index{integration on locally compact metric space}
$(S,d)$ is a real-valued linear function $I$ on the linear space
$C(S)$ such that (i) $I(X)>0$ for some $X\in C(S)$, and (ii) for
each $X\in C(S)$ with $I(X)>0$ there exists a point $x$ in $S$
for which $X(x)>0.$ $\square$
\end{defn}
The Riemann-Stieljes integration defined for a distribution function
$F$ on $R$ is an integration on $(R,d)$ where $d$ is the Euclidean
metric, and is denoted by $\int\cdot dF$, with $I(X)$ written as
$\int X(x)dF(x)$ for each $X\in C(S)$. Riemann-Stieljes integrals
provide an abundance of examples for integration on locally compact
metric spaces. 

It follows from the linearity of $I$ that if $X,Y\in C(S)$ are such
that $I(X)>I(Y)$, then there exists a point $x$ in $S$ for which
$X(x)>Y(x).$ The positivity condition (ii), extended in the next
proposition, is a powerful tool in proving existence theorems. It
translates a condition on integrals into the existence of a point
in $S$ with certain properties. To prove the next proposition, we
need the following lemma which will be used again in a later chapter.
This lemma, from \cite{Chan75}  is a pleasant surprise because, in
general, the convergence of a series of non-negative real numbers
does not follow constructively from the boundlessness of partial sums.
\begin{lem}
\label{Lem. Chan's lemma}\emph{ }\textbf{\emph{(Positivity of a linear
function on a linear space of functions).}} Suppose $I$ is a linear
function on a linear space $L$ of functions on a set $S$. Suppose
$I$ satisfies the following condition: for each $X_{0}\in L$ there
exists a non-negative function $Z\in L$ such that, for each sequence
$(X_{i})_{i=1,2,\cdots}$ of non-negative functions in $L$ with $\sum_{i=1}^{\infty}I(X_{i})<I(X_{0})$,
there exists $x\in S$ with (i) $Z(x)=1$ and (ii) $\sum_{i=1}^{p}X_{i}(x)\leq X_{0}(x)$
for each $p>0$. Then, for each $X_{0}\in L$ and for each sequence
$(X_{i})_{i=1,2,\cdots}$ of non-negative functions in $L$ with $\sum_{i=1}^{\infty}I(X_{i})<I(X_{0})$,
there exists $x\in S$ such that $\sum_{i=1}^{\infty}X_{i}(x)$ converges
and is less than $X_{0}(x)$.
\end{lem}
\begin{proof}
Classically, the convergence of $\sum_{i=1}^{\infty}X_{i}(x)$ follows
trivially from the boundlessness of the partial sums. Note that if
the constant function $1$ is a member of $L$, then the lemma can
be simplified with $Z\equiv1$, or with $Z$ altogether omitted.

Suppose $X_{0}\in L$ and $(X_{i})_{i=1,2,\cdots}$ is a sequence
of non-negative functions in $L$ with $\sum_{i=1}^{\infty}I(X_{i})<I(X_{0})$.
Let $Z$ be as given in the hypothesis. Choose a positive real number
$\alpha$ so small that 

\[
\alpha I(Z)+\sum_{i=1}^{\infty}I(X_{i})+\alpha<I(X_{0})
\]
Choose an increasing sequence $(n_{k})_{k=1,2,\cdots}$ of integers
such that 

\[
\sum_{i=n(k)}^{\infty}I(X_{i})<2^{-2k}\alpha
\]
for each $k\geq1$. 

Consider the sequence of functions 

\[
(\alpha Z,X_{1},2\sum_{i=n(1)}^{n(2)}X_{i},X_{2},2^{2}\sum_{i=n(2)}^{n(3)}X_{i},X_{3},\cdots)
\]

It can easily be verified that the series of the corresponding values
for the function $I$ then converges to a sum less than $\alpha I(Z)+\sum_{i=1}^{\infty}I(X_{i})+\alpha$,
which is in turn less than $I(X_{0})$ by the choice of the number
$\alpha$. 

By the hypothesis, there exists a point $x\in S$ with $Z(x)=1$ such
that

\[
\alpha Z(x)+X_{1}(x)+\cdots+X_{k}(x)+2^{k}\sum_{i=n(k)}^{n(k+1)}X_{i}(x)\leq X_{0}(x)
\]
for each $k\geq1$. In particular $\sum_{i=n(k)}^{n(k+1)}X_{i}(x)\leq2^{-k}X_{0}(x)$
so $\sum_{i=1}^{\infty}X_{i}(x)<\infty$. The last displayed inequality
implies also that

\[
\alpha Z(x)+\sum_{i=1}^{\infty}X_{i}(x)\leq X_{0}(x)
\]
Because $Z(x)=1$, we have 

\[
\alpha+\sum_{i=1}^{\infty}X_{i}(x)\leq X_{0}(x)
\]
as desired. 
\end{proof}
\begin{prop}
\label{Prop. Positivity of I on (S,C)} \textbf{\emph{(Positivity
of an integration on a locally compact metric space).}} Let $I$ be
an integration on a locally compact metric space $(S,d)$. Let $(X_{i})_{i=0,1,2,\cdots}$
be a sequence in $C(S)$ such that $X_{i}$ is non-negative for $i\geq1$,
and such that $\sum_{i=1}^{\infty}I(X_{i})<I(X_{0})$. Then there
exists $x\in S$ such that $\sum_{i=1}^{\infty}X_{i}(x)<X_{0}(x)$. 
\end{prop}
\begin{proof}
Let $K$ be a compact support of $X_{0}$. The set $B\equiv\{x\in S:d(x,K)\leq1\}$
is bounded. Hence, since $S$ is locally compact, there exists a compact
subset $K'$ such that $B\subset K'$. Define $Z\equiv(1-d(\cdot,K'))_{+}$. 

Let $\varepsilon\in(0,1)$ be arbitrary. By Lemma \ref{Lem. Existence of metrically discrete epsilon approximations},
there exists a metrically discrete and enumerated finite set $A\equiv\{y_{1},\cdots,y_{n}\}$
which is an $\varepsilon$-approximation of $K$. Let $\{Y_{y(1)},\cdots,Y_{y(n)}\}$
be the $\varepsilon$\emph{-}partition of unity determined by $A$,
as in Definition \ref{Def.. Epsilon-parittion of unity}. For short,
abuse notations and write $Y_{k}\equiv Y_{y(k)}$ for each $k=1,\cdots n$.
By Proposition \ref{Prop. Properties of  epsilon partition of unity },
we have $\sum_{k=1}^{n}Y_{k}\leq1$, with equality prevailing on $K\subset\bigcup_{x\in A}(d(\cdot,x)\leq\varepsilon)$.
It follows that $\sum_{k=1}^{n}I(X_{i}Y_{k})\leq I(X_{i})$ for each
$i\geq0$, with equality in the case $i=0$. Therefore

\[
\sum_{k=1}^{n}\sum_{i=1}^{\infty}I(X_{i}Y_{k})\leq\sum_{i=1}^{\infty}I(X_{i})<I(X_{0})=\sum_{k=1}^{n}I(X_{0}Y_{k})
\]
Hence there exists some $k=1,\cdots n$ for which

\[
\sum_{i=1}^{\infty}I(X_{i}Y_{k})<I(X_{0}Y_{k})
\]
Again by Proposition \ref{Prop. Properties of  epsilon partition of unity },
for each $x\in S$ with $Y_{k}(x)>0$, we have $d(x,y_{k})<2\varepsilon$.
Hence 
\[
(Y_{k}(x)>0\mbox{ and }Y_{k}(x')>0)\Rightarrow(d(x,x')\leq4\varepsilon\mbox{ and }x\in B\subset K')
\]
for $x,x'\in S$.

Define $Z_{1}\equiv Y_{k}$. By repeating the above argument with
$\varepsilon_{m}=\frac{1}{4m}$ $(m=1,2,\cdots)$, we can construct
inductively a sequence of non-negative continuous functions $(Z_{m})_{m=1,2,\cdots}$
such that, for each $m\geq1$ and for each $x,x'\in S$, we have
\begin{equation}
(Z_{m}(x)>0\mbox{ and }Z_{m}(x')>0)\Rightarrow(d(x,x')\leq\frac{1}{m}\mbox{ and }x,x'\in K')\label{eq:temp2-1}
\end{equation}
and such that 

\begin{equation}
\sum_{i=1}^{\infty}I(X_{i}Z_{1}\cdots Z_{m})<I(X_{0}Z_{1}\cdots Z_{m})\label{eq:temp-137}
\end{equation}

Since all terms in \ref{eq:temp-137} are non-negative, the same inequality
holds if the infinite sum is replaced by the partial sum of the first
$m$ terms. By the positivity of $I$, this implies for each $m\geq1$
 the existence of a point $x_{m}$ such that

\begin{equation}
\sum_{i=1}^{m}X_{i}Z_{1}\cdots Z_{m}(x_{m})<X_{0}Z_{1}\cdots Z_{m}(x_{m})\label{eq:temp-16}
\end{equation}

In particular $Z_{p}(x_{m})>0$ for each $p\leq m$. Therefore the
inference \ref{eq:temp2-1} yields $x_{p}\in K'$ and $d(x_{p},x_{m})\leq\frac{1}{p}$
for each $p\leq m$. Hence $(x_{m})_{m=1,2,\cdots}$ is a Cauchy sequence
in $K'$ and converges to some point $x\in K'$. By the definition
of the function $Z$ at the beginning of this proof, we have $Z(x)=1$. 

Canceling positive common factors on both sides of inequality \ref{eq:temp-16},
we obtain $\sum_{i=1}^{p}X_{i}(x_{m})<X_{0}(x_{m})$ for each $p\leq m$.
Letting $m\rightarrow\infty$ yields $\sum_{i=1}^{p}X_{i}(x)\leq X_{0}(x)$
for each $p\geq1$. 

The conditions in Lemma \ref{Lem. Chan's lemma} have been established.
The conclusion of the present proposition follows.
\end{proof}

\section{Integration Space \textemdash{} the Daniell Integral}

Integration on a locally compact space is a special case of Daniell
integration, introduced next. 
\begin{defn}
\label{Def. Integration Space} \textbf{(integration Space).} An \emph{integration
space\index{integration space}} is a triple $(\Omega,L,I)$ where
$\Omega$ is a non-empty set, $L$ is a set of real-valued functions
on $\Omega$, and $I$ is a non-zero real-valued function with $domain(I)=L$,
satisfying the following conditions.
\end{defn}
\begin{enumerate}
\item If $X,Y\in L$ and $a,b\in R$, then $aX+bY,|X|$, and $X\wedge1$
belong to $L$, and $I(aX+bY)=aI(X)+bI(Y)$. In particular, if $X,Y\in L$,
then there exists $\omega\in domain(X)\cap domain(Y)$.
\item If a sequence $(X_{i})_{i=0,1,2,\cdots}$ of functions in $L$ is
such that $X_{i}$ is non-negative for each $i\geq1$ and such that
$\sum_{i=1}^{\infty}I(X_{i})<I(X_{0})$, then there exists a point
$\omega\in\bigcap_{i=0}^{\infty}domain(X_{i})$ such that $\sum_{i=1}^{\infty}X_{i}(\omega)<X_{0}(\omega)$.
This condition will be referred to as the \emph{positivity condition}\index{positivity condition for integration}
for $I.$
\item For each $X\in L$, we have $I(X\wedge n)\rightarrow I(X)$ and $I(|X|\wedge n^{-1})\rightarrow0$
as $n\rightarrow\infty$.
\end{enumerate}
$I$ is then called an $integration$\index{integration} or \emph{integral}
on $(\Omega,L)$, and $I(X)$ called the \index{integral}\emph{integral}
of $X$. A function $X\in L$ is said to \emph{integrable relative
to} $I$. $\square$

Note that given $X\in L$, the function $X\wedge n=n((\frac{1}{n}X)\wedge1)$
belongs to $L$ by condition 1 of Definition \ref{Def. Integration Space}.
Similarly the function $|X|\wedge n^{-1}$ belongs to $L$. Hence
$I(X\wedge n)$ and $I(|X|\wedge n^{-1})$ in Condition 3 are defined

In the following, in order to minimize clutter, we will write $IX$
for $I(X)$, and $IXY$ etc for $I(XY)$ etc, when there is no risk
of confusion. 

Note that, in general, there is no assumption that two functions should
have a point in the intersection of their domains. The positivity
condition is an existence condition useful in many constructions. 

One trivial example of an integration space is the triple $(\Omega,L,\delta_{\omega})$
where $\omega$ is a given point in a given set $\Omega$, where  is
the set of all functions $X$ on $\Omega$ whose domains contain $\omega$,
and where $\delta_{\omega}$ is defined on $L$ by $\delta_{\omega}(X)=X(\omega)$.
The integration $\delta_{\omega}$ is called the \emph{\index{point mass}point
mass} at $\omega$. 
\begin{prop}
\label{Prop. Integration on S is  an integration space} \textbf{\emph{(An
integration on a locally compact space entails an integration space).
}}Let $I$ be an integration on the locally compact metric space $(S,d)$
as defined in Definition \ref{Def. integration on loc compact space}.
Then $(S,C(S,d),I)$ is an integration space.
\end{prop}
\begin{proof}
The positivity condition in Definition \ref{Def. Integration Space}
has been proved for $(S,C(S,d),I)$ in Proposition \ref{Prop. Positivity of I on (S,C)}.
The other conditions are trivial.
\end{proof}
The next proposition collects some simple properties of integration
spaces.
\begin{prop}
\label{Prop. Properties of  Integration Space} \textbf{\emph{(Basic
properties of an integration space). }}Let $(\Omega,L,I)$ be an integration
space. Then the following holds.
\end{prop}
\begin{enumerate}
\item \emph{If $X,Y\in L$, then $X\vee Y,X\wedge Y\in L$. If in addition
$a>0$, then $X\wedge a\in L$ and $I(X\wedge a)$ is continuous in
$a$.}
\item \emph{If $X\in L$, then $X_{+},X_{-}\in L$ and $IX=IX_{+}+IX_{-}$.}
\item \emph{For any $X\in L$ with $IX>0$, there exists $\omega$ such
that $X(\omega)>0.$ }
\item \emph{Suppose $X(\omega)\geq0$ for each $\omega\in domain(X)$. Then
we have $IX\geq0.$ }
\item \emph{There exists a non-negative $X\in L$ such that $IX=1$. }
\item \emph{For any sequence $(X_{i})_{i=1,2\cdots}$ in $L$, there exists
a point $\omega\in\cap_{i=0}^{\infty}domain(X_{i})$.}
\end{enumerate}
\begin{proof}

1. The first part follows from $X\vee Y=(X+Y+|X-Y|)/2$ and $X\wedge Y=(X+Y-|X-Y|)/2$.
The second part follows from $X\wedge a=a(\frac{X}{a}\wedge1)$ and
,in view of Condition 3 in Definition \ref{Def. Integration Space},
from $I|X\wedge a-X\wedge b|\leq I(|b-a|\wedge|X|)$ for $a,b>0$.

2. The conclusion follows from $X_{+}=X\vee0X$, $X_{-}=X\wedge0X$,
and $X=X_{+}+X_{-}$. 

3. Suppose $X\in L$ has integral $IX>0.$ The positivity condition
in Definition \ref{Def. Integration Space}, applied to the sequence
$(X,0X,0X,...)$, guarantees an $\omega$ such that $X(\omega)>0$.

4. Suppose $IX<0$. Then $I(-X)>0,$ and part 3 of this proposition
would give an $\omega\in domain(X)$ with $X(\omega)<0$, a contradiction.
Hence $IX\geq0$.

5. Since $I$ is nonzero and linear, there exists $X$ such that $I(X)>0$.
By part 4 of this proposition, and by the linearity of $I$, we see
that $IX\leq I|X|$ and so $I|X|>0$. Let $X_{0}$ denote the function
$|X|/I|X|$. Then $X_{0}$ is non-negative and $IX_{0}=1$.

6. Let a non-negative $X_{0}\in L$ be such that $IX_{0}=1$. Suppose
$(X_{i})_{i=1,2\cdots}$ is a sequence of functions in $L$. Then
the sequence $(X_{0},0X_{1},0X_{2},\cdots)$ trivially satisfies the
requirements in the positivity condition in Definition \ref{Def. Integration Space},
which therefore guarantees a point in the intersection of the domains. 
\end{proof}
\begin{defn}
\label{Def. integration subspace} \textbf{(Integration subspace).}
Let $(\Omega,L,I)$ be an integration space. Let $L'$ be a subfamily
of $L$ such that $(\Omega,L',I)$ is an integration space. We will
then call $(\Omega,L',I)$ an \index{integration subspace}integration
subspace of $(\Omega,L,I)$. When confusion is unlikely, we will abuse
terminology and simply call $L'$ an integration subspace of $L$,
with $\Omega$ and $I$ understood. 
\end{defn}
\begin{prop}
\label{Prop. A closed subspace L' of L yields an integration subspace}
\textbf{\emph{(Linear subspace of integrable functions closed to absolute
values and minimum with constants is an integrations subspace).}}
Let $(\Omega,L,I)$ be an integration space. Let $L'$ be a linear
subspace of $L$ such that if $X,Y\in L'$ then $|X|,X\wedge1\in L'$.
Then $(\Omega,L',I)$ is an integration subspace of $(\Omega,L,I)$. 
\end{prop}
\begin{proof}
By hypothesis, $L'$ is closed to linear operations, absolute values,
and the operation of taking minimum with the constant 1. Condition
1 in Definition \ref{Def. Integration Space} for an integration space
is thus satisfied by $L'$. Conditions 2 and 3  are inherited by $(\Omega,L',I)$
from $(\Omega,L,I)$. 
\end{proof}
\begin{prop}
\label{Prop. Integration induced by a surjection} \textbf{\emph{(Integration
induced by a surjection).}} Let $(\Omega,L,I)$ be an integration
space. Let $\pi:\bar{\Omega}\rightarrow\Omega$ be a function from
some set $\bar{\Omega}$ onto $\Omega$. For each $f\in L$ write
$f(\pi)\equiv f\circ\pi$. Define $\overline{L}\equiv\{f(\pi):f\in L\}$
and define $\bar{I}:\overline{L}\rightarrow R$ by $\bar{I}X\equiv I(f)$
for each $f\in L$ and $X=f(\pi)\in\overline{L}$. Then $(\bar{\Omega},\overline{L},\bar{I})$
is an integration space.
\end{prop}
\begin{proof}
Suppose $X=f(\pi)=g(\pi)$ for some $f,g\in L$. Let $\omega\in domain(f)$
be arbitrary. Since $\pi$ is an on-to function, there exists $\varpi\in\bar{\Omega}$
such that $\pi(\varpi)=\omega\in domain(f)$. It follows that $\varpi\in domain(f(\pi))=domain(g(\pi))$
and so $\omega=\pi(\varpi)\in domain(g)$. Since $\omega\in domain(f)$
is arbitrary, we see that $domain(f)\subset domain(g)$ and, by symmetry,
$domain(f)=domain(g)$. Moreover $f(\omega)=f(\pi(\varpi))=g(\pi(\varpi))=g(\omega)$.
We conclude that $f=g$. 

Next let  and $a,b\in R$ be arbitrary, where $f,h\in L$. Then $af+bh\in L$
and so $aX+bY=(af+bh)(\pi)\in\overline{L}$. Furthermore $\bar{I}(aX+bY)\equiv I(af+bh)=aI(f)+bI(h)\equiv a\bar{I}X+b\bar{I}Y$.
Thus $\overline{L}$ is a linear space and $\bar{I}$ is a linear
function. Similarly, $|X|=|f|(\pi)\in\overline{L}$ and $a\wedge X=(a\wedge f)(\pi)\in\overline{L}$.
Furthermore $\bar{I}(a\wedge X)\equiv I(a\wedge f)\rightarrow I(f)\equiv\bar{I}X$
as $a\rightarrow\infty$, while $\bar{I}(a\wedge|X|)\equiv I(a\wedge|f|)\rightarrow0$
as $a\rightarrow0$. Thus Conditions 1 and 3 in Definition \ref{Def. Integration Space}
for an integration space are verified for the triple $(\bar{\Omega},\overline{L},\bar{I})$.

It remains to prove the positivity condition, Condition 2 in Definition
\ref{Def. Integration Space}. To that end, let $(X_{i})_{i=0,1,2,\cdots}$
be a sequence in $\overline{L}$ such that $X_{i}$ is non-negative
for each $i\geq1$ and such that $\sum_{i=1}^{\infty}\bar{I}X_{i}<\bar{I}X_{0}$.
For each $i\geq0$ let $f_{i}\in L$ be such that $X_{i}=f_{i}(\pi)$.
Then, since $\pi$ is an on-to function, $f_{i}\geq0$ for each $i\geq1$.
Moreover $\sum_{i=1}^{\infty}I(f_{i})\equiv\sum_{i=1}^{\infty}\bar{I}X_{i}<\bar{I}X_{0}\equiv I(f_{0})$.
Since $I$ is an integration, there exists $\omega\in\bigcap_{i=0}^{\infty}domain(f_{i})$
such that $\sum_{i=1}^{\infty}f_{i}(\omega)<f_{0}(\omega)$. Let $\varpi\in\bar{\Omega}$
be such that $\pi(\varpi)=\omega$. Then 
\[
\varpi\in\bigcap_{i=0}^{\infty}domain(f_{i}(\pi))=\bigcap_{i=0}^{\infty}domain(X_{i}).
\]
By hypothesis$\sum_{i=1}^{\infty}X_{i}(\varpi)=\sum_{i=1}^{\infty}f_{i}(\omega)<f_{0}(\omega)=X_{0}(\varpi)$.
All the conditions in Definition \ref{Def. Integration Space} have
been established. Accordingly, $(\bar{\Omega},\overline{L},\bar{I})$
is an integration space.
\end{proof}

\section{Complete Extension of Integrations }

Because discontinuous real random variables will be of interest, integration
spaces like $(S,C(S,d),I)$ will need to be expanded. More generally,
given an integration space $(\Omega,L,I)$, we can expand the set
$L$ to a larger set $L_{1},$ and extend the integration $I$ to
$L_{1}$, by summing a series of small pieces in $L$, small in the
sense that the integrals of the absolute values of these pieces sum
to a finite number. This is analogous to of the usual extension of
rational numbers to reals by representing a real number as the sum
of an absolutely convergent series of rational numbers.
\begin{defn}
\label{Def. Integrablee functions and Completion of integration space}
\textbf{(Integrable functions, and Completion of an integration space).}
Let $(\Omega,L,I)$ be an integration space. A function $X$ on a
subset of $\Omega$ is called an \emph{\index{integrable function}integrable
function} if there exists a sequence $(X_{n})_{n=1,2,\cdots}$ in
$L$ such that (i) $\sum_{i=1}^{\infty}I|X_{i}|<\infty$, (ii) $domain(X)$
contains the set 
\[
D\equiv\{\omega\in\cap_{i=1}^{\infty}domain(X_{i}):\sum_{i=1}^{\infty}|X_{i}(\omega)|<\infty\},
\]
and (iii) $X(\omega)=\sum_{i=1}^{\infty}X_{i}(\omega)$ for each $\omega\in D$.
The sequence $(X_{n})_{n=1,2,\cdots}$ is called a \emph{\index{representation of an integrable function}representation}
of $X$ by elements of $L$ relative to $I$. The set of integrable
functions will be denoted by $L_{1}$. Define the sum
\[
I_{1}(X)\equiv\sum_{i=1}^{\infty}IX_{i}
\]
and call it the \emph{integral} of $X$. Then $I_{1}$ is called the
complete extension\index{complete extension of an integration} of
$I$. Likewise, $L_{1}$ and $(\Omega,L_{1},I_{1})$ are called the
\emph{complet}e \emph{extensions} \index{complete extension of an integration space},
or simply \emph{completion,} of $L$ and $(\Omega,L,I)$ respectively.

The next proposition and theorem prove that (i') $I$ is well defined
on $L$, and (ii) $(\Omega,L_{1},I_{1})$ is indeed an integration
space with $L\subset L_{1}$ and $I=I_{1}|L$. Henceforth we can use
the same symbol $I$ to denote the given integration and its complete
extension, and write $I$ also for $I_{1}$.

An integration space $(\Omega,L,I)$ is said to be \index{complete integration space}complete
if $(\Omega,L,I)=(\Omega,L_{1},I)$. $\square$
\end{defn}
Suppose $(X_{n})_{n=1,2,\cdots}$ is a representation of $X$. If
we define $X_{0}\equiv\sum_{i=1}^{\infty}X_{i}$, then $(X_{n})_{n=1,2,\cdots}$
is a\emph{ }representation also of $X_{0}$, with $domain(X_{0})=D\subset domain(X)$
and $X=X_{0}$ on $domain(X_{0})$. 

The next proposition shows that $I_{1}X$ is well-defined; in other
words, it is independent of the representation.
\begin{prop}
\label{Prop. Complete extension of integration is well-defined} \textbf{\emph{(Complete
extension of integration is well defined). }}If $(X_{n})_{n=1,2,\cdots}$
and $(Y_{n})_{n=1,2,\cdots}$ are two representations of the integrable
function $X$, then $\sum_{i=1}^{\infty}IX_{i}=\sum_{i=1}^{\infty}IY_{i}$.
\end{prop}
\begin{proof}
By the definition of a representation, the series $\sum_{i=1}^{\infty}I|X_{i}|$
and $\sum_{i=1}^{\infty}I|Y_{i}|$ converge. Suppose $\sum_{i=1}^{\infty}IX_{i}<\sum_{i=1}^{\infty}IY_{i}$.
Then for some large number $m$ we have

\[
\sum_{i=m+1}^{\infty}I|X_{i}|+\sum_{i=m+1}^{\infty}I|Y_{i}|+\sum_{i=m+1}^{\infty}I|Y_{i}-X_{i}|<\sum_{i=1}^{m}IY_{i}-\sum_{i=1}^{m}IX_{i}=I\sum_{i=1}^{m}(Y_{i}-X_{i})
\]
The conditions in Definition \ref{Def. Integration Space} of integration
space then implies the existence of a point $\omega\in\cap_{i=1}^{\infty}(domain(X_{i})\cap domain(Y_{i}))$
such that

\[
\sum_{i=m+1}^{\infty}|X_{i}(\omega)|+\sum_{i=m+1}^{\infty}|Y_{i}(\omega)|+\sum_{i=m+1}^{\infty}|Y_{i}(\omega)-X_{i}(\omega)|<\sum_{i=1}^{m}Y_{i}(\omega)-\sum_{i=1}^{m}X_{i}(\omega)
\]
Hence, applying the triangle inequality, we have 

\[
0<\sum_{i=1}^{m}Y_{i}(\omega)-\sum_{i=1}^{m}X_{i}(\omega)-\sum_{i=m+1}^{\infty}|Y_{i}(\omega)-X_{i}(\omega)|
\]

\[
\leq|\sum_{i=1}^{m}Y_{i}(\omega)-\sum_{i=1}^{m}X_{i}(\omega)|-|\sum_{i=m+1}^{\infty}Y_{i}(\omega)-\sum_{i=m+1}^{\infty}X_{i}(\omega)|
\]

\[
\leq|\sum_{i=1}^{\infty}Y_{i}(\omega)-\sum_{i=1}^{\infty}X_{i}(\omega)|=|X(\omega)-X(\omega)|=0
\]
The next to last equality is because both $(Y_{n})_{n=1,2,\cdots}$
and $(X_{n})_{n=1,2,\cdots}$ are, by hypothesis, representations
of $X$. Thus the assumption $\sum_{i=1}^{\infty}IX_{i}<\sum_{i=1}^{\infty}IY_{i}$
leads to a contradiction. Therefore $\sum_{i=1}^{\infty}IX_{i}\geq\sum_{i=1}^{\infty}IY_{i}$.
Similarly $\sum_{i=1}^{\infty}IX_{i}\leq\sum_{i=1}^{\infty}IY_{i}$,
and the equality follows.
\end{proof}
\begin{thm}
\label{Thm. (Omega,L1,I) is integration space} \textbf{\emph{(Complete
extensiion of an integration space is an integration space). }}Let
$(\Omega,L,I)$ be an integration space. Then $(\Omega,L_{1},I_{1})$
is an integration space. Moreover, $L\subset L_{1}$, and $I_{1}X=IX$
for each $X\in L$.
\end{thm}
\begin{proof}
Let $X\in L$ be arbitrary. Then $(X,0X,0X,...)$ is a representation
of $X$. Hence $X\in L_{1}$ and $I_{1}X=IX$. It remains to verify,
for the triple $(\Omega,L_{1},I_{1})$, the conditions in Definition
\ref{Def. Integration Space} of integration spaces. Proposition \ref{Prop. Properties of  Integration Space}
and Condition 2 in Definition \ref{Def. Integration Space} guarantees
that $domain(X)$ is non-empty for each $X\in L_{1}$.

First let $X,Y\in L_{1}$ be arbitrary, with representations $(X_{n})_{n=1,2,\cdots}$
and $(Y_{n})_{n=1,2,\cdots}$ respectively. Let $a,b\in R$. Then
clearly the sequence 
\[
(aX_{n}+bY_{n},X_{n},-X_{n},Y_{n},-Y_{n})_{n=1,2,\cdots}
\]
is a representation of $aX+bY$. The seemingly redundant terms $X_{n},-X_{n},$
$Y_{n},-Y_{n},$ $\cdots$ are included to ensure that the absolute
convergence of the resulting series at some $\omega$ implies that
$\omega\in domain(a+X+bY)$. Similar tricks will be used several times
later without further comments. Thus we see that $aX+bY\in L_{1}$
and $I_{1}(aX+bY)=aI_{1}X+bI_{1}Y$. Similarly, let $a>0$ be arbitrary.
Because 
\[
I|a\wedge\sum_{i=1}^{n}X_{i}-a\wedge\sum_{i=1}^{n-1}X_{i}|\leq I|X_{n}|
\]
 for each $n\geq1$, the sequence
\[
(a\wedge\sum_{i=1}^{n}X_{i}-a\wedge\sum_{i=1}^{n-1}X_{i},X_{n},-X_{n})_{n=1,2,\cdots}
\]
 is a representation of $X\wedge a$. Hence $X\wedge a$ belongs to
$L_{1}$, and 
\begin{equation}
I_{1}(X\wedge a)=\lim_{n\rightarrow\infty}I(a\wedge\sum_{i=1}^{n}X_{i})\leq\lim_{n\rightarrow\infty}I(\sum_{i=1}^{n}X_{i})=I_{1}X\label{eq:temp-1}
\end{equation}
In particular $1\wedge X\in L_{1}$. Again, because 
\[
I|(|\sum_{i=1}^{n}X_{i}|-|\sum_{i=1}^{n-1}X_{i}|)|\leq I|X_{n}|
\]
for each $n\geq1$, the sequence
\[
(|\sum_{i=1}^{n}X_{i}|-|\sum_{i=1}^{n-1}X_{i}|,X_{n},-X_{n})_{n=1,2,\cdots}
\]
 is a representation of $|X|$. Hence $|X|$ belongs to to $L_{1}$,
with 
\begin{equation}
I_{1}|X|=\lim_{n\rightarrow\infty}I\:|\sum_{i=1}^{n}X_{i}|.\label{eq:temp0}
\end{equation}
It follows that $I_{1}$ is a nonnegative function on $L_{1}$. In
other words, if $X,Y\in L_{1}$ are such that $X\leq Y$ on $domain(X)\cap domain(Y)$,
then $I_{1}X\leq I_{1}Y$.

We next verify the positivity condition in Definition \ref{Def. Integration Space}.
To that end, let $(X_{i})_{i=0,1,2,\cdots}$ be a sequence of functions
in $L_{1}$ such that $X_{i}$ is non-negative for $i\geq1$ and such
that 
\[
\sum_{i=1}^{\infty}I_{1}X_{i}<I_{1}X_{0}
\]
For each $i\geq0$, let $(X_{i,k})_{k=1,2,\cdots}$ be a representation
of $X_{i}$. Then $\sum_{k=1}^{\infty}I|X_{i,k}|$ $<\infty$ for
each $i\geq0$. Since $X_{i}$ is non-negative for $i\geq1$, we see
from part 1 of this proof that
\[
\lim_{n\rightarrow\infty}I|\sum_{k=1}^{n}X_{i,k}|=I_{1}|X_{i}|=I_{1}X_{i}
\]
 for $i\geq1$. Therefore there exists a sequence $(m_{i})_{i=0,1,2,\cdots}$
of integers such that
\[
\sum_{i=1}^{\infty}I|\sum_{k=1}^{m(i)}X_{i,k}|+\sum_{i=0}^{\infty}\sum_{k=m(i)+1}^{\infty}I|X_{i,k}|<\sum_{k=1}^{m(0)}IX_{0,k}
\]
The positivity condition of \ref{Def. Integration Space} on $(\Omega,L,I)$
therefore guarantees the existence of $\omega\in\bigcap_{i=0}^{\infty}\bigcap_{k=1}^{m(i)}domain(X_{i,k})$
such that 
\begin{equation}
\sum_{i=1}^{\infty}|\sum_{k=1}^{m(i)}X_{i,k}(\omega)|+\sum_{i=0}^{\infty}\sum_{k=m(i)+1}^{\infty}|X_{i,k}(\omega)|<\sum_{k=1}^{m(0)}X_{0,k}(\omega)\label{eq:temp}
\end{equation}
It follows from the definition of a representation that $X_{i}(\omega)$
is defined and $\sum_{k=1}^{\infty}X_{i,k}(\omega)=X_{i}(\omega)$
for each $i\geq0$. Moreover 
\[
\sum_{i=1}^{\infty}X_{i}(\omega)=\sum_{i=1}^{\infty}\sum_{k=1}^{\infty}X_{i,k}(\omega)=\sum_{i=1}^{\infty}\sum_{k=1}^{m(i)}X_{i,k}(\omega)+\sum_{i=1}^{\infty}\sum_{k=m(i)+1}^{\infty}X_{i,k}(\omega)
\]
\[
\leq\sum_{i=1}^{\infty}|\sum_{k=1}^{m(i)}X_{i,k}(\omega)|+\sum_{i=0}^{\infty}\sum_{k=m(i)+1}^{\infty}|X_{i,k}(\omega)|-\sum_{k=m(0)+1}^{\infty}|X_{0,k}(\omega)|
\]
\[
<\sum_{k=1}^{m(0)}X_{0,k}(\omega)-\sum_{k=m(0)+1}^{\infty}|X_{0,k}(\omega)|\leq\sum_{k=1}^{\infty}X_{0,k}(\omega)=X_{0}(\omega)
\]
where the next to last inequality follows form inequality \ref{eq:temp-1}
above. This proves condition 2 of Definition \ref{Def. Integration Space}
for $(\Omega,L_{1},I_{1})$. 

Now let $X\in L_{1}$, with a representation $(X_{i})_{i=1,2,\cdots}$
in $L$. Then, for every $m>0$, the sequence 
\[
(X_{1},-X_{1},X_{2},-X_{2},\cdots,X_{m},-X_{m},X_{m+1},X_{m+2},X_{m+3},\cdots)
\]
is a representation of the function $X-\sum_{i=1}^{m}X_{i}\in L_{1}$.
Therefore, applying equation \ref{eq:temp0} in above, we have
\begin{equation}
I_{1}|X-\sum_{i=1}^{m}X_{i}|=\lim_{n\rightarrow\infty}I|\sum_{i=m+1}^{n}X_{i}|\leq\sum_{i=m+1}^{\infty}I|X_{i}|\rightarrow0\label{eq:temp1}
\end{equation}
Hence, for any given $\varepsilon>0$, there exists $m$ so large
that $I_{1}|X-\sum_{i=1}^{m}X_{i}|<\varepsilon$. Write $Y\equiv\sum_{i=1}^{m}X_{i}$.
Then 
\[
I_{1}|X-Y|<\varepsilon
\]
and so 
\[
I_{1}|X\wedge n-Y\wedge n|\leq I_{1}|X-Y|<\varepsilon
\]
for each $n\geq1$. In view of Condition 3 of Definition\ref{Def. Integration Space}
for $(\Omega,L,I)$, there exists $p\geq1$ so large that 
\[
|I_{1}(Y\wedge n)-I_{1}Y|=|I(Y\wedge n)-IY|<\varepsilon
\]
for each $n\geq p$. Hence 
\[
I_{1}X\geq I_{1}(X\wedge n)>I_{1}(Y\wedge n)-\varepsilon>I_{1}Y-2\varepsilon>I_{1}X-3\varepsilon
\]
for each $n\geq p$. Since $\varepsilon$ is arbitrary, we have $I_{1}(X\wedge n)\rightarrow I_{1}X$
as $n\rightarrow\infty$. Separately, again by Condition 3 of Definition
\ref{Def. Integration Space} for $(\Omega,L,I)$, there exists $p\geq1$
so large that $I(|Y|\wedge n^{-1})<\varepsilon$ for each $n\geq p$.
Therefore
\[
I_{1}(|X|\wedge n^{-1})\leq I_{1}(|Y|\wedge n^{-1})+I_{1}|X-Y|<2\varepsilon
\]
Hence $I_{1}(|X|\wedge n^{-1})\rightarrow0$ as $n\rightarrow\infty$.
All three conditions in Definition \ref{Def. Integration Space} have
been verified for $(\Omega,L_{1},I_{1})$ to be an integration space.
\end{proof}
\begin{cor}
\label{Cor.  L is dense in L1} \textbf{\emph{($L$ is dense in its
complete extension).}} If $X\in L_{1}$ has representation $(X_{i})_{i=1,2,\cdots}$,
then $\lim_{n\rightarrow\infty}$ $I|X-\sum_{i=1}^{m}X_{i}|=0.$
\end{cor}
\begin{proof}
See expression \ref{eq:temp1} in the proof of \ref{Thm. (Omega,L1,I) is integration space}. 
\end{proof}
Henceforth we will, write $I$ also for $I_{1}$. In words, we will
use the symbol for a given integration also for its complete extension.
\begin{prop}
\label{Prop. Completeness of L1} \textbf{\emph{(Complete extension
of an integration spacce is a complete metric space). }} Let $X,Y\in L_{1}$
be arbitrary. Define $\rho_{I}(X,Y)=I|X-Y|.$ Then $(L_{1},\rho_{I})$
is a complete metric space, and $L$ is a dense subset of $L_{1}$. 
\end{prop}
\begin{proof}
The proof that $(L_{1},\rho_{I})$ is a metric space is trivial. Corollary
\ref{Cor.  L is dense in L1} implies that $L$ is a dense subset
of $(L_{1},\rho_{I})$. It remains to prove that the latter is complete.

Let $(X_{n})_{n=1,2,\cdots}$ be any Cauchy sequence in $(L_{1},\rho_{I})$.
To prove completeness we need to find some $Z\in L_{1}$ such that
$\lim_{n\rightarrow\infty}I|X_{n}-Z|=0$. Let $(n_{i})_{i=1,2,\cdots}$
be an increasing sequence such that $I|X_{n(i)}-X_{n}|<2^{-i}$ for
each $n>n_{i}$ and $i\geq1$. Because $L$ is dense in $L_{1},$
there exists for each $i\geq1$ some $Z_{i}\in L$ such that $I|X_{n(i)}-Z_{i}|<2^{-i}$.
Then $I|Z_{i}-Z_{i+1}|<2^{-i+1}+2^{-i-1}$ for each $i\geq1.$ Hence
the sequence ($Z_{1},Z_{2}-Z_{1},Z_{3}-Z_{2},\cdots)$ is the representation
of some $Z\equiv Z_{1}+\sum_{i=1}^{\infty}(Z_{i+1}-Z_{i})\in L_{1}$.
Corollary \ref{Cor.  L is dense in L1} therefore implies that $\lim_{i\rightarrow\infty}I|Z_{i}-Z|=0$.
At the same time, for each $n>n_{i}$ and $i\geq1$ we have 
\[
I|X_{n}-Z|\leq I|X_{n}-X_{n(i)}|+I|X_{n(i)}-Z|<2^{-i}+I|X_{n(i)}-Z|
\]
Combining, we see that $\lim_{n\rightarrow\infty}I|X_{n}-Z|=0$. 
\end{proof}
\begin{cor}
\label{Cor. (L1)1=00003DL1} \textbf{\emph{(Nothing is gained from
further complete extension). }}Let $L_{1}$ be space of integrable
functions on $(\Omega,L,I)$. Let $(L_{1})_{1}$ be space of integrable
functions on $(\Omega,L_{1},I)$. Then $(L_{1})_{1}=L_{1}$. 
\end{cor}
\begin{proof}
Let $Z\in(L_{1})_{1}.$ By Proposition \ref{Prop. Completeness of L1},
for any $\varepsilon>0$ there exists a function $Y\in L_{1}$ with
$I|Z-Y|<\varepsilon$, and there exists a function $X\in L$ with
$I|Y-X|<\varepsilon$, and so $I|Z-X|<2\varepsilon$. Thus we can
construct a sequence in $L$ which converges to $Z$. Proposition
\ref{Prop. Completeness of L1} implies that $Z\in L_{1}$.
\end{proof}
\begin{cor}
\label{Cor. I=00003DI' on L0, L0 dense in complete L=00003D>I=00003DI' on L}
\textbf{\emph{(Two integrations on the same space of integrable functions
are equal if they agree on some dense subset).}} Let $(\Omega,L,I)$
and $(\Omega,L,I')$ be complete integration spaces. Suppose $I=I'$
on some subset $L_{0}$ of $L$ which is dense in $L$ relative to
the metric defined by $\rho_{I}(X,Y)=I|X-Y|$ for all $X,Y\in L$.
Then $I=I'$ on $L$.
\end{cor}
\begin{proof}
Let $X\in L$ be arbitrary. Let $(X_{n})_{n=1,2,\cdots}$ be a sequence
in $L_{0}$ which converges to $X$ relative to the metric $\rho_{I}$.
As in the proof of Proposition \ref{Prop. Completeness of L1}, we
can then construct a sequence $(Y_{n})_{n=1,2,\cdots}$ in $L_{0}$
which is a representation of $X$ relative to $I$. Since $I=I'$
on $L_{0}$, it follows immediately that $(Y_{n})_{n=1,2,\cdots}$
is also a representation of $X$ relative to $I'$, with 
\[
I'X=\sum_{i=1}^{\infty}I'Y_{n}=\sum_{i=1}^{\infty}IY_{n}=IX.
\]
Since $X\in L$ is arbitrary, we conclude that $I=I'$.
\end{proof}
For monotone sequences, we have a very useful theorem for establishing
convergence in $L_{1}$.
\begin{thm}
\textbf{\emph{(Monotone Convergence Theorem).}} \label{Thm. Monotone Convergence}
Let $(\Omega,L_{1},I)$ be a complete integration space. Suppose $(X_{i})_{i=1,2,\cdots}$
is a sequence in $L_{1}$ such that $X_{i-1}\leq X_{i}$ on $\cap_{i=1}^{\infty}domain(X_{i})$,
and such that $\lim_{i\rightarrow\infty}I(X_{i})$ exists. Then $X\equiv\lim_{i\rightarrow\infty}X_{i}\in L_{1}$.
Moreover $\lim_{i\rightarrow\infty}I|X-X_{i}|=0$. Similarly, suppose
$(Y_{i})_{i=1,2,\cdots}$ is a sequence in $L_{1}$ such that $Y_{i-1}\geq Y_{i}$
on $\cap_{i=1}^{\infty}domain(Y_{i})$, and such that $\lim_{i\rightarrow\infty}I(Y_{i})$
exists. Then $Y\equiv\lim_{i\rightarrow\infty}Y_{i}\in L_{1}$. Moreover
$\lim_{i\rightarrow\infty}I|Y-Y_{i}|=0$.
\end{thm}
\begin{proof}
The sequence $(X_{1},X_{2}-X_{1},X_{3}-X_{2},\cdots)$ is obviously
a representation of $X$, with $IX=\lim_{i\rightarrow\infty}I(X_{1}+(X_{2}-X_{1})+\cdots+(X_{i}-X_{i}))=\lim_{i\rightarrow\infty}IX_{i}$.
Corollary \ref{Cor.  L is dense in L1} implies $\lim_{i\rightarrow\infty}I|X-X_{i}|=0$.
The second part of the theorem follows by applying the first part
to the sequence $(-Y_{i})_{i=1,2,\cdots}$m. 
\end{proof}

\section{Integrable Sets}

To model an event in a chance experiment which may or may not have
occurred after the outcome is observed, we can use a function with
only two possible values, $1$ or $0$. Equivalently we can specify
the subset of those outcomes that realize the event. We make these
notions precise in the present section.
\begin{defn}
\label{Def. Indicators and Mutually Eclusive} \textbf{(Indicators
and mutally excllusive subsets). }Subsets $A_{1},\cdots,A_{n}$ of
a set $\Omega$ are said to be \emph{mutually exclusive} if $A_{i}A_{j}=\phi$
for all $i,j$ with $i\neq j$. A function $X$ on a set $\Omega$
with only two possible values, 1 or $0$, is called an \emph{indicator}.\index{indicator}
Indicators $X_{1},\cdots,X_{n}$ are said to be \index{mutually exclusive}
\emph{mutually exclusive} if the sets $\{\omega\in domain(X_{i}):X_{i}(\omega)=1\}$
$(i=1,\cdots,n)$ are mutually exclusive. $\square$
\end{defn}
In the remainder of this section, let $(\Omega,L,I)$ be a complete
integration space. Recall that an integrable function need not be
defined everywhere. However, they are defined almost everywhere in
the sense of the next definition.
\begin{defn}
\label{Def. Full Set} \textbf{(Full set and almost everywhere properties).}
A subset $D$ of $\Omega$ is called a \emph{full set}\index{full set}
if $D\supset domain(X)$ for some integrable function $X\in L$. Two
functions $Y,Z$ on $\Omega$ are said to be equal \emph{almost everywhere}\index{almost everywhere},
with abbreviation $Y=Z$ a.e., if $Y=Z$ on a full set $D$. In other
words, $Y=Z$ a.e. if there exists a full set $D$ such that (i) $D\cap domain(Y)=D\cap domain(Z)$
and (ii) $Y(\omega)=Z(\omega)$ for each $\omega\in D\cap domain(Y)$.
In general, a statement about a general element $\omega$ of $\Omega$
is said to hold\emph{ almost everywhere}\index{almost everywhere},
a.e. for short, if it holds for each $\omega$ in a full set. $\square$
\end{defn}
For example, since according to the terminology established in the
Introduction of this book, the statement $Y\leq Z$ means that for
each $\omega\in\Omega$ we have (i) $\omega\in domain(Y)\Leftrightarrow\omega\in domain(Z)$
and (ii) $Y(\omega)\leq Z(\omega)$ if $\omega\in domain(Y)$. So
the statement $Y\leq Z$ a.e. means that there exists some full set
$D$ such that for each $\omega\in D$ the conditions (i) and (ii)
hold. Equivalently, $Y\leq Z$ a.e. iff (i') $\omega\in D\cap domain(Y)\Leftrightarrow\omega\in D\cap domain(Z)$
and (ii') $Y(\omega)\leq Z(\omega)$ for each $\omega\in D\cap domain(Y)$.
Thus $Y\leq Z$ a.e. iff $Y\leq Z$ on some full set $D$.

Similarly, if $A,B$ are subsets of $\Omega$, then $A\subset B$
a.e. iff $AD\subset BD$ for some full set $D$. 

Every integrable function is defined a.e. The last sentence, however,
does not tell us anything until we explore the properties of full
sets, in the next proposition.
\begin{prop}
\label{Prop. Property of Full Sets} \textbf{\emph{(Properties of
full sets). }}Let $X,Y,Z\in L$ denote integrable functions.
\end{prop}
\begin{enumerate}
\item \emph{A subset which contains a full set is a full set. The intersection
of a sequence of full sets is again a full set. }
\item \emph{Suppose $W$ is a function on $\Omega$ and $W=X$ a.e. Then
$W$ is an integrable function with $IW=IX$.}
\item \emph{If $D$ is a full set then $D=domain(X)$ for some $X\in L$}
\item \emph{$X=Y$ a.e. if and only if $I|X-Y|=0$. }
\item \emph{If $X\leq Y$ a.e. then $IX\leq IY$.}
\item \emph{If $X\leq Y$ a.e. and $Y\leq Z$ a.e. then $X\leq Z$ a.e.
Moreover, if $X\leq Y$ a.e. and $X\geq Y$ a.e. then $X=Y$ a.e. }
\item \emph{Almost everywhere equality is an equality relation in $L$.
In other words, for all $X,Y,Z\in L$ we have (i) $X=X$ a.e. (ii)
if $X=Y$ a.e. then $Y=X$ a.e. and (iii) if $X=Y$ a.e. and $Y=Z$
a.e. then $X=Z$ a.e. }
\end{enumerate}
\begin{proof}

1. Suppose $D_{n}\supset domain(X_{n})$ where $X_{n}\in L$ for $n\geq1$.
Define $X\equiv\sum_{n=1}^{\infty}0X_{n}$. Then $X\in L$ since $L$
is complete. Moreover 
\[
\bigcap_{n=1}^{\infty}D_{n}\supset\bigcap_{n=1}^{\infty}domain(X_{n})=\bigcap_{n=1}^{\infty}domain(0X_{n})=domain(X)
\]

2. By the definition of a.e. equality, there exists a full set $D$
such that $D\cap domain(W)=D\cap domain(X)$ and $W(\omega)=X(\omega)$
for each $\omega\in D\cap domain(X)$. By the definition of a full
set, $D\supset domain(Z)$ for some $Z\in L$. It is then easily verified
that the sequence $(X,0Z,0Z,\cdots)$ is a representation of the function
$W$. Therefore $W\in L_{1}=L$ with $IW=IX$.

3. Suppose $D$ is a full set. By definition $D\supset domain(X)$
for some $X\in L$. Define a function $W$ by $domain(W)\equiv D$
and $W(\omega)\equiv0$ for each $\omega\in D$. Then $W=0X$ on the
full set $domain(X)$. Hence by assertion 2 above, $W$ is an integrable
function, with $D=domain(W)$. 

4. Suppose $X=Y$ a.e. Then $|X-Y|=0(X-Y)$ a.e. Hence $I|X-Y|=0$
according to assertion 2. Suppose conversely that $I|X-Y|=0$. Then
the function defined by $Z\equiv\sum_{n=1}^{\infty}|X-Y|$ is integrable.
By definition 
\[
domain(Z)\equiv\{\omega\in domain(X-Y):\sum_{n=1}^{\infty}|X(\omega)-Y(\omega)|<\infty\}
\]
\[
=\{\omega\in domain(X)\cap domain(Y):X(\omega)=Y(\omega)\}
\]
Thus we see that $X=Y$ on the full set $domain(Z)$. 

5. Because $Y-X=|Y-X|$ a.e., we have, by assertion 4, $I(Y-X)=I|Y-X|\geq0$.

6. Suppose  $X\leq Y$ a.e. and $Y\leq Z$ a.e. Then there exists
a full set $D$ such that $D\cap domain(X)=D\cap domain(Y)$ and $X(\omega)\leq Y(\omega)$
for each $\omega\in D\cap domain(X)$. Similarly, there exists a full
set $D'$ such that $D'\cap domain(Y)=D'\cap domain(Z)$ and $Y(\omega)\leq Z(\omega)$
for each $\omega\in D'\cap domain(Y)$. By assertion 1, the set $DD'$
is a full set. Furthermore, $DD'\cap domain(X)=DD'\cap domain(Y)=DD'\cap domain(Z)$
and $X(\omega)\leq Y(\omega)\leq Z(\omega)$ for each $\omega\in DD'\cap domain(X)$
. It follows that $X\leq Z$ a.e. The remainder of the assertion is
similarly proved.

7. Trivial consequence of assertion 4. 
\end{proof}
\begin{defn}
\label{Def. Intgrable Set, Indicator, Measure, Complement} \textbf{(Integrable
set, measure of integrable set, complement of integrable set, and
null set).} A subset $A$ of $\Omega$ is called an \emph{integrable
set} \index{integrable set} if there exists an indicator $X$ which
is an integrable function such that $A=(X=1)$. In this case call
$X$ an indicator of $A$. We then define the meas\emph{ure} \index{measure}
of $A$ to be $\mu(A)\equiv IX$, and call the set $(X=0)$ a meas\emph{ure-theoretic
complement} \index{measure-theoretic complement} of $A$. We write
$1_{A}$ for an indicator of $A$\index{indicator of an integrable set},
and write $A^{c}$ for a measure-theoretic complement of $A$. An
integrable set with measure $\mu(A)=0$ is called a \emph{null set}.\index{null set}
$\square$
\end{defn}
Two distinct integrable indicators $X$ and $Y$ can be indicators
of the same integrable set $A$; hence $1_{A}$ is not uniquely defined
relative to the set-theoretic equality for functions. However, as
shown in the next proposition, given an integrable set, its indicator,
measure, and measure-theoretic complement are all uniquely defined
relative to a.e. equality.
\begin{prop}
\label{Prop. Property of Integrable Sets} \textbf{\emph{(Properties
of integrable sets).}} Let $A$ and $B$ be integrable sets. Let $X,Y$
be integrable indicators of $A,B$ respectively.
\end{prop}
\begin{enumerate}
\item \emph{$A=B$ a.e. iff $X=Y$ a.e. In particular, $1_{A}$ is well-defined
relative to a.e. equality, and the }meas\emph{ure $\mu(A)$ is well-defined.}
\item \emph{If $A=B$ a.e., then $(X=0)=(Y=0)$ a.e. In particular, $A^{c}$
is well-defined relative to equality a.e.}
\item \emph{The empty set $\phi$ is a null set, and $\Omega$ is a full
set.}
\item \emph{Any full set is a }meas\emph{ure-theoretic complement of a null
set.}
\item \emph{Any }meas\emph{ure-theoretic complement of a null set is a full
set.}
\item \emph{If $C$ is a subset of $\Omega$ such that $C=A$ a.e., then
$C$ is integrable with $\mu(A)=\mu(C)$.}
\end{enumerate}
\begin{proof}
By the definition of an indicator for an integrable set, we have $A=(X=1)$
and $B=(Y=1)$. Let $D$ be an arbitrary full set. Then the intersection
$D'\equiv D\cap domain(X)\cap domain(Y)$ is a full set. Since $D'D=D',$
we have $D'DA=D'D(X=1)$ and $D'DB=D'D(Y=1)$. 

1. Suppose $A=B$ on the full set $D$. Then $DA=DB$. It follows
from the previous paragraph that $D'(X=1)=D'D(X=1)=D'D(Y=1)=D'(Y=1)$.
By the remark following Definition \ref{Def. Indicators and Mutually Eclusive}we
see that for each $\omega\in D'$, $X(\omega)$ and $Y(\omega)$ are
defined and equal. Hence $X=Y$ a.e. Moreover, it follows from Proposition
\ref{Prop. Property of Full Sets} that $\mu(A)\equiv IX=IY\equiv\mu(B)$.
Conversely, suppose $X=Y$ a.e. with $D\cap domain(X)=D\cap domain(Y)$
and $X(\omega)=Y(\omega)$ for each $\omega\in D\cap domain(X)$.
Then $D'A=D'DA=D'D(X=1)$ $=D'D(Y=1)=D'(Y=1)=D'B$. Hence $A=B$ a.e.

2. Suppose $A=B$ a.e. In the above proof for assertion 1, we see
that for each $\omega$ in the full set $D'$, we have $X(\omega)=0$
iff $Y(\omega)=0$. 

3. Let $X$ be any integrable function. Then $0X$ is an indicator
for $\phi$, with $\mu(\phi)=I(0X)=0$. Hence $\phi$ is a null set.
Trivially $\Omega\supset domain(X)$, and so $\Omega$ is a full set.

4. Suppose $D$ is a full set. By Proposition \ref{Prop. Property of Full Sets},
we have $D=domain(X)$ for some integrable function $X$. Since $\phi=(0X=1)$
we see that $0X$ is an indicator for $\phi$. Hence $\phi^{c}\equiv(0X=0)=domain(X)=D$
is a measure-theoretic complement of the null set $\phi$.

5. Suppose $A$ is a null set, with $Z\in L$ as an indicator and
$A^{c}\equiv(Z=0)$. Since $IZ=0,$ the function $X\equiv\sum_{i=1}^{\infty}Z$
is integrable. Moreover,
\[
domain(X)=\{\omega\in domain(Z):\sum_{i=1}^{\infty}Z(\omega)=0\}=(Z=0)=A^{c}
\]
 Hence $A^{c}$ is a full set. 

6. Suppose $C=A$ on the full set $D$. Define a function $W$ by
$domain(W)=C\cup(X=0)$ and $W(\omega)=1$ or $0$ according as $\omega\in C$
or $(X=0)$. Then $W=X$ on the full set $D\cap domain(X)$. By Proposition
\ref{Prop. Property of Full Sets}, the function $W$ is integrable.
Hence $C=(W=1)$ has an integrable indicator. Moreover $\mu(C)=IW=IX=\mu(A)$. 
\end{proof}
Suppose each of a sequence of statements is valid a.e. Then in view
of Proposition \ref{Prop. Property of Full Sets}, there exists a
full set on which all of these statements are valid; in other words,
a.e. we have the validity of all the statements. For example if $(A_{n})_{n=1,2,\cdots}$
is a sequence of integrable sets with $A_{n}\subset A_{n+1}$ a.e.
for each $n>0$, then $A_{1}\subset A_{2}\subset\cdots$ a.e.
\begin{prop}
\label{Prop. Measure of Boolean Integrable Sets} \textbf{\emph{(Basics
of meaures of integrable sets)}}. Let $A,B$ be integrable sets, with
indicators $1_{A},1_{B}$ respectively, and with $A^{c}\equiv(1_{A}=0)$
and $B^{c}\equiv(1_{B}=0)$ . Then the following holds.
\end{prop}
\begin{enumerate}
\item \emph{$A\cup A^{c}$ is a full set, and $AA^{c}=\phi$.}
\item \emph{$A\cup B$ is an integrable set, with $1_{A\cup B}=1_{A}\vee1_{B}$
a.e.}
\item \emph{$AB$ is an integrable set, with $1_{AB}=1_{A}\wedge1_{B}$
a.e. Moreover $AB^{c}$ is an integrable set, with $1_{AB^{c}}=1_{A}-1_{A}\wedge1_{B}$
a.e. Furthermore $A(AB^{c})^{c}=AB$.}
\item \emph{$\mu(A\cup B)+\mu(AB)=\mu(A)+\mu(B)$. }
\item \emph{If $A\supset B$ a.e. then $\mu(AB^{c})=\mu(A)-\mu(B)$.}
\end{enumerate}
\begin{proof}
$\,$

1. We have $A=(1_{A}=1)$ and $A^{c}=(1_{A}=0)$. Hence $AA^{c}=\phi$.
Moreover $A\cup A^{c}=domain(1_{A})$, a full set.

2. Define the function $X$ by $domain(X)\equiv(A\cup B)\cup(A^{c}B^{c})$
and $X(\omega)\equiv1$ or $0$ according as $\omega\in A\cup B$
or $\omega\in A^{c}B^{c}$. Then $X=1_{A}\vee1_{B}$ on the full set
$domain(1_{A}\vee1_{B})$. Hence $X$ is an integrable function according
to \ref{Prop. Property of Full Sets}. Since $A\cup B=(X=1)$, the
function $X$ is an indicator of $A\cup B$. In other words $1_{A\cup B}=X=1_{A}\vee1_{B}$
a.e. 

3. Obviously $AB=(1_{A}\wedge1_{B}=1)$. Hence $1_{AB}=1_{A}\wedge1_{B}$.
Next define the function $X$ by $domain(X)\equiv(AB^{c})\cup(A^{c}\cup B)$
and $X(\omega)\equiv1$ or $0$ according as $\omega\in AB^{c}$ or
$\omega\in A^{c}\cup B$. Then $X=1_{A}-1_{A}\wedge1_{B}$ on the
full set $domain(1_{A}\wedge1_{B})$. Hence $X$ is an integrable
function according to \ref{Prop. Property of Full Sets}. Since $AB^{c}=(X=1)$,
the function $X$ is an indicator of $AB^{c}$. In other words $1_{AB^{c}}=X=1_{A}-1_{A}\wedge1_{B}$
a.e. Furthermore, $A(AB^{c})^{c}=A(X=0)=A(1_{A}\wedge1_{B}=1)=AB$.

4. Since $1_{A}\vee1_{B}+1_{A}\wedge1_{B}=1_{A}+1_{B}$, the conclusion
follows from linearity of $I$.

5. Suppose $AD\supset BD$ for a full set $D$. Write $A'\equiv AD$
and $B'\equiv BD$ and define $B'^{c}\equiv B^{c}D$. Then $B'^{c}$
is a measure-theoretic complement of $B'$. We have $A'\supset B'$
and so $B'=B'A'$. According to Proposition \ref{Prop. Property of Integrable Sets},
the sets $A'$ and $B'$ are integrable, with $\mu(A')=\mu(A)$ and
$\mu(B')=\mu(B)$. By the same token, since $A'B'^{c}=AB^{c}$ on
the full set $D$, the set $A'B'^{c}$ is integrable with $\mu(A'B'^{c})=\mu(AB^{c})$.
On the other hand, from assertion 4, we have
\[
\mu(B')+\mu(A'B'^{c})=\mu(B'A')+\mu(B'^{c}A')
\]
\[
=\mu((B'A')\cup(B'^{c}A'))=\mu((B'\cup B'^{c})A')=\mu(A')
\]
 where the next to last equality is because $(B'^{c}\cup B')A'=A'$
on the full set $B'^{c}\cup B'$. The assertion is proved. 
\end{proof}
\begin{prop}
\emph{\label{Prop. Sequences of intgrable sets}} \textbf{\emph{(Sequence
of integrable sets).}} For each $n\geq1$ let $A_{n}$ be an integrable
set with a measure-theoretic complement $A_{n}^{c}$. Then the following
holds.
\end{prop}
\begin{enumerate}
\item \emph{If $A_{n}\subset A_{n+1}$ a.e. for each $n\geq1$, and if $\mu(A_{n})$
converges, then $\bigcup_{n=1}^{\infty}A_{n}$ is an integrable set
with $\mu(\bigcup_{n=1}^{\infty}A_{n})=\lim_{n\rightarrow\infty}\mu(A_{n})$
and $(\bigcup_{n=1}^{\infty}A_{n})^{c}=\bigcap_{n=1}^{\infty}A_{n}^{c}$.}
\item \emph{If $A_{n}\supset A_{n+1}$ a.e. for each $n\geq1$, and if $\mu(A_{n})$
converges, then $\bigcap_{n=1}^{\infty}A_{n}$ is an integrable set
with $\mu(\bigcap_{n=1}^{\infty}A_{n})=\lim_{n\rightarrow\infty}\mu(A_{n})$
and $(\bigcap_{n=1}^{\infty}A_{n})^{c}=\bigcup_{n=1}^{\infty}A_{n}^{c}$.}
\item \emph{If $A_{n}A_{m}=\phi$ a.e. for each $n>m\geq1,$ and if $\sum_{n=1}^{\infty}\mu(A_{n})$
converges, then $\bigcup_{n=1}^{\infty}A_{n}$ is an integrable set
with $\mu(\bigcup_{n=1}^{\infty}A_{n})=\sum_{n=1}^{\infty}\mu(A_{n})$.}
\item \emph{If $\sum_{n=1}^{\infty}\mu(A_{n})$ converges, then $\bigcup_{n=1}^{\infty}A_{n}$
is an integrable set with $\mu(\bigcup_{n=1}^{\infty}A_{n})\leq\sum_{n=1}^{\infty}\mu(A_{n})$
.}
\end{enumerate}
\begin{proof}
For each $n\geq1$ let $1_{A_{n}}$be the integrable indicator of
$A_{n}$ such that $A_{n}^{c}=(1_{A_{n}}=0)$.

1. Define a function $Y$ by
\[
domain(Y)\equiv(\bigcup_{n=1}^{\infty}A_{n})\cup(\bigcap_{n=1}^{\infty}A_{n}^{c})
\]
with $Y(\omega)\equiv1$ or $0$ according as $\omega\in\bigcup_{n=1}^{\infty}A_{n}$
or $\omega\in\bigcap_{n=1}^{\infty}A_{n}^{c}$. Then $(Y=1)=\bigcup_{n=1}^{\infty}A_{n}$
and $(Y=0)=\bigcap_{n=1}^{\infty}A_{n}^{c}$. For each $n\geq1$,
we have $A_{n}\subset A_{n+1}$ a.e. and so $1_{A_{n+1}}\geq1_{A_{n}}$
a.e. By assumption we have the convergence of 
\[
I(1_{A_{1}})+I(1_{A_{2}}-1_{A_{1}})+\cdots+I(1_{A_{n}}-1_{A_{n-1}})=I1_{A_{n}}=\mu(A_{n})
\]
as $n\rightarrow\infty$. Hence $X\equiv1_{A_{1}}+(1_{A_{2}}-1_{A_{1}})+(1_{A_{3}}-1_{A_{2}})+\cdots$
is an integrable function. Consider an arbitrary $\omega\in domain(X)$.
The limit
\[
\lim_{n\rightarrow\infty}1_{A_{n}}(\omega)=\lim_{n\rightarrow\infty}(1_{A_{1}}(\omega)+(1_{A_{2}}-1_{A_{1}})(\omega)+\cdots+(1_{A_{n}}-1_{A_{n-1}})(\omega))=X(\omega)
\]
exists, and is either 0 or 1 since it is the limit of a sequence in
$\{0,1\}$. Suppose $X(\omega)=1$. Then $1_{A_{n}}(\omega)=1$ for
some $n\geq1$. Hence $\omega\in\bigcup_{n=1}^{\infty}A_{n}$ and
so $Y(\omega)\equiv1=X(\omega)$. Suppose $X(\omega)=0$. Then $1_{A_{n}}(\omega)=0$
for each $n\geq1$. Hence $\omega\in\bigcap_{n=1}^{\infty}A_{n}^{c}$
and so $Y(\omega)\equiv0=X(\omega)$. Combining, we see that $Y=X$
on the full set $domain(X)$. According to Proposition \ref{Prop. Property of Full Sets},
we therefore have $Y\in L$. Thus $\bigcup_{n=1}^{\infty}A_{n}=(Y=1)$
is an integrable set with $Y$ as its indicator, and has measure equal
to
\[
IY=IX=\lim_{n\rightarrow\infty}I1_{A_{n}}=\lim_{n\rightarrow\infty}\mu(A_{n})
\]
Moreover $(\bigcup_{n=1}^{\infty}A_{n})^{c}=(Y=0)=\bigcap_{n=1}^{\infty}A_{n}^{c}$. 

2. Similar.

3. Write $B_{n}=\bigcup_{i=1}^{n}A_{i}$. Repeated application of
Proposition \ref{Prop. Measure of Boolean Integrable Sets} leads
to $\mu(B_{n})=\sum_{i=1}^{n}\mu(A_{i})$. From assertion 1 we see
that $\bigcup_{n=1}^{\infty}A_{n}=\bigcup_{n=1}^{\infty}B_{n}$ is
an integrable set with $\mu(\bigcup_{n=1}^{\infty}A_{n})=\mu(\bigcup_{n=1}^{\infty}B_{n})=\lim_{n\rightarrow\infty}\mu(B_{n})=\sum_{i=1}^{\infty}\mu(A_{i}).$

4. Define $B_{1}=A_{1}$ and $B_{n}=(\bigcup_{k=1}^{n}A_{k})(\bigcup_{k=1}^{n-1}A_{k})^{c}$
for $n>1$. Let $D$ denote the full set $\bigcap_{k=1}^{\infty}(A_{k}\cup A_{k}^{c})(B_{k}\cup B_{k}^{c})$.
Clearly $B_{n}B_{k}=\phi$ on $D$ for each positive integer $k<n$.
This implies $\mu(B_{n}B_{k})=0$ for each positive integer $k<n$.
Furthermore, for every $\omega\in D$, we have $\omega\in\bigcup_{k=1}^{\infty}A_{k}$
iff there is a smallest $n>0$ such that $\omega\in\bigcup_{k=1}^{n}A_{k}$.
Since for every $\omega\in D$ either $\omega\in A_{k}$ or $\omega\in A_{k}^{c},$
we have $\omega\in\bigcup_{k=1}^{\infty}A_{k}$ iff there is an $n>0$
such that $\omega\in B_{n}$. In other words $\bigcup_{k=1}^{\infty}A_{k}=\bigcup_{k=1}^{\infty}B_{k}$
a.e. Moreover $\mu(B_{n})=\mu(\bigcup_{k=1}^{n}A_{k})-\mu(\bigcup_{k=1}^{n-1}A_{k})$.
Hence the sequence $(B_{n})$ of integrable sets satisfies the hypothesis
in assertion 3. Therefore $\bigcup_{k=1}^{\infty}B_{k}$ is an integrable
set, with
\[
\mu(\bigcup_{k=1}^{\infty}A_{k})=\mu(\bigcup_{k=1}^{\infty}B_{k})=\lim_{n\rightarrow\infty}\sum_{k=1}^{n}\mu(B_{k})\leq\lim_{n\rightarrow\infty}\sum_{k=1}^{n}\mu(A_{k})=\sum_{n=1}^{\infty}\mu(A_{n}).
\]
\end{proof}
\begin{prop}
\label{Prop. Xn->X in L1 implies subseq -> X a.s} \textbf{\emph{(Convergence
in $L$ implies an a.e. convergent subsequence). }}Let $X\in L$ and
let $(X_{n})_{n=1,2,\cdots}$ be a sequence in $L$. If $I|X_{n}-X|\rightarrow0$
then there exists a subsequence $(Y_{n})_{n=1,2,\cdots}$ such that
$Y_{n}\rightarrow X$ a.e.
\end{prop}
\begin{proof}
Let $(Y_{n})_{n=1,2,\cdots}$ be a subsequence such that $I|Y_{n}-X|<2^{-n}$.
Then the sequence $(Z_{n})_{n=1,2,\cdots}$ defined as $(X,-X+Y_{1},X-Y_{1},-X+Y_{2},X-Y_{2},\cdots)$
is a representation of $X$. Define $Z\equiv\sum_{n=1}^{\infty}Z_{n}\in L$
On the full set $domain(Z)$, we then have $Y_{n}=(Z_{1}+\cdots Z_{2n})\rightarrow X$. 
\end{proof}
We will use the next theorem many times to construct integrable functions.
\begin{thm}
\label{Thm. |X-Yn|<Zn where Yn,Zn integrable =000026 IZn->0 implies X integrable}
\textbf{\emph{(A sufficient condition for a functiion to be integrable).
}}Suppose $X$ is a function defined a.e. on $\Omega$. Suppose there
exist two sequences $(Y_{n})_{n=1,2,\cdots}$ and $(Z_{n})_{n=1,2,\cdots}$
in $L$ such that $|X-Y_{n}|\leq Z_{n}$ a.e. for each $n\geq1$ and
such that $IZ_{n}\rightarrow0$. Then $X\in L$. Moreover $I|X-Y_{n}|\rightarrow0$.
\end{thm}
\begin{proof}
According to Proposition \ref{Prop. Xn->X in L1 implies subseq -> X a.s},
by passing to a subsequence, we can assume that $Z_{n}\rightarrow0$
a.e. Since, by assumption, $|X-Y_{n}|\leq Z_{n}$ a.e. for each $n\geq1$,
it follows that $Y_{n}\rightarrow X$ a.e. On the other hand, we have
$|Y_{n}-Y_{m}|\leq|Y_{n}-X|+|Y_{m}-X|\leq Z_{n}+Z_{m}$. Consequently
$I|Y_{n}-Y_{m}|\leq IZ_{n}+IZ_{m}\rightarrow0$ as $n,m\rightarrow\infty$.
By the completeness of $L$, there exists $Y\in L$ such that $I|Y_{n}-Y|\rightarrow0$.
By passing again to a subsequence, we may assume that $Y_{n}\rightarrow Y$
a.e. Combining, we see that $X=Y$ a.e. According to Proposition \ref{Prop. Property of Full Sets},
we therefore have $X\in L$. Moreover, $I|X-Y_{n}|\leq IZ_{n}\rightarrow0$
\end{proof}

\section{Abundance of Integrable Sets }

In this section let $(\Omega,L,I)$ be a complete integration space.

Let $X$ be any function defined on a subset of $\Omega$ and let
$t$ be a real number. Recall from the Notations and Conventions in
the Introduction that we use the abbreviation $(t\leq X)$ for the
subset $\{\omega\in domain(X):t\leq X(\omega)\}$ Similar notations
are used for $(X<t)$, $(X\leq t)$ and $(X<t)$. We will also write
$(t<X\leq u)$ etc for the intersection $(t<X)(X\leq u)$ etc. If
$J$ is a subset of $R$, let $J_{c}$ denote the metric complement
of $J$ in $R$. 

We will show in this section that if $X$ is an integrable function,
then $(t\leq X)$ and $(t<X)$ are integrable sets for each positive
$t$ in the metric complement of some countable subset of $R$. 

Define some functions which will serve as approximations for step
functions on $R$. For any real numbers $0<s<t$ define $g_{s,t}(x)\equiv\frac{x\wedge t-x\wedge s}{t-s}$
. Then the function $g_{s,t}(X)\equiv\frac{X\wedge t-X\wedge s}{t-s}$
is integrable for all $s,t\in R$ with $0<s<t$. Clearly $1\geq g_{t',t}\geq g_{s,s'}\geq0$
for all $t',t,s,s'\in R$ with $t'<t\leq s<s'$,. If we can prove
that $\lim_{s\uparrow t}Ig_{s,t}(X)$ exists, then we can use the
Monotone Convergence Theorem to show that $\lim_{s\uparrow t}g_{s,t}(X)$
is integrable and is an indicator of $(t\leq X)$, proving that the
latter set is integrable. Classically the existence of $\lim_{s\uparrow t}Ig_{s,t}(X)$
is trivial since for fixed $t$ the integral $Ig_{s,t}(X)$ is nonincreasing
in $s$ and bounded from below by $0$. A constructive proof that
the limit exists for all but countably many $t$'s is given below.
The proof is in terms of a general theory of profiles which finds
applications also outside measure or integration theory.
\begin{defn}
\label{Def. Profile}\textbf{(Profile). }Let $K$ be a non-empty open
interval in $R$. Let $G$ be a family of continuous functions on
$R$, such that $0\leq g\leq1$ for each $g\in G$. Let $t\in K$
and $g\in G$ be arbitrary. We say $t$ \emph{precedes} $g$ and write
$t\diamondsuit g$ if $g=0$ on $(-\infty,t]\cap K$. We say $g$
\emph{precedes} $t$ and write $g\diamondsuit t$ if $g=1$ on $[t,\infty)\cap K$.
We write $t\diamondsuit g\diamondsuit s$ and say $g$ \emph{separates}\index{separation}
$t$ and $s$ if both $t\diamondsuit g$ and $g\diamondsuit s$. We
say $G$ separates points in $K$ if for all $t,s\in K$ with $t<s$
 there exists $g\in G$ such that $t\diamondsuit g\diamondsuit s$.
A function $\lambda$ on $G$ is said to be nondecreasing if for each
$g,g'$ with $g\leq g'$ on $K$ we have $\lambda(g)\leq\lambda(g'$).
We say $(G,\lambda)$ is a \index{profile}\emph{profile} on the interval
$K$ if $G$ separates points in $K$ and if $\lambda$ is a nondecreasing
function on $G$. We say that a closed interval $[t,s]\subset K$
has a positive real number $\alpha$ as a \emph{profile bound}\index{profile bound},
and write $[t,s]\ll\alpha$, if there exist $t',s'\in K$ and $f,g\in G$
such that (i) $f\diamondsuit t'$, $t'<t\leq s<s'$, $s'\diamondsuit g$,
and (ii) $\lambda(f)-\lambda(g)<\alpha$. Suppose $a,b\in R$ and
$a\leq b$. We say that the open interval $(a,b)\subset K$ has a
positive real number $\alpha$ as a profile bound, and write $(a,b)\ll\alpha$
if $[t,s]\ll\alpha$ for each closed subinterval $[t,s]$ of $(a,b)$.
Note that the open interval $(a,b)$, defined as the set $\{x\in R:a<x<b\}$,
can be empty. $\square$
\end{defn}
Note that $t\diamondsuit g$ is merely an abbreviation for $1_{[t,\infty)}\geq g$;
and $g\diamondsuit t$ is an abbreviation for $g\geq1_{[t,\infty)}$. 

The motivating example of a profile is when $K\equiv(0,\infty)$,$G\equiv\{g_{s,t}:s,t\in K\mbox{ and }0<s<t\}$,
and the function $\lambda$ is defined on $G$ by $\lambda(g)\equiv Ig(X)$
for each $g\in G$. It can easily be verified that $(G,\lambda)$
is a profile on $K$. 

In the following let $(G,\lambda)$ be a general profile on an open
interval $K$ in $R$. The next lemma lists some basic properties.
\begin{lem}
\label{Lem. Profile Basics} \textbf{\emph{(Basics of profiles)}}. 
\end{lem}
\begin{enumerate}
\item \emph{If $f\diamondsuit t$, $t\leq s$, and $s\diamondsuit g$ then
$f\geq g$ and $\lambda(f)\geq\lambda(g)$.}
\item \emph{If $t\leq s$ and $s\diamondsuit g$ then $t\diamondsuit g$. }
\item \emph{If $g\diamondsuit t$ and $t\leq s$ then $g\diamondsuit s$.}
\item \emph{In view of the transitivity in assertions 2 and 3 above, we
can rewrite, without ambiguity, condition (i) in Definition \ref{Def. Profile}
as $f\diamondsuit t'<t\leq s<s'\diamondsuit g$.}
\item \emph{Suppose $[t,s]\ll\alpha$ and $t_{0}<t\leq s<s_{0}$. Let $\varepsilon>0$
be arbitrary. Then there exist $t_{1},s_{1}\in K$ and $f_{1},g_{1}\in G$
such that (i) $t_{0}\diamondsuit f_{1}\diamondsuit t_{1}<t\leq s<s_{1}\diamondsuit g_{1}\diamondsuit s_{0}$,
(ii) $\lambda(f_{1})-\lambda(g_{1})<\alpha$, and (iii) $t-\varepsilon<t_{1}<t$
and $s<s_{1}<s+\varepsilon$.}
\item \emph{Every closed sub-interval of $K$ has a finite profile bound. }
\end{enumerate}
\begin{proof}
We will prove assertions 5 and 6, the rest being trivial. 

Suppose $[t,s]\ll\alpha$ and $t_{0}<t\leq s<s_{0}$. Let $\varepsilon>0$
be arbitrary. Then there exist $t',s'\in K$ and $f,g\in G$ such
that (i') $f\diamondsuit t'<t\leq s<s'\diamondsuit g$, and (ii')
$\lambda(f)-\lambda(g)<\alpha$. Pick real numbers $t'',t_{1}$ such
that $t_{0}\vee t'\vee(t-\varepsilon)<t''<t_{1}<t$. Since $G$ separates
points in $K$, there exists $f_{1}\in G$ such that $t_{0}<t''\diamondsuit f_{1}\diamondsuit t_{1}<t$.
Since $f\diamondsuit t'<t''\diamondsuit f_{1}$, we have, in view
of assertion 1, $\lambda(f)\geq\lambda(f_{1})$. Similarly we obtain
$s'',s_{1}\in K$ and $g_{1}\in G$ such that $s<s_{1}<s''<s_{0}\wedge s'\wedge(s+\varepsilon)$
and $s<s_{1}\diamondsuit g_{1}\diamondsuit s''<s_{0}$ with $\lambda(g_{1})\geq\lambda(g)$.
Hence $\lambda(f_{1})-\lambda(g_{1})\leq\lambda(f)-\lambda(g)<\alpha$.
Conditions (i) and (iii) are obviously also satisfied. Assertion 5
is proved.

Given any interval $[t,s]\subset K$, let $t'',t',s',s''$ be members
of $K$ such that $t''<t'<t\leq s<s'<s''$. Since G separates points
in $K$, there exist $f,g\in G$ such that $t''\diamondsuit f\diamondsuit t'<t\leq s<s'\diamondsuit g\diamondsuit s''$.
Hence $[t,s]\ll\alpha$ for any real number $\alpha$ such that $\lambda(f)-\lambda(g)<\alpha$.
Assertion 6 is proved. 
\end{proof}
\begin{lem}
\label{Lem. Low profile except q points} \textbf{\emph{(Bound for
the number of intervald with significant profiles). }}Let $(G,\lambda)$
be a profile on a proper open interval $K$ in $R$. Let $[a,b]$
be a closed sub-interval of $K$ with $[a,b]\ll\alpha$. Let $\varepsilon>0$
be arbitrary. Let $q$ be any integer with $q\geq\alpha/\varepsilon$.
Then there exists a sequence $s_{0}=a\leq s_{1}\leq\cdots\leq s_{q}=b$
of $(q+1)$ points in $K$ such that $(s_{k-1},s_{k})\ll\varepsilon$
for each $k=1,\cdots q$.
\end{lem}
\begin{proof}
For abbreviation write $d_{n}\equiv2^{-n}(b-a)$ for each $n\geq1$.
By hypothesis $[a,b]\ll\alpha$. Hence there exist $a',b'\in K$ and
$f',f''\in G$ such that (i) $f'\diamondsuit a'<a\leq b<b'\diamondsuit f''$,
and (ii) $\lambda(f')-\lambda(f'')<\alpha\leq q\varepsilon$. Define$F'\equiv\lambda(f')$
and $F''\equiv\lambda(f'')$. Then $0\leq F'-F''<q\varepsilon$. 

Let $n\geq1$ be arbitrary. For $i=0,\cdots,2^{n}$ define $t_{n,i}\equiv a+id_{n}$.
Clearly $t_{n,0}=a$ and $t_{n,2^{n}}=b$. Define $D_{n}\equiv\{t_{n,i}:0\leq i\leq2^{n}\}$.
The set $D_{n}$ is a result of binary subdivisions of the interval
$[a,b]$. Specifically, consider any $n\geq1$ and $t=t_{n,i}\in D_{n}$.
Then we have 
\[
t_{n,i}\equiv c+id_{n}=c+2id_{n+1}\equiv t_{n+1,2i}\in D_{n+1}
\]
Hence $D_{n}\subset D_{n+1}.$ For each $1\leq i\leq2^{n}$, we have
$t_{n,i-1},t_{n,i}\in[a,b]\subset K$ and so there exists a function
$f_{n,i}\in G$ with $t_{n,i-1}\diamondsuit f_{n,i}\diamondsuit t_{n,i}$.
In addition, define $f_{n,0}\equiv f'$ and $f_{n,2^{n}+1}\equiv f''$.
Then $f_{n,0}\diamondsuit t_{n,0}\equiv a'$ and $b\equiv t_{n,2^{n}}\diamondsuit f_{n,2^{n}+1}$.
Combining, we have
\begin{equation}
f_{n,0}\diamondsuit t_{n,0}\diamondsuit f_{n,1}\diamondsuit t_{n,1}\cdots\diamondsuit f_{n,2^{n}}\diamondsuit t_{n,2^{n}}\diamondsuit f_{n,2^{n}+1}\label{eq:temp-15}
\end{equation}
Next, for each $t=t_{n,i}\in D_{n}$ define $F_{n}(t)\equiv F_{n}(t_{n,i})\equiv\lambda(f_{n,i})$.
By the relation \ref{eq:temp-15} and Lemma \ref{Lem. Profile Basics},
we see that $F_{n}$ is a nonincreasing function on $D_{n}$: 
\[
\mbox{ if s,\ensuremath{t\in D_{n}}are such that }s\leq t\mbox{, then }F_{n}(s)\geq F_{n}(t)
\]
.

Next let $t=t_{n,i}\in D_{n}$ and let $\ensuremath{s=t_{n+1,j}\in D_{n+1}}$.
Suppose $s\leq t-d_{n}$. Then $t>a$. Consequently $i\geq1$ and
$t_{n,i-1}=t-d_{n}\geq s$. Hence $f_{n+1,j}\diamondsuit t_{n+1,j}=s\leq t_{n,i-1}\diamondsuit f_{n,i}$.
Hence, by Lemma \ref{Lem. Profile Basics}, we have $F_{n+1}(s)\equiv\lambda(f_{n+1,j})\geq\lambda(f_{n,i})\equiv F_{n}(t$).
We have thus proved that 
\begin{equation}
\mbox{ if \ensuremath{t\in D_{n}}and \ensuremath{s\in D_{n+1}}are such that }s\leq t-d_{n}\mbox{, then }F_{n+1}(s)\geq F_{n}(t)\label{eq:temp-21}
\end{equation}
Similarly we can prove that 
\begin{equation}
\mbox{ if \ensuremath{t\in D_{n}}and \ensuremath{s\in D_{n+1}}are such that }t+d_{n+1}\leq s,\mbox{ then }F_{n}(t)\geq F_{n+1}(s)\label{eq:temp-22}
\end{equation}

For each $n\geq1$ and integer $i$ with $0\leq i\leq2^{n}$ define
$S_{n,i}\equiv F_{n}(t_{n,0})-F_{n}(t_{n,i})=F'-F_{n}(t_{n,i})$.
Pick an $\varepsilon'>0$ such that 
\begin{equation}
q\varepsilon>q\varepsilon'>F'-F''\label{eq:temp-19}
\end{equation}
and such that $|\varepsilon'-S_{n,i}/k|>0$ for each $i$ with $0\leq i\leq2^{n}$,
for each $k=1,\cdots,q$, and for each $n\geq1$. Now let $n\geq1$
be arbitrary. Define also $S_{n,2^{n}+1}\equiv q\varepsilon'$. Then
we have
\begin{equation}
0=S_{n,0}\leq S_{n,1}\leq\cdots\leq S_{n,2^{n}+1}\equiv q\varepsilon'\label{eq:temp-14}
\end{equation}
Consider each $k=1,\cdots,q$. From inequality \ref{eq:temp-14} we
see that $S_{n,0}<k\varepsilon'\leq S_{n,2^{n}+1}$. Hence there exists
an integer$i_{n,k}$ with $0\leq i_{n,k}\leq2^{n}$ such that 
\begin{equation}
\mbox{ }S_{n,i_{n,k}}<k\varepsilon'\leq S_{n,i_{n,k}+1}\label{eq:temp-13}
\end{equation}
Define $s_{n,k}\equiv t_{n,i_{n,k}}\in[a,b]$. Clearly $s_{n,k}\leq s_{n,k+1}$
for $k<q$. Moreover $S_{n,2^{n}}<q\varepsilon'\equiv S_{n,2^{n}+1}$.
Hence $i_{n,q}=2^{n}$ and
\begin{equation}
s_{n,q}\equiv t_{n,i_{n,q}}=t_{n,2^{n}}=b\label{eq:temp-8}
\end{equation}

Fix any $k=1,\cdots,q$. We will show that the sequence $(s_{n,k})_{n=1,2,\cdots}$
converges. Consider the terms $s_{n,k}$ and $s_{n+1,k}$. For ease
of notations, write $i=i_{n,k}$ and $j=i_{n+1,k}$.

Suppose $s_{n,k}<s_{n+1,k}$. Then $t_{n,i}=s_{n,k}<b$ and so $i\leq2^{n}-1$.
It follows that $S_{n,i+1}\equiv F'-F_{n}(t_{n,i+1})$. By the definition
for $i_{n,k}$ and $i_{n+1,k}$ we have 
\[
F'-F_{n+1}(t_{n+1,j})\equiv S_{n+1,j}<k\varepsilon'\leq S_{n,i+1}\equiv F'-F_{n}(t_{n,i+1})
\]
whence
\[
F_{n}(t_{n,i+1})<F_{n+1}(t_{n+1,j})
\]
This implies, in view of inference \ref{eq:temp-22}, that $t_{n,i+1}+d_{n+1}\geq t_{n+1,j}$.
Hence
\[
s_{n+1,k}=t_{n+1,j}\leq t_{n,i+1}+d_{n+1}=t_{n,i}+d_{n}+d_{n+1}=s_{n,k}+d_{n}+d_{n+1}
\]
Therefore $s_{n+1,k}-s_{n,k}\leq d_{n}+d_{n+1}$.

On the other hand, suppose $s_{n,k}>s_{n+1,k}$. Then $t_{n+1,j}=s_{n+1,k}<b$
and so $j\leq2^{n+1}-1$. It follows that $S_{n+1,j+1}\equiv F'-F_{n+1}(t_{n+1,j+1})$.
By the definition for $i_{n,k}$ and $i_{n+1,k}$ we have
\[
F'-F_{n}(t_{n,i})\equiv S_{n,i}<k\varepsilon'\leq S_{n+1,j+1}\equiv F'-F_{n+1}(t_{n+1,j+1})
\]
whence
\[
F_{n+1}(t_{n+1,j+1})<F_{n}(t_{n,i})
\]
This implies, in view of inference \ref{eq:temp-21}, that $t_{n+1,j+1}\geq t_{n,i}-d_{n}.$
Hence
\[
s_{n+1,k}=t_{n+1,j}=t_{n+1,j+1}-d_{n+1}\geq t_{n,j}-d_{n}-d_{n+1}=s_{n,k}-d_{n}-d_{n+1}
\]
Therefore $s_{n,k}-s_{n+1,k}\leq d_{n}+d_{n+1}$.

Combining, we obtain $|s_{n,k}-s_{n+1,k}|\leq d_{n}+d_{n+1}\equiv3\cdot2^{-n-1}(b-a).$
Thus we see that the sequence $(s_{n,k})_{n=1,2,\cdots}$ is Cauchy,
and converges to some $s_{k}\in[a,b]$. Furthermore, for $k=1,\cdots,q$,
\[
|s_{n,k}-s_{k}|=\lim_{p\rightarrow\infty}|(s_{n,k}-s_{n+1,k})+(s_{n+1,k}-s_{n+2,k})+\cdots+(s_{p-1,k}-s_{p,k})|
\]
\begin{equation}
\leq3(2^{-n-1}+2^{-n}+\cdots)(b-a)\leq3\cdot2^{-n}(b-a)=3d_{n}\label{eq:temp-18}
\end{equation}
For ease of notations, we will also define $s_{n,0}\equiv s_{0}\equiv a$
and $s_{n,q+1}\equiv s_{q+1}\equiv b$. Then $s_{k}\leq s_{k+1}$
for $0\leq k\leq q$.

Now let $k=0,\cdots,q$ be arbitrary. Suppose $[t,s]\subset(s_{k},s_{k+1})$
for some real numbers $t\leq s$. We will show that $[t,s]\ll\varepsilon$.
To this end, let $n$ be so large that $s_{k}+4d_{n}<t\leq s<s_{k+1}-4d_{n}.$
This implies, in view of inequality \ref{eq:temp-18}, that $s_{n,k}+d_{n}<t\leq s<s_{n,k+1}-d_{n}$.
For abbreviation write $i\equiv i_{n,k}$ and $j\equiv i_{n,k+1}$.
According to the definition of $i_{n,k}$ and $i_{n,k+1}$ we then
have $k\varepsilon'\leq S_{n,i+1}$ and $S_{n,j}<(k+1)\varepsilon'$.
Moreover $t_{n,i}=s_{n,k}<s_{n,k+1}\leq b$. In view of equality \ref{eq:temp-8},
we have $k<q$ and $i<2^{n}$. Hence $t_{n,i+1}=t_{n,i}+d_{n}=s_{n,k}+d_{n}<t$.
Similarly, we have $S_{n,j}<(k+1)\varepsilon'\leq S_{n,j+1}$. Moreover
$t_{n,j}=s_{n,k+1}>s_{n,k}\geq a$. We therefore have $j>0$. Hence
$t_{n,j-1}=t_{n,j}-d_{n}=s_{n,k+1}-d_{n}>s$. Combining, we have $f_{n,i+1}\diamondsuit t_{n,i+1}<t\leq s<t_{n,j-1}\diamondsuit f_{n,j}$.
Furthermore, $S_{n,j}-S_{n,i+1}<(k+1)\varepsilon'-k\varepsilon'=\varepsilon'$.
Equivalently $F_{n}(t_{n,i+1})-F_{n}(t_{n,j})<\varepsilon'$, which
is in turn equivalent to $\lambda(f_{n,i+1})-\lambda(f_{n,j})<\varepsilon'$.
Therefore $[t,s]\ll\varepsilon'<\varepsilon$. Since $[t,s]$ is an
arbitrary closed sub-interval of $(s_{k},s_{k+1})$, we have proved
that $(s_{k},s_{k+1})\ll\varepsilon$. 
\end{proof}
\begin{thm}
\label{Thm. Profile Smooth}\textbf{\emph{ (All but countably many
points have arbitrarily low profile). }}Let $(G,\lambda)$ be a profile
on a proper open interval $K$. Then there exists a countable subset
$J$ of $K$ such that for each $t\in K\cap J_{c}$ we have $[t,t]\ll\varepsilon$
for arbitrarily small $\varepsilon>0$.
\end{thm}
\begin{proof}
Let $[a,b]\subset[a_{2},b_{2}]\subset\cdots$ be a sequence of subintervals
of $K$ such that $K=\bigcup_{p=1}^{\infty}[a_{p},b_{p}]$. According
to Lemma \ref{Lem. Low profile except q points}, there exists for
each $p\geq1$, a finite sequence $s_{0}^{(p)}=a_{p}\leq s_{1}^{(p)}\leq\cdots\leq s_{q_{p}}^{(p)}=b_{p}$
such that $(s_{k-1}^{(p)},s_{k}^{(p)})\ll\frac{1}{p}$ for $k=1,\cdots,q_{p}$.
Define $J\equiv\{s_{k}^{(p)}:1\leq k\leq q_{p};p\geq1\}$. Suppose
$t\in K\cap J_{c}$. Let $\varepsilon>0$ be arbitrary. Let $p\geq1$
be so large that $t\in[a_{p},b_{p}]$ and that $\frac{1}{p}<\varepsilon$.
By the definition of the metric complement $J_{c}$ we have $|t-s_{k}^{(p)}|>0$
for each $k=1,\cdots,q_{p}$. Hence $t\in(s_{k-1}^{(p)},s_{k}^{(p)})$
for some $k=1,\cdots,q_{p}$. But $(s_{k-1}^{(p)},s_{k}^{(p)})\ll\frac{1}{p}$.
We conclude that $[t,t]\ll\frac{1}{p}<\varepsilon$. 
\end{proof}
An immediate application the preceding Theorem \ref{Lem. Low profile except q points}
is to establish the abundance of integrable sets, in the next theorem.
\begin{thm}
\label{Thm. (t<X) is integrable} \textbf{\emph{(Abundance of integrable
sets).}} Given an integrable function $X$ on the complete integration
space $(\Omega,L,I)$, there exists a countable subset $J$ of $(0,\infty)$
such that for each positive real number $t$ in the metric complement
$J_{c}$ of $J$, the sets $(t\leq X)$ and $(t<X)$ are integrable
sets, with $(t\leq X)^{c}=(X<t)$ and $(t<X)^{c}=(X\leq t)$. Furthermore,
the measures $\mu(t\leq X)$ and $\mu(t<X)$ are equal and are continuous
at each $t>0$ with $t\in J_{c}$.
\end{thm}
\begin{proof}
Recall the previously defined profile $(G,\lambda)$ on the interval
$K\equiv(0,\infty)$, where $G\equiv\{g_{s,t}:s,t\in K\mbox{ and }0<s<t\}$,
and $\lambda(g)\equiv Ig(X)$ for each $g\in G$. Here $g_{s,t}$
denotes the function defined on $R$ by $g_{s,t}(x)\equiv\frac{x\wedge t-x\wedge s}{t-s}$
for each $x\in R$. Let the countable subset $J$ of $K$ be constructed
as in Theorem \ref{Thm. Profile Smooth}. 

Suppose $t\in K\cap J_{c}$. We have $[t,t]\ll\frac{1}{p}$ for each
$p\geq1$. Recursively applying Lemma \ref{Lem. Profile Basics},
we can construct two sequences $(u_{p})_{p=0,1,\cdots}$ and $(v_{p})_{p=0,1,\cdots}$
in $K$, and two sequences $(f_{p})_{p=1,2,\cdots}$ and $(g_{p})_{p=1,2,\cdots}$
in $G$ such that for each $p\geq1$ we have (i) $u_{p-1}\diamondsuit f_{p}\diamondsuit u_{p}<t<v_{p}\diamondsuit g_{p}\diamondsuit v_{p-1}$,
(ii) $\lambda(f_{p})-\lambda(g_{p})<\frac{1}{p}$, and (iii) $t-\frac{1}{p}<u_{p}<v_{p}<t+\frac{1}{p}$. 

Consider $p,q\geq1$. We have $f_{q}\diamondsuit u_{q}<t<v_{p}\diamondsuit g_{p}$.
Hence $\lambda(f_{q})\geq\lambda(g_{p})>\lambda(f_{p})-\frac{1}{p}$.
By symmetry we also have $\lambda(f_{p})>\lambda(f_{q})-\frac{1}{q}$.
Combining, we see that $|\lambda(f_{q})-\lambda(f_{p})|<\frac{1}{p}+\frac{1}{q}$.
Hence $(\lambda(f_{p}))_{p=1,2,\cdots}$ is a Cauchy sequence and
converges. Similarly $(\lambda(g_{p}))_{p=1,2,\cdots}$ converges.
In view of condition (ii), the two limits are equal. 

By the definition of $\lambda$, we see that $\lim_{p\rightarrow\infty}If_{p}(X)$
exists. Since $f_{p-1}\geq f_{p}$ for each $p>1$, the Monotone Convergence
Theorem \ref{Thm. Monotone Convergence} implies that $Y\equiv\lim_{p\rightarrow\infty}f_{p}(X)$
is an integrable function, with $\lim_{p\rightarrow\infty}I|f_{p}(X)-Y|=0$.
Likewise $Z\equiv\lim_{p\rightarrow\infty}g_{p}(X)$ is an integrable
function, with $\lim_{p\rightarrow\infty}I|g_{p}(X)-Z|=0$. Furthermore
\[
I|Y-Z|=\lim_{p\rightarrow\infty}I|f_{p}(X)-g_{p}(X)|=\lim_{p\rightarrow\infty}(\lambda(f_{p})-\lambda(g_{p}))=0
\]
According to Proposition \ref{Prop. Property of Full Sets}, we have
$Y=Z$ a.e.

We next show that $Y$ is an indicator with $(Y=1)=(t\leq X)$. Consider
$\omega\in domain(Y)$. Suppose $Y(\omega)>0$. Then $\omega\in domain(X)$
and $f_{p}(X(\omega))\geq Y(\omega)>0$ for each $p\geq1$. It follows,
in view of condition (i) above, that $u_{p-1}\leq X(\omega)$ and
so $f_{p-1}(X(\omega))=1$ for each $p>1$. Passing to limit as $p\rightarrow\infty$,
we conclude that $t\leq X(\omega)$ and so $Y(\omega)=1$. In particular
$Y$ can have only two possible values, 0 or 1. Thus $Y$ is an indicator.
We have also seen that $(Y=1)\subset(t\leq X)$. Conversely, suppose
$t\leq X(\omega)$. Then, in view of $f_{p}\diamondsuit u_{p}<t$
in condition (i) above, we have $f_{p}(X(\omega))=1$ for each $p\geq1$.
It follows trivially that $\lim_{p\rightarrow\infty}f_{p}(X(\omega))=1$
and so $\omega\in domain(Y)$ and $Y(\omega)=1$. Summing up, the
set $(X\geq t)$ has $Y$ as an indicator. 

We will now prove that $(Y=0)=(X<t)$. Let $\omega\in(Y=0)$. Then
by the definition of $Y$ we have $\omega\in domain(X)$ and $0\equiv\lim_{p\rightarrow\infty}f_{p}(X(\omega))$.
Then there exists $p\geq1$ such that $1>f_{p}(X(\omega))$. This
implies, in view of $f_{p}\diamondsuit u_{p}$ in condition (i) above,
that $X(\omega)\leq u_{p}<t$. In other words $\omega\in(X<t)$. Conversely,
suppose $\omega\in(X<t)$. Then $X(\omega)<t$. Since $u_{p}\uparrow t$
as $p\uparrow\infty$ there exists $q$ so large that $X(\omega)<u_{p-1}$
for each $p>q$. In view of $u_{p-1}\diamondsuit f_{p}$ in condition
(i) above, we have $f_{p}(X(\omega))=0$ for each $p\geq q$. It follows
that $\lim_{p\rightarrow\infty}f_{p}(X(\omega))=0$ and so $\omega\in domain(Y)$
and $Y(\omega)=0$. Summing up, we have $(Y=0)=(X<t)$. Thus $(t\leq X)^{c}=(X<t)$.

Similarly we can prove that $(t<X)$ has $Z$ as an indicator, and
that $(t<X)^{c}=(X\leq t)$. It follows that $\mu(t\leq X)=IY=IZ=\mu(t<X)$.

It remains to show that $\mu(t\leq X)$ is continuous at $t$. Let
$p>1$ be arbitrary. Recall that $u_{p-1}\diamondsuit f_{p}\diamondsuit u_{p}<t<v_{p}\diamondsuit g_{p}\diamondsuit v_{p-1}$
where $u_{p}$ and $v_{p}$ are arbitrarily close to $t$ if $p$
is sufficiently large. From the previous paragraphs we see that $|\lambda(f_{p})-\mu(t\leq X)|=\lambda(f_{p})-\lim_{q\rightarrow\infty}\lambda(f_{q})\leq\frac{1}{p}$.
Now consider any $t'\in K\cap J_{c}$ such that $t'\in(u_{p},v_{p})$.
We can similarly construct an arbitrarily large $q$, points $t'_{q-1},t'_{q},s'_{q}$
and $s'_{q-1}$ in $K$ that are arbitrarily close to $t'\in(t_{p},s_{p})$,
and functions $f'_{q},g'_{q}\in G$ such that 
\[
t_{p}<t'_{q-1}\diamondsuit f'_{q}\diamondsuit t'_{q}<t<s'_{q}\diamondsuit g'_{q}\diamondsuit s'_{q-1}<s_{p}
\]
and such that $|\lambda(f'_{q})-\mu(t'\leq X)|\leq\frac{1}{q}$. It
follows that $f_{p}\diamondsuit t_{p}<t'_{q-1}\diamondsuit f'_{q}$
and so $\lambda(f_{p})\geq\lambda(f'_{q})$. Similarly $\lambda(g'_{q})\geq\lambda(g_{p})$.
Hence
\[
0\leq\lambda(f_{p})-\lambda(f'_{q})<\lambda(g_{p})+\frac{1}{p}-\lambda(g'_{q})\leq\frac{1}{p}
\]
Using the triangle inequality twice, we obtain $|\mu(t'\leq X)-\mu(t\leq X)|<\frac{1}{q}+\frac{2}{p}$.
Since $q$ is arbitrarily large, we see that $|\mu(t'\leq X)-\mu(t\leq X)|\leq\frac{2}{p}$
for each $t'$ in the neighborhood $(t_{p},s_{p})$ of $t$. Continuity
of $\mu(t\leq X)$ at $t$ has thus been established. 
\end{proof}
\begin{cor}
\label{Cor. Countable Exceptional Points for Integrable X} \textbf{\emph{(Abundance
of integrable sets). }}Let $X$ be an integrable function. There exists
a countable subset $J$ of $R$  such that for each $t$ in the metric
complement $J_{c}$ of $J$ the following conditions hold.
\end{cor}
\begin{enumerate}
\item \emph{If $t>0$ then $(t<X)$ and $(t\leq X)$ are integrable, with
equal }meas\emph{ures that are continuous at $t$.}
\item \emph{If $t<0$ then $(X<t)$ and $(X\leq t)$ are integrable, with
equal }meas\emph{ures that are continuous at $t$.}
\end{enumerate}
\begin{proof}
Apply Theorem \ref{Thm. (t<X) is integrable} to $X$ and $-X$ and
let $J$ be the union of the two corresponding countable exceptional
sets. 
\end{proof}
\begin{defn}
\label{Def. Regular=000026 Continuity Pts of Integrable Func} \textbf{(Regular
and continuity points of an integrable function relative to an integrable
set).} Let $X$ be an integrable function, let $A$ be an integrable
set, and let $t\in R$. We say that $t$ is a \index{regular point of integrable function}regular
point of $X$ relative to $A$ if (i) there exists a sequence $(s_{n})_{n=1,2,\cdots}$
of real numbers decreasing to $t$ such that $(s_{n}<X)A$ is integrable
for each $n\geq1$ and such that $\lim_{n\rightarrow\infty}\mu(s_{n}<X)A$
exists, and (ii) there exists a sequence $(r_{n})_{n=1,2,\cdots}$
of real numbers increasing to $t$ such that $(r_{n}<X)A$ is integrable
for each $n\geq1$ and such that $\lim_{n\rightarrow\infty}\mu(r_{n}<X)A$
exists. If in addition the two limits in (i) and (ii) are equal, then
we call $t$ a \index{continuity point of integrable function}continuity
point of $X$ relative to $A$. We say that a positive real number
$t>0$ is a regular point of $X$ if conditions (i) and (ii), with
$A$ omitted, are satisfied. We say that a negative real number $t<0$
is a regular point of $X$ if $-t$ is a regular point of $-X$. $\square$
\end{defn}
\begin{cor}
\label{Cor. Properties of Regular Pts of Integrable Func} \textbf{\emph{(Simple
properties of regular and continuity points). }}Let $X$ be an integrable
function, let $A$ be an integrable set, and let $t$ be a regular
point of $X$ relative to $A$. Then the following holds.
\end{cor}
\begin{enumerate}
\item \emph{If $u$ is a regular point of $X$, then $u$ is a regular point
for $X$ relative to any integrable set $B$. If $u$ is a continuity
point of $X$, then $u$ is a continuity point for $X$ relative to
any integrable set $B$.}
\item \emph{All but countably many real numbers are continuity points of
$X$.}
\item \emph{All but countably many real numbers are continuity points of
$X$ relative to $A$. Hence all but countably many real numbers are
regular points of $X$ relative to $A$. }
\item \emph{The sets $A(t<X),A(t\leq X),A(X<t)$ , $A(X\leq t)$, and $A(X=t)$
are integrable sets.}
\item \emph{$(X\leq t)A=A((t<X)A)^{c}$ a.e. and $(t<X)A=A((X\leq t)A)^{c}$
a.e.}
\item \emph{$(X<t)A=A((t\leq X)A)^{c}$ a.e. and $(t\leq X)A=A((X<t)A)^{c}$
a.e.}
\item \emph{For a.e. $\omega\in A$, we have $t<X(\omega)$, $t=X(\omega)$,
or $t>X(\omega)$. Thus we have a limited, but very useful, version
of the principle of excluded middle.}
\item \emph{Let $\varepsilon>0$ be arbitrary. There exists $\delta>0$
such that if $r\in(t-\delta,t]$ and $A(X<r)$ is integrable, then
$\mu(A(X<t))-\mu(A(X<r))<\varepsilon$. There exists $\delta>0$ such
that if $s\in[t,t+\delta)$ and $A(X\leq s)$ is integrable, then
$\mu(A(X\leq s))-\mu(A(X\leq t))<\varepsilon$.}
\item \emph{If $t$ is a continuity point of $X$ relative to $A$, then
$\mu((t<X)A)=$ $\mu((t\leq X)A)$.}
\end{enumerate}
\begin{proof}
1. Suppose $u>0$ and $u$ is a regular point of $X$. Then by Definition
\ref{Def. Regular=000026 Continuity Pts of Integrable Func}, (i$'$)
there exists a sequence $(s_{n})_{n=1,2,\cdots}$ of real numbers
decreasing to $u$ such that $(s_{n}<X)$ is integrable and $\lim_{n\rightarrow\infty}\mu(s_{n}<X)$
exists, and (ii$'$) there exists a sequence $(r_{n})_{n=1,2,\cdots}$
of positive real numbers increasing to $u$ such that $(r_{n}<X)$
is integrable and $\lim_{n\rightarrow\infty}\mu(r_{n}<X)$ exists.
Now let $B$ be any integrable set. Then for all $m>n\geq1$, we have
\[
0\leq\mu((s_{m}<X)B)-\mu((s_{n}<X)B)=\mu(B(s_{m}<X)(s_{n}<X)^{c})
\]
 
\[
\leq\mu((s_{m}<X)(s_{n}<X)^{c})=\mu(s_{m}<X)-\mu(s_{n}<X)\downarrow0
\]
Therefore $(\mu((s_{n}<X)B))_{n=1,2,\cdots}$ is a Cauchy sequence
and converges, verifying condition (i) in Definition \ref{Def. Regular=000026 Continuity Pts of Integrable Func}.
Condition (ii) is similarly verified. Hence $u$ is a regular point
of $X$ relative to $B$. Suppose in addition that $u$ is a continuity
point of $X$. Then for each $n\geq1$, we have
\[
0\leq\mu((r_{n}<X)B)-\mu((s_{n}<X)B)=\mu(B(r_{n}<X)(s_{n}<X)^{c})
\]
\[
\leq\mu((r_{n}<X)(s_{n}<X)^{c})=\mu(r_{n}<X)-\mu(s_{n}<X)\downarrow0
\]
Therefore the two sequences $(\mu((r_{n}<X)B))_{n=1,2,\cdots}$ and
$(\mu((s_{n}<X)B))_{n=1,2,\cdots}$ have the same limit. Thus assertion
1 is proved for the case $t>0.$ The case $t<0$ is similar.

2. Assertion 2 is an immediate consequence of Corollary \ref{Cor. Countable Exceptional Points for Integrable X}. 

3. Assertion 3 is follows from assertions 1 and 2. 

4. Let $(s_{n})_{n=1,2,\cdots}$ and $(r_{n})_{n=1,2,\cdots}$ be
sequences of real numbers satisfying conditions (i) and (ii) in Definition
\ref{Def. Regular=000026 Continuity Pts of Integrable Func}. Then
$(s_{n}<X)A$ is an integrable set for $n\geq1$, and $\lim_{n\rightarrow\infty}\mu(s_{n}<X)A$
exists. Since $(s_{n})$ decreases to $t$, we have $(t<X)A=\bigcup_{n=1}^{\infty}(s_{n}<X)A$
a.e. This union is integrable, and $\mu(s_{n}<X)A\uparrow\mu((t<X)A)$,
according to Proposition \ref{Prop. Sequences of intgrable sets}.
Define $B\equiv(X\leq t)A$ and $C\equiv(t<X)A$. Because the set
$C$ has just been proved integrable, the set $D\equiv C\cup C^{c}$
is a full set. Consider $\omega\in AD$. Then either $\omega\in C$
or $\omega\in C^{c}$. Consider $\omega\in B$. Then $X(\omega)\leq t$.
If $\omega\in C$ then we would have $t<X(\omega)$, a contradiction.
Hence $\omega\in C^{c}$. Conversely, consider $\omega\in C^{c}$.
If $t<X(\omega)$ we would have $\omega\in CC^{c}=\phi$, a contradiction.
Hence $X(\omega)\leq t$ and so $\omega\in B$. Summing up, we have
$ADB=ADC^{c}$. In other words $(X\leq t)A\equiv B=AB=AC^{c}=A((t<X)A)^{c}$
a.e. It follows from Proposition \ref{Prop. Measure of Boolean Integrable Sets}
that $(X\leq t)A\equiv B$ is integrable. Similarly, since $(r_{n})$
increases to $t$, we have $(X<t)A=\bigcup_{n=1}^{\infty}A((r_{n}<X)A)^{c}$
a.e. This union is integrable, with $\mu(A((r_{n}<X)A)^{c})=\mu(A)-\mu((r_{n}<X)A)\uparrow\mu((X<t)A)$.
Since we can show, in a proof similar to the one above, that $(t\leq X)A=A((X<t)A)^{c}$
a.e., it follows from Proposition \ref{Prop. Measure of Boolean Integrable Sets}
that $(X\leq t)A$ is integrable. Since $A(X=t)=A(t\leq X)((t<X)A)^{c}$
a.e. the set $A(X=t)$ is integrable. Assertion 4 is proved. 

5. We have seen in the proof of assertion 4 that $B=AC^{c}$ a.e.
where $B\equiv(X\leq t)A$ and $C\equiv(t<X)A$. Using Proposition
\ref{Prop. Measure of Boolean Integrable Sets}, we have $AB^{c}=A(AC^{c})^{c}=AC$
a.e. 

6. Similar.

7. With the notations in the above proof for assertion 4, we have
$B=ADC^{c}$ where $D\equiv C\cup C^{c}$ is a full set. Consider
any $\omega\in AD$. Then, we have either $t<X(\omega)$ or $\omega\in C^{c}$.
Hence, we have either $t<X(\omega)$ or $\omega\in B$. Therefore
$t<X(\omega)$ or $X(\omega)\leq t$. In other words, for a.e. $\omega\in A$,
we have $t<X(\omega)$ or $X(\omega)\leq t$. Similarly, for a.e.
$\omega\in A$, we have $X(\omega)<t$ or $t\leq X(\omega)$. Combining,
for a.e. $\omega\in A$, we have $t<X(\omega)$, $X(\omega)<t$, or
$t\leq X(\omega)\leq t$. Assertion 7 is proved.

8. Use the notations in the above proof for assertion 4. Let $\varepsilon>0$
be arbitrary. Let $n\geq1$ be so large that $\mu((t<X)A)-\mu(s_{n}<X)A<\varepsilon$.
Define $\delta=s_{n}-t$. Suppose $s\in[t,t+\delta)$ and $A(X\leq s)$
is integrable. Then $s<s_{n}$ and so 
\[
\mu(A(X\leq s))-\mu(A(X\leq t))=(\mu(A)-\mu(s<X)A)-(\mu(A)-\mu((t<X)A))
\]
\[
=\mu((t<X)A)-\mu(s<X)A\leq\mu((t<X)A)-\mu(s_{n}<X)A<\varepsilon
\]
This proves the second half of assertion 8, the first half having
a similar proof.

9. Suppose $t$ is a continuity point of $X$ relative to $A$. Then
the limits $\lim_{n\rightarrow\infty}\mu(s_{n}<X)A$ and $\lim_{n\rightarrow\infty}\mu(r_{n}<X)A$
are equal. The proof of assertion 4 therefore shows $\mu((t<X)A)=\mu(A)-\mu((X<t)A)$
which is in turn equal to $\mu((t\leq X)A)$ in view of assertion
6. 
\end{proof}
\begin{defn}
\label{Convention. Only regular Pts of Integrable X used} \textbf{(Convention
of implicit assumption of regular points of integrable functions)}
Let $X$ be an integrable function, and let $A$ be an integrable
set. Henceforth, if the integrability of the set $(X<t)A$ or $(X\leq t)A$,
for some $t\in R$, is required in a discussion, then it is understood
that the real number $t$ has been chosen from the regular points
of the integrable function $X$ relative to the integrable set $A$. 

Likewise, if the integrability of the set $(t<X)$ or $(t\leq X)$,
for some $t>0$, is required in a discussion, then it is understood
that the number $t>0$  has been chosen from the regular points of
the integrable function $X$. 

Separately, we will sometimes  write $(X<t;Y\leq s;\cdots)$ for $(X<t)(Y\leq s)\cdots$
for brevity.$\square$
\end{defn}
Recall that $C_{ub}(R)$ is the space of bounded and uniformly continuous
functions on $R$. 
\begin{prop}
\label{Prop. X, A in L and f in C_u(R) =00003D>  f(X)1_A in L} \textbf{\emph{(Product
of bounded continuous function of an integrable function and an integrable
indicator is integrable).}} Suppose $X\in L$, $A$ is an integrable
set, and $f\in C_{ub}(R)$. Then $f(X)1_{A}\in L$. In particular,
if $X\in L$ is bounded, then $X1_{A}$ is integrable.
\end{prop}
\begin{proof}
Let $c>0$ be so large that $|f|\leq c$ on $R$. Let $\varepsilon>0$
be arbitrary. Since $X$ is integrable, there exists $a>0$ so large
that $I|X|-I|X|\wedge(a-1)<\varepsilon$. Since $f$ is uniformly
continuous, there exists a sequence $-a=t_{0}<t_{1}<\cdots<t_{n}=a$
whose mesh is so small that $|f(t_{i})-f(x)|\leq\varepsilon$ for
each $x\in(t_{i-1},t_{i}]$, for each $i=1,\cdots,n$. Then
\[
Y\equiv\sum_{i=1}^{n}f(t_{i})1_{(t(i-1)<X\leq t(i))A}
\]
is an integrable function. Moreover, since $1_{(|X|>a)A}\leq|X|-|X|\wedge(a-1)$,
we have
\[
|f(X)1_{A}-Y|\leq|\sum_{i=1}^{n}(f(X)-f(t_{i}))1_{(t(i-1)<X\leq t(i))A}|+c1_{(|X|>a)A}
\]
\[
\leq\varepsilon1_{A}+c(|X|-|X|\wedge(a-1))\quad a.e.,
\]
where
\[
I(\varepsilon1_{A}+c(|X|-|X|\wedge(a-1)))
\]
\[
=\varepsilon\mu(A)+c(I|X|-I|X|\wedge(a-1))<\varepsilon\mu(A)+c\varepsilon\rightarrow0
\]
as $\varepsilon\rightarrow0$. Hence, by Theorem \ref{Thm. |X-Yn|<Zn where Yn,Zn integrable =000026 IZn->0 implies X integrable},
$f(X)1_{A}\in L$. 

Now suppose $X\in L$ is bounded. Let $b>0$ be such that $|X|\leq b$.
Define $f\in C_{ub}(R)$ by $f(r)\equiv b\wedge r\vee-b$. Then $X=f(X)$
and so, according to the first part of this proposition, $X1_{A}\in L$.
\end{proof}

\section{Uniform Integrability}

In this section, let $(\Omega,L,I)$ be a complete integration space.
We will give some useful propositions on bounds of integrals and measures.
\begin{prop}
\textbf{\emph{(Chebychev's Inequality).}} \label{Prop. Chebychev's inequalty}Let
$X\in L$ be arbitrary. Then he following holds. 

1. (First and common version). If $t>0$ is a regular point of $X$,
then we have $\mu(|X|>t)\leq t^{-1}I|X|$. 

2. (Second version). If $I|X|<b$ for some $b>0$, then for each $s>0$,
we have $(|X|>s)\subset B$ for some integrable set $B$ with $\mu(B)<s^{-1}b$.
This second version of \emph{Chebychev's inequality\index{Chebychev's inequality}}
is useful when a real number $s>0$ is given without any assurance
that the set $(|X|>s)$ is integrable.
\end{prop}
\begin{proof}
1. $1_{(|X|>t)}\leq t^{-1}|X|$. 

2. Take an arbitrary regular point $t$ of the integrable function
$X$ in the open interval $(b^{-1}I|X|s,s)$. Let $B\equiv(|X|>t)$.
By Assertion 1, we then have $\mu(B)\leq t^{-1}I|X|<s^{-1}b$. Moreover,
$(|X|>s)\subset(|X|>t)\equiv B$.
\end{proof}
\begin{prop}
\label{Prop. Existence of modulus of integrability} \textbf{\emph{(Bounds
related to integrable functions).}} Let$X\in L$ be arbitrary. Let
$A$ be an arbitrary integrable set. Then the following holds.
\end{prop}
\begin{enumerate}
\item \emph{$X1_{A}\in L$.}
\item \emph{$I(|X|1_{A})\rightarrow0$ as $\mu(A)\rightarrow0$. Specifically,
for each $\varepsilon>0$ there exists $\delta=\delta(\varepsilon)>0$
such that $I(|X|1_{A})\leq\varepsilon$ if $\mu(A)<\delta(\varepsilon)$.}
\item \emph{$I(|X|1_{(|X|>a)})\rightarrow0$ as $a\rightarrow\infty$. Specifically,
suppose $I|X|\leq b$ for some $b>0$, and suppose the operation $\delta$
is as in assertion 2. For each $\varepsilon>0$, if we define $\eta(\varepsilon)\equiv b/\delta(\varepsilon)$,
then $I(|X|1_{(|X|>a)})\leq\varepsilon$ for each $a>\eta(\varepsilon)$.}
\item \emph{Suppose an operation $\eta>0$ is such that $I(|X|1_{(|X|>a)})\leq\varepsilon$
for each $a>\eta(\varepsilon)$. Then the operation $\delta$ defined
by $\delta(\varepsilon)\equiv\frac{\varepsilon}{2}/\eta(\frac{\varepsilon}{2})$
satisfies the conditions in assertion 2. }
\end{enumerate}
\begin{proof}
1. Let $n>0$ be arbitrary. Then $|X1_{A}|\wedge n=|X|\wedge(n1_{A})$
is integrable. Moreover for $n>p$ we have $I(|X1_{A}|\wedge n-|X1_{A}|\wedge p)\leq I(|X|\wedge n-|X|\wedge p)\rightarrow0$
as $p\rightarrow\infty$ since $|X|\in L$. Hence $\lim_{n\rightarrow\infty}I(|X1_{A}|\wedge n)$
exists. By the Monotone Convergence Theorem, the limit $|X1_{A}|=\lim_{n\rightarrow\infty}|X1_{A}|\wedge n$
is integrable. Similarly, $|X_{+}1_{A}|$ is integrable and so also
is $X1_{A}=2|X_{+}1_{A}|-|X1_{A}|$.

2. Suppose $a>0$. Since $|X|1_{A}\leq(|X|-|X|\wedge a)1_{A}+a1_{A}$
we have $I(|X|1_{A})\leq I|X|-I(|X|\wedge a)+a\mu(A)$. Given any
$\varepsilon>0$, since $X$ is integrable, there exists $a>0$ so
large that $I|X|-I(|X|\wedge a)<\varepsilon$. Then for each $A$
with $\mu(A)<\varepsilon/a$ we have $I(|X|1_{A})<2\varepsilon$.

3. Suppose $a>\eta(\varepsilon)\equiv b/\delta(\varepsilon)$ where
$\delta$ is an operation as in assertion 2. Chebychev's inequality
gives $\mu(|X|>a)\leq I|X|/a\leq b/a<\delta(\varepsilon)$. Hence
$I(1_{(|X|>a)}|X|)<\varepsilon$. 

4. Suppose $\mu(A)<\varepsilon/\eta(\varepsilon)$. For each $a>\eta(\varepsilon)$
we have $I(|X|1_{A})\leq I(a1_{A(X\leq a)}+1_{(|X|>a)}|X|)\leq a\mu(A)+\varepsilon\leq a\varepsilon/\eta(\varepsilon)+\varepsilon$.
By taking $a$ arbitrarily close to $\eta(\varepsilon)$ we see that
$I(|X|1_{A})\leq\eta(\varepsilon)\varepsilon/\eta(\varepsilon)+\varepsilon=2\varepsilon$.
Replace $\varepsilon$ by $\frac{\varepsilon}{2}$ and the assertion
is proved. 
\end{proof}
Note that in the proof for assertion 4 of Proposition \ref{Prop. Existence of modulus of integrability},
we use a real number $a>\eta(\varepsilon)$ arbitrarily close to $\eta(\varepsilon)$
rather than simply $a=\eta(\varepsilon)$. This ensures that $a$
can be a regular point of $|X|$, as required in Convention \ref{Convention. Only regular Pts of Integrable X used}. 
\begin{defn}
\label{Def. Uniform Integrability}\label{Def. Modulus of integrability}
\textbf{(Uniform integrability and simple modulus of integrability).
}A family $G$ of integrable functions is said to be \emph{uniformly
integrable} \index{uniform integrability} if for each $\varepsilon>0$,
there exists $\eta(\varepsilon)$ such that $E(|X|1_{(|X|>a)})\leq\varepsilon$
for each $a>\eta(\varepsilon)$, for each $X\in G$. The operation
$\eta$ is then called a \emph{simple modulus of integrability}\index{simple modulus of integrability}
of $G$. 

Proposition \ref{Prop. Existence of modulus of integrability} ensures
that each family $G$ consisting finitely many integrable functions
is uniformly integrable. $\square$
\end{defn}
\begin{prop}
\emph{\label{Prop. Alternative def of uniform integrabitty, and simple modulus of integrability}
}\textbf{\emph{(Alternative definition of uniform integrability, and
modulus of integrability).}}\emph{ }Suppose the integration space
$(\Omega,L,I)$ is such that $1\in L$ and $I1=1$. Then a family
$G$ of integrable r.r.v.'s is uniformly integrable \index{uniform integrable r.r.v.'s}
iff \emph{(i)} there exists $b\geq0$ such that $I|X|\leq b$ for
each $X\in G$, and \emph{(ii)} for each $\varepsilon>0$, there exists
$\delta(\varepsilon)$ such that $I|X|1_{A}\leq\varepsilon$ for each
integrable set $A$ with $\mu(A)<\delta(\varepsilon)$, and for each
$X\in G$. The operation $\delta$ is then called a \emph{modulus
of integrability}\index{modulus of integrability} of $G$. 
\end{prop}
\begin{proof}
First suppose the family $G$ is uniformly integrable. In other words,
for each $\varepsilon>0$, there exists $\eta(\varepsilon)$ such
that 
\[
I(|X|1_{(|X|>a)})\leq\varepsilon
\]
for each $a>\eta(\varepsilon)$, and for each $X\in G$. Define $b\equiv\eta(1)+2$.
Let $X\in G$ be arbitrary. Take any $a\in(\eta(1),\eta(1)+1)$. Then
\[
I|X|=I(1_{(|X|>a)}|X|)+I(1_{(|X|\leq a)}|X|)\leq1+aI1=1+a<1+\eta(1)+1=b,
\]
where the second equality follows from the hypothesis that $I1=1$.
This verifies Condition (i) \ref{Def. Uniform Integrability}. Now
let $\varepsilon>0$ be arbitrary. Define the operation $\delta$
by\emph{ $\delta(\varepsilon)\equiv\frac{\varepsilon}{2}/\eta(\frac{\varepsilon}{2})$.
}Then Assertion 4 of Proposition \ref{Prop. Existence of modulus of integrability}
implies that\emph{ $I(|X|1_{A})\leq\varepsilon$} for each integrable
set $A$ with $\mu(A)<\delta(\varepsilon)$, for each $X\in G$. This
verifies Condition (ii).

Conversely, suppose the Conditions (i) and (ii) hold. For each $\varepsilon>0$,
define $\eta(\varepsilon)\equiv b/\delta(\varepsilon)$.\emph{ }Then,
according to Assertion 3 of Proposition \ref{Prop. Existence of modulus of integrability}\emph{,}
we have $I(|X|1_{(|X|>a)})\leq\varepsilon$ for each $a>\eta(\varepsilon)$.
Thus $G$ is uniformly integrable according to Definition \ref{Def. Uniform Integrability}. 
\end{proof}
\begin{prop}
\label{Prop |X|<integrable Y for all X in G =00003D> G unif integrable}\textbf{\emph{
}}\textbf{\textup{(Dominated}}\textbf{\emph{ uniform integrability).}}
If there is an integrable functions $Y$ such that $|X|\leq Y$ for
each $X$ in a family $G$ of integrable functions, then $G$ is uniformly
integrable. 
\end{prop}
\begin{proof}
Note that $b\equiv I|Y|$ satisfies Conditions (i) in Definition \ref{Def. Uniform Integrability}.
Let $\varepsilon>0$ be arbitrary. Then Assertion 3 of Proposition
\ref{Prop. Existence of modulus of integrability} guarantees an operation
$\eta$ such that, for each $\varepsilon>0$, we have $I(1_{(|Y|>a)}|Y|)\leq\varepsilon$
for each $a>\eta(\varepsilon)$. Hence, for each $X\in G$, and for
each $a>\eta(\varepsilon)$, we have 
\[
I(1_{(|X|>a)}|X|)\leq I(1_{(Y>a)}Y)\leq\varepsilon.
\]
Thus $\eta$ is a common simple modulus of integrability for members
$X$ of $G$. The conditions in Definition \ref{Def. Uniform Integrability}
have been verified for the family $G$ to be uniformly integrable. 
\end{proof}
\begin{prop}
\label{Prop. Integrable func is limit in L1 of combo of indicators}
\textbf{\emph{(Each integrable function is the $L_{1}$ limit of some
sequence of linear combinations of integrable indicators).}}
\end{prop}
\begin{enumerate}
\item \emph{Suppose $X$ is an integrable function with $X$$\geq0$ a.e.
Then there exists a sequence $(Y_{k})_{k=1,2,\cdots}$ such that for
each $k\geq1$ we have (i) $Y_{k}\equiv\sum_{i=1}^{n_{k}-1}t_{k,i}1_{(t_{k,i}<X\leq t_{k,i+1})}\in L$
for some sequence $0<t_{k,1}<\cdots<t_{k,n_{k}}$, (ii) $I|Y_{k}-X|\rightarrow0$,
and (iii) $Y_{k}\uparrow X$ on $D\equiv\bigcap_{k=1}^{\infty}domain(Y_{k})$.
Moreover, we can take }$n_{k}\equiv2^{2k}$\emph{ and $t_{k,i}\equiv2^{-k}ia$
for $i=1,\cdots,n_{k}$, where $a$ is some positive real number.}
\item \emph{Suppose $X$ is an integrable function. Then there exists a
sequence $(Z_{k})_{k=1,2,\cdots}$ of linear combinations of mutually
exclusive integrable indicators such that $I|X-Z_{k}|\leq2^{-k}$
and such that $Z_{k}\rightarrow X$ on $\bigcap_{k=1}^{\infty}domain(Z_{k})$.
Furthermore, there exists a sequence $(U_{k})_{k=1,2,\cdots}$ of
linear combinations of integrable indicators which is a representation
of $X$ in $L$.}
\item \emph{Suppose $X$ and $X'$ are bounded integrable functions. Then
$XX'$ is integrable.}
\end{enumerate}
\begin{proof}
1. Let $a>0$ be such that $(a<i^{-1}2^{k}X)$ is integrable for all
$k,i\geq1$. For $k\geq1$ define $Y_{k}\equiv\sum_{i=1}^{n_{k}-1}t_{k,i}1_{(t_{k,i}<X\leq t_{k,i+1})}\in L$
where $n_{k}\equiv2^{2k}$ and $t_{k,i}\equiv2^{-k}ia$ for $i=1,\cdots,n_{k}$.
For all $k,i\geq1$, the set $(t_{k,i}<X)\equiv(2^{-k}ia<X)=(a<i^{-1}2^{k}X)$
is integrable. Hence $Y_{k}\in L$ for $k\geq1$. By definition, $t_{k,n_{k}}\equiv2^{k}a\rightarrow\infty$
and $t_{k,i}-t_{k,i-1}=2^{-k}a\rightarrow0$ for $i=1,\cdots,n_{k}$.
Let $h>k\geq1$ be arbitrary. Consider any $\omega\in D$. Suppose
$Y_{k}(\omega)>0$. Then $Y_{k}(\omega)=t_{k,i}1_{(t_{k,i}<X(\omega)\leq t_{k,i+1})}=t_{k,i}$
for some $i=1,\cdots,n_{k}-1$. Write $p\equiv2^{h-k}i$ and $q\equiv2^{h-k}(i+1)\leq2^{h-k}n_{k}\equiv2^{h+k}\leq n_{h}$.
Then
\[
t_{h,p}\equiv2^{-h}pa\equiv2^{-h}(2^{h-k}i)a=t_{k,i}<X(\omega)
\]
\[
\leq t_{k,i+1}=2^{-h}(2^{h-k}(i+1))a\equiv2^{-h}qa\equiv t_{h,q}
\]
Therefore there exists $j$ with $p\leq j<q$ such that $t_{h,j}<X(\omega)\leq t_{h,j+1}$.
Consequently 
\begin{equation}
Y_{k}(\omega)=t_{k,i}=t_{h,p}\leq t_{h,j}=Y_{h}(\omega)<X(\omega)\leq t_{h,j+1}\leq t_{h,q}=t_{k,i+1}\label{eq:temp-62}
\end{equation}
Thus we see that $0\leq Y_{k}\leq Y_{h}\leq X$ on $D$ for $h>k\geq1$.
Next, let $\varepsilon>0$ be arbitrary. Then either $X(\omega)>0$
or $X(\omega)<\varepsilon$. In the first case, $Y_{m}(\omega)>0$
for some $m\geq1$, whence $Y_{k}(\omega)>0$ for each $k\geq m$,
and so, for each $k\geq k_{0}\equiv m\vee\log_{2}(a\varepsilon^{-1})$,
we see from inequality \ref{eq:temp-62} that 
\[
X(\omega)-Y_{k}(\omega)<t_{k,i+1}-t_{k,i}=2^{-k}a<\varepsilon
\]
In the second case, we have, trivially, $X(\omega)-Y_{k}(\omega)<\varepsilon$
for each $k\geq1$. Combining, we have $Y_{k}\uparrow X$ on $D$.
We will show next that $IY_{k}\uparrow IX$. By Proposition \ref{Prop. Existence of modulus of integrability},
 there exists $k_{1}\geq1$ so large that $I(1_{(2^{k}a<X)}X)\leq\varepsilon$
for each $k\geq k_{1}$. At the same time, since $X\geq0$ is integrable,
there exists $k_{2}\geq1$ so large that $I(2^{-k}a\wedge X)<\varepsilon$
for each $k\geq k_{2}$. Hence, for each $k\geq k_{0}\vee k_{1}\vee k_{2}$,
we have 
\[
I(X-Y_{k})=I(X-Y_{k})1_{(X\leq2^{-k}a)}+I(X-Y_{k})1_{(2^{-k}a<X\leq2^{k}a)}+I(X-Y_{k})1_{(2^{k}a<X)}
\]
\[
\leq I(2^{-k}a\wedge X)1_{(X\leq2^{-k}a)}+I(2^{-k}a\wedge X)1_{(2^{-k}a<X\leq2^{k}a)}+IX1_{(2^{k}a<X)}
\]
\[
\leq I(2^{-k}a\wedge X)1_{(X\leq2^{k}a)}+\varepsilon
\]
\[
\leq I(2^{-k}a\wedge X)+\varepsilon<2\varepsilon
\]
Since $\varepsilon>0$ is arbitrary, we conclude that $I|Y_{k}-X|\rightarrow0$.
Assertion 1 is proved.

2. By assertion 1, we see that there exists a sequence \emph{$(Y_{k}^{+})_{k=1,2,\cdots}$
}of linear combinations of mutually exclusive indicators such that
$I|X_{+}-Y_{k}^{+}|<2^{-k-1}$ for each $k\geq1$ and such that \emph{$Y_{k}^{+}\uparrow X_{+}$
on $D^{+}\equiv\bigcap_{k=1}^{\infty}domain(Y_{k}^{+})$.} By the
same token, there exists a sequence \emph{$(Y_{k}^{-})_{k=1,2,\cdots}$
}of linear combinations of mutually exclusive indicators such that\emph{
}$I|X_{-}-Y_{k}^{-}|<2^{-k-1}$ for each $k\geq1$ and such that \emph{$Y_{k}^{-}\uparrow X_{-}$
on $D^{-}\equiv\bigcap_{k=1}^{\infty}domain(Y_{k}^{-})$. }For each
$k\geq1$ define $Z_{k}\equiv Y_{k}^{+}-Y_{k}^{-}$ whence $I|X-Z_{k}|\leq I|X_{+}-Y_{k}^{+}|+I|X_{-}-Y_{k}^{-}|<2^{-k}$.
Moreover, we see from the proof of assertion 1 that, for each $k\geq1$,
$Y_{k}^{+}$ can be taken to be a linear combination of indicators
of subsets of $(X_{+}>0)$, and, by the same token, $Y_{k}^{-}$ can
be taken to be a linear combination of indicators of subsets of $(X_{-}>0)$.
Since $(X_{+}>0)$ and $(X_{-}>0)$ are disjoint, so $Z_{k}\equiv Y_{k}^{+}-Y_{k}^{-}$
is a linear combination of mutually exclusive indicators. Since \emph{$Y_{k}^{+}\uparrow X_{+}$}
on $D^{+}$ and \emph{$Y_{k}^{-}\uparrow X_{-}$ on $D^{-}$}, we
have $Z_{k}\rightarrow X=X_{+}-X_{-}$ on $\bigcap_{k=1}^{\infty}domain(Z_{k})=D^{+}\cap D^{-}$.
Next, define $Z_{0}\equiv0$ and define $U_{k}\equiv Z_{k}-Z_{k-1}$
for each $k\geq1$. Then $\sum_{k=1}^{\infty}I(U_{k})<\infty$ and
$\sum_{k=1}^{\infty}U_{k}=X$ on $\bigcap_{k=1}^{\infty}domain(U_{k})$.
Hence $(U_{k})_{k=1,2,\cdots}$ is a representation of $X$ in $L$. 

3. The assertion is trivial if $X$ and $X'$ are integrable indicators.
Hence it is also valid if $X$ and $X'$ are linear combinations of
integrable indicators. Now suppose $X$ and $X'$ are integrable functions
bounded in absolute value by some $a>0$. By assertion 2, there exists
sequences $(Z_{n})_{n=1,2,\cdots}$ and $(Z'_{n})_{n=1,2,\cdots}$
of linear combinations of integrable indicators such that $I|X-Z_{n}|\rightarrow0$
and $I|X'-Z'_{n}|\rightarrow0$. Then, for each $n\geq1$, $Z_{n}Z_{n}'$
is integrable by the previous remarks, and $|XX'-Z_{n}Z_{n}'|\leq a|X-Z_{n}|+a|X'-Z'_{n}|$.
Therefore, by Theorem \ref{Thm. |X-Yn|<Zn where Yn,Zn integrable =000026 IZn->0 implies X integrable},
$XX'$ is integrable. 
\end{proof}

\section{Measurable Functions and Measurable Sets\label{sec, Measurable functions} }

In this section, let $(\Omega,L,I)$ be a complete integration space,
and let $(S,d)$ be a complete metric space with a fixed reference
point $x_{\circ}\in S$. In the case where $S=R$, it is understood
that $d$ is the Euclidean metric and that $x_{\circ}=0$. 

We will write $\mu A\equiv\mu(A)$ for the measure of an integrable
set $A$, and similarly write $\mu AB\equiv\mu(AB)$ for integrable
sets $A$ and $B$. Recall that $C_{ub}(S)$ is the space of bounded
and uniformly continuous real-valued functions on $S$. Recall from
the Notations and Conventions in the Introduction that if $X$ is
a real-valued function on $\Omega$ and if $t\in R$, then we use
the abbreviation $(t\leq X)$ for the subset $\{\omega\in domain(X):t\leq X(\omega)\}$.
Similarly with ``$\leq$'' replaced by ``$\geq$'', ``$<$'',
or ``$=$''. As usual we write $a_{b}$ interchangeably with $a(b)$
to lessen the burden on subscripts.
\begin{defn}
\label{Def. Measurable Function} \textbf{(Measurable functions).}
A function $X$ from $(\Omega,L,I)$ to the complete metric space
$(S,d)$ is called a \emph{measurable function} if, for each integrable
set $A$ and each $f\in C_{ub}(S)$, we have (i) $f(X)1_{A}\in L$,
and (ii)  $\mu(d(x_{\circ},X)>a)A\rightarrow0$ as $a\rightarrow\infty$.
A subset $B$ of $\Omega$ is said to be a \emph{measurable set} \index{measurable set}
if $B=(X=1)$ for some real-valued measurable indicator function $X$.
The set $(X=0)$ is then called a \emph{measure-theoretic complement}
\index{measure-theoretic complement} of $B$. We write $1_{B}$ for
a measurable indicator function of $B$\index{indicator of a measurable set},
and write $B^{c}$ for a measure-theoretic complement of $B$. If
the constant function $1$ is integrable, then Conditions (i) and
(ii) reduce to (i') $f(X)\in L$, and (ii')  $\mu(d(x_{\circ},X)$
as $a\rightarrow\infty$. $\square$
\end{defn}
It is obvious that if condition (ii) holds for one point $x_{\circ}\in S$,
then it holds for any point $x'_{\circ}\in S$. The next lemma shows
that, given condition (i), the measure in condition (ii) is well-defined
for all but countably many $a\in R$. Thus condition (ii) makes sense.
\begin{lem}
\label{Lem. f(X)1_A in L for all f in C(S) =00003D> (d(x0,X)>a)A integrable}
\textbf{\emph{(Integrability of some basic sets). }}Let $X$ be a
function from $\Omega$ to $S$. Suppose $f(X)1_{A}\in L$ for each
$f\in C_{ub}(S)$ and for each integrable set $A$. Let $A$ be an
arbitrary integrable set. Then the set $(d(x_{\circ},X)>a)A$ is integrable
for all but countably many $a\in R$. Thus $\mu(d(x_{\circ},X)>a)A$
is well-defined for all but countably many $a\in R$. 
\end{lem}
\begin{proof}
Let $n\geq0$ be arbitrary. Then $h_{n}\equiv1\wedge(n+1-d(x_{\circ},\cdot))_{+}\in C_{ub}(S)$
and so $h_{n}(X)1_{A}\in L$ by hypothesis. Hence all but countably
many $b\in(0,1)$ are regular points of $h_{n}(X)1_{A}$. Therefore
the set 
\[
(d(x_{\circ},X)>n+1-b)A=(h_{n}(X)1_{A}<b)A
\]
is integrable for all but countably many $b\in(0,1)$. Equivalently,
$(d(x_{\circ},X)>a)A$ is integrable for all but countably many $a\in(n,n+1)$.
Since $n\geq0$ is arbitrary, we see that $(d(x_{\circ},X)>a)A$ is
integrable for all but countably many points $a>0$. For each $a\leq0$,
the set $(d(x_{\circ},X)>a)A=A$ is integrable by hypothesis. 
\end{proof}
The next proposition gives an obviously equivalent condition to (ii)
in Definition \ref{Def. Measurable Function}.
\begin{prop}
\label{Prop. Alternative definition of measurablility} \textbf{\emph{(Alternative
definition of measurable functions).}} For each $n\geq0$, define
the function $h_{n}\equiv1\wedge(n+1-d(x_{\circ},\cdot))_{+}\in C_{ub}(S)$.
A function $X$ from $(\Omega,L,I)$ to the complete metric space
$(S,d)$ is a \emph{measurable function} iff, for each integrable
set $A$ and each $f\in C_{ub}(S)$, we have \emph{(i)} $f(X)1_{A}\in L$,
and \emph{(ii)}  $Ih_{n}(X)1_{A}\rightarrow\mu(A)$ as $n\rightarrow\infty$. 
\end{prop}
\begin{proof}
Suppose Conditions (i) and (ii) hold. Let $n\geq0$ be arbitrary.
We need to verify that the function $X$ is measurable. Then, since
$h_{n}\in C_{ub}(S)$, we have $h_{n}(X)1_{A}\in L$ by Condition
(i). Let $A$ be an arbitrary integrable set. Then, for each $n\geq1$
and $a>n+1$, 
\[
Ih_{n}(X)1_{A}\leq\mu(d(x_{\circ},X)\leq a)A\leq\mu(A).
\]
Letting $n\rightarrow\infty$, Condition (ii) and the last displayed
inequality imply that $\mu(d(x_{\circ},X)\leq a)A\rightarrow\mu(A)$
as $a\rightarrow\infty$. Equivalently $\mu(d(x_{\circ},X)>a)A\rightarrow0$
as $a\rightarrow\infty$. The conditions in Definition \ref{Def. Measurable Function}
are satisfied for $X$ to be measurable.

Conversely, suppose $X$ is measurable. Then Definition \ref{Def. Measurable Function}
of measurability implies Condition (i) in the present lemma. It implies
also that $\mu(d(x_{\circ},X)\leq a)A\rightarrow\mu(A)$ as $a\rightarrow\infty$.
At the same time, for each $a>0$ and $n>a$,
\[
\mu(d(x_{\circ},X)\leq a)A\leq Ih_{n}(X)1_{A}\leq\mu(A).
\]
Letting $a\rightarrow\infty$, we see that $Ih_{n}(X)1_{A}\rightarrow\mu(A)$.
\end{proof}
\begin{prop}
\label{Prop. Every Measurable X is defined a.e.} \textbf{\emph{(Basic
properties of measurable functions).}}

1. The domain of each measurable function is a full set. In particular
if $A$ is a measurable set, then $A\cup A^{c}$ is a full set. 

2. Each function that is equal a.e. to a measurable function is itself
measurable. 

3. Each integrable function is a real-valued measurable function.
Each integrable set is measurable. 
\end{prop}
\begin{proof}
1. Suppose $X$ is a measurable function. Let $A$ be an integrable
set. Let $f\equiv0$ be the constant $0$ function. Then $f\in C_{ub}(S)$.
Hence, by condition (i) in Definition \ref{Def. Measurable Function},
we have $f(X)1_{A}\in L$. Consequently $D\equiv domain(f(X)1_{A})$
is a full set. Since $domain(X)=domain(f(X))\supset D$, we see that
$domain(X)$ is a full set. In other words, $X$ is defined a.e. Now
let $A$ be an arbitrary measurable set. In other words, $1_{A}$
is measurable. Then e $A\cup A^{c}=domain(1_{A})$ is a full set according
to the previous argument.

2. Now suppose $Y$ is a function on $\Omega$, with values in $S$,
such that $Y=X$ a.e. where $X$ is a measurable function. Let $A$
be any integrable set. Let $f\in C_{ub}(S)$ be arbitrary. Then, by
condition (i) in Definition \ref{Def. Measurable Function}, we have
$f(X)1_{A}\in L$. Moreover, because $Y=X$ a.e., we have $f(Y)1_{A}=f(X)1_{A}$
a.e. Consequently $f(X)1_{A}\in L$. Again because $Y=X$ a.e., 
\[
\mu(d(x_{\circ},Y)>a)A=\mu(d(x_{\circ},X)>a)A\rightarrow0
\]
as $a\rightarrow\infty$. Thus the conditions in Definition \ref{Def. Measurable Function}
are verified for $Y$ to be measurable.

3. Next, let $X$ be any integrable function. Let $f\in C_{ub}(R)$
be arbitrary and let $A$ be an arbitrary integrable set. By Proposition
\ref{Prop. X, A in L and f in C_u(R) =00003D>  f(X)1_A in L}, we
have $f(X)1_{A}\in L$, which establishes condition (i) of Definition
\ref{Def. Measurable Function}. By Chebychev's inequality, $\mu(|X|>a)A\leq a^{-1}I|X|\rightarrow0$
as $a\rightarrow\infty$. Condition (ii) of Definition \ref{Def. Measurable Function}
follows. Hence $X$ is measurable. In particular, $1_{A}$ and $A$
are measurable.
\end{proof}
Suppose two real-valued measurable functions $X$ and $Y$ are indicators
to the same measurable set $A$. Then $X=Y$ on $D\equiv domain(X)\cap domain(Y)$
and so $X=Y$ a.e. Therefore the indicator $1_{A}$ is well-defined
relative to a.e. equality. Moreover $(X=0)D$ $=$ $(Y=0)D$ and so
$(X=0)=(Y=0)$ a.e. Hence the measure-theoretic complement is also
well-defined relative to a.e. equality.

The next proposition will be used repeatedly to construct measurable
functions from given ones.
\begin{prop}
\label{Prop. Basing seq of measurable functions on measurable partition}
\textbf{\emph{(Construction of a measurable function from pieces of
given measurable functions on measurable sets in a disjoint union).}}
Let $(S,d)$ be a complete metric space. Let $(X_{i},A_{i})_{i=1,2,\cdots}$
be a sequence where, for each $i,j\geq1$, $X_{i}$ is a measurable
function on $(\Omega,L,I)$ with values in $S$, \emph{and (i}) $A_{i}$
is a measurable subset of $\Omega$, \emph{(ii)} if $i\neq j$ then
$A_{i}A_{j}=\phi$, \emph{(iii)} $\bigcup_{k=1}^{\infty}A_{k}$ is
a full set, and \emph{(iv)} $\sum_{k=1}^{\infty}\mu A_{k}A=\mu A$
for each integrable set $A$. 

Define a function $X$ on $domain(X)\equiv\bigcup_{i=1}^{\infty}domain(X_{i})A_{i}$
by $X\equiv X_{i}$ on $domain(X_{i})A_{i}$, for each $i\geq1$.
Then $X$ is a measurable function on $\Omega$ with values in $S$. 

The same conclusion holds for a finite sequence $(X_{i},A_{i})_{i=1,\cdots,n}$.
\end{prop}
\begin{proof}
We will give the proof for the infinite case only. For each $n\geq1$
define $h_{n}\equiv1\wedge(n+1-d(x_{\circ},\cdot))_{+}\in C_{ub}(S)$. 

Let $f\in C_{ub}(S)$ be arbitrary, with $|f|\leq c$ on $S$ for
some $c>0$. Let $A$ be an arbitrary integrable set. Since
\[
\sum_{i=1}^{\infty}If(X_{i})1_{A(i)A}\leq c\sum_{i=1}^{\infty}\mu A_{i}A<\infty,
\]
the function $Y\equiv\sum_{i=1}^{\infty}f(X_{i})1_{A(i)A}$ is integrable.
At the same time $f(X)1_{A}=Y$ on the full set
\[
(\bigcup_{i=1}^{\infty}A_{i})(\bigcap_{i=1}^{\infty}domain(X_{i})).
\]
Hence $f(X)1_{A}$ is integrable. In particular $h_{n}(X)1_{A}$ is
integrable for each $n\geq1$. Moreover
\[
Ih_{n}(X)1_{A}=\sum_{i=1}^{\infty}Ih_{n}(X_{i})1_{A(i)A}\uparrow\sum_{i=1}^{\infty}\mu(A_{i}A)=\mu(A).
\]
Hence, by Lemma \ref{Prop. Alternative definition of measurablility},
$X$ is a measurable function.
\end{proof}
Next is a metric space lemma.
\begin{lem}
\label{Lem. f unif continuous on A and on B, and A,B well covers S =00003D> f unif cont}
\textbf{\emph{(Sufficient condition for uniform continuity on a metric
space).}} Let (S,d) be an arbitrary metric space. Let $A,B$ be subset
of $S$ and let $a>0$ be such that, for each $x\in S$ we have either
\emph{(i)} $(d(\cdot,x)<a)\subset A$, or \emph{(ii)} $(d(\cdot,x)<a)\subset B$.
Suppose $\lambda:S\rightarrow R$ is a function with $domain(\lambda)=S$
such that $\lambda$ is uniformly continuous on each of $A$ and $B$.
Then $\lambda$ is uniformly continuous on $S$.
\end{lem}
\begin{proof}
Let $\varepsilon>0$ be arbitrary. Since $\lambda$ is uniformly continuous
on each of $A$ and $B$, there exists $\delta_{0}>0$ so small that
$|\lambda(x)-\lambda(y)|<\varepsilon$ for each $x,y$ with $d(x,y)<\delta_{0}$,
provided that either $x,y\in A$ or $x,y\in B$. 

Let $\delta\equiv a\wedge\delta_{0}$. Consider each $x,y\in S$ with
$d(x,y)<\delta$. By hypothesis, either condition (i) or condition
(ii) holds. Assume that condition (i) holds. Then since $d(x,x)=0<a$
and $d(y,x)<\delta\leq a$ we have $x,y\in A$. Hence, since $d(y,x)<\delta\leq\delta_{0}$,
we have $|\lambda(x)-\lambda(y)|<\varepsilon$. Similarly, if condition
(ii) holds, then $|\lambda(x)-\lambda(y)|<\varepsilon$. Combining,
we see that $\lambda$ is uniformly continuous on $S$. 
\end{proof}
\begin{prop}
\label{Prop. X meas,  f unif continuous and bd on bd subsets =00003D> f(X) meas}
\textbf{\emph{(A continuous function of a measurable function is  measurable).}}
Let $(S,d)$ and $(S'd')$ be complete metric spaces. Let $X$ be
a measurable function on $(\Omega,L,I)$, with values in $S$. Suppose
a function $f:(S,d)\rightarrow(S',d')$ with $domain(f)=S$ is uniformly
continuous on each bounded subset of $S$, and bounded on each bounded
subset of $S$. Then the composite function $f(X)\equiv f\circ X$
is measurable. In particular, $d(x,X)$ is a real-valued measurable
function for each $x\in S$. 
\end{prop}
\begin{proof}
We need to prove that $Y\equiv f(X)$ is measurable. To that end,
let $g\in C_{ub}(S')$ be arbitrary, with $|g|\leq b$ for some $b>0$.
Consider arbitrary integrable set $A$ and $\varepsilon>0$. Since
$X$ is measurable by hypothesis, there exists $a>0$ so large that
$\mu(B)<\varepsilon$ where $B\equiv(d(x_{\circ},X)>a)A$. Define
$h\equiv1\wedge(a-d(x_{\circ},\cdot))_{+}\in C(S)$.

The function $f$ is, by hypothesis, uniformly continuous on the bounded
set $G\equiv(d(\cdot,x_{\circ})<2+a)$. By assumption, $g$ is uniformly
continuous. Therefore $(g\circ f)$ and $(g\circ f)h$ are uniformly
continuous on $G$. At the same time $h=0$ on $H\equiv(d(\cdot,x_{\circ})>a)$.
Hence $(g\circ f)h=0$ on $H$. Thus $(g\circ f)h$ is uniformly continuous
on $H$.

Now consider each $x\in S$. Either (i) $d(x,x_{\circ})<a+\frac{3}{2}$,
or (ii) $d(x,x_{\circ})>a+\frac{1}{2}$. In Case (i), we have $(d(\cdot,x)<\frac{1}{2})\subset(d(\cdot,x_{\circ})<2+a)\equiv G$.
In Case (ii), we have $(d(\cdot,x)<\frac{1}{2})\subset(d(\cdot,x_{\circ})>a)\equiv H$.
Combining, Lemma \ref{Lem. f unif continuous on A and on B, and A,B well covers S =00003D> f unif cont}
implies that $(g\circ f)h$ is uniformly continuous on S. Moreover,
since $(g\circ f)h$ is bounded on $G$ by hypothesis, and is equal
to $0$ on $H$, it is bounded on $S$. In short $(g\circ f)h\in C_{ub}(S)$.
Since $X$ is measurable, the function $g(Y)h(X)1_{A}=(g\circ f)(X)h(X)1_{A}$
is integrable. At the same time,
\[
|g(Y)1_{A}-g(Y)h(X)1_{A}|\leq b(1-h(X))1_{A}
\]
where 
\[
I(1-h(X))1_{A}\leq\mu(d(x_{\circ},X)>a)A=\mu(B)<\varepsilon.
\]
Since $\varepsilon>0$ is arbitrary, Theorem \ref{Thm. |X-Yn|<Zn where Yn,Zn integrable =000026 IZn->0 implies X integrable}
implies that the function $g(Y)1_{A}$ is integrable.

Now let $c>a$ be arbitrary. By hypothesis, there exists $c'>0$ so
large that $|f|\leq c'$ on $(d(x_{\circ},\cdot)<c)$. Combining,
we see that
\[
\mu(d(x_{\circ},Y)>c')A\equiv\mu(d(x_{\circ},f(X))>c')A
\]
\[
\leq\mu(d(x_{\circ},X)\geq c)A\leq\mu(d(x_{\circ},X)>a)A<\varepsilon.
\]
Since $\varepsilon>0$ is arbitrary, we conclude that $\mu(d(x_{\circ},Y)>c')A\rightarrow0$
as $c'\rightarrow\infty$. Thus we have verified the conditions of
Definition \ref{Def. Measurable Function} for $Y$ to be measurable.
In other words, $f(X)$ is measurable.
\end{proof}
\begin{cor}
\label{Cor. Identity function X(x)=00003Dx  is meas on completion of (S,C(S),I)}
(\textbf{\emph{Condition for measurability of identity function, and
of continuous function of a measurable function).}} Let (S,d) be a
complete metric space. Suppose $(S,C_{ub}(S),I)$ is an integration
space, with completion $(S,L,I)$, such that $Ih_{k}1_{A}\uparrow\mu(A)$
for each integrable set $A$, where $h_{k}\equiv1\wedge(1+k-d(\cdot,x_{\circ}))_{+}$
for each $k\geq1$, Then the following holds.

1. The identity function $X:(S,L,I)\rightarrow(S,d)$, defined by
$X(x)\equiv x$ for each $x\in S$, is a measurable function on $(S,L,I)$.

2. Let $(S',d')$ be a second complete metric space. Suppose a function
$f:(S,d)\rightarrow(S',d')$ with $domain(f)=S$ is uniformly continuous
and bounded on each bounded subset of $S$. Then the function $f:(S,L,I)\rightarrow(S',d')$
is measurable. In particular, $d(x,\cdot)$ is a real-valued measurable
function for each $x\in S$. 
\end{cor}
\begin{proof}
1. Let $f\in C_{ub}(S)$ be arbitrary, and let $A$ be an arbitrary
integrable set. Then $f(X)\equiv f\in L$. Hence $f(X)1_{A}\in L$
. Moreover, $Ih_{k}(X)1_{A}=Ih_{k}1_{A}\uparrow\mu(A)$ by hypothesis.
Hence $X$ is measurable according to Lemma \ref{Prop. Alternative definition of measurablility}.

2. The conditions in the hypothesis of Proposition \ref{Prop. X meas,  f unif continuous and bd on bd subsets =00003D> f(X) meas}
are satisfied by the functions $X:(S,L,I)\rightarrow(S,d)$ and $f:(S,d)\rightarrow(S',d')$.
Accordingly, the function $f\equiv f(X)$ is measurable. 
\end{proof}
The next proposition says that, in the case where $(S,d)$ is locally
compact, the conditions for measurability in Definition \ref{Def. Measurable Function}
can be weakened somewhat, by replacing $C_{ub}(S)$ with the subset
$C(S)$.
\begin{prop}
\label{Prop. Sufficient condition of measurability in case S locally compact. }
\textbf{\emph{(Sufficient condition for measurability in case $S$
is locally compact).}} Let $(S,d)$ be a locally compact metric space.
Define $h_{n}\equiv1\wedge(n+1-d(x_{\circ},\cdot))_{+}\in C(S)$ for
each $n\geq1$. Let $X$ be a function from $(\Omega,L,I)$ to $(S,d)$
such that $f(X)1_{A}\in L$ for each integrable set $A$ and each
$f\in C(S)$. Then the following holds.
\end{prop}
\begin{enumerate}
\item If $Ih_{n}(X)1_{A}\uparrow\mu(A)$ for each integrable set $A$, then
$X$ is measurable.
\item If $\mu(d(x_{\circ},X)>a)A\rightarrow0$ as $a\rightarrow\infty$,
for each integrable set $A$, then $X$ is measurable.
\end{enumerate}
\begin{proof}
Let $A$ be an arbitrary integrable set. Let $g\in C_{ub}(S)$ be
arbitrary, with $|g|\leq c$ on $S$. For each $n\geq1$, since $S$
is locally compact, we have $h_{n},h_{n}g\in C(S)$, Hence, for each
$n\geq1$, we have $h_{n}(X)g(X)1_{A}\in L,$and, by hypothesis, 
\begin{equation}
|g(X)1_{A}-h_{n}(X)g(X)1_{A}|\leq c(1-h_{n}(X))1_{A}.\label{eq:temp-281}
\end{equation}

1. By hypothesis, $I(1-h_{n}(X))1_{A}\rightarrow0$. Hence, in view
of inequality \ref{eq:temp-281}, Theorem \ref{Thm. |X-Yn|<Zn where Yn,Zn integrable =000026 IZn->0 implies X integrable}
is applicable, to yield $g(X)1_{A}\in L,$ where $g\in C_{ub}(S)$
is arbitrary. Thus the conditions in Proposition \ref{Prop. Alternative definition of measurablility}
are satisfied for $X$ to be measurable. 

2. For each $a>0$ and for each $n>a$, 
\[
0\leq\mu(A)-Ih_{n}(X)1_{A}=I(1-h_{n}(X))1_{A}\leq\mu(d(x_{\circ},X)>a)A,
\]
which, by hypothesis, converges to $0$ as $a\rightarrow\infty$.
Hence, by Assertion 1, $X$ is measurable.
\end{proof}
\begin{defn}
\label{Def. Regular =000026 Continuity Pts of Measurable Func} \textbf{(Regular-
and Continuity points of a measurable function relative to each integrable
set).} Suppose $X$ is a real-valued measurable function on $(\Omega,L,I)$.
We say that $t\in R$ is a \index{regular point of measurable function}\emph{regular
point} \emph{of $X$ relative to an integrable set} $A$ if (i) there
exists a sequence $(s_{n})_{n=1,2,\cdots}$ of real numbers decreasing
to $t$ such that $(s_{n}<X)A$ is integrable for each $n\geq1$ and
such that $\lim_{n\rightarrow\infty}\mu(s_{n}<X)A$ exists, and (ii)
there exists a sequence $(r_{n})_{n=1,2,\cdots}$ of real numbers
increasing to $t$ such that $(r_{n}<X)A$ is integrable for each
$n\geq1$ and such that $\lim_{n\rightarrow\infty}\mu(r_{n}<X)A$
exists. If, in addition, the two limits in (i) and (ii) are equal,
then we call $t$ a \index{continuity point of measurable function}\emph{continuity
point} of $X$ \emph{relative to} $A$. 

If a real number $t$ is a regular point of $X$ relative to each
integrable set $A$, then we call $t$ a \emph{regular point} of $X$.
If a real number $t$ is a continuity point of $X$ relative to each
integrable set $A$, then we say $t$ is a \emph{continuity point}
of $X$. $\square$
\end{defn}
The next proposition shows that regular points and continuity points
of a real-valued measurable function are abundant, and that they inherit
the properties of regular points and continuity points of integrable
functions.
\begin{prop}
\label{Prop.  Countable Exceptional Pts for Meas X int A} \textbf{\emph{(All
but countably many points are continuous points of a real measurable
function, relative to each given integrable set $A$).}} Let $X$
be a real-valued measurable function on $(\Omega,L,I)$. Let $A$
be an integrable set and let $t$ be a regular point of $X$ relative
to $A$.
\end{prop}
\begin{enumerate}
\item \emph{All but countably many $u\in R$ are continuity points of $X$
relative to $A.$ Hence all but countably many $u\in R$ are regular
points of $X$ relative to $A.$}
\item \emph{$A^{c}$ is a measurable set.}
\item \emph{The sets $A(t<X),A(t\leq X),A(X<t)$,$A(X\leq t)$, and $A(X=t)$
are integrable sets.}
\item \emph{$(X\leq t)A=A((t<X)A)^{c}$ a.e., and $(t<X)A=A((X\leq t)A)^{c}$
a.e.}
\item \emph{$(X<t)A=A((t\leq X)A)^{c}$ a.e., and $(t\leq X)A=A((X<t)A)^{c}$
a.e.}
\item \emph{For a.e. $\omega\in A$, we have $t<X(\omega)$, $t=X(\omega)$,
or $t>X(\omega)$. Thus we have a limited, but useful, version of
the principle of excluded middle.}
\item \emph{Let $\varepsilon>0$ be arbitrary. There exists $\delta>0$
such that if $r\in(t-\delta,t]$ and $A(X<r)$ is integrable, then
$\mu(A(X<t))-\mu(A(X<r))<\varepsilon$. There exists $\delta>0$ such
that if $s\in[t,t+\delta)$ and $A(X\leq s)$ is integrable, then
$\mu(A(X\leq s))-\mu(A(X\leq t))<\varepsilon$.}
\item \emph{If $t$ is a continuity point of $X$ relative to $A$, then
$\mu((t<X)A)=$ $\mu((t\leq X)A)$.}
\end{enumerate}
\begin{proof}
In the special case where $X$ is an integrable function, the assertions
have been proved in Corollary \ref{Cor. Properties of Regular Pts of Integrable Func}.
In general, suppose $X$ is a real-valued measurable function. Let
$n\geq1$ be arbitrary. Then, by Definition \ref{Def. Measurable Function},
$((-n)\vee X\wedge n)1_{A}\in L$. Hence all but countably many $u\in(-n,n)$
are continuity points of the integrable function $((-n)\vee X\wedge n)1_{A}$
relative to $A$. On the other hand, for each $t\in(-n,n)$, we have
$(X<t)A=(((-n)\vee X\wedge n)1_{A}<t)A$. Hence, a point $u\in(-n,n)$
is a continuity point of $((-n)\vee X\wedge n)1_{A}$ relative to
$A$ iff it is a continuity point of $X$ relative to $A$. Combining,
we see that all but countably many points in the interval $(-n,n)$
are continuity points of $X$ relative to $A$. Therefore all but
countably many points in $R$ are continuity points of $X$ relative
to $A$. This proves assertion 1. The remaining assertions are proved
by similarly reducing to assertions about the integrable functions
$((-n)\vee X\wedge n)1_{A}$. 
\end{proof}
Suppose $X$ is a measurable function. Note that we defined the regular
points and continuity points of $X$ relative to a specific integrable
set $A$. In the case of a $\sigma$-finite integration space, to
be defined next, all but countably many real numbers $t$ are regular
points of $X$ relative to each integrable set $A$.
\begin{defn}
\label{Def. Sigma Finite; I-basis} \textbf{($\sigma$-finiteness
and $I$-basis).}\emph{ The complete integration space $(\Omega,L,I)$
is said to be \index{finite integration space}finite i}f the constant
function $1$ is integrable, and $(\Omega,L,I)$ is said to be sigma
finite, or $\sigma$\emph{-finite} \index{sigma-finite}, if there
exists a sequence $(A_{k})_{k=1,2,\cdots}$ of integrable sets with
positive measures such that (i) $A_{k}\subset A_{k+1}$ for $k=1,2,\cdots$
, (ii) $\bigcup_{k=1}^{\infty}A_{k}$ is a full set, and (iii) for
any integrable set $A$ we have $\mu(A_{k}A)\rightarrow\mu(A)$. The
sequence $(A_{k})_{k=1,2,\cdots}$ is then called an $I$-\emph{basis}
\index{I-basis} for $(\Omega,L,I)$. $\square$
\end{defn}
If $(\Omega,L,I)$ is finite, then it is $\sigma$-finite, with an
\emph{$I$-}basis\emph{ }given by $(A_{k})_{k=1,2,\cdots}$ where
$A_{k}\equiv\Omega$ for each $k\geq1$. In particular, if $(S,C_{ub}(S),I)$
is an integration space, with completion $(S,L,I)$, then the constant
function 1 is integrable, and so $(S,L,I)$ is finite. 
\begin{lem}
\textbf{\emph{\label{Lem. Completion of integration on locally compact S results in sigma finite}
(Completion of an integration on a locally compact metric space results
in a$\sigma$-finite integration space).}} Suppose $(S,d)$ is locally
compact. Let $(S,L,I)$ be the completion of some integration space
$(S,C(S),I)$. Then $(\Omega,L,I)\equiv(S,L,I)$ is $\sigma$-finite.
Specifically, there exists an increasing sequence $(a_{k})_{k=1,2,\cdots}$
of positive real numbers with $a_{k}\uparrow\infty$ such that\emph{
}$(A_{k})_{k=1,2,\cdots}$ is an $I$-\emph{basis} for $(S,L,I)$,
where $A_{k}\equiv(d(x_{\circ},\cdot)\leq a_{k})$ for each $k\geq1$. 
\end{lem}
\begin{proof}
Consider each $k\geq1$. Define $X_{k}\equiv1\wedge(k+1-d(x_{\circ},\cdot))_{+}\in C(S)\subset L$.
Let $c\in(0,1)$ be arbitrary and let $a_{k}\equiv k+1-c$. Then the
set $A_{k}\equiv(d(x_{\circ},\cdot)\leq a_{k})=(X_{k}\geq c)$ is
integrable. Conditions (i) and (ii) in Definition \ref{Def. Sigma Finite; I-basis}
are easily verified. For condition (iii), consider any integrable
set $A$. According to Assertion 2 of Corollary \ref{Cor. Identity function X(x)=00003Dx  is meas on completion of (S,C(S),I)},
the real-valued function $d(x_{\circ},\cdot)$ is measurable on $(\Omega,L,I)\equiv(S,L,I)$.
Hence 
\[
\mu A-\mu A_{k}A=\mu(d(x_{\circ},\cdot)>a_{k})A\rightarrow0
\]
Thus condition (iii) in Definition \ref{Def. Sigma Finite; I-basis}
is also verified for $(A_{k})_{k=1,2,\cdots}$ to be an $I$-basis. 
\end{proof}
\begin{prop}
\label{Prop. Continuity Pts for Measurable Function} \textbf{\emph{(In
the case of a $\sigma$-finite integration space, all but countably
many points are continuous points of a real measurable function).}}
Suppose $X$ is a real-valued measurable function on a $\sigma$-finite
integration space $(\Omega,L,I)$. Then all but countably many real
numbers $t$ are continuity points, hence regular points, of $X$.
\end{prop}
\begin{proof}
Let $(A_{k})_{k=1,2,\cdots}$ be an $I$-basis for $(\Omega,L,I)$.
According to Proposition \ref{Prop.  Countable Exceptional Pts for Meas X int A},
for each $k$ there exists a countable subset $J_{k}$ of $R$ such
that if $t\in(J_{k})_{c}$, where $(J_{k})_{c}$ stands for the metric
complement of $J_{k}$ in $R$, then $t$ is a continuity point of
$X$ relative to $A_{k}$. Define $J\equiv\bigcup_{k=1}^{\infty}J_{k}$. 

Consider each $t\in J_{c}$. Let the integrable set $A$ be arbitrary.
According to condition (iii) in Definition \ref{Def. Sigma Finite; I-basis},
we can select a subsequence $(A_{k(n)})_{n=1,2,\cdots}$ of $(A_{k})$
such that $\mu(A)-\mu(AA_{k(n)})<\frac{1}{n}$ for each $n\geq1$.
Let $n\geq1$ be arbitrary. Write $B_{n}\equiv A_{k(n)}$. Then $t\in(J_{k(n)})_{c}$
and so $t$ is a continuity point of $X$ relative to $B_{n}$. Consequently,
according to Proposition \ref{Prop.  Countable Exceptional Pts for Meas X int A},
the sets $(X<t)B_{n}$ and $(X\leq t)B_{n}$ are integrable, with
$\mu(X<t)B_{n}=\mu(X\leq t)B_{n}$. Furthermore, according to the
same proposition, there exists $\delta_{n}>0$ such that (i) if $r\in(t-\delta_{n},t]$
and $(X<r)B_{n}$ is integrable, then $\mu(X<t)B_{n}-\mu(X<r)B_{n}<\frac{1}{n}$,
and (ii) if $s\in[t,t+\delta_{n})$ and $(X\leq s)B_{n}$ is integrable,
then $\mu(X\leq s)B_{n}-\mu(X\leq t)B_{n}<\frac{1}{n}$. Let $r_{0}\equiv t-1$
and $s_{0}\equiv t+1$. Inductively we can select  $r_{n}\in(t-\delta_{n},t)\cap(r_{n-1},t)\cap(t-\frac{1}{n},t)$
such that $r_{n}$ is a regular point of $X$ relative to both $B_{n}$
and $A$. Similarly we can select $s_{n}\in(t,t+\delta_{n})\cap(t,s_{n-1})\cap(t,t+\frac{1}{n})$
such that $s_{n}$ is a regular point of $X$ relative to both $B_{n}$
and $A$. Then, for each $n\geq1$, we have 
\[
\mu(r_{n}<X)A-\mu(s_{n}<X)A
\]
\[
=\mu A(r_{n}<X)(X\leq s_{n})A\leq\mu(r_{n}<X)(X\leq s_{n})B_{n}+\mu(AB_{n}^{c})
\]
\[
=\mu(X\leq s_{n})B_{n}-\mu(X\leq r_{n})B_{n}+\mu(A)-\mu(AB_{n})
\]
\[
\leq(\mu(X\leq t)B_{n}+\frac{1}{n})-(\mu(X<t)B_{n}-\frac{1}{n})+(\mu(A)-\mu(AA_{k(n)}))
\]
\begin{equation}
\leq\frac{1}{n}+\frac{1}{n}+\frac{1}{n}\label{eq:temp-5}
\end{equation}
Since the sequence $(\mu(r_{n}<X)A$) is nonincreasing and the sequence
$(\mu(s_{n}<X)A)$ is nondecreasing, inequality \ref{eq:temp-5} implies
that both sequences converge, and to the same limit. By Definition
\ref{Def. Regular =000026 Continuity Pts of Measurable Func}, $t$
is a continuity point of $X$ relative to $A$. 
\end{proof}
We now expand Convention \ref{Convention. Only regular Pts of Integrable X used}
to cover also measurable functions.
\begin{defn}
\label{Convention. Only Regular Pts of Measurble Funcs used} \textbf{(Convention
regarding regular points of measurable functions).} Let $X$ be a
real-valued measurable function, and let $A$  be an integrable set.
Henceforth, when the integrability of the set $(X<t)A$ or $(X\leq t)A$
is required in a discussion, for some $t\in R$, then it is understood
that the real number $t$ has been chosen from the regular points
of the measurable function $X$ relative to the given integrable set
$A$. 

Furthermore, if $(\Omega,L,I)$ is a $\sigma$-finite integration
space, when the measurability of the set $(X<t)$ or $(X\leq t)$
is required in a discussion, for some $t\in R$, then it is understood
that the real number $t$ has been chosen from the regular points
of the measurable functions $X$.$\square$
\end{defn}
\begin{cor}
\label{Cor. In sigma-finite space, (X<=00003Dt)^c=00003D(t<X) etc}
\textbf{\emph{(Properties of regular points).}} Let $X$ be a real-valued
measurable function on a $\sigma$-finite integration space $(\Omega,L,I)$,
and let $t$ be a regular point of $X$. Then \emph{$(X\leq t)=(t<X){}^{c}$
}and $(t<X)=(X\leq t){}^{c}$ are measurable sets. Similarly, $(X<t)=(t\leq X){}^{c}$,
and $(t\leq X)=(X<t){}^{c}$ are measurable sets. The equalities here
are understood to be a.e. equalities, and the measure theoretic complement
of a measurable set is defined a.e.
\end{cor}
\begin{proof}
We will prove only the first alleged equality, the rest being similar.
Define an indicator function $Y$, with $domain(Y)=(X\leq t)\cup(t<X)$,
by $Y=1$ on $(X\leq t)$ and $Y=0$ on $(t<X)$. It suffices to show
that $Y$ satisfies conditions (i) and (ii) of Definition \ref{Def. Measurable Function}
for a measurable function. To that end, consider an arbitrary $f\in C_{ub}(R)$
and an arbitrary integrable subset $A$. By hypothesis, and by Definition
\ref{Convention. Only Regular Pts of Measurble Funcs used}, $t$
is a regular point of $X$ relative to $A$. Moreover
\[
f(Y)1_{A}=f(1)1_{(X\leq t)A}+f(0)1_{(t<X)A},
\]
which is integrable according to Proposition \ref{Prop.  Countable Exceptional Pts for Meas X int A}.
Thus condition (i) of Definition \ref{Def. Measurable Function} is
verified. Moreover, since $Y$ has only the possible values of $0$
and 1, the set $(|Y|>a)$ is empty for each $a>1$. Hence, trivially,
$\mu(|Y|>a)A\rightarrow0$ as $a\rightarrow\infty$, thereby establishing
also condition (ii) of Definition \ref{Def. Measurable Function}.
Consequently $Y$ is a measurable indicator for $(X\leq t)$, and
$(t<X)=(Y=0)=(X\leq t)^{c}$. 
\end{proof}
\begin{prop}
\label{Prop. vector of measurable func (S)   is meas func (S)} \textbf{\emph{(A
vector of measurable functions constitutes a measurable function).}}
Let $(\Omega,L,I)$ be a complete integration space. 

1. Let $(S',d'),(S'',d'')$ be complete metric spaces. Let $(\widetilde{S},\widetilde{d})\equiv(S'\times S'',d'\otimes d'')$
be their product metric space. Let $X':\varOmega\rightarrow S'$ and
$X'':\varOmega\rightarrow S''$ be measurable functions. Define $X:\Omega\rightarrow\widetilde{S}$
by $X(\omega)\equiv(X'(\omega),X''(\omega))$ for each $\omega\in domain(X)\equiv domain(X')\cap domain(X'')$.
Then $X$ is a measurable function. 

2. Suppose $(S,d)$ is a complete metric space and a function $g:\widetilde{S}\rightarrow S$
is \emph{(i)} uniformly continuous on bounded subsets, and \emph{(ii)}
bounded on bounded subsets. Then $g(X)\equiv g(X',X'')$ is a measurable
function with values in $S$. 

The above assertions generalize trivially to any finite number of
measurable functions $X',X'',\cdots,X^{(n)}$.

3. As a special case of Assertion 2 above, suppose $(S',d')=(S'',d'')$.
Then for arbitrary measurable functions $X',X'':\varOmega\rightarrow S'$
, the function $d'(X',X''):\Omega\rightarrow R$ is measurable.
\end{prop}
\begin{proof}
1. Let $x'_{\circ},x''_{\circ}$ be arbitrary, but fixed, reference
points of $(S',d'),(S'',d'')$ respectively. Designate the point $x_{\circ}\equiv(x'_{\circ},x''_{\circ})$
as the reference point in the product $(\widetilde{S},\widetilde{d})$.
Consider each integrable set $A$, and consider each $f\in C_{ub}(\widetilde{S})$
with a bound $b>0$ and with a modulus of continuity $\delta_{f}$.
Let $\varepsilon>0$ be arbitrary. Since $X'$ and $X''$ are, by
hypothesis, measurable, there exists a continuity point $a>0$ of
the measurable functions $d'(x'_{\circ},X'),$$d''(x''_{\circ},X'')$
so large that $\mu(d'(x'_{\circ},X')>a)A<\varepsilon$ and $\mu(d''(x''_{\circ},X'')>a)A<\varepsilon.$ 

Now take a sequence $-a\equiv a_{0}<a_{1}<\cdots<a_{n}\equiv a$ such
that $a_{i}-a_{i-1}<\delta_{f}(\varepsilon)$ for each $i=1,\cdots,n$.
For each $i,j=1,\cdots,n$, define the integrable set 
\[
\Delta_{i,j}\equiv(a_{i-1}<d'(x'_{\circ},X')\leq a_{i})A\cup(a_{j-1}<d''(x''_{\circ},X'')\leq a_{j})A.
\]
Partition the set $H\equiv\{(i,j):1\leq i\leq n;1\leq j\leq n\}$
of double subscripts into two disjoint subsets $H_{1}$ and $H_{2}$
such that (i) $\mu\Delta_{i,j}>0$ for each $(i,j)\in H_{1},$ and
(ii) $\mu\Delta_{i,j}<n^{-2}\varepsilon$ for each $(i,j)\in H_{2}$.
Then, for each $(i,j)\in H_{1},$ there exists $\omega_{i,j}\in\Delta_{i,j}$.
Note that, for each $(i,j)\in H_{1},$ we have 
\[
\widetilde{d}(X(\omega_{i,j}),X)\equiv d'(X'(\omega_{i,j}),X')\vee d''(X''(\omega_{i,j}),X'')\leq|a_{i}-a_{i-1}|\vee|a_{j}-a_{j-1}|<\delta_{f}(\varepsilon)
\]
on $\Delta_{i,j}$, whence
\[
|f(X(\omega_{i,j})-f(X)|<\varepsilon
\]
on $\Delta_{i,j}$. Note also that $\Delta_{i,j}\Delta_{i',j'}=\phi$
for each $(i,j),(i',j')\in H$ with $(i,j)\neq(i',j')$. Moreover,
the set $\Delta\equiv(\bigcup_{(i,j)\in H}\Delta_{i,j})^{c}$ is measurable,
with 
\[
D\equiv\Delta\cup\bigcup_{(i,j)\in H}\Delta_{i,j}
\]
equal to a full set. At the same time $\Delta_{i,j}\Delta=\phi$ for
each $(i,j)\in H$. Consequetnly, for each integrable set $B$, we
have 
\[
\mu\Delta B+\bigcup_{(i,j)\in H}\mu\Delta_{i,j}B=\mu(\Delta\cup\bigcup_{(i,j)\in H}\Delta_{i,j})B=\mu B.
\]
Thus the finite family $\{\Delta_{i,j}:(i,j)\in H\}\cup\{\Delta\}$
of measurable sets satisfies Conditions (i-iv) in Proposition \ref{Prop. Basing seq of measurable functions on measurable partition}.
Accordingly, we can define two integrable functions $Y,Z:\Omega\rightarrow R$
by $domain(Y)\equiv domain(Z)\equiv D$ and by (i') $Y\equiv f(X(\omega_{i,j}))$
and $Z\equiv\varepsilon1_{A}$ on $\Delta_{i,j}$, for each $(i,j)\in H_{1},$
(ii') $Y\equiv0$ and $Z\equiv b1_{A}$ on $\Delta_{i,j}$, for each
$(i,j)\in H_{2},.$ and (iii') $Y\equiv0$ and $Z\equiv b1_{A}$ ,
on $\Delta$. Then, for each $(i,j)\in H_{1},$ we have 
\[
|Y-f(X)|1_{A}=|f(X(\omega_{i,j}))-f(X)|1_{A}<\varepsilon1_{A}=Z
\]
on $\Delta_{i,j}$. At the same time, for each $(i,j)\in H_{2},$
we have 
\[
|Y-f(X)|1_{A}=|0-f(X)|1_{A}\leq b1_{A}=Z
\]
on $\Delta_{i,j}$. Likewise, 
\[
|Y-f(X)|1_{A}\leq|0-f(X)|1_{A}\leq b1_{A}=Z
\]
on $\Delta$. Summing up, we obtain
\[
|Y1_{A}-f(X)1_{A}|\leq Z
\]
 on the full set $D\equiv\Delta\cup\bigcup_{(i,j)\in H}\Delta_{i,j}$.
Now estimate
\[
IZ=\sum_{(i,j)\in H(1)}IZ1_{\Delta(i,j)}+\sum_{(i,j)\in H(2)}IZ1_{\Delta(i,j)}+IZ1_{\Delta}
\]
\[
=\sum_{(i,j)\in H(1)}\varepsilon I1_{A\Delta(i,j)}+\sum_{(i,j)\in H(2)}bI1_{A\Delta(i,j)}+bI1_{A\Delta}
\]
\[
=\varepsilon\sum_{(i,j)\in H(1)}\mu(A\Delta_{i,j})+b\sum_{(i,j)\in H(2)}\mu(A\Delta_{i,j})+b\mu(A\Delta)
\]
\[
\leq\varepsilon\mu(\bigcup_{(i,j)\in H(1)}A\Delta_{i,j})+b\sum_{(i,j)\in H(2)}\mu(\Delta_{i,j})+b\mu((d'(x'_{\circ},X')>a)A\cup(d''(x''_{\circ},X'')>a)A)
\]
\[
\leq\varepsilon\mu(A)+bn^{2}n^{-2}\varepsilon+b(\varepsilon+\varepsilon)=\varepsilon\mu(A)+3b\varepsilon,
\]
where $\varepsilon>0$ is arbitrary. By repeating the above argument
with a sequence $(\varepsilon_{k})_{k=1,2,\cdots}$with $\varepsilon_{k}\downarrow0$,
we can construct a two sequences of integrable functions $(Y_{k}1_{A})_{k=1,2,\cdots}$
and $(Z_{k})_{k=1,2,\cdots}$ such that 
\[
|Y_{k}1_{A}-f(X)1_{A}|\leq Z_{k}\qquad\mathrm{a.s}.
\]
and such that
\[
IZ_{k}\leq\varepsilon_{k}\mu(A)+3b\varepsilon_{k}\downarrow0
\]
as $k\rightarrow\infty$. The conditions in Theorem \ref{Thm. |X-Yn|<Zn where Yn,Zn integrable =000026 IZn->0 implies X integrable}
are satisfied, to yield $f(X)1_{A}\in L$, where $f\in C_{ub}(S)$
and the integrable set $A$ are arbitrary. At the same time,
\[
(\widetilde{d}(x_{\circ},X)>a)A\equiv(d'(x'_{\circ},X')\vee d''(x''_{\circ},X'')>a)A
\]
is integrable, with 
\[
\mu(\widetilde{d}(x_{\circ},X)>a)A\leq\mu(d'(x'_{\circ},X')>a)A+\mu(d(x''_{\circ},X'')>a)A<2\varepsilon,
\]
where $\varepsilon>0$ is arbitrarily small. All the conditions in
Proposition \ref{Prop. Sufficient condition of measurability in case S locally compact. }
for $X:\Omega\rightarrow\widetilde{S}$ to be measurable have thus
been verified. Assertion 1 is proved.

2. Next, suppose $(S,d)$ is a complete metric space and a function
$g:\widetilde{S}\rightarrow S$ is (i) uniformly continuous on bounded
subsets, and (ii) bounded on bounded subsets. Then the composite function
$g(X):\Omega\rightarrow S$ is a measurable function according to
Proposition \ref{Prop. X meas,  f unif continuous and bd on bd subsets =00003D> f(X) meas}.
Assertion 2 is proved.

3. Now suppose $(S',d')=(S'',d'')$, and suppose $S=R$. Then the
distance function $d':\widetilde{S}\equiv S'^{2}\rightarrow R$ is
uniformly continuous, and is bounded on bounded subsets. Hence $d'(X',X'')$
is measurable by Assertion 2.
\end{proof}
\begin{cor}
\label{Cor. X,Y measurable =00003D> aX+bY, |X|, X**a measurable}
\textbf{\emph{(Operations preserving measurability).}} Let $X,Y$
be real-valued measurable functions. Then $aX+bY$, $1$, $X\vee Y$,
$X\wedge Y$, $|X|$, and $X^{\alpha}$ are measurable functions for
any real numbers $a,b,\alpha$ with $\alpha\geq0$. Let $A,B$ be
measurable sets. Then $A\cup B$ and $AB$ are measurable. Moreover
$(A\cup B)^{c}=A^{c}B^{c}$ and $(AB)^{c}=A^{c}\cup B^{c}$.
\end{cor}
\begin{proof}
$aX+bY$, $1$, $X\vee Y$, $X\wedge Y$, $|X|$, and $X^{\alpha}$
are uniformly continuous real-valued functions of $X$ and $Y$ which
are bounded on bounded subsets of $R$, and so are measurable by Assertion
2 of Proposition \ref{Prop. vector of measurable func (S)   is meas func (S)}.
Let $A,B$ be measurable sets with indicators $U,V$ respectively.
Then $U\vee V$ is a measurable indicator, with
\[
(U\vee V=1)=(U=1)\cup(V=1)=A\cup B.
\]
Hence $A\cup B$ is a measurable set, with $U\vee V$ as indicator.
Moreover, $(A\cup B)^{c}=(U\vee V=0)=(U=0)(V=0)=A^{c}B^{c}$. Similarly
$AB$ is measurable, with $(AB)^{c}=A^{c}\cup B^{c}$.
\end{proof}
\begin{prop}
\label{Prop. measurable dominated by integrable is Integrable} \textbf{\emph{(Measurable
function dominated by integrable function is integrable).}} If $X$
is a real-valued measurable function such that $|X|\leq Y$ a.e. for
some non-negative integrable function $Y$, then $X$ is integrable.
In particular, if $A$ is a measurable set and $Z$ is an integrable
function, then $Z1_{A}$ is integrable.
\end{prop}
\begin{proof}
Let $(a_{n})_{n=1,2,\cdots}$ be an increasing sequence of positive
real numbers with $a_{n}\rightarrow\infty$. Let $(b_{n})_{n=1,2,\cdots}$
be a decreasing sequence of positive real numbers with $b_{n}\rightarrow0$.
Let $n\geq1$ be arbitrary. Then the function $X_{n}\equiv((-a_{n})\wedge X\vee a_{n})1_{(b(n)<Y)}$
is integrable. Moreover, $X=X_{n}$ on $(b_{n}<Y)(|X|\leq a_{n})$.
Hence
\[
|X-X_{n}|\leq|X-X_{n}|(1_{(Y\leq b(n))}+1_{(|X|>a(n))})
\]
\begin{equation}
\leq2Y(1_{(Y\leq b(n))}+1_{(|X|>a(n))})\leq2Y\wedge b_{n}+2Y1_{(Y>a(n))}\label{eq:temp-277}
\end{equation}
By Proposition \ref{Prop. Existence of modulus of integrability}
we have $I(Y1_{(Y>a(n))})\rightarrow0$. At the same time $I(Y\wedge b_{n})\rightarrow0$
since $Y$ is integrable. Therefore $I(2Y\wedge b_{n}+2Y1_{(Y>a(n)})\rightarrow0$.
Hence, by Theorem \ref{Thm. |X-Yn|<Zn where Yn,Zn integrable =000026 IZn->0 implies X integrable},
inequality \ref{eq:temp-277} implies that $X$ is integrable. 
\end{proof}

\section{Convergence of Measurable Functions}

In this section, let $(\Omega,L,I)$ be a complete integration space,
and let $(S,d)$ be a complete metric space, with a fixed reference
point $x_{\circ}\in S$. In the case where $S=R$, it is understood
that $d$ is the Euclidean metric and that $x_{\circ}=0$. We will
introduce several notions of convergence of measurable functions on
$(\Omega,L,I)$ with values in $(S,d)$.

First recall some notations and definitions. We will write $(a_{i})$
for short for a given sequence $(a_{i})_{i=1,2,\cdots}$ . If $(Y_{i})_{i=1,2,\cdots}$
is a sequence of functions from a set $\Omega'$ to the metric space
$S$, and if 
\[
D\equiv\{\omega\in\cup_{n=1}^{\infty}\cap_{i=n}^{\infty}domain(Y_{i}):\lim_{i\rightarrow\infty}Y_{i}(\omega)\:exists\;in\;S\}
\]
is non-empty, then the function $\lim_{i\rightarrow\infty}Y_{i}$
is defined by $domain(\lim_{i\rightarrow\infty}$$Y_{i})\equiv D$
and by $\lim_{i\rightarrow\infty}Y_{i}(\omega)$ for each $\omega\in D$.
We write $\mu A\equiv\mu(A)\equiv I1_{A}$ for the measure of an integrable
set $A$, and similarly write $\mu AB\equiv\mu(A;B)\equiv\mu(AB)$
for integrable sets $A$ and $B$. If $Z$ is a real-valued function
on $\Omega$ and if $a\in R$ then
\[
(Z\leq a)\equiv\{\omega\in domain(Z):Z(\omega)\leq a\}\subset domain(Z).
\]
Similarly when $\leq$ is replaced by $<,\geq,>,$ or $=$. 
\begin{defn}
\label{Def Convergence in meas, a.u.,L1,a.e.} \textbf{(Convergence
in measure, a.u., a.e., and in $L_{1}$).} For each $n\geq1$, let
$X,X_{n}$ be functions on the complete integration space \emph{$(\Omega,L,I)$,}
with values in the complete metric space $(S,d)$.
\end{defn}
\begin{enumerate}
\item The sequence $(X_{n})$ is said to converge to $X$ uniform\emph{ly
\index{convergence uniformly} on a subset} $A$ of $\Omega$ if,
for each $\varepsilon>0$, there exists $p\geq1$ so large that $A\subset\bigcap_{n=p}^{\infty}(d(X_{n},X)\leq\varepsilon)$. 
\item The sequence $(X_{n})$ is said to converge to $X$ \emph{almost }uniform\emph{ly}\index{convergence almost uniformly}
(\emph{a.u.}) if, for each integrable set $A$ and real number $\varepsilon>0$,
there exists an integrable set $B$ with $\mu(B)<\varepsilon$ such
that $X_{n}$ converges to $X$ uniformly on $AB^{c}$. 
\item The sequence $(X_{n})$ is said to \emph{converge to} $X$ \emph{in}
meas\emph{ure} \index{convergence in measure} if, for each integrable
set $A$ and each $\varepsilon>0$, there exists $p\geq1$ so large
that, for each $n\geq p$, there exists an integrable set $B_{n}$
with $\mu(B_{n})<\varepsilon$ and $AB_{n}^{c}\subset(d(X_{n},X)\leq\varepsilon)$. 
\item The sequence $(X_{n})$ is said to be \index{Cauchy in measure}\emph{Cauchy
in }meas\emph{ure} if for each integrable set $A$ and each $\varepsilon>0$,
there exists $p\geq1$ so large that for each $m,n\geq p$, there
exists an integrable set $B_{m,n}$ with $\mu(B_{m,n})<\varepsilon$
and $AB_{m,n}^{c}\subset(d(X_{n},X_{m})\leq\varepsilon)$. 
\item Suppose $S=R$ and $X,X_{n}\in L$ for $n\geq1$. The sequence $(X_{n})$
is said to \emph{converge to $X$ in $L_{1}$} \index{convergence in L_1@convergence in L\_1}
if $I|X_{n}-X|\rightarrow0$.
\item The sequence $(X_{n})$ is said to converge to $X$ on a subset $A$
if $A\subset domain(\lim_{n\rightarrow\infty}X_{n})$ and\emph{ }if
$X=\lim_{n\rightarrow\infty}X_{n}$ on $A$. The sequence $(X_{n})$
is said to converge to $X$ \index{convergence almost everywhere}\emph{almost
everywhere }(a.e.)\emph{ }if $(X_{n})$ converges to $X$ on some
full subset $DA$ of each integrable set $A$. $\square$
\end{enumerate}
We will use the abbreviation $X_{n}\rightarrow X$ to stand for ``$(X_{n})$
converges to $X$'', in whichever sense specified. 
\begin{prop}
\label{Prop. a.u. Convergence =00003D> Convergence in measure} \textbf{\emph{(a.u.
Convergence implies convergence in measure, etc).}} For each $n\geq1$,
let $X,X_{n}$ be functions on the complete integration space \emph{$(\Omega,L,I)$,}
with values in the complete metric space $(S,d)$. Then the following
holds.

1. If\emph{ $X_{n}\rightarrow X$} a.u. then \emph{(i) }$X$ is defined
a.e.\emph{ }on each integrable set\emph{ $A$,} \emph{(ii)} $X_{n}\rightarrow X$
in measure, and (iii) $X_{n}\rightarrow X$ a.e.

2. If \emph{(i)} $X_{n}$ is measurable for each $n\geq1$, \emph{and
(ii)} \emph{$X_{n}\rightarrow X$} in measure, then $X$ is measurable.

3. If \emph{(i)} $X_{n}$ is measurable for each $n\geq1$, and \emph{(ii)}
$X_{n}\rightarrow X$ a.u\emph{.}, then $X$ is measurable.
\end{prop}
\begin{proof}
1. Suppose $X_{n}\rightarrow X$ a.u. Let the integrable set $A$
and $n\geq1$ be arbitrary. Then, by Definition \ref{Def Convergence in meas, a.u.,L1,a.e.},
there exists an integrable set $B_{n}$ with $\mu(B_{n})<2^{-n}$
such that $X_{n}\rightarrow X$ uniformly on $AB_{n}^{c}.$ Hence
there exists $p\equiv p_{n}\geq1$ so large that 
\begin{equation}
AB_{n}^{c}\subset\bigcap_{h=p(n)}^{\infty}(d(X_{k},X)\leq2^{-h}).\label{eq:temp-183}
\end{equation}
In particular, $AB_{n}^{c}\subset domain(X)$. Define the integrable
set $B\equiv\bigcap_{k=1}^{\infty}\bigcup_{n=k}^{\infty}B_{n}$. Then
\[
\mu(B)\leq\sum_{n=k}^{\infty}\mu(B_{n})\leq\sum_{n=k}^{\infty}2^{-n}=2^{-k+1}
\]
for each $k\geq1$. Hence $B$ is a null set, and $D\equiv B^{c}$
is a full set. Moreover,
\[
AD\equiv AB^{c}=\bigcup_{k=1}^{\infty}\bigcap_{n=k}AB_{n}^{c}\subset\bigcup_{k=1}^{\infty}\bigcap_{n=k}domain(X)\subset domain(X).
\]
In other words, $X$ is defined a.e. on the integrable set $A.$ Part
(i) of Assertion 1 is proved. 

Now let $\varepsilon>0$ be arbitrary. Let $m\geq1$ be so large that
$2^{-p(m)}<\varepsilon$. Then, for each $n\geq p_{m},$ we have $\mu(B_{n})<2^{-n}<\varepsilon.$
Moreover,
\[
AB_{n}^{c}\subset(d(X_{n},X)\leq2^{-n})\subset(d(X_{n},X)<2^{-p(m)})\subset(d(X_{n},X)<\varepsilon).
\]
Thus the condition in Definition \ref{Def Convergence in meas, a.u.,L1,a.e.}
is verified for $X_{n}\rightarrow X$ in measure. Part (ii) of Assertion
1 is proved. Furthermore, since $X_{n}\rightarrow X$ uniformly on
$X_{n}\rightarrow X$ uniformly on $AB_{n}^{c},$ it follows that
$X_{n}\rightarrow X$ at each point in $AD=\bigcup_{k=1}^{\infty}\bigcap_{n=k}AB_{n}^{c}$,
where $D$ is a full set. Thus $X_{n}\rightarrow X$ a.e. on $A$,
where $A$ is an arbitrary integrable set. In other words, $X_{n}\rightarrow X$
a.e. Part (iii) of Assertion 1 is also proved. 

2. Suppose $X_{n}$ is measurable for each $n\geq1$, and suppose
\emph{$X_{n}\rightarrow X$} in measure\emph{. }We need to prove that
$X$ is measurable. To that end, let $f\in C_{ub}(S)$ be arbitrary.
Then $|f|\leq c$ on $S$ for some $c>0$, and $f$ has a modulus
of continuity $\delta_{f}$. Let the integrable set $A$ and $m\geq1$
be arbitrary. Take any $\alpha_{m}>0$ with $\alpha_{m}<2^{-m}\wedge\delta_{f}(2^{-m})$.
By hypothesis, $X_{n}\rightarrow X$ in measure. Hence, there exists
$p_{m}\geq1$ such that, for each $n\geq p_{m}$, there exists an
integrable set $B_{n}$ with $\mu(B_{n})<\alpha_{m}$ and with 
\[
AB_{n}^{c}\subset(d(X_{n},X)\leq\alpha_{m}).
\]
In particular, $AB_{n}^{c}\subset domain(X)$. Define the integrable
set $B\equiv\bigcap_{k=1}^{\infty}\bigcup_{n=k}^{\infty}B_{n}$. Then
$B$ is a null set, and $D\equiv B^{c}$ is a full set. Moreover,
\[
AD\equiv AB^{c}=\bigcup_{k=1}^{\infty}\bigcap_{n=k}AB_{n}^{c}\subset domain(X).
\]
In other words, $X$ is defined a.e. on the integrable set $A.$ Define
a function $Y:\Omega\rightarrow S$ by $domain(Y)\equiv AD\cup A^{c}D$,
and by $Y\equiv X$ and $Y\equiv x_{\circ}$ on $AD$ and $A^{c}D$
respectively. Then $f(Y)$ is defined a.e., and $f(Y)1_{A}=f(X)$
on $AD$. We will show that $Y$ is a measurable function.

To that end, let $m\geq1$ be arbitrary, and consider each $n\geq p_{m}$.
Then, since $\alpha_{m}<\delta_{f}(2^{-m})$, we have
\[
AB_{n}^{c}\subset A(d(X_{n},X)\leq\alpha_{m})\subset A(|f(X_{n})-f(X)|<2^{-m})=A(|f(X_{n})-f(Y)|<2^{-m}).
\]
Hence
\[
|f(Y)1_{A}-f(X_{n})1_{A}|\leq|f(Y)-f(X_{n})|(1_{AB(n)^{c}}+1_{B(n)})
\]
\begin{equation}
\leq2^{-m}1_{A}+2c1_{B(n)}.\label{eq:temp-283}
\end{equation}
Write $Z_{n}\equiv|f(Y)1_{A}-f(X_{n})1_{A}|.$ Then
\begin{equation}
I(Z_{n})\leq\varepsilon\mu(A)+2c\mu(B_{n})<2^{-m}\mu(A)+2c\alpha_{m}<2^{-m}\mu(A)+2c2^{-m},\label{eq:temp-284}
\end{equation}
Furthermore, $X_{n}$ is measurable for each $n\geq1$, by hypothesis.
Hence $f(X_{n})1_{A}\in L$ for each $n\geq1$. Therefore Theorem
\ref{Thm. |X-Yn|<Zn where Yn,Zn integrable =000026 IZn->0 implies X integrable}
implies that $f(X)1_{A}\in L$. Condition (i) in Definition \ref{Def. Measurable Function}
has been proved for $X.$ Hence, by Lemma \ref{Lem. f(X)1_A in L for all f in C(S) =00003D> (d(x0,X)>a)A integrable},
the set $(d(X,x_{\circ})>a)A$ is integrable for all but countably
many $a\in R$.

It remains to verify Condition (ii) in Definition \ref{Def. Measurable Function}
for $X$ to be measurable. For that purpose, consider each $m\geq1$
and $n\geq p_{m}$. Since $X_{n}$ is a measurable function, there
exists $a>0$ so large that 
\[
\mu(d(X_{n},x_{\circ})>a)A<2^{-m}.
\]
Take any $b>a+1>a+\alpha_{m}>a+2^{-m}$. Then
\[
(d(X,x_{\circ})>b)A\subset(d(X_{n},x_{\circ})>a)A\cup(d(X_{n},X)>\alpha_{m})A\subset(d(X_{n},x_{\circ})>a)A\cup AB_{n}.
\]
Hence
\[
\mu((d(X,x_{\circ})>b)A)\leq\mu((d(X_{n},x_{\circ})>a)A)+\mu(AB_{n})<2^{-m}+2^{-m}=2^{-m+1},
\]
where $2^{-m}>0$ is arbitrarily small. We conclude that $\mu(d(X,x_{\circ})>b)A\rightarrow0$
as $b\rightarrow\infty$. Thus Condition (ii) in Definition \ref{Def. Measurable Function}
is also verified. Accordingly, the function $X$ is measurable.

3. Suppose (i) $X_{n}$ is measurable for each $n\geq1$, and (ii)
$X_{n}\rightarrow X$ a.u\emph{.} Then, by Assertion 1, we have $X_{n}\rightarrow X$
in measure. Hence, by Assertion 2, the function $X$ is measurable.
Assertion 3 and the proposition is proved.
\end{proof}
\begin{prop}
\label{Prop.Sigma-finite:  Seq Cauchy in measure =00003D> a.u convergent subsequence}
\textbf{\emph{(In case of $\sigma$-finite $\Omega$, each sequence
Cauchy in measure converges in measure, and contains an a.u. convergent
subsequence) }}Suppose $(\Omega,L,I)$ is $\sigma$-finite\emph{.}
For each $n,m\geq1$, let $X_{n}$ be a function on\emph{ }$(\Omega,L,I)$,
with values in the complete metric space $(S,d)$, such that $d(X_{n},X_{m})$
is measurable.\emph{ }Suppose $(X_{n})_{n=1,2,\cdots}$ is Cauchy
in measure. Then there exists a subsequence $(X_{n(k)})_{k=1,2,\cdots}$
such that $X\equiv\lim_{k\rightarrow\infty}X_{n(k)}$ is a measurable
function, with $X_{n(k)}\rightarrow X$ a.u. and $X_{n(k)}\rightarrow X$
a.e. Moreover, $X_{n}\rightarrow X$ in measure.
\end{prop}
\begin{proof}
Let $(A_{k})_{k=1,2,\cdots}$ be a sequence of integrable sets that
is an $I$-basis of $(\Omega,L,I)$. Thus (i) $A_{k}\subset A_{k+1}$
for each $k\geq1$, and (ii) $\bigcup_{k=1}^{\infty}A_{k}$ is a full
set, and (iii) for any integrable set $A$ we have $\mu(A_{k}A)\rightarrow\mu(A)$. 

By hypothesis, $(X_{n})$ is Cauchy in measure. By Definition \ref{Def Convergence in meas, a.u.,L1,a.e.},
for each $k\geq1$ there exists $n_{k}\geq1$ such that, for each
$m,n\geq n_{k}$, there exists an integrable set $B_{m,n,k}$ with
\begin{equation}
\mu(B_{m,n,k})<2^{-k}\label{eq:temp-280}
\end{equation}
and with
\begin{equation}
A_{k}B_{m,n,k}^{c}\subset(d(X_{n},X_{m})\leq2^{-k}).\label{eq:temp-278}
\end{equation}
By inductively replacing $n_{k}$ with $n_{1}\vee\cdots\vee n_{k}$
we may assume that $n_{k+1}\geq n_{k}$ for each $k\geq1$. Define
\[
B_{k}\equiv B_{n(k+1),n(k),k}
\]
 for each $k\geq1$. Then $\mu(B_{k})<2^{-k}$ and
\begin{equation}
A_{k}B_{k}^{c}\subset(d(X_{n(k+1)},X_{n(k)})\leq2^{-k}).\label{eq:temp-128}
\end{equation}
For each $i\geq1$ let 
\begin{equation}
C_{i}\equiv\bigcup_{k=i}^{\infty}B_{k}.\label{eq:temp-282}
\end{equation}
Then
\begin{equation}
\mu(A_{i}^{c}C_{i})\leq\mu(C_{i})\leq\sum_{k=i}^{\infty}2^{-k}=2^{-i+1}\label{eq:temp-279}
\end{equation}
for each $i\geq1$, whence $\bigcap_{i=1}^{\infty}A_{i}^{c}C_{i}$
is a null set. Hence $D\equiv\bigcup_{i=1}^{\infty}A_{i}C_{i}^{c}$
is a full set. 

Let $i\geq1$ be arbitrary. Then 
\begin{equation}
A_{i}C_{i}^{c}\subset\bigcap_{k=i}^{\infty}A_{i}B_{k}^{c}\subset\bigcap_{k=i}^{\infty}A_{k}B_{k}^{c}\subset\bigcap_{k=i}^{\infty}(d(X_{n(k+1)},X_{n(k)})\leq2^{-k}),\label{eq:temp-266}
\end{equation}
in view of relation \ref{eq:temp-128}. Note that the second inclusion
is because $A_{j}\subset A_{k}$ for each $k\geq i$. Therefore, since
$(S,d)$ is complete, $(X_{n(k)})_{k=1,2,\cdots}$ converges uniformly
on $A_{i}C_{i}^{c}$. In other words, $X_{n(k)}\rightarrow X$ uniformly
on $A_{i}C_{i}^{c}$, where $X\equiv\lim_{k\rightarrow\infty}X_{n(k)}$. 

Next let $A$ be an arbitrary integrable set. Let $\varepsilon>0$
be arbitrary. In view of in view of Condition (iii) above, there exists
$i\geq1$ be so large that $2^{-i+1}<\varepsilon$ and $\mu(AA_{i}^{c})<\varepsilon$.
Such an $i$ exists . Let $B\equiv AA_{i}^{c}\cup C_{i}$. Then $\mu(B)<2\varepsilon$.
Moreover, $AB^{c}=AA_{i}C_{i}^{c}\subset A_{i}C_{i}^{c}$, whence
$X_{n(k)}\rightarrow X$ uniformly on $AB^{c}$. Since $\varepsilon>0$
is arbitrary, we conclude that $X_{n(k)}\rightarrow X$ a.u. It then
follows from Proposition \ref{Prop. a.u. Convergence =00003D> Convergence in measure}
that $X_{n(k)}\rightarrow X$ in measure, $X_{n(k)}\rightarrow X$
a.e., and $X$ is measurable. Now define, for each $m\geq n_{i}$,
\[
\bar{B}_{m}\equiv AA_{i}^{c}\cup B_{m,n(i),i}\cup C_{i}.
\]
Then, in view of expressions \ref{eq:temp-280} and \ref{eq:temp-279},
we have, for each $m\geq n_{i}$, 
\begin{equation}
\mu(\bar{B}_{m})=\mu(AA_{i}^{c}\cup B_{m,n(i),i}\cup C_{i})<\varepsilon+2^{-i}+2^{-i+1}<3\varepsilon.\label{eq:temp-280-1}
\end{equation}
Moreover, 
\[
A\bar{B}_{m}^{c}=A(A^{c}\cup A_{i})B_{m,n(i),i}^{c}C_{i}^{c}
\]
\[
=AA_{i}B_{m,n(i),i}^{c}C_{i}^{c}=(A_{i}B_{m,n(i),i}^{c})(AC_{i}^{c})
\]
\[
\subset(d(X_{m},X_{n(i)})\leq2^{-i})(d(X,X_{n(i)})\leq2^{-i+1})
\]
\[
\subset(d(X,X_{m})\leq2^{-i}+2^{-i+1})
\]
\begin{equation}
\subset(d(X,X_{m})<2\varepsilon)\label{eq:temp-279-1}
\end{equation}
for each $m\geq n(i)$, where the first inclusion is because of expressions
\ref{eq:temp-278} and \ref{eq:temp-279}. Since $\varepsilon>0$
is arbitrarily, we have verified the condition in Definition \ref{Def Convergence in meas, a.u.,L1,a.e.}
for $X_{m}\rightarrow X$ in measure.
\end{proof}
\begin{prop}
\emph{\label{Prop. Xn values in complete S, Then Xn->X in meas iff  Imin(1,d(Xn,X))1_A ->0 -1}
}\textbf{\emph{(Convergence in measure in terms of convergence of
integrals).}}\emph{ }For each $n\geq1$, let $X,X_{n}$ be functions
on\emph{ }$(\Omega,L,I)$, with values in the complete metric space
$(S,d)$, such that $d(X_{n},X_{m})$ and $d(X_{n},X)$ are measurable,
for each $n,m\geq1$.\emph{ }Then the following holds.

1. If $I(1\wedge d(X_{n},X))1_{A}\rightarrow0$ for each integrable
set $A$, then $X_{n}\rightarrow X$ in measure. 

2. Conversely, if $X_{n}\rightarrow X$ in measure, then $I(1\wedge d(X_{n},X))1_{A}\rightarrow0$
for each integrable set $A$.

3. The sequence $(X_{n})$ is Cauchy in measure iff $I(1\wedge d(X_{n},X_{m}))1_{A}\rightarrow0$
as $n,m\rightarrow\infty$ for each integrable set $A$. 
\end{prop}
\begin{proof}
Let the integrable set $A$ and the positive real number $\varepsilon\in(0,1)$
be arbitrary. 

1. Suppose $I(1\wedge d(X_{n},X))1_{A}\rightarrow0$. Let $\varepsilon>0$
be arbitrary. Then, by Chebychev's inequality,
\[
\mu(d(X_{n},X)>\varepsilon)A)\leq\mu(1\wedge d(X_{n},X))1_{A}>\varepsilon)\leq\varepsilon^{-1}I(1\wedge d(X_{n},X))1_{A}\rightarrow0
\]
as $n\rightarrow\infty$. In particular, there exists $p\geq1$ so
large that $\mu(1\wedge d(X_{n},X)>\varepsilon)A<\varepsilon$ for
each $n\geq p$. Now consider each $n\geq p$. Define the integrable
set $B_{n}\equiv(1\wedge d(X_{n},X)>\varepsilon)A$. Then $\mu(B_{n})<\varepsilon$
and $AB_{n}^{c}\subset(d(X_{n},X)\leq\varepsilon)$. Thus \emph{$X_{n}\rightarrow X$
}in measure.

2. Conversely, suppose \emph{$X_{n}\rightarrow X$ }in measure. Let
$\varepsilon>0$ be arbitrary. Then there exists $p\geq1$ so large
that, for each $n\geq p$, there exists an integrable set $B_{n}$
with $\mu(B_{n})<\varepsilon$ and $AB_{n}^{c}\subset(d(X_{n},X)\leq\varepsilon)$.
Hence
\[
I(1\wedge d(X_{n},X))1_{A}=I(1\wedge d(X_{n},X))1_{AB(n)}+I(1\wedge d(X_{n},X))1_{AB(n)^{c}}
\]
\[
\leq I1_{B(n)}+I\varepsilon1_{A}<\varepsilon+\varepsilon\mu(A),
\]
where $\varepsilon>0$. Thus $I(1\wedge d(X_{n},X))1_{A}\rightarrow0$.

3. The proof of Assertion 4 is similar to that of Assertions 1 and
2.
\end{proof}
The next Proposition will  be convenient for establishing a.u. convergence. 
\begin{prop}
\emph{\label{Prop. Sufficiency for a.u. convergence} }\textbf{\emph{(Sufficient
condition for a.u. convergence).}}\emph{ Suppose $(\Omega,L,I)$ is
$\sigma$-finite, with an $I$-basis $(A_{i})_{i=1,2,\cdots}.$} For
each $n,m\geq1$, let $X_{n}$ be a function on\emph{ }$(\Omega,L,I)$,
with values in the complete metric space $(S,d)$, such that $d(X_{n},X_{m})$
is measurable.\emph{ Suppose that, for each $i\geq1$, there exists
a sequence $(\varepsilon_{n})_{n=1,2,\cdots}$ of positive real numbers
such that $\sum_{n=1}^{\infty}\varepsilon_{n}<\infty$ and such that
$I(1\wedge d(X_{n},X_{n+1}))1_{A(i)}<\varepsilon_{n}^{2}$ for each
$n\geq1$. Then $X\equiv\lim_{n\rightarrow\infty}X_{n}$ exists on
a full set, and $X_{n}\rightarrow X$ }a.u. If, in addition, $X_{n}$
is measurable for each $n\geq1$, then the limit $X$ is measurable.
\end{prop}
\begin{proof}
For abbreviation write $Z_{n}\equiv1\wedge d(X_{n+1},X_{n})$ for
each $n\geq1$. Let $A$ be an arbitrary integrable set and let $\varepsilon>0$.
Let \emph{$i\geq1$} be so large that $\mu(AA_{i}^{c})<\varepsilon$.
By hypothesis, there exists a sequence $(\varepsilon_{n})_{n=1,2,\cdots}$
of positive real numbers such that $\sum_{n=1}^{\infty}\varepsilon_{n}<\infty$
and such that $IZ_{n}1_{A(i)}<\varepsilon_{n}^{2}$ for each $n\geq1$.
Chebychev's inequality then implies that 
\[
\mu(Z_{n}>\varepsilon_{n})A_{i}\leq I(\varepsilon_{n}^{-1}Z_{n}1_{A(i)})<\varepsilon_{n}
\]
for each $n\geq1$. Let $p\geq1$ be so large that $\sum_{n=p}^{\infty}\varepsilon_{n}<1\wedge\varepsilon$.
Let $C\equiv\bigcup_{n=p}^{\infty}(Z_{n}>\varepsilon_{n})A_{i}$ and
let $B\equiv AA_{i}^{c}\cup C$. Then $\mu(B)<\varepsilon$. Moreover,
\[
AB^{c}\subset AA_{i}C^{c}=AA_{i}\bigcap_{n=p}^{\infty}(Z_{n}\leq\varepsilon_{n})\subset\bigcap_{n=p}^{\infty}(Z_{n}\leq\varepsilon_{n})
\]
\[
\subset\bigcap_{n=p}^{\infty}(d(X_{n+1},X_{n})\leq\varepsilon_{n}).
\]
Since $\sum_{n=p}^{\infty}\varepsilon_{n}<\infty$, it follows that
$X_{n}\rightarrow X$ uniformly on $AB^{c}$, where $X\equiv\lim_{n\rightarrow\infty}X_{n}$.
Since $A$ and $\varepsilon>0$ are arbitrary, we see that $X_{n}\rightarrow X$
a.u. 

If, in addition, $X_{n}$ is measurable for each $n\geq1$, then $X$
is measurable by Proposition \ref{Prop. a.u. Convergence =00003D> Convergence in measure}.
\end{proof}
\begin{prop}
\label{Prop. If f continuous on R2 then f(X,Y) continuous wrt convergence in measure}
\textbf{\emph{(A continuous function preserves convergence in measure).}}
Let $(\Omega,L,I)$ be a complete integration space. Let $(S',d'),(S'',d'')$
be locally compact metric spaces and let $(\widetilde{S},\widetilde{d})\equiv(S',d')\otimes(S'',d'')$
denote the product metrics space. Let $X',X'_{1}X'_{2},\cdots$ be
a sequence of measurable functions with values in $S'$ such that
$X'_{n}\rightarrow X'$ in measure. Similarly, let $X'',X''_{1}X''_{2},\cdots$
be a sequence of measurable functions with values in $S''$ such that
$X''_{n}\rightarrow X''$ in measure.

Let $f:(\widetilde{S},\widetilde{d})\rightarrow S$ be a continuous
function with values in a complete metric space $(S,d)$ which is
\emph{(i)} uniformly continuous on bounded subsets, and \emph{(ii)}
bounded on bounded subsets. Then $f(X'_{n},X''_{n})\rightarrow f(X',X'')$
in measure as $n\rightarrow\infty$.

Generalization to $m\geq2$ sequences of measurable functions is similar.
\end{prop}
\begin{proof}
Let $x'_{\circ}$ and $x''_{\circ}$ be fixed reference points in
$S',S''$ respectively. Write $x_{\circ}\equiv(x'_{\circ},x''_{\circ})$.
For each $x,y\in\widetilde{S}$, write $x\equiv(x',x'')$ and $y\equiv(y',y'')$.
Likewise, write $X\equiv(X',X'')$ and $X_{n}\equiv(X'_{n},X''_{n})$
for each $n\geq1$. Note that, for each $n\geq1$, the functions $f(X),f(X_{n})$
are measurable functions with values in $S$, thanks to Assertion
2 of Proposition \ref{Prop. vector of measurable func (S)   is meas func (S)} 

Let $A$ be an arbitrary integrable set. Let $\varepsilon>0$ be arbitrary.
By condition (ii) in Definition \ref{Def. Measurable Function}, there
exists $a>0$ so large that $\mu(B')<\varepsilon$ and $\mu(B'')<\varepsilon$,
where $B'\equiv(d'(x'_{\circ},X')>a)A$ and $B''\equiv(d''(x''_{\circ},X'')>a)A$.
Since $(\widetilde{S},\widetilde{d})$ is locally compact, the bounded
subset $(d'(x'_{\circ},\cdot)\leq a)\times(d''(x''_{\circ},\cdot)\leq a)$
is contained in some compact subset. On the other hand, by hypothesis,
the function $f:\widetilde{S}\rightarrow S$ is uniformly continuous
on each compact subset of $\widetilde{S}$. Hence there exists $\delta_{1}>0$
be so small that, for each 
\[
x,y\in(d'(x'_{\circ},\cdot)\leq a)\times(d''(x''_{\circ},\cdot)\leq a)\subset\widetilde{S}
\]
with $\widetilde{d}(x,y)<\delta_{1}$, we have $d(f(x),f(y))<\varepsilon$.
Take any $\delta\in(0,\frac{1}{2}\delta_{1})$. For each $n\geq1$,
define $C'_{n}\equiv(d'(X'_{n},X')\geq\delta)A$ and $C''_{n}\equiv(d''(X_{n}'',X'')\geq\delta)A$.
By hypothesis, $X'_{n}\rightarrow X',$ and $X''_{n}\rightarrow X''$
in measure as $n\rightarrow\infty$. Hence there exists $p\geq1$
so large that $\mu(C'_{n})<\varepsilon$ and $\mu(C''_{n})<\varepsilon$
for each $n\geq p$. Consider any $n\geq p$. We have
\[
\mu(B'\cup B''\cup C'_{n}\cup C''_{n})<4\varepsilon.
\]
Moreover,
\[
A(B'\cup B''\cup C_{n}'\cup C_{n}'')^{c}=AB'^{c}B''{}^{c}C_{n}'^{c}C_{n}''^{c}
\]
\[
=A(d'(x'_{\circ},X')\leq a;d''(x''_{\circ},X'')\leq a;d'(X'_{n},X')<\delta;d''(X_{n}'',X'')<\delta)
\]
\[
\subset A((d'(x'_{\circ},X')\leq a)\times(d''(x''_{\circ},X'')\leq a))(\widetilde{d}(X_{n},X)<\delta_{1})
\]
\[
\subset(d(f(X),f(X{}_{n}))<\varepsilon).
\]
Since $\varepsilon>0$ and $A$ are arbitrary, the condition in Definition
\ref{Def Convergence in meas, a.u.,L1,a.e.} is verified for $f(X_{n})\rightarrow f(X)$
in measure. 
\end{proof}
\begin{thm}
\textbf{\emph{(Dominated Convergence Theorem)}}. \label{Thm. Dominated Convergence}\index{Dominated Convergence Theorem}
Let $(X_{n})_{n=1,2,\cdots}$ be a sequence of real-valued measurable
functions on the complete integration space $(\Omega,L,I)$, and let
$X$ be a real-valued function defined a.e. on $\Omega$, with $X_{n}\rightarrow X$
in measure. Suppose there exists an integrable function $Y$ such
that $|X|\leq Y$ a.e. and $|X_{n}|\leq Y$ a.e. for each $n\geq1$.
Then $X,X_{n}$ are integrable for each $n\geq1$, and $I|X_{n}-X|\rightarrow0$.
\end{thm}
\begin{proof}
By Proposition \ref{Prop. measurable dominated by integrable is Integrable},
$X_{n}$ is integrable for each $n\geq1$.

Let $\varepsilon>0$ be arbitrary. Since $Y$ is integrable and is
non-negative a.e., there exists $a>0$ so small that $I(Y\wedge a)<\varepsilon$.
Define $A\equiv(Y>a)$. Then $|X_{n}-X|\leq2Y=2(Y\wedge a)$ a.e.
on $A^{c}$ for each $n\geq1$. By Proposition \ref{Prop. Existence of modulus of integrability},
there exists $\delta\in(0,\varepsilon/(1+\mu A))$ so small that $IY1_{B}<\varepsilon$
for each integrable set $B$ with $\mu(B)<\delta$. On the other hand,
by hypothesis, $X_{n}\rightarrow X$ in measure. Hence there exists
$m>0$ so large that, for each $n\geq m$, we have $AB_{n}^{c}\subset(|X_{n}-X|\leq\delta)$
for some integrable set $B_{n}$ with $\mu(B_{n})<\delta$. Combining,
for each $n\geq m$, we have
\[
|X-X_{n}|\leq|X-X_{n}|1_{AB(n)}+|X-X_{n}|1_{AB(n)^{c}}+|X-X_{n}|1_{A^{c}}
\]
\begin{equation}
\leq2Y1_{B(n)}+\delta1_{A}+2(Y\wedge a)\quad a.e.,\label{eq:temp-285}
\end{equation}
where
\begin{equation}
I(2Y1_{B(n)}+\delta1_{A}+2(Y\wedge a))\leq2\varepsilon+\delta\mu(A)+2\varepsilon\leq2\varepsilon+\varepsilon+2\varepsilon.\label{eq:temp-286}
\end{equation}
Since $\varepsilon>0$ is arbitrary, inequalities \ref{eq:temp-285}
and \ref{eq:temp-286}, together with Theorem \ref{Thm. |X-Yn|<Zn where Yn,Zn integrable =000026 IZn->0 implies X integrable},
imply that $X$ is integrable and that $I|X_{n}-X|\rightarrow0$. 
\end{proof}
The next definition introduces Newton's notation for the Riemann-Stieljes
integration relative to a distribution function. 
\begin{defn}
\label{Def.  Newton's integral notation} \textbf{(Newton's notation).}
Suppose $F$ is a distribution function on $R$. Let $I$ be the Riemann-Stieljes
integration with respect to $F$, and let $(R,L,I)$ be the completion
of $(R,C(R),I)$. We will use the notation $\int\cdot dF$ for $I$.
For each $X\in L$, we write $\int XdF$ or $\int X(x)dF(x)$ for
$IX$. An integrable function in $L$ is then to be integrable relative
to $F$, and a measurable function on $(R,L,I)$ said to be measurable
relative to $F$. 

Suppose $X$ is a measurable function relative to $F$, and suppose
$s,t\in R$ such that the functions $1_{(s\wedge t,s]}X$ and $1_{(s\wedge t,t]}X$
are integrable relative to $F$. Then we write 
\[
\int_{s}^{t}XdF\equiv\int_{s}^{t}X(x)dF(x)\equiv\int X1_{(s\wedge t,t]}dF-\int X1_{(s\wedge t,s]}dF.
\]
Thus
\[
\int_{s}^{t}XdF=-\int_{t}^{s}XdF.
\]
If $A$ is a measurable set relative to $F$ such that $X1_{A}$ is
integrable, then we write
\[
\int_{A}XdF\equiv\int_{x\in A}X(x)dF(x)\equiv\int X1_{A}dF.
\]

In the special case where $F(x)\equiv x$ for $x\in R$, we write
$\int\cdot dx$ for $\int\cdot dF$. Let $s<t$ in $R$ be arbitrary.
The integration spaces $(R,L,\int\cdot dx)$ and $([s,t],L_{[s,t]},\int_{s}^{t}\cdot dx)$
are called the \index{Lebesgue integration space}\emph{Lebesgue integration
spaces }on $R$ and\emph{ $[s,t]$ }respectively\emph{,} and $\int\cdot dx$
and $\int_{s}^{t}\cdot dx$ are called the \emph{Lebesgue integration}.
Then an integrable function in $L$ or $L_{[s,t]}$ is said to be
Lebesgue integrable; and a measurable function is said to be  \emph{Lebesgue
}meas\emph{urable}\index{Lebesgue measurable function}. \emph{$\square$}
\end{defn}
Since the identity function $Z$, defined by $Z(x)\equiv x$ for each
$x\in R$, is continuous and is therefore a measurable function on
$(R,L,\int\cdot dF),$ all but countably many $t\in R$ are regular
points of $Z$. Hence $(s,t]=(s<Z\leq t)$ is a measurable set in
$(R,L,\int\cdot dF)$ for all but countably many $s,t\in R$. In other
words $1_{(s,t]}$ is measurable relative to $F$ for all but countably
many $s,t\in R$. Therefore the definition of $\int_{s}^{t}X(x)dF(x)$
is not vacuous.
\begin{prop}
\label{Prop. Intervals are Lebesgue integrable} \textbf{\emph{(Intervals
are Lebesgue integrable). }}Let $s,t\in R$ be arbitrary with $s\leq t$.
Then each of the intervals $[s,t]$, $(s,t)$, $(s,t]$, and $[s,t)$
is Lebesgue integrable, with Lebesgue measure equal to $t-s$, and
with measure-theoretic complements$(-\infty,s)\cup(t,\infty)$,$(-\infty,s]\cup[t,\infty)$,$(-\infty,s]\cup(t,\infty)$,$(-\infty,s)\cup[t,\infty)$,
respectively. Each of the intervals $(-\infty,s)$, $(-\infty,s]$,
$(s,\infty)$, and $[s,\infty)$ is Lebesgue measurable. 
\end{prop}
\begin{proof}
Consider the Lebesgue integration $\int\cdot dx$ and the Lebesgue
measure $\mu$.

Let $a,b\in R$ be such that $a<s\leq t<b$. Define $f\equiv f_{a,s,t,b}\in C(R)$
such that $f\equiv1$ on $[s,t]$, $f\equiv0$ on $(-\infty,a]\cup[b,\infty)$,
and $f$ is linear on $[a,s]$ and on $[t,b]$. Let $t_{0}<\cdots<t_{n}$
be any partition in the definition of a Riemann-Stieljes sum $S\equiv\sum_{i=1}^{n}f(t_{i})(t_{i}-t_{i-1})$
such that $a=t_{j}$ and $b=t_{k}$ for some $j,k=1,\cdots,n$ with
$j\leq k$. Then $S=\sum_{i=j+1}^{k}f(t_{i})(t_{i}-t_{i-1})$ since
$f$ has $[a,b]$ as support. Hence $0\leq S\leq\sum_{i=j+1}^{k}(t_{i}-t_{i-1})=t_{k}-t_{j}=b-a$.
Now let $n\rightarrow\infty,$ $t_{0}\rightarrow-\infty$, $t_{n}\rightarrow\infty$,
and let the mesh of the partition approach $0$. It follows from the
last inequality that $\int f(x)dx\leq b-a$. Similarly $t-s\leq\int f(x)dx$.
Now, with $s,t$ fixed, let $(a_{k})_{k=1,2,\cdots}$ and $(b_{k})_{k=1,2,\cdots}$
be sequences in $R$ such that $a_{k}\uparrow s$ and $b_{k}\downarrow t$,
and let $g_{k}\equiv f_{a_{k},s,t,b_{k}}$ Then, by the previous argument,
we have $t-s\leq\int g_{k}(x)dx\leq b_{k}-a_{k}\downarrow t-s$ Hence,
by the Monotone Convergence Theorem, the limit $g\equiv\lim_{k\rightarrow\infty}g_{k}$
is integrable, with integral $t-s$. It is obvious that $g=1$ or
$0$ on $domain(g).$ In other words, $g$ is an indicator function.
Moreover, $[s,t]=(g=1)$. Hence $[s,t]$ is an integrable set, with
$1_{[s,t]}=g$, with measure $\mu([s,t])=\int g(x)dx=t-s$, and with
measure-theoretic complement $[s,t]^{c}=(-\infty,s)\cup(t,\infty)$. 

Next consider the half open interval $(s,t]$. Since $(s,t]=\bigcup_{k=1}^{\infty}[s+\frac{1}{k},t]$
where $[s+\frac{1}{k},t]$ is integrable for each $k\geq1$, and where
$\mu([s+\frac{1}{k},t])=t-s-\frac{1}{k}\uparrow t-s$ as $k\rightarrow\infty$,
we have the integrability of $(s,t]$, and $\mu([s,t])=\lim_{k\rightarrow\infty}\mu([s+\frac{1}{k},t])=t-s$.
Moreover
\[
(s,t]^{c}=\bigcap_{k=1}^{\infty}[s+\frac{1}{k},t]^{c}
\]
\[
=\bigcap_{k=1}^{\infty}((-\infty,s+\frac{1}{k})\cup(t,\infty))=(-\infty,s]\cup(t,\infty).
\]
 The proofs for the intervals $(s,t)$ and $[s,t)$ are similar.

Now consider the interval $(-\infty,s)$. Define the function $X$
on the full set $D$ by $X(x)=1$ or 0 according as $x\in(-\infty,s)$
or $x\in[s,\infty)$. Let $A$ be any integrable subset of $R$. Then,
for each $n\geq1$ with $n>-s$, we have $|X1_{A}-1_{[-n,s)}1_{A}|\leq1_{(-\infty,-n)}1_{A}$
on the full set $D(A\cup A^{c})([-n,s)\cup[-n,s)^{c})$. At the same
time, $\int1_{(-\infty,-n)}(x)1_{A}(x)dx\rightarrow0$. Therefore,
by Theorem \ref{Thm. |X-Yn|<Zn where Yn,Zn integrable =000026 IZn->0 implies X integrable},
The function $X1_{A}$ is integrable. It follows that, for any $f\in C(R)$,
the function $f(X)1_{A}=f(1)X1_{A}+f(0)(1_{A}-X1_{A})$ is integrable.
We have thus verified condition (i) in Definition \ref{Def. Measurable Function}
for $X$ to be measurable. At the same time, since $|X|\leq1$, we
have trivially $\mu(|X|>a)=0$ for each $a>1$. Thus condition (ii)
in Definition \ref{Def. Measurable Function} is also verified. We
conclude that $(-\infty,s)$ is measurable. Similarly we can prove
that each of $(-\infty,s]$, $(s,\infty)$, and $[s,\infty)$ is measurable.
\end{proof}

\section{Product Integrals and Fubini's Theorem}

In the next definition and the following lemma, let $\Omega'$ and
$\Omega''$ be two sets, and let $L'$ and $L''$ be linear spaces
of real-valued functions on $\Omega'$ and $\Omega''$ respectively,
such that if $X',Y'$ are indicators in $L'$ then $X'Y'\in L'$,
and such that if $X'',Y''$ are indicators in $L''$ then $X''Y''\in L''$.
It follows that if $X',Y'$ are indicators in $L'$ then $Y'(1-X')=Y'-Y'X'\in L'$.
Similarly for $L''$. 
\begin{defn}
\label{Def. Direct-product-of functions}\textbf{(Direct product of
functions). }Let $X'\in L'$ and $X''\in L''$ be arbitrary. Define
a function $X'\otimes X'':\Omega'\times\Omega''\rightarrow R$ by
$domain(X'\otimes X'')\equiv domain(X')\times domain(X'')$ and by
$(X'\otimes X'')(\omega',\omega'')\equiv X'(\omega')X''(\omega'')$.
The function $X'\otimes X''$ is then called the direct product of
the functions $X'$ and $X''$. When the risk of confusion is low,
we will write $X'\otimes X''$ and $X'X''$ interchangeably. Generalization
to direct products $X'\otimes\cdots\otimes X^{(n)}$ with $n\geq2$
factors is similar. Further generalization to direct products $\otimes_{i=1}^{\infty}X^{(i)}$
of countably many factors is similar, provided that we restrict its
domain by
\[
domain(\otimes_{i=1}^{\infty}X^{(i)})
\]
\[
\equiv\{(\omega',\omega'',\cdots)\in\prod_{i=1}^{\infty}domain(X^{(i)}):\quad X^{(i)}(\omega^{(i)})\rightarrow1\;\mathrm{and}\;\prod_{i=1}^{\infty}X^{(i)}(\omega^{(i)})\;\mathrm{converges}\}.
\]
\end{defn}
$\square$
\begin{defn}
\label{Def. Simple functions} \textbf{(Simple functions).} A real
valued function $X$ on $\Omega'\times\Omega''$ is called a \emph{simple
function} \emph{relative to} $L',L''$ \index{simple function} if
$X=\sum_{i=0}^{n}\sum_{j=0}^{m}c_{i,j}X'{}_{i}X''_{j}$ where (i)
$n,m\geq1$, (ii) $X'_{1},\cdots,X'_{n}\in L'$ are mutually exclusive
indicators, (iii) $X''_{1},\cdots,X''_{m}\in L''$ are mutually exclusive
indicators, (iv) $X'_{0}=1-\sum_{i=1}^{n}X'_{i}$, (v) $X''_{0}=1-\sum_{j=1}^{m}X''_{j}$,
(vi) $c_{i,j}\in R$ with $c_{i,j}=0$ or $|c_{i,j}|>0$ for each
$i=0,\cdots,n$ and $j=0,\cdots,m$, and (vii) $c_{i,0}=0=c_{0,j}$
for each $i=0,\cdots,n$ and $j=0,\cdots,m$. Let $L_{0}$ denote
the set of simple functions on $\Omega'\times\Omega''$. Two simple
functions are said to be equal if they have the same domain and the
same values on the common domain. In other words, equality in $L_{0}$
is the set-theoretic equality. Note that, in the above notations,
$X'_{0},\cdots,X'_{n}$ are mutually exclusive indicators on $\Omega'$
that sum to $1$ on the intersection of their domains. Similarly $X''_{0},\cdots,X''_{m}$
are mutually exclusive indicators on $\Omega''$ that sum to $1$
on the intersection of their domains. The definition can be extended
in a straightforward manner to \emph{simple functions} relative to
linear spaces $L^{(1)},\cdots,L^{(k)}$ of functions on any $k\geq1$
sets $\Omega^{(1)},\cdots,\Omega^{(k)}$ respectively\emph{.} $\square$

If $1\in L'$ and $1\in L''$ then $X'_{0}\in L'$ and $X''_{0}\in L''$,
and the definition can be simplified. To be precise, if $1\in L'$
and $1\in L''$, then $X$ is a\emph{ }simple function\emph{ }relative\emph{
to} $L',L''$ iff $X=\sum_{i=0}^{n}\sum_{j=0}^{m}c_{i,j}X'{}_{i}X''_{j}$
where (i$'$) $n,m\geq1$, (ii$'$) $X'_{0},\cdots,X'_{n}\in L'$
are mutually exclusive indicators, (iii$'$) $X''_{0},\cdots,X''_{m}\in L''$
are mutually exclusive indicators, (iv$'$) $\sum_{i=0}^{n}X'_{i}=1$,
(v$'$) $\sum_{j=0}^{m}X''_{j}=1$, and (vi') $c_{i,j}\in R$ with
$c_{i,j}=0$ or $|c_{i,j}|>0$ for each $i=0,\cdots,n$ and $j=0,\cdots,m$.

In the notations of Definition \ref{Def. Simple functions}, if $Y'$
is an indicator with $Y'\in L'$, then we have $Y'X'_{0}=Y'(1-\sum_{i=1}^{n}X'_{i})\in L'$
by the remark preceding the definition, whence $Y'X'_{i}\in L'$ for
each $i=0,\cdots,n$. Similarly for $L''$.
\end{defn}
\begin{lem}
\label{Lem. Simple Fcts form linear space, closed to |.| and min(a,.)}
\textbf{\emph{(Simple functions constitute a linear space).}} Let
$L_{0}$ be the set of simple functions on $\Omega'\times\Omega''$
\emph{relative to} $L',L''$. Then the following holds.
\end{lem}
\begin{enumerate}
\item \emph{If $X\in L_{0}$ then $|X|,a\wedge X\in L_{0}$ for each $a>0$.
Specifically, let $X$ be a simple function with }$X=\sum_{i=0}^{n}\sum_{j=0}^{m}c_{i,j}X'{}_{i}X''_{j}$\emph{
as in Definition \ref{Def. Simple functions}. Then (i) $|X|=\sum_{i=0}^{n}\sum_{j=0}^{m}|c_{i,j}|X'{}_{i}X''_{j}\in L_{0}$
and (ii) $a\wedge X=\sum_{i=0}^{n}\sum_{j=0}^{m}(a\wedge c_{i,j})X'{}_{i}X''_{j}\in L_{0}$
for each $a>0$.}
\item \emph{$L_{0}$ is a linear space. }
\item \emph{Assertions 1 and 2 generalize to the case of simple functions
relative to linear spaces }$L^{(1)},\cdots,L^{(k)}$\emph{ of functions
on sets $\Omega^{(1)},\cdots,\Omega^{(k)}$ respectively}, \emph{for
any $k\geq1$ }.
\end{enumerate}
\begin{proof}

1. Consider any $(\omega',\omega'')\in domain(X)$. Either $X'_{i}(\omega')X''_{j}(\omega'')=0$
for each $i=1,\cdots,n$ and $j=1,\cdots,m$, or $X'_{k}(\omega')X''_{h}(\omega'')=1$
for exactly one pair of $k,h$ with $k=1,\cdots,n$ and $h=1,\cdots,m$.
In the first case, 
\[
|X(\omega',\omega'')|=0=\sum_{i=1}^{n}\sum_{j=1}^{m}|c_{i,j}|X'_{i}(\omega')X''_{j}(\omega'').
\]
In the second case,
\[
|X(\omega',\omega'')|=|c_{k,h}|=\sum_{i=1}^{n}\sum_{j=1}^{m}|c_{i,j}|X'_{i}(\omega')X''_{j}(\omega'').
\]
 $ $. Condition (i) is thus established. Condition (ii) is similarly
proved.

2. Obviously $L_{0}$ is closed under scalar multiplication. To show
that $L_{0}$ is closed also under addition, let
\[
X=\sum_{i=0}^{n}\sum_{j=0}^{m}c_{i,j}X'{}_{i}X''_{j}
\]
 and \emph{
\[
Y=\sum_{k=0}^{p}\sum_{h=0}^{q}b_{k,h}Y'{}_{k}Y''_{h}
\]
}be simple functions, as in Definition \ref{Def. Simple functions}.
Then
\[
\sum_{i=0}^{n}\sum_{j=0}^{m}X'{}_{i}X''_{j}=(\sum_{i=0}^{n}X'{}_{i})(\sum_{j=0}^{m}X''_{j})=1
\]
 on $domain(X)$ and similarly $\sum_{k=0}^{p}\sum_{h=0}^{q}Y'{}_{k}Y''_{h}=1$
on $domain(Y)$. Hence 
\[
X+Y=\sum_{i=0}^{n}\sum_{j=0}^{m}c_{i,j}X'{}_{i}X''_{j}+\sum_{k=0}^{p}\sum_{h=0}^{q}b_{k,h}Y'{}_{k}Y''_{h}
\]
\[
=\sum_{i=0}^{n}\sum_{j=0}^{m}c_{i,j}X'{}_{i}X''_{j}(\sum_{k=0}^{p}\sum_{h=0}^{q}Y'{}_{k}Y''_{h})
\]
\[
+\sum_{k=0}^{p}\sum_{h=0}^{q}b_{k,h}Y'{}_{k}Y''_{h}(\sum_{i=0}^{n}\sum_{j=0}^{m}X'{}_{i}X''_{j})
\]
\begin{equation}
=\sum_{i=0}^{n}\sum_{k=0}^{p}\sum_{j=0}^{m}\sum_{h=0}^{q}(c_{i,j}+b_{k,h})(X'{}_{i}Y'{}_{k})(X''_{j}Y''_{h}).\label{eq:temp-114}
\end{equation}
For each $i=0,\cdots,n$ and $k=0,\cdots,p$, the function $X'_{i}Y'_{k}$
is an indicator, and, by the remark preceding this lemma, belongs
to $L'$ if $(i,k)\neq(0,0)$. Moreover
\[
\sum_{i=0}^{n}\sum_{k=0}^{p}X'{}_{i}Y'_{k}=(\sum_{i=0}^{n}X'{}_{i})(\sum_{k=0}^{p}Y'_{k})=1.
\]
Suppose, $(i_{1},k_{1})\neq(i_{2},k_{2})$. Then either $i_{1}\neq i_{2}$
or $k_{1}\neq k_{2}$. In the first case, we have $X'_{i_{1}}X'{}_{i_{2}}=0$
and so $(X'_{i_{1}}Y'_{k_{1}})(X'{}_{i_{2}}Y'_{k_{2}})=0$. Similarly
$(X'_{i_{1}}Y'_{k_{1}})(X'{}_{i_{2}}Y'_{k_{2}})=0$ in the second
case. Summing up, we see that the sequence
\[
(X'_{0}Y'_{0},X'_{1}Y'_{0},X'_{0}Y'_{1},X'_{1}Y'_{1},\cdots,X'_{n}Y'_{p})
\]
satisfies conditions (ii-iv) in Definition \ref{Def. Simple functions}
Similarly, the sequence
\[
(X''_{0}Y''_{0},X''_{1}Y''_{0},X''_{0}Y''_{1},X''_{1}Y''_{1},\cdots,X''_{m}Y''_{q})
\]
satisfies conditions (v-vii) in Definition \ref{Def. Simple functions}.
Moreover, $c_{i,0}+b_{k,0}=0+0=0$ for each $i=0,\cdots,n$ and $k=0,\cdots,p$.
Similarly $c_{0,j}+b_{0,h}=0+0=0$ for each $j=0,\cdots,m$ and $h=0,\cdots,q$.
The right-hand side of equality \ref{eq:temp-114} is thus seen to
be a simple function. Therefore $X+Y$ is a simple function, proving
that $L_{0}$ is closed relative to addition and that it is a linear
space.

3. The proof for the general case is similar to the above proof for
$k=2$, and is omitted.
\end{proof}
In the remainder of this section, let $(\Omega',L',I')$ and $(\Omega'',L'',I'')$
be complete integration spaces. If $X',Y'\in L'$ are indicators,
then $X'Y'\in L'$. Similarly for $L''$. Let $\Omega\equiv\Omega'\times\Omega''$
and let $L_{0}$ denote the linear space of simple functions on $\Omega$.
If \emph{$X=\sum_{i=0}^{n}\sum_{j=0}^{m}c_{i,j}X'{}_{i}X''_{j}\in L_{0}$}
as in Definition \ref{Def. Simple functions}, define
\begin{equation}
I(X)=\sum_{i=0}^{n}\sum_{j=0}^{m}c_{i,j}I'(X'{}_{i})I''(X''_{j})\label{eq:temp-115}
\end{equation}

$\:$
\begin{lem}
\label{Lem. product integral of simple funcs is well def} \textbf{\emph{(Product
integral of simple functions is well defined). }}The function $I$
defined by equality \ref{eq:temp-115} is well defined, and is a linear
function on $L_{0}$. 
\end{lem}
\begin{proof}
It is obvious that $I(aX)=aI(X)$ for each $X\in L_{0}$ and $a\in R$.
Suppose $X,Y\in L_{0}$. Using the notations in equality\ref{eq:temp-114},
we have 
\[
I(X+Y)=\sum_{i=0}^{n}\sum_{k=0}^{p}\sum_{j=0}^{m}\sum_{h=0}^{q}(c_{i,j}+b_{k,h})I'(X'{}_{i}Y'{}_{k})I''(X''_{j}Y''_{h})
\]
\[
=\sum_{i=0}^{n}\sum_{j=0}^{m}c_{i,j}I'(X'{}_{i}\sum_{k=0}^{p}Y'{}_{k})I''(X''_{j}\sum_{h=0}^{q}Y''_{h})
\]
\[
+\sum_{k=0}^{p}\sum_{h=0}^{q}b_{k,h}I'(\sum_{i=0}^{n}X'{}_{i}Y'{}_{k})I''(\sum_{j=0}^{m}X''_{j}Y''_{h})
\]
\[
=\sum_{i=0}^{n}\sum_{j=0}^{m}c_{i,j}I'(X'{}_{i})I''(X''_{j})+\sum_{k=0}^{p}\sum_{h=0}^{q}b_{k,h}I'(Y'{}_{k})I''(Y''_{h})
\]
\[
=I(X)+I(Y).
\]
Thus $I$ is a linear operation. Next, suppose a simple function
\[
X=\sum_{i=0}^{n}\sum_{j=0}^{m}c_{i,j}X'{}_{i}X''_{j}
\]
is such that $X=0$. Then $c_{i,j}X'{}_{i}X''_{j}=0$ for each $i=0,\cdots,n$
and $j=0,\cdots,m$. It follows that the right-hand side of equality
\ref{eq:temp-115} vanishes, whence $I(X)=0$. Now suppose simple
functions $X$ and $Y$ are such that $X=Y$. Then $X-Y$ is a simple
function and $X-Y=0$. Hence $I(X)-I(Y)=I(X-Y)=0$, or $I(X)=I(Y)$.
Thus $I$ is a well-defined function.
\end{proof}
\begin{thm}
\label{Thm.  Product integration space} \textbf{\emph{(Integration
on space of simple functions).}} Let I be defined as in equality \ref{eq:temp-115}.
Then the triple $(\Omega,L_{0},I)$ is an integration space.
\end{thm}
\begin{proof}
We need to verify the three conditions in Definition \ref{Def. Integration Space}.

The linearity of $L_{0}$ has been proved in Lemma \ref{Lem. Simple Fcts form linear space, closed to |.| and min(a,.)}.
The linearity of the function $I$ has been proved in Lemma \ref{Lem. product integral of simple funcs is well def}. 

Next consider any $X\in L_{0}$, with $X=\sum_{i=0}^{n}\sum_{j=0}^{m}c_{i,j}X'{}_{i}X''_{j}$
in the notations of Definition \ref{Def. Simple functions}. By Lemma
\ref{Lem. Simple Fcts form linear space, closed to |.| and min(a,.)},
\[
|X|=\sum_{i=0}^{n}\sum_{j=0}^{m}|c_{i,j}|X'{}_{i}X''_{j}\in L_{0}
\]
and
\[
a\wedge X=\sum_{i=0}^{n}\sum_{j=0}^{m}(a\wedge c_{i,j})X'{}_{i}X''_{j}\in L_{0}
\]
for each $a>0$. Hence
\[
I(X\wedge a)=\sum_{i=0}^{n}\sum_{j=0}^{m}(a\wedge c_{i,j})I'(X'{}_{i})I''(X''_{j})
\]
 
\[
\rightarrow\sum_{i=0}^{n}\sum_{j=0}^{m}c_{i,j}I'(X'_{i})I''(X''_{j})\equiv I(X)
\]
as $a\rightarrow\infty$. Likewise
\[
I(|X|\wedge a)=\sum_{i=0}^{n}\sum_{j=0}^{m}(a\wedge|c_{i,j}|)I'(X'{}_{i})I''(X''_{j})\rightarrow0
\]
as $a\rightarrow0$. Conditions 1 and 3 in Definition \ref{Def. Integration Space}
are thus satisfied. 

It remains to prove condition 2 in Definition \ref{Def. Integration Space},
the positivity condition. To that end, suppose $(X_{k})_{k=0,1,2,\cdots}$
is a sequence of functions in $L_{0}$ such that $X_{k}\geq0$ for
$k\geq1$ and such that $\sum_{k=1}^{\infty}I(X_{k})<I(X_{0})$. For
$k\geq0$, we have $X_{k}=\sum_{i=0}^{n_{k}}\sum_{j=0}^{m_{k}}c_{k,i,j}X'{}_{k,i}X''_{k,j}$
as in Definition \emph{\ref{Def. Simple functions}. }It follows that
\[
\sum_{k=1}^{\infty}\sum_{i=0}^{n_{k}}\sum_{j=0}^{m_{k}}c_{k,i,j}I'(X'{}_{k,i})I''(X''_{k,j})=\sum_{k=1}^{\infty}I(X_{k})
\]
\[
<I(X_{0})=\sum_{i=0}^{n_{k}}\sum_{j=0}^{m_{k}}c_{0,i,j}I'(X'{}_{0,i})I''(X''_{0,j}).
\]
In view of the positivity condition on the integration $I'$, there
exists $\omega'\in\Omega'$ such that
\[
\sum_{k=1}^{\infty}\sum_{i=0}^{n_{k}}\sum_{j=0}^{m_{k}}c_{k,i,j}X'{}_{k,i}(\omega')I''(X''_{k,j})<\sum_{i=0}^{n_{k}}\sum_{j=0}^{m_{k}}c_{0,i,j}X'{}_{0,i}(\omega')I''(X''_{0,j}).
\]
In view of the positivity condition of the integration $I''$, the
last inequality in turn yields some $\omega''\in\Omega''$ such that
\[
\sum_{k=1}^{\infty}\sum_{i=0}^{n_{k}}\sum_{j=0}^{m_{k}}c_{k,i,j}X'{}_{k,i}(\omega')X''_{k,j}(\omega'')<\sum_{i=0}^{n_{k}}\sum_{j=0}^{m_{k}}c_{0,i,j}X'{}_{0,i}(\omega')X''_{0,j}(\omega'').
\]
Equivalently $\sum_{k=1}^{\infty}X_{k}(\omega',\omega'')<X_{0}(\omega',\omega'')$.
The positivity condition for $I$ has thus also been verified. We
conclude that $(\Omega,L_{0},I)$ is an integration space. 
\end{proof}
\begin{defn}
\label{Def.  Product integration space} \textbf{(Product of two integration
spaces).} The completion of the integration space $(\Omega,L_{0},I)$
is denoted by, 
\[
(\Omega,L,I)\equiv(\Omega'\times\Omega'',L'\otimes L'',I'\otimes I''),
\]
and is called the \emph{product integration space} \index{product integration space}
of $(\Omega',L',I')$ and $(\Omega'',L'',I'')$. The integration $I$
is called the \emph{product integration}. \index{product integration}$\square$
\end{defn}
\begin{prop}
\label{Prop. Cartesian product of integrable funcs is integrable}
\textbf{\emph{(Product of integrable functions is integrable relative
to the product integration).}} Let $(\Omega,L,I)$ denote the product
integration space of $(\Omega',L',I')$ and $(\Omega'',L'',I'')$. 

1. Suppose $X'\in L'$ and $X''\in L''$. Then $X'\otimes X''\in L$
and $I(X'\otimes X'')=(I'X')(I''X'')$. 

2. Moreover, if $D'$ and $D''$ are full subsets of $\Omega'$ and
$\Omega''$ respectively, then $D'\times D''$ is a full subset of
$\Omega$. 
\end{prop}
\begin{proof}
First suppose $X'=\sum_{i=1}^{n}a'_{i}X'_{i}$ and $X''=\sum_{i=1}^{m}a''_{i}X''_{i}$
where (i) $X'_{1},\cdots,X'_{n}$ are mutually exclusive integrable
indicator in $(\Omega',L',I')$, (ii) $X''_{1},\cdots,X''_{m}$ are
mutually exclusive integrable integrator in $(\Omega'',L'',I'')$,
and (iii) $a'_{1},\cdots,a'_{n},a''_{1},\cdots,a''_{m}\in R$. Define
$X'_{0}\equiv1-\sum_{i=1}^{n}X'_{i}$, define $X''_{0}\equiv1-\sum_{i=1}^{m}X''_{i}$,
and define $a'_{0}\equiv0\equiv a''_{0}$. Then $X'X''=\sum_{i=0}^{n}\sum_{j=0}^{m}a'_{i}a''_{j}X'_{i}X''_{j}\in L_{0}\subset L$.
Moreover
\[
I(X'X'')=\sum_{i=0}^{n}\sum_{j=0}^{m}a'_{i}a''_{j}I'(X'_{i})I''(X''_{j})
\]
\[
=(\sum_{i=0}^{n}a'_{i}I'X'_{i})(\sum_{j=0}^{m}a''_{j}I''X''_{j})=(I'X')(I''X'').
\]

Next let $D'$ be a full subset of $\Omega'$ and define the measurable
function $Y'$ on $\Omega'$ by $domain(Y')\equiv D'$ and $Y'\equiv1$
on $D'$. Thus $Y'$ is an indicator of $D'$, and $X'\equiv1-Y'=0Y'$
is an integrable indicator. Define $Y''\equiv1$ and $X''\equiv0$
on $\Omega''$. Then $X''$ is an integrable indicator on $\Omega''$.
Hence $X\equiv0(Y'Y''+Y'X''+X'Y'')+X'X''$ is a simple function, with
$IX=(I'X')(I''X'')=0$. On the other hand $domain(X)=domain(Y')\times domain(Y'')\equiv D'\times\Omega''$.
Therefore $D'\times\Omega''$ is a full subset of $\Omega$. Suppose,
$D''$ is a full subset of $\Omega''$. Then similarly $\Omega'\times D''$
is also a full subset of $\Omega$. It follows that $D'\times D''=(D'\times\Omega'')\cap(\Omega'\times D'')$
is a full subset of $\Omega$, proving the last assertion of the proposition.

Now consider arbitrary $X'\in L'$ and $X''\in L''$. We need to show
that $X'X''\in L$. By linearity, there is no loss of generality in
assuming that $X'\geq0$ and $X''\geq0$. By Proposition \ref{Prop. Integrable func is limit in L1 of combo of indicators},
there exist sequences $(X'_{k})_{k=1,2,\cdots}$ and $(X''_{k})_{k=1,2,\cdots}$
where (i) for each $n\geq1$, the functions $X'_{k}$ and $X''_{k}$
are linear combinations of mutually exclusive integrable indicators
in $(\Omega',L',I')$ and $(\Omega'',L'',I'')$ respectively, (ii)
$0\leq X'_{k}\uparrow X'$ and $0\leq X''_{k}\uparrow X''$ on $D'\equiv\bigcap_{k=1}^{\infty}domain(X'_{k})$
and $D''\equiv\bigcap_{k=1}^{\infty}domain(X''_{k})$ respectively,
and (iii) $I'X'_{k}\uparrow IX'$ and $I'X''_{k}\uparrow IX''$. Let
$k\geq1$ be arbitrary. By the first paragraph of this proof, we have
$I(X'_{k}X''_{k})=I'(X'_{k})I''(X''_{k})\uparrow I'(X')I''(X'')$.
Therefore, by the Monotone Convergence Theorem, $X'_{k}X''_{k}\uparrow X$
a.e. relative to $I$ for some $X\in L$, and $IX=(I'X')(I''X'')$.
On the other hand $X'_{k}X''_{k}\uparrow X'X''$ on the set $D'\times D''$,
which is a full set as seen in the previous paragraph. Thus $X'X''=X$
a.e. Hence $X'X''\in L$ and $I(X'X'')=IX=(I'X')(I''X'')$. 
\end{proof}
Next is Fubini's Theorem which enables the calculation of the product
integral as iterated integrals.
\begin{thm}
\textbf{\emph{(Fubini's Theorem for product of two integration spaces)}}\label{Thm. Fubini for product of two integraton spacces}
\index{Fubini's Theorem} Let $(\Omega,L,I)\equiv(\Omega,L'\otimes L'',I'\otimes I'')$
be the product integration space of $(\Omega',L',I')$ and $(\Omega'',L'',I'')$.
Let $X\in L'\otimes L''$ be arbitrary. 

Then there exists a full subset $D'$ of $\Omega'$ such that (i$'$)
for each $\omega'\in D'$, the function $X(\omega',\cdot)$ is a member
of $L''$, (ii$'$) the function $I''X$ defined by $domain(I''X)\equiv D'$
and $(I''X)(\omega')\equiv I''(X(\omega',\cdot))$ for each $\omega'\in D'$
is a member of $L'$, and (iii$'$) $IX=I'(I''X)$. 

Similarly, there exists a full subset $D''$ of $\Omega''$ such that
(i$''$) for each $\omega''\in D''$, the function $X(\cdot,\omega'')$
is a member of $L'$, (ii$''$) the function $I'X$ defined by $domain(I'X)\equiv D''$
and $(I'X)(\omega'')\equiv I'(X(\cdot,\omega''))$ for each $\omega''\in D''$
is a member of $L''$, and (iii$'$) $IX=I''(I'X)$. 
\end{thm}
\begin{proof}
First consider a simple function $X=\sum_{i=0}^{n}\sum_{j=0}^{m}c_{i,j}X'{}_{i}X''_{j}$,
in the notations of Definition \ref{Def. Simple functions}. Define
$D'\equiv\bigcap_{i=0}^{n}domain(X'_{i})$. Then $D'$ is a full subset
of $\Omega'$. Let $\omega'\in D'$ be arbitrary. Then $X(\omega',\cdot)=\sum_{i=0}^{n}\sum_{j=0}^{m}c_{i,j}X'{}_{i}(\omega')X''_{j}\in L''$,
verifying condition (i'). Define the function $I''X$ as in condition
(ii'). Then
\[
(I''X)(\omega')\equiv I''(X(\omega',\cdot))=\sum_{i=0}^{n}\sum_{j=0}^{m}c_{i,j}X'{}_{i}(\omega')I''X''_{j}
\]
for each $\omega'\in D'$. Thus $I''X=\sum_{i=0}^{n}\sum_{j=0}^{m}(c_{i,j}I''X''_{j})X'{}_{i}\in L'$,
which verifies condition (ii'). It follows from the last equality
that
\[
I'(I''X)=\sum_{i=0}^{n}\sum_{j=0}^{m}(c_{i,j}I''X''_{j})(I'X'{}_{i})=IX,
\]
 proving also condition (iii'). Thus conditions (i'-iii') are proved
in the case of a simple function $X$. 

Next let $X\in L'\otimes L''$ be arbitrary. Then there exists a sequence
$(X_{k})_{k=1,2,\cdots}$ of simple functions which is a representation
of $X$ relative to the integration $I$. Let $k\geq1$ be arbitrary.
We have
\[
X_{k}=\sum_{i=0}^{n_{k}}\sum_{j=0}^{m_{k}}c_{k,i,j}X'{}_{k,i}X''_{k,j},
\]
in the notations of Definition \ref{Def. Simple functions}. By Lemma
\ref{Lem. Simple Fcts form linear space, closed to |.| and min(a,.)},
we have
\[
|X_{k}|=\sum_{i=0}^{n_{k}}\sum_{j=0}^{m_{k}}|c_{k,i,j}|X'{}_{k,i}X''_{k,j}.
\]

For each $k\geq1$ we have $IX_{k}=I'(I''X_{k})$ and $I|X_{k}|=I'(I''|X_{k}|)$
by the first part of this proof. Therefore
\[
\sum_{k=1}^{\infty}I'|I''X_{k}|\leq\sum_{k=1}^{\infty}I'(I''|X_{k}|)=\sum_{k=1}^{\infty}I|X_{k}|<\infty.
\]
 Hence the functions $Y\equiv\sum_{k=1}^{\infty}I''X_{k}$ and $Z\equiv\sum_{k=1}^{\infty}I''|X_{k}|$
are in $L'$, with 
\[
I'Y=\sum_{k=1}^{\infty}I'(I''X_{k})=\sum_{k=1}^{\infty}IX_{k}=IX.
\]
Consider any $\omega'\in D'\equiv domain(Z)$. Then $\sum_{k=1}^{\infty}I''|X_{k}(\omega',\cdot)|<\infty$.
Moreover, if $\omega''\in\Omega''$ is such that $\sum_{k=1}^{\infty}|X_{k}(\omega',\cdot)|(\omega'')<\infty$,
then $\sum_{k=1}^{\infty}|X_{k}(\omega',\omega'')|<\infty$ and
\[
X(\omega',\cdot)(\omega'')\equiv X(\omega',\omega'')=\sum_{k=1}^{\infty}X_{k}(\omega',\omega'')\equiv\sum_{k=1}^{\infty}X_{k}(\omega',\cdot)(\omega'')
\]
In other words, for each $\omega'\in D'$, the sequence $(X_{k}(\omega',\cdot))_{k=1,2,\cdots}$
is a representation of $X(\omega',\cdot)$ in $L''$, and so $X(\omega',\cdot)\in L''$
with 
\[
(I''X)(\omega')\equiv I''X(\omega',\cdot)=\sum_{k=1}^{\infty}I''X_{k}(\omega',\cdot)=Y(\omega')
\]
Thus we see that $I''X=Y$ on the full set $D'$. Since $Y\in L'$,
so also $I''X\in L'$. Moreover, $I'(I''X)=I'(Y')=IX$. Conditions
(i'-iii') have thus been verified for an arbitrary $X\in L$. 

Conditions (i''-iii''), where the roles of $I'$ and $I''$ are reversed,
is proved similarly. 
\end{proof}
Following is the straightforward generalization of Fubini's theorem
to product integration to many factors.
\begin{defn}
\label{Def.  Product of many iintegration space} \textbf{(Product
of several integration spaces).} Let $n\geq1$ be arbitrary. Let $(\Omega^{(1)},L^{(1)},I^{(1)}),\cdots,(\Omega^{(n)},L^{(n)},I^{(n)})$
be complete integration spaces. If $n=1$, let
\[
(\prod_{i=1}^{n}\Omega^{(i)},\bigotimes_{i=1}^{n}L^{(i)},\bigotimes_{i=1}^{n}I^{(i)})\equiv(\Omega^{(1)},L^{(1)},I^{(1)}).
\]
Inductively for $n\geq2,$ define
\[
(\prod_{i=1}^{n}\Omega^{(i)},\bigotimes_{i=1}^{n}L^{(i)},\bigotimes_{i=1}^{n}I^{(i)})\equiv(\prod_{i=1}^{n-1}\Omega^{(i)},\bigotimes_{i=1}^{n-1}L^{(i)},\bigotimes_{i=1}^{n-1}I^{(i)})\bigotimes(\Omega^{(n)},L^{(n)},I^{(n)})
\]
where the product of the two integration spaces $(\prod_{i=1}^{n-1}\Omega^{(i)},\bigotimes_{i=1}^{n-1}L^{(i)},\bigotimes_{i=1}^{n-1}I^{(i)})$
and $(\Omega^{(n)},L^{(n)},I^{(n)})$ on the right-hand side is as
in Definition \ref{Def.  Product integration space}. Then, 
\[
(\prod_{i=1}^{n}\Omega^{(i)},\bigotimes_{i=1}^{n}L^{(i)},\bigotimes_{i=1}^{n}I^{(i)})
\]
is called the \index{product integration space}\emph{product integration
space} of the given integration spaces. 

In the special case where $(\Omega^{(1)},L^{(1)},I^{(1)})=\cdots=(\Omega^{(n)},L^{(n)},I^{(n)})$
are all equal to the same integration space $(\Omega_{0},L_{0},I_{0})$
we write
\[
(\Omega_{0}^{n},L_{0}^{\otimes n},I_{0}^{\otimes n})\equiv(\prod_{i=1}^{n}\Omega^{(i)},\bigotimes_{i=1}^{n}L^{(i)},\bigotimes_{i=1}^{n}I^{(i)})
\]
and call it the $n$-th \emph{power of the integration space}\index{power integration space}
$(\Omega_{0},L_{0},I_{0})$. 
\end{defn}
$\square$
\begin{thm}
\textbf{\emph{\label{Thm. Fubini's theorem for product of many integration spaces.}(Fubini's
Theorem for product of several integration spaces).}}\index{Fubini's Theorem}\textbf{\emph{
}}Let $n\geq1$ be arbitrary. Let $(\Omega^{(1)},L^{(1)},I^{(1)}),\cdots,(\Omega^{(n)},L^{(n)},I^{(n)})$
be complete integration spaces. Let $(\Omega,L,I)\equiv(\prod_{i=1}^{n}\Omega^{(i)},\bigotimes_{i=1}^{n}L^{(i)},\bigotimes_{i=1}^{n}I^{(i)})$
be their product space. Let $X\in\bigotimes_{i=1}^{n}L^{(i)}$ be
arbitrary. Let $k\in\{1,\cdots n\}$ be arbitrary. Then there exists
a full subset $D^{(k)}$ of $\prod_{i=1;i\neq k}^{n}\Omega^{(i)}$
such that \emph{(i)} for each $\widehat{\omega}_{k}\equiv(\omega_{1},\cdots,\omega_{k-1},\omega_{k+1},\cdots,\omega_{n})\in D^{(k)}$,
the function $X_{\widehat{\omega}(k)}:\Omega^{(k)}\rightarrow R$,
defined by
\[
X_{\widehat{\omega}(k)}(\omega_{k})\equiv X(\omega_{1},\cdots,\omega_{k-1},\omega_{k},\omega_{k+1},\cdots,\omega_{n}),
\]
is a member of $L^{(k)}$, \emph{(ii)} the function $\widehat{X}_{k}:\prod_{i=1;i\neq k}^{n}\Omega^{(i)}\rightarrow R$,
defined by $domain(\widehat{X}_{k})\equiv D^{(k)}$ and by $\widehat{X}_{k}(\widehat{\omega}_{k})\equiv I^{(k)}(X_{\widehat{\omega}(k)})$
for each $\widehat{\omega}_{k}\in D^{(k)},$ is a member of $\bigotimes_{i=1;i\neq k}^{n}L^{(i)}$,
and \emph{(iii) }$(\bigotimes_{i=1}^{n}I^{(i)})X=(\bigotimes_{i=1;i\neq k}^{n}I^{(i)})\widehat{X}_{k}$. 

A special example of is where $Y_{k}\in L^{(k)}$ is given for each
$k=1,\cdots,n$, and where we define the function $X:\Omega\rightarrow R$
by 
\[
X(\omega_{1},\cdots,\omega_{n})\equiv Y_{1}(\omega_{1})\cdots Y_{n}(\omega_{n})
\]
for each $(\omega_{1},\cdots,\omega_{n})\in\Omega$ such that $\omega_{k}\in domain(Y_{k})$
for each $k=1,\cdots,n$. Then $X\in\bigotimes_{i=1}^{n}L^{(i)}$
and 
\[
IX=I_{1}(Y_{1})\cdots I_{n}(Y_{n}).
\]
\end{thm}
\begin{proof}
The case where $n=1$ is trivial. The case where $n=2$ is proved
in Theorem \ref{Thm. Fubini for product of two integraton spacces}
and Proposition \ref{Thm. measurable function on one factor of product integration space can be regarded as measurable on product}.
The proof of the general case is by induction on $n$ , and is straightforward
and omitted.
\end{proof}
\begin{thm}
\label{Thm. measurable function on one factor of product integration space can be regarded as measurable on product}
\textbf{\emph{(A measurable function on one factor of a product integration
space can be regarded as }}$\mathbf{measurable}$\textbf{\emph{ on
the product).}} 

1. Let $(\Omega,L,I)$ be the product integration space of $(\Omega',L',I')$
and $(\Omega'',L'',I'')$. Let $X'$ be an arbitrary measurable function
on $(\Omega',L',I')$ with values in some complete metric space $(S,d)$.
Define $X:\Omega\rightarrow S$ by $X(\omega)\equiv X'(\omega')$
for each $(\omega',\omega'')\in\Omega$ such that $\omega'\in domain(X')$.
Then $X$ is measurable on $(\Omega L,I)$ with values in $(S,d)$.
Moreover $If(X)1_{A'\times\Omega''}=I'f(X')1_{A'}$ for each $f\in C_{ub}(S)$
and each integrable subset $A'$ of $\Omega'$. 

2. Similarly, with an arbitrary measurable function $X''$ on $(\Omega'',L'',I'')$
and with $X:\Omega\rightarrow S$ by $X(\omega)\equiv X''(\omega'')$
for each $(\omega',\omega'')\in\Omega$ such that $\omega''\in domain(X'')$,
the function $X$ is measurable on $(\Omega L,I)$ with values in
$(S,d)$. Moreover $If(X)1_{\Omega''}=I'f(X')1_{A'}$ for each $f\in C_{ub}(S)$
and each integrable subset $A'$ of $\Omega'$. 

3. More generally, let $n\geq1$ be arbitrary. Let $(\Omega^{(1)},L^{(1)},I^{(1)}),\cdots,(\Omega^{(n)},L^{(n)},I^{(n)})$
be complete integration spaces. Let $(\Omega,L,I)\equiv(\prod_{i=1}^{n}\Omega^{(i)},\bigotimes_{i=1}^{n}L^{(i)},\bigotimes_{i=1}^{n}I^{(i)})$
be their product space. Let $i=1,\cdots,n$ be arbitrary, and suppose
$X^{(i)}$ is a measurable function on $(\Omega^{(i)},L^{(i)},I^{(i)})$
with values in some complete metric space $(S,d)$. Define the function
$X:\Omega\rightarrow S$ by $X(\omega)\equiv X^{(i)}(\omega_{i})$
for each $\omega\equiv(\omega_{1},\cdots,\omega_{n})\in\Omega$ such
that $\omega_{i}\in domain(X^{(i)})$. Then $X$ is a measurable function
on $(\Omega L,I)$ with values in $(S,d)$. Moreover $If(X)1_{A}=I^{(i)}f(X^{(i)})1_{A(i)}$
for each $f\in C_{ub}(S)$ and each integrable subsets $A_{i}$ of
$\Omega^{(i)}$, where $A\equiv\prod_{k=1}^{n}A_{k},$ where $A_{k}\equiv\Omega^{(k)}$
for each $k=1,\cdots,n$ with $k\neq i$. 

4. Suppose, in addition, $\Omega^{(k)}$ is an integrable set with
$I^{(k)}\Omega^{(k)}=1$, for each $k=1,\cdots,n$. Then $X$ is a
measurable function with values in $(S,d)$, such that  $If(X)=I^{(i)}f(X^{(i)})$
for each $f\in C_{ub}(S)$. Anticipating a definition later, we say
that the measurable function $X$ has the same distribution as $X^{(i)}$. 
\end{thm}
\begin{proof}
Let $x_{\circ}\in S$ be an arbitrary, but fixed, reference point.
For each $n\geq0$, define $h_{n}\equiv1\wedge(n+1-d(x_{\circ},\cdot))_{+}\in C_{ub}(S)$. 

Let $f\in C_{ub}(S)$ and $g\in L$ be arbitrary. Then $|f|\leq b$
for some $b>0$. First assume that $g\equiv1_{A'\times A''}$ where
$A',A''$ are integrable subsets of $\Omega',\Omega''$ respectively.
Then $f(X)g=(f(X')1_{A'})1_{A''}$ is integrable according to \ref{Prop. Cartesian product of integrable funcs is integrable}.
Moreover,
\[
Ih_{n}(X)1_{A'\times A''}=(I'h_{n}(X')1_{A'})(I''1_{A''})\uparrow(I'1_{A'})(I''1_{A''})=I_{A}
\]
as $n\rightarrow\infty$. Hence, by linearity, if $g$ is a simple
function $\Omega$ relative to $L',L''$, then we have (i) $f(X)g\in L$,
and (ii) $Ih_{n}(X)g\rightarrow Ig$. Now let $(g_{k})_{k=1,2,\cdots}$
be a sequence of simple functions which is a representation of $g\in L$.
Then 
\begin{equation}
\sum_{k=1}^{\infty}I|f(X)g_{k}|\leq b\sum_{k=1}^{\infty}I|g_{k}|<\infty\label{eq:temp-238}
\end{equation}
where $b>0$ is any bound for $f\in C_{ub}(S)$. Hence 
\[
f(X)g=\sum_{k=1}^{\infty}f(X)g_{k}\in L.
\]
Similarly, $h_{n}(X)g\in L$ and
\[
Ih_{n}(X)g=\sum_{k=1}^{\infty}Ih_{n}(X)g_{k}
\]
for each $n\geq0$. Now $I|h_{n}(X)g_{k}|\leq I|g_{k}|$, and, by
Condition (ii) above, $Ih_{n}(X)g_{k}\rightarrow Ig_{k}$ as $n\rightarrow\infty$,
for each $k\geq1$. Hence $Ih_{n}(X)g\rightarrow\sum_{k=1}^{\infty}Ig_{k}=Ig$
as $n\rightarrow\infty$. In particular, if $A$ is an arbitrary integrable
subset of $\Omega$, then $Ih_{n}(X)1_{A}\uparrow I1_{A}\equiv\mu(A)$,
where $\mu$ is the measure relative to $I$. We have verified the
conditions in Proposition \ref{Prop. Alternative definition of measurablility}
for $X$ to be $\mathrm{measurable}$. Assertion 1 is proved. Assertion
2 is proved similarly. Assertion 3 follows from Assertion1 1 and 2,
by induction. Assertion 4 is a special case of Assertion 3, where
$A_{k}\equiv\Omega^{(k)}$ for each $k=1,\cdots,n$ . 
\end{proof}
For products of integrations based on locally compact spaces, the
following proposition will be convenient.
\begin{prop}
\label{Prop. Completion of C(S1xS2) =00003D compC(S1) x CompC(S2)}
\textbf{\emph{(Product of integration spaces based on locally compact
metric spaces). }}For each $i=1,\cdots,n$, let $(S_{i},d_{i})$ be
a locally compact metric space, and let $(S_{i},C(S_{i}),I^{(i)})$
be an integration space, with completion $(S_{i},L^{(i)},I^{(i)})$.
Let $(S,d)$ be their product metric space, and let
\[
(S,L,I)\equiv(\prod_{i=1}^{n}S_{i},\bigotimes_{i=1}^{n}L^{(i)},\bigotimes_{i=1}^{n}I^{(i)})
\]
be the product integration space. Then $C(S)\subset L$, and $(S,C(S),I)$
is an integration space with $(S,L,I)$ as completion.
\end{prop}
\begin{proof}
Consider only the case $n=2$, the general case being similar. For
arbitrary real valued functions $V_{1},V_{2}$ on $S_{1},S_{2}$ respectively,
we will abuse notations and write $V_{1}V_{2}$ for the function whose
value at $x\equiv(x_{1},x_{2})$ is $V_{1}(x_{1})V_{2}(x_{2})$ for
each $x\equiv(x_{1},x_{2})\in S$. By Definition \ref{Def.  Product integration space},
the product integration space $(S,L,I)$ is the completion of the
subspace $(S,L_{0},I)$ of simple functions.

Let $X\in C(S)$ be arbitrary. Since $X$ has compact support, there
exists $V_{i}\in C(S_{i})$ for $i=1,2$ such that (i) $0\leq V_{i}\leq1$
for $i=1,2,$ (ii) if $x\equiv(x_{1},x_{2})\in S$ is such that $|X(x)|>0$
then $V_{1}(x_{1})=1=V_{2}(x_{2})$, and (iii) $I^{(1)}V_{1}>0$ and
$I^{(2)}V_{2}>0$. Let $\varepsilon>0$ be arbitrary. By Proposition
\ref{Prop. Approx  by  interpolation}, there exist $U_{i,1},\cdots,U_{i,m}\in C(S_{i})$
for $i=1,2$ such that
\[
|X-\sum_{k=1}^{m}U_{1,k}U_{2,k}|<\varepsilon.
\]
Multiplication by $V_{1}V_{2}$ yields, in view of Condition (ii),
\begin{equation}
|X-\sum_{k=1}^{m}(V_{1}U_{1,k})(V_{2}U_{2,k})|<\varepsilon V_{1}V_{2}.\label{eq:temp-2}
\end{equation}
Since $C(S_{i})\subset L^{(i)}$ for each $i=1,2$, we have $V_{1}V_{2}\in L$
and $(V_{1}U_{1,k})(V_{2}U_{2,k})\in L$ for each $k=1,\cdots,m$,
according to Proposition \ref{Prop. Cartesian product of integrable funcs is integrable}.
Since $I(\varepsilon V_{1}V_{2})>0$ is arbitrarily small, inequality
\ref{eq:temp-2} implies that $X\in L$, thanks to Theorem \ref{Thm. |X-Yn|<Zn where Yn,Zn integrable =000026 IZn->0 implies X integrable}.
Since $X\in C(S)$ is arbitrary, we conclude that $C(S)\subset L$. 

Since $I$ is a linear function on $L$, and since $C(S)$ is a linear
subspace of $L$, it is a linear function on $C(S)$. Since $I(V_{1}V_{2})=(I^{(1)}V_{1})(I^{(2)}V_{2})>0$,
the triple $(S,C(S),I)$ satisfies condition (i) of Definition \ref{Def. integration on loc compact space}.
Condition (ii) of Definition \ref{Def. integration on loc compact space},
the positivity condition, follows trivially from the positivity condition
of $(S,L,I$). Hence $(S,C(S),I)$ is an integration space. Since
$C(S)\subset L$ and since $(S,L,I)$ is complete, the completion
$\overline{L}$ of $C(S)$ relative to $I$ is such that $\overline{L}\subset L$.

We will show that, conversely, $L\subset\overline{L}$. To that end,
consider any $Y_{1}\in L^{(1)}$ and $Y_{2}\in L^{(2)}$. Let $\varepsilon>0$
be arbitrary. Then there exists $U_{i}\in C(S_{i})$ such that $I^{(i)}|U_{i}-Y_{i}|<\varepsilon$
for each $i=1,2$. Consequently
\[
I|Y_{1}Y_{2}-U_{1}U_{2}|\leq I|Y_{1}(Y_{2}-U_{2})|+I|(Y_{1}-U_{1})U_{2}|
\]
\[
=I^{(1)}|Y_{1}|\cdot I^{(2)}|Y_{2}-U_{2}|+I^{(1)}|Y_{1}-U_{1}|\cdot I_{2}|U_{2}|
\]
\[
\leq I^{(1)}|Y_{1}|\varepsilon+\varepsilon(I^{(2)}|Y_{2}|+\varepsilon).
\]
Since $\varepsilon>0$ is arbitrary while $U_{1}U_{2}\in C(S)$, we
see that $Y_{1}Y_{2}\in\overline{L}$. Since every simple function
on $S$, as in Definition \ref{Def. Simple functions}, is a linear
combination of functions of the form $Y_{1}Y_{2}$ where $Y_{1}\in L^{(1)}$
and $Y_{2}\in L^{(2)}$, we see that $L_{0}\subset\overline{L}$.
On the other hand $\overline{L}$ is complete relative to $I$. Hence
the completion $L$ of $L_{0}$ is contained in $\overline{L}$.

Summing up, we have $\overline{L}=L$. In other words, the completion
of $(S,C(S),I)$ is $(S,L,I)$.
\end{proof}
\begin{prop}
\label{Prop. Product of sigma finite integration spaces is sigma finite}
\textbf{\emph{(Product of $\sigma$-finite integration spaces is $\sigma$-finite).}}
Let $(\Omega',L',I')$ and $(\Omega'',L'',I'')$ be arbitrary integration
spaces which are $\sigma$-finite, with $I$-bases $(A'_{k})_{k=1,2,\cdots}$
and $(A''_{k})_{k=1,2,\cdots}$ respectively. Then the product integration
space 
\[
(\Omega,L,I)\equiv(\Omega'\times\Omega'',L'\otimes L'',I'\otimes I'')
\]
is $\sigma$-finite, with an $I$-basis $(A_{k})_{k=1,2,\cdots}\equiv(A'_{k}\times A''_{k})_{k=1,2,\cdots}$. 
\end{prop}
\begin{proof}
By the definition of an $I$-basis, we have $A'_{k}\subset A'_{k+1}$
and $A''_{k}\subset A''_{k+1}$ for each $k\geq1$. Hence, $A_{k}\equiv A'_{k}\times A''_{k}\subset A'_{k+1}\times A''_{k+1}$
for each $k\geq1$. Consequently, 
\[
\bigcup_{k=1}^{\infty}(A'_{k}\times A''_{k})=(\bigcup_{k=1}^{\infty}A'_{k})\times(\bigcup_{k=1}^{\infty}A''_{k}).
\]
Again by the definition of an $I$-basis, the two unions on the right-hand
side are full subsets in $\Omega',\Omega''$ respectively. Hence the
union on the left-hand side is, according to Proposition \ref{Prop. Cartesian product of integrable funcs is integrable}
a full set in $\Omega$. 

Now let $f\equiv1_{B'}1_{B''}$ where let $B',B''$ are arbitrary
integrable subsets in $\Omega',\Omega''$ respectively. Then 
\[
I(1_{A(k)}f)=I'(1_{A'(k)}1_{B'})I''(1_{A''(k)}1_{B''})\rightarrow I'(1_{B'})I''(1_{B''})=If.
\]
By linearity, it follows that 
\[
I(1_{A(k)}g)\rightarrow Ig
\]
for each simple function $g$ on $\Omega'\times\Omega''$ relative
to $L',L''$. Consider each $h\in L$. Let $\varepsilon>0$ be arbitrary.
Since $(\Omega,L,I)$ is the completion of $(\Omega,L_{0},I)$, where
$L_{0}$ is the space of simple functions on $\Omega'\times\Omega''$
relative to $L',L''$, it follows that $I|h-g|<\varepsilon$ for some
$g\in L_{0}$. Hence
\[
|I1_{A(k)}f-If|\leq|I1_{A(k)}f-I1_{A(k)}g|+|I1_{A(k)}g-Ig|+|Ig-If|<\varepsilon
\]
\[
\leq I|f-g|+|I1_{A(k)}g-Ig|+|Ig-If|<3\varepsilon
\]
for sufficiently large $k\geq1$. Since $\varepsilon>0$ is arbitrary,
we conclude that $I(1_{A(k)}h)\rightarrow Ih$. In particular, if
$A$ is an arbitrary integrable subset of $\Omega$, we have $I(1_{A(k)}1_{A})\rightarrow I1_{A}.$
In other words, $\mu(A_{k}A)\rightarrow\mu(A)$ for each integrable
set $A\subset\Omega$. We have verified the conditions in Definition
\ref{Def. Sigma Finite; I-basis} for $(\Omega,L,I)$ to be $\sigma$-finite,
with $(A_{k})_{k=1,2,\cdots}\equiv(A'_{k}\times A''_{k})_{k=1,2,\cdots}$
as an $I$-basis. 
\end{proof}
The next definition establishes some familiar notations for the special
cases of the Lebesgue integration space on $R^{n}$.
\begin{defn}
\label{Def. Lebesgue integral in n dimension} \textbf{(Lebesgue integration
on $R^{n}$).} The product integration space 
\[
(R^{n},\overline{L},\int\cdots\int\cdot dx_{1}\cdots dx_{n})\equiv(\prod_{i=1}^{n}R,\bigotimes_{i=1}^{n}L,\bigotimes_{i=1}^{n}\int\cdot dx)
\]
is called the \index{Lebesgue integration space}\emph{Lebesgue integration
space of dimension $n$}. Similarly when $R^{n}$ is replaced by an
interval $\prod_{i=1}^{n}[s_{i},t_{i}]\subset R^{n}$. When confusion
is unlikely we will also abbreviate $\int\cdots\int\cdot dx_{1}\cdots dx_{n}$
to $\int\cdot dx$, with the understanding that the dummy variable
$x$ is now a member of $R^{n}$. An integrable function relative
to $\int\cdots\int\cdot dx_{1}\cdots dx_{n}$ will be called \emph{Lebesgue
integrable}. $\square$
\end{defn}
\begin{cor}
\label{Cor. n-dim Lebesgue completion of C(Rn) =00003D n-dim Lesbesgue integrbles}
\textbf{\emph{(Power of Lebesgue integration space based on $R^{n}$
is equal to the completion of Lebesgue integration on the locally
compact metric space $R^{n}$). }}Let $n\geq1$ be arbitrary. Then,
in the notations of Definition \ref{Def. Lebesgue integral in n dimension},
we have $C(R^{n})\subset\overline{L}$. Moreover, $(R^{n},C(R^{n}),\int\cdots\int\cdot dx_{1}\cdots dx_{n})$
is an integration space, and \textup{\emph{its completion}}\emph{
}is equal to the Lebesgue integration space \textup{$(R^{n},\overline{L},\int\cdots\int\cdot dx_{1}\cdots dx_{n})$.}
\end{cor}
\begin{proof}
Let $S_{i}\equiv R$ for each $i=1,\cdots,n$. Let $S\equiv R^{n}.$
Proposition \ref{Prop. Completion of C(S1xS2) =00003D compC(S1) x CompC(S2)}
then applies and yields the desired conclusions.
\end{proof}
\begin{defn}
\label{Def.  Product of countably many complete integration spacess}
\textbf{(Product of countably many complete integration spaces).}
For each $n\geq1$, let $(\Omega^{(n)},L^{(n)},I^{(n)})$ be a complete
integration space. Consider the Cartesian product $\overline{\Omega}\equiv\prod_{i=1}^{\infty}\Omega^{(i)}$.
Let $n\geq1$ be arbitrary. Let $(\prod_{i=1}^{n}\Omega^{(i)},\bigotimes_{i=1}^{n}L^{(i)},\bigotimes_{i=1}^{n}I^{(i)})$
be the product of the first $n$ complete integration spaces. For
each $g\in\bigotimes_{i=1}^{n}L^{(i)}$, define a function $\overline{g}$
on $\overline{\Omega}$ by $doamin(\overline{g})\equiv doamin(g)\times\prod_{i=n+1}^{\infty}\Omega^{(i)}$,
and by $\overline{g}(\omega_{1},\omega_{2},\cdots)\equiv g(\omega_{1},\cdots,\omega_{n})$
for each $(\omega_{1},\omega_{2},\cdots)\in doamin(\overline{g})$.
Let
\[
G_{n}\equiv\{\overline{g}:g\in\bigotimes_{i=1}^{n}L^{(i)}\}.
\]
Then $G_{n}\subset G_{n+1}$. Let $\overline{L}\equiv\bigcup_{n=1}^{\infty}G_{n}$
and define a function $\overline{I}:\overline{L}\rightarrow R$ by
$\overline{I}(\overline{g})\equiv(\bigotimes_{i=1}^{n}I^{(i)})(g)$
if $\overline{g}\in G_{n}$, for each $\overline{g}\in\overline{L}$.
The next theorem says that $\overline{L}$ is a linear space, that
$\overline{I}$ is a well defined linear function, and that $(\overline{\Omega},\overline{L},\overline{I})$
is an integration space. 

Let $(\prod_{i=1}^{\infty}\Omega^{(i)},\bigotimes_{i=1}^{\infty}L^{(i)},\bigotimes_{i=1}^{\infty}I^{(i)})$
denote the completion of $(\overline{\Omega},\overline{L},\overline{I})$,
and call it the \emph{product of the given sequence of complete integration
spaces}\index{product of a sequence of complete integration spaces}.

In the special case where $(\Omega^{(1)},L^{(1)},I^{(1)})=(\Omega^{(2n)},L^{(2)},I^{(2)})=\cdots$
are all equal to the same integration space $(\Omega_{0},L_{0},I_{0})$,
then we write
\[
(\Omega_{0}^{\infty},L_{0}^{\otimes\infty},I_{0}^{\otimes\infty})\equiv(\prod_{i=1}^{\infty}\Omega^{(i)},\bigotimes_{i=1}^{\infty}L^{(i)},\bigotimes_{i=1}^{\infty}I^{(i)})
\]
and call it the countable \emph{power of the integration space}\index{power integration space}
$(\Omega_{0},L_{0},I_{0})$. 
\end{defn}
$\square$
\begin{thm}
\textbf{\emph{\label{Thm. Countable product of complete integrations is well defined.}
(Countable product of complete integration spaces is well defined).}}
Assume the same terms and notations in Definition \ref{Def.  Product of countably many complete integration spacess}.
Then the following holds.

1. The set $\overline{L}$ of functions is a linear space. Moreover,
$\overline{I}$ is a well-defined linear function, and $(\overline{\Omega},\overline{L},\overline{I})$
is an integration space. 

2. Let $N\geq1$ be arbitrary,. Let $Z^{(N)}$ be a measurable function
on $(\Omega^{(N)},L^{(N)},I^{(N)})$ with values in some complete
metric space $(S,d)$. Define the function $\overline{Z}^{(N)}:\Omega\rightarrow S$
by $\overline{Z}^{(N)}(\omega)\equiv Z^{(N)}(\omega_{i})$ for each
$\omega\equiv(\omega_{1},\omega_{2},\cdots)\in\overline{\Omega}$
such that $\omega_{N}\in domain(Z^{(N)})$. Let $M\geq1$ be arbitrary.
Let $f_{j}\in C_{ub}(S,d)$ be arbitrary for each $j\leq M$. Then
\[
\overline{I}(\prod_{j=1}^{M}f_{j}(\overline{Z}^{(j)})=\prod_{j=1}^{M}I^{(j)}f_{j}(\overline{Z}^{(j)})
\]

3. For each $N\geq1$, the function $\overline{Z}^{(N)}$ is measurable
 on the countable product space $(\prod_{i=1}^{\infty}\Omega^{(i)},\bigotimes_{i=1}^{\infty}L^{(i)},\overline{I})$.
\end{thm}
\begin{proof}
1. Obviously $G_{n}$ and $\overline{L}$ are linear spaces. Suppose
$\overline{g}=\overline{h}$ for some $\overline{g}\in G_{n}$ and
$\overline{h}\in G_{m}$ with $n\leq m$. Then 
\[
h(\omega_{1},\cdots,\omega_{m})\equiv\overline{h}(\omega_{1},\omega_{2},\cdots)=\overline{g}(\omega_{1},\omega_{2},\cdots)
\]
\[
\equiv g(\omega_{1},\cdots,\omega_{n})=g(\omega_{1},\cdots,\omega_{n})1_{A}(\omega_{n+1},\cdots,\omega_{m}),
\]
where $A\equiv\prod_{i=n+1}^{m}\Omega^{(i)}$. Hence
\[
\overline{I}(\overline{h})\equiv(\bigotimes_{i=1}^{m}I^{(i)})(h)=(\bigotimes_{i=1}^{m}I^{(i)})(g\otimes1_{A})=(\bigotimes_{i=1}^{n}I^{(i)}\otimes\bigotimes_{i=n+1}^{m}I^{(i)})(g\otimes1_{A})
\]
\[
=((\bigotimes_{i=1}^{n}I^{(i)})(g))\cdot(\bigotimes_{i=n+1}^{m}I^{(i)})(1_{A})=((\bigotimes_{i=1}^{n}I^{(i)})(g))\cdot1=\overline{I}(\overline{g}).
\]
Thus the function $\overline{I}$ is well defined. Linearity of $\overline{I}$
is obvious. The verification of the other conditions in Definition
\ref{Def. Integration Space} is straightforward. Accordingly, $(\overline{\Omega},\overline{L},\overline{I})$
is an integration space. 

2. In view of Fubini's Theorem \ref{Thm. Fubini for product of two integraton spacces},
the proof of Assertions 2 and 3 are straightforward and omitted.
\end{proof}
Following are two results which will be convenient for future reference.
\begin{prop}
\label{Prop. Region below graph of positive integrable func is integrable}
\textbf{\emph{(Region below graph of }}$\mathbf{measurable}$\textbf{\emph{
function is }}$\mathbf{measurable}$\textbf{\emph{ in product space).}}
Let $(Q,L,I)$ be a complete integration space which is $\sigma$-finite.
Let $(\Theta,\Lambda,I_{0})\equiv(\Theta,\Lambda,\int\cdot d\theta)$
be the Lebesgue integration space based on $\Theta\equiv R$ or $\Theta\equiv[0,1]$.
Let $\lambda:Q\rightarrow R$ be an arbitrary $measurable$ function
on $(Q,L,I)$. Then the sets \textup{
\[
A_{\lambda}\equiv\{(t,\theta)\in Q\times\Theta:\theta\leq\lambda(t)\}
\]
and} 
\[
A'_{\lambda}\equiv\{(t,\theta)\in Q\times\Theta:\theta<\lambda(t)\}
\]
are $measurable$ on $(Q,L,I)\otimes(\Theta,\Lambda,I_{0})$. Suppose,
in addition, that $\lambda$ is a non-negative integrable function.
Then the sets \textup{
\[
B_{\lambda}\equiv\{(t,\theta)\in Q\times\Theta:0\leq\theta\leq\lambda(t)\}
\]
and} 
\[
B'_{\lambda}\equiv\{(t,\theta)\in Q\times\Theta:0\leq\theta<\lambda(t)\}
\]
are integrable, with 
\begin{equation}
(I\otimes I_{0})B_{\lambda}=I\lambda=(I\otimes I_{0})B'_{\lambda}.\label{eq:temp-190}
\end{equation}
\end{prop}
\begin{proof}
Let $g$ be the identity function on $\Theta$, with $g(\theta)\equiv\theta$
for each $\theta\in\Theta$. By Proposition \ref{Thm. measurable function on one factor of product integration space can be regarded as measurable on product},
$g$ and $\lambda$ can be regarded as $\mathrm{measurable}$ functions
on $Q\times\Theta$. Define the function $f:Q\times\Theta\rightarrow R$
by $f(t,\theta)\equiv g(\theta)-\lambda(t)\equiv\theta-\lambda(t)$
for each $(t,\theta)\in Q\times\Theta.$ Then $f$ is the difference
of two real valued $\mathrm{measurable}$ functions on $Q\times\Theta$.
Hence $f$ is $\mathrm{measurable}$. Therefore there exists a sequence
$(a_{n})_{n=1,2,\cdots}$ in $(0,\infty)$ with $a_{n}\downarrow0$
such that $(f\leq a_{n})$ is $\mathrm{measurable}$ for each $n\geq1$.
We will write $a_{n}$ and $a(n)$ interchangeably. 

Let $A\subset Q$ and $B\subset\Theta$ be arbitrary integrable subsets
of $Q$ and $\Theta$ respectively. Let $h:\Theta\rightarrow\Theta$
be the identity function, with $h(\theta)\equiv\theta$ for each $\theta\in\Theta$.
Let $m\geq n$ be arbitrary. Then 
\[
I\otimes I_{0}(1_{(f\leq a(n))(A\times B)}-1_{(f\leq a(m))(A\times B)})
\]
\[
=I\otimes I_{0}(1_{(a(m)<f\leq a(n))(A\times B)})
\]
\[
=I\otimes I_{0}(1_{(\lambda-a(n)\leq h<\lambda-a(m))(A\times B)})
\]
\begin{equation}
=I(1_{A}(\int_{[\lambda-a(n),\lambda-a(m))}1_{B}(\theta)d\theta)).\label{eq:temp-237}
\end{equation}
Since, for each $t\in Q$, the Lebesgue measure 
\[
I_{0}[\lambda(t)-a_{n},\lambda(t)-a_{m})=(a_{n}-a_{m})\downarrow0
\]
as $n\rightarrow\infty$, and since $1_{B}$ is integrable, Proposition
\ref{Prop. Existence of modulus of integrability} implies that 
\[
\int_{[\lambda-a(n),\lambda-a(m))}1_{B}(\theta)d\theta\downarrow0
\]
uniformly on $Q$. Since $1_{A}$ is integrable, the Dominated Convergence
Theorem implies that the right-hand side of equality \ref{eq:temp-237}
converges to $0$ as $n\rightarrow\infty$. Consequently, $I\otimes I_{0}(1_{(f\leq a(n))(A\times B)}$
converges as $n\rightarrow\infty$. Therefore, by the Monotone Convergence
Theorem, the limit $1_{(f\leq0)(A\times B)}$ is integrable on $Q\times\Theta$. 

Now let $C$ be an arbitrary integrable subset of $Q\times\Theta$.
Let $(A_{i})_{i=1,2,\cdots}$ and $(B_{i})_{i=1,2,\cdots}$ be $I$-bases
of the $\sigma$-finite integration spaces $(Q,L,I)$ and $(\Theta,\Lambda,I_{0})$
respectively. Then, by Proposition \ref{Prop. Product of sigma finite integration spaces is sigma finite},
$Q\times\Theta$ is $\sigma$-finite with an $I$-basis $(A_{i}\times B_{i})_{i=1,2,\cdots}$.
By the previous paragraph, $1_{(f\leq0)(A(i)\times B(i))C}=1_{(f\leq0)(A(i)\times B(i))}1_{C}$
is integrable on $Q\times\Theta$, for each $i\geq1$. Moreover, as
$i,j\rightarrow\infty$ with $j\geq i$, we have 
\[
0\leq I1_{(f\leq0)C(A(j)\times B(j))}-I1_{(f\leq0)C(A(i)\times B(i))}\leq I1_{C(A(j)\times B(j))}-I1_{C(A(i)\times B(i))}\rightarrow0
\]
Hence, by the Monotone Convergence Theorem, $1_{(f\leq0)CD}$ is integrable,
where $D\equiv\bigcup_{i=1}^{\infty}(A_{i}\times B_{i})$ is a full
set. Consequently, $1_{(f\leq0)C}$ is integrable. In other words,
$1_{(f\leq0)}1_{C}$ is integrable. Since the integrable subset $C$
of $Q\times\Theta$ is arbitrary, we conclude that $1_{(f\leq0)}$
is $\mathrm{measurable}$. Equivalently, $(f\leq0)$ is $\mathrm{measurable}$.
Recalling the definition of $f$ at the beginning of this proof, we
obtain 
\[
A_{\lambda}\equiv\{(t,\theta)\in Q\times\Theta:\theta-\lambda(t)\leq0\}=\{(t,\theta)\in Q\times\Theta:0\leq f(t,\theta)\}\equiv(f\leq0),
\]
whence $A_{\lambda}$ is $\mathrm{measurable}$. Similarly $A'_{\lambda}$
is $\mathrm{measurable}$. 

Suppose, in addition, that $\lambda$ is a non-negative integrable
function. Then, for each $t\in Q$, we have 
\[
\int1_{B(\lambda)}(t,\theta)d\theta=\int(1_{A(\lambda)}(t,\theta)-1_{A(0)}(t,\theta)d\theta=\int_{(0,\lambda(t)]}d\theta=\lambda(t).
\]
Fubini's Theorem therefore yields $(I\otimes I_{0})B_{\lambda}=I\lambda$,
the first half of equality \ref{eq:temp-190}. The second half is
similarly proved.
\end{proof}
$\square$
\begin{prop}
\label{Prop. Region between non-neg functions in C(Q) are integrable.}
\textbf{\emph{(Regions between graphs of integrable functions).}}\textbf{
}Let $(Q,L,I)$ be a complete integration space which is $\sigma$-finite.
Suppose $\lambda_{0}\equiv0\leq\lambda_{1}\leq\lambda_{2}\leq\cdots\leq\lambda_{n}$
are integrable functions with $\lambda_{n}\leq1$. Define $\lambda_{n+1}\equiv1$.
For each $k=1,\cdots,n+1$, define
\[
\Delta_{k}\equiv\{(t,\theta)\in Q\times R:\theta\in(\lambda_{k-1}(t),\lambda_{k}(t))\}.
\]
Then $\Delta_{1},\cdots,\Delta_{n+1}$ are mutually exclusive $\mathrm{measurable}$
subsets in $(Q,L,I)\otimes(R,M,J)$ whose union is a full set. Moreover,
$\Delta_{1},\cdots,\Delta_{n}$ are integrable in $(Q,L,I)\otimes(R,M,J)$,
with integrals equal to $I\lambda_{1}$, $I\lambda_{2}-I\lambda_{1}$,
$\cdots,I\lambda_{n}-I\lambda_{n-1}$ respectively.
\end{prop}
\begin{proof}
Use Proposition \ref{Prop. Region below graph of positive integrable func is integrable}
above. $\square$
\end{proof}

\section{Supplements and Exercises}

$\,$
\begin{xca}
\label{Ex. B=000026B's Def of measurable function} Let $X$ be a
function defined a.e. on a complete integration space $(\Omega,L,I)$.
Show that $X$ is $\mathrm{measurable}$ iff it satisfies the condition
({*}) that for each integrable set $A$ and $\varepsilon>0$ there
exist an integrable set $B$ and an integrable function $Y$ such
that $B\subset A$, $\mu(AB^{c})<\varepsilon$, and $|X-Y|<\varepsilon$
on $B$. Said condition ({*}) is used as the definition of a $\mathrm{measurable}$
function in \cite{BishopBridges85}. Thus our definition of a $\mathrm{measurable}$
function, which we find more convenient, is equivalent to the one
in \cite{BishopBridges85}.
\end{xca}
\emph{Hint.} Suppose $X$ is $\mathrm{measurable}$. Let $A$ be any
integrable set, and let $\varepsilon>0$ be arbitrary. By condition
(ii) in Definition \ref{Def. Measurable Function}, there exists $a>0$
so large that $\mu(|X|\geq a)A<\varepsilon$. Define $B\equiv(|X|<a)A$.
Then $B$ is an integrable set, and $\mu(AB^{c})<\varepsilon$. Define
$Y\equiv(-a)\vee X1_{B}\wedge a$. Then $Y$ is integrable by condition
(i) in Definition \ref{Def. Measurable Function}. Moreover, $Y-X=0$
on $B$. This verifies condition ({*}).

Conversely, suppose condition ({*}) holds. Let $\varepsilon>0$ be
arbitrary. Let the integrable set $B$ and the integrable function
$Y$ satisfy condition $(*)$ for $\varepsilon$. Let $a>0$ be arbitrary.
Write $X_{a}\equiv(-a)\vee X\wedge a$ and $Y_{a}\equiv(-a)\vee Y\wedge a$.
Then
\[
|X_{a}1_{A}-Y_{a}1_{A}|
\]
\[
\leq|X_{a}1_{A}-X_{a}1_{B}|+|X_{a}1_{B}-Y_{a}1_{B}|+|Y_{a}1_{B}-Y_{a}1_{A}|
\]
\[
\leq a1_{AB^{c}}+\varepsilon1_{B}+a1_{AB^{c}}
\]
where
\[
I(a1_{AB^{c}}+\varepsilon1_{B}+a1_{AB^{c}})\leq a\varepsilon+\varepsilon\mu(A)+a\varepsilon\rightarrow0
\]
as $\varepsilon\rightarrow0$ with $a$ fixed. Hence $(-a)\vee X1_{A}\wedge a=X_{a}1_{A}\in L$
by Theorem \ref{Thm. |X-Yn|<Zn where Yn,Zn integrable =000026 IZn->0 implies X integrable}.
This verifies condition (i) in Definition \ref{Def. Measurable Function}.
Next, let $\varepsilon\in(0,1)$ be arbitrary but fixed. Let $a>0$
be so large that $a>I|Y|/\varepsilon$. By Chebychev's inequality,
we have $\mu(|Y|>a)\leq I|Y|/a<\varepsilon$. Hence
\[
\mu(|X|>a+1)A\leq\mu(|X|>a+1)B+\mu(AB^{c})
\]
\[
\leq\mu(|Y|>a)B+\mu(AB^{c})<I|Y|/a+\varepsilon<2\varepsilon
\]
Since $\varepsilon>0$ is arbitrary, we have $\mu(|X|>a)A\rightarrow0$
as $a\rightarrow\infty$, verifying condition (ii) in Definition \ref{Def. Measurable Function}.
Hence $X$ is a $\mathrm{measurable}$ function. $\square$
\begin{xca}
\label{Ex. Subset A of B with same measure =00003D> AUB'  is full set}
Let $A\subset B$ be integrable sets in a complete integration space
$(\Omega,L,I)$ such that $\mu A=\mu B$. Then $A\cup B^{c}$ is a
full set. 
\end{xca}
\emph{{[}Hint{]}$\mu(BA^{c})=\mu B-\mu A=0$. }Hence $BA^{c}$ is
a null set, and its complement $A\cup B^{c}$ is a full set. $\square$

\chapter{Probability Space}

In this chapter, we specialize the study of complete integration spaces
to the case where the constant function 1
 is integrable and has integral equal to $1$. An integrable function
can then be interpreted as an observable in a probabilistic experiment
which, on repeated observations, has an expected value given by its
integral. Likewise, an integrable set can be interpreted as an event,
and its measure as the probability for said event to occur. We will
transition from terms used in measure theory to commonly used terms
in probability theory. Then we will introduce and study more concepts
and tools common in probability theory.

In this chapter, unless otherwise specified, $(S,d)$ will denote
a complete metric space, not necessarily locally compact. Let $x_{\circ}\in S$
be an arbitrary, but fixed, reference point. Recall that $C_{ub}(S)\equiv C_{ub}(S,d)$
stands for the space of bounded and uniformly continuous functions
on $S$, and that $C(S)\equiv C(S,d)$ stands for the space of continuous
functions on $S$ with compact support.

Let $n\geq1$ be arbitrary. Define the auxiliary function $h_{n}\equiv1\wedge(1+n-d(\cdot,x_{\circ}))_{+}\in C_{ub}(S)$.
Note that the function $h_{n}$ has bounded support. Hence $h_{n}\in C(S)$
if $(S,d)$ is locally compact. 

Separately, for each integration space $(\Omega,L,J)$, we will let
$(\Omega,\overline{L},J)$ denote its complete extension. 

\section{Random Variables}
\begin{defn}
\textbf{(Probability Space and r.v.'s)} \label{Def. Probability Space, r.v.s,  Events, Distribution}
Henceforth, unless otherwise specified, $(\Omega,L,E)$ will denote
a complete integration space in which the constant function $1$ is
integrable, with $E1=1$. Then $(\Omega,L,E)$ is called a \emph{probability
space}\index{probability space}. The integration $E$ is called an
\emph{expectation}\index{expectation}, and the integral $EX$ of
each $X\in L$ is called the \index{expected value}\emph{expected
value} of $X$.

A $\mathrm{measurable}$ function $X$ on $(\Omega,L,E)$ with values
in a complete metric space $(S,d)$ is called a\emph{ random variable}\index{random variable},
or r.v. for abbreviation. Two r.v.'s are considered equal if they
have equal values on a full subset of $\Omega$. A real-valued $\mathrm{measurable}$
function $X$ on $(\Omega,L,E)$ is then called a\emph{ real random
variable}\index{real random variable}, or r.r.v. for abbreviation.
An integrable function $X$ is called an \emph{integrable real random
variable}\index{integrable real random variable}, its integral $EX$
called its\emph{ expected value}\index{expected value}. 

A $\mathrm{measurable}$ set is sometimes called an \emph{event}\index{event}.
It is then integrable because $1_{A}\leq1$, and its measure $\mu(A)$
is called its \emph{probability}\index{probability of an event} and
denoted by $P(A)$ or $PA$. The function $P$ on the set of $\mathrm{measurable}$
sets is called the \emph{probability function}\index{probability function}
corresponding to the expectation $E$. Sometimes we will write $E(A)$
for $P(A)$. The set $\Omega$ is called the \emph{sample space}\index{sample space},
a point $\omega\in\Omega$ called a \emph{sample}\index{sample} or
an \emph{outcome}\index{outcome}. If an outcome $\omega$ belongs
to an event $A$, the event $A$ is said to \emph{occur} for $\omega$,
and $\omega$ is said to \emph{realize} $A$. 

The phrases ``\emph{almost surely}'', ``\emph{almost sure}'',
and the abbreviation ``a.s.'' \index{almost sure}will stand for
``almost everywhere'' or its abbreviation ``a.e.''. Henceforth,
unless otherwise specified, equality of r.v.'s and equality of events
will mean a.s. equality, and the term ``complement'' for events
will stand for ``measure-theoretic complement''. If $X$ is an integrable
r.r.v. and $A,B,\cdots$ are events, we will sometimes write $E(X;A,B,\cdots)$
for $EX1_{AB\cdots}$. 

Let $X\in L$ be arbitrary. We will sometimes use the more suggestive
notation
\[
\int E(d\omega)X(\omega)\equiv EX,
\]
where $\omega$ is a dummy variable. For example, if $Y\in L\otimes L\otimes L$,
we can define a function $Z\in L\otimes L$ by the formula
\[
Z(\omega_{1},\omega_{3})\equiv\int E(d\omega_{2})Y(\omega_{1},\omega_{2},\omega_{3})\equiv EY(\omega_{1},\cdot,\omega_{3})
\]
for each $(\omega_{1},\omega_{3})\in\Omega^{2}$ for which the right-hand
side is defined.

$\square$

In \cite{Billingsley68}, a r.v. is called a random element, and a
r.r.v. is called a random variable. Our usage of the two terms follows
\cite{Neveu 65}, for the benefit of both acronyms.
\end{defn}
Being a $\mathrm{measurable}$ function, a r.v. inherits all the definitions
and properties for $\mathrm{measurable}$ functions developed in preceding
chapters. In particular, since the constant function $1$ is integrable,
$\Omega$ is an integrable set, with $1$ as probability, a probability
space is trivially $\sigma$-finite. Therefore r.v.'s inherit the
theorems on $\mathrm{measurable}$ functions which require a $\sigma$-finite
integration space. 

First we restate Definition \ref{Def. Regular =000026 Continuity Pts of Measurable Func}
of regular points in a simpler form, in the context of a probability
space. The reader can verify that, in the present context, the restated
definition below is equivalent to the 
\begin{defn}
\label{Def. Regular=000026 Continuity Pts r.r.v.} \textbf{(Regular
and continuity points of a r.r.v.).} Let $(\Omega,L,E)$ be an arbitrary
probability space. Let $X$ be a r.r.v. on $(\Omega,L,E)$. Then a
point $t\in R$ is a \index{regular point of a r.r.v.}regular point
of $X$ if (i) there exists a sequence $(s_{n})_{n=1,2,\cdots}$ of
real numbers decreasing to $t$ such that $(X\leq s_{n})$ is a measurable
set for each $n\geq1$ and such that $\lim_{n\rightarrow\infty}P(X\leq s_{n})$
exists, and (ii) there exists a sequence $(r_{n})_{n=1,2,\cdots}$
of real numbers increasing to $t$ such that $(X\leq r_{n})$ is measurable
for each $n\geq1$ and such that $\lim_{n\rightarrow\infty}P(X\leq r_{n})$
exists. If in addition the two limits in (i) and (ii) are equal, then
we call $t$ a \index{continuity point of a r.r.v.}continuity point
of the r.r.v. $X$. 

Note that Condition (i) implies that, if $t\in R$ is a regular point
of a r.r.v. $X$, then $(X\leq t)$ is measurable, with 
\[
P(X\leq s_{n})\downarrow P(X\leq t)
\]
for each sequence $(s_{n})_{n=1,2,\cdots}$ satisfying Condition (i),
thanks to the Monotone Convergence Theorem. Similarly, in that case
$(X<t)$ is measurable , with 
\[
P(X\leq r_{n})\uparrow P(X<t)
\]
for each sequence $(r_{n})_{n=1,2,\cdots}$ satisfying Condition (ii).
Consequently, $(X=t)$ is measurable, and $P(X=t)=0$ if $t$ is a
continuity point. $\square$
\end{defn}
We re-iterate Convention \ref{Convention. Only Regular Pts of Measurble Funcs used}
regarding regular points, now in the context of a probability space
and r.r.v.'s.
\begin{defn}
\label{Convention. Only Regular Pts r.r.v.'s} \textbf{(Convention
regarding regular points of r.r.v.'s).} Let $X$ be an arbitrary r.r.v.
When the measurability of the set $(X<t)$ or $(X\leq t)$ is required
in a discussion for some $t\in R$, it is understood that the real
number $t$ has been chosen from the regular points of the r.r.v.
$X$. 

For example, a sequence of statements like ``Let $t\in R$ be arbitrary.
$\cdots$ Then $P(X_{i}>t)<a$ for each $i\geq1$'' means ``Let
$t\in R$ be arbitrary, such that $t$ is a regular point for $X_{i}$
for each $i\geq1$. $\cdots$ Then $P(X_{i}>t)<a$ for each $i=1,2,\cdots$''.
The purpose of this convention is to obviate unnecessary distraction
from the main arguments.

If, for another example, the measurability of the set $(X\leq0)$
is required in a discussion, we would need to first supply a proof
that $0$ is a regular point of $X$, or, instead of $(X\leq0)$,
use $(X\leq a)$ as a substitute, where $a$ is some regular point
near $0$. Unless the exact value $0$ is essential to the discussion,
the latter, usually effortless, alternative will be used. The implicit
assumption of regularity of the point $a$ is clearly possible, for
example, when we have the freedom to pick the number $a$ from some
open interval, thanks to Proposition \ref{Prop.  Countable Exceptional Pts for Meas X int A},
which says that all but countably many real numbers are regular points
of $X$. 

Classically, all $t\in R$ are regular points for each r.r.v. $X$,
and so this convention would be redundant classically.

$\square$
\end{defn}
In the case of a $\mathrm{measurable}$ indicator $X$, it is easily
seen that $0$ and $1$ are regular points. We recall that the indicator
$1_{A}$ and the complement $A^{c}$ of an event are uniquely defined
relative to a.s. equality.
\begin{prop}
\label{Prop. Basic Properties of r.v.} \textbf{\emph{(Basic Properties
of r.v.'s).}} Let $(\Omega,L,E)$ be a probability space.
\end{prop}
\begin{enumerate}
\item \emph{Suppose $A$ is an event. Then $A^{c}$ is an event. Moreover
$(A^{c})^{c}=A$ and $P(A^{c})=1-P(A)$.}
\item \emph{A subset $A$ of $\Omega$ is a \index{full set} full set iff
it is an event with probability 1.}
\item \emph{Let $(S,d)$ be a complete metric space. A function $X:\Omega\rightarrow S$
is a r.v. with values in $(S,d)$ iff }(i)\emph{ $f(X)\in L$ for
each $f\in C_{ub}(S,d)$, and }(ii)\emph{  $P(d(X,x_{\circ})\geq a)\rightarrow0$
as $a\rightarrow\infty$. Note that if $d$ is bounded, then Condition
}(ii)\emph{ is automatically satisfied.}
\item \emph{Let $(S,d)$ be a complete metric space, with a reference point
$x_{\circ}$. }For each $n\geq1$, define $h_{n}\equiv1\wedge(1+n-d(\cdot,x_{\circ}))_{+}\in C_{ub}(S)$.\emph{
Then a function $X:\Omega\rightarrow S$ is a r.v. iff }(i)\emph{
$f(X)\in L$ for each $f\in C_{ub}(S)$ and }(iii)\emph{  $Eh_{n}(X)\uparrow1$
as $n\rightarrow\infty$.} In that case, we have $E|f(X)-f(X)h_{n}(X)|\rightarrow0$,
where $fh_{n}\in C(S)$
\item \emph{Let $(S,d)$ be a locally compact metric space, with a reference
point $x_{\circ}$. For each $n\geq1$, define the function $h_{n}$
as} above.\emph{ Then $h_{n}\in C(S)$. A function $X:\Omega\rightarrow S$
is a r.v. iff }(iv)\emph{ $f(X)\in L$ for each $f\in C(S)$ and }(iii)\emph{
 $Eh_{n}(X)\uparrow1$ as $n\rightarrow\infty$. In that case, for
each} $f\in C_{ub}(S)$\emph{, there exists a sequence $(g_{n})_{n=1,2,\cdots}$
in $C(S)$ such that $E|f(X)-g_{n}(X)|\rightarrow0$. }
\item \emph{If $X$ is an integrable r.r.v. and $A$ is an event, then $EX=E(X;A)+E(X;A^{c})$.}
\item \emph{A point $t\in R$ is a \index{regular point of integrable function}regular
point of a r.r.v. $X$ iff it is a regular point relative to $\Omega$. }
\item \emph{If $X$ is a r.r.v. such that }$(t-\varepsilon<X<t)\cup(t<X<t+\varepsilon)$
\emph{is a null set  for some $t\in R$ and $\varepsilon>0$, then
the point $t\in R$ is a regular point of $X$.}
\end{enumerate}
\begin{proof}

1. Suppose $A$ is an event with indicator $1_{A}$ and complement
$A^{c}=(1_{A}=0)$. Because $1$ is integrable, so is $1-1_{A}$.
At the same time $A^{c}=(1_{A}=0)=(1-1_{A}=1)$. Hence $A^{c}$ is
an event with indicator $1-1_{A}$. Moreover $P(A^{c})=E(1-1_{A})=1-P(A)$.
Repeating the argument with the event $A^{c}$, we see that
\[
(A^{c})^{c}=(1-(1-1_{A})=1)=(1_{A}=1)=A
\]

2. Suppose $A$ is a full set. Since any two full sets are equal a.s.,
we have $A=\Omega$ a.s. Hence $P(A)=P(\Omega)=1$. Conversely, if
$A$ is an event with $P(A)=1$ then, according to Assertion 1, $A^{c}$
is a null set with $A=(A^{c})^{c}$. Hence by Proposition \ref{Prop. Property of Integrable Sets},
$A$ is a full set.

3. Suppose $X$ is a r.v. Since $\Omega$ is an integrable set, Conditions
(i) and (ii) hold as special cases of Conditions (i) and (ii) in Definition
\ref{Def. Measurable Function} when we take $A=\Omega$.

Conversely, suppose conditions (i) and (ii) hold. Let $f\in C_{ub}(S)$
be arbitrary and let $A$ be an arbitrary integrable set. Then $f(X)\in L$
by condition (i), and so $f(X)1_{A}\in L$. Moreover $P(d(x_{\circ},X)\geq a)A\leq P(d(x_{\circ},X)\geq a)\rightarrow0$
as $a\rightarrow\infty$. Thus conditions (i) and (ii) in Definition
\ref{Def. Measurable Function} are established for $X$ to be a $\mathrm{measurable}$
function. In other words, $X$ is a r.v. Assertion 3 is proved.

4. Given Condition (i), the Conditions (ii) and (iii) are equivalent
to each other, thanks to \ref{Prop. Alternative definition of measurablility},
Thus Assertion 4 follows from Assertion 3.

\emph{5. }Suppose $(S,d)$ is locally compact.\emph{ }Assume that
Condition (iii) holds. In view of Assertion 4, we need only verify
that Conditions (i) and (iv) are then equivalent. Trivially Condition
(i) implies Condition (iv). Conversely, suppose Condition (iv) holds.
Let $f\in C_{ub}(S)$ be arbitrary. We need to prove that $f(X)\in L$.
There is no loss of generality in assuming that $0\leq f\leq b$ for
some $b>0$. Then 
\[
E(f(X)h_{m}(X)-f(X)h_{n}(X))\leq bE(h_{m}(X)-h_{n}(X))\rightarrow0
\]
as $m\geq n\rightarrow\infty$, thanks to Condition (iii). Thus $Ef(X)h_{n}(X)$
converges as $n\rightarrow\infty$. Hence the Monotone Convergence
Theorem implies that $\lim_{n\rightarrow\infty}f(X)h_{n}(X)$ is integrable.
Since $\lim_{n\rightarrow\infty}fh_{n}=f$ on $S$, so $f(X)=\lim_{n\rightarrow\infty}f(X)h_{n}(X)\in L$.
Thus Condition (i) holds. Summing up, given Condition (iii), the Conditions
(i) and (iv) are equivalent to each other, as alleged. The Monotone
Convergence Theorem implies also that $E|f(X)h_{n}(X)-f(X)|\rightarrow0$,
where $fh_{n}\in C(S)$ for each $n\geq1$. Assertion 5 is proved. 

6. $EX=EX(1_{A}+1_{A^{c}})=EX1_{A}+EX1_{A^{c}}\equiv E(X;A)+E(X;A^{c})$.

7. Trivial.

8. Suppose $X$ is a r.r.v. such that $B\equiv(t-\varepsilon<X<t)\cup(t<X<t+\varepsilon)$
is a null set\emph{ }for some $t\in R$ and $\varepsilon>0$. Let
$(s_{n})_{n=1,2,\cdots}$ be a sequence of regular points of $X$
in $(t,t+\varepsilon)$ which decreases to $t$. Then $(s_{n}<X)=(s_{n+1}<X)$
a.s., because $(s_{n+1}<X\leq s_{n})\subset B$ is a null set. Hence
$\lim_{n\rightarrow\infty}P(s_{n}<X)$ exists. Similarly, there exists
a sequence $(r_{n})_{n=1,2,\cdots}$ of regular points of $X$ in
$(t-\varepsilon,t)$ which increases to $t$ such that $\lim_{n\rightarrow\infty}P(r_{n}<X)$.
The conditions in Definition \ref{Def. Regular =000026 Continuity Pts of Measurable Func}
have been proved for $t$ to be a regular point of $X$.
\end{proof}
We will make heavy use of the following Borel-Cantelli Lemma, so much
so we will not bother mentioning its name.
\begin{prop}
\textbf{\emph{(First Borel-Cantelli Lemma)}}\label{Prop. First Borel-Cantelli}
\index{first Borel-Cantelli lemma} Suppose $(A_{n})_{n=1,2,\cdots}$
is a sequence of events such that $\sum_{n=1}^{\infty}P(A_{n})$ converges.
Then a.s. only a finite number of the events $A_{n}$ occur. More
precisely, we have  $P(\bigcup_{k=1}^{\infty}\bigcap_{n=k}^{\infty}A_{n}^{c})=1$.
\end{prop}
\begin{proof}
By Proposition \ref{Prop. Sequences of intgrable sets}, for each
$k\geq1$, the union $B_{k}\equiv\bigcup_{n=k}^{\infty}A_{n}$ is
an event, with $P(B_{k})\leq\sum_{n=k}^{\infty}P(A_{n})\rightarrow0$.
Hence $\lim_{k\rightarrow\infty}P(B_{k}^{c})=1$. Therefore, again
by Proposition \ref{Prop. Sequences of intgrable sets}, the union
$B\equiv\bigcup_{k=1}^{\infty}B_{k}^{c}$ is an event, with 
\[
1=P(B)\equiv P(\bigcup_{k=1}^{\infty}B_{k}^{c})=P(\bigcup_{k=1}^{\infty}\bigcap_{n=k}^{\infty}A_{n}^{c}).
\]
\end{proof}
\begin{defn}
\label{Def. Moments} \textbf{($L_{p}$ space).} Let $X,Y$ be arbitrary
r.r.v.'s Let $p\in[1,\infty)$ be arbitrary. If $X^{p}$ is integrable,
define $\left\Vert X\right\Vert _{p}\equiv(E|X|^{p})^{1/p}$. Define
$L_{p}$ to be the family of all r.r.v. $X$ such that $X^{p}$ is
integrable. We will refer to $\left\Vert X\right\Vert _{p}$ as the
\emph{$L_{p}$-norm\index{L$_{p}$-norm}} of $X$. Let $n\geq1$ be
an integer. If $X\in L_{n}$, then $E|X|^{n}$ is called the $n\textrm{th}$
\emph{absolute moment}\index{absolute moment}, and $EX^{n}$ the
$n\textrm{th}$ \emph{moment}\index{moment}, of $X$. If $X\in L_{1}$,
then $EX$ is also called the \emph{mean}\index{mean} of $X$. 

If $X,Y\in L_{2}$, then, according to Proposition \ref{Prop. Hoelder, Minkowski, Lyapunov}
below, $X,Y,$ and $(X-EX)(Y-EY)$ are integrable. Then $E(X-EX)^{2}$
and $E(X-EX)(Y-EY)$ are respectively called the \emph{variance}\index{variance}
of $X$ and the \emph{covariance}\index{covariance} of $X$ and $Y$.
The square root of the variance of $X$ is called the \emph{standard
deviation}\index{standard deviation} of $X$. $\square$
\end{defn}
Next are several basic inequalities for $L_{p}$.
\begin{prop}
\label{Prop. Hoelder, Minkowski, Lyapunov}\textbf{\emph{ (Basic inequalities
in $L_{p}).$}} Let $p,q\in[1,\infty)$ be arbitrary.
\end{prop}
\begin{enumerate}
\item (\textbf{Hoelder's inequality}) \index{Hoelder's inequality} \emph{Suppose
$p,q>1$ and $\frac{1}{p}+\frac{1}{q}=1$. If $X\in L_{p}$ and $Y\in L_{q}$,
then $XY\in L_{1}$ and $E|XY|\leq\left\Vert X\right\Vert _{p}\left\Vert Y\right\Vert _{q}$.
The special case where $p=q=2$ is referred to as the }\index{Cauchy-Schwarz inequality}\emph{Cauchy-Schwarz
inequality.}
\item \textbf{(Minkowski's inequality)} \index{Minkowski's inequality}
\emph{If $X,Y\in L_{p}$, then $X+Y\in L_{p}$ and $\left\Vert X+Y\right\Vert _{p}\leq\left\Vert X\right\Vert _{p}+\left\Vert Y\right\Vert _{p}$}
\item \textbf{(Lyapunov's inequality)} \index{Lyapunov's inequality}\emph{
If $p\leq q$ and $X\in L_{q}$, then $X\in L_{p}$ and $\left\Vert X\right\Vert _{p}\leq\left\Vert X\right\Vert _{q}$}
\end{enumerate}
\begin{proof}
1. Write $\alpha,\beta$ for $\frac{1}{p},\frac{1}{q}$ respectively.
Then $x^{\alpha}y^{\beta}\leq\alpha x+\beta y$ for non-negative $x,y$.
This can be seen by noting that, with $y$ fixed, the function $f$
defined by $f(x)\equiv\alpha x+\beta y-x^{\alpha}y^{\beta}$ is equal
to $0$ at $x=y$, is decreasing for $x<y$, and is increasing for
$x>y$. Let $a,b\in R$ be arbitrary with $a>\left\Vert X\right\Vert _{p}$
and $b>\left\Vert Y\right\Vert _{q}$. Replacing $x,y$ by $|X/a|^{p},|Y/b|^{q}$
respectively, we see that
\[
|XY|\leq(\alpha|X/a|^{p}+\beta|Y/b|^{q})ab
\]
 It follows that $|XY|$ is integrable, with integral bounded by
\[
E|XY|\leq(\alpha\left\Vert X\right\Vert _{p}^{p}/a^{p}+\beta\left\Vert Y\right\Vert _{q}^{q}/b^{q})ab
\]
 As $a\rightarrow\left\Vert X\right\Vert _{p}$ and $b\rightarrow\left\Vert Y\right\Vert _{q}$,
the last bound approaches $\left\Vert X\right\Vert _{p}\left\Vert Y\right\Vert _{q}$.

2. Suppose first that $p>1$. Let $q\equiv\frac{p}{p-1}$. Then $\frac{1}{p}+\frac{1}{q}=1$.
Because $|X+Y|^{p}\leq(2(|X|\vee|Y|))^{p}\leq2^{p}(|X|^{p}+|Y|^{p})$,
we have $X+Y\in L_{p}$. It follows trivially that $|X+Y|^{p-1}\in L_{q}$.
Applying Hoelder's inequality, we estimate
\[
E|X+Y|^{p}\leq E|X+Y|^{p-1}|X|+E|X+Y|^{p-1}|Y|
\]
\[
\leq(E|X+Y|^{(p-1)q})^{1/q}(\left\Vert X\right\Vert _{p}+\left\Vert Y\right\Vert _{p})
\]
\begin{equation}
=(E|X+Y|^{p})^{1/q}(\left\Vert X\right\Vert _{p}+\left\Vert Y\right\Vert _{p})\label{eq:temp-4}
\end{equation}
Suppose $\left\Vert X+Y\right\Vert _{p}>\left\Vert X\right\Vert _{p}+\left\Vert Y\right\Vert _{p}$.
Then inequality \ref{eq:temp-4}, when divided by $(E|X+Y|^{p})^{1/q}$,
would imply $\left\Vert X+Y\right\Vert _{p}=(E|X+Y|^{p})^{1-1/q}\leq\left\Vert X\right\Vert _{p}+\left\Vert Y\right\Vert _{p}$,
a contradiction. This proves Minkowski's inequality for $p>1$. Suppose
now $p\geq1$. Then $|X|^{r},|Y|^{r}\in L_{p/r}$ for any $r<1$.
The preceding proof of the special case of Minkowski's inequality
for the exponent $\frac{p}{r}>1$ therefore implies
\begin{equation}
(E(|X|^{r}+|Y|^{r})^{p/r})^{r/p}\leq(E(|X|^{p})^{r/p}+(E(|Y|^{p})^{r/p}\label{eq:temp-6}
\end{equation}
Since
\[
(|X|^{r}+|Y|^{r})^{p/r}\leq2^{p/r}(|X|^{r}\vee|Y|^{r})^{p/r}
\]
\[
=2^{p/r}(|X|^{p}\vee|Y|^{p})\leq2^{p/r}(|X|^{p}+|Y|^{p})\in L
\]
 we can let $r\rightarrow1$ and apply the Dominated Convergence Theorem
to the left-hand side of inequality \ref{eq:temp-6}. Thus we conclude
that $(|X|+|Y|)^{p}\in L$, and that $(E(|X|+|Y|)^{p})^{1/p}\leq(E(|X|^{p})^{1/p}+(E(|X|^{p})^{1/p}$.
Minkowski's inequality is proved. 

3. Since $|X|^{p}\leq1\vee|X|^{q}\in L$, we have $X\in L_{p}$. Suppose
$E|X|^{p}>(E|X|^{q})^{p/q}$. Let $r\in(0,p)$ be arbitrary. Clearly
$|X|^{r}\in L_{q/r}$. Applying Hoelder's inequality to $|X|^{r}$
and $1$, we obtain
\[
E|X|^{r}\leq(E|X|^{q})^{r/q}
\]
At the same time $|X|^{r}\leq1\vee|X|^{q}\in L$. As $r\rightarrow p$
the Dominated Convergence Theorem yields $E|X|^{p}\leq(E|X|^{q})^{p/q}$,
establishing Lyapunov's inequality. 
\end{proof}
Next we restate and simplify some definitions and theorems of convergence
of $\mathrm{measurable}$ functions, in terms of r.v.'s
\begin{defn}
\label{Def Convergence of r.v.'s in probability etc} \textbf{(Convergence
in probability, a.u., a.s., and in $L_{1}$).} For each $n\geq1$,
let $X_{n},X$ be a functions on the probability space \emph{$(\Omega,L,E)$,}
with values in the complete metric space $(S,d)$. 
\end{defn}
\begin{enumerate}
\item The sequence $(X_{n})$ is said to \emph{converge to} $X$ \emph{almost
uniformly}\index{convergence almost uniformly} (\emph{a.u.}) on the
probability space $(\Omega,L,E)$ if $X_{n}\rightarrow X$ a.u. on
the integration space $(\Omega,L,E)$. In that case we write $X=\mathrm{a.u.}\lim_{n\rightarrow\infty}X_{n}$.
Since \emph{$(\Omega,L,E)$ }is a probability space, $\Omega$ is
a full set. It can therefore be easily verified that $X_{n}\rightarrow X$
a.u. iff for each $\varepsilon>0$, there exists a $\mathrm{measurable}$
set $B$ with $P(B)<\varepsilon$ such that $X_{n}$ converges to
$X$ uniformly on $B^{c}$. 
\item The sequence $(X_{n})$ is said to \emph{converge to} $X$ \emph{in}
\emph{probability} \index{convergence in probability} on the probability
space $(\Omega,L,E)$ if $X_{n}\rightarrow X$ in measure. Then we
write $X_{n}\rightarrow X$ in probability. It can easily be verified
that $X_{n}\rightarrow X$ in probability iff for each $\varepsilon>0$,
there exists $p\geq1$ so large that, for each $n\geq p$, there exists
a $\mathrm{measurable}$ set $B_{n}$ with $P(B_{n})<\varepsilon$
such that $B_{n}^{c}\subset(d(X_{n},X)\leq\varepsilon)$. 
\item The sequence $(X_{n})$ is said to be \index{Cauchy in probability}\emph{Cauchy
in probability }if it is Cauchy in measure\emph{. }It can easily be
verified that $(X_{n})$ is Cauchy in probability\emph{ }iff for each
$\varepsilon>0$, there exists $p\geq1$ so large that for each $m,n\geq p$,
there exists a $\mathrm{measurable}$ set $B_{m,n}$ with $P(B_{m,n})<\varepsilon$
such that $B_{m,n}^{c}\subset(d(X_{n},X_{m})\leq\varepsilon)$. 
\item The sequence $(X_{n})$ is said to converge to $X$ \index{convergence almost surely}\emph{almost
surely }(a.s.)\emph{ }if $X_{n}\rightarrow X$ a.e. 
\end{enumerate}
\begin{prop}
\label{Prop.  a.u. Convergence =00003D> convergence in prob, etc}
\textbf{\emph{(a.u. Convergence implies convergence in probability,
etc).}} For each $n\geq1$, let $X,X_{n}$ be functions on the probability
space $(\Omega,L,E)$, with values in the complete metric space $(S,d)$.
Then the following holds.

1. If\emph{ $X_{n}\rightarrow X$} a.u. then \emph{(i) }$X$ is defined
a.e.\emph{,} \emph{(ii)} $X_{n}\rightarrow X$ in probability, and
(iii) $X_{n}\rightarrow X$ a.s.

2. If \emph{(i)} $X_{n}$ is a r.v. for each $n\geq1$, \emph{and
(ii)} \emph{$X_{n}\rightarrow X$} in probability, then $X$ is a
r.v.

3. If \emph{(i)} $X_{n}$ is a r.v. for each $n\geq1$, and \emph{(ii)}
$X_{n}\rightarrow X$ a.u\emph{.}, then $X$ is a r.v.

4. If \emph{(i)} $X_{n}$ is a r.v. for each $n\geq1$, and \emph{(ii)}
$(X_{n})_{n=1,2,\cdots}$ is Cauchy in probability, then there exists
a subsequence $(X_{n(k)})_{k=1,2,\cdots}$ such that $X\equiv\lim_{k\rightarrow\infty}X_{n(k)}$
is a r.v., with $X_{n(k)}\rightarrow X$ a.u. and $X_{n(k)}\rightarrow X$
a.s. Moreover, $X_{n}\rightarrow X$ in probability.

5. Suppose \emph{(i)} $X_{n},X$ are r.r.v.'s for each $n\geq1$,
\emph{(ii)} $X_{n}\uparrow X$ in probability, and \emph{(iii)} $a\in R$
is a regular point of $X_{n},X$ for each $n\geq0$. Then $P((X_{n}>a)B)\uparrow P((X>a)B)$
for each measurable set $B$.
\end{prop}
\begin{proof}
Assertions 1-3 are trivial consequences of the corresponding assertions
in Proposition \ref{Prop. a.u. Convergence =00003D> Convergence in measure}.
Assertion 4 is a trivial consequence of Proposition \ref{Prop.Sigma-finite:  Seq Cauchy in measure =00003D> a.u convergent subsequence}.
It remains to prove Assertion 5. To that end, let $\varepsilon>0$
be arbitrary. Then, because $a$ is a regular point of the r.r.v.
$X$, there exists $a'>a$ such that $P(a'\geq X>a)<\varepsilon$.
Since, by hypothesis, $X_{n}\uparrow X$ in probability, there exists
$m\geq1$ so large that $P(X-X_{n}>a'-a)<\varepsilon$ for each $n\geq m$.
Now let $n\geq m$ be arbitrary. Let $A\equiv(a'\geq X>a)\cup(X-X_{n}>a'-a)$.
Then $P(A)<2\varepsilon$. Moreover,
\[
P((X>a)B)-P((X_{n}>a)B)\leq P(X>a;X_{n}\leq a)
\]
\[
=P((X>a;X_{n}\leq a)A^{c})+P(A)
\]
\[
<P((X>a)\cap(X_{n}\leq a)\cap((a'<X)\cup(X\leq a))\cap(X-X_{n}\leq a'-a))+2\varepsilon
\]
\[
=P((X_{n}\leq a)\cap(a'<X)\cap(X-X_{n}\leq a'-a))+2\varepsilon
\]
\[
=0+2\varepsilon=2\varepsilon.
\]
Since $P(A)<2\varepsilon$ is arbitrarily small, we see that $P((X_{n}>a)B)\uparrow P((X>a)B)$,
as alleged in Assertion 5..
\end{proof}
The next definition and proposition shows that convergence in probability
can be metrized.
\begin{defn}
\label{Def. Metric rho_prob on space of r.v.'s} \textbf{(Probability
metric on the space of r.v.'s). }Let \emph{$(\Omega,L,E)$} be a probability
space. Let $(S,d)$ be a complete metric space. We will let $M(\Omega,S)$
denote the space of r.v.'s on $(\Omega,L,E)$ with values in $(S,d)$,
where two r.v.'s are considered equal if they are equal a.s. Define
the metric \emph{
\begin{equation}
\rho_{Prob}(X,Y)\equiv E(1\wedge d(X,Y))\label{eq:temp-445-2}
\end{equation}
for each $X,Y\in M(\Omega,S)$. The next proposition proves that $\rho_{Prob}$
is indeed a metric. We will call $\rho_{Prob}$ the probability metric
on the space} $M(\Omega,S)$ of r.v.'s.
\end{defn}
\begin{prop}
\label{Prop. Basics of the probability metric}\textbf{\emph{ (Basics
of the probability metric $\rho_{Prob}$ on the space}} $M(\Omega,S)$
\textbf{\emph{of r.v.'s).}} Let \emph{$(\Omega,L,E)$} be a probability
space. Let $X,X_{1},X_{2},\cdots$ be r.v.'s with values in the complete
metric space $(S,d)$. Then the following holds.
\end{prop}
\begin{enumerate}
\item \emph{The pair $(M(\Omega,S),\rho_{Prob})$ is a metric space. Note
that $\rho_{Prob}\leq1$.}
\item $X_{n}\rightarrow X$ \emph{in} \emph{probability} \emph{iff, for
each $\varepsilon>0$, there exists $p\geq1$ so large that $P(d(X_{n},X)>\varepsilon)<\varepsilon$
for each $n\geq p$. }
\item \emph{Sequential convergence relative to $\rho_{Prob}$ is equivalent
to convergence in probability. }
\item \emph{The metric space $(M(\Omega,S),\rho_{Prob})$ is complete. }
\item \emph{Suppose there exists a sequence $(\varepsilon_{n})_{n=1,2,\cdots}$
of positive real numbers such that $\sum_{n=1}^{\infty}\varepsilon_{n}<\infty$
and such that $\rho_{Prob}(X_{n},X_{n+1})\equiv E(1\wedge d(X_{n},X_{n+1}))<\varepsilon_{n}^{2}$
for each $n\geq1$. Then $Y\equiv\lim_{n\rightarrow\infty}X_{n}$
is a r.v., and $X_{n}\rightarrow Y$ }a.u.
\end{enumerate}
\begin{proof}
1. Let\emph{ $X,Y\in M(\Omega,S)$} be arbitrary. Then $d(X,Y)$ is
a r.r.v according to Proposition \ref{Prop. vector of measurable func (S)   is meas func (S)}.
Hence $1\wedge d(X,Y)$ is an integrable function, and $\rho_{Prob}$
is well defined in equality \ref{eq:temp-445-2}. Symmetry and triangle
inequality for the function $\rho_{Prob}$ are obvious from its definition.
Suppose $\rho_{Prob}(X,Y)\equiv E(1\wedge d(X,Y))=0$. Let $(\varepsilon_{n})_{n=1,2,\cdots}$
be a sequence in $(0,1)$ with $\varepsilon_{n}\downarrow0$. The
Chebychev's inequality implies 
\[
P(d(X,Y)>\varepsilon_{n})=P(1\wedge d(X,Y)>\varepsilon_{n})\leq\varepsilon_{n}^{-1}E(1\wedge d(X,Y))=0
\]
for each $n\geq1$. Hence $A\equiv\bigcup_{n=1}^{\infty}(d(X,Y)>\varepsilon_{n})$
is a null set. On the full set $A^{c}$ , we have $d(X,Y)\leq\varepsilon_{n}$
for each $n\geq1$. Therefore $d(X,Y)=0$ on the full set $A^{c}$.
Thus $X=Y$ in $M(\Omega,S)$. Summing up, $\rho_{Prob}$ is a metric.

2. Suppose $X_{n}\rightarrow X$ in probability. Let $\varepsilon>0$
be arbitrary. Then, according to Definition \ref{Def Convergence of r.v.'s in probability etc},
there exists $p\geq1$ so large that, for each $n\geq p$, there exists
an integrable set $B_{n}$ with $P(B_{n})<\varepsilon$ and $B_{n}^{c}\subset(d(X_{n},X)\leq\varepsilon)$.
Now consider each $n\geq p$. Then\emph{ }$P(d(X_{n},X)>\varepsilon)\leq P(B_{n})<\varepsilon$
for each $n\geq p$. Conversely, suppose, for each $\varepsilon>0$,
there exists $p\geq1$ so large that $P(d(X_{n},X)>\varepsilon)<\varepsilon$
for each $n\geq p$. Let $\varepsilon>0$ be arbitrary and define
the integrable set $B_{n}\equiv(d(X_{n},X)>\varepsilon)$ for each
$n\geq1$. Then $P(B_{n})<\varepsilon$ and $B_{n}^{c}\subset(d(X_{n},X)\leq\varepsilon)$.
Hence $X_{n}\rightarrow X$ in probability according to Definition
\ref{Def Convergence of r.v.'s in probability etc}. 

3. Suppose\emph{ $\rho_{Prob}(X_{n},X)\equiv E(1\wedge d(X_{n},X))\rightarrow0$.}
Let $\varepsilon>0$ be arbitrary. Take $p\geq1$ so large that $E(1\wedge d(X_{n},X))<\varepsilon(1\wedge\varepsilon)$
for each $n\geq p$. Then Chebychev's inequality implies that
\[
P(d(X_{n},X)>\varepsilon)\leq P(1\wedge d(X_{n},X)\geq1\wedge\varepsilon)\leq(1\wedge\varepsilon)^{-1}E(1\wedge d(X_{n},X))<\varepsilon
\]
for each $n\geq p$. Thus $X_{n}\rightarrow X$ in probability, by
Assertion 2. Conversely, suppose $X_{n}\rightarrow X$ in probability.
Then, by Assertion 2, there exists $p\geq1$ so large that $P(d(X_{n},X)>\varepsilon)<\varepsilon$
for each $n\geq p$. Hence
\[
E(1\wedge d(X_{n},X))=E(1\wedge d(X_{n},X))1_{(d(X(n),X)>\varepsilon)}+E(1\wedge d(X_{n},X))1_{(d(X(n),X)\leq\varepsilon)}
\]
\[
\leq E1_{(d(X(n),X)>\varepsilon)}+\varepsilon=P(d(X_{n},X)>\varepsilon)+\varepsilon<2\varepsilon
\]
for each $n\geq p$. Thus $\rho_{Prob}(X_{n},X)\equiv E(1\wedge d(X_{n},X))\rightarrow0$.
Assertion 3 is proved.

4. Suppose\emph{ $\rho_{Prob}(X_{n},X_{m})\equiv E(1\wedge d(X_{n},X_{m}))\rightarrow0$
as $n,m\rightarrow\infty$. }Let $\varepsilon>0$ be arbitrary. Take
$p\geq1$ so large that $E(1\wedge d(X_{n},X_{m}))<\varepsilon(1\wedge\varepsilon)$
for each $n,m\geq p$. Then Chebychev's inequality implies that
\[
P(d(X_{n},X_{m})>\varepsilon)\leq P(1\wedge d(X_{n},X_{m})\geq1\wedge\varepsilon)\leq(1\wedge\varepsilon)^{-1}E(1\wedge d(X_{n},X_{m}))<\varepsilon
\]
for each $n,m\geq p$. Thus the sequence $(X_{n})_{n=1,2,\cdots}$
of functions is Cauchy in probability. Hence Proposition \ref{Prop.Sigma-finite:  Seq Cauchy in measure =00003D> a.u convergent subsequence}
implies that $X\equiv\lim_{k\rightarrow\infty}X_{n(k)}$ is a r.v.
for some subsequence $(X_{n(k)})_{k=1,2,\cdots}$ of $(X_{n})_{n=1,2,\cdots}$
, and that $X_{n}\rightarrow X$ in probability. By Assertion 3, it
then follows that $\rho_{Prob}(X_{n},X)\rightarrow0$. Thus The metric
space $(M(\Omega,S),\rho_{Prob})$ is complete\emph{,} and Assertion
4 is proved. 

5. Assertion 5 is a trivial special case of Proposition \ref{Prop. Sufficiency for a.u. convergence}.
\end{proof}
\begin{cor}
\textbf{\emph{\label{Cor. Reciprocal of an a.s. postive r.r.v.} (Reciprocal
of an a.s. positive r.r.v.)}}\emph{ Let $X$ be a nonnegative r.r.v.
such that $P(X<a)\rightarrow0$ as $a\rightarrow0$. }Define the function
$X^{-1}$ by $domain(X^{-1})\equiv D\equiv(X>0)$ and $X^{-1}(\omega)\equiv(X(\omega))^{-1}$
for each $\omega\in D$. \emph{Then $X^{-1}$ is a r.r.v.}
\end{cor}
\begin{proof}
Let $a_{1}>a_{2}>\cdots>0$ be a sequence such that $PD_{k}^{c}\rightarrow0$
where $D_{k}\equiv P(X\geq a_{k})$ for each $k\geq1$. Then $D=\bigcup_{k=1}^{\infty}D_{k}$,
whence $D$ is a full set. Let $j\geq k\geq1$ be arbitrary. Define
the r.r.v. $Y_{k}\equiv(X\vee a_{k})^{-1}1_{D(k)}$. Then $X^{-1}1_{D(k)}=Y_{k}$.
Moreover $Y_{j}\geq Y_{k}$, and 
\[
(Y_{j}-Y_{k}>0)\subset(1_{D(j)}-1_{D(k)}>0)=D_{j}D_{k}^{c}\subset D_{k}^{c}.
\]
Consequently, since $PD_{k}^{c}\rightarrow0$ as $k\rightarrow\infty$,
the sequence $(Y_{k})_{k=1,2,\cdots}$ converges a.u. Hence, according
to Proposition \ref{Prop.  a.u. Convergence =00003D> convergence in prob, etc},
$Y\equiv\lim_{k\rightarrow\infty}Y_{k}$ is a r.r.v. Since $X^{-1}=Y$
on the full set $D$, so $X^{-1}$ is a r.r.v.
\end{proof}
We see in Proposition \ref{Prop. Basics of the probability metric}
that convergence in $L_{1}$ of r.r.v.'s implies convergence in probability.
The next proposition gives the converse in the case of uniform integrability. 

.
\begin{prop}
\label{Prop. Unif integrable+convgnce in pr-> L1 convg} \textbf{\emph{(Uniform
integrability of sequence of r.r.v.'s and convergence in probability
implies convergence in $L_{1}$).}} Suppose $(X_{n})_{n=1,2,\cdots}$
is a uniformly integrable sequence of r.r.v.'s. If $(X_{n})$ converges
in probability to some r.r.v. $X$, then $X$ is integrable and $X_{n}\rightarrow X$
in $L_{1}$.
\end{prop}
\begin{proof}
Let $\varepsilon>0$. By Proposition \ref{Prop. Alternative def of uniform integrabitty, and simple modulus of integrability},
there exists $\delta>0$ so small that $E(|X_{n}|;A)<\varepsilon$
for each $n\geq1$ and for each event $A$ with $PA<\delta$. By hypothesis
$X_{n}\rightarrow X$ in probability. Hence there exists an integer
$p\geq1$ so large that $PA_{n}<\delta\wedge\varepsilon$ where $A_{n}\equiv(|X_{n}-X|>\varepsilon)$
for each $n\geq p$. Therefore
\[
E|X_{n}-X_{m}|\leq E(|X_{n}-X_{m}|;\:|X_{n}-X_{m}|>2\varepsilon)+2\varepsilon
\]
\[
\leq E(|X_{n}|+|X_{m}|;\:A_{n}\cup A_{m})+2\varepsilon
\]
\[
\leq E(|X_{n}|;A_{n})+E(|X_{n}|;A_{m})+E(|X_{m}|;A_{n})+E(|X_{m}|;A_{m})+2\varepsilon\leq6\varepsilon
\]
for all $m,n\geq p$. Since $\varepsilon>0$ is arbitrary, we have
$E|X_{n}-X_{m}|\rightarrow0$ as $n,m\rightarrow\infty$. Hence $E|X_{n}-Y|\rightarrow0$
for some integrable r.r.v. $Y$, thanks to the completeness of $(\Omega,L,E)$.
Moreover, $P(|Y-X_{n}|>\varepsilon)\rightarrow0$ by Chebychev's inequality.
It follows that
\[
P(|Y-X|>2\varepsilon)\leq P(|Y-X_{n}|>\varepsilon)+P(|X_{n}-X|>\varepsilon)\rightarrow0.
\]
Since $\varepsilon>0$ is arbitrary, we see that $X=Y$ a.s. Since
$Y$ is integrable, so is $X$, with
\[
\lim_{n\rightarrow\infty}E|X_{n}-X|=\lim_{n\rightarrow\infty}E|X_{n}-Y|=0.
\]
\end{proof}
\begin{prop}
\label{Prop. Alternative definition of a.u. convergence} \textbf{\emph{(Necessary
and sufficient condition for a.u. convergence).}} For $n\geq1$, let
$X,X_{n}$ be r.v.'s with values in the locally compact metric space
$(S,d)$. Then the following two conditions are equivalent: \noun{(}\emph{i}\noun{)}
for each $\varepsilon>0$, there exist an integrable set $B$ with
$P(B)<\varepsilon$ and an integer $m\geq1$ such that for each $n\geq m$
we have $d(X,X_{n})\leq\varepsilon$ on $B^{c}$, and \emph{(ii)}
$X_{n}\rightarrow X$ a.u. 
\end{prop}
\begin{proof}
Suppose Condition (i) holds. Let $(\varepsilon_{k})_{k=1,2,\cdots}$
be a sequence of positive real numbers with $\sum_{k=1}^{\infty}\varepsilon_{k}<\infty$.
By hypothesis, for each $k\geq1$ there exist an integrable set \emph{$B_{k}$
}with $P(B_{k})<\varepsilon_{k}$ and an integer $m_{k}\geq1$ such
that for each $n\geq m_{k}$ we have $d(X,X_{n})\leq\varepsilon_{k}$
on $D_{n}B_{k}^{c}$ for some full set $D_{n}$\emph{. }Let $\varepsilon>0$
be arbitrary. Let $p\geq1$ be so large that $\sum_{k=p}^{\infty}\varepsilon_{k}<\varepsilon$
and define $A\equiv\bigcup_{k=p}^{\infty}(B_{k}\cup\bigcup_{n=m_{k}}^{\infty}D_{n}^{c})$.
Then $P(A)\leq\sum_{k=p}^{\infty}\varepsilon_{k}<\varepsilon$. Moreover,
on $A^{c}=\bigcap_{k=p}^{\infty}(\bigcap_{n=m_{k}}^{\infty}D_{n}B_{k}^{c})$
we have \emph{$d(X,X_{n})\leq\varepsilon_{k}$ }for each $n\geq m_{k}$
and each $k\geq p$. Therefore $X_{n}\rightarrow X$ uniformly on
$A^{c}$. Since $P(A)$ is arbitrarily small, $X_{n}\rightarrow X$
a.u. Thus Condition (ii) is verified.

Conversely, suppose Condition (ii) holds. Let $\varepsilon>0$ be
arbitrary. Then, by \ref{Def Convergence in meas, a.u.,L1,a.e.},
there exists there exists a $\mathrm{measurable}$ set $B$ with $P(B)<\varepsilon$
such that $X_{n}$ converges to $X$ uniformly on $B^{c}$. Hence
there exists $m\geq1$ so large that $\bigcup_{n=m}^{\infty}(d(X_{n},X)>{\normalcolor \varepsilon)\subset B}$.
In particular, for each $n\geq m$, we have $(d(X_{n},X)>{\normalcolor \varepsilon)\subset B}$,
whence $d(X,X_{n})\leq\varepsilon$ on $B^{c}$. Condition (i) is
established.
\end{proof}
\begin{defn}
\label{Def. Probability Subspace} \textbf{(Probability subspace).}
Let $(\Omega,L,E)$ be a probability space and let $L'$ be a subset
of $L$. If $(\Omega,L',E)$ is a probability space, then we call
$(\Omega,L',E)$ a \emph{\index{probability subspace}} \emph{probability
subspace} of $(\Omega,L,E)$. When confusion is unlikely, we will
abuse terminology and simply call $L'$ a probability subspace of
$L$, with $\Omega$ and $E$ understood. 

Let $G$ be a non-empty family of r.v.'s with values in a complete
metric space $(S,d)$. Define 
\[
L_{C(ub)}(G)\equiv\{f(X_{1},\cdots,X_{n}):n\geq1;f\in C_{ub}(S^{n});X_{1},\cdots,X_{n}\in G\}.
\]
Then $(\Omega,L_{C(ub)}(G),E)$ is an integration subspace \emph{of}
$(\Omega,L,E)$. Its completion 
\[
L(G)\equiv L(X:X\in G)\equiv\overline{L_{C(ub)}(G)}
\]
will be called the \emph{\index{probability subspace generated by family of r.v.'s}probability
subspace of} $L$ \emph{generated by the family} $G$. 

If $G$ is a finite or countably infinite set $\{X_{1},X_{2},\cdots\}$,
we will write write $L(X_{1},X_{2},\cdots)$ for $L(G)$. $\square$
\end{defn}
$\,$

Note that $L_{C(ub)}(G)$ is a linear subspace of $L$ containing
constants and is closed to the operation of maximum and absolute values.
Hence $(\Omega,L_{C(ub)}(G),E)$ is indeed an integration space, according
to Proposition \ref{Prop. A closed subspace L' of L yields an integration subspace}.
Since $1\in L_{C(ub)}(G)$ with $E1=1$, the completion $(\Omega,L(G),E)$
is a probability space. Any r.r.v. in $L(G)$ has its value determined
once all the values of the r.v.'s in the generating family $G$ have
been observed. Intuitively, $L(G)$ contains all the information obtainable
by observing the values of all $X\in G$. 
\begin{prop}
Let $(\Omega,L,E)$ be a probability space. Let $G$ be a non-empty
family of r.v.'s with values in a locally compact metric space $(S,d)$.
Let
\[
L_{C}(G)\equiv\{f(X_{1},\cdots,X_{n}):n\geq1;f\in C(S^{n});X_{1},\cdots,X_{n}\in G\}.
\]
Then $(\Omega,L_{C}(G),E)$ is an integration subspace \emph{of} $(\Omega,L,E)$.
Moreover its completion $\overline{L_{C}(G)}$ is equal to $L(G)\equiv\overline{L_{C(ub)}(G)}$. 
\end{prop}
\begin{proof}
Note first that $L_{C}(G)\subset L_{C(ub)}(G),$ and $L_{C}(G)$ is
a linear subspace of $\overline{L_{C(ub)}(G)}$ such that if $U,V\in L_{C}(G)$
then $|U|,U\wedge1\in L_{C}(G)$. Hence $L_{C}(G)$ is an integration
subspace \emph{of} $(\Omega,\overline{L_{C(ub)}(G)},E)$ according
to Proposition \ref{Prop. A closed subspace L' of L yields an integration subspace}.
Consequently $\overline{L_{C}(G)}\subset\overline{L_{C(ub)}(G)}$.

Conversely, let $U\in L_{C(ub)}(G)$ be arbitrary. Then $U=f(X_{1},\cdots,X_{n})$
for some $f\in C_{ub}(S^{n})$ and some $X_{1},\cdots,X_{n}\in G$.
Then, by Proposition \ref{Prop. Basic Properties of r.v.}, there
exists a sequence $(g_{k})_{k=1,2,\cdots}$ in $C(S)$ such that $E|f(X)-g_{k}(X)|\rightarrow0$,
where we write $X\equiv(X_{1},\cdots,X_{n})$. Since $g_{k}(X)\in L_{C(G)}\subset\overline{L_{C}(G)}$
for each $k\geq1$, and since $\overline{L_{C}(G)}$ is complete,
we see that $U=f(X)\in\overline{L_{C}(G)}$. Since $U\in L_{C(ub)}(G)$
is arbitrary, we obtain $L_{C(ub)}(G)\subset\overline{L_{C}(G)}$.
Consequently $\overline{L_{C(ub)}(G)}\subset\overline{L_{C}(G)}$.

Summing up, $\overline{L_{C}(G)}=\overline{L_{C(ub)}(G)}\equiv L(G)$,
as alleged.
\end{proof}
The next lemma sometimes comes in handy.
\begin{lem}
\label{Lem. Intgersection of probability subspaces is a prob subspace}
\textbf{\emph{(Intersection of probability subspaces is a probability
subspace).}} Let $(\Omega,L,E)$ be a probability space. Let $\widehat{L}$
be a nonempty family of probability subspaces $L'$ of $L$. Then
$L''\equiv\bigcap_{L'\in\widehat{L}}L'$ is a probability subspace
of $L$. 
\end{lem}
\begin{proof}
Clearly the intersection $L''$ is a linear subspace of $L$, contains
the constant function $1$ with $E1=1$, and is such that if $X,Y\in L''$
then $|X|,X\wedge1\in L''$. Hence it is an integration subspace of
$L$, according to Proposition \ref{Prop. A closed subspace L' of L yields an integration subspace}.
At the same time, since the sets $L'$ in the family $\widehat{L}$
are closed in the space $L$ relative to the norm $E|\cdot|$, so
is their intersection $L''$. Since $L$ is complete relative to $E$,
so is the closed subspace $L''$. Summing up, $(\Omega,L'',I)$ is
a probability subspace of $(\Omega,L,I)$. 
\end{proof}

\section{Probability Distributions on a Metric Space}
\begin{defn}
\label{Def. distributions on complete metric space} \textbf{(Distribution
on a complete metric space).} Suppose $(S,d)$ is a complete metric
space. Let $n\geq1$ be arbitrary. Recall the function $h_{n}\equiv1\wedge(1+n-d(\cdot,x_{\circ}))_{+}\in C_{ub}(S,d),$
where $x_{\circ}\in S$ is an arbitrary, but fixed, reference point.
Note that the function $h_{n}$ has bounded support. Hence $h_{n}\in C(S)$
if $(S,d)$ is locally compact. Let $J$ be an integration on $(S,C_{ub}(S,d)),$
in the sense of Definition \ref{Def. Integration Space}. Suppose
$Jh_{n}\uparrow1$ as $n\rightarrow\infty$. Then the integration
$J$ is called a \index{distribution on a metric space}\emph{ probability
distribution,} or simply a \emph{distribution, on $(S,d)$}. We will
let $\widehat{J}(S,d)$ denote the set of distributions on the complete
metric space\index{set of distributions on  complete metric space}
\emph{$(S,d)$.} $\square$
\end{defn}
\begin{lem}
\label{Lem. Distribution basics} \textbf{\emph{(Distribution basics).
}}Suppose $(S,d)$ is a complete metric space. Then the following
holds.

1. Let $J$ be an arbitrary distribution on $(S,d)$. Then $1\in L$
and $J1=1$, where $(S,L,J)\equiv(S,\overline{C_{ub}(S)},J)$. Thus
$(S,L,J)$ is a probability space.

2. Suppose the metric space $(S,d)$ is bounded. Let $J$ be an integration
on $(S,C_{ub}(S))$ such that $J1=1$. Then the integration $J$ is
a distribution on $(S,d)$. 

3. Suppose $(S,d)$ is locally compact. Let $J$ be an integration
on $(S,C(S))$ in the sense of Definition \ref{Def. integration on loc compact space}.
Suppose $Jh_{n}\uparrow1$ as $n\rightarrow\infty$. Then $J$ is
a distribution on $(S,d)$. 
\end{lem}
\begin{proof}
1. By Definition \ref{Def. distributions on complete metric space},
$Jh_{n}\uparrow1$ as $n\rightarrow\infty$. At the same time $h_{n}\uparrow1$
on $S$. The Monotone Convergence Theorem therefore implies that $1\in L$
and $J1=1$.

2. Suppose $(S,d)$ is bounded. Then $h_{n}=1$ for sufficiently large
$n\geq1$. Hence, trivially $Jh_{n}\uparrow J1=1$, where the equality
is by assumption. Therefore the integration $J$ ion $(S,C_{ub}(S))$
satisfies the conditions in Definition \ref{Def. distributions on complete metric space}
to be a distribution.

3. Since $(S,d)$ is locally compact. Then $h_{n}\in C(S)$ for each
$n\geq1$.Moreover $Jh_{n}\uparrow1$ by hypothesis. Let $(S,L,J)$
denote the completion of $(S,C(S),J)$. Let $f\in C_{ub}(S)$ be arbitrary,
with some bound $b\geq0$ for $|f|$. Then
\[
J|h_{m}f-h_{n}f|\leq bJ|h_{m}-h_{n}|=bJ(h_{m}-h_{n})\rightarrow0
\]
as $m\geq n\rightarrow\infty$. Hence the sequence $(h_{n}f)_{n=1,2,\cdots}$is
Cauchy in the complete integration space $L$ relative to $J$. Therefore
$g\in L$ and $Jg=\lim_{n\rightarrow\infty}Jh_{n}f$, where $g\equiv\lim_{n\rightarrow\infty}(h_{n}f)$.
At the same time, $\lim_{n\rightarrow\infty}(h_{n}f)=f$ on $S$.
Hence $f=g\in L$, with $Jf=Jg=\lim_{n\rightarrow\infty}Jh_{n}f$.
Since $f\in C_{ub}(S)$ is arbitrary, we conclude that $C_{ub}(S)\subset L$.
Consequently $(S,C_{ub}(S),J)$ is an integration subspace of $(S,L,J)$.
Moreover, in the special case $f\equiv1$, we obtain $1\in L$ with
$J1=\lim_{n\rightarrow\infty}Jh_{n}=1$. Thus the integration $J$
on $C_{ub}(S)$ satisfies the conditions in Definition \ref{Def. distributions on complete metric space}
to be a distribution.
\end{proof}
\begin{defn}
\label{Def. Distribution induced by r.v. w/ values in complete metric space}
\textbf{(Distribution induced by a r.v.)} Let $X$ be a r.v. on a
probability space $(\Omega,L,E)$ with values in the complete metric
space $(S,d)$. For each $f\in C_{ub}(S)$, define $E_{X}f\equiv Ef(X)$.
Lemma \ref{Lem. Distribution induced by r.v is indeed a disstributon }
below proves that $E_{X}$ is a distribution on $(S,d)$. We will
call $E_{X}$ the\emph{ distribution}\index{distribution induced by a r.v.}
\emph{on $S$ induced by the r.v. $X$.} The completion $(S,L_{X},E_{X})\equiv(S,\overline{C_{ub}(S)},E_{X})$
of $(S,C_{ub}(S),E_{X})$ is a probability space, called the \index{probability space induced by a r.v.}\emph{probability
space induced on the complete metric space} $(S,d)$\emph{ by the
r.v.} $X$.\emph{ }$\square$
\end{defn}
\begin{lem}
\label{Lem. Distribution induced by r.v is indeed a disstributon }
\textbf{\emph{(Distribution induced by a r.v. is indeed a distribution).}}
Let $X$ be an arbitrary r.v. on a probability space $(\Omega,L,E)$
with values in the complete metric space $(S,d)$. Then \emph{the
function} $E_{X}$ introduced in Definition \ref{Def. Distribution induced by r.v. w/ values in complete metric space}
is indeed a distribution. 
\end{lem}
\begin{proof}
Let $f\in C_{ub}(S)$ be arbitrary. By Proposition \ref{Prop. Basic Properties of r.v.},
we have $f(X)\in L$. Hence $E_{X}f\equiv Ef(X)$ is well-defined.
The space $C_{ub}(S)$ is linear, contains constants, and is closed
to absolute values and taking minimums. The remaining conditions in
Definition \ref{Def. Integration Space} for $E_{X}$ to be an integration
on $(S,C_{ub}(S))$ follow from the corresponding conditions for $E$.
Moreover, $E_{X}h_{n}\equiv Eh_{n}(X)\uparrow1$ as $n\rightarrow\infty$,
where the convergence is by by Assertion 4 in Proposition \ref{Prop. Basic Properties of r.v.}.
All the conditions in Definition \ref{Def. distributions on complete metric space}
have been verified for $E_{X}$ to be a distribution. 
\end{proof}
\begin{prop}
\label{Prop. Each distribution on complete (S,d) is I_X} \textbf{\emph{(Each
Distribution is induced by some r.v.)}} Suppose $J$ is a distribution
on a complete metric space $(S,d)$. Let $(S,L,J)$ denote the completion
of the integration space $(S,C_{ub}(S),J)$. Then the following holds.

1. The identity function $X:(S,L,J)\rightarrow(S,d)$, defined by
$X(x)=x$ for each $x\in S$, is a r.v.

2. The function $d(\cdot,x_{\circ})$ is a r.r.v. on $(S,L,J)$.

3. $J=E_{X}$. Thus each distribution is induced by some r.v.
\end{prop}
\begin{proof}
By Lemma \ref{Lem. Distribution basics}, $(S,L,J)$ is a probability
space, and $Jh_{n}\uparrow1$ as $n\rightarrow\infty$. Hence the
hypothesis in Corollary \ref{Cor. Identity function X(x)=00003Dx  is meas on completion of (S,C(S),I)}
is satisfied. Accordingly, $X$ is a r.v. on $(\Omega,L,E)\equiv(S,L,J)$,
and $d(\cdot,x_{\circ})$ is a r.r.v. Moreover, for each $f\in C_{ub}(S)$,
we have $Jf\equiv Ef\equiv Ef(X)\equiv E_{X}f$. Hence $J=E_{X}$
on $C_{ub}(S)$. Therefore the completion of $(S,C_{ub}(S,J),J)$
and that of $(S,C_{ub}(S,J),E_{X})$ are the same. In other words,
$(S,L,J)=(S,L_{X},E_{X})$.
\end{proof}
\begin{prop}
\label{Prop. L(X)=00003D=00007Bf(X): f in L_X=00007D} \textbf{\emph{(Relation
between probability spaces generated and induced by a r.v.)}} Suppose
$X$ is a r.v. on the probability space $(\Omega,L,E)$ with values
in a complete metric space $(S,d)$. Let $(\Omega,L(X),E)$ be the
probability subspace generated by $\{X\}$. Let $(S,L_{X},E_{X})$
be the probability space induced on $(S,d)$ by $X$. Let $f:S\rightarrow R$
be an arbitrary function. Then the following holds.
\end{prop}
\begin{enumerate}
\item $f\in L_{X}$ iff $f(X)\in L(X)$, in which case $E_{X}f=Ef(X)$.
\item \emph{$f$} is a r.r.v. on $(S,L_{X},E_{X})$ iff $f(X)$ is a r.r.v.
on $(\Omega,L(X),E)$. 
\end{enumerate}
\begin{proof}
1. Suppose $f\in L_{X}$. Then there exists sequence  $(f_{n})_{n=1,2,\cdots}$
in $C_{ub}(S)$ such that $E_{X}|f_{n}-f|\rightarrow0$ and $f=\lim_{n\rightarrow\infty}f_{n}$.
Consequently
\[
E|f_{n}(X)-f_{m}(X)|\equiv E_{X}|f_{n}-f_{m}|\rightarrow0.
\]
Thus $(f_{n}(X))_{n=1,2,\cdots}$ is a Cauchy sequence in $L(X)$
relative to the expectation $E$. Since $L(X)$ is complete, we have
$Y\equiv\lim_{n\rightarrow\infty}f_{n}(X)\in L(X)$ with
\[
E|f_{n}(X)-Y|\rightarrow0,
\]
whence
\[
EY=\lim_{n\rightarrow\infty}Ef_{n}(X)\equiv\lim_{n\rightarrow\infty}E_{X}f_{n}=E_{X}f.
\]
Since $f(X)=\lim_{n\rightarrow\infty}f_{n}(X)=Y$ on the full set
$domain(Y)$, it follows that $f(X)\in L(X)$, with $Ef(X)=EY=E_{X}f$. 

Conversely, suppose $Z\in L(X)$. We will show that Z$=f(X)$ for
some integrable function $f$ relative to $E_{X}$. Since $L(X)$
is, by definition, the completion of $L_{C(ub)}(X)\equiv\{f(X):f\in C_{ub}(S)\}$,
the latter is dense in the former, relative to the norm $E|\cdot|$.
Hence there exists a sequence $(f_{n})_{n=1,2,\cdots}$ in $C_{ub}(S)$
such that 
\begin{equation}
E|Z-f_{n}(X)|\rightarrow0.\label{eq:temp-133}
\end{equation}
Consequently
\[
E_{X}|f_{n}-f_{m}|\equiv E|f_{n}(X)-f_{m}(X)|\rightarrow0.
\]
Hence $E_{X}|f_{n}-f|\rightarrow0$ where $f\equiv\lim_{n\rightarrow\infty}f_{n}\in\overline{C_{ub}(S)}\equiv L_{X}$.
. By the first part of this proof in the previous paragraph, we have
\begin{equation}
E|f_{n}(X)-Y|\rightarrow0,\label{eq:temp-185}
\end{equation}
where
\begin{equation}
Y=f(X)\qquad\mathrm{a.s.}\label{eq:temp-186}
\end{equation}
Convergence expressions \ref{eq:temp-133} and \ref{eq:temp-185}
together implies that $Z=Y$ a.s., which, together with equality \ref{eq:temp-186},
in turn yields $Z=f(X)$, where $f\in L_{X}$. Assertion 1 is proved.

2. For each $n\geq1$, define $g_{n}\equiv1\wedge(1+n-|\cdot|)_{+}\in C_{ub}(R)$.
Suppose the function\emph{ $f$} is a r.r.v. on $(S,L_{X},E_{X})$.
Then, by Proposition \ref{Prop. Basic Properties of r.v.}, we have
(i) $g\circ f\in L_{X}$ for each $g\in C_{ub}(R)$, and (ii)  $E_{X}g_{n}\circ f\uparrow1$
as $n\rightarrow\infty$. 

In view of Condition (i), we have $g(f(X))\equiv g\circ f(X)\in L(X)$
for each $g\in C_{ub}(R)$ by Assertion 1. Moreover, $Eg_{n}(f(X))=E_{X}g_{n}\circ f\uparrow1$
as $n\rightarrow\infty$. Combining, we can apply Assertion 4 of Proposition
\ref{Prop. Basic Properties of r.v.} to the function $f(X):\Omega\rightarrow R$
in the place of $X:\Omega\rightarrow S$, and conclude that $f(X)$
is a r.r.v. on $(\Omega,L(X),E)$.

Conversely, suppose $f(X)$ is a r.r.v. on $(\Omega,L(X),E)$. Then,
again by Assertion 4 of Proposition \ref{Prop. Basic Properties of r.v.},
we have (i') $g(f(X))\in L(X)$ for each $g\in C_{ub}(R)$, and (ii')
 $E(g_{n}(f(X))\uparrow1$ as $n\rightarrow\infty$. In view of Condition
(i'), we have $g\circ f\in L_{X}$ for each $g\in C_{ub}(R)$ by Assertion
1 of the present proposition. Moreover, $E_{X}g_{n}\circ f=Eg_{n}(f(X))=\uparrow1$
as $n\rightarrow\infty$. Combining, we see that $f$ is a r.r.v.
on $(S,L_{X},E_{X}),$ again by Assertion 4 of Proposition \ref{Prop. Basic Properties of r.v.}. 
\end{proof}
\begin{prop}
\label{Prop.  t is continuity pt of f  iff t is continuity pt of f(X)}
\textbf{\emph{(Regular points of a r.r.v. $f$ relative to induced
distribution by a r.v. X are same as regular points of $f(X)).$}}
Suppose $X$ is a r.v. on the probability space $(\Omega,L,E)$ with
values in a complete metric space $(S,d)$. Suppose $f$ is a r.r.v.
on $(S,L_{X},E_{X})$.  Then $t\in R$ is a regular point of $f$
iff it is a regular point of $f(X)$. Similarly, $t\in R$ is a continuity
point of $f$ iff it is a continuity point of $f(X)$. 
\end{prop}
\begin{proof}
Suppose $f$ is a r.r.v. on $(S,L_{X},E_{X})$. By Definition \ref{Def. Regular =000026 Continuity Pts of Measurable Func},
$t$ is a regular point of $f$ iff (i) there exists a sequence $(s_{n})_{n=1,2,\cdots}$
of real numbers decreasing to $t$ such that $(s_{n}<f$) is integrable
relative to $E_{X}$ for each $n\geq1$, and $\lim_{n\rightarrow\infty}P_{X}(s_{n}<f)$
exists, and (ii) there exists a sequence $(r_{n})_{n=1,2,\cdots}$
of real numbers increasing to $t$ such that $(r_{n}<f)$ is integrable
relative to $E_{X}$ for each $n\geq1$, and $\lim_{n\rightarrow\infty}P_{X}(r_{n}<f)$
exists. In view of Proposition \ref{Prop. L(X)=00003D=00007Bf(X): f in L_X=00007D},
conditions (i) and (ii) are equivalent to: (i') there exists a sequence
$(s_{n})_{n=1,2,\cdots}$ of real numbers decreasing to $t$ such
that $(s_{n}<f(X)$) is integrable relative to $E$ for each $n\geq1$,
and $\lim_{n\rightarrow\infty}P(s_{n}<f(X))$ exists, and (ii') there
exists a sequence $(r_{n})_{n=1,2,\cdots}$ of real numbers increasing
to $t$ such that $(r_{n}<f(X))$ is integrable relative to $E$ for
each $n\geq1$, and $\lim_{n\rightarrow\infty}P(r_{n}<f(X))$ exists.
In other words, $t$ is a regular point of $f$ iff $t$ is a regular
point of $f(X)$. 

Moreover, a regular point $t$ of $f$ is a continuity point of $f$
iff the two limits in conditions (i) and (ii) exist and are equal.
Equivalently, $t$ is a continuity point of $f$ iff the two limits
in conditions (i') and (ii') exist and are equal. Combining, we conclude
that $t$ is a continuity point of $f$ iff it is a continuity point
of $f(X)$. 
\end{proof}

\section{Weak Convergence of Distributions}

Recall that, if $X$ is a r.v. on a probability space $(\Omega,L,E)$
with values in $S$, then $E_{X}$ denotes the distribution induced
on $S$ by $X$.
\begin{defn}
\label{Def. Weak Convergence of Distributions on Metric Space} \textbf{(Weak
convergence of distributions on a complete metric space).} Recall
that $\widehat{J}(S,d)$ denotes the set of distributions on the complete
metric space $(S,d)$. A sequence $(J_{n})_{n=1,2,\cdots}$ in $\widehat{J}(S,d)$
is said to \emph{converge weakly}\index{weak convergence of distributions}
to $J\in\widehat{J}(S,d)$ if $J_{n}f\rightarrow Jf$ for each $f\in C_{ub}(S)$.
We then write $J_{n}\Rightarrow J$. Suppose $X,X_{1},X_{2},\cdots$
are r.v.'s with values in $S$, not necessarily on the same probability
space. The sequence $(X_{n})_{n=1,2,\cdots}$ is said to \emph{converge
weakly}\index{weak convergence of r.v.'s}, or to \index{converge in distribution}\emph{converge
in distribution\index{convergence of r.v.'s in distribution}}, to
$X$ if $E_{X(n)}\Rightarrow E_{X}$. We then write $X_{n}\Rightarrow X$.
$\square$
\end{defn}
\begin{prop}
\label{Prop. metric space valued Xn -> X in prob  implies  Xn =00003D>X  }
\textbf{\emph{(Convergence in probability implies weak convergence).}}
Let $(X_{n})_{n=0,1,\cdots}$ be a sequence of r.v.'s on the same
probability space $(\Omega,L,E)$, with values in a complete metric
space $(S,d)$. If$X_{n}\rightarrow X_{0}$ in probability, then $X_{n}\Rightarrow X_{0}$.
\end{prop}
\begin{proof}
Suppose $X_{n}\rightarrow X_{0}$\emph{ }in\emph{ }probability.\emph{
}Let\emph{ }$f\in C_{ub}(S)$ be arbitrary, with $|f|\leq c$ for
some $c>0$, and with a modulus of continuity $\delta_{f}$. Let $\varepsilon>0$
be arbitrary.\emph{ }By Definition \ref{Def Convergence of r.v.'s in probability etc}
of convergence in probability,\emph{ }there exists $p\geq1$ so large
that, for each $n\geq p$, there exists an integrable set $B_{n}$
with $P(B_{n})<\varepsilon$ and
\[
B_{n}^{c}\subset(d(X_{n},X_{0})<\delta_{f}(\varepsilon))\subset(|f(X_{n})-f(X_{0})|<\varepsilon).
\]
Consider each $n\geq p$. Then
\[
|Ef(X_{n})-Ef(X_{0})|=E|f(X_{n})-f(X_{0})|1_{B(n)}+E|f(X_{n})-f(X_{0})|1_{B(n)^{c}}
\]
\[
\leq2cP(B_{n})+\varepsilon<2c\varepsilon+\varepsilon.
\]
Since $\varepsilon>0$ is arbitrarily small, we conclude that $Ef(X_{n})\rightarrow Ef(X_{0})$.
Equivalently, $J_{X(n)}f\rightarrow J_{X(0)}f$ . Since $f\in C_{ub}(S)$
is arbitrary, we have $J_{X(n)}\Rightarrow J_{X(0)}$. In other words,
$X_{n}\Rightarrow X_{0}$.
\end{proof}
\begin{lem}
\label{Lem. locally compact S, Ip=00003D>I iff Ip(f)->I(f) for each f in C(S)}
\textbf{\emph{(Weak convergence of distributions on a locally compact
metric space).}} Suppose (S,d) is locally compact. Suppose $J,J',J_{p}\in\widehat{J}(S,d)$
for each $p\geq1$. Then $J_{p}\Rightarrow J$ iff $J_{p}f\rightarrow Jf$
for each $f\in C(S)$. Moreover, $J=J'$ if $Jf=J'f$ for each $f\in C(S)$.
Consequently, a distribution on a locally compact metric space is
uniquely determined by the expectation of continuous functions with
compact supports.
\end{lem}
\begin{proof}
Since $C(S)\subset C_{ub}(S)$, it suffices to prove the ``if''
part. To that end, suppose $J_{p}f\rightarrow Jf$ for each $f\in C(S)$.
Let $g\in C_{ub}(S)$ be arbitrary. We need to prove that $J_{p}g\rightarrow Jg$
. Let $\varepsilon>0$ be arbitrary. We assume, without loss of generality,
that $0\leq g\leq1$. Since $J$ is a distribution, there exists $n\geq1$
so large that $J(1-h_{n})<\varepsilon$, where $h_{n}\in C_{ub}(S)$
is defined at the beginning of this chapter. Since $h_{n},gh_{n}\in C(S)$,
we have, by hypothesis, $J_{m}h_{n}\rightarrow Jh_{n}$ and $J_{m}gh_{n}\rightarrow Jgh_{n}$
as $m\rightarrow\infty$. Hence 
\[
|J_{m}g-Jg|\leq|J_{m}g-J_{m}gh_{n}|+|J_{m}gh_{n}-Jgh_{n}|+|Jgh_{n}-Jg|
\]
\[
\leq|1-J_{m}h_{n}|+|J_{m}gh_{n}-Jgh_{n}|+|Jh_{n}-1|<\varepsilon+\varepsilon+\varepsilon
\]
for sufficiently large $m\geq1$. Since $\varepsilon>0$ is arbitrary,
we conclude that $J_{m}g\rightarrow Jg$, where $g\in C_{ub}(S)$
is arbitrary. Thus $J_{p}\Rightarrow J$.

Now suppose $Jf=J'f$ for each $f\in C(S)$. Define $J_{p}\equiv J'$
for each $p\geq1$. Then $J_{p}f\equiv J'f=Jf$ for each $f\in C(S)$.
Hence by the previous paragraphs, $J'g\equiv J_{p}g\rightarrow Jg$
for each $g\in C_{ub}(S)$. Thus $J'g=Jg$ for each $g\in C_{ub}(S)$.
In other words, $J=J'$ on $C_{ub}(S)$. We conclude that $J=J'$
as distributions.
\end{proof}
\begin{defn}
\label{Def. Distribution metric} \textbf{(Distribution metric for
distributions on a locally compact metric space).} Suppose the metric
space $(S,d)$ is locally compact, with the reference point $x_{\circ}\in S$.
Let $\xi\equiv(A_{n})_{n=1,2,\cdots}$ be a binary approximation of
$(S,d)$ relative to $x_{\circ}$. Let 
\[
\pi\equiv(\{g_{n,x}:x\in A_{n}\})_{n=1,2,\cdots}
\]
be the partition of unity of $(S,d)$ determined by $\xi$, as in
Definition \ref{Def. Partition of unity for locally compact (S,d)}.

Let $\widehat{J}(S,d)$ denote the set of distributions on the locally
compact metric space $(S,d)$. Let $J,J'\in\widehat{J}(S,d)$ be arbitrary.
Define
\end{defn}
\begin{equation}
\rho_{Dist,\xi}(J,J')\equiv\sum_{n=1}^{\infty}2^{-n}|A_{n}|^{-1}\sum_{x\in A(n)}|Jg_{n,x}-J'g_{n,x}|\label{eq:temp-119}
\end{equation}
and call $\rho_{Dist,\xi}$ the \index{distribution metric for a locally compact space}\emph{distribution
metric on $\widehat{J}(S,d)$ relative to the binary approximation
}$\xi$. The next proposition shows that $\rho_{Dist,\xi}$ is indeed
a metric, and that sequential convergence relative to $\rho_{Dist,\xi}$
is equivalent to weak convergence. Note that $\rho_{Dist,\xi}\leq1$.
$\square$
\begin{prop}
\label{Prop. rho_xi convergence=00003D Weak Convergence} \textbf{\emph{(Sequential
metrical convergence implies weak convergence, on a locally compact
metric space).}} Suppose the metric space $(S,d)$ is locally compact,
with the reference point $x_{\circ}\in S$. Let $\xi\equiv(A_{n})_{n=1,2,\cdots}$
be a binary approximation of $(S,d)$ relative to $x_{\circ}$, with
a corresponding modulus of local compactness $\left\Vert \xi\right\Vert \equiv(|A_{n}|)_{n=1,2,\cdots}$
of $(S,d)$. Let $\rho_{Dist,\xi}$ be the function introduced in
Definition \ref{Def. Distribution metric}.

Let $J_{p}\in\widehat{J}(S,d)$ for $p\geq1$. Let\emph{ $f\in C(S)$
}be arbitrary, with a modulus of continuity $\delta_{f}$, with $|f|\leq1$,
and with $(d(\cdot,x_{\circ})\leq b)$ as support for some $b>0$.
Then the following holds.
\end{prop}
\begin{enumerate}
\item \emph{Let $\varepsilon>0$ be arbitrary. Then there exists $\delta_{\widehat{J}}(\varepsilon)\equiv\delta_{\widehat{J}}(\varepsilon,\delta_{f},b,\left\Vert \xi\right\Vert )>0$
such that, for each }$J,J'\in\widehat{J}(S,d)$\emph{ with $\rho_{Dist,\xi}(J,J')<\delta_{\widehat{J}}(\varepsilon)$
we have $|Jf-J'f|<\varepsilon$. }
\item \emph{Let 
\[
\pi\equiv(\{g_{n,x}:x\in A_{n}\})_{n=1,2,\cdots}
\]
be the partition of unity of $(S,d)$ de}term\emph{ined by $\xi$.}
\emph{Suppose $J_{p}g_{n,x}\rightarrow Jg_{n,x}$ }as $p\rightarrow\infty$,\emph{
for each $x\in A_{n}$, for each $n\geq1$. Then $\rho_{Dist,\xi}(J_{p},J)\rightarrow0$
.}
\item \emph{$J_{p}f\rightarrow Jf$ for each $f\in C(S)$ iff $\rho_{Dist,\xi}(J_{p},J)\rightarrow0$.
Thus $J_{p}\Rightarrow J$ iff $\rho_{Dist,\xi}(J_{p},J)\rightarrow0$.}
\item \emph{$\rho_{Dist,\xi}$ is a metric.}
\end{enumerate}
\begin{proof}
1.\emph{ }Let $\varepsilon>0$ be arbitrary.\emph{ }Let $n\equiv[0\vee(1-\log_{2}\delta_{f}(\frac{\varepsilon}{3}))\vee\log_{2}b]_{1}$.
We will show that 
\[
\delta_{\widehat{J}}(\varepsilon)\equiv\delta_{\widehat{J}}(\varepsilon,\delta_{f},b,\left\Vert \xi\right\Vert )\equiv\frac{1}{3}2^{-n}|A_{n}|^{-1}\varepsilon
\]
has the desired property. To that end, suppose $J,J'\in\widehat{J}(S,d)$
are such that \emph{$\rho_{Dist,\xi}(J,J')<\delta_{\widehat{J}}(\varepsilon)$}.
By Definition \ref{Def. Partition of unity for locally compact (S,d)}
of $\pi$, the sequence $\{g_{n,x}:x\in A_{n}\}$ is a $2^{-n}$-partition
of unity determined by $A_{n}$. Separately, by hypothesis, the function
$f$ has support
\[
(d(\cdot,x_{\circ})\leq b)\subset(d(\cdot,x_{\circ})\leq2^{n})\subset\bigcup_{x\in A(n)}(d(\cdot,x)\leq2^{-n}),
\]
where the first inclusion is because $b<2^{n}$, and the second inclusion
is by Definition \ref{Def. Binary approximationt and Modulus of local compactness}.
Since $2^{-n}<\frac{1}{2}\delta_{f}(\frac{1}{3}\varepsilon)$, Proposition
\ref{Prop. Approx  by  interpolation} then implies that 
\begin{equation}
\left\Vert f-g\right\Vert \leq\frac{\varepsilon}{3}\label{eq:temp-218}
\end{equation}
where 
\[
g\equiv\sum_{x\in A(n)}f(x)g_{n,x}.
\]

By the definition of $\rho_{Dist,\xi}$, we have 
\[
2^{-n}|A_{n}|^{-1}\sum_{x\in A(n)}|Jg_{n,x}-J'g_{n,x}|\leq\rho_{Dist,\xi}(J,J')<\delta_{\widehat{J}}(\varepsilon).
\]
Therefore
\[
|Jg-J'g|\equiv|\sum_{x\in A(n)}f(x)(Jg_{n,x}-J'g_{n,x})|
\]
\[
\leq\sum_{x\in A(n)}|Jg_{n,x}-J'g_{n,x}|<2^{n}|A_{n}|\delta_{\widehat{J}}(\varepsilon)\equiv\frac{1}{3}\varepsilon.
\]
Combining with inequality \ref{eq:temp-218}, we obtain 
\[
|Jf-J'f|\leq|Jg-J'g|+\frac{2}{3}\varepsilon<\frac{1}{3}\varepsilon+\frac{2}{3}\varepsilon=\varepsilon.
\]
Assertion 1 is proved.

2. Suppose \emph{$J_{p}g_{n,x}\rightarrow Jg_{n,x}$ }as $p\rightarrow\infty$,
for each\emph{ $x\in A_{n}$, }for each\emph{ $n\geq1$}. Let $\varepsilon>0$
be arbitrary. Note that
\[
\rho_{Dist,\xi}(J,J_{p})\equiv\sum_{n=1}^{\infty}2^{-n}|A_{n}|^{-1}\sum_{x\in A(n)}|Jg_{n,x}-J_{p}g_{n,x}|
\]
\[
\leq\sum_{n=1}^{m}2^{-n}\sum_{x\in A(n)}|Jg_{n,x}-J_{p}g_{n,x}|+2^{-m}.
\]
We can first fix $m\geq1$ so large that $2^{-m}<\frac{1}{2}\varepsilon$.
Then, for sufficiently large $p\geq1,$ the  last sum is also less
than $\frac{1}{2}\varepsilon$, whence $\rho_{Dist,\xi}(J,J_{p})<\varepsilon$.
Since $\varepsilon>0$ is arbitrary, we have $\rho_{Dist,\xi}(J,J_{p})\rightarrow0$.

3. Suppose $\rho_{Dist,\xi}(J_{p},J)\rightarrow0$. Then Assertion
1 implies that $J_{p}f\rightarrow Jf$ for each $f\in C(S)$. Hence
$J_{p}\Rightarrow J$, thanks to Lemma \ref{Lem. locally compact S, Ip=00003D>I iff Ip(f)->I(f) for each f in C(S)}.
Conversely, suppose $J_{p}f\rightarrow Jf$ for each $f\in C(S)$.
Then, in particular,\emph{ $J_{p}g_{n,x}\rightarrow Jg_{n,x}$ }as
$p\rightarrow\infty$, for each\emph{ $x\in A_{n}$, }for each\emph{
$n\geq1$}. Hence $\rho_{Dist,\xi}(J_{p},J)\rightarrow0$ by Assertion
2. Applying to the special case where $J_{p}=J'$ for each $p\geq1$,
we obtain $\rho_{Dist,\xi}(J',J)=0$ iff $J=J'$.

4. Symmetry and the triangle inequality required for a metric follow
trivially from the defining equality \ref{eq:temp-119}. Hence $\rho_{Dist,\xi}$
is a metric.
\end{proof}
From the defining equality \ref{eq:temp-119}, we have $\rho_{Dist,\xi}(J,J')\leq1$
for each $J,J'\in\widehat{J}(S,d)$. Hence the metric space $(\widehat{J}(S,d),\rho_{Dist,\xi})$
is bounded. It is not necessarily complete. An easy counterexample
is by taking $S\equiv R$ with the Euclidean metric, and taking $J_{p}$
to be the point mass at $p$ for each $p\geq0$. In other words $J_{p}f\equiv f(p)$
for each $f\in C(R)$. Then $\rho_{Dist,\xi}(J_{p},J_{q})\rightarrow0$
as $p,q\rightarrow\infty$. On the other hand $J_{p}f\rightarrow0$
for each $f\in C(R)$. Hence if $\rho_{Dist,\xi}(J_{p},J)\rightarrow0$
for some $J\in\widehat{J}(S,d)$, then $Jf=0$ for each $f\in C(R)$,
and so $J=0$, contradicting the condition for $J$ to be a distribution
and an integration. The obvious problem here is that the mass of the
distributions $J_{p}$ escapes to infinity as $p\rightarrow\infty$.
The notion of tightness, defined next for a subfamily of $\widehat{J}(S,d)$,
is to prevent this from happening.
\begin{defn}
\label{Def. Modulus of Tightness} \textbf{(Tightness).} Suppose the
metric space $(S,d)$ is locally compact. Let $\beta:(0,\infty)\rightarrow[0,\infty)$
be an operation. Let $\overline{J}$ be a subfamily of $\widehat{J}(S,d)$,
such that, for each $\varepsilon>0$ and for each $J\in\overline{J}$,
we have $P_{J}(d(\cdot,x_{\circ})>a)<\varepsilon$ for each $a>\beta(\varepsilon)$,
where $P_{J}$ is the probability function of the distribution $J$.
Then we say the subfamily $\overline{J}$ is \index{tight family of distributions}\emph{tight},
with $\beta$ as a \emph{modulus of tightness}\index{modulus of tightness}
relative to the reference point $x_{\circ}$. We say that a distribution
$J$ has modulus of tightness $\beta$ if the singleton family $\{J\}$
has modulus of tightness $\beta$.

A family $M$ of r.v.'s with values in the locally compact metric
space $(S,d)$, not necessarily on the same probability space, is
said to be \index{tight family of r.v.'s}\emph{tight,} with modulus
of tightness $\beta$, if the family $\{E_{X}:X\in M\}$ is tight
with modulus of tightness $\beta$. We will say that a r.v. $X$ has
modulus of tightness $\beta$ if the singleton$\{X\}$ family has
modulus of tightness $\beta$. $\square$
\end{defn}
We emphasize that we have defined tightness of a subfamily $\overline{J}$
of $\widehat{J}(S,d)$ only when the metric space $(S,d)$ is locally
compact, even as weak convergence in $\widehat{J}(S,d)$ is defined
for the more general case of any complete metric space $(S,d)$.

Note that, according to Proposition \ref{Prop. Each distribution on complete (S,d) is I_X},
$d(\cdot,x_{\circ})$ is a r.r.v. relative to each distribution $J$.
Hence, given each $J\in\overline{J}$, the set $(d(\cdot,x_{\circ})>a)$
is integrable relative to $J$ for all but countably many $a>0$.
Therefore the probability $P_{J}(d(\cdot,x_{\circ})>a)$ makes sense
for all but countably many $a>0$. However, the countable exceptional
set of values of $a$ depends on $J$.

A modulus of tightness for a family $M$ of r.v.'s gives the uniform
rate of convergence $P(d(x_{\circ},X)>a))\rightarrow0$ as $a\rightarrow\infty$,
independent of $X\in M$, where the probability function $P$ and
the corresponding expectation $E$ are specific to $X$. This is analogous
to a modulus of uniform integrability for a family $G$ of integrable
r.r.v.'s, which gives the rate of convergence $E(|X|;|X|>a)\rightarrow0$
as $a\rightarrow\infty$, independent of $X\in G$.

The next lemma will be convenient.
\begin{lem}
\label{Lem. Tightness of rrv from bound of E|X|**p} \textbf{\emph{(A
family of r.r.v.'s bounded in $L_{p}$ is tight).}} Let $p>0$ be
arbitrary. Let $M$ be a family of r.r.v.'s such that $E|X|^{p}\leq b$
for each $X\in M$, for some $b\geq0$. Then the family $M$ is tight,
with a modulus of tightness $\beta$ relative to $0\in R$ defined
by $\beta(\varepsilon)\equiv b^{\frac{1}{p}}\varepsilon^{-\frac{1}{p}}$
for each $\varepsilon>0$.
\end{lem}
\begin{proof}
Let $X\in M$ be arbitrary. Let $\varepsilon>0$ be arbitrary. Then,
for each $a>\beta(\varepsilon)\equiv b^{\frac{1}{p}}\varepsilon^{-\frac{1}{p}},$
we have
\[
P(|X|>a)=P(|X|^{p}>a^{p})\leq a^{-p}E|X|^{p}\leq a^{-p}b<\varepsilon,
\]
where the first inequality is Chebychev's, and the second is by the
definition of the constant $b$ in the hypothesis. Thus $X$ has the
operation $\beta$ as a modulus of tightness relative to $0\in R$. 
\end{proof}
If a family $\overline{J}$ of distributions is tight relative to
a reference point $x_{\circ}$, then it is tight relative to any other
reference point $x'_{0}$, thanks to the triangle inequality. Intuitively,
tightness limits the escape of mass to infinity as we go through distributions
in $\overline{J}$. Therefore a tight family of distributions remains
so after a finite-distance shift of the reference point.
\begin{prop}
\label{Prop.Weak convergence with tightness is complete} \textbf{\emph{(Tightness
and convergence of a sequence of distributions at each member of $C(S)$
implies weak convergence to some distribution).}} Suppose the metric
space $(S,d)$ is locally compact. Let $\{J_{n}:n\geq1\}$ be a tight
family of distributions, with a modulus of tightness $\beta$ relative
to the reference point $x_{\circ}$. 

Suppose $J(f)\equiv\lim_{n\rightarrow\infty}J_{n}(f)$ exists for
each $f\in C(S)$. Then $J$ is a distribution, and $J_{n}\Rightarrow J$.
Moreover, $J$ has the modulus of tightness $\beta+2$. 
\end{prop}
\begin{proof}
Clearly $J$ is a linear function on $C(S)$. Suppose $f\in C(S)$
is such that $Jf>0$. Then, in view of the convergence in the hypothesis,
there exists $n\geq1$ such that $J_{n}f>0$. Since $J_{n}$ is an
integration, there exists $x\in S$ such that $f(x)>0$. We have thus
verified condition (ii) in Definition \ref{Def. integration on loc compact space}
for $J$. 

Next let $\varepsilon\in(0,1)$ be arbitrary, and take any $a>\beta(\varepsilon)$.
Then $P_{n}(d(\cdot,x_{\circ})>a)<\varepsilon$ for each $n\geq1$,
where $P_{n}\equiv P_{J(n)}$ is the probability function for $J_{n}$.
Define $h_{k}\equiv1\wedge(1+k-d(\cdot,x_{\circ}))_{+}\in C(S)$ for
each $k\geq1$. Take any $m\equiv m(\varepsilon,\beta)\in(a,a+2)$.
Then $h_{m}\geq1_{(d(\cdot,x(\circ))\leq a)}$, whence
\begin{equation}
J_{n}h_{m}\geq P_{n}(d(\cdot,x_{\circ})\leq a)>1-\varepsilon\label{eq:temp-486}
\end{equation}
for each each $n\geq1$. By hypothesis, $J_{n}h_{m}\rightarrow Jh_{m}$
as $n\rightarrow\infty$. Inequality \ref{eq:temp-486} therefore
yields
\begin{equation}
Jh_{m}\geq1-\varepsilon>0.\label{eq:temp-281-1}
\end{equation}
We have thus verified also condition (i) in Definition \ref{Def. integration on loc compact space}
for $J$ to be an integration on $(S,d)$. Therefore, by Proposition
\ref{Prop. Integration on S is  an integration space}, $(S,C(S),J)$
is an integration space. At the same time, inequality \ref{eq:temp-281-1}
implies that $Jh_{m}\uparrow1$. We conclude that $J$ is a distribution.
Since $J_{n}f\rightarrow Jf$ for each $f\in C(S)$ by hypothesis,
Lemma \ref{Lem. locally compact S, Ip=00003D>I iff Ip(f)->I(f) for each f in C(S)}
implies that $J_{n}\Rightarrow J$. 

Now note that inequality \ref{eq:temp-281-1} implies that 
\begin{equation}
P_{J}(d(\cdot,x_{\circ})\leq a+2)=J1_{(d(\cdot,x_{\circ})\leq a+2)}\geq Jh_{m}\geq1-\varepsilon>0,\label{eq:temp-281-1-1}
\end{equation}
where $a>\beta(\varepsilon)$ is arbitrary. Thus $J$ is tight with
the modulus of tightness $\beta+2$.
\end{proof}
\begin{cor}
\label{Cor. Tight Cauchy sequence of distribtuions converges} \textbf{\emph{(A
tight $\rho_{Dist,\xi}$-Cauchy sequence of distributions converges).}}
Let $\xi$ be a binary approximation of a locally compact metric space
$(S,d)$ relative to a reference point $x_{\circ}\in S$. Let $\rho_{Dist,\xi}$
be the distribution metric on the space $\widehat{J}(S,d)$ of distributions,
determined by $\xi$. Suppose the subfamily $\{J_{n}:n\geq1\}\subset\widehat{J}(S,d)$
of distributions is tight, with a modulus of tightness $\beta$ relative
to $x_{\circ}$. 

If $\rho_{Dist,\xi}(J_{n},J_{m})\rightarrow0$ as $n,m\rightarrow\infty$.
Then $J_{n}\Rightarrow J$ and $\rho_{Dist,\xi}(J_{n},J)\rightarrow0$,
for some $J\in\widehat{J}(S,d)$ with the modulus of tightness $\beta+2$.
\end{cor}
\begin{proof}
Suppose $\rho_{Dist,\xi}(J_{n},J_{m})\rightarrow0$ as $n,m\rightarrow\infty$.
Let $f\in C(S)$ be arbitrary. We will prove that $J(f)\equiv\lim_{n\rightarrow\infty}J_{n}(f)$
exists. Let $\varepsilon>0$ be arbitrary. Then there exists $a>\beta(\varepsilon)$
such that $P_{n}(d(\cdot,x_{\circ})>a)<\varepsilon$ for each $n\geq1$,
where $P_{n}\equiv P_{J(n)}$ is the probability function for $J_{n}$.
Let $k\geq1$ be so large that $k\geq a$, and recall that
\[
h_{k}\equiv1\wedge(1+k-d(\cdot,x_{\circ}))_{+}.
\]
Then 
\[
J_{n}h_{k}\geq P_{n}(d(\cdot,x_{\circ})\leq a)>1-\varepsilon
\]
for each each $n\geq1$. At the same time $fh_{k}\in C(S)$. Hence,
since $\rho_{Dist,\xi}(J_{n},J_{m})\rightarrow0$, implies that $(J_{n}fh_{k})_{n=1,2,\cdots}$
is a Cauchy sequence of real numbers, according to Assertion 1 of
Proposition \ref{Prop. rho_xi convergence=00003D Weak Convergence}.
Hence $Jfh_{k}\equiv\lim_{n\rightarrow\infty}J_{n}(fh_{k})$ exists.
Consequently, 
\[
|J_{n}f-J_{m}f|\leq|J_{n}f-J_{n}fh_{k}|+|J_{n}fh_{k}-J_{m}fh_{k}|+|J_{m}fh_{k}-J_{m}f|
\]
\[
\leq|1-J_{n}h_{k}|+|J_{n}fh_{k}-J_{m}fh_{k}|+|J_{m}h_{k}-1|
\]
\[
\leq\varepsilon+|J_{n}fh_{k}-J_{m}fh_{k}|+\varepsilon<\varepsilon+\varepsilon+\varepsilon
\]
for sufficiently large $n,m\geq1$. Since $\varepsilon>0$ is arbitrary,
we conclude that $J(f)\equiv\lim_{n\rightarrow\infty}J_{n}f$ exists
for each $f\in C(S)$. By Proposition \ref{Prop.Weak convergence with tightness is complete},
$J$ is a distribution with the modulus of tightness $\beta+2$, and
$J_{n}\Rightarrow J$. Proposition \ref{Prop. rho_xi convergence=00003D Weak Convergence}
then implies that $\rho_{Dist,\xi}(J_{n},J)\rightarrow0$.
\end{proof}
\begin{prop}
\label{Prop. Sequence converging in distribution is tight} \textbf{\emph{(A
weakly convergent sequence of distributions on a locally compact metric
space is tight).}} Suppose the metric space $(S,d)$ is locally compact.
Let $J,J_{n}$ be distributions for each $n\geq1$. Suppose $J_{n}\Rightarrow J$.
Then the family $\{J,J_{1},J_{2}\cdots\}$ is tight. In particular,
any finite family of distributions on $S$ is tight, and any finite
family of r.v.'s with values in $S$ is tight.
\end{prop}
\begin{proof}
For each $n\geq1$ write $P$ and $P_{n}$ for $P_{J}$ and $P_{J(n)}$
respectively. Since $J$ is a distribution, we have $P(d(\cdot,x_{\circ})>a)\rightarrow0$
as $a\rightarrow\infty$. Thus any family consisting of a single distribution
$J$ is tight. Let $\beta_{0}$ be a modulus of tightness of $\{J\}$
with reference to $x_{\circ}$, and, for each $k\geq1$, let $\beta_{k}$
be a modulus of tightness of $\{J_{k}\}$ with reference to $x_{\circ}$.
Let $\varepsilon>0$ be arbitrary. Let $a>\beta_{0}(\frac{\varepsilon}{2})$
and define $f\equiv1\wedge(a+1-d(\cdot,x_{\circ}))_{+}$. Then $f\in C(S)$
with 
\[
1_{(d(\cdot,x(\circ))>a+1)}\leq1-f\leq1_{(d(\cdot,x(\circ))>a)}.
\]
Hence $1-Jf\leq P(d(\cdot,x_{\circ})>a)<\frac{\varepsilon}{2}$. By
hypothesis, we have $J_{n}\Rightarrow J$. Hence there exists $m\geq1$
so large that $|J_{n}f-Jf|<\frac{\varepsilon}{2}$ for each $n>m$.
Consequently
\[
P_{n}(d(\cdot,x_{\circ})>a+1)\leq1-J_{n}f<1-Jf+\frac{\varepsilon}{2}<\varepsilon
\]
 for each $n>m$. Define $\beta(\varepsilon)\equiv(a+1)\vee\beta_{1}(\varepsilon)\vee\cdots\vee\beta_{m}(\varepsilon)$.
Then, for each $a'>\beta(\varepsilon)$ we have 

(i) $P(d(\cdot,x_{\circ})>a')\leq P(d(\cdot,x_{\circ})>a)<\frac{\varepsilon}{2}$, 

(ii) $P_{n}(d(\cdot,x_{\circ})>a')\leq P_{n}(d(\cdot,x_{\circ})>a+1)<\varepsilon$
for each $n>m$, and 

(iii) $a'>\beta_{n}(\varepsilon)$ and so $P_{n}(d(\cdot,x_{\circ})>a')\leq\varepsilon$
for each $n=1,\cdots,m$. 

Since $\varepsilon>0$ is arbitrary,  the family $\{J,J_{1},J_{2}\cdots\}$
is tight. 
\end{proof}
The next proposition provides some alternative characterization of
weak convergence in the case of locally compact $(S,d)$. 
\begin{prop}
\label{Prop. Modulus of continuity of J->Jf for fixed Lipshitz f}
\textbf{\emph{(Modulus of continuity of the function $J\rightarrow Jf$
for functions $f$ with fixed Lipschitz constant).}} Suppose $(S,d)$
is locally compact, with a reference point $x_{\circ}$. Let $\xi\equiv(A_{n})_{n=1,2,\cdots}$
be a binary approximation of $(S,d)$ relative to $x_{\circ}$, with
a corresponding modulus of local compactness $\left\Vert \xi\right\Vert \equiv(|A_{n}|)_{n=1,2,\cdots}$
of $(S,d)$. Let $\rho_{Dist,\xi}$ be the distribution metric on
the space $\widehat{J}(S,d)$ of distributions on $(S,d),$ determined
by $\xi$, as introduced in Definition \ref{Def. Distribution metric}.. 

Let $J,J',J_{p}$ be distributions on $(S,d)$, for each $p\geq1$.
Let $\beta$ be a modulus of tightness of $\{J,J'\}$ relative to
$x_{\circ}$. Then the following holds. 
\end{prop}
\begin{enumerate}
\item \emph{Let $f\in C(S,d)$ be arbitrary with $|f|\leq1$ and with modulus
of continuity $\delta_{f}$. Then, for each $\varepsilon>0$, there
exists $\widetilde{\Delta}(\varepsilon,\delta_{f},\beta,\left\Vert \xi\right\Vert )>0$
such that if $\rho_{Dist,\xi}(J,J')<\widetilde{\Delta}(\varepsilon,\delta_{f},\beta,\left\Vert \xi\right\Vert )$
then $|Jf-J'f|<\varepsilon$.}
\item \emph{The following three conditions are equivalent: }\emph{\noun{(i)}}\emph{
$J_{p}f\rightarrow Jf$ for each Lipschitz continuous $f\in C(S)$,
}(ii)\emph{ $J_{p}\Rightarrow J$, and }(iii)\emph{ $J_{p}f\rightarrow Jf$
for each Lipschitz continuous $f$ which is bounded.}
\end{enumerate}
\begin{proof}
By Definition \ref{Def. Binary approximationt and Modulus of local compactness},
we have
\begin{equation}
(d(\cdot,x_{\circ})\leq2^{n})\subset\bigcup_{x\in A(n)}(d(\cdot,x)\leq2^{-n})\label{eq:temp-505-1}
\end{equation}
and
\begin{equation}
\bigcup_{x\in A(n)}(d(\cdot,x)\leq2^{-n+1})\subset(d(\cdot,x_{\circ})\leq2^{n+1})\label{eq:temp-502-1}
\end{equation}
for each $n\geq1$. 

1. Let 
\[
b\equiv[1+\beta(\frac{\varepsilon}{4})]_{1}.
\]
Write $h\equiv1\wedge(b-d(\cdot,x_{\circ}))_{+}\in C(S)$. Then $h$
and $fh$ have support $(d(\cdot,x_{\circ})\leq b)$, and $h=1$ on
$(d(\cdot,x_{\circ})\leq b-1)$. 

Moreover, since $h$ has Lipschitz constants $1$, the function $fh$
has a modulus of continuity $\delta_{fh}$ defined by $\delta_{fh}(\alpha)\equiv\frac{\alpha}{2}\wedge\delta_{f}(\frac{\alpha}{2})$.
Hence, by Proposition \ref{Prop. rho_xi convergence=00003D Weak Convergence},
there exists $\delta_{\widehat{J}}(\frac{\varepsilon}{2})\equiv\delta_{\widehat{J}}(\frac{\varepsilon}{2},\delta_{fh},b,\left\Vert \xi\right\Vert )>0$
such that if $\rho_{Dist,\xi}(J,J')<\delta_{\widehat{J}}(\frac{\varepsilon}{2})$
then
\begin{equation}
|Jfh-J'fh|<\frac{\varepsilon}{2}.\label{eq:temp-146}
\end{equation}
More precisely, according to said proposition, we can let
\[
n\equiv[0\vee(1-\log_{2}\delta_{fh}(\frac{\varepsilon}{3}))\vee\log_{2}b]_{1}
\]
\[
\equiv[0\vee\log_{2}(\frac{\varepsilon}{6}\wedge\delta_{f}(\frac{\varepsilon}{6}))\vee\log_{2}b]_{1}
\]
and
\[
\delta_{\widehat{J}}(\frac{\varepsilon}{2},\delta_{fh},b,\left\Vert \xi\right\Vert )\equiv\frac{1}{6}2^{-n}|A_{n}|^{-1}\varepsilon.
\]

Now define 
\[
\widetilde{\Delta}(\varepsilon)\equiv\widetilde{\Delta}(\varepsilon,\delta_{f},\beta,\left\Vert \xi\right\Vert )\equiv\delta_{\widehat{J}}(\frac{\varepsilon}{2},\delta_{fh},b,\left\Vert \xi\right\Vert )\equiv\frac{1}{6}2^{-n}|A_{n}|^{-1}\varepsilon.
\]
Suppose $\rho_{Dist,\xi}(J,J')<\widetilde{\Delta}(\varepsilon)$.
We need to prove that\emph{ $|Jf-J'f|<\varepsilon$.} To that end,
note that, since $J,J'$ have tightness modulus $\beta$, and since
$1-h=0$ on $(d(\cdot,x_{\circ})\leq b-1)$ where $b-1>\beta(\frac{\varepsilon}{4})$,
we have $J(1-h)\leq\frac{\varepsilon}{4}$ and $J'(1-h)\leq\frac{\varepsilon}{4}$.
Consequently,
\begin{equation}
|Jf-Jfh|=|Jf(1-h)|\leq J(1-h)\leq\frac{\varepsilon}{4}.\label{eq:temp-160}
\end{equation}
Similarly,
\begin{equation}
|J'f-J'fh|\leq\frac{\varepsilon}{4}.\label{eq:temp-161}
\end{equation}
Combining inequalities \ref{eq:temp-146}, \ref{eq:temp-160}, and
\ref{eq:temp-161}, we obtain
\[
|Jf-J'f|<\frac{\varepsilon}{4}+\frac{\varepsilon}{2}+\frac{\varepsilon}{4}\varepsilon=\varepsilon.
\]
Assertion 1 is thus proved.

2. We need to prove that Conditions (i-iii) are equivalent. To that
end, first suppose (i) $J_{p}f\rightarrow Jf$ for each Lipschitz
continuous $f\in C(S)$. Let 
\[
\pi\equiv(\{g_{n,x}:x\in A_{n}\})_{n=1,2,\cdots}
\]
be the partition of unity of $(S,d)$ determined by $\xi$. Then,
for each $n\geq1$ and each $x\in A_{n}$, we have $J_{p}g_{n,x}\rightarrow Jg_{n,x}$
as $p\rightarrow\infty$, because $g_{n,x}\in C(S)$ is Lipschitz
continuous by Proposition \ref{Prop. Properties of parittion of unity-1}.
Hence $\rho_{Dist,\xi}(J_{p},J)\rightarrow0$ and $J_{p}\Rightarrow J$
by Proposition \ref{Prop. rho_xi convergence=00003D Weak Convergence}.
Thus we have proved that Condition (i) implies Condition (ii). 

Suppose next that $J_{p}\Rightarrow J$. Then $J_{p}f\rightarrow Jf$
for each $f\in C(S)$. Hence, since $(S,d)$ is locally compact, we
have $\rho_{Dist,\xi}(J_{p},J)\rightarrow0$ by Proposition \ref{Prop. rho_xi convergence=00003D Weak Convergence}.
Separately, in view of Proposition \ref{Prop. Sequence converging in distribution is tight},
the family $\{J,J_{1},J_{2},\cdots\}$ is tight, with some modulus
of tightness $\beta$. Let $f\in C(S)$ be Lipschitz continuous. We
need to prove that $J_{p}f\rightarrow Jf$. By linearity, we may assume
that $|f|\leq1$, whence $J_{p}f\rightarrow Jf$ by Assertion 1. Thus
Condition (ii) implies Condition (iii). 

Finally, Condition (iii) trivially implies Condition (i) . Assertion
2 is proved.
\end{proof}

\section{Probability Density Functions and Distribution Functions}

Useful distributions can be obtained by using integrable functions
as density functions. The Riemann-Stieljes integration gives rise
to other examples of distributions on $R$. This section makes these
terms precise for later reference.
\begin{defn}
\label{Def. p.d.f.  on completion of (S,C(S),I)} \textbf{(probability
density function).} Let $I$ be an integration on a locally compact
metric space $(S,d)$, and let $(S,\Lambda,I)$ denote the completion
of the integration space $(S,C(S),I)$. Let $g\in\Lambda$ be arbitrary,
with $g\geq0$ and $Ig=1$. Then $g$ will be called a \emph{probability
density function}, \index{probability density function} or $p.d.f.$
for short, on the integration space $(S,\Lambda,I)$. Define $I_{g}h\equiv Igh$
for each $h\in C(S)$. Then $(S,C(S),I_{g})$ is an integration space,
with a completion $(\Omega,\Lambda_{g},I_{g})$ which is a probability
space.

Suppose $X$ is a r.v. with values in $S$ such that $X$ induces
the distribution $I_{g}$. In other words, Suppose $E_{X}=I_{g}$.
Then $X$ is said to have the p.d.f. $g$. $\square$

Frequently used p.d.f.'s are defined on $(S,\Lambda,I)\equiv(R^{n},\Lambda,\int\cdot dx)$,
the $n$-dimensional Euclidean space equipped with the Lebesgue integral,
and on $(S,\Lambda,I)\equiv(\{1,2,\cdots\},\Lambda,I)$ with the counting
measure $I$ defined by $Ig\equiv\sum_{n=1}^{\infty}g(n)$ for each
$g\in C(S)$. 
\end{defn}
\begin{prop}
\label{Prop. p.d.f.  basics} \textbf{\emph{(Integrable functions
relative to a p.d.f.) }}Use the notations of Definition \ref{Def. p.d.f.  on completion of (S,C(S),I)}.
Let $g$ be a p.d.f. on $(S,\Lambda,I)$. Let $f,h$ be an arbitrary
$\mathrm{measurable}$ function on $(S,\Lambda,I)$ such that $hg\in\Lambda$.
Then (i) $h\in\Lambda_{g}$ and $I_{g}h=Ihg$, and (ii) $f$ is $\mathrm{measurable}$
on $(S,\Lambda_{g},I_{g})$. 
\end{prop}
\begin{proof}
1. First suppose that $h\in\Lambda$ and $|h|\leq a$ for some $a>0$.
Let $k\geq1$ be arbitrary. Then there exists $h_{k}\in C(S)$ such
that $I|h_{k}-h|<\frac{1}{k}$. By replacing $h_{k}$ with $-a\vee h_{k}\wedge a$,
we may assume that $|h_{k}|\leq a$. It follows that $\bar{h}\equiv\lim_{k\rightarrow\infty}h_{k}\in\Lambda$,
that $h=\bar{h}$ on the full subset $D\equiv domain(\bar{h})$ of
\emph{$(S,\Lambda,I)$.} Hence $I_{g}|h_{k}-h_{j}|\equiv I|h_{k}-h_{j}|g\rightarrow0$
as $k,j\rightarrow\infty$ by the Dominated Convergence Theorem. Since
$h_{k}\in C(S)\subset\Lambda_{g}$ for each $k\geq0$, we conclude
that $\bar{h}\in\Lambda_{g}$ and that 
\[
I_{g}\bar{h}=\lim_{k\rightarrow\infty}I_{g}h_{k}=\lim_{k\rightarrow\infty}Ih_{k}g=I\bar{h}g.
\]
Hence $D\equiv domain(\bar{h})$ is a full set also of $(S,\Lambda_{g},I_{g})$.
Since $h=\bar{h}$ on $D$, we conclude that \emph{$h\in\Lambda_{g}$
and $I_{g}h=I_{g}\bar{h}=Ihg$. }

2. Next, suppose \emph{$h\geq0$.} Let $x_{\circ}\in S$ be an arbitrary,
but fixed, reference point. Let $m\geq1$ be arbitrary. Then $(m\wedge h)f_{m}\in\Lambda$
where $f_{m}\equiv1\wedge(m-d(\cdot,x_{\circ}))_{+}\in C(S)$. Hence,
by Step 1 of this proof, we have $(m\wedge h)f_{m}\in\Lambda_{g}$
and 
\[
I_{g}(m\wedge h)f_{m}\equiv I(m\wedge h)f_{m}g.
\]
Since $(m\wedge h)f_{m}\uparrow h$, the Dominated Convergence Theorem
implies that $I(m\wedge h)f_{m}g\uparrow Ihg$. In turn, the Monotone
Convergence Theorem then implies that $h\in\Lambda_{g}$ with $I_{g}h=Ihg$. 

3. Finally, let $h$ be an arbitrary nonnegative $\mathrm{measurable}$
function on $(S,\Lambda,I)$ such that $hg\in\Lambda$. Step 2 above
implies $h_{+},|h|\in\Lambda_{g}$. Hence, since $h=2h_{+}-|h|$,
we have \emph{$h\in\Lambda_{g}$} by linearity, with
\[
I_{g}h=2I_{g}h_{+}-I_{g}|h|=I(2h_{+}-|h|)g=Ihg.
\]
Assertion (i) is proved.

4. By Assertion 1, we have $f_{k}\equiv-k\vee f\wedge k\in\Lambda_{g}$
and $1\wedge|f_{k}-f|\in\Lambda_{g}$ for each $k\geq1$. Moreover,
$I_{g}1_{(1\wedge|f_{k}-f|)}=I1_{(1\wedge|f_{k}-f|)}g\rightarrow0$
by the Monotone Convergence Theorem. Hence $f_{k}\rightarrow f$ in
probability on $(S,\Lambda_{g},I_{g})$, whence $f$ is $\mathrm{measurable}$
on $(S,\Lambda_{g},I_{g})$. 
\end{proof}
\begin{prop}
\textbf{\emph{\label{Prop. X has pdf g =00003D> g(X)>0  a.s.}(p.d.f.
of a r.v. $X$ is a.s. positive at $X$).}} Let $X:(\Omega,L,E)\rightarrow(S,d)$
be a r.v. with the p.d.f. $g$ on $(S,\Lambda,I)$. In other words,
$(S,\Lambda_{g},E_{X})=(S,\Lambda_{g},I_{g})$. Then $P_{g}(g\leq\varepsilon)\equiv I_{g}1_{(g\leq\varepsilon)}\rightarrow0$
as $\varepsilon\rightarrow0$. Moreover, $P(g(X)\leq\varepsilon)\rightarrow0$
as $\varepsilon\rightarrow0$. Consequently, \textup{$g(X)>0$ a.s.
on $(\Omega,L,E)$, and $g>0$ a.s. on $(S,\Lambda_{g},I_{g})$.}
\end{prop}
\begin{proof}
By Proposition \ref{Prop. p.d.f.  basics}, $g(X)$ is $\mathrm{measurable}$
on $\Omega$ and $g$ is $\mathrm{measurable}$ on $(S,\Lambda_{g},I_{g})$.
Then
\[
P_{g}(g\leq\varepsilon)\equiv I_{g}1_{(g\leq\varepsilon)}=Ig1_{(g\leq\varepsilon)}\leq Ig\wedge\varepsilon\rightarrow0
\]
as $\varepsilon\rightarrow0$, because $g\in\Lambda$.Consequently,
\[
P(g(X)\leq\varepsilon)=E1_{(g(X)\leq\varepsilon)}=E1_{(g\leq\varepsilon)}(X)\equiv E_{X}1_{(g\leq\varepsilon)}=I_{g}1_{(g\leq\varepsilon)}\rightarrow0
\]
as $\varepsilon\rightarrow0$, where the last inequality is because
$(S,\Lambda_{g},E_{X})=(S,\Lambda_{g},I_{g})$ by hypothesis.
\end{proof}
Distributions on $R$ can be studied in terms of their corresponding
distribution functions, as introduced earlier in Definition \ref{Def. Distribution Func}
and specialized to probability distribution functions. 

Recall the convention that if $F$ is a function, then we write $F(t)$
only with the implicit assumption that $t\in domain(F)$.
\begin{defn}
\label{Def. PDF} \textbf{(Probability Distribution Functions).} Suppose
$F$ is a distribution function on $R$ satisfying the following conditions:
(i) $F(t)\rightarrow0$ as $t\rightarrow-\infty$, and $F(t)\rightarrow1$
as $t\rightarrow\infty$, (ii) for each $t\in domain(F)$, the left
limit $\lim_{r<t;r\rightarrow t}F(s)$ exists, (iii) for each $t\in domain(F)$,
the right limit $\lim_{s>t;s\rightarrow t}F(s)$ exists and is equal
to $F(t)$, (iv) $domain(F)$ contains the metric complement $A_{c}$
of some countable subset $A$ of $R$, and (v) if $t\in R$ is such
that both the above-defined left- and right limits exist, then $t\in domain(F)$.
Then $F$ is called a \index{probability distribution function}\emph{probability
distribution function,} or a \index{P.D.F.}P.D.F. for abbreviation\emph{.
}A point $t\in domain(F)$ is called a regular point of $F$. A point
$t\in domain(F)$ at which the above-defined left- and right limits
are equal is called a \emph{continuity point}\index{continuity point of a P.D.F.}
of $F$.

Suppose $X$ is a r.r.v. on a probability space $(\Omega,L,E)$. Let
$F_{X}$ be the function defined by (i) $domain(F_{X})\equiv\{t\in R:t\mbox{ is a regular point }$
$\mbox{of }X\}$, and (ii) $F_{X}(t)\equiv P(X\leq t)$ for each $t\in domain(F_{X})$.
Then $F_{X}$ is called the \index{P.D.F. of a r.r.v.} P.D.F. of
$X$. $\square$
\end{defn}
Recall in the following that $\int\cdot dF$ denotes the Riemann-Stieljes
integration relative to a distribution function $F$ on $R$\emph{.}
\begin{prop}
\label{Prop.  F_X  is indeed a PDF} \textbf{\emph{($F_{X}$ is indeed
a P.D.F.)}} Let $X$ be a r.r.v. on a probability space $(\Omega,L,E)$
with $F_{X}$ as in Definition \ref{Def. PDF}. Let $E_{X}$ denote
the distribution induced on $R$ by $X$. Then the following holds.
\end{prop}
\begin{enumerate}
\item \emph{$F_{X}$ is a P.D.F.}
\item \emph{$\int\cdot dF_{X}=E_{X}$. }
\end{enumerate}
\begin{proof}
For abbreviation, write $J\equiv E_{X}$ and $F\equiv F_{X}$, and
write $P$ for the probability function associated to $E$.

1. We are to verify conditions (i) through (v) in Definition \ref{Def. PDF}
for $F$. Condition (i) holds because $P(X\leq t)\rightarrow0$ as
$t\rightarrow-\infty$ and $P(X\leq t)=1-P(X>t)\rightarrow1$ as $t\rightarrow\infty$,
by the definition of a $\mathrm{measurable}$ function. Next consider
any $t\in domain(F)$. Then $t$ is a regular point of $X$, by the
definition of $F_{X}$. Hence there exists a sequence $(s_{n})_{n=1,2,\cdots}$
of real numbers decreasing to $t$ such that $(X\leq s_{n})$ is integrable
for each $n\geq1$ and such that $\lim_{n\rightarrow\infty}P(X\leq s_{n})$
exists. Since $P(X\leq s_{n+2})\leq F(s)\leq P(X\leq s_{n})$ for
each $s\in(s_{n+2},s_{n})$ and for each $n\geq1$, we see that $\lim_{s>t;s\rightarrow t}F(s)$
exists. Similarly $\lim_{r<t;r\rightarrow t}F(s)$ exists. Moreover,
by Proposition \ref{Prop. Sequences of intgrable sets}, we have $\lim_{s>t;s\rightarrow t}F(s)=\lim_{n\rightarrow\infty}P(X\leq s_{n})=P(X\leq t)\equiv F(t)$.
Conditions (ii) and (iii) in Definition \ref{Def. PDF} have thus
been verified. Condition (iv) in Definition \ref{Def. PDF} follows
from Assertion 1 of Proposition \ref{Prop.  Countable Exceptional Pts for Meas X int A}.
Condition (v) remains. Suppose $t\in R$ is such that both $\lim_{r<t;r\rightarrow t}F(r)$
and $\lim_{s>t;s\rightarrow t}F(s)$ exist. Then there exists a sequence
$(s_{n})_{n=1,2,\cdots}$ in $domain(F)$ decreasing to $t$ such
that $F(s_{n})$ converges. This implies that $(X\leq s_{n})$ is
an integrable set, and that $P(X\leq s_{n})$ converges. Hence $(X>s_{n})$
is an integrable set, and $P(X>s_{n})$ converges. Similarly, there
exists a sequence $(r_{n})_{n=1,2,\cdots}$ increasing to $t$ such
that $(X>r_{n})$ is an integrable set and $P(X>r_{n})$ converges.
We have thus verified the conditions in Definition \ref{Def. Regular =000026 Continuity Pts of Measurable Func}
for $t$ to be a regular point of $X$. In other words, $t\in domain(F)$.
Condition (v) in Definition \ref{Def. PDF} have thus also been verified.
Summing up, $F\equiv F_{X}$ is a P.D.F.

2. Note that both $\int\cdot dF_{X}$ and $E_{X}$ are complete extensions
of integrations defined on $(R,C(R))$. Hence it suffices to prove
that they are equal on $C(R)$. Let $f\in C(R)$ be arbitrary. We
need to show that $\int f(x)dF_{X}(x)=E_{X}f$. Let $\varepsilon>0$
be arbitrary, and let $\delta_{f}$ be a modulus of continuity for
$f$. The Riemann-Stieljes integral $\int f(t)dF_{X}(t)$ is, by definition,
the limit of Riemann-Stieljes sums $S(t_{1},\cdots,t_{n})=\sum_{i=1}^{n}f(t_{i})(F_{X}(t_{i})-F_{X}(t_{i-1}))$
as $t_{1}\rightarrow-\infty$ and $t_{n}\rightarrow\infty$ with the
mesh of the partition $t_{1}<\cdots<t_{n}$ approaching $0$. Consider
such a Riemann-Stieljes sum where the mesh is smaller than $\delta_{f}(\varepsilon)$,
and where $[t_{1},t_{n}]$ contains a support of $f$. Then 
\[
|S(t_{1},\cdots,t_{n})-Ef(X)|=|E\sum_{i=1}^{n}(f(t_{i})-f(X))1_{(t_{i-1}<X\leq t_{i})}|\leq\varepsilon
\]
Passing to the limit, we have $\int f(t)dF_{X}(t)=Ef(X)\equiv E_{X}f$. 
\end{proof}
Proposition \ref{Prop.  F_X  is indeed a PDF} says that $F_{X}$
is a P.D.F. for each r.r.v. $X$. The next proposition gives the converse.
\begin{prop}
\label{Prop.  PDF =00003D FX} \textbf{\emph{(Basics of P.D.F.) }}\emph{The
following holds.}
\end{prop}
\begin{enumerate}
\item \emph{Let $J$ be any distribution on $R$, and let $(R,L,J)$ denote
the completion of $(R,C_{ub}(R),J)$. Then $J=\int\cdot dF_{X}$ where
$F_{X}$ is the P.D.F. of the r.r.v. $X$ on $(R,L,J)$ defined by
$X(x)\equiv x$ for $x\in R$. }
\item \emph{Let $F$ be a P.D.F. For each $t\in domain(F)$, the interval
$(-\infty,t]$ is integrable relative to $\int\cdot dF$, and $\int1_{(-\infty,t]}dF=F(t)$. }
\item \emph{If two P.D.F.'s $F$ and $F'$ equal on some dense subset $D$
of $domain(F)\cap domain(F')$, then $F=F'$.}
\item \emph{If two P.D.F.'s $F$ and $F'$ are such that $\int\cdot dF=\int\cdot dF'$,
then $F=F'$.}
\item \emph{Let $F$ be any P.D.F. Then $F=F_{X}$ for some r.r.v. $X$.}
\item \emph{Let $F$ be any P.D.F. Then all but countably many $t\in R$
are continuity points of $F$.}
\item \emph{Let} \emph{$J$ be any distribution on $R$. Then there exists
a unique P.D.F. F such that $J=\int\cdot dF$. Thus there is a bijection
between distributions on $R$ and P.D.F.'s. For that reason, we will
often abuse }term\emph{inology and refer to  $F$ as a distribution,
and write $F$ for $J$.}
\end{enumerate}
\begin{proof}
1. According to Proposition \ref{Prop. Each distribution on complete (S,d) is I_X},
$X$ is a r.v. on \emph{$(R,L,J)$}. Moreover, for each $f\in C_{ub}(R)$,
we have $f(X)=f\in C_{ub}(R)$. Hence, in view of Proposition \ref{Prop.  F_X  is indeed a PDF},
we have $Jf=Jf(X)\equiv E_{X}f=\int f(x)dF_{X}(x)$. Assertion 1 is
validated.

2. Define $J\equiv\int\cdot dF$. Consider any $t,s\in domain(F)$
with $t<s$, and any $f\in C(R)$ with $0\leq f\leq1$ such that $[t,s]$
is a support of $f$. Then $Jf\equiv\int f(x)dF(x)$ is the limit
of Riemann-Stieljes sums $\sum_{i=1}^{n}f(t_{i})(F(t_{i})-F(t_{i-1}))$,
where the sequence $t_{0}<\cdots<t_{n}$ includes the points $t,s$. 

Consider any such Riemann-Stieljes sum. If $i=0,\cdots,n$ is such
that $t_{i}<t$ or $t_{i}>s$, then $f(t_{i})=0$. We can, by excluding
such indices $i$, assume that $t_{0}=t$ and $t_{n}=s$. It follows
that the Riemann-Stieljes sums in question are bounded by $\sum_{i=1}^{n}(F(t_{i})-F(t_{i-1}))=F(s)-F(t).$
Passing to the limit, we see that $Jf\leq F(s)-F(t)$ for each $f\in C(R)$
with $0\leq f\leq1$ such that $[t,s]$ is a support of $f$. A similar
argument shows that $Jf\geq F(s)-F(t)$ for each $f\in C(R)$ with
$0\leq f\leq1$ such that $f=1$ on $[t,s]$. 

By condition (i) in Definition \ref{Def. PDF}, there exists a decreasing
sequence $(r_{k})_{k=1,2,\cdots}$ in $domain(F)$ such that $r_{1}<t$,
$r_{k}\rightarrow-\infty$, and $F(r_{k})\rightarrow0$. By condition
(iii) in Definition \ref{Def. PDF}, we have $F(s_{n})\rightarrow F(t)$
for some sequence $(s_{n})_{n=1,2,\cdots}$ such that $s_{n}\downarrow t$. 

For each $k,n\geq1$, let $f_{k,n}\in C(R)$ be defined by $f_{k,n}=1$
on $[r_{k},s_{n+1}]$, $f_{k,n}=0$ on $(-\infty,r_{k+1}]\cup[s_{n},\infty)$,
and $f_{k,n}$ is linear on $[r_{k+1},r_{k}]$ and on $[s_{n+1},s_{n}]$.
Consider any $n\geq1$ and $j>k\geq1$. Then $0\leq f_{k,n}\leq f_{j,n}\leq1$,
and $f_{j,n}-f_{k,n}$ has $[r_{j+1},r_{k}]$ as support. Therefore,
as seen earlier, $Jf_{j,n}-Jf_{k,n}\leq F(r_{k})-F(r_{j+1})\rightarrow0$
as $j\geq k\rightarrow\infty$. Hence the Monotone Convergence Theorem
implies that $f_{n}\equiv\lim_{k\rightarrow\infty}f_{k,n}$ is integrable,
with $Jf_{n}=\lim_{k\rightarrow\infty}Jf_{k,n}$. Moreover, $f_{n}=1$
on $(-\infty,s_{n+1}]$, $f_{n}=0$ on $[s_{n},\infty)$, and $f_{n}$
is linear on $[s_{n+1},s_{n}]$. 

Now consider any $m\geq n\geq1$. Then $0\leq f_{m}\leq f_{n}\leq1$,
and $f_{n}-f_{m}$ has $[t,s_{n}]$ as support. Therefore, as seen
earlier, $Jf_{n}-Jf_{m}\leq F(s_{n})-F(t)\rightarrow0$ as $m\geq n\rightarrow\infty$.
Hence, the Monotone Convergence Theorem implies that $g\equiv\lim_{n\rightarrow\infty}f_{n}$
is integrable, with $Jg=\lim_{n\rightarrow\infty}Jf_{n}$. It is evident
that on $(-\infty,t]$ we have $f_{n}=1$ for each $n\geq1$. Hence
$g$ is defined and equal to $1$ on $(-\infty,t]$. Similarly, $g$
is defined and equal to $0$ on $(t,\infty)$. 

Consider any $x\in domain(g).$ Then either $g(x)>0$ or $g(x)<1$.
Suppose $g(x)>0$. Then the assumption $x>t$ would imply $g(x)=0$,
a contradiction. Hence $x\in(-\infty,t]$ and so $g(x)=1$. On the
other hand, suppose $g(x)<1$. Then $f_{n}(x)<1$ for some $n\geq1$,
whence $x\geq s_{n+1}$ for some $n\geq1$. Hence $x\in(t,\infty)$
and so $g(x)=0$. Combining, we see that $1$ and $0$ are the only
possible values of $g$. In other words, $g$ is an integrable indicator.
Moreover $(g=1)=(-\infty,t]$ and $(g=0)=(t,\infty)$. Thus the interval
$(-\infty,t]$ is an integrable set with $1_{(-\infty,t]}=g$. 

Finally, for any $k,n\geq1$, we have $F(s_{n+1})-F(r_{k})\leq Jf_{k,n}\leq F(s_{n})-F(r_{k+1})$.
Letting $k\rightarrow\infty$, we obtain $F(s_{n+1})\leq Jf_{n}\leq F(s_{n})$
for $n\geq1$. Letting $n\rightarrow\infty$, we obtain in turn $Jg=F(t)$.
In other words $J1_{(-\infty,t]}=F(t)$. Assertion 2 is proved.

3. Consider any $t\in domain(F)$. Let $(s_{n})_{n=1,2,\cdots}$ be
a decreasing sequence in $D$ converging to $t$. By hypothesis, $F'(s_{n})=F(s_{n})$
for each $n\geq1$. At the same time $F(s_{n})\rightarrow F(t)$ since
$t\in domain(F)$. Therefore$F'(s_{n})\rightarrow F(t)$. By the monotonicity
of $F'$, it follows that $\lim_{s>t;s\rightarrow t}F'(s)=F(t)$.
Similarly, $\lim_{r<t;r\rightarrow t}F'(r)$ exists. Therefore, according
to Definition \ref{Def. PDF}, we have $t\in domain(F')$, and $F'(t)=\lim_{s>t;s\rightarrow t}F'(s)=F(t)$.
We have thus proved that $domain(F)\subset domain(F')$ and $F'=F$
on $domain(F)$. By symmetry $domain(F)=domain(F').$

4. Write $J\equiv\int\cdot dF=\int\cdot dF'$. Consider any $t\in D\equiv domain(F)\cap domain(F')$.
By assertion 2, the interval $(-\infty,t]$ is integrable relative
to $J$, with $F(t)=J1_{(-\infty,t]}=F'(t)$. Since $D$ is a dense
subset of $R$, we have $F=F'$ by Assertion 3. This proves Assertion
4.

5. Let $F$ be any P.D.F. By assertion 1, we have $\int\cdot dF=\int\cdot dF_{X}$
for some r.r.v. $X$. Therefore $F=F_{X}$ according to assertion
4. Assertion 5 is proved.

6. Let $F$ be any P.D.F. By assertion 5, $F=F_{X}$ for some r.r.v.
$X$. Hence $F(t)=F_{X}(t)\equiv P(X\leq t)$ for each regular point
$t$ of $X$. Consider any continuity point $t$ of $X$. Then, by
Definition \ref{Def. Regular =000026 Continuity Pts of Measurable Func},
we have $\lim_{n\rightarrow\infty}P(X\leq s_{n})=\lim_{n\rightarrow\infty}P(X\leq r_{n})$
for some decreasing sequence $(s_{n})$ with $s_{n}\rightarrow t$
and some increasing sequence $(r_{n})$ with $r_{n}\rightarrow t$.
Since $P(X\leq r_{n})\leq F(x)\leq P(X\leq s_{n})$ for all $x\in(r_{n},s_{n})$,
it follows that $\lim_{s>t;s\rightarrow t}F(s)$ $=\lim_{r<t;r\rightarrow t}F(r)$.
Summing up, every continuity point of $X$ is a continuity point of
$F$. By Proposition \ref{Prop.  Countable Exceptional Pts for Meas X int A},
all but countably many $t\in R$ are continuity points of $X$. Hence
all but countably many $t\in R$ are continuity points of $F$. This
validates Assertion 6.

7. Let $J$ be arbitrary. By Assertion 1, there exists a P.D.F. F
such that $J=\int\cdot dF$. Uniqueness of $F$ follows from Assertion
4. The proposition is proved.
\end{proof}

\section{The Skorokhod Representation}

In  this section, let $(S,d)$ be a locally compact metric space with
an arbitrary, but fixed, reference point $x_{\circ}\in S$. Let 
\[
(\Theta_{0},L_{0},I)\equiv([0,1],L_{0},\int\cdot dx)
\]
denote the Lebesgue integration space based on the unit interval $[0,1]$,
and let $\mu$ the corresponding Lebesgue measure. 

Given two distributions $E$ and $E'$ on the locally compact metric
space $(S,d)$, we saw in Proposition \ref{Prop. Each distribution on complete (S,d) is I_X}
that they are equal to the distributions induced by $X$ and $X'$
respectively, where $X$ and $X'$ are some r.v.'s with values in
$S$. The underlying probability spaces on which $X$ and $X'$ are
respectively defined are in general different. Therefore functions
of both $X$ and $X'$, e.g. $d(X,X')$, and their associated probabilities
need not make sense. Additional conditions on joint probabilities
are needed to construct one probability space on which both $X$ and
$X'$ are defined. 

One such condition is independence, to be made precise in a later
section, where knowledge on the value of $X$ has no effect whatsoever
on the probabilities concerning $X'$. 

In some other situations, it is desirable to have models where $X=X'$
if $E=E'$, and more generally where $d(X,X')$ is small when $E$
is close to $E'$. In this section, we construct the Skorokhod representation
which, to each distribution $E$ on $S$, assigns a unique r.v. $X:[0,1]\rightarrow S$
which induces $E$. In the context of applications to random fields,
Theorem 3.1.1 of \cite{Skorohod56} introduced said representation
and proves that it is continuous relative to weak convergence of $E$
and a.u. convergence of $X$. We will prove this result, for applications
in the next chapter.

In addition, we will prove that, when restricted to a tight subset
of distributions, Skorokhod's representation is uniformly continuous
relative to the distribution metric $\rho_{Dist,\xi}$, and the metric
$\rho_{Prob}$ on r.v.'s. The metrics $\rho_{Dist,\xi}$ and $\rho_{Prob}$
were introduced in Definition \ref{Def. Distribution metric} and
in Proposition \ref{Prop. Basics of the probability metric} respectively. 

The Skorokhod representation is a generalization of the \index{quantile mapping}\emph{quantile
mapping} which, to each P.D.F. $F$, assigns the r.r.v. $X\equiv F^{-1}:[0,1]\rightarrow R$
on the probability space $[0,1]$ relative to the uniform distribution,
where $X$ can easily shown to induce the P.D.F. $F$.  

Skorokhod's proof, in terms of Borel sets, is recast here in terms
of a given partition of unity $\pi$. The use of a partition of unity
facilitates the proof of the aforementioned metrical continuity.

Recall that $[\cdot]_{1}$ is an operation which assigns to each $r\in(0,\infty)$
an integer $[r]_{1}\in(r,r+2)$. 
\begin{thm}
\label{Thm.Construction of Skorood  Representation} \textbf{\emph{(Construction
of the Skorokhod Representation)}} Let $\xi\equiv(A_{n})_{n=1,2,\cdots}$
be a binary approximation of the locally compact metric space $(S,d)$,
relative to the reference point $x_{\circ}\in S$. Let $\widehat{J}(S,d)$
be the set of distributions on $(S,d)$. Recall that $M(\Theta_{0},S)$
stands for the space of r.v.'s on the probability space $(\Theta_{0},L_{0},I)$
with values in $(S$,d).

Then there exists a function
\[
\Phi_{Sk,\xi}:\widehat{J}(S,d)\rightarrow M(\Theta_{0},S)
\]
such that, for each $E\in\widehat{J}(S,d)$, the r.v. $X\equiv\Phi_{Sk,\xi}(E):\Theta_{0}\rightarrow S$
induces the distribution $E$, or $I_{X}=E$ in symbols. 

The function $\Phi_{Sk,\xi}$ is called the\emph{ Skorokhod representation}\index{Skorokhod representation}
of distributions on $(S,d)$ determined by $\xi$.
\end{thm}
\begin{proof}
Let $E\in\widehat{J}(S,d)$ be arbitrary. Let 
\[
\pi\equiv(\{g_{n,x}:x\in A_{n}\})_{n=1,2,\cdots}
\]
be the partition of unity of $(S,d)$ determined by $\xi$, as in
Definition \ref{Def. Partition of unity for locally compact (S,d)}.

1. Let $n\geq1$ be arbitrary. By Definition of \ref{Def. Binary approximationt and Modulus of local compactness},
the enumerated finite set $A_{n}\equiv\{x_{n,1},\cdots,x_{n,\kappa(n)}\}$
is a $2^{-n}-$approximation of $(d(\cdot,x_{\circ})\leq2^{n})$.
In other words, 
\begin{equation}
A_{n}\subset(d(\cdot,x_{\circ})\leq2^{n})\label{eq:temp-435}
\end{equation}
and 
\begin{equation}
(d(\cdot,x_{\circ})\leq2^{n})\subset\bigcup_{x\in A(n)}(d(\cdot,x)\leq2^{-n}).\label{eq:temp-505-4}
\end{equation}
Recall from Proposition \ref{Prop. Properties of parittion of unity-1}
that $0\leq g_{n,x}\leq\sum_{x\in A(n)}g_{n,x}\leq1$, 
\begin{equation}
(g_{n,x}>0)\subset(d(\cdot,x)\leq2^{-n+1}),\label{eq:temp-433}
\end{equation}
for each $x\in A_{n}$, and that
\begin{equation}
\bigcup_{x\in A(n)}(d(\cdot,x)\leq2^{-n})\subset(\sum_{x\in A(n)}g_{n,x}=1).\label{eq:temp-434}
\end{equation}
Define $K_{n}\equiv\kappa_{n}+1$, and define the sequence
\[
(f_{n,1},\cdots,f_{n,K(n)})
\]
\begin{equation}
\equiv(g_{n,x(n,1)},,\cdots,g_{n,x(n,\kappa(n))},(1-\sum_{x\in A(n)}g_{n,x}))\label{eq:temp-200}
\end{equation}
of nonnegative continuous functions on $S$. Then
\begin{equation}
\sum_{k=1}^{K(n)}f_{n,k}=1\label{eq:temp-432}
\end{equation}
on $S$. 

2. For the purpose of this proof, an open interval is defined by the
pair of its end points $a,b$, where $0\leq a\leq b\leq1$. Two open
intervals $(a,b),(a',b')$ are considered equal if $a=a'$ and $b=b'$.
For arbitrary open intervals $(a,b),(a',b')\subset[0,1]$ we will
write $(a,b)<(a',b')$ if $b\leq a'$.

3. Let $n\geq1$ be arbitrary. Define the product set 
\[
B_{n}\equiv\{1,\cdots,K_{1}\}\times\cdots\times\{1,\cdots,K_{n}\}.
\]
Let $\mu$ denote the Lebesgue measure on $[0,1]$. Define the open
interval $\Theta\equiv(0,1)$. Then, since
\[
\sum_{k=1}^{K(1)}Ef_{n,k}=E\sum_{k=1}^{K(1)}f_{n,k}=E1=1,
\]
we can subdivide the open interval $\Theta$ into mutually exclusive
open subintervals $\Theta_{1},\cdots,$ $\Theta_{K(1)},$ such that
\[
\mu\Theta_{k}=Ef_{n,k}
\]
for each $k=1,\cdots,K_{1}$, and such that $\Theta_{k}<\Theta_{j}$
for each $k=1,\cdots,\kappa_{1}$ with $k<j$.

4. We will construct, for each $n\geq1$, a family of mutually exclusive
open subintervals
\[
\{\Theta_{k(1),\cdots,k(n)}:(k_{1},\cdots,k_{n})\in B_{n}\}
\]
of $(0,1)$ such that, for each $(k_{1},\cdots,k_{n})\in B_{n}$,
we have

(i) $\mu\Theta_{k(1),\cdots,k(n)}=Ef_{1,k(1)}\cdots f_{n,k(n)},$ 

(ii) $\Theta_{k(1),\cdots,k(n)}\subset\Theta_{k(1),\cdots,k(n-1)}$
if $n\geq2$, 

(iii) $\Theta_{k(1),\cdots,k(n-1),k}<\Theta_{k(1),\cdots,k(n-1),j}$
for each $k,j=1,\cdots,\kappa_{n}$ with $k<j$.

5. Proceed inductively. Step 3 above gave the construction for $n=1$.
Now suppose the construction has been carried out for some $n\geq1$
such that Conditions (i-iii) are satisfied. Consider each $(k_{1},\cdots,k_{n})\in B_{n}$.
Then

\[
\sum_{k=1}^{K(n+1)}Ef_{1,k(1)}\cdots f_{n,k(n)}f_{n+1,k}=Ef_{1,k(1)}\cdots f_{n,k(n)}\sum_{k=1}^{K(n+1)}f_{n+1,k}
\]
\[
=Ef_{1,k(1)}\cdots f_{n,k(n)}=\mu\Theta_{k(1),\cdots,k(n)},
\]
where the last equality is because of Condition (i) in the induction
hypothesis. Hence we can subdivide $\Theta_{k(1),\cdots,k(n)}$ into
$K_{n+1}$ mutually exclusive open subintervals 
\[
\Theta_{k(1),\cdots,k(n),1},\cdots,\Theta_{k(1),\cdots,k(n),K(n+1)}
\]
such that 
\begin{equation}
\mu\Theta_{k(1),\cdots,k(n),k(n+1)}=Ef_{1,k(1)}\cdots f_{n,k(n)}f_{n+1,k(n+1)}\label{eq:temp-147-2-1}
\end{equation}
for each $k_{n+1}=1,\cdots,K_{n+1},.$ Thus Condition (i) holds for
$n+1$. In addition, we can arrange these open subintervals such that
\[
\Theta_{k(1),\cdots,k(n),k}<\Theta_{x(1),\cdots,x(n),j}
\]
for each $k,j=1,\cdots,K_{n+1}$ with $k<j$. This establishes Condition
(iii) for $n+1$. Condition (ii) also holds for $n+1$ since, by construction,
$\Theta_{k(1),\cdots,k(n),k(n+1)}$ is a subinterval of $\Theta_{k(1),\cdots,k(n)}$
for each \textbf{$(k_{1},\cdots,k_{n+1})\in B_{n+1}$}. Induction
is completed. 

6. Note that Condition (i) implies that
\[
\mu\sum_{(k(1),\cdots,k(n))\in B(n)}\Theta_{k(1),\cdots,k(n)}
\]
\begin{equation}
=E(\sum_{k(1)=1}^{K(1)}f_{1,K(1)})\cdots(\sum_{k(n)=1}^{K(n)}f_{n,K(n)})=E1=1\label{eq:temp-487}
\end{equation}
for each $n\geq1$. Hence
\begin{equation}
D\equiv\bigcap_{n=1}^{\infty}\bigcup_{(k(1),\cdots,k(n))\in B(n)}\Theta_{k(1),\cdots,k(n)}\label{eq:temp-145}
\end{equation}
is a full subset of $[0,1]$.

7. Let $\theta\in D$ be arbitrary. Consider each $n\geq1$. Then
$\theta\in\Theta_{k(1),\cdots,k(n)}$ for some unique sequence $(k_{1},\cdots,k_{n})\in B_{n}$
since the intervals in each union in equality \ref{eq:temp-145} are
mutually exclusive. By the same token, $\theta\in\Theta_{j(1),\cdots,j(n+1)}$
for some unique $(j_{1},\cdots,j_{n+1})\in B_{n+1}$. Then $\theta\in\Theta_{j(1),\cdots,j(n)}$
in view of Condition (ii) in Step 4 above. Hence, by uniqueness of
the sequence $(k_{1},\cdots,k_{n})$, we have $(j_{1},\cdots,j_{n})=(k_{1},\cdots,k_{n})$.
Now define $k_{n+1}\equiv j{}_{n+1}$. It follows that $\theta\in\Theta_{k(1),\cdots,k(n+1)}$.
Thus we obtain inductively a unique sequence $(k_{p})_{p=1,2,\cdots}$
such that $k_{p}\in\{1,\cdots,K_{p}\}$ and $\theta\in\Theta_{k(1),\cdots,k(p)}$
for each $p\geq1$. 

Since the open interval $\Theta_{k(1),\cdots,k(n)}$ contains the
given point $\theta$, it has positive Lebesgue measure. In view of
Condition (i) in Step 4 above, it follows that
\begin{equation}
Ef_{1,k(1)}\cdots f_{n,k(n)}>0,\label{eq:temp-136}
\end{equation}
where $n\geq1$ is arbitrary. 

8. Define the function $X_{n}:[0,1]\rightarrow(S,d)$ by 
\[
domain(X_{n})\equiv D\subset\bigcup_{(k(1),\cdots,k(n))\in B(n)}\Theta_{k(1),\cdots,k(n)},
\]
and by
\begin{equation}
X_{n}\equiv x_{n,k(n)}\;\mathrm{or}\;x_{\circ}\quad\mathrm{on}\;D\Theta_{k(1),\cdots,k(n)},\;\mathrm{according\;as\;\mathit{k_{n}}\leq\kappa_{n}\;or\;}k_{n}=\kappa_{n}+1,\label{eq:temp-149}
\end{equation}
for each $(k_{1},\cdots,k_{n})\in domain(X_{n})$. Then, according
to Proposition \ref{Prop. Basing seq of measurable functions on measurable partition},
$X_{n}\in M(\Theta_{0},S)$. In other words, $X_{n}$ is a r.v with
values in the metric space $(S,d)$. Now define the function $X:[0,1]\rightarrow(S,d)$
by 
\[
domain(X)\equiv\{\theta\in D:\lim_{n\rightarrow\infty}X_{n}(\theta)\quad exists\},
\]
and by
\[
X(\theta)\equiv\lim_{n\rightarrow\infty}X_{n}(\theta)
\]
for each $\theta\in domain(X)$. We proceed to prove that the function
$X$ is a r.v. by showing that $X_{n}\rightarrow X$ a.u.

9. To that end, let $n\geq1$ be arbitrary. Define 
\begin{equation}
m\equiv m_{n}\equiv n\vee[\log_{2}(1\vee\beta(2^{-n}))]_{1},\label{eq:temp-441}
\end{equation}
where $\beta$ is the given modulus of tightness of the distribution
$E$ relative to the reference point $x_{\circ}\in S$. Then 
\[
2^{m}>\beta(2^{-n}).
\]
Take an arbitrary $\alpha_{n}\in(\beta(2^{-n}),2^{m})$. Then 
\begin{equation}
E(d(\cdot,x_{\circ})>\alpha_{n})\leq2^{-n}\label{eq:temp-148}
\end{equation}
because $\beta$ is a modulus of tightness of $E$. At the same time,
\[
(d(\cdot,x_{\circ})\leq\alpha_{n})\subset(d(\cdot,x_{\circ})\leq2^{m})
\]
\begin{equation}
\subset\bigcup_{x\in A(m)}(d(\cdot,x)\leq2^{-m})\subset(\sum_{x\in A(m)}g_{m,x}=1),\label{eq:temp-152}
\end{equation}
where the second and third inclusion are by relations \ref{eq:temp-505-4}
and \ref{eq:temp-434} respectively. Define the Lebesgue $\mathrm{measurable}$
set 
\begin{equation}
D_{n}\equiv\bigcup_{(k(1),\cdots,k(m))\in B(m);k(m)\leq\kappa(m)}\Theta_{k(1),\cdots,k(m)}\subset[0,1].\label{eq:temp-382}
\end{equation}
Then
\[
\mu(D_{n}^{c})=\sum_{(k(1),\cdots,k(m))\in B(m);k(m)=K(m)}\mu\Theta_{k(1),\cdots,k(m)}
\]
\[
=\sum_{(k(1),\cdots,k(m))\in B(m);k(m)=K(m)}Ef_{1,k(1)}\cdots f_{m,k(m)}
\]
\[
=\sum_{k(1)=1}^{K(1)}\cdots\sum_{k(m-1)=1}^{K(m-1)}Ef_{1,k(1)}\cdots f_{m-1,k(m)-1}f_{m,K(m)}
\]
\[
=Ef_{m,K(m)}=E(1-\sum_{x\in A(m)}g_{m,x})
\]
\begin{equation}
\leq E(d(\cdot,x_{\circ})>\alpha_{n})\leq2^{-n},\label{eq:temp-430}
\end{equation}
where the first inequality is thanks to relation \ref{eq:temp-152},
and the second is inequality \ref{eq:temp-148}. 

10. Consider each $\theta\in D.$ By Step 7, there exists a unique
sequence $(k_{p})_{p=1,2,\cdots}$ such that $k_{p}\in\{1,\cdots,K_{p}\}$
and $\theta\in\Theta_{k(1),\cdots,k(p)}$ for each $p\geq1$. In particular,
$\theta\in\Theta_{k(1),\cdots,k(m)}$. In view of the defining equality
\ref{eq:temp-382} for the set $D_{n}$, it follows that $k_{m}\leq\kappa_{m}$,
whence 
\[
f_{m,k(m)}\equiv g_{m,x(m,k(m))},
\]
according to the defining equality \ref{eq:temp-200}. Moreover, by
Condition (ii) in Step 4 above, we have $\theta\in\Theta_{k(1),\cdots,k(q)}\subset\Theta_{k(1),\cdots,k(m)}$
for each $q\geq m$. Suppose, for the sake of a contradiction, that
$k_{m+1}=K_{m+1}$. Then, by the defining equality \ref{eq:temp-200},
we have 
\[
f_{m+1,k(m+1)}\equiv1-\sum_{x\in A(m+1)}g_{m+1,x}.
\]
Hence inequality \ref{eq:temp-136}, applied to $m+1$, yields
\begin{equation}
0<Ef_{1,k(1)}\cdots f_{m-1,k(m-1)}g_{m,x(m,k(m))}f_{m+1,k(m+1)}(1-\sum_{x\in A(m+1)}g_{m+1,x})\label{eq:temp-46}
\end{equation}
On the other hand, using successively the relations \ref{eq:temp-433}
\ref{eq:temp-435}, \ref{eq:temp-505-4}, and \ref{eq:temp-434},
we obtain 
\[
(g_{m,x(m,k(m))}>0)\subset(d(\cdot,x_{m,k(m)})\leq2^{-m+1})\subset(d(\cdot,x_{\circ})\leq2^{m}+2^{-m+1})
\]
\[
\subset(d(\cdot,x_{\circ})\leq2^{m+1})\subset\bigcup_{x\in A(m+1)}(d(\cdot,x)\leq2^{-m-1})
\]
\[
\subset(\sum_{x\in A(m+1)}g_{m+1,x}=1).
\]
Hence the right-hand side of the strict inequality \ref{eq:temp-46}
vanishes, while the left-hand side is $0$, a contradiction. We conclude
that $k_{m+1}\leq\kappa_{m+1}.$ Repeating these steps, we obtain,
for each $q\geq m$, the inequality
\[
k_{q}\leq\kappa_{q},
\]
whence
\begin{equation}
f_{q,k(q)}\equiv g_{q,x(q,k(q))}.\label{eq:temp-440}
\end{equation}
It follows from the defining equality \ref{eq:temp-149} that 
\[
X_{q}(\theta)=x_{q,k(q)}
\]
for each $q\geq m\equiv m_{n}$, where $\theta\in DD_{n}$ is arbitrary.
Summing up, we have
\begin{equation}
DD_{n}\subset\bigcap_{q=m}^{\infty}(X_{q}=x_{q,k(q)}).\label{eq:temp-438}
\end{equation}

11. Continue with $\theta\in DD_{n}$ and the corresponding unique
sequence $(k_{p})_{p=1,2,\cdots}$ in the previous step. Let $q\geq p\geq m$
be arbitrary. For abbreviation, write $y_{p}\equiv x_{p,k(p)}$. Then
inequality \ref{eq:temp-136} and equality \ref{eq:temp-440} together
imply that 
\[
Ef_{1,k(1)}\cdots f_{m-1,k(m-1)}g_{m,y(m)}\cdots g_{p,y(p)}\cdots g_{q,y(q)}>0.
\]
Hence there exists $z\in S$ such that
\[
(f_{1,k(1)}\cdots f_{m-1,k(m-1)}g_{m,y(m)}\cdots g_{p,y(p)}\cdots g_{q,y(q)})(z)>0,
\]
whence $g_{p,y(p)}(z)>0$ and $g_{q,y(q)}(z)>0$. Consequently, by
relation \ref{eq:temp-433}, we obtain 
\begin{equation}
d(y_{p},y_{q})\leq d(y_{p},z)+d(z,y_{q})\leq2^{-p+1}+2^{-q+1}\rightarrow0\label{eq:temp-439}
\end{equation}
as $p,q\rightarrow\infty$. Since $(S,d)$ is complete, we have 
\[
X_{p}(\theta)\equiv x_{p,k(p)}\equiv y_{p}\rightarrow y
\]
as $p\rightarrow\infty$, for some $y\in S$. Hence $\theta\in domain(X)$,
with $X(\theta)\equiv y$. Moreover, with $q\rightarrow\infty$ in
inequality \ref{eq:temp-439}, we obtain 
\begin{equation}
d(X_{p}(\theta),X(\theta))\leq2^{-p+1},\label{eq:temp-154}
\end{equation}
where $p\geq m\equiv m_{n}$ and $\theta\in DD_{n}$ are arbitrary.
Since $\mu(DD_{n})^{c}=\mu D_{n}^{c}\leq2^{-n}$ is arbitrarily small
when $n\geq1$ is sufficiently large, we conclude that $X_{n}\rightarrow X$
a.u. relative to the Lebesgue measure $I$, as $n\rightarrow\infty$.
It follows that the function $X:[0,1]\rightarrow S$ is $\mathrm{measurable}$.
In other words, $X$ is a r.v.

12. It remains to verify that $I_{X}=E$, where $I_{X}$ is the distribution
induced by $X$ on $S$. For that purpose, let $h\in C(S)$ be arbitrary.
We need to prove that $Ih(X)=Eh$, where, without loss of generality,
we assume that $|h|\leq1$ on $S$. Let $\delta_{h}$ be a modulus
of continuity of the function h. Let $\varepsilon>0$ be arbitrary.
Let $n\geq1$ be so large that (i') 
\begin{equation}
2^{-n}<\varepsilon\wedge\frac{1}{2}\delta_{h}(\frac{\varepsilon}{3}),\label{eq:temp-436}
\end{equation}
and (ii') $f$ is supported by $(d(\cdot,x_{\circ})\leq2^{n})$. Then
relation \ref{eq:temp-505-4} implies that $f$ is supported by $\bigcup_{x\in A(n)}(d(\cdot,x)\leq2^{-n})$.
At the same time, by the defining equality \ref{eq:temp-149} of the
simple r.v. $X_{n}:[0,1]\rightarrow S$, we have
\[
Ih(X_{n})=\sum_{(k(1),\cdots,k(n))\in B(n);k(n)\leq\kappa(n)}h(x_{n,k(n)})\mu\Theta_{k(1),\cdots,k(n)}+h(x_{\circ})\mu(D_{n}^{c})
\]
\[
=\sum_{k(1)=1}^{K(1)}\cdots\sum_{k(n-1)=1}^{K(n-1)}\sum_{k(n)=1}^{\kappa(n)}h(x_{n,k(n)})Ef_{1,k(1)}\cdots f_{n-1,k(n-1)}f_{n,k(n)}+h(x_{\circ})\mu(D_{n}^{c})
\]
\[
=\sum_{k=1}^{\kappa(n)}h(x_{n,k})Ef_{n,k}+h(x_{\circ})\mu(D_{n}^{c})
\]
\[
=\sum_{k=1}^{\kappa(n)}h(x_{n,k(n)})Eg_{n,x(n,k)}+h(x_{\circ})\mu(D_{n}^{c})
\]
\[
=\sum_{x\in A(n)}h(x)Eg_{n,x}+h(x_{\circ})\mu(D_{n}^{c})
\]
where the third equality is thanks to equality \ref{eq:temp-432}.
Hence
\[
|Ih(X_{n})-E\sum_{x\in A(n)}h(x)g_{n,x}|
\]
\begin{equation}
\leq|h(x_{\circ})\mu(D_{n}^{c})|\leq\mu(D_{n}^{c})\leq2^{-n}<\varepsilon.\label{eq:temp-431}
\end{equation}
At the same time, since $A_{n}$ is a $2^{-n}$\emph{-}partition of
unity of $(S,d)$, with $2^{-n}<\frac{1}{2}\delta_{h}(\frac{\varepsilon}{3})$,
Proposition \ref{Prop. Approx  by  interpolation} implies that$\left\Vert \sum_{x\in A(n)}h(x)g_{n,x}-h\right\Vert \leq\varepsilon$.
Hence
\[
|E\sum_{x\in A(n)}h(x)g_{n,x}-Eh|\leq\varepsilon.
\]
Inequality \ref{eq:temp-431} therefore yields 
\[
|Ih(X_{n})-Eh|<2\varepsilon.
\]
Since $\varepsilon>0$ is arbitrarily small, we have $Ih(X_{n})\rightarrow Eh$
as $n\rightarrow\infty$. On the other hand, $h(X_{n})\rightarrow h(X)$
a.u. relative to $I$. Hence the Dominated Convergence Theorem implies
that $Ih(X_{n})\rightarrow Ih(X)$. It follows that $Eh=Ih(X)$, where
$h\in C(S)$ is arbitrary. We conclude that $E=I_{X}$. 

Define $\Phi_{Sk,\xi}(E)\equiv X$, and the theorem is proved. 
\end{proof}
\begin{thm}
\label{Thm. Continuity of Skorohod Representation}\textbf{\emph{
(Metrical Continuity of Skorokhod Representation) }}Let $\xi\equiv(A_{n})_{n=1,2,\cdots}$
be a binary approximation of the locally compact metric space $(S,d)$
relative to the reference point $x_{\circ}\in S$. Let $\left\Vert \xi\right\Vert \equiv(\kappa_{n})_{n=1,2,\cdots}$
be the modulus of local compactness of $(S,d)$ corresponding to $\xi$.
In other words, $\kappa_{n}\equiv|A_{n}|$ is the number of elements
in the enumerated finite set $A_{n}$, for each $n\geq1$.

Let $\widehat{J}(S,d)$ be the set of distributions on $(S,d)$. Let
$\widehat{J}^{\beta}(S,d)$ be a tight subset of $\widehat{J}(S,d)$,
with a modulus of tightness $\beta$ relative to $x_{\circ}$. Recall
the probability metric $\rho_{Prob}$ on $M(\Theta_{0},S)$ defined
in Definition \ref{Def. Metric rho_prob on space of r.v.'s}. Then
the Skorokhod representation 
\[
\Phi_{Sk,\xi}:(\widehat{J}(S,d),\rho_{Dist,\xi})\rightarrow(M(\Theta_{0},S),\rho_{Prob})
\]
constructed in Theorem \ref{Thm.Construction of Skorood  Representation}
is uniformly continuous on the subset $\widehat{J}^{\beta}(S,d)$,
with a modulus of continuity $\delta_{Sk}(\cdot,\left\Vert \xi\right\Vert ,\beta)$
depending only on $\left\Vert \xi\right\Vert $ and $\beta$.
\end{thm}
\begin{proof}
Refer to the proof of Theorem \ref{Thm.Construction of Skorood  Representation}
for notations. In particular, let 
\[
\pi\equiv(\{g_{n,x}:x\in A_{n}\})_{n=1,2,\cdots}
\]
denote the partition of unity of $(S,d)$ determined by $\xi$. 

1. Let $n\geq1$ be arbitrary. Recall from Proposition \ref{Prop. Properties of parittion of unity-1}
that, for each $x\in A_{n}$, the functions $g_{n,x}$ and $\sum_{y\in A(n)}g_{n,y}$
in $C(S)$ have Lipschitz constant $2^{n+1}$ and have values in $[0,1]$.
Consequently, each of the functions $f_{n,1},\ldots,f_{n,K(n)}$ defined
in formula \ref{eq:temp-200}has Lipschitz constant $2^{n+1}$. 

2. Let 
\[
(k_{1},\cdots,k_{n})\in B_{n}\equiv\{1,\cdots,K_{1}\}\times\cdots\times\{1,\cdots,K_{n}\}
\]
be arbitrary. Then the function
\[
h_{k(1),\cdots,k(n)}\equiv n^{-1}\sum_{p=1}^{n}\sum_{k=1}^{k(p)-1}f_{1,k(1)}\cdots f_{p-1,k(p-1)}f_{p,k}\in C_{ub}(S)
\]
has values in $[0,1]$. By elementary operations of Lipschitz constants,
as in Exercise \ref{Ex.  Lipshitz Constants}, the function $h_{k(1),\cdots,k(n)}$
has Lipschitz constant given by
\[
n^{-1}\sum_{p=1}^{n}\sum_{k=1}^{k(p)-1}(2^{1+1}+2^{2+1}+\cdots+2^{p+1})
\]
\[
=2^{2}n^{-1}\sum_{p=1}^{n}(\kappa_{p}-1)(1+2+\cdots+2^{p-1})<2^{2}n^{-1}\sum_{p=1}^{n}\kappa_{n}2^{p}
\]
\[
<2^{2}n^{-1}\kappa_{n}2^{n+1}=n^{-1}2^{n+3}\kappa_{n}.
\]

3. Now let $E,E'\in\widehat{J}^{\beta}(S,d)$ be arbitrary. Let the
objects $\{\Theta_{k(1),\cdots,k(n)}:(k_{1},\cdots,k_{n})\in B_{n};n\geq1\},$
$D$, $(X_{n})_{n=1,2,\cdots},$ and $X$ be constructed as in Theorem
\ref{Thm.Construction of Skorood  Representation} relative to $E$.
Let the objects $\{\Theta'_{k(1),\cdots,k(n)}:(k_{1},\cdots,k_{n})\in B_{n};n\geq1\},$
$D'$, $(X'_{n})_{n=1,2,\cdots},$ and $X'$ be similarly constructed
relative to $E'$. 

4. Let $\varepsilon>0$ be arbitrary. Fix
\[
n\equiv[3-\log_{2}\varepsilon]_{1}.
\]
Thus $2^{-n+3}<\varepsilon$. As in the proof of Theorem \ref{Thm.Construction of Skorood  Representation},
let 
\begin{equation}
m\equiv m_{n}\equiv n\vee[\log_{2}(1\vee\beta(2^{-n}))]_{1},\label{eq:temp-441-1}
\end{equation}
and let
\begin{equation}
c\equiv m^{-1}2^{m+3}\kappa_{m},\label{eq:temp-157}
\end{equation}
where $\beta$ is the given modulus of tightness of the distributions
$E$ in $\widehat{J}^{\beta}(S,d)$ relative to the reference point
$x_{\circ}\in S$. Let
\begin{equation}
\alpha\equiv2^{-n}\prod_{p=1}^{m}K_{p}{}^{-1}=2^{-n}|B_{m}|^{-1}.\label{eq:temp-155}
\end{equation}
By Proposition \ref{Prop. Modulus of continuity of J->Jf for fixed Lipshitz f},
there exists $\widetilde{\Delta}(m^{-1}\alpha,c,\beta,\left\Vert \xi\right\Vert )>0$
such that if 
\[
\rho_{Dist,\xi}(E,E')<\delta_{Sk}(\varepsilon,\left\Vert \xi\right\Vert ,\beta)\equiv\widetilde{\Delta}(m^{-1}\alpha,c,\beta,\left\Vert \xi\right\Vert ),
\]
then 
\begin{equation}
|Ef-E'f|<m^{-1}\alpha\label{eq:temp-163}
\end{equation}
for each $f\in C_{ub}(S)$ with Lipschitz constant $c>0$ and with
$|f|\leq1$. 

5. Suppose
\begin{equation}
\rho_{Dist,\xi}(E,E')<\delta_{Sk}(\varepsilon,\left\Vert \xi\right\Vert ,\beta).\label{eq:temp-162}
\end{equation}
We will prove that 
\[
\rho_{Prob}(X,X')<\varepsilon.
\]
To that end, let $(k_{1},\cdots,k_{m})\in B_{m}$ be arbitrary. We
will calculate the endpoints of the open interval 
\[
(a_{k(1),\cdots,k(m)},b_{k(1),\cdots,k(m)})\equiv\Theta_{k(1),\cdots,k(m)}.
\]
Recall that, by construction, $\{\Theta_{k(1),\cdots,k(m-1),k}:1\leq k\leq K_{m}\}$
is the set of subintervals in a partition of the open interval $\Theta_{k(1),\cdots,k(m-1)}$
into mutually exclusive open subintervals, with 
\[
\Theta_{k(1),\cdots,k(m-1),k}<\Theta_{k(1),\cdots,k(m-1),j}
\]
if $1\leq k<j\leq K_{m}$. Hence the left endpoint of $\Theta_{k(1),\cdots,k(m)}$
is 
\[
a_{k(1),\cdots,k(m)}=a_{k(1),\cdots,k(m-1)}+\sum_{k=1}^{k(m)-1}\mu\Theta_{k(1),\cdots,k(m-1),k}
\]
\begin{equation}
=a_{k(1),\cdots,k(m-1)}+\sum_{k=1}^{k(m)-1}Ef_{1,k(1)}\cdots f_{m-1,k(m-1)}f_{m,k},\label{eq:temp-437}
\end{equation}
where the second equality is due to Condition (i) in Step 4 of the
proof of Theorem \ref{Thm.Construction of Skorood  Representation}.
Recursively, we then obtain
\[
a_{k(1),\cdots,k(m)}
\]
\[
=a_{k(1),\cdots,k(m-2)}+\sum_{k=1}^{k(m-1)-1}Ef_{1,k(1)}\cdots f_{m-2,k(m-2)}f_{m-1,k}+\sum_{k=1}^{k(m)-1}Ef_{1,k(1)}\cdots f_{m-1,k(m-1)}f_{m,k}
\]
\[
=\cdots=\sum_{p=1}^{m}\sum_{k=1}^{k(p)-1}Ef_{1,k(1)}\cdots f_{p-1,k(p-1)}f_{p,k}\equiv mEh_{k(1),\cdots,k(m)}.
\]

6. Similarly, write 
\[
(a'_{k(1),\cdots,k(m)},b'_{k(1),\cdots,k(m)})\equiv\Theta'_{k(1),\cdots,k(m)}.
\]
Then 
\[
a'_{k(1),\cdots,k(m)}=mE'h_{k(1),\cdots,k(m)}.
\]
Therefore
\[
|a_{k(1),\cdots,k(m)}-a'_{k(1),\cdots,k(m)}|
\]
\begin{equation}
=m|Eh_{k(1),\cdots,k(m)}-E'h_{k(1),\cdots,k(m)}|<mm^{-1}\alpha=\alpha.\label{eq:temp-170}
\end{equation}
where the inequality is by applying inequality \ref{eq:temp-163}
to the function $f\equiv h_{x(1),\cdots,x(m)}$, which has been observed
in Step 2 above to have values in $[0,1]$ and to have Lipschitz constant
$c\equiv m^{-1}2^{m+3}\kappa_{m}$. By symmetry, we can similarly
prove that 
\begin{equation}
|b_{k(1),\cdots,k(m)}-b'_{k(1),\cdots,k(m)}|<\alpha,\label{eq:temp-143}
\end{equation}

7. Inequality \ref{eq:temp-154} in Step 11 of the proof of Theorem
\ref{Thm.Construction of Skorood  Representation} gives
\begin{equation}
d(X_{m},X)\leq2^{-m+1},\label{eq:temp-154-2}
\end{equation}
on $DD_{n}$. Now partition the set $B_{m}\equiv B_{m,0}\cup B_{m,1}\cup B_{m,2}$
into three disjoint subsets, where 
\[
B_{m,0}\equiv\{(k_{1},\cdots,k_{m})\in B_{m}:k_{m}=K_{m}\},
\]
\[
B_{m,1}\equiv\{(k_{1},\cdots,k_{m})\in B_{m}:k_{m}\leq\kappa_{m};\mu\Theta_{k(1),\cdots,k(m)}>2\alpha\},
\]
\[
B_{m,2}\equiv\{(k_{1},\cdots,k_{m})\in B_{m}:k_{m}\leq\kappa_{m};\mu\Theta_{k(1),\cdots,k(m)}<3\alpha\},
\]
Define the set
\[
H\equiv\bigcup_{(k(1),\cdots,k(m))\in B(m,1)}\tilde{\Theta}_{k(1),\cdots,k(m)}\subset[0,1],
\]
where the open interval 
\[
\tilde{\Theta}_{k(1),\cdots,k(m)}\equiv(a_{k(1),\cdots,k(m)}+\alpha,b_{k(1),\cdots,k(m)}-\alpha)
\]
is the $\alpha$-interior of $\Theta_{k(1),\cdots,k(m)}$ for each
$(k_{1},\cdots,k_{m})\in B_{m,1}$. Then
\[
H^{c}=\bigcup_{(k(1),\cdots,k(m))\in B(m,1)}\Theta_{k(1),\cdots,k(m)}\tilde{\Theta}_{k(1),\cdots,k(m)}^{c}\cup
\]
\[
\bigcup_{(k(1),\cdots,k(m))\in B(m,2)}\Theta_{k(1),\cdots,k(m)}\cup\bigcup_{(k(1),\cdots,k(m))\in B(m,0)}\Theta_{k(1),\cdots,k(m)}.
\]
Hence
\[
\mu H^{c}=\sum_{(k(1),\cdots,k(m))\in B(m,1)}2\alpha+\sum_{(k(1),\cdots,k(m))\in B(m,2)}\mu\Theta_{k(1),\cdots,k(m)}
\]
\[
+\mu\bigcup_{(k(1),\cdots,k(m))\in B(m);k(m)=K(m)}\Theta_{k(1),\cdots,k(m)}
\]
\[
<\sum_{(k(1),\cdots,k(m))\in B(m,1)}2\alpha+\sum_{(k(1),\cdots,k(m))\in B(m,2)}3\alpha+\mu D_{n}^{c}
\]
\[
<|B_{m,1}|2\alpha+|B_{m,2}|3\alpha+2^{-n}
\]
\begin{equation}
\leq|B_{m}|3\alpha+2^{-n}=3\cdot2^{-n}+2^{-n}=2^{-n+2},\label{eq:temp-165}
\end{equation}
where the next-to-last equality is from the defining equality \ref{eq:temp-155}.
Note for later reference that the set $H$ depends only on $E$, and
not on $E'$.

8. Now let $\theta\in HDD'$ be arbitrary. Then, according to the
definitions of $H_{n},D$, and $D'$, we have 
\[
\theta\in\tilde{\Theta}_{k(1),\cdots,k(m)}\Theta'_{j(1),\cdots,j(m)}
\]
for some $(k_{1},\cdots,k_{m})\in B_{m,1}$ and some $(j_{1},\cdots,j{}_{m})\in B_{m}$.
Hence, in view of inequalities \ref{eq:temp-170} and \ref{eq:temp-143},
we have
\[
\theta\in(a_{k(1),\cdots,k(m)}+\alpha,b_{k(1),\cdots,k(m)}-\alpha)
\]
\[
\subset(a'_{k(1),\cdots,k(m)},b'_{k(1),\cdots,k(m)})\equiv\Theta'_{k(1),\cdots,k(m)}.
\]
Consequently, 
\begin{equation}
\theta\in\Theta'_{k(1),\cdots,k(m)}\Theta'_{j(1),\cdots,j(m)}.\label{eq:temp-157-1-1}
\end{equation}
The intersection of the two open intervals on the right-hand side
would however be empty unless their subscripts are identical. Hence
$(k_{1},\cdots,k_{m})=(j{}_{1},\cdots,j{}_{m})$. In particular, $j_{m}=k_{m}\leq\kappa_{m}$,
where the inequality is because $(k_{1},\cdots,k_{m})\in B_{m,1}$.
Hence $\theta\in D_{n}D'_{n}$ by the defining equality \ref{eq:temp-382}
for the sets $D_{n}$ and $D'_{n}$. At the same time, by the defining
equality \ref{eq:temp-149} for the r.v.'s $X_{m}$ and $X'_{m}$,
we have 
\[
X_{m}(\theta)=x_{m,k(m)}=x_{m,j(m)}=X'_{m}(\theta).
\]
Since $\theta\in HDD'$ is arbitrary, we have proved that (i) $HDD'\subset D_{n}D'_{n}DD'$,
and (ii) $X_{m}=X'_{m}$ on $HDD'$. 

9. By inequality \ref{eq:temp-154} in the proof of Theorem \ref{Thm.Construction of Skorood  Representation},
we have
\begin{equation}
d(X_{m},X)\vee d(X'_{m},X')\leq2^{-m+1}\label{eq:temp-154-1}
\end{equation}
on $D_{n}D'_{n}DD'$. Combining with Conditions (i) and (ii) in the
previous step, we obtain 
\[
HDD'\subset(X_{m}=X'_{m})\cap(d(X_{m},X)\leq2^{-m+1})\cap(d(X'_{m},X')\leq2^{-m+1})
\]
\begin{equation}
\subset(d(X,X')\leq2^{-m+2})\subset(d(X,X')\leq2^{-n+2})\subset(d(X,X')<\varepsilon),\label{eq:temp-164}
\end{equation}
where $D,D'$ are full sets. 

10. From relation \ref{eq:temp-164} and inequality \ref{eq:temp-165},
we deduce 
\[
\rho_{Prob}(X,X')=I(1\wedge d(X,X');H)+I(1\wedge d(X,X');H^{c})\leq2^{-m+2}+\mu H^{c}
\]
\[
<2^{-m+2}+2^{-n+2}\leq2^{-n+2}+2^{-n+2}=2^{-n+3}<\varepsilon,
\]
where $E,E'\in\widehat{J}^{\beta}(S,d)$ are arbitrary such that $\rho_{Dist,\xi}(E,E')<\delta_{Sk}(\varepsilon,\left\Vert \xi\right\Vert ,\beta)$,
where $X\equiv\Phi_{Sk,\xi}(E)$ and $X'\equiv\Phi_{Sk,\xi}(E')$,
and where $\varepsilon>0$ is arbitrary. Thus the mapping $\Phi_{Sk,\xi}:(\widehat{J}(S,d),\rho_{Dist,\xi})\rightarrow(M(\Theta_{0},S),\rho_{Prob})$
is uniformly continuous on the subspace $\widehat{J}^{\beta}(S,d)$,
with $\delta_{Sk}(\cdot,\left\Vert \xi\right\Vert ,\beta)$ as a modulus
of continuity. The theorem is proved.
\end{proof}
Skorokhod's continuity theorem in \cite{Skorohod56}, in terms of
a.u. convergence, is a consequence of the preceding proof.
\begin{thm}
\label{Thm. Skorood's Continuity in terms off weak convergence and a.u. convergence}
\textbf{\emph{(Continuity of Skorokhod representation in }}term\textbf{\emph{s
of weak convergence and a.u. convergence).}} Let $\xi$ be a binary
approximation of the locally compact metric space $(S,d)$, relative
to the reference point $x_{\circ}\in S$. 

Let $E,E^{(1)},E^{(2)},\cdots$ be a sequence of  distributions on
$(S,d)$ such that $E^{(n)}\Rightarrow E$. Let \textup{$X\equiv\Phi_{Sk,\xi}(E)$}
and \textup{$X^{(n)}\equiv\Phi_{Sk,\xi}(E^{(n)})$} for each $n\geq1$.
Then $X^{(n)}\rightarrow X$ a.u\textup{.}
\end{thm}
\begin{proof}
Let $\left\Vert \xi\right\Vert $ be the modulus of local compactness
of $(S,d)$ corresponding to $\xi$. By hypothesis, $E^{(n)}\Rightarrow E$.
Hence $\rho_{Dist,\xi}(E,E^{(n)})\rightarrow0$ by Proposition \ref{Prop. rho_xi convergence=00003D Weak Convergence}.
By Proposition \ref{Prop. Sequence converging in distribution is tight},
the family $\widehat{J}^{\beta}(S,d)\equiv\{E,E^{(1)},E^{(2)},\cdots\}$
is tight, with some modulus of tightness $\beta$. Let $\varepsilon>0$
be arbitrary. Let $\delta_{Sk}(\varepsilon,\left\Vert \xi\right\Vert ,\beta)>0$
be defined as in Theorem \ref{Thm. Continuity of Skorohod Representation}.
In Step 9 of the proof of Theorem \ref{Thm. Continuity of Skorohod Representation},
we see that there exists a Lebesgue $\mathrm{measurable}$ subset
$H$ of $[0,1]$ which depends only on $E$, with $\mu H^{c}<\varepsilon$,
such that for each $E'\in\widehat{J}^{\beta}(S,d)$ we have\emph{
}
\begin{equation}
H\subset(d(X,X')<\varepsilon)\quad a.s.,\label{eq:temp-147-1-1}
\end{equation}
where $X'\equiv\Phi_{Sk,\xi}(E')$, provided that $\rho_{Dist,\xi}(E,E')<\delta_{Sk}(\varepsilon,\left\Vert \xi\right\Vert ,\beta)$.
Hence, if $p\geq1$ is so large that $\rho_{Dist,\xi}(E,E^{(n)})<\delta_{Sk}(\varepsilon,\left\Vert \xi\right\Vert ,\beta)$
for each $n\geq p$, then 
\[
d(X,X^{(n)})\leq\varepsilon
\]
a.s. on $H^{c}$, for each $n\geq p$. Consequently, $X_{n}\rightarrow X$
a.u. according to Proposition \ref{Prop. Alternative definition of a.u. convergence}. 
\end{proof}

\section{Independence and Conditional Expectation}

The product space introduced in \ref{Def.  Product integration space},
gives a model for compounding two independent experiments into one.
This section introduces the notion of conditional expectations, which
will later be used for a more general method of compounding probability
spaces, or, in the other direction, analyzing an r.v. in terms of
simpler ones. 
\begin{defn}
\label{subsec:Def  Independent sets of r.v.-1}\index{independent events}\index{independent r.v.'s}
\textbf{(Independent set of r.v.'s).} Let $(\Omega,L,E)$ be a probability
space. A finite set $\{X_{1},\cdots,X_{n}\}$ of r.v.'s where $X_{i}$
has values in a complete metric space $(S_{i},d_{i})$, for each $i=1,\cdots,n$
is said to be \emph{independent} if \emph{
\begin{equation}
Ef_{1}(X_{1})\cdots f_{n}(X_{n})=Ef_{1}(X_{1})\cdots Ef_{n}(X_{n})\label{eq:temp-55-1-1}
\end{equation}
}for each\emph{ $f_{1}\in C_{ub}(S_{1}),\cdots,f_{n}\in C_{ub}(S_{n})$.
}In that case we will also simply say that $X_{1},\cdots,X_{n}$ are
independent. A sequence of events $A_{1},\cdots,A_{n}$ is said to
be independent if $1_{A(1)},\cdots,1_{A(n)}$ are independent r.r.v.'s. 

An arbitrary set of r.v.'s is said to be independent if every finite
subset is independent. $\square$
\end{defn}
\begin{prop}
\label{Prop. Independent r.v.s from product space} \textbf{\emph{(Independent
r.v.'s from product space).}} Let $F_{1},\cdots,F_{n}$ be distributions
on the locally compact metric spaces $(S_{1},d_{1}),\cdots,(S_{n},d_{n})$
respectively. \textup{\emph{Let }}\textup{$(S,d)\equiv(S_{1}\times\cdots,S_{n},d_{1}\otimes\cdots\otimes d_{n})$
}\textup{\emph{be the product metric space}}\textup{.} \textup{\emph{Consider
the product integration space}}\textup{ 
\[
(\Omega,L,E)\equiv(S,L,F_{1}\otimes\cdots\otimes F_{n})\equiv\bigotimes_{j=1}^{n}(S_{j},L_{j},F_{j}),
\]
}\textup{\emph{where $(S_{i},L_{i},F_{i})$ is the probability space
that is the completion of }}\textup{$(S_{i},C_{ub}(S_{i}),F_{i})$,
}\textup{\emph{for each }}\textup{$i=1,\cdots,n$. }\textup{\emph{Then
the following holds.}}

\textup{\emph{1. Let}}\textup{ }$i=1,\cdots,n$ be arbitrary. Define\textup{
}\textup{\emph{the coordinate r.v.}}\textup{ $X_{i}:\Omega\rightarrow S_{i}$
by $X_{i}(\omega)\equiv\omega_{i}$} for each $\omega\equiv(\omega_{1},\cdots,\omega_{n})\in\Omega$.
Then the r.v.'s \textup{$X_{1},\cdots,X_{n}$ are }\textup{\emph{independent}}\textup{.
}Moreover,\textup{ $X_{i}$ }\textup{\emph{induces the distribution
$F_{i}$ on}} $(S_{i},d_{i})$ for each $i=1,\cdots,n$. 

\textup{2. $F_{1}\otimes\cdots\otimes F_{n}$ }\textup{\emph{is a
distribution on}}\textup{ $(S,d)$.}\textup{\emph{ Specifically it
is}}\textup{ }the distribution\textup{ $F$ }\textup{\emph{induced
on $(S,d)$ by}}\textup{ }the r.v. $X\equiv(X_{1},\cdots,X_{n})$.
\end{prop}
\begin{proof}
1. By Proposition \ref{Prop. X meas,  f unif continuous and bd on bd subsets =00003D> f(X) meas},
the continuous functions $X_{1},\cdots,X_{n}$ on $(S,L,E)$ are $\mathrm{measurable}$.
Let $f_{i}\in C_{ub}(S_{i})$ be arbitrary, for each $i=1,\cdots,n$.
Then
\begin{equation}
Ef_{1}(X_{1})\cdots f_{n}(X_{n})=F_{1}f_{1}\cdots F_{n}f_{n}\label{eq:temp-268-1}
\end{equation}
by Fubini's Theorem. Let $i=1,\cdots,n$ be arbitrary. In the special
case where $f_{j}\equiv1$ for each $j=1,\cdots,n$ with $j\neq i$,
we obtain, from equality \ref{eq:temp-268-1}, 
\begin{equation}
Ef_{i}(X_{i})=F_{i}f_{i}.\label{eq:temp-383}
\end{equation}
Hence equality \ref{eq:temp-268-1} yields 
\[
Ef_{1}(X_{1})\cdots f_{n}(X_{n})=Ef_{1}(X_{1})\cdots Ef_{n}(X_{n})
\]
where $f_{i}\in C_{ub}(S_{i})$ is arbitrary for each $i=1,\cdots,n$.
Thus the r.v.'s $X_{1},\cdots,X_{n}$ are independent. Moreover equality
\ref{eq:temp-383} shows that the r.v. $X_{i}$ induces the distribution
$F_{i}$ on\emph{ $(S_{i},d_{i})$} for each $i=1,\cdots,n$. 

2. Since $X$ is an r.v. with values in $S$, it induces a distribution
$E_{X}$ on $(S,d)$. Hence
\[
E_{X}f\equiv Ef(X)=(F_{1}\otimes\cdots\otimes F_{n})f
\]
for each $f\in C_{ub}(S)$. Thus $F_{1}\otimes\cdots\otimes F_{n}=E_{X}$
is a distribution $F$ on $(S,d)$.
\end{proof}
\begin{prop}
\label{Prop. Basics of independence} \textbf{\emph{(Basics of independence).}}
Let $(\Omega,L,E)$ be a probability space. Suppose the each $i=1,\cdots,n$,
let $X_{i}$ be a r.v. with values in a complete metric space $(S_{i},d_{i})$,
and let $(S_{i},L_{X(i)},E_{X(i)})$ be the the probability space
it induces on $(S_{i},d_{i})$. Suppose the r.v.'s \emph{$X_{1},\cdots,X_{n}$
}are independent. Then, for arbitrary $f_{1}\in L_{X(1)},\cdots,f_{n}\in L_{X(n)}$,
we have 

\emph{
\begin{equation}
E\prod_{i=1}^{n}f_{i}(X_{i})=\prod_{i=1}^{n}Ef_{i}(X_{i}).\label{eq:temp-55-2}
\end{equation}
}
\end{prop}
\begin{proof}
Consider each $i=1,\cdots,n$. Let $f_{i}\in L_{X(i)}$ be arbitrary.
By Definition \ref{Def. Distribution induced by r.v. w/ values in complete metric space},
$L_{X(i)}$ is the completion of $(\Omega,C_{ub}(S_{i}),E_{X(i)})$.
The r.r.v. $f_{i}\in L_{X(i)}$ is therefore the $L_{1}-$limit relative
to $E_{X(i)}$ of a sequence $(f_{i,h})_{h=1,2,\cdots}$ in $C_{ub}(S_{i})$
as $h\rightarrow\infty$. Moreover, according to \ref{Prop. L(X)=00003D=00007Bf(X): f in L_X=00007D},
we have $f_{i}(X_{i})\in L(X_{i})$ with $Ef_{i}(X_{i})=E_{X(i)}f_{i}$,
. Hence
\[
E|f_{i,h}(X_{i})-f_{i}(X_{i})|=E_{X}|f_{i,h}-f_{i}|\rightarrow0
\]
as $h\rightarrow\infty$. By passing to subsequences if necessary,
we may assume that
\begin{equation}
f_{i,h}(X_{i})\rightarrow f_{i}(X_{i})\quad a.u.\label{eq:temp-384}
\end{equation}
as $h\rightarrow\infty$, for each $i=1,\cdots,n$. 

First consider the case where $f_{i}\geq0$ for each $i=1,\cdots,n$.
Let $a>0$ be arbitrary. In view of the independence of the r.v.'s
\emph{$X_{1},\cdots,X_{n}$, }we have 
\[
E\prod_{i=1}^{n}(0\vee f_{i,h}(X_{i})\wedge a)=\prod_{i=1}^{n}E(0\vee f_{i,h}(X_{i})\wedge a)\equiv\prod_{i=1}^{n}E_{X(i)}(0\vee f_{i,h}\wedge a).
\]
In view of the a.u. convergence \ref{eq:temp-384}, we can let $h\rightarrow\infty$
and apply the Dominated Convergence Theorem to obtain
\[
E\prod_{i=1}^{n}(f_{i}(X_{i})\wedge a)=\prod_{i=1}^{n}E_{X(i)}(f_{i}\wedge a).
\]
Now let $a\rightarrow\infty$ and apply the Monotone Convergence Theorem
to obtain
\[
E\prod_{i=1}^{n}f_{i}(X_{i})=\prod_{i=1}^{n}E_{X(i)}(f_{i})=\prod_{i=1}^{n}Ef_{i}(X_{i}).
\]
The same equality for arbitrary $f_{1}\in L_{X(1)},\cdots,f_{n}\in L_{X(n)}$
follows by linearity.
\end{proof}
We next define the conditional expectation of a r.r.v., as the revised
expectation given the observed values of all the r.v.'s in a family
$G$. 
\begin{defn}
\emph{\label{Def: Conditional Expectation}} \textbf{(Conditional
expectation).} Let $(\Omega,L,E)$ be a probability space, and let
$L'$ be a probability subspace of $L$. Let $Y\in L$ be arbitrary.
If there exists $X\in L'$ such that $EZY=EZX$ for each indicator
$Z\in L'$, then we say that $X$ is the \emph{conditional expectation}
\index{conditional expectation} of $Y$ given $L'$, and define $E(Y|L')\equiv X$.
We will call $L_{|L'}\equiv\{Y\in L:E(Y|L')\quad exists\}$ the subspace
of \emph{conditionally integrable r.r.v.'s}\index{conditionally integrable r.r.v.'s}
given the subspace $L'$. 

In the special case where $L'\equiv L(G)$ is the probability subspace
generated by a given family of r.v.'s with values in some complete
metric space $(S,d)$, we will simply write $E(Y|G)\equiv E(Y|L')$
and say that $L_{|G}\equiv L_{|L'}$ is the subspace of \emph{conditionally
integrable r.r.v.'s}\index{conditionally integrable r.r.v.'s} given
the family $G$. In the case where $G\equiv\{V_{1},\cdots,V_{m}\}$
for some $m\geq1$, we write also $E(Y|V_{1},\cdots,V_{m})\equiv E(Y|G)\equiv E(Y|L')$. 

In the case where $m=1$, and where $V_{1}=1_{A}$ for some measurable
set $A$ with $P(A)>0$, it can easily be verified that, for arbitrary
$Y\in L$ , the conditional $E(Y|1_{A})$ exists is given by $E(Y|1_{A})=P(A)^{-1}E(Y1_{A})1_{A}.$
In that case, we will write 
\[
E_{A}(Y)\equiv P(A)^{-1}E(Y1_{A})
\]
for each $Y\in L$, and write $P_{A}(B)\equiv E_{A}(1_{B})$ for each
measurable set $B$. The next lemma proves that $(\Omega,L,E_{A})$
is a probability space, called the \emph{conditional probability space
given the event}\index{conditional probability space given an event}
$A$. 

More generally, if $Y_{1},\cdots,Y_{n}\in L_{|L'}$ then we define
the vector 
\[
E((Y_{1},\cdots,Y_{n})|L')\equiv(E(Y_{1}|L'),\cdots,E(Y_{n}|L'))
\]
of integrable r.r.v.'s in $L'$.

Let $A$ be an arbitrary $\mathrm{measurable}$ subset of $(\Omega,L,E)$.
If $1_{A}\in L_{|L'}$ we will write $P(A|L')\equiv E(1_{A}|L')$
and call $P(A|L')$ the \emph{conditional probability}\index{conditional probability}
of the event $A$ given the probability subspace $L'$. If $1_{A}\in L_{|G}$
for some given family of r.v.'s with values in some complete metric
space $(S,d)$, we will simply write $P(A|G)\equiv E(1_{A}|G)$. In
the case where $G\equiv\{V_{1},\cdots,V_{m}\}$, we write also $P(A|V_{1},\cdots,V_{m})\equiv E(1_{A}|V_{1},\cdots,V_{m})$. 

$\square$ 
\end{defn}
Before proceeding, note that the statement $E(Y|L')=X$ asserts two
things: that $E(Y|L')$ exists, and that it is equal to $X$. We have
defined the conditional expectation without the sweeping classical
assertion of its existence. Before we use a particular conditional
expectation, we will first supply a proof of its existence. 
\begin{lem}
\label{Lem. Conditional prob space is indeed a probability space-1}
\textbf{\emph{(Conditional probability space given an event is indeed
a probability space). }}Let the measurable set $A$ be arbitrary,
with $P(A)>0$. Then\emph{ the triple $(\Omega,L,E_{A})$ is indeed
a a probability space.}
\end{lem}
\begin{proof}
We need to verify the conditions in Definition \ref{Def. Integration Space}
for an integration space.

1. Clearly $E_{A}$ is a linear function on the linear $L$.

2. Let $(Y_{i})_{i=0,1,2,\cdots}$ be an arbitrary sequence of functions
in $L$ such that $Y_{i}$ is non-negative for each $i\geq1$ and
such that $\sum_{i=1}^{\infty}E_{A}(Y_{i})<E_{A}(Y_{0})$. Then $\sum_{i=1}^{\infty}E(Y_{i}1_{A})<E(Y_{0}1_{A})$
by the definition of the function $E_{A}$. Hence, since $E$ is an
integration, there exists $\omega\in\bigcap_{i=0}^{\infty}domain(Y_{i}1_{A})$
such that $\sum_{i=1}^{\infty}Y_{i}(\omega)1_{A}(\omega)<Y_{0}(\omega)1_{A}(\omega)$.
It follows that $1_{A}(\omega)>0$. Dividing by $1_{A}(\omega)$,
we obtain $\sum_{i=1}^{\infty}Y_{i}(\omega)<Y_{0}(\omega)$.

3. Let $Y\in L$ be arbitrary. Then $E_{A}(Y\wedge n)\equiv E(Y\wedge n)1_{A}\rightarrow E(Y)1_{A}\equiv E_{A}(Y)$
as $n\rightarrow\infty$. Similarly , $E_{A}(|Y|\wedge n^{-1})\equiv E(|Y|\wedge n^{-1})1_{A}\rightarrow0$
as $n\rightarrow\infty$. 

Summing up, all three conditions in Definition \ref{Def. Integration Space}
are satisfied by the triple $(\Omega,L,E_{A})$. Because $L$ is complete
relative to the integration $E$, in the sense of Definition \ref{Def. Integrablee functions and Completion of integration space},
so it can trivially be verified that $L$ is complete relative to
the integration $E_{A}$. Because $1\in L$ with $E_{A}(1)=1$, the
complete integration space\emph{ $(\Omega,L,E_{A})$ }is a probability
space.
\end{proof}
$\,$

We will show that the conditional expectation is unique if it exists,
two r.v.'s considered equal if they are equal a.s. The next two proposition
proves basic properties of conditional expectations. They would be
trivial classically, because the principle of infinite search would
imply, via the Radon-Nikodym Theorem, that $L_{|L'}=L$.
\begin{prop}
\label{Prop. Basics of Conditional expectations} \textbf{\emph{(Basics
of conditional expectation).}} Let $(\Omega,L',E)$ be a probability
subspace of a probability space $(\Omega,L,E)$. Then the following
holds. 
\end{prop}
\begin{enumerate}
\item \emph{Suppose $Y_{1}=Y_{2}$ a.s. in $L$, and suppose $X_{1},X_{2}\in L'$
are such that $EZY_{j}=EZX_{j}$ for each $j=1,2,$ for each indicator
$Z\in L'$. Then $X_{1}=X_{2}$ a.s. Consequently, the conditional
expectation, if it exists, is uniquely defined.}
\item \emph{Suppose $X,Y\in L_{|L'}$. Then $aX+bY\in L_{|L'}$, and 
\[
E(aX+bY|L')=aE(X|L')+bE(Y|L'),
\]
for each $a,b\in R$. If, in addition, $X\leq Y$ a.s., then $E(X|L')\leq E(Y|L')$
a.s. In particular, if }$|X|\in L_{|L'}$ also, then \emph{$|E(X|L')|\leq E(|X||L')$
a.s.}
\item \emph{$E(E(Y|L'))=E(Y)$ for each $Y\in L_{|L'}$. Moreover, $L'\subset L_{|L'}$,
and $E(X|L')=X$ for each $X\in L'$.}
\item \emph{Suppose $Y\in L_{|L'}$. In other words, suppose the conditional
expectation $E(Y|L')$ exists. Then $ZY\in L_{|L'}$, and $E(ZY|L')=ZE(Y|L')$,
for each bounded $Z\in L'$. }
\item \emph{Let $Y\in L$ be arbitrary. Let $G$ be an arbitrary set of
r.v.'s with values in some complete metric space $(S,d)$. Suppose
there exists $X\in L(G)$ such that 
\begin{equation}
EYh(V_{1},\cdots,V_{k})=EXh(V_{1},\cdots,V_{k})\label{eq:temp-403}
\end{equation}
for each $h\in C_{ub}(S^{k})$, for each finite subset $\{V_{1},\cdots,V_{k}\}\subset G,$
and for each $k\geq1$, then $E(Y|G)=X$. }
\item Let $L''$ be a probability subspace with $L''\subset L'$.\emph{
Suppose $X\in L_{|L'}$ with }$Y\equiv E(X|L')$.\emph{ Then $Y\in L_{L''}$
iff $X\in L_{L''}$, in which case $E(X|L'')=E(E(X|L')|L'')$. }
\item \emph{If $Y\in L$, and if $Y,Z$ are independent for each indicator
$Z\in L'$, then $E(Y|L')=EY$.}
\item Let $Y$ be a r.r.v.. with $Y^{2}\in L$. Suppose $X\equiv E(Y|L')$
exists. Then $X^{2},(Y-X)^{2}\in L$ and $EY^{2}=EX^{2}+E(Y-X)^{2}.$
Consequently, $EX^{2}\leq EY^{2}$ and $E(Y-X)^{2}\leq EY^{2}.$
\end{enumerate}
\begin{proof}
1. Let $t>0$ be arbitrary and let $Z\equiv1_{(t<X_{1}-X_{2})}\in L'$.
Then 
\[
tP(t<X_{1}-X_{2})\leq EZ(X_{1}-X_{2})=EZ(Y_{1}-Y_{2})=0,
\]
whence $P(t<X_{1}-X_{2})=0$. It follows that that $(0<X_{1}-X_{2})$
is a null set, and that \emph{$X_{1}\geq X_{2}$} a.s. By symmetry,
$X_{1}=X_{2}$ a.s.\emph{ }

2. Suppose $X,Y\in L_{|L'}$. Let $Z\in L'$ be an arbitrary indicator.
Then
\[
E(Z(aX+bY))=aE(ZX)+bE(ZY)
\]
\[
=aE(ZE(X|L'))+bE(ZE(Y|L'))
\]
\[
=E(Z(aE(X|L')+bE(Y|L'))).
\]
Hence 
\[
E(aX+bY|L')=aE(X|L')+bE(Y|L').
\]
The remainder of Assertion 2 follows from assertion 1.

3. If $Y\in L_{|L'}$, then 
\[
E(E(Y|L'))=E(1E(Y|L'))=E(1Y)=E(Y).
\]
Separately, if $Y\in L'$, then, for each indicator $Z\in L'$, we
have trivially $E(ZY)=E(ZY)$, whence $E(Y|L')=Y$.

4. Suppose\emph{ $Y\in L_{|L'}$, }with\emph{ $X\equiv E(Y|L')$.
}Then\emph{, }by definition, $EZY=EZX$ for each indicator $Z\in L'$.
The equality extends to all linear combinations of integrable indicators.
Since such linear combinations are dense in $L$, the equality extends,
by the Dominated Convergence Theorem, to all each bounded integrable
r.r.v.'s $Z\in L'$. Moreover if $U,Z\in L'$ are bounded and integrable
r.r.v.'s, so is $UZ$, and the previous statement implies that $E(UZY)=E(UZX)$,
whence $E(UY|L')=UX=UE(Y|L')$.

5. Let $Y\in L$ be arbitrary. Suppose equality \ref{eq:temp-403}
holds. Let $Z$ be an arbitrary indicator in $L'\equiv L(G)$. Then
$Z$ is the $L_{1}-$limit of some sequence $(h_{n}(V_{n,1},\cdots,V_{n,k(n)}))_{n=1,2,\cdots}$
where $h_{n}\in C_{ub}(S^{k(n)})$ for each $n\geq1$. Hence, by the
Dominated Convergence Theorem 
\[
EYZ=\lim_{n\rightarrow\infty}EYh_{n}(V_{n,1},\cdots,V_{n,k(n)})=\lim_{n\rightarrow\infty}EXh_{n}(V_{n,1},\cdots,V_{n,k(n)})=EXZ,
\]
where the second equality is due to equality \ref{eq:temp-403}. Thus
$E(Y|G)\equiv E(Y|L')=X$.

6. Let $U\in L''$ an arbitrary indicator. Then $U\in L'$. First,
suppose $X\in L_{L''}$, with \emph{$Z\equiv E(X|L'')\in L''$}. Then,
by assertions 4 and 3 above, we have
\[
EUY\equiv E(UE(X|L'))=E(UE(E(X|L'')|L'))=E(UE(X|L''))\equiv EUZ,
\]
whence $E(Y|L'')=Z$. Consequently, 
\[
E(X|L'')\equiv Z=E(Y|L'')\equiv E(E(X|L')|L'').
\]
Conversely, suppose \emph{$Y\in L_{L''}$, with} \emph{$Z\equiv E(Y|L'')$}.
Then, since $U\in L'$ and $U\in L''$, we have $EUX=EUY=EUZ.$ Hence
$E(X|L'')=Z$ and \emph{$X\in L_{L''}$.}

7.\emph{ }Suppose $Y\in L$, and suppose $Y,Z$ are independent for
each indicator $Z\in L'$. Then, for each indicator $Z\in L'$ we
have
\[
E(ZY)=(EZ)(EY)=E(ZEY).
\]
Since trivially $EY\in L'$, it follows that $E(Y|L')=EY$.

8. Let $Y$ be a r.r.v.. with $Y^{2}\in L$. Suppose $X\equiv E(Y|L')$
exists. Since $Y\in L_{2}$, there exists a deceasing sequence $\varepsilon_{1}>\varepsilon_{2}>\cdots$
of positive real numbers such that $EY^{2}1_{A}<2^{-k}$ for each
$\mathrm{measurable}$ set $A$ with $P(A)<\varepsilon_{k}$, for
each $k\geq1$. Since $X$ is a r.r.v., there exists a sequence $0\equiv a_{0}<a_{1}<a_{2}<\cdots$
of positive real numbers with $a_{k}\rightarrow\infty$ such that
$P(|X|\geq a_{k})<\varepsilon_{k}$. Let $k\geq1$ be arbitrary. Then
\[
EY^{2}1_{(|X|\geq a(k))}<2^{-k}
\]
Write $Z_{k}\equiv1_{(a(k+1)>|X|\geq a(k))}$. Then $Z_{k},XZ_{k},X^{2}Z_{k}\in L'$
are bounded in absolute value by $1,a_{k+1},a_{k+1}^{2}$ respectively.
Hence 
\[
EY^{2}1_{(a(k+1)>|X|\geq a(k))}=E((Y-X)+X)^{2}Z_{k}
\]
\[
=E(Y-X)^{2}Z_{k}+2E(Y-X)XZ_{k}+EX^{2}Z_{k}
\]
\[
=E(Y-X)^{2}Z_{k}+2E(YXZ_{k})-2E(XXZ_{k})+EX^{2}Z_{k}
\]
\[
=E(Y-X)^{2}Z_{k}+2E(E(Y|L')XZ_{k})-2E(XXZ_{k})+EX^{2}Z_{k}
\]
\[
\equiv E(Y-X)^{2}Z_{k}+2E(XXZ_{k})-2E(XXZ_{k})+EX^{2}Z_{k}
\]
\begin{equation}
=E(Y-X)^{2}Z_{k}+EX^{2}Z_{k},\label{eq:temp-407}
\end{equation}
where the fourth equality is by applying Assertion 4 to $Y$ and to
the bounded r.r.v.\emph{ $XZ_{k}\in L'$.} Since $Y^{2}\in L$ by
assumption, we have 
\[
\sum_{k=0}^{\infty}EX^{2}1_{(a(k+1)>|X|\geq a(k))}\leq\sum_{k=0}^{\infty}EY^{2}1_{(a(k+1)>|X|\geq a(k))}=EY^{2}.
\]
Consequently 
\[
X^{2}=\sum_{k=0}^{\infty}X^{2}1_{(a(k+1)>|X|\geq a(k))}\in L.
\]
Similarly, $(Y-X)\in L^{2}.$ Moreover, summing equality \ref{eq:temp-407}
over $k=0,1,\cdots,$ we obtain 
\[
EY^{2}=E(Y-X)^{2}+EX^{2}.
\]
Assertion 8 and the proposition are proved.
\end{proof}
\begin{prop}
\label{Prop. .Space of onditionally integrable functions is closed relative to L1}
\textbf{\emph{(Space of conditionally integrable functions given a
probability subspace is closed).}} Let $(\Omega,L,E)$ be a probability
space. Let $(\Omega,L',E)$ be a probability subspace of $(\Omega,L,E)$.
Let $L_{|L'}$ be the space of r.r.v.'s conditionally integrable given
$L'$. Then the following holds.

1. Let $X,Y\in L$ be arbitrary. Suppose  \emph{$EUX\leq EUY$ }for
each indicator $U\in L'$. Then \emph{$EZX\leq EZY$ }for each bounded
nonnegative r.r.v. \emph{$Z\in L'$. }

2. Suppose $Y\in L_{|L'}$. Then $E|E(Y|L')|\leq E|Y|$. 

3. The linear subspace $L_{|L'}$ of $L$ is closed relative to the
$L_{1}$-norm. 
\end{prop}
\begin{proof}
1. Suppose  \emph{$EUX\leq EUY$ }for each indicator $U\in L'$. Then,
by linearity, \emph{$EVX\leq EVY$ }for each nonnegative linear combination
$Y$ of indicators in \emph{$L'$. }Now consider each bounded nonnegative
r.r.v. \emph{$Z\in L'$. }We may assume, without loss of generality,
that $Z$ has values in $[0,1]$. Then $E|Z-V_{n}|\rightarrow0$ for
some sequence $(V_{n})_{n=1,2,\cdots}$of nonnegative linear combinations
of indicators in $L')$, with values in \emph{$[0,1]$.} By passing
to a subsequence, we may assume that $V_{n}\rightarrow Z$ a.s. Hence,
by the Dominated Convergence Theorem, we have 
\[
EZX=\lim_{n\rightarrow\infty}EV_{n}X\leq\lim_{n\rightarrow\infty}EV_{n}Y=EZY.
\]

2. Suppose $Y\in L_{|L'}$, with $E(Y|L')=X\in L'$. Let $\varepsilon>0$
be arbitrary. Then, since $X$ is integrable, there exists $a>0$
such that $EX1_{(|X|\leq a)}<\varepsilon$. Then 
\[
E|X|=EX1_{(X>a)}-EX1_{(X<-a)}+EY1_{(|X|\leq a)}
\]
\[
<EX1_{(X>a)}-EX1_{(X<-a)}+\varepsilon
\]
\[
=EY1_{(Y>a)}-EY1_{(Y<-a)}+\varepsilon
\]
\[
\leq E|Y|1_{(X>a)}+E|Y|1_{(X<-a)}+\varepsilon\leq E|Y|+\varepsilon,
\]
where the second equality is because $Y=E(X|L')$. Since $\varepsilon>0$
is arbitrarily small, we conclude that $E|X|\leq E|Y|$, as alleged.

3. Let $(Y_{n})_{n=1,2,\cdots}$be a sequence in $L_{|L'}$ such that
$E|Y_{n}-Y|\rightarrow0$ for some $Y\in L$. For each $n\geq1$,
let $X_{n}\equiv E(Y_{n}|L'$). Then, By Assertion 2 above, we have
\[
E|X_{n}-X_{m}|=E|E(Y_{n}-Y_{m}|L')|\leq E|Y_{n}-Y_{m}|\rightarrow0
\]
as $n,m\rightarrow\infty$. Thus $(X_{n})_{n=1,2,\cdots}$ is a Cauchy
sequence in the complete metric space $L'$ relative to the $L_{1}$-norm.
It follows that $E|X_{n}-X|\rightarrow0$ for some $X\in L'$, as
$n\rightarrow\infty$. Hence, for each indicator $Z\in L'$, we have
\[
EXZ=\lim_{n\rightarrow\infty}EX_{n}Z=\lim_{n\rightarrow\infty}EY_{n}Z=EYZ.
\]
It follows that $E(Y|L')=X$ and $Y\in L_{|L'}$.
\end{proof}

\section{Normal Distributions \label{sec:Normal-Distribution} }

The classical development of the topics in the remainder of this chapter
is exemplar of constructive mathematics. However, some tools in said
development have been given many proofs, some constructive and others
not. An example is the spectral theorem for symmetric matrices below.
For ease of reference, we therefore present some such topics here,
using only constructive proofs. 

Recall some notations and basic theorems from Matrix Algebra. 
\begin{defn}
\label{Def. Matrix notations} \textbf{(Matrix notations). }For an
arbitrary $m\times n$ matrix 
\[
\theta\equiv[\theta_{i,j}]_{i=1,\cdots,m;j=1,\cdots,n}\equiv\left[\begin{array}{ccc}
\theta_{1,1}, & \cdots, & \theta_{1,n}\\
\cdot & \cdots & \cdot\\
\cdot & \cdots & \cdot\\
\cdot & \cdots & \cdot\\
\theta_{m,1}, & \cdots, & \theta_{m,n}
\end{array}\right],
\]
 of real or complex elements $\theta_{i,j}$, we will let
\[
\theta^{T}\equiv[\theta_{j,i}]_{j=1,\cdots,n;1=1,\cdots m}=\left[\begin{array}{ccc}
\theta_{1,1}, & \cdots, & \theta_{m,1}\\
\cdot & \cdots & \cdot\\
\cdot & \cdots & \cdot\\
\cdot & \cdots & \cdot\\
\theta_{1,n} & ,\cdots, & \theta_{m,n}
\end{array}\right]
\]
denote the transpose, which is an $n\times m$ matrix. If $n=m$ and
$\theta=\theta^{T}$, then $\theta$ is said to be symmetric. If $\theta_{i,j}=0$
for each $i,j=1,\cdots,n$ with $i\neq j$, then $\theta$ is called
a diagonal matrix. For each sequence of complex numbers $(\lambda_{1},\cdots,\lambda_{n})$,
write $diag(\lambda_{1},\cdots,\lambda_{n})$ for the diagonal matrix
$\theta$ with $\theta_{i,i}=\lambda_{i}$ for each $i=1,\cdots,n$.
A matrix $\theta$ is said to be real if all its elements $\theta_{i,j}$
are real numbers. Unless otherwise specified, all matrices in the
following are assumed to be real. 

For an arbitrary sequence $\bar{\mu}\equiv(\mu_{1},\cdots,\mu_{n})\in R^{n}$,
we will abuse notations and let $\bar{\mu}$ denote also the column
vector 
\[
\bar{\mu}\equiv(\mu_{1},\cdots,\mu_{n})\equiv\left[\begin{array}{c}
\mu_{1}\\
.\\
.\\
.\\
\mu_{n}
\end{array}\right].
\]
Thus $\bar{\mu}^{T}=[\mu_{1},\cdots,\mu_{n}]$. A $1\times1$ matrix
is identified with its only entry. Hence, if $\bar{\mu}\in R^{n}$,
then 
\[
|\mu|\equiv\left\Vert \bar{\mu}\right\Vert \equiv\sqrt{\bar{\mu}^{T}\bar{\mu}}=\sqrt{\sum_{i=1}^{n}\mu_{i}^{2}}.
\]
We will let $I_{n}$ denote the $n\times n$ diagonal matrix $diag(1,\cdots,1)$.
When the dimension $n$ is understood, we write simply $I\equiv I_{n}$.
Likewise, we will write $0$ for any matrix whose entries are all
equal to the real number $0$, with dimensions understood from the
context. 

The determinant of an $n\times n$ matrix $\theta$ is denoted by
$\det\theta$. The $n$ complex roots $\lambda_{1},\cdots,\lambda_{n}$
of the polynomial $\det(\theta-\lambda I)$ of degree $n$ are called
the\emph{ eigenvalues} of $\theta$. Then $\det\theta=\lambda_{1}\cdots\lambda_{n}$.
Let $j=1,\cdots,n$ be arbitrary. Then there exists a nonzero column
vector $x_{j}$, whose elements are in general complex, such that
$\theta x_{j}=\lambda_{j}x_{j}$. The vector $x_{j}$ is called an
\emph{eigenvector}\index{eigenvector} for the eigenvalue $\lambda_{j}$.
If $\theta$ is real and symmetric, then the $n$ eigenvalues $\lambda_{1},\cdots,\lambda_{n}$
are real.

Let $\overline{\sigma}$ be a symmetric $n\times n$ matrix whose
elements are real. Then $\overline{\sigma}$ is said to be \index{nonnegative definite matrix}\emph{nonnegative
definite} if $x^{T}\overline{\sigma}x\geq0$ for each $x\in R^{n}$.
In that case all its eigenvalues are nonnegative, and, for each eigenvalue,
there exists a real eigenvector whose elements are real. It is said
to be \index{positive definite matrix}\emph{positive definite} if
$x^{T}\overline{\sigma}x>0$ for each nonzero $x\in R^{n}$. In that
case all its eigenvalues are positive, whence $\overline{\sigma}$
is nonsingular, with an inverse $\overline{\sigma}^{-1}$. An $n\times n$
real matrix $U$ is said to be orthogonal if $U^{T}U=I$. This is
equivalent to saying that the column vectors of $U$ form an orthonormal
basis of $R^{n}$. $\square$ 
\end{defn}
\begin{thm}
\label{Thm. Spectral Theorem for symmetric matrices} \textbf{\emph{(Spectral
Theorem for Symmetric Matrices).}} Let $\theta$ be an arbitrary $n\times n$
symmetric matrix. Then the following holds.

1. There exists an orthogonal matrix $U$ such that $U^{T}\theta U=\Lambda$,
where
\[
\varLambda\equiv diag(\lambda_{1},\cdots,\lambda_{n})
\]
and $\lambda_{1},\cdots,\lambda_{n}$ are eigenvalues of $\theta$. 

2. Suppose, in addition, that $\lambda_{1},\cdots,\lambda_{n}$ are
nonnegative. Define the symmetric matrix $A\equiv U\varLambda^{\frac{1}{2}}U^{T}$,
where $\varLambda^{\frac{1}{2}}=diag(\lambda_{1}^{\frac{1}{2}},\cdots,\lambda_{n}^{\frac{1}{2}})$.
Then $\theta=AA^{T}$.
\end{thm}
\begin{proof}
1. Proceed by induction on $n$. The assertion is trivial if $n=1$.
Suppose the assertion has been proved for $n-1$. Recall that, for
an arbitrary unit vector $v_{n}$, there exist $v_{1},\cdots,v_{n-1}\in R^{n}$
such that $v_{1},\cdots,v_{n-1},v_{n}$ form an orthonormal basis
of $R^{n}$. Now let $v_{n}$ be an eigenvector of $\theta$ corresponding
to $\lambda_{n}$. Let $V$ be the $n\times n$ matrix whose $i$-th
column is $v_{i}$ for each $i=1,\cdots,n$. Then $V$ is an orthogonal
matrix. Define an $(n-1)\times(n-1)$ symmetric matrix $\eta$ by
$\eta_{i,j}\equiv v_{i}^{T}\theta v_{j}$ for each $i,j=1,\cdots,n-1$.
By the induction hypothesis, there exists an $(n-1)\times(n-1)$ orthogonal
matrix
\[
W\equiv\left[\begin{array}{ccc}
w_{1,1}, & \cdots, & w_{1,n-1}\\
\cdot & \cdots & \cdot\\
\cdot & \cdots & \cdot\\
\cdot & \cdots & \cdot\\
w_{n-1,1}, & \cdots, & w_{n-1,n-1}
\end{array}\right]
\]
such that 
\begin{equation}
W^{T}\eta W=\Lambda_{n-1}=diag(\lambda_{1},\cdots,\lambda_{n-1})\label{eq:temp-483}
\end{equation}
for some $\lambda_{1},\cdots,\lambda_{n-1}\in R$. Define the $n\times n$
matrices
\[
W'\equiv\left[\begin{array}{cccc}
w_{1,1}, & \cdots, & w_{1,n-1}, & 0\\
\cdot & \cdots & \cdot & \cdot\\
\cdot & \cdots & \cdot & \cdot\\
\cdot & \cdots & \cdot & \cdot\\
w_{n-1,1}, & \cdots, & w_{n-1,n-1}, & 0\\
0, & \cdots, & 0, & 1
\end{array}\right]
\]
and $U\equiv VW'$. Then it is easily verified that $U$ is orthogonal.
Moreover, 
\[
U^{T}\theta U=W'^{T}V^{T}\theta VW'=W'^{T}\left[\begin{array}{cccc}
v_{1}^{T}\theta v_{1}, & \cdots, & v_{1}^{T}\theta v_{n-1}, & v_{1}^{T}\theta v_{n}\\
\cdot & \cdots & \cdot & \cdot\\
\cdot & \cdots & \cdot & \cdot\\
\cdot & \cdots & \cdot & \cdot\\
v_{n-1}^{T}\theta v_{1}, & \cdots, & v_{n-1}^{T}\theta v_{n-1}, & v_{n-1}^{T}\theta v_{n}\\
v_{n}^{T}\theta v_{1}, & \cdots, & v_{n}^{T}\theta v_{n-1}, & v_{n}^{T}\theta v_{n}
\end{array}\right]W'
\]
\[
=\left[\begin{array}{cccc}
w_{1,1}, & \cdots, & w_{1,n-1}, & 0\\
\cdot & \cdots & \cdot & \cdot\\
\cdot & \cdots & \cdot & \cdot\\
\cdot & \cdots & \cdot & \cdot\\
w_{n-1,1}, & \cdots, & w_{n-1,n-1}, & 0\\
0, & \cdots, & 0, & 1
\end{array}\right]^{T}\left[\begin{array}{cccc}
\eta_{1,1}, & \cdots, & \eta_{1,n-1}, & 0\\
\cdot & \cdots & \cdot & \cdot\\
\cdot & \cdots & \cdot & \cdot\\
\cdot & \cdots & \cdot & \cdot\\
\eta_{n-1,1}, & \cdots, & \eta_{n-1,n-1}, & 0\\
0, & \cdots, & 0, & \lambda_{n}
\end{array}\right]
\]
\[
\left[\begin{array}{cccc}
w_{1,1}, & \cdots, & w_{1,n-1}, & 0\\
\cdot & \cdots & \cdot & \cdot\\
\cdot & \cdots & \cdot & \cdot\\
\cdot & \cdots & \cdot & \cdot\\
w_{n-1,1}, & \cdots, & w_{n-1,n-1}, & 0\\
0, & \cdots, & 0, & 1
\end{array}\right]
\]
\[
=\left[\begin{array}{cccc}
\lambda_{1}, & \cdots, & 0, & 0\\
\cdot & \cdots & \cdot & \cdot\\
\cdot & \cdots & \cdot & \cdot\\
\cdot & \cdots & \cdot & \cdot\\
0 & \cdots, & \lambda_{n-1}, & 0\\
0 & \cdots, & 0, & \lambda_{n}
\end{array}\right]\equiv\varLambda\equiv diag(\lambda_{1},\cdots,\lambda_{n}),
\]
where the fourth equality is thanks to equality \ref{eq:temp-483}.
Induction is completed. The equality $U^{T}\theta U=\Lambda$ implies
that $\theta U=U\Lambda$ and that $\lambda_{i}$ is an eigenvalue
of $\theta$ with an eigenvector given by the $i$-th column of $U$.
Assertion 1 is thus proved.

Since 
\[
\theta=U\varLambda U^{T}=U\varLambda^{\frac{1}{2}}\varLambda^{\frac{1}{2}}U^{T}=U\varLambda^{\frac{1}{2}}U^{T}U\varLambda^{\frac{1}{2}}U^{T}=AA^{T},
\]
Assertion 2 is proved.
\end{proof}
\begin{defn}
\label{Def. n-dim Normal Distribution, positive definite } \textbf{(Normal
distribution with positive definite covariance).} Let $n\geq1$ and
$\bar{\mu}\in R^{n}$ be arbitrary. Let $\overline{\sigma}$ be an
arbitrary positive definite $n\times n$ matrix. Then the function
defined on $R^{n}$ by
\begin{equation}
\varphi_{\bar{\mu},\overline{\sigma}}(y)\equiv(2\pi)^{-\frac{n}{2}}(\det\overline{\sigma})^{-\frac{1}{2}}\exp(-\frac{1}{2}(y-\bar{\mu})^{T}\overline{\sigma}^{-1}(y-\bar{\mu}))\label{eq:temp-36}
\end{equation}
for each $y\in R^{n}$ is a p.d.f. Let $\Phi_{\bar{\mu},\overline{\sigma}}$
be the corresponding distribution on $R^{n}$, and let $Y\equiv(Y_{1},\cdots,Y_{n})$
be any r.v. with values in $R^{n}$ and with $\Phi_{\bar{\mu},\overline{\sigma}}$
as its distribution. Then $\varphi_{\bar{\mu},\overline{\sigma}},\Phi_{\bar{\mu},\overline{\sigma}},Y$,
and $Y_{1},\cdots,Y_{n}$ are respectively said to be the \emph{normal
p.d.f \index{normal p.d.f}, normal distribution \index{normal distribution},
normally distributed \index{normally distributed}, and jointly normal}
\index{jointly normal}, with mean $\bar{\mu}$ and covariance matrix
$\overline{\sigma}$. Proposition \ref{Prop. Basics of Normal distributions}
below justifies the terminology. The p.d.f. $\varphi_{0,I}$ and the
distribution $\Phi_{0,I}$ are said to be \index{standard normal p.d.f.}\emph{
standard normal}, where $I$ is the identity matrix. 

In the case where $n=1$, define $\sigma\equiv\sqrt{\overline{\sigma}}$
and write $\Phi_{\mu,\sigma^{2}}$ also for the P.D.F. associated
with the distribution $\Phi_{\mu,\sigma^{2}}$, and call it a \index{normal P.D.F.}\emph{normal
P.D.F.} Thus $\Phi_{0,1}(x)=\int_{-\infty}^{x}\varphi_{0,1}(u)du$
for each $x\in R$. 

In Definition \ref{Def. n-dim Normal Distribution, nonegative definite},
we will generalize the definition of normal distributions to an arbitrary
nonnegative definite matrix $\overline{\sigma}$. $\square$
\end{defn}
\begin{prop}
\label{Prop. Basics of std normal pdf} \textbf{\emph{(Basics of standard
normal distribution).}} Consider the case $n=1$. Then the following
holds.

1. The function $\varphi_{0,1}$ on $R$ defined by 
\[
\varphi_{0,1}(x)\equiv\frac{1}{\sqrt{2\pi}}\exp(-\frac{1}{2}x^{2})
\]
is a p.d.f. on $R$ relative to the Lebesgue measure. Thus $\Phi_{0,1}$
is a P.D.F. on R.

2. Write $\Phi\equiv\Phi_{0,1}$. We will call $\Psi\equiv1-\Phi:[0,\infty)\rightarrow(0,\frac{1}{2}]$
the tail of $\Phi$. Then $\Phi(-x)=1-\Phi(x)$ for each $x\in R$.
Moreover, \emph{
\[
\Psi(x)\leq e^{-x^{2}/2}
\]
}for each $x\geq0$.

3. The inverse $\bar{\Psi}:(0,\frac{1}{2}]\rightarrow[0,\infty)$
of $\Psi$ is a decreasing function from $(0,1)$ to $R$ such that
$\bar{\Psi}(\varepsilon)\rightarrow\infty$ as $\varepsilon\rightarrow0$.
Moreover $\bar{\Psi}(\varepsilon)\leq\sqrt{-2\log\varepsilon}$ for
$\varepsilon\in(0,\frac{1}{2}]$.
\end{prop}
\begin{proof}
1. We calculate 
\[
(\frac{1}{\sqrt{2\pi}}\int_{-\infty}^{+\infty}e^{-x^{2}/2}dx)^{2}=\frac{1}{2\pi}\int_{-\infty}^{+\infty}\int_{-\infty}^{+\infty}e^{-(x^{2}+y^{2})/2}dxdy
\]
\[
=\frac{1}{2\pi}\int_{0}^{2\pi}\int_{0}^{+\infty}e^{-r^{2}/2}rdrd\theta=\frac{1}{2\pi}\int_{0}^{2\pi}(-e^{-r^{2}/2})|_{0}^{+\infty}d\theta
\]
\begin{equation}
=\frac{1}{2\pi}\int_{0}^{2\pi}d\theta=1,\label{eq:temp-29-1}
\end{equation}
where the change of variables from $(x,y)$ to $(r,\theta)$ is defined
by $x=r\cos\theta$ and $y=r\sin\theta$. Thus $\varphi_{0,1}$ is
Lebesgue integrable, with integral equal to 1, hence a p.d.f. on $R$. 

In the above proof, we used a series of steps: (i) the function $e^{-x^{2}/2}$
is integrable relative to the Lebesgue integration $J\equiv\int\cdot dx$,
(ii) the function $e^{-(x^{2}+y^{2})/2}$ is integrable relative to
$\int\int\cdot dxdy\equiv J\otimes J$, (iii) Fubini's Theorem equates
the double integral to successive integrals in either order, (iv)
a disk $D_{a}$ with center $0$ and radius $a>0$ is integrable relative
to $J\otimes J$, (v) the double integral $\int\int e^{-(x^{2}+y^{2})/2}dxdy$
is equal to the limit of $\int\int1_{D(a)}(x,y)e^{-(x^{2}+y^{2})/2}dxdy$
as $a\rightarrow\infty$, and (vi) we make a change of integration
variables from $(x,y)$ to $(r,\theta)$ in the last double integral.
Step (i) follows from an estimate of $|c_{a}-c_{a'}|\rightarrow0$
as $a,a'\rightarrow\infty$, where  $c_{a}\equiv\int_{0}^{a}e^{-x^{2}/2}dx$.
Step (ii) is justified by Corollary \ref{Prop. Cartesian product of integrable funcs is integrable}.
The use of Fubini's Theorem in step (iii) is justified by the conclusion
of step (ii). The integrability of $D_{a}$ in step (iv) follows because
$D_{a}=(Z\leq a)1_{[-a,a]\times[-a,a]}$ where $Z$ is the continuous
function defined by $Z(x,y)\equiv\sqrt{x^{2}+y^{2}}$. Step (v) is
an application of the Monotone Convergence Theorem. Step (vi), the
change of integration variables from rectangular- to polar coordinates,
is by Corollary \ref{Cor. Integrating in Polar Coordinates} in the
Appendix. In the remainder of this book, such slow motion, blow-by-blow
justifications will mostly be left to the reader.

2. Note first that 
\[
\Phi(-x)\equiv\int_{-\infty}^{-x}\varphi(u)du
\]
\[
=\int_{x}^{\infty}\varphi(-v)dv=\int_{x}^{\infty}\varphi(v)dv=1-\Phi(x),
\]
where we made a change of integration variables $v=-u$ and noted
that $\varphi(-v)=\varphi(v)$.

Next, if $x\in[\frac{1}{\sqrt{2\pi}},\infty)$, then
\[
\Psi(x)\equiv\frac{1}{\sqrt{2\pi}}\int_{x}^{\infty}e^{-u^{2}/2}du
\]
\[
\leq\frac{1}{\sqrt{2\pi}}\int_{x}^{\infty}\frac{u}{x}e^{-u^{2}/2}du=\frac{1}{\sqrt{2\pi}}\frac{1}{x}e^{-x^{2}/2}\leq e^{-x^{2}/2}.
\]
 On the other hand, if $x\in[0,\frac{1}{\sqrt{2\pi}})$, then
\[
\Psi(x)\leq\Psi(0)=\frac{1}{2}<\exp(-(\frac{1}{\sqrt{2\pi}})^{2}/2)\leq e^{-x^{2}/2}.
\]
Therefore, by continuity, $\Psi(x)\leq e^{-x^{2}/2}$ for each $x\in[0,\infty)$. 

3. Consider any $\varepsilon\in(0,1)$. Define $x\equiv\sqrt{-2\log\varepsilon}$.
Then $\Psi(x)\leq e^{-x^{2}/2}=\varepsilon$ by Assertion 2. Since
$\bar{\Psi}$ is a decreasing function, it follows that 
\[
\sqrt{-2\log\varepsilon}\equiv x=\bar{\Psi}(\Psi(x))\geq\bar{\Psi}(\varepsilon).
\]
\end{proof}
\begin{prop}
\label{Prop. Moments of standard normal rrv}\textbf{\emph{ (Moments
of standard normal r.r.v.)}} Suppose a r.r.v. $X$ has the standard
normal distribution $\Phi_{0,1}$, with p.d.f. $\varphi_{0,1}(x)\equiv\frac{1}{\sqrt{2\pi}}e^{-x^{2}/2}$.
Then $X^{m}$ is integrable for each $m\geq0$, Moreover, for each
even integer $m\equiv2k\geq0$ we have 
\[
EX^{m}=EX^{2k}=(2k-1)(2k-3)\cdots3\cdot1=(2k)!2^{-k}/k!,
\]
while $EX^{m}=0$ for each odd integer $m>0$. 
\end{prop}
\begin{proof}
Let $m\geq0$ be any even integer. Let $a>0$ b arbitrary. Then, integrating
by parts, we have 

\begin{equation}
\frac{1}{\sqrt{2\pi}}\int_{-a}^{a}x^{m+2}e^{-x^{2}/2}dx=\frac{1}{\sqrt{2\pi}}(-x^{m+1}e^{-x^{2}/2})|_{-a}^{a}+(m+1)\frac{1}{\sqrt{2\pi}}\int_{-a}^{a}x^{m}e^{-x^{2}/2}dx\label{eq:temp-9}
\end{equation}
Since the function $g_{m}$ defined by $g_{m}(x)\equiv1_{[-a,a]}(x)x^{m}$
for each $x\in R$ is Lebesgue $\mathrm{measurable}$ and is bounded,
Proposition \ref{Prop. p.d.f.  basics} implies that $g_{m}$ is integrable
relative to the P.D.F. $\Phi_{0,1}$, which has $\varphi_{0,1}$ as
p.d.f.. Moreover, according to Proposition \ref{Prop. p.d.f.  basics},
equality \ref{eq:temp-9} can be re-written as
\[
\int g_{m+2}(x)d\Phi_{0,1}(x)=\frac{1}{\sqrt{2\pi}}(-x^{m+1}e^{-x^{2}/2})|_{-a}^{a}+(m+1)\int g_{m}(x)d\Phi_{0,1}(x)
\]
or, in view of Proposition \ref{Prop.  F_X  is indeed a PDF}, as
\begin{equation}
Eg_{m+2}(X)=\frac{1}{\sqrt{2\pi}}(-x^{m+1}e^{-x^{2}/2})|_{-a}^{a}+(m+1)Eg_{m}(X)\label{eq:temp-20}
\end{equation}
The Lemma is trivial for $m=0$. Suppose the Lemma has been prove
for integers up to and including the even integer $m\equiv2k-2$.
By the induction hypothesis, $X^{m}$ is integrable. At the same time,
$g_{m}(X)\rightarrow X^{m}$ in probability as $a\rightarrow\infty$.
Hence, by the Dominated Convergence Theorem, we have $Eg_{m}(X)\uparrow EX^{m}$
as $a\rightarrow\infty$. Since $|a|^{m+1}e^{-a^{2}/2}\rightarrow0$
as $a\rightarrow\infty$, equality \ref{eq:temp-20} yields $Eg_{m+2}(X)\uparrow(m+1)EX^{m}$
as $a\rightarrow\infty$. The Monotone Convergence Theorem therefore
implies that $X^{m+2}$ is integrable, with $EX^{m+2}=(m+1)EX^{m}$,
or 
\[
EX^{2k}=(2k-1)EX^{2k-2}=\cdots=(2k-1)(2k-3)\cdots1=(2k)!2^{-k}/k!
\]
Since $X^{m+2}$ is integrable, so is $X^{m+1}$, according to Lyapunov's
inequality. Moreover, 
\[
EX^{m+1}=\int x^{m+1}d\Phi_{0,1}(x)=\int x^{m+1}\varphi_{0,1}(x)dx=0
\]
since $x^{m+1}\varphi_{0,1}(x)$ is an odd function of $x\in R$.
Induction is completed. 
\end{proof}
The next proposition shows that $\varphi_{\bar{\mu},\overline{\sigma}}$
and $\Phi_{\bar{\mu},\overline{\sigma}}$ in Definition \ref{Def. n-dim Normal Distribution, positive definite }
are well defined. 
\begin{prop}
\label{Prop. Basics of Normal distributions} \textbf{\emph{(Basics
of normal distributions with positive definite covariance).}} Let
$n\geq1$ and $\bar{\mu}\in R^{n}$ be arbitrary. Let $\overline{\sigma}$
be an arbitrary positive definite $n\times n$ matrix. Use the notations
in Definition \ref{Def. n-dim Normal Distribution, positive definite }.
Then the following holds.

1. $\varphi_{\bar{\mu},\overline{\sigma}}$ is indeed a p.d.f. on
$R^{n}$, i.e. $\int\varphi_{\bar{\mu},\overline{\sigma}}(x)dx=1$,
where $\int\cdot dx$ stands for the Lebesgue integration on $R^{n}.$
Thus the corresponding distribution $\Phi_{\bar{\mu},\overline{\sigma}}$
on $R^{n}$ is well defined. Moreover, $\Phi_{\bar{\mu},\overline{\sigma}}$
is equal to the distribution of the r.v. $Y\equiv\bar{\mu}+AX$ where
$A$ is an arbitrary $n\times n$ matrix with $\overline{\sigma}=AA^{T}$
and where $X$ is an arbitrary r.v. with values in $R^{n}$ and with
the standard normal distribution $\Phi_{0,I}$. In short, linear combinations
of a finite set of standard normal r.r.v.'s and the constant $1$
are jointly normal. More generally, linear combinations of a finite
set of jointly normal r.r.v.'s are jointly normal. 

2. Let\textup{\emph{ $Z\equiv(Z_{1},\cdots,Z_{n})$ be a r.v. }}with
values in $R^{n}$\textup{\emph{ with distribution $\Phi_{\bar{\mu},\overline{\sigma}}$.
Then }}$EZ=\bar{\mu}$ and $E(Z-\bar{\mu})(Z-\bar{\mu})^{T}=\overline{\sigma}$.

3. Let $Z_{1},\cdots,Z_{n}$ be jointly normal r.r.v.'s. Then $Z_{1},\cdots,Z_{n}$
are independent iff they are pairwise uncorrelated. In particular,
if $Z_{1},\cdots,Z_{n}$ are jointly standard normal, then they are
independent.
\end{prop}
\begin{proof}
For each $x\equiv(x_{1},\cdots,x_{n})\in R^{n}$, we have, by Definition
\ref{Def. n-dim Normal Distribution, positive definite },
\[
\varphi_{0,I}(x_{1},\cdots,x_{n})\equiv\varphi_{0,I}(x)\equiv(2\pi)^{-\frac{n}{2}}\exp(-\frac{1}{2}x^{T}x)
\]
\[
=\prod_{i=1}^{n}(\frac{1}{\sqrt{2\pi}}\exp(-\frac{1}{2}x_{i}^{2}))=\varphi_{0,1}(x_{1})\cdots\varphi_{0,1}(x_{n}).
\]
Since $\varphi_{0,1}$ is a p.d.f. on $R$ according to Proposition
\ref{Prop. Basics of std normal pdf} above, Proposition \ref{Prop. Cartesian product of integrable funcs is integrable}
implies that the Cartesian product $\varphi_{0,I}$ is a p.d.f. on
$R^{n}$. Let $X\equiv(X_{1},\cdots,X_{n})$ be an arbitrary r.v.
with values in $R^{n}$ and with p.d.f. $\varphi_{0,I}$. Then
\[
Ef_{1}(X_{1})\cdots f_{n}(X_{n})
\]
\[
=\int\cdots\int f_{1}(x_{1})\cdots f_{n}(x_{n})\varphi_{0,1}(x_{1})\cdots\varphi_{0,1}(x_{n})dx_{1}\cdots dx_{n}
\]
\begin{equation}
=\prod_{i=1}^{n}\int f_{i}(x_{i})\varphi_{0,1}(x_{i})dx_{i}=Ef_{1}(X_{1})\cdots Ef_{n}(X_{n})\label{eq:temp-29}
\end{equation}
for each $f_{1},\cdots,f_{n}\in C(R)$. Separately, for each $i=1,\cdots,n$,
the r.r.v. $X_{i}$ has distribution $\varphi_{0,1}$, whence $X_{i}$
has $m$-th moment for each $m\geq0$, with $EX_{i}^{m}=0$ if $m$
is odd, according to Proposition \ref{Prop. Moments of standard normal rrv}.

1. Next let $\overline{\sigma},\bar{\mu}$ be as given. Let $A$ be
an arbitrary $n\times n$ matrix such that $\overline{\sigma}=AA^{T}$.
By \ref{Thm. Spectral Theorem for symmetric matrices}, such a matrix
$A$ exists. Then $\det(\overline{\sigma})=\det(A)^{2}$. Since $\overline{\sigma}$
is positive definite, it is nonsingular and so is $A$. Let $X$ be
an arbitrary r.v. with values in $R^{n}$ and with the standard normal
distribution $\Phi_{0,I}$. Define the r.v. $Y\equiv\bar{\mu}+AX$.
Then, for arbitrary $f\in C(R^{n})$, we have
\[
Ef(Y)=Ef(\bar{\mu}+AX)=\int f(\bar{\mu}+Ax)\varphi_{0,I}(x)dx
\]
\[
\equiv(2\pi)^{-\frac{n}{2}}\int f(\bar{\mu}+Ax)\exp(-\frac{1}{2}x^{T}x)dx
\]
\[
=(2\pi)^{-\frac{n}{2}}\det(A)^{-1}\int f(y)\exp(-\frac{1}{2}(y-\bar{\mu})^{T}(A^{-1})^{T}A^{-1}(y-\bar{\mu}))dy
\]
\[
=(2\pi)^{-\frac{n}{2}}\det(\overline{\sigma})^{-\frac{1}{2}}\int f(y)\exp(-\frac{1}{2}(y-\bar{\mu})^{T}\overline{\sigma}^{-1}(y-\bar{\mu}))dy
\]
\[
\equiv\int f(y)\varphi_{\bar{\mu},\overline{\sigma}}(y)dy,
\]
where the fourth equality is by the change of integration variables
$y=\bar{\mu}+Ax$. Thus $\varphi_{\bar{\mu},\overline{\sigma}}$ is
the p.d.f. on $R^{n}$ of the r.v. $Y$, and $\Phi_{\bar{\mu},\overline{\sigma}}$
is the distribution of $Y$. 

2. Next, let\emph{ $Z_{1},\cdots,Z_{n}$ }be jointly normal r.r.v.'s
with distribution $\Phi_{\bar{\mu},\overline{\sigma}}$. By Assertion
1, there exist a standard normal r.v. $X\equiv(X_{1},\cdots,X_{n})$
on some probability space $(\Omega',L',E')$, and an $n\times n$
matrix $AA^{T}=\overline{\sigma}$, such that $E'f(\bar{\mu}+AX)=\Phi_{\bar{\mu},\overline{\sigma}}(f)=Ef(Z)$
for each $f\in C(R^{n})$. Thus $Z$ and $Y\equiv\bar{\mu}+AX$ induce
the same distribution on $R^{n}$. Let $i,j=1,\cdots,n$ be arbitrary.
Since \emph{$X_{i},X_{j},X_{i}X_{j}$} and, therefore, $Y_{i},Y_{j},Y_{i}Y_{j}$
are integrable, so are $Z_{i},Z_{j},Z_{i}Z_{j}$, with ,\emph{ }
\[
EZ=E'Y=\bar{\mu}+AE'X=\bar{\mu},
\]
and
\[
E(Z-\bar{\mu})(Z-\bar{\mu})^{T}=E'(Y-\bar{\mu})(Y-\bar{\mu})^{T}=AE'XX^{T}A^{T}=AA^{T}=\overline{\sigma}.
\]

3. Suppose $Z_{1},\cdots,Z_{n}$ are pairwise uncorrelated. Then $\overline{\sigma}_{i,j}=E(Z_{i}-\bar{\mu}_{i})(Z_{j}-\bar{\mu}_{j})=0$
for each $i,j=1,\cdots,n$ with $i\neq j$. Thus $\overline{\sigma}$
and $\overline{\sigma}^{-1}$ are diagonal matrices, with $(\overline{\sigma}^{-1})_{i,j}=\overline{\sigma}_{i,i}$or
$0$ according as $i=j$ or not. Hence, for each $f_{1},\cdots,f_{n}\in C(R)$,
we have
\[
Ef_{1}(Z_{1})\cdots f_{n}(Z_{n})
\]
\[
=(2\pi)^{-\frac{n}{2}}(\det\overline{\sigma})^{-\frac{1}{2}}\int\cdots\int f(z_{1})\cdots f(z_{n})\exp(-\frac{1}{2}(z-\bar{\mu})^{T}\overline{\sigma}^{-1}(z-\bar{\mu}))dz_{1}\cdots dz_{n}
\]
\[
=(2\pi)^{-\frac{n}{2}}(\overline{\sigma}_{1,1}\cdots\overline{\sigma}_{n,n})^{-\frac{1}{2}}\int\cdots\int f(z_{1})\cdots f(z_{n})\exp\sum_{i=1}^{n}(-\frac{1}{2}(z_{i}-\bar{\mu}_{i})\overline{\sigma}_{i,i}^{-1}(z_{i}-\bar{\mu}_{i}))dz_{1}\cdots dz_{n}
\]
\[
=(2\pi\overline{\sigma}_{i,i})^{-\frac{1}{2}}\int f(z_{i})\exp(-\frac{1}{2}(z_{i}-\bar{\mu}_{i})\overline{\sigma}_{i,i}^{-1}(z_{i}-\bar{\mu}_{i}))dz_{i}
\]
\[
=Ef_{1}(Z_{1})\cdots Ef_{n}(Z_{n}).
\]
We conclude that $Z_{1},\cdots,Z_{n}$ are independent if they are
pairwise uncorrelated. The converse is trivial.
\end{proof}
Next we generalize the definition of normal distribution to include
the case where the covariance matrix nonnegative definite.
\begin{defn}
\label{Def. n-dim Normal Distribution, nonegative definite} \textbf{(Normal
distribution with nonnegative definite covariance).} Let $n\geq1$
and $\bar{\mu}\in R^{n}$ be arbitrary. Let $\overline{\sigma}$ be
an arbitrary nonnegative definite $n\times n$. Define the \index{normal distribution}\emph{normal
distribution} $\Phi_{\bar{\mu},\overline{\sigma}}$on $R^{n}$ by
\begin{equation}
\Phi_{\bar{\mu},\overline{\sigma}}(f)\equiv\lim_{\varepsilon\rightarrow0}\Phi_{\bar{\mu},\overline{\sigma}+\varepsilon I}(f)\label{eq:temp-295}
\end{equation}
for each $f\in C(R^{n})$, where, for each $\varepsilon>0$, the function
$\Phi_{\bar{\mu},\overline{\sigma}+\varepsilon I}$ is the normal
distribution on $R^{n}$ introduced in Definition \ref{Def. n-dim Normal Distribution, positive definite }
for the positive definite matrix $\overline{\sigma}+\varepsilon I$.
Lemma \ref{Lem. Normal distribution is well defined for nonneg definite matrix }
below proves that $\Phi_{\bar{\mu},\overline{\sigma}}$ well defined
and is indeed a distribution. 

A sequence $Z_{1},\cdots,Z_{n}$ of r.r.v.'s is said to be \emph{jointly
normal,} with $\Phi_{\bar{\mu},\overline{\sigma}}$ as distribution,
if $Z\equiv(Z_{1},\cdots,Z_{n})$ has the distribution $\Phi_{\bar{\mu},\overline{\sigma}}$
on $R^{n}$. 
\end{defn}
\begin{lem}
\label{Lem. Normal distribution is well defined for nonneg definite matrix }
\textbf{\emph{(Normal distribution with nonnegative definite covariance
is well defined).}} Use the notations and assumptions in Definition
\ref{Def. n-dim Normal Distribution, nonegative definite}. Then the
following holds.

1. The the limit $\lim_{\varepsilon\rightarrow0}\Phi_{\bar{\mu},\overline{\sigma}+\varepsilon I}(f)$
in equality \ref{eq:temp-295} exists for each $f\in C(R^{n})$. Moreover,
$\Phi_{\bar{\mu},\overline{\sigma}}$ is the distribution of $Y\equiv\bar{\mu}+AX$
for some standard normal $X\equiv(X_{1},\cdots,X_{n})$ and some $n\times n$
matrix $A$ with $AA^{T}=\overline{\sigma}$.

2. If $\overline{\sigma}$ is positive definite, then $\Phi_{\bar{\mu},\overline{\sigma}}(f)=\int f(y)\varphi_{\bar{\mu},\overline{\sigma}}(y)dy$,
where $\varphi_{\bar{\mu},\overline{\sigma}}$ was defined in Definition
\ref{Def. n-dim Normal Distribution, positive definite }. Thus Definition
\ref{Def. n-dim Normal Distribution, nonegative definite} of $\Phi_{\bar{\mu},\overline{\sigma}}$
for a nonnegative definite $\overline{\sigma}$ is consistent with
the previous Definition \ref{Def. n-dim Normal Distribution, positive definite }
for a positive definite $\overline{\sigma}$.

3. Let \textup{$Z\equiv(Z_{1},\cdots,Z_{n})$} be an arbitrary r.v.
with values in $R^{n}$ and with distribution $\Phi_{\bar{\mu},\overline{\sigma}}$.
Then $Z_{1}^{k(1)}\cdots Z_{n}^{k(n)}$ is integrable for each $k_{1},\cdots,k_{n}\geq0$.
In particular, $Z$ has mean $\bar{\mu}$ and covariance matrix $\overline{\sigma}$.
\end{lem}
\begin{proof}
1. Let $\varepsilon>0$ be arbitrary. Then $\overline{\sigma}+\varepsilon I$
is positive definite. Hence, the normal distribution $\Phi_{\bar{\mu},\overline{\sigma}+\varepsilon I}$
has been defined. Separately, Theorem \ref{Thm. Spectral Theorem for symmetric matrices}
implies that there exists an orthogonal matrix $U$ such that $U^{T}\overline{\sigma}U=\Lambda$,
where $\varLambda\equiv diag(\lambda_{1},\cdots,\lambda_{n})$ is
a diagonal matrix whose diagonal elements consist of the eigenvalues
$\lambda_{1},\cdots,\lambda_{n}$ of $\overline{\sigma}$. These eigenvalues
are nonnegative since $\overline{\sigma}$ is nonnegative definite.
Hence, again by Theorem \ref{Thm. Spectral Theorem for symmetric matrices},
we have
\begin{equation}
\overline{\sigma}+\varepsilon I=A_{\varepsilon}A_{\varepsilon}^{T},\label{eq:temp-303}
\end{equation}
where 
\begin{equation}
A_{\varepsilon}\equiv U\varLambda_{\varepsilon}^{\frac{1}{2}}U^{T},\label{eq:temp-41}
\end{equation}
where $\varLambda_{\varepsilon}^{\frac{1}{2}}=diag(\sqrt{\lambda_{1}+\varepsilon},\cdots,\sqrt{\lambda_{n}+\varepsilon})$. 

Now let $X$ be an arbitrary r.v. on $R^{n}$ with the standard normal
distribution $\Phi_{0,I}$. In view of equality \ref{eq:temp-303},
Proposition \ref{Prop. Basics of Normal distributions} implies that
$\Phi_{\bar{\mu},\overline{\sigma}+\varepsilon I}$ is equal to the
distribution of the r.v.
\[
Y^{(\varepsilon)}\equiv\bar{\mu}+A_{\varepsilon}X.
\]
Define $A\equiv U\varLambda^{\frac{1}{2}}U^{T}$, where $\varLambda^{\frac{1}{2}}=diag(\sqrt{\lambda_{1}},\cdots,\sqrt{\lambda_{n}})$
and define $Y\equiv\bar{\mu}+AX$. Then
\[
E|A_{\varepsilon}X-AX|^{2}=EX^{T}(A_{\varepsilon}-A)^{T}(A_{\varepsilon}-A)X
\]
\[
=\sum_{i=1}^{n}\sum_{j=1}^{n}\sum_{k=1}^{n}EX_{i}U_{i,j}(\sqrt{\lambda_{j}+\varepsilon}-\sqrt{\lambda_{j}})^{2}U_{j,k}X_{k}
\]
\[
=\sum_{i=1}^{n}\sum_{j=1}^{n}U_{i,j}(\sqrt{\lambda_{j}+\varepsilon}-\sqrt{\lambda_{j}})^{2}U_{j,i}
\]
\[
=\sum_{j=1}^{n}(\sqrt{\lambda_{j}+\varepsilon}-\sqrt{\lambda_{j}})^{2}\sum_{i=1}^{n}U_{i,j}U_{j,i}
\]
\[
=\sum_{j=1}^{n}(\sqrt{\lambda_{j}+\varepsilon}-\sqrt{\lambda_{j}})^{2}\rightarrow0
\]
as $\varepsilon\rightarrow0$. Lyapunov's inequality then implies
that
\[
E|Y^{(\varepsilon)}-Y|=E|A_{\varepsilon}X-AX|\leq(E|A_{\varepsilon}X-AX|^{2})^{\frac{1}{2}}\rightarrow0
\]
as $\varepsilon\rightarrow0$. In other words, $Y^{(\varepsilon)}\rightarrow Y$
in probability. Consequently, the distribution $\Phi_{\bar{\mu},\overline{\sigma}+\varepsilon I}$
converges to the distribution $F_{Y}$ of $Y$. We conclude that the
limit $\Phi_{\bar{\mu},\overline{\sigma}}(f)$ in equality \ref{eq:temp-295}
exists and is equal to $EF(Y)$. In other words, $\Phi_{\bar{\mu},\overline{\sigma}}$
is the distribution of the r.v. $Y\equiv\bar{\mu}+AX$. Moreover,
\[
AA^{T}=U\varLambda^{\frac{1}{2}}U^{T}U\varLambda^{\frac{1}{2}}U^{T}=U\varLambda U^{T}=\overline{\sigma}.
\]
Assertion 1 is proved.

2. Next suppose $\overline{\sigma}$ is positive definite. Then $\varphi_{\bar{\mu},\overline{\sigma}+\varepsilon I}\rightarrow\varphi_{\bar{\mu},\overline{\sigma}}$
uniformly on compact subsets of $R^{n}$. Hence 
\[
\lim_{\varepsilon\rightarrow0}\Phi_{\bar{\mu},\overline{\sigma}+\varepsilon I}(f)=\lim_{\varepsilon\rightarrow0}\int f(y)\varphi_{\bar{\mu},\overline{\sigma}+\varepsilon I}(y)dy=\int f(y)\varphi_{\bar{\mu},\overline{\sigma}}(y)dy
\]
for each $f\in C(R^{n})$. Therefore Definition \ref{Def. n-dim Normal Distribution, nonegative definite}
is consistent with Definition \ref{Def. n-dim Normal Distribution, positive definite },
proving Assertion 2.

3. Now let $Z\equiv(Z_{1},\cdots,Z_{n})$ be any r.v. with values
in $R^{n}$ and with distribution $\Phi_{\bar{\mu},\overline{\sigma}}$.
By Assertion 1, $\Phi_{\bar{\mu},\overline{\sigma}}$ is the distribution
of $Y\equiv\bar{\mu}+AX$ for some standard normal $X\equiv(X_{1},\cdots,X_{n})$
and some $n\times n$ matrix $A$ with $AA^{T}=\overline{\sigma}$.
Thus $Z$ and $Y$ has the same distribution. Let $k_{1},\cdots,k_{n}\geq0$
be arbitrary. Then the r.r.v. $Y_{1}^{k(1)}\cdots Y_{n}^{k(n)}$ is
a linear combination of products $X_{1}^{j(1)}\cdots X_{n}^{j(n)}$
integrable where $j_{1},\cdots,j_{n}\geq0$, each of which is integrable
in view of Proposition \ref{Prop. Moments of standard normal rrv}
and Proposition \ref{Prop. Cartesian product of integrable funcs is integrable}.
Hence $Y_{1}^{k(1)}\cdots Y_{n}^{k(n)}$ is integrable. It follows
that $Z_{1}^{k(1)}\cdots Z_{n}^{k(n)}$ is integrable. $EZ=EY=\bar{\mu}$
and 
\[
E(Z-\bar{\mu})(Z-\bar{\mu})^{T}=E(Y-\bar{\mu})(Y-\bar{\mu})^{T}=EAXX^{T}A^{T}=AA^{T}=\overline{\sigma}.
\]
In other words, $Z$ has mean $\bar{\mu}$ and covariance matrix $\overline{\sigma}$,
proving Assertion 3.
\end{proof}
We will need some bounds related to the normal p.d.f. in later sections.

Recall from Proposition \ref{Prop. Basics of std normal pdf} the
standard normal P.D.F. $\Phi$ on $R$, its tail $\Psi$, and the
inverse function $\bar{\Psi}$ of the latter.
\begin{lem}
\label{Lem. Bounds for Integrals  relative to Normal PDF} \textbf{\emph{(Some
bounds for normal probabilities). }}
\end{lem}
\begin{enumerate}
\item \emph{Suppose $h$ is a measurable function on $R$ relative to the
Lebesgue integration. If $|h|\leq a$ on $[-\alpha,\alpha]$ and $|h|\leq b$
on $[-\alpha,\alpha]^{c}$ for some $a,b,\alpha>0$, then
\[
\int|h(x)|\varphi_{0,\sigma}(x)dx\leq a+2b\Psi(\frac{\alpha}{\sigma})
\]
for each $\sigma>0$.}
\item \emph{In general, let }$n\geq1$\emph{ be arbitrary. Let $I$ denote
the $n\times n$ identity matrix. Suppose $f$ is a Lebesgue integrable
function on $R^{n}$, with $|f|\leq1$. Let $\sigma>0$ be arbitrary.
Define a function $f_{\sigma}$ on $R^{n}$ by
\[
f_{\sigma}(x)\equiv\int_{y\in R^{n}}f(x-y)\varphi_{0,\sigma I}(y)dy
\]
for each $x\in R^{n}$. Suppose $f$ is continuous at some $t\in R^{n}$.
In other words, suppose, for arbitrary $\varepsilon>0$, there exists
$\delta_{f}(\varepsilon,t)>0$ such that $|f(t)-f(r)|<\varepsilon$
for each $r\in R^{n}$ with $|r-t|<\delta_{f}(\varepsilon,t)$. Let
$\varepsilon>0$ be arbitrary. Let }$\alpha\equiv\delta_{f}(\frac{\varepsilon}{2},t)>0$\emph{
and let
\begin{equation}
\sigma<\alpha/\bar{\Psi}(\frac{1}{2}(1-(1-\frac{\varepsilon}{4})^{\frac{1}{n}})).\label{eq:temp-10}
\end{equation}
Then
\[
|f_{\sigma}(t)-f(t)|\leq\varepsilon.
\]
}
\item \emph{Again consider the case $n=1$. Let $\varepsilon>0$ be arbitrary.
Suppose $\sigma>0$ is so small that $\sigma<\varepsilon/\bar{\Psi}(\frac{\varepsilon}{8})$.
Let $r,s\in R$ be arbitrary with $r+2\varepsilon<s$. Let $f\equiv1_{(r,s]}$.
Then $1_{(r+\varepsilon,s-\varepsilon]}-\varepsilon\leq f_{\sigma}\leq1_{(r-\varepsilon,s+\varepsilon]}+\varepsilon$.}
\end{enumerate}
\begin{proof}
1. We estimate
\[
\frac{1}{\sqrt{2\pi}\sigma}\int|h(x)|e^{-x^{2}/(2\sigma^{2})}dx
\]
\[
\leq\frac{1}{\sqrt{2\pi}\sigma}\int_{-\alpha}^{\alpha}ae^{-x^{2}/(2\sigma^{2})}dx+\frac{1}{\sqrt{2\pi}\sigma}\int_{|x|>\alpha}be^{-x^{2}/(2\sigma^{2})}dx
\]
\[
\leq a+\frac{1}{\sqrt{2\pi}}\int_{|u|>\frac{\alpha}{\sigma}}be^{-u^{2}/2}du=a+2b\Psi(\frac{\alpha}{\sigma})
\]

2. Let $f,t,\varepsilon,\delta_{f},$ $\alpha$, and $\sigma$ be
as given. Then inequality \ref{eq:temp-10} implies that
\[
(1-(1-\frac{\varepsilon}{4})^{\frac{1}{n}})>2\Psi(\frac{\alpha}{\sigma}),
\]
whence
\begin{equation}
2(1-(1-2\Psi(\frac{\alpha}{\sigma}))^{n})<\frac{\varepsilon}{2}.\label{eq:temp-55}
\end{equation}
Then, $|f(t-u)-f(t)|<\frac{\varepsilon}{2}$ for $u\in R^{n}$ with
$\left\Vert u\right\Vert \equiv|u_{1}|\vee\cdots\vee|u_{n}|<\alpha$.
By hypothesis $\sigma\leq\alpha/\bar{\Psi}(\frac{\varepsilon}{8})$.
Hence $\frac{\alpha}{\sigma}\geq\bar{\Psi}(\frac{\varepsilon}{8})$
and so $\Psi(\frac{\alpha}{\sigma})\leq\frac{\varepsilon}{8}$. Hence,
by Assertion 1, we have
\[
|f_{\sigma}(t)-f(t)|=|\int(f(t-u)-f(t))\varphi_{0,\sigma I}(u)du|
\]
\[
\leq\int_{u:\left\Vert u\right\Vert <\alpha}|f(t-u)-f(t)|\varphi_{0,\sigma I}(u)du+\int_{u:\left\Vert u\right\Vert \geq\alpha}|f(t-u)-f(t)|\varphi_{0,\sigma I}(u)du
\]
\[
\leq\frac{\varepsilon}{2}+2(1-\int_{u:\left\Vert u\right\Vert <\alpha}\varphi_{0,\sigma I}(u)du)
\]
\[
=\frac{\varepsilon}{2}+2(1-(\Phi(\frac{\alpha}{\sigma})-\Phi(\frac{\alpha}{\sigma}))^{n})
\]
\[
=\frac{\varepsilon}{2}+2(1-(1-2\Psi(\frac{\alpha}{\sigma}))^{n})
\]
\[
<\frac{\varepsilon}{2}+\frac{\varepsilon}{2}=\varepsilon,
\]
as desired, where the last inequality is from inequality \ref{eq:temp-55}. 

3. Define $f_{\sigma}(x)\equiv\int f(x-y)\varphi_{0,\sigma}(y)dy$
for each $x\in R$. Consider each $t\in(r+\varepsilon,s-\varepsilon]$.
Then $f$ is constant in a neighborhood of $t$, hence continuous
at $t$. More precisely, let $\delta_{f}(\theta,t)\equiv\varepsilon$
for each $\theta>0$. Then 
\[
(t-\delta_{f}(\theta),t+\delta_{f}(\theta))=(t-\varepsilon,t+\varepsilon)\subset(r,s]\subset(f=1)
\]
for each $\theta>0$. Let $\alpha\equiv\delta_{f}(\frac{\varepsilon}{2},t)\equiv\varepsilon$.
Then, by hypothesis\emph{
\begin{equation}
\sigma<\varepsilon\bar{\Psi}(\frac{\varepsilon}{8})^{-1}=\alpha/\bar{\Psi}(\frac{1}{2}(1-(1-\frac{\varepsilon}{4}))).\label{eq:temp-10-1}
\end{equation}
}Hence, by Assertion 2, we have
\[
|f(t)-f_{\sigma}(t)|\leq\varepsilon,
\]
where $t\in(r+\varepsilon,s-\varepsilon]$ is arbitrary. Since $1_{(r+\varepsilon,s-\varepsilon]}(t)\leq f(t)$,
it follows that 
\begin{equation}
1_{(r+\varepsilon,s-\varepsilon]}(t)-\varepsilon\leq f_{\sigma}(t)\label{eq:temp-34}
\end{equation}
for each $t\in(r+\varepsilon,s-\varepsilon]$. Since $f_{\sigma}\geq0$,
inequality \ref{eq:temp-34} is trivially satisfied for $t\in(-\infty,r+\varepsilon]\cup(s-\varepsilon,\infty)$.
We have thus proved that inequality \ref{eq:temp-34} holds on $domain(1_{(r+\varepsilon,s-\varepsilon]})$.
Next consider any $t\in(-\infty,r-\varepsilon]\cup(s+\varepsilon,\infty)$.
Again, for arbitrary $\theta>0$ we have $|f(t)-f(u)|=0<\theta$ for
each $u\in(t-\delta_{f}(\theta),t+\delta_{f}(\theta))$. Hence, by
Assertion 2, we have $f_{\sigma}(t)=f_{\sigma}(t)-f(t)<\varepsilon$.
It follows that
\begin{equation}
f_{\sigma}(t)\leq1_{(r-\varepsilon,s+\varepsilon]}(t)+\varepsilon\label{eq:temp-35}
\end{equation}
for each $t\in(-\infty,r-\varepsilon]\cup(s+\varepsilon,\infty)$.
Since $f_{\sigma}\leq1$, inequality \ref{eq:temp-35} is trivially
satisfied for $t\in(r-\varepsilon,s+\varepsilon]$. We have thus proved
that inequality \ref{eq:temp-35} holds on $domain(1_{(r-\varepsilon,s+\varepsilon]})$.
Assertion 3 is proved. 
\end{proof}

\section{Characteristic Functions}

In previous sections we analyzed distributions $J$ on a locally compact
metric space $(S,d)$ in terms of their values $Jg$ at basis functions
$g$ in a partition of unity. In the special case where $(S,d)$ is
the Euclidean space $R,$ the basis functions can be replaced by the
exponential functions $h_{\lambda}$with $\lambda\in R$, where $h_{\lambda}(x)\equiv e^{i\lambda x}$,
where $i\equiv\sqrt{-1}$. The result is characteristic functions,
a most useful in the study of distributions of r.r.v.'s. 

The classical development of this tool in usual texts, e.g.  \cite{Chung68}
or \cite{Loeve60} is constructive, except for infrequent and superficial
appeals to the principle of infinite search. The bare essentials of
this material is presented here for completeness and for ease of reference.
The reader who is familiar with the topic and is comfortable that
the classical treatment is constructive, or easily made so, can skip
over this and the next section and come back for reference.

We will be working with complex-valued $\mathrm{measurable}$ functions.
Let $\mathbb{C}$ denote the complex plane equipped with the usual
metric.
\begin{defn}
\textbf{\emph{\label{Def. Complex valued integrable function}}}\textbf{(Complex
valued integrable function)}\textbf{\emph{.}} Let $I$ be an integration
on a locally compact metric space $(S,d)$, and let $(S,\Lambda,I)$
denote the completion of the integration space $(S,C(S),I)$. A function
$X\equiv IU+iIV:S\rightarrow\mathbb{C}$ whose real part $U$ and
imaginary part $V$ are $\mathrm{measurable}$ on $(S,\Lambda,I)$
is said to be $\mathrm{measurable}$\index{mathrm{measurable} function, complex valued@$\mathrm{measurable}$ function, complex valued}
on $(S,\Lambda,I)$. If both $U,V$ are integrable, then $X$ is said
to be \index{integrable function, complex valued}\emph{integrable},
with integral $IX\equiv IU+iIV$. $\square$
\end{defn}
By separation into real and imaginary parts, the complex-valued functions
immediately inherit the bulk of the theory of integration developed
hitherto in this book for real-valued functions. One exception is
the very basic inequality $|IX|\leq I|X|$ when $|X|$ is integrable.
Its trivial proof in the case of real valued integrable functions
relies on the linear ordering of $R$, which is absent in $\mathbb{C}$.
The next lemma gives a proof for complex valued integrable functions.
\begin{lem}
\label{Lem. Complex |IX|<=00003DI|X|} \textbf{\emph{$(|IX|\leq I|X|$
for complex valued integrable function $X$).}} Use the notations
in Definition \ref{Def. Complex valued integrable function}. Let
$X:S\rightarrow\mathbb{C}$ be an arbitrary complex valued function.
Then the function $X$ is \textup{\emph{measurable}}\emph{ }in the
sense of Definition \ref{Def. Complex valued integrable function}
iff it is \textup{\emph{measurable}} in the sense of Definition \ref{Def. Complex valued integrable function}.
In other words, the former is consistent with the latter. Moreover,
if $X$ is \textup{\emph{measurable}} and if $|X|\in L$, then $X$
is integrable with $|IX|\leq I|X|$.
\end{lem}
\begin{proof}
Write $X\equiv IU+iIV$, where $U,V$ are the real and imaginary parts
of $X$ respectively. 

1. Suppose $X$ is $\mathrm{measurable}$ in the sense of Definition
\ref{Def. Complex valued integrable function}. Then $U,V:(S,\Lambda,I)\rightarrow R$
are $\mathrm{measurable}$ functions. Therefore the function $(U,V):(S,\Lambda,I)\rightarrow R^{2}$
is $\mathrm{measurable}$. At the same time, we have $X=f(U,V)$,
where the continuous function $f:R^{2}\rightarrow\mathbb{C}$ is defined
by $f(u,v)\equiv u+iv$. Hence $X$ is $\mathrm{measurable}$ in the
sense of Definition \ref{Def. Measurable Function}, according to
Proposition \ref{Prop. X meas,  f unif continuous and bd on bd subsets =00003D> f(X) meas}. 

Conversely, suppose $X(S,\Lambda,I)\rightarrow\mathbb{C}$ is $\mathrm{measurable}$
in the sense of Definition \ref{Def. Measurable Function}. Note that
$U,V$ are continuous functions of $X$. Hence, again by Proposition
\ref{Prop. X meas,  f unif continuous and bd on bd subsets =00003D> f(X) meas},
both $U,V$ are $\mathrm{measurable}$. Thus $X$ is $\mathrm{measurable}$
in the sense of Definition \ref{Def. Complex valued integrable function}.

2. Suppose $X$ is $\mathrm{measurable}$ and $|X|\in L$. Then, by
Definition \ref{Def. Complex valued integrable function}, both $U$
and $V$ are $\mathrm{measurable}$, with $|U|\vee|V|\leq|X|\in L$,
it follows that $U,V\in L$. Thus $X$ is integrable according to
Definition \ref{Def. Complex valued integrable function}. 

Let $\varepsilon>0$ be arbitrary. Then either (i) $I|X|<3\varepsilon$,
or (ii) $I|X|>2\varepsilon$.

First consider Case (i). Then 
\[
|IX|=|IU+iIV|\leq|IU|+|iIV|\leq I|U|+I|V|\leq2I|X|<I|X|+3\varepsilon.
\]

Now consider Case (ii). By the Dominated Convergence Theorem, there
exists $a>0$ so small that $I(|X|\wedge a)<\varepsilon$. Then 
\begin{equation}
I|X|1_{(|X|\leq a)}\leq I(|X|\wedge a)<\varepsilon.\label{eq:temp-463}
\end{equation}
Write $A\equiv(a<|X|)$. Then 
\[
|IU1_{A}-IU|=|IU1_{(|X|\leq a)}|\leq I|U|1_{(|X|\leq a)}\leq I|X|1_{(|X|\leq a)}<\varepsilon.
\]
Similarly, $|IV1_{A}-IV|<\varepsilon$. Hence
\begin{equation}
|I(X1_{A})-IX|=|I(U1_{A}-IU)+i(IV1_{A}-IV)|<2\varepsilon.\label{eq:temp-461}
\end{equation}
Write $c\equiv I|X|1_{A}$. Then it follows that 
\[
c\equiv I|X|1_{A}=I|X|-I|X|1_{(|X|\leq a)}>2\varepsilon-2\varepsilon=0,
\]
where the inequality is on account of Condition (ii) and inequality
\ref{eq:temp-463}. Now define a probability integration space $(S,L,E)$
using $g\equiv c^{-1}|X|1_{A}$ as a probability density function
on the integration space $(S,\Lambda,I)$. Thus 
\[
E(Y)\equiv c^{-1}I(Y|X|1_{A})
\]
for each $Y\in L$. Then 
\[
|c^{-1}I(X1_{A})|\equiv|E(\frac{X}{|X|\vee a}1_{A})|=|E(\frac{U}{|X|\vee a}1_{A})+iE(\frac{V}{|X|\vee a}1_{A})|
\]
\[
=((E(\frac{U}{|X|\vee a}1_{A}))^{2}+(E(\frac{V}{|X|\vee a}1_{A}))^{2})^{\frac{1}{2}}
\]
\[
\leq(E(\frac{U^{2}}{(|X|\vee a)^{2}}1_{A})+E(\frac{V^{2}}{(|X|\vee a)^{2}}1_{A}))^{\frac{1}{2}}=(E(\frac{|X|^{2}}{(|X|\vee a)^{2}}1_{A}))^{\frac{1}{2}}\leq1,
\]
where the inequality is thanks to Lyapunov. Hence $|I(X1_{A})|\leq c\equiv I|X|1_{A}$.
Inequality \ref{eq:temp-461} therefore yields 
\[
|IX|<|I(X1_{A})|+2\varepsilon\leq I|X|1_{A}+2\varepsilon<I|X|+3\varepsilon.
\]
Summing up, we have $|IX|<I|X|+3\varepsilon$ regardless of Case (i)
or Case (ii), where $\varepsilon>0$ is arbitrary. We conclude that
$I|X|\leq I|X|$.
\end{proof}
\begin{lem}
\label{Lem. bounds for |1-exp(ix)| etc} \textbf{\emph{(Basic inequalities
for exponentials).}} Let $x,y,x',y'\in R$ be arbitrary, with $y\leq0$
and $y'\leq0$. Then 
\begin{equation}
|e^{ix}-1|\leq2\wedge|x|\label{eq:temp-478}
\end{equation}
and
\[
|e^{ix+y}-e^{ix'+y'}|\leq2\wedge|x-x'|+1\wedge|y-y'|.
\]
\end{lem}
\begin{proof}
If $x\geq0$, then
\[
|e^{ix}-1|^{2}=|\cos x-1+i\sin x|^{2}=2(1-\cos x)
\]
\begin{equation}
=2\int_{0}^{x}\sin udu\leq2\int_{0}^{x}udu\leq x^{2}.\label{eq:temp-478-1}
\end{equation}
Hence, by symmetry and continuity, $|e^{ix}-1|\leq|x|$ for arbitrary
$x\in R$. At the same time, $|e^{ix}-1|\leq2$. Equality \ref{eq:temp-478}
follows. 

Now assume $y\geq y'$. 
\[
|e^{ix+y}-e^{ix'+y'}|\leq|e^{ix+y}-e^{ix'+y}|+|e^{ix'+y}-e^{ix'+y'}|
\]
\[
\leq|e^{ix}-e^{ix'}|e^{y}+|e^{y}-e^{y'}|\leq|e^{i(x-x')}-1|e^{y}+e^{y}(1-e^{-(y-y')})
\]
\[
\leq(2\wedge|x-x'|)e^{y}+(1-e^{-(y-y')})\leq2\wedge|x-x'|+1\wedge|y-y'|
\]
\textbackslash{}where the last inequality is because $y'\leq y\leq0$
by assumption. Hence, by symmetry and continuity, the same inequality
holds for arbitrary $y,y'\leq0$.
\end{proof}
Recall the matrix notations and basics from Definition \ref{Def. Matrix notations}.
Moreover, we will write $|x|\equiv(x_{1}^{2}+\cdots+x_{n}^{2})^{\frac{1}{2}}$
and write $\left\Vert x\right\Vert \equiv|x_{1}|\vee\cdots\vee|x_{n}|$
for each $x\equiv(x_{1},\cdots,x_{n})\in R^{n}$.
\begin{defn}
\label{Def. Characteristic functions and Fourier Transforms} \textbf{(Characteristic
function, Fourier transform, and convolution).} Let $n\geq1$ be arbitrary.

1. Let $X\equiv(X_{1},\cdots,X_{n}$) be a r.v. with values in $R^{n}$.
The \emph{\index{characteristic function}characteristic function}
of $X$ is the complex-valued function $\psi_{X}$ on $R^{n}$ defined
by
\[
\psi_{X}(\lambda)\equiv E\exp i\lambda^{T}X\equiv E\cos(\lambda^{T}X)+iE\sin(\lambda^{T}X).
\]
for each $\lambda\in R^{n}$. 

2. Let $J$ be an arbitrary distribution on $R^{n}$. The \emph{characteristic
function} of $J$ is defined to be $\psi_{J}\equiv\psi_{X}$, where
$X$ is any r.v. with values in $R^{n}$ such that $E_{X}=J$. Thus
$\psi_{J}(\lambda)\equiv Jh_{\lambda}$ where $h_{\lambda}(x)\equiv\exp i\lambda^{T}X$
for each $\lambda,x\in R^{n}$. 

3. If $g$ is a complex-valued integrable function on $R^{n}$ relative
to the Lebesgue integration, the \emph{Fourier transform} \index{Fourier transform}
of $g$ is defined to be the complex valued function $\hat{g}$ on
$R^{n}$ with
\[
\hat{g}(\lambda)\equiv\int_{x\in R^{n}}(\exp i\lambda^{T}x)g(x)dx
\]
for $\lambda\in R^{n},$ where $\int\cdot dx$ signifies the Lebesgue
integration on $R^{n}$, and where $x\in R^{n}$ is the integration
variable. The \index{convolution}\emph{convolution} of two complex-valued
Lebesgue integrable functions $f,g$ on $R^{n}$ is the complex valued
function $f\star g$ defined by $(f\star g)(x)\equiv\int_{y\in R^{n}}f(x-y)g(y)dy$
for each $x\in R^{n}$.

4. Suppose $n=1$. Let $F$ be an P.D.F. on $R$. The \emph{characteristic
function} of $F$ is defined as $\psi_{F}\equiv\psi_{J}$ , where
$J\equiv\int\cdot dF$. If, in addition, $F$ has a p.d.f. $f$ on
$R$, then the characteristic function of $f$ is defined as $\psi_{f}\equiv\psi_{F}$.
In that case, $\psi_{F}(\lambda)=\int e^{i\lambda t}f(t)dt\equiv\hat{f}(\lambda)$
for each $\lambda\in R$.

$\square$

We can choose to express the characteristic function in terms of the
r.v. $X$, or in terms of the distribution $J$, or, in the case $n=1$,
the P.D.F., as a matter of convenience. A theorem proved in one set
of notations will be used in another set without further comment.
\end{defn}
\begin{lem}
\label{Lem. Basics of Convolutions} \textbf{\emph{(Basics of convolution).}}
Let $f,g,h$ be complex-valued Lebesgue integrable functions on $R^{n}$.
Then the following holds.
\end{lem}
\begin{enumerate}
\item \emph{$f\star g$ is Lebesgue integrable.}
\item \emph{$f\star g=g\star f$}
\item \emph{$(f\star g)\star h=f\star(g\star h)$ }
\item \emph{$(af+bg)\star h=a(f\star h)+b(g\star h)$ for all complex numbers
$a,b$.}
\item \emph{Suppose $n=1$, and suppose $g$ is a p.d.f}.\emph{ If $|f|\leq a$
for some $a\in R$ then $|f\star g|\leq a$. If $f$ is real-valued
with $a\leq f\leq b$ for some $a,b\in R$, then $a\leq f\star g\leq b$.}
\item \emph{$\widehat{f\star g}=\hat{f}\hat{g}$}
\item $|\hat{f}|\leq\left\Vert f\right\Vert _{1}\equiv\int_{x\in R^{n}}|f(x)|dx$
\end{enumerate}
\begin{proof}
If $f$ and $g$ are real-valued, then the integrability of $f\star g$
follows from Corollary \ref{Cor. If X,Y integrable, then  X(u-v)Y(v) integrable}
in the appendices. Assertion 1 then follows by linearity. We will
prove Assertions 6 and 7, the remaining assertions  left as an exercise.
For Assertion 6, note that, for each $\lambda\in R^{n}$, we have
\[
\widehat{f\star g}(\lambda)\equiv\int(\exp i\lambda^{T}x)(\int f(x-y)g(y)dy)dx
\]
\[
=\int_{-\infty}^{\infty}(\int_{-\infty}^{\infty}(\exp i\lambda^{T}(x-y))f(x-y)dx)(\exp i\lambda^{T}y)g(y)dy
\]
\[
=\int_{-\infty}^{\infty}(\int_{-\infty}^{\infty}(\exp i\lambda^{T}u)f(u))du)(\exp i\lambda^{T}y)g(y)dy
\]
\[
=\int\hat{f}(\lambda)(\exp i\lambda^{T}y)g(y)dy=\hat{f}(\lambda)\hat{g}(\lambda),
\]
as asserted. At the same time, for each $\lambda\in R^{n}$, we have
\[
|\hat{f}(\lambda)|\equiv|\int_{x\in R^{n}}(\exp i\lambda^{T}x)f(x)dx|
\]
\[
\leq\int_{x\in R^{n}}|(\exp i\lambda^{T}x)f(x)|dx|=\int_{x\in R^{n}}|f(x)|dx,
\]
where the inequality is by Lemma \ref{Lem. Complex |IX|<=00003DI|X|}.
Assertion 7 is verified.
\end{proof}
\begin{prop}
\label{Prop. Unifrom continuity of Char fcts} \textbf{\emph{(Uniform
continuity of characteristic functions).}} Let $X$ be a r.v. with
values in $R^{m}$. Let $\beta_{X}$ be a modulus of tightness of
$X$. Then the following holds.
\end{prop}
\begin{enumerate}
\item \emph{$|\psi_{X}(\lambda)|\leq1$ and $\psi_{a+BX}(\lambda)=\exp(i\lambda^{T}a)\psi_{X}(\lambda^{T}B)$
for each $a,\lambda\in R^{n}$ and for each $n\times m$ matrix $B$.}
\item \emph{$\psi_{X}$ is }uniform\emph{ly continuous. More precisely,
$\psi_{X}$ has a modulus of continuity given by $\delta(\varepsilon)\equiv\frac{\varepsilon}{3}/\beta(\frac{\varepsilon}{3})$
for $\varepsilon>0$.}
\item \emph{If $g$ is a Lebesgue integrable function on $R^{n}$, then
$\hat{g}$ is }uniform\emph{ly continuous. More precisely, for each
$\varepsilon>0$ there exists $\gamma\equiv\gamma_{g}(\varepsilon)>0$
so large that }$\int1_{(|x|>\gamma)}|g(x)|dx<\varepsilon$.\emph{
Then a modulus of continuity of $\hat{g}$ is given by $\delta(\varepsilon)\equiv\frac{\varepsilon}{\left\Vert g\right\Vert +2}/\gamma_{g}(\frac{\varepsilon}{\left\Vert g\right\Vert +2})$
for $\varepsilon>0$, where $\left\Vert g\right\Vert \equiv\int|g(t)|dt$.}
\end{enumerate}
\begin{proof}
1. For each $\lambda\in R^{n}$ we have $|\psi_{X}(\lambda)|\equiv|E\exp(i\lambda^{T}X)|\leq E|\exp(i\lambda^{T}X)|=E1=1.$
Moreover 
\[
\psi_{a+BX}(\lambda)=\exp(i\lambda^{T}a)E(i\lambda^{T}BX)=\exp(i\lambda^{T}a)\psi_{X}(\lambda^{T}B)
\]
\emph{for each $a,\lambda\in R^{n}$ and for each $n\times m$ matrix
$B$}.

2. Let $\varepsilon>0$. Let $\delta(\varepsilon)\equiv\frac{\varepsilon}{3}/\beta(\frac{\varepsilon}{3})$.
Suppose $h\in R^{n}$ is such that $|h|<\delta(\varepsilon)$. Then
$\beta(\frac{\varepsilon}{3})<\frac{\varepsilon}{3|h|}$. Pick $a\in(\beta(\frac{\varepsilon}{3}),\frac{\varepsilon}{3|h|})$.
Then $P(|X|>a)<\frac{\varepsilon}{3}$ by the definition of $\beta$.
On the other hand, for each $x\in R^{n}$ with $|x|\leq a$, we have
$|\exp(ih^{T}x)-1|\leq|h^{T}x|\leq|h|a<\frac{\varepsilon}{3}$. Hence,
for each $\lambda\in R^{n}$,
\[
|\psi_{X}(\lambda+h)-\psi_{X}(\lambda)|\leq E|\exp(i\lambda^{T}X)(\exp(ih^{T}X)-1)|
\]
\[
\leq E(|\exp(ih^{T}X)-1|;|X|\leq a)+2P(|X|>a)<\frac{\varepsilon}{3}+2\frac{\varepsilon}{3}=\varepsilon
\]

3. Proceed in the same manner as above. Let $\varepsilon>0$. Write
$\varepsilon'\equiv\frac{\varepsilon}{\left\Vert g\right\Vert +2}$.
Let $\delta(\varepsilon)\equiv\frac{\varepsilon'}{\gamma_{g}(\varepsilon')}$.
Suppose $h\in R^{n}$ is such that $|h|<\delta(\varepsilon)$. Then
$\gamma_{g}(\varepsilon')<\frac{\varepsilon'}{|h|}$. Pick $a\in(\gamma_{g}(\varepsilon'),\frac{\varepsilon'}{|h|})$.
Then $\int_{|x|>a}|g(x)|dx<\varepsilon'$ by the definition of $\gamma_{g}$.
Moreover, for each $x\in R^{n}$ with $|x|\leq a$, we have $|\exp(ih^{T}x)-1|\leq|h^{T}x|\leq|h|a<\varepsilon'$.
Hence, for each $\lambda\in R^{n}$,
\[
|\hat{g}(\lambda+h)-\hat{g}(\lambda)|\leq\int|\exp(i\lambda^{T}x)(\exp(ih^{T}x)-1)g(x)|dx
\]
\[
\leq\int_{|x|\leq a}|(\exp(ih^{T}x)-1)g(x)|dx+\int_{|x|>a}|2g(x)|dx
\]
\[
\leq\varepsilon'\int|g(x)|dx+2\varepsilon'=\varepsilon'(\left\Vert g\right\Vert +2)=\varepsilon.
\]
\end{proof}
\begin{lem}
\label{Lem. Char Fct of  Normal} \textbf{\emph{(Characteristic function
of normal distribution).}} Let $\Phi_{\bar{\mu},\overline{\sigma}}$
be an arbitrary normal distribution on $R^{n}$, with mean $\bar{\mu}$
and covariance matrix $\overline{\sigma}$. Then the characteristic
function of $\Phi_{\bar{\mu},\overline{\sigma}}$is given by
\[
\psi_{\bar{\mu},\overline{\sigma}}(\lambda)\equiv\exp(i\bar{\mu}^{T}\lambda-\frac{1}{2}\lambda^{T}\overline{\sigma}\lambda)
\]
for each $\lambda\in R^{n}$. 
\end{lem}
\begin{proof}
1. Consider the special case where $n=1$, $\bar{\mu}=0$, and $\overline{\sigma}=1$.
Let $X$ be a r.r.v. with the standard normal distribution $\Phi_{0,1}$.
By \ref{Prop. Moments of standard normal rrv}, $X^{p}$ is integrable
for each $p\geq0$, with $m_{p}\equiv EX^{p}=(2k)!2^{-k}/k!$ if $p$
is equal to some even integer $2k$, and with $m_{p}\equiv EX^{p}=0$
otherwise. Using these moment formulas, we compute the characteristic
function
\[
\psi_{0,1}(\lambda)=\frac{1}{\sqrt{2\pi}}\int e^{i\lambda x}e^{-x^{2}/2}dx=\frac{1}{\sqrt{2\pi}}\int\sum_{p=0}^{\infty}\frac{(i\lambda x)^{p}}{p!}e^{-x^{2}/2}dx
\]
\[
=\sum_{p=0}^{\infty}\frac{(i\lambda)^{p}}{p!}m_{p}=\sum_{k=0}^{\infty}\frac{(-1)^{k}\lambda^{2k}}{(2k)!}m_{2k}
\]
\[
=\sum_{k=0}^{\infty}\frac{(-1)^{k}\lambda^{2k}}{(2k)!}(2k)!2^{-k}/k!=\sum_{k=0}^{\infty}\frac{(-\lambda^{2}/2)^{k}}{k!}=e^{-\lambda^{2}/2}
\]
where Fubini's Theorem justifies any change in the order of integration
and summation.

2. Now consider the general case. By Lemma \ref{Lem. Normal distribution is well defined for nonneg definite matrix },
$\Phi_{\bar{\mu},\overline{\sigma}}$ is the distribution of a r.v.
$Y=\bar{\mu}+AX$ for some matrix $A$ with $\overline{\sigma}\equiv AA^{T}$
and for some r.v. $X$ with the standard normal p.d.f. $\varphi_{0,I}$
on $R^{n}$, where $I$ is the $n\times n$ identity matrix. Let $\lambda\in R^{n}$
be arbitrary. Write $\theta\equiv A^{T}\lambda$. Then
\[
\psi_{\bar{\mu},\overline{\sigma}}(\lambda)\equiv E\exp(i\lambda^{T}Y)\equiv E\exp(i\lambda^{T}\bar{\mu}+i\lambda^{T}AX)
\]
\[
=\int_{x\in R^{n}}\exp(i\lambda^{T}\bar{\mu}+i\theta^{T}x)\varphi_{0,I}(x)dx
\]
\[
=\exp(i\lambda^{T}\bar{\mu})\int\cdots\int\exp(i\sum_{j=1}^{n}\theta_{j}x_{j})\varphi_{0,1}(x_{1})\cdots\varphi_{0,1}(x_{n})dx_{1}\cdots dx_{n},
\]
where we used Theorem \ref{Theorem. Change of Variables for integrable integrands}
for the change of integration variables. By Fubini's Theorem and by
the first part of this proof, this reduces to 
\[
\psi_{\bar{\mu},\overline{\sigma}}(\lambda)=\exp(i\lambda^{T}\bar{\mu})\prod_{j=1}^{n}(\int\exp(i\theta_{j}x_{j})\varphi_{0,1}(x_{j})dx_{j})
\]
\[
=\exp(i\lambda^{T}\bar{\mu})\prod_{j=1}^{n}\exp(-\frac{1}{2}\theta_{j}^{2})=\exp(i\lambda^{T}\bar{\mu})\exp(-\frac{1}{2}\theta^{T}\theta)
\]
\[
=\exp(i\lambda^{T}\bar{\mu}-\frac{1}{2}\lambda^{T}AA^{T}\lambda)\equiv\exp(i\lambda^{T}\bar{\mu}-\frac{1}{2}\lambda^{T}\overline{\sigma}\lambda).
\]
\end{proof}
\begin{cor}
\label{Cor. Convolution with normal density} \textbf{\emph{(Convolution
with normal density).}}\emph{ }Suppose $f$ is a Lebesgue integrable
function on $R^{n}$. Let $\sigma>0$ be arbitrary. Write $\overline{\sigma}\equiv\sigma^{2}I$,
where $I$ is the $n\times n$ identity matrix. Define $f_{\sigma}\equiv f\star\varphi_{0,\overline{\sigma}}$.
Then
\[
f_{\sigma}(t)=(2\pi)^{-n}\int\exp(i\lambda^{T}t-\frac{1}{2}\sigma^{2}\lambda^{T}\lambda)\hat{f}(-\lambda)d\lambda
\]
for each $t\in R^{n}$.
\end{cor}
\begin{proof}
In view of Lemma \ref{Lem. Char Fct of  Normal}, we have, for each
$t\in R^{n}$,
\[
f_{\sigma}(t)\equiv\int\varphi_{0,\overline{\sigma}}(t-x)f(x)dx=(2\pi\sigma^{2})^{-\frac{n}{2}}\int\exp(-\frac{1}{2\sigma^{2}}(t-x)^{T}(t-x))f(x)dx
\]
\[
=(2\pi\sigma^{2})^{-\frac{n}{2}}\int\psi_{0,I}(\sigma^{-1}(t-x))f(x)dx
\]
\[
=(2\pi\sigma^{2})^{-\frac{n}{2}}\int\int(2\pi)^{-\frac{n}{2}}\exp(i\sigma^{-1}(t-x)^{T}y-\frac{1}{2}y^{T}y)f(x)dxdy
\]
\[
=(2\pi\sigma)^{-n}\int\exp(-i\sigma^{-1}y^{T}t-\frac{1}{2}y^{T}y)\hat{f}(-\sigma^{-1}y)dy
\]
\[
=(2\pi)^{-n}\int\exp(-i\lambda^{T}t-\frac{1}{2}\sigma^{2}\lambda^{T}\lambda)\hat{f}(-\lambda)d\lambda
\]
Note that in the double integral, the integrand is a continuous function
in $(x,y)$ and is bounded in absolute value by a constant multiple
of $\exp(-\frac{1}{2}y^{T}y)f(x)$ which is, by Proposition \ref{Prop. Cartesian product of integrable funcs is integrable},
Lebesgue integrable on $R^{2n}$. This justifies the changes in order
of integration, thanks to Fubini. 
\end{proof}
The next theorem recovers a distribution on $R^{n}$ from its characteristic
function.
\begin{thm}
\label{Thm. Inversion  Formula  for Ch fcts} \textbf{\emph{(Inversion
formula for characteristic functions).}} Let $J,J'$ be distribution
on $R^{n}$, with characteristic functions $\psi_{J},\psi_{J'}$ respectively.
Let $f$ be an arbitrary Lebesgue integrable function on $R^{n}$.
Let $\hat{f}$ denote the Fourier Transform of $f$. Let $\sigma>0$
be arbitrary. Write $\overline{\sigma}\equiv\sigma^{2}I$, where $I$
is the $n\times n$ identity matrix. Define $f_{\sigma}\equiv f\star\varphi_{0,\overline{\sigma}}$.
Then the following holds. 
\end{thm}
\begin{enumerate}
\item \emph{We have 
\begin{equation}
Jf_{\sigma}=(2\pi)^{-n}\int\exp(-\frac{1}{2}\sigma^{2}\lambda^{T}\lambda)\hat{f}(-\lambda)\psi_{J}(\lambda)d\lambda.\label{eq:temp-466}
\end{equation}
 }
\item \emph{Suppose $f\in C_{ub}(R^{n})$ and $|f|\leq1$. Let $\varepsilon>0$
be arbitrary. Suppose $\sigma>0$ is so small that 
\[
\sigma\leq\delta_{f}(\frac{\varepsilon}{2})/\sqrt{-2\log(\frac{1}{2}(1-(1-\frac{\varepsilon}{4})^{\frac{1}{n}}))}.
\]
Then
\[
|Jf-Jf_{\sigma}|\leq\varepsilon.
\]
Consequently $Jf=\lim_{\sigma\rightarrow0}Jf_{\sigma}$.}
\item \emph{Suppose $f\in C_{ub}(R)$ is arbitrary such that $\hat{f}$
is Lebesgue integrable on }$R^{n}$. \emph{Then} 
\[
Jf=(2\pi)^{-n}\int\hat{f}(-\lambda)\psi_{J}(\lambda)d\lambda.
\]
\item \emph{If $\psi_{J}$ is Lebesgue integrable on $R^{n}$, then $J$
has a p.d.f. Specifically, then
\[
Jf=(2\pi)^{-n}\int f(x)\hat{\psi_{J}}(-x)dx
\]
 for each $f\in C(R^{n})$.}
\item $J=J'$ iff $\psi_{J}=\psi_{J'}$.
\end{enumerate}
\begin{proof}
Write $\left\Vert f\right\Vert \equiv\int|f(t)|dt<\infty$. Then $|\hat{f}|\leq\left\Vert f\right\Vert $. 

1. Consider the function $Z(\lambda,x)\equiv\exp(i\lambda^{T}x-\frac{1}{2}\sigma^{2}\lambda^{T}\lambda)\hat{f}(-\lambda)$
on the product space $(R^{n},L_{0},I_{0})\otimes(R^{n},L,J)$ where
$(R^{n},L_{0},I_{0})\equiv(R^{n},L_{0},\int\cdot dx)$ is the Lebesgue
integration space and where $(R^{n},L,J)$ is the probability space
that is the completion of $(R^{n},C_{ub}(R^{n}),J)$. The function
$Z$ is a continuous function of $(\lambda,x)$. Hence $Z$ is $\mathrm{measurable}$.
Moreover, $|Z|\leq U$ where $U(\lambda,x)\equiv\left\Vert f\right\Vert e^{-\sigma^{2}\lambda^{2}/2}$
is integrable. Hence $Z$ is integrable by the Dominated Convergence
Theorem. 

Define $h_{\lambda}(t)\equiv\exp(it^{T}\lambda)$ for each $t,\lambda\in R^{n}$.
Then $J(h_{\lambda})=\psi_{J}(\lambda)$ for each $\lambda\in R^{n}$.
Corollary \ref{Cor. Convolution with normal density} then implies
that 
\[
Jf_{\sigma}=(2\pi)^{-n}\int J(h_{\lambda})\exp(-\frac{1}{2}\sigma^{2}\lambda^{T}\lambda)\hat{f}(-\lambda)d\lambda
\]
\begin{equation}
\equiv(2\pi)^{-n}\int\psi_{J}(\lambda)\exp(-\frac{1}{2}\sigma^{2}\lambda^{T}\lambda)\hat{f}(-\lambda)d\lambda,\label{eq:temp-25}
\end{equation}
proving Assertion 1.

2. Now suppose $f\in C_{ub}(R^{n})$ with modulus of continuity $\delta_{f}$
with $|f|\leq1$. Recall that\emph{ $\bar{\Psi}:(0,\frac{1}{2}]\rightarrow[0,\infty)$}
denotes the inverse of the tail function $\Psi\equiv1-\Phi$ of the
standard normal P.D.F.\emph{ $\Phi$. }Proposition \ref{Prop. Basics of std normal pdf}
says that $\bar{\Psi}(\alpha)\leq\sqrt{-2\log\alpha}$ for each $\alpha\in(0,\frac{1}{2}]$.
Hence\emph{
\[
\sigma\leq\delta_{f}(\frac{\varepsilon}{2})/\sqrt{-2\log(\frac{1}{2}(1-(1-\frac{\varepsilon}{4})^{\frac{1}{n}}))}\leq\delta_{f}(\frac{\varepsilon}{2})/\bar{\Psi}(\frac{1}{2}(1-(1-\frac{\varepsilon}{4})^{\frac{1}{n}})),
\]
}where the first inequality is by hypothesis. Therefore Lemma \ref{Lem. Bounds for Integrals  relative to Normal PDF}
implies that $|f_{\sigma}-f|\leq\varepsilon$. Consequently $|Jf-Jf_{\sigma}|\leq\varepsilon$.
Hence \emph{$Jf=\lim_{\sigma\rightarrow0}Jf_{\sigma}$. }This proves
Assertion 2.

3. Now let $f\in C_{ub}(R^{n})$ be arbitrary. Then, by linearity,
Assertion 2 implies that
\begin{equation}
\lim_{\sigma\rightarrow0}Jf_{\sigma}=Jf.\label{eq:temp-484}
\end{equation}
\emph{
\begin{equation}
Jf_{\sigma}=(2\pi)^{-n}\int\exp(-\frac{1}{2}\sigma^{2}\lambda^{T}\lambda)\hat{f}(-\lambda)\psi_{J}(\lambda)d\lambda.\label{eq:temp-466-2}
\end{equation}
}Suppose $\hat{f}$ is Lebesgue integrable on $R^{n}$. Then the integrand
in equality \ref{eq:temp-466} is dominated in absolute value by the
integrable function $|\hat{f}|$, and converges a.u. on $R^{n}$ to
the function $\hat{f}(-\lambda)\psi_{J}(\lambda)$ as $\sigma\rightarrow0$.
Hence the Dominated Convergence Theorem implies that 
\[
\lim_{\sigma\rightarrow0}Jf_{\sigma}=(2\pi)^{-n}\int\hat{f}(-\lambda)\psi_{J}(\lambda)d\lambda.
\]
Combining with equality \ref{eq:temp-484}, Assertion 3 is proved.

4. Next consider the case where $\psi_{J}$ is Lebesgue integrable.
Suppose $f\in C(R^{n})$ with $|f|\leq1$. Then the function $U_{\sigma}(x,\lambda)\equiv f(x)\psi_{J}(\lambda)e^{-i\lambda x-\sigma^{2}\lambda^{2}/2}$
is an integrable function relative to the product Lebesgue integration
on $R^{2n}$, and is dominated in absolute value by the integrable
function $f(x)\psi_{J}(\lambda)$. Moreover, $U_{\sigma}\rightarrow U_{0}$
uniformly on compact subsets of $R^{2n}$ where $U_{0}(x,\lambda)\equiv f(x)\psi_{J}(\lambda)e^{-i\lambda x}$.
Hence $U_{\sigma}\rightarrow U_{0}$ in measure relative to $I_{0}\otimes I_{0}$.
The Dominated Convergence Theorem therefore yields, as $\sigma\rightarrow0$,
\[
Jf_{\sigma}=(2\pi)^{-n}\int\exp(-\frac{1}{2}\sigma^{2}\lambda^{T}\lambda)\hat{f}(-\lambda)\psi_{J}(\lambda)d\lambda
\]
\[
=(2\pi)^{-n}\int\exp(-\frac{1}{2}\sigma^{2}\lambda^{T}\lambda)\psi_{J}(\lambda)\int\exp(-i\lambda x)f(x)dxd\lambda
\]
\[
\rightarrow(2\pi)^{-n}\int\psi_{J}(\lambda)\int\exp(-i\lambda x)f(x)dxd\lambda
\]
\[
=(2\pi)^{-n}\int\hat{\psi}(-x)f(x)dx.
\]
On the other hand, by Assertion 2, we have $If_{\sigma}\rightarrow If$
as $\sigma\rightarrow0$. Assertion 4 is proved.

5. Assertion 5 follows from Assertion 4.
\end{proof}
\begin{defn}
\label{Def. Metric on characteristic functions} \textbf{(Metric of
characteristic functions). }Let $n\geq1$ be arbitrary. Let $\psi,\psi'$
be arbitrary characteristic functions on $R^{n}$. Define 
\begin{equation}
\rho_{char}(\psi,\psi')\equiv\rho_{ch,n}(\psi,\psi')\equiv\sum_{j=1}^{\infty}2^{-j}\sup_{|\lambda|\leq j}|\psi(\lambda)-\psi'(\lambda)|.\label{eq:temp-12}
\end{equation}
Then $\rho_{char}$ is a metric. $\square$

We have seen earlier that characteristic functions are continuous
and bounded in absolute values by $1$. Hence the supremum inside
the parentheses in equality \ref{eq:temp-12} exists and is bounded
by $2$. Thus $\rho_{char}$ is well-defined. In view of Theorem \ref{Thm. Inversion  Formula  for Ch fcts},
it is easily seen that $\rho_{char}$ is a metric. Convergence relative
to $\rho_{char}$ is equivalent to uniform convergence on each compact
subset of $R^{n}$.

The next theorem shows that the correspondence between distributions
on $R^{n}$ and their characteristic functions is uniformly continuous
when restricted to a tight subset. 
\end{defn}
\begin{thm}
\textbf{\emph{(Continuity Theorem for characteristic functions).}}
\label{Thm. Continuity Theorem for ch functions} Let $\xi$ be an
arbitrary binary approximation of $R$. Let $n\geq1$ be arbitrary,
but fixed. Let $\xi^{n}\equiv(A_{p})_{p=1,2,\cdots}$ be the  binary
approximation of $R^{n}$ which is the $n$-th power of $\xi$,. Let
$\left\Vert \xi^{n}\right\Vert $ be the modulus of local compactness
of $R^{n}$ associated with $\xi^{n}$. Let $\rho_{Dist,\xi^{n}}$
be the corresponding distribution metric on the space of distributions
on $R^{n}$, as in Definition \ref{Def. Distribution metric}. Let
$\widehat{J}_{0}$ be a family of distributions on $R^{n}$. 

Let $J,J'\in\widehat{J}_{0}$ be arbitrary, with corresponding characteristic
functions $\psi,\psi'$. Then the following holds.

1. For each $\varepsilon>0$, there exists $\delta_{ch,dstr}(\varepsilon,n)>0$
such that if $\rho_{char}(\psi,\psi')<\delta_{ch,dstr}(\varepsilon,n)$
then $\rho_{Dist,\xi^{n}}(J,J')<\varepsilon$.

2. Suppose $\widehat{J}_{0}$ is tight, with some modulus of tightness
$\beta$. Then, for each $\varepsilon>0$, there exists $\delta_{dstr,ch}(\varepsilon,\beta,\left\Vert \xi^{n}\right\Vert )>0$
such that if $\rho_{Dist,\xi^{n}}(J,J')<\delta_{dstr,ch}(\varepsilon,\beta)$
then $\rho_{char}(\psi,\psi')<\varepsilon$.

3. If $(J_{m})_{m=0,1,\cdots}$ is a sequence of distributions on
$R^{n}$ with a corresponding sequence $(\psi_{m})_{m=0,1,\cdots}$
of characteristic functions such that $\rho_{char}(\psi_{m},\psi_{0})\rightarrow0$,
then $J_{m}\Rightarrow J_{0}$.
\end{thm}
\begin{proof}
Let 
\begin{equation}
\pi_{R^{n}}\equiv(\{g_{p,x}:x\in A_{p}\})_{p=1,2,\cdots}\label{eq:temp-153-1-1}
\end{equation}
be the  partition of unity of $R^{n}$ determined by $\xi^{n}$, as
introduced in Definition \ref{Def. Partition of unity for locally compact (S,d)}.
Thus $\left\Vert \xi^{n}\right\Vert \equiv(|A_{p}|)_{p=1,2,\cdots}$.
Let 
\[
v_{n}\equiv\int_{|y|\leq1}dy,
\]
the volume of the unit sphere $\{y\in R^{n}:|y|\leq1\}$ in $R^{n}$.

1. Let $\varepsilon>0$ be arbitrary. For abbreviation, write 
\[
\alpha\equiv\frac{1}{8}\varepsilon.
\]
Let $p\equiv[0\vee(1-\log_{2}\varepsilon)]_{1}$. Thus
\[
2^{-p}<\frac{\varepsilon}{2}.
\]
For each\emph{ $\theta>0$} define 
\begin{equation}
\delta_{p}(\theta)\equiv2^{-p-1}\theta>0.\label{eq:temp-241}
\end{equation}
Recall from Proposition \ref{Prop. Basics of std normal pdf} the
standard normal P.D.F. $\Phi$ on $R$, its decreasing tail function$\Psi:[0,\infty)\rightarrow(0,\frac{1}{2}]$,
and the inverse function $\bar{\Psi}:(0,\frac{1}{2}]\rightarrow[0,\infty)$
of the latter. Define 
\begin{equation}
\sigma\equiv\delta_{p}(\frac{\alpha}{2})/\sqrt{-2\log(\frac{1}{2}(1-(1-\frac{\alpha}{4})^{\frac{1}{n}}))}>0.\label{eq:temp-42}
\end{equation}
Define
\[
m\equiv[\sigma^{-1}n^{\frac{1}{2}}\bar{\Psi}(\frac{1}{2}\wedge v_{n}^{-1}\varepsilon2^{-5}(2\pi)^{\frac{n}{2}}\sigma^{n}n^{-1})]_{1}.
\]
Thus $m\geq1$ is so large that 
\begin{equation}
v_{n}2^{2}(2\pi)^{-\frac{n}{2}}\sigma^{-n}n\Psi(\sigma n^{-\frac{1}{2}}m)<\frac{1}{8}\varepsilon.\label{eq:temp-27}
\end{equation}
Finally, define
\begin{equation}
\delta_{ch,dstr}(\varepsilon)\equiv\delta_{ch,dstr}(\varepsilon,n)\equiv v_{n}^{-1}\varepsilon2^{-m-3}(2\pi)^{\frac{n}{2}}\sigma^{n}>0.\label{eq:temp-57}
\end{equation}

Now suppose the characteristic functions $\psi,\psi'$ on $R^{n}$
are such that
\begin{equation}
\rho_{char}(\psi,\psi')\equiv\sum_{j=1}^{\infty}2^{-j}\sup_{\left|\lambda\right|\leq j}|\psi(\lambda)-\psi'(\lambda)|<\delta_{ch,dstr}(\varepsilon).\label{eq:temp-12-1}
\end{equation}
We will prove that $\rho_{Dist,\xi^{n}}(J,J')<\varepsilon$. To that
end, first note that, with $m\geq1$ as defined above, the last displayed
inequality implies
\begin{equation}
\sup_{\left|\lambda\right|\leq m}|\psi(\lambda)-\psi'(\lambda)|<2^{m}\delta_{ch,dstr}(\varepsilon).\label{eq:temp-26}
\end{equation}
Next, let $k=1,\cdots,p$ and $x\in A_{k}$ be arbitrary. Write $f\equiv g_{k,x}$
for abbreviation. Then, by Proposition \ref{Prop. Properties of parittion of unity-1},
$f$ has values in $[0,1]$ and has Lipschitz constant $2^{k+1}\leq2^{p+1}$.
Consequently, $f$ has the modulus of continuity $\delta_{p}$\emph{
}defined in equality \ref{eq:temp-241} above\emph{.} Hence, in view
of equality \ref{eq:temp-42}, Theorem \ref{Thm. Inversion  Formula  for Ch fcts}
implies that
\begin{equation}
|Jf-Jf_{\sigma}|\leq\alpha\equiv\frac{1}{8}\varepsilon,\label{eq:temp-31}
\end{equation}
and that
\begin{equation}
Jf_{\sigma}=(2\pi)^{-n}\int\exp(-\frac{1}{2}\sigma^{2}\lambda^{T}\lambda)\hat{f}(-\lambda)\psi(\lambda)d\lambda,\label{eq:temp-43}
\end{equation}
where $f_{\sigma}\equiv f\star\varphi_{0,\sigma^{2}I}$, where $I$
is the $n\times n$ identity matrix, and where $\hat{f}$ stands for
the Fourier transform of $f$. Moreover, by Proposition \ref{Prop. Properties of parittion of unity-1},
the function $f\equiv g_{k,x}$ has the sphere $\{y\in R^{n}:|y-x|\leq2^{-k+1}\}$
as support. Therefore 
\begin{equation}
|\hat{f}|\leq\int f(y)dy\leq v_{n}(2^{-k+1})^{n}<v_{n}.\label{eq:temp-56}
\end{equation}
where $v_{n}$ is the volume of the unit $n$-sphere in $R^{n}$,
as defined previously. 

By equality \ref{eq:temp-43} for $J$ and a similar equality for
$J'$, we have
\[
|Jf_{\sigma}-J'f_{\sigma}|=(2\pi)^{-n}|\int\exp(-\frac{1}{2}\sigma^{2}\lambda^{T}\lambda)\hat{f}(-\lambda)(\psi(\lambda)-\psi'(\lambda))d\lambda|
\]
\[
\leq(2\pi)^{-n}\int_{|\lambda|\leq m}\exp(-\frac{1}{2}\sigma^{2}\lambda^{T}\lambda)|\hat{f}(-\lambda)(\psi(\lambda)-\psi'(\lambda)|d\lambda
\]
\begin{equation}
+(2\pi)^{-n}\int_{|\lambda|>m}\exp(-\frac{1}{2}\sigma^{2}\lambda^{T}\lambda)|\hat{f}(-\lambda)(\psi(\lambda)-\psi'(\lambda)|d\lambda).\label{eq:temp-268}
\end{equation}
In view of inequalities \ref{eq:temp-56} and \ref{eq:temp-26}, the
first summand in the last sum is bounded by
\[
(2\pi)^{-n}v_{n}2^{m}\delta_{ch,dstr}(\varepsilon)\int\exp(-\frac{1}{2}\sigma^{2}\lambda^{T}\lambda)d\lambda
\]
\[
\leq(2\pi)^{-n}v_{n}2^{m}\delta_{ch,dstr}(\varepsilon)(2\pi)^{\frac{n}{2}}\sigma^{-n}=\frac{1}{8}\varepsilon,
\]
where the last equality is from the defining equality \ref{eq:temp-57}.
The second summand is bounded by
\[
(2\pi)^{-n}v_{n}2\int_{|\lambda|>m}\exp(-\frac{1}{2}\sigma^{2}\lambda^{T}\lambda)d\lambda
\]
\[
\leq(2\pi)^{-n}v_{n}2\int\cdots\int_{|\lambda_{1}|\vee\cdots\vee|\lambda_{n}|>m/\sqrt{n}}\exp(-\frac{1}{2}\sigma^{2}(\lambda_{1}^{2}+\cdots+\lambda_{n}^{2}))d\lambda_{1}\cdots d\lambda_{n}
\]
\[
\leq(2\pi)^{-n}v_{n}2(2\pi)^{\frac{n}{2}}\sigma^{-n}\int\cdots\int_{|\lambda_{1}|\vee\cdots\vee|\lambda_{n}|>\sigma m/\sqrt{n}}\varphi_{0,1}(\lambda_{1})\cdots\varphi_{0,1}(\lambda_{n})d\lambda_{1}\cdots d\lambda_{n}
\]
\[
\leq v_{n}2(2\pi)^{-\frac{n}{2}}\sigma^{-n}\sum_{j=1}^{n}\int\cdots\int_{|\lambda_{j}|>\sigma m/\sqrt{n}}\varphi_{0,1}(\lambda_{1})\cdots\varphi_{0,1}(\lambda_{n})d\lambda_{1}\cdots d\lambda_{n}
\]
\[
=v_{n}2(2\pi)^{-\frac{n}{2}}\sigma^{-n}\sum_{j=1}^{n}\int_{|\lambda_{j}|>\sigma m/\sqrt{n}}\varphi_{0,1}(\lambda_{j})d\lambda_{j}
\]
\[
=v_{n}2^{2}(2\pi)^{-\frac{n}{2}}\sigma^{-n}n\Psi(\sigma n^{-\frac{1}{2}}m)<\frac{1}{8}\varepsilon,
\]
where the last inequality follows from inequality \ref{eq:temp-27}.
Hence inequality \ref{eq:temp-268} yields
\[
|Jf_{\sigma}-J'f_{\sigma}|\leq\frac{1}{8}\varepsilon+\frac{1}{8}\varepsilon=\frac{1}{4}\varepsilon.
\]
Combining with inequality \ref{eq:temp-31} for $J$ and a similar
inequality for $J'$, we obtain
\[
|Jf-J'f|\leq|Jf-Jf_{\sigma}|+|Jf_{\sigma}-J'f_{\sigma}|+|J'f-J'f_{\sigma}|
\]
\[
\leq\frac{1}{8}\varepsilon+\frac{1}{4}\varepsilon+\frac{1}{8}\varepsilon=\frac{\varepsilon}{2},
\]
where $f\equiv g_{k,x}$, where $k=1,\cdots,p$ and $x\in A_{p}$
are arbitrary. Hence\emph{
\[
\rho_{Dist,\xi^{n}}(J,J')\equiv\sum_{k=1}^{\infty}2^{-k}|A_{k}|^{-1}\sum_{x\in A(k)}|Jg_{k,x}-J'g_{k,x}|
\]
\[
\leq\sum_{k=1}^{p}2^{-k}\frac{\varepsilon}{2}+\sum_{k=p+1}^{\infty}2^{-k}
\]
\[
\leq\frac{\varepsilon}{2}+2^{-p}<\frac{\varepsilon}{2}+\frac{\varepsilon}{2}=\varepsilon.
\]
}Assertion 1 has been proved.

2. Conversely\emph{, }let $\varepsilon>0$ be arbitrary. Write $p\equiv[0\vee(2-\log_{2}\varepsilon)]_{1}$.
For each $\theta>0$ define $\delta_{p}(\theta)\equiv p^{-1}\theta$.
By Proposition \ref{Prop. Modulus of continuity of J->Jf for fixed Lipshitz f},
there exists $\widetilde{\Delta}(\frac{\varepsilon}{4},\delta_{p},\beta,\left\Vert \xi_{R^{n}}\right\Vert )>0$
such that if 
\[
\rho_{Dist,\xi^{n}}(J,J')<\widetilde{\Delta}(\frac{\varepsilon}{4},\delta_{p},\beta,\left\Vert \xi^{n}\right\Vert )
\]
then, for each $f\in C_{ub}(R^{n})$ with modulus of continuity $\delta_{p}$
and with $|f|\leq1$, we have 
\begin{equation}
|Jf-J'f|<\frac{\varepsilon}{4}.\label{eq:temp-44}
\end{equation}
Define  
\begin{equation}
\delta_{dstr,ch}(\varepsilon,\beta,\left\Vert \xi^{n}\right\Vert )\equiv\widetilde{\Delta}(\frac{\varepsilon}{4},\delta_{p},\beta,\left\Vert \xi^{n}\right\Vert ).\label{eq:temp-485}
\end{equation}
We will prove that $\delta_{dstr,ch}(\varepsilon,\beta,\left\Vert \xi^{n}\right\Vert )$
has the desired properties. To that end, suppose 
\[
\rho_{Dist,\xi^{n}}(J,J')<\delta_{dstr,ch}(\varepsilon,\beta,\left\Vert \xi^{n}\right\Vert )\equiv\widetilde{\Delta}(\frac{\varepsilon}{4},\delta_{p},\beta,\left\Vert \xi^{n}\right\Vert ).
\]
Let $\lambda\in R^{n}$ be arbitrary with $|\lambda|\leq p$. Define
the function 
\[
h_{\lambda}(x)\equiv\exp(i\lambda^{T}x)\equiv\cos\lambda^{T}x+i\sin\lambda^{T}x
\]
for each $x\in R^{n}$. Then, using inequality \ref{eq:temp-478},
we obtain
\[
|\cos\lambda^{T}x-\cos\lambda^{T}y|\leq|\exp(i\lambda^{T}x)-\exp(i\lambda^{T}y)|
\]
\[
=|\exp(i\lambda^{T}(x-y))-1|\leq|\lambda^{T}(x-y)|\leq p|x-y|,
\]
for each $x,y\in R^{n}$. Hence the function $\cos(\lambda^{T}\cdot)$
on $R^{n}$ has modulus of continuity $\delta_{p}$. Moreover, $|\cos(\lambda^{T}\cdot)|\leq1$.
Hence, inequality \ref{eq:temp-44} is applicable and yields
\[
|J\cos(\lambda^{T}\cdot)-J'\cos(\lambda^{T}\cdot)|<\frac{\varepsilon}{4}.
\]
Similarly
\[
|J\sin(\lambda^{T}\cdot)-J'\sin(\lambda^{T}\cdot)|<\frac{\varepsilon}{4}.
\]
Combining,
\[
|\psi(\lambda)-\psi'(\lambda)|=|Jh_{\lambda}-J'h_{\lambda}|
\]
\begin{equation}
\leq|J\cos(\lambda^{T}\cdot)-J'\cos(\lambda^{T}\cdot)|+|J\sin(\lambda^{T}\cdot)-J'\sin(\lambda^{T}\cdot)|<\frac{\varepsilon}{2},\label{eq:temp-28}
\end{equation}
where $\lambda\in R^{n}$ is arbitrary with $|\lambda|\leq p$. We
conclude that
\[
\rho_{char}(\psi,\psi')\equiv\sum_{j=1}^{\infty}2^{-j}\sup_{\left|\lambda\right|\leq j}|\psi(\lambda)-\psi'(\lambda)|
\]
\emph{
\[
\leq\sum_{j=1}^{p}2^{-j}\sup_{|\lambda|\leq p}|\psi(\lambda)-\psi'(\lambda)|+\sum_{j=p+1}^{\infty}2^{-j}2
\]
\[
\leq\frac{\varepsilon}{2}+2^{-p+1}\leq\frac{\varepsilon}{2}+\frac{\varepsilon}{2}=\varepsilon.
\]
}

3. Finally, suppose $\rho_{char}(\psi_{m},\psi_{0})\rightarrow0$.
Then, by Assertion 1, we have $\rho_{Dist,\xi^{n}}(J_{m},J_{0})\rightarrow0$
as $m\rightarrow\infty$. Hence Proposition \ref{Prop. rho_xi convergence=00003D Weak Convergence}
implies that $J_{m}\Rightarrow J_{0}$.

The theorem is proved.
\end{proof}
The following propositions relate the moments of a r.r.v. $X$ to
the derivatives of its characteristic function. 
\begin{prop}
\label{Prop Taylor of ch fct} \textbf{\emph{(Taylor expansion of
characteristic functions).}} Let $n\geq1$ be arbitrary, and let $X$
be an arbitrary r.r.v. Suppose $X^{n}$ is integrable, with a simple
modulus of integrability $\eta_{intg}$ in the sense of Definition
\ref{Def. Modulus of integrability}, and with $E|X|^{n}\leq b$ for
some $b>0$. Let $\psi$ denote the characteristic function of $X$.
Define the remainder $r_{n}(\lambda)$ by
\[
\psi(\lambda)\equiv\sum_{k=0}^{n}(i\lambda)^{k}EX^{k}/k!+r_{n}(\lambda)
\]
for each \emph{$\lambda\in R$.} Then the following holds.
\end{prop}
\begin{enumerate}
\item \emph{The characteristic function $\psi$ has continuous derivative
of order $n$ on $R$, with
\begin{equation}
\psi^{(k)}(\lambda)=i^{k}EX^{k}e^{i\lambda X}\label{eq:temp-464}
\end{equation}
for each $\lambda\in R$, for each $k=0,\cdots,n$. In particular
the $k$-th moment of $X$ is given by $EX^{k}=(-i)^{k}\psi^{(k)}(0)$,
for each $k=0,\cdots,n$. Moreover, $\psi^{(n)}$ is }uniform\emph{ly
continuous on $R$, with a modulus of continuity $\delta_{\psi,n}$
on $R$ defined by
\[
\delta_{\psi,n}(\varepsilon)\equiv\frac{\varepsilon}{2b}(\eta_{intg}(\frac{\varepsilon}{2}))^{-\frac{1}{n}}
\]
for each $\varepsilon>0$.}
\item \emph{For each $\lambda\in R$ with
\[
|\lambda|<\frac{n!\varepsilon}{2b}(\eta_{intg}(\frac{n!\varepsilon}{2}))^{-\frac{1}{n}},
\]
we have $|r_{n}(\lambda)|<\varepsilon|\lambda|^{n}$.}
\item \emph{Suppose $X^{n+1}$ is integrable. Then, for each $t\in R$,
we have 
\[
\psi(t)\equiv\sum_{k=0}^{n}\psi^{(k)}(t_{0})(t-t_{0})^{k}/k!+\bar{r}_{n}(t)
\]
where
\[
|\bar{r}_{n}(t)|\leq|t-t_{0}|^{n+1}E|X|^{n+1}/(n+1)!.
\]
}
\end{enumerate}
\begin{proof}
We first observe that $(e^{iax}-1)/a\rightarrow ix$ uniformly for
$x$ in any compact interval $[-t,t]$ as $a\rightarrow0$. This can
be shown by first noting that, for arbitrary $\varepsilon>0$, Taylor's
Theorem in the Appendix implies $|e^{iax}-1-iax|\leq a^{2}x^{2}/2$
and so $|a^{-1}(e^{iax}-1)-ix|\leq ax^{2}/2<\varepsilon$ for each
$x\in[-t,t]$, provided that $a<2\varepsilon/t^{2}$.

1. Let $\lambda\in R$ be arbitrary. Proceed inductively. The assertion
is trivial if $n=0$. Suppose the assertion has been proved for $k=0,\cdots n-1$.
Let $\varepsilon>0$ be arbitrary, and let $t$ be so large that $P(|X|>t)<\varepsilon$.
For $a>0$, define $D_{a}\equiv i^{k}X^{k}e^{i\lambda X}(e^{iaX}-1)/a$.
By the observation at the beginning of this proof, $D_{a}$ converges
uniformly to $i^{k+1}X^{k+1}e^{i\lambda X}$ on $(|X|\leq t)$ as
$a\rightarrow0$. Thus we see that $D_{a}$ converges a.u. to $i^{k+1}X^{k+1}e^{i\lambda X}$.
At the same time $|D_{a}|\leq|X|^{k}|(e^{iaX}-1)/a|\leq|X|^{k+1}$
where $|X|^{k+1}$ is integrable. The Dominated Convergence Theorem
 applies, yielding $\lim_{a\rightarrow0}ED_{a}=i^{k+1}EX^{k+1}e^{i\lambda X}$.
On the other hand, by the induction hypothesis $ED_{a}\equiv a^{-1}(i^{k}EX^{k}e^{i(\lambda+a)X}-i^{k}EX^{k}e^{i\lambda X})=a^{-1}(\psi^{(k)}(\lambda+a)-\psi^{(k)}(\lambda))$.
Combining, we see that $\frac{d}{d\lambda}\psi^{(k)}(\lambda)$ exists
and is equal to $i^{k+1}EX^{k+1}e^{i\lambda X}$. Induction is completed. 

We next prove the continuity of $\psi^{(n)}$. To that end, let $\varepsilon>0$
be arbitrary. Let $\lambda,a\in R$ be arbitrary with 
\[
|a|<\delta_{\psi,n}(\varepsilon)\equiv\frac{\varepsilon}{2b}(\eta_{intg}(\frac{\varepsilon}{2}))^{-\frac{1}{n}}.
\]
Then
\[
|\psi^{(n)}(\lambda+a)-\psi^{(n)}(\lambda)|=|EX^{n}e^{i(\lambda+a)X}-EX^{n}e^{i\lambda X}|
\]
\[
\leq E|X|^{n}|e^{iaX}-1|\leq E|X|^{n}(2\wedge|aX|)
\]
\[
\leq2E|X|^{n}1_{(|X|^{n}>\eta_{intg}(\frac{\varepsilon}{2}))}+E|X|^{n}(|aX|)1_{(|X|^{n}\leq\eta_{intg}(\frac{\varepsilon}{2}))}
\]
\[
\leq\frac{\varepsilon}{2}+|a|(\eta_{intg}(\frac{\varepsilon}{2}))^{\frac{1}{n}}E|X|^{n}<\frac{\varepsilon}{2}+\frac{\varepsilon}{2b}E|X|^{n}\leq\frac{\varepsilon}{2}+\frac{\varepsilon}{2}=\varepsilon.
\]
Thus $\delta_{\psi,n}$ is the modulus of continuity of $\psi^{(n)}$
on $R$. Assertion 1 is verified.

2. Assertion 2 is an immediate consequence of Assertion 1 above and
Assertion 2 of Theorem \ref{Thm. Tayor's Theorem}\emph{ }when we
set $f\equiv\psi$, $t=\lambda$ and $t_{0}=0$ in the latter.

3. Suppose $X^{n+1}$ is integrable\emph{. }Then $\psi^{(n+1)}$ exists
on $R$, with $|\psi^{(n+1)}|\leq E|X|^{n+1}$ by equality \ref{eq:temp-464}.
Hence $\bar{r}_{n}(t)\leq E|X|^{n+1}|t-t_{0}|^{n+1}/(n+1)!$ according
to Assertion 3 of Theorem \ref{Thm. Tayor's Theorem}.
\end{proof}
For the proof of a partial converse, we need some basic equalities
for binomial coefficients. 
\begin{lem}
\label{subsec:binomial coefs} \textbf{\emph{(Binomial coefficients).
}}For each  $n\geq1$ the sum 
\[
\sum_{k=0}^{n}(_{k}^{n})(-1)^{k}k^{j}=0
\]
for $j=0,\cdots,n-1$, and 
\[
\sum_{k=0}^{n}(_{k}^{n})(-1)^{k}k^{n}=(-1)^{n}n!.
\]
\end{lem}
\begin{proof}
Differentiate $j$ times the binomial expansion
\[
(1-e^{t})^{n}=\sum_{k=0}^{n}(_{k}^{n})(-1)^{k}e^{kt}
\]
to get
\[
n(n-1)\cdots(n-j+1)(-1)^{j}(1-e^{t})^{n-j}=\sum_{k=0}^{n}(_{k}^{n})(-1)^{k}k^{j}e^{kt},
\]
and then set $t$ to $0$.
\end{proof}
Classical proofs for the next theorem in familiar texts rely on Fatou's
Lemma, which is not constructive because it trivially implies the
principle of infinite search. The following proof contains an easy
fix.
\begin{prop}
\label{Prop If Char fct has 2n derivative, then EX**2n exists} \textbf{\emph{(Moments
of r.r.v. and derivatives of its characteristic function).}} Let $\psi$
denote the characteristic function of $X$. Let $n\geq1$ be arbitrary.
If $\psi$ has a continuous derivative of order $2n$ in some neighborhood
of $\lambda=0$, then $X^{2n}$ is integrable.
\end{prop}
\begin{proof}
Write $\lambda_{k}\equiv2^{-k}$ for each $k\geq1$. Then
\[
\frac{\sin^{2}(\lambda_{k}X)}{\lambda_{k}^{2}}=\frac{\sin^{2}(2\lambda_{k+1}X)}{(2\lambda_{k+1})^{2}}=\frac{(2\sin(\lambda_{k+1}X)\cos(\lambda_{k+1}X))^{2}}{(2\lambda_{k+1})^{2}}
\]
\[
=\frac{\sin^{2}(\lambda_{k+1}X)}{\lambda_{k+1}^{2}}\cos^{2}(\lambda_{k+1}X)\leq\frac{\sin^{2}(\lambda_{k+1}X)}{\lambda_{k+1}^{2}}
\]
for each $k\geq1$. Thus we see that the sequence $((\frac{\sin^{2}(\lambda_{k}X)}{\lambda_{k}^{2}})^{n})_{k=1,2,\cdots}$
of integrable r.r.v.'s is nondecreasing. Since $\psi^{(2n)}$ exists,
we have, by Taylor's Theorem, Theorem \ref{Thm. Tayor's Theorem}
in the Appendix, 
\[
\psi(\lambda)=\sum_{j=0}^{2n}\frac{\psi^{(j)}(0)}{j!}\lambda^{j}+o(\lambda^{2n})
\]
as $\lambda\rightarrow0$. Hence for any $\lambda\in R$ we have
\[
E(\frac{\sin\lambda X}{\lambda})^{2n}=E(\frac{e^{i\lambda X}-e^{-i\lambda X}}{2i\lambda})^{2n}=(2i\lambda)^{-2n}E\sum_{k=0}^{2n}(_{k}^{2n})(-1)^{k}e^{i(2n-2k)\lambda X}
\]
\[
=(2i\lambda)^{-2n}\sum_{k=0}^{2n}(_{k}^{2n})(-1)^{k}\psi((2n-2k)\lambda)
\]
\[
=(2i\lambda)^{-2n}\sum_{k=0}^{2n}(_{k}^{2n})(-1)^{k}\{\sum_{j=0}^{2n}\frac{\psi^{(j)}(0)}{j!}(2n-2k)^{j}\lambda{}^{j}+o(\lambda^{2n})\}
\]
\[
=(2i\lambda)^{-2n}\{o(\lambda^{2n})+\sum_{j=0}^{2n}\frac{\psi^{(j)}(0)\lambda^{j}}{j!}\sum_{k=0}^{2n}(_{k}^{2n})(-1)^{k}(2n-2k)^{j}\}
\]
\[
=o(1)+(2i\lambda)^{-2n}\{\psi^{(2n)}(0)\lambda^{2n}2^{2n}\}=(-1)^{n}\psi^{(2n)}(0)
\]
in view of Lemma \ref{subsec:binomial coefs}. Consequently $E(\frac{\sin\lambda_{k}X}{\lambda_{k}})^{2n}\rightarrow(-1)^{n}\psi^{(2n)}(0)$.
At the same time $\frac{\sin\lambda_{k}t}{\lambda_{k}}\rightarrow t$
uniformly for $t$ in any compact interval. Hence $(\frac{\sin\lambda_{k}X}{\lambda_{k}})^{2n}\uparrow X^{2n}$
a.u. as $k\rightarrow\infty$. Therefore, by the Monotone Convergence
Theorem, the limit r.r.v. $X^{2n}=\lim_{k\rightarrow\infty}(\frac{\sin\lambda_{k}X}{\lambda_{k}})^{2n}$
is integrable. 
\end{proof}
\begin{prop}
\label{Prop. Product distribution and direct product of characteristic functions}
\textbf{\emph{(Product distribution and direct product of characteristic
function).}} Let $F_{1},F_{2}$ be distributions on $R^{n}$ and $R^{m}$
respectively, with the characteristic functions $\psi_{1},\psi_{2}$
respectively. Let the function $\psi_{1}\otimes\psi_{2}$ be defined
by 
\[
(\psi_{1}\otimes\psi_{2})(\lambda)\equiv\psi_{1}(\lambda_{1})\psi_{2}(\lambda_{2})
\]
for each $\lambda\equiv(\lambda_{1},\lambda_{2})\in R^{n+m}$ be arbitrary,
where$\lambda_{1}\in R^{n}$ and $\lambda_{2}\in R^{m}$. Let $F$
be a distribution on $R^{n+m}$ with characteristic function $\psi$.
Then $F=F_{1}\otimes F_{2}$ iff $\psi=\psi_{1}\otimes\psi_{2}$.
\end{prop}
\begin{proof}
Suppose $F=F_{1}\otimes F_{2}$. Let $\lambda\equiv(\lambda_{1},\lambda_{2})\in R^{n+m}$
be arbitrary, where$\lambda_{1}\in R^{n}$ and $\lambda_{2}\in R^{m}$.
Let $\exp(i\lambda^{T}\cdot)$ be the function on $R^{n+m}$ whose
value at arbitrary $x\equiv(x_{1},x_{2})\in R^{n+m}$, where $x_{1}\in R^{n}$
and $x_{2}\in R^{m}$, is $\exp(i\lambda^{T}x)$. Similarly let $\exp(i\lambda_{1}^{T}\cdot),\exp(i\lambda_{2}^{T}\cdot)$
be the functions whose values at $(x_{1},x_{2})\in R^{n+m}$ are $\exp(i\lambda_{1}^{T}x_{1}),\exp(i\lambda_{2}^{T}x_{2})$
respectively. Then
\[
\psi(\lambda)\equiv F\exp(i\lambda^{T}\cdot)=F\exp(i\lambda_{1}^{T}\cdot)\exp(i\lambda_{2}^{T}\cdot)
\]
\[
=(F_{1}\otimes F_{2})\exp(i\lambda_{1}^{T}\cdot)\exp(i\lambda_{2}^{T}\cdot)
\]
\[
=(F_{1}\exp(i\lambda_{1}^{T}\cdot))(F_{2}\exp(i\lambda_{2}^{T}\cdot))
\]
\[
=\psi_{1}(\lambda_{1})\psi_{2}(\lambda_{2})=(\psi_{1}\otimes\psi_{2})(\lambda).
\]
Thus $\psi=\psi_{1}\otimes\psi_{2}$.

Conversely, suppose $\psi=\psi_{1}\otimes\psi_{2}$. Let $G\equiv F_{1}\otimes F_{2}$.
Then $G$ has characteristic function $\psi_{1}\otimes\psi_{2}$ by
the previous paragraph. Thus the distributions $F$ and $G$ have
the same characteristic function $\psi$. By Theorem \ref{Thm. Inversion  Formula  for Ch fcts},
it follows that $F=G\equiv F_{1}\otimes F_{2}$.
\end{proof}
\begin{cor}
\label{Cor. Independence in terms of characteristic functions} \textbf{\emph{(Independence
in }}term\textbf{\emph{s of characteristic functions).}} Let $X_{1}:\Omega\rightarrow R^{n}$
and $X_{2}:\Omega\rightarrow R^{m}$ be r.v.'s on a probability space
$(\Omega,L,E)$, with characteristic functions $\psi_{1},\psi_{2}$
respectively. Let $\psi$ be the characteristic function of the r.v.
$X\equiv(X_{1},X_{2}):\Omega\rightarrow R^{n+m}$. Then \emph{$X_{1},X_{2}$
}are independent iff $\psi=\psi_{1}\otimes\psi_{2}$.
\end{cor}
\begin{proof}
Let $F,F_{1},F_{2}$ be the distributions induced by $X,X_{1},X_{2}$
on $R^{n+m},R^{n},R^{m}$ respectively. Then \emph{$X_{1},X_{2}$
}are independent iff $F=F_{1}\otimes F_{2}$, by Definition \ref{subsec:Def  Independent sets of r.v.-1}.
Since $F=F_{1}\otimes F_{2}$ iff $\psi=\psi_{1}\otimes\psi_{2}$,
according to Proposition \ref{Prop. Product distribution and direct product of characteristic functions},
the corollary is proved.
\end{proof}
\begin{prop}
\textup{\label{Prop. Conditonal expectation of jointly normal r.r.v.'s}
}\textbf{\textup{(Conditional expectation of jointly normal r.r.v.'s).}}\textup{
Let $Z_{1},\cdots,Z_{n},$ $Y_{1},\cdots,Y_{m}$} be arbitrary jointly
normal r.r.v.'s with mean $0$. Suppose the covariance matrix $\overline{\sigma}_{Z}\equiv EZZ^{T}$
of $Z\equiv(Z_{1},\cdots,Z_{n})$ is positive definite. Let $\overline{\sigma}_{Y}\equiv EYY^{T}$
be the covariance matrix of $Y\equiv(Y_{1},\cdots,Y_{m})$. Define
the $n\times m$ cross-covariance matrix $c_{Z,Y}\equiv EZY^{T}$,
and define the $n\times m$ matrix $b_{Y}\equiv\overline{\sigma}_{Z}^{-1}c_{Z,Y}$.
Then the following holds.
\end{prop}
\begin{enumerate}
\item \emph{The $m\times m$ matrix $\sigma_{Y|Z}\equiv\overline{\sigma}_{Y}-c_{Z,Y}^{T}\overline{\sigma}_{Z}^{-1}c_{Z,Y}$
is nonnegative definite. }
\item \emph{For each $f\in L_{Y}$, we have
\[
E(f(Y)|Z)=\Phi_{b_{Y}^{T}Z,\sigma_{Y|Z}}f.
\]
Heuristically, given $Z$, the conditional distribution of $Y$ is
normal with mean $b_{Y}^{T}Z$ and covariance matrix $\sigma_{Y|Z}$.
In particular, $E(Y|Z)=b_{Y}^{T}Z=c_{Z,Y}^{T}\overline{\sigma}_{Z}^{-1}Z$. }
\item \emph{The r.v.'s $V\equiv E(Y|Z)$ and $X\equiv Y-E(Y|Z)$ are independent
normal r.v.'s with values in $R^{m}$.}
\item \emph{$EY^{T}Y=EV^{T}V+EX^{T}X.$}
\end{enumerate}
\begin{proof}
1. Let $X\equiv(X_{1},\cdots,X_{m})\equiv Y-b_{Y}^{T}Z.$ Thus $Y=b_{Y}^{T}Z+X.$
Then $Z_{1},\cdots,Z_{n},X_{1},\cdots,X_{m}$ are jointly normal according
to Proposition \ref{Prop. Basics of Normal distributions}. Furthermore,
\[
EZX^{T}=EZY^{T}-EZZ^{T}b_{Y}\equiv c_{Z,Y}-\overline{\sigma}_{Z}b_{Y}=0,
\]
while the covariance matrix of $X$ is given by
\[
\overline{\sigma}_{X}\equiv EXX^{T}=EYY^{T}-EYZ^{T}b_{Y}-Eb_{Y}^{T}ZY^{T}+Eb_{Y}^{T}ZZ^{T}b_{Y}
\]
\[
=\overline{\sigma}_{Y}-c_{Z,Y}^{T}b_{Y}-b_{Y}^{T}c_{Z,Y}+b_{Y}^{T}\overline{\sigma}_{Z}b_{Y}
\]
\[
=\overline{\sigma}_{Y}-c_{Z,Y}^{T}\overline{\sigma}_{Z}^{-1}c_{Z,Y}-b_{Y}^{T}c_{Z,Y}+b_{Y}^{T}c_{Z,Y}
\]
\[
=\overline{\sigma}_{Y}-c_{Z,Y}^{T}\overline{\sigma}_{Z}^{-1}c_{Z,Y}\equiv\sigma_{Y|Z},
\]
whence $\sigma_{Y|Z}$ is nonnegative definite.

2. Hence the r.v. $U\equiv(Z,X)$ in $R^{n+m}$ has mean $0$ and
covariance matrix
\[
\overline{\sigma}_{U}\equiv\left[\begin{array}{cc}
\overline{\sigma}_{Z} & 0\\
0 & \overline{\sigma}_{X}
\end{array}\right]\equiv\left[\begin{array}{cc}
\overline{\sigma}_{Z} & 0\\
0 & \sigma_{Y|Z}
\end{array}\right].
\]
Accordingly, $U$ has the characteristic function
\[
E\exp(i\lambda^{T}U)=\psi_{0,\overline{\sigma}_{U}}(\lambda)\equiv\exp(-\frac{1}{2}\lambda^{T}\overline{\sigma}_{U}\lambda)
\]
\[
=\exp(-\frac{1}{2}\theta^{T}\overline{\sigma}_{Z}\theta)\exp(-\frac{1}{2}\gamma^{T}\overline{\sigma}_{X}\gamma)
\]
\[
=E\exp(i\theta^{T}Z)E\exp(i\gamma^{T}X),
\]
for each $\lambda\equiv(\theta_{1},\cdots,\theta_{n},\gamma_{1},\cdots,\gamma_{m})\in R^{n+m}$.
It follows from Corollary \ref{Cor. Independence in terms of characteristic functions}
that $Z,X$ are independent. In other words, the distribution $E_{(Z,X)}$
induced by $(Z,X)$ on $R^{n+m}$ is given by the product distribution,
\[
E_{(Z,X)}=E_{Z}\otimes E_{X}
\]
of $E_{Z},E_{X}$ induced on $R^{n},R^{m}$ respectively by $Z,X$
respectively. 

Now let $f\in L_{Y}$ be arbitrary. Thus $f(Y)\in L$. Let $z\equiv(z_{1},\cdots,z_{n})\in R^{n}$
and $z\equiv(x_{1},\cdots,x_{m})\in R^{m}$ be arbitrary. Define
\[
\tilde{f}(z,x)\equiv f(b_{Y}^{T}z+x)
\]
and
\[
\bar{f}(z)\equiv E_{X}\tilde{f}(z,\cdot)\equiv E\tilde{f}(z,X)\equiv Ef(b_{Y}^{T}z+X)=\Phi_{b_{Y}^{T}z,\sigma_{Y|Z}}f.
\]
We will prove that the r.r.v. $\bar{f}(Z)$ is the condition expectation
of $f(Y)$ given $Z$. To that end, let $g\in C(R^{n})$ be arbitrary.
Then, by Fubini's Theorem
\[
Ef(Y)g(Z)=E\tilde{f}(Z,X)g(Z)=E_{(Z,X)}(\tilde{f}g)=E_{Z}\otimes E_{X}(\tilde{f}g)
\]
\[
=E_{Z}(E_{X}(\tilde{f}g))=E_{Z}\bar{f}g=E\bar{f}(Z)g(Z).
\]
It follows that $E(f(Y)|Z)=\bar{f}(Z)\equiv\Phi_{b_{Y}^{T}Z,\sigma_{Y|Z}}f$.
In particular, $E(Y|Z)=b_{Y}^{T}Z=c_{Z,Y}^{T}\overline{\sigma}_{Z}^{-1}Z$. 

3. By Step 2, the r.v.'s $Z,X$ are independent normal. Hence the
r.v.'s $V\equiv E(Y|Z)=b_{Y}^{T}Z$ and $X\equiv Y-E(Y|Z)$ are independent
normal.

4. Hence $EV^{T}X=(EV^{T})(EX)=0$. It follows that 
\[
EY^{T}Y=E(V+X)^{T}(V+X)=EV^{T}V+EX^{T}X.
\]
\end{proof}

\section{The Central Limit Theorem}

Let $X_{1},\cdots,X_{n}$ be independent r.r.v.'s with mean $0$ and
standard deviations $\sigma_{1},\cdots,\sigma_{n}$ respectively.
Define $\sigma$ by $\sigma^{2}=\sigma_{1}^{2}+\cdots+\sigma_{n}^{2}$
and consider the distribution $F$ of the scaled sum $X=(X_{1}+\cdots+X_{n})/\sigma$.
By replacing $X_{i}$ with $X_{i}/\sigma$ we may assume that $\sigma=1$.
The Central Limit Theorem says that, if each individual summand $X_{i}$
is small relative to the sum $X$, then $F$ is close to the standard
normal distribution $\Phi_{0,1}$. 

One criterion, due to Lindberg and Feller, for the summands $X_{k}$
$(k=1,\cdots,n$) to be individually small relative to the sum, is
for
\[
\theta(r)\equiv\sum_{k=1}^{n}(E1_{|X_{k}|>r}X_{k}^{2}+E1_{|X_{k}|\leq r}|X_{k}|^{3})
\]
to be small for some $r\geq0$. 
\begin{lem}
\label{Lem. bound for sum of sigmak**3} Suppose $r\geq0$ is such
that $\theta(r)<\frac{1}{8}$. Then 
\begin{equation}
\sum_{k=1}^{n}\sigma_{k}^{3}\leq\theta(r)\label{eq:temp-60}
\end{equation}
\end{lem}
\begin{proof}
Consider each $k=1,\cdots,n$. Then, since $\theta(r)<\frac{1}{8}$
by hypothesis, we have $z\equiv E1_{|X_{k}|>r}X_{k}^{2}<\frac{1}{8}$
and $a\equiv E1_{|X_{k}|\leq r}|X_{k}|^{3}<\frac{1}{8}$. A consequence
is that $(z+a^{2/3})^{3/2}\leq z+a$ which can be seen by noting that
the two sides are equal at $z=0$ and by comparing first derivatives
relative to $z$ on $[0,\frac{1}{8}]$. Lyapunov's inequality then
implies that 
\[
\sigma_{k}^{3}=(EX_{k}^{2}1_{(|X_{k}|>r)}+EX_{k}^{2}1_{(|X_{k}|\leq r)})^{3/2}
\]
\[
\leq(EX_{k}^{2}1_{(|X_{k}|>r)}+(E|X_{k}|^{3}1_{(|X_{k}|\leq r)})^{2/3})^{3/2}
\]
\[
\equiv(z+a^{2/3})^{3/2}\leq z+a\equiv EX_{k}^{2}1_{(|X_{k}|>r)}+E|X_{k}|^{3}1_{(|X_{k}|\leq r)}.
\]
Summing over $k$, we obtain inequality \ref{eq:temp-60}. 
\end{proof}
\begin{thm}
\label{Thm.Central-Limit Theorem}\index{Central Limit Theorem} \textbf{\emph{(Central
Limit Theorem)}}. Let $f\in C(R)$ and $\varepsilon>0$ be arbitrary.
Then there exists $\delta>0$ such that, if $\theta(r)<\delta$ for
some $r\geq0$, then 
\begin{equation}
|\int f(x)dF(x)-\int f(x)d\Phi_{0,1}(x)|<\varepsilon.\label{eq:temp-117}
\end{equation}
\end{thm}
\begin{proof}
Let $\xi_{R}$ be an arbitrary, but fixed, binary approximation of
$R$ relative to the reference point $0$. We assume, without loss
of generality, that $|f(x)|\leq1$. Let $\delta_{f}$ be a modulus
of continuity of $f$, and let $b>0$ be so large that $f$ has $[-b,b]$
as support. Let $\varepsilon>0$ be arbitrary. By Proposition \ref{Prop. rho_xi convergence=00003D Weak Convergence},
there exists $\delta_{\widehat{J}}(\varepsilon,\delta_{f},b,\left\Vert \xi_{R}\right\Vert )>0$
such that, if the distributions $F,\Phi_{0,1}$ satisfy 
\begin{equation}
\rho_{\xi(R)}(F,\Phi_{0,1})<\varepsilon'\equiv\delta_{\widehat{J}}(\varepsilon,\delta_{f},b,\left\Vert \xi_{R}\right\Vert ),\label{eq:temp-116}
\end{equation}
then inequality \ref{eq:temp-117} holds. Separately, according to
Corollary \ref{Thm. Continuity Theorem for ch functions}, there exists
$\delta_{ch,dstr}(\varepsilon')>0$ such that, if the characteristic
functions $\psi_{F},\psi_{0,1}$ of $F,\Phi_{0,1}$ respectively satisfy
\begin{equation}
\rho_{char}(\psi_{F},\psi_{0,1})\equiv\sum_{j=1}^{\infty}2^{-j}\sup_{|\lambda|\leq j}|\psi_{F}(\lambda)-\psi_{0,1}(\lambda)|<\varepsilon''\equiv\delta_{ch,dstr}(\varepsilon'),\label{eq:temp-121}
\end{equation}
then inequality \ref{eq:temp-116} holds. 

Now take $m\geq1$ be so large that $2^{-m+2}<\varepsilon''$, and
define 
\[
\delta\equiv\frac{1}{8}\wedge\frac{1}{6}m^{-3}\varepsilon''.
\]
Suppose 
\[
\theta(r)<\delta
\]
for some $r\geq0$. Then $\theta(r)<\frac{1}{8}$. We will show that
inequality \ref{eq:temp-117} holds.

To that end, let $\lambda\in[-m,m]$ and $k=1,\cdots,n$ be arbitrary.
Let $\varphi_{k}$ denote the characteristic function of $X_{k}$,
and let $Y_{k}$ be a normal r.r.v. with mean $0$, variance $\sigma_{k}^{2}$,
and characteristic function $e^{-\sigma_{k}^{2}\lambda^{2}/2}$. Then
\[
E|Y_{k}|^{3}=\frac{2}{\sqrt{2\pi}\sigma_{k}}\int_{0}^{\infty}y^{3}\exp(-\frac{1}{2\sigma_{k}^{2}}y^{2})dy.
\]
\[
=\frac{4\sigma_{k}^{3}}{\sqrt{2\pi}}\int_{0}^{\infty}u\exp(-u)du=\frac{4\sigma_{k}^{3}}{\sqrt{2\pi}}=2\sqrt{\frac{2}{\pi}}\sigma_{k}^{3},
\]
where we made a change of integration variables $u\equiv-\frac{1}{2\sigma_{k}^{2}}y^{2}$.
Moreover, since $\sum_{k=1}^{n}\sigma_{k}^{2}=1$ by assumption, and
since all characteristic functions have absolute value bounded by
1, we have
\[
|\varphi_{F}(\lambda)-e^{-\lambda^{2}/2}|=|\prod_{k=1}^{n}\varphi_{k}(\lambda)-\prod_{k=1}^{n}e^{-\sigma_{k}^{2}\lambda^{2}/2}|
\]
\begin{equation}
\leq\sum_{k=1}^{n}|\varphi_{k}(\lambda)-e^{-\sigma_{k}^{2}\lambda^{2}/2}|.\label{eq:temp-54}
\end{equation}
By Proposition \ref{Prop Taylor of ch fct}, the Taylor expansions
up to degree 2 for the characteristic functions $\varphi_{k}(\lambda)$
and $e^{-\sigma_{k}^{2}\lambda^{2}/2}$ are equal because the two
corresponding distributions have equal first and second moments. Hence
the difference of the two functions is equal to the difference of
the two remainders in their respective Taylor expansions. Again by
Proposition \ref{Prop Taylor of ch fct}, the remainder for $\varphi_{k}(\lambda)$
is bounded by
\[
\lambda^{2}EX_{k}^{2}1_{(|X_{k}|>r)}+\frac{|\lambda|^{3}}{3!}E|X_{k}|^{3}1_{(|X_{k}|\leq r)}
\]
\[
\leq m^{3}(EX_{k}^{2}1_{(|X_{k}|>r)}+E|X_{k}|^{3}1_{(|X_{k}|\leq r)}).
\]
By the same token, the remainder for $e^{-\sigma_{k}^{2}\lambda^{2}/2}$
is bounded by a similar expression, where $X_{k}$ is replaced by
$Y_{k}$ and where $r\geq0$ is replaced by $s\geq0$, which becomes,
as $s\rightarrow\infty$, 
\[
m^{3}E|Y_{k}|^{3}=2\sqrt{\frac{2}{\pi}}m^{3}\sigma_{k}^{3}<2m^{3}\sigma_{k}^{3}.
\]
Combining, inequality \ref{eq:temp-54} yields, for each $\lambda\in[-m,m]$,
\[
|\varphi_{F}(\lambda)-e^{-\lambda^{2}/2}|\leq m^{3}\sum_{k=1}^{n}(EX_{k}^{2}1_{(|X_{k}|>r)}+E|X_{k}|^{3}1_{(|X_{k}|\leq r)})+2m^{3}\sum_{k=1}^{n}\sigma_{k}^{3}
\]
\[
\leq3m^{3}\theta(r)\leq3m^{3}\delta\leq\frac{\varepsilon''}{2},
\]
where the second inequality follows from the definition of $\theta(r)$
and from Lemma \ref{Lem. bound for sum of sigmak**3}. Hence, since
$|\psi_{F}-\psi_{0,1}|\leq2$, we obtain 
\[
\rho_{char}(\psi_{F},\psi_{0,1})\leq\sum_{j=1}^{m}2^{-j}\sup_{|\lambda|\leq j}|\psi_{F}(\lambda)-\psi_{0,1}(\lambda)|+2^{-m+1}
\]
\[
\leq\frac{\varepsilon''}{2}+\frac{\varepsilon''}{2}=\varepsilon''\equiv\delta_{ch,dstr}(\varepsilon'),
\]
establishing inequality \ref{eq:temp-121}. Consequently, inequality
\ref{eq:temp-116}, and, in turn, inequality \ref{eq:temp-117} follow.
The theorem is proved.
\end{proof}
\begin{cor}
\textbf{\emph{(Lindberg's Central Limit Theorem)}} For each $p=1,2,\cdots$,
let $n_{p}\geq1$ be arbitrary, and let $(X_{p,1},\cdots,X_{p,n(p)})$
be an independent sequence of r.r.v.'s with mean $0$ and variance
$\sigma_{p,k}^{2}$ such that $\sum_{k=1}^{n(p)}\sigma_{p,k}^{2}=1$.
Suppose for each $r>0$ we have
\begin{equation}
\lim_{p\rightarrow\infty}\sum_{k=1}^{n(p)}EX_{p,k}^{2}1_{(|X(p,k)|>r)}=0\label{eq:temp-123}
\end{equation}
Then $\sum_{k=1}^{n(p)}X_{p,k}$ converges in distribution to the
standard normal distribution $\Phi_{0,1}$ as $p\rightarrow\infty$.
\end{cor}
\begin{proof}
Let $\delta>0$ be arbitrary. According to Theorem \ref{Thm.Central-Limit Theorem},
it suffices to show that there exists $r>0$ such that, for sufficiently
large $p$, we have
\[
\sum_{k=1}^{n(p)}EX_{p,k}^{2}1_{(|X(p,k)|>r)}<\frac{\delta}{2}
\]
 and
\[
\sum_{k=1}^{n(p)}E|X_{p,k}|^{3}1_{(|X(p,k)|<r)}<\frac{\delta}{2}.
\]
For that purpose, take any $r\in(0,\frac{\delta}{2})$. Then the first
of the last two inequalities holds for sufficiently large $p$, in
view of inequality \ref{eq:temp-123} in the hypothesis. The second
follows from
\[
\sum_{k=1}^{n(p)}E|X_{p,k}|^{3}1_{(|X(p,k)|<r)}\leq r\sum_{k=1}^{n(p)}EX_{p,k}^{2}=r<\frac{\delta}{2}.
\]
\end{proof}
Because of the importance of the Central Limit Theorem, much work
since the early development of probability theory has been dedicated
to an optimal rate of convergence, culminating in the Feller's bound:
$\sup_{t\in R}|F(t)-\Phi(t)|\leq6\theta(r)$ for each $r\geq0$ .
The proof on pages 544-546 of \cite{FellerII71}, which is a careful
analysis of the difference $\varphi_{k}(\lambda)-e^{-\sigma_{k}^{2}\lambda^{2}/2}$,
contains a few typos and omitted steps which serve to keep the reader
on the toes. That proof contains also a superfluous assumption that
the maximum distance between two P.D.F.'s is always attained at some
point in $R$. There is no constructive proof for the general validity
of that assumption. There is however an easy constructive substitute
which says that, if one of the two P.D.F.s is continuously differentiable,
then the supremum distance exists: there is a sequence $(t_{n})$
in $R$ such that $\lim_{n\rightarrow\infty}|F(t_{n})-\Phi_{0,1}(t_{n})|$
exists and bounds any $|F(t)-\Phi_{0,1}(t)|$. This is sufficient
for  Feller's proof. 

\section{Supplements and Exercises}
\begin{xca}
\label{Ex. E(Y|X) may not exists even assuming X continuous} \textbf{(General
existence of $E(Y|X)$ implies the principle of infinite search).} 
\end{xca}
Hint. We will give a counter example where $X$ has a p.d.f. To that
end, let $\Omega=[0,1]\times\{-1,1\}$ be equipped with the Euclidean
metric. Define the distribution $E_{0}\equiv I_{0}\otimes I_{1}$
on $\Omega$, where $I_{0}$ is the Lebesgue integration on $[0,1]$,
and where $I_{1}$ is the probability integration on $\{-1,1\}$ which
assigns equal probabilities to each of $-1$ and 1. Let $(r_{n})_{n=1,2,\cdots}$
be an arbitrary 0-1 sequence with at most a 1.

For each $n\geq1$, define a distribution $E_{n}$on $\Omega$ by
$E_{n}Z\equiv\int_{0}^{1}Z(x,(-1)^{[2^{n}x]})dx$ for each $Z\in C(\Omega)$,
where $[a]$ stands for the integer part of a real number $a$, whence
$[2^{n}x]$ is defined for a.e. $x\in[0,1]$. Then $|E_{n}Z-E_{0}Z|\rightarrow0$
for each $Z\in C(\Omega)$. Hence
\[
EZ\equiv E_{0}Z+\sum_{n=1}^{\infty}r_{n}(E_{n}Z-E_{0}Z)
\]
exists. Then $(\Omega,C(\Omega),E)$ is an integration space with
$E1=1$. Therefore its completion $(\Omega,L,E)$ is a probability
space. We have also 
\[
|E_{m}Z-EZ|\leq|E_{m}Z-E_{0}Z|+\sum_{n=1}^{\infty}r_{n}(E_{n}Z-E_{0}Z)\rightarrow0\textrm{\quad as }m\rightarrow\infty.
\]
Let $X,Y$ be the first and second coordinate functions: $X(x,y)=x$
and $Y(x,y)=y$ for each $(x,y)\in\Omega$. It can easily be shown
that $E_{0}(Y|X)=0$ and $E_{n}(Y|X)=(-1)^{[2^{n}X]}$ for each $n\geq1$.
Moreover $X$ has a p.d.f.: $Ef(X)=\int_{0}^{1}f(t)dt$. 

Assume that $E(Y|X)$ exists. Let $a\in(0,1)$ be such that $A=(E(Y|X)>a)$
is $\mathrm{measurable}$. Either $P(A)>0$ or $P(A)<\frac{1}{2}$.

Suppose $P(A)>0$. Then
\[
EY1_{A}=E(E(Y|X)1_{A})=E(E(Y|X)1_{(E(Y|X)>a)})>aP(A)>0
\]
At the same time, $E_{0}Y1_{A}=E_{0}(E_{0}(Y|X)1_{A})=0$. Therefore
$\sum_{n=1}^{\infty}r_{n}(E_{n}Y1_{A}-E_{0}Y1_{A})=EY1_{A}>0$, and
so there exists an $n$-th summand in the series which is positive.
Hence that $r_{n}>0$ for some $n\geq1$. 

Suppose, on the other hand, that $P(A)<\frac{1}{2}$. Suppose $r_{n}>0$
for some $n$. Then $E=E_{n}$. Hence $A=(E_{n}(Y|X)>a)$ and so $P(A)=P((-1)^{[2^{n}X]})=\frac{1}{2}$,
a contradiction. Hence $r_{n}=0$ for each $n$.

Thus from the general existence of $E(Y|X)$ we have deduced the principle
of infinite search. $\square$
\begin{xca}
\label{Ex. E(Y|X) may exist w/o  E(|Y||X) existing} \textbf{(Constructively,
existence of $E(Z|X)$ does not imply that of $E(|Z||X)$).} Let $X,Y$
be the r.r.v.'s as constructed in the preceding exercise. Let $V$
be a r.r.v. independent of $X,Y$ such that $P(V=1)=P(V=-1)=\frac{1}{2}$.
Define $Z=V(1+Y)$. Then $E(Z|X)=0$. At the same time, $|Z|=(1+Y)$.
Hence the existence of $E(|Z||X)$ would imply that of $E(Y|X)$.
Combining with Exercise \ref{Ex. E(Y|X) may not exists even assuming X continuous},
we see that a proof that existence of \textbf{$E(Z|X)$} implies that
of $E(|Z||X)$ would also prove the principle of infinite search. 
\end{xca}
$\,$
\begin{xca}
\label{Ex. Absolute Val of Integral of complex function} (Integration
of complex valued functions). Let $f$ be a complex-valued integrable
function on and integration space $(\Omega,L,I)$. Then $|If|\leq I|f|$. 
\end{xca}
Proof. Let $g$ and $h$ denote the real and imaginary parts of $f$
respectively. Thus $f\equiv g+ih$. By hypothesis both $g$ and $h$
are integrable relative to $I$. Hence, by definition, $f$ is integrable.
Moreover, 
\[
|If|^{2}\equiv|Ig+iIh|^{2}\leq|Ig|^{2}+|Ih|^{2}\leq Ig^{2}+Ih^{2}=I|f|^{2}.
\]
$\square$
\begin{xca}
\label{Ex. Weak Convergence =00003D convergence at all bounded continuouts fcts}
(Weak Convergence is equivalent to convergence at all bounded continuous
functions). Let $(S,d)$ be a locally compact metric space. Let $I,I_{n}\in Q_{S}$
for each $n\geq1$. Prove that $I_{n}\Rightarrow I$ iff $I_{n}f\rightarrow If$
for each $f\in C_{b}(S)$. 
\end{xca}
Hint. Suppose $I_{n}\Rightarrow I$. Equivalently $I_{n}g\rightarrow Ig$
for each $g\in C(S)$. Moreover $\{I,I_{1},I_{2},\cdots\}$ is tight.
Consider $f\in C_{b}(S)$. Without loss of generality, assume that
$0\leq f\leq1$. Let $\varepsilon>0$ be arbitrary. Let $x_{\circ}\in S$
be fixed and let $a>0$ be so large that $g\equiv1\wedge(a-d(\cdot,x_{\circ}))_{+}$
has $0\leq1-Ig<\varepsilon$. Then $0\leq If-Ifg=If(1-g)\leq1-Ig\leq\varepsilon$.
Since $g,fg\in C(S)$, there exists $m$ so large that $|I_{n}g-Ig|<\varepsilon$
and $|I_{n}fg-Ifg|<\varepsilon$ for each $n\geq m$. Consider any
$n\geq m$. We then have $0\leq1-I_{n}g\leq1-Ig+|I_{n}g-Ig|<2\varepsilon$.
Consequently $0\leq I_{n}f-I_{n}fg\leq I_{n}(1-g)<2\varepsilon$.
Combining, we have $|I_{n}f-If|\leq|I_{n}f-I_{n}fg|+|I_{n}fg-Ifg|+|Ifg-If|<4\varepsilon$.
Since $\varepsilon>0$ is arbitrary, we have proved that $I_{n}g\rightarrow Ig$
for each $g\in C(S)$  implies $I_{n}f\rightarrow If$ for each $f\in C_{b}(S)$.
The converse is trivial since $C(S)\subset C_{b}(S)$. $\square$

\part{Stochastic Process}

\chapter{Random Fields and Stochastic Processes}

In this chapter, unless otherwise specified, $(S,d)$ will denote
a locally compact metric space, with an arbitrary, but fixed, reference
point $x_{\circ}.$

\section{Random Field and Finite Joint Distributions}

In this section, we introduce random fields, their marginal distributions,
and some notions of their continuity. and let $Q$ be a set. Let $(\Omega,L,E)$
be a probability space. 
\begin{defn}
\label{Def. Random Fields and Stochastic Processes} \textbf{(Random
Fields).} Suppose a function 
\[
X:Q\times\Omega\rightarrow S
\]
is such that, for each $t\in Q$, the function $X_{t}\equiv X(t,\cdot)$
is a r.v. on $\Omega$ with values in $S$. Then $X$ is called a
\emph{random field}\index{random field}, or r.f. for abbreviation,
with \index{sample space}\emph{sample space} $\Omega$, with \index{parameter set}\emph{parameter
set} $Q$, and with \emph{state space}\index{state space} $S$. To
be precise, we will sometimes write
\[
X:Q\times(\Omega,L,E)\rightarrow(S,d).
\]
We will let $\widehat{R}(Q\times\Omega,S)$ denote the set of such
r.f.'s. Two r.f.'s $X,Y\in\widehat{R}(Q\times\Omega,S)$ are considered
equal if $X_{t}=Y_{t}$ a.s. on $\Omega$, for each $t\in Q$.

Let $X\in\widehat{R}(Q\times\Omega,S)$ be arbitrary. For each $\omega\in\Omega$
such that $domain(X(\cdot,\omega))$ is nonempty, the function $X(\cdot,\omega)$
is called a \index{sample function}\emph{sample function}. If $K$
is a subset of $Q$, then we write
\[
X|K\equiv X|(K\times\Omega):K\times\Omega\rightarrow S,
\]
and call the r.f. $X|K$ the \index{restriction of a r.f.}\emph{restriction}
of $X$ to $K$.

In the special case where the parameter set $Q$ is a subset of $R$,
the r.f. $X$ is called a \emph{stochastic process}\index{stochastic process},
or simply a \emph{process}. In that case, the variable $t\in Q$ is
often called the \index{time parameter} \emph{time parameter.}$\square$ 

$\,$

When the parameter set $Q$ is countably infinite, we can view a r.f.
with state space $(S,d)$ as a r.v. with values in $(S^{\infty},d^{\infty})$.
\end{defn}
\begin{lem}
\label{Lem. Random field with countable parametr set is a r.v.} \textbf{\emph{(Random
field with countable parameter set can be regarded as a r.v. with
values in the path space, and conversely).}} Let $X:Q\times(\Omega,L,E)\rightarrow(S,d)$
be an arbitrary r.f. where the parameter set $Q=\{t_{1},t_{2},\cdots\}$
is countably infinite. Then  the function $(X_{t(1)},X_{t(2)},\cdots):\Omega\rightarrow S^{Q}$
is a r.v. on $(\Omega,L,E)$ with values in the complete metric space
$(S^{Q},d^{Q})$. The converse also holds.
\end{lem}
\begin{proof}
1. Suppose $X:Q\times(\Omega,L,E)\rightarrow(S,d)$ is a r.f. Let
$x_{\circ}$ be an arbitrary, but fixed, reference point in $S$.
Let $f\in C_{ub}(S^{\infty},d^{\infty})$ be arbitrary, with a modulus
of continuity $\delta_{f}$. Let $\varepsilon>0$ be arbitrary. Let
$m\geq1$ be so large that $2^{-m}<\delta_{f}(\varepsilon)$. Define
a function $f_{m}$ on $S^{m}$ by 
\[
f_{m}(x_{1},\cdots,x_{m})\equiv f(x_{1},\cdots,x_{m},x_{\circ},x_{\circ},\cdots)
\]
for each $(x_{1},\cdots,x_{m})\in S^{m}$. Then it is easily verified
that $f_{m}\in C_{ub}(S^{m},d^{m})$. Hence $f_{m}(X_{t(1)},\cdots,X_{t(m)})\in L$
. At the same time, 
\[
d^{\infty}((X_{t(1)},X_{t(2)},\cdots),(X_{t(1)},\cdots,X_{t(m)},x_{\circ},x_{\circ},\cdots))
\]
\[
=\sum_{k=m+1}^{\infty}2^{-k}\widehat{d}(X_{t(k)},x_{\circ})\leq\sum_{k=m+1}^{\infty}2^{-k}=2^{-m}<\delta_{f}(\varepsilon).
\]
Consequently 
\[
|f(X_{t(1)},X_{t(2)},\cdots)-f_{m}(X_{t(1)},\cdots,X_{t(m)})|<\varepsilon.
\]
Thus $f(X_{t(1)},X_{t(2)},\cdots)$ is the uniform limit of a sequence
in $L$, hence is itself a member of $L$. Since the complete metric
space $(S^{\infty},d^{\infty})$ is bounded, Proposition \ref{Prop. Basic Properties of r.v.}
implies that the function $(X_{t(1)},X_{t(2)},\cdots)$ is a r.v. 

2. Conversely, suppose $X\equiv(X_{t(1)},X_{t(2)},\cdots)$ is a r.v.
on $(\Omega,L,E)$ with values in the complete metric space $(S^{Q},d^{Q})$.
Let $n\geq1$ be arbitrary. Define the function $f:(S^{Q},d^{Q})\rightarrow(S,d)$
by $f(x_{1},x_{2},\cdots,x_{n},\cdots)\equiv x_{n}$ for each $(x_{1},x_{2},\cdots,x_{n},\cdots)\in(S^{Q},d^{Q})$.
Then it can easily be verified that the function $f$ is uniformly
continuous and is bounded on bounded subsets of $(S^{Q},d^{Q})$.
Hence Proposition \ref{Prop. X meas,  f unif continuous and bd on bd subsets =00003D> f(X) meas}
implies that $f\circ X\equiv X_{t(n)}$ is a r.v. on $(\Omega,L,E)$
with values in $(S,d)$. Since $n\geq1$ is arbitrary, we conclude
that $X:Q\times(\Omega,L,E)\rightarrow(S,d)$ is a r.f. 
\end{proof}
In general, when the parameter set $Q$ is a metric space, we introduce
three notions of continuity of a r.f. They correspond to the terminology
in \cite{Neveu 65}. For ease of presentation, we restrict our attention
to the special case where $Q$ is bounded. The generalization to a
locally compact metric space is straightforward.
\begin{defn}
\label{Def. continuity in prob, continuity a.u., and a.u. continuity}
\textbf{(Continuity of r.f. on a bounded metric parameter space).}
Let $X:Q\times\Omega\rightarrow S$ be a r.f., where $(S,d)$ is a
locally compact metric space and where $(Q,d_{Q})$ is a bounded metric
space. Thus $d_{Q}\leq b$ for some $b\geq0$. 
\end{defn}
\begin{enumerate}
\item Suppose, for each $\varepsilon>0$, there exists $\delta_{Cp}(\varepsilon)>0$
such that
\[
E(1\wedge d(X_{t},X_{s}))\leq\varepsilon
\]
for each $s,t\in Q$ with $d_{Q}(t,s)<\delta_{Cp}(\varepsilon)$.
Then the r.f. $X$ is said to be \index{r.f., continuous in probability}\emph{continuous
in probability}, with the operation $\delta_{Cp}$ as a \emph{modulus
of continuity in probability}\index{modulus of continuity in probability}.
We will let $\widehat{R}_{Cp}(Q\times\Omega,S)$ denote the set of
r.f.'s which are continuous in probability, with the given bounded
metric space as parameter space. 
\item Suppose $domain(X(\cdot,\omega))$ is dense in $Q$ for a.e. $\omega\in\Omega.$
Suppose, in addition, that for each $\varepsilon>0$, there exists
$\delta_{cau}(\varepsilon)>0$ such that, for each $s\in Q$, there
exists a $\mathrm{measurable}$ set $D_{s}\subset domain(X_{s})$
with $P(D_{s}^{c})<\varepsilon$ such that for each $\omega\in D_{s}$
and for each $t\in domain(X(\cdot,\omega))$ with $d_{Q}(t,s)<\delta_{cau}(\varepsilon)$,
we have
\[
d(X(t,\omega),X(s,\omega))\leq\varepsilon.
\]
Then the r.f. $X$ is said to be \index{r.f., continuity a.u.}\emph{continuous
a.u.}, with the operation $\delta_{cau}$ as a \emph{modulus of continuity
a.u.}\index{modulus of continuity a.u. of r.f.} on $Q$. 
\item Suppose, for each $\varepsilon>0$, there exist $\delta_{auc}(\varepsilon)>0$
and a $\mathrm{measurable}$ set $D$ with $P(D^{c})<\varepsilon$
such that
\[
d(X(t,\omega),X(s,\omega))\leq\varepsilon,
\]
and $s,t\in domain(X(\cdot,\omega)$ with $d_{Q}(t,s)<\delta_{auc}(\varepsilon)$,
for each $\omega\in D$. Then the r.f. $X$ is said to be \emph{a.u.
continuous}\index{r.f., a.u. continuous}, with the operation $\delta_{auc}$
as a \emph{modulus of a.u. continuity.}\index{modulus of a.u. continuity of r.f.}.
$\square$
\end{enumerate}
The reader can give simple examples of stochastic processes which
are continuous in probability but not continuous a.u., and of processes
which are continuous a.u. but not a.u. continuous.
\begin{defn}
\label{Def. continuity in prob etc for arbitrary paameter metric space}
\textbf{(Continuity of r.f. on an arbitrary metric parameter space).}
Let $X:Q\times\Omega\rightarrow S$ be a r.f., where $(S,d)$ is a
locally compact metric space and where $(Q,d_{Q})$ is an arbitrary
metric space. 

The r.f. $X$ is said to be \emph{continuous in probability }if, for
each bounded subset $K$ of $Q$, the restricted r.f. $X|K:K\times\Omega\rightarrow S$
is continuous in probability. 

The r.f. $X$ is said to be \emph{continuous a.u. }if, for each bounded
subset $K$ of $Q$, the restricted r.f. $X|K:K\times\Omega\rightarrow S$
is continuous a.u. 

The r.f. $X$ is said to be\emph{ a.u.} \emph{continuous} if, for
each bounded subset $K$ of $Q$, the restricted r.f. $X|K:K\times\Omega\rightarrow S$
is a.u. continuous. $\square$
\end{defn}
\begin{prop}
\label{Prop. Alternative characterization of r.f. continuity}. \textbf{\emph{(Alternative
definitions of r.f. continuity). }}Let $X:Q\times\Omega\rightarrow S$
be a r.f., where $(S,d)$ is a locally compact metric space and where
$(Q,d_{Q})$ is a bounded metric space.  Then the following holds. 
\end{prop}
\begin{enumerate}
\item Suppose $X$ is continuous in probability, with a modulus of continuity
in probability $\delta{}_{Cp}$.\emph{ \index{continuity in probability of r.f.}
Let $\varepsilon>0$ be arbitrary. Define $\delta{}_{cp}(\varepsilon)\equiv\delta{}_{Cp}(2^{-2}(1\wedge\varepsilon)^{2})>0$.
Then, for each $s,t\in Q$ with $d_{Q}(t,s)<\delta{}_{cp}(\varepsilon)$,
there exists measurable set $D_{t,s}$ with $PD_{t,s}^{c}\leq\varepsilon$
such that
\[
d(X(t,\omega),X(s,\omega))\leq\varepsilon
\]
for each $\omega\in D_{t,s}$. Conversely, if there exists an operation
$\delta{}_{cp}$ with the above described properties, then the r.f.
$X$ is continuous in probability, with a modulus of continuity in
probability $\delta{}_{Cp}$ defined by $\delta{}_{Cp}(\varepsilon)\equiv\delta{}_{cp}(\frac{1}{2}\varepsilon)$
for each $\varepsilon>0$.}
\item \emph{$X$ is }continuous a.u.\emph{ \index{continuity a.u. of r.f.}
iff, for each $\varepsilon>0$ and $s\in Q$ there exists a measurable
 set $D_{s}\subset domain(X_{s})$ with $P(D_{s}^{c})<\varepsilon$,
such that, for each $\alpha>0$, there exists }$\delta'_{cau}(\alpha,\varepsilon)>0$\emph{
such that
\[
d(X(t,\omega),X(s,\omega))\leq\alpha
\]
for each $t\in Q$ with $d_{Q}(t,s)<\delta'_{cau}(\alpha,\varepsilon)$,
for each $\omega\in D_{s}$. }
\item \emph{$X$ is }a.u. continuous\emph{ \index{a.u. continuity of r.f.}
iff, for each $\varepsilon>0$, there exists a measurable set $D$
with $P(D^{c})<\varepsilon$, such that, for each $\alpha>0$, there
exists }$\delta'_{auc}(\alpha,\varepsilon)>0$\emph{ such that
\[
d(X(t,\omega),X(s,\omega))\leq\alpha
\]
for each $s,t\in domain(X_{s})\cap domain(X_{t})$ with $d_{Q}(t,s)<\delta'_{auc}(\alpha,\varepsilon)$,
for each $\omega\in D$. Moreover, if such an operation $\delta'_{auc}$
exists, then $X$ has a modulus of a.u. continuity given by }$\delta{}_{auc}(\varepsilon)\equiv\delta'_{auc}(\varepsilon,\varepsilon)$
for each $\varepsilon>0$.
\end{enumerate}
\begin{proof}
As usual, write $\widehat{d}\equiv1\wedge d$. 

1. Suppose $X$ is continuous in probability, with a modulus of continuity
in probability $\delta{}_{Cp}$. Let $\varepsilon>0$ be arbitrary.
Write $\varepsilon'\equiv1\wedge\varepsilon$. Suppose $s,t\in Q$
are arbitrary with 
\[
d_{Q}(t,s)<\delta{}_{cp}(\varepsilon)\equiv\delta{}_{Cp}(2^{-2}(1\wedge\varepsilon)^{2})\equiv\delta{}_{Cp}(2^{-2}\varepsilon'{}^{2}).
\]
Then, by Definition \ref{Def. continuity in prob, continuity a.u., and a.u. continuity}
of $\delta{}_{Cp}$ as a modulus of continuity in probability, we
have $E\widehat{d}(X_{t},X_{s})\leq2^{-2}\varepsilon'{}^{2}<\varepsilon'^{2}$.
Take any $\alpha\in(2^{-1}\varepsilon',\varepsilon')$ such that the
set $D_{t,s}\equiv(\widehat{d}(X_{t},X_{s})\leq\alpha)$ is $\mathrm{measurable}$.
Then Chebychev's inequality implies that $P(D_{t,s}^{c})<\alpha<\varepsilon'\leq\varepsilon$.
Moreover, for each $\omega\in D_{t,s}$, we have $\widehat{d}(X(t,\omega),X(s,\omega))\leq\alpha<\varepsilon'\leq\varepsilon$.
Thus the operation $\delta{}_{cp}$ has the properties described in
Assertion 1.

Conversely, suppose $\delta{}_{cp}$ is an operation with the properties
described in Assertion 1. Let $\varepsilon>0$ be arbitrary. Let $s,t\in Q$
be arbitrary with $d_{Q}(t,s)<\delta{}_{cp}(\frac{1}{2}\varepsilon)$.
Then, by hypothesis, there exists a $\mathrm{measurable}$ subset
$D_{t,s}$ with $PD_{t,s}^{c}\leq\frac{1}{2}\varepsilon$ such that,
for each $\omega\in D_{t,s}$, we have $d(X_{t}(\omega),X_{s}(\omega))\leq\frac{1}{2}\varepsilon$.
It follows that $E(1\wedge d(X_{t},X_{s}))\leq\frac{1}{2}\varepsilon+PD_{t,s}^{c}\leq\varepsilon$.
Thus $X$ is continuous in probability.

2. Suppose $X$ is continuous a.u., with $\delta_{cau}$ as a modulus
of continuity a.u. Let $\varepsilon>0$ and $s\in Q$ be arbitrary.
Then there exists, for each $k\geq1$, a $\mathrm{measurable}$ set
$D_{s,k}$ with $P(D_{s,k}^{c})<2^{-k}\varepsilon$ such that, for
each $\omega\in D_{s,k}$, we have
\[
d(X(t,\omega),X(s,\omega))\leq2^{-k}\varepsilon
\]
for each $t\in Q$ with $d_{Q}(t,s)<\delta_{cau}(2^{-k}\varepsilon)$.
Let $D_{s}\equiv\bigcap_{k=1}^{\infty}D_{s,k}$. Then $P(D_{s}^{c})<\sum_{k=1}^{\infty}2^{-k}\varepsilon=\varepsilon.$
Now let $\alpha>0$ be arbitrary. Let $k\geq1$ be so large that $2^{-k}<\alpha$,
and let 
\[
\delta'_{cau}(\alpha,\varepsilon)\equiv\delta_{cau}(2^{-k}\varepsilon).
\]
Consider each $\omega\in D_{s}$ and $t\in Q$ with $d_{Q}(t,s)<\delta'_{cau}(\alpha,\varepsilon)$.
Then $\omega\in D_{s,k}$ and $d_{Q}(t,s)<\delta_{cau}(2^{-k}\varepsilon)$.
Hence 
\[
d(X(t,\omega),X(s,\omega))\leq2^{-k}\varepsilon<\alpha.
\]
Thus the operation $\delta'_{cau}$ has the described properties in
Assertion 2.

Conversely, let $\delta'_{cau}$ be an operation with the properties
described in Assertion 2. Let $\varepsilon>0$ be arbitrary. Let $s\in Q$
be arbitrary. Then there exists a $\mathrm{measurable}$ set $D_{s}$
with $P(D_{s}^{c})<\varepsilon$ such that, for each $\omega\in D_{s}$,
and $t\in Q$ with $d_{Q}(t,s)<\delta{}_{cau}(\varepsilon)\equiv\delta'_{cau}(\varepsilon,\varepsilon)$,
we have $d(X(t,\omega),X(s,\omega))\leq\varepsilon$. Thus the r.f.
$X$ is continuous a.u., with the operation $\delta{}_{cau}$ as a
modulus of continuity a.u. Assertion 2 is proved.

3. For Assertion 3, proceed almost verbatim as in the above proof
of Assertion 2. Suppose the r.f. $X$ is a.u. continuous, with $\delta_{auc}$
as a modulus of a.u. continuity. Let $\varepsilon>0$ be arbitrary.
Then there exists, for each $k\geq1$, a $\mathrm{measurable}$ set
$D_{k}$ with $P(D_{k}^{c})<2^{-k}\varepsilon$ such that, for each
$\omega\in D_{k}$, we have
\[
d(X(t,\omega),X(s,\omega))\leq2^{-k}\varepsilon
\]
for each $s,t\in Q$ with $d_{Q}(t,s)<\delta_{auc}(2^{-k}\varepsilon)$.
Let $D\equiv\bigcap_{k=1}^{\infty}D_{k}$. Then $P(D^{c})<\sum_{k=1}^{\infty}2^{-k}\varepsilon=\varepsilon.$
Now let $\alpha>0$ be arbitrary. Let $k\geq1$ be so large that $2^{-k}<\alpha$,
and let 
\[
\delta'_{auc}(\alpha,\varepsilon)\equiv\delta_{auc}(2^{-k}\varepsilon).
\]
Consider each $\omega\in D$ and $s,t\in Q$ with $d_{Q}(t,s)<\delta'_{auc}(\alpha,\varepsilon)$.
Then $\omega\in D_{k}$ and $d_{Q}(t,s)<\delta_{auc}(2^{-k}\varepsilon)$.
Hence 
\[
d(X(t,\omega),X(s,\omega))\leq2^{-k}\varepsilon<\alpha.
\]
Thus the operation $\delta'_{auc}$ has the properties described in
Assertion 3.

Conversely, let $\delta'_{auc}$ be an operation with the properties
described in Assertion 3. Let $\varepsilon>0$ be arbitrary. Then
there exists a $\mathrm{measurable}$ set $D$ with $P(D^{c})<\varepsilon$
such that, for each $\omega\in D$, and $s,t\in Q$ with $d_{Q}(t,s)<\delta{}_{auc}(\varepsilon)\equiv\delta'_{auc}(\varepsilon,\varepsilon)$,
we have $d(X(t,\omega),X(s,\omega))\leq\varepsilon$. Thus the r.f.
$X$ is a.u. continuous, with the operation $\delta{}_{auc}$ as a
modulus of a.u. continuity. Assertion 3 is proved.
\end{proof}
\begin{prop}
\label{Prop. R.F. continuous a.u. =00003D> continuous in prob} \textbf{\emph{(a.u.
Continuity implies continuity a.u., etc.)}} Let $X:Q\times\Omega\rightarrow S$
be a r.f., where $(S,d)$ is a locally compact metric space and where
$(Q,d_{Q})$ is a bounded metric space. Then a.u. continuity of $X$
implies continuity a.u. which in turn implies continuity in probability.
\end{prop}
\begin{proof}
Let $\varepsilon>0$ be arbitrary. Suppose $X$ is a.u. continuous,
with modulus of a.u. continuity given by $\delta_{auc}$. Let $D$
be a $\mathrm{measurable}$ set satisfying Condition 3 in Definition
\ref{Def. continuity in prob, continuity a.u., and a.u. continuity}.
Then $D_{s}\equiv D$ satisfies Condition 2 in Definition \ref{Def. continuity in prob, continuity a.u., and a.u. continuity}.
Accordingly, $X$ is continuous a.u.

Now suppose $X$ is continuous a.u., with modulus of continuity a.u.
given by $\delta_{cau}$. Let $D_{s}$ be a $\mathrm{measurable}$
set satisfying Condition 2 in Definition \ref{Def. continuity in prob, continuity a.u., and a.u. continuity}.
Then $D_{t,s}\equiv D_{s}$, satisfies the conditions in Assertion
1 of Proposition \ref{Prop. Alternative characterization of r.f. continuity},
provided that we define $\delta{}_{cp}\equiv\delta_{cau}$. Accordingly,
$X$ is continuous in probability.
\end{proof}
\begin{defn}
\label{Def. Marginal Distributions of a rnadom field} \textbf{(Marginal
distributions of a r.f.).} Let $X:Q\times\Omega\rightarrow S$ be
a r.f. Let $n\geq1$ be arbitrary, and let $t\equiv(t_{1},\cdots,t_{n})$
be an arbitrary sequence in $Q$. Let $F_{t(1),\cdots,t(n)}$ denote
the distribution induced on $(S^{n},d^{(n)})$ by the r.v. $(X_{t(1)},\cdots,X_{t(n)})$.Then
\begin{equation}
F_{t(1),\cdots,t(n)}f\equiv Ef(X_{t(1)},\cdots,X_{t(n)})\label{eq:temp-197}
\end{equation}
for each $f\in C_{ub}(S^{n})$. We call the indexed family
\[
F\equiv\{F_{t(1),\cdots,t(n)}:n\geq1\:\mbox{and}\:t_{1},\cdots,t_{n}\in Q\}
\]
the family of \index{marginal distributions}\emph{ marginal distributions}
\emph{of} $X$. We will say that the r.f. $X$ extends the family
$F$ of finite joint distributions, and that $X$ is an \emph{extension}\index{extension of family of f.j.d.'s}
of $F$.

Let $X':Q\times\Omega'\rightarrow S$ be a r.f. with sample space
$(\Omega',L',E')$. Then $X$ and $X'$ are said to be \index{equivalent stochastic processes}\emph{equivalent}
if their marginal distributions at each finite sequence in $Q$ are
the same. In other words, $X$ and $X'$ are said to be equivalent
if
\[
Ef(X_{t(1)},\cdots,X_{t(n)})=E'f(X'_{t(1)},\cdots,X'_{t(n)})
\]
for each $f\in C_{ub}(S^{n})$, for each sequence $(t_{1},\cdots,t_{n})$
in $Q$, for each $n\geq1$. In short, two r.f.'s are equivalent if
they extend the same family of finite joint distributions. $\square$ 
\end{defn}

\section{Consistent Family of Finite Joint Distributions}

In the last section, we saw that each r.f. gives rise to a family
of marginal distributions. Conversely, we seek conditions for a family
$F$ of finite joint distributions to be the family of marginal distributions
of some r.f. We will presently show that a necessary condition is
consistency, to be defined next. In the following chapters we will
present various sufficient conditions on $F$ for the construction
of r.f.'s with $F$ as the family of marginal distributions and with
desired properties of sample functions. 
\begin{defn}
\label{Def.Consistent Finite Joint  Distributions-1} \textbf{(Consistent
family of f.j.d.'s).} Let $Q$ be a set. Suppose, for each $n\geq1$
and for each finite sequence $t_{1},\cdots,t_{n}$ in $Q$, a distribution
$F_{t(1),\cdots,t(n)}$ is given on the locally compact metric space
$(S^{n},d^{(n)})$, which will be called a \emph{finite joint distribution}\index{finite joint distribution},
or \index{f.j.d.}\emph{f.j.d.} for short. Then the indexed family
\[
F\equiv\{F_{t(1),\cdots,t(n)}:n\geq1\:\mbox{and}\:t_{1},\cdots,t_{n}\in Q\}
\]
is said to be a \index{consistent family of f.j.d.'s} \emph{consistent
family }of\emph{ f.j.d.'s }with parameter set\emph{ $Q$} and state
space $S$, if the following \emph{\index{consistency condition}
Kolmogorov consistency condition} is satisfied. 

Let $n,m\geq1$ be arbitrary. Let $t\equiv(t_{1},\cdots,t_{m})$ be
an arbitrary sequence in $Q$, and let $i\equiv(i_{1},\cdots,i_{n})$
be an arbitrary sequence in $\{1,\cdots,m\}$. Define the continuous
function $i^{*}:S^{m}\rightarrow S^{n}$ by 
\begin{equation}
i^{*}(x_{1},\cdots,x_{m})\equiv(x_{i(1)},\cdots,x_{i(n)})\label{eq:temp-171}
\end{equation}
for each $(x_{1},\cdots,x_{m})\in S^{m}$, and call $i^{*}$ the \emph{dual
function of the sequence}\index{dual function of a sequence} $i$.
Then, for each $f\in C_{ub}(S^{n})$, we have
\begin{equation}
F_{t(1),\cdots,t(m)}(f\circ i^{*})=F_{t(i(1)),\cdots,t(i(n))}f,\label{eq:temp-196}
\end{equation}
or, in short,
\[
F_{t}(f\circ i^{*})=F_{t\circ i}(f).
\]

We will let $\widehat{F}(Q,S)$ denote the set of consistent families
of \emph{f.j.d.'s }with parameter set\emph{ $Q$} and state space
$S$. When there is little risk of confusion, we will call a consistent
family of f.j.d.'s simply a \emph{consistent family}\index{consistent family}.

$\square$
\end{defn}
Note that for an arbitrary $f\in C_{ub}(S^{n})$ we have $f\circ i^{*}\in C_{ub}(S^{m})$
and so $f\circ i^{*}$ is integrable relative to $F_{t(1),\cdots,t(m)}$.
Hence the left-hand side of equality \ref{eq:temp-196} makes sense. 

When the parameter set is a countable discrete subset of $R$, we
have the following proposition with a simple sufficient condition
for the construction of a consistent family of f.j.d.'s.

First some notations.
\begin{defn}
\textbf{\label{Def. Notations for sequences} (Notations for sequences).
}Given any sequence $(a_{1},\cdots,a_{m})$ of objects, we will use
the shorter notation $a$ for the sequence. When there is little risk
of confusion, we will write $\kappa\sigma\equiv\kappa\circ\sigma$
for the composite of two functions $\sigma:A\rightarrow B$ and $\kappa:B\rightarrow C$.
Separately, for each $m\geq n\geq1$, define the sequence 
\[
\kappa\equiv\kappa_{n,m}:\{1,\cdots,m-1\}\rightarrow\{1,\cdots,m\}
\]
by 
\[
(\kappa_{1},\cdots,\kappa_{m-1})\equiv(1,\cdots,\widehat{n},\cdots,m)\equiv(1,\cdots,n-1,n+1,\cdots,m),
\]
where the caret on the top of an element in a sequence signifies the
omission of that element. Let $\kappa^{*}\equiv\kappa_{n,m}^{*}$
denote the dual function of sequence $\kappa$. Thus 
\[
\kappa^{*}x=\kappa^{*}(x_{1},\cdots,x_{m})\equiv x\kappa=(x_{\kappa(1)},\cdots,x_{\kappa(m-1)})=(x_{1},\cdots,\widehat{x_{n}},\cdots,x_{m})
\]
for each $x\equiv(x_{1},\cdots,x_{m})\in S^{m}$. In words, the function
$\kappa_{n,m}^{*}$ deletes the $n$-th entry of the sequence $(x_{1},\cdots,x_{m})$.
\end{defn}
\begin{lem}
\label{Lem. Consistency when parameter set is discrete subset of R}
\textbf{\emph{(Consistency when parameter set is discrete subset of
$R$).}} Let $(S,d)$ be a locally compact metric space. Let $Q$
be an arbitrary metrically discrete subset of $R$. Suppose, for each
$m\geq1$ and nonincreasing sequence $r\equiv(r_{1},\cdots,r_{m})$
in $Q$, a distribution $F_{r(1),\cdots,r(m)}$ on $(S^{m},d^{m})$
is given, such that
\begin{equation}
F_{r(1),\cdots,\widehat{r(n)}\cdots,r(m)}f=F_{r(1),\cdots,r(m)}(f\circ\kappa_{n,m}^{*}),\label{eq:temp-429-1}
\end{equation}
or, equivalently, 
\begin{equation}
F_{r\circ\kappa(n,m)}f=F_{r}(f\circ\kappa_{n,m}^{*}),\label{eq:temp-429-1-1}
\end{equation}
for each $f\in C_{ub}(S^{m-1})$, for each $n=1,\cdots,m$. Then the
family 
\[
F\equiv\{F_{r(1),\cdots,r(n)}:n\geq1;r_{1}\leq\cdots\leq r_{n}\;in\;Q\}
\]
of f.j.d.'s can be uniquely extended to a consistent family of f.j.d.'s
\begin{equation}
F=\{F_{s(1),\cdots,s(m)}:m\geq1;s_{1},\cdots,s_{m}\in Q\}\label{eq:temp-399}
\end{equation}
with parameter $Q$. 
\end{lem}
\begin{proof}
1. Let the integers $m,n$, with $m\geq n\geq1$, and the increasing
sequence $r\equiv(r_{1},\cdots,r_{m})$ in $Q$ be arbitrary. Let
$r'\equiv(r'_{1},\cdots,r'_{n})$ be an arbitrary subsequence of $r$.
Then $m>h\equiv m-n\geq0$. Moreover, $r'$ can be obtained by deleting
$h$ elements in the sequence $r$. Specifically, $r'=\kappa^{*}r=r\kappa$,
where 
\[
\kappa\equiv\kappa_{n(h),m}\kappa_{n(h-1),m-1}\cdots\kappa_{n(1),m-h}:\{1,\cdots,m-h\}\rightarrow\{1,\cdots,m\}
\]
if $h>0$, and where $\kappa$ is the identity function if $h=0$.
Hence, by repeated application of equality \ref{eq:temp-429-1} in
the hypothesis, we obtain
\[
F_{r'}f=F_{r\kappa}f=F_{r}(f\circ\kappa^{*})
\]
for each $f\in C(\overline{S}^{m-h})$. 

2. Let the sequence $s\equiv(s_{1},\cdots,s_{p})$ in $Q$ be arbitrary.
Let $r\equiv(r_{1},\cdots,r_{m})$ be an arbitrary increasing sequence
in $Q$ such that $s$ is a sequence in $\{r_{1},\cdots,r_{m}\}$.
Then, because the sequence $r$ is increasing, there exists a unique
function $\sigma:\{1,\cdots,p\}\rightarrow\{1,\cdots,m\}$ such that
$s=r\sigma$. Let $f\in C_{ub}(S^{p},d^{p})$ be arbitrary. Define
\begin{equation}
\overline{F}_{s(1),\cdots,s(p)}f\equiv\overline{F}_{s}f\equiv F_{r}(f\circ\sigma*)\label{eq:temp-403-2-1-1}
\end{equation}

We will verify that $\overline{F}_{s}f$ is well defined. To that
end, let $r'\equiv(r'_{1},\cdots,r'_{m'})$ be a second increasing
sequence in $Q$ such that $s$ is a sequence in $\{r'_{1},\cdots,r'_{m'}\}$,
and let $\sigma':\{1,\cdots,p\}\rightarrow\{1,\cdots,m'\}$ be the
corresponding function such that $s=r'\sigma'$. We need to verify
that $F_{r}(f\circ\sigma*)=F_{r'}(f\circ\sigma'*)$. To that end,
let $\bar{r}\equiv(\bar{r}_{1},\cdots,\bar{r}_{\overline{m}})$ be
an arbitrary supersequence of $r$ and $r'$. Then, $s=r\sigma=\bar{r}\kappa\sigma$,
while $s=r'\sigma'=\bar{r}\kappa'\sigma'$. Hence, by uniqueness,
we have $\kappa\sigma=\kappa'\sigma'$. Consequently, 
\[
F_{r}(f\circ\sigma*)=F_{\bar{r}}(f\circ\sigma*\circ\kappa^{*})=F_{\bar{r}}(f\circ\sigma'*\circ\kappa'^{*})=F_{r'}(f\circ\sigma'*).
\]
This shows that $\overline{F}_{s}f$ is well defined in equality \ref{eq:temp-403-2-1-1}.
The same equality says that $\overline{F}_{s}$ is the distribution
induced by the r.v. $\sigma*:(S^{m},\overline{C_{ub}}(S^{m},d^{m}),F_{r})\rightarrow(S^{p},d^{p})$,
where $\overline{C_{ub}}(S^{m},d^{m})^{-}$ stands for the completion
of $C(S^{m},d^{m})$ relative to the distribution $F_{r}$. In particular,
$\overline{F}_{s(1),\cdots,s(p)}\equiv\overline{F}_{s}$ is a distribution.

3. Next, let $s\equiv(s_{1},\cdots,s_{q})$ be arbitrary sequence
in $Q$, and let $(s_{i(1)},\cdots,s_{i(p)})$ be an arbitrary subsequence
of $s$. Write $i\equiv(i_{1},\cdots,i_{p})$. Let the increasing
sequence $r\equiv(r_{1},\cdots,r_{m})$ be arbitrary such that $s$
is a sequence in $\{r_{1},\cdots,r_{m}\}$, and let $\sigma:\{1,\cdots,q\}\rightarrow\{1,\cdots,m\}$
such that $s=r\sigma$. Then $si=r\sigma i$. Hence, for each $f\in C_{ub}(S^{p},d^{p})$,
we have 
\[
\overline{F}_{si}f=\overline{F}_{r\sigma i}f\equiv F_{r}(f\circ i^{*}\sigma^{*})\equiv\overline{F}_{s}(f\circ i^{*}).
\]
Thus the family 
\[
\overline{F}\equiv\{\overline{F}_{s(1),\cdots,s(p)}:p\geq1;s_{1},\cdots,s_{p}\in Q\}
\]
of f.j.d.'s is consistent. 

4. Lastly, let $s\equiv(s_{1},\cdots,s_{q})$ be arbitrary increasing
sequence in $Q$. Write $r\equiv s$. Then $s=r\sigma$ where $\sigma:\{1,\cdots,q\}\rightarrow\{1,\cdots,q\}$
is the identity function. Hence 
\[
\overline{F}_{s(1),\cdots,s(q)}f\equiv F_{r(1),\cdots,r(q)}f\circ\sigma^{*}=F_{s(1),\cdots,s(q)}f
\]
for each $f\in\overline{C_{ub}}(S^{q},d^{q})$. In other words, $\overline{F}_{s(1),\cdots,s(q)}\equiv F_{s(1),\cdots,s(q)}$.
Thus the family $\overline{F}$ is an extension of the family $F$,
and we can simply write $F$ for $\overline{F}$. The lemma is proved.
\end{proof}
The next lemma extends the consistency condition \ref{eq:temp-196}
to integrable functions.
\begin{prop}
\label{Prop. consistency of f.j.d. on L} \textbf{\emph{(Consistency
condition extends to integrable functions).}} Suppose the consistency
condition \ref{eq:temp-196} holds for each $f\in C_{ub}(S^{n})$
and for the family $F$ of f.j.d.'s. Then a real valued function $f$
on $S^{n}$ is integrable relative to \textup{$F_{t(i(1)),\cdots,t(i(n))}$}\textup{\emph{
iff $f\circ i^{*}$ is integrable relative to $F_{t(1),\cdots,t(m)}$,
in which case}} condition \ref{eq:temp-196} holds for $f$.
\end{prop}
\begin{proof}
Since $i^{*}:(S^{m},d^{(m)})\rightarrow(S^{n},d^{(n)})$ is uniformly
continuous, $i^{*}$ is a r.v. on the completion of $(S^{m},C(S^{m}),F_{t(1),\cdots,t(m)})$
and has values in $S^{n}$, whence it induces a distribution on $S^{n}$.
Equality \ref{eq:temp-196} then implies that the distribution thus
induced is equal to $F_{t(i(1)),\cdots,t(i(n))}$. Therefore, according
to Proposition \ref{Prop. L(X)=00003D=00007Bf(X): f in L_X=00007D},
a function $f:S^{n}\rightarrow R$ is integrable relative to $F_{t(i(1)),\cdots,t(i(n))}$
iff $f(i^{*})$ is integrable relative to $F_{t(1),\cdots,t(m)}$,
in which case
\[
F_{t(i(1)),\cdots,t(i(n))}f=F_{t(1),\cdots,t(m)}f(i^{*})\equiv F_{t(1),\cdots,t(m)}f\circ i^{*}.
\]
\end{proof}
\begin{prop}
\label{Prop. Marginal Distributions  are consistent}\textbf{\emph{
(Marginal distributions are consistent).}} Let $X:Q\times\Omega\rightarrow S$
be a r.f. Then the family $F$ of marginal distributions of $X$ is
consistent.
\end{prop}
\begin{proof}
Let $n,m\geq1$ and $f\in C_{ub}(S^{n})$ be arbitrary. Let $t\equiv(t_{1},\cdots,t_{m})$
be an arbitrary sequence in $Q$, and let $i\equiv(i_{1},\cdots,i_{n})$
be an arbitrary sequence in $\{1,\cdots,m\}$. Using the defining
equalities \ref{eq:temp-197} and \ref{eq:temp-171}, we obtain 
\[
F_{t(1),\cdots,t(m)}(f\circ i^{*})\equiv E((f\circ i^{*})(X_{t(1)},\cdots,X_{t(m)}))
\]
\[
\equiv Ef(X_{t(i(1))},\cdots,X_{t(i(n))})\equiv F{}_{t(i(1)),\cdots,t(i(n))}f.
\]
Thus the consistency condition \ref{eq:temp-196} holds.
\end{proof}
\begin{defn}
\label{Def. Restriction of set of Consistent  Families to subset of paramter set}
\textbf{(Restriction to a subset of the parameter set).} Let $(S,d)$
be a locally compact metric space. Recall that $\widehat{F}(Q,S)$
is the set of consistent families of f.j.d.'s with parameter set $Q$
and state space $S$. Let $Q'$ be any subset of $Q$. For each $F\in\widehat{F}(Q,S)$
define 
\begin{equation}
F|Q'\equiv\Phi_{Q|Q'}(F)\equiv\{F_{s(1),\cdots,s(n)}:n\geq1;s_{1},\cdots,s_{n}\in Q'\}\label{eq:temp-134-3}
\end{equation}
and call $F|Q'$ the \emph{restriction} \index{restriction of a consistent family}
of the consistent family $F$ to $Q'$. The function 
\[
\Phi_{Q|Q'}:\widehat{F}(Q,S)\rightarrow\widehat{F}(Q',S)
\]
will be called the restriction mapping of consistent families with
parameter set $Q$ to consistent families with parameter set $Q'$. 

Let $\widehat{F}_{0}\subset\widehat{F}(Q,S)$ be arbitrary. Denote
its image under the mapping $\Phi_{Q|Q'}$ by 
\begin{equation}
\widehat{F}_{0}|Q'\equiv\Phi_{Q|Q'}(\widehat{F}_{0})=\{F|Q':F\in\widehat{F}_{0}\},\label{eq:temp-194}
\end{equation}
and call $\widehat{F}_{0}|Q'$ the restriction of the set $\widehat{F}_{0}$
of consistent families to $Q'$. 
\end{defn}
$\square$

We next introduce a metric on the set $\widehat{F}(Q,S)$ when $Q$
is countably infinite. 
\begin{defn}
\label{Def. Marginal metric} \textbf{(Marginal metric on set of consistent
families of f.j.d.'s with countably infinite parameter set). }Let
$(S,d)$ be a locally compact metric space, with a binary approximation
$\xi$ relative to some fixed reference point $x_{\circ}$. Let $n\geq1$
be arbitrary. Recall that $\xi^{n}$ is the $n$-th power of $\xi$,
and is a binary approximation of $(S^{n},d^{n})$ relative to $x_{\circ}^{(n)}\equiv(x_{\circ},\cdots,x_{\circ})\in S^{n}$,
as in Definition \ref{Def. Finite product and power of binary approx of locally compact S}.
Recall from Definition \ref{Def. Distribution metric}, the distribution
metric $\rho_{Dist,\xi^{n}}$ on the set of distributions on $(S^{n},d^{n})$,
and, from Proposition \ref{Prop. rho_xi convergence=00003D Weak Convergence},
that sequential convergence relative to $\rho_{Dist,\xi^{n}}$ is
equivalent to weak convergence. 

Let $Q\equiv\{t_{1},t_{2},\cdots\}$ be an enumerated, countably infinite,
parameter set. Recall that $\widehat{F}(Q,S)$ is the set of consistent
families of f.j.d.'s with parameter set $Q$ and  state space $S$.
 Define a metric $\widehat{\rho}_{Marg,\xi,Q}$ on $\widehat{F}(Q,S)$
by 
\begin{equation}
\widehat{\rho}_{Marg,\xi}(F,F')\equiv\widehat{\rho}_{Marg,\xi,Q}(F,F')\equiv\sum_{n=1}^{\infty}2^{-n}\rho_{Dist,\xi^{n}}(F_{t(1),\cdots,t(n)},F'_{t(1),\cdots,t(n)})\label{eq:temp-208}
\end{equation}
for each $F,F'\in\widehat{F}(Q,S)$. The next lemma proves that $\widehat{\rho}_{Marg,\xi}$
metric on families of f.j.d.'s with countable parameters is indeed
a metric. We will call $\widehat{\rho}_{Marg,\xi}$ the \index{marginal metric}\emph{marginal
metric} for the set $\widehat{F}(Q,S)$ of consistent families of
f.j.d.'s, relative to the binary approximation $\xi$ of the locally
compact state space $(S,d)$. Note that $\widehat{\rho}_{Marg,\xi}\leq1$
because $\rho_{Dist,\xi^{n}}\leq1$ for each $n\geq1$. We emphasize
that the metric $\widehat{\rho}_{Marg,\xi,Q}$ depends on the ordering
in the enumerated set $Q$. Two different enumeration leads to two
different metrics, which are however equivalent. We drop the subscript
$Q$ when it is understood from context.

As observed above, sequential convergence relative to $\rho_{Dist,\xi^{n}}$
is equivalent to weak convergence of distributions on $(S^{n},d^{n})$,
for each $n\geq1$. Hence, for each sequence $(F^{(m)})_{m=0,1,2,\cdots}$in
$\widehat{F}(Q,S)$, we have $\widehat{\rho}_{Marg,\xi}(F^{(m)},F^{(0)})\rightarrow0$
iff $F_{t(1),\cdots,t(n)}^{(m)}\Rightarrow F_{t(1),\cdots,t(n)}^{(0)}$
as $m\rightarrow\infty$, for each $n\geq1$. $\square$
\end{defn}
\begin{lem}
\label{Lem. rho_F_tilde is a metric if Q is countable.} The marginal
metric $\widehat{\rho}_{Marg,\xi}\equiv\widehat{\rho}_{Marg,\xi,Q}$
defined in Definition \ref{Def. Metric on random fields w/ countable parameters}
is indeed a metric. 
\end{lem}
\begin{proof}
1. Symmetry and triangle inequality for $\widehat{\rho}_{Marg,\xi}$
follow from their respective counterparts for $\rho_{Dist,\xi^{n}}$
 for each $n\geq1$ in the defining equality \ref{eq:temp-208}.

2. Suppose $F=F'$. Then each summand in the right-hand side of equality
\ref{eq:temp-208} vanishes. Consequently $\widehat{\rho}_{Marg,\xi}(F,F')=0$.
Conversely, suppose $F,F'\in\widehat{F}(Q,S)$ are such that $\widehat{\rho}_{Marg,\xi}(F,F')=0$.
For each $n\geq1$, the defining equality \ref{eq:temp-208} implies
that $\rho_{Dist,\xi^{n}}(F_{t(1),\cdots,t(n)},F'_{t(1),\cdots,t(n)})=0$.
Hence, since $\rho_{Dist,\xi^{n}}$ is a metric, we have $F_{t(1),\cdots,t(n)}=F'_{t(1),\cdots,t(n)}$,
for each $n\geq1$. Now let $m\geq1$ and $s_{1},\cdots,s_{m}\in Q$
be arbitrary. Then there exists $n\geq1$ so large that $s_{k}=t_{i(k)}$
for some $i_{k}\in\{1,\cdots,n$\}, for each $k=1,\cdots,m$. By the
consistency condition \ref{eq:temp-196}
\[
F_{s(1),\cdots,s(m)}f=F_{t(i(1)),\cdots,t(i(m))}f=F_{t(1),\cdots,t(n)}(f\circ i^{*})
\]
\[
=F'_{t(1),\cdots,t(n)}(f\circ i^{*})=F'_{t(i(1)),\cdots,t(i(m))}f=F'_{s(1),\cdots,s(m)}f
\]
for each $f\in C_{ub}(S^{n})$. We thus see that $F_{s(1),\cdots,s(m)}=F'_{s(1),\cdots,s(m)}$
as distributions on $S^{m}$ for each $s_{1},\cdots,s_{m}\in Q$.
In other words, $F=F'$. Summing up, $\widehat{\rho}_{Marg,\xi}$
is a metric.
\end{proof}
\begin{defn}
\label{Def. Continuous in probabillity of Finite Joint Distributions}
\textbf{(Continuity in probability of consistent families).} Let $(S,d)$
be a locally compact metric space. Write $\widehat{d}\equiv1\wedge d$.
Let $(Q,d_{Q})$ be a metric space. Recall that $\widehat{F}(Q,S)$
is the set of consistent families of f.j.d.'s with parameter space
$Q$ and state space $(S,d)$. Let $F\in\widehat{F}(Q,S)$ be arbitrary.

1. Suppose $(Q,d_{Q})$ is bounded. Suppose, for each $\varepsilon>0$,
there exists $\delta_{Cp}(\varepsilon)>0$ such that
\[
F_{s,t}\widehat{d}\leq\varepsilon
\]
for each $s,t\in Q$ with $d_{Q}(s,t)<\delta_{Cp}(\varepsilon)$.
Then the consistent family $F$ of f.j.d.'s is said to be \index{f.j.d.'s, continuous in probability}\emph{
continuous in probability}, with $\delta_{Cp}$ as a \index{modulus of continuity of f.j.d.'s}\emph{modulus
of continuity in probability}.

2. More generally, let the metric space $(Q,d_{Q})$ be arbitrary,
not necessarily bounded. Then the consistent family $F$ of f.j.d.'s
is said to be \index{f.j.d.'s, continuous in probability}\emph{ continuous
in probability} if, for each bounded subset $K$ of $Q$, the restricted
consistent family $F|K$ is continuous in probability.  We will let
$\widehat{F}_{Cp}(Q,S)$ denote the subset of $\widehat{F}(Q,S)$
whose members are continuous in probability. $\square$ 
\end{defn}
\begin{lem}
\label{Lem. Continuity of finite joint distributions} \textbf{\emph{(Continuity
in probability extends to f.j.d.'s of higher dimensions).}} Let $(S,d)$
be a locally compact metric space. Let $(Q,d_{Q})$ be a bounded metric
space. Suppose the consistent family $F$ of f.j.d.'s with state space
$S$ and parameter space $Q$ is continuous in probability, with a
modulus of continuity in probability $\delta_{Cp}$. 

Let $m\geq1$ be arbitrary. Let $f\in C_{ub}(S^{m},d^{m})$ be arbitrary
with a modulus of continuity $\delta_{f}$ and with $|f|\leq1$. Let
and $\varepsilon>0$ be arbitrary. Then there exists $\delta_{fjd}(\varepsilon,m,\delta_{f},\delta_{Cp})>0$
such that, for each $s_{1},\cdots,s_{m},t_{1},\cdots,t_{m}\in Q$
with 
\begin{equation}
\bigvee_{k=1}^{m}d_{Q}(s_{k},t_{k})<\delta_{fjd}(\varepsilon,m,\delta_{f},\delta_{Cp}),\label{eq:temp-141}
\end{equation}
we have
\begin{equation}
|F_{s(1),\cdots,s(m)}f-F_{t(1),\cdots,t(m)}f|\leq\varepsilon.\label{eq:temp-192}
\end{equation}
\end{lem}
\begin{proof}
Let $m\geq1$ and $f\in C_{ub}(S^{m},d^{m})$ be as given. Write 
\[
\alpha\equiv\frac{1}{8}m^{-1}\varepsilon(1\wedge\delta_{f}(\frac{\varepsilon}{2}))
\]
and define 
\[
\delta_{fjd}(\varepsilon,m,\delta_{f},\delta_{Cp})\equiv\delta_{Cp}(\alpha).
\]
Suppose $s_{1},\cdots,s_{m},t_{1},\cdots,t_{m}\in Q$ satisfy inequality
\ref{eq:temp-141}. Then
\begin{equation}
\bigvee_{k=1}^{m}d_{Q}(s_{k},t_{k})<\delta_{Cp}(\alpha).\label{eq:temp-182}
\end{equation}
Let $i\equiv(1,\cdots,m)$ and $j\equiv(m+1,\cdots,2m)$. Thus $i$
and $j$ are sequences in $\{1,\cdots,2m\}$. Let $x\in S^{2m}$ be
arbitrary. Then
\[
(f\circ i^{*})(x_{1},\ldots,x_{2m})\equiv f(x_{i(1)},\cdots,x_{i(m)})=f(x_{1},\ldots,x_{m})
\]
and
\[
(f\circ j^{*})(x_{1},\ldots,x_{2m})\equiv f(x_{j(1)},\cdots,x_{j(m)})=f(x_{m+1},\ldots,x_{2m}),
\]
where $i^{*},j^{*}$ are as defined in Definition \ref{Def.Consistent Finite Joint  Distributions-1}
relative to $i,j$ respectively. Consider each $k=1,\cdots,m$. Let
$h\equiv(k,m+k)$. Thus $h$ is a sequence in $\{1,\cdots,2m\}$.
Let 
\[
(r_{1},\cdots,r_{2m})\equiv(s_{1},\cdots,s_{m},t_{1},\cdots,t_{m}).
\]
Then
\[
F_{r(1),\cdots,r(2m)}(\widehat{d}\circ h^{*})=F_{r(h(1)),r(h(2))}\widehat{d}
\]
\begin{equation}
=F_{r(k),r(m+k)}\widehat{d}=F_{s(k),t(k)}\widehat{d}<\alpha,\label{eq:temp-195}
\end{equation}
where the inequality follows from inequality \ref{eq:temp-182} in
view of the definition of $\delta_{Cp}$ as a modulus of continuity
in probability of the family $F$. Now take any 
\[
\delta_{0}\in(\frac{1}{2}(1\wedge\delta_{f}(\frac{\varepsilon}{2})),1\wedge\delta_{f}(\frac{\varepsilon}{2})).
\]
Let
\[
A_{k}\equiv\{x\in S^{2m}:\widehat{d}(x_{k},x_{m+k})>\delta_{0}\}=(\widehat{d}\circ h^{*}>\delta_{0})\subset S^{2m}.
\]
In view of inequality \ref{eq:temp-195}, Chebychev\noun{'s} inequality
yields
\[
F_{r(1),\cdots,r(2m)}(A_{k})=F_{r(1),\cdots,r(2m)}(\widehat{d}\circ h^{*}>\delta_{0})<\delta_{0}^{-1}\alpha
\]
\[
<2(1\wedge\delta_{f}(\frac{\varepsilon}{2}))^{-1}\alpha=\frac{1}{4}m^{-1}\varepsilon\alpha^{-1}\alpha=\frac{1}{4}m^{-1}\varepsilon.
\]
Let $A\equiv\bigcup_{k=1}^{m}A_{k}\subset S^{2m}$. Then
\[
F_{r(1),\cdots,r(2m)}(A)\leq\sum_{k=1}^{m}F_{r(1),\cdots,r(2m)}(A_{k})\leq\frac{1}{4}\varepsilon.
\]
Now consider each $x\in A^{c}$. We have 
\[
1\wedge d(x_{k},x_{m+k})\equiv\widehat{d}(x_{k},x_{m+k})\leq\delta_{0}<1
\]
for each $k=1,\cdots,m$, whence
\[
d^{m}((x_{1},\ldots,x_{m}),(x_{m+1},\ldots,x_{2m}))\equiv\bigvee_{k=1}^{m}d(x_{k},x_{m+k})\leq\delta_{0}<\delta_{f}(\frac{\varepsilon}{2}).
\]
Consequently
\[
|(f\circ i^{*})(x)-(f\circ j^{*})(x)|=|f(x_{1},\ldots,x_{m})-f(x_{m+1},\ldots,x_{2m})|<\frac{\varepsilon}{2}
\]
for each $x\in A^{c}$. By hypothesis, $|f|\leq1$. Hence
\[
|F_{s(1),\cdots,s(m)}f-F_{t(1),\cdots,t(m)}f|=|F_{r(i(1)),\cdots,r(i(n))}f-F_{r(j(1)),\cdots,r(j(n))}f|
\]
\[
=|F_{r(1),\cdots,r(2m)}(f\circ i^{*})-F_{r(1),\cdots,r(2m)}(f\circ j^{*})|
\]
\[
=|F_{r(1),\cdots,r(2m)}(f\circ i^{*}-f\circ j^{*})|
\]
\[
\leq\frac{\varepsilon}{2}F_{r(1),\cdots,r(2m)}(A^{c})+2F_{r(1),\cdots,r(2m)}(A)
\]
\begin{equation}
\leq\frac{\varepsilon}{2}+\frac{\varepsilon}{2}=\varepsilon.\label{eq:temp-193}
\end{equation}
as desired. 
\end{proof}
$\,$
\begin{defn}
\label{Def. Metric on  of continuous in prob families of finite joint distributions.}
\textbf{(Metric space of consistent families which are continuous
in probability).} Let $(S,d)$ be a locally compact metric space,
with a reference point $x_{\circ}\in S$ and a binary approximation
$\xi\equiv(A_{n})_{n=1,2,\cdots}$ relative to $x_{\circ}$. Let $(Q,d_{Q})$
be a locally compact metric space. Let $Q_{\infty}\equiv\{q_{1},q_{2},\cdots\}$
be an arbitrary enumerated, countably infinite, and  dense subset
of $Q$. 

Recall that $\widehat{F}(Q,S)$ is be the set of consistent families
of f.j.d.'s with parameter set $Q$ and state space $S$. Let $\widehat{F}_{Cp}(Q,S)$
denote the subset of $\widehat{F}(Q,S)$ whose members are continuous
in probability. 

Relative to the countably infinite parameter subset $Q_{\infty}$
and the binary approximation $\xi$, define a metric $\widehat{\rho}_{Cp,\xi,Q|Q(\infty)}$
on $\widehat{F}_{Cp}(Q,S)$ by
\[
\widehat{\rho}_{Cp,\xi,Q|Q(\infty)}(F,F')\equiv\widehat{\rho}_{Marg,\xi,Q(\infty)}(F|Q_{\infty},F'|Q_{\infty})
\]
\begin{equation}
\equiv\sum_{n=1}^{\infty}2^{-n}\rho_{Dist,\xi^{n}}(F_{q(1),\cdots,q(n)},F'_{q(1),\cdots,q(n)})\label{eq:temp-208-1}
\end{equation}
for each $F,F'\in\widehat{F}_{Cp}(Q,S)$, where $\widehat{\rho}_{Marg,\xi,Q(\infty)}$
is the marginal metric on $\widehat{F}(Q_{\infty},S)$ introduced
in Definition \ref{Def. Marginal metric}. In other words,
\[
\widehat{\rho}_{Cp,\xi,Q|Q(\infty)}(F,F')\equiv\widehat{\rho}_{Marg,\xi,Q(\infty)}(\Phi_{Q|Q(\infty)}(F),\Phi_{Q|Q(\infty)}(F'))
\]
for each $F,F'\in\widehat{F}_{Cp}(Q,S)$. The next lemma shows that
$\widehat{\rho}_{Cp,\xi,Q|Q(\infty)}$ is  indeed a metric. Then,
trivially, the mapping
\[
\Phi_{Q|Q(\infty)}:(\widehat{F}_{Cp}(Q,S),\widehat{\rho}_{Cp,\xi,Q|Q(\infty)})\rightarrow(\widehat{F}_{Cp}(Q_{\infty},S),\widehat{\rho}_{Marg,\xi,Q(\infty)})
\]
is an isometry. Note that $0\leq\widehat{\rho}_{Cp,\xi,Q|Q(\infty)}\leq1$.
 $\square$
\end{defn}
\begin{lem}
\textup{\emph{\label{Lem.Each F in  F_hat is equicontiuous in prob =00003D> Rho_Hat is a metric}
The function $\widehat{\rho}_{Cp,\xi,Q|Q(\infty)}$ defined in Definition
\ref{Def. Metric on  of continuous in prob families of finite joint distributions.}
is a metric on }}\emph{$\widehat{F}_{Cp}(Q,S)$}. 
\end{lem}
\begin{proof}
Suppose $F,F'\in\widehat{F}_{Cp}(Q,S)$ are such that $\widehat{\rho}_{Cp,\xi,Q|Q(\infty)}(F,F')=0$.
By the defining equality \ref{eq:temp-208-1}, we have $\widehat{\rho}_{Marg,\xi,Q(\infty)}(F|Q_{\infty},F'|Q_{\infty})=0$.
Hence, since $\widehat{\rho}_{Marg,\xi,Q(\infty)}$ is a metric on
$\widehat{F}(Q_{\infty},S)$, we have $F|Q_{\infty}=F'|Q_{\infty}$.
In other words,
\[
F_{q(1),\cdots,q(n)}=F'_{q(1),\cdots,q(n)}
\]
for each $n\geq1$. Hence, for each $m\geq1$ and each $s_{1},\cdots,s_{m}\in Q_{\infty}$,
we can let $n\geq1$ be so large that $\{s_{1},\cdots,s_{m}\}\subset\{q_{1},\cdots,q_{n}\}$
and obtain, consistency of $F$, 
\begin{equation}
F_{s(1),\cdots,s(m)}=F'_{s(1),\cdots,s(m)}\label{eq:temp-228}
\end{equation}
Now let $m\geq1$ and $t_{1},\cdots,t_{m}\in Q$ be arbitrary. Let
$f\in C(S^{m})$ be arbitrary. For each $i=1,\ldots,m$, let $(s_{i}^{(p)})_{p=1,2,\cdots}$
be a sequence in $Q_{\infty}$ with $d_{Q}(s_{i}^{(p)},r_{i})\rightarrow0$
as $p\rightarrow\infty$. Then there exists a bounded subset $K\subset Q$
such that ($t_{1},\cdots,t_{m})$ and $(s_{i}^{(p)})_{p=1,2,\cdots;i=1,\cdots,m}$
are in $K$. Since $F|K$ is continuous in probability, we have, by
Lemma \ref{Lem. Continuity of finite joint distributions}, 
\[
F_{s(p,1),\cdots,s(p,m)}f\rightarrow F_{t(1),\cdots,t(m)}f,
\]
where we write $s(p,i)\equiv s_{i}^{(p)}$ to lessen the burden on
subscripts. Similarly,
\[
F'_{s(p,1),\cdots,s(p,m)}f\rightarrow F'_{t(1),\cdots,t(m)}f,
\]
On the other hand
\[
F_{s(p,1),\cdots,s(p,m)}f=F'_{s(p,1),\cdots,s(p,m)}f
\]
in view of equality \ref{eq:temp-228}. Combining, 
\[
F'_{t(1),\cdots,t(m)}f=F_{t(1),\cdots,t(m)}f.
\]
We conclude that $F=F'$.

Conversely, suppose $F=F'$. Then trivially $\widehat{\rho}_{Cp,\xi,Q|Q(\infty)}(F,F')=0$
from equality \ref{eq:temp-208-1}. The triangle inequality and symmetry
of $\widehat{\rho}_{Cp,\xi,Q|Q(\infty)}$  follow from equality \ref{eq:temp-208-1}
and from the fact that $\rho_{Dist,\xi^{n}}$ is a metric for each
$n\geq1$. Summing up, $\widehat{\rho}_{Cp,\xi,Q|Q(\infty)}$ is a
metric.
\end{proof}

\section{Daniell-Kolmogorov Extension}

In this and the next section, let $Q\equiv\{t_{1},t_{2},\cdots\}$
denote a countable parameter set. For simplicity of presentation,
and without loss of generality, we will assume that $t_{n}=n$ for
each $n\geq1$. Thus 
\[
Q\equiv\{t_{1},t_{2},\cdots\}\equiv\{1,2,\cdots\}.
\]
However we state the theorems in terms of a more general countable
set $Q\equiv\{t_{1},t_{2},\cdots\},$ for ease of later reference
when more structure on the set $Q$ is introduced, when, for example,
the set $Q$ is the set of dyadic rationals in $[0,\infty)$.

Recall that $\widehat{F}(Q,S)$ is the set of consistent families
of f.j.d.'s with parameter set $Q$ and the locally compact state
space $(S,d)$. We will prove the \index{Daniell-Kolmogorov Extension}Daniell-Kolmogorov
Extension Theorem, which constructs, for each member $F\in\widehat{F}(Q,S)$,
a probability space $(S^{\infty},L,E)$ and a r.f. $U:Q\times(S^{\infty},L,E)\rightarrow S$
with marginal distributions given by $F$. 

Furthermore, we will prove the uniform metrical continuity of the
Daniell-Kolmogorov Extension, with a modulus of continuity dependent
only on a modulus of local compactness $\left\Vert \xi\right\Vert $
of $(S,d)$. Said metrical continuity implies continuity relative
to weak convergence.

Recall that $[\cdot]_{1}$ is an operation which assigns to each $c\geq0$
an integer $[c]_{1}$ in the interval $(c,c+2)$. As usual, for arbitrary
symbols $a$ and $b$, we will write $a_{b}$ and $a(b)$ interchangeable.
\begin{defn}
\label{Def. Path space, and coordinate function} \textbf{(Path space,
coordinate function, and distributions on path space).} Let $S^{Q}\equiv\prod_{t\in Q}S$
denote the space of functions from $Q$ to $S$, called the \emph{path
space}\index{path space}. Relative to the enumerated set $Q$, define
a complete metric $d^{Q}$ on $S^{Q}$, by
\[
d^{Q}(x,y)\equiv\sum_{i=0}^{\infty}2^{-i}(1\wedge d(x_{t(i)},y_{t(i)}))
\]
for arbitrary $x,y\in S^{Q}$. Define the function $U:Q\times S^{Q}\rightarrow S$
by $U(r,v)\equiv v_{r}$ for each $(r,v)\in Q\times S^{Q}$. The function
$U$ is called the \index{coordinate function}\emph{coordinate function}
$Q\times S^{Q}$. Note that $d^{Q}\leq1$ and that $(S^{Q},d^{Q})$
is compact if $(S,d)$ is compact. 

Conforming to usage in Definition \ref{Def. distributions on complete metric space},
we will let $\widehat{J}(S^{Q},d^{Q})$ denote the set of distributions
on the complete path space $(S^{Q},d^{Q})$. 
\end{defn}
$\square$ 
\begin{thm}
\label{Thm. Compact Daniell-Kolmogorov Extension}\textbf{\emph{ (Compact
Daniell-Kolmogorov Extension).}} Suppose the metric space $(S,d)$
is compact. Then there exists a function 
\[
\overline{\Phi}_{DK}:\widehat{F}(Q,S)\rightarrow\widehat{J}(S^{Q},d^{Q})
\]
such that, for each consistent family of f.j.d.'s $F\in\widehat{F}(Q,S)$,
the distribution $E\equiv\overline{\Phi}_{DK}(F)$ satisfies the conditions
\emph{(i)} the coordinate function 
\[
U:Q\times(S^{Q},L,E)\rightarrow(S,d)
\]
is a r.f., where $L$ is the completion of $C(S^{Q},d^{Q})$ relative
to the distribution $E$, and \emph{(ii) the r.f.} $U$ has marginal
distributions given by the family $F$. 

The function $\overline{\Phi}_{DK}$ will be called the \index{Compact Daniell-Kolmogorov Extension}
\emph{Compact Daniell-Kolmogorov Extension}.
\end{thm}
\begin{proof}
Note that, since $(S,d)$ is compact by hypothesis, its countably
infinite power $(S^{Q},d^{Q})=(S^{\infty},d^{\infty})$ is compact. 

1. Consider each $F\equiv\{F_{1,\cdots,k}:k\geq1\}\in\widehat{F}(Q,S)$.
Let $f\in C(S^{\infty},d^{\infty})$ be arbitrary, with a modulus
of continuity $\delta_{f}$. For each $n\geq1$, define the function
$f_{n}\in C(S^{n},d^{n})$ by
\begin{equation}
f_{n}(x_{1},x_{2},\cdots,x_{n})\equiv f(x_{1},x_{2},\cdots,x_{n},x_{\circ},x_{\circ},\cdots)\label{eq:temp-49-2}
\end{equation}
for each $(x_{1},x_{2},\cdots,x_{n})\in S^{n}$. Consider each $m\geq n\geq1$
so large that $2^{-n}<\delta_{f}(\varepsilon)$. Define the function
$f_{n,m}\in C(S^{m},d^{m})$ by
\begin{equation}
f_{n,m}(x_{1},x_{2},\cdots,x_{m})\equiv f_{n}(x_{1},x_{2},\cdots,x_{n})\equiv f(x_{1},x_{2},\cdots,x_{n},x_{\circ},x_{\circ},\cdots)\label{eq:temp-49-1-1}
\end{equation}
for each $(x_{1},x_{2},\cdots,x_{m})\in S^{m}$. Consider the initial-section
subsequence $i\equiv(i_{1},\cdots,i_{n})\equiv(1,\cdots,n)$ of the
sequence $(1,\cdots,m)$. Let $i^{*}:S^{m}\rightarrow S^{n}$ be the
dual of the sequence $i$, as in Definition \ref{Def.Consistent Finite Joint  Distributions-1}.
Then, for each $(x_{1},\cdots,x_{m})\in S^{m}$, we have
\[
f_{n,m}(x_{1},x_{2},\cdots,x_{m})\equiv f_{n}(x_{1},x_{2},\cdots,x_{n})=f_{n}(x_{i(1)},x_{i(2)},\cdots,x_{i(n)})=f_{n}\circ i*(x_{1},x_{2},\cdots,x_{m}).
\]
In short,
\[
f_{n,m}=f_{n}\circ i*,
\]
whence, by the consistency of the family $F$ of f.j.d.'s, we obtain
\begin{equation}
F_{1,\cdots,m}f_{n,m}=F_{1,\cdots,m}f_{n}\circ i*=F_{i(1),\cdots,i(n)}f_{n}=F_{1,\cdots,n}f_{n}.\label{eq:temp-203-1}
\end{equation}
At the same time,
\[
d^{\infty}((x_{1},x_{2},\cdots,x_{m},x_{\circ},x_{\circ},\cdots),(x_{1},x_{2},\cdots,x_{n},x_{\circ},x_{\circ},\cdots))
\]
\[
\equiv\sum_{i=1}^{n}2^{i}(1\wedge d((x_{i},x_{i}))+\sum_{i=n+1}^{m}2^{i}(1\wedge d((x_{i},x_{\circ}))+\sum_{i=m+1}^{\infty}2^{i}(1\wedge d((x_{\circ},x_{\circ}))
\]
\[
=0+\sum_{i=n+1}^{m}2^{i}(1\wedge d((x_{i},x_{\circ}))+0\leq2^{-n}<\delta_{f}(\varepsilon)
\]
for each $(x_{1},x_{2},\cdots,x_{m})\in S^{m}$. Hence
\[
|f_{m}(x_{1},x_{2},\cdots,x_{m})-f_{n,m}(x_{1},x_{2},\cdots,x_{m})|
\]
\[
=|f(x_{1},x_{2},\cdots,x_{m},x_{\circ},x_{\circ},\cdots)-f(x_{1},x_{2},\cdots,x_{n},x_{\circ},x_{\circ},\cdots)|<\varepsilon
\]
for each $(x_{1},x_{2},\cdots,x_{m})\in S^{n}$. Consequently, $|F_{1,\cdots,m}f_{m}-F_{1,\cdots,m}f_{n,m}|\leq\varepsilon.$
Combined with equality \ref{eq:temp-203-1}, this yields
\begin{equation}
|F_{1,\cdots,m}f_{m}-F_{1,\cdots,n}f_{n,}|\leq\varepsilon,\label{eq:temp-48}
\end{equation}
where $m\geq n\geq1$ are arbitrary with $2^{-n}<\delta_{f}(\varepsilon)$.
Thus we see that the sequence $(F_{1,\cdots,n}f_{n,})_{n=1,2,\cdots}$
of real numbers is Cauchy, and has a limit. Define
\begin{equation}
Ef\equiv\lim_{n\rightarrow\infty}F_{1,\cdots,n}f_{n,}.\label{eq:temp-269}
\end{equation}
Letting $m\rightarrow\infty$ in inequality \ref{eq:temp-48}, we
obtain
\begin{equation}
|Ef-F_{1,\cdots,n}f_{n,}|\leq\varepsilon,\label{eq:temp-48-1}
\end{equation}
where $n\geq1$ is arbitrary with $2^{-n}<\delta_{f}(\varepsilon)$. 

3. We proceed to prove that $E$ is an integration on the compact
metric space $(S^{\infty},d^{\infty})$ in the sense of Definition
\ref{Def. integration on loc compact space}. We will first verify
that the function $E$ is linear. To that end, let $f,g\in C(S^{\infty},d^{\infty})$
and $a,b\in R$ be arbitrary. For each $n\geq1$, define the function
$f_{n}$ relative to $f$ as in equality \ref{eq:temp-49-2}. Similarly
define the functions $g_{n},(af+bg)_{n}$ relative to the functions
$g,af+bg$ respectively. Then the defining equality \ref{eq:temp-49-2}
implies that $(af+bg)_{n}=af_{n}+bg_{n}$ for each $n\geq1$. Hence
\[
E(af+bg)\equiv\lim_{n\rightarrow\infty}F_{1,\cdots,n}(af+bg)_{n,}=\lim_{n\rightarrow\infty}F_{1,\cdots,n}(af_{n}+bg_{n})
\]
\[
=a\lim_{n\rightarrow\infty}F_{1,\cdots,n}f_{n,}+b\lim_{n\rightarrow\infty}F_{1,\cdots,n}g_{n,}\equiv aEf+bEg._{,}
\]
Thus $E$ is a linear function. Moreover, in the special case where
$f\equiv1$, we have
\begin{equation}
E1\equiv Ef\equiv\lim_{n\rightarrow\infty}F_{1,\cdots,n}f_{n,}=\lim_{n\rightarrow\infty}F_{1,\cdots,n}1=1>0.\label{eq:temp-299}
\end{equation}
Inequality \ref{eq:temp-299} immediately shows that the triple $(S^{\infty},C(S^{\infty},d^{\infty}),E)$
satisfies Condition (i) of Definition \ref{Def. integration on loc compact space}.
It remains to verify also Condition (ii), the positivity condition,
of Definition \ref{Def. integration on loc compact space}. To that
end, let $f\in C(S^{\infty},d^{\infty})$ be arbitrary with $Ef>0$.
Then, by equality \ref{eq:temp-269}, we have $F_{1,\cdots,n}f_{n,}>0$
for some $n\geq1$. Hence, since $F_{1,\cdots,n}$ is a distribution,
there exists $(x_{1},x_{2},\cdots,x_{n})\in S^{n}$ such that $f_{n,}(x_{1},x_{2},\cdots,x_{n})>0$.
Therefore 
\[
f(x_{1},x_{2},\cdots,x_{n},x_{\circ},x_{\circ},\cdots)\equiv f_{n,}(x_{1},x_{2},\cdots,x_{n})>0.
\]
Thus the positivity condition is also verified. Accordingly, $E$
is an integration on the compact metric space $(S^{\infty},d^{\infty})$. 

4. Since the compact metric space $(S^{\infty},d^{\infty})$ is bounded,
and since $E1=1,$ Lemma \ref{Lem. Distribution basics} implies that
$E$ is a distribution on $(S^{\infty},d^{\infty})$, and that the
completion $(S^{\infty},L,E)$ of the integration space $(S^{\infty},C(S^{\infty},d^{\infty}),E)$
is a probability space. In symbols, $E\in\widehat{J}(S^{Q},d^{Q})$.
Define $\overline{\Phi}_{DK}(F)\equiv E$. Thus we have constructed
the function $\overline{\Phi}_{DK}:\widehat{F}(Q,S)\rightarrow\widehat{J}(S^{Q},d^{Q})$.

5. It remains to show that the coordinate function $U:Q\times(S^{\infty},L,E)\rightarrow S$
is a r.f. with marginal distributions given by the family $F$. To
that end, let $m\geq n\geq1$ and $g\in C(S^{n},d^{n})$ be arbitrary.
Define a function $f\in C(S^{\infty},d^{\infty})$ by 
\begin{equation}
f(x_{1},x_{2},\cdots)\equiv g(x_{1},\cdots,x_{n})\label{eq:temp-300}
\end{equation}
for each $x\equiv(x_{1},x_{2},\cdots)\in S^{\infty}$, and define
the function $f_{m}$ relative to $f$ as in equality \ref{eq:temp-49-2}.
Thus
\begin{equation}
f_{m}(x_{1},x_{2},\cdots,x_{m})\equiv f(x_{1},x_{2},\cdots,x_{m},x_{\circ},x_{\circ},\cdots)\equiv g(x_{1},\cdots,x_{n}),\label{eq:temp-49-1-1-1-1}
\end{equation}
for each $(x_{1},x_{2},\cdots,x_{m})\in S^{m}$. Consequently, 
\[
f_{m}(x_{1},x_{2},\cdots,x_{m})=g(x_{1},\cdots,x_{n})=f_{n}(x_{1},x_{2},\cdots,x_{n})=f_{n}\circ i*(x_{1},x_{2},\cdots,x_{m})
\]
for each $(x_{1},x_{2},\cdots,x_{m})\in S^{m}$, where $i\equiv(i_{1},\cdots,i_{n})\equiv(1,\cdots,n)$
is the initial-section subsequence of the sequence $(1,\cdots,m)$,
and where $i^{*}:S^{m}\rightarrow S^{n}$ be is the dual of the sequence
$i$,. In short, $f_{m}=f_{n}\circ i*$. At the same time,
\[
g(U_{1},\cdots,U_{n})(x)\equiv g(U_{1}(x),\cdots,U_{n}(x))=g(x_{1},\cdots,x_{n})=f(x)
\]
for each $x\equiv(x_{1},x_{2},\cdots)\in S^{\infty}$. In short, $g(U_{1},\cdots,U_{n})=f$.
Combining,
\begin{equation}
Eg(U_{1},\cdots,U_{n})=Ef\equiv\lim_{m\rightarrow\infty}F_{1,\cdots,m}f_{m}=\lim_{m\rightarrow\infty}F_{1,\cdots,m}f_{n}\circ i*=F_{1,\cdots,n}f_{n}=F_{1,\cdots,n}g,\label{eq:temp-301}
\end{equation}
where $n\geq1$ and $g\in C(S^{n},d^{n})$ are arbitrary, and where
the fourth equality is by the consistency of the family $F$ of f.j.d.'s.
Equality \ref{eq:temp-301} implies that $(U_{1},\cdots,U_{n})$ is
a r.v. on the sample space $(S^{\infty},L,E)$, with values in $(S^{n},d^{n})$
and with distribution $F_{1,\cdots,n}$. It follows that $U_{n}$
is a r.v. on the sample space $(S^{\infty},L,E)$, with values in
$S$, where $n\geq1$ is arbitrary. Summing up, we conclude that the
coordinate function 
\[
U:Q\times(S^{\infty},L,E)\rightarrow S
\]
is a r.f. Equality \ref{eq:temp-301} says that $U$ has marginal
distributions given by the family $F$. The theorem is proved.
\end{proof}
We proceed to prove the continuity of the Compact Daniell-Kolmogorov
Extension relative to the two metrics specified next. 
\begin{defn}
\textbf{\label{Def. Specification of binary approximation. Recall marginal metric}(Specification
of binary approximation of state space, and related marginal metric
on the set of consistent families of f.j.d.'s).} Let $\xi\equiv(A_{k})_{k=1,2,}$
be an arbitrary binary approximation of the locally compact state
space $(S,d)$ relative to the reference point $x_{\circ}$, in the
sense of Definition \ref{Def. Binary approximationt and Modulus of local compactness}.
Recall that $\widehat{F}(Q,S)$ is then equipped with the marginal
metric $\widehat{\rho}_{Marg,\xi,Q}$ defined relative to $\xi$ in
Definition \ref{Def. Marginal metric}, and that sequential convergence
relative to this metric $\widehat{\rho}_{Marg,\xi,Q},$ is equivalent
to weak convergence of corresponding sequences of f.j.d.'s. 
\end{defn}
$\square$
\begin{defn}
\textbf{\label{Specification of binary approximation of comact path space, and distribution on said path space}
(Specification of binary approximation of compact path space, and
distribution metric on the set of distributions on said path space).}
Suppose the state space $(S,d)$ is compact. Let $\xi\equiv(A_{k})_{k=1,2,}$
be an arbitrary binary approximation of $(S,d)$ relative to the reference
point $x_{\circ}$. Recall that, since the metric space $(S^{\infty},d^{\infty})$
is compact, the countable power $\xi^{\infty}\equiv(B_{k})_{k=1,2,}$
of $\xi$ is defined and is a binary approximation of $(S^{\infty},d^{\infty})$,
according to Definition \ref{Def. Countable product of binary approxximations; compact}
and Lemma \ref{Lem. Countable Product of binary approximation is a binary approxt; compacts}.
Recall that, since $(S,d)$ is compact by assumption, the\emph{ set
$\widehat{J}(S^{\infty},d^{\infty})$} of distributions is equipped
with the distribution metric $\rho_{Dist,\xi^{\infty}}$, defined
relative to $\xi^{\infty}$ in Definition \ref{Def. Distribution metric},
and that sequential convergence relative to this metric $\rho_{Dist,\xi^{\infty}}$
is equivalent to weak convergence. Note that the metric $\rho_{Dist,\xi^{\infty}}$
is defined only when the state space $(S,d)$ is compact. Write $\rho_{Dist,\xi^{Q}}\equiv\rho_{Dist,\xi^{\infty}}$.
$\square$
\end{defn}
\begin{thm}
\textbf{\noun{\label{Thm. Continuity of the Compact DK Extension}
(}}\textbf{\emph{Continuity of the Compact Daniell-Kolmogorov Extension).}}
Suppose $(S,d)$ is compact. Let $\xi\equiv(A_{k})_{k=1,2,}$ be an
arbitrary binary approximation of $(S,d)$ relative to the reference
point $x_{\circ}$. Then the Compact Daniell-Kolmogorov Extension
\[
\overline{\Phi}_{DK}:(\widehat{F}(Q,S),\widehat{\rho}_{Marg,\xi,Q})\rightarrow(\widehat{J}(S^{Q},d^{Q}),\rho_{Dist,\xi^{Q}})
\]
constructed in Theorem \ref{Thm. Compact Daniell-Kolmogorov Extension}
is uniformly continuous, with modulus of continuity $\hspace*{1em}$
$\overline{\delta}_{DK}(\cdot,\left\Vert \xi\right\Vert )$ dependent
only on the modulus of local compactness $\left\Vert \xi\right\Vert \equiv(|A_{k}|)_{k=1,2,}$
of the compact metric space $(S,d)$.
\end{thm}
\begin{proof}
1. Let $\varepsilon\in(0,1)$ be arbitrary. For abbreviation, write
$c\equiv2^{4}\varepsilon{}^{-1}$ and $\alpha\equiv2^{-1}\varepsilon$.
Let $m\equiv[\log_{2}2\varepsilon{}^{-1}]_{1}$. Define the operation
$\delta_{c}$ by 
\begin{equation}
\delta_{c}(\varepsilon')\equiv c^{-1}\varepsilon'\label{eq:temp-209-1-2}
\end{equation}
for each $\varepsilon'>0$. Take $n\geq1$ so large that 
\[
2^{-n}<\frac{1}{3}c^{-1}\alpha=\delta_{c}(\frac{\alpha}{3}).
\]
Note that, by the definition of the operation $[\cdot]_{1}$, we have
$2\varepsilon{}^{-1}<2^{m}<2\varepsilon{}^{-1}\cdot2^{2}=\frac{1}{2}c.$
Hence
\[
2^{m+1}<c
\]
and 
\[
2^{-m}<\frac{1}{2}\varepsilon.
\]

2. Let $F,F'\in\widehat{F}(Q,S)$ be arbitrary, with $F\equiv\{F_{1,\cdots,k}:k\geq1\}$
and $F'\equiv\{F'_{1,\cdots,k}:k\geq1\}.$ Consider the distributions
$F_{1,\cdots,n}\in F$ and $F'_{1,\cdots,n}\in F$'. Since, by hypothesis
$d\leq1$, we have the product metric $d^{n}\leq1$ also. Hence, trivially,
the distributions $F_{1,\cdots,n},F'_{1,\cdots,n}$ on $(S^{n},d^{n})$
have modulus of tightness equal to $1$. Let $\xi^{n}$ be the power
$n$-th power of $\xi$ , as in Definition \ref{Def. Finite product and power of binary approx of locally compact S}.
Thus $\xi^{n}$ is a binary approximation for $(S^{n},d^{n})$. Recall
the distribution metric $\rho_{Dist,\xi^{n}}$ relative to $\xi^{n}$
on the set of distributions on $(S^{n},d^{n})$, as introduced in
Definition \ref{Def. Distribution metric}. Then Assertion 1 of Proposition
\ref{Prop. Modulus of continuity of J->Jf for fixed Lipshitz f} applies
to the compact metric space $(S^{n},d^{n})$ and the distribution
metric $\rho_{Dist,\xi^{n}}$, to yield 
\begin{equation}
\widetilde{\Delta}\equiv\widetilde{\Delta}(\frac{\alpha}{3},\delta_{c},1,\left\Vert \xi^{n}\right\Vert )>0\label{eq:temp-252-2-1}
\end{equation}
such that, if
\[
\rho_{Dist,\xi^{n}}(F_{1,\cdots,n},F'_{1,\cdots,n})<\widetilde{\Delta},
\]
then
\begin{equation}
|F_{1,\cdots,n}g-F'_{1,\cdots,n}g|<\frac{\alpha}{3}\label{eq:temp-186-1-1}
\end{equation}
for each $g\in C(S^{n},d^{n})$ with $|g|\leq1$ and with modulus
of continuity $\delta_{c}$. Recall from Lemma \ref{Lem. Product separant is a separant}
that the modulus of local compactness $\left\Vert \xi^{n}\right\Vert $
of $(S^{n},d^{n})$ is determined by the modulus of local compactness
$\left\Vert \xi\right\Vert $ of $(S,d)$. Hence we can define
\begin{equation}
\overline{\delta}_{DK}(\varepsilon)\equiv\overline{\delta}_{DK}(\varepsilon,\left\Vert \xi\right\Vert )\equiv2^{-n}\widetilde{\Delta}>0.\label{eq:temp-253-2-1}
\end{equation}
We will prove that $\overline{\delta}_{DK}$ is a modulus of continuity
of the Compact Daniell-Kolmogorov Extension $\overline{\Phi}_{DK}$. 

3. Suppose, for that purpose, that
\begin{equation}
\sum_{k=1}^{\infty}2^{-k}\rho_{Dist,\xi^{k}}(F_{1,\cdots,k},F'_{1,\cdots,k})\equiv\widehat{\rho}_{Marg,\xi,Q}(F,F')<\overline{\delta}_{DK}(\varepsilon).\label{eq:temp-260-2-1}
\end{equation}
We need to show that, then, 
\[
\rho_{Dist,\xi^{\infty}}(E,E')<\varepsilon
\]
where $E\equiv\overline{\Phi}_{DK}(F)$ and $E'\equiv\overline{\Phi}_{DK}(F')$.

4. To that end, let
\[
\pi\equiv(\{g_{k,x}:x\in B_{k}\})_{k=1,2,\cdots}
\]
be the partition of unity of the compact metric space $(S^{\infty},d^{\infty})$
determined by its binary approximation $\xi^{\infty}\equiv(B_{k})_{k=1,2,}$,
as in Definition \ref{Def. Partition of unity for locally compact (S,d)}.
In other words, the family $\{g_{k,x}:x\in B_{k}\}$ of basis functions
is the $2^{-k}$\emph{-}partition of unity of $(S^{\infty},d^{\infty})$
determined by the enumerated finite subset $B_{k}$, for each $k\geq1$.
Moreover, according to Definition \ref{Def. Distribution metric},
we have

\begin{equation}
\rho_{Dist,\xi^{\infty}}(E,E')\equiv\sum_{k=1}^{\infty}2^{-k}|B_{k}|^{-1}\sum_{x\in B(k)}|Eg_{k,x}-E'g_{k,x}|\label{eq:temp-119-3}
\end{equation}

5. Next note that inequality \ref{eq:temp-260-2-1} immediately yields
\begin{equation}
\rho_{Dist,\xi^{n}}(F_{1,\cdots,n},F'_{1,\cdots,n})<2^{n}\overline{\delta}_{DK}(\varepsilon)=\widetilde{\Delta}.\label{eq:temp-220-2-1}
\end{equation}
Consider each $k=1,\cdots,m$. Let $x\in B_{k}$ be arbitrary. Proposition
\ref{Prop. Properties of  epsilon partition of unity } says that
the basis function $g_{k,x}$ has values in $[0,1]$, and has Lipschitz
constant $2^{k+1}$ on $(S^{\infty},d^{\infty})$, where $2^{k+1}\leq2^{m+1}<c$.
Hence the function $g_{k,x}$ has Lipschitz constant $c$, and, equivalently,
has the modulus of continuity $\delta_{c}$. Now define the function
$g_{k,x,n}\in C(S^{n},d^{n})$ by 
\[
g_{k,x,n}(y_{1},\cdots,y_{n})\equiv g_{k,x}(y_{1},\cdots,y_{n},x_{\circ},x_{\circ},\cdots)
\]
for each $(y_{1},\cdots,y_{n})\in S^{n}$. Then, for each $(z_{1},z_{2},\cdots,z_{n}),(y_{1},y_{2},\cdots,y_{n})\in(S^{n},d^{n}),$
we have
\[
|g_{k,x,n}(z_{1},z_{2},\cdots,z_{n})-g_{k,x,n}(y_{1},y_{2},\cdots,y_{n})|
\]
\[
\equiv|g_{k,x}(z_{1},z_{2},\cdots,z_{n},x_{\circ},x_{\circ},\cdots)-g_{k,x}(y_{1},y_{2},\cdots,y_{n},x_{\circ},x_{\circ},\cdots)|
\]
\[
\leq2^{k+1}d^{\infty}((z_{1},z_{2},\cdots,z_{n},x_{\circ},x_{\circ},\cdots),(y_{1},y_{2},\cdots,y_{n},x_{\circ},x_{\circ},\cdots))
\]
\[
<c\sum_{k=1}^{n}2^{-k}d(z_{k},y_{k})\leq c\bigvee_{k=1}^{n}d(z_{k},y_{k})
\]
\[
\equiv cd^{n}((z_{1},z_{2},\cdots,z_{n}),(y_{1},y_{2},\cdots,y_{n})).
\]
Thus the function $g_{k,x,n}$ also has Lipschitz constant $c$, and,
equivalently, has the modulus of continuity $\delta_{c}$. In addition,
$|g_{k,x}|\leq1,$ whence $|g_{k,x,n}|\leq1$. In view of inequality
\ref{eq:temp-220-2-1}, all the conditions for inequality \ref{eq:temp-186-1-1}
have now been verified for the function $g_{k,x,n}$. Accordingly,
\[
|F_{1,\cdots,n}g_{k,x,n}-F'_{1,\cdots,n}g_{k,x,n}|<\frac{\alpha}{3}.
\]
At the same time, since $2^{-n}<\delta_{c}(\frac{\alpha}{3})$, where
$\delta_{c}$ is a modulus of continuity of the function $g_{k,x,n}$,
inequality \ref{eq:temp-48-1} in the proof of Theorem \ref{Thm. Compact Daniell-Kolmogorov Extension}
applies to the functions $g_{k,x},g_{k,x,n}$ in the place of $f,f_{n}$,
and to the constant $\frac{\alpha}{3}$ in the place of $\varepsilon$,
to yield 
\begin{equation}
|Eg_{k,x}-F_{1,\cdots,n}g_{k,x,n}|\leq\frac{\alpha}{3},\label{eq:temp-48-1-1}
\end{equation}
with a similar inequality when $E,F$ are replaced by $E',F'$ respectively.
The triangle inequality therefore leads to
\begin{equation}
|Eg_{k,x}-E'g_{k,x}|\leq|F_{1,\cdots,n}g_{k,x,n}-F'_{1,\cdots,n}g_{k,x,n}|+\frac{2}{3}\alpha<\frac{\alpha}{3}+\frac{2}{3}\alpha=\alpha,\label{eq:temp-219-1-1}
\end{equation}
where $k=1,\cdots,m$ and $x\in B_{k}$ are arbitrary. It follows
that
\[
\rho_{Dist,\xi^{\infty}}(\overline{\Phi}_{DK}(F),\overline{\Phi}_{DK}(F'))=\rho_{Dist\xi^{\infty}}(E,E')\equiv\sum_{k=1}^{\infty}2^{-k}|B_{k}|^{-1}\sum_{x\in B(k)}|Eg_{k,x}-E'g_{k,x}|
\]
\[
\leq\sum_{k=1}^{m}2^{-k}\alpha+\sum_{k=m+1}^{\infty}2^{-k}<\alpha+2^{-m}<\alpha+\frac{1}{2}\varepsilon=\varepsilon,
\]
where $F,F'\in\widehat{F}(Q,S)$ are arbitrary with $\widehat{\rho}_{Marg,\xi,Q}(F,F')<\overline{\delta}_{DK}(\varepsilon,\left\Vert \xi\right\Vert )$,
where $\varepsilon>0$ is arbitrary. Thus the Compact Daniell-Kolmogorov
Extension
\[
\overline{\Phi}_{DK}:(\widehat{F}(Q,S),\widehat{\rho}_{Marg,\xi,Q})\rightarrow(\widehat{J}(S^{Q},d^{Q}),\rho_{Dist,\xi^{Q}})
\]
is uniformly continuous on $\widehat{F}(Q,S)$, with modulus of continuity
$\overline{\delta}_{DK}(\cdot,\left\Vert \xi\right\Vert )$. The theorem
is proved.
\end{proof}
To generalize Theorems \ref{Thm. Compact Daniell-Kolmogorov Extension}
and \ref{Thm. Continuity of the Compact DK Extension} to a locally
compact, but not necessarily compact, state space $(S,d)$, we (i)
identify each consistent family of f.j.d.'s on the latter with one
on the one-point  compactification $(\overline{S},\overline{d})\equiv(S\cup\{\Delta\},\overline{d})$
whose f.j.d.'s assign probability 1 to powers of $S$, (ii) apply
Theorems \ref{Thm. Compact Daniell-Kolmogorov Extension} and \ref{Thm. Continuity of the Compact DK Extension}
to the compact state space $(\overline{S},\overline{d})$, resulting
in distributions on the path space $(\overline{S}^{Q}\overline{d}^{Q})$,
and (iii) prove that these distributions assign probability 1 to the
path subspace $(S^{Q},d^{Q})$, and can therefore be regarded as distributions
on the latter.

The remainder of this section makes this precise.
\begin{lem}
\textbf{\emph{\label{Lem. Identify family of f.j.d.s on S with family of f.j.d.'s on S_bar}
(Identifying each consistent family of f.j.d.'s with state space $S$
with a consistent family of f.j.d.'s with state space $\overline{S}\equiv S\cup\{\Delta\}$).}}
Suppose $(S,d)$ is locally compact, not necessarily compact. There
exists an injection
\[
\psi:\widehat{F}(Q,S)\rightarrow\widehat{F}(Q,\overline{S})
\]
such that, for each $F\equiv\{F_{1,\cdots.n}:n\geq1\}\in\widehat{F}(Q,S)$,
with $\overline{F}\equiv\{\overline{F}_{1,\cdots,n}:n\geq1\}\equiv\psi(F)$,
we have 
\begin{equation}
\overline{F}_{1,\cdots,n}\overline{f}\equiv F_{1,\cdots.n}(\overline{f}|S^{n})\label{eq:temp-52}
\end{equation}
for each $\overline{f}\in C(\overline{S}^{n},\overline{d}^{n})$,
for each $n\geq1$. Moreover, for each $F\in\widehat{F}(Q,S)$, with
$\overline{F}\equiv\psi(F)$, and for each $n\geq1$, the set $S^{n}$
is a full subset of $\overline{S}^{n}$ relative to the distribution
$\overline{F}_{1,\cdots,n}$ on $(\overline{S}^{n},\overline{d}^{n})$. 

Henceforth, we will identify $F$ with $\overline{F}\equiv\psi(F)$.
In words, each consistent family of f.j.d.'s with state space $S$
is regarded as a consistent family of f.j.d.'s with state space $\overline{S}$
which assign probability $1$ to powers of $S.$ Thus 
\[
\widehat{F}(Q,S)\subset\widehat{F}(Q,\overline{S}).
\]
\end{lem}
\begin{proof}
Consider each $F\equiv\{F_{1,\cdots.n}:n\geq1\}\in\widehat{F}(Q,S)$.
Let $n\geq1$ be arbitrary. Let $\overline{f}\in C(\overline{S}^{n},\overline{d}^{n})$
be arbitrary. Then $\overline{f}|S^{n}\in C_{ub}(S^{n},d^{n})$ by
Corollary \ref{Cor. Extension of f in C(S**n) to point of infinity}.
Hence $\overline{f}|S^{n}$ is integrable relative to the distribution
$F_{1,\cdots.n}$ on $(S^{n},d^{n})$, according to Definition \ref{Def. distributions on complete metric space}.
Therefore we can define 
\begin{equation}
\overline{F}_{1,\cdots,n}\overline{f}\equiv F_{1,\cdots.n}(\overline{f}|S^{n}).\label{eq:temp-442}
\end{equation}
Since $F_{1,\cdots.n}$ is a distribution, the right-hand side is
a linear function of $\overline{f}$. Hence $\overline{F}_{1,\cdots,n}$
is a linear function on $C(\overline{S}^{n},\overline{d}^{n})$. Suppose
$\overline{F}_{1,\cdots,n}\overline{f}>0$. Then $F_{1,\cdots.n}(\overline{f}|S^{n})>0$.
Again, since $F_{1,\cdots.n}$ is a distribution, it follows that
there exists $x\in S$ such $\overline{f}(x)=(\overline{f}|S^{n})(x)>0$.
Thus $\overline{F}_{1,\cdots,n}$ is an integration on the compact
metric space $(\overline{S}^{n},\overline{d}^{n})$. Moreover $\overline{F}_{1,\cdots,n}1\equiv F_{1,\cdots.n}(1)=1$.
Therefore $\overline{F}_{1,\cdots,n}$ is a distribution. 

2. Next, we need to verify that the family $\overline{F}\equiv\{\overline{F}_{1,\cdots,n}:n\geq1\}$
is consistent. To that end, let $m\geq n\geq1$ be arbitrary. Let
$i\equiv(i_{1},\cdots,i_{n})\equiv(1,\cdots,n)$ be the initial subsequence
$(1,\cdots,m)$. Let $i^{*}:\overline{S}^{m}\rightarrow\overline{S}^{n}$
be its dual function. Then trivially $(\overline{f}\circ i^{*})|S^{m}=(\overline{f}|S^{m})\circ i^{*}$
on $S^{n}$. Hence
\[
\overline{F}_{1,\cdots,m}(\overline{f}\circ i^{*})\equiv F_{1,\cdots.m}(\overline{f}\circ i^{*}|S^{m})
\]
\begin{equation}
=F_{1,\cdots.m}((\overline{f}|S^{m})\circ i^{*})=F_{i(1),\cdots,i(m)}(\overline{f}|S^{n})\equiv\overline{F}_{1,\cdots,n}\overline{f},\label{eq:temp-196-1}
\end{equation}
where the third equality follows from the consistency of the family
$F$. Thus the family 
\[
\overline{F}\equiv\psi(F)\equiv\{\overline{F}_{1,\cdots,n}:n\geq1\}
\]
of f.j.d.'s with state space $(\overline{S},\overline{d})$ is consistent.
In other words, $\psi(F)\in\widehat{F}(Q,\overline{S})$. From the
defining equality \ref{eq:temp-442}, we see that, if $F=F'\in\widehat{F}(Q,S)$
then $\psi(F)=\psi(F')$. We conclude that $\psi$ is a well defined
function.

3. Now suppose $F,F'\in\widehat{F}(Q,S)$ are such that $\overline{F}\equiv\psi(F)=\psi(F')\equiv\overline{F}'$.
Let $n\geq1$ be arbitrary. Consider each $f\in C(S^{n},d^{n})$.
Then, by Corollary \ref{Cor. Extension of f in C(S**n) to point of infinity},
there exists $\overline{f}\in C(\overline{S}^{n},\overline{d}^{n})$
such that $f=\overline{f}|S^{n}$. The defining equality \ref{eq:temp-442}
therefore implies that 
\[
F_{1,\cdots.n}(f)=F_{1,\cdots.n}(\overline{f}|S^{n})=\overline{F}_{1,\cdots,n}\overline{f}
\]
\[
=\overline{F}'_{1,\cdots,n}\overline{f}=F'_{1,\cdots,n}(\overline{f}|S^{n})=F'_{1,\cdots.n}(f).
\]
By Lemma \ref{Lem. locally compact S, Ip=00003D>I iff Ip(f)->I(f) for each f in C(S)},
it follows that $F_{1,\cdots.n}=F'_{1,\cdots.n}$ as distributions,
where $n\geq1$ is arbitrary. Hence $F=F'$. Thus $\psi$ is an injection.
The lemma is proved.
\end{proof}
\begin{lem}
\textbf{\emph{\label{Lem. Identifying some distribution on Sbar^inf with distribution on S^inf}
(Identifying each distribution in $\overline{\Phi}_{DK}(\widehat{F}(Q,S))$
with a distribution on the path space $(S^{Q},d^{Q})$). }}Since $(\overline{S},\overline{d})$
is a compact metric space, Theorem \ref{Thm. Compact Daniell-Kolmogorov Extension}
yields the Compact Daniell-Kolmogorov Extension 
\[
\overline{\Phi}_{DK}:\widehat{F}(Q,\overline{S})\rightarrow\widehat{J}(\overline{S}^{Q},\overline{d}^{Q})
\]
\[
.
\]
Since $\widehat{F}(Q,S)\subset\widehat{F}(Q,\overline{S})$ according
to Lemma \ref{Lem. Identify family of f.j.d.s on S with family of f.j.d.'s on S_bar},
the image $\overline{\Phi}_{DK}(\widehat{F}(Q,S))$ is well defined
and is a subset of $\widehat{J}(\overline{S}^{Q},\overline{d}^{Q})$.
Define
\begin{equation}
\widehat{J}_{DK}(S^{Q},d^{Q})\equiv\overline{\Phi}_{DK}(\widehat{F}(Q,S))\subset\widehat{J}(\overline{S}^{Q},\overline{d}^{Q}).\label{eq:temp-53}
\end{equation}
Let $\overline{E}\in\widehat{J}_{DK}(S^{Q},d^{Q})$ be arbitrary.
In other words, $\overline{E}=\overline{\Phi}_{DK}(F)$ for some $F\in\widehat{F}(Q,S)$.
Then $S^{Q}$ is a full subset relative to the distribution $\overline{E}$.
Define
\begin{equation}
\varphi(\overline{E})\equiv E\equiv\overline{E}|C_{ub}(S^{Q},d^{Q}).\label{eq:temp-203-2}
\end{equation}
Then the following holds.

1. $\varphi(\overline{E})\equiv E\in\widehat{J}(S^{Q},d^{Q}).$ 

2. The coordinate function $U:Q\times(S^{Q},L,E)\rightarrow(S,d)$
is a r.f. with marginal distributions given by the family $F$, where
$(S^{Q},L,E)$ is the completion of $(S^{Q},C_{ub}(S^{Q},d^{Q}),E)$. 

3. The function 
\begin{equation}
\varphi:\widehat{J}_{DK}(S^{Q},d^{Q})\rightarrow\widehat{J}(S^{Q},d^{Q})\label{eq:temp-135-2}
\end{equation}
thus defined is an injection. 

Henceforth, we will identify $\overline{E}$ with $E\equiv\varphi(\overline{E})$.
In words, each distribution on the path space $(\overline{S}^{Q},\overline{d}^{Q})$
which is the image under the mapping $\overline{\Phi}_{DK}$ of some
family of f.j.d.'s with state space $(S,d)$ will be identified with
a distribution on the path space $(S^{Q},d^{Q})$. Thus 
\[
\widehat{J}_{DK}(S^{Q},d^{Q})\subset\widehat{J}(S^{Q},d^{Q})\cap\widehat{J}(\overline{S}^{Q},\overline{d}^{Q}).
\]
\end{lem}
\begin{proof}
1. Let $\overline{E}\in\widehat{J}_{DK}(S^{Q},d^{Q})$ be arbitrary.
Then, by the defining equality \ref{eq:temp-53}, there exists $F\in\widehat{F}(Q,S)$
such that $\overline{E}=\overline{\Phi}_{DK}(\overline{F})$, where
$\overline{F}\equiv\psi(F)$. Theorem \ref{Thm. Compact Daniell-Kolmogorov Extension},
applied to the compact metric space $(\overline{S},\overline{d})$
and the consistent family $\overline{F}$ of f.j.d.s with state space
$(\overline{S},\overline{d})$, says that the coordinate function
\[
\overline{U}:Q\times(\overline{S}^{Q},\overline{L},\overline{E})\rightarrow(\overline{S},\overline{d})
\]
is a r.f. with marginal distributions given by the family $\overline{F}$,
where $(\overline{S}^{Q},\overline{L},\overline{E})$ is the completion
of $(\overline{S}^{Q},C(\overline{S}^{Q},\overline{d}^{Q}),\overline{E})$. 

Note that each function $f\in C_{ub}(S^{Q},d^{Q})$ can be regarded
as a function on $\overline{S}^{Q}$ with $domain(f)\equiv S^{Q}\subset\overline{S}^{Q}$.
We will prove that $S^{Q}$ is a full set in $(\overline{S}^{Q},\overline{L},\overline{E})$,
and that $C_{ub}(S^{Q},d^{Q})\subset\overline{L}$. We will then show
that the restricted function $E\equiv\overline{E}|C_{ub}(S^{Q},d^{Q})$
is a distribution on $(S^{Q},d^{Q})$. 

2. To that end, let $n,m\geq1$ be arbitrary. Define the function
\[
h_{m}\equiv1\wedge(1+m-d(\cdot,x_{\circ}))_{+}\in C(S,d).
\]
Define the function $\overline{h}_{m}\in C(\overline{S},\overline{d})$
by $\overline{h}_{m}(y)\equiv h_{m}(y)$ or $\overline{h}_{m}(y)\equiv0$
according as $y\in S$ or $y=\Delta$. Then $h_{m}=\overline{h}_{m}|S$.
Moreover, for each point $y\in(S,d)$ we have $\overline{h}_{m}(y)=1$
if $m\geq1$ is sufficiently large. At the same time, $\overline{h}_{m}(\Delta)=0$
for each $m\geq1$. It follows that $\lim_{m\rightarrow\infty}\overline{h}_{m}(\overline{U}_{n})(y)=1$
or $0$ according as $\overline{U}_{n}(y)\in S$ or $\overline{U}_{n}(y)=\Delta$,
for each $y\in domain(\overline{U}_{n})$. In short, 
\[
\lim_{m\rightarrow\infty}\overline{h}_{m}(\overline{U}_{n})=1_{(\overline{U}(n)\in S)}
\]
on $domain(\overline{U}_{n})$. Moreover, 
\[
\overline{E}\overline{h}_{m}(\overline{U}_{n})=\overline{F}_{n}\overline{h}_{m}\equiv F_{n}(\overline{h}_{m}|S)=F_{n}h_{m}\uparrow1
\]
as $m\rightarrow\infty$, where the first equality is because the
r.f. $\overline{U}$ has marginal distributions given by the family
$\overline{F}$, where the second equality is by the defining formula
\ref{eq:temp-52} of the family $\overline{F}\equiv\psi(F)$, and
where the convergence is because $F_{n}$ is a distribution on the
locally compact metric space $(S,d)$. The Monotone Convergence Theorem
therefore implies that the indicator $1_{(\overline{U}(n)\in S)}$
is integrable on $(\overline{S}^{\infty},\overline{L},\overline{E})$,
with integral $1$. Thus $(\overline{U}_{n}\in S)$ is a full subset
of the probability space $(\overline{S}^{\infty},\overline{L},\overline{E})$,
where $n\geq1$ is arbitrary. Since, by the definition of the coordinate
function $\overline{U}$, we have
\[
S^{\infty}=\bigcap_{m=1}^{\infty}\{(x_{1},x_{2},\cdots)\in\overline{S}^{\infty}:x_{m}\in S\}
\]
\[
=\bigcap_{m=1}^{\infty}\{(x_{1},x_{2},\cdots)\in\overline{S}^{\infty}:\overline{U}_{m}(x_{1},x_{2},\cdots)\in S\}=\bigcap_{m=1}^{\infty}(\overline{U}_{m}\in S),
\]
it follows that $S^{\infty}$is the intersection of a sequence of
full subsets, and is itself a full subset of $(\overline{S}^{\infty},\overline{L},\overline{E})$. 

3. Now note that $U_{n}=\overline{U}_{n}$ on the full subset $S^{\infty}$
of $(\overline{S}^{\infty},\overline{L},\overline{E})$, where 
\[
U:Q\times S^{\infty}\rightarrow S
\]
is the coordinate function as in Definition \ref{Def. Path space, and coordinate function},
where $n\geq1$ is arbitrary. It follows that $U$ is a r.f. on $(\overline{S}^{\infty},\overline{L},\overline{E})$
with values in the metric space $(S,d)$. 

4. Next consider each $n\geq1$, $\overline{f}\in C_{ub}(\overline{S}^{n},\overline{d}^{n})$
and $f\in C_{ub}(S^{n},d^{n})$, such that $\overline{f}|S^{n}=f$.
Then
\begin{equation}
f(U_{1},\cdots,U_{n})\in C_{ub}(S^{Q},d^{Q}).\label{eq:temp-47}
\end{equation}
Separately, according to Proposition \ref{Prop. vector of measurable func (S)   is meas func (S)},
the function $(U_{1},\cdots,U_{n})$ on $(\overline{S}^{\infty},\overline{L},\overline{E})$
is a r.v. with values in $(S^{n},d^{n})$. Hence $f(U_{1},\cdots,U_{n})\in\overline{L}$
by Proposition \ref{Prop. X meas,  f unif continuous and bd on bd subsets =00003D> f(X) meas}.
Therefore
\begin{equation}
\overline{f}(\overline{U}_{1},\cdots,\overline{U}_{n})|S^{\infty}=f(U_{1},\cdots,U_{n})\in\overline{L}.\label{eq:temp-135}
\end{equation}
Combining,
\[
Ef(U_{1},\cdots,U_{n})=\overline{E}f(U_{1},\cdots,U_{n})=\overline{E}\overline{f}(\overline{U}_{1},\cdots,\overline{U}_{n})
\]
\begin{equation}
=\overline{F}_{1,\cdots,n}\overline{f}\equiv F_{1,\cdots,n}(\overline{f}|S^{n})=F_{1,\cdots,n}(f),\label{eq:temp-49-3}
\end{equation}
where the first equality is by the equality \ref{eq:temp-203-2} and
relation \ref{eq:temp-47}, where the second equality is from equality
\ref{eq:temp-135}, where the third equality is because the r.f. $\overline{U}:Q\times(\overline{S}^{Q},\overline{L},\overline{E})\rightarrow(\overline{S},\overline{d})$
has marginal distributions given by the family $\overline{F},$ and
where the fourth equality is by the defining equality \ref{eq:temp-52}
in Lemma \ref{Lem. Identify family of f.j.d.s on S with family of f.j.d.'s on S_bar}.

5. We proceed to prove that the function $E$ is a distribution on
the complete metric space $(S^{\infty},d^{\infty})$. To that end,
let $f\in C_{ub}(S^{\infty},d^{\infty})$ be arbitrary, with a modulus
of continuity $\delta_{f}$., and with $|f|\leq b$ for some $b>0$.
We will prove that $f\in\overline{L}$. Define the function $\overline{f}:\overline{S}^{\infty}\rightarrow R$
by $domain(\overline{f})\equiv S^{\infty}$ and by $\overline{f}(x)\equiv f(x)$
for each $x\in S^{\infty}$. 

Let $\varepsilon>0$ be arbitrary. Let $n\geq1$ be so large that
$2^{-n}<\delta_{f}(\varepsilon).$ Define $f_{n}\in C_{ub}(S^{n},d^{n})$
by 
\begin{equation}
f_{n}(x_{1},x_{2},\cdots,x_{n})\equiv f(x_{1},x_{2},\cdots,x_{n},x_{\circ},x_{\circ},\cdots)\label{eq:temp-49-2-1}
\end{equation}
for each $(x_{1},x_{2},\cdots,x_{n})\in S^{n}$. Then
\[
d^{\infty}((x_{1},x_{2},\cdots,x_{n},x_{\circ},x_{\circ},\cdots),(x_{1},x_{2},\cdots))\leq\sum_{i=n+1}^{\infty}2^{-i}=2^{-n}<\delta_{f}(\varepsilon)
\]
 for each $x=(x_{1},x_{2},\cdots)\in S^{\infty}$. Hence 
\[
|f_{n}(U_{1},\cdots,U_{n})-f(x_{1},x_{2},\cdots)|=|f_{n}(x_{1},x_{2},\cdots,x_{n})-f(x_{1},x_{2},\cdots)|
\]
\[
=|f(x_{1},x_{2},\cdots,x_{n},x_{\circ},x_{\circ},\cdots)-f(x_{1},x_{2},\cdots)|<\varepsilon
\]
for each $x=(x_{1},x_{2},\cdots)\in S^{\infty}$. In other words,
\begin{equation}
|f_{n}(U_{1},\cdots,U_{n})-f|<\varepsilon\label{eq:temp-50}
\end{equation}
on the full set $S^{\infty}$, where $n\geq1$ is arbitrary with $2^{-n}<\delta_{f}(\varepsilon).$
Thus
\begin{equation}
f_{n}(U_{1},\cdots,U_{n})\rightarrow f\label{eq:temp-184}
\end{equation}
in probability on $(\overline{S}^{Q},\overline{L},\overline{E})$.
At the same time, since $f_{n}\in C_{ub}(S^{n},d^{n}),$ we have $f_{n}(U_{1},\cdots,U_{n})\in\overline{L}$
according to relation \ref{eq:temp-135}, while $|f|\leq b$ and $|f_{n}|\leq b$
for each $n\geq1$. Therefore, in view of the convergence relation
\ref{eq:temp-184}, the Dominated Convergence Theorem implies that
$f\in\overline{L}$, and that 
\begin{equation}
Ef=\lim_{n\rightarrow\infty}Ef_{n}(U_{1},\cdots,U_{n})=\lim_{n\rightarrow\infty}\overline{E}f_{n}(\overline{U}_{1},\cdots,\overline{U}_{n})=\overline{E}f,\label{eq:temp-187}
\end{equation}
where the second equality is by applying equality \ref{eq:temp-49-3}
to $f_{n}$ for each $n\geq1$. Since $f\in C_{ub}(S^{\infty},d^{\infty})$
is arbitrary, we see that $C_{ub}(S^{\infty},d^{\infty})\subset\overline{L}$.

6. We will now verify that $(S^{\infty},C_{ub}(S^{\infty},d^{\infty}),E)$
is an integration space. First note that the space $C_{ub}(S^{\infty},d^{\infty})$
is linear, contains constants, and is closed to the operations of
absolute values and minimums. Linearity of the function $E$ follows
from that of $\overline{E}$, in view of equality \ref{eq:temp-187}.
Now suppose a sequence $(g_{i})_{i=0,1,2,\cdots}$ of functions in
$C_{ub}(S^{\infty},d^{\infty})$ is such that $g_{i}$ is non-negative
for each $i\geq1$ and such that $\sum_{i=1}^{\infty}Eg_{i}<Eg_{0}$.
Then, by equality \ref{eq:temp-187}, we have $\sum_{i=1}^{\infty}\overline{E}g_{i}<\overline{E}g_{0}$.
Hence, since $\overline{E}$ is an integration, there exists a point
$x\in\bigcap_{i=0}^{\infty}domain(g_{i})$ such that $\sum_{i=1}^{\infty}g_{i}(x)<g_{0}(x)$.
Thus the positivity condition in Definition \ref{Def. Integration Space}
has been verified for the function $E$. Now let $g\in C_{ub}(S^{\infty},d^{\infty})$
be arbitrary. Then $E(g\wedge k)=\overline{E}(g\wedge k)\rightarrow\overline{E}(g)=Eg$
as $k\rightarrow\infty$, where the convergence is because $\overline{E}$
is an integration. Similarly, $E(|g|\wedge k^{-1})\rightarrow0$ as
$k\rightarrow\infty$. Thus all the conditions in Definition \ref{Def. Integration Space}
have been verified for $(S^{\infty},C_{ub}(S^{\infty},d^{\infty}),E)$
to be an integration space. 

7. Since $d^{\infty}\leq1$ and $E1=\overline{E}1=1$, Assertion 2
of Lemma \ref{Lem. Distribution basics} implies that the integration
$E$ is a distribution on the complete metric space $(S^{\infty},d^{\infty})$.
In other words, $\varphi(\overline{E})\equiv E\in\widehat{J}(S^{Q},d^{Q})$.
We see that the function $\varphi$ is well defined. Moreover, let
$(S^{\infty},L,E)$ be the completion of the integration space $(S^{\infty},C_{ub}(S^{\infty},d^{\infty}),E)$.
Then equality \ref{eq:temp-49-3} implies that the coordinate function
$U$ is a r.f. with sample space $(S^{\infty},L,E)$ and with marginal
distributions given by the family $F$.

8. It remains to prove that $\varphi$ is an injection. To that end,
let a second distribution $\overline{E}'\in\widehat{J}_{DK}(S^{Q},d^{Q})$
be arbitrary. Suppose 
\[
E\equiv\varphi(\overline{E})=\varphi(\overline{E}')\equiv E'.
\]
Let $\overline{f}\in(\overline{S}^{Q},\overline{d}^{Q})$ be arbitrary.
Then $f\equiv\overline{f}|S^{\infty}\in C_{ub}(S^{\infty},d^{\infty})$.
Then equality \ref{eq:temp-187} yields
\begin{equation}
\overline{E}\overline{f}=\overline{E}f=Ef=E'f=\overline{E}'f=\overline{E}'\overline{f},\label{eq:temp-187-1}
\end{equation}
where the first and last equality are because $\overline{f}=f$ on
the full subset $S^{\infty}$relative to $\overline{E}$ and to $\overline{E}'$.
Thus $\overline{E}=\overline{E}'$ as distributions on the compact
metric space $(\overline{S}^{Q},\overline{d}^{Q})$. We conclude that
the function $\varphi$ is an injection. The lemma is proved.
\end{proof}
We are now ready to prove the Daniell-Kolmogorov Extension Theorem,
where the state space is required only to be locally compact.
\begin{thm}
\label{Thm. Daniell-Kolmogorov Extension}\textbf{\emph{ (Daniell-Kolmogorov
Extension, and its continuity).}} Suppose $(S,d)$ is locally compact,
not necessarily compact, with a binary approximation $\xi$. Recall
the set
\[
\widehat{J}_{DK}(S^{Q},d^{Q})\subset\widehat{J}(S^{Q},d^{Q})\cap\widehat{J}(\overline{S}^{Q},\overline{d}^{Q})
\]
of distributions, defined in Lemma \ref{Lem. Identifying some distribution on Sbar^inf with distribution on S^inf}.
Then the following holds.

1. \textbf{\emph{(Existence).}} There exists a function 
\[
\Phi_{DK}:\widehat{F}(Q,S)\rightarrow\widehat{J}_{DK}(S^{Q},d^{Q})
\]
such that, for each consistent family $F\in\widehat{F}(Q,S)$ of f.j.d.'s,
the distribution $E\equiv\Phi_{DK}(F)$ satisfies the conditions \emph{(i)}
the coordinate function 
\[
U:Q\times(S^{Q},L,E)\rightarrow S
\]
is a r.f., where $L$ is the completion of $C_{ub}(S^{Q},d^{Q})$
relative to the distribution $E$, and \emph{(ii) }the r.f. $U$ has
marginal distributions given by the family $F$. The function $\Phi_{DK}$
will be called the \index{Daniell-Kolmogorov Extension} \emph{Daniell-Kolmogorov
Extension}.

2. \textbf{\emph{(Continuity). }}Let $\overline{\xi}$ be the compactification
of the given binary approximation $\xi$, as constructed in Corollary
\ref{Cor. Compactification of binary approximation}. Thus $\overline{\xi}$
is a binary approximation of $(\overline{S},\overline{d})$ relative
to the fixed reference point $x_{\circ}\in S$. Since the metric space
$(\overline{S},\overline{d})$ is compact, the countable power $\overline{\xi}^{\infty}\equiv(B_{k})_{k=1,2,}$
of $\overline{\xi}$ is defined and is a binary approximation of $(\overline{S}^{\infty},\overline{d}^{\infty})$,
according to Definition \ref{Def. Countable product of binary approxximations; compact}.
Recall that the\emph{ set $\widehat{J}(\overline{S}^{\infty},\overline{d}^{\infty})$}
is then equipped with the distribution metric $\rho_{Dist,\overline{\xi}^{\infty}}$
defined relative to $\overline{\xi}^{\infty}$, according to Definition
\ref{Def. Distribution metric}, and that convergence of a sequence
of distributions on $(\overline{S}^{\infty},\overline{d}^{\infty})$
relative to the metric $\rho_{Dist,\overline{\xi}^{\infty}}$ is equivalent
to weak convergence. Write $\rho{}_{Dist,\overline{\xi}^{Q}}\equiv\rho_{Dist,\overline{\xi}^{\infty}}$.
Since $\widehat{J}_{DK}(S^{Q},d^{Q})$ is a subset of $\widehat{J}(\overline{S}^{\infty},\overline{d}^{\infty})$,
we have a metric subspace $(\widehat{J}_{DK}(S^{Q},d^{Q}),\rho{}_{Dist,\overline{\xi}^{Q}})$.
Then the \emph{Daniell-Kolmogorov Extension}
\[
\Phi_{DK}:(\widehat{F}(Q,S),\widehat{\rho}_{Marg,\xi,Q})\rightarrow(\widehat{J}_{DK}(S^{Q},d^{Q}),\rho{}_{Dist,\overline{\xi}^{Q}})
\]
is uniformly continuous, with a modulus of continuity $\delta_{DK}(\cdot,\left\Vert \xi\right\Vert )$
dependent only on the modulus of local compactness $\left\Vert \xi\right\Vert \equiv(|A_{k}|)_{k=1,2,}$
of the locally compact state space $(S,d)$.
\end{thm}
\begin{proof}
1. Apply Theorems \ref{Thm. Compact Daniell-Kolmogorov Extension}
and \ref{Thm. Continuity of the Compact DK Extension} to the compact
metric space $(\overline{S},\overline{d})$ to obtain the Compact
Daniell-Kolmogorov Extension
\[
\overline{\Phi}_{DK}:(\widehat{F}(Q,\overline{S}),\widehat{\rho}_{Marg,\xi,Q})\rightarrow(\widehat{J}(\overline{S}^{Q},\overline{d}^{Q}),\rho_{Dist,\overline{\xi}^{Q}}),
\]
which is uniformly continuous with modulus of continuity $\overline{\delta}_{DK}(\cdot,\left\Vert \overline{\xi}\right\Vert )$
dependent only on the modulus of local compactness $\left\Vert \overline{\xi}\right\Vert $
of the compact metric space $(\overline{S},\overline{d})$. 

2. By the defining equality \ref{eq:temp-53} in Lemma \ref{Lem. Identifying some distribution on Sbar^inf with distribution on S^inf},
we have 
\begin{equation}
\overline{\Phi}_{DK}(\widehat{F}(Q,S))\equiv\widehat{J}_{DK}(S^{Q},d^{Q}),\label{eq:temp-53-1}
\end{equation}
where $\widehat{F}(Q,S)$ is a subset of $\widehat{F}(Q,\overline{S})$.
Hence we can define the restricted mapping 
\begin{equation}
\Phi_{DK}\equiv\overline{\Phi}_{DK}|\widehat{F}(Q,S):(\widehat{F}(Q,S),\widehat{\rho}_{Marg,\xi,Q})\rightarrow(\widehat{J}_{DK}(S^{Q},d^{Q}),\rho_{Dist,\overline{\xi}^{Q}}),\label{eq:temp-49}
\end{equation}
which inherits the continuity and modulus of continuity $\overline{\delta}_{DK}(\cdot,\left\Vert \overline{\xi}\right\Vert )$
from $\overline{\Phi}_{DK}$. Thus the mapping $\Phi_{DK}$ is uniformly
continuous. According to Corollary \ref{Cor. Compactification of binary approximation},
$\left\Vert \overline{\xi}\right\Vert $ in turn depends only on the
modulus of local compactness $\left\Vert \xi\right\Vert $ of the
locally compact metric space $(S,d)$. Hence we can define
\[
\delta_{DK}(\cdot,\left\Vert \xi\right\Vert )\equiv\overline{\delta}_{DK}(\cdot,\left\Vert \overline{\xi}\right\Vert ),
\]
and $\Phi_{DK}$ has the modulus of continuity $\delta_{DK}(\cdot,\left\Vert \xi\right\Vert )$
which depends only on $\left\Vert \xi\right\Vert $.

3. Let $F\in\widehat{F}(Q,S)$ be arbitrary. Let $E\equiv\Phi_{DK}(F)$.
Then $E=\overline{\Phi}_{DK}(F)$ by the defining equality \ref{eq:temp-49}
of the function $\Phi_{DK}$. Hence Assertion 2 of Lemma \ref{Lem. Identifying some distribution on Sbar^inf with distribution on S^inf}
is applicable to $\overline{E}\equiv E=\overline{\Phi}_{DK}(F)$ and
$F\in\widehat{F}(Q,S)$, and says that the coordinate function $U:Q\times(S^{Q},L,E)\rightarrow(S,d)$
is a r.f. with marginal distributions given by the family $F$, where
$(S^{Q},L,E)$ is the completion of $(S^{Q},C_{ub}(S^{Q},d^{Q}),E)$.

The theorem is proved.
\end{proof}

\section{Daniell-Kolmogorov-Skorokhod Extension}

We use the notations as in the previous section. In particular, $Q\equiv\{t_{1},t_{2},\cdots\}$
denotes a countable parameter set, and, for simplicity of presentation,
and without loss of generality, we assume that $t_{n}=n$ for each
$n\geq1$. Thus 
\[
Q\equiv\{t_{1},t_{2},\cdots\}\equiv\{1,2,\cdots\}.
\]

For two consistent families $F$ and $F'$ of f.j.d.'s with the parameter
set $Q$ and the locally compact state space $(S,d)$, the Daniell-Kolmogorov
Extension in the previous section produces two corresponding distributions
$E$ and $E'$ on the path space $(S^{\infty},d^{\infty})$, such
that the families $F$ and $F'$ of f.j.d.'s are the marginal distributions
of the r.f.'s $U:Q\times(S^{Q},L,E)\rightarrow(S,d)$ and $U:Q\times(S^{Q},L',E')\rightarrow S$
receptively, even as the underlying coordinate function $U$ remains
the same. 

In contrast, Theorem 3.1.1 in \cite{Skorohod56} combines the Daniell-Kolmogorov
Extension with Skorokhod's Representation Theorem, presented as Theorem
\ref{Thm.Construction of Skorood  Representation} in the present
work, and produces (i) as the sample space, the fixed probability
space 
\[
(\Theta_{0},L_{0},I_{0})\equiv([0,1],L_{0},\int\cdot dx)
\]
based on the uniform distribution $\int\cdot dx$ on the unit interval
$[0,1]$, and (ii) for each $F\in\widehat{F}(Q,S)$, a r.f. $Z:Q\times(\Theta_{0},L_{0},I_{0})\rightarrow(S,d)$
with marginal distributions given by $F$. The sample space is fixed,
but two different families $F$ and $F'$ of f.j.d.'s result in two
different r.f.'s $Z$ and $Z'$. Theorem 3.1.1 in \cite{Skorohod56}
shows that the Daniell-Kolmogorov-Skorokhod Extension thus obtained
is continuous relative to weak convergence in $\widehat{F}(Q,S)$.
Because the r.f.'s produced can be regarded as r.v.'s on the same
probability space $(\Theta_{0},L_{0},I_{0})$ with values in the path
space $(S^{Q},d^{Q})$, we will have at out disposal the familiar
tools of making new r.v.'s, including the taking of continuous function
of given r.v.'s and the taking of limits in various senses. These
operations on such r.v.'s would be clumsy or impossible in terms of
distributions on the path space, . This will be clear as we go along. 

Note that, in Theorem \ref{Thm.Construction of Skorood  Representation}
of the present work, we recast the aforementioned Skorokhod's Representation
Theorem in terms of partitions of unity in the sense of Definition
\ref{Def. Partition of unity for locally compact (S,d)}; namely,
where Borel sets are used in \cite{Skorohod56}, we use continuous
basis functions with compact support. This will facilitate the subsequent
proof of metrical continuity of the following Daniell-Kolmogorov-Skorokhod
Extension, and the derivation of an accompanying modulus of continuity. 

Recall from Definition \ref{Def. Random Fields and Stochastic Processes}
that $\widehat{R}(Q\times\Theta_{0},S)$ denotes the set of r.f.'s
with parameter set $Q$, sample space $(\Theta_{0},L_{0},I_{0})$,
and state space $(S,d)$. We first identify each r.f. in $\widehat{R}(Q\times\Theta_{0},S)$
with a r.v. in $M(\Omega,S^{Q})$ 
\begin{defn}
\label{Def. Metric on random fields w/ countable parameters} \textbf{(Metric
space of r.f.'s with countable parameter set).} Suppose the state
space $(S,d)$ is locally compact, not necessarily compact. Let $(\Omega,L,E)$
be an arbitrary probability space. Recall that $\widehat{R}(Q\times\Omega,S)$
denotes the set of r.f.'s $Z:Q\times\Omega\rightarrow(S,d)$. By Lemma
\ref{Lem. Random field with countable parametr set is a r.v.}, each
r.f. $Z:Q\times\Omega\rightarrow(S,d)$ can be regarded as a r.v.
$Z:\Omega\rightarrow(S^{Q},d^{Q})$. Thus the set $\widehat{R}(Q\times\Omega,S)$
can be identified with a the set $M(\Omega,S^{Q})$ of r.v.'s with
values in the path space. Then the set $\widehat{R}(Q\times\Omega,S)$
inherits from $M(\Omega,S^{Q})$ the probability metric $\rho_{Prob}$,
defined in Definition \ref{Def. Metric rho_prob on space of r.v.'s}.
More precisely, define
\[
\widehat{\rho}_{Prob,Q}(Z,Z')\equiv\rho_{Prob}(Z,Z')\equiv E(1\wedge d^{Q}(Z,Z'))
\]
\emph{
\begin{equation}
=E(1\wedge\sum_{n=1}^{\infty}2^{-n}(1\wedge d(Z_{t(n)},Z'_{t(n)})))=E\sum_{n=1}^{\infty}2^{-n}(1\wedge d(Z_{t(n)},Z'_{t(n)})).\label{eq:temp-445-2-1}
\end{equation}
}for each $Z,Z'\in\widehat{R}(Q\times\Omega,S)=M(\Omega,S^{Q})$.
Note that $\widehat{\rho}_{Prob,Q}\leq1$.

`In view of the right-hand side of the defining equality \ref{eq:temp-445-2-1},
the metric $\widehat{\rho}_{Prob,Q}$ is determined by the enumeration
$(t_{1},t_{2},\cdots)$ of the countably infinite set $Q$. A different
enumeration would produce a different, albeit equivalent,  metric.
We omit the subscript $Q$ for the enumerated set only when it is
understood from context.

Note that equality \ref{eq:temp-445-2-1} implies that sequential
convergence of r.f.'s in $\widehat{R}(Q\times\Omega,S)$ relative
to the metric $\rho_{Prob}$ is equivalent to convergence in probability
and, therefore, to the weak convergence of the sequence, when the
r.f.'s are viewed as r.v.'s. $\square$
\end{defn}
\begin{thm}
\label{Thm. Compact Daniell-Kolmogorov-Skorohod} \textbf{\emph{(Compact
Daniell-Kolmogorov-Skorokhod Extension).}} Suppose $(S,d)$ is compact.
Let $\xi\equiv(A_{p})_{p=1,2,\cdots}$ be an arbitrary binary approximation
of state space $(S,d)$. Then there exists a function
\[
\overline{\Phi}_{DKS,\xi}:\widehat{F}(Q,S)\rightarrow\widehat{R}(Q\times\Theta_{0},S)
\]
such that, for each $F\in\widehat{F}(Q,S)$, the r.f. $Z\equiv\overline{\Phi}_{DKS,\xi}(F):Q\times\Theta_{0}\rightarrow S$
has marginal distributions given by the family $F$. The function
$\Phi_{DKS,\xi}$ constructed in the proof below will be called the
\emph{Daniell-Kolmogorov-Skorokhod Extension} \index{Daniell-Kolmogorov-Skorokhod Extension}
relative to the binary approximation $\xi$ of $(S,d)$. 
\end{thm}
\begin{proof}
1. Consider the Compact Daniell-Kolmogorov Extension
\[
\overline{\Phi}_{DK}:\widehat{F}(Q,S)\rightarrow\widehat{J}(S^{Q},d^{Q}),
\]
which maps each consistent family $F\in\widehat{F}(Q,S)$ of f.j.d.'s
to a distribution $E\equiv\overline{\Phi}_{DK}(F)$ on $(S^{Q},d^{Q})$,
such that (i) the coordinate function $U:Q\times(S^{Q},L,E)\rightarrow S$
is a r.f., where $L$ is the completion of $C(S^{Q},d^{Q})$ relative
to the distribution $E$, and (ii) $U$ has marginal distributions
given by the family $F$. 

2. Since the path space $(S,d)$ is compact by hypothesis, the countable
power $\xi^{\infty}\equiv(B_{k})_{k=1,2,}$ of $\xi$ is defined and
is a binary approximation of $(S^{Q},d^{Q})\equiv(S^{\infty},d^{\infty})$,
according to Definition \ref{Def. Countable product of binary approxximations; compact}.
Recall that the\emph{ set $\widehat{J}(S^{Q},d^{Q})$} is then equipped
with the distribution metric $\rho_{Dist,\xi^{\infty}}$ defined relative
to $\xi^{\infty}$, according to Definition \ref{Def. Distribution metric},
and that convergence of a sequence of distributions on $(S^{Q},d^{Q})$
relative to the metric $\rho_{Dist,\xi^{\infty}}$ is equivalent to
weak convergence. Write $\rho_{Dist,\xi^{Q}}\equiv\rho_{Dist,\xi^{\infty}}$. 

2. Recall from \ref{Def. Metric rho_prob on space of r.v.'s} that\emph{
}$M(\Theta_{0},S^{Q})$ denotes the set of r.v.'s $Z\equiv(Z_{t(1)},Z_{t(2)},\cdots)\equiv(Z_{1},Z_{2},\cdots)$
on $(\Theta_{0},L_{0},I_{0})$, with values in the compact path space
$(S^{Q},d^{Q})$. Theorem \ref{Thm.Construction of Skorood  Representation}
constructed the Skorokhod representation
\[
\Phi_{Sk,\xi^{\infty}}:\widehat{J}(S^{Q},d^{Q})\rightarrow M(\Theta_{0},S^{Q})
\]
such that, for each distribution $E\in\widehat{J}(S^{Q},d^{Q})$,
with $Z\equiv\Phi_{Sk,\xi^{\infty}}(E)$, we have 
\begin{equation}
E=I_{0,Z},\label{eq:temp-51}
\end{equation}
where $I_{0,Z}$ is the distribution induced on the compact metric
space $(S^{Q},d^{Q})$ by the r.v. $Z$, in the sense of \ref{Def. Distribution induced by r.v. w/ values in complete metric space}. 

3. We will now verify that the composite function
\begin{equation}
\overline{\Phi}_{DKS,\xi}\equiv\Phi_{Sk,\xi^{\infty}}\circ\overline{\Phi}_{DK}:\widehat{F}(Q,S)\rightarrow M(\Theta_{0},S^{Q})\label{eq:temp-188}
\end{equation}
has the desired properties. To that end, let the consistent family
$F\in\widehat{F}(Q,S)$ of f.j.d.'s be arbitrary. Let $Z\equiv\overline{\Phi}_{DKS,\xi}(F).$
Then $Z\equiv\Phi_{Sk,\xi^{\infty}}(E)$ for some $E\equiv\overline{\Phi}_{DK}(F)$.
We need only verify that the r.v. $Z=(Z_{1},Z_{2},\cdots)\in M(\Theta_{0},S^{Q})$,
when viewed as a r.f. $Z:Q\times(\Theta_{0},L_{0},I_{0})\rightarrow(S,d)$,
has marginal distributions given by the family $F$. 

4. To that end, let $n\geq1$ and $g\in C(S^{n},d^{n})$ be arbitrary.
Define the function $f\in C(S^{\infty},d^{\infty})$ by
\[
f(x_{1},x_{2},\cdots)\equiv g(x_{1},\cdots,x_{n})
\]
for each $(x_{1},x_{2},\cdots)\in S^{\infty}$. Then, for each $x\equiv(x_{1},x_{2},\cdots)\in S^{\infty}$,
we have, by the definition of the coordinate function $U$,
\begin{equation}
f(x)=f(x_{1},x_{2},\cdots)=f(U_{1}(x),U_{2}(x),\cdots)=g(U_{1},\cdots,U_{n})(x).\label{eq:temp-166}
\end{equation}
Therefore
\[
I_{0}g(Z_{1},\cdots,Z_{n})=I_{0}f(Z_{1},Z_{2}\cdots)=I_{0}f(Z)=I_{0,Z}f
\]
\begin{equation}
=Ef=Eg(U_{1},\cdots,U_{n})=F_{1,\cdots,n}g,\label{eq:temp-215}
\end{equation}
where the third equality is by the definition of the induced distribution
$I_{0,Z}$, the fourth follows from equality \ref{eq:temp-51}, the
fifth is by equality \ref{eq:temp-166}, and the last is by Condition
(ii) in Step 1. Since $n\geq1$ and $g\in C(S^{n},d^{n})$ are arbitrary,
we conclude that the r.f. $\overline{\Phi}_{DKS,\xi}(F)=Z$ has marginal
distributions given by the family $F$. The theorem is proved.
\end{proof}
\begin{thm}
\emph{\label{Thm. Continuity of Compact Daniell-Kolmogorov-Skorohod Extension}
}\textbf{\emph{(Continuity of Compact Daniell-Kolmogorov-Skorokhod
Extension).}}\emph{ }Use the same assumptions and notations as in
Theorem\emph{ \ref{Thm. Compact Daniell-Kolmogorov-Skorohod}}. In
particular, suppose the state space $(S,d)$ is compact. Recall that
the modulus of local compactness of $(S,d)$ corresponding to the
binary approximation $\xi\equiv(A_{p})_{p=1,2,\cdots}$ is defined
as the sequence
\[
\left\Vert \xi\right\Vert \equiv(|A_{p}|)_{p=1,2,\cdots}.
\]
of integers. Then the Compact Daniell-Kolmogorov-Skorokhod Extension
\begin{equation}
\overline{\Phi}_{DKS,\xi}:(\widehat{F}(Q,S),\widehat{\rho}_{Marg,\xi,Q})\rightarrow(\widehat{R}(Q\times\Theta_{0},S),\widehat{\rho}_{Prob,Q})\label{eq:temp-159}
\end{equation}
is uniformly continuous with a modulus of continuity $\overline{\delta}_{DKS}(\cdot,\left\Vert \xi\right\Vert )$
dependent only on $\left\Vert \xi\right\Vert $. The marginal metric
$\widehat{\rho}_{Marg,\xi,Q}$ and the probability metric $\widehat{\rho}_{Prob,Q}$
were introduced in Definitions \ref{Def. Marginal metric} and \ref{Def. Metric on random fields w/ countable parameters}
respectively.
\end{thm}
\begin{proof}
1. By the defining equality \ref{eq:temp-188} in Theorem \ref{Thm. Compact Daniell-Kolmogorov-Skorohod},
we have
\begin{equation}
\overline{\Phi}_{DKS,\xi}\equiv\Phi_{Sk,\xi^{\infty}}\circ\overline{\Phi}_{DK},\label{eq:temp-188-1}
\end{equation}
where the Compact Daniell-Kolmogorov Extension
\[
\overline{\Phi}_{DK}:(\widehat{F}(Q,S),\widehat{\rho}_{Marg,\xi,Q})\rightarrow(\widehat{J}(S^{Q},d^{Q}),\rho_{Dist,\xi^{Q}})
\]
is uniformly continuous according to Theorem \ref{Thm. Continuity of the Compact DK Extension},
with modulus of continuity 
\[
\overline{\delta}_{DK}(\cdot,\left\Vert \xi\right\Vert )
\]
dependent only on the modulus of local compactness $\left\Vert \xi\right\Vert \equiv(|A_{k}|)_{k=1,2,}$
of the compact metric space $(S,d)$. 

2. Separately, the metric space $(S,d)$ is compact by hypothesis.
Hence its countable power $(S^{Q},d^{Q})$ is compact. Moreover, the
countable power $\xi^{Q}$ is defined and is a binary approximation
of $(S^{Q},d^{Q})$. Moreover, since $d^{Q}\leq1$, the set $\widehat{J}(S^{Q},d^{Q})$
of distributions on $(S^{Q},d^{Q})$ is trivially tight, with the
modulus of tightness $\beta\equiv1$. Hence Theorem \ref{Thm. Continuity of Skorohod Representation}
is applicable to the metric space $(S^{Q},d^{Q})$ along with its
binary approximation $\xi^{Q}$ , and implies that the Skorokhod representation
\[
\Phi_{Sk,\xi^{\infty}}:(\widehat{J}(S^{Q},d^{Q}),\rho_{Dist,\xi^{Q}})\rightarrow(M(\Theta_{0},S),\rho_{Prob})
\]
is uniformly continuous, with a modulus of continuity $\delta_{Sk}(\cdot,\left\Vert \xi^{Q}\right\Vert ,1)$
depending only on $\left\Vert \xi^{Q}\right\Vert $. Equivalently,
\[
\Phi_{Sk,\xi^{\infty}}:(\widehat{J}(S^{Q},d^{Q}),\rho_{Dist,\xi^{Q}})\rightarrow(\widehat{R}(Q\times\Theta_{0},S),\rho_{Prob})
\]
is uniformly continuous, with a modulus of continuity $\delta_{Sk}(\cdot,\left\Vert \xi^{Q}\right\Vert ,1)$.

3. Combining, we see that the composite function $\overline{\Phi}_{DKS,\xi}$
in equality \ref{eq:temp-159} is uniformly continuous, with a modulus
of continuity given by the composite operation
\[
\overline{\delta}_{DKS}(\cdot,\left\Vert \xi\right\Vert )\equiv\overline{\delta}_{DK}(\delta_{Sk}(\cdot,\left\Vert \xi^{Q}\right\Vert ,1),\left\Vert \xi\right\Vert ),
\]
where we observe that the modulus of local compactness $\left\Vert \xi^{Q}\right\Vert $of
the countable power $(S^{Q},d^{Q})$ is determined by the modulus
of local compactness $\left\Vert \xi\right\Vert $ of the compact
metric space $(S,d)$, according to Lemma \ref{Lem. Countable Product of binary approximation is a binary approxt; compacts}.
The theorem is proved. 
\end{proof}
Now the Daniell-Kolmogorov-Skororkhod Extension Theorem, where the
state space is required only to be locally compact.
\begin{thm}
\label{Thm. DKS Extension, construction and continuity}\textbf{\emph{
(Daniell-Kolmogorov-Skorokhod Extension, and its continuity).}} Suppose
$(S,d)$ is locally compact, not necessarily compact, with a binary
approximation $\xi$. Then the following holds.

1. \textbf{\emph{(Existence).}} There exists a function 
\[
\Phi_{DKS,\xi}:\widehat{F}(Q,S)\rightarrow\widehat{R}(Q\times\Theta_{0},S)
\]
such that, for each $F\in\widehat{F}(Q,S)$, the r.f. $Z\equiv\Phi_{DKS,\xi}(F):Q\times\Theta_{0}\rightarrow S$
has marginal distributions given by the family $F$. The function
$\Phi_{DKS,\xi}$ will be called the Daniell\emph{-}Kolmogorov-Skorokhod
Extension relative to the binary approximation $\xi$ of $(S,d)$. 

2. \textbf{\emph{(Continuity). }}\emph{The Daniell-Kolmogorov-Skorokhod
Extension}
\[
\Phi_{DKS,\xi}:(\widehat{F}(Q,S),\widehat{\rho}_{Marg,\xi,Q})\rightarrow(\widehat{R}(Q\times\Theta_{0},S),\widehat{\rho}_{Prob,Q})
\]
is uniformly continuous, with a modulus of continuity $\delta_{DKS}(\cdot,\left\Vert \xi\right\Vert )$
dependent only on the modulus of local compactness $\left\Vert \xi\right\Vert \equiv(|A_{k}|)_{k=1,2,}$
of the locally compact state space $(S,d)$.
\end{thm}
\begin{proof}
1. Let $\overline{\xi}$ be the compactification of the given binary
approximation $\xi$, as constructed in Corollary \ref{Cor. Compactification of binary approximation}.
Thus $\overline{\xi}$ is a binary approximation of $(\overline{S},\overline{d})$
relative to the fixed reference point $x_{\circ}\in S$. Since the
metric space $(\overline{S},\overline{d})$ is compact, the countable
power $\overline{\xi}^{\infty}$ of $\overline{\xi}$ is defined and
is a binary approximation of $(\overline{S}^{\infty},\overline{d}^{\infty})$,
according to Definition \ref{Def. Countable product of binary approxximations; compact}.
Apply Theorems \ref{Thm. Compact Daniell-Kolmogorov-Skorohod} and
\ref{Thm. Continuity of Compact Daniell-Kolmogorov-Skorohod Extension}
to the compact metric space $(\overline{S},\overline{d})$ to obtain
the Compact Daniell-Kolmogorov-Skorokhod Extension
\begin{equation}
\overline{\Phi}_{DKS,\overline{\xi}}:(\widehat{F}(Q,\overline{S}),\widehat{\rho}_{Marg,\xi,Q})\rightarrow(\widehat{R}(Q\times\Theta_{0},\overline{S}),\rho_{Prob}),\label{eq:temp-159-1}
\end{equation}
which is uniformly continuous with a modulus of continuity $\overline{\delta}_{DKS}(\cdot,\left\Vert \overline{\xi}\right\Vert )$
dependent only on $\left\Vert \overline{\xi}\right\Vert $. Specifically,
\[
\overline{\Phi}_{DKS,\overline{\xi}}\equiv\Phi_{Sk,\overline{\xi}^{\infty}}\circ\overline{\Phi}_{DK},
\]
where 
\[
\overline{\Phi}_{DK}:(\widehat{F}(Q,\overline{S}),\widehat{\rho}_{Marg,\xi,Q})\rightarrow(\widehat{J}(\overline{S}^{Q},\overline{d}^{Q}),\rho_{Dist,\overline{\xi}^{Q}})
\]
and 
\[
\Phi_{Sk,\overline{\xi}^{Q}}:\widehat{J}(\overline{S}^{Q},\overline{d}^{Q})\rightarrow(M(\Theta_{0},\overline{S}^{Q}),\rho_{Prob})
\]

2. Since $\widehat{F}(Q,S)$ is a subset of $\widehat{F}(Q,\overline{S})$,
we can define the restricted mapping 
\begin{equation}
\Phi_{DKS,\xi}\equiv\overline{\Phi}_{DKS,\overline{\xi}}|\widehat{F}(Q,S):(\widehat{F}(Q,S),\widehat{\rho}_{Marg,\xi,Q})\rightarrow(\widehat{R}(Q\times\Theta_{0},\overline{S}),\rho_{Prob}),\label{eq:temp-49-1}
\end{equation}
which inherits the continuity and modulus of continuity $\overline{\delta}_{DKS}(\cdot,\left\Vert \overline{\xi}\right\Vert )$
from $\overline{\Phi}_{DKS,\overline{\xi}}$. Thus the mapping $\Phi_{DKS,\xi}$
is uniformly continuous. According to Corollary \ref{Cor. Compactification of binary approximation},
$\left\Vert \overline{\xi}\right\Vert $ in turn depends only on the
modulus of local compactness $\left\Vert \xi\right\Vert $ of the
locally compact metric space $(S,d)$. Hence we can define
\[
\delta_{DKS}(\cdot,\left\Vert \xi\right\Vert )\equiv\overline{\delta}_{DKS}(\cdot,\left\Vert \overline{\xi}\right\Vert ).
\]
Then $\Phi_{DKS,\xi}$ has the modulus of continuity $\delta_{DKS}(\cdot,\left\Vert \xi\right\Vert )$
which depends only on $\left\Vert \xi\right\Vert $.

3. Let $F\in\widehat{F}(Q,S)$ be arbitrary. Write $Z\equiv\Phi_{DKS,\xi}(F)$.
It remains to prove that $Z\in\widehat{R}(Q\times\Theta_{0},S)$ and
that it has marginal distributions given by the family $F$. To that
end, note that, since $\Phi_{DKS,\xi}\equiv\overline{\Phi}_{DKS,\overline{\xi}}|\widehat{F}(Q,S)$,
we have $Z=\overline{\Phi}_{DKS,\overline{\xi}}(F)$. Moreover, write
$E\equiv\overline{\Phi}_{DK}(F)\in\widehat{J}(\overline{S}^{Q},\overline{d}^{Q})$.
Then, by the defining equality \ref{eq:temp-188} for the function
$\overline{\Phi}_{DKS,\overline{\xi}}$, we have 
\[
Z\equiv\Phi_{DKS,\xi}(F)=\Phi_{Sk,\overline{\xi}^{Q}}(E)\in M(\Theta_{0},\overline{S}^{Q}).
\]
Furthermore, Theorem \ref{Thm.Construction of Skorood  Representation}
defined and constructed the Skorokhod representation $\Phi_{Sk,\overline{\xi}^{\infty}}$,
such that the r.v. 
\[
Z\equiv\Phi_{Sk,\overline{\xi}}(E):(\Theta_{0},L_{0},I_{0})\rightarrow(\overline{S}^{Q},\overline{d}^{Q})
\]
induces the distribution $E$ on the metric space $(\overline{S}^{Q},\overline{d}^{Q})$.
In other words, $Ef=I_{0}f(Z)$ for each $f\in C(\overline{S}^{Q},\overline{d}^{Q})$.
In particular,
\begin{equation}
Eg(U_{1},\cdots,U_{n})=I_{0}g(Z_{1},\cdots,Z_{n})\label{eq:temp-167}
\end{equation}
 for each $g\in C_{ub}(S^{n},d^{n})$, for each $n\geq1$. Here $U:Q\times S^{Q}\rightarrow S$
is the coordinate function. 

On the other hand, Theorem \ref{Thm. Daniell-Kolmogorov Extension}
says that $\Phi_{DK}\equiv\overline{\Phi}_{DK}|\widehat{F}(Q,S)$,
whence 
\[
E=\overline{\Phi}_{DK}(F)=\Phi_{DK}(F)\in\widehat{J}(S^{Q},d^{Q}),
\]
and that (i) the coordinate function 
\[
U:Q\times(S^{Q},L,E)\rightarrow(S,d)
\]
is a r.f., where $L$ is the completion of $C_{ub}(S^{Q},d^{Q})$
relative to the distribution $E$, and (ii)\emph{ }the r.f. $U$ has
marginal distributions given by the family $F$. Hence 
\begin{equation}
Eg(U_{1},\cdots,U_{n})=F_{1,\cdots,n}g\label{eq:temp-191}
\end{equation}
 for each $g\in C_{ub}(S^{n},d^{n})$, for each $n\geq1$. Combining
equalities \ref{eq:temp-167} and \ref{eq:temp-191}, we obtain
\[
I_{0}g(Z_{1},\cdots,Z_{n})=F_{1,\cdots,n}g
\]
for each $g\in C_{ub}(S^{n},d^{n})$, for each $n\geq1$. Summing
up, we conclude that 
\[
Z:Q\times(\Theta_{0},L_{0},I_{0})\rightarrow(S,d)
\]
is a r.f. with marginal distributions given by the consistent family
$F$ of f.j.d.'s. The theorem is proved.
\end{proof}
As a corollary, we prove Skorokhod's sequential continuity theorem,
essentially Theorem 3.1.1 in \cite{Skorohod56}.
\begin{thm}
\textbf{\emph{\label{Thm:Continuuity of DKS extension relative to a.u. convergence}
(Sequential continuity of the Daniell-Komogorov-Skorokhod Extension).}}
Let $F^{(0)},F^{(1)},F^{(2)},\cdots$ be an arbitrary sequence in
$\widehat{F}(Q,S)$ such that 
\begin{equation}
\widehat{\rho}_{Marg,\xi,Q}(F^{(p)},F^{(0)})\rightarrow0.\label{eq:temp-203-3}
\end{equation}
For each $p\geq0$, write 
\[
Z^{(p)}\equiv\Phi_{DKS,\xi}(F^{(p)})\in(\widehat{R}(Q\times\Theta_{0},S),\rho_{Prob})=M(\Theta_{0},S^{Q})
\]
Then 
\[
Z^{(p)}\rightarrow Z^{(0)}\qquad\mathrm{a.u.}
\]
as r.v.'s on $(\Theta_{0},L_{0},I_{0})$ with values in the path space
$(S^{Q},d^{Q})$.
\end{thm}
\begin{proof}
Let $p\geq0$ be arbitrary. Write $Z^{(p)}\equiv\Phi_{DKS,\xi}(F^{(p)})$
and write $\overline{E}^{(p)}\equiv\overline{\Phi}_{DK}(F^{(p)})\in\widehat{J}(\overline{S}^{Q},\overline{d}^{Q})$.
Then
\begin{equation}
Z^{(p)}\equiv\overline{\Phi}_{DKS,\overline{\xi}}(F^{(p)})\equiv\Phi_{Sk,\overline{\xi}^{\infty}}(\overline{\Phi}_{DK}(F^{(p)}))\equiv\Phi_{Sk,\overline{\xi}^{\infty}}(\overline{E}^{(p)}).\label{eq:temp-49-1-2}
\end{equation}
Since 
\[
\overline{\Phi}_{DK}:(\widehat{F}(Q,S),\widehat{\rho}_{Marg,\xi,Q})\rightarrow(\widehat{J}(\overline{S}^{Q},\overline{d}^{Q}),\rho_{Dist,\overline{\xi}^{Q}})
\]
is uniformly continuous, convergence relation \ref{eq:temp-203-3}
implies $\rho_{Dist,\overline{\xi}^{Q}}(\overline{E}^{(p)},\overline{E}^{(0)})\rightarrow0$.
At the same time, the metric space $(\overline{S}^{Q},\overline{d}^{Q})$
is compact. Hence Theorem \ref{Thm. Skorood's Continuity in terms off weak convergence and a.u. convergence}
is applicable, and implies that 
\[
\Phi_{Sk,\overline{\xi}^{Q}}(\overline{E}^{(p)})\rightarrow\Phi_{Sk,\overline{\xi}^{Q}}(\overline{E}^{(0)})\qquad\mathrm{a.u.}
\]
on $(\Theta_{0},L_{0},I_{0})$. In view of equality \ref{eq:temp-49-1-2},
this can be rewritten as 
\[
Z^{(p)}\rightarrow Z^{(0)}\qquad\mathrm{a.u.}
\]
as r.v.'s on $(\Theta_{0},L_{0},I_{0})$ with values in $(\overline{S}^{Q},\overline{d}^{Q})$.
In other words,
\begin{equation}
\overline{d}^{\infty}(Z^{(p)},Z^{(0)})\rightarrow0\quad\mathrm{a.u}.\label{eq:temp-178}
\end{equation}

5. We proceed to show that 
\begin{equation}
d^{\infty}(Z^{(p)},Z^{(0)})\rightarrow0\quad\mathrm{a.u}.\label{eq:temp-178-1}
\end{equation}
To that end, let $n\geq1$ be arbitrary. Then there exists $b>0$
so large that 
\[
I_{0}B^{c}<2^{-n}
\]
where 
\begin{equation}
B\equiv(\bigvee_{k=1}^{n}d(Z_{t(k)}^{(0)},x_{\circ})\leq b).\label{eq:temp-212}
\end{equation}
Let $h\equiv n\vee[\log_{2}b]_{1}$. In view of the a.u. convergence
\ref{eq:temp-178}, there exist $m\geq h$ and a $\mathrm{measurable}$
subset $A$ of $(\Theta_{0},L_{0},I_{0})$ with 
\[
I_{0}A^{c}<2^{-n}
\]
such that 
\begin{equation}
1_{A}\sum_{k=1}^{\infty}2^{-k}\overline{d}(Z_{t(k)}^{(p)},Z{}_{t(k)}^{(0)})\leq2^{-2h-1}|A_{h}|^{-2}\label{eq:temp-211}
\end{equation}
for each $p\geq m$. 

6. Now consider each $\theta\in AB$ and each $k=1,\cdots,n$. Then,
by inequality \ref{eq:temp-211}, we have
\begin{equation}
\overline{d}(Z_{t(k)}^{(p)}(\theta),Z{}_{t(k)}^{(0)}(\theta))\leq2^{k}2^{-2h-1}|A_{h}|^{-2}\leq2^{-h-1}|A_{h}|^{-2},\label{eq:temp-213}
\end{equation}
while equality \ref{eq:temp-212} yields 
\begin{equation}
d(Z_{t(k)}^{(0)}(\theta),x_{\circ})\leq b<2^{h}.\label{eq:temp-214}
\end{equation}
According to Assertion (ii) of Theorem \ref{Thmf. Construction of one-point compactificaton from separant},
regarding the one-point compactification $(\overline{S},\overline{d})$
of $(S,d)$ relative to the binary approximations $\xi\equiv(A_{p})_{p=1,2,\cdots}$,
if $y\in(d(\cdot,x_{\circ})\leq2^{h})$ and $z\in S$ are such that
\[
\overline{d}(y,z)<2^{-h-1}|A_{h}|^{-2},
\]
then 
\[
d(y,z)<2^{-h+2}.
\]
Hence inequalities \ref{eq:temp-213} and \ref{eq:temp-214} together
yield
\[
d(Z_{t(k)}^{(p)}(\theta),Z{}_{t(k)}^{(0)}(\theta))<2^{-h+2}\leq2^{-n+2},
\]
where $\theta\in AB$, and $k=1,\cdots,n$ are arbitrary. Consequently,
recalling the notation $\widehat{d}\equiv1\wedge d$, we obtain
\[
1_{AB}\sum_{k=1}^{\infty}2^{-k}\widehat{d}(Z_{t(k)}^{(p)},Z{}_{t(k)}^{(0)})\leq1_{AB}\sum_{k=1}^{n}2^{-k}\widehat{d}(Z_{t(k)}^{(p)},Z{}_{t(k)}^{(0)})+\sum_{k=n+1}^{\infty}2^{-k}
\]
\[
\leq1_{AB}\sum_{k=1}^{n}2^{-k}2^{-n+2}+\sum_{k=n+1}^{\infty}2^{-k}
\]
\[
<2^{-n+2}+2^{-n}<2^{-n+1},
\]
where $2^{-n+1}$ and $I_{0}(AB)^{c}<2^{-n+1}$ are arbitrarily small.
Hence, by Proposition \ref{Prop. Alternative definition of a.u. convergence},
we have 
\[
\sum_{k=1}^{\infty}2^{-k}\widehat{d}(Z_{t(k)}^{(p)},Z{}_{t(k)}^{(0)})\rightarrow0\quad a.u.
\]
Equivalently, 
\begin{equation}
d^{\infty}(Z^{(p)},Z^{(0)})\rightarrow0\quad\mathrm{a.u}.\label{eq:temp-178-1-1}
\end{equation}
In other words, $Z^{(p)}\rightarrow Z^{(0)}$ a.u. in $M(\Theta_{0},S^{n})$,
as alleged.
\end{proof}

\chapter{Measurable Random Field}

\section{Measurable R.F.'s which are Continuous in Probability}

In this chapter, let $(S,d)$ be a locally compact metric space, not
necessarily a linear space or ordered. Let $(Q,d_{Q})$ be a compact
metric space endowed with an arbitrary, but fixed, integration. 

Consider each consistent family $F$ of f.j.d.'s with parameter space
$Q$ and state space $S$ which is continuous in probability, we will
construct a $\mathrm{measurable}$ r.f. $X:Q\times\Omega\rightarrow S$
which extends $F$, in a sense to be defined presently. We will also
prove that the construction is metrically continuous.

In the special case where $S\equiv[-\infty,\infty]$ and where $Q$
is a subinterval of $[-\infty,\infty]$, where the symbol $[-\infty,\infty]$
stands for a properly defined and metrized two-point compactification
of the real line, the main theorem in Section III.4 of \cite{Neveu 65},
gives a classical construction. The construction in \cite{Neveu 65}
uses a sequence of step processes on half open intervals, and then
uses the limit supremum of this sequence as the desired $\mathrm{measurable}$
process $X$. Existence of a limit supremum is, however, by invoking
principle of infinite search, and not constructive. 

\cite{Potthoff09} gives a constructive proof of existence in t the
case where $S\equiv[-\infty,\infty]$ and where $(Q,d_{Q})$ is a
metric space, by using linear combinations, with stochastic coefficients,
of certain deterministic basis functions, as successive $L_{1}-$approximations
to the desired $\mathrm{measurable}$ random field $X,$ obviating
the use of any limit supremum. These deterministic basis functions
are continuous on $(Q,d_{Q})$ with values in the state space $[0,1]$,
and are from a partition of unity of $(Q,d_{Q})$.

In the general case, where neither a linear structure nor an ordering
is available on the state space $(S,d)$, the aforementioned limit
supremum or linear combinations of basis functions would not be available.
We will go around these difficulties by replacing the linear combinations
with \index{stochastic interpolation}\emph{stochastic interpolations}
of the basis functions, essentially by continuously varying the probability
weighting of basis functions. This method of construction of $\mathrm{measurable}$
processes, and the subsequent theorem of metrical continuity of the
construction, in epsilon-delta terms, seem hitherto unknown.

In the rest of this section, we will make the last paragraph precise.
\begin{defn}
\label{Def. Specification of state space and parameter space} \textbf{(Specification
of locally compact state space and compact parameter space, and their
binary approximations).} In this section, let $(S,d)$ be a locally
compact metric space, with a binary approximation $\xi\equiv(A_{n})_{n=1,2,\cdots}$
relative to some arbitrary, but fixed, reference point $x_{\circ}\in S$. 

Let $(Q,d_{Q})$ be a compact metric space, with $d_{Q}\leq1$, and
with a binary approximation $\xi_{Q}\equiv(B_{n})_{n=1,2,\cdots}$
relative to some arbitrary, but fixed, reference point $q_{\circ}\in Q$.
Let $I$ be an arbitrary but fixed, distribution on $(Q,d_{Q})$,
and let $(Q,\Lambda,I)$ denote the probability space which is the
completion of $(Q,C(Q,d_{Q}),I)$. This distribution provides measurable
sets and $\mathrm{measurable}$ functions, thereby facilitating the
definition of measurability. It is otherwise unimportant, and will
be called a reference distribution.

The assumption of compactness of $(Q,d_{Q})$ simplifies presentation.
The generalization of the results to a locally compact parameter space
$(Q,d_{Q})$ is easy, by considering each member in a sequence $(Q_{i})_{i=1,2,\cdots}$of
compact and integrable subsets which forms an $I$-basis of $Q$.
This generalization is straightforward and left to the reader.
\end{defn}
$\square$
\begin{defn}
\label{Def. Metric space of measurable random fields} \textbf{(Metric
space of measurable}\emph{ }\textbf{r.f.'s)} Let $(\Omega,L,E)$ be
an arbitrary probability space. Recall from Definition \ref{Def. Random Fields and Stochastic Processes}
the space $\widehat{R}(Q\times\Omega,S)$ of r.f.'s $X:Q\times(\Omega,L,E)\rightarrow(S,d)$,
with sample space $(\Omega,L,E)$, the compact parameter space $Q$,
and the locally compact state space $(S,d)$. Recall from Definition
\ref{Def. Metric rho_prob on space of r.v.'s} the metric space $(M(Q\times\Omega,S),\rho_{Prob})$
of r.v.'s on the product probability space $(Q\times\Omega,\Lambda\otimes L,I\otimes E)$\emph{
}with values in the state space\emph{ $(S,d)$.} We will say that
a r.f. $X\in\widehat{R}(Q\times\Omega,S)$ is\emph{ measurable }\index{measurable r.f.}
if \emph{$X\in M(Q\times\Omega,S)$. }We will write 
\[
\widehat{R}_{Meas}(Q\times\Omega,S)\equiv\widehat{R}(Q\times\Omega,S)\cap M(Q\times\Omega,S)
\]
for the space of $\mathrm{measurable}$ r.f.'s. Note that $\widehat{R}_{Meas}(Q\times\Omega,S)$
then inherits the probability metric $\rho_{Prob}$ on $M(Q\times\Omega,S)$,
which is defined, according to Definition \ref{Def. Metric rho_prob on space of r.v.'s},
by

\emph{
\begin{equation}
\rho_{Prob}(X,Y)\equiv I\otimes E(1\wedge d(X,Y))\label{eq:temp-445-2-2}
\end{equation}
}for each\emph{ $X,Y\in M(Q\times\Omega,S)$. }
\end{defn}
$\square$ 
\begin{defn}
\label{Def. Metric on measurable r.f. continuous in probability}
\textbf{(Two metrics on the space of }$\mathbf{measurable}$\textbf{
r.f.'s which are continuous in probability).} Recall from Definition
\ref{Def. continuity in prob, continuity a.u., and a.u. continuity}
the set $\widehat{R}_{Cp}(Q\times\Omega,S)$ of r.f.'s which are continuous
in probability. Let 
\[
\widehat{R}_{Meas,Cp}(Q\times\Omega,S)\equiv\widehat{R}_{Meas}(Q\times\Omega,S)\cap\widehat{R}_{Cp}(Q\times\Omega,S)
\]
be the subset of the metric space $(\widehat{R}_{Meas}(Q\times\Omega,S),\rho_{Prob})$
whose members are continuous in probability. As a subset, it inherits
the probability metric $\rho_{Prob}$ from the latter. Define a second
metric $\rho_{Sup,Prob}$ on this set $\widehat{R}_{Meas,Cp}(Q\times\Omega,S)$
by
\begin{equation}
\rho_{Sup,Prob}(X,Y)\equiv\sup_{t\in Q}E(1\wedge(X_{t},Y_{t}))\label{eq:temp-227}
\end{equation}
for each $X,Y\in\widehat{R}_{Meas,Cp}(Q\times\Omega,S)$. Note in
the above definition that $E(1\wedge(X_{t},Y_{t}))$ is a continuous
function on the compact metric space $(Q,d_{Q})$, on account of continuity
in probability, whence the supremum exists. Note also that defining
formulas \ref{eq:temp-445-2-2} and \ref{eq:temp-227} implies that
$\rho_{Prob}\leq\rho_{Sup,Prob}$. In words, $\rho_{Sup,Prob}$ is
a stronger metric than $\rho_{Prob}$ on the space of $\mathrm{measurable}$
r.f. which are continuous in probability.
\end{defn}
$\square$
\begin{defn}
\label{Def. Specification of  Q_inf, and partion of unity of Q} \textbf{(Specification
of a countable dense subset of the parameter space, and a partition
of unity of Q). }By Definition \ref{Def. Specification of state space and parameter space},
$\xi_{Q}\equiv(B_{n})_{n=1,2,\cdots}$ is an arbitrary, but fixed,
binary approximation of the compact metric space $(Q,d_{Q})$ relative
to the reference point $q_{\circ}\in Q$. Thus $B_{1}\subset B_{2}\subset\cdots$
is a sequence of metrically discrete and enumerated finite subsets
of $Q$, with $B_{n}\equiv\{q_{n,1},\cdots,q_{n,\gamma(n)}\}$ for
each $n\geq1$. 

1. Define the set
\begin{equation}
Q_{\infty}\equiv\{t_{1},t_{2},\cdots\}\equiv\bigcup_{n=1}^{\infty}B_{n}.\label{eq:temp-238-1-1}
\end{equation}
Note that, by assumption, $d_{Q}\leq1$. Hence, for each $n\geq1$,
we have, by Definition \ref{Def. Binary approximationt and Modulus of local compactness}
of a binary approximation, 
\begin{equation}
Q=(d_{Q}(\cdot,q_{\circ})\leq2^{n})\subset\bigcup_{q\in B(n)}(d_{Q}(\cdot,q)\leq2^{-n}).\label{eq:temp-505-2}
\end{equation}
Hence $Q_{\infty}\equiv\bigcup_{n=1}^{\infty}B_{n}$ is a metrically
discrete, countably infinite, and dense subset of $(Q,d_{Q})$. Moreover,
we can fix an enumeration of $Q_{\infty}$ in such a manner that 
\[
\{t_{1},t_{2},\cdots,t_{\gamma(n)}\}=B_{n}\equiv\{q_{n,1},\cdots,q_{n,\gamma(n)}\},
\]
where $\gamma_{n}\equiv|B_{n}|$, for each $n\geq1$. 

2. Let 
\[
\pi_{Q}\equiv(\{\lambda_{n,q}:q\in B_{n}\})_{n=1,2,\cdots}
\]
be the partition of unity\index{partition of unity} of $(Q,d_{Q})$
determined by $\xi_{Q}$. Let $n\geq1$ be arbitrary. Then, for each
$q\in B_{n}$, the basis function $\lambda_{n,q}\in C(Q,d_{Q})$ has
values in $[0,1]$ and has support $(d_{Q}(\cdot,q)\leq2^{-n+1})$.
Moreover 
\begin{equation}
Q\subset\bigcup_{q\in B(n)}(d_{Q}(\cdot,q)\leq2^{-n})\subset(\sum_{q\in B(n)}\lambda_{n,q}=1),\label{eq:temp-130}
\end{equation}
where the second inclusion is according to Proposition \ref{Prop. Properties of parittion of unity-1}.
Define the auxiliary continuous functions $\lambda_{n,0}^{+}\equiv0$,
and 
\[
\lambda_{n,k}^{+}\equiv\sum_{i=1}^{k}\lambda_{n,q(n,i)}
\]
for each $k=1,\cdots,\gamma_{n}$. Then
\[
0\equiv\lambda_{n,0}^{+}\leq\lambda_{n,1}^{+}\leq\cdots\leq\lambda_{n,\gamma(n)}^{+}=1.
\]
$\square$
\end{defn}
In the following, recall some miscellaneous short-hand notations.
For an arbitrary integrable set $A$ in a complete integration space
$(\Omega,L,E)$, we write $EA,$ $E(A)$, and $E1_{A}$ interchangeably.
Thus, if $(\Omega,L,E)$ is a probability space, then $EA\equiv P(A)$
is the probability of $A$. Recall also that $[\cdot]_{1}$ is the
operation which assigns to each $a\in R$ an integer $[a]_{1}\in(a,a+2)$.
As usual, we write $\widehat{d}\equiv1\wedge d$, and write a subscripted
expression $x_{y}$ interchangeably with $x(y)$.
\begin{thm}
\textbf{\emph{\label{Thm. Extension of measurable r.f. from Q_inf to Q, given continuity in prob}
(Extension of measurable r.f. with parameter set $Q_{\infty}$ to
the full parameter set $Q$, given continuity in probability).}} Consider
the locally compact metric space $(S,d)$, without necessarily any
linear structure or ordering. Let $(\Omega_{0},L_{0},E_{0})$ be an
arbitrary probability space. Recall the space $\widehat{R}_{Cp}(Q_{\infty}\times\Omega_{0},S)$
of r.f.'s which are continuous in probability over the parameter subspace
$(Q_{\infty},d_{Q})$. Recall the space $\widehat{R}_{Meas,Cp}(Q\times\Omega,S)$
of r.f.'s which are defined and continuous in probability on the full
parameter space $(Q,d_{Q})$. 

Then there exists a probability space $(\Omega,L,E),$ and a function
\[
\Phi_{meas,\xi(Q)}:\widehat{R}_{Cp}(Q_{\infty}\times\Omega_{0},S)\rightarrow\widehat{R}_{Meas,Cp}(Q\times\Omega,S)
\]
such that, for each $Z\in\widehat{R}_{Cp}(Q_{\infty}\times\Omega_{0},S)$
with a modulus of continuity in probability $\delta_{Cp}$, the r.f.
\[
X\equiv\Phi_{meas,\xi(Q)}(Z):Q\times(\Omega,L,E)\rightarrow S
\]
satisfies the following conditions.

1. For a.e. $\theta\in\Theta_{1}$, we have $X_{s}(\theta,\cdot)=Z_{s}$
a.s. on $\Omega_{0}$ for each $s\in Q_{\infty}$.

2. The r.f. $X|Q_{\infty}$ is equivalent to $Z$. 

3. The r.f. $X$ is \textup{\emph{measurable}}\emph{ }and continuous
in probability, with the same modulus of continuity in probability
$\delta_{Cp}$ as $Z$.

4. There exists a full subset $D$ of $(\Omega,L,E)$ such that, for
each $\omega\in D,$ and for each $t\in Q$, there exists a sequence
$(s_{j})$ in $Q_{\infty}$ with $d_{Q}(t,s_{j})\rightarrow0$ and
$d(X(t,\omega),X(s_{j},\omega))\rightarrow0$ as $j\rightarrow\infty$.
\end{thm}
\begin{proof}
1. Let 
\[
(\Theta_{1},L_{1},I_{1})\equiv([0,1],L_{1},\int\cdot d\theta)
\]
denote the Lebesgue integration space based on the interval $[0,1]$.
Define the product sample space 
\[
(\Omega,L,E)\equiv(\Theta_{1},L_{1},I_{1})\otimes(\Omega_{0},L_{0},E_{0}).
\]

2. Consider each r.f. $Z\in\widehat{R}_{Cp}(Q_{\infty}\times\Omega_{0},S)$,
with a modulus of continuity in probability $\delta_{Cp}$. Define
the full subset 
\[
D_{0}\equiv\bigcap_{q\in Q(\infty)}domain(Z_{q})\subset\Omega_{0}.
\]

3. Augment each sample from $\Omega_{0}$ with a secondary sample
from $\Theta_{1}$. More precisely, define a function 
\[
\widetilde{Z}:Q_{\infty}\times\Omega\rightarrow S
\]
by 
\[
domain(\widetilde{Z})\equiv Q_{\infty}\times\Theta_{1}\times D_{0}
\]
and 
\[
\widetilde{Z}(q,\theta,\omega_{0})\equiv Z(q,\omega_{0})
\]
for each $(q,\theta,\omega_{0})\in domain(\widetilde{Z})$. Then,
for each $q\in Q_{\infty}$, the function $\widetilde{Z}_{q}$ is
a .r.v. on $\Omega$, according to Propositions \ref{Prop. vector of measurable func (S)   is meas func (S)}
and \ref{Prop. X meas,  f unif continuous and bd on bd subsets =00003D> f(X) meas}.
Thus $\widetilde{Z}$ is a r.f. We proceed to extend the r.f. $\widetilde{Z}$,
by a sequence of \index{stochastic interpolation}stochastic interpolations,
to a $\mathrm{measurable}$ r.f. $X:Q\times\Omega\rightarrow S$.

4. To that end, let $m\geq1$ and $k=1,\cdots,\gamma_{m}$ be arbitrary,
where $\gamma_{m}$ is as in Definition \ref{Def. Specification of  Q_inf, and partion of unity of Q}.
 Define 
\[
\Delta_{m,k}\equiv\{(t,\theta)\in Q\times\Theta_{1}:\theta\in(\lambda_{m,k-1}^{+}(t),\lambda_{m,k}^{+}(t))\}.
\]
Relation \ref{eq:temp-130} says that $\lambda_{m,\gamma(m)}^{+}=1$.
Hence Proposition \ref{Prop. Region between non-neg functions in C(Q) are integrable.}
implies that the sets $\Delta_{m,1},\cdots,\Delta_{m,\gamma(m)}$
are mutually disjoint integrable subsets of $Q\times\Theta_{1}$,
and that their union $\bigcup_{i=1}^{\gamma(m)}\Delta_{m,i}$ is a
full subset. Define a function $X^{(m)}:Q\times\Omega\rightarrow S$
by 
\[
domain(X^{(m)})\equiv\bigcup_{i=1}^{\gamma(m)}\Delta_{m,i}\times D_{0},
\]
and by 
\begin{equation}
X^{(m)}(t,\omega)\equiv\widetilde{Z}(q_{m,i},\omega)\equiv Z(q_{m,i},\omega_{0})\label{eq:temp-120-2}
\end{equation}
for each $(t,\omega)\equiv(t,\theta,\omega_{0})\in\Delta_{m,i}\times D_{0}$,
for each $i=1,\cdots,\gamma_{m}$. Since the $\mathrm{measurable}$
sets $\Delta_{m,1}\times D_{0},\cdots,\Delta_{m,\gamma(m)}\times D_{0}$
are mutually exclusive in $Q\times\Omega$, with union equal to a
full subset, the function $X^{(m)}:Q\times\Omega\rightarrow S$ is
$\mathrm{measurable}$ on $(Q\times\Omega,\Lambda\otimes L,I\otimes E)$,
by Proposition \ref{Prop. Basing seq of measurable functions on measurable partition}.

5. Now let $t\in Q$ be arbitrary. Define the open interval 
\begin{equation}
\Delta_{m,k,t}\equiv(\lambda_{m,k-1}^{+}(t),\lambda_{m,k}^{+}(t))\equiv\{\theta\in\Theta_{1}:\theta\in(\lambda_{m,k-1}^{+}(t),\lambda_{m,k}^{+}(t))\}.\label{eq:temp-118}
\end{equation}
Then $\bigcup_{i=1}^{\gamma(m)}\Delta_{m,i,t}$ is a full subset of
$\Theta_{1}\equiv[0,1]$. Hence $\Delta_{m,1,t}\times D_{0},$ $\cdots,$
$\Delta_{m,\gamma(m),t}\times D_{0}$ are mutually exclusive $\mathrm{measurable}$
subsets of $\Omega$ whose union is a full subset of $\Omega$. Furthermore,
by the definition of $X^{(m)}$ in Step 2, we have 
\begin{equation}
X_{t}^{(m)}(\omega)\equiv\widetilde{Z}(q_{m,i},\omega)\equiv Z(q_{m,i},\omega_{0})\label{eq:temp-120}
\end{equation}
for each $(\theta,\omega_{0})\in\Delta_{m,i,t}\times D_{0}$, for
each $i=1,\cdots,\gamma_{m}$. Hence $X_{t}^{(m)}:\Omega\rightarrow S$
is a r.v. by Proposition \ref{Prop. Basing seq of measurable functions on measurable partition}.
Thus we see that $X^{(m)}:Q\times\Omega\rightarrow S$ is a r.f. By
Step 2, $X^{(m)}$ is a $\mathrm{measurable}$ function. Therefore
$X^{(m)}$ is a $\mathrm{measurable}$ r.f. 

Intuitively, for each $t\in Q$, the r.v. $X_{t}^{(m)}$ is set to
the r.v. $Z_{q(m,i)}$ with probability $|\Delta_{m,i,t}|=\lambda_{m,q(m,i),t}$,
for each $i=1,\cdots,\gamma_{m}$. In this sense, $X^{(m)}$ is a
stochastic interpolation \index{stochastic interpolation} of $Z_{q(m,1)},\cdots,Z_{q(m,\gamma(m))}$.
 Note that the probabilities $\lambda_{m,q(m,1),t},\cdots,$ $\lambda_{m,q(m,\gamma(m)),t}$
are continuous functions of $t$. We will later prove that the r.f.
$X^{(m)}$ is continuous in probability, even though its sample functions
are piecewise constant.

6. We will first construct an a.u. convergent subsequence of $(X^{(m)})_{m=1,2,\cdots}$.
Let $m_{0}\equiv0$. Let $j\geq1$ be arbitrary. Write, for abbreviation,
\begin{equation}
\varepsilon_{j}\equiv2^{-j}.\label{eq:temp-238-3}
\end{equation}
Let
\begin{equation}
n_{j}\equiv j\vee[(2-\log_{2}\delta_{Cp}(\varepsilon_{j}))]_{1}.\label{eq:temp-507}
\end{equation}
Recursively on $j\geq1$, define 
\begin{equation}
m_{j}\equiv m_{j-1}\vee n_{j}.\label{eq:temp-223-2}
\end{equation}
Define 
\begin{equation}
X\equiv\lim_{j\rightarrow\infty}X^{(m(j))}.\label{eq:temp-176}
\end{equation}
A priori, the limit need not exist anywhere. We will show that actually
$X_{t}^{(m(j))}\rightarrow X_{t}$ a.u. for each $t\in Q$, and that
therefore $X$ is a well defined function and is a r.f. 

7. To that end, let $j\geq1$ be arbitrary. Then $m_{j}\geq n_{j}\geq j$.
Consider each $t\in Q$ and $s\in Q_{\infty}$ with
\begin{equation}
d_{Q}(s,t)<2^{-1}\delta_{Cp}(\varepsilon_{j}).\label{eq:temp-172}
\end{equation}
Consider each $\theta\in\bigcup_{i=1}^{\gamma(m(j))}\Delta_{m(j),i,t}$.
Then, for each $\omega_{0}\in D_{0}$, we have 
\[
(\theta,\omega_{0})\in\bigcup_{i=1}^{\gamma(m(j))}\Delta_{m(j),i,t}\times D_{0}\equiv domain(X_{t}^{(m(j))}),
\]
and 
\[
\hat{d(}\widetilde{Z}_{s}(\theta,\omega_{0}),X_{t}^{(m(j))}(\theta,\omega_{0}))
\]
\[
=\sum_{i=1}^{\gamma(m(j))}1_{\Delta(m(j),i,t)}(\theta)\hat{d(}Z_{s}(\omega_{0}),X_{t}^{(m(j))}(\theta,\omega_{0}))
\]
\begin{equation}
=\sum_{i=1}^{\gamma(m(j))}1_{\Delta(m(j),i,t)}(\theta)\hat{d(}Z_{s}(\omega_{0}),Z_{q(m(j),i)}(\omega_{0})),\label{eq:temp-230}
\end{equation}
where the last inequality follows from the defining formula \ref{eq:temp-120}.
Hence, since $D_{0}$ is a full set, we obtain 
\begin{equation}
E_{0}\hat{d(}\widetilde{Z}_{s}(\theta,\cdot),X_{t}^{(m(j))}(\theta,\cdot))=\sum_{i=1}^{\gamma(m(j))}1_{\Delta(m(j),i,t)}(\theta)E_{0}\hat{d(}Z_{s},Z_{q(m(j),i)})\label{eq:temp-504}
\end{equation}
Suppose the summand with index $i$ on the right-hand side is positive.
Then $\Delta_{m(j),i,t}\equiv(\lambda_{m(j),i-1}^{+}(t),\lambda_{m(j),i}^{+}(t))$
is a non-empty open interval. Hence
\[
\lambda_{m(j),i-1}^{+}(t)<\lambda_{m(j),i}^{+}(t).
\]
Equivalently, $\lambda_{m(j),q(m(j),i)}(t)>0.$ At the same time,
the continuous function $\lambda_{m(j),q(m(j),i)}$ on $Q$ has support
$(d_{Q}(\cdot,q_{m(j),i})\leq2^{-m(j)+1})$, as observed in the remarks
preceding this theorem. Consequently, 
\begin{equation}
d_{Q}(t,q_{m(j),i})\leq2^{-m(j)+1}.\label{eq:temp-261}
\end{equation}
Inequalities \ref{eq:temp-172} and \ref{eq:temp-261} together imply
that 
\[
d_{Q}(s,q_{m(j),i})\leq d_{Q}(s,t)+d_{Q}(t,q_{m(j),i})<\frac{1}{2}\delta_{Cp}(\varepsilon_{j})+2^{-m(j)+1}
\]
\begin{equation}
<\frac{1}{2}\delta_{Cp}(\varepsilon_{j})+\frac{1}{2}\delta_{Cp}(\varepsilon_{j})<\delta_{Cp}(\varepsilon_{j}),\label{eq:temp-236}
\end{equation}
where the third inequality follows from defining formulas \ref{eq:temp-507}
and \ref{eq:temp-223-2}. By the definition of $\delta_{Cp}$ as a
modulus of continuity in probability of $Z$ on $Q_{\infty}$, inequality
\ref{eq:temp-236} yields 
\[
E_{0}\hat{d(}Z_{s},Z_{q(m(j),i)})\leq\varepsilon_{j}.
\]
Summing up, the above inequality holds for the $i$-th summand in
the right-hand side of equality \ref{eq:temp-504} if said $i$-th
summand is positive. Equality \ref{eq:temp-504} therefore results
in
\begin{equation}
E_{0}\hat{d(}\widetilde{Z}_{s}(\theta,\cdot),X_{t}^{(m(j))}(\theta,\cdot))\leq\varepsilon_{j}\sum_{i=1}^{\gamma(m(j))}1_{\Delta(m(j),i,t)}(\theta)=\varepsilon_{j},\label{eq:temp-506}
\end{equation}
where $\theta$ is an arbitrary member of the full set $\bigcup_{i=1}^{\gamma(m(j))}\Delta_{m(j),i,t}$.
Therefore, by Fubini's Theorem,
\begin{equation}
E\hat{d(}\widetilde{Z}_{s},X_{t}^{(m(j))})\equiv I_{1}\otimes E_{0}\hat{d(}\widetilde{Z}_{s},X_{t}^{(m(j))})\leq\varepsilon_{j},\label{eq:temp-131-3-1}
\end{equation}
where $t\in Q$ and $s\in Q_{\infty}$ are arbitrary with
\[
d_{Q}(s,t)<\frac{1}{2}\delta_{Cp}(\varepsilon_{j}).
\]

8. Now let $j\geq1$ and $t\in Q$ be arbitrary. Take any $s\in Q_{\infty}$
such that 
\[
d_{Q}(s,t)<\frac{1}{2}\delta_{Cp}(\varepsilon_{j})\wedge\frac{1}{2}\delta_{Cp}(\varepsilon_{j+1}).
\]
Then it follows from inequality \ref{eq:temp-131-3-1} that
\[
E\hat{d(}X_{t}^{(m(j))},X_{t}^{(m(j+1))})
\]
\[
\leq E\hat{d(}\widetilde{Z}_{s},X_{t}^{(m(j))})+E\hat{d(}\widetilde{Z}_{s},X_{t}^{(m(j+1))})
\]
\begin{equation}
\leq\varepsilon_{j}+\varepsilon_{j+1}\equiv2^{-j}+2^{-j-1}<2^{-j+1},\label{eq:temp-229}
\end{equation}
where $t\in Q$ and $j\geq1$ are arbitrary. Hence, by Assertion 5
of Proposition \ref{Prop. Basics of the probability metric}, the
function $X_{t}\equiv\lim_{j\rightarrow\infty}X_{t}^{(m(j))}$ is
a r.v., and 
\[
X_{t}^{(m(j))}\rightarrow X_{t}\quad a.u.
\]
We conclude that 
\[
X:Q\times\Omega\rightarrow S
\]
is a r.f. 

9. We will now show that $X$ is a $\mathrm{measurable}$ r.f. Note
that, by Fubini's Theorem, inequality \ref{eq:temp-229} implies 
\begin{equation}
(I\otimes E)\hat{d(}X^{(m(j))},X^{(m(j+1))})\leq2^{-j+1}\label{eq:temp-173}
\end{equation}
for each $j\geq1$. Since $X^{(m(j))}$ is a $\mathrm{measurable}$
function on $(Q\times\Omega,\Lambda\otimes L,I\otimes E)$, for each
$j\geq1$, Assertion 5 of Proposition \ref{Prop. Basics of the probability metric}
implies that $X_{t}^{(m(j))}\rightarrow X$ a.u. on $(Q\times\Omega,\Lambda\otimes L,I\otimes E)$,
and that $X\equiv\lim_{j\rightarrow\infty}X_{t}^{(m(j))}$ is a $\mathrm{measurable}$
function on $(Q\times\Omega,\Lambda\otimes L,I\otimes E)$. Thus $X$
is a $\mathrm{measurable}$ r.f. 

10. Define the full subset 
\[
D_{1}\equiv\bigcap_{t\in Q(\infty)}\bigcap_{n=1}^{\infty}\bigcup_{i=1}^{\gamma(n)}\Delta_{n,i,t}
\]
of $\Theta_{1}\equiv[0,1]$. Let $\theta\in D_{1}$ be arbitrary.
For each $s\in Q_{\infty}$, letting $j\rightarrow\infty$ in inequality
\ref{eq:temp-131-3-1} with $t=s$, we obtain $E_{0}\hat{d(}\widetilde{Z}_{s},X_{s}(\theta,\cdot))=0,$
whence
\begin{equation}
Z_{s}=\widetilde{Z}_{s}(\theta,\cdot)=X_{s}(\theta,\cdot)\quad a.s.\label{eq:temp-226}
\end{equation}
on $(\Omega_{0},L_{0},E_{0})$. Condition 1 in the conclusion of the
theorem is proved. 

11. Now let $K\geq1$, $f\in C_{ub}(S^{K})$, and $s_{1},\cdots,s_{K}\in Q_{\infty}$
be arbitrary. Then, in view of equality \ref{eq:temp-226}, Fubini's
Theorem implies
\[
Ef(X_{s(1)},\cdots,X_{s(K)})=(I_{1}\otimes E_{0})f(X_{s(1)},\cdots,X_{s(K)})=E_{0}f(Z_{s(1)},\cdots,Z_{s(K)}).
\]
Thus the r.f.'s $X|Q_{\infty}$ and $Z$ are equivalent, establishing
Condition 2 in the conclusion of the theorem.

12. We will prove that $X$ is continuous in probability. For that
purpose, let $\varepsilon>0$ be arbitrary. Let $j\geq1$ and $t,t'\in Q$
be arbitrary with
\[
d_{Q}(t,t')<\delta_{Cp}(\varepsilon).
\]
Since $Q_{\infty}$ is dense in $Q$, there exist $s,s'\in Q_{\infty}$
with $d_{Q}(s,s')<\delta_{Cp}(\varepsilon)$ and 
\[
d_{Q}(t,s)\vee d_{Q}(t',s')<\frac{1}{2}\delta_{Cp}(\varepsilon_{j}).
\]
It follows that $E_{0}\hat{d(}Z_{s},Z_{s'})\leq\varepsilon$. We can
then apply inequality \ref{eq:temp-131-3-1} to obtain 
\[
E\hat{d(}X_{t}^{(m(j))},X_{t'}^{(m(j))})\leq E\widehat{d}(X_{t}^{(m(j))},\widetilde{Z}_{s})+E\hat{d(}\widetilde{Z}_{s},\widetilde{Z}_{s'})+E\hat{d(}\widetilde{Z}_{s'},X_{t'}^{(m(j))})
\]
\[
\leq\varepsilon_{j}+(I\otimes E_{0})\hat{d(}1_{\Theta(1)}\otimes Z_{s},1_{\Theta(1)}\otimes Z_{s'})+\varepsilon_{j}.
\]
\[
=\varepsilon_{j}+E_{0}\hat{d(}Z_{s},Z_{s'})+\varepsilon_{j}\leq\varepsilon_{j}+\varepsilon+\varepsilon_{j},
\]
where the equality is thanks to Fubini's Theorem. Letting $j\rightarrow\infty$
yields
\[
E\hat{d(}X_{t},X_{t'})\leq\varepsilon,
\]
where $\varepsilon>0$ is arbitrarily small. Summing up, the r.f.
$X$ is continuous in probability on $Q$, with $\delta_{Cp}$ as
a modulus of continuity in probability. Condition 3 has been established.

13. For each $s\in Q_{\infty}$, letting $j\rightarrow\infty$ in
inequality \ref{eq:temp-131-3-1} with $t=s$, we obtain
\begin{equation}
E\hat{d(}\widetilde{Z}_{s},X_{s})\equiv I_{1}\otimes E_{0}\hat{d(}\widetilde{Z}_{s},X_{s})=0.\label{eq:temp-131-3-1-1}
\end{equation}
Hence 
\begin{equation}
D_{2}\equiv\bigcap_{s\in Q(\infty)}(\widetilde{Z}_{s}=X_{s})\label{eq:temp-226-1}
\end{equation}
is a full subset of $(\Omega,L,E)$. Define the full subset 
\[
D\equiv D_{2}\sqcap(D_{1}\times D_{0})\cap\bigcap_{s\in Q(\infty)}\bigcap_{j=1}^{\infty}domain(X_{s}^{(m(j))})
\]
of the sample space $(\Omega,L,E)$. 

Consider each $\omega\equiv(\theta,\omega_{0})\in D\subset D_{2}$.
Then $\theta\in D_{1}$ and $\omega_{0}\in D_{0}$. Let
\[
t\in domain(X(\cdot,\omega))
\]
be arbitrary. In other words, $(t,\omega)\in domain(X)$. Hence, by
the defining equality \ref{eq:temp-176}, we have
\begin{equation}
X(t,\omega)=\lim_{j\rightarrow\infty}X^{(m(j))}(t,\omega).\label{eq:temp-176-1}
\end{equation}
Let $\varepsilon>0$ be arbitrary. Let $J\geq1$ be so large that
\begin{equation}
(t,\omega)\in domain(X^{(m(j))})\equiv\bigcup_{i=1}^{\gamma(m(j))}\Delta_{m(j),i}\times D_{0}\label{eq:temp-223}
\end{equation}
and 
\[
d(X(t,\omega),X^{(m(j))}(t,\omega))<\varepsilon
\]
for each $j\geq J$. Now consider each $j\geq J$. By relation \ref{eq:temp-223},
there exists $i_{j}=1,\cdots,\gamma_{m(j)}$ such that 
\[
(t,\theta,\omega_{0})\equiv(t,\omega)\in\Delta_{m(j),i(j)}\times D_{0}.
\]
Hence 
\begin{equation}
X^{(m(j))}(t,\omega)=\widetilde{Z}(q_{m(j),i(j)},\omega)=X(q_{m(j),i(j)},\omega),\label{eq:temp-120-2-1-1}
\end{equation}
where the first equality follows from equality \ref{eq:temp-120-2},
and the second from equality \ref{eq:temp-226-1} and from the membership
$\omega\in D_{2}$. At the same time, since $t\in\Delta_{m(j),i(j)}$,
we have
\begin{equation}
d_{Q}(t,q_{m(j),i(j)})\leq2^{-m(j)+1}\label{eq:temp-261-1}
\end{equation}
according to inequality \ref{eq:temp-261}. Summing up, for each $\omega$
in the full set $D$, and for each $t\in Q$, the sequence $(s_{j})\equiv(q_{m(j),i(j)})$
in $Q_{\infty}$ is such that
\[
d_{Q}(t,s_{j})\rightarrow0
\]
and
\[
d(X(t,\omega),X(s_{j},\omega))\rightarrow0.
\]
Thus Condition 4 of the conclusion of the theorem has also been proved.
\end{proof}
Recall the set $\widehat{F}(Q,S)$ of consistent families of f.j.d.'s
with the parameter set $Q$ and state space $S$. Recall the subset
$\widehat{F}_{Cp}(Q,S)$ whose members are continuous in probability,
equipped with the metric $\widehat{\rho}_{Cp,\xi,Q|Q(\infty)}$ defined
in Definition \ref{Def. Metric on  of continuous in prob families of finite joint distributions.}
by
\[
\widehat{\rho}_{Cp,\xi,Q|Q(\infty)}(F,F')\equiv\widehat{\rho}_{Marg,\xi,Q}(F|Q_{\infty},F'|Q_{\infty})
\]
\begin{equation}
\equiv\sum_{n=1}^{\infty}2^{-n}\rho_{Dist,\xi^{n}}(F_{q(1),\cdots,q(n)},F'_{q(1),\cdots,q(n)})\label{eq:temp-208-1-1}
\end{equation}
for each $F,F'\in\widehat{F}_{Cp}(Q,S)$. 

The $\mathrm{measurable}$ extension of a consistent family of f.j.d.'s
which is continuous in probability is another immediate corollary
of Theorem \ref{Thm. Compact Daniell-Kolmogorov-Skorohod}. 
\begin{thm}
\textbf{\emph{\label{Thm. Measurable Extension of consistent family of fjd.s continuous in prob}(Construction
of measurable r.f. from family of consistent f.j.d.'s which is continuous
in probability).}} Consider the locally compact metric space $(S,d)$,
without necessarily any linear structure or ordering. Let 
\[
(\Theta_{0},L_{0},I_{0})\equiv([0,1],L_{0},\int\cdot d\theta)
\]
denote the Lebesgue integration space based on the interval $[0,1]$.
Then there exists a function
\[
\Phi_{meas,\xi,\xi(Q)}:\widehat{F}_{Cp}(Q,S)\rightarrow\widehat{R}_{Meas,Cp}(\Theta_{0}\times\Omega,S)
\]
such that, for each $F\in\widehat{F}_{Cp}(Q,S)$, the $\mathrm{measurable}$
r.f. 
\[
X\equiv\Phi_{meas,\xi,\xi(Q)}(F):Q\times(\Theta_{0},L_{0},I_{0})\rightarrow S
\]
has marginal distributions given by the family $F$. 

We will refer to the function $\Phi_{meas,\xi,\xi(Q)}$ as the \index{measurable extension}
$\mathrm{measurable}$\emph{ extension} relative to the binary approximations
$\xi$ and $\xi_{Q}$ of $(S,d)$ and $(Q,d_{q})$ respectively.
\end{thm}
\begin{proof}
1. Let $\Phi_{Q,Q(\infty)}:$$\widehat{F}_{Cp}(Q,S)\rightarrow\widehat{F}_{Cp}(Q_{\infty},S)$
be the function defined by
\[
\Phi_{Q,Q(\infty)}(F)\equiv F|Q_{\infty}
\]
for each $F\in\widehat{F}_{Cp}(Q,S)$. Let $\Phi_{DKS,\xi}$ be the
Daniell-Kolmogorov-Skorokhod extension as constructed in Theorem \ref{Thm. DKS Extension, construction and continuity}.
Let $\Phi_{meas,\xi(Q)}$ be the function constructed in Theorem \ref{Thm. Extension of measurable r.f. from Q_inf to Q, given continuity in prob}.

2. We will prove that the composite function 
\begin{equation}
\Phi_{meas,\xi,\xi(Q)}\equiv\Phi_{meas,\xi(Q)}\circ\Phi_{DKS,\xi}\circ\Phi_{Q,Q(\infty)}\label{eq:temp-232}
\end{equation}
has the desired properties. To that end, let $F\in\widehat{F}_{Cp}(Q,S)$
be arbitrary, with a modulus of continuity of probability $\delta_{Cp}$.
Then $F|Q_{\infty}$ is also continuous in probability, with the same
modulus of continuity of probability $\delta_{Cp}$. Let 
\[
Z\equiv\Phi_{DKS,\xi}(F|Q_{\infty})
\]
be the Daniell-Kolmogorov-Skorokhod extension of $F|Q_{\infty}$ relative
to the binary approximation $\xi$ of $(S,d)$. Thus 
\[
Z:Q_{\infty}\times(\Theta_{0},L_{0},I_{0})\rightarrow S
\]
is a r.f. with marginal distributions given by $F|Q_{\infty}$. It
follows that the r.f. $Z$ has modulus of continuity of probability
$\delta_{Cp}$. Hence Theorem \ref{Thm. Extension of measurable r.f. from Q_inf to Q, given continuity in prob}
applies, and yields the $\mathrm{measurable}$ r.f. 
\[
X\equiv\Phi_{meas,\xi(Q)}(Z):Q\times(\Theta_{0},L_{0},I_{0})\rightarrow S.
\]
Now define 
\[
\Phi_{meas,\xi,\xi(Q)}(F)\equiv X=\Phi_{meas,\xi(Q)}\circ\Phi_{DKS,\xi}\circ\Phi_{Q,Q(\infty)}(F).
\]
According to Theorem \ref{Thm. Extension of measurable r.f. from Q_inf to Q, given continuity in prob},
$X$ is continuous in probability with the same modulus of continuity
of probability $\delta_{Cp}$. Moreover, $X|Q_{\infty}$ is equivalent
to $Z$. Hence $X|Q_{\infty}$ has marginal distributions given by
$F|Q_{\infty}$. 

4. It remains to prove that the r.f. $X$ has marginal distributions
given by $F$. To that end, consider each $k\geq1$, $r_{1},\cdots,r_{k}\in Q$,
$s_{1},\cdots,s_{k}\in Q_{\infty}$, and $f\in C_{ub}(S^{k},d^{k})$.
Then
\[
I_{0}f(X_{s(1)},\cdots,X_{s(k)})=I_{0}f(Z_{s(1)},\cdots,Z_{s(k)})=F_{s(1),\cdots,s(k)}f.
\]
Now let $s_{k}\rightarrow r_{k}$ in $Q$ for each $k=1,\cdots,k$.
Then the left-hand side converges to $I_{0}f(X_{r(1)},\cdots,X_{r(k)})$,
on account of the continuity in probability of $X$. The right-hand
side converges to $F_{r(1),\cdots,r(k)}f$ by the continuity in probability
of $F$, according to Lemma \ref{Lem. Continuity of finite joint distributions}.
Hence
\[
I_{0}f(X_{r(1)},\cdots,X_{r(k)})=F_{r(1),\cdots,r(k)}f.
\]
The theorem is proved.
\end{proof}
We will next prove the metrical continuity of the mapping $\Phi_{meas,\xi,\xi(Q)}$.
Recall from Definition \ref{Def. Metric on  of continuous in prob families of finite joint distributions.}
the metric $\widehat{\rho}_{Cp,\xi,Q|Q(\infty)}$ on the space $\widehat{F}_{Cp}(Q,S)$.
\begin{thm}
\textbf{\emph{\label{Thm. Continuity-of-Measurable Extension} (Continuity
of construction of measurable r.f.'s)}} Let $\widehat{F}_{Cp,\delta(Cp)}(Q,S)$
be an arbitrary subset of $\widehat{F}_{Cp}(Q,S)$ whose members have
a common modulus of continuity in probability $\delta_{Cp}$. Recall
that 
\[
(\Theta_{0},L_{0},I_{0})\equiv([0,1],L_{0},\int\cdot d\theta)
\]
denote the Lebesgue integration space based on the interval $[0,1]$. 

Then the onstruction 
\begin{equation}
\Phi_{meas,\xi,\xi(Q)}:(\widehat{F}_{Cp}(Q,S),\widehat{\rho}_{Cp,\xi,Q|Q(\infty)})\rightarrow(\widehat{R}_{Meas,Cp}(Q\times\Theta_{0},S),\rho_{Sup,Prob})\label{eq:temp-122}
\end{equation}
in Theorem \ref{Thm. Measurable Extension of consistent family of fjd.s continuous in prob}
is uniformly continuous on the subset $\widehat{F}_{Cp,\delta(Cp)}(Q,S)$,
with a modulus of continuity $\delta_{fjd,meas}(\cdot,\delta_{Cp},\left\Vert \xi\right\Vert ,\bigl\Vert\xi_{Q}\bigr\Vert)$
dependent only on $\delta_{Cp},$ $\left\Vert \xi\right\Vert $, and
$\bigl\Vert\xi_{Q}\bigr\Vert$.
\end{thm}
\begin{proof}
Refer to the proofs of  Theorem \ref{Thm. Extension of measurable r.f. from Q_inf to Q, given continuity in prob}
and Theorem \ref{Thm. Measurable Extension of consistent family of fjd.s continuous in prob},
where the defining equality \ref{eq:temp-232} leads to 
\begin{equation}
\Phi_{meas,\xi,\xi(Q)}|\widehat{F}_{Cp,\delta(Cp)}(Q,S)=\Phi_{meas,\xi(Q)}\circ\Phi_{DKS,\xi}\circ(\Phi_{Q,Q(\infty)}|\widehat{F}_{Cp,\delta(Cp)}(Q,S)).\label{eq:temp-175}
\end{equation}
For the uniform continuity of $\Phi_{meas,\xi,\xi(Q)}|\widehat{F}_{Cp,\delta(Cp)}(Q,S)$,
we need only verify the continuity of the three functions on the right-hand
side, and compound their moduli of continuity.

1. By Definition \ref{Def. Metric on  of continuous in prob families of finite joint distributions.},
the function
\[
\Phi_{Q,Q(\infty)}:(\widehat{F}_{Cp}(Q,S),\widehat{\rho}_{Cp,\xi,Q|Q(\infty)})\rightarrow(\widehat{F}(Q_{\infty},S),\widehat{\rho}_{Marg,\xi,Q}),
\]
defined by $\Phi_{Q,Q(\infty)}(F)\equiv F|Q_{\infty}$ for each $F\in\widehat{F}_{Cp}(Q,S)$,
is metric preserving. It is therefore uniformly continuous, with a
trivial modulus of continuity $\delta_{Q,Q(\infty)}$ given by $\delta_{Q,Q(\infty)}(\varepsilon)\equiv\varepsilon$
for each $\varepsilon>0$.

2. Let $F\in\widehat{F}_{Cp,\delta(Cp)}(Q,S)$ be arbitrary. By hypothesis,
$F$ is continuous in probability on $Q$, with a modulus of continuity
in probability $\delta_{Cp}$.\emph{ }Hence its restriction\emph{
}$F|Q_{\infty}$ is trivially continuous in probability on $Q_{\infty}$,
with the same modulus of continuity in probability $\delta_{Cp}$.
According to Theorem \ref{Thm. DKS Extension, construction and continuity},
the Daniell-Kolmogorov-Skorokhod Extension
\[
\Phi_{DKS,\xi}:(\widehat{F}(Q_{\infty},S),\widehat{\rho}_{Marg,\xi,Q(\infty)})\rightarrow(\widehat{R}(Q\times\Theta_{0},S),\rho_{Q\times\Theta(0),S})
\]
is uniformly continuous on the subset 
\[
\widehat{F}_{Cp,\delta(Cp)}(Q,S)|Q_{\infty}\equiv\Phi_{Q,Q(\infty)}(\widehat{F}_{Cp,\delta(Cp)}(Q,S)),
\]
with a modulus of continuity $\delta_{DKS}(\cdot,\left\Vert \xi\right\Vert )$
dependent only on the modulus of local compactness $\left\Vert \xi\right\Vert \equiv(|A_{k}|)_{k=1,2,}$
of the locally compact state space $(S,d)$.

3. It remains to verify that the function 
\[
\Phi_{meas,\xi(Q)}:(\widehat{R}(Q_{\infty}\times\Theta_{0},S),\rho_{Q\times\Theta(0),S})\rightarrow(\widehat{R}_{Meas,Cp}(Q\times\Theta_{0},S),\rho_{Sup,Prob})
\]
is uniformly continuous on the subset 
\[
\widehat{R}_{0}\equiv\Phi_{DKS,\xi}(\widehat{F}_{Cp,\delta(Cp)}(Q,S)|Q_{\infty})\equiv\Phi_{DKS,\xi}(\Phi_{Q,Q(\infty)}(\widehat{F}_{Cp,\delta(Cp)}(Q,S))).
\]

4. As in the proof of  Theorem \ref{Thm. Extension of measurable r.f. from Q_inf to Q, given continuity in prob}
define, for each $j\geq1$, 
\begin{equation}
\varepsilon_{j}\equiv2^{-j},\label{eq:temp-238-3-1}
\end{equation}
\begin{equation}
n_{j}\equiv j\vee[(2-\log_{2}\delta_{Cp}(\varepsilon_{j}))]_{1},\label{eq:temp-507-1}
\end{equation}
and 
\begin{equation}
m_{j}\equiv m_{j-1}\vee n_{j},\label{eq:temp-223-2-1}
\end{equation}
where $m_{0}\equiv0$. 

5. Now let $\varepsilon>0$ be arbitrary. Let
\[
j\equiv[0\vee(4-\log_{2}\varepsilon)]_{1}.
\]
Then $2^{-j}<2^{-4}\varepsilon$. Define
\[
\delta_{meas}(\varepsilon,\delta_{Cp},\bigl\Vert\xi_{Q}\bigr\Vert)\equiv2^{-\gamma(m(j))-4}\varepsilon^{2}.
\]
Let $Z,Z'\in\widehat{R}_{0}$ be arbitrary such that
\begin{equation}
\rho_{Q\times\Theta(0),S}(Z,Z')<\delta_{meas}(\varepsilon,\delta_{Cp},\bigl\Vert\xi_{Q}\bigr\Vert).\label{eq:temp-168}
\end{equation}
We will verify that $\rho_{Sup,Prob}(\Phi_{meas,\xi(Q)}(Z),\Phi_{meas,\xi(Q)}(Z'))<\varepsilon$. 

6. Write $X\equiv\Phi_{meas,\xi(Q)}(Z)$, and $X'\equiv\Phi_{meas,\xi(Q)}(Z')$.
Thus $Z,Z':Q_{\infty}\times\Theta_{0}\rightarrow S$ are r.f.'s, and
$X,X':Q\times\Theta_{0}\rightarrow S$ are $\mathrm{measurable}$
r.f.'s. As in the proof of  Theorem \ref{Thm. Extension of measurable r.f. from Q_inf to Q, given continuity in prob},
define the full subset 
\[
D_{0}\equiv\bigcap_{q\in Q(\infty)}domain(Z_{q})
\]
of $\Omega_{0}\equiv\Theta_{0}\equiv[0,1]$. Similarly, define the
full subset $D_{0}'\equiv\bigcap_{q\in Q(\infty)}domain(Z'_{q})$
\[
D_{0}\equiv\bigcap_{q\in Q(\infty)}domain(Z_{q})
\]
of $\Theta_{0}$. Then $D_{0}D_{0}'$ is a full subset of $\Theta_{0}$.
Note that inequality \ref{eq:temp-168} is equivalent to 
\begin{equation}
I_{0}\sum_{i=1}^{\infty}2^{-i}\widehat{d}(Z_{t(i)},Z'_{t(i)})<2^{-\gamma(m(j))-4}\varepsilon^{2}.\label{eq:temp-131}
\end{equation}
Hence, by Chebychev's inequality, there exists a $\mathrm{measurable}$
set $A\subset\Theta_{0}$ with 
\[
I_{0}A^{c}<2^{-2}\varepsilon,
\]
such that 
\begin{equation}
\sum_{i=1}^{\infty}2^{-i}\hat{d(}Z(t_{i},\omega_{0}),Z'(t_{i},\omega_{0}))\leq2^{-\gamma(m(j))-2}\varepsilon,\label{eq:temp-256}
\end{equation}
for each $\omega_{0}\in A$. 

7. Let $\omega_{0}\in AD_{0}D'_{0}$ be arbitrary. Then inequality
\ref{eq:temp-256} trivially implies that 
\begin{equation}
\bigvee_{i=1}^{\gamma(m(j))}\hat{d(}Z(t_{i},\omega_{0}),Z'(t_{i},\omega_{0}))\leq2^{-2}\varepsilon.\label{eq:temp-257}
\end{equation}
Let $t\in Q$ be arbitrary. Recall that

\[
\widetilde{\Delta}\equiv\bigcup_{i=1}^{\gamma(m(j))}\Delta_{n,i,t}\equiv\bigcup_{i=1}^{\gamma(m(j))}(\lambda_{n,i-1}^{+}(t),\lambda_{n,i}^{+}(t))
\]
is a full subset of $\Theta_{1}\equiv[0,1]$. Let $\theta\in$$\widetilde{\Delta}$
be arbitrary. Then $\theta\in$$\Delta_{n,i,t}$ for some $i=1,\cdots,\gamma_{m(j)}$.
Hence $(\theta,\omega_{0})\in\Delta_{n,i,t}\times D_{0}D'_{0}$. Therefore
the defining equality \ref{eq:temp-120} in the proof of  Theorem
\ref{Thm. Extension of measurable r.f. from Q_inf to Q, given continuity in prob}
says that 
\begin{equation}
X_{t}^{(m(j))}(\theta,\omega_{0})=Z(q_{m,i},\omega_{0})\label{eq:temp-120-1}
\end{equation}
and
\begin{equation}
X_{t}'^{(m(j))}(\theta,\omega_{0})=Z'(q_{m,i},\omega_{0}).\label{eq:temp-120-1-1}
\end{equation}
Consequently, in view of inequality \ref{eq:temp-257}, we have 
\[
\widehat{d}(X_{t}^{(m(j))}(\theta,\omega_{0}),X_{t}'^{(m(j))}(\theta,\omega_{0}))
\]
\begin{equation}
=\widehat{d}(Z(q_{m,i},\omega_{0}),Z'(q_{m,i},\omega_{0}))\leq2^{-2}\varepsilon,\label{eq:temp-259-1}
\end{equation}
where $(\theta,\omega_{0})\in\widetilde{\Delta}\times AD_{0}D'_{0}$
is arbitrary. Since $\widehat{d}\leq1$, it follows that
\[
I_{0}\widehat{d}(X_{t}^{(m(j))},X_{t}'^{(m(j))})\leq I_{0}\widehat{d}(X_{t}^{(m(j))},X_{t}'^{(m(j))})1_{\widetilde{\Delta}\times AD(0)D'(0)}+I_{0}1_{\widetilde{\Delta}\times(AD(0)D'(0))^{c}}
\]
\[
\leq2^{-2}\varepsilon+(I_{1}\otimes I_{0})1_{\widetilde{\Delta}\times A^{c}}=2^{-2}\varepsilon+I_{0}1_{A^{c}}<2^{-2}\varepsilon+2^{-2}\varepsilon=2^{-1}\varepsilon.
\]

8. Separately, inequality \ref{eq:temp-229} in the proof of  Theorem
\ref{Thm. Extension of measurable r.f. from Q_inf to Q, given continuity in prob}
implies that 
\[
I_{0}\hat{d(}X_{t}{}^{(m(j))},X{}_{t})\leq\sum_{i=j}^{\infty}2^{-i+1}=2^{-j+2},
\]
and, similarly, that
\[
I_{0}\hat{d(}X_{t}'^{(m(j))},X'_{t})<2^{-j+2}.
\]
Combining the last three displayed inequalities, we obtain 
\[
I_{0}\hat{d(}X_{t},X'_{t})<2^{-1}\varepsilon+2^{-j+2}+2^{-j+2}=2^{-1}\varepsilon+2^{-j+3},
\]
where $t\in Q$ is arbitrary. Therefore
\[
\rho_{Sup,Prob}(X,X')\equiv\sup_{t\in Q}I_{0}\widehat{d}(X_{t},X'_{t})\leq2^{-1}\varepsilon+2^{-j+3}<2^{-1}\varepsilon+2^{-1}\varepsilon=\varepsilon.
\]
In other words, 
\[
\rho_{Sup,Prob}(\Phi_{meas,\xi(Q)}(Z),\Phi_{meas,\xi(Q)}(Z'))<\varepsilon,
\]
as alleged. Thus $\delta_{meas}(\cdot,\delta_{Cp},\bigl\Vert\xi_{Q}\bigr\Vert)$
is a modulus of continuity of $\Phi_{meas,\xi(Q)}$ on 
\[
\widehat{R}_{0}\equiv\Phi_{DKS,\xi}(\widehat{F}_{Cp,\delta(Cp)}(Q,S)|Q_{\infty})=\Phi_{DKS,\xi}(\Phi_{Q,Q(\infty)}(\widehat{F}_{Cp,\delta(Cp)}(Q,S))).
\]

9. Combining, we conclude that the composite function $\Phi_{meas,\xi,\xi(Q)}|\widehat{F}_{0}$
in equality \ref{eq:temp-175} is uniformly continuous, with the composite
modulus of continuity
\[
\delta_{fjd,meas}(\cdot,\delta_{Cp},\left\Vert \xi\right\Vert ,\bigl\Vert\xi_{Q}\bigr\Vert)
\]
\[
=\delta_{Q,Q(\infty)}(\delta_{DKS}(\delta_{meas}(\varepsilon,\delta_{Cp},\bigl\Vert\xi_{Q}\bigr\Vert),\left\Vert \xi\right\Vert ))
\]
\[
=\delta_{DKS}(\delta_{meas}(\varepsilon,\delta_{Cp},\bigl\Vert\xi_{Q}\bigr\Vert),\left\Vert \xi\right\Vert ),
\]
as desired.
\end{proof}

\section{Measurable Gaussian Random Fields}

Let $(Q,d_{Q})$ be a compact metric space, with $d_{Q}\leq1$ and
with an arbitrary, but fixed, distribution $I$. As an application
of the Theorems \ref{Thm. Measurable Extension of consistent family of fjd.s continuous in prob},
we will construct a $\mathrm{measurable}$ Gaussian r.f. $X:Q\times\Omega\rightarrow R$
from its continuous mean and covariance functions, and will prove
the continuity of this construction. For that purpose we need only
prove that from the mean and covariance functions we can construct
a consistent family of normal f.j.d.'s which is continuous in probability,
and that the construction is continuous.
\begin{defn}
\label{Def. Gaussian r.f.} \textbf{(Gaussian r.f.)} A r.f. $X:Q\times\Omega\rightarrow R$
is said to be \emph{Gaussian}\index{Gaussian r.f.} if it has marginal
distributions which are normal. The functions $\mu(t)\equiv EX_{t}$
and $\sigma(t,s)\equiv E(X_{t}-EX_{t})(X_{s}-EX_{s})$ are called
the mean and covariance functions, respectively, of the r.f. $X$.
 $\square$ 
\end{defn}
Without loss of generality, we will treat only the case where the
r.f. is centered, with $EX_{t}=0$ for each $t\in Q$. The more general
case where the mean is a non-trivial continuous function $\mu$ follows
by adding $\mu$ to a centered r.f. $X$. 

Recall the matrix terminologies in Definition \ref{Def. Matrix notations}. 
\begin{defn}
\label{Def. nonnegative definite functions} \textbf{(Nonnegative
definite functions).} Let $D$ be an arbitrary nonempty set. Write
$D^{2}\equiv D\times D$.. Let $\sigma:D^{2}\rightarrow[0,\infty)$
be an arbitrary symmetric function. If, for each $m\geq1$ and for
each $r_{1},\cdots,r_{m}\in D$, the square matrix $[\sigma(r_{k},r_{h})]_{k=1,\cdots,m;h=1,\cdots,m}$
is nonnegative definite, then $\sigma$ is said to be a \index{nonnegative definite function}\emph{nonnegative
definite function} on the set $D^{2}$. If, for each $m\geq1$ and
for each $r_{1},\cdots,r_{m}\in D$, the matrix $[\sigma(r_{k},r_{h})]_{k=1,\cdots,m;h=1,\cdots,m}$
is positive definite, then $\sigma$ is said to be a \index{positive definite function}\emph{
positive definite function} on the set $D^{2}$. $\square$
\end{defn}
In the rest of this chapter, let $\xi$ and  $\xi_{Q}\equiv(B_{n})_{n=1,2,\cdots}$
be arbitrary, but fixed, binary approximations of the Euclidean state
space $(S,d)\equiv(R,d)$ and and the compact parameter space $(Q,d_{Q})$
respectively, as specified in Definitions \ref{Def. Specification of state space and parameter space}
and \ref{Def. Specification of  Q_inf, and partion of unity of Q}
respectively. Recall the enumerated, countably infinite, dense subset
\begin{equation}
Q_{\infty}\equiv\{t_{1},t_{2},\cdots\}\equiv\bigcup_{n=1}^{\infty}B_{n}\label{eq:temp-238-1-1-2}
\end{equation}
of $Q$, where $B_{n}\equiv\{q_{n,1},\cdots,q_{n,\gamma(n)}\}=\{t_{1},\cdots,t_{\gamma(n)}\}$
as sets, for each $n\geq1$. 
\begin{prop}
\label{Prop. Consistency of f.j.d.'s generated by covariance function}
\textbf{\emph{(Consistency of family of normal f.j.d.'s generated
by covariance function).}} Let $\sigma:Q\times Q\rightarrow[0,\infty)$
be a continuous nonnegative definite function. For each $m\geq1$
and each $r_{1},\cdots,r_{m}\in Q$, write the nonnegative definite
matrix 
\begin{equation}
\overline{\sigma}\equiv[\sigma(r_{k},r_{h})]_{k=1,\cdots,m;h=1,\cdots,m},\label{eq:temp-138}
\end{equation}
and define 
\begin{equation}
F_{r(1),\cdots,r(m)}^{\sigma}\equiv\Phi_{0,\overline{\sigma}},\label{eq:temp-3-1-1}
\end{equation}
where $\Phi_{0,\overline{\sigma}}$ is the normal distribution with
mean $0$ and covariance matrix $\overline{\sigma}$. Then the following
holds.

1. The family
\begin{equation}
F^{\sigma}\equiv\Phi_{covar,fjd}(\sigma)\equiv\{F_{r(1),\cdots,r(m)}^{\sigma}:m\geq1;r_{1},\cdots,r_{m}\in Q\}\label{eq:temp-132}
\end{equation}
of f.j.d.'s is consistent. 

2. The consistent family  $F^{\sigma}$ is continuous in probability.
In symbols, $F^{\sigma}\in\widehat{F}_{Cp}(Q,R)$. Specifically, suppose
$\delta_{0}$ is a modulus of continuity of $\sigma$ on the compact
metric space $(Q^{2},d_{Q}^{2})$. Then $F^{\sigma}$ has a modulus
of continuity in probability defined by
\[
\delta_{Cp}(\varepsilon)\equiv\delta_{Cp}(\varepsilon,\delta_{0})\equiv\delta_{0}(\frac{1}{2}\varepsilon^{2})
\]
 for each $\varepsilon>0$. 
\end{prop}
\begin{proof}
1. Let $n,m\geq1$ be arbitrary. Let $r\equiv(r_{1},\cdots,r_{m})$
be an arbitrary sequence in $Q$, and let $j\equiv(j_{1},\cdots,j_{n})$
be an arbitrary sequence in $\{1,\cdots,m\}$. Let the matrix $\overline{\sigma}$
be defined as in equality \ref{eq:temp-138} above. By Lemma \ref{Lem. Normal distribution is well defined for nonneg definite matrix },
$F_{r(1),\cdots,r(m)}^{\sigma}\equiv\Phi_{0,\overline{\sigma}}$ is
the distribution of a r.v. $Y\equiv AZ$, where $A$ is an $m\times m$
matrix such that $\overline{\sigma}\equiv AA^{T}$, and where $Z$
is a standard normal r.v. on some probability space $(\Omega,L,E)$,
with values in $R^{m}$. Let the dual function $j^{*}:R^{m}\rightarrow R^{n}$
be defined by
\[
j^{*}(x_{1},\cdots,x_{m})\equiv(x_{j(1)},\cdots,x_{j(n)})
\]
for each $x\equiv(x_{1},\cdots,x_{m})\in R^{m}$. Then $j^{*}(x)=Bx$
for each $x\in R^{m}$, where the $n\times m$ matrix 
\[
B\equiv[b_{k,h}]_{k=1,\cdots,n;h=1,\cdots,m},
\]
is defined by $b_{k,h}\equiv1$ or 0 according as $h=j_{k}$ or $h\neq j_{k}$.
Let $\widetilde{A}\equiv BA$. Define the $n\times n$ matrix
\[
\tilde{\sigma}\equiv\widetilde{A}\widetilde{A}^{T}=BAA^{T}B^{T}=B\overline{\sigma}B^{T}=[\sigma(r_{j(k)},r_{j(h)})]_{k=1,\cdots,n;h=1,\cdots,n}.
\]
Then, by the defining formula \ref{eq:temp-3-1-1}, 
\[
F_{r(j(1)),\cdots,r(j(n))}^{\sigma}\equiv\Phi_{0,\tilde{\sigma}}.
\]
At the same time, the r.v.
\[
\tilde{Y}\equiv j^{*}(Y)=BY=BAZ\equiv\widetilde{A}Z
\]
has the normal characteristic function defined by 
\[
E(\exp i\lambda^{T}\widetilde{A}Z)=\exp(-\frac{1}{2}\lambda^{T}\widetilde{A}\widetilde{A}^{T}\lambda)=\exp(-\frac{1}{2}\lambda^{T}\tilde{\sigma}\lambda)
\]
for each $\lambda\in R^{n}$. Hence $\tilde{Y}$ has the normal distribution
$\Phi_{0,\tilde{\sigma}}$. Combining, we see that, for each $f\in C(R^{n})$,
\[
F_{r(1),\cdots,r(m)}^{\sigma}(f\circ j^{*})=E(f\circ j^{*}(Y))=Ef(\tilde{Y})=\Phi_{0,\tilde{\sigma}}(f)=F_{r(j(1)),\cdots,r(j(n))}^{\sigma}f.
\]
We conclude that the family $F^{\sigma}$ of f.j.d.'s is consistent. 

2. Now consider the case where $m=2$. Consider each $r_{1},r_{2}\in Q$
with
\begin{equation}
d_{Q}(r_{1},r_{2})<\delta_{Cp}(\varepsilon)\equiv\delta_{0}(\frac{1}{2}\varepsilon^{2}).\label{eq:temp-125}
\end{equation}
Let $d$ denote the Euclidean metric for $R$. As in Step 1, there
exists a r.v. $Y\equiv(Y_{1},Y_{2})$ with values in $R^{2}$ with
the normal distribution $\Phi_{0,\overline{\sigma}}$, where $\overline{\sigma}\equiv[\sigma(r_{k},r_{h})]_{k=1,2;h=1,2}$.
Then
\[
F_{r(1),r(2)}^{\sigma}(1\wedge d)=\Phi_{0,\overline{\sigma}}(1\wedge d)\leq\Phi_{0,\overline{\sigma}}d=E|Y_{1}-Y_{2}|\leq\sqrt{E(Y_{1}-Y_{2})^{2}}
\]
\[
=\sqrt{\sigma(r_{1},r_{1})-2\sigma(r_{1},r_{2})+\sigma(r_{2},r_{2})}
\]
\[
\leq\sqrt{|\sigma(r_{1},r_{1})-\sigma(r_{1},r_{2})|+|\sigma(r_{2},r_{2})-\sigma(r_{1},r_{2})|}
\]
\[
\leq\sqrt{\frac{1}{2}\varepsilon^{2}+\frac{1}{2}\varepsilon^{2}}=\varepsilon,
\]
where the second inequality is Lyapunov's inequality, and the last
is due to inequality \ref{eq:temp-125}. Thus $F^{\sigma}$ is continuous
in probability, with $\delta_{Cp}(\cdot,\delta_{0})$ as a modulus
of continuity in probability.
\end{proof}
Recall from Definition \ref{Def. Metric on  of continuous in prob families of finite joint distributions.}
the metric space $(\widehat{F}_{Cp}(Q,R),\widehat{\rho}_{Cp,\xi,Q|Q(\infty)})$
of consistent families of f.j.d.'s with parameter space $(Q,d_{Q})$
and state space $R$.
\begin{prop}
\label{Prop. covariance function to normal fjd  is continuous} \textbf{\emph{(Normal
f.j.d.'s depend continuously on covariance function).}} Let $G$ denote
the set of continuous nonnegative definite functions $\sigma:Q\times Q\rightarrow R$.
Equip $G$ with the metric $d_{G}$ defined by
\[
d_{G}(\sigma,\sigma')\equiv\sup_{(t,s)\in Q\times Q}|\sigma(t,s)-\sigma'(t,s)|
\]
for each $\sigma,\sigma'\in G$. Then the function 
\[
\Phi_{covar,fjd}:(G,d_{G})\rightarrow(\widehat{F}_{Cp}(Q,R),\widehat{\rho}_{Cp,\xi,Q|Q(\infty)})
\]
in Proposition \ref{Prop. Consistency of f.j.d.'s generated by covariance function}
is uniformly continuous, with a modulus of continuity $\delta_{covar,fjd}$
defined in equality \ref{eq:temp-180} in the proof below.
\end{prop}
\begin{proof}
1. Let $\varepsilon>0$ be arbitrary. Let $n\geq1$ be arbitrary.
By Theorem \ref{Thm. Continuity Theorem for ch functions}, there
exists $\delta_{ch,dstr}(\varepsilon,n)>0$ such that, for arbitrary
distributions $J,J'$ on $R^{n}$ whose respective characteristic
functions $\psi,\psi'$ satisfy 
\begin{equation}
\rho_{ch,n}(\psi,\psi')\equiv\sum_{j=1}^{\infty}2^{-j}\sup_{|\lambda|\leq j}|\psi(\lambda)-\psi'(\lambda)|<\delta_{ch,dstr}(\varepsilon,n),\label{eq:temp-37-1}
\end{equation}
we have 
\begin{equation}
\rho_{Dist,\xi^{n}}(J,J')<\varepsilon,\label{eq:temp-125-1}
\end{equation}
where $\rho_{Dist,\xi^{n}}$ is the metric on the space of distributions
on $R^{n}$, as in Definitions \ref{Def. Distribution metric}

2. Let $\varepsilon>0$ be arbitrary. Let $m\geq1$ be so large that
$2^{-m+1}<\varepsilon$. Let $K\geq1$ be so large that 
\[
2^{-K+1}<\alpha\equiv1\wedge e^{-1}\bigwedge_{n=1}^{m}\delta_{ch,dstr}(\frac{\varepsilon}{2},n).
\]
Then, since $0\leq\alpha\leq1$, basic calculus shows that $e^{\alpha}-1\leq\alpha(e-1)$.
Define
\begin{equation}
\delta_{covar,fjd}(\varepsilon)\equiv2K^{-2}m^{-2}\alpha.\label{eq:temp-180}
\end{equation}
We will verify that $\delta_{covar,fjd}$ is the desired modulus of
continuity of $\Phi_{covar,fjd}$.

3. To that end, let $\sigma,\sigma'\in G$ be arbitrary such that
\begin{equation}
d_{G}(\sigma,\sigma')<\delta_{covar,fjd}(\varepsilon).\label{eq:temp-37}
\end{equation}
Let $F^{\sigma}\equiv\Phi_{covar,fjd}(\sigma)$ and $F^{\sigma'}\equiv\Phi_{covar,fjd}(\sigma')$
be constructed as in Theorem \ref{Thm. Extension to Gaussian r.f.}.
We will show that $\widehat{\rho}_{Cp,\xi,Q|Q(\infty)}(F^{\sigma},F^{\sigma'})<\varepsilon$. 

4. First note that inequality \ref{eq:temp-37} is equivalent to 
\begin{equation}
\sup_{(t,s)\in Q\times Q}|\sigma(t,s)-\sigma'(t,s)|<2K^{-2}m^{-2}\alpha.\label{eq:temp-11}
\end{equation}
Next, let $n=1,\cdots,m$ be arbitrary. The joint normal distribution
$F_{t(1),\cdots,t(n)}^{\sigma}$ has characteristic function defined
by 
\[
\chi_{t(1),\cdots,t(n)}^{\sigma}(\lambda)\equiv\exp(-\frac{1}{2}\sum_{k=1}^{n}\sum_{h=1}^{n}\lambda_{k}\sigma(t_{k},t_{h})\lambda_{h})
\]
for each $\lambda\equiv(\lambda_{1},\cdots,\lambda_{n})\in R^{n}$,
with a similar equality for $\sigma'$. It follows that
\[
\rho_{char}(\chi_{t(1),\cdots,t(n)}^{\sigma},\chi_{t(1),\cdots,t(n)}^{\sigma'})\equiv\sum_{j=1}^{\infty}2^{-j}\sup_{|\lambda|\leq j}|\chi_{t(1),\cdots,t(n)}^{\sigma}(\lambda)-\chi_{t(1),\cdots,t(n)}^{\sigma'}(\lambda)|
\]
\[
<\sup_{|\lambda|\leq K}|\chi_{t(1),\cdots,t(n)}^{\sigma}(\lambda)-\chi_{t(1),\cdots,t(n)}^{\sigma'}(\lambda)|+\sum_{j=K+1}^{\infty}2^{-j}\cdot2
\]
\[
=\sup_{|\lambda|\leq K}|\exp(-\frac{1}{2}\sum_{k=1}^{n}\sum_{h=1}^{n}\lambda_{k}\sigma(t_{k},t_{h})\lambda_{h})-\exp(-\frac{1}{2}\sum_{k=1}^{n}\sum_{h=1}^{n}\lambda_{k}\sigma'(t_{k},t_{h})\lambda_{h})|+2^{-K+1}.
\]
By the real variable inequality $|e^{-x}-e^{-y}|=e^{-y}|e^{-x+y}-1|\leq e^{|x-y|}-1$
for arbitrary $x,y\geq0$, the last displayed expression is bounded
by
\[
\sup_{|\lambda|\leq K}(\exp(\frac{1}{2}\sum_{k=1}^{n}\sum_{h=1}^{n}|\lambda_{k}(\sigma(t_{k},t_{h})-\sigma'(t_{k},t_{h})\lambda_{h}|)-1)+2^{-K+1}
\]
\[
\leq(\exp(\frac{K^{2}}{2}\sum_{k=1}^{n}\sum_{h=1}^{n}|\sigma(t_{k},t_{h})-\sigma'(t_{k},t_{h})|)-1)+2^{-K+1}
\]
\[
\leq(\exp(\frac{K^{2}}{2}\sum_{k=1}^{n}\sum_{h=1}^{n}2K^{-2}m^{-2}\alpha)-1)+2^{-K+1}
\]
\[
\leq(e^{\alpha}-1)+2^{-K+1}\leq\alpha(e-1)+\alpha=\alpha e\leq\delta_{ch,dstr}(\frac{\varepsilon}{2},n),
\]
where the second inequality is from inequality \ref{eq:temp-11} above.
Hence, according to inequality \ref{eq:temp-125-1}, we have
\[
\rho_{Dist,\xi^{n}}(F_{t(1),\cdots,t(n)}^{\sigma},F{}_{t(1),\cdots,t(n)}^{\sigma'})<\frac{\varepsilon}{2},
\]
where $n=1,\cdots,m$ is arbitrary. Therefore, according to Definition
\ref{Def. Metric on  of continuous in prob families of finite joint distributions.},
\[
\widehat{\rho}_{Cp,\xi,Q|Q(\infty)}(F^{\sigma},F^{\sigma'})\equiv\sum_{n=1}^{\infty}2^{-n}\rho_{Dist,\xi^{n}}(F_{t(1),\cdots,t(n)}^{\sigma},F{}_{t(1),\cdots,t(n)}^{\sigma})
\]
\[
\leq\sum_{n=1}^{m}2^{-n}\rho_{Dist,\xi^{n}}(F_{t(1),\cdots,t(n)}^{\sigma},F{}_{t(1),\cdots,t(n)}^{\sigma})+\sum_{n=m+1}^{\infty}2^{-n}
\]
\begin{equation}
\leq\sum_{n=1}^{m}2^{-n}\frac{\varepsilon}{2}+2^{-m}<\frac{\varepsilon}{2}+\frac{\varepsilon}{2}<\varepsilon,\label{eq:temp-208-1-1-1}
\end{equation}
where we used the bounds $0\leq\rho_{Dist,\xi^{n}}\leq1$ for each
$n\geq1$.

Since $\varepsilon>0$ is arbitrarily small, we conclude that the
function $\Phi_{covar,fjd}$ is uniformly continuous, with modulus
of continuity $\delta_{covar,fjd}$.
\end{proof}
Now we can mechanically apply the theorems in the previous section.
As in the previous section, let
\[
(\Theta_{0},L_{0},I_{0})\equiv(\Theta_{1},L_{1},I_{1})\equiv([0,1],L_{1},\int\cdot d\theta)
\]
be the Lebesgue integration space based on the interval $[0,1]$,
and let 
\[
(\Omega,L,E)\equiv(\Theta_{1},L_{1},E_{1})\otimes(\Omega_{0},L_{0},E_{0}).
\]

\begin{thm}
\label{Thm. Extension to Gaussian r.f.} \textbf{\emph{(Construction
of }}$\mathbf{measurable}$\textbf{\emph{ Gaussian r.f. from continuous
covariance function).}} Let $\sigma:Q\times Q\rightarrow R$ be a
continuous nonnegative definite function. Then there exists a \textup{\emph{measurable}}\emph{
}Gaussian r.f. 
\[
X\equiv\Phi_{cov,gauss,\xi,\xi(Q)}(\sigma)\equiv\Phi_{meas,\xi,\xi(Q)}\circ\Phi_{covar,fjd}(\sigma):Q\times\Omega\rightarrow R
\]
which is continuous in probability, and which is such that $EX_{t}=0$
and $EX_{t}X_{s}=\sigma(t,s)$ for each $t,s\in Q$. We will call
the function $\Phi_{cov,gauss,\xi,\xi(Q)}$ the \textup{\emph{measurable}}\emph{
} Gaussian extension relative to the binary approximations $\xi$
and $\xi_{Q}$.
\end{thm}
\begin{proof}
By Proposition \ref{Prop. Consistency of f.j.d.'s generated by covariance function},
the family $F^{\sigma}\equiv\Phi_{covar,fjd}(\sigma)$ of normal f.j.d.'s
is consistent and continuous in probability. Hence Theorem \ref{Thm. Measurable Extension of consistent family of fjd.s continuous in prob}
is applicable to $F^{\sigma}$ and yields the $\mathrm{measurable}$
r.f. 
\[
X\equiv\Phi_{meas,\xi,\xi(Q)}(F^{\sigma})\equiv\Phi_{meas,\xi,\xi(Q)}\circ\Phi_{covar,fjd}(\sigma),
\]
with marginal distributions given by $F^{\sigma}$. Since $F^{\sigma}$
is continuous in probability by Proposition \ref{Prop. Consistency of f.j.d.'s generated by covariance function},
so is $X$. 
\end{proof}
Recall from Definition \ref{Def. Metric on measurable r.f. continuous in probability}
the metric space  $(\widehat{R}_{Meas,Cp}(Q\times\Omega,R),\rho_{Sup,Prob})$
of $\mathrm{measurable}$ r.f.'s $X:Q\times\Omega\rightarrow R$ which
are continuous in probability. Thus 
\[
\rho_{Sup,Prob}(X,Y)\equiv\sup_{t\in Q}E1\wedge|X_{t}-Y_{t}|
\]
for each $X,Y\in\widehat{R}_{Meas,Cp}(Q\times\Omega,R)$. 
\begin{thm}
\label{Thm. Continuity  of Gaussian extension} \textbf{\emph{(Continuity
of the construction of }}$\mathbf{measurable}$\textbf{\emph{ Gaussian
r.f.'s).}} Use the same assumptions and notations as in Proposition
\ref{Prop. covariance function to normal fjd  is continuous} and
Theorem \ref{Thm. Extension to Gaussian r.f.}. Suppose $G_{0}$ is
a subset of the set $G$ of continuous nonnegative definite functions
$\sigma:Q\times Q\rightarrow R$ whose members share a common modulus
of continuity $\delta_{0}$ on $Q\times Q$. Suppose, in addition,
that there exists $b_{0}\geq0$ such that $\sigma(t,t)\leq b_{0}$
for each $t\in Q$, for each $\sigma\in G_{0}$. 

Then the \textup{\emph{measurable}}\emph{ } Gaussian extension 
\[
\Phi_{cov,gauss,\xi,\xi(Q)}:(G,d_{G})\rightarrow(\widehat{R}_{Meas,Cp}(Q\times\Omega,R),\rho_{Sup,Prob})
\]
constructed in Theorem \ref{Thm. Extension to Gaussian r.f.} is uniformly
continuous on the subset $G_{0}$ of $G$, with a modulus of continuity
$\delta_{cov,gauss}(\cdot,\delta_{0},b_{0},\left\Vert \xi\right\Vert ,\bigl\Vert\xi_{Q}\bigr\Vert)$.
\end{thm}
\begin{proof}
1. By the construction in Theorem \ref{Thm. Extension to Gaussian r.f.},
\[
\Phi_{cov,gauss,\xi,\xi(Q)}\equiv\Phi_{meas,\xi,\xi(Q)}\circ\Phi_{covar,fjd}.
\]

2. By Lemma \ref{Prop. covariance function to normal fjd  is continuous},
the function 
\[
\Phi_{covar,fjd}:(G,d_{G})\rightarrow(\widehat{F}_{Cp}(Q,R),\widehat{\rho}_{Cp,\xi,Q|Q(\infty)})
\]
is uniformly continuous, with a modulus of continuity $\delta_{covar,fjd}$. 

3. By Proposition \ref{Prop. Consistency of f.j.d.'s generated by covariance function},
the members of $\widehat{F}_{Cp,\delta(Cp)}(Q,S)\equiv\Phi_{covar,fjd}(G_{0})$
share the same modulus of continuity in probability defined by $\delta_{Cp}(\varepsilon,\delta_{0})\equiv\delta_{0}(\frac{1}{2}\varepsilon^{2})$
for each $\varepsilon>0$. 

4. Hence $\widehat{F}_{Cp,\delta(Cp)}(Q,S)$ satisfies the conditions
in the hypothesis of Theorem \ref{Thm. Continuity-of-Measurable Extension}.
According to Assertion 2 of Theorem \ref{Thm. Continuity-of-Measurable Extension},
the function 
\[
\Phi_{meas,\xi,\xi(Q)}:(\widehat{F}_{Cp,\delta(Cp)}(Q,S),\widehat{\rho}_{Cp,\xi,Q|Q(\infty)})\rightarrow(\widehat{R}_{Meas,Cp}(Q\times\Omega,R),\rho_{Sup,Prob})
\]
is uniformly continuous on $\widehat{F}_{Cp,\delta(Cp)}(Q,S)$, with
a modulus of continuity
\[
\delta_{fjd,meas}(\cdot,\delta_{Cp},\beta,\left\Vert \xi\right\Vert ,\bigl\Vert\xi_{Q}\bigr\Vert).
\]

5. Combining, the composite function $\Phi_{cov,gauss,\xi,\xi(Q)}$
is uniformly continuous, with a modulus of continuity defined by the
composite operation
\[
\delta_{cov,gauss}(\varepsilon,\delta_{0},b_{0},\left\Vert \xi\right\Vert ,\bigl\Vert\xi_{Q}\bigr\Vert)\equiv\delta_{covar,fjd}(\delta_{fjd,meas}(\varepsilon,\delta_{Cp}(\cdot,\delta_{0}),\beta,\left\Vert \xi\right\Vert ,\bigl\Vert\xi_{Q}\bigr\Vert))
\]
for each $\varepsilon>0$.
\end{proof}

\chapter{Martingales}

In this chapter, we define a martingale $X\equiv\{X_{t}:t=1,2,\ldots\}$
for modeling one's fortune in a fair game of chance. Then we will
prove the basic theorems on martingales which have wide-ranging applications.
Among these is the a.u. convergence of $X_{t}$ as $t\rightarrow\infty$.
Our proof is constructive and quantifies rates of convergence by means
of a maximal inequality. There are proofs in traditional texts which
also are constructive and quantify rates similarly by means of maximal
inequalities. These traditional maximal inequalities, however, require
the integrability of $|X_{t}|^{p}$ for some $p>1$, or at least the
integrability of $|X_{t}|\log|X_{t}|$. For the separate case of $p=1$,
the classical proof of a.u. convergence is by a separate inference
from certain upcrossing inequalities. Such inference is essentially
equivalent to the principle of infinite search, and is not constructive. 

In contrast, the maximal inequality we present requires only the integrability
of $|X_{t}$|. Therefore, thanks to Lyapunov's inequality, it is at
once applicable to the case of integrable $|X_{t}|^{p}$ for any given
$p\geq1$, without having to first determine whether $p>1$ or $p=1$. 

For readers who are uninitiated in the subject, the previous paragraphs
are perhaps confusing, but will become clear as we proceed. For the
rich body of classical results on, and applications of, martingales,
see e.g. \cite{Doob53,Chung68,Durret84}.
\begin{defn}
\textbf{\label{Def. Assumptions and Notations}(Assumptions and Notations).
}In this chapter, let $(S,d)$ be a locally compact metric space with
an arbitrary, but fixed, reference point $x_{\circ}$. Let $(\Omega,L,E)$
be an arbitrary probability space. Unless otherwise specified, a r.v.
refers to a $\mathrm{measurable}$ function with values in $S$. 

If $(\Omega,L',E)$ is a probability subspace of $(\Omega,L,E)$,
we will simply call $L'$ a probability subspace of $L$ when $\Omega$
and $E$ are understood. Let $Q$ denote an arbitrary nonempty subset
of $R$, called the time parameter set. 

For abbreviation, we will write $A\in L$ if $A$ is a $\mathrm{measurable}$
subset of $(\Omega,L,E)$. Thus $A\in L$ iff $1_{A}\in L$, in which
case we will write $P(A)$, $PA$, $E1_{A}$, and $EA$ interchangeably,
and write $E(X;A)\equiv EX1_{A}$ for each $X\in L$. As usual, we
write a subscripted expression $x_{y}$ interchangeably with $x(y)$. 
\end{defn}
$\square$

\section{\label{Section :Filtration}Filtrations}

Let $Q$ denote an arbitrary nonempty subset of $R$.
\begin{defn}
\label{Def. Filtration} \textbf{(Filtration). }Suppose that, for
each $t\in Q$, there exists a probability subspace $(\Omega,L^{(t)},E)$
of $(\Omega,L,E)$, such that $L^{(t)}\subset L^{(s)}$ for each $t,s\in Q$
with $t\leq s$. Then the family $\mathcal{L}\equiv\{L^{(t)}:t\in Q\}$
is called a \emph{filtration}\index{filtration} in $(\Omega,L,E)$
with time parameter set $Q$. The filtration $\mathcal{L}$ is said
to be \index{right continuous filtration}\emph{right continuous }if,
for each $t\in Q$, we have
\[
L^{(t)}=\bigcap_{s\in Q;s>t}L^{(s)}.
\]

Suppose, in addition, that $\widetilde{Q}$ is a subset of $Q$. Then
a stochastic process $X:\widetilde{Q}\times\Omega\rightarrow S$ is
 said to be \index{adapted process} \emph{adapted} to the filtration
$\mathcal{L}$ if $X_{t}$ is a r.v. on $(\Omega,L^{(t)},E)$ for
each $t\in Q$. 
\end{defn}
$\square$

The probability space $L^{(t)}$ can be regarded as the observable
history up to the time $t$. Thus a process $X$ adapted to $\mathcal{L}$
is such that $X_{t}$ is observable at the time $t$, for each $t\in Q$.
Note that if all points in the set $Q$ are isolated points in $Q$,
then each filtration $\mathcal{L}$ with time parameter set $Q$ is
right continuous.
\begin{defn}
\label{Def. Natutal filtration of a process} \textbf{(Natural filtration
of a stochastic process).} Let $X:Q\times\Omega\rightarrow S$ be
an arbitrary stochastic process. For each $t\in Q$, define the set
\[
G^{(X,t)}\equiv\{X_{r}:r\in Q;\;r\leq t\},
\]
and let 
\[
L^{(X,t)}\equiv L(X_{r}:r\in Q;\;r\leq t)\equiv L(G^{(X,t)})
\]
be the\emph{ }probability subspace of\emph{ }$L$\emph{ }generated
by the set $G^{(X,t)}$of r.v.'s. Then the family $\mathcal{L_{X}}\equiv\{L^{(X,t)}:t\in Q\}$
is called the \emph{natural filtration of the process\index{filtration generated by a process}}
$X$. $\square$ 
\end{defn}
\begin{lem}
\label{Lem. Filtration genrerated by process is indeed a filtration}
\textbf{\emph{(A natural filtration is indeed a filtration).}} Let
$X:Q\times\Omega\rightarrow S$ be an arbitrary stochastic process.
Then the natural filtration $\mathcal{L}_{X}$ of \textbf{$X$ }is
a filtration to which the process $X$ is adapted.
\end{lem}
\begin{proof}
For each $t\leq s$ in $Q$ we have $G^{(X,t)}\subset G^{(X,s)}$
whence $L^{(X,t)}\subset L^{(X,s)}$. Thus $\mathcal{L}_{X}$ is a
filtration. Let $t\in Q$ be arbitrary. Then $f(X_{t})\in L(G^{(X,t)})\equiv L^{(X,t)}$
for each $f\in C_{ub}(S,d)$. At the same time, because $X_{t}$ is
a r.v. on $(\Omega,L,E)$\emph{, }we have\emph{ $P(d(X_{t},x_{\circ})\geq a)\rightarrow0$
as $a\rightarrow\infty$}. Hence $X_{t}$ is a r.v. on $(\Omega,L^{(X,t)},E)$
according to Proposition \ref{Prop. Basic Properties of r.v.}. Thus
the process $X$ is adapted to the its natural filtration $\mathcal{L}_{X}$.
\end{proof}
\begin{defn}
\label{Def. Right -limit extsion and right continuity of filtration}
\textbf{(Right-limit extension and right continuity of a filtration).}
Suppose (i) $\overline{Q}=[0,\infty)$ or (ii) $\overline{Q}\equiv[0,a]$
for some $a>0$. Suppose $Q$ is a subset which is dense in $\overline{Q}$
and which, in Case (ii), contains the end point $a$. Let $\mathcal{L}\equiv\{L^{(t)}:t\in Q\}$
be an arbitrary filtration of a given probability space $(\Omega,L,E)$. 

In Case (i) define, for each $t\in\overline{Q}$, the probability
subspace 
\begin{equation}
L^{(t+)}\equiv\bigcap\{L^{(s)}:s\in Q\cap(t,\infty)\}.\label{eq:temp-330}
\end{equation}
of $L$. In Case (ii) define, for each $t\in\overline{Q}$, the probability
subspace 
\begin{equation}
L^{(t+)}\equiv\bigcap\{L^{(s)}:s\in Q\cap(t,a]\cup\{a\})\}.\label{eq:temp-330-2}
\end{equation}
Then the filtration $\mathcal{L}^{+}\equiv\{L^{(t+)}:t\in\overline{Q}\}$
is called the \emph{right-limit extension}\index{right-limit extension of a filtration}
of the filtration $\mathcal{L}$. 

If $Q=\overline{Q}$ and $L^{(t)}=L^{(t+)}$for each $t\in\overline{Q}$,
then $\mathcal{L}$ is said to be \emph{a right continuous} \emph{filtration}
\index{right continuous filtration}.
\end{defn}
$\square$
\begin{lem}
\label{Lem. Right-limit extesnion of filtration is right continuous}
\textbf{\emph{(Right-limit extension of a filtration is right continuous).}}
In the notations of Definition \ref{Def. Right -limit extsion and right continuity of filtration},
we have $(\mathcal{L}^{+})^{+}=\mathcal{L}^{+}$. In words, the right-limit
extension of the filtration $\mathcal{L}\equiv\{L^{(t)}:t\in Q\}$
is right continuous.
\end{lem}
\begin{proof}
We will give the proof only for the case where $\overline{Q}=[0,\infty)$,
the proof for the case where $\overline{Q}\equiv[0,a]$ being similar.
To that end, let $t\in\overline{Q}$ be arbitrary. Then 
\[
(L^{(t+)+})\equiv\bigcap\{L^{(s+)}:s\in\overline{Q}\cap(t,\infty)\}
\]
\[
\equiv\bigcap\{\bigcap\{L^{(u)}:u\in Q\cap(s,\infty)\}:s\in\overline{Q}\cap(t,\infty)\}
\]
\[
=\bigcap\{L^{(u)}:u\in Q\cap(t,\infty)\}\equiv L^{(t+)},
\]
where the third equality is because $u\in Q\cap(t,\infty)$ iff $u\in Q\cap(s,\infty)$
for some $s\in Q\cap(t,\infty)$, thanks to the assumption that $Q$
is dense in $\overline{Q}$. 
\end{proof}

\section{\label{Sec. Stopping time}Stopping Times}
\begin{defn}
\label{Def. r.r.v.-with-values in a subset of R} \textbf{(r.r.v.
with values in a subset of $R$). }Let $A$ denote an arbitrary nonempty
subset of $R$. We say that a r.r.v. $\eta$ has values in the subset
$A$ if $(\eta\in A)$ is a full set. $\square$
\end{defn}
\begin{lem}
\label{Lem. r.r.v. wth values in an increasing sequence in R} \textbf{\emph{(r.r.v.
with values in an increasing sequence in $R$).}} Let the subset $A\equiv\{t_{0},t_{1},\cdots\}\subset R$
be arbitrary such that $t_{n-1}<t_{n}$ for each $n\geq1$. Then a
r.r.v. $\eta$ has values in A iff \emph{(i)} $(\eta=t_{n})$ is measurable
for each $n\geq0$, and \emph{(ii)} $\sum_{n=1}^{\infty}P(\eta=t_{n})=1$. 
\end{lem}
\begin{proof}
Recall Definition \ref{Def. Regular =000026 Continuity Pts of Measurable Func}
of regular points of a real-valued $\mathrm{measurable}$ function.
Suppose the r.r.v. $\eta$ has values in $A$. For convenience, write
$t_{-1}\equiv t_{0}-(t_{1}-t_{0})$. Consider each $t_{n}\in A$ with
$n\geq0$. Write $\Delta_{n}\equiv(t_{n}-t_{n-1})\wedge(t_{n+1}-t_{n})>0$.
Then there exist regular points $t,s$ of the r.r.v. $\eta$ such
that 
\[
t_{n-1}<t_{n}-\Delta_{n}<s<t_{n}<t<t+\Delta_{n}<t_{n+1}.
\]
 Then $(\eta=t_{n})=(\eta\leq t)(\eta\leq s)^{c}(\eta\in A)$. Since
$(\eta\leq t)$, $(\eta\leq s)$ , and $(\eta\in A)$ are measurable
subsets, it follows that the set $(\eta=t_{n})$ is measurable. At
the same time $P(\eta\leq t_{m})\uparrow1$ as $m\rightarrow\infty$
since $\eta$ is a r.r.v. Hence
\[
\sum_{n=1}^{m}P(\eta=t_{n})=P(\eta\leq t_{m})\uparrow1
\]
as $m\rightarrow\infty$. In other words, $\sum_{n=1}^{\infty}P(\eta=t_{n})=1$.
Thus we have proved that if the r.r.v. $\eta$ has values in $A$
then Conditions (i) and (ii) holds. The converse is trivial.
\end{proof}
\begin{defn}
\textbf{\label{Def. Stopping times and simple stoppingg times.} (Stopping
time, space of integrable observables at a stopping time, and simple
stopping time). }Let $Q$ denote an arbitrary nonempty subset of $R$.
Let $\mathcal{L}$ be an arbitrary right continuous filtration with
time parameter set $Q$. Then a r.r.v. $\tau$ with values in $Q$
is called a \index{stopping time}\emph{stopping time} relative to
the filtration\emph{ $\mathcal{L}$} if 
\begin{equation}
(\tau\leq t)\in L^{(t)}\label{eq:temp-221}
\end{equation}
for each regular point $t\in Q$ of the r.r.v. $\tau$. We will omit
the reference to $\mathcal{L}$ when it is understood from context,
and simply say that $\tau$ is a stopping time. Each r.v. relative
to the probability subspace 
\[
L^{(\tau)}\equiv\{Y\in L:Y1_{(\tau\leq t)}\in L^{(t)}\mbox{ for each regular point \ensuremath{t\in Q} of }\tau\}
\]
is said to be \index{r.v. observable at stopping time}\emph{ observable
at the stopping time} $\tau$. Each member of $L^{(\tau)}$ is called
integrable observable \emph{at the stopping time} $\tau$. 

Let $X:Q\times\Omega\rightarrow S$ be an arbitrary stochastic process
adapted to the filtration $\mathcal{L}$. Define the function $X_{\tau}$
by 
\[
domain(X_{\tau})\equiv\{\omega\in domain(\tau):(\tau(\omega),\omega)\in domain(X)\},
\]
and by
\begin{equation}
X_{\tau}(\omega)\equiv X(\tau(\omega),\omega)\label{eq:temp-516-1-1}
\end{equation}
for each $\omega\in domain(X_{\tau})$. Then the function $X_{\tau}$
is called the \emph{observable of the process $X$ at the stopping
time $\tau$\index{observable of a process at a simple stopping time}}.
In general, $X_{\tau}$ need not be a well defined r.v. We will need
to prove that $X_{\tau}$ is a well defined r.v. in each application
before using it as such.

A stopping time $\tau$ with values in some discrete finite subset
of $Q$, is called a\emph{ simple stopping time}\index{simple stopping time}.
$\square$
\end{defn}
We leave it as an exercise to verify that $L^{(\tau)}$ is indeed
a probability subspace. A trivial example of a stopping time is a
deterministic time $\tau\equiv s$, where $s\in Q$ is arbitrary.

The next lemma generalizes the defining equality \ref{eq:temp-221}
and will be convenient.
\begin{lem}
\textbf{\emph{\label{Lem. basics of stopping times}(Basic properties
of stopping times). }}Suppose $Q=[0,1]$ or $Q=[0,\infty)$. $\mathcal{L}$
be an arbitrary right continuous filtration with time parameter set
$Q$. Let $\tau$ is a stopping time relative to the filtration $\mathcal{L}$.
Let $t\in Q$ be an arbitrary regular point of the r.r.v. $\tau$.
Then $(\tau<t),(\tau=t)\in L^{(t)}$.
\end{lem}
\begin{proof}
Let $(s_{k})_{k=1,2,\cdots}$ be an increasing sequence of regular
points in $Q$ of $\tau$ such that $s_{k}\uparrow t$ and such that
$P(\tau\leq s_{k})\uparrow P(\tau\leq t)$. In other words $E|1_{\tau\leq s(k)}-1_{\tau<t}|\rightarrow0$.
Since $\tau$ is a stopping time relative to a filtration $\mathcal{L}$,
we have $1_{\tau\leq s(k)}\in L^{(s(k))}\subset L^{(t)}$ for each
$k\geq1$. Hence $1_{\tau<t}\in L^{(t)}$ and $1_{\tau=t}=1_{\tau\leq t}-1_{\tau<t}\in L^{(t)}$.
Equivalently, $(\tau<t),(\tau=t)\in L^{(t)}$.
\end{proof}
\begin{defn}
\textbf{\label{Def. Specialization to a discrete parameter set }(Specialization
to a discrete parameter set). }In the remainder of this section, assume
that the parameter set $Q\equiv\{0,\Delta,2\Delta,\cdots\}$ is equally
spaced, with some fixed $\Delta>0$, and let $\mathcal{L}\equiv\{L^{(t)}:t\in Q\}$
be an arbitrary, but fixed, filtration in $(\Omega,L,E)$ with parameter
$Q$. Note that the filtration is then trivially right continuous.
\end{defn}
$\square$
\begin{prop}
\label{Prop. simple stopping times}\textbf{\emph{ (Basic properties
of stopping times, discrete case).}} Let $\tau$ and $\tau'$ be stopping
times with values in $Q\equiv\{0,\Delta,2\Delta,\cdots\}$, relative
to the filtration \emph{$\mathcal{L}$}. For each $n\geq-1$, write
$t_{n}\equiv n\Delta$ for convenience. Then the following holds.

1. Let $\eta$ be a r.r.v. with values in $Q$ . Then $\eta$ is a
stopping time iff $(\eta=t_{n})\in L^{(t(n))}$ for each $n\geq0$. 

2. $\tau\wedge\tau'$, $\tau\vee\tau'$ are stopping times.

3. If $\tau\leq\tau'$ then $L^{(\tau)}\subset L^{(\tau')}$.

4. Let $X:Q\times\Omega\rightarrow S$ be an arbitrary stochastic
process adapted to the filtration $\mathcal{L}$. Then $X_{\tau}$
is a well defined r.v. on the probability space $(\Omega,L^{(\tau)},E)$. 
\end{prop}
\begin{proof}
1. By Lemma \ref{Lem. r.r.v. wth values in an increasing sequence in R},
the set $(\eta=t_{n})$ is measurable for each $n\geq0$, and $\sum_{n=1}^{\infty}P(\eta=t_{n})=1$.
Suppose $\eta$ is a stopping time. Let $n\geq0$ be arbitrary. Then
$(\eta\leq t_{n})\in L^{(t(n))}$. Moreover, if $n\geq1$, then $(\eta\leq t_{n-1})^{c}\in L^{(t(n-1))}\subset L^{(t(n))}$.
If $n=0$, then $(\eta\leq t_{n-1})^{c}=(\eta\geq0)$ is a full set,
whence $(\eta\geq0)\in L^{(t(n))}$. Combining, we see that $(\eta=t_{n})=(\eta\leq t_{n})(\eta\leq t_{n-1})^{c}\in L^{(t(n))}.$
We have proved the ``only if'' part of Assertion 1. 

Conversely, suppose $(\eta=t_{n})\in L^{(t(n))}$ for each $n\geq0$.
Let $t\in Q$ be arbitrary. Then $t=t_{m}$ for some $m\geq0$. Hence
$(\eta\leq t)=\bigcup_{n=0}^{m}(\eta=t_{n})$, where, by assumption,
$(\eta=t_{n})\in L^{(t(n))}\subset L^{(t(m))}$ for each $n=0,\cdots,m$.
Thus we see that $(\eta\leq t)$, where $t\in Q$ is arbitrary. We
conclude that $\eta$ is a stopping time.

2. Let $t=\in Q$ be arbitrary. Then
\[
(\tau\wedge\tau'\leq t)=(\tau\leq t)\cup(\tau'\leq t)\in L^{(t)},
\]
and
\[
(\tau\vee\tau'\leq t)=(\tau\leq t)(\tau'\leq t)\in L^{(t)}.
\]
Thus $\tau\wedge\tau'$ and $\tau\vee\tau'$ are stopping times.

3. Let $Y\in L^{(\tau)}$ be arbitrary. Consider each $t\in Q$. Then,
since $\tau\leq\tau'$, 
\[
Y1_{(\tau'\leq t)}=\sum_{s\in Q}Y1_{(\tau=s)}1_{(\tau'\leq t)}=\sum_{s\in[0,t]Q}Y1_{(\tau=s)}1_{(\tau'\leq t)}\in L^{(t)}.
\]
Thus $Y\in L^{(\tau')},$ where $Y\in L^{(\tau)}$ is arbitrary. We
conclude that $L^{(\tau)}\subset L^{(\tau')}$. 

4. Let $X:Q\times\Omega\rightarrow S$ be an arbitrary stochastic
process adapted to the filtration $\mathcal{L}$. Define the full
sets $A\equiv\bigcap_{n=0}^{\infty}domain(X_{t})$ and $B\equiv\bigcup_{n=0}^{\infty}(\tau=t_{n})$.
Consider each $\omega\in AB$. Then $(\tau(\omega),\omega)=(t_{n},\omega)\in domain(X)$
on $(\tau=t_{n})$ for each $n\geq0$. In short, $X_{\tau}$ is defined
and is equal to the r.v. $X_{t(n)}$ on $(\tau=t_{n})$, for each
$n\geq0$. Since $\bigcup_{n=0}^{\infty}(\tau=t_{n})$ is a full set,
the function $X_{\tau}$ is therefore a r.v. according to Proposition
\ref{Prop. Basing seq of measurable functions on measurable partition}. 
\end{proof}
Simple first exit times from a time-varying neighborhood, introduced
next, are examples of simple stopping times. 
\begin{defn}
\label{Def. Simple First Exit time} \textbf{(Simple first exit time).}
Let $Q'\equiv\{s_{0},\cdots,s_{n}\}$ be a finite subset of $Q\equiv\{0,\Delta,2\Delta,\cdots\}$,
where $(s_{0},\cdots,s_{n})$ is an increasing sequence. Let $\mathcal{L}\equiv\{L^{(t)}:t\in Q\}$
be a filtration.

1. Let $x:Q'\rightarrow S$ be an arbitrary function. Let $b:Q'\rightarrow(0,\infty)$
be an arbitrary function such that, for each $t,r,s\in Q$' , we have
$b(s)\leq d(x_{t},x_{r})$ or $b(s)>d(x_{t},x_{r})$. Let $t\in Q'$
be arbitrary. Define 
\[
\eta_{t,b,Q'}(x)\equiv\sum_{r\in Q';t\leq r}r1_{(d(x(t),x(r))>b(r))}\prod_{s\in Q';t\leq s<r}1_{(d(x(t),x(s))\leq b(s))}
\]
\begin{equation}
+s_{n}\prod_{s\in Q';t\leq s}1_{(d(x(s),x(t))\leq b(s))}.\label{eq:temp-142-2-1}
\end{equation}
In words, $\eta_{t,b,Q'}(x)$ is the first time $r\in[t,s_{n}]Q'$
such that $x_{r}$ is at a distance greater than $b(r)$ from the
initial position $x_{t}$, with $\eta_{t,b,Q'}(x)$ set to the final
time $s_{n}\in Q'$ if no such $r$ exists. Then $\eta_{t,b,Q'}(x)$
is called the \emph{simple first exit tim}e\index{simple first exit time}
for the function $x|[t,s_{n}]Q'$ to exit the time-varying $b$-neighborhood
of $x_{t}$. In the special case where $b(r)=\alpha$ for each $r\in Q'$
for some constant $\alpha>0$, we will write simply $\eta_{t,\alpha,Q'}(x)$
for $\eta_{t,b,Q'}(x)$.

2. More generally, let $X:Q'\times\Omega\rightarrow S$ be an arbitrary
process adapted to the filtration $\mathcal{L}$. Let $b:Q'\rightarrow(0,\infty)$
be an arbitrary function such that, for each $t,r,s\in Q$' , the
real number $b(s)$ is a regular point for the r.r.v. $d(X_{t},X_{r})$.
Let $t\in Q'$ be arbitrary. Define the r.r.v. $\eta_{t,b,Q'}(X)$
on $\Omega$ defined by
\[
\eta_{t,b,Q'}(X)\equiv\sum_{r\in Q';t\leq r}r1_{(d(X(t),X(r))>b(r))}\prod_{s\in Q';t\leq s<r}1_{(d(X(t),X(s))\leq b(s))}
\]
\begin{equation}
+s_{n}\prod_{s\in Q';t\leq s}1_{(d(X(s),X(t))\leq b(s))}\label{eq:temp-142-2}
\end{equation}
is a r.r.v. called the \emph{simple first exit tim}e\index{simple first exit time}
for the process $X|[t,s_{n}]Q'$ to exit the time-varying $b$-neighborhood
of $X_{t}$. When there is little risk of confusion as to the identity
of the process $X$, we will omit the reference to $X$, write $\eta_{t,b,Q'}$
for $\eta_{t,b,Q'}(X)$, and abuse notations by writing $\eta_{t,b,Q'}(\omega)$
for $\eta_{t,b,Q'}(X(\omega))$, for each $\omega\in\bigcap_{u\in Q'}domain(X_{u})$.
\end{defn}
$\square$

The next proposition verifies that $\eta_{t,b,Q'}(X)$ is a simple
stopping time relative to the filtration $\mathcal{L}$. It also proves
some simple properties that are intuitively obvious when described
in words.
\begin{prop}
\label{Prop. Basics of simple Exit times} \textbf{\emph{(Basic properties
of simple first exit times). }}Let $Q'\equiv\{s_{0},\cdots,s_{n}\}$
be a finite subset of $Q$, where $(s_{0},\cdots,s_{n})$ is an increasing
sequence. Use the assumptions and notations in Part 2 of Definition
\ref{Def. Simple First Exit time}. Let $t,s\in Q'\equiv\{s_{0},\cdots,s_{n}\}$
be arbitrary. Let $\omega\in domain(\eta_{t,b,Q'})$ be arbitrary.
Then the following holds.

1. $t\leq\eta_{t,b,Q'}(\omega)\leq s_{n}$. 

2. The r.r.v. $\eta_{t,b,Q'}$ is a simple stopping time relative
to the filtration $\mathcal{L}$.

3. If $\eta_{t,b,Q'}(\omega)<t_{n}$ then $d(X(t,\omega),X(\eta_{t,b,Q'}(\omega),\omega))>b(\eta_{t,b,Q'}(\omega)).$
In words, if the simple first exit time occurs before the final time,
then the sample path exits successfully at the simple first exit time.

4. If $t\leq s<\eta_{t,b,Q'}(\omega)$, then $d(X(t,\omega),X(s,\omega))\leq b(s)$.
In words, before the simple first exit time, the sample path remains
in the $b$-neighborhood. Moreover, if 
\[
d(X(\eta_{t,b,Q'}(\omega),\omega),X(t,\omega))\leq b(\eta_{t,b,Q'}(\omega)),
\]
then
\[
d(X(s,\omega),X(t,\omega))\leq b(s)
\]
for each $s\in Q'$ with $t\leq s$. In words, if the sample path
is in the $b$-neighborhood at the simple first exit time, then it
is in the $b$-neighborhood at any time prior to the simple first
exit time.

Conversely, if $r\in[t,s_{n})Q'$ is such that $d(X(t,\omega),X(s,\omega))\leq b(s)$
for each $s\in(t,r]Q'$, then $r<\eta_{t,b,Q'}(\omega)$. In words,
if the sample path stays within the the $b$-neighborhood up to and
including a certain time, then the simple first exit time can come
only after that time.

5. Suppose $s_{0}=s_{k(0)}<s_{k(1)}<\cdots<s{}_{k(p)}=s_{n}$ is a
subsequence of $s_{0}<s_{1}<\cdots<s_{n}$. Define $Q''\equiv\{s_{k(1)},\cdots,s{}_{k(p)}\}$.
Let $t\in Q''\subset Q'$ be arbitrary. Then $\eta_{t,b,Q'}\leq\eta_{t,b,Q''}$.
In other words, if the process $X$ is sampled at more time points,
then the simple first exit time can occur no later. 
\end{prop}
\begin{proof}
By hypothesis, the process $X$ is adapted to the filtration $\mathcal{L}$. 

1. Assertion 1 is obvious from the defining equality \ref{eq:temp-142-2}. 

2. By equality \ref{eq:temp-142-2}, for each $r\in\{t,\cdots,s_{n-1}\},$
we have
\[
(\eta_{t,b,Q'}=r)=(d(X_{r},X_{t})>b(r))\bigcap_{s\in Q';t\leq s<r}(d(X_{s},X_{t})\leq b(s))
\]
\begin{equation}
\in L^{(r)}\subset L^{(s(n))}.\label{eq:temp-142}
\end{equation}
Consequently,
\[
(\eta_{t,b,Q'}=s_{n})=\bigcap_{r\in Q';r<t(n)}(\eta_{t,b,Q'}=r)^{c}\in L^{(s(n))}.
\]
Hence $\eta_{t,b,Q'}$ is a simple stopping time relative to $\mathcal{L}'$
with values in $Q'$, according to Proposition \ref{Prop. simple stopping times}. 

3. Assertion 3 is obvious from the defining equality \ref{eq:temp-142-2}. 

4. Suppose $t<s<r\equiv\eta_{t,b,Q'}(\omega)$. Then
\begin{equation}
d(X(t,\omega),X(s,\omega))\leq b(s)\label{eq:temp-521}
\end{equation}
by equality \ref{eq:temp-142}. The last inequality is trivially satisfied
if $t=s$. Hence if $r\equiv\eta_{t,b,Q'}(\omega)=s_{n}$ with $d(X(t,\omega),X(r,\omega))\leq b(r)$
then inequality \ref{eq:temp-521} holds for each $s\in Q'$ with
$t\leq s\leq r$.

Conversely, suppose $r\in Q'$ is such that $t\leq r<s_{n}$ and such
that
\[
d(X(t,\omega),X(s,\omega))\leq b(s)
\]
for each $s\in(t,r]Q'$. Suppose $s\equiv\eta_{t,b,Q'}(\omega)\leq r<s_{n}$.
Then $d(X(t,\omega),X(s,\omega))>b(s)$ by Assertion 3, a contradiction.
Hence $\eta_{t,b,Q'}(\omega)>r$. Assertion 4 is verified.

5.Let $t\in Q''\subset Q'$ be arbitrary. Suppose, for the sake of
a contradiction, that $s\equiv\eta_{t,b,Q''}(\omega)<\eta_{t,b,Q'}(\omega)\leq s_{n}$.
Then $t<s$ and $s\in Q''\subset Q'$. Hence, by Assertion 4 applied
to the time $s$ and to the simple first exit time $\eta_{t,b,Q'}$,
we have 
\[
d(X(t,\omega),X(s,\omega))\leq b(s).
\]
On the other hand, by Assertion 3 applied to the time $s$ and to
the simple first exit time $\eta_{t,b,Q''}$, we have 
\[
d(X(t,\omega),X(s,\omega))>b(s),
\]
 a contradiction. Hence $\eta_{t,b,Q''}(\omega)\geq\eta_{t,b,Q'}(\omega)$.
Assertion 5 is proved.
\end{proof}

\section{Martingales}
\begin{defn}
\label{Def. Martingale and Suibmartingale} \textbf{(Martingale and
Submartingale). }Let $Q$ be an arbitrary nonempty subset of $R$.
Let $\mathcal{L}\equiv\{L^{(t)}:t\in Q\}$ be an arbitrary right continuous
filtration in $(\Omega,L,E)$. Let $X:Q\times\Omega\rightarrow R$
be a stochastic process such that $X_{t}\in L^{(t)}$ for each $t\in Q$. 

1. The process $X$ is called a \index{martingale}\emph{martingale}
relative to $\mathcal{L}$ if, for each $t,s\in Q$ with $t\leq s$,
we have \emph{$EZX_{t}=EZX_{s}$} for each indicator \emph{$Z\in L^{(t)}$}.
Accordingly to Definition \ref{Def: Conditional Expectation}, the
last condition is equivalent to $E(X_{s}|L^{(t)})=X_{t}$\emph{ }for
each $t,s\in Q$ with $t\leq s$\emph{. }

2. The process $X$ is called a \emph{\index{wide-sense submartingale}
wide-sense submartingale} relative to $\mathcal{L}$ if, for each
$t,s\in Q$ with $t\leq s$, we have \emph{$EZX_{t}\leq EZX_{s}$}
for each indicator \emph{$Z\in L^{(t)}$}. If, in addition, $E(X_{s}|L^{(t)})$
exists for each $t,s\in Q$ with $t\leq s$, then $X$ is called a
\emph{\index{submartingale}submartingale} relative to $\mathcal{L}$.

3. The process $X$ is called a \emph{\index{wide-sense supermartingale}
wide-sense supermartingale} relative to $\mathcal{L}$ if, for each
$t,s\in Q$ with $t\leq s$, we have \emph{$EZX_{t}\geq EZX_{s}$}
for each indicator \emph{$Z\in L^{(t)}$}. If, in addition, $E(X_{s}|L^{(t)})$
exists for each $t,s\in Q$ with $t\leq s$, then $X$ is called a
\emph{\index{submartingale}supermartingale} relative to $\mathcal{L}$.

When there is little risk of confusion, we will omit the explicit
reference to the given filtration $\mathcal{L}$.
\end{defn}
$\square$

Clearly a submartingale is also a wide-sense submartingale. The two
notions are classically equivalent because, classically, with the
benefit of the principle of infinite search, the conditional expectation
always exists. Hence, any result that we prove for wide-sense submartingales
holds classically for submartingales. 
\begin{prop}
\label{Prop. Martingale basics}\textbf{\emph{(Martingale basics).}}
Let $X:Q\times\Omega\rightarrow R$ be an arbitrary process adapted
to the right continuous filtration $\mathcal{L}\equiv\{L^{(t)}:t\in Q\}$.
Unless otherwise specified, all martingales and wide-sense submartingales
are relative to the filtration $\mathcal{L}$. Then the following
holds.

1. The process $X$ is a martingale iff it is both a wide-sense submartingale
and a wide-sense supermartingale. 

2. The process $X$ is a wide-sense supermartingale iff $-X$ is a
wide-sense submartingale.

3. The expectation $EX_{t}$ is constant for $t\in Q$ if $X$ is
a martingale. Moreover, $EX_{t}$ is nondecreasing in $t$ if $X$
is a wide-sense submartingale. 

4. Suppose $X$ is a martingale. Then $|X|$ is a wide-sense submartingale.
In particular, $E|X_{t}|$ is nondecreasing in $t\in Q$.

5. Suppose $X$ is a martingale. Let $a\in Q$ be arbitrary. Then
the family $\{X_{t}:t\in(-\infty,a]Q\}$ is uniformly integrable.

6. Let $\mathcal{\overline{L}}\equiv\{\overline{L}^{(t)}:t\in Q\}$
be an arbitrary filtration such that $L^{(t)}\subset\overline{L}^{(t)}$
for each $t\in Q$. Suppose $X$ is a wide-sense submartingale relative
to the filtration $\mathcal{\overline{L}}$. Then $X$ is a wide-sense
submartingale relative to the filtration $\mathcal{\mathcal{L}}$.
The same assertion holds for martingales.

7. Suppose $X$ is a wide-sense submartingale relative to the filtration
$\mathcal{L}$. Then it is a a wide-sense submartingale relative to
the natural filtration $\mathcal{L}_{X}\equiv\{L^{(X,t)}:t\in Q\}$
of the process $X$.
\end{prop}
\begin{proof}
1. Assertions 1-3 being trivial, we will prove Assertions 4-7 only. 

2. To that end, let $t,s\in Q$ with $t\leq s$ be arbitrary. Let
the indicator \emph{$Z\in L^{(t)}$ }and the real number $\varepsilon>0$
be arbitrary\emph{. }Then\emph{
\[
E(|X_{s}|Z;X_{t}>\varepsilon)\geq E(X_{s}Z;X_{t}>\varepsilon)=E(X_{t}Z;X_{t}>\varepsilon)=E(|X_{t}|Z;X_{t}>\varepsilon),
\]
}where the first equality is from the definition of a martingale.
Since $-X$ is also a martingale, we have similarly\emph{
\[
E(|X_{s}|Z;X_{t}<-\varepsilon)\geq E(-X_{s}Z;X_{t}<-\varepsilon)=E(-X_{t}Z;X_{t}<-\varepsilon)=E(|X_{t}|Z;X_{t}<-\varepsilon).
\]
}Adding the last two displayed inequalities, we obtain\emph{
\[
E(|X_{s}|;Z)\geq E(|X_{s}|Z;X_{t}>\varepsilon)+E(|X_{s}|Z;X_{t}<-\varepsilon)
\]
\[
\geq E(|X_{t}|Z;X_{t}>\varepsilon)+E(|X_{t}|Z;X_{t}<-\varepsilon)=E(|X_{t}|Z)-E(|X_{t}|Z;|X_{t}|\leq\varepsilon).
\]
}Since 
\[
E(|X_{t}|Z;|X_{t}|\leq\varepsilon)\leq E(|X_{t}|;|X_{t}|\leq\varepsilon)\rightarrow0
\]
as $\varepsilon\rightarrow0$, we conclude that 
\[
E(|X_{s}|;Z)\geq E(|X_{t}|;Z),
\]
where $t,s\in Q$ with $t\leq s$ and the indicator \emph{$Z\in L^{(t)}$}
are arbitrary. Thus the process $|X|$ is a wide-sense submartingale.
Assertion 4 is proved.

2. Suppose $X$ is a martingale. Consider each $a\in Q$ be arbitrary.
Let $t\in Q$ be arbitrary with $t\leq a$, and let $\varepsilon>0$
be arbitrary. Then, since $X_{a}$ is integrable, there exists $\delta\equiv\delta_{X(a)}(\varepsilon)>0$
so small that $E|X_{a}|1_{A}<\varepsilon$ for each $\mathrm{measurable}$
set $A$ with $P(A)<\delta$. Now let $\gamma>\beta(\varepsilon)\equiv E|X_{a}|\delta^{-1}$
be arbitrary. Then, by Chebychev's inequality, 
\[
P(|X_{t}|>\gamma)\leq E|X_{t}|\gamma^{-1}\leq E|X_{a}|\gamma^{-1}<\delta,
\]
where the second inequality is because $|X|$ is a wide-sense submartingale
by Assertion 4. Hence
\[
E|X_{t}|1_{(X(t)>\gamma)}\leq E|X_{a}|1_{(X(t)>\gamma)}<\varepsilon,
\]
where the first inequality is because $|X|$ is a wide-sense submartingale.
Since $t\in Q$ is arbitrary with $t\leq a$, we conclude that the
family $\{X_{t}:t\in(-\infty,a]Q\}$ is uniformly integrable, with
a simple modulus of uniform integrability $\beta$. Assertion 5 has
been verified.

3. To prove Assertion 6, assume that the process $X$ is a wide-sense
submartingale relative to the filtration $\mathcal{\overline{L}}$.
Let $t,s\in Q$ with $t\leq s$ be arbitrary. Consider each indicator
\emph{$Z\in L^{(t)}$. }Then $Z\in\overline{L}^{(t)}$ by the assumption
on $\mathcal{\overline{L}}$. Hence \emph{ $EZX_{t}\leq EZX_{s}$,
}where the indicator $Z\in L^{(t)}$ is arbitrary. Thus $X$ is a
wide-sense submartingale relative to $\mathcal{L}$. The proof for
martingales is similar. Assertion 6 is proved.

4. It remains to prove Assertion 7. To that end, suppose $X$ is a
wide-sense submartingale relative to $\mathcal{L}$. Note that, for
each $t\in Q$, we have $X_{r}\in L^{(t)}$ for each $r\in[0,t]Q$.
Hence $L^{(X,t)}\equiv L(X_{r}:r\in[0,t]Q)\subset L^{(t)}$. Hence
Assertion 6 implies that $X$ is a wide-sense submartingale relative
to $\mathcal{L}_{X}$. The proof for martingales is similar. Assertion
7 and the proposition are proved.
\end{proof}
\begin{defn}
\textbf{\label{Def. Specialization to a discrete parameter set -1}(Specialization
to a discrete parameter set). }In the remainder of this section, unless
otherwise specified, assume the parameter set $Q\equiv\{0,\Delta,2\Delta,\cdots\}$
with some fixed $\Delta>0$, and let $\mathcal{L}\equiv\{L^{(t)}:t\in Q\}$
be an arbitrary, but fixed, filtration in $(\Omega,L,E)$ with parameter
$Q$. Note that the filtration $\mathcal{L}$ is then trivially right
continuous. For ease of notations, we will assume, without loss of
generality, by a change of units if necessary, that $\Delta=1$. 
\end{defn}
$\square$

If a martingale $X_{t}$ is used to model a gambler's fortune at the
current time $t$, then the conditional expectation of said fortune
at a later time $s$, given all information up to and including the
current time $t$, is exactly his or her current fortune. Thus a martingale
$X$ is a model for a fair game of chance. Similarly, a submartingale
can be used to model a favorable game. 
\begin{thm}
\label{Thm. Doob decomposition assuming conditional expectation}
\textbf{\emph{(Doob decomposition).}} Let $Y:\{0,1,2,\cdots\}\times(\Omega,L,E)\rightarrow R$
be a process which is adapted to the filtration $\mathcal{L}\equiv\{L^{(n)}:n\geq0\}$.
Suppose the conditional expectation $E(Y_{m}|L^{(n)})$ exists for
each $m,n\geq0$ with $n\leq m$. For each $n\geq0$, define 
\begin{equation}
X_{n}\equiv Y_{0}+\sum_{k=1}^{n}(Y_{k}-E(Y_{k}|L^{(k-1)})).\label{eq:temp-527}
\end{equation}
and
\begin{equation}
A_{n}\equiv\sum_{k=1}^{n}(E(Y_{k}|L^{(k-1)})-Y_{k-1}),\label{eq:temp-144}
\end{equation}
where an empty sum is by convention equal to $0$. Then $X:\{0,1,2,\cdots\}\times\Omega\rightarrow R$
is a martingale relative to the filtration $\mathcal{L}$. Moreover,
$A_{n}\in L^{(n-1)}$ and $Y_{n}=X_{n}+A_{n}$ for each $n\geq1$.
\end{thm}
\begin{proof}
From the defining equality \ref{eq:temp-527}, we see that $X_{n}\in L^{(n)}$
for each $n\geq1$. Hence the process $X:\{0,1,2,\cdots\}\times\Omega\rightarrow R$
is adapted to the filtration $\mathcal{L}$. Let $m>n\geq1$ be arbitrary.
Then 
\[
E(X_{m}|L^{(n)})=E(\{X_{n}+\sum_{k=n+1}^{m}(Y_{k}-E(Y_{k}|L^{(k-1)}))\}|L^{(n)})
\]
\[
=X_{n}+\sum_{k=n+1}^{m}\{E(Y_{k}|L^{(n)})-E(E(Y_{k}|L^{(k-1)})|L^{(n)})\}
\]
\[
=X_{n}+\sum_{k=n+1}^{m}\{E(Y_{k}|L^{(n)})-E(Y_{k}|L^{(n)})\}=X_{n},
\]
where we used basic properties of conditional expectations in Proposition
\ref{Prop. Basics of Conditional expectations}. Thus the process
$X$ is a martingale relative to the filtration $\mathcal{L}$. Moreover,
$A_{n}\in L^{(n-1)}$ because all the summands in the defining equality
\ref{eq:temp-144} are members of $L^{(n-1)}$. 
\end{proof}
Intuitively, Theorem \ref{Thm. Doob decomposition assuming conditional expectation}
says that a multi-round game $Y$ can be turned into a fair game $X$
by charging a fair price determined at each round as the conditional
expectation of payoff at the next round, with the cumulative cost
of entry equal to $A_{n}$ by the time $n$. 

The next theorem of Doob and its corollary are key to the analysis
of martingales. It proves that, under reasonable conditions, a fair
game can never be turned to a favorable one by sampling at a sequence
of stopping times, or by stopping at some stopping time which which
cannot see the future. The reader can look up ``gambler's ruin''
in the literature for a counterexample where said reasonable conditions
is not assumed, where a fair coin tossing game can be turned into
an almost sure win by stopping when and only when the gambler is ahead
by one dollar. This latter strategy sounds intriguing except for the
lamentable fact that, to achieve almost sure winning against a house
with infinite capital, the strategy would require the gambler to stay
in the game for unbounded number of rounds and to have infinite capital
to avoid bankruptcy.

The next theorem and its proof are essentially restatements of parts
of Theorems 9.3.3 and 9.3.4 in \cite{Chung68}, except that, for the
case of wide-sense submartingales, we add a condition to make the
theorem constructive. 
\begin{thm}
\label{Thm.Doob's  Optional sstopping theorem}\textbf{\emph{ (Doob's
optional sampling theorem).}} Let 
\[
X:\{0,1,2,\cdots\}\times(\Omega,L,E)\rightarrow R
\]
be a wide-sense submartingale relative to a filtration $\mathcal{L}\equiv\{L^{(k)}:k\geq0\}$.
Let $\overline{\tau}\equiv(\tau_{n})_{n=1,2,\cdots}$ be a nondecreasing
sequence of stopping times with values in $\{0,1,2,\cdots\}$ relative
to the filtration $\mathcal{L}$. Define the function $X_{\overline{\tau}}:\{0,1,2,\cdots\}\times(\Omega,L,E)\rightarrow R$
by $X_{\overline{\tau},n}\equiv X_{\tau(n)}$ for each $n\geq0$.
Suppose one of the following three conditions holds.

\emph{(i)} The function $X_{\tau(n)}$ is an integrable r.r.v. for
each $n\geq0$, and the family $\{X_{n}:n\geq0\}$ of r.r.v.'s is
uniformly integrable. 

\emph{(ii)} For each $m\geq1$, there exists some $M_{m}\geq0$ such
that $\tau_{m}\leq M_{m}$.

\emph{(iii)} The given process $X$ is a martingale, and the family
$\{X_{n}:n\geq0\}$ of r.r.v.'s is uniformly integrable. 

Then $X_{\overline{\tau}}$ is a wide-sense submartingale relative
to the filtration $\mathcal{L}^{\overline{\tau}}\equiv\{L^{(\tau(n))}:n\geq0\}.$
If the given process $X$ is a martingale, then $X_{\overline{\tau}}$
is a martingale relative to the filtration $\mathcal{L}^{\overline{\tau}}$.
\end{thm}
\begin{proof}
Recall that $Q\equiv\{0,1,2,\cdots\}$. Let $m\geq n$ and the indicator
$Z\in L^{(\tau(n))}$ be arbitrary. We need to prove that the function
$X_{\tau(n)}$ is integrable, and that
\begin{equation}
E(X_{\tau(m)}Z)\geq E(X_{\tau(n)}Z).\label{eq:temp-599-2}
\end{equation}
First we will prove that $X_{\tau(n)}$ is integrable. 

1. Suppose Condition (i) holds.Then the function $X_{\tau(n)}$ is
integrable by assumption.

2. Suppose Condition (ii) holds. Then the function
\[
X_{\tau(n)}=\sum_{u=0}^{M(n)}X_{\tau(n)}1_{(\tau(n)=u)}=\sum_{u=0}^{M(n)}X_{u}1_{(\tau(n)=u)}
\]
is a finite sum of integrable r.r.v., hence $X_{\tau(n)}$ is itself
an integrable r.r.v.

3. Suppose Condition (iii) holds. Then $X$ is a martingale. Hence
$|X|$ is a wide-sense submartingale and $E|X_{t}|$ is nondecreasing
in $t\in Q$, according to Assertion 4 of Proposition \ref{Prop. Martingale basics}.
Consider each $v,v'\in Q$ with $v\leq v'$. Then it follows that
\[
\sum_{u=v}^{v'}E|X_{\tau(n)}|1_{(\tau(n)=u)}=\sum_{u=v}^{v'}E|X_{u}|1_{(\tau(n)=u)}\leq\sum_{u=v}^{v'}E|X_{v'}|1_{(\tau(n)=u)}
\]
\begin{equation}
=E|X_{v'}|1_{(v\leq\tau(n)\leq v')}\leq\alpha_{v,v'}\equiv E|X_{v'}|1_{(v\leq\tau(n))}.\label{eq:temp-603}
\end{equation}
Let $v\rightarrow\infty$. Since $\tau_{n}$ is a nonnegative r.r.v.,
we have $P(v\leq\tau_{n})\rightarrow0$. Therefore $\alpha_{v,v'}\rightarrow0$,
thanks to the uniform integrability of the family $\{X_{t}:t\in Q\}$
of r.r.v.'s. under Condition (iii). Summing up, we conclude that $\sum_{u=v}^{v'}E|X_{\tau(n)}|1_{(\tau(n)=u)}\rightarrow0$
as $v\rightarrow\infty$. Thus $\sum_{u=0}^{\infty}E|X_{\tau(n)}|1_{(\tau(n)=u)}<\infty$.
Consequently, the function $X_{\tau(n)}=\sum_{u=0}^{\infty}X_{\tau(n)}1_{(\tau(n)=u)}$
is an integrable r.r.v. 

4. Thus we see that $X_{\tau(n)}$ is an integrable r.r.v. under any
one of the three Conditions (i-iii). It remains to prove \ref{eq:temp-599-2}.
To that end, let $u,v\in Q$ be arbitrary with $u\leq v$. Then $Z1_{\tau(n)=u}\in L^{(u)}\subset L^{(v)}$.
Hence
\[
Y_{u,v}\equiv X_{v}Z1_{\tau(n)=u}\in L^{(v)}.
\]
Moreover, 
\[
EY_{u,v}1_{\tau(m)\geq v}\equiv EX_{v}Z1_{\tau(n)=u}1_{\tau(m)\geq v}
\]
\[
=EX_{v}Z1_{\tau(n)=u}1_{\tau(m)=v}+EX_{v}Z1_{\tau(n)=u}1_{\tau(m)\geq v+1}
\]
\[
=EX_{\tau(m)}Z1_{\tau(n)=u}1_{\tau(m)=v}+EX_{v}Z1_{\tau(n)=u}1_{\tau(m)\geq v+1}
\]
\[
\leq EX_{\tau(m)}Z1_{\tau(n)=u}1_{\tau(m)=v}+EX_{v+1}Z1_{\tau(n)=u}1_{\tau(m)\geq v+1},
\]
where the inequality is because the indicator 
\[
Z1_{\tau(n)=u}1_{\tau(m)\geq v+1}=Z1_{\tau(n)=u}(1-1_{\tau(m)\leq v})\in L^{(v)}
\]
 and because $X$ is, by hypothesis, a wide-sense submartingale. In
short
\begin{equation}
EY_{u,v}1_{\tau(m)\geq v}\leq EX_{\tau(m)}Z1_{\tau(n)=u}1_{\tau(m)=v}+EY_{u,v+1}1_{\tau(m)\geq v+1},\label{eq:temp-600-1}
\end{equation}
where $v\in[u,\infty)Q$ is arbitrary. Let $\kappa\geq0$ be arbitrary.
Applying inequality \ref{eq:temp-600-1} successively to $v=u,u+1,u+2,\cdots,u+\kappa$,
we obtain
\[
EY_{u,u}1_{\tau(m)\geq u}\leq EX_{\tau(m)}Z1_{\tau(n)=u}1_{\tau(m)=u}+EY_{u,u+1}1_{\tau(m)\geq u+1}
\]
\[
\leq EX_{\tau(m)}Z1_{\tau(n)=u}1_{\tau(m)=u}+EX_{\tau(m)}Z1_{\tau(n)=u}1_{\tau(m)=u+1}+EY_{u,u+2}1_{\tau(m)\geq u+2}
\]
\[
\leq\cdots
\]
\[
\leq EX_{\tau(m)}Z1_{\tau(n)=u}\sum_{v\in[u,u+\kappa]Q}1_{\tau(m)=v}+EY_{u,u+(\kappa+1)}1_{\tau(m)\geq u+(\kappa+1)}
\]
\[
=EX_{\tau(m)}Z1_{\tau(n)=u}1_{u\leq\tau(m)\leq u+\kappa}+EX_{u+(\kappa+1)}Z1_{\tau(n)=u}1_{\tau(m)\geq u+(\kappa+1)}
\]
\[
=EX_{\tau(m)}Z1_{\tau(n)=u}1_{\tau(m)\leq u+\kappa}+EX_{u+(\kappa+1)}Z1_{\tau(n)=u}1_{\tau(m)\geq u+(\kappa+1).}
\]
\[
=EZX_{\tau(m)}Z1_{\tau(n)=u}-EX_{\tau(m)}Z1_{\tau(n)=u}1_{\tau(m)\geq u+(\kappa+1)}+EX_{u+(\kappa+1)}Z1_{\tau(n)=u}1_{\tau(m)\geq u+(\kappa+1).}.
\]
\begin{equation}
\equiv EX_{\tau(m)}Z1_{\tau(n)=u}-EX_{\tau(m)}1_{A(\kappa)}+EX_{u+(\kappa+1)}1_{A(\kappa)},\label{eq:temp-601-1}
\end{equation}
where $A_{\kappa}$ is the measurable set whose indicator is $1_{A(\kappa)}\equiv Z1_{\tau(n)=u}1_{\tau(m)\geq u+(\kappa+1)}$
and whose probability is therefore bounded by 
\[
P(A_{\kappa})\leq P(\tau_{m}\geq u+(\kappa+1)).
\]
Now let $\kappa\rightarrow\infty$. Then $P(A_{\kappa})\rightarrow0$
because $X_{\tau(m)}$ is an integrable r.r.v., as proved in Steps
1-3. Consequently, the second summand on the right-hand side of inequality
\ref{eq:temp-601-1} tends to $0$. Now consider the third summand
on the right-hand side of inequality \ref{eq:temp-601-1}. Suppose
Condition (ii) holds. Then, as soon as $\kappa$ is so large that
$u+(\kappa+1)\geq M_{m}$, we have $P(A_{\kappa})=0$ , whence said
two summands vanish as $\kappa\rightarrow\infty$. Suppose, alternatively,
Condition (i) or (iii) holds, then the last summand tends to $0$,
thanks to the uniform integrability of the family $\{X_{t}:t\in[0,\infty)\}$
of r.r.v.'s guaranteed by Condition (i) or (iii). Summing up, the
second and third summand both tend to $0$ as $\kappa\rightarrow\infty$,
with only the first summand on the right-hand side of inequality \ref{eq:temp-601-1}
surviving, to yield
\[
EY_{u,u}1_{\tau(m)\geq u}\leq EX_{\tau(m)}Z1_{\tau(n)=u}.
\]
Equivalently,
\[
EX_{u}Z1_{\tau(n)=u}1_{\tau(m)\geq u}\leq EX_{\tau(m)}Z1_{\tau(n)=u}.
\]
Since $(\tau_{n}=u)\subset(\tau_{m}\geq u)$, this last inequality
simplifies to
\[
EX_{\tau(n)}Z1_{\tau(n)=u}\leq EX_{\tau(m)}Z1_{\tau(n)=u},
\]
where $u\in Q\equiv\{0,1,2,\cdots\}$ is arbitrary. Summation over
$u\in Q$ then yields the desired equality \ref{eq:temp-599-2}. Thus
$X_{\overline{\tau}}$ is a wide-sense submartingale relative to the
filtration $\mathcal{L}^{\overline{\tau}}\equiv\{L^{(\tau(n))}:n=0,1,\cdots\}.$
The first part of the conclusion of the theorem, regarding wide-sense
submartingales, has been proved.

5. Finally, suppose the given wide-sense submartingale $X$ is actually
a martingale. Then $-X$ is a wide-sense submartingale, and so by
the preceding arguments, both the processes $X_{\overline{\tau}}$
and $-X_{\overline{\tau}}$ are a wide-sense submartingale relative
to the filtration $\mathcal{L}^{\overline{\tau}}.$ Combining, we
conclude that $X_{\overline{\tau}}$ is a martingale if $X$ is a
martingale, provided that one of the three Conditions (i-iii) holds.
The theorem is proved.
\end{proof}
\begin{cor}
\label{Cor.Doob's  Optional sampling of martingales} \textbf{\emph{(Doob's
optional stopping theorem for a finite game). }}Let $n\geq1$ be arbitrary.
Write $Q'\equiv\{0,1,\cdots,n\}\equiv\{t_{0},t_{1},\cdots,t_{n}\}\subset Q$.
Let $X:Q'\times\Omega\rightarrow R$ be a process adapted to the filtration
$\mathcal{L}\equiv\{L^{(t)}:t\in Q\}$. Let $\tau$ be an arbitrary
simple stopping time relative to $\mathcal{L}$ with values in Q'.
Define the r.r.v.
\[
X_{\tau}\equiv\sum_{t\in Q'}X_{t}1_{(\tau=t)}\in L'^{(\tau)}.
\]
Define the process $X':\{0,1,2\}\times\Omega\rightarrow R$ by
\[
(X'_{0},X'_{1},X'_{2})\equiv(X_{t(0)}X_{\tau},X_{t(n)}).
\]
Define the filtration $\mathcal{L}'\equiv\{L'^{(i)}:i=0,1,2\}$ by
\[
(L'^{(0)},L'^{(1)},L'^{(2)})\equiv(L^{(t(0))},L^{(\tau)},L^{(t(n))}).
\]
Then the following holds. 

1. If the process $X$ is a wide-sense submartingale relative to $\mathcal{L}$,
then the process $X'$ is a wide-sense submartingale relative to the
filtration $\mathcal{L}'$.

2. If the process $X$ is a martingale relative to $\mathcal{L}$,
then the process $X'$ is a martingale relative to $\mathcal{L}'$.
\end{cor}
\begin{proof}
Extend the process $X$ to the parameter set $Q\equiv\{0,1,\cdots\}$
by $X_{t}\equiv X_{t\wedge n}$ for each $t\in\{0,1,\cdots\}$. Likewise
extend the filtration $\mathcal{L}$ by defining $L^{(t)}\equiv L^{(t\wedge n)}$
for each $t\in\{0,1,\cdots\}$. Trivially we can verify that the extended
process $X:\{0,1,\cdots\}\times\Omega\rightarrow R$ retains the same
property of being a martingale or wide-sense submartingale, respectively,
as the given process being a martingale or wide-sense submartingale,
relative to the extended filtration $\mathcal{L}.$ Now define a sequence
$\overline{\tau}\equiv(\tau_{0},\tau_{1},\cdots)$ of stopping times
by $\tau_{0}\equiv t_{0}$, $\tau_{1}\equiv\tau,$ and $\tau_{m}\equiv t_{n}$
for each $m\geq2$. Then it can easily be verified that the sequence
$\overline{\tau}$ satisfies Condition (ii) of Theorem \ref{Thm.Doob's  Optional sstopping theorem}.
Hence the process $X_{\overline{\tau}}$ defined in Theorem \ref{Thm.Doob's  Optional sstopping theorem}
is a martingale if $X$ is a martingale, and is a wide-sense submartingale
if $X$ is a wide-sense submartingale. Since 
\[
(X'_{0},X'_{1},X'_{2})\equiv(X_{t(0)}X_{\tau},X_{t(n)})=(X_{\tau(0)}X_{\tau(1)},X_{\tau(2)})
\]
and
\[
(L'^{(0)},L'^{(1)},L'^{(2)})\equiv(L^{(t(0))},L^{(\tau)},L^{(t(n))})=(L^{(\tau(0))},L^{(\tau(1))},L^{(\tau(2))}),
\]
the conclusion of the corollary follows.
\end{proof}

\section{Convexity and Martingale Convergence}

Next consider the a.u. convergence of martingales, and of wide-sense
submartingales in general. Suppose $X:\{1,2,\cdots\}\times\Omega\rightarrow R$
is a martingale relative to some filtration $\mathcal{L}\equiv\{L^{(n)}:n=1,2,\cdots\}$.
A classical theorem says that, if $E|X_{n}|$ is bounded as $n\rightarrow\infty$,
then $X_{n}$ converges a.u. as $n\rightarrow\infty$. The theorem
can be proved, classically, by the celebrated upcrossing inequality
of Doob, thanks to the principle of infinite search. See, for example,
\cite{Durret84}. While the upcrossing inequality is constructive,
the inference of a.u. convergence from it is not. 

As a matter of fact, the following example shows that the martingale
convergence theorem, as stated above, implies the principle of infinite
search. Let $(a_{n})_{n=1,2,\cdots}$ be a arbitrary nondecreasing
sequence in $\{0,1\}$. Let $Y$ be an arbitrary r.r.v. which takes
the value $-1$ or $+1$ with equal probabilities. For each $n\geq1$
define $X_{n}\equiv1+a_{n}Y$. Then the process $X:\{1,2,\cdots\}\times\Omega\rightarrow\{0,1,2\}$
is a martingale relative to its natural filtration, with $E|X_{n}|=EX_{n}=1$
for each $n\geq1$. Suppose $X_{n}\rightarrow X$ a.u. Then there
exists $b\in(0,1)$ such that the set $(X<b)$ is $\mathrm{measurable}$.
Then either (i) $P(X<b)<\frac{1}{2}$ or (ii) $P(X<b)>0$. In Case
(i), we must have $a_{n}=0$ for each $n\geq1$. In Case (ii), because
of a.u. convergence, there exists $b'\in(b,1)$ such that $P(X_{n}<b')>0$
for some $n\geq1$, whence $a_{n}=1$ for some $n\geq1$. Since the
nondecreasing sequence $(a_{n})_{n=1,2,\cdots}$ is arbitrary, we
have deduced the principle of infinite search from said classical
theorem of martingale convergence .

Thus the boundlessness $X_{n}$ together with the constancy of $E|X_{n}|$
is not sufficient for the constructive a.u. convergence. Boundedness
is not the issue. Convexity is. The function $|x|$ simply does not
have any positive convexity away from $x=0$. 

With strictly convex functions $\lambda(x)$, to be presently defined,
which have positive and continuous second derivatives, as a natural
alternative to the function $|x|$, we will generalize Bishop's maximal
inequality for martingales, Theorem 3 in Chapter 8 of \cite{Bishop67},
to wide-sense submartingales. We then use the convergence of $E\lambda(X_{n}),$
as the criterion for a.u. convergence, obviating the use of upcrossing
inequalities. We will actually use a specific strictly convex function
$\overline{\lambda}$ such that $|\overline{\lambda}(x)|\leq3|x|$
for each $x\in R$. Then the boundlessness and convergence of $E\overline{\lambda}(X_{n})$
follows, classically, from the boundlessness of $E|X_{n}|$. Thus
we will have a criterion for constructive a.u. convergence which,
from the classical view point, imposes no additional condition beyond
the boundlessness of $E|X_{n}|$. The proof, being constructive, produces
rates of a.u. convergence.
\begin{defn}
\label{Def. Strictly Convex function} \textbf{(Strictly convex function).}
A continuous function $\lambda:R\rightarrow R$ is said to be\emph{
strictly convex}\index{strictly convex function} if it has a positive
and continuous second derivative $\lambda''$ on $R$. 
\end{defn}
$\square$

This definition generalizes the admissible functions in Chapter 8
of \cite{Bishop67}. The conditions of symmetry and nonnegativity
of $\lambda$ are dropped, so that we can admit an increasing function
$\lambda$. Correspondingly, we need to generalize Bishop's version
of Jensen's inequality, Lemma 2 in Chapter 8 of \cite{Bishop67},
to the following theorem.
\begin{thm}
\label{Thm. Bishop-Jensen inequality} \textbf{\emph{(Bishop-Jensen
inequality). }}Let $\lambda:R\rightarrow R$ be a strictly convex
function. Let the continuous function $\theta$ be defined by $\theta(x)\equiv\inf_{y\in[-x,x]}\lambda''(y)>0$
for each $x>0$. Define the continuous function $g:R^{2}\rightarrow[0,\infty)$
by
\begin{equation}
g(x_{0},x_{1})\equiv\frac{1}{2}(x_{1}-x_{0})^{2}\theta(|x_{0}|\vee|x_{1}|)\label{eq:temp-189}
\end{equation}
for each $(x_{0},x_{1})\in R^{2}$. 

Let $X_{0}$ and $X_{1}$ be integrable r.r.v.'s on $(\Omega,L,E)$
such that $\lambda(X_{0}),\lambda(X_{1})$ are integrable. Suppose
either \emph{(i)} $E(X_{1}|X_{0})=X_{0}$ , or \emph{(ii) }the strictly
convex function $\lambda:R\rightarrow R$ is nondecreasing and\emph{
$EUX_{0}\leq EUX_{1}$ }for each indicator $U\in L(X_{0})$. Then
the r.r.v. $g(X_{0},X_{1})$ is integrable, with 
\begin{equation}
0\leq Eg(X_{0},X_{1})\leq E\lambda(X_{1})-E\lambda(X_{0}).\label{eq:temp-124}
\end{equation}
\end{thm}
\begin{proof}
1. Let $(x_{0},x_{1})\in R^{2}$ be arbitrary. Then
\[
0\leq g(x_{0},x_{1})=\theta(|x_{0}|\vee|x_{1}|)\int_{v=x(0)}^{x(1)}\int_{u=x(0)}^{v}dudv
\]
\[
\leq\int_{v=x(0)}^{x(1)}(\int_{u=x(0)}^{v}\lambda''(u)du)dv=\int_{v=x(0)}^{x(1)}(\lambda'(v)-\lambda'(x_{0}))dv
\]
\begin{equation}
=\lambda(x_{1})-\lambda(x_{0})-\lambda'(x_{0})(x_{1}-x_{0}),\label{eq:temp-259}
\end{equation}
where the last equality is a direct application of the Fundamental
Theorem of Calculus. 

2. Let $(\Omega,L(X_{0}),E)$ denote the probability subspace of $(\Omega,L,E)$
generated by the r.r.v. $X_{0}$. Let $V\in L(X_{0})$ be an arbitrary
indicator such that $\lambda'(X_{0})V$ is bounded. Suppose Condition
(i) in the hypothesis holds, i.e. $E(X_{1}|X_{0})=X_{0}.$ Then, by
the properties of conditional expectations and the boundlessness of
the r.r.v. $\lambda'(X_{0})V\in L(X_{0})$, we obtain
\[
E(X_{1}-X_{0})\lambda'(X_{0})V=E(X_{0}-X_{0})\lambda'(X_{0})V=0.
\]
Suppose Condition (ii) holds. Then \emph{$EUX_{0}\leq EUX_{1}$ }for
each indicator $U\in L(X_{0})$, and the function $\lambda:R\rightarrow R$
is nondecreasing. Hence the bounded r.r.v. \emph{$\lambda'(X_{0})V\in L(X_{0})$
}is nonnegative. Therefore, by Proposition \ref{Prop. .Space of onditionally integrable functions is closed relative to L1},
we have
\[
E\lambda'(X_{0})VX_{0}\leq E\lambda'(X_{0})VX_{1}.
\]
Summing up, in either case, we have\emph{
\begin{equation}
E(X_{1}-X_{0})\lambda'(X_{0})V\geq0,\label{eq:temp-17}
\end{equation}
}for each indicator $V\in L(X_{0})$ such that $\lambda'(X_{0})V$
is bounded. 

Now the r.r.v.'s $\lambda(X_{0}),\lambda(X_{1}),X_{0},X_{1}$ are
integrable by hypothesis. Let $b>a>0$ be arbitrary. Since the function
$\lambda'$ is continuous, it is bounded on $[-b,b]$. Hence the r.r.v.
$\lambda'(X_{0})1_{(a\geq|X(0|)}$ is bounded. Therefore inequality
\ref{eq:temp-259} implies that $g(X_{0},X_{1})1_{(a\geq|X(0|)}$
is integrable. At the same time, inequality \ref{eq:temp-17} holds
with $V\equiv1_{(b\geq|X(0)|>a)}$. Consequently,
\[
0\leq Eg(X_{0},X_{1})1_{(b\geq|X(0|)}-Eg(X_{0},X_{1})1_{(a\geq|X(0|)}
\]
\[
=Eg(X_{0},X_{1})V\leq E(\lambda(X_{1})-\lambda(X_{0})-\lambda'(X_{0})(X_{1}-X_{0}))V
\]
\[
\leq E(\lambda(X_{1})-\lambda(X_{0}))V=E(\lambda(X_{1})-\lambda(X_{0}))1_{(b\geq|X(0)|>a)}\rightarrow0
\]
as $b>a\rightarrow\infty$, where the second inequality is due to
inequality \ref{eq:temp-259}. Hence the integral $g(X_{0},X_{1})1_{(a\geq|X(0|)}$
converges as $a\rightarrow\infty$. It follows from the Monotone Convergence
Theorem that the r.r.v. 
\[
g(X_{0},X_{1})=\lim_{a\rightarrow\infty}g(X_{0},X_{1})1_{(a\geq|X(0|)}
\]
is integrable, with 
\[
Eg(X_{0},X_{1})=\lim_{a\rightarrow\infty}Eg(X_{0},X_{1})1_{(a\geq|X(0|)}
\]
\[
\leq\lim_{a\rightarrow\infty}E(\lambda(X_{1})-\lambda(X_{0})-\lambda'(X_{0})(X_{1}-X_{0}))1_{(a\geq|X(0|)}
\]
\[
\leq\lim_{a\rightarrow\infty}E(\lambda(X_{1})-\lambda(X_{0}))1_{(a\geq|X(0|)}=E\lambda(X_{1})-E\lambda(X_{0}),
\]
where the first inequality follows from inequality \ref{eq:temp-259},
and the second inequality follows from inequality \ref{eq:temp-17}.
The theorem is proved.
\end{proof}
In the following, keep in mind the convention in Definition \ref{Convention. Only Regular Pts r.r.v.'s}
regarding regular points of r.r.v.'s. Now we are ready to prove the
advertised maximal inequality.
\begin{defn}
\textbf{\label{Def. Special convex function}(The special convex funtion).
}Define the continuous function $\overline{\lambda}:R\rightarrow R$
by
\begin{equation}
\overline{\lambda}(x)\equiv2x+(e^{-|x|}-1+|x|)\label{eq:temp-524-1}
\end{equation}
for each $x\in R$. We will call $\overline{\lambda}$ the \index{special convex function}\emph{special
convex function}.
\end{defn}
\begin{thm}
\label{Thm. Bishop's maximal inequality for w.s. submartingales.}
\textbf{\emph{(A maximal inequality for wide-sense submartingales).}}\textbf{
}Then the following holds.

1. The special convex function $\overline{\lambda}$ is increasing
and strictly convex, with 
\begin{equation}
|x|\leq|\overline{\lambda}(x)|\leq3|x|\label{eq:temp-549}
\end{equation}
for each $x\in R$. 

2. Let $Q\equiv\{t_{0},t_{1},\cdots,t_{n}\}$ be an arbitrary finite
subset of $R$ , with $t_{0}<t_{1}<\cdots<t_{n}$. Let $X:Q\times\Omega\rightarrow R$
be an arbitrary wide-sense submartingale relative to the filtration
$\mathcal{L}\equiv\{L^{(t(i))}:i=1,\cdots,n\}$. Let $\varepsilon>0$
be arbitrary. Suppose
\begin{equation}
E\overline{\lambda}(X_{t(n)})-E\overline{\lambda}(X_{t(0)})<\frac{1}{6}\varepsilon^{3}\exp(-3(E|X_{t(0)}|\vee E|X_{t(n)}|)\varepsilon^{-1}).\label{eq:temp-523}
\end{equation}
Then 
\begin{equation}
P(\bigvee_{k=0}^{n}|X_{t(k)}-X_{t(0)}|>\varepsilon)<\varepsilon.\label{eq:temp-509}
\end{equation}
We emphasize that the last two displayed inequalities are regardless
of how large $n\geq0$ is. We also note that, in view of inequality
\ref{eq:temp-549}, the r.r.v. $\overline{\lambda}(Y)$ is integrable
for each integrable r.r.v. $Y$. Thus inequality \ref{eq:temp-523}
is in contrast to the classical counterpart which requires either
$X_{t(n)}^{p}$ is integrable for some $p>1$ or $|X_{t(n)}|\log|X_{t(n)}|$
is integrable. 
\end{thm}
\begin{proof}
1. First note that $\overline{\lambda}(0)=0$. Elementary calculus
yields a continuous first derivative $\overline{\lambda}'$ on $R$
such that
\begin{equation}
\overline{\lambda}'(x)=2+(-e^{-x}+1)\geq2\label{eq:temp-454-2}
\end{equation}
for each $x\geq0$, and such that
\begin{equation}
\overline{\lambda}'(x)=2+(e^{x}-1)\geq1\label{eq:temp-454-1-1}
\end{equation}
for each $x\leq0$. Thus the function $\overline{\lambda}$ is increasing.
Moreover $\overline{\lambda}$ has a positive and continuous second
derivative
\[
\overline{\lambda}''(x)=e^{-|x|}
\]
for each $x\in R$. Hence
\[
\theta(x)\equiv\inf_{y\in[-x,x]}\overline{\lambda}''(y)=e^{-x}>0
\]
for each $x>0$. Thus the function $\overline{\lambda}$ is strictly
convex. Furthermore, since $0\leq e^{-r}-1+r\leq r$ for each $r\geq0$,
the triangle inequality yields
\[
|x|=2|x|-|x|\leq2|x|-(e^{-|x|}-1+|x|)
\]
\[
\leq|\overline{\lambda}(x)|\leq2|x|+(e^{-|x|}-1+|x|)\leq2|x|+|x|=3|x|.
\]
This establishes the desired inequality \ref{eq:temp-549}. 

As in Theorem \ref{Thm. Bishop-Jensen inequality}, define the continuous
function $g:R^{2}\rightarrow[0,\infty)$ by
\begin{equation}
g(x_{0},x_{1})\equiv\frac{1}{2}(x_{1}-x_{0})^{2}\theta(|x_{0}|\vee|x_{1}|)\equiv\frac{1}{2}(x_{1}-x_{0})^{2}\exp(-(|x_{0}|\vee|x_{1}|))\label{eq:temp-189-2}
\end{equation}
for each $(x_{0},x_{1})\in R^{2}$. 

2. By relabeling if necessary, we assume, without loss of generality,
that 
\[
(t_{0},t_{1},\cdots,t_{n})=(0,1,\cdots,n).
\]
Let $\varepsilon>0$ be as given. Write $K\equiv E|X_{0}|\vee E|X_{n}|$,
$b\equiv3K\varepsilon^{-1}$, and 
\[
\gamma\equiv\frac{1}{2}\varepsilon^{2}\theta(b)=\frac{1}{2}\varepsilon^{2}e^{-b}\equiv\frac{1}{2}\varepsilon^{2}\exp(-3K\varepsilon^{-1}).
\]
Then inequality \ref{eq:temp-523} in the hypothesis can be rewritten
as
\begin{equation}
E\overline{\lambda}(X_{n})-E\overline{\lambda}(X_{0})<\frac{1}{3}\varepsilon\gamma.\label{eq:temp-523-1}
\end{equation}
Let $\tau\equiv\eta_{0,\varepsilon,Q}$ be the simple first exit time
of the process $X$ after $t=0$ from the $\varepsilon$-neighborhood
of $X_{0}$, in the sense of Definition \ref{Def. Simple First Exit time}.
Define the probability subspace $L^{(\tau)}$ relative to the simple
stopping time $\tau$, as in Definition \ref{Def. Stopping times and simple stoppingg times.}.
Define the r.r.v.
\[
X_{\tau}\equiv\sum_{t\in Q}X_{t}1_{(\tau=t)}\in L^{(\tau)}.
\]
As in Corollary \ref{Cor.Doob's  Optional sampling of martingales},
define the process $X':\{0,1,2\}\times\Omega\rightarrow R$ by
\[
(X'_{0},X'_{1},X'_{2})\equiv(X_{0},X_{\tau},X_{n}),
\]
and define the filtration $\mathcal{L}'\equiv\{L'^{(i)}:i=0,1,2\}$
by 
\[
(L'^{(0)},L'^{(1)},L'^{(2)})\equiv(L^{(0)},L^{(\tau)},L^{(n)}).
\]
Then, by Corollary \ref{Cor.Doob's  Optional sampling of martingales},
the process $X'$ is a wide-sense submartingale relative to the filtration
$\mathcal{L}'$. In other words, \emph{$EUX'_{i-1}\leq EUX'_{i}$
}for each indicator $U\in L(X'_{i-1})\subset L'^{(i-1)}$, for each
$i=1,2$. 

Thus the conditions in Theorem \ref{Thm. Bishop-Jensen inequality}
are satisfied by the pair $X'_{i-1},X'_{i}$ of r.r.v.'s and by the
strictly convex function $\overline{\lambda}$, for each $i=1,2$.
Accordingly, for each $i=1,2$, Theorem \ref{Thm. Bishop-Jensen inequality}
implies that the nonnegative r.r.v. $Y_{i}\equiv g(X'_{i-1},X'_{i})$
is integrable, with 
\begin{equation}
0\leq EY_{i}\leq E\overline{\lambda}(X'_{i})-E\overline{\lambda}(X'_{i-1}).\label{eq:temp-124-2}
\end{equation}
Consequently.
\[
0\leq EY_{1}\leq E\overline{\lambda}(X'_{1})-E\overline{\lambda}(X'_{0})\leq E\overline{\lambda}(X'_{2})-E\overline{\lambda}(X'_{0}).
\]
\[
\equiv E\overline{\lambda}(X_{n})-E\overline{\lambda}(X_{0})<\frac{1}{3}\varepsilon\gamma,
\]
where the last inequality is inequality \ref{eq:temp-523-1}. Chebychev's
inequality therefore yields a $\mathrm{measurable}$ set $A$ with
$P(A^{c})<\frac{1}{3}\varepsilon$ such that 
\[
A\subset(Y_{1}\leq\gamma).
\]
Next, consider each $i=0,1$. Then, by equality \ref{eq:temp-524-1}
and equality \ref{eq:temp-124-2}, we obtain
\[
E|X'_{i}|\leq E\overline{\lambda}(X'_{i})\leq E\overline{\lambda}(X'_{2})\equiv E\overline{\lambda}(X_{n})\leq K\equiv\frac{1}{3}\varepsilon b,
\]
Chebychev's inequality therefore yields $P(B_{i}^{c})<\frac{1}{3}\varepsilon$
where $B_{i}\equiv(|X'_{i}|\leq b)$. 

Now consider each $\omega\in AB_{0}B_{1}$. Then $|X'_{0}(\omega)|\vee|X'_{1}(\omega)|\leq b$.
Hence
\[
\frac{1}{2}(X_{\tau}(\omega)-X_{0}(\omega))^{2}\theta(b)\equiv\frac{1}{2}(X'_{1}(\omega)-X'_{0}(\omega))^{2}\theta(b)
\]
\[
\leq\frac{1}{2}(X'_{1}(\omega)-X'_{0}(\omega))^{2}\theta(|X'_{0}(\omega)|\vee|X'_{1}(\omega)|)\equiv g(X'_{0}(\omega),X'_{1}(\omega))
\]
\[
\equiv Y_{1}(\omega)\leq\gamma\equiv\frac{1}{2}\varepsilon^{2}\theta(b),
\]
where the first inequality is because $\theta$ is a decreasing function,
and where the last inequality is because $\omega\in A\subset(Y_{1}\leq\gamma)$.
Dividing by $\frac{1}{2}\theta(b)$ and taking square roots, we obtain
\[
|X_{\eta(0,\varepsilon,Q)}(\omega)-X_{0}(\omega)|\equiv|X_{\tau}(\omega)-X_{0}(\omega)|\leq\varepsilon
\]
for each $\omega\in AB_{0}B_{1}$. It follows from the basic properties
of simple first exit time that, on $AB_{0}B_{1}$, we have
\[
\bigvee_{t\in Q}|X_{t}-X_{0}|\leq\varepsilon.
\]
Summing up, 
\[
P(\bigvee_{k=0}^{n}|X_{t(k)}-X_{t(0)}|>\varepsilon)\leq P(AB_{1}B_{2})^{c}<\frac{1}{3}\varepsilon+\frac{1}{3}\varepsilon+\frac{1}{3}\varepsilon=\varepsilon,
\]
as alleged. 
\end{proof}
Theorem \ref{Thm. Bishop's maximal inequality for w.s. submartingales.}
easily leads to the following a.u. convergence theorem for wide-sense
submartingales. We emphasize that, while constructive proofs of a.u.
convergence for a martingale $X:\{1,2,\cdots\}\times\Omega\rightarrow R$
are well known if $\lim_{k\rightarrow}E|X_{k}|^{p}$ exists for some
$p>1$, the following theorem requires no $L_{p}$-integrability for
$p>1$. 
\begin{thm}
\label{Thm. Specific strictly convex function and a.u. convergence of submartingles}
\textbf{\emph{(a.u. Convergence of wide-sense submartingales).}} Let
$X:Q\times\Omega\rightarrow R$ be an arbitrary wide-sense submartingale
relative to its natural filtration $\mathcal{L}$, where the parameter
set $Q$ is either $Q\equiv\{1,2,\cdots\}$ or $Q\equiv\{\cdots,-2,-1\}$.
Suppose \emph{(i)} there exists $b_{0}>0$ such that $b_{0}\geq E|X_{t}|$
for each $t\in Q$. As in the previous theorem, define the increasing
strictly convex function $\overline{\lambda}:R\rightarrow R$ by
\begin{equation}
\overline{\lambda}(x)\equiv2x+(e^{-|x|}-1+|x|)\label{eq:temp-524}
\end{equation}
for each $x\in R$. Then, for each $t\in Q$, we have $|X_{t}|\leq|\overline{\lambda}(X_{t})|\leq3|X_{t}|$,
whence the r.r.v. $\overline{\lambda}(X_{t})$ is integrable, with
$|E\overline{\lambda}(X_{t})|\leq3b_{0}$. 

Suppose, in addition,\emph{ (ii)} $\lim_{|t|\rightarrow\infty}E\overline{\lambda}(X_{t})$
exists. Then the following holds.

1. $X_{t}\rightarrow Y$ a.u. as $|t|\rightarrow\infty$ in $Q$,
for some r.r.v. $Y$. 

2. A rate of the above a.u. convergence can be obtained as follows.
Let $(b_{h})_{h=1,2,\cdots}$ be an arbitrary sequence of positive
real numbers such that $b_{h}\leq$$b_{0}$ and such that $b_{h}\geq E|X_{t}|$
for each $t\in Q$ with $|t|\geq h$, for each $h\geq1$. Let $k_{0}\equiv0$,
and let $(k_{m})_{m=1,2,\cdots}$ be a nondecreasing sequence of nonnegative
integers such that 
\begin{equation}
|E\overline{\lambda}(X_{t})-E\overline{\lambda}(X_{s})|<\frac{1}{6}2^{-3m}\exp(-2^{m}3b_{k(m-1)})\label{eq:temp-523-2}
\end{equation}
for each $t,s\in Q$ with $s\leq t$ and $|t|,|s|\geq k_{m}$, for
each $m\geq1$. Then, for each $m\geq1$, there exists a \textup{\emph{measurable}}\emph{
}set $A_{m}$ with $P(A_{m}^{c})<2^{-m+2}$ such that
\begin{equation}
A_{m}\subset\bigcap_{p=m}^{\infty}\bigcap_{t\in Q;|t|\geq k(p)}^{\infty}(|X_{t}-Y|\leq2^{-p+3}).\label{eq:temp-233-1-4}
\end{equation}

Note that, if $X$ is a martingale with $Q\equiv\{1,2,\cdots\}$,
and if both $\lim_{t\rightarrow\infty}E|X_{t}|$ and $\lim_{t\rightarrow\infty}Ee^{-|X(t)|}$
exist, then both Conditions \emph{(i)} and \emph{(ii)} hold because
$EX_{t}$ is a constant. In general, note also that if Condition\emph{
(i)} holds, then, since $E\overline{\lambda}(X_{t})$ is a nondecreasing
function of $t\in Q$ according to Theorem \ref{Thm. Bishop-Jensen inequality},
and since $|E\overline{\lambda}(X_{t})|\leq3E|X_{t}|\leq3b_{0}$ for
each $t\in Q$ , Condition \emph{(ii)} is, classically, automatically
satisfied.
\end{thm}
\begin{proof}
1. First note that Assertion 1 of Theorem \ref{Thm. Bishop's maximal inequality for w.s. submartingales.}
implies that $|X_{t}|\leq|\overline{\lambda}(X_{t})|\leq3|X_{t}|$
whence the r.r.v. $\overline{\lambda}(X_{t})$ is integrable, with
$|E\overline{\lambda}(X_{t})|\leq E|\overline{\lambda}(X_{t})|\leq3|X_{t}|\leq3b_{0}$,
for each $t\in Q$. 

2. Let $h\geq1$ be arbitrary. Condition (i) in the hypothesis guarantees
that there exists $b_{h}\in(0,$$b_{0}]$ such that $b_{h}>E|X_{t}|$
for each $t\in Q$ with $|t|\geq h$. If necessary, we can always
take $b_{h}\equiv b_{0}$. 

3. Let $m\geq1$ be arbitrary. Take any $\varepsilon_{m}\in(2^{-m},2^{-m+1})$.
Then inequality \ref{eq:temp-523-2} implies that
\begin{equation}
E\overline{\lambda}(X_{t})-E\overline{\lambda}(X_{s})<\frac{1}{6}2^{-3m}\exp(-3b_{k(m-1)}2^{m})<\frac{1}{6}\varepsilon_{m}^{3}\exp(-3b_{k(m-1)}\varepsilon_{m}^{-1})\label{eq:temp-523-2-3}
\end{equation}
for each $t,s\in Q$ with $s\leq t$ and $|t|,|s|\geq k_{m}$.

4. We will first prove the theorem for the case where $Q\equiv\{1,2,\cdots\}$.
Let $m\geq1$ be arbitrary. Then $E|X_{t}|\leq b_{k(m-1)}$ for each
$t\in Q$ with $t\geq k_{m-1}$. In particular, $E|X_{k(m)}|\vee E|X_{k(m+1)}|\leq b_{k(m-1)}$.
Hence, since $e^{-r}$ is a decreasing function of $r\in R$, inequality
\ref{eq:temp-523-2-3} yields
\begin{equation}
E\overline{\lambda}(X_{k(m+1)})-E\overline{\lambda}(X_{k(m)})<\frac{1}{6}\varepsilon_{m}^{3}\exp(-3(E|X_{k(m)}|\vee E|X_{k(m+1)}|)\varepsilon_{m}^{-1}).\label{eq:temp-523-2-1}
\end{equation}
Therefore Theorem \ref{Thm. Bishop's maximal inequality for w.s. submartingales.}
implies that $P(B_{m})<\varepsilon_{m}$, where we define
\begin{equation}
B_{m}\equiv(\bigvee_{k=k(m)}^{k(m+1)}|X_{k}-X_{k(m)}|>\varepsilon_{m}).\label{eq:temp-509-2}
\end{equation}
Now define $A_{m}\equiv\bigcap_{h=m}^{\infty}B_{h}^{c}$. Then 
\[
P(A_{m}^{c})\leq\sum_{h=m}^{\infty}P(B_{h})<\sum_{h=m}^{\infty}\varepsilon_{h}<\sum_{h=m}^{\infty}2^{-h+1}=2^{-m+2}.
\]
Consider each $\omega\in A_{m}$. Let $p\geq m$ and $j\geq i\geq k_{p}$
be arbitrary. Then $k_{h}\leq i\leq k_{h+1}$ and $k_{n}\leq j\leq k_{n+1}$
for some $n\geq h\geq p$. Consequently, 
\[
|X_{i}(\omega)-X_{j}(\omega)|
\]
\[
\leq|X_{i}(\omega)-X_{k(h)}(\omega)|+|X_{k(h)}(\omega)-X_{k(h+1)}(\omega)|+\cdots+|X_{k(n)}(\omega)-X_{j}(\omega)|
\]
\[
\leq\varepsilon_{h}+(\varepsilon_{h}+\cdots+\varepsilon_{n})
\]
\begin{equation}
<2^{-h+1}+\sum_{\kappa=h}^{\infty}2^{-\kappa+1}=2^{-h+1}+2^{-h+2}<2^{-p+3},\label{eq:temp-233}
\end{equation}
where the second inequality is because $\omega\in A_{m}\subset B_{h}^{c}B_{h+1}^{c}\cdots B_{n}^{c}$.
Since $2^{-p+3}$ is arbitrarily small for sufficiently large $p\geq m$,
we see that the sequence $(X_{\kappa}(\omega))_{\kappa=1,2,\cdots}$
of real numbers is Cauchy, and so $Y(\omega)\equiv\lim_{\kappa\rightarrow\infty}X_{\kappa}(\omega)$
exists. Fixing $i$ and letting $j\rightarrow\infty$ in inequality
\ref{eq:temp-233}, we obtain 
\[
|X_{i}(\omega)-Y(\omega)|\leq2^{-p+3}
\]
for each $i\geq k_{p}$, where $p\geq m$ is arbitrary, for each $\omega\in A_{m}$.
Hence
\[
A_{m}\subset\bigcap_{p=m}^{\infty}\bigcap_{i=k(p)}^{\infty}(|X_{i}-Y|\leq2^{-p+3})=\bigcap_{p=m}^{\infty}\bigcap_{t\in Q;|t|\geq k(p)}^{\infty}(|X_{t}-Y|\leq2^{-p+3}).
\]
Thus $X_{i}\rightarrow Y$ uniformly on the $\mathrm{measurable}$
set $A_{m}$, where $P(A_{m}^{c})<2^{-m+2}$ is arbitrarily small
when $m\geq1$ is sufficiently large. In other words, $X_{i}\rightarrow Y$
a.u. By Proposition \ref{Prop.  a.u. Convergence =00003D> convergence in prob, etc},
the function $Y$ is a r.r.v. The theorem has been proved for the
case where $Q\equiv\{1,2,\cdots\}$.

5. The proof for the case where $Q\equiv\{\cdots,-2,-1\}$ is almost
a mirror image of the preceding paragraph. Let $m\geq1$ be arbitrary.
Then $E|X_{t}|\leq b_{k(m-1)}$ for each $t\in Q$ with $t\leq-k_{m-1}$.
In particular, $E|X_{-k(m+1)}|\vee E|X_{-k(m)}|\leq b_{k(m-1)}$.
Hence, since $e^{-r}$ is a decreasing function of $r\in R$, inequality
\ref{eq:temp-523-2-3} yields
\begin{equation}
E\overline{\lambda}(X_{-k(m)})-E\overline{\lambda}(X_{-k(m+1)})<\frac{1}{6}\varepsilon_{m}^{3}\exp(-3(E|X_{-k(m+1)}|\vee E|X_{-k(m)}|)\varepsilon_{m}^{-1}).\label{eq:temp-523-2-1-1}
\end{equation}
Therefore Theorem \ref{Thm. Bishop's maximal inequality for w.s. submartingales.}
implies that $P(B_{m})<\varepsilon_{m}$, where we define
\begin{equation}
B_{m}\equiv(\bigvee_{k=k(m)}^{k(m+1)}|X_{-k}-X_{-k(m+1)}|>\varepsilon_{m}).\label{eq:temp-509-2-1}
\end{equation}
Now define $A_{m}\equiv\bigcap_{h=m}^{\infty}B_{h}^{c}$. Then 
\[
P(A_{m}^{c})\leq\sum_{h=m}^{\infty}P(B_{h})<\sum_{h=m}^{\infty}\varepsilon_{h}<\sum_{h=m}^{\infty}2^{-h+1}=2^{-m+2}.
\]
Consider each $\omega\in A_{m}$. Let $p\geq m$ and $j\geq i\geq k_{p}$
be arbitrary. Then $k_{h}\leq i\leq k_{h+1}$ and $k_{n}\leq j\leq k_{n+1}$
for some $n\geq h\geq p$. Consequently, 
\[
|X_{-j}(\omega)-X_{-i}(\omega)|
\]
\[
\leq|X_{-j}(\omega)-X_{-k(n+1)}(\omega)|+|X_{-k(n+1)}(\omega)-X_{-k(n)}(\omega)|+\cdots+|X_{-k(h+1)}(\omega)-X_{i}(\omega)|
\]
\[
\leq\varepsilon_{n}+(\varepsilon_{n}+\cdots+\varepsilon_{h})
\]
\begin{equation}
<2^{-n+1}+\sum_{\kappa=h}^{\infty}2^{-\kappa+1}=2^{-n+1}+2^{-h+2}<2^{-p+3},\label{eq:temp-233-1}
\end{equation}
where the second inequality is because $\omega\in A_{m}\subset B_{n}^{c}B_{n+1}^{c}\cdots B_{h}^{c}$.
Since $2^{-p+3}$ is arbitrarily small for sufficiently large $p\geq m$,
we see that the sequence $(X_{-\kappa}(\omega))_{\kappa=1,2,\cdots}$
of real numbers is Cauchy, and so $Y(\omega)\equiv\lim_{\kappa\rightarrow\infty}X_{-\kappa}(\omega)$
exists. Fixing $i$ and letting $j\rightarrow\infty$ in inequality
\ref{eq:temp-233-1}, we obtain 
\[
|X_{-i}(\omega)-Y(\omega)|\leq2^{-p+3}
\]
for each $i\geq k_{p}$, where $p\geq m$ is arbitrary, for each $\omega\in A_{m}$.
Hence
\[
A_{m}\subset\bigcap_{p=m}^{\infty}\bigcap_{i=k(p)}^{\infty}(|X_{-i}-Y|\leq2^{-p+3})=\bigcap_{p=m}^{\infty}\bigcap_{t\in Q;|t|\geq k(p)}^{\infty}(|X_{t}-Y|\leq2^{-p+3}).
\]
Thus $X_{-i}\rightarrow Y$ uniformly on the $\mathrm{measurable}$
set $A_{m}$, with $P(A_{m}^{c})<2^{-m+2}$ arbitrarily small when
$m\geq1$ is sufficiently large. In other words, $X_{t}\rightarrow Y$
a.u. as $|t|\rightarrow\infty$ with $t\in Q\equiv\{\cdots,-2,-1\}$.
By Proposition \ref{Prop.  a.u. Convergence =00003D> convergence in prob, etc},
the function $Y$ is a r.r.v. The theorem has been proved also for
the case where $Q\equiv\{\cdots,-2,-1\}$.
\end{proof}

\section{The Law of Large Numbers}

Applications of martingales are numerous. One application is to prove
the Strong Law of Large numbers (SLLN). This theorem says that if
$Z_{1},Z_{2},\cdots$ is a sequence of integrable independent and
identically distributed r.r.v's, with mean $0$, then $n^{-1}(Z_{1}+\cdots+Z_{n})\rightarrow0$
a.u. Historically the first proof of this theorem in its generality,
due to Kolmogorov, is constructive, complete with rates of convergence.
See, for example, Theorem 5.4.2 of \cite{Chung68}. Subsequently,
remarkable proofs are also given in terms of a.u. martingale convergence
via Doob's upcrossing inequality. See, for example, Theorem 9.4.1
of \cite{Chung68}. As observed earlier, the theorem that deduces
a.u. convergence from upcrossing inequalities actually implies the
principle of infinite search, and cannot be made constructive. For
that reason, we present below a constructive proof by a simple application
of Theorem \ref{Thm. Specific strictly convex function and a.u. convergence of submartingles}
in the previous section. A similar constructive proof is most likely
known in the literature.

First the weak law, with a well known proof by characteristic functions. 
\begin{thm}
\textbf{\emph{\label{Thm. Weak Law of Large Numbers.} (Weak Law of
Large Numbers). }}Suppose $Z_{1},Z_{2},\cdots$ is a sequence of integrable,
independent, and identically distributed r.r.v's with mean $0$, on
some probability space $(\Omega,L,E)$. Let $\eta_{intg}$ be a simple
modulus of integrability of $Z_{1}$, in the sense of Definition \ref{Def. Modulus of integrability}.
For each $m\geq1$, let \textbf{\emph{$S_{m}\equiv m^{-1}(Z_{1}+\cdots+Z_{m})$}}.
Then
\[
E|S_{m}|\rightarrow0
\]
as\emph{ }\textbf{\emph{$m\rightarrow\infty$}}. More precisely, for
each $m\geq1$, there exists an integer $q_{m}\equiv q_{m,\eta(intg)}\geq1$
such that $E|S_{k}|\leq2^{-m}$ for each $k\geq q_{m}$, .
\end{thm}
\begin{proof}
1. By hypothesis, the independent r.r.v.'s $Z_{1},Z_{2},\cdots$ have
a common distribution $J$ on $R$. Hence they share a common characteristic
function $\psi$. Therefore, for each $n\geq1$, the characteristic
function of the r.r.v. $S_{n}$ is given by $\psi_{n}\equiv\psi^{n}(\frac{\cdot}{n})$.
Let $J_{n}$ denote the distribution of $S_{n}.$ Let $J_{0}$ denote
the distribution on $R$ which assigns probability $1$ to the point
$0\in R$. Then the characteristic function of $J_{0}$ is the constant
function $\psi_{0}\equiv1$ on $R$. Define the remainder function
$r_{1}$ of the first degree Taylor expansion of the characteristic
function $\psi$ by
\[
\psi(u)\equiv1+iuEZ_{1}+r_{1}(u)=1+r_{1}(u)
\]
for each $u\in R$, where the mean $EZ_{1}$ vanishes by hypothesis. 

2. Separately, take an arbitrary $a>\eta_{intg}(1)$. Then $E|Z_{1}|\leq E|Z_{1}|1_{(Z(1)>a)}+a\leq1+a$,
by the Definition \ref{Def. Modulus of integrability} of a simple
modulus of integrability. Letting $a\downarrow\eta_{intg}(1)$ then
yields
\begin{equation}
E|Z_{k}|=E|Z_{1}|\leq b\equiv1+\eta_{intg}(1)\label{eq:temp-456}
\end{equation}
for each $k\geq1$.

3. Let $n\geq1$ be arbitrary. Define the positive real number\emph{
}
\begin{equation}
c\equiv\pi^{-1}2^{2n+4},\label{eq:temp-468}
\end{equation}
and integer 
\begin{equation}
p_{n}\equiv p_{n,\eta(intg)}\equiv[\frac{8bc^{3}}{\pi}\eta_{intg}(\frac{\pi}{8c^{2}})]_{1}.\label{eq:temp-469}
\end{equation}
Consider each $k\geq p_{n}$ and each $u\in[-c,c]$. Write $\alpha\equiv\frac{\pi}{4c^{2}}$
for short. Then\emph{ 
\[
|\frac{u}{k}|\leq\frac{c}{k}\leq\frac{c}{p_{n}}<\frac{\pi c}{8bc^{3}}(\eta_{intg}(\frac{\pi}{8c^{2}}))^{-1}=\frac{\alpha}{2b}(\eta_{intg}(\frac{\alpha}{2}))^{-1}.
\]
}Hence, by Assertion 2 of Proposition\emph{ \ref{Prop Taylor of ch fct},
}where the dimension is set to \emph{$1$, }and where $X,\lambda,\varepsilon,b$
are replaced by $Z,\frac{u}{n},\alpha,b$ respectively, we obtain
\[
|r_{1}(\frac{u}{k})|<\alpha|\frac{u}{k}|\leq\alpha\frac{c}{k}=\frac{\pi}{4kc}.
\]
For abbreviation, write $z\equiv r_{1}(\frac{u}{k})$. Then $|z|<\frac{\pi}{4kc}$.
Therefore the binomial expansion yields 
\[
|(1+r_{1}(\frac{u}{k}))^{k}-1|\equiv|(1+z)^{k}-1|
\]
\[
=|(1+C_{1}^{k}z+C_{2}^{k}z^{2}+\cdots+C_{k}^{k}z^{k})-1|
\]
\[
\leq C_{1}^{k}\frac{\pi}{4kc}+C_{2}^{k}(\frac{\pi}{4kc})^{2}+\cdots+C_{k}^{k}(\frac{\pi}{4kc})^{k}
\]
\[
\leq(\frac{\pi}{4c})^{1}+(\frac{\pi}{4c})^{2}\cdots+(\frac{\pi}{4c})^{k}<\frac{\pi}{4c}(1-\frac{\pi}{4c})^{-1}
\]
\[
\equiv\frac{\pi}{4c}(1-\frac{\pi}{4\cdot\pi^{-1}2^{2n+4}})^{-1}=\frac{\pi}{4c}(1-\pi^{2}2^{-2n-6})^{-1}<\frac{\pi}{4c}\cdot2=\frac{\pi}{2c},
\]
where 
\[
C_{j}^{k}\equiv\frac{k(k-1)\cdots(k-j+1)}{j!}\leq k^{j}
\]
for each $j=1,\cdots,k.$ Consequently, 
\[
|\psi_{k}(u)-\psi_{0}(u)|=|\psi^{k}(\frac{u}{k})-1|
\]
\begin{equation}
=|(1+r_{1}(\frac{u}{k}))^{k}-1|<\frac{\pi}{2c},\label{eq:temp-476}
\end{equation}
where $u\in[-c,c]$, $k\geq p_{n}$, and $n\geq1$ are arbitrary.

4. Now let $a\in(2^{-n},2^{-n+1})$ be arbitrary. Define the function
$f\in C(R)$ by 
\begin{equation}
f(x)\equiv a^{-1}(1-a^{-1}|x|)_{+}\label{eq:temp-489}
\end{equation}
for each $x\in R$. Then the Fourier Transform $\hat{f}$ of the function
$f$ satisfies
\[
|\hat{f}(u)|\equiv|\int_{x\in R}e^{iux}f(x)dx|=\frac{2}{a}\int_{0}^{a}(\cos ux)(1-\frac{x}{a})dx
\]
\[
=\frac{2}{a}\int_{0}^{a}\frac{\sin ux}{u}\frac{1}{a}dx=2\frac{(1-\cos au)}{a^{2}u^{2}}\leq1\wedge\frac{4}{a^{2}u^{2}}
\]
for each $u\in R$, where the third equality is by integration by
parts in Calculus. Therefore $\hat{f}$ is Lebesgue integrable on
$R$. Hence Assertion 3 of Theorem \ref{Thm. Inversion  Formula  for Ch fcts}
implies that
\[
J_{k}f=(2\pi)^{-1}\int\hat{f}(-u)\psi_{k}(u)du,
\]
with a similar equality when $k$ is replaced by $0$. Consequently,
\[
2\pi|J_{k}f-J_{0}f|=|\int\hat{f}(-u)(\psi_{k}(u)-\psi_{0}(u))du|
\]
\[
=|\int_{|u|\leq c}\hat{f}(-u)((1+r_{1}(\frac{u}{k}))^{k}-1)du+\int_{|u|>c}\hat{f}(-u)(\psi_{k}(u)-\psi_{0}(u))du|
\]
\[
\leq\int_{|u|\leq c}|(1+r_{1}(\frac{u}{k}))^{k}-1|du+\int_{|u|>c}\frac{4}{a^{2}u^{2}}\cdot2du
\]
\[
<\int_{|u|\leq c}\frac{\pi}{2c}du+\frac{16}{a^{2}c}\leq\pi+\frac{16}{2^{-2n}c}=\pi+\pi=2\pi,
\]
whence
\[
|J_{k}f-J_{0}f|\leq1.
\]
At the same time, from the defining formula \ref{eq:temp-489}, we
see that $1_{[-a,a]}\geq af$. Hence 
\[
P(|S_{k}|>a)=J_{k}(1-1_{[-a,a]})
\]
\[
\leq1-aJ_{k}f\leq1-a(J_{0}f-1)=1-a(f(0)-1)=1-a(a^{-1}-1)=a,
\]
where $a\in(2^{-n},2^{-n+1})$ , $k\geq p_{n}$, and $n\geq1$ are
arbitrary. Since $a$ is arbitrarily small for sufficiently large
$n\geq1$, we conclude that $S_{k}\rightarrow0$ in probability.

5.\emph{ }To prove $E|S_{k}|\rightarrow0$, first note that, by Proposition
\ref{Prop. Existence of modulus of integrability}, for each $k\geq1$,
the r.r.v. $Z_{k}$ has a modulus of integrability $\delta$ defined
by 
\[
\delta(\varepsilon)\equiv\frac{\varepsilon}{2}/\eta_{intg}(\frac{\varepsilon}{2})
\]
for each $\varepsilon>0$. Let $m\geq1$ be arbitrary. Let $n\equiv(m+2)\vee[1-\log_{2}(\delta(2^{-m-1}))]_{1}$.
Take an arbitrary $a\in(2^{-n},2^{-n+1})$. Define 
\begin{equation}
q_{m,\eta(intg)}\equiv p_{n,\eta(intg)}.\label{eq:temp-467}
\end{equation}
Consider each $k\geq q_{m,\eta(intg)}\equiv p_{n,\eta(intg)}$. Then,
by Step 4, we have 
\[
P(|S_{k}|>a)<a<2^{-n+1}<\delta(2^{-m-1}).
\]
Hence, since $\delta$ is a modulus of integrability of $Z_{\kappa}$
for each $\kappa=1,\cdots,k$, it follows that 
\[
E|S_{k}|1_{(|S(k)|>a)}\leq k^{-1}\sum_{\kappa=1}^{k}E|Z_{\kappa}|1_{(|S(k)|>a)}\leq k^{-1}\sum_{\kappa=1}^{k}2^{-m-1}=2^{-m-1}.
\]
Consequently,
\[
E|S_{k}|\leq E|S_{k}|1_{(|S(k)|\leq a)}+E|S_{k}|1_{(|S(k)|>a)}\leq a+2^{-m-1}
\]
\[
<2^{-n+1}+2^{-m-1}\leq2^{-m-1}+2^{-m-1}=2^{-m},
\]
where $m\geq1$ and $k\geq q_{m,\eta(intg)}$ are arbitrary. We conclude
that $E|S_{k}|\rightarrow0$ as $k\rightarrow\infty$.
\end{proof}
\begin{thm}
\textbf{\emph{\label{Thm. Strong Law of Large Numbers} (Strong Law
of Large Numbers). }}Suppose $Z_{1},Z_{2},\cdots$ is a sequence of
integrable, independent, and identically distributed r.r.v's with
mean $0$, on some probability space $(\Omega,L,E)$. Let $\eta_{intg}$
be a simple modulus of integrability of $Z_{1}$, in the sense of
Definition \ref{Def. Modulus of integrability}. Then 
\[
S_{k}\equiv k^{-1}(Z_{1}+\cdots+Z_{k})\rightarrow0\qquad a.u.
\]
as $k\rightarrow\infty$. More precisely, for each $m\geq1$ there
exists an integer $k_{m}\equiv k_{m,\eta(intg)}$ and a a \textup{\emph{measurable}}\emph{
}set $A_{m}$, with $P(A_{m}^{c})<2^{-m+2}$ and with
\begin{equation}
A_{m}\subset\bigcap_{p=m}^{\infty}\bigcap_{k=k(p)}^{\infty}(|S_{k}|\leq2^{-p+3}).\label{eq:temp-233-1-1-1}
\end{equation}
\end{thm}
\begin{proof}
1. Let $m\geq j\geq1$ be arbitrary, and let $I_{j}$ denote the distribution
of $Z_{j}$ on $R$. Then, in view of the the hypothesis of independence
and identical distribution, the r.v. $(Z_{1},\cdots,Z_{j},\cdots,Z_{m})$
with values in $R^{m}$ has the same distribution as the r.v $(Z_{j},\cdots,Z_{1},\cdots,Z_{m})$,
where, for brevity, the latter stands for the sequence obtained from
$(Z_{1},\cdots,Z_{j},\cdots,Z_{m})$ by swapping the first- and the
$j$-th members. Now let $h$ be an arbitrary Thus 
\begin{equation}
Eh(Z_{1},\cdots,Z_{j},\cdots,Z_{m})=Eh(Z_{j},\cdots,Z_{1},\cdots,Z_{m})\label{eq:temp-488}
\end{equation}
for each integrable function $h$ on $R^{m}$ relative to the joint
distribution $E_{Z(1),\cdots,Z(m)}$. 

2. Let $Q\equiv\{\cdots,-2,-1\}$. For each $k\geq1$, define $X_{-k}\equiv S_{k}$.
Let $\mathcal{\mathcal{L}\equiv\mathcal{L}_{X}}\equiv\{L^{(X,t)}:t\in Q\}$
be the natural filtration of the process $X:Q\times\Omega\rightarrow R$.
Let $t\in Q$ be arbitrary. Then $t=-n$ for some $n\geq1$. Hence
$(\Omega,L^{(X,t)},E)$ is the\emph{ }probability subspace of\emph{
}$(\Omega,L,E)$\emph{ }generated by the family
\[
G^{(X,t)}\equiv\{X_{r}:r\in Q;\;r\leq t\}=\{S_{m}:m\geq n\}.
\]
In other words, $(\Omega,L^{(X,t)},E)$ is the completion of the integration
space
\[
(\Omega,L_{C(ub)}(G^{(X,t)}),E),
\]
where
\begin{equation}
L_{C}(G^{(X,t)})=\{f(S_{n},\cdots,S_{m}):m\geq n;f\in C(R^{m-n+1})\}.\label{eq:temp-490}
\end{equation}
By Lemma \ref{Lem. Filtration genrerated by process is indeed a filtration},
the process $X$ is adapted to its natural filtration $\mathcal{L}$. 

3. We will prove that the process $X$ is a martingale relative to
the filtration $\mathcal{L}$. To that end, let $s,t\in Q$ be arbitrary
with $t\leq s$. Then $t=-n$ and $s=-k$ for some $n\geq k\geq1$.
Let $Y\in L_{C}(G^{(X,t)})$ be arbitrary. Then, in view of equality
\ref{eq:temp-490}, we have 
\[
Y=f(S_{n},\cdots,S_{m})
\]
\[
\equiv f((Z_{1}+\cdots+Z_{n})n^{-1}\cdots,(Z_{1}+\cdots+Z_{m})m^{-1})
\]
for some $f\in C_{ub}(R^{m-n+1})$, for some $m\geq n$. Let $j=1,\cdots,n$
be arbitrary. Then, since the r.r.v. $Y$ is bounded, the r.r.v. $YZ_{j}$
is integrable. Hence, by equality \ref{eq:temp-488}, we have 
\[
Ef((Z_{1}+\cdots+Z_{j}+\cdots+Z_{n})n^{-1},\cdots,(Z_{1}+\cdots+Z_{j}+\cdots+Z_{m})m^{-1})Z_{j}
\]
\[
=Ef((Z_{j}+\cdots+Z_{1}+\cdots+Z_{n})n^{-1},\cdots,(Z_{j}+\cdots+Z_{1}+\cdots+Z_{m})m^{-1})Z_{1}
\]
\[
=Ef((Z_{1}+\cdots+Z_{j}+\cdots+Z_{n})n^{-1},\cdots,(Z_{1}+\cdots+Z_{j}+\cdots+Z_{m})m^{-1})Z_{1}.
\]
In short, 
\[
EYZ_{j}=EYZ_{1}
\]
for each $j=1,\cdots,n$. Therefore, since $k\leq n$, we have
\[
EYS_{k}=k^{-1}E(YZ_{1}+\cdots+YZ_{k})=EYZ_{1}.
\]
In particular, $EYS_{n}=EYZ_{1}=EYS_{k}.$ In other words, $EYX_{t}=EYX_{s}$,
where $X_{t}\in L^{(X,t)}$, and where $Y\in L_{C}(G^{(X,t)})$ is
arbitrary. Hence 
\[
E(X_{s}|L^{(X,t)})=X_{t},
\]
according to Assertion 5 of Proposition \ref{Prop. Basics of Conditional expectations}.
Since $s,t\in Q$ are arbitrary with $t\leq s$, the process $X$
is a martingale relative to its natural filtration $\mathcal{L}\equiv\{L^{(X,t)}:t\in Q\}$,.

4. Let $m\geq1$ be arbitrary and write $h\equiv3m+8$. By Theorem
\ref{Thm. Weak Law of Large Numbers.}, there exists an integer $q_{h,\eta(intg)}\geq1$
so large that $E|S_{k}|\leq2^{-h}$ for each $k\geq q_{h,\eta(intg)}$.
Define 
\[
k_{m}\equiv k_{m,\eta(intg)}\equiv q_{h,\eta(intg)}\equiv q_{3m+8,\eta(intg)}\geq1.
\]
Thus
\begin{equation}
E|X_{-k}|\equiv E|S_{k}|\leq b_{m}\equiv2^{-3m+8}\label{eq:temp-453}
\end{equation}
for each $k\geq k_{m}$. Since $2^{-3m-8}$ is arbitrarily small for
sufficiently large $m\geq1$, we see that $E|X_{t}|\rightarrow0$
as $k\equiv|t|\rightarrow\infty$, with $t\in Q$. Consequently $X_{t}\rightarrow0$
in probability. Hence $e^{-|X(t)||}\rightarrow1$ in probability.
The Dominated Convergence Theorem then implies that $Ee^{-|X(t)|}\rightarrow1$
as $|t|\rightarrow\infty$ with $t\in Q$. 

5. With arbitrary $m\geq1$, consider each $t,s\in Q$ with $|t|,|s|\geq k_{m}$.
Recall the increasing and strictly convex function $\overline{\lambda}:R\rightarrow R$
defined in Theorem \ref{Thm. Specific strictly convex function and a.u. convergence of submartingles}
by
\begin{equation}
\overline{\lambda}(x)\equiv2x+(e^{-|x|}-1+|x|)\label{eq:temp-524-1-1-1}
\end{equation}
for each $x\in R$. Then
\[
|E\overline{\lambda}(X_{t})-E\overline{\lambda}(X_{s})|
\]
\[
\equiv|E(2X_{t}+(e^{-|X(t)|}-1+|X_{t}|))-E(2X_{s}+(e^{-|X(s)|}-1+|X_{s}|))|
\]
\[
=|E(e^{-|X(t)|}-1+|X_{t}|)-E(e^{-|X(s)|}-1+|X_{s}|)|
\]
\[
\leq E(|X_{t}|+|X_{t}|)+E(|X_{s}|+|X_{s}|))|<2^{2}2^{-3m-8}=2^{-3m-6}
\]
\begin{equation}
=2^{-3m-1}3^{-2}(2^{-5}3^{2})<2^{-3m-1}3^{-2}e^{-1},\label{eq:temp-452}
\end{equation}
where the second equality is because $X$ is a martingale, where the
first inequality is thanks to the elementary inequality $|e^{-r}-1|\leq r$
for each $r\geq0$, and where the second inequality is from inequality
\ref{eq:temp-453}. Thus $\lim_{|t|\rightarrow\infty}E\overline{\lambda}(X_{t})$
exists. 

6. Continuing with arbitrary $m\geq1$, take any $\varepsilon_{m}\in(2^{-m},2^{-m+1})$.
Since
\[
2^{m}3b_{m}<2^{m}3\cdot2^{-3m-7}=2^{-2m-7}3\leq2^{-9}3<1.
\]
 inequality \ref{eq:temp-452} implies 
\[
|E\overline{\lambda}(X_{t})-E\overline{\lambda}(X_{s})|<2^{-3m-1}3^{-2}e^{-1}
\]
\begin{equation}
<2^{-3m-1}3^{-2}\exp(-2^{m}3b_{m})<\frac{1}{6}\varepsilon_{m}^{3}\exp(-3b_{m}\varepsilon_{m}^{-1})\label{eq:temp-451}
\end{equation}
for each $t,s\in Q$ with $|t|,|s|\geq k_{m}$. In view of inequalities
\ref{eq:temp-453} and \ref{eq:temp-451}, all the conditions in the
hypothesis of Theorem \ref{Thm. Specific strictly convex function and a.u. convergence of submartingles}
are satisfied. Accordingly, $X_{t}\rightarrow Y$ a.u. as $|t|\rightarrow\infty$,
for some r.r.v. $Y$. At the same time, since $E|X_{t}|\rightarrow0$
as $|t|\rightarrow\infty$, some subsequence of $(X_{t})_{t=-1,-2,\cdots}$
converges to $0$ a.u. as $|t|\rightarrow\infty$. Hence $Y=0$ a.u.
Summing up, $S_{k}\equiv X_{-k}\rightarrow0$ a.u. as $k\rightarrow\infty$.
Moreover, in view of inequality \ref{eq:temp-451}, Theorem \ref{Thm. Specific strictly convex function and a.u. convergence of submartingles}
implies that, for each $m\geq1$, there exists a $\mathrm{measurable}$
set $A_{m}$ with $P(A_{m}^{c})<2^{-m+2}$ such that
\begin{equation}
A_{m}\subset\bigcap_{p=m}^{\infty}\bigcap_{t\in Q;|t|\geq k(p)}^{\infty}(|X_{t}-Y|\leq2^{-p+3})\equiv\bigcap_{p=m}^{\infty}\bigcap_{k=k(p)}^{\infty}(|S_{k}|\leq2^{-p+3}),\label{eq:temp-233-1-4-1}
\end{equation}
as desired.
\end{proof}

\chapter{a.u. Continuous Processes on $[0,1]$}

In this chapter \label{Chapter. Processes with continuous sample paths},
let $(S,d)$ be a locally compact metric space. Unless otherwise specified,
this will serve as the state space for the processes in this chapter.
We consider an arbitrary consistent family $F$ of f.j.d.'s which
is continuous in probability, with state space $S$ and parameter
set $[0,1]$. We will find conditions on the f.j.d.'s in $F$ under
which an a.u. continuous process $X$ can be constructed with marginal
distributions given by the family $F$. 

The classical approach to the existence of such processes $X$, as
elaborated in $\;$ $\;$ $\;$ \cite{Billingsley74}, uses the following
theorem. 
\begin{thm}
\label{Thm. Prokhorov's Relative Compactness Theorem} \textbf{\emph{(Prokhorov's
Relative Compactness Theorem). }}Each tight family $\overline{J}$
of distributions on a locally compact metric space $(H,d_{H})$ is
relative compact, in the sense that each sequence in $\overline{J}$
contains a subsequence which converges weakly to some distribution
on $(H,d_{H})$.
\end{thm}
Prokhorov's theorem however implies the principle of infinite search,
and is therefore not constructive. This can be seen as follows. Let
$(r_{n})_{n=1,2,\cdots}$ be an arbitrary nondecreasing sequence in
$H\equiv\{0,1\}$. Let the doubleton $H$ be endowed with the Euclidean
metric $d_{H}$ defined by $d_{H}(x,y)=|x-y|$ for each $x,y\in H$.
For each $n\geq1$, let $J_{n}$ be the distribution on $(H,d_{H})$
which assigns unit mass to $r_{n}$; in other words, $J_{n}(f)\equiv f(r_{n})$
for each $f\in C(H,d_{H})$. Then the family $\overline{J}\equiv\{J_{1},J_{2},\cdots\}$
is tight, and Prokhorov's theorem implies that $J_{n}$ converges
weakly to some distribution $J$ on $(H,d_{H})$. It follows that
$J_{n}g$ converges as $n\rightarrow0$, where $g\in C(H,d_{H})$
is defined by $g(x)=x$ for each $x\in H$. Thus $r_{n}\equiv g(r_{n})\equiv J_{n}g$
converges as $n\rightarrow0$. Since $(r_{n})_{n=1,2,\cdots}$ is
an arbitrary nondecreasing sequence in $\{0,1\}$, the principle of
infinite search follows from Prokhorov's theorem.

In our constructions, we will bypass any use of Prokhorov's theorem
or to any unjustified supremums, in favor of direct proofs using Borel-Cantelli
estimates. We will give a necessary and sufficient condition on the
f.j.d.'s in the family $F$, for $F$ to be extendable to an a.u.
continuous process $X$. We will call this condition \emph{$C$-}regularity.
We will derive a modulus of a.u. continuity of the process $X$ from
a given modulus of continuity in probability and a given modulus of
$C$-regularity of the consistent family $F$, to be defined presently.
We will also prove that the extension is uniformly metrically continuous
on an arbitrary set $\widehat{F}_{0}$ of such consistent families
$F$ which share a common modulus of $C$-regularity. 

In essence, the material presented in Sections 1 and 2 of the present
work is a constructive and more general version of, materials from
Section 7 of Chapter 2 of \cite{Billingsley74}, the latter treating
only the special case where $S=R$. We remark that the generalization
to the arbitrary locally compact state space $(S,d)$ is not entirely
trivial, because we forego the convenience of linear interpolation
in $R$.

A subsequent chapter in the present work will introduce a condition,
analogous to $C$-regularity, for the treatment of processes which
are, almost uniformly, right continuous with left limits, again with
a general locally compact metric space as state space.

In Section 3, we will prove a generalization of Kolmogorov's theorem
for a.u. locally Hoelder continuity, in a sense to be made precise
in Section 3, with state space $(S,d)$. 

Separately, in Section 4, in the case of Gaussian processes, we will
present the sufficient condition and the proof in \cite{Garsia 70}
for the construction of an a.u. continuous process given the modulus
of continuity of the covariance function. A minor modification of
their proof makes it strictly constructive.

We note that, for a more general parameter space which is a subset
of $R^{m}$ for some $m\geq0$, with some restriction on its local
$\varepsilon$-entropy, \cite{Potthoff09-1} gives sufficient conditions
on the pair distributions $F_{t,s}$ to guarantee the construction
of an a.u. continuous or an a.u. locally Hoelder, real-valued, random
field. 

In this and later chapters we will use the following notations for
the dyadic rationals. 
\begin{defn}
\textbf{\label{Def. Notations for dyadic rationals} (Notations for
dyadic rationals).} For each $m\geq0$, define $p_{m}\equiv2^{m}$,
$\Delta_{m}\equiv2^{-m}$, and define the enumerated set of dyadic
rationals
\[
Q_{m}\equiv\{t_{0},t_{1},\cdots,t_{p(m)}\}=\{q_{m,0},\cdots,q_{m,p(m)}\}
\]
\[
\equiv\{0,\Delta_{m},2\Delta_{m},\cdots,1\}\subset[0,1],
\]
and 
\[
Q_{\infty}\equiv\bigcup_{m=0}^{\infty}Q_{m}\equiv\{t_{0},t_{1},\cdots\},
\]
where the second equality is equality of sets without the enumeration.
Let $m\geq0$ be arbitrary. Then the enumerated set $Q_{m}$ is a
$2^{-m}$-approximation of $[0,1]$, with $Q_{m}\subset Q_{m+1}$.
Conditions in Definition \ref{Def. Binary approximationt and Modulus of local compactness}
can easily be verified for the sequence 
\[
\xi_{[0,1]}\equiv(Q_{m})_{m=1,2,\cdots}
\]
to be a binary approximation of $[0,1]$ relative to the reference
point $q_{\circ}\equiv0$. 

In addition, for each $m\geq0$, define the enumerated set of dyadic
rationals
\[
\overline{Q}_{m}\equiv\{u_{0},u_{1},\cdots,u_{p(2m)}\}\equiv\{0,2^{-m},2\cdot2^{-m},\cdots,2^{m}\}\subset[0,2^{m}],
\]
and 
\[
\overline{Q}_{\infty}\equiv\bigcup_{m=0}^{\infty}\overline{Q}_{m}\equiv\{u_{0},u_{1},\cdots\},
\]
where the second equality is equality of sets without the enumeration. 
\end{defn}
$\square$
\begin{defn}
\textbf{\label{Def. Misc notations}(Miscellaneous Notations and conventions).
}As usual, to lighten notational burden, we will write an arbitrary
subscripted symbol $a_{b}$ interchangeably with $a(b)$. We will
write $TU$ for a composite function $T\circ U$. If $f:A\rightarrow B$
is a function from a set $A$ to a set $B$, and if $A'$ is a nonempty
subset $A$, then the restricted function $f|A':A'\rightarrow B$
will also be denoted simply by $f$ when there is little risk of confusion.
If $A$ is a $\mathrm{measurable}$ subset on a probability space
$(\Omega,L,E)$, then we will write $PA$, $P(A)$, $E(A)$, $EA$,
or $E1_{A}$, interchangeably. For arbitrary r.r.v. $Y\in L$ and
$\mathrm{measurable}$ subsets $A,B$, we will write $(A;B)\equiv AB$,
$E(Y;A)\equiv E(Y1_{A})$, and $A\in L$. For further abbreviations,
we drop parentheses when there is little risks of confusion. For example,
we write $1_{Y\leq\beta;Z>\alpha}\equiv1_{(Y\leq\beta)(Z>\alpha)}$.
.For an arbitrary integrable function $X\in L$,  will sometimes use
the more suggestive notation $\int E(d\omega)X(\omega)$ for $EX$,
where $\omega$ is a dummy variable.

Let $Y$ be an arbitrary r.r.v. Recall from Definition \ref{Convention. Only Regular Pts r.r.v.'s}
the convention that if measurability of the set $(Y<\beta)$ or $(Y\leq\beta)$
is required in a discussion, for some $\beta\in R$, then it is understood
that the real number $\beta$ has been chosen from the regular points
of the r.r.v. $Y$. 

Recall that $[\cdot]_{1}$ is an operation which assigns to each $a\in R$
an integer $[a]_{1}\in(a,a+2)$.$\square$
\end{defn}

\section{Extension of a.u. Continuous processes with dyadic rational parameters
to parameters in $[0,1]$}

Our approach to extend a given family $F$ of f.j.d.'s which is continuous
in probability on the parameter set $[0,1]$ is as follows. First
note that $F$ carries no more useful information than its restriction
$F|Q_{\infty}$, where $Q_{\infty}$ is the dense subset of dyadic
rationals in $[0,1]$, because the family can be recovered from the
$F|Q_{\infty}$, thanks to continuity in probability. Hence we can
first extend the family $F|Q_{\infty}$ to a process $Z:Q_{\infty}\times\Omega\rightarrow S,$
by the Daniell-Kolmogorov Theorem or the Daniell-Kolmogorov-Skorokhod
Theorem. Then any condition of the family $F$ is equivalent to a
condition to $Z$.

In particular, in the current context, any condition on f.j.d.'s to
make $F$ extendable to an a..u. continuous process $X:[0,1]\times\Omega\rightarrow S$
can be discussed in the more general terms of a process $Z:Q_{\infty}\times\Omega\rightarrow S,$
the latter to be extended by limit to a process $X$. It is intuitively
obvious that any a.u. continuous process $Z:Q_{\infty}\times\Omega\rightarrow S$
is extendable to an a.u. continuous process $X:[0,1]\times\Omega\rightarrow S$,
because $Q_{\infty}$ is dense $[0,1]$. In this section, we will
make this precise, and prove that the extension construction is itself
a metrically continuous construction. 
\begin{defn}
\label{Def. C_hat and rho_C_hat} \textbf{(Metric Space of a.u. Continuous
Processes). }Let $C[0,1]$ be the space of continuous functions $x:[0,1]\rightarrow S$,
endowed with the uniform\emph{ metric}\index{uniform metric} defined
by 
\begin{equation}
d_{C[0,1]}(x,y)\equiv\sup_{t\in[0,1]}d(x(t),y(t))\label{eq:temp-520}
\end{equation}
for each $x,y\in C[0,1]$. Write $\widehat{d}_{C[0,1]}\equiv1\wedge d_{C[0,1]}$.Let
$\widehat{C}[0,1]$ denote the set of stochastic processes $X:[0,1]\times\Omega\rightarrow S$
which are a.u. continuous on $[0,1]$. Define a metric $\rho_{\widehat{C}[0,1]}$
on $\widehat{C}[0,1]$ by
\begin{equation}
\rho_{\widehat{C}[0,1]}(X,Y)\equiv E\sup_{t\in[0,1]}\widehat{d}(X_{t},Y_{t})\equiv E\widehat{d}_{C[0,1]}(X,Y)\label{eq:temp-202}
\end{equation}
for each $X,Y\in\widehat{C}[0,1]$. The next lemma says that $(\widehat{C}[0,1],\rho_{\widehat{C}[0,1]})$
is a well-defined metric space. $\square$
\end{defn}
\begin{lem}
\label{Lem. rho_C_hat is a metric} \textbf{\emph{($\rho_{\widehat{C}[0,1]}$
is a metric).}} The function $\sup_{t\in[0,1]}\widehat{d}(X_{t},Y_{t})$
is a r.r.v.. The function $\rho_{\widehat{C}[0,1]}$ is well-defined
and is a metric.
\end{lem}
\begin{proof}
Let $X,Y\in\widehat{C}[0,1]$ be arbitrary, with moduli of a.u. continuity
$\delta_{auc}^{X},\delta_{auc}^{Y}$ respectively. First note that
the function $\sup_{t\in[0,1]}\widehat{d}(X_{t},Y_{t})$ is defined
a.s., on account of continuity on $[0,1]$ of $X,Y$ on a full subset
of $\Omega$. We need to prove that it is $\mathrm{measurable}$,
so that the expectation in the defining formula \ref{eq:temp-202}
makes sense.

To that end, let $\varepsilon>0$ be arbitrary. Then there exist $\mathrm{measurable}$
sets $D_{X},D_{Y}$ with $P(D_{X}^{c})\vee P(D_{Y}^{c})<\varepsilon$
such that 
\[
d(X_{t},X_{s})\vee d(Y_{t},Y_{s})\leq\varepsilon
\]
for each $t,s\in[0,1]$ with $|t-s|<\delta\equiv\delta_{auc}^{X}(\varepsilon)\wedge\delta_{auc}^{Y}(\varepsilon)$.
Now let the sequence $s_{0},\cdots,s_{n}$ be an arbitrary $\delta$-approximation
of $[0,1]$, for some $n\geq1$. Then, for each $t\in[0,1]$, we have
$|t-s_{k}|<\delta$ for some $k=0,\cdots,n$, whence
\begin{equation}
d(X_{t},X_{s(k)})\vee d(Y_{t},Y_{s(k)})\leq\varepsilon\label{eq:temp-274}
\end{equation}
on $D\equiv D_{X}D_{Y}$, which in turn implies
\[
|d(X_{t},Y_{t})-d(X_{s(k)},Y_{s(k)})|\leq2\varepsilon
\]
on $D$. It follows that
\[
|\sup_{t\in[0,1]}\widehat{d}(X_{t},Y_{t})-\bigvee_{k=1}^{n}\widehat{d}(X_{s(k)},Y_{s(k)})|\leq2\varepsilon
\]
on $D$, where $Z\equiv\bigvee_{k=1}^{n}\widehat{d}(X_{t(k)},Y_{t(k)})$
is a r.r.v., and where $P(D^{c})<2\varepsilon$. For each $p\geq1$,
we can repeat this argument with $\varepsilon\equiv\frac{1}{p}$.
Thus we obtain a sequence $Z_{p}$ of r.r.v.'s with $Z_{p}\rightarrow\sup_{t\in[0,1]}\widehat{d}(X_{t},Y_{t})$
in probability as $p\rightarrow\infty$. The function $\sup_{t\in[0,1]}\widehat{d}(X_{t},Y_{t})$
is accordingly a r.r.v., and, being bounded by $1$, integrable. Summing
up, the expectation in equality \ref{eq:temp-202} exists, and $\rho_{\widehat{C}[0,1]}$
is well-defined.

Verification of the conditions for the function $\rho_{\widehat{C}[0,1]}$
to be a metric is straightforward and omitted.
\end{proof}
$\square$
\begin{defn}
\label{Def. Extension by limit from Q_inf to =00005B0,1=00005D}(\textbf{Extension
by limit of a process with parameter set $Q_{\infty}$).} Let $Z:Q_{\infty}\times\Omega\rightarrow S$
be an arbitrary process. Define a function by 
\[
domain(X)\equiv\Phi_{Lim}(Z)\equiv\{(r,\omega)\in[0,1]\times\Omega:\lim_{s\rightarrow r;s\in Q(\infty)}Z(s,\omega)\;\mathrm{exists\}}
\]
and 
\[
X(r,\omega)\equiv\lim_{s\rightarrow r;s\in Q(\infty)}Z(s,\omega)
\]
for each $(r,\omega)\in domain(X)$. We will call $X$ the \emph{extension-by-limit
of the process $Z$ to the parameter set $[0,1]$.} A similar definition
is made where the interval $[0,1]$ is replaced by the interval $[0,\infty)$,
and where the set $Q_{\infty}$ of dyadic rationals in $[0,1]$ is
replaced by the set $\overline{Q}_{\infty}$ of dyadic rationals in
$[0,\infty)$. 

We emphasize that, absent any additional conditions on the process
$Z$, the function $X$ need not be a process; it need not even be
a well-defined function.
\end{defn}
$\square$
\begin{thm}
\textbf{\emph{\label{Thm. Completion of process on dense parameter subset}
(Extension by limit of a.u. continuous process on $Q_{\infty}$ to
a.u. continuous process on $[0,1]$; and metrical continuity of said
extension).}} Let $\widehat{R}_{0}$ be a subset of $\widehat{R}(Q_{\infty}\times\Omega,S)$
whose members are a.u. continuous with a common modulus of a.u. continuity
$\delta_{auc}$. Then the following holds. 

1. Let $Z\in\widehat{R}_{0}$ be arbitrary\emph{.} Then its extension-by-limit
$X\equiv\Phi_{Lim}(Z)$ is an a.u. continuous process such that $X_{t}=Z_{t}$
on $domain(X_{t})$ for each $t\in Q_{\infty}$. Moreover, the process
$X$ has the same modulus of a.u. continuity $\delta_{auc}$ as $Z$. 

2. Recall that $(\widehat{R}(Q_{\infty}\times\Omega,S),\rho_{Prob})$
is the metric space of processes $Z:Q_{\infty}\times\Omega\rightarrow S$.
The extension-by-limit
\[
\Phi_{Lim}:(\widehat{R}_{0},\rho_{Prob})\rightarrow(\widehat{C}[0,1],\rho_{\widehat{C}[0,1]})
\]
is uniformly continuous, with a modulus of continuity $\delta_{Lim}(\cdot,\delta_{auc})$.
\end{thm}
\begin{proof}
1. Let $Z\in\widehat{R}_{0}$ be arbitrary\emph{.} Let $\varepsilon>0$
be arbitrary. Then, by hypothesis, there exists $\delta_{auc}(\varepsilon)>0$
and a $\mathrm{measurable}$ set $D\subset\Omega$ with $P(D^{c})<\varepsilon$
such that, for each $\omega\in D$ and for each $s,s'\in Q_{\infty}$
with $|s-s'|<\delta_{auc}(\varepsilon)$, we have 
\begin{equation}
d(Z(s,\omega),Z(s',\omega))\leq\varepsilon.\label{eq:temp-511}
\end{equation}

Next let $\omega\in D$ and $r,r'\in[0,1]$ be arbitrary with $|r-r'|<\delta_{auc}(\varepsilon)$.
Letting $s.s'\rightarrow r$ with $s,s'\in Q_{\infty}$, we have $|s-s'|\rightarrow0$,
and so $d(Z(s,\omega),Z(s',\omega))\rightarrow0$ as $s,s'\rightarrow r$.
Since $(S,d)$ is complete, we conclude that the limit 
\[
X(r,\omega)\equiv\lim_{s\rightarrow r;s\in Q(\infty)}Z(s,\omega)
\]
exists. Moreover, letting $s'\rightarrow r$ with $s,s'\in Q_{\infty}$,
inequality \ref{eq:temp-511} yields 
\begin{equation}
d(Z(s,\omega),X(r,\omega))\leq\varepsilon.\label{eq:temp-511-1}
\end{equation}
Since $\varepsilon>0$ is arbitrary, we see that $Z_{s}\rightarrow X_{r}$
a.u. as $s\rightarrow r$. Hence $X_{r}$ is a r.v. Thus $X:[0,1]\times\Omega_{0}\rightarrow S$
is a stochastic process.

Now let $s\rightarrow r$ and $s'\rightarrow r'$ with $s,s'\in Q_{\infty}$
in inequality \ref{eq:temp-511}. Then we obtain 
\begin{equation}
d(X(r,\omega),X(r',\omega))\leq\varepsilon,\label{eq:temp-498}
\end{equation}
where $\omega\in D$ and $r,r'\in[0,1]$ are arbitrary with $|r-r'|<\delta_{auc}(\varepsilon)$.
Thus $X$ has the same modulus of a.u. continuity $\delta_{auc}$
as $Z$. 

2. It remains to verify that the mapping $\Phi_{Lim}$ is a continuous
function. To that end, let $\varepsilon>0$ be arbitrary. Write $\alpha\equiv\frac{1}{3}\varepsilon$.
Let $m\equiv m(\varepsilon,\delta_{auc})\geq1$ be so large that $2^{-m}<\delta_{auc}(\alpha)$.
Define 
\[
\delta_{Lim}(\varepsilon,\delta_{auc})\equiv2^{-p(m)-1}\alpha^{2}.
\]
Let $Z,Z'\in\widehat{R}_{0}$ be arbitrary such that
\[
\rho_{Prob}(Z,Z')<\delta_{Lim}(\varepsilon,\delta_{auc}).
\]
Equivalently,
\begin{equation}
E\sum_{n=0}^{\infty}2^{-n-1}\widehat{d}(Z_{t(n)},Z'_{t(n)})<2^{-p(m)-1}\alpha^{2}.\label{eq:temp-512}
\end{equation}
Then, by Chebychev's inequality, there exists a $\mathrm{measurable}$
set $A$ with $P(A^{c})<\alpha$ such that, for each $\omega\in A$,
we have
\[
\sum_{n=0}^{\infty}2^{-n-1}\widehat{d}(Z(t_{n},\omega),Z'(t_{n},\omega))<2^{-p(m)-1}\alpha,
\]
whence
\begin{equation}
\widehat{d}(Z(t_{n},\omega),Z'(t_{n},\omega))<\alpha\label{eq:temp-513}
\end{equation}
for each $n=0,\cdots,p_{m}$. 

Now let $X\equiv\Phi_{Lim}(Z)$ and $X'\equiv\Phi_{Lim}(Z')$. By
Assertion 1, the processes $X$ and $X'$ have the same modulus of
a.u. continuity $\delta_{auc}$ as $Z$ and $Z'$. Hence, there exists
$\mathrm{measurable}$ sets $D,D'$ with $P(D^{c})\vee P(D'^{c})<\alpha$
such that, for each $\omega\in DD'$, we have 
\begin{equation}
\widehat{d}(X(r,\omega),X(s,\omega))\vee\widehat{d}(X'(r,\omega),X'(s,\omega))\leq\alpha\label{eq:temp-498-1}
\end{equation}
for each $r,s\in[0,1]$ with $|r-s|<\delta_{auc}(\alpha)$. 

Now consider each $\omega\in ADD'.$ Let $r\in[0,1]$ be arbitrary.
Since $t_{0},\cdots,t_{p(m)}$ is a $2^{-m}$-approximation of $[0,1]$,
there exists $n=0,\cdots,p_{m}$ such that $|r-t_{n}|<2^{-m}<\delta_{auc}(\alpha)$.
Then inequality \ref{eq:temp-498-1} holds with $s\equiv t_{n}$.
Combining inequalities \ref{eq:temp-513} and \ref{eq:temp-498-1},
we obtain
\[
\widehat{d}(X(r,\omega),X'(r,\omega))
\]
\[
\leq\widehat{d}(X(r,\omega),X(s,\omega))+\widehat{d}(X(s,\omega),X'(s,\omega))+\widehat{d}(X'(r,\omega),X'(s,\omega))
\]
\[
=\widehat{d}(X(r,\omega),X(s,\omega))+\widehat{d}(Z(s,\omega),Z'(s,\omega))+\widehat{d}(X'(r,\omega),X'(s,\omega))<3\alpha,
\]
where $\omega\in ADD'$ and $r\in[0,1]$ are arbitrary. It follows
that 
\[
\rho_{\widehat{C}[0,1]}(X,X')\equiv E\sup_{r\in[0,1]}\widehat{d}(X_{r},X'_{r})
\]
\[
\leq E\sup_{r\in[0,1]}\widehat{d}(X_{r},X'_{r})1_{ADD'}+P(ADD')^{c}<3\alpha+3\alpha=6\alpha\equiv\varepsilon.
\]
We conclude that $\delta_{Lim}(\cdot,\delta_{auc})$ is a modulus
of continuity of $\Phi_{Lim}$.
\end{proof}

\section{$C$-regular families of f.j.d.'s and $C$-regular processes }
\begin{defn}
\label{Def.  C-regularity of F} \textbf{(C-regularity).} Let $(\Omega,L,E)$
be an arbitrary sample space. Let $Z:Q_{\infty}\times\Omega\rightarrow S$
be an arbitrary process. We will say that $Z$ is a \emph{C-regular}
\index{C-regularity} process if there exists an increasing sequence
$\overline{m}\equiv(m_{n})_{n=0,1,\cdots}$ of positive integers,
called the \index{modulus of C-regularity}\emph{ modulus of C-regularity
of the process $Z$, such that}, for each $n\geq0$ and for each $\beta_{n}>2^{-n}$
such that the set 
\[
A_{t,s}^{(n)}\equiv(d(Z_{t},Z_{s})>\beta_{n})
\]
is $\mathrm{measurable}$ for each $s,t\in Q_{\infty}$, we have 
\begin{equation}
P(C_{n})<2^{-n}\label{eq:temp-204}
\end{equation}
where 
\begin{equation}
C_{n}\equiv\bigcup_{t\in Q(m(n))}\bigcup_{s\in[t,t']Q(m(n+1))}A_{t,s}^{(n)}\cup A_{s,t'}^{(n)},\label{eq:temp-201-4}
\end{equation}
where, for each $t\in Q_{m(n)}$, we abuse notations and write $t'\equiv1\wedge(t+2^{-m(n)})\in Q_{m(n)}$. 

Let $F$ be a consistent family of f.j.d.'s which is continuous in
probability on $[0,1]$. Then the family $F$ of consistent f.j.d.'s
is said to be \emph{C-regular} \index{C-regularity}, with the sequence
$\overline{m}\equiv(m_{n})_{n=0,1,\cdots}$ as a \index{modulus of C-regularity}\emph{
modulus of C-regularity if $F|Q_{\infty}$ is family of marginal distributions
of some} \emph{C-regular process $Z$.}
\end{defn}
$\square$

We will prove that a process on $[0,1]$ is a.u. continuous iff it
is $C$-regular. Note that $C$-regularity is a condition on the f.j.d.'s
while a.u. continuity is a condition on sample paths. 
\begin{thm}
\label{Thm.  a.u. continuoity on Q0 =00003D> C-regularity} \textbf{\emph{(a.u.
Continuity implies $C$-regularity).}} Let $(\Omega,L,E)$ be an arbitrary
sample space. Let $X:[0,1]\times\Omega\rightarrow S$ be an a.u. continuous
process, with a modulus of a.u. continuity $\delta_{auc}$. Then the
process $X$ is C-regular, with a modulus of $C$-regularity given
by $\overline{m}\equiv(m_{n})_{n=0,1,\cdots}$, where $m_{0}\equiv1$
and
\[
m_{n}\equiv[m_{n-1}\vee(-\log_{2}\delta_{auc}(2^{-n-1}))]_{1}
\]
for each $n\geq1$. 
\end{thm}
\begin{proof}
First note that $X$ is continuous in probability. Separately, let
$n\geq0$ be arbitrary. By Definition \ref{Def. continuity in prob, continuity a.u., and a.u. continuity}
of a.u. continuity, there exists a $\mathrm{measurable}$ set $D_{n}$
with $P(D_{n}^{c})<2^{-n}$ such that, for each $\omega\in D_{n}$
and for each $s,t\in Q_{\infty}$ with $|t-s|<\delta_{auc}(2^{-n})$,
we have
\begin{equation}
d(X(t,\omega),X(s,\omega))\leq2^{-n}.\label{eq:temp-134}
\end{equation}
Let $\beta_{n}>2^{-n}$ be arbitrary and let 
\[
A_{t,s}^{(n)}\equiv(d(X_{t},X_{s})>\beta_{n})
\]
for each $s,t\in Q_{\infty}$. Define 
\begin{equation}
C_{n}\equiv\bigcup_{t\in Q(m(n))}\bigcup_{s\in[t,t']Q(m(n+1))}A_{t,s}^{(n)}\cup A_{s,t'}^{(n)},\label{eq:temp-201-4-2}
\end{equation}
where, as before, for each $t\in Q_{m(n)}$ we abuse notations and
write $t'\equiv1\wedge(t+2^{-m(n)})\in Q_{m}$. Suppose, for the sake
of a contradiction, that $P(D_{n}C_{n})>0$. Then there exists some
$\omega\in D_{n}C_{n}$. Hence, by equality \ref{eq:temp-201-4-2},
there exists $t\in Q_{m(n)}$ and $s\in[t,t']Q_{m(n+1)}$ with 
\begin{equation}
d(X(t,\omega),X(s,\omega))\vee d(X(s,\omega),X(t',\omega))>\beta_{n}.\label{eq:temp-302}
\end{equation}
It follows that 
\begin{equation}
|s-t|\vee|s-t'|\leq2^{-m(n)}<\delta_{auc}(2^{-n})\label{eq:temp-294}
\end{equation}
Inequalities \ref{eq:temp-294} and \ref{eq:temp-134} together imply
that
\[
d(X(t,\omega),X(s,\omega))\vee d(X(s,\omega),X(t',\omega))\leq2^{-n}<\beta{}_{n},
\]
contradicting inequality \ref{eq:temp-302}. We conclude that $P(D_{n}C_{n})=0$.
Consequently, 
\[
P(C_{n})=P(D_{n}\cup D_{n}^{c})C_{n}=P(D_{n}^{c}C_{n})\leq P(D_{n}^{c})<2^{-n}.
\]
Thus the conditions in Definition \ref{Def.  C-regularity of F} are
satisfied for the family $F$ of marginal distributions of $X$ to
be $C$-regular, with modulus of $C$-regularity given by $\overline{m}$.
\end{proof}
The next theorem is the converse of Theorem \ref{Thm.  a.u. continuoity on Q0 =00003D> C-regularity},
and is the main theorem in this section. 
\begin{thm}
\label{Thm. C-regular =00003D> a.u. continuity} \textbf{\emph{($C$-regularity
implies a.u. continuity).}} Let $(\Omega,L,E)$ be an arbitrary sample
space. Let $F$ be a $C$-regular family of consistent f.j.d.'s. Then
there exists an a.u. continuous process $X:[0,1]\times\Omega\rightarrow S$
with marginal distributions given by $F$.

Specifically, let $\overline{m}\equiv(m_{n})_{n=0,1,\cdots}$ be a
modulus of C-regularity of $F$. Let $Z:Q_{\infty}\times\Omega\rightarrow S$
be an arbitrary process with marginal distributions given by $F|Q_{\infty}$.
Let $\varepsilon>0$ be arbitrary. Define $h\equiv[0\vee(4-\log_{2}\varepsilon)]_{1}$
and $\delta_{auc}(\varepsilon,\overline{m})\equiv2^{-m(h)}$. Then
$\delta_{auc}(\cdot,\overline{m})$ is a modulus of a.u. continuity
of $Z$.

Moreover, the extension-by-limit $X\equiv\Phi_{Lim}(Z):[0,1]\times\Omega\rightarrow S$
of the process $Z$ to the full parameter set $[0,1]$ is a.u. continuous,
with the same modulus of a.u. continuity $\delta_{auc}(\cdot,\overline{m})$,
and with marginal distributions given by $F$.
\end{thm}
\begin{proof}
1. First let $n\geq0$ be arbitrary. Take any $\beta_{n}\in(2^{-n},2^{-n+1})$.
Then, by Definition \ref{Def.  C-regularity of F}, 
\begin{equation}
P(C_{n})<2^{-n}\label{eq:temp-204-1}
\end{equation}
where 
\begin{equation}
C_{n}\equiv\bigcup_{t\in Q(m(n))}\bigcup_{s\in[t,t']Q(m(n+1))}(d(Z_{t},Z_{s})>\beta_{n})\cup(d(Z_{s},Z_{t'})>\beta_{n}),\label{eq:temp-201-4-4}
\end{equation}
where, as before, for each $t\in Q_{m(n)}$ we abuse notations and
write $t'\equiv1\wedge(t+2^{-m(n)})\in Q_{m(n)}$.

2. Now  define $D_{n}\equiv(\bigcup_{j=n}^{\infty}C_{j})^{c}$. Then
\[
P(D_{n}^{c})\leq\sum_{j=n}^{\infty}P(C_{j})<\sum_{j=n}^{\infty}2^{-j}=2^{-n+1}.
\]
Consider each $\omega\in D_{n}$. Consider each $t\in Q_{m(n)}$.
For each $s\in[t,t']Q_{m(n+1)}$, since $\omega\in C_{n}^{c}$, we
have 
\[
d(Z(t,\omega),Z(s,\omega))\leq\beta_{n}<2^{-n+1}.
\]
In short
\begin{equation}
[t,t']Q_{m(n+1)}\subset(d(Z(\cdot,\omega),Z(t,\omega))<2^{-n+1}),\label{eq:temp-201}
\end{equation}
where $t\in Q_{m(n)}$ is arbitrary. Repeating the above argument
with $n$ replaced by $n+1$ and with $t$ replaced by each $s\in[t,t']Q_{m(n+1)}$,
we obtain
\begin{equation}
[s,s']Q_{m(n+2)}\subset(d(Z(\cdot,\omega),Z(s,\omega))<2^{-n}),\label{eq:temp-206}
\end{equation}
where $s'\equiv1\wedge(s+2^{-m(n+1)})\in Q_{m(n+1)}$. Since 
\[
[t,t']Q_{m(n+2)}=[t,t']Q_{m(n+1)}\cap\bigcup_{s\in[t,t']Q(m(n+1))}[s,s']Q_{m(n+2)},
\]
relations \ref{eq:temp-201} and \ref{eq:temp-206} together yield
\[
[t,t']Q_{m(n+2)}\subset(d(Z(\cdot,\omega),Z(t,\omega))<2^{-n+1}+2^{-n}).
\]
Inductively with $k=n+1,n+2,\cdots$, we obtain
\[
[t,t']Q_{m(k)}\subset(d(Z(\cdot,\omega),Z(t,\omega))<2^{-n+1}+2^{-n}+\cdots+2^{-k+2})
\]
\[
\subset(d(Z(\cdot,\omega),Z(t,\omega))<2^{-n+2})
\]
for each $k\geq n+1$. Therefore 
\[
[t,t']Q_{\infty}=[t,t']\bigcup_{k=n+1}^{\infty}Q_{m(k)}\subset(d(Z(\cdot,\omega),Z(t,\omega))<2^{-n+2}).
\]
In particular $d(Z(t',\omega),Z(t,\omega))<2^{-n+2}$, and so the
last displayed condition implies
\begin{equation}
[t,t']Q_{\infty}\subset(d(Z(\cdot,\omega),Z(t,\omega))\vee d(Z(\cdot,\omega),Z(t',\omega))<2^{-n+3}),\label{eq:temp-126}
\end{equation}
where $n\geq1$, $\omega\in D_{n}$, and $t\in Q_{m(n)}$ are arbitrary.

3. Continuing with arbitrary $n\geq1$, suppose $r,s\in Q_{\infty}$
are arbitrary such that $0<s-r<2^{-m(n)}$. Then there exist $t,u\in Q_{m(n)}$
with $t\leq u$ such that $r\in[t,t']$ and $s\in[u,u']$, where $u'\equiv1\wedge(u+2^{-m(n)})\in Q_{m(n)}$.
If $t'<u$ then $s\geq u\geq t'+2^{-m(n)}\geq r+2^{-m(n)}$, a contradiction.
Hence $u\leq t'$. On the other hand, $t\leq u$ by the choice of
$t,u$. Consequently, $u=t$ or $u=t'$. At the same time, according
to relation \ref{eq:temp-126}, we have 
\[
d(Z(r,\omega),Z(t,\omega))\vee d(Z(r,\omega),Z(t',\omega))<2^{-n+3}.
\]
Similarly,
\[
d(Z(s,\omega),Z(u,\omega))\vee d(Z(s,\omega),Z(u',\omega))<2^{-n+3}.
\]
If $u=t$, then it follows that
\[
d(Z(r,\omega),Z(s,\omega))\leq d(Z(r,\omega),Z(t,\omega))+d(Z(u,\omega),Z(s,\omega))<2^{-n+4}.
\]
If $u=t'$, then similarly
\[
d(Z(r,\omega),Z(s,\omega))\leq d(Z(r,\omega),Z(t',\omega))+d(Z(u,\omega),Z(s,\omega))<2^{-n+4}
\]

Summing up, for each $\omega\in D_{n}$, for each $n\geq1$, and for
each $r,s\in Q_{\infty}$ with $0<s-r<2^{-m(n)}$, we have 
\begin{equation}
d(Z(r,\omega),Z(s,\omega))<2^{-n+4}.\label{eq:temp-231}
\end{equation}
By symmetry, the last inequality therefore holds for each $\omega\in D_{n}$,
for each $n\geq1$, and for each $r,s\in Q_{\infty}$ with $|s-r|<2^{-m(n).}$. 

4. Now let $\varepsilon>0$ be arbitrary. Let $n\equiv[0\vee(4-\log_{2}\varepsilon)]_{1}$
and $\delta_{auc}(\varepsilon,\overline{m})\equiv2^{-m(n)}$, as in
the hypothesis. By the previous paragraphs, we see that the $\mathrm{measurable}$
set $D_{n}\equiv(\bigcup_{j=n}^{\infty}C_{j})^{c}$ is such that $P(D_{n}^{c})\leq2^{-n+1}<\varepsilon$
and such that $d(Z(r,\omega),Z(s,\omega))<2^{-n+4}<\varepsilon$ for
each $r,s\in Q_{\infty}$ with $|s-r|<\delta_{auc}(\varepsilon,\overline{m})\equiv2^{-m(n).}$.
Thus the process $Z$ is a.u. continuous, with $\delta_{auc}(\cdot,\overline{m})$
as a modulus of a.u. continuity of $Z$. 

5. By Proposition \ref{Thm. Completion of process on dense parameter subset},
the complete extension $X\equiv\Phi_{Lim}(Z)$ of the process $Z$
to the full parameter set $[0,1]$ is a.u. continuous with the same
modulus of a.u. continuity $\delta_{auc}(\cdot,\overline{m})$. 

The theorem is proved.
\end{proof}
\begin{thm}
\textbf{\emph{\label{Thm. Continuity of a.u. continuous extension}
(Continuity of extension-by-limit of C-regular processes).}} Recall
the metric space $(\widehat{R}(Q_{\infty}\times\Omega,S),\rho_{Prob})$
of stochastic processes with parameter set $Q_{\infty}\equiv\{t_{0},t_{1},\cdots\}$,
sample space $(\Omega,L,E)$, and state space $(S,d)$. Let $\widehat{R}_{0}$
be a $C$-equiregular subset of $\widehat{R}(Q_{\infty}\times\Omega,S)$
with a modulus of $C$-regularity $\overline{m}\equiv(m_{n})_{n=0,1,\cdots}$.
Let $(\widehat{C}[0,1],\rho_{\widehat{C}[0,1]})$ be the metric space
of a.u. continuous processes on $[0,1]$, as in Definition \ref{Def. C_hat and rho_C_hat}. 

Then the extension-by-limit
\[
\Phi_{Lim}:(\widehat{R}_{0},\rho_{Prob})\rightarrow(\widehat{C}[0,1],\rho_{\widehat{C}[0,1]})
\]
as in Definition \ref{Def. Extension by limit from Q_inf to =00005B0,1=00005D},
is uniformly continuous, with a modulus of continuity $\delta_{reg,auc}(\cdot,\overline{m})$
\end{thm}
\begin{proof}
Let $\varepsilon>0$ be arbitrary. Define
\[
j\equiv[0\vee(4-\log_{2}\frac{\varepsilon}{4})]_{1},
\]
\[
h_{j}\equiv2^{m(j)},
\]
and 
\[
\delta_{reg,auc}(\varepsilon,\overline{m})\equiv2^{-h(j)-2j-j}.
\]
We will prove that $\delta_{reg,auc}(\cdot,\overline{m})$ is a modulus
of continuity of $\Phi_{Lim}$ on $\widehat{R}_{0}$. 

1. Let $Z,Z'\in\widehat{R}_{0}$ be arbitrary and let $X\equiv\Phi_{Lim}(Z)$,
$X'\equiv\Phi_{Lim}(Z')$. Suppose 
\begin{equation}
\rho_{prob,Q(\infty)}(Z,Z')\equiv E\sum_{n=0}^{\infty}2^{-n-1}\widehat{d}(Z_{t(n)},Z_{t(n)}')<\delta_{reg,auc}(\varepsilon,\overline{m}).\label{eq:temp-207}
\end{equation}
We need to prove that $\rho_{\widehat{C}[0,1]}(X,X')<\varepsilon$. 

To that end, first note that, by Step 4 in the proof of Theorem \ref{Thm. C-regular =00003D> a.u. continuity},
there exist $\mathrm{measurable}$ sets $D_{j},D'_{j}$ with $P(D_{j}^{c})\vee P(D_{j}'^{c})<2^{-j+1}$
such that 
\begin{equation}
d(X_{r},X_{s})\vee d(X'_{r},X'_{s})=d(Z_{r},Z_{s})\vee d(Z'_{r},Z'_{s})\leq2^{-j+4}\label{eq:temp-126-4-1}
\end{equation}
on $D_{j}D'_{j}$, for each $r,s\in Q_{\infty}$ with $|r-s|<2^{-m(j)}$.
Consider each $\omega\in D_{j}D'_{j}$ and $t\in[0,1]$. Then there
exists $s\in Q_{m(j)}$ such that $|t-s|<2^{-m(j)}$. Letting $r\rightarrow t$
with $r\in Q_{\infty}$ and $|r-s|<2^{-m(j)}$, inequality \ref{eq:temp-126-4-1}
yields 
\[
d(X_{t}(\omega),X_{s}(\omega))\vee d(X'_{t}(\omega),X'_{s}(\omega))\leq2^{-j+4}.
\]
Consequently,
\[
|d(X_{t}(\omega),X'_{t}(\omega))-d(X_{s}(\omega),X'_{s}(\omega))|<2^{-j+5},
\]
where $\omega\in D_{j}D'_{j}$ and $t\in[0,1]$ are arbitrary. Therefore
\begin{equation}
|\sup_{t\in[0,1]}d(X_{t},X'_{t})-\bigvee_{s\in Q(m(j))}d(X_{s},X'_{s})|\leq2^{-j+5}\label{eq:temp-239}
\end{equation}
on $D_{j}D'_{j}$. Note here that Lemma \ref{Lem. rho_C_hat is a metric}
earlier proved that the supremum is a r.r.v. 

2. Separately, take any $\alpha\in(2^{-j},2^{-j+1})$ and define 
\begin{equation}
A_{j}\equiv(\bigvee_{s\in Q(m(j))}d(X_{s},X'_{s})\leq\alpha).\label{eq:temp-224}
\end{equation}
Then inequality \ref{eq:temp-239} and equality \ref{eq:temp-224}
together yield 
\begin{equation}
G_{j}\equiv D_{j}D'_{j}A_{j}\subset(\sup_{t\in[0,1]}d(X_{t},X'_{t})\leq2^{-j+5}+2^{-j+1}).\label{eq:temp-262}
\end{equation}

3. By inequality \ref{eq:temp-207}, we have 
\begin{equation}
\rho_{prob,Q(\infty)}(Z,Z')<\delta_{reg,auc}(\varepsilon,\overline{m})\equiv2^{-h(j)-2j-1},\label{eq:temp-267}
\end{equation}
where $h_{j}\equiv2^{m(j)}$. Hence 
\[
E\bigvee_{s\in Q(m(j))}\widehat{d}(X_{s},X'_{s})=E\bigvee_{k=0}^{h(j)}\widehat{d}(X_{t(k)},X'_{t(k)})\leq2^{h(j)+1}E\sum_{k=0}^{h(j)}2^{-k-1}\widehat{d}(X_{t(k)},X'_{t(k)})
\]
\begin{equation}
\leq2^{h(j)+1}E\sum_{k=0}^{h(j)}2^{-k-1}\widehat{d}(Z_{t(k)},Z'_{t(k)})\leq2^{h(j)+1}\rho_{prob,Q(\infty)}(Z,Z')<2^{-2j}<\alpha^{2}.\label{eq:temp-258}
\end{equation}
Chebychev's inequality therefore implies that 
\[
P(A_{j}^{c})<\alpha<2^{-j+1}.
\]
Hence
\[
\rho_{\widehat{C}[0,1]}(X,X')\equiv E\sup_{t\in[0,1]}\widehat{d}(X_{t},X'_{t})\leq E1_{G(j)}\sup_{t\in[0,1]}\widehat{d}(X_{t},X'_{t})+P(G_{j}^{c})
\]
\[
\leq(2^{-j+5}+2^{-j+1})+P(A_{j}^{c})+P(D_{j}^{c})+P(D_{j}'^{c})
\]
\[
<(2^{-j+5}+2^{-j+1})+2^{-j+1}+2^{-j+1}+2^{-j+1}<2^{-j+6}<\varepsilon.
\]
Since $\varepsilon>0$ is arbitrary, we see that $\delta_{reg,auc}(\cdot,\overline{m})$
is a modulus of continuity of $\Phi_{Lim}$.
\end{proof}
Theorems \ref{Thm. C-regular =00003D> a.u. continuity} and \ref{Thm. Continuity of a.u. continuous extension}
can now be restated in terms of $C$-regular consistent families of
f.j.d.'s. 
\begin{cor}
\label{Cor. Extension to a.u. continuous processes-1} \textbf{\emph{(Construction
of a.u. continuous processes from $C$-regular families of f.j.d.'s
to )}} Let 
\[
(\Theta_{0},L_{0},I_{0})\equiv([0,1],L_{0},\int\cdot dx)
\]
denote the Lebesgue integration space based on the interval $\Theta_{0}\equiv[0,1]$.
Let $\xi$ be a fixed binary approximation of $(S,d)$ relative to
a reference point $x_{\circ}\in S$. As usual, write $\widehat{d}\equiv1\wedge d$.
Recall from Definition \ref{Def. Metric on  of continuous in prob families of finite joint distributions.}
the metric space $(\widehat{F}_{Cp}([0,1],S),\widehat{\rho}_{Cp,\xi,[0,1]|Q(\infty)})$
of consistent families of f.j.d.'s which are continuous in probability,
with parameter set $[0,1]$ and state space $(S,d)$. Let $\widehat{F}_{0}$
be a subset of $\widehat{F}_{Cp}([0,1],S)$ whose members are $C$-regular
and share a common modulus of C-regularity $\overline{m}\equiv(m_{n})_{n=0,1,2,\cdots}$.
Define the restriction function $\Phi_{[0,1],Q(\infty)}:\widehat{F}_{0}\rightarrow\widehat{F}_{0}|Q_{\infty}$
by $\Phi_{[0,1],Q(\infty)}(F)\equiv F|Q_{\infty}$ for each $F\in\widehat{F}_{0}$.
Then the following holds.

1. The function 
\begin{equation}
\Phi_{fjd,auc,\xi}\equiv\Phi_{Lim}\circ\Phi_{DKS,\xi}\circ\Phi_{[0,1],Q(\infty)}:(\widehat{F}_{0},\widehat{\rho}_{Cp,\xi,[0,1]|Q(\infty)})\rightarrow(\widehat{C}[0,1],\rho_{\widehat{C}[0,1]})\label{eq:temp-272}
\end{equation}
is well defined, where $\Phi_{DKS,\xi}$ is the Daniell-Kolmogorov-Skorokhod
extension constructed in Theorem \ref{Thm. Compact Daniell-Kolmogorov-Skorohod},
and where $\Phi_{Lim}$ is the extension-by-limit constructed in Theorem
\ref{Thm. C-regular =00003D> a.u. continuity}. 

2. For each consistent family $F\in\widehat{F}_{0}$, the a.u. continuous
process $X\equiv\Phi_{fjd,auc,\xi}(F)$ has marginal distributions
given by $F$. 

3. The  construction  $\Phi_{fjd,auc,\xi}$ is uniformly continuous.
\end{cor}
\begin{proof}
1. Let $F\in\widehat{F}_{0}$ be arbitrary. By hypothesis, $F$ is
$C$-regular, with $\overline{m}$ as a modulus of C-regularity. Since
the process $Z\equiv\Phi_{DKS,\xi}(F|Q_{\infty}):Q_{\infty}\times\Theta_{0}\rightarrow S$
extends $F|Q_{\infty}$, so is $Z$. In other words, $Z\in\widehat{R}_{0}$,
where $\widehat{R}_{0}$ is the set of $C$-regular processes on $Q_{\infty}$,
with sample space $(\Theta_{0},L_{0},I_{0})$, and with $\overline{m}$
as a modulus of $C$-regularity. In other words, the set of processes
$\widehat{R}_{0}$ is $C$-equiregular. Hence the a.u. continuous
process $X\equiv\Phi_{Lim}(Z)$ is well defined by Theorem \ref{Thm. C-regular =00003D> a.u. continuity},
with $X|Q_{\infty}=Z$. Thus the composite mapping in equality \ref{eq:temp-272}
is well defined. Assertion 1 is verified.

2. Being $C$-regular, the family $F\in\widehat{F}_{0}$ is continuous
in probability. Hence, for each $r_{1},\cdots,r_{n}\in[0,1]$, and
$f\in C_{ub}(S^{n})$ we have 
\[
F_{r(1),\cdots,r(n)}f=\lim_{s(i)\rightarrow r(i);s(i)\in Q(\infty);i=1,\cdots,n}F_{s(1),\cdots,s(n)}f
\]
\[
=\lim_{s(i)\rightarrow r(i);s(i)\in Q(\infty);i=1,\cdots,n}I_{0}f(Z_{s(1)},\cdots,Z_{s(n)})
\]
\[
=\lim_{s(i)\rightarrow r(i);s(i)\in Q(\infty);i=1,\cdots,n}I_{0}f(X_{s(1)},\cdots,X_{s(n)})
\]
\[
=I_{0}f(X_{r(1)},\cdots,X_{r(n)}),
\]
where the last equality follows from the a.u. continuity of $X$.
We conclude that $F$ is the family of marginal distributions of $X$,
proving Assertion 2.

3. Recall the metric space $(\widehat{R}(Q_{\infty}\times\Theta_{0},S),\rho_{prob,Q(\infty)})$
of processes $Z:Q_{\infty}\times\Theta_{0}\rightarrow S$. Then the
uniform continuity of 
\[
\Phi_{[0,1],Q(\infty)}:(\widehat{F}_{0},\widehat{\rho}_{Cp,\xi,[0,1]|Q(\infty)})\rightarrow(\widehat{F}_{0}|Q_{\infty},\widehat{\rho}_{Marg,\xi,Q(\infty)})
\]
is trivial from Definition \ref{Def. Metric on  of continuous in prob families of finite joint distributions.}.
The Daniel-Kolmogorv-Skorokhod Extension 
\[
\Phi_{DKS,\xi}:(\widehat{F}_{0}|Q_{\infty},\widehat{\rho}_{Marg,\xi,Q(\infty)})\rightarrow(\widehat{R}(Q_{\infty}\times\Theta_{0},S),\rho_{prob,Q(\infty)})
\]
is uniformly continuous by Theorem \ref{Thm. DKS Extension, construction and continuity}.
Moreover, by Step 1, we have 
\[
\Phi_{DKS,\xi}(\widehat{F}_{0}|Q_{\infty})\subset\widehat{R}_{0}\subset\widehat{R}(Q_{\infty}\times\Omega,S),
\]
where the set of processes $\widehat{R}_{0}$ is $C$-equiregular.
Therefore, finally, Theorem \ref{Thm. Continuity of a.u. continuous extension}
says that 
\[
\Phi_{Lim}:(\widehat{R}_{0},\rho_{prob,Q(\infty)})\rightarrow(\widehat{C}[0,1],\rho_{\widehat{C}[0,1]})
\]
is uniformly continuous. Combining, the composite function 
\[
\Phi_{fjd,auc,\xi}\equiv\Phi_{Lim}\circ\Phi_{DKS,\xi}\circ\Phi_{[0,1],Q(\infty)}
\]
is uniformly continuous. Assertion 3 is proved.
\end{proof}

\section{Sufficient Condition for a.u. locally Hoelder Continuity}

Let $(S,d)$ be a locally compact metric space. In Theorems \ref{Thm.  a.u. continuoity on Q0 =00003D> C-regularity}
and \ref{Thm. C-regular =00003D> a.u. continuity}, we saw that the
$C$-regularity of a process\emph{ $Z:Q_{\infty}\times\Omega\rightarrow(S,d)$
}is necessary and sufficient for the a.u continuity of its extension-by-limit
$X\equiv\Phi_{Lim}(Z):[0,1]\times\Omega\rightarrow S$. In this section,
we will prove a sufficient condition, on pairwise joint distributions,
for a.u. continuity of $\Phi_{Lim}(Z)$. 

Refer to the last two sections for notations. For a $\mathrm{measurable}$
set $A$ relative to an arbitrary probability subspace $(\Omega,L,E)$,
we will write $A\in L$ and $1_{A}\in A$ interchangeably, and write
$P(A),E(A),$ and $E1_{A}$ interchangeably. As usual, we will write
the symbols $a_{b}$ and $a(b)$ interchangeably. Recall also the
convention that, for an arbitrary r.r.v. $U$ and for any $a\in R$,
we write $P(U\leq a)$ or $P(U<a)$ only with the explicit or implicit
condition that the real number $a$ has been so chosen that the sets
$(U\leq a)$ or $(U<a)$ are $\mathrm{measurable}$.
\begin{thm}
\label{Thm. Sufficient Condition on pair disributions for a.u. exension}
\textbf{\emph{(A sufficient condition on pair distributions for a.u.
continuous extension).}} Let $\kappa\geq1$ be arbitrary. Let $\overline{\gamma}\equiv(\gamma_{k})_{k=\kappa,\kappa+1,\cdots}$
and $\overline{\varepsilon}\equiv(\varepsilon_{k})_{k=\kappa,\kappa+1,\cdots}$
be two sequences of positive real numbers with $\sum_{k=\kappa}^{\infty}\gamma_{k}<\infty$
and $\sum_{k=1}^{\infty}\varepsilon_{k}<\infty$. 

Let $Z:Q_{\infty}\times(\Omega,L,E)\rightarrow(S,d)$ be an arbitrary
process such that, for each $k\geq\kappa$ and for each $\alpha{}_{k}\geq\gamma_{k}$,
we have 
\begin{equation}
\sum_{t\in[0,1)Q(k)}P(d(Z_{t},Z_{t+\Delta(k+1)})>\alpha_{k})+\sum_{t\in[0,1)Q(k)}P(d(Z_{t+\Delta(k+1)},Z_{t+\Delta(k)})>\alpha_{k})\leq2\varepsilon_{k}.\label{eq:temp-225-4}
\end{equation}
Then the extension-by-limit $X\equiv\Phi_{Lim}(Z):[0,1]\times\Omega\rightarrow(S,d)$
is an a.u. continuous process. 

Specifically, there exists a sequence $(D_{n})_{n=\kappa,\kappa+1,\cdots}$
of \textup{\emph{measurable}} sets such that (\emph{i)} $D_{\kappa}\subset D_{\kappa+1}\subset\cdots,$
\emph{(ii)} for each $n\geq\kappa$, we have 
\begin{equation}
P(D_{n}^{c})\leq2\sum_{k=n}^{\infty}\varepsilon_{k},\label{eq:temp-40-1}
\end{equation}
and \emph{(iii)} for each $n\geq\kappa,$ and for each $\omega\in D_{n}$,
we have 
\begin{equation}
d(X_{r}(\omega),X_{s}(\omega))<8\sum_{k=n}^{\infty}\gamma_{k}\label{eq:temp-127-1-1-1}
\end{equation}
for each $r,s\in[0,1]$ with $|r-s|\leq2^{-n}$. 

Consequently, the process $X$ has a modulus of a.u. continuity defined
as follows. Let $\varepsilon>0$ be arbitrary. Take $n\geq\kappa$
so large that $8\sum_{k=n}^{\infty}\gamma_{k}\vee2\sum_{k=n}^{\infty}\varepsilon_{k}<\varepsilon$.
Define $\delta{}_{auc}(\varepsilon,\overline{\gamma},\overline{\varepsilon})\equiv2^{-n}$.
Then the operation $\delta{}_{auc}(\cdot,\overline{\gamma},\overline{\varepsilon})$
is a modulus of a.u. continuity for the process $X$.
\end{thm}
\begin{proof}
1. Let $Z:Q_{\infty}\times\Omega\rightarrow S$ be as given in the
hypothesis. Let $k\geq\kappa$ be arbitrary, and take any $\alpha_{k}\in[\gamma_{k},2\gamma_{k}).$
Define 
\begin{equation}
C_{k}\equiv\bigcup_{t\in[0,1)Q(k)}(d(Z_{t},Z_{t+\Delta(k+1)})>\alpha_{k})\cup(d(Z_{t+\Delta(k+1)},Z_{t+\Delta(k)})>\alpha_{k}).\label{eq:temp-201-4-5}
\end{equation}
Then $P(C_{k})\leq2\varepsilon_{k}$, thanks to inequality \ref{eq:temp-225-4}
in the hypothesis. Moreover, for each $\omega\in C_{k}^{c}$, and
for each $t\in Q_{k}$ and $s\in Q_{k+1}$ with $|t-s|\leq\Delta_{k+1}$,
we have 
\begin{equation}
d(Z_{t},Z_{s})\leq\alpha_{k}.\label{eq:temp-179}
\end{equation}

2. Let $n\geq\kappa$ be arbitrary, but fixed till further notice.
Define the $\mathrm{measurable}$ set
\[
D_{n}\equiv(\bigcup_{k=n}^{\infty}C_{k})^{c}.
\]
Then
\begin{equation}
PD_{n}^{c}\leq\sum_{k=n}^{\infty}2\varepsilon_{k}.\label{eq:temp-481}
\end{equation}
Moreover, $D_{n}\subset D_{m}$ for each $m\geq n$. 

Now let $r,s\in Q_{\infty}$ be arbitrary with $|s-r|<2^{-n}$. First
assume that $0<s-r<2^{-n}$. Then there exists $r_{n},s_{n}\in Q_{n}$
such that $r\in[r_{n},(r_{n}+2^{-n})\wedge1]$ and $s\in[s_{n},(s_{n}+2^{-n})\wedge1]$.
It follows that 
\[
r_{n}\leq r<s\leq(s_{n}+2^{-n}),
\]
which implies $r_{n}\leq s_{n}$ , and that
\[
s_{n}\leq s<r+2^{-n}\leq(r_{n}+2^{-n})+2^{-n},
\]
which implies $s_{n}\leq(r_{n}+2^{-n})$. Combining, we obtain 
\begin{equation}
s_{n}=r_{n}\quad\mathrm{o}\mathrm{r}\quad s_{n}=(r_{n}+2^{-n}).\label{eq:temp-38}
\end{equation}

Separately. we have $r\in[r_{n},r_{n}+2^{-n}]Q_{\infty}$. Inductively,
consider each $k=n,n+1,\cdots$. Then we have either (i) $r\in[r_{k},r_{k}+2^{-k-1}]Q_{\infty}$,
or (ii) $r\in[r_{k}+2^{-k-1},r_{k}+2^{-k}]Q_{\infty}$. In Case (i)
let $r_{k+1}\equiv r_{k}$. In Case (ii) let $r_{k+1}\equiv r_{k}+2^{-k-1}$.
Then we have $r_{k+1}\in Q_{k+1}$, and, in either case, $r\in[r_{k+1},r_{k+1}+2^{-k-1}]Q_{\infty}$.
Moreover, 
\begin{equation}
|r_{k}-r_{k+1}|\leq2^{-k-1}\equiv\Delta_{k+1}\label{eq:temp-298}
\end{equation}
for each $k\geq n$, and so $r_{k}\rightarrow r$. We can construct
a similar sequence $(s_{k})_{k=n,n+1,\cdots}$ relative to $s$ such
that
\begin{equation}
|s_{k}-s_{k+1}|\leq2^{-k-1}\equiv\Delta_{k+1}\label{eq:temp-298-1}
\end{equation}
for each $k\geq n$. 

4. Now consider each $\omega\in D_{n}$, and consider each $k\geq n$.
Then $\omega\in D_{n}\subset C_{k}^{c}$. Hence equalities \ref{eq:temp-38}
and inequality \ref{eq:temp-179} imply that 
\begin{equation}
d(Z_{r(n)}(\omega),Z_{s(n)}(\omega))\leq2\alpha_{n}.\label{eq:temp-275}
\end{equation}
Similarly, inequalities \ref{eq:temp-298} and \ref{eq:temp-179}
imply that 
\[
d(Z_{r(k)}(\omega),Z_{r(k+1)}(\omega))\leq\alpha_{k}.
\]
Since $r_{k}\rightarrow r$, it follows that
\begin{equation}
d(Z_{r(n)}(\omega),Z_{r}(\omega))\leq\sum_{k=n}^{\infty}d(Z_{r(k)}(\omega),Z_{r(k+1)}(\omega))\leq\sum_{k=n}^{\infty}\alpha_{k}.\label{eq:temp-288}
\end{equation}
Similarly
\begin{equation}
d(Z_{s(n)}(\omega),Z_{s}(\omega))\leq\sum_{k=n}^{\infty}\alpha_{k}.\label{eq:temp-289}
\end{equation}
Combining inequalities \ref{eq:temp-275}, \ref{eq:temp-288}, and
\ref{eq:temp-289}, we obtain
\[
d(Z_{r}(\omega),Z_{s}(\omega))\leq d(Z_{r(n)}(\omega),Z_{s(n)}(\omega))+d(Z_{r}(\omega),Z_{r(n)}(\omega))+d(Z_{s}(\omega),Z_{s(n)}(\omega))
\]
\begin{equation}
\leq2\alpha_{n}+2\sum_{k=n}^{\infty}\alpha_{k}<4\sum_{k=n}^{\infty}\alpha_{k}<8\sum_{k=n}^{\infty}\gamma_{k}\label{eq:temp-297}
\end{equation}
where $r,s\in Q_{\infty}$ are arbitrary with $0<s-r<2^{-n}$. The
same inequality \ref{eq:temp-297} holds, by symmetry, for arbitrary
$r,s\in Q_{\infty}$ with $0<r-s<2^{-n}$. It holds trivially in the
case where $r=s$. Summing up, we have
\begin{equation}
d(Z_{r}(\omega),Z_{s}(\omega))<8\sum_{k=n}^{\infty}\gamma_{k},\label{eq:temp-127}
\end{equation}
for arbitrary $r,s\in Q_{\infty}$ with $|r-s|<2^{-n}$, for arbitrary
$\omega\in D_{n}$, where $P(D_{n}^{c})\leq\sum_{k=n}^{\infty}2\varepsilon_{k}$.
It follows that the process $Z$ is a.u. continuous, with a modulus
of a.u. continuity $\delta{}_{auc}$ defined as follows. Let $\varepsilon>0$
be arbitrary. Let $n\geq\kappa$ be so large that $8\sum_{k=n}^{\infty}\gamma_{k}\vee2\sum_{k=n}^{\infty}\varepsilon_{k}<\varepsilon$.
Define $\delta{}_{auc}(\varepsilon)\equiv\delta{}_{auc}(\varepsilon,\overline{\gamma},\varepsilon_{k})\equiv2^{-n}$.
Proposition \ref{Thm. Completion of process on dense parameter subset}
then says that the extension-by-limit $X\equiv\Phi_{Lim}(Z)$ is an
a.u. continuous process, with the same modulus of a.u. continuity
$\delta{}_{auc}$. By continuity, inequality \ref{eq:temp-127} immediately
extends to the process $X$ for arbitrary $r,s\in[0,1]$ with $|r-s|\leq2^{-n}$,
yielding the desired inequality \ref{eq:temp-127-1-1-1}. 
\end{proof}
As a corollary, we will prove a theorem due to Kolmogorov, which gives
a sufficient condition for the construction of an a.u. locally Hoelder
process, in the following sense.
\begin{defn}
\label{Def. a.u.. lcoally hoelder process} \textbf{(a.u. locally
Hoelder process).} Let $a>0$ be arbitrary. A process $X:[0,a]\times\Omega\rightarrow(S,d)$
is said to be \emph{a.u. locally Hoelder }\index{a.u. locally Hoelder process}
\emph{continuous}, or \emph{a.u. locally Hoelder} for short, if there
exist constants $c,\theta>0$ such that, for each $\varepsilon>0$
there exists some $\delta_{LocHldr}(\varepsilon)>0$ and some $\mathrm{measurable}$
set $D$ with $P(D^{c})<\varepsilon$ such that, for each $\omega\in D$,
we have 
\begin{equation}
d(X_{r}(\omega),X_{s}(\omega))<c|r-s|^{\theta}\label{eq:temp-514-2-2}
\end{equation}
for each $r,s\in[0,a]$ with $|r-s|<\delta_{LocHldr}(\varepsilon)$.
The process is then said to have\emph{ a.u. locally Hoelder exponent}
\index{a..u. locally Hoelder exponent}$\theta$,\emph{ a.u. locally
Hoelder coefficient }\index{a.u. locally Hoelder coefficient}\emph{
$c$,} and \emph{modulus of a..u. locally Hoelder continuity}\index{modulus of a..u. locally Hoelder continuity}.
\end{defn}
\begin{thm}
\label{Thm. Sufficient condition for  a.u. locallyHoelder continuity}
\textbf{\emph{(A sufficient condition on pair distributions for a.u.
locally Hoelder continuity).}} Let $(S,d)$ be a locally compact metric
space. Let $c_{0},u,w>0$ be arbitrary. Let $\theta$ be arbitrary
such that $\theta<u^{-1}w$. 

Suppose $Z:Q_{\infty}\times(\Omega,L,E)\rightarrow(S,d)$ is an arbitrary
process such that
\begin{equation}
P(d(Z_{r},Z_{s})>b)\leq c_{0}b^{-u}|r-s|^{1+w}\label{eq:temp-225-3-1}
\end{equation}
for each $b>0,$ for each $r,s\in Q_{\infty}$. Then the extension-by-limit
$X\equiv\Phi_{Lim}(Z):[0,1]\times\Omega\rightarrow(S,d)$ is a.u.
locally Hoelder with exponent $\theta$. 

Note that inequality \ref{eq:temp-225-3-1} is satisfied if 
\begin{equation}
Ed(Z_{r},Z_{s})^{u}\leq c_{0}|r-s|^{1+w}\label{eq:temp-514}
\end{equation}
\end{thm}
\begin{proof}
Let $Z:Q_{\infty}\times(\Omega,L,E)\rightarrow(S,d)$ be an arbitrary
process such that inequality \ref{eq:temp-225-3-1} holds. 

1. For abbreviation, define the constants $a\equiv(w-\theta u)>0$,
$c_{1}\equiv8(1-2^{-\theta})^{-1},$ and $c\equiv c_{1}2^{\theta}$.
Let $\kappa\geq2$ be so large that $2^{-ka}k^{2}\leq1$ for each
$k\geq\kappa$. Let $k\geq\kappa$ be arbitrary. Define
\begin{equation}
\varepsilon_{k}\equiv c_{0}2^{-w-1}k^{-2}\label{eq:temp-225}
\end{equation}
and
\begin{equation}
\gamma_{k}\equiv2^{-kw/u}k^{2/u}.\label{eq:temp-577}
\end{equation}
Take any $\alpha_{k}\geq\gamma_{k}$ such that the set $(d(Z_{t},Z_{s})>\alpha_{k})$
is $\mathrm{measurable}$ for each $t,s\in Q_{\infty}$. Let $t\in[0,1)Q_{k}$
be arbitrary. We estimate 
\[
P(d(Z_{t},Z_{t+\Delta(k+1)})>\alpha_{k})\leq c_{0}\alpha_{k}^{-u}\Delta_{k+1}^{1+w}
\]
\[
\leq c_{0}\gamma_{k}^{-u}2^{-(k+1)w}2^{-(k+1)}
\]
\[
=c_{0}2^{kw}k^{-2}2^{-(k+1)w}2^{-(k+1)}
\]
\[
=c_{0}2^{-w-1}k^{-2}2^{-k}
\]
where the first inequality is thanks to inequality \ref{eq:temp-225-3-1}.
Similarly
\[
P(d(Z_{t+\Delta(k+1)},Z_{t+\Delta(k)})>\alpha_{k})\leq c_{0}2^{-w-1}k^{-2}2^{-k}.
\]
Combining, we obtain
\[
\sum_{t\in[0,1)Q(k)}P(d(Z_{t},Z_{t+\Delta(k+1)})>\alpha_{k})+\sum_{t\in[0,1)Q(k)}P(d(Z_{t+\Delta(k+1)},Z_{t+\Delta(k)})>\alpha_{k})
\]
\[
\leq2\cdot2^{k}(c_{0}2^{-w-1}k^{-2}2^{-k})=2c_{0}2^{-w-1}k^{-2}\equiv2\varepsilon_{k},
\]
where $k\geq\kappa$ is arbitrary. Since $\sum_{k=\kappa}^{\infty}\gamma_{k}<\infty$
and $\sum_{k=\kappa}^{\infty}\varepsilon_{k}<\infty$, the conditions
in the hypothesis of Theorem \ref{Thm. Sufficient Condition on pair disributions for a.u. exension}
are satisfied by the objects $Z,(\gamma_{k})_{k=\kappa,,\kappa+1\cdots},(\varepsilon_{k})_{k=\kappa,,\kappa+1\cdots}$.
Accordingly, the extension-by-limit $X\equiv\Phi_{Lim}(Z):[0,1]\times\Omega\rightarrow(S,d)$
is a.u. continuous, such that there exists a sequence $(D_{n})_{n=\kappa,,\kappa+1,\cdots}$
of $\mathrm{measurable}$ sets such that (i) $D_{\kappa}\subset D_{\kappa+1}\subset\cdots,$
(ii) for each $n\geq\kappa$, we have 
\begin{equation}
P(D_{n}^{c})\leq2\sum_{k=n}^{\infty}\varepsilon_{k},\label{eq:temp-40}
\end{equation}
and (iii) for each $n\geq\kappa,$ and for each $\omega\in D_{n}$,
we have 
\begin{equation}
d(X_{r}(\omega),X_{s}(\omega))<8\sum_{k=n}^{\infty}\gamma_{k}\label{eq:temp-127-1-1}
\end{equation}
for each $r,s\in[0,1]$ with $|r-s|\leq2^{-n}$. 

2. We will now estimate bounds for the partial sum on the right-hand
side of each of the inequalities \ref{eq:temp-40} and \ref{eq:temp-127-1-1}.
To that end, consider each $n\geq\kappa.$ Then
\begin{equation}
2\sum_{k=n}^{\infty}\varepsilon_{k}\equiv2\sum_{k=n}^{\infty}c_{0}2^{-w-1}k^{-2}\leq2c_{0}2^{-w-1}\int_{y=n-1}^{\infty}y^{-2}dy\leq2c_{0}2^{-w-1}(n-1)^{-1}.\label{eq:temp-199}
\end{equation}
At the same time, 
\[
8\sum_{k=n}^{\infty}\gamma_{k}\equiv8\sum_{k=n}^{\infty}2^{-kw/u}k^{2/u}\equiv8\sum_{k=n}^{\infty}2^{-k\theta u/u}2^{-ka/u}k^{2/u}
\]
\[
=8\sum_{k=n}^{\infty}2^{-k\theta}2^{-ka/u}k^{2/u}\leq8\sum_{k=n}^{\infty}2^{-k\theta}
\]
\begin{equation}
=8\cdot2^{-n\theta}(1-2^{-\theta})^{-1}\equiv c_{1}2^{-n\theta}\label{eq:temp-240}
\end{equation}

3. Let $\varepsilon>0$ be arbitrary. Take $m\geq0$ so large that
$2c_{0}2^{-w-1}(m-1)^{-1}<\varepsilon$. Inequality \ref{eq:temp-199}
then implies that $P(D_{m}^{c})\leq\varepsilon.$ Define $\delta_{LocHldr}(\varepsilon)\equiv2^{-m}.$
Consider each $\omega\in D_{m}\subset D_{m+1}\subset\cdots$. Then,
for each $n\geq m$, we have $\omega\in D_{n}$, whence inequalities
\ref{eq:temp-127-1-1} and \ref{eq:temp-240} together imply that

\begin{equation}
d(X_{r}(\omega),X_{s}(\omega))<c_{1}2^{-n\theta}\label{eq:temp-127-1-1-3}
\end{equation}
for each $(r,s)\in G_{n}\equiv\{(r,s)\in[0,1]^{2}:$ $2^{-n-1}\leq|r-s|\leq2^{-n}\}$.
Hence 

\begin{equation}
d(X_{r}(\omega),X_{s}(\omega))<c_{1}2^{\theta}2^{-(n+1)\theta}\leq c_{1}2^{\theta}|r-s|^{\theta}\equiv c|r-s|^{\theta}\label{eq:temp-127-1-1-3-1}
\end{equation}
for each $(r,s)\in G_{n}$. Therefore 

\begin{equation}
d(X_{r}(\omega),X_{s}(\omega))\leq c|r-s|^{\theta}\label{eq:temp-127-1-1-3-1-1}
\end{equation}
for each $(r,s)\in\bigcup_{n=m}^{\infty}G_{n}$. Since $\bigcup_{n=m}^{\infty}G_{n}$
is dense in $G\equiv\{(r,s)\in[0,1]^{2}:|r-s|\leq2^{-m}\}$, and since
$X(\cdot,\omega)$ is a continuous function, it follows that inequality\ref{eq:temp-127-1-1-3-1-1}
holds for each $(r,s)\in G$. In other words, it holds for each $r,s\in[0,1]$
with $|r-s|\leq2^{-m}\equiv\delta_{LocHldr}(\varepsilon)$. Thus the
process $X$ is a.u. locally Hoelder with exponent $\theta$, as alleged.
\end{proof}
The next corollary implies Theorem 12.4 of \cite{Billingsley68}.
The latter asserts only a.u. continuity and is only for real-valued
processes. 
\begin{cor}
\label{Cor. Sufficient condition for Hoelder continuity-1} \textbf{\emph{(A
sufficient condition on pair distributions for time-scaled a.u. locally
Hoelder continuity).}} Let $(S,d)$ be a locally compact metric space.
Let $c_{0},u,w>0$ be arbitrary. Let $\theta$ be arbitrary such that
$\theta<u^{-1}w$. Let $K:[0,1]\rightarrow[0,1]$ be an arbitrary
continuous and nondecreasing function, with $K(0)=0$ and $K(1)=1,.$

Suppose $Z:Q_{\infty}\times(\Omega,L,E)\rightarrow(S,d)$ is an arbitrary
process such that
\begin{equation}
P(d(Z_{r},Z_{s})>b)\leq c_{0}b^{-u}|K(r)-K(s)|^{1+w}\label{eq:temp-225-3-1-2}
\end{equation}
for each $b>0,$ for each $r,s\in Q_{\infty}$. Then the extension-by-limit
$X\equiv\Phi_{Lim}(Z):[0,1]\times\Omega\rightarrow(S,d)$, subject
to a deterministic time scaling, is a.u. locally Hoelder. More precisely,
there exists a continuous and strictly increasing function, with $G(0)=0$
and $G(1)=1$, such that the process $X_{G}:[0,1]\times\Omega\rightarrow(S,d)$,
defined by $X_{G}(t)\equiv X_{G(t)}$ for each $t\in[0,1]$, is a.u.
locally Hoelder with exponent $\theta$.
\end{cor}
\begin{proof}
1. Fix any $\alpha\geq0$ such that the continuous function $H:[0,1]\rightarrow[0,1]$,
defined by 
\begin{equation}
H(t)\equiv(K(t)+\alpha t)(1+\alpha)^{-1}\label{eq:temp-199-1-1}
\end{equation}
for each $t\in[0,1]$, is strictly increasing. Clearly $H(0)=0$ and
$H(1)=1$. Let $G\equiv H^{-1}$ be the inverse function of $H$,
which is also a continuous increasing function, with $G(0)=0$ and
$G(1)=1$. Write $a\equiv1+\alpha$. Then, equality \ref{eq:temp-199-1-1}
implies that, for each $t\leq s\in[0,1]$, we have
\[
(s-t)\equiv H\circ G(s)-H\circ G(t)\geq a^{-1}\alpha(G(s)-G(t)),
\]
whence
\[
G(s)-G(t)\leq\alpha^{-1}a(s-t).
\]

2. For each $k\geq1$, define 
\begin{equation}
\zeta_{k}\equiv\bigvee_{t\in Q(k)}(K(t+\Delta_{k})-K(t)),\label{eq:temp-275-3-2}
\end{equation}
\begin{equation}
\varepsilon_{k}\equiv2^{-u}c_{0}\zeta_{k}^{w/2},\label{eq:temp-225-3}
\end{equation}
and
\begin{equation}
\gamma_{k}\equiv\zeta_{k}^{w/2u}.\label{eq:temp-577-4}
\end{equation}
Then $(\zeta_{k})_{k=1,2,\cdots}$ is a non increasing sequence in
{[}0,1{]}. Let $\varepsilon>0$ be arbitrary. Let $k\geq1$ be so
large that $\Delta_{k}\equiv2^{-k}<\delta_{K}(\varepsilon)$. Then
$0\leq K(t+\Delta_{k})-K(t)<\varepsilon$ for each $t\in Q_{k}$.
Hence $\zeta_{k}\rightarrow0$. Therefore, for each $\lambda>0$,
we have $\sum_{k=1}^{\infty}\zeta_{k}^{\lambda}<\infty.$ Consequently,
$\sum_{k=1}^{\infty}\gamma_{k}<\infty$ and $\sum_{k=1}^{\infty}\varepsilon_{k}<\infty.$

3. Now let $Z:Q_{\infty}\times(\Omega,L,E)\rightarrow(S,d)$ be an
arbitrary process such that inequality \ref{eq:temp-225-3-1-2} holds.
Let $k\geq1$ be arbitrary, and take any $\alpha_{k}\geq\gamma_{k}.$
We estimate 
\[
P(d(Z_{t},Z_{t+\Delta(k+1)})>\alpha_{k})\leq c_{0}\alpha_{k}^{-u}|K(t+\Delta_{k+1})-K(t)|^{1+w}
\]
\[
=c_{0}\alpha_{k}^{-u}|K(t+\Delta_{k+1})-K(t)|^{w}(K(t+\Delta_{k+1})-K(t))
\]
\[
\leq c_{0}\alpha_{k}^{-u}\zeta_{k}^{w}(K(t+\Delta_{k})-K(t))
\]
\[
<c_{0}\gamma_{k}^{-u}\zeta_{k}^{w}(K(t+\Delta_{k})-K(t))
\]
\[
\equiv c_{0}\zeta_{k}^{-w/2}\zeta_{k}^{w}(K(t+\Delta_{k})-K(t))
\]
\[
=c_{0}\zeta_{k}^{w/2}(K(t+\Delta_{k})-K(t)).
\]
where the first inequality is thanks to inequality \ref{eq:temp-225-3-1-2},
and where the second inequality is by the defining formula \ref{eq:temp-275-3-2}
for $\zeta_{k}$. Similarly
\[
P(d(Z_{t+\Delta(k+1)},Z_{t+\Delta(k)})>\alpha_{k})\leq c_{0}\zeta_{k}^{w/2}(K(t+\Delta_{k})-K(t))
\]
Combining, we obtain
\[
\sum_{t\in[0,1)Q(k)}P(d(Z_{t},Z_{t+\Delta(k+1)})>\alpha_{k})+\sum_{t\in[0,1)Q(k)}P(d(Z_{t+\Delta(k+1)},Z_{t+\Delta(k)})>\alpha_{k})
\]
\[
\leq2c_{0}\zeta_{k}^{w/2}\sum_{t\in Q(k)'}(K(t+\Delta_{k})-K(t))
\]
\[
\leq2c_{0}\zeta_{k}^{w/2}(K(1)-K(0))=2c_{0}\zeta_{k}^{w/2}\equiv2\varepsilon_{k},
\]
where $k\geq0$ is arbitrary. Thus the conditions in the hypothesis
of Theorem \ref{Thm. Sufficient Condition on pair disributions for a.u. exension}
are satisfied by the process $Z$. Accordingly, the extension-by-limit
$X\equiv\Phi_{Lim}(Z):[0,1]\times\Omega\rightarrow(S,d)$ is an a.u.
continuous process. Inequality \ref{eq:temp-225-3-1-2} extends, by
continuity, to
\begin{equation}
P(d(X_{r},X_{s})>b)\leq c_{0}b^{-u}|K(r)-K(s)|^{1+w}\label{eq:temp-225-3-1-1-1}
\end{equation}
for each $b>0,$ for each $r,s\in[0,1]$. 

4. Define the process $Y:Q_{\infty}\times(\Omega,L,E)\rightarrow(S,d)$
by $Y_{t}\equiv X_{G(t)}$ for each $t\in Q_{\infty}$. Let $b>0$
and $r,s\in Q_{\infty}$ be arbitrary. Then
\[
P(d(Y_{r},Y_{s})>b)\equiv P(d(X_{G(r)},X_{G(s)})>b)
\]
\[
\leq c_{0}b^{-u}|K(G(r))-K(G((s))|^{1+w}\leq c_{0}b^{-u}|H(G(r))-H(G((s))|^{1+w}
\]
\begin{equation}
=c_{0}b^{-u}|r-s|^{1+w},\label{eq:temp-40-1-1}
\end{equation}
where the first inequality is thanks to inequality \ref{eq:temp-225-3-1-1-1},
and where the second inequality is due to the defining equality \ref{eq:temp-199-1-1}.
Thus the process $Y$ satisfies the conditions in the hypothesis of
Theorem \ref{Thm. Sufficient condition for  a.u. locallyHoelder continuity}
Accordingly, the extension-by-limit $\overline{Y}_{t}\equiv\Phi_{Lim}(Y):[0,1]\times\Omega\rightarrow(S,d)$
is an a.u. locally Hoelder process with exponent $\theta$, hence
a.u. continuous. Therefore, a.s., we have, for each $t\in[0,1]$,
the equality 
\[
\overline{Y}_{t}\equiv\lim_{r\rightarrow t;r\in Q(\infty)}Y_{r}\equiv\lim_{r\rightarrow t;r\in Q(\infty)}X_{G(r)}
\]
\[
=\lim_{G(r)\rightarrow G(t);G(r)\in G(Q(\infty))}X_{G(r)}=\lim_{s\rightarrow G(t);s\in G(Q(\infty))}X_{s}=X_{G(t)}.
\]
It follows that the process $X_{G}:[0,1]\times\Omega\rightarrow(S,d)$
is a.u. locally Hoelder, as asserted.
\end{proof}

\section{The Brownian Motion}

An application of Theorem \ref{Thm. Sufficient condition for  a.u. locallyHoelder continuity}
is in the construction of the all important Brownian motion. 
\begin{defn}
\label{Def. Brownian Motion in R}\textbf{ (Brownian Motion in $R$).}
An a.u. continuous process $B:[0,\infty)\times(\varOmega,L,E)\rightarrow R$
is called a \emph{Brownian Motion}\index{Brownian Motion} if (i)
$B_{0}=0$, (ii) for each sequence $0\equiv t_{0}\leq t_{1}\leq\cdots\leq t_{n-1}\leq t_{n}$
in $[0,\infty)$, the r.r.v.'s $B_{t(1)}-B_{t(0)},\cdots,B_{t(n)}-B_{t(n-1)}$
are independent, and (iii) for each $s,t\in[0,\infty)$, the r.r.v.
$B_{t}-B_{s}$ is normal with mean $0$ and variance $|t-s|$. Recall
here Definition \ref{Def. continuity in prob etc for arbitrary paameter metric space}
of an a.u. continuous process with a metric parameter space. $\square$
\end{defn}
In the following, let $\overline{Q}_{\infty}$ stand for the set of
dyadic rationals in $[0,\infty)$. 
\begin{thm}
\label{Thm.Construction of  Brownian Motion on =00005B0,inf)} \textbf{\emph{(Construction
of Brownian motion in $R$).}} Brownian Motions in $R$ exist. Specifically,
the following holds.

1. Let $Z:\overline{Q}_{\infty}\times(\Omega,L,E)\rightarrow R$ be
an arbitrary process such that \emph{(i)} $Z_{0}=0$, \emph{(ii)}
for each sequence $0\equiv t_{0}\leq t_{1}\leq\cdots\leq t_{n-1}\leq t_{n}$
in $\overline{Q}_{\infty}$, the r.r.v.'s $Z_{t(1)}-Z_{t(0)},\cdots,Z_{t(n)}-Z_{t(n-1)}$
are independent, and \emph{(iii)} for each $s,t\in\overline{Q}_{\infty}$,
the r.r.v. $Z_{t}-Z_{s}$ is normal with mean $0$ and variance $|t-s|$.
Then the extension-by-limit 
\[
B\equiv\Phi_{Lim}(Z):[0,\infty)\times\Omega\rightarrow R
\]
 is a Brownian motion.

2. For each $n\geq1$ and for each $t_{1},\cdots,t_{n}\in\overline{Q}_{\infty}$,
define the f.j.d.
\[
F_{t(1),\cdots,t(m)}\equiv\Phi_{0,\overline{\sigma}}
\]
 where 
\[
\overline{\sigma}\equiv[\sigma(t_{k},t_{j})]_{k=1,\cdots,n;j=1,\cdots,n}\equiv[t_{k}\wedge t]_{jk=1,\cdots,n;j=1,\cdots,n}.
\]
Then the family $F\equiv\{F_{t(1),\cdots,t(m)}:m\geq1;t_{1},\cdots,t_{m}\in[0,\infty)\}$
of f.j.d.'s is consistent and is continuous in probability.

3. Let $Z:\overline{Q}_{\infty}\times(\Omega,L,E)\rightarrow R$ be
an arbitrary process with marginal distributions given by the family
$F|\overline{Q}_{\infty}$, where $F$ is defined in Assertion 2 above.
Then the extension-by-limit 
\[
B\equiv\Phi_{Lim}(Z):[0,\infty)\times\Omega\rightarrow R
\]
 is a Brownian motion.
\end{thm}
\begin{proof}
For convenience, let $U,U_{1},U_{2},\cdots$ be an independent sequence
of standard normal r.r.v.'s. on some probability space $(\widetilde{\Omega},\widetilde{L},\widetilde{E})$. 

1. Let $Z:\overline{Q}_{\infty}\times(\Omega,L,E)\rightarrow R$ be
an arbitrary process such that Conditions (i-iii) hold. Let $b>0$,
and $s_{1},s_{2}\in\overline{Q}_{\infty}$be arbitrary. Then, by Condition
(iii), the r.r.v. $Z_{s(1)}-Z_{s(2)}$ is normal with mean $0$ and
variance $|s_{1}-s_{2}|$. Consequently, by the formulas in Proposition
\ref{Prop. Moments of standard normal rrv} for moments of standard
normal r.r.v.'s, we obtain $E(Z_{s(1)}-Z_{s(2)})^{4}=3|s_{1}-s_{2}|^{2}.$
Chebychev's inequality then implies that, for each $b>0$, we have
\[
P(|Z_{s(1)}-Z_{s(2)}|>b)=P((Z_{s(1)}-Z_{s(2)})^{4}>b^{4})
\]
\begin{equation}
\leq b^{-4}E(Z_{s(1)}-Z_{s(2)})^{4}=3b^{-4}|s_{1}-s_{2}|^{2},\label{eq:temp-581}
\end{equation}
where $s_{1},s_{2}\in\overline{Q}_{\infty}$ are arbitrary. 

2. Let $N\geq0$ be arbitrary and consider the shifted process $Z^{N}:Q_{\infty}\times(\Omega,L,E)\rightarrow R$
defined by $Z_{s}^{N}\equiv Z_{N+s}$ for each $s\in Q_{\infty}$.
Then, for each $b>0$, and $s_{1},s_{2}\in Q_{\infty}$, we have
\[
P(|Z_{s(1)}^{N}-Z_{s(2)}^{N}|>b)\equiv P(|Z_{N+s(1)}-Z_{N+s(2)}|>b)
\]
\[
\leq3b^{-4}(|(N+s_{1})-(N+s_{2})|^{2}=3b^{-4}|s_{1}-s_{2}|^{2},
\]
where the inequality follows from inequality \ref{eq:temp-581}. Thus
the process $Y$ satisfies the hypothesis of Theorem \ref{Thm. Sufficient condition for  a.u. locallyHoelder continuity},
with $c_{o}=3$, $u=4$, and $w=1$. Accordingly, the extension-by-limit
$W\equiv\Phi_{Lim}(Y):[0,1]\times\Omega\rightarrow R$ is a.u. locally
Hoelder, hence a.u. continuous. In particular, for each $t\in[N,N+1]$,
the limit 
\[
B_{t}\equiv\lim_{r\rightarrow t;r\in\overline{Q}(\infty)[N,N+1]}Z_{r}\equiv\lim_{N+s\rightarrow t;s\in Q(\infty)}Z_{N+s}\equiv\lim_{s\rightarrow t-N;s\in Q(\infty)}Z_{s}^{N}\equiv W_{t-N}.
\]
exist and are equal as r.r.v.'s. In other words, $B|[N,N+1]:[N,N+1]\times\Omega\rightarrow R$
is a well defined process. Moreover, since the process $W\equiv\Phi_{Lim}(Y):[0,1]\times\Omega\rightarrow R$
is a.u. continuous, we see that $B|[N,N+1]$ is a.u. continuous, where
$N\geq0$ is arbitrary. Combining, it follows that the process $B:[0,\infty)\times\Omega\rightarrow R$
is an a.u. continuous process, in the sense of Definition \ref{Def. continuity in prob etc for arbitrary paameter metric space}.
Note that $B_{0}=Z_{0}=0$, in view of Condition (i).

3. Let the sequence $0\equiv t_{0}\leq t_{1}\leq\cdots\leq t_{n-1}\leq t_{n}$
in $[0,\infty)$ and the sequence $0\equiv s_{0}\leq s_{1}\leq\cdots\leq s_{n-1}\leq s_{n}$
in $\overline{Q}_{\infty}$ be arbitrary. Let $f_{i}\in C_{ub}(R)$
be arbitrary for each $i=1,\cdots,n$. Then, using Conditions (ii)
and (iii) in the hypothesis, we obtain
\[
E\prod_{i=1}^{n}f_{i}(B_{s(i)}-B_{s(i-1)})=E\prod_{i=1}^{n}f_{i}(Z_{s(i)}-Z_{s(i-1)})
\]
\begin{equation}
=\prod_{i=1}^{n}Ef_{i}(Z_{s(i)}-Z_{s(i-1)})=\prod_{i=1}^{n}\intop_{R^{m}}\Phi_{0,s(i)-s(i-1)}(du)f_{i}(u).\label{eq:temp-501-1-1}
\end{equation}
Now let $s_{i}\rightarrow t_{i}$ for each $i=1,\cdots,n$. Since
the process $B$ is a.u. continuous, the left-hand side of equality
\ref{eq:temp-501-1-1} converges to $E\prod_{i=1}^{n}f_{i}(B_{t(i)}-B_{t(i-1)})$.
At the same time, since 
\[
\intop_{R^{m}}\Phi_{0,t}(du)f_{i}(u)=\widetilde{E}f(\sqrt{t}U)
\]
is a continuous function of $t$, the right-hand side of equality
\ref{eq:temp-501-1-1} converges to
\[
\prod_{i=1}^{n}\intop_{R^{m}}\Phi_{0,t(i)-t(i-1)}(du)f_{i}(u).
\]
Combining, we see that 
\[
E\prod_{i=1}^{n}f_{i}(B_{t(i)}-B_{t(i-1)})=\prod_{i=1}^{n}\intop_{R^{m}}\Phi_{0,t(i)-t(i-1)}(du)f_{i}(u).
\]
Consequently, the r.r.v.'s $B_{t(1)}-B_{t(0)},\cdots,B_{t(n)}-B_{t(n-1)}$
are independent, with normal distributions with mean 0 and variances
given by $t_{1}-t_{0},\cdots,t_{n}-t_{n-1}$ respectively.

All the conditions in Definition \ref{Def. Brownian Motion in R}
have been verified for the process $B$ to be a Brownian motion. Assertion
1 is proved.

4. To prove Assertion 2, define the function $\sigma:[0,\infty)^{2}\rightarrow[0,\infty)$
by $\sigma(s,t)\equiv s\wedge t$ for each $(s,t)\in[0,\infty)^{2}$.
The function $\sigma$ is clearly symmetric and continuous. We will
verify that it is nonnegative definite in the sense of Definition
\ref{Def. nonnegative definite functions}. To that end, let $n\geq1$
and $t_{1},\cdots,t_{n}\in[0,\infty)$ be arbitrary. We need only
show that the square matrix 
\[
\overline{\sigma}\equiv[\sigma(t_{k},t_{j})]_{k=1,\cdots,n;j=1,\cdots,n}\equiv[t_{k}\wedge t]_{jk=1,\cdots,n;j=1,\cdots,n}
\]
 is nonnegative definite. Let $(\lambda_{k},\cdots,\lambda_{k})\in R^{n}$
be arbitrary. We wish to prove that 
\begin{equation}
\sum_{k=1}^{n}\sum_{j=1}^{n}\lambda_{k}(t_{k}\wedge t_{j})\lambda_{j}\geq0.\label{eq:temp-40-2}
\end{equation}
First assume that $|t_{k}-t_{j}|>0$ if $k\neq j$. Then there exists
a permutation $\pi$ of the indices $1,\cdots,n$ such that $t_{\pi(k)}\leq t_{\pi(j)}$
iff $k\leq j$. It follows that 
\begin{equation}
\sum_{k=1}^{n}\sum_{j=1}^{n}\lambda_{k}(t_{k}\wedge t_{j})\lambda_{j}=\sum_{k=1}^{n}\sum_{j=1}^{n}\lambda_{\pi(k)}(t_{\pi(k)}\wedge t_{\pi(j)})\lambda_{\pi(j)}\equiv\sum_{k=1}^{n}\sum_{j=1}^{n}\theta_{k}(s_{k}\wedge s_{j})\theta_{j},\label{eq:temp-38-1}
\end{equation}
where we write $s_{k}\equiv t_{\pi(k)}$ and $\theta_{k}\equiv\lambda_{\pi(k)}$
for each $k=1,\cdots,n$. Recall the independent standard normal r.r.v.'s.$U_{1},\cdots,U_{n}$
non the probability space $(\widetilde{\Omega},\widetilde{L},\widetilde{E})$.
Thus $\widetilde{E}U_{k}U_{j}=1$ or $0$ according as $k=j$ or not.
Define $V_{k}\equiv\sum_{i=1}^{k}\sqrt{s_{i}-s_{i-1}}U_{i}$ for each
$k=1,\cdots,n$, where $s_{0}\equiv0$. Then $\widetilde{E}V_{k}=0$
and 
\begin{equation}
\widetilde{E}V_{k}V_{j}=\sum_{i=1}^{k\wedge j}(s_{i}-s_{i-1})\widetilde{E}U_{i}^{2}=\sum_{i=1}^{k\wedge j}(s_{i}-s_{i-1})=s_{k\wedge j}=s_{k}\wedge s_{j}\label{eq:temp-580}
\end{equation}
for each $k,j=1,\cdots,n$. Consequently,
\[
\sum_{k=1}^{n}\sum_{j=1}^{n}\theta_{k}(s_{k}\wedge s_{j})\theta_{j}=\widetilde{E}\sum_{k=1}^{n}\sum_{j=1}^{n}\theta_{k}V_{k}V_{j}\theta_{j}=\widetilde{E}(\sum_{k=1}^{n}\theta_{k}V_{k})^{2}\geq0.
\]
Hence the sum on the left-hand side of equality \ref{eq:temp-38-1}
is non-negative. In other words, inequality \ref{eq:temp-40-2} is
valid if the point $(t_{1},\cdots,t_{n})\in[0,\infty)^{n}$ is such
that $|t_{k}-t_{j}|>0$ if $k\neq j$ . Since the set of such points
is dense in $[0,\infty)^{n}$, inequality \ref{eq:temp-40-2} holds,
by continuity, for each $(t_{1},\cdots,t_{n})\in[0,\infty)^{n}.$
In other words, the function $\sigma:[0,\infty)^{2}\rightarrow[0,\infty)$
is nonnegative definite according to Definition \ref{Def. nonnegative definite functions}.

5. For each $m\geq1$ and each sequence $t_{1},\cdots,t_{m}\in[0,\infty)$,
write the nonnegative definite matrix 
\begin{equation}
\overline{\sigma}\equiv[\sigma(t_{k},t_{h})]_{k=1,\cdots,m;h=1,\cdots,m},\label{eq:temp-138-2}
\end{equation}
and define 
\begin{equation}
F_{t(1),\cdots,t(m)}\equiv\Phi_{0,\overline{\sigma}},\label{eq:temp-3-1-1-2}
\end{equation}
where $\Phi_{0,\overline{\sigma}}$ is the normal distribution with
mean $0$ and covariance matrix $\overline{\sigma}$. Take any $M\geq1$
so large that $t_{1},\cdots,t_{m}\in[0,M]$. Proposition \ref{Prop. Consistency of f.j.d.'s generated by covariance function}
says that the family 
\[
F^{(M)}\equiv\{F_{r(1),\cdots,r(m)}:m\geq1;r_{1},\cdots,r_{m}\in[0,M]\}
\]
is consistent and is continuous in probability. Hence, for each $f\in C(R^{n})$,
and for each sequence mapping $i:\{1,\cdots,n\}\rightarrow\{1,\cdots,m\}$,
we have 
\begin{equation}
F_{t(1),\cdots,t(m)}(f\circ i^{*})=F_{t(i(1)),\cdots,t(i(n))}f,\label{eq:temp-196-2-1}
\end{equation}
where the dual function $i^{*}:R^{m}\rightarrow R^{n}$ is defined
by 
\begin{equation}
i^{*}(x_{1},\cdots,x_{m})\equiv(x_{i(1)},\cdots,x_{i(n)})\label{eq:temp-171-1}
\end{equation}
for each $(x_{1},\cdots,x_{m})\in R^{m}$. Thus the family
\begin{equation}
F\equiv\{F_{t(1),\cdots,t(m)}:m\geq1;t_{1},\cdots,t_{m}\in[0,\infty)\}\label{eq:temp-132-2}
\end{equation}
of f.j.d.'s is consistent and is continuous in probability. Assertion
2 is proved.

6. To prove the remaining Assertion 3, let $Z:\overline{Q}_{\infty}\times(\Omega,L,E)\rightarrow R$
be an arbitrary process with marginal distributions given by the family
$F|\overline{Q}_{\infty}$, where $F$ is defined in Assertion 2 above.
Such a process $Z$ exists by the Daniel-Kolmogorv- or the Daniel-Kolmogorv-Skorokhod
Extension Theorem.

7. Let $t_{1},t_{2}\in\overline{Q}_{\infty}$ be arbitrary. Then,
according to Steps 5 and 6 above, the r.r.v.'s $Z_{t(1)},Z_{t(2)}$
have a jointly normal distribution given by $F_{t(1),t(2)}\equiv\Phi_{0,\overline{\sigma}}$
where $\overline{\sigma}\equiv[t_{r}\wedge t_{h}]_{k=1,2;;h=1,2}.$
Hence 
\[
EZ_{t(1)}Z_{t(2)}=t_{1}\wedge t_{2},
\]
It follows that $Z_{t(1)}-Z_{t(2)}$ is a normal r.r.v. with mean
$0$, and with variance given by
\[
E(Z_{t(1)}-Z_{t(2)})^{2}=EZ_{t(1)}^{2}+ZB_{t(2)}^{2}-2EZ_{t(1)}Z_{t(2)}=t_{1}+t_{2}-2t_{1}\wedge t_{2}=|t_{1}-t_{2}|.
\]

8. Now let $0\equiv t_{0}\leq t_{1}\leq\cdots\leq t_{n-1}\leq t_{n}$
be arbitrary in $\overline{Q}_{\infty}$. Then the r.r.v.'s $Z_{t(1)},\cdots,Z_{t(n)}$
have joint distribution $F_{t(1),\cdots,t(n)}$ according to Steps
5 and 6. Hence $Z_{t(1)},\cdots,Z_{t(n)}$ are jointly normal. Therefore
the r.r.v.'s $Z_{t(1)}-Z_{t(0)},\cdots,Z_{t(n)}-Z_{t(n-1)}$ are jointly
normal. Moreover, for each $i,k=1,\cdots,n$ with $i<k$, we have
\[
E(Z_{t(i)}-Z_{t(i-1)})(Z_{t(k)}-Z_{t(k-1)})
\]
\[
=EZ_{t(i)}Z_{t(k)}-EZ_{t(i)}Z_{t(k-1)}-EZ_{t(i-1)}Z_{t(k)}+EZ_{t(i-1)}Z_{t(k-1)}
\]
\begin{equation}
=t_{t}-t_{i}-t_{i-1}+t_{i-1}=0.\label{eq:temp-290}
\end{equation}
Thus the jointly normal r.r.v.'s $Z_{t(1)}-Z_{t(0)},\cdots,Z_{t(n)}-Z_{t(n-1)}$
are pairwise uncorrelated. Hence, by Assertion 3 of Proposition \ref{Prop. Basics of Normal distributions},
they are mutually independent. Summing up Steps 7 and 8, all of Conditions
(i-iii) of Assertion 1 have been verified for the process $Z$. Accordingly,
the extension-by-limit 
\[
B\equiv\Phi_{Lim}(Z):[0,\infty)\times\Omega\rightarrow R
\]
 is a Brownian motion. Assertion 3 and the Theorem is proved.
\end{proof}
The following corollary is Levy's well know result on the a.u. Hoelder
continuity of a Brownian motion. A stronger theorem by Levy gives
the best modulus of a.u. continuity of a Brownian motion, and shows
that the $\theta\in(0,\frac{1}{2})$ is the best Hoelder exponent
that can be hoped for; a.u. local Hoelder continuity for Brownian
motion with exponent $\frac{1}{2}$ fails. 
\begin{cor}
\textbf{\emph{\label{Cor. (Brownian-Motion-on finite interval is a.u. Hoelder with exponent <1/2}
(Brownian Motion on a finite interval is a.u. locally Hoelder with
exponent less than $\frac{1}{2}$).}} Let $B:[0,\infty)\times(\varOmega,L,E)\rightarrow R$
be a Brownian motion. Let $\theta\in(0,\frac{1}{2})$ and $a>0$ be
arbitrary. Then $B|[0,a]$ is a.u. locally Hoelder with exponent $\theta$.
\end{cor}
\begin{proof}
Since $\theta<\frac{1}{2}$, there exists $m\geq0$ be so large that
$\theta<(2+2m)^{-1}m$. Consider the process $X:[0,1]\times(\varOmega,L,E)\rightarrow R$
defined by $X_{t}\equiv B_{at}$ for each $t\in[0,1]$. Consider each
$b>0$ and each $r,s\in[0,1]$. Then the r.r.v. $X_{s}-X_{r}\equiv B_{ar}-B_{at}$
is normally distributed with mean $0$ and variance $a|r-s|$. Therefore
\[
E|X_{s}-X_{r}|^{2+2m}=\overline{c}_{0}a^{1+m}|r-s|^{1+m},
\]
where $\overline{c}_{0}\equiv EU^{2+2m}$ is the $(2+m)$-th moment
of a standard normal r.r.v. $U$. Thus the process $X|Q_{\infty}$satisfies
inequality \ref{eq:temp-514} of Theorem \ref{Thm. Sufficient condition for  a.u. locallyHoelder continuity}
with $u\equiv2+2m$, $c_{0}\equiv\overline{c}_{0}a^{1+m}$, and $w\equiv m$.
Note that $\theta<u^{-1}w$ by the choice of $m$. Hence, accordingly
to Theorem \ref{Thm. Sufficient condition for  a.u. locallyHoelder continuity},
the process $X$ is a.u. locally Hoelder with exponent $\theta$,
with some a.u. locally Hoelder coefficient\emph{ }$c$, and with some
modulus of a..u. locally Hoelder continuity $\delta_{LocHldr}.$

Let $\varepsilon>0$ be arbitrary. Then, according to Definition \ref{Def. a.u.. lcoally hoelder process},
there exists $\mathrm{measurable}$ set $D$ with $P(D^{c})<\varepsilon$
such that, for each $\omega\in D$, we have 
\begin{equation}
|X_{r}(\omega)-X_{s}(\omega)|<c|r-s|^{\theta}\label{eq:temp-514-2-2-1-1}
\end{equation}
for each $r,s\in[0,1]$ with $|r-s|<\delta_{LocHldr}(\varepsilon)$.Now
consider each $\omega\in D$ and each $t,u\in[0,a]$ with $|t-u|<a\delta_{LocHldr}(\varepsilon)$.
Then inequality \ref{eq:temp-514-2-2-1-1} yields
\begin{equation}
|B_{t}(\omega)-B_{u}(\omega)|\equiv|X_{t/a}(\omega)-X_{u/a}(\omega)|<c|a^{-1}t-a^{-s}s|^{\theta}=ca^{-1}|t-s|^{\theta}.\label{eq:temp-514-2-2-1-1-1}
\end{equation}
Thus we see that the process $B|[0,a]$ is a.u.. locally Hoelder with
exponent $\theta$, according to Definition \ref{Def. a.u.. lcoally hoelder process},
as alleged.
\end{proof}

\section{The Garsia-Rodemich-Rumsey Theorem}

In this section, we will restrict our attention to real-valued Gaussian
processes with parameter set $[0,1]$. We will let $\sigma:[0,1]\times[0,1]\rightarrow R$
be an arbitrary continuous symmetric positive definite function. 

\cite{Garsia 70} gives a condition on the modulus of continuity of
$\sigma$ under which there exists an a.u. continuous Gaussian process
$X$ with $\sigma$ as covariance functions. The Garsia-Rodemich-Rumsey
proof is by showing that the partial sums of the Karhunen-Loeve expansion
relative to $\sigma$ are, under said condition, a.u. convergent to
an a.u. continuous process. We will quote the key real-variable lemma
in \cite{Garsia 70}. We will then present a proof of the main theorem
which is, in essence, the proof in the cited paper except that we
dispense with an unnecessary appeal to a version of the submartingale
convergence theorem which states that every submartingale with bounded
expectations converges a.u. to some r.r.v. This version implies the
principle of infinite search. The constructive version, Theorem \ref{Thm. Specific strictly convex function and a.u. convergence of submartingles},
requires additional information on the convergence of some sequences
of expectations, in order to yield the a.u. convergence and the measurability
of the limiting r.r.v.

Instead of supplying a proof of the convergence of the submartingale
derived from the Karhunen-Loeve expansion, we bypass both the Karhunen-Loeve
expansion and the use of submartingales. We will derive Borel-Cantelli
styled estimates on conditional expectations, thus sticking to elementary
time-domain analysis and obviating the need, for the present purpose,
of more ground work of spectral analysis of the covariance function.
We note that the use of conditional expectations in relation to the
Karhunen-Loeve expansion is mentioned in \cite{Garsia 70} for a related
result. 

First some notations. Define 
\begin{equation}
\Delta\sigma(s,t)\equiv\sigma(s,s)+\sigma(t,t)-2\sigma(s,t)\label{eq:temp-500}
\end{equation}
for each $t,s\in[0,1]$. It follows from the continuity of $\sigma$
that $\Delta\sigma(s,t)\rightarrow0$ as $|s-t|\rightarrow0$.

In the following, recall that $\int_{0}^{1}\cdot dp$ denotes the
Riemann-Stieljes integration relative to an arbitrary distribution
function $p$ on $[0,1]$. 
\begin{defn}
\textbf{\label{Def. Two auxilliary functions} (Two auxiliary functions).}
Introduce the auxiliary function
\begin{equation}
\Psi(v)\equiv\exp(\frac{1}{4}v^{2})\label{eq:temp-526}
\end{equation}
for each $v\in[0,\infty)$, with its inverse
\begin{equation}
\Psi^{-1}(u)\equiv2\sqrt{\log u}\label{eq:temp-7-2-1}
\end{equation}
for each $u\in[1,\infty)$.
\end{defn}
$\square$

Next we cite, without the proof from \cite{Garsia 70}, a remarkable
real variable lemma . It derives the global modulus of continuity,
inequality \ref{eq:temp-531} below, of a function from a condition
on its local properties, inequality \ref{eq:temp-7} below. It is
key to the main theorem. 
\begin{lem}
\label{Lem:GRR lemma} \textbf{\emph{(Garsia-Rodemich-Rumsey Real
Variable Lemma).}} Let the function $\Psi$ and its inverse $\Psi^{-1}$
be as in Definition \ref{Def. Two auxilliary functions}. Let $\overline{p}:[0,1]\rightarrow[0,\infty)$
be an arbitrary continuous nondecreasing function with $\overline{p}(0)=0$.
Let $f$ be an arbitrary continuous function on $[0,1],$ and let
$B>0$ be such that the function 
\[
\Psi(\frac{|f(t)-f(s)|}{\overline{p}(|t-s|)})
\]
of $(t,s)\in[0,1]^{2}$ is integrable, with
\begin{equation}
\int_{0}^{1}\int_{0}^{1}\Psi(\frac{|f(t)-f(s)|}{\overline{p}(|t-s|)})dtds\leq B,\label{eq:temp-7}
\end{equation}
 Then
\begin{equation}
|f(t)-f(s)|\leq8\int_{0}^{|t-s|}\Psi^{-1}(\frac{4B}{u^{2}})d\overline{p}(u)\label{eq:temp-531}
\end{equation}
for each $(t,s)\in[0,1]^{2}:$
\end{lem}
\begin{proof}
See \cite{Garsia 70}. 
\end{proof}
Recall from Definition \ref{Def. Notations for dyadic rationals}
some more notations for dyadic rationals in $[0,1]$. For each $N\geq0$,
we have $p_{N}\equiv2^{N}$,$\overline{p}_{N}\equiv2^{2N},$ $\Delta_{N}\equiv2^{-N}$,
and the enumerated finite sets of dyadic rationals
\[
Q_{N}\equiv\{t_{0},t_{1},\cdots,t_{p(N)}\}=\{q_{N.0},\cdots,q_{N,p(N)}\}\equiv\{0,\Delta_{N},2\Delta_{N},\cdots,1\},
\]
where the second equality is equality of sets without the enumeration,
and
\[
Q_{\infty}\equiv\bigcup_{N=0}^{\infty}Q_{N}\equiv\{t_{0},t_{1},\cdots\}.
\]

Recall that $[\cdot]_{1}$ is the operation which assigns to each
$a\in R$ an integer $[a]_{1}\in(a,a+2)$. Recall also the matrix
notations in Definition \ref{Def. Matrix notations}, and the basic
properties of conditional distributions established in Propositions
\ref{Prop. Basics of Conditional expectations} and \ref{Prop. Conditonal expectation of jointly normal r.r.v.'s}.
As usual, to lessen the burden on subscripts, we write the symbols
$x_{y}$ and $x(y)$ interchangeably for any expressions $x$ and
$y$. 
\begin{lem}
\label{Lem. Interpolation of Gaussian procee} \textbf{\emph{(Interpolation
of Gaussian process with conditional expectations).}} Let $Y:Q_{\infty}\times\Omega\rightarrow R$
be an arbitrary centered Gaussian process with a continuous positive
definite covariance function $\sigma$. Thus $E(Y_{t}-Y_{s})^{2}=\Delta\sigma(t,s)$
for each $t,s\in Q_{\infty}$. Let $\overline{p}:[0,1]\rightarrow[0,\infty)$
be an arbitrary continuous nondecreasing function such that 
\[
\bigvee_{0\leq s,t\leq1;|s-t|\leq u}(\Delta\sigma(s,t))^{1/2}\leq\overline{p}(u)
\]
for each $u\in[0,1]$. Then there exists $Y^{(n)}:[0,1]\times\Omega\rightarrow R$
such that the following holds.

1. Let $n\geq0$ and $t\in[0,1]$ be arbitrary. Define the r.r.v.
\[
Y_{t}^{(n)}\equiv Y_{t}(n)\equiv E(Y_{t}|Y_{t(0)},\cdots,Y_{t(n)}).
\]
Then, for each fixed $n\geq0$, the process $Y^{(n)}:[0,1]\times\Omega\rightarrow R$
is an a.u. continuous centered Gaussian process. Moreover, $Y_{r}^{(n)}=Y_{r}$
for each $r\in\{t_{0},\cdots,t_{n}\}$. We will call the process $Y^{(n)}$
the \emph{interpolated approximation of $Y$ by conditional expectations}\index{interpolated Gaussian process by conditional expectations}
on $\{t_{0},\cdots,t_{n}\}$.

2. For each fixed $t\in[0,1]$, the process $Y_{t}:\{0,1,\cdots\}\times\Omega\rightarrow R$
is a martingale relative to the filtration $\mathcal{L}\equiv\{L(Y_{t(0)},\cdots,Y_{t(n)}):n=0,1,\cdots\}$. 

3. If $m>n\geq1$, define
\[
Z_{t}^{(m,n)}\equiv Y_{t}^{(m)}-Y_{t}^{(n)}\in L(Y_{t(0)},\cdots,Y_{t(m)})
\]
for each $t\in[0,1]$. Let $\overline{\Delta}>0$ be arbitrary. Suppose
$n$ is so large that the subset $\{t_{0},\cdots,t_{n}\}$ is a $\overline{\Delta}$-approximation
of $[0,1]$. Define the continuous nondecreasing function $\overline{p}_{\overline{\Delta}}:[0,1]\rightarrow[0,\infty)$
by 
\[
\overline{p}_{\overline{\Delta}}(u)\equiv\overline{p}(u)\wedge2\overline{p}(\overline{\Delta})
\]
for each $u\in[0,1]$. Then 
\[
E(Z_{t}^{(m,n)}-Z_{s}^{(m,n)})^{2}\leq\overline{p}_{\overline{\Delta}}^{2}(|t-s|).
\]
\end{lem}
\begin{proof}
First note that since $Y:Q_{\infty}\times\Omega\rightarrow R$ is
centered Gaussian with covariance function $\sigma$, we have 
\[
E(Y_{t}-Y_{s})^{2}=EY_{t}^{2}-2EY_{t}Y_{s}+EY_{s}^{2}=\sigma(t,s)-2\sigma(t,s)+\sigma(t,s)\equiv\Delta\sigma(t,s)
\]
for each $t,s\in Q_{\infty}$

1. Let $n\geq0$ be arbitrary. Then the r.v. $U_{n}\equiv(Y_{t(0)},\cdots,Y_{t(n)})$
with values in $R^{n+1}$ is normal, with mean $0$ and has the positive
definite covariance matrix
\[
\overline{\sigma}_{n}\equiv EU_{n}U_{n}^{T}=[\sigma(t_{h},t_{j})]_{h=0,\cdots;n;j=0,\cdots,n}.
\]
For each $t\in[0,1]$, define
\[
c_{n,t}\equiv(\sigma(t,t_{0}),\cdots,\sigma(t,t_{n}))\in R^{n+1}
\]
and define the Gaussian process $\overline{Y}^{(n)}:[0,1]\times\Omega\rightarrow R$
by
\[
\overline{Y}_{t}^{(n)}\equiv c_{n,t}^{T}\overline{\sigma}_{n}^{-1}U_{n}.
\]
Then, since $c_{n,t}$ is continuous in $t$, the process $\overline{Y}^{(n)}$
is a.u. continuous.

Moreover, for each $t\in Q_{\infty}$, the conditional expectation
of $Y_{t}$ given $U_{n}$ is, according to Proposition \ref{Prop. Conditonal expectation of jointly normal r.r.v.'s},
given by 
\begin{equation}
E(Y_{t}|Y_{t(0)},\cdots,Y_{t(n)})=E(Y_{t}|U_{n})=c_{n,t}^{T}\overline{\sigma}_{n}^{-1}U_{n}\equiv\overline{Y}_{t}^{(n)}.\label{eq:temp-525}
\end{equation}
Thus $Y_{t}^{(n)}\equiv Y_{t}(n)=\overline{Y}_{t}^{(n)}$. Hence,
since $\overline{Y}^{(n)}$ is an a.u. continuous and centered Gaussian
process, so is $Y_{t}^{(n)}$. Assertion 1 is proved. Note that, for
each $r\in\{t_{0},\cdots,t_{n}\}$ and for each $m\geq n$, we have
$r\in\{t_{0},\cdots,t_{m}\}$, whence $Y_{r}\in L(U_{m})$
\begin{equation}
Y_{r}^{(m)}=E(Y_{r}|U_{m})=Y_{r},\label{eq:temp-17-1}
\end{equation}
where the second equality is a trivial consequence of the conditional
expectation.

2. Let $m>n\geq1$ be arbitrary. Then, for each $t\in Q_{\infty}$,
we have 
\[
E(Y_{t}^{(m)}|U_{n})=E(E(Y_{t}|U_{m})|U_{n})=E(Y_{t}|U_{n})=Y_{t}^{(n)},
\]
where the first and third equality are by equality \ref{eq:temp-525},
and where the second equality is because $L(U_{n})\subset L(U_{m})$.
Hence, for each $V\in L(U_{n})$, we have
\[
EY_{t}^{(m)}V=EY_{t}^{(n)}V
\]
for each $t\in Q_{\infty}$, and, by continuity, also for each $t\in[0,1]$.
Thus $E(Y_{t}^{(m)}|U_{n})=Y_{t}^{(n)}$ for each $t\in[0,1]$. We
conclude that, for each fixed $t\in[0,1]$, the process $Y_{t}:\{0,1,\cdots\}\times\Omega\rightarrow R$
a martingale relative to the filtration $\{L(U_{n}):n=0,1,\cdots\}$.
Assertion 2 is proved.

3. Let $m>n\geq1$ and $t,s\in[0,1]$ be arbitrary. Then
\[
Z_{t}^{(m,n)}-Z_{s}^{(m,n)}\equiv Y_{t}^{(m)}-Y_{t}^{(n)}-Y_{s}^{(m)}+Y_{s}^{(n)}=Y_{t}^{(m)}-Y_{s}^{(m)}-E(Y_{t}^{(m)}-Y_{s}^{(m)}|U_{n}).
\]
Hence, by Proposition \ref{Prop. Basics of Conditional expectations}
\[
E(Z_{t}^{(m,n)}-Z_{s}^{(m,n)})^{2}\leq E(Y_{t}^{(m)}-Y_{s}^{(m)})^{2}
\]
Suppose $t,s\in Q_{\infty}.$ Then equality \ref{eq:temp-525} implies
that
\begin{equation}
Y_{t}^{(m)}-Y_{s}^{(m)}=E(Y_{t}-Y_{s}|U_{m}).\label{eq:temp-445}
\end{equation}
Hence, 
\[
E(Z_{t}^{(m,n)}-Z_{s}^{(m,n)})^{2}\leq E(Y_{t}^{(m)}-Y_{s}^{(m)})^{2}\leq E(Y_{t}-Y_{s})^{2}=\Delta\sigma(t,s),
\]
where $t,s\in Q_{\infty}.$ are arbitrary, and where the second inequality
is by equality \ref{eq:temp-445} and by Proposition \ref{Prop. Basics of Conditional expectations}.
By continuity, we therefore have
\begin{equation}
E(Z_{t}^{(n,m)}-Z_{s}^{(n,m)})^{2}\leq E(Y_{t}^{(m)}-Y_{s}^{(m)})^{2}\leq\Delta\sigma(t,s)\leq\overline{p}^{2}(|t-s|)\label{eq:temp-144-3}
\end{equation}
where $t,s\in[0,1]$ are arbitrary. 

Now let $\overline{\Delta}>0$ be arbitrary. Suppose $n\geq1$ is
so large that the subset $\{t_{0},\cdots,t_{n}\}$ is a $\overline{\Delta}$-approximation
of $[0,1]$. Let $t',s'\in\{t_{0},\cdots,t_{n}\}$ be such that $|t-t'|\vee|s-s'|<\overline{\Delta}$.
Then equality \ref{eq:temp-17-1} implies that 
\begin{equation}
Z_{t'}^{(m,n)}\equiv Y_{t'}^{(m)}-Y_{t'}^{(n)}=Y_{t'}-Y_{t'}=0,\label{eq:temp-17-3}
\end{equation}
with a similar inequality for $s'$. Applying inequality \ref{eq:temp-144-3}
to $t,t'$ in place of $t,s$, we obtain
\[
E(Z_{t}^{(m,n)}-Z_{t'}^{(m,n)})^{2}\leq\overline{p}^{2}(|t-t'|)\leq\overline{p}^{2}(\overline{\Delta}),
\]
and a similar inequality for the pair $s,s'$ in place of $t,t'$.
In addition, equality \ref{eq:temp-17-3} implies 
\[
Z_{t}^{(m,n)}-Z_{s}^{(m,n)}=(Z_{t}^{(m,n)}-Z_{t'}^{(m,n)})-(Z_{s}^{(m,n)}-Z_{s'}^{(m,n)}).
\]
Hence Minkowski's inequality yields
\begin{equation}
\sqrt{E(Z_{t}^{(m,n)}-Z_{s}^{(m,n)})^{2}}\leq\sqrt{E(Z_{t}^{(m,n)}-Z_{t'}^{(m,n)})^{2}}+\sqrt{E(Z_{s}^{(m,n)}-Z_{s'}^{(m,n)})^{2}}\leq2\overline{p}(\overline{\Delta}).\label{eq:temp-189-3}
\end{equation}
Combining inequalities \ref{eq:temp-144-3} and \ref{eq:temp-189-3},
we obtain
\[
E(Z_{t}^{(m,n)}-Z_{s}^{(m,n)})^{2}\leq(\overline{p}(|t-s|)\wedge2\overline{p}(\overline{\Delta}))^{2}\equiv\overline{p}_{\overline{\Delta}}^{2}(|t-s|).
\]
Assertion 3 is proved.
\end{proof}
The next lemma prepares for the proof of the main theorem. It contains
a redundant assumption of a.u. continuity, which will be stripped
off in the main theorem.
\begin{lem}
\label{Lem. GRR modulus of continuity assuming a.u.contintuiy} \textbf{\emph{(Modulus
of a.u. continuity with the redundant assumption of a.u.continuity).}}
Let $V:[0,1]\times\Omega\rightarrow R$ be an arbitrary a.u. continuous
and centered Gaussian process, with a continuous positive definite
covariance function $\sigma$. Thus 
\[
\xi_{t,s}\equiv E(V_{t}-V_{s})^{2}=\Delta\sigma(t,s)
\]
 for each $t,s\in[0,1]$, where the operator $\Delta$ is defined
in equality \ref{eq:temp-500}. Let $\overline{p}:[0,1]\rightarrow[0,\infty)$
be a continuous increasing function, with $\overline{p}(0)=0$, such
that $\sqrt{-\log u}$ is integrable relative to the distribution
function $\overline{p}$ on $[0,1]$. Thus 
\begin{equation}
\int_{0}^{1}\sqrt{-\log u}d\overline{p}(u)<\infty.\label{eq:temp-144-1-1}
\end{equation}

Suppose
\begin{equation}
\bigvee_{0\leq s,t\leq1;|s-t|\leq u}\xi_{t,s}\leq\overline{p}(u)^{2}\label{eq:temp-499}
\end{equation}
for each $u\in[0,1]$. Then there exists an integrable r.r.v. $B$
with $EB\leq\sqrt{2}$ such that 
\begin{equation}
|V(t,\omega)-V(s,\omega)|\leq16\int_{0}^{|t-s|}\sqrt{\log(\frac{4B(\omega)}{u^{2}})}d\overline{p}(u)\label{eq:temp-3-1-2-1}
\end{equation}
 for each $t,s\in[0,1]$, for each $\omega\in domain(B)$.
\end{lem}
\begin{proof}
With positive definiteness of the function $\sigma$, the defining
equality \ref{eq:temp-500} implies that $\Delta\sigma(s,t)>0$ for
each $s,t\in[0,1]$ with $|s-t|>0$. Hence, in view of inequality
\ref{eq:temp-58}, we have $\overline{p}(u)>0$ for each $u\in(0,1]$. 

Define the full subset 
\[
D\equiv\{(t,s)\in[0,1]^{2}:|t-s|>0\}
\]
of $[0,1]^{2}$. Because the process $V$ is a.u. continuous, there
exists a full set $A\subset\Omega$ such that $V(\cdot,\omega)$ is
continuous on $[0,1]$. Moreover, $V$ is a $\mathrm{measurable}$
function on $[0,1]\times\Omega$. Define the function $U:[0,1]^{2}\times\Omega\rightarrow R$
by 
\[
domain(U)\equiv D\times A
\]
and
\[
U(t,s,\omega)\equiv\Psi(\frac{|V(t,\omega)-V(s,\omega)|}{\overline{p}(|t-s|)})
\]
\begin{equation}
\equiv\exp(\frac{1}{4}\frac{(V(t,\omega)-V(s,\omega))^{2}}{\overline{p}(|t-s|)^{2}})\label{eq:temp-532}
\end{equation}
for each $(t,s,\omega)\in domain(U)$. Then $U$ is a $\mathrm{measurable}$
on $[0,1]^{2}\times\Omega$. 

Let $t,s\in[0,1]$. be arbitrary. Then $\xi_{t,s}\leq\overline{p}(|t-s|)^{2}$.
Hence 
\begin{equation}
\frac{1}{\sqrt{2\pi\xi_{t,s}}}\exp(\frac{1}{4}\frac{u^{2}}{\overline{p}(|t-s|)^{2}})\exp(-\frac{1}{2}\frac{u^{2}}{\xi_{t,s}})\leq\frac{1}{\sqrt{2\pi\xi_{t,s}}}\exp(-\frac{1}{4}\frac{u^{2}}{\xi_{t,s}}).\label{eq:temp-518}
\end{equation}
The $\mathrm{measurable}$ function of $(t,s,u)$ on the right-hand
side is integrable on $[0,1]^{2}\times R$ relative to the product
Lebesgue integration, with 
\[
\int_{0}^{1}\int_{0}^{1}\int_{-\infty}^{\infty}\frac{1}{\sqrt{2\pi\xi_{t,s}}}\exp(-\frac{1}{4}\frac{u^{2}}{\xi_{t,s}})dudtds=\int_{0}^{1}\int_{0}^{1}\sqrt{2}dtds=\sqrt{2}.
\]
Hence the $\mathrm{measurable}$ function of $(t,s,u)$ on the left-hand
side of \ref{eq:temp-518} is integrable on $[0,1]^{2}\times R$,
with integral bounded by $\sqrt{2}$. 

Now let $b>0$ be arbitrary. Note that the r.r.v. $V(t,\cdot)-V(s,\cdot)$
is Gaussian, with mean $0$ and variance $\xi_{t,s}$. Hence
\[
E\int_{0}^{1}\int_{0}^{1}b\wedge Udtds
\]
\[
=\int_{0}^{1}\int_{0}^{1}\int_{-\infty}^{\infty}b\wedge\exp(\frac{1}{4}\frac{u^{2}}{\overline{p}(|t-s|)^{2}})\varphi_{0,\xi(t,s)}(u)dudtds
\]
\[
=\int_{0}^{1}\int_{0}^{1}\int_{-\infty}^{\infty}b\wedge\frac{1}{\sqrt{2\pi\xi_{t,s}}}\exp(\frac{1}{4}\frac{u^{2}}{\overline{p}(|t-s|)^{2}})\exp(-\frac{1}{2}\frac{u^{2}}{\xi_{t,s}})dudtds
\]
\begin{equation}
\rightarrow\int_{0}^{1}\int_{0}^{1}\int_{-\infty}^{\infty}\frac{1}{\sqrt{2\pi\xi_{t,s}}}\exp(\frac{1}{4}\frac{u^{2}}{\overline{p}(|t-s|)^{2}})\exp(-\frac{1}{2}\frac{u^{2}}{\xi_{t,s}})dudtds\leq\sqrt{2},\label{eq:temp-537}
\end{equation}
as $b\rightarrow\infty$. Therefore The Monotone Convergence Theorem
\ref{Thm. Monotone Convergence} implies that the r.r.v. $U$ is integrable
on $[0,1]^{2}\times\Omega$. Hence, by Fubini's Theorem, the function
\[
B\equiv\int_{0}^{1}\int_{0}^{1}Udtds\equiv\int_{0}^{1}\int_{0}^{1}\Psi(\frac{|V_{t}-V_{s}|}{\overline{p}(|t-s|)})dtds
\]
 is an integrable r.r.v, with expectation given by 
\[
E(B)\equiv E\int_{0}^{1}\int_{0}^{1}Udtds\leq\sqrt{2}.
\]
Consider each $\omega\in domain(B)$. Then
\begin{equation}
B(\omega)\equiv\int_{0}^{1}\int_{0}^{1}\Psi(\frac{|V(t,\omega)-V(s,\omega)|}{\overline{p}(|t-s|)})dtds.\label{eq:temp-499-1}
\end{equation}
In view of equality \ref{eq:temp-499-1}, Lemma \ref{Lem:GRR lemma}
implies that 
\begin{equation}
|V(t,\omega)-V(s,\omega)|\leq8\int_{0}^{|t-s|}\Psi^{-1}(\frac{4B(\omega)}{u^{2}})d\overline{p}(u)\equiv16\int_{0}^{|t-s|}\sqrt{\log(\frac{4B(\omega)}{u^{2}})}d\overline{p}(u).\label{eq:temp-3-1-2}
\end{equation}
The lemma is proved.

.
\end{proof}
\begin{thm}
\label{Thm.  Garsia}.\textbf{\emph{ (Garsia-Rodemich-Rumsey Theorem).}}
Let $\overline{p}:[0,1]\rightarrow[0,\infty)$ be a continuous increasing
function, with $\overline{p}(0)=0$, such that $\sqrt{-\log u}$ is
integrable relative to the distribution function $\overline{p}$ on
$[0,1]$. Thus 
\begin{equation}
\int_{0}^{1}\sqrt{-\log u}d\overline{p}(u)<\infty.\label{eq:temp-144-1}
\end{equation}
Let $\sigma:[0,1]\times[0,1]\rightarrow R$ be an arbitrary symmetric
positive definite function such that 
\begin{equation}
\bigvee_{0\leq s,t\leq1;|s-t|\leq u}(\Delta\sigma(s,t))^{1/2}\leq\overline{p}(u).\label{eq:temp-58}
\end{equation}
Then there exists an a.u. continuous centered Gaussian process $X:[0,1]\times\Omega\rightarrow R$
with $\sigma$ as covariance function. Moreover, there exists an integrable
r.r.v. $B$ with $EB\leq\sqrt{2}$ such that 
\[
|X(t,\omega)-X(s,\omega)|\leq16\int_{0}^{|t-s|}\sqrt{\log(\frac{4B(\omega)}{u^{2}})}d\overline{p}(u)
\]
for each $t,s\in[0,1]$, for each $\omega\in domain(B)$.
\end{thm}
\begin{proof}
1. As observed in the beginning of the proof of Lemma \ref{Lem. GRR modulus of continuity assuming a.u.contintuiy},
the positive definiteness of the function $\sigma$ implies that $\overline{p}(u)>0$
for each $u\in(0,1]$. 

2. Let 
\begin{equation}
F^{\sigma}\equiv\Phi_{covar,fjd}(\sigma)\label{eq:temp-515}
\end{equation}
be the consistent family of normal f.j.d.'s on the parameter set $[0,1]$
associated with mean function $0$ and the given covariance function
$\sigma$, as defined in equalities \ref{eq:temp-132} and \ref{eq:temp-3-1-1}
of Theorem \ref{Thm. Extension to Gaussian r.f.}. Let $Y:Q_{\infty}\times\Omega\rightarrow R$
be an arbitrary process with marginal distributions given by $F^{\sigma}|Q_{\infty}$,
the restriction of the family $F^{\sigma}$ of the normal f.j.d.'s
to the countable parameter subset $Q_{\infty}$.

3. By hypothesis, the function $\overline{p}$ is continuous at $0$,
with $\overline{p}(0)=0$. Hence there is a modulus of continuity
$\delta_{\overline{p},0}:(0,\infty)\rightarrow(0,\infty)$ such that
$\overline{p}(u)<\varepsilon$ for each $u$ with $0\leq u<\delta_{\overline{p},0}(\varepsilon)$,
for each $\varepsilon>0$.

4. Also by hypothesis, the function $\sqrt{-\log u}$ is integrable
relative to $\int_{0}^{1}\cdot d\overline{p}$. Hence there exists
a modulus of integrability $\delta_{\overline{p},1}:(0,\infty)\rightarrow(0,\infty)$
such that 
\[
\int_{(0,c]}\sqrt{-\log u}d\overline{p}(u)<\varepsilon
\]
for each $c\in[0,1]$ with $\int_{(0,c]}d\overline{p}=\overline{p}(c)<\delta_{\overline{p},1}(\varepsilon)$.

5. Let $b>0$ be arbitrary. Then 
\begin{equation}
\Psi^{-1}(\frac{4b}{u^{2}})\equiv2\sqrt{2\log b-2\log u}\leq2\sqrt{2}(\sqrt{\log b}+\sqrt{-\log u})\label{eq:temp-276}
\end{equation}
where the functions of $u$ on both ends have domain $(0,1]$ and
are continuous on $(0,1]$. Hence these functions are $\mathrm{measurable}$
relative to $\int_{0}^{1}\cdot d\overline{p}$. Since the right-hand
side of inequality \ref{eq:temp-276} is an integrable function of
$u$ relative to the integration $\int_{0}^{1}\cdot d\overline{p}$,
so is the function on the left-hand side. 

6. Define $m_{0}\equiv0$. Let $k\geq1$ be arbitrary, but fixed till
further notice. In view of Steps 3-5, there exists $m_{k}\equiv m_{k}(\delta_{\overline{p},0},\delta_{\overline{p},1})\geq m_{k-1}$
so large that 
\begin{equation}
16\int_{0}^{\Delta(m(k))}\Psi^{-1}(\frac{k^{2}}{u^{2}})d\overline{p}(u)<k^{-2},\label{eq:temp-45-2}
\end{equation}
where $\Delta_{m(k)}\equiv2^{-m(k)}$. Write $p_{0}\equiv p_{m(0)}=1$
and $p_{m(k)}\equiv2^{m(k)}$. 

7. Let $n\geq0$ be arbitrary. Define, as in Lemma \ref{Lem. Interpolation of Gaussian procee},
the interpolated process $Y^{(n)}:[0,1]\times\Omega\rightarrow R$
of $Y:Q_{\infty}\times\Omega\rightarrow R$ by conditional expectations
on $\{t_{0},\cdots,t_{n}\}$. Then Lemma \ref{Lem. Interpolation of Gaussian procee}
implies that (i) $Y^{(n)}$ is a centered Gaussian process, (ii) $Y^{(n)}$
is a.u. continuous, and (iii) $Y_{r}^{(n)}=Y_{r}$ for each $r\in\{t_{0},\cdots,t_{n}\}$.
Consequently the difference process 
\[
Z^{(k)}\equiv Z^{(m(k),m(k-1))}\equiv Y^{(m(k))}-Y^{(m(k-1))}
\]
is a.u. continuous. Note that from Condition (iii), we have
\[
Z_{0}^{(k)}\equiv Y_{0}^{(m(k))}-Y_{0}^{(m(k-1))}=Y_{0}-Y_{0}=0.
\]
For convenience define the trivial process $Z^{(0)}\equiv0$.

8. Let $t,s\in[0,1]$ be arbitrary with $|t-s|>0$. Then, since $\{t_{0},t_{1},\cdots,t_{j(m(k-1))}\}$
is a $\Delta_{m(k-1)}$-approximation of $[0,1]$, we have, by Lemma
\ref{Lem. Interpolation of Gaussian procee}, 
\begin{equation}
\xi_{k,t,s}\equiv E(Z_{t}^{(k)}-Z_{s}^{(k)})^{2}\equiv E(Z_{t}^{(m(k),m(k-1))}-Z_{s}^{(m(k),m(k-1))})^{2}\leq\overline{p}_{k}^{2}(|t-s|)\label{eq:temp-538}
\end{equation}
where
\[
\overline{p}_{k}(u)\equiv2(\overline{p}(u)\wedge\overline{p}(\Delta_{m(k-1)}))
\]
for each $u\geq0$. Note that $\overline{p}_{k}(u)$ is constant for
$u>\Delta_{m(k-1)}$. Hence the definition of Riemann-Stieljes integrals
implies that, for each nonnegative function $f$ on $[0,1]$ which
is integrable relative to the distribution function $\overline{p}$,
we have 
\[
\int_{0}^{1}f(u)d\overline{p}_{k}(u)=\int_{0}^{\Delta(m(k-1))}f(u)d\overline{p}_{k}(u)
\]
\begin{equation}
=2\int_{0}^{\Delta(m(k-1))}f(u)d\overline{p}(u)<\infty.\label{eq:temp-249-1}
\end{equation}
In particular
\begin{equation}
\int_{0}^{1}\sqrt{-\log u}d\overline{p}_{k}(u)=2\int_{0}^{\Delta(m(k-1))}\sqrt{-\log u}d\overline{p}(u)<\infty\label{eq:temp-249-1-1}
\end{equation}

9. Inequalities \ref{eq:temp-538} and \ref{eq:temp-249-1-1} show
that the a.u. continuous process $Z^{(k)}$ and the function $\overline{p}_{k}$
satisfy the conditions in the hypothesis of Lemma \ref{Lem. GRR modulus of continuity assuming a.u.contintuiy}.
Accordingly, there exists an integrable r.r.v. $B_{k}$ with $EB_{k}\leq\sqrt{2}$
such that 
\[
|Z^{(k)}(t,\omega)-Z^{(k)}(s,\omega)|\leq16\int_{0}^{|t-s|}\sqrt{\log(\frac{4B_{k}(\omega)}{u^{2}})}d\overline{p}_{k}(u)
\]
\begin{equation}
=16\int_{0}^{|t-s|\wedge\Delta(m(k-1))}\sqrt{\log(\frac{4B_{k}(\omega)}{u^{2}})}d\overline{p}(u)\label{eq:temp-3-1-2-1-1}
\end{equation}
for each $t,s\in[0,1]$, for each $\omega\in domain(B_{k})$. 

10. Let $\alpha_{k}\in(2^{-3}k^{2},2^{-2}k^{2})$ be arbitrary, and
define $A_{k}\equiv(B_{k}\leq\alpha_{k}).$ Chebychev's inequality
then implies that 
\[
P(A_{k}^{c})\equiv P(B_{k}>\alpha_{k})\leq\alpha_{k}^{-1}\sqrt{2}<2^{3}\sqrt{2}k^{-2}.
\]
Consider each $\omega\in A_{k}$. Then inequality \ref{eq:temp-3-1-2-1-1}
implies that, for each $t,s\in[0,1]$, we have
\[
|Z_{t}^{(k)}(\omega)-Z_{s}^{(k)}(\omega)|\leq16\int_{0}^{\Delta(m(k-1))}\sqrt{\log(\frac{4\alpha_{k}}{u^{2}})}d\overline{p}(u)
\]
\begin{equation}
\leq16\int_{0}^{\Delta(m(k-1))}\sqrt{\log(\frac{k^{2}}{u^{2}})}d\overline{p}(u)<(k-1)^{-2},\label{eq:temp-59}
\end{equation}
where the last inequality is by inequality \ref{eq:temp-45-2}. In
particular, if we set $s=0$ and recall that $Z_{0}^{(k)}=0$, we
obtain
\[
|Y_{t}^{(m(k))}(\omega)-Y_{t}^{(m(k-1))}(\omega)|\equiv|Z_{t}^{(k)}(\omega)|<(k-1)^{-2},
\]
where $\omega\in A_{k}$ is arbitrary, if $k>1$. Since $P(A_{k}^{c})<2^{3}\sqrt{2}k^{-2}$,
we conclude that $Y_{t}^{(m(k))}$ converges a.u. to the limit r.r.v.
$X_{t}\equiv\lim_{k\rightarrow\infty}Y_{t}^{(m(k))}$. Thus we obtain
the limiting process $X:[0,1]\times\Omega\rightarrow R$.

11. We will next prove that the process $X$ is a.u. continuous. To
that end, note that, since $Y^{(m(k))}$ is an a.u. continuous process
according to Condition (ii) in Step 7, there exist a $\mathrm{measurable}$
set $D_{k}$ with $P(D_{k}^{c})<k^{-1}$ and some $\delta_{k}>0$,
such that, for each $\omega\in D_{k}$, we have 
\[
|Y_{t}^{(m(k))}(\omega)-Y_{s}^{(m(k))}(\omega)|<k^{-1}
\]
for each $t,s\in[0,1]$ with $|s-t|<\delta_{k}$. Separately, define
the $\mathrm{measurable}$ set $C_{k}\equiv\bigcap_{h=k+1}^{\infty}A_{h}$.
Then $P(C_{k}^{c})\leq\sum_{h=k+1}^{\infty}2^{3}\sqrt{2}k^{-2}<2^{3}\sqrt{2}k^{-1}$. 

Now consider each $\omega\in D_{k}C_{k}$, and each $t,s\in[0,1]$
with $|s-t|<\delta_{k}$. Then
\[
X_{t}(\omega)=Y_{t}^{(m(k))}(\omega)+\sum_{h=k+1}^{\infty}(Y_{t}^{(m(h))}(\omega)-Y_{t}^{(m(h-1))}(\omega))\equiv Y_{t}^{(m(k))}(\omega)+\sum_{h=k+1}^{\infty}Z_{t}^{(h)}(\omega),
\]
with a similar equality when $t$ is replaced by $s$. Hence
\[
|X_{t}(\omega)-X_{s}(\omega)|\leq|Y_{t}^{(m(k))}(\omega)-Y_{s}^{(m(k))}(\omega)|+|\sum_{h=k+1}^{\infty}(Z_{t}^{(h)}(\omega)-Z_{s}^{(h)}(\omega))|
\]
\[
\leq|Y_{t}^{(m(k))}(\omega)-Y_{s}^{(m(k))}(\omega)|+\sum_{h=k+1}^{\infty}h^{-2}
\]
\[
\leq k^{-1}+k^{-1}=2k^{-1},
\]
where $\omega\in D_{k}C_{k}$ and $t,s\in[0,1]$ with $|s-t|<\delta_{k}$
are arbitrary. Since $P(D_{k}C_{k})^{c}<k^{-1}+2^{3}\sqrt{2}k^{-1}$
and $2k^{-1}$ are arbitrarily small if $k\geq1$ is sufficiently
large, we see that $X:[0,1]\times\Omega\rightarrow R$ is an a.u.
continuous process. Consequently, the process $X$ is continuous in
probability.

12. Now we will verify that the process $X$ is Gaussian, centered,
and has covariance function $\sigma$. Note that $X|Q_{\infty}=Y$.
Hence $X|Q_{\infty}$ has marginal distributions given by the family
$F^{\sigma}|Q_{\infty}$ of f.j.d.'s. Since the process $X$ and the
family $F^{\sigma}$ are continuous in probability, and since the
subset $Q_{\infty}$ is dense in the parameter set $[0,1]$, it follows
that $X$ has marginal distributions given by the family $F^{\sigma}$.
Thus $X$ is Gaussian, centered, and has covariance function $\sigma$.

13. Therefore
\[
\xi_{t,s}\equiv E(X_{t}-X_{s})^{2}=\Delta\sigma(t,s)\leq\overline{p}^{2}(|t-s|).
\]
In view of inequalities \ref{eq:temp-144-1} and \ref{eq:temp-58}
in the hypothesis, the conditions in Lemma \ref{Lem. GRR modulus of continuity assuming a.u.contintuiy}
are satisfied by the process $X$ and the function $\overline{p}$.
Hence Lemma \ref{Lem. GRR modulus of continuity assuming a.u.contintuiy}
implies the existence of an integrable r.r.v. $B$ with $EB\leq\sqrt{2}$
such that 
\begin{equation}
|X(t,\omega)-X(s,\omega)|\leq16\int_{0}^{|t-s|}\sqrt{\log(\frac{4B(\omega)}{u^{2}})}d\overline{p}(u)\label{eq:temp-3-1-2-1-2-1}
\end{equation}
 for each $t,s\in[0,1]$, for each $\omega\in domain(B)$, as desired.
\end{proof}

\chapter{a.u. Càdlàg Processes}

In this chapter, let $(S,d)$ be a locally compact metric space, with
a fixed reference point $x_{\circ}$. As usual, write $\widehat{d}\equiv1\wedge d$.
We will study processes $X:[0,\infty)\times\Omega\rightarrow S$ whose
sample paths are right continuous with left limits, or càdlàg (the
commonly used French acronym \textquotedbl{}continue à droite, limite
à gauche\textquotedbl{}). 

Classically, the proof of existence of such processes relies on Prokhorov's
Relative Compactness Theorem. As discussed in the beginning of Chapter
\ref{Chapter. Processes with continuous sample paths} of the present
book, this theorem implies the principle of infinite search. We will
therefore bypass Prokhorov's theorem, in favor of direct proofs using
Borel-Cantelli estimates.

In Section 1 a version of Skorokhod's definition of càdlàg functions
from $[0,\infty)$ to $S$. Each càdlàg function will come with a
modulus of càdlàg, much as a continuous function comes with a modulus
of continuity. In Section 2 we study a Skorokhod metric $d_{D}$ on
the space $D$ of càdlàg functions. 

In Section 3 we define an a.u. càdlàg process $X:[0,1]\times\Omega\rightarrow S$
as a process which is continuous in probability and which has, almost
uniformly, càdlàg sample functions. In Section 4, we introduce a $D$-regular
process $Z:Q_{\infty}\times\Omega\rightarrow S$, in terms of the
marginal distributions of $Z$, where $Q_{\infty}$ is the set of
dyadic rationals in $[0,1]$. We then prove, in Sections 4 and 5,
that a process $X:[0,1]\times\Omega\rightarrow S$ is a.u. càdlàg
iff its restriction $X|Q_{\infty}$ is $D$-regular, or equivalently,
iff $X$ is the extension, by right limit, of a $D$-regular process
$Z$. Thus we obtain a characterization of an a.u. càdlàg processes
in terms of conditions on its marginal distributions. Equivalently,
we have a procedure to construct an a.u. càdlàg process $X$ from
a consistent family $F$ of f.j.d.'s which is $D$-regular. We will
derive the modulus of a.u. càdlàg of $X$ from the given modulus of
$D$-regularity of $F$. 

In Section 6, we will prove that this construction is metrically continuous,
in epsilon-delta terms. Such continuity of construction also seems
to be hitherto unknown. In Sections 7 we apply the construction to
obtain a.u. càdlàg processes with strongly right continuous marginal
distributions; in Section 8, to a.u. càdlàg Martingales; in Section
9, to processes which are right Hoelder in a sense to be made precise
there. In Section 10, we state the generalization of definitions and
results in Sections 1-9, to the parameter interval $[0,\infty)$,
without giving the straightforward proofs.

Before proceeding, we remark that our constructive method for a.u.
càdlàg processes is by using certain accordion functions, defined
in Definition \ref{Def. Accordian function}, as time-varying boundaries
for hitting times. This will be clarified as we go along. This method
was first used in \cite{Chan74} to construct an a.u. càdlàg Markov
process from a given strongly continuous semigroup. 
\begin{defn}
\textbf{\label{Def. Specification of binaary approx of =00005B0,1=00005D-1}
(Notations for dyadic rationals).} For ease of reference, we restate
he following notations in Definition \ref{Def. Notations for dyadic rationals}
related to dyadic rationals. For each $m\geq0$, define $p_{m}\equiv2^{m}$,
$\Delta_{m}\equiv2^{-m}$, and recall the enumerated set of dyadic
rationals
\[
Q_{m}\equiv\{t_{0},t_{1},\cdots,t_{p(m)}\}=\{q_{m,0},\cdots,q_{m,p(m)}\}\equiv\{0,\Delta_{m},2\Delta_{m},\cdots,1\}\subset[0,1],
\]
where the second equality is equality of sets without the enumeration,
and recall the enumerated set
\[
Q_{\infty}\equiv\bigcup_{m=0}^{\infty}Q_{m}\equiv\{t_{0},t_{1},\cdots\},
\]
where the second equality is equality of sets without the enumeration.
Thust
\[
Q_{m}\equiv\{q_{m,0},\cdots,q_{m,p(m)}\}\equiv\{0,2^{-m},2\cdot2^{-m},\cdots,1\}
\]
is a $2^{-m}$-approximation of $[0,1]$, with $Q_{m}\subset Q_{m+1}$,
for each $m\geq0$.

Moreover, for each $m\geq0$, recall the enumerated set of dyadic
rationals
\[
\overline{Q}_{m}\equiv\{u_{0},u_{1},\cdots,u_{p(2m)}\}\equiv\{0,2^{-m},2\cdot2^{-m},\cdots,2^{m}\}\subset[0,2^{m}],
\]
and 
\[
\overline{Q}_{\infty}\equiv\bigcup_{m=0}^{\infty}\overline{Q}_{m}\equiv\{u_{0},u_{1},\cdots\},
\]
where the second equality is equality of sets without the enumeration.$\square$
\end{defn}

\section{Càdlàg Functions}

Recall some notations and  conventions. To minimize clutter, a subscripted
expression $a_{b}$ will be written interchangeably with $a(b)$.
For an arbitrary function $x$, we write $x(t)$ only with the explicit
or implicit condition that $t\in domain(x)$. If $X:A\times\Omega\rightarrow S$
is a random field, and if $B$ is a subset of $A$, then $X|B\equiv X|(B\times\Omega)$
denotes the random field obtained by restricting the parameter set
to $B$. 
\begin{defn}
\label{Def. Pointwise continuity and right continuity on =00005B0,a=00005D-1}
\textbf{(Pointwise left- and right limits).} Let $Q$ be an arbitrary
subset of $[0,\infty)$. Let the function $x:Q\rightarrow S$ be arbitrary
such that $domain(x)$ is dense in $Q$. Let the point $t\in Q$ be
arbitrary.

The function $x$ is said to be \emph{right continuous} \index{right continuity}
at a point $t\in domain(x)$ if $\lim_{r\rightarrow t;r\geq t}x(r)=x(t)$.
The function $x$ is said to have a \emph{left limit at a point }$t\in Q$
\index{left limit} if $\lim_{r\rightarrow t;r<t}x(r)$ exists. 

Suppose, for each $t\in Q$ such that $\lim_{r\rightarrow t;r\geq t}x(r)$
exists, we have $t\in domain(x)$. Then we say that the function $x$
is \index{right complete}\emph{right complete}.
\end{defn}
$\square$

Recall that the function $x$ is said to be \emph{continuous at $t$}\index{pointwise continuity}
if $t\in domain(x)$ and if $\lim_{r\rightarrow t}x(r)=x(t)$. Trivially,
if $x$ is continuous at $t$ then it is  right continuous and has
left limit at $t$.

The next definition is essentially Skorokhod's characterization of
càdlàg functions. 
\begin{defn}
\label{Def. Cadlag functions, D=00005B0,1=00005D} \textbf{(Càdlàg
function on $[0,1]$).} Let $(S,d)$ be a locally compact metric space.
Let $x:[0,1]\rightarrow S$ be a function such that $domain(x)$ contains
the enumerated set $Q_{\infty}$ of dyadic rationals in $[0,1]$.
Suppose the following conditions are satisfied.

1. (Right continuity). The function $x$ is right continuous at each
$t\in domain(x)$, and is continuous at $t=1$. 

2. (Right completeness).  Let $t\in[0,1]$ be arbitrary. If $\lim_{r\rightarrow t;r>t}x(r)$
exists, then $t\in domain(x)$.

3. (Approximation by step functions). For each $\varepsilon>0$, there
exist $\delta_{cdlg}(\varepsilon)>0$, $p\geq1$, and a sequence
\begin{equation}
0=\tau_{0}<\tau_{1}<\cdots<\tau_{p-1}<\tau_{p}=1\label{eq:temp-331-3}
\end{equation}
in $domain(x)$, such that (i) for each $i=1,\cdots,p$, we have
\[
\tau_{i}-\tau_{i-1}\geq\delta_{cdlg}(\varepsilon)
\]
and (ii) for each $i=0,\cdots,p-1$, we have
\begin{equation}
d(x,x(\tau_{i}))\leq\varepsilon,\label{eq:temp-307-4}
\end{equation}
on the interval $\theta_{i}\equiv[\tau_{i},\tau_{i+1})$ or $\theta_{i}\equiv[\tau_{i},\tau_{i+1}]$
according as $i\leq p-2$ or $i=p-1$. We will call $(\tau_{i})_{i=0,\cdots,p}$
a sequence of \emph{$\varepsilon$-division points\index{varepsilon-division points@$\varepsilon$-division points}}
of $x$ with separation at least $\delta_{cdlg}(\varepsilon)$.

Then $x$ said to be a \emph{càdlàg} \index{càdlàg functions} function
on $[0,1]$ with values in $S$, with the operation $\delta_{cdlg}$
as a \emph{modulus of càdlàg}\index{modulus of càdlàg}. Here we let
brevity supersede grammar.

We will let $D[0,1]$ denote the set of  càdlàg functions. Two members
of $D[0,1]$ are considered equal if they are equal as functions,
i.e. if they have the same domain and have equal values in the common
domain.

$\square$
\end{defn}
Note that Condition 3 implies that the end points $0,1$ are in $domain(x)$.
Condition 3 implies also that $p\leq\delta_{cdlg}(\varepsilon)^{-1}$.
Let $x,y\in D[0,1]$ be arbitrary, with moduli of càdlàg $\delta_{cdlg},\delta'_{cdlg}$
respectively. Then the operation $\delta_{cdlg}\wedge\delta'_{cdlg}$
is obviously a common modulus of càdlàg of $x,y$. The next lemma
is a simple consequence of right continuity, and generalizes its counterpart
for $C[0,1]$.
\begin{lem}
\label{Lem.  If x=00003Dy on dense A then x=00003Dy. If f(x,y)<=00003Dc on A then...-1}
\textbf{\emph{(A càdlàg function is uniquely determined by its values
on a dense subset of its domain).}} Let $x,y\in D[0,1]$ be arbitrary.
Suppose $B\equiv domain(x)\cap domain(y)$ contains a dense subset
$A$ of $[0,1]$. Then the following holds.

1. Let $t\in B$ and $\alpha>0$ be arbitrary. Then there exists $r\in[t,t+\alpha)\cap A$
such that 
\begin{equation}
d(x(t),x(r))\vee d(y(t),y(r))\leq\alpha.\label{eq:temp-45}
\end{equation}

2. Let $f:S^{2}\rightarrow R$ be a uniformly continuous function.
Let $c\in R$ be arbitrary such that $f(x(r),y(r))\leq c$ for each
$r\in A$. Then $f(x,y)\leq c$. In other words $f(x(t),y(t))\leq c$
for each $t\in domain(x)\cap domain(y)$. The same assertion holds
when $"\leq"$ is replaced by $"\geq"$ or by $"="$. In particular,
if $d(x(r),y(r))\leq\varepsilon$ for each $r\in A$, for some $\varepsilon>0$,
then $f(x,y)\leq\varepsilon$ on $domain(x)\cap domain(y)$.

3. Suppose $x(r)=y(r)$ for each $r\in A$. Then $x=y$. In other
words, $domain(x)=domain(y)$ and $x(t)=y(t)$ for each $t\in domain(x)$. 

4. Let $\lambda:[0,1]\rightarrow[0,1]$ be an arbitrary continuous
and increasing function with $\lambda(0)=0$ and $\lambda(1)=1$.
Then $x\circ\lambda\in D[0,1]$.
\end{lem}
\begin{proof}
Let $\delta_{cdlg}$ be a common modulus of càdlàg of $x$ and $y$. 

1. Let $t\in B$ and $\alpha>0$ be arbitrary. Let $(\tau_{i})_{i=0,\cdots,p}$
and $(\tau'_{i})_{i=0,\cdots,p'}$ be sequences of $\frac{\alpha}{2}$-division
points of $x$ and $y$ respectively, with separation at least $\delta_{cdlg}(\frac{\alpha}{2})$.
Then $\tau_{p-1}\vee\tau'_{p'-1}<1$. Hence either (i) $t<1$, or
(ii) $\tau_{p-1}\vee\tau'_{p'-1}<t$. Consider Case (i). Since $x,y$
are right continuous at $t$, according to Definition \ref{Def. Cadlag functions, D=00005B0,1=00005D},
and since $A$ is dense in $[0,1]$, there exists $r$ in $A\cap[t,1\wedge(t+\alpha))$
such that 
\[
d(x(t),x(r))\vee d(y(t),y(r))\leq\alpha,
\]
as desired. Consider Case (ii). Take $r\in(\tau_{p-1}\vee\tau'_{p'-1},t)\cap A$.
Then $t,r\in[\tau_{p-1},1]\cap[\tau'_{p'-1},1].$ Hence Condition
3 in Definition \ref{Def. Cadlag functions, D=00005B0,1=00005D} implies
that
\[
d(x(t),x(r))\leq d(x(\tau_{p-1}),x(t))\vee d(x(\tau_{p-1}),x(r))\leq\frac{\alpha}{2}+\frac{\alpha}{2}=\alpha.
\]
Similarly, $d(y(t),y(r))\leq\alpha.$ Assertion 1 is proved.

2. Let $t\in B$ be arbitrary. By Assertion 1, for each $k\geq1$,
there exists $r_{k}\in[t,t+\frac{1}{k})A$ such that inequality 
\[
d(x(t),x(r_{k}))\vee d(y(t),y(r_{k}))\leq\frac{1}{k}.
\]
holds. Hence, by right continuity of $x,y$ and continuity of $f$,
we have 
\[
f(x(t),y(t))=\lim_{k\rightarrow\infty}f(x(r_{k}),y(r_{k}))\leq c.
\]

3. By hypothesis, $d(x(r),y(r))=0$ for each $r\in A$. Let $t\in domain(x)$
and $\varepsilon>0$ be arbitrary. Since $x$ is right continuous
at $t$, there exists $c>0$ such that 
\begin{equation}
d(x(r),x(t))<\varepsilon\label{eq:temp-3-2}
\end{equation}
for each $r\in[t,t+c)\cap domain(x)$. Consider each $s\in[t,t+c)\cap domain(y)$.
Let $\alpha\equiv(t+c-s)\wedge\varepsilon$. By Assertion 1 applied
to the pair $y,y$ in $D[0,1]$, there exists $r\in[s,s+\alpha)\cap A$
such that $d(y(s),y(r))\leq\alpha\leq\varepsilon$. Then $r\in[t,t+c)\cap A$,
whence inequality \ref{eq:temp-3-2} holds. Combining,
\[
d(y(s),x(t))\leq d(y(s),y(r))+d(y(r),x(r))+d(x(r),x(t))<\varepsilon+0+\varepsilon=2\varepsilon.
\]
Since $\varepsilon>0$ is arbitrary, we see that $\lim_{s\rightarrow t;s>t}y(s)$
exists and is equal to $x(t)$. Hence the right completeness Condition
2 in Definition \ref{Def. Cadlag functions, D=00005B0,1=00005D} implies
that $t\in domain(y)$. Condition 1 in Definition \ref{Def. Cadlag functions, D=00005B0,1=00005D}
then implies that $y(t)=\lim_{s\rightarrow t;s>t}y(s)=x(t)$. Since
$t\in domain(x)$ is arbitrary, we conclude that $domain(x)\subset domain(y)$
and $x=y$ on $domain(x)$. By symmetry, $domain(x)=domain(y)$.

4. Since $\lambda$ is continuous and increasing, it has an inverse
$\lambda^{-1}$ which is also continuous and increasing, with some
modulus of continuity $\bar{\delta}$. Let $\delta_{cdlg}$ be a modulus
of càdlàg of $x$. We will prove that $x\circ\lambda$ is càdlàg,
with $\delta_{1}\equiv\bar{\delta}\circ\delta_{cdlg}$ as a modulus
of càdlàg. To that end, let $\varepsilon>0$ be arbitrary. Let 
\begin{equation}
0\equiv\tau_{0}<\tau_{1}<\cdots<\tau_{p-1}<\tau_{p}=1\label{eq:temp-331-2-2}
\end{equation}
be a sequence of \emph{$\varepsilon$-}division points\emph{\index{varepsilon-division points@$\varepsilon$-division points}}
of $x$ with separation at least $\delta_{cdlg}(\varepsilon)$. Thus,
for each $i=1,\cdots,p$, we have $\tau_{i}-\tau_{i-1}\geq\delta_{cdlg}(\varepsilon).$
Suppose $\lambda\tau_{i}-\lambda\tau_{i-1}<\bar{\delta}(\delta_{cdlg}(\varepsilon))$
for some $i=1,\cdots,p$. Then, since $\bar{\delta}$ is a modulus
of continuity of the inverse function $\lambda^{-1}$, it follows
that 
\[
\tau_{i}-\tau_{i-1}=\lambda^{-1}\lambda\tau_{i}-\lambda^{-1}\lambda\tau_{i-1}<\delta_{cdlg}(\varepsilon),
\]
 a contradiction. Hence 
\[
\lambda\tau_{i}-\lambda\tau_{i-1}\geq\bar{\delta}(\delta_{cdlg}(\varepsilon))\equiv\delta_{1}(\varepsilon)
\]
for each $i=1,\cdots,p$. Moreover, for each $i=0,\cdots,p-1$, we
have 
\begin{equation}
d(x,x(\tau_{i}))\leq\varepsilon\label{eq:temp-307-2-2}
\end{equation}
on the interval $\theta_{i}\equiv[\tau_{i},\tau_{i+1})$ or $\theta_{i}\equiv[\tau_{i},\tau_{i+1}]$
according as $i\leq p-2$ or $i=p-1$. Since the function $\lambda$
is increasing, it follows that 
\begin{equation}
d(x\circ\lambda,x\circ\lambda(\tau_{i}))\leq\varepsilon\label{eq:temp-307-2-1-1}
\end{equation}
on the interval $\theta'_{i}\equiv[\lambda^{-1}\tau_{i},\lambda^{-1}\tau_{i+1})$
or $\theta'_{i}\equiv[\lambda^{-1}\tau_{i},\lambda^{-1}\tau_{i+1}]$
according as $i\leq p-2$ or $i=p-1$. Thus the sequence 
\[
0=\lambda^{-1}\tau_{0}<\lambda^{-1}\tau_{1}<\cdots<\lambda^{-1}\tau_{p-1}<\lambda^{-1}\tau_{p}=1
\]
is a sequence of \emph{$\varepsilon$-}division points\emph{\index{varepsilon-division points@$\varepsilon$-division points}}
of $x$ with separation at least $\delta_{1}(\varepsilon)$. Condition
3 in Definition \ref{Def. Cadlag functions, D=00005B0,1=00005D} has
been proved for the function $x\circ\lambda$. In view of the monotonicity
and continuity of the function $\lambda$, the other conditions can
also be easily verified. Accordingly, the function $x\circ\lambda$
is càdlàg, with a modulus of càdlàg $\delta_{1}$. 
\end{proof}
\begin{prop}
\label{Prop. Points of continuitys of cadlag functions-1} \textbf{\emph{(Points
of continuity of càdlàg function)}}. Let $x\in D[0,1]$ be arbitrary
with a modulus of càdlàg $\delta_{cdlg}$. Then the function $x$
on $[0,1]$ is continuous at the end points $0$ and $1$. Moreover,
$domain(x)$ contains all but countably many points in $[0,1]$. 

More precisely, for each $k\geq1$, let $(\tau_{k,i})_{i=0,\cdots,p(k)}$
be a sequence of $\frac{1}{k}$-division points of $x$ with separation
at least $\delta_{cdlg}(\frac{1}{k})$. Then the following holds.

1. Define the set $A\equiv\bigcap_{k=1}^{\infty}\bigcup_{i=0}^{p(k)-1}\theta{}_{k,i}$,
where $\theta{}_{k,i}\equiv[\tau_{k,i},\tau_{k,i+1})$ or $\theta{}_{k,i}\equiv[\tau_{k,i},\tau_{k,i+1}]$
according as $i=0,\cdots,p_{k}-2$ or $i=p_{k}-1$. Then the set $A$
contains all but countably many points in $[0,1]$, and is a subset
of $domain(x)$. 

2. Define the set $A'\equiv\bigcap_{k=1}^{\infty}\bigcup_{i=0}^{p(k)-1}\theta'_{k,i}$,
where $\theta'_{k,i}\equiv[0,\tau_{k,1}]$ or $\theta'_{k,i}\equiv(\tau_{k,i},\tau_{k,i+1}]$
according as $i=0$ or $i=1,\cdots,p_{k}-1$. Then the set $A'$ contains
all but countably many points in $[0,1]$, and the function $x$ has
a left limit at each $t\in A'$.

3. Define the set $A''\equiv\bigcap_{k=1}^{\infty}\bigcup_{i=0}^{p(k)-1}(\tau_{k,i},\tau_{k,i+1})$.
Then the set $A''$ contains all but countably many points in $[0,1]$,
and the function $x$ is continuous at each $t\in A''$. 

4. The function $x$ is bounded on $domain(x)$. Specifically, 
\[
d(x_{\circ},x(t))\leq b\equiv\bigvee_{i=0}^{p(1)-1}d(x_{\circ},x(\tau_{1,i}))+1
\]
for each $t\in domain(x)$, where $x_{\circ}$ is an arbitrary, but
fixed, reference point in $S$.
\end{prop}
\begin{proof}
By Definition \ref{Def. Cadlag functions, D=00005B0,1=00005D}, we
have $1\in domain(x)$. Condition 3 in Definition \ref{Def. Cadlag functions, D=00005B0,1=00005D}
implies that $0=\tau_{1,0}\in domain(x)$ and that $x$ is continuous
at $0$ and $1$.

1. Let $t\in A$ be arbitrary. Let $k\geq1$ be arbitrary. Then $t\in\theta{}_{k,i}$
for some $i_{.}=0,\cdots,p_{k}-1$. Let $\delta_{0}\equiv\tau_{k,i+1}-t$
or $\delta_{0}\equiv2$ according as $i=0,\cdots,p_{k}-2$ or $i=p_{k}-1$.
Then $domain(x)\cap[t,t+\delta_{0})$ is a nonempty subset of $\theta{}_{k,i}$.
Moreover, by Condition 3 of Definition \ref{Def. Cadlag functions, D=00005B0,1=00005D},
we have $d(x(r),x(\tau_{k,i}))\leq\frac{1}{k}$ for each $r\in domain(x)\cap[t,t+\delta_{0})$.
Hence $d(x(r),x(s))\leq\frac{2}{k}$ for each $r,s\in domain(x)\cap[t,t+\delta_{0})$.
Since $\frac{2}{k}$ is arbitrarily small, and since the metric space
$(S,d)$ is complete, we see that $\lim_{r\rightarrow t;r\geq t}x(r)$
exists. The right completeness Condition 2 of Definition \ref{Def. Cadlag functions, D=00005B0,1=00005D}
therefore implies that $t\in domain(x)$. We conclude that $A\subset domain(x)$. 

2. Let $t\in A'$ be arbitrary. Let $k\geq1$ be arbitrary. Then $t\in\theta'_{k,i}$
for some $i_{.}=0,\cdots,p_{k}-1$. Let $\delta_{0}\equiv2$ or $\delta_{0}\equiv t-\tau_{k,i}$
according as $i_{.}=0$ or $i_{.}=1,\cdots,p_{k}-1$. Then $domain(x)\cap(t-\delta_{0},t)$
is a nonempty subset of $\theta{}_{k,i}$. Moreover, by Condition
3 of Definition \ref{Def. Cadlag functions, D=00005B0,1=00005D},
we have $d(x,x(\tau_{k,i}))\leq\frac{1}{k}$ for each $r\in domain(x)\cap(t-\delta_{0},t)$.
Argument similar to the previous paragraph then shows that $\lim_{r\rightarrow t;r<t}x(r)$
exists. 

3. Since $A''\subset A$, we have $A'\subset domain(x)$, thanks to
Assertion 1 above. Let $t\in A''$ be arbitrary. Let $k\geq1$ be
arbitrary. Then $t\in(\tau_{k,i},\tau_{k,i+1})$ for some $i_{.}=0,\cdots,p_{k}-1$.
Hence, by Condition 3 of Definition \ref{Def. Cadlag functions, D=00005B0,1=00005D},
we have $d(x(r),x(t))\leq\frac{2}{k}$ for each $r\in domain(x)\cap(\tau_{k,i},\tau_{k,i+1})$.
We conclude that the function $x$ is continuous at $t$ . 

4. Finally, observe that each of the sets $A$, $A'$, and $A''$
contains the metric complement of the countable subset $\{\tau_{k,i}\}$.
Thus each contains all but countably many points in $[0,1]$, and
is dense in $[0,1]$. Now let $t\in A\subset domain(x)$ be arbitrary.
Then $t\in\theta_{1,i}$ for some $i=0,\cdots,p_{1}-1$. Hence
\[
d(x_{\circ},x(t))\leq d(x_{\circ},x(\tau_{1,i}))+d(x(\tau_{1,i}),x(t)))\leq b\equiv\bigvee_{j=0}^{p(1)-1}d(x_{\circ},x(\tau_{1,j}))+1.
\]
Since $A$ is dense in $[0,1]$ and since the function $d:S^{2}\rightarrow R$
is uniformly continuous, Lemma \ref{Lem.  If x=00003Dy on dense A then x=00003Dy. If f(x,y)<=00003Dc on A then...-1}
implies that $d(x_{\circ},x(r))\leq b$ for each $r\in domain(x)$.
\end{proof}
\begin{prop}
\label{Prop. modulus of right continuity a.e.} \textbf{\emph{(a.u.
right continuity of càdlàg function)}}. Let $x\in D[0,1]$ be arbitrary,
with a modulus of càdlàg $\delta_{cdlg}$. For each $\alpha>0$, let
$h\equiv[2+0\vee-\log_{2}\alpha]_{1}$, and define
\[
\delta_{rc}(\alpha,\delta_{cdlg})\equiv2^{-h}\delta_{cdlg}(2^{-h})>0.
\]
Let $\varepsilon>0$ be arbitrary. Then there exist a Lebesgue measurable
subset $A$ of $domain(x)$ with Lebesgue measure $\mu(A)<\varepsilon$,
such that for each $\alpha\in(0,\varepsilon)$, we have
\[
d(x(t),x(s))<\alpha
\]
for each $t\in A\cap domain(x)$ and $s\in[t,t+\delta_{rc}(\alpha,\delta_{cdlg}))\cap domain(x)$. 
\end{prop}
\begin{proof}
1. Let $h\geq0$ be arbitrary. Write $\alpha_{h}\equiv2^{-h}$. Then
there exist an integer $p_{h}\geq1$ and a sequence
\begin{equation}
0=\tau_{h,0}<\tau_{h,1}<\cdots<\tau_{h,p-1}<\tau_{h,p(h)}=1\label{eq:temp-331-3-1}
\end{equation}
in $domain(x)$, such that (i) for each $i=1,\cdots,p_{h}$, we have
\begin{equation}
\tau_{h,i}-\tau_{h,i-1}\geq\delta_{cdlg}(\alpha_{h})\label{eq:temp-368}
\end{equation}
and (ii) for each $i=0,\cdots,p_{h}-1$, we have
\begin{equation}
d(x,x(\tau_{h,i}))\leq\alpha_{h},\label{eq:temp-307-4-2}
\end{equation}
on the interval $\theta_{h,i}\equiv[\tau_{h,i},\tau_{h,i+1})$ or
$\theta_{h,i}\equiv[\tau_{h,i},\tau_{h,i+1}]$ according as $i\leq p_{h}-2$
or $i=p_{h}-1$.

2. Let $i=0,\cdots,p_{h}-1$ be arbitrary. Define 
\begin{equation}
\overline{\theta}_{h,i}\equiv[\tau_{h,i},\tau_{h,i+1}-\alpha_{h}(\tau_{h,i+1}-\tau_{h,i}))\subset\theta_{h,i}.\label{eq:temp-479-1}
\end{equation}
Define $\overline{\theta}_{h}\equiv\bigcup_{i=0}^{p(h)-1}\overline{\theta}_{h,i}$.
Then 
\[
\mu(\overline{\theta}_{h})=\sum_{i=0}^{p(h)-1}\mu(\overline{\theta}_{h,i})=\sum_{i=0}^{p(h)-1}(\tau_{h,i+1}-\tau_{h,i})(1-\alpha_{h})=1-\alpha_{h},
\]
whence $\mu\overline{\theta}_{h}{}^{c}=\alpha_{h}\equiv2^{-h}$, where
$h\geq0$ is arbitrary. 

3. Now let $\varepsilon>0$ be arbitrary, and let $k\equiv[1+0\vee-\log_{2}\varepsilon]_{1}$.
Define $A\equiv\bigcup_{h=k}^{\infty}\overline{\theta}_{h}{}^{c}$.
Then 
\[
\mu(A)\leq\sum_{h=k+1}^{\infty}2^{-h}=2^{-k}<\varepsilon.
\]
Consider each $t\in A^{c}\cap domain(x)$. Let $\alpha\in(0,\varepsilon)$
be arbitrary, and let $h\equiv[2+0\vee-\log_{2}\alpha]_{1}$. Then
\[
h>2+0\vee-\log_{2}\alpha>2+(1+0\vee-\log_{2}\varepsilon)>k,
\]
Hence $h\geq k+1$ and so $t\in A^{c}\subset\overline{\theta}_{h}$.
Therefore $t\in\overline{\theta}_{h,i}\equiv[\tau_{h,i},\tau_{h,i+1}-\alpha_{h}(\tau_{h,i+1}-\tau_{h,i}))$
for some $i=0,\cdots,p_{h}-1$. Moreover, 
\[
c
\]
Now let $s\in[t,t+\delta_{rc}(\alpha,\delta_{cdlg}))\cap domain(x)$
be arbitrary. Then
\[
\tau_{h,i}\leq s\leq t+\delta_{rc}(\alpha,\delta_{cdlg})
\]
\[
<\tau_{h,i+1}-\alpha_{h}(\tau_{h,i+1}-\tau_{h,i})+\delta_{rc}(\alpha,\delta_{cdlg})
\]
\[
\leq\tau_{h,i+1}-2^{-h}\delta_{cdlg}(2^{-h})+2^{-h}\delta_{cdlg}(2^{-h})=\tau_{h,i+1}
\]
Hence $s,t\in[\tau_{h,i},\tau_{h,i+1})$. It follows that 
\[
d(x(s),x(\tau_{h,i}))\vee d(x(t),x(\tau_{h,i}))\leq\alpha_{h}
\]
and therefore that 
\[
d(x(s),x(t))\leq2\alpha_{h}=2^{-h+1}<\alpha.
\]
\end{proof}
The next proposition shows that if a function satisfies all the conditions
in Definition \ref{Def. Cadlag functions, D=00005B0,1=00005D} except
perhaps the right completeness Condition 2, then it can be right completed
and extended to a càdlàg function. This is analogous to the completion
of a uniformly continuous function on a dense subset of $[0,1]$.
\begin{prop}
\label{Prop. Cadlag completion}\textbf{\emph{ (Right-limit extension
and càdlàg completion).}} Let $(S,d)$ be a locally compact metric
space. Suppose $Q=[0,1]$ or $Q=[0,\infty)$. Let $x:Q\rightarrow S$
be a function whose domain is dense in $Q$, and which is right continuous
at each $t\in domain(x)$. Define its right-limit extension $\overline{x}:Q\rightarrow S$
by 
\begin{equation}
domain(\overline{x})\equiv\{t\in Q;\lim_{r\rightarrow t;r\geq t}x(r)\;exists\},\label{eq:temp-417-1}
\end{equation}
and by 
\begin{equation}
\overline{x}(t)\equiv\lim_{r\rightarrow t;r\geq t}x(r)\label{eq:temp-395-1-1}
\end{equation}
for each $t\in domain(\overline{x})$. Then the following holds.

1. The function $\overline{x}$ is right continuous at each $t\in domain(\overline{x})$.

2. Suppose $t\in Q$ is such that $\lim_{r\rightarrow t;r\geq t}\overline{x}(r)$
exists. Then $t\in domain(\overline{x})$. 

3. Suppose $Q=[0,1]$. Suppose, in addition, that $\delta_{cdlg}:(0,\infty)\rightarrow(0,\infty)$
is an operation such that $x$ and $\delta_{cdlg}$ satisfy Conditions
3 in Definition \ref{Def. Cadlag functions, D=00005B0,1=00005D}.
Then $\overline{x}\in D[0,1]$. Moreover, $\overline{x}$ has $\delta_{cdlg}$
as a modulus of càdlàg. Furthermore, $x=\overline{x}|domain(x)$\textup{.}
We will then call $\overline{x}$ the \index{càdlàg completion}\emph{càdlàg
completion} of $x$. 
\end{prop}
\begin{proof}
1. Since, by hypothesis, $x$ is right continuous at each $t\in domain(x)$,
it follows from the definition of $\overline{x}$ that $domain(x)\subset domain(\overline{x})$
and that $\overline{x}=x$ on $domain(x)$. In other words, $x=\overline{x}|domain(x).$
Since $domain(x)$ is, by hypothesis, dense in $Q$, so is $domain(\overline{x})$.
Now let $t\in domain(\overline{x})$ and $\varepsilon>0$ be arbitrary.
Then, by the defining equality \ref{eq:temp-395-1-1}, 
\[
\overline{x}(t)\equiv\lim_{r\rightarrow t;r\geq t}x(r).
\]
Hence there exists $\delta_{0}>0$ such that 
\begin{equation}
d(\overline{x}(t),x(r))\leq\varepsilon\label{eq:temp-396-1-1}
\end{equation}
for each $r\in domain(x)\cap[t,t+\delta_{0})$. Let $s\in domain(\overline{x})\cap[t,t+\delta_{0})$
be arbitrary. Then, again by the defining equalities \ref{eq:temp-417-1}
and \ref{eq:temp-395-1-1}, there exists a sequence $(r_{j})_{j=1,2,\cdots}$
in $domain(x)\cap[s,t+\delta_{0})$ such that $r_{j}\rightarrow s$
and $\overline{x}(s)=\lim_{j\rightarrow\infty}x(r_{j}).$ For each
$j\geq1$, we then have $r_{j}\in domain(x)\cap[t,t+\delta_{0})$.
Hence inequality \ref{eq:temp-396-1-1} holds for $r=r_{j}$, for
each $j\geq1$. Letting $j\rightarrow\infty$, we therefore obtain
$d(\overline{x}(t),\overline{x}(s))\leq\varepsilon$. Since $\varepsilon>0$
is arbitrary, we conclude that $\overline{x}$ is right continuous
at $t$. Assertion 1 has thus been verified.

2. Next suppose $\lim_{r\rightarrow t;r\geq t}\overline{x}(r)$ exists.
Then, since $x=\overline{x}|domain(x)$, the right limit 
\[
\lim_{r\rightarrow t;r\geq t}x(r)=\lim_{r\rightarrow t;r\geq t}\overline{x}(r)
\]
exists. Hence $t\in domain(\overline{x})$ by the defining equality
\ref{eq:temp-417-1}. Condition 2 of Definition \ref{Def. Cadlag functions, D=00005B0,1=00005D}
has been proved for $\overline{x}$. 

3. Now let $\varepsilon>0$ be arbitrary. Because $x=\overline{x}|domain(x)$,
each sequence $(\tau_{i})_{i=0,\cdots,p}$ of \emph{$\varepsilon$-}division
points of $x$, with separation at least $\delta_{cdlg}(\varepsilon)$,
is also a sequence of \emph{$\varepsilon$-}division points of $\overline{x}$
with separation at least $\delta_{cdlg}(\varepsilon)$. Therefore
Condition 3 in Definition \ref{Def. Cadlag functions, D=00005B0,1=00005D}
holds for $\overline{x}$ and the operation $\delta_{cdlg}$. Summing
up, the function $\overline{x}$ is càdlàg, with $\delta_{cdlg}$
as a modulus of càdlàg.
\end{proof}
The next definition introduces simple càdlàg functions as càdlàg completion
of step functions.
\begin{defn}
\label{Def. simple Cadlag functions} \textbf{(Simple càdlàg function).}
Let $0=\tau_{0}<\cdots<\tau_{p-1}<\tau_{p}=1$ be an arbitrary sequence
in $[0,1]$ such that 
\[
\bigwedge_{i=1}^{p`}(\tau_{i}-\tau_{i-1})\geq\delta_{0}
\]
for some $\delta_{0}>0$. Let $x_{0},\cdots,x_{p-1}$ be an arbitrary
sequence in $S$ . 

Define a function $z:[0,1]\rightarrow S$ by 
\begin{equation}
domain(z)\equiv\bigcup_{i=0}^{p-1}\theta_{i},\label{eq:temp-409-1}
\end{equation}
where $\theta_{i}\equiv[\tau_{i},\tau_{i+1})$ or $\theta_{i}\equiv[\tau_{i},\tau_{i+1}]$
according as $i=0,\cdots,p-2$ or $i=p-1$, and by
\[
z(r)\equiv x_{i}
\]
for each $r\in\theta_{i}$, for each $i=0,\cdots,p-1$. Let $x\equiv\overline{Z}\in D[0,1]$
be the càdlàg completion  of $z$. Then $x$ is called the \emph{simple}
\emph{càdlàg} \index{simple càdlàg functions} \emph{function} determined
by the pair of sequences $((\tau_{i})_{i=0,\cdots,p-1},(x_{i})_{i=0,\cdots,p-1})$.
In symbols, we then write
\[
x\equiv\Phi_{smpl}((\tau_{i})_{i=0,\cdots,p-1},(x_{i})_{i=0,\cdots,p-1}),
\]
or simply $x\equiv\Phi_{smpl}((\tau_{i}),(x_{i}))$ when the range
of subscripts is understood. The sequence $(\tau_{i})_{i=0,\cdots,p}$
is then called the sequence of \emph{division points}\index{division points of simple càdlàg function}
of the simple càdlàg function $x$. The next lemma verifies that $x$
is a well-defined càdlàg function, with the constant operation $\delta_{cdlg}(\cdot)\equiv\delta_{0}$
as a modulus of càdlàg.

$\square$
\end{defn}
\begin{lem}
\label{Lem. Simple cadlag functions are cadlag-1} \textbf{\emph{(Simple
càdlàg functions are well defined).}} Use the notations and assumptions
in Definition \ref{Def. simple Cadlag functions}. Then $z$ and $\delta_{cdlg}$
satisfy the conditions in Proposition \ref{Prop. Cadlag completion}.
Accordingly, the càdlàg completion $\overline{Z}\in D[0,1]$ of $z$
is well-defined. 
\end{lem}
\begin{proof}
First note that $domain(z)$ contains the metric complement of $\{\tau_{1},\cdots,\tau_{p}\}$
in $[0,1]$. Hence $domain(z)$ is dense in $[0,1]$. Let $t\in domain(z)$
be arbitrary. Then $t\in\theta_{i}$ for some $i=0,\cdots,p-1$. Hence,
for each $r\in\theta_{i}$, we have $z(r)\equiv x_{i}\equiv z(t)$.
Therefore $z$ is right continuous at $t$. Conditions 1 in Definition
\ref{Def. Cadlag functions, D=00005B0,1=00005D} has been verified
for $z$. The proof of Condition 3 in Definition \ref{Def. Cadlag functions, D=00005B0,1=00005D}
for $z$ and $\delta_{cdlg}$ being trivial, the conditions in Proposition
\ref{Prop. Cadlag completion} are satisfied.
\end{proof}
\begin{lem}
\label{Lem. Inserting division points preserves a simple cadlag function }
\textbf{\emph{(Insertion of division points leave a simple càdlàg
function unchanged). }}Let $p\geq1$ be arbitrary. Let $0\equiv q_{0}<q_{1}<\cdots<q_{p}\equiv1$
be an arbitrary sequence in $[0,1]$, with an arbitrary subsequence
$0\equiv q_{i(0)}<q_{i(1)}<\cdots<q_{i(\kappa)}\equiv1$. Let $(w_{0},\cdots,w_{\kappa-1})$
be an arbitrary sequence in $S$. Let 
\[
x\equiv\Phi_{smpl}((q_{i(k)})_{k=0,\cdots,\kappa-1},(w_{k})_{k=0,\cdots,\kappa-1}).
\]
Let
\begin{equation}
y=\Phi_{smpl}((q_{j})_{j=0,\cdots,p-1},(x(q_{j}))_{j=0,\cdots,p-1}).\label{eq:temp-569}
\end{equation}
Then $x=y$. 
\end{lem}
\begin{proof}
Let $j=0,\cdots,p-1$ and $t\in[q_{j},q_{j+1})$ be arbitrary. Then
$t\in[q_{j},q_{j+1})\subset[q_{i(k)},q_{i(k-1)})$ for some unique
$k=0,\cdots,\kappa-1$. Hence $y(t)=x(q_{j})=w_{k}=x(t)$. Thus $y=x$
on the dense subset $\bigcup_{j=0}^{p-1}[q_{j},q_{j+1})$ of $domain(y\cap domain(x)$.
Hence, by Lemma \ref{Lem.  If x=00003Dy on dense A then x=00003Dy. If f(x,y)<=00003Dc on A then...-1},
we have $y=x.$
\end{proof}

\section{Skorokhod Space $D[0,1]$ of Càdlàg Functions}

Following Skorokhod, via \cite{Billingsley68}, we proceed to define
a metric on the space $D[0,1]$ of càdlàg functions. This metric is
similar to the supremum metric in $C[0,1]$, except that it allows
a small distortion of the time scale in $[0,1]$ by some continuous
and increasing function $\lambda:[0,1]\rightarrow[0,1]$ with $\lambda(0)=0$
and $\lambda(1)=1$. 

Let $\lambda,\lambda'$ be any such continuous and increasing functions.
We will write, for abbreviation, $\lambda t$ for $\lambda(t)$, for
each $t\in[0,1]$. We will write $\lambda^{-1}$ for the inverse of
$\lambda$, and write $\lambda'\lambda\equiv\lambda'\circ\lambda$
for the composite function. 
\begin{defn}
\label{Def. Skorokhod Metric} \textbf{(Skorokhod metric).} Let $\Lambda$
denote the set of continuous and increasing functions $\lambda:[0,1]\rightarrow[0,1]$
with $\lambda0=0$ and $\lambda1=1$, such that there exists $c>0$
with 
\begin{equation}
|\log\frac{\lambda t-\lambda s}{t-s}|\leq c,\label{eq:temp-244-2}
\end{equation}
or, equivalently,
\[
e^{-c}(t-s)\leq\lambda t-\lambda s\leq e^{c}(t-s),
\]
for each $0\leq s<t\leq1$. We will call $\Lambda$ the set of \index{admissible functions on {[}0,1{]}}
\emph{admissible functions} on $[0,1]$.

Let $x,y\in D[0,1]$ be arbitrary. Let $A_{x,y}$ denote the set consisting
of all pairs $(c,\lambda)\in[0,\infty)\times\Lambda$ such that inequality
\ref{eq:temp-244-2} holds for each $0\leq s<t\leq1$, and such that
\begin{equation}
d(x,y\circ\lambda)\leq c.\label{eq:temp-245}
\end{equation}
Recall that the last inequality means $d(x(t),y(\lambda t))\leq c$
for each $t\in domain(t)\cap\lambda^{-1}domain(y)$.

Let
\[
B_{x,y}\equiv\{c\in[0,\infty):(c,\lambda)\in A_{x,y}\;for\;some\;\lambda\}
\]
Define the metric $d_{D[0,1]}$ on $D[0,1]$ by

\begin{equation}
d_{D[0,1]}(x,y)\equiv\inf B_{x,y}.\label{eq:temp-242}
\end{equation}
We will presently prove that $d_{D[0,1]}$ is well-defined and is
indeed a metric, called the \emph{Skorokhod metric }\index{Skorokhod metric}
on $D[0,1]$. When the interval $[0,1]$ is understood, we write $d_{D}$
for $d_{D[0,1]}$.
\end{defn}
$\square$

Intuitively, the number $c$ bounds both the (i) error in the time
measurement, represented by the distortion $\lambda$, and (ii) the
supremum distance between the functions $x$ and $y$ when allowance
is made for said error. Existence of the infimum in equality \ref{eq:temp-242}
would follow easily from the principle of infinite search. We will
supply such a a constructive proof, in the following Lemmas \ref{Lem. B_xy basics}
through \ref{Lem. Prep Lemma Arzela-Ascoli for D=00005B0,1=00005D}
and Proposition \ref{Lem. Skorokod metric is well defined}. Proposition
\ref{Prop. Skorokhod metric is indeed a metric} will complete the
proof that $d_{D[0,1]}$ is a metric. Then we will prove that the
Skorokhod metric space $(D[0,1],d_{D})$ is complete.

First two elementary lemmas. 
\begin{lem}
\label{Lem. Condition for existence of infimum} \textbf{\emph{(A
condition for existence of infimum or supremum).}} Let $B$ be an
arbitrary nonempty subset of $R$. 

1. Suppose, for each $k\geq0$, there exists $\alpha_{k}\in R$ such
that \emph{(i)} $\alpha_{k}\leq c+2^{-k}$ for each $c\in B$, and
\emph{(ii)} $c\leq\alpha_{k}+2^{-k}$ for some $c\in B$. Then $\inf B$
exists, and $\inf B=\lim_{k\rightarrow\infty}\alpha_{k}$, with $\alpha_{k}-2^{-k}\leq\inf B\leq\alpha_{k}+2^{-k}$
for each $k\geq0$.

2. Suppose, for each $k\geq0$, there exists $\alpha_{k}\in R$ such
that \emph{(iii)} $\alpha_{k}\geq c-2^{-k}$ for each $c\in B$, and
\emph{(iv)} $c\geq\alpha_{k}-2^{-k}$ for some $c\in B$. Then $\sup B$
exists, and $\sup B=\lim_{k\rightarrow\infty}\alpha_{k}$ , with $\alpha_{k}-2^{-k}\leq\sup B\leq\alpha_{k}+2^{-k}$
for each $k\geq0$.
\end{lem}
\begin{proof}
Let $h,k\geq0$ be arbitrary. Then, by Condition (ii), there exists
$c\in B$ such that $c\leq\alpha_{k}+2^{-k}$. At the same time, by
Condition (i), we have $\alpha_{h}\leq c+2^{-h}\leq\alpha_{k}+2^{-k}+2^{-h}$.
Similarly, $\alpha_{k}\leq\alpha_{h}+2^{-h}+2^{-k}.$ Thus $|\alpha_{h}-\alpha_{k}|\leq2^{-h}+2^{-k}$.
We conclude that the limit $\alpha\equiv\lim_{k\rightarrow\infty}\alpha_{k}$
exists. Let $c\in B$ be arbitrary. Letting $k\rightarrow\infty$
in Condition (i), we see that $\alpha\leq c$. Thus $\alpha$ is a
lower bound for the set $B$. Suppose $\beta$ is a second lower bound
for $B$. By condition (ii) there exists $c\in B$ be such that $c\leq\alpha_{k}+2^{-k}$.
Then $\beta\leq c\leq\alpha_{k}+2^{-k}$. Letting $k\rightarrow\infty$,
we obtain $\beta\leq\alpha$. Thus $\alpha$ is the greatest lower
bound of the set $B$. In other words, $\inf B$ exists and is equal
to $\alpha$, as alleged. Assertion 1 is proved. The proof of Assertion
2 is similar.
\end{proof}
\begin{lem}
\label{Lem. log((gt-gs)/(t-s))} \textbf{\emph{(Logarithm of certain
difference quotients).}} Let 
\[
0=\tau_{0}<\cdots<\tau_{p-1}<\tau_{p}\equiv1
\]
be an arbitrary sequence in $[0,1]$. Suppose the function $\lambda\in\Lambda$
is linear on $[\tau_{i},\tau_{i+1}]$ for each $i=0,\cdots,p-1$.
Then 
\begin{equation}
\sup_{0\leq s<t\leq1}|\log\frac{\lambda t-\lambda s}{t-s}|=\alpha\equiv\bigvee_{i=0}^{p-1}|\log\frac{\lambda\tau_{i+1}-\lambda\tau_{i}}{\tau_{i+1}-\tau_{i}}|.\label{eq:temp-305}
\end{equation}
\end{lem}
\begin{proof}
Let $s,t\in A\equiv\bigcup_{i=0}^{p-1}(\tau_{i},\tau_{i+1})$ be arbitrary,
with $s<t$. Then $\tau_{i}<s<\tau_{i+1}$ and $\tau_{j}<t<\tau_{j+1}$
for some $i,j=0,\cdots,p-1$ with $i\leq j$. Hence, in view of the
linearity of $\lambda$ on each of the intervals $[\tau_{i},\tau_{i+1})$
and $[\tau_{j},\tau_{j+1})$, we obtain 
\[
\lambda t-\lambda s=(\lambda t-\lambda\tau_{j})+(\lambda\tau_{j}-\lambda\tau_{j-1})+\cdots+(\lambda\tau_{i+2}-\lambda\tau_{i+1})+(\lambda\tau_{i+1}-\lambda s)
\]
\[
\leq e^{\alpha}(t-\tau_{j})+e^{\alpha}(\tau_{j}-\tau_{j-1})+\cdots+e^{\alpha}(\tau_{j}-\tau_{j-1})+e^{\alpha}(\tau_{i+1}-s)
\]
\[
=e^{\alpha}(t-s).
\]
Similarly $\lambda t-\lambda s\geq e^{-\alpha}(t-s)$. Thus
\[
e^{-\alpha}(t-s)\leq\lambda t-\lambda s\leq e^{\alpha}(t-s),
\]
where $s,t\in A$ with $s<t$ are arbitrary. Since $A$ is dense in
$[0,1]$, the last displayed inequality holds, by continuity, for
each $s,t\in[0,1]$ with $s<t$. Equivalently, 
\begin{equation}
|\log\frac{\lambda t-\lambda s}{t-s}|\leq\alpha\label{eq:temp-317}
\end{equation}
for each $s,t\in[0,1]$ with $s<t$. At the same time, for each $\varepsilon>0$,
there exists  $j=0,\cdots,m$ such that 
\[
|\log\frac{\lambda\tau{}_{j+1}-\lambda\tau{}_{j}}{\tau_{j+1}-\tau_{j}}|>\alpha-\varepsilon.
\]
Thus
\begin{equation}
\alpha<|\log\frac{\lambda u-\lambda v}{u-v}|+\varepsilon\label{eq:temp-318}
\end{equation}
where $u\equiv\tau{}_{j+1}$ and $v\equiv\lambda\tau_{j}$. Since
$\varepsilon>0$ is arbitrary, inequalities \ref{eq:temp-317} and
\ref{eq:temp-318} together imply the desired equality \ref{eq:temp-305},
thanks to Lemma \ref{Lem. Condition for existence of infimum}. 
\end{proof}
Now some metric-like properties of the sets $B_{x,y}$ introduced
in Definition \ref{Def. Skorokhod Metric}.
\begin{lem}
\label{Lem. B_xy basics} \textbf{\emph{(Metric-like properties of
the sets $B_{x,y}$).}} Let $x,y,z\in D[0,1]$ be arbitrary. Then
$B_{x,y}$ is nonempty. Moreover, the following holds.

1. $0\in B_{x,x}$. More generally, if $d(x,y)\leq b$ for some $b\geq0$,
then $b\in B_{x,y}$.

2. \textup{$B_{x,y}=B_{x,y}$.}

3. Let $c\in B_{x,y}$ and $c'\in B_{y,z}$ be arbitrary. Then $c+c'\in B_{x,z}$.
Specifically, suppose $(c,\lambda)\in A_{x,y}$ and $(c',\lambda')\in A_{y,z}$
for some $\lambda,\lambda'\in\Lambda$. Then $(c+c',\lambda,\lambda')\in A_{x,z}$. 
\end{lem}
\begin{proof}
1. Let $\lambda_{0}:[0,1]\rightarrow[0,1]$ be the identity function.
Then, trivially, $\lambda_{0}$ is admissible and $d(x,x\circ\lambda_{0})=0$.
Hence $(0,\lambda_{0})\in A_{x,x}$. Consequently, $0\in B_{x,x}$.
More generally, if $d(x,y)\leq b$ then for some $b\geq0$, then $d(x,y\circ\lambda_{0})=d(x,y)\leq b=0$,
whence $b\in B_{x,y}$.

2. Next consider each $c\in B_{x,y}$. Then there exists $(c,\lambda)\in A_{x,y}$,
satisfying inequalities \ref{eq:temp-244-2} and \ref{eq:temp-245}.
For each $0\leq s<t\leq1$, if we write $u\equiv\lambda^{-1}t$ and
$v\equiv\lambda^{-1}s$, then 
\begin{equation}
|\log\frac{\lambda^{-1}t-\lambda^{-1}s}{t-s}|=|\log\frac{u-v}{\lambda u-\lambda v}|=|\log\frac{\lambda u-\lambda v}{u-v}|\leq c.\label{eq:temp-244-1-1}
\end{equation}
Consider each $t\in domain(y)\cap(\lambda^{-1})^{-1}domain(x)$. Then
$u\equiv\lambda^{-1}t\in domain(x)\cap\lambda^{-1}domain(y)$. Hence
\begin{equation}
d(y(t),x(\lambda^{-1}t))=d(y(\lambda u),x(u))\leq c.\label{eq:temp-245-1-1}
\end{equation}
Thus $(c,\lambda^{-1})\in A_{y,x}$. Consequently $c\in B_{y,x}$.
Since $c\in B_{x,y}$ is arbitrary, we conclude that $B_{x,y}\subset B_{y,x}$,
and, by symmetry, that $B_{x,y}=B_{y,x}$. 

3. Consider arbitrary $c\in B_{x,y}$ and $c'\in B_{y,z}$. Then $(c,\lambda)\in A_{x,y}$
and $(c',\lambda')\in A_{y,z}$ for some $\lambda,\lambda'\in\Lambda$.
The composite function $\lambda'\lambda$ on $[0,1]$ then satisfies
\[
|\log\frac{\lambda'\lambda t-\lambda'\lambda s}{t-s}|=|\log\frac{\lambda'\lambda t-\lambda'\lambda s}{\lambda t-\lambda s}+\log\frac{\lambda t-\lambda s}{t-s}|\leq c+c'.
\]
Let
\[
r\in A\equiv domain(x)\cap\lambda^{-1}domain(y)\cap(\lambda'\lambda)^{-1}domain(z)
\]
be arbitrary. Then
\[
d(x(r),z(\lambda'\lambda r))\leq d(x(r),y(\lambda r))+d(y(\lambda r),z(\lambda'\lambda r))\leq c+c'.
\]
By Proposition \ref{Prop. Points of continuitys of cadlag functions-1},
the set $A$ is dense in $[0,1]$. It therefore follows from Lemma
\ref{Lem.  If x=00003Dy on dense A then x=00003Dy. If f(x,y)<=00003Dc on A then...-1}
that 
\[
d(x,z\circ(\lambda'\lambda))\leq c+c'.
\]
Combining, we see that $(c+c',\lambda'\lambda)\in B_{x,z}$. 
\end{proof}
\begin{defn}
\label{Def notations for Azela-Ascoli} \textbf{(Notations).} We will
use the following notations.

1. Recall, from Definition \ref{Def. Skorokhod Metric}, the set $\Lambda$
of admissible functions on $[0,1]$. Let $\lambda_{0}\in\Lambda$
be the identity function, i.e. $\lambda_{0}t=t$ for each $t\in[0,1]$.
Let $m'\geq m$ be arbitrary. Then $\Lambda_{m,m'}$ will denote the
finite subset of $\Lambda$ consisting of functions $\lambda$ such
that (i) $\lambda Q_{m}\subset Q_{m'}$, and (ii) $\lambda$ is linear
on $[q_{m,i},q_{m,i+1}]$ for each $i=0,\cdots,p_{m}-1$.

2. Let $B$ be an arbitrary compact subset of $(S,d)$, and let $\delta:(0,\infty)\rightarrow(0,\infty)$
be an arbitrary operation. Then $D_{B,\delta}[0,1]$ will denote the
subset of $D[0,1]$ consisting of càdlàg functions $x$ with values
in the compact set $B$ and with $\delta_{cdlg}\equiv\delta$ as a
modulus of càdlàg. 

3. Let $U\equiv\{u_{1},\cdots,u_{M}\}$ be an arbitrary finite subset
of $(S,d)$. Then $D_{simple,m,U}[0,1]$ will denote the finite subset
of $D[0,1]$ consisting of simple càdlàg functions with values in
$U$ and with $q_{m}$ as a sequence of division points.

4. Let $\delta_{\log@1}$ be a modulus of continuity at 1 of the natural
logarithm function $\log$. Specifically, let $\delta_{\log@1}(\varepsilon)\equiv1-e^{-\varepsilon}$
for each $\varepsilon>0$. Note that $0<\delta_{\log@1}<1$. $\square$

The next lemma proves that $\inf B_{x,y}$ exists for arbitrary simple
càdlàg functions $x,y$. With only finite searches, we need to give
a constructive proof, along with a method of approximating the alleged
infimum to arbitrary precision. 
\end{defn}
\begin{lem}
\label{Lem. d_D(x,y) well defined for certain simple cadlags} \textbf{\emph{($d_{D[0,1]}$
is well defined on $D_{simple,m,U}[0,1]$).}} Let $M,m\geq1$ be arbitrary.
Let $U\equiv\{u_{1},\cdots,u_{M}\}$ be an arbitrary finite subset
of $(S,d)$. Let $x,y\in D_{simple,m,U}[0,1]$ be arbitrary. Then
$d_{D}(x,y)\equiv\inf B_{x,y}$ exists. 

Specifically, take any $b\geq\bigvee_{i,j=0}^{M}d(u_{i},u_{j})$.
Let $(x_{i})_{i=0,\cdots,p(m)-1}$ be an arbitrary sequence in $U$,
and consider the simple càdlàg function
\[
x\equiv\Phi_{smpl}((q_{m,i})_{i=0,\cdots,p(m)-1},(x_{i})_{i=0,\cdots,p(m)-1}).
\]
Let $k\geq0$ be arbitrary. Take $m'\geq m$ so large that
\begin{equation}
2^{-m'}\leq2^{-m-2}e^{-b}\delta_{\log@1}(2^{-k})<2^{-m-2}e^{-b}.\label{eq:temp-325}
\end{equation}
For each $\nu\in\Lambda_{m,m'}$, define 

\begin{equation}
\beta_{\nu}\equiv\bigvee_{i=0}^{p(m)-1}(|\log\frac{\nu q_{m,i+1}-\nu q_{m,i}}{q_{m,i+1}-q_{m,i}}|\vee d(y(\nu q_{m,i}),x(q_{m,i}))\vee d(y(q_{m,i}),x(\nu^{-1}q_{m,i}))).\label{eq:temp-331}
\end{equation}
Then there exists $\nu\in\Lambda_{m,m'}$ with $(\beta_{\nu},\nu)\in A_{x,y}$
such that $d_{D}(x,y)\in[\beta_{\nu}-2^{-k+1},\beta_{\nu}]$. 

Thus $\beta_{\nu}$ is a $2^{-k+1}$-approximation of $d_{D}(x,y)$,
and $(d_{D}(x,y)+2^{-k+1},\nu)\in A_{x,y}$.
\end{lem}
\begin{proof}
1. Let $m,k,m'$ be as given. For abbreviation, write $\varepsilon\equiv2^{-k}$,
$p\equiv p_{m}\equiv2^{m}$, and 
\[
\tau_{i}\equiv\xi_{i}\equiv q_{m,i}\equiv i2^{-m}\in Q_{m}
\]
for each $i=0,\cdots,p$. Similarly, write $n\equiv p_{m'}\equiv2^{m'}$,
$\Delta\equiv2^{-m'}$, and 
\[
\eta_{j}\equiv q_{m',j}\equiv j2^{-m'}\in Q_{m'}
\]
for each $j=0,\cdots,n$. Then, by hypothesis, $x\equiv\Phi_{smpl}((\tau_{i})_{i=0,\cdots,p},(x_{i})_{i=0,\cdots,p-1})$.
Similarly, 
\[
y\equiv\Phi_{smpl}((\xi_{i})_{i=0,\cdots,p},(y_{i})_{i=0,\cdots,p-1})
\]
for some sequence $(y_{i})$ in $U$. By the definition of simple
càdlàg functions, we have $x=x_{i}$ and $y=y_{i}$ on $[\tau_{i},\tau_{i+1})\equiv[\xi_{i},\xi_{i+1})$,
for each $i=0,\cdots,p-1$. By hypothesis, 
\[
b\geq\bigvee_{i,j=0}^{M}d(u_{i},u_{j})\geq\bigvee_{i,j=0}^{p-1}d(x_{i},y_{j})\geq d(x(t),y(t))
\]
for each $t\in domain(x)\cap domain(y)$. Hence $b\in B_{x,y}$ by
Lemma \ref{Lem. B_xy basics}. Define 
\[
\alpha_{k}\equiv\bigwedge_{\nu\in\Lambda(m,m')}\beta_{\nu}.
\]
We will prove that (i) $\alpha_{k}\leq c+2^{-k}$ for each $c\in B_{x,y}$,
and (ii) $c\leq\alpha_{k}+2^{-k}$ for some $c\in B_{x,y}$. It will
then follow from Lemma \ref{Lem. Condition for existence of infimum}
that $\inf B_{x,y}$ exists, and that $\inf B_{x,y}=\lim_{k\rightarrow\infty}\alpha_{k}$. 

2. We will first prove Condition (i). To that end, let $c\in B_{x,y}$
be arbitrary. Since $b\in B_{x,y}$, there is no loss of generality
in assuming that $c\leq b$. Since $c\in B_{x,y}$ by assumption,
there exists $\lambda\in\Lambda$ such that $(c,\lambda)\in A_{x,y}$.
In other words, 
\begin{equation}
|\log\frac{\lambda t-\lambda s}{t-s}|\leq c\label{eq:temp-244-2-1}
\end{equation}
for each $0\leq s<t\leq1$, and 
\begin{equation}
d(x,y\circ\lambda)\leq c.\label{eq:temp-245-1}
\end{equation}
Consider each $i=1,\cdots,p-1$. There exists $j_{i}=1,\cdots,n-1$
such that 
\begin{equation}
\eta_{j(i)-1}<\lambda\tau_{i}<\eta_{j(i)+1}.\label{eq:temp-385}
\end{equation}
Either (i') 
\[
d(y(\eta_{j(i)-1}),x(\tau_{i}))<d(y(\eta_{j(i)}),x(\tau_{i}))+\varepsilon,
\]
or (ii') 
\[
d(y(\eta_{j(i)}),x(\tau_{i}))<d(y(\eta_{j(i)-1}),x(\tau_{i}))+\varepsilon.
\]
In case (i'), define $\zeta_{i}\equiv\eta_{j(i)-1}$. In case (ii'),
define $\zeta_{i}\equiv\eta_{j(i)}$. Then, in both Cases (i') and
(ii'), we have
\begin{equation}
\zeta_{i}-\Delta<\lambda\tau_{i}<\zeta_{i}+2\Delta,\label{eq:temp-389}
\end{equation}
and 
\begin{equation}
d(y(\zeta_{i}),x(\tau_{i}))\leq d(y(\eta_{j(i)-1}),x(\tau_{i}))\wedge d(y(\eta_{j(i)}),x(\tau_{i}))+\varepsilon.\label{eq:temp-327}
\end{equation}
At the same time, in view of inequality \ref{eq:temp-385}, there
exists a point 
\[
s\in(\lambda\tau_{i},\lambda\tau_{i+1})\cap((\eta_{j(i)-1},\eta_{j(i)})\cup(\eta_{j(i)},\eta_{j(i)+1})).
\]
Then $t\equiv\lambda^{-1}s\in(\tau_{i},\tau_{i+1})$, whence $x(\tau_{i})=x(t)$.
Moreover, either $s\in(\eta_{j(i)-1},\eta_{j(i)})$, in which case
$y(\eta_{j(i)-1})=y(s)$, or $s\in(\eta_{j(i)},\eta_{j(i)+1})$, in
which case $y(\eta_{j(i)})=y(s)$. In either case, inequality \ref{eq:temp-327}
yields 
\begin{equation}
d(y(\zeta_{i}),x(\tau_{i}))\leq d(y(s),x(\tau_{i}))+\varepsilon=d(y(\lambda t),x(t))+\varepsilon\leq c+\varepsilon,\label{eq:temp-328}
\end{equation}
where the last inequality is from inequality \ref{eq:temp-245-1}.
Now let $\zeta_{0}\equiv0$ and $\zeta_{p}\equiv1$. Then, for each
$i=0,\cdots,p-1$, inequality \ref{eq:temp-389} implies that
\begin{equation}
\zeta_{i}-\Delta<\lambda\tau_{i}<\zeta_{i}+2\Delta,\label{eq:temp-389-1}
\end{equation}
\[
\zeta_{i+1}-\zeta_{i}>(\lambda\tau_{i+1}-2\Delta)-(\lambda\tau_{i}+\Delta)>\lambda\tau_{i+1}-\lambda\tau_{i}-4\Delta
\]
\[
\geq e^{-c}(\tau_{i+1}-\tau_{i})-4\Delta\geq e^{-b}2^{-m}-2^{-m'+2}\geq0,
\]
where the last inequality is from inequality \ref{eq:temp-325}. Thus
$(\zeta_{i})_{i=0,\cdots,p}$ is an increasing sequence in $Q_{m'}$.
As such, it determines a unique function $\mu\in\Lambda_{m,m'}$ which
is linear on $[\tau_{i},\tau_{i+1}]$ for each $i=0,\cdots,p-1$,
with $\mu\tau_{i}\equiv\zeta_{i}$ for each $i=0,\cdots,p$. By Lemma
\ref{Lem. log((gt-gs)/(t-s))}, we then have 
\[
\sup_{0\leq s<t\leq1}|\log\frac{\mu t-\mu s}{t-s}|=\bigvee_{i=0}^{p-1}|\log\frac{\mu\tau_{i+1}-\mu\tau_{i}}{\tau_{i+1}-\tau_{i}}|
\]
\[
=\bigvee_{i=0}^{p-1}|\log(\frac{\lambda\tau_{i+1}-\lambda\tau_{i}}{\tau_{i+1}-\tau_{i}}\cdot\frac{\mu\tau_{i+1}-\mu\tau_{i}}{\lambda\tau_{i+1}-\lambda\tau_{i}})|
\]
\[
\leq\bigvee_{i=0}^{p-1}|\log\frac{\lambda\tau_{i+1}-\lambda\tau_{i}}{\tau_{i+1}-\tau_{i}}|+\bigvee_{i=0}^{p-1}|\log(\frac{\zeta_{i+1}-\zeta_{i}}{\lambda\tau_{i+1}-\lambda\tau_{i}})|
\]
\[
\leq c+\bigvee_{i=0}^{p-1}|\log(1+\frac{(\zeta_{i+1}-\lambda\tau_{i+1})-(\zeta_{i}-\lambda\tau_{i})}{\lambda\tau_{i+1}-\lambda\tau_{i}})|
\]
\begin{equation}
\equiv c+\bigvee_{i=0}^{p-1}|\log(1+a_{i})|,\label{eq:temp-250}
\end{equation}
where, for each $i=0,\cdots,p-1$, 
\[
a_{i}\equiv\frac{(\zeta_{i+1}-\lambda\tau_{i+1})-(\zeta_{i}-\lambda\tau_{i})}{\lambda\tau_{i+1}-\lambda\tau_{i}},
\]
with 
\[
|a_{i}|\leq\frac{2^{-m'+1}+2^{-m'+1}}{\exp(-c)(\tau_{i+1}-\tau_{i})}=e^{c}2^{-m'+2+m}\leq e^{b}2^{-m'+2+m}<\delta_{\log@1}(\varepsilon),
\]
thanks again to inequality \ref{eq:temp-325}. Hence $\bigvee_{i=0}^{p-1}|\log(1+a_{i})|<\varepsilon$.
Inequality \ref{eq:temp-250} therefore yields 
\begin{equation}
\sup_{0\leq s<t\leq1}|\log\frac{\mu t-\mu s}{t-s}|=\bigvee_{i=0}^{p-1}|\log\frac{\mu\tau_{i+1}-\mu\tau_{i}}{\tau_{i+1}-\tau_{i}}|<c+\varepsilon.\label{eq:temp-306}
\end{equation}
Separately, by the definition of $\mu$, inequality \ref{eq:temp-328}
implies
\begin{equation}
\bigvee_{i=0}^{p-1}d(y(\mu\tau_{i}),x(\tau_{i}))\equiv\bigvee_{i=0}^{p-1}d(y(\zeta_{i}),x(\tau_{i}))\leq c+\varepsilon.\label{eq:temp-307}
\end{equation}
We will proceed to verify also that 
\begin{equation}
\bigvee_{j=0}^{p-1}d(y(\xi_{j}),x(\mu^{-1}\xi_{j}))\leq c+\varepsilon,\label{eq:temp-307-2}
\end{equation}
where we recall that $\xi_{i}\equiv\tau_{i}\equiv q_{m,i}\equiv i2^{-m}\in Q_{m}$,
for each $i=0,\cdots,p$. To that end, consider each $j=1,\cdots,p-1$.
Then, for each $i=0,\cdots,p-1$, we have $\xi_{j},\mu\tau_{i}\equiv\zeta_{i}\in Q_{m'}$.
Hence 
\begin{equation}
\xi_{j}\in[\mu\tau_{i},\mu\tau_{i+1})\equiv[\zeta_{i},\zeta_{i+1})\label{eq:temp-247}
\end{equation}
for some $i=0,\cdots,p-1$. Either (i'') $\xi_{j}=\zeta_{i}$, or
(ii'') $\xi_{j}\geq\zeta_{i}+\Delta$. First consider Case (i'').
Then 
\[
d(y(\xi_{j}),x(\mu^{-1}\xi_{j}))=d(y(\zeta_{i}),x(\mu^{-1}\zeta_{i}))=d(y(\zeta_{i}),x(\tau_{i}))\leq c+\varepsilon
\]
by inequality \ref{eq:temp-328}. Now consider Case (ii''). Then,
in view of inequality \ref{eq:temp-389}, we obtain 
\[
\lambda\tau_{i}<\zeta_{i}+2\Delta\leq\xi_{j}+\Delta.
\]
At the same time, relation \ref{eq:temp-247} implies that $\xi_{j}<\zeta_{i+1}$.
Hence
\[
\xi_{j}\leq\zeta_{i+1}-\Delta<\lambda\tau_{i+1},
\]
where the inequality on the right-hand side is from inequality \ref{eq:temp-389}
applied to $i+1$. Combining the last two displayed inequalities,
we see that there exists a point 
\[
s\in[\lambda\tau_{i}\vee\xi_{j},\lambda\tau_{i+1}\wedge(\xi_{j}+\Delta)).
\]
It follows that $s\in[\xi_{j},\xi_{j}+\Delta)\subset[\zeta_{i},\zeta_{i+1})$,
and that $r\equiv\lambda^{-1}s\in[\tau_{i},\tau_{i+1})$. Therefore
$y(s)=y(\xi_{j})$ and $x(r)=x(\tau_{i})$ by the definition of the
simple càdlàg functions $x$ and $y$. Moreover, 
\[
\mu^{-1}\xi_{j}\in[\mu^{-1}\zeta_{i},\mu^{-1}\zeta_{i+1})\equiv[\tau_{i},\tau_{i+1}),
\]
whence $x(\mu^{-1}\xi_{j})=x(\tau_{i})$. Combining, we obtain
\[
d(y(\xi_{j}),x(\mu^{-1}\xi_{j}))=d(y(s),x(\tau_{i}))=d(y(s),x(r))\equiv d(y(\lambda r),x(r))\leq c+\varepsilon,
\]
where $j=1,\cdots,p-1$ is arbitrary. Thus we have verified inequality
\ref{eq:temp-307-2}. 

Inequalities \ref{eq:temp-306}, \ref{eq:temp-307}, and  \ref{eq:temp-307-2}
together imply that
\[
\beta_{\mu}\equiv\bigvee_{i=0}^{p-1}(|\log\frac{\mu\tau_{i+1}-\mu\tau_{i}}{\tau_{i+1}-\tau_{i}}|\vee d(y(\mu\tau_{i}),x(\tau_{i}))\vee d(y(\xi_{i}),x(\mu^{-1}\xi_{i})))\leq c+\varepsilon.
\]
Hence $\alpha_{k}\equiv\bigwedge_{\nu\in\Lambda(m,m')}\beta_{\nu}\leq c+\varepsilon$,
where $c\in B_{x,y}$ is arbitrary. The desired Condition (i) is established.

3. Since $\Lambda_{m,m'}$ is a finite set, there exists $\nu\in\Lambda_{m,m'}$
such that $\beta_{\nu}<\alpha_{k}+\varepsilon.$ Thus 
\[
\bigvee_{i=0}^{p-1}(|\log\frac{\nu\tau_{i+1}-\nu\tau_{i}}{\tau_{i+1}-\tau_{i}}|\vee d(y(\nu\tau_{i}),x(\tau_{i}))\vee d(y(\xi_{i}),x(\nu^{-1}\xi_{i})))\equiv\beta_{\nu}<\alpha_{k}+\varepsilon.
\]
Then, by Lemma \ref{Lem. log((gt-gs)/(t-s))}, we have 
\begin{equation}
\sup_{0\leq s<t\leq1}|\log\frac{\nu t-\nu s}{t-s}|=\bigvee_{i=0}^{p-1}|\log\frac{\nu\tau_{i+1}-\nu\tau_{i}}{\tau_{i+1}-\tau_{i}}|\leq\beta_{\nu}.\label{eq:temp-388}
\end{equation}
We will show that $d(y\circ\nu,x)\leq\beta_{\nu}$. To that end, consider
each point $t$ in the dense subset $A\equiv Q_{m'}^{c}\cup\nu^{-1}Q_{m'}^{c}$
of $[0,1]$. Then 
\begin{equation}
t\in[\tau_{i},\tau_{i+1})\cap\nu^{-1}[\xi_{j},\xi_{j+1})\label{eq:temp-248}
\end{equation}
for some $i,j=0,\cdots,p-1$. It follows that either (i''') $\nu^{-1}\xi_{j}\leq\tau_{i}$,
or (ii''') $\tau_{i}<\nu^{-1}\xi_{j}$. In Case (i'''), we have $\xi_{j}\leq\nu\tau_{i}\leq\nu t<\xi_{j+1}$,
whence 
\[
d(y(\nu t),x(t))=d(y(\nu\tau_{i}),x(\tau_{i}))\leq\beta_{\nu},
\]
where we used the definition of simple càdlàg functions. In Case (ii'''),
relation \ref{eq:temp-248} implies that $\tau_{i}<\nu^{-1}\xi_{j}\leq t<\tau_{i+1}$,
whence 
\[
d(y(\nu t),x(t))=d(y(\xi_{j}),x(\nu^{-1}\xi_{j}))\leq\beta_{\nu},
\]
where we used, once more, the definition of simple càdlàg functions.
Since $t$ is an arbitrary member of the dense subset $A$ of $[0,1]$,
we conclude that 
\begin{equation}
d(y\circ\nu,x)\leq\beta_{\nu}\label{eq:temp-386}
\end{equation}
Combining with inequality \ref{eq:temp-388}, we conclude that $(\beta_{\nu},\nu)\in A_{x,y}$.
Since $\beta_{\nu}<\alpha_{k}+2^{-k}$, it follows that $c\equiv\alpha_{k}+2^{-k}\in B_{x,y}$,
thus proving Condition (ii).

4. Since $k\geq0$ is arbitrary, it follows from Lemma \ref{Lem. Condition for existence of infimum}
that both $\lim_{k\rightarrow\infty}\alpha_{k}$ and $d_{D}(x,y)\equiv\inf B_{x,y}$
exist, and are equal to each other.

5. Finally, let $c\in B_{x,y}$ be arbitrary. Then, by Step 2, there
exists $\mu\in\Lambda_{m,m'}$ such that $\beta_{\mu}\leq c+\varepsilon\equiv c+2^{-k}$.
By Step 3, there exists $\nu\in\Lambda_{m,m'}$ with $\beta_{\nu}<\alpha_{k}+2^{-k}\leq\beta_{\mu}+2^{-k}$
such that $(\beta_{\nu},\nu)\in A_{x,y}.$ Combining, $\beta_{\nu}<c+2^{-k+1}$.
Since $c\in B_{x,y}$ is arbitrary, it follows that $\beta_{\nu}\leq\inf B_{x,y}+2^{-k+1}$.
In other words, $\beta_{\nu}\leq d_{D}(x,y)+2^{-k+1}$. Summing up,
we have $\nu\in\Lambda_{m,m'}$ with $(\beta_{\nu},\nu)\in A_{x,y}$,
such that $d_{D}(x,y)\in[\beta_{\nu}-2^{-k+1},\beta_{\nu}]$. The
Lemma is proved. 
\end{proof}
The next lemma prepares for a subsequent generalization of Arzela-Ascoli
Theorem to càdlàg functions. 
\begin{lem}
\label{Lem. Prep Lemma Arzela-Ascoli for D=00005B0,1=00005D} \textbf{\emph{(Preparatory
lemma for Arzela-Ascoli Theorem for $D[0,1]$). }}Let $B$ be an arbitrary
compact subset of $(S,d)$. Let $k\geq0$ be arbitrary. Let $U\equiv\{v_{1},\cdots,v_{M}\}$
be a $2^{-k-1}$-approximation of of $B$. Let $x\in D[0,1]$ be arbitrary
with values in the compact set $B$ and with a modulus of \textup{\emph{càdlàg}}
$\delta_{cdlg}$. Let $m\geq1$ be so large that
\[
m\equiv m(k,\delta_{cdlg})\equiv[0\vee(1-\log_{2}(\delta_{\log@1}(2^{-k})\delta_{cdlg}(2^{-k-1})))]_{1}
\]
Then there exist an increasing sequence $(\eta_{i})_{i=0,\cdots,n-1}$
in $Q_{m}$, and a sequence $(u_{i})_{i=0,\cdots,n-1}$ in $U$, such
that
\[
2^{-k}\in B_{x,\overline{x}},
\]
\emph{where 
\[
\overline{x}\equiv\Phi_{smpl}((\eta_{i})_{i=0,\cdots,n-1},(u_{i})_{i=0,\cdots,n-1})\in D_{simple,m,U}[0,1].
\]
}
\end{lem}
\begin{proof}
Let $k,m$ be as given. Then 
\begin{equation}
2^{-m}<2^{-1}\delta_{\log@1}(2^{-k})\delta_{cdlg}(2^{-k-1}).\label{eq:temp-243}
\end{equation}
For abbreviation, write $\varepsilon\equiv2^{-k}$, $p\equiv p_{m}\equiv2^{m}$,
and $(q_{0},\cdots,q_{p})\equiv(0,2^{-m},2^{-m+1},\cdots,1)$. By
Definition \ref{Def. Cadlag functions, D=00005B0,1=00005D}, there
exists a sequence of $\frac{\varepsilon}{2}$-division points $(\tau_{i})_{i=0,\cdots,n}$
of the càdlàg function  $x$, with separation at least $\delta_{cdlg}(\frac{\varepsilon}{2})$.
Thus 
\begin{equation}
\bigwedge_{i=0}^{n-1}(\tau_{i+1}-\tau_{i})\geq\delta_{cdlg}(\frac{\varepsilon}{2}).\label{eq:temp-324}
\end{equation}
Define $\eta_{0}\equiv0$ and $\eta_{n}\equiv1$. Consider each $i=1,\cdots,n-1$.
Then there exists $j_{i}=1,\cdots,p-1$ such that $\tau_{i}\in(q_{j(i)-1},q_{j(i)+1})$.
Define $\eta_{i}\equiv q_{j(i)}\in Q_{m}$. Then 
\[
|\eta_{i}-\tau_{i}|<2^{-m}<2^{-1}\delta_{\log@1}(2^{-k})\delta_{cdlg}(2^{-k-1})<2^{-1}\delta_{cdlg}(2^{-k-1})\equiv\frac{1}{2}\delta_{cdlg}(\frac{\varepsilon}{2}),
\]
In view of inequality \ref{eq:temp-324}, it follows that 
\[
{\normalcolor \bigwedge_{i=0}^{n-1}(\eta_{i+1}-\eta_{i})>\bigwedge_{i=0}^{n-1}(\tau_{i+1}-\tau_{i})-}\delta_{cdlg}(\frac{\varepsilon}{2})\geq0.
\]
Thus $(\eta_{i})_{i=0,\cdots,n}$ is an increasing sequence in $[0,1]$.
Therefore we can define the increasing function $\nu\in\Lambda$ by
(i) $\nu\tau_{i}=\eta_{i}$ for each $i=0,\cdots,n$, and (ii) $\nu$
is linear on $[\tau_{i},\tau_{i+1}]$ for each $i=1,\cdots,n-1$. 

By Lemma \ref{Lem. log((gt-gs)/(t-s))}, 
\[
\sup_{0\leq s<t\leq1}|\log\frac{\nu t-\nu s}{t-s}|=\bigvee_{j=0}^{n-1}|\log\frac{\nu\tau_{i+1}-\nu\tau_{i}}{\tau_{i+1}-\tau_{i}}|
\]
\begin{equation}
=\bigvee_{i=0}^{n-1}|\log\frac{\eta_{i+1}-\eta_{i}}{\tau_{i+1}-\tau_{i}}|=\bigvee_{i=0}^{n-1}|\log(1+a_{i})|\label{eq:temp-153}
\end{equation}
where, for each $i=0,\cdots,n-1$, 
\[
a_{i}\equiv\frac{(\eta_{i+1}-\tau_{i+1})-(\eta_{i}-\tau_{i})}{\tau_{i+1}-\tau_{i}},
\]
with
\[
|a_{i}|\leq\frac{2^{-m+1}}{|\tau_{i+1}-\tau_{i}|}<2^{-m+1}\delta_{cdlg}(\frac{\varepsilon}{2})^{-1}<\delta_{\log@1}(\varepsilon),
\]
where the inequality is from the hypothesis. Hence $|\log(1+a_{i})|\leq\varepsilon$
for each $i=0,\cdots,n-1$ by the definition of the definition of
the modulus of continuity $\delta_{\log@1}$ in Definition \ref{Def notations for Azela-Ascoli}.
Inequality \ref{eq:temp-153} therefore yields
\begin{equation}
\sup_{0\leq s<t\leq1}|\log\frac{\nu t-\nu s}{t-s}|\leq\varepsilon.\label{eq:temp-3}
\end{equation}

Next, by hypothesis, the càdlàg function $x$ has values in the compact
set $B$, and that $U\equiv\{v_{1},\cdots,v_{M}\}$ is an $2^{-k-1}$-approximation
of $B$. Hence, for each $i=0,\cdots,n-1$, there exists $u_{i}\in U$
such that $d(u_{i},x(\tau_{i}))<2^{-k-1}\equiv\frac{\varepsilon}{2}$.
Also by hypothesis, we have \emph{
\[
\overline{x}\equiv\Phi_{smpl}((\eta_{i})_{i=0,\cdots,n-1},(u_{i})_{i=0,\cdots,n-1})\in D_{simple,m,U}[0,1].
\]
}We will prove that $(\varepsilon,\nu)\in A_{x,\overline{x}}$, and
therefore that\emph{ $\varepsilon\in B_{x,\overline{x}}$. }To that
end, consider each $i=0,\cdots,n-1$ and each $t\in[\tau_{i},\tau_{i+1})$.
Then $\nu t\in[\eta_{i},\eta_{i+1}).$ Hence
\[
d(x(t),\overline{x}(\nu t))
\]
\[
\leq d(x(t),x(\tau_{i}))+d(x(\tau_{i}),u_{i})+d(u_{i},\overline{x}(\eta_{i}))+d(\overline{x}(\eta_{i}),\overline{x}(\nu t))<\frac{\varepsilon}{2}+\frac{\varepsilon}{2}+0+0=\varepsilon.
\]
Since the set $A\equiv\bigcup_{i=0}^{m}[\tau_{i},\tau_{i+1})$ is
dense in $[0,1]$, Lemma \ref{Lem.  If x=00003Dy on dense A then x=00003Dy. If f(x,y)<=00003Dc on A then...-1}
implies that
\[
d(x,\overline{x}\circ\nu)\leq\varepsilon.
\]
Combining with inequality \ref{eq:temp-3}, we see that $(\varepsilon,\nu)\in A_{x,\overline{x}}$.
It follows that\emph{ $\varepsilon\in B_{x,\overline{x}}$. }The lemma
is proved.
\end{proof}
The next lemma proves the existence of $\inf B_{x,y}$ for each $x,y\in D[0,1]$.
The next proposition then verifies that the Skorokhod metric is indeed
a metric.
\begin{lem}
\label{Lem. Skorokod metric is well defined} \textbf{\emph{(Skorokhod
metric is well defined).}} Let $x,y\in D[0,1]$ be arbitrary. Then
the infimum $d_{D[0,1]}(x,y)\equiv\inf B_{x,y}$ in Definition \ref{Def. Skorokhod Metric}
exists.
\end{lem}
\begin{proof}
Let $\delta\equiv\delta_{cdlg}$ be a common modulus of càdlàg of
$x$ and $y$. Let $k\geq0$ be arbitrary and write $\varepsilon\equiv2^{-k}$.
Take $m\geq1$ large enough to satisfy inequality \ref{eq:temp-243}
in Lemma \ref{Lem. Prep Lemma Arzela-Ascoli for D=00005B0,1=00005D}.
By Lemma \ref{Prop. Points of continuitys of cadlag functions-1},
there exists $b\geq0$ such that $d(x,x_{\circ})\vee d(y,x_{\circ})\leq b$.
Since $(S,d)$ is locally compact, the bounded set $(d(\cdot,x_{\circ})\leq b)$
is contained in some compact subset $B$ of $S$. Hence $x,y\in D_{B,\delta}[0,1]$.
Let $U\equiv\{v_{1},\cdots,v_{M}\}$ be an $\frac{\varepsilon}{2}$-approximation
of of $B$. Then, according to Lemma \ref{Lem. Prep Lemma Arzela-Ascoli for D=00005B0,1=00005D},
there exists $\overline{x},\bar{y}\in D_{simple,m,U}[0,1]$ such that
$\varepsilon\in B_{x,\overline{x}}$ and $\varepsilon\in B_{y,\bar{y}}$.
At the same time, by Lemma \ref{Lem. d_D(x,y) well defined for certain simple cadlags},
the infimum $\alpha\equiv\inf B_{\overline{x},\bar{y}}$ exists.

Now let $c\in B_{x,y}$ be arbitrary. Then $c+\varepsilon+\varepsilon\in B_{\overline{x},\bar{y}}$
by Lemma \ref{Lem. B_xy basics}. Hence $\alpha\leq c+2\varepsilon$.
Conversely, take any $\overline{c}\in B_{\overline{x},\bar{y}}$ such
that $\overline{c}<\inf B_{\overline{x},\bar{y}}+\varepsilon\equiv\alpha+\varepsilon$.
Then $c\equiv\overline{c}+\varepsilon+\varepsilon\in B_{x,y}$ by
Lemma \ref{Lem. B_xy basics}. Hence $c\equiv\overline{c}+2\varepsilon<\alpha+3\varepsilon$.
Since $\varepsilon\equiv2^{-k}>0$ is arbitrarily small, \ref{Lem. Condition for existence of infimum}
implies that $\inf B_{x,y}$ exists.
\end{proof}
\begin{prop}
\label{Prop. Skorokhod metric is indeed a metric} \textbf{\emph{(Skorokhod
metric is indeed a metric).}} $(D[0,1],d_{D})$ is a metric space.
Moreover, if $x,y\in D[0,1]$ are such that $d(x,y)\leq c$ on a dense
subset of $domain(x)\cap domain(y)$, then $d_{D}(x,y)\leq c$.
\end{prop}
\begin{proof}
Let $x,y,z\in D[0,1]$ be arbitrary. Let $\delta_{cdlg}$ be a common
modulus of càdlàg of $x,y$ and $z$. By Lemma \ref{Lem. Skorokod metric is well defined},
$d_{D}(x,y)\equiv\inf B_{x,y}$ exists. By Lemma \ref{Lem. B_xy basics},
we have $0\in B_{x,x}$ and $B_{x,y}=B_{x,y}$. It follows immediately
that $d_{D}(x,x)=0$ and $d_{D}(x,y)=d_{D}(y,x)$. Now let $\varepsilon>0$
be arbitrary. By the definition of infimums, there exist $c\in B_{x,y}$
and $c'\in B_{y,z}$ such that $c<\inf B_{x,y}+\varepsilon\equiv d_{D}(x,y)+\varepsilon$
and $c'<\inf B_{y,z}+\varepsilon\equiv d_{D}(y,z)+\varepsilon$. Hence,
again by Lemma \ref{Lem. B_xy basics}, we have 
\[
d_{D}(x,z)\equiv\inf B_{x,z}\leq c+c'<d_{D}(x,y)+\varepsilon+d_{D}(y,z)+\varepsilon.
\]
Since $\varepsilon>0$ is arbitrary, it follows that $d_{D}(x,z)\leq d_{D}(x,y)+d_{D}(y,z)$.

It remains to prove that if $d_{D}(x,y)=0$ then $x=y$. To that end,
suppose $d_{D}(x,y)=0$. Let $\varepsilon>0$ be arbitrary. Let $(\tau_{i})_{i=1,\cdots,p}$
and $(\eta_{j})_{j=1,\cdots,n}$ be sequences of $\varepsilon$-division
points of $x,y$ respectively. Let $m\geq\varepsilon^{-1/2}$ be arbitrary.
Consider each $k\geq m\vee2$. Then $\varepsilon_{k}\equiv k^{-2}\leq m^{-2}\leq\varepsilon$.
Moreover, since $d_{D}(x,y)=0<\varepsilon_{k}$, we have $\varepsilon_{k}\in B_{x,y}$.
Therefore there exists, by Definition \ref{Def. Skorokhod Metric},
some $\lambda_{k}\in\Lambda$ such that
\begin{equation}
|\log\frac{\lambda_{k}r-\lambda_{k}s}{r-s}|\leq\varepsilon_{k}\label{eq:temp-317-1}
\end{equation}
for each $r,s\in[0,1]$ with $s\leq r$, and such that 
\begin{equation}
d(x(t),y(\lambda_{k}t))\leq\varepsilon_{k}\label{eq:temp-318-1}
\end{equation}
for each $t\in domain(x)\cap\lambda_{k}^{-1}domain(y)$. Then inequality
\ref{eq:temp-317-1} implies that
\begin{equation}
e^{-\varepsilon(k)}r\leq\lambda_{k}r\leq e^{\varepsilon(k)}r,\label{eq:temp-320}
\end{equation}
for each $r\in[0,1]$. Define
\[
C_{k}\equiv(\bigcup_{j=0}^{n}[e^{-\varepsilon(k)}\eta_{j},e^{\varepsilon(k)}\eta_{j}])^{c}.
\]
where the superscript $c$ signifies the measure-theoretic complement
of a Lebesgue $\mathrm{measurable}$ set. in {[}0,1{]}. Let $\mu$
denote the Lebesgue measure on $[0,1]$. Then
\[
\mu(C_{k})\geq1-\sum_{i=0}^{n}(e^{\varepsilon(k)}\eta_{j}-e^{-\varepsilon(k)}\eta_{j})
\]
\[
\geq1-\sum_{j=0}^{n}(e^{\varepsilon(k)}-e^{-\varepsilon(k)})\geq1-(n+1)(e^{2\varepsilon(k)}-1)\geq1-2(n+1)e\varepsilon_{k}\equiv1-2(n+1)ek^{-2},
\]
where $k\geq m$ is arbitrary, and where we used the elementary inequality
$e^{r}-1\leq er$ for each $r\in(0,1)$. Define 
\[
C\equiv\bigcup_{h=m}^{\infty}\bigcap_{k=h+1}^{\infty}C_{k}.
\]
Then 
\[
\mu(C^{c})\leq\sum_{k=h+1}^{\infty}\mu(C_{k}^{c})\leq\sum_{k=h+1}^{\infty}2(n+1)ek^{-2}=2(n+1)eh^{-1}
\]
for each $h\geq m$. Hence $\mu(C^{c})=0$ and $C$ is a full subset
of $[0,1]$. Consequently 
\[
A\equiv C\cap domain(x)\cap domain(y)\cap\bigcap_{k=m}^{\infty}\lambda_{k}^{-1}domain(y)
\]
is a full subset of $[0,1]$ and, as such, is dense in $[0,1]$. Now
let $t\in A$ be arbitrary. Let $h\geq m$ be arbitrary. Then there
exists $k\geq h$ such that $t\in C_{k}$. Hence there exists $j=0,\cdots,n-1$
such that $t\in(e^{\varepsilon(k)}\eta_{j},e^{-\varepsilon(k)}\eta_{j+1})$.
It then follows from inequality \ref{eq:temp-320} that
\[
\eta_{j}<e^{-\varepsilon(k)}t\leq\lambda_{k}t\leq e^{\varepsilon(k)}t<\eta_{j+1},
\]
whence 
\[
\lambda_{k}t,t\in(\eta_{j},\eta_{j+1}).
\]
Therefore 
\[
d(x(t),y(t))=d(x(t),y(\lambda_{k}t))+d(y(\lambda_{k}t),y(t))\leq\varepsilon_{k}+2\varepsilon\leq3\varepsilon,
\]
where $t\in A$ is arbitrary. Since $\varepsilon>0$ is arbitrary,
we see that $x=y$ on the dense subset $A$. It follows from Lemma
\ref{Lem.  If x=00003Dy on dense A then x=00003Dy. If f(x,y)<=00003Dc on A then...-1}
that $x=y$. Summing up, we conclude that $d_{D}$ is a metric.

Finally, suppose $x,y\in D[0,1]$ are such that $d(x(t),y(t))\leq c$
for each $t$ in a dense subset of $domain(x)\cap domain(y)$. Then
$c\in B_{x,y}$ by Lemma \ref{Lem. B_xy basics}. Hence $d_{D}(x,y)\equiv\inf B_{x,y}\leq c$.
The Proposition is proved.
\end{proof}
The next theorem is now a trivial consequence of Lemma \ref{Lem. Prep Lemma Arzela-Ascoli for D=00005B0,1=00005D}.
\begin{thm}
\label{Th,m Arzela-Ascoli for d_D} \textbf{\emph{(Arzela-Ascoli Theorem
for ($D[0,1],d_{D})$). }}Let $B$ be an arbitrary compact subset
of $(S,d)$. Let $k\geq0$ be arbitrary. Let $U\equiv\{v_{1},\cdots,v_{M}\}$
be a $2^{-k-1}$-approximation of of $B$. Let $x\in D[0,1]$ be arbitrary
with values in the compact set $B$ and with a modulus of \textup{\emph{càdlàg}}
$\delta_{cdlg}$. Let $m\geq1$ be so large that
\[
m\equiv m(k,\delta_{cdlg})\equiv[0\vee(1-\log_{2}(\delta_{\log@1}(2^{-k})\delta_{cdlg}(2^{-k-1})))]_{1}
\]
Then there exist an increasing sequence $(\eta_{i})_{i=0,\cdots,n-1}$
in $Q_{m}$, and a sequence $(u_{i})_{i=0,\cdots,n-1}$ in $U$, such
that
\[
d_{D[0,1]}(x,\overline{x})<2^{-k},
\]
where\emph{ 
\[
\overline{x}\equiv\Phi_{smpl}((\eta_{i})_{i=0,\cdots,n-1},(u_{i})_{i=0,\cdots,n-1})\in D_{simple,m,U}[0,1].
\]
}
\end{thm}
\begin{proof}
By Lemma \ref{Lem. Prep Lemma Arzela-Ascoli for D=00005B0,1=00005D},
there exist an increasing sequence $(\eta_{i})_{i=0,\cdots,n-1}$
in $Q_{m}$, and a sequence $(u_{i})_{i=0,\cdots,n-1}$ in $U$, such
that $2^{-k}\in B_{x,\overline{x}}.$ At the same time $d_{D[0,1]}(x,\overline{x})\equiv\inf B_{x,y}$
according to Definition \ref{Def. Skorokhod Metric}. Hence $d_{D[0,1]}(x,\overline{x})<2^{-k},$
as desired.
\end{proof}
\begin{thm}
\label{Thm. Skorokod metric is complete} \textbf{\emph{(Skorokhod
space is Complete).}} The Skorokhod space $(D[0,1],d_{D})$ is complete. 
\end{thm}
\begin{proof}
1. Let $(y_{k})_{k=1,2,\cdots}$ be an arbitrary Cauchy sequence in
$(D[0,1],d_{D})$. We need to prove that $d_{D}(y_{k},y)\rightarrow0$
for some $y\in D[0,1]$. Since $(y_{k})_{k=1,2,\cdots}$ is Cauchy,
it suffices to show that some subsequence converges. Hence, by passing
to a subsequence if necessary, there is no loss in generality in assuming
that
\begin{equation}
d_{D}(y_{k},y_{k+1})<2^{-k}\label{eq:temp-355}
\end{equation}
for each $k\geq1$. Let $\delta_{k}\equiv\delta_{cdlg,k}$ be a modulus
of càdlàg of $y_{k}$, for each $k\geq1$. We may assume, again without
loss of generality, that $\delta_{k}\leq1$ for each $k\geq1$. For
convenience, let $y_{0}\equiv x_{\circ}$ denote the constant càdlàg
function on $[0,1]$. 

2. Let $k\geq1$ be arbitrary. Then 
\[
d_{D}(y_{0},y_{k})\leq d_{D}(y_{0},y_{1})+2^{-1}+\cdots+2^{-k+1}<b_{0}\equiv d_{D}(y_{0},y_{1})+1.
\]
Hence $d(x_{\circ},y_{k}(t))\leq b_{0}$ for each $t\in domain(y_{k})$.
Thus the values of $y_{k}$ are contained in some compact subset $B$
of $(S,d)$ which contains $(d(x_{\circ},\cdot)\leq b_{0})$. Define
$b\equiv2b_{0}$. Then $d(y_{k}(t),y_{h}(s))\leq2b_{0}\equiv b$,
for each $t\in domain(y_{k})$, $s\in domain(y_{h})$, for each $h,k\geq0$. 

3. Next, refer to Definition \ref{Def. Notations for dyadic rationals}
for notations related to dyadic rationals, and refer to Definition
\ref{Def. Skorokhod Metric} for the notations related to the Skorokhod
metric. In particular, for each $x,y\in D[0,1]$, recall the sets
$A_{x,y}$ and $B_{x,y}$. Thus $d_{D}(x,y)\equiv\inf B_{x,y}$. Let
$\lambda_{0}\in\Lambda$ denote the identity function on $[0,1]$. 

4. The next two steps will replace the sequence $(y_{k})_{k=0,1,\cdots}$
with a sequence $(x_{k})_{k=0,1,\cdots}$ of simple càdlàg functions
with division points which are dyadic rationals in $[0,1]$. To that
end, take an arbitrary $m_{0}\geq0$, and, inductively for each $k\geq1$,
take $m_{k}\geq1$ so large that 
\begin{equation}
2^{-m(k)}<2^{-m(k-1)-2}e^{-b}\delta_{\log@1}(2^{-k})\delta_{k}(2^{-k-1}).\label{eq:temp-243-2-1}
\end{equation}

5. Define $x_{0}\equiv y_{0}\equiv x_{\circ}$. Let $k\geq1$ be arbitrary,
and write $n_{k}\equiv p_{m(k)}$. Then 
\[
2^{-m(k)}<2^{-1}\delta_{\log@1}(2^{-k})\delta_{k}(2^{-k-1}).
\]
Hence Theorem \ref{Th,m Arzela-Ascoli for d_D}, the Arzela-Ascoli
Theorem for the Skorokhod space, applies to the càdlàg function $y_{k}$,
to yield a simple càdlàg function $x_{k}$ with the sequence 
\[
\tau_{k}\equiv(\tau_{k,0},\tau_{k,1},\cdots,\tau_{k,n(k)-1})\equiv q_{m(k)}\equiv(q_{m(k),0},q_{m(k),1},\cdots,q_{m(k),n(k)-1})
\]
of dyadic rationals as division points, such that 
\begin{equation}
d_{D}(y_{k},x_{k})\leq2^{-k}.\label{eq:temp-323}
\end{equation}
It follows that
\[
d_{D}(x_{k},x_{k+1})\leq d_{D}(x_{k},y_{k})+d_{D}(y_{k},y_{k+1})+d_{D}(y_{k+1},x_{k+1})
\]
\begin{equation}
<2^{-k}+2^{-k}+2^{-k-1}<2^{-k+2}.\label{eq:temp-322}
\end{equation}
Moreover, since $x_{0}\equiv x_{\circ}\equiv y_{0}$, we have
\[
d_{D}(x_{\circ},x_{k})\leq d_{D}(y_{0},y_{k})+2^{-k}
\]
\[
<d_{D}(y_{0},y_{1})+d_{D}(y_{1},y_{2})+\cdots++d_{D}(y_{k-1},y_{k})+2^{-k}
\]
\[
<d_{D}(y_{0},y_{1})+2^{-1}+\cdots+2^{-(k-1)}+2^{-k}<b_{0}\equiv d_{D}(y_{0},y_{1})+1.
\]
Thus $d(x_{\circ},x_{k}(t))\leq b_{0}$ for each $t\in domain(x_{k})$.
Consequently,
\begin{equation}
d(x_{k}(t),x_{h}(s))\leq b\equiv2b_{0}\label{eq:temp-244}
\end{equation}
for each $t\in domain(x_{k})$, $s\in domain(x_{h})$, for each $h,k\geq0$. 

6. Now let $k\geq0$ be arbitrary. We will construct $\lambda_{k+1}\in\Lambda_{m(k+1),m(k+2)}$
such that 
\begin{equation}
(2^{-k+2},\lambda_{k+1})\in A_{x(k),x(k+1)}\label{eq:temp-394}
\end{equation}
To that end, note that $x_{k},x_{k+1}\in D_{simple,m(k+1),U}[0,1]$,
where $U$ is the finite subset of $S$ consisting of the values of
the two simple càdlàg functions $x_{k},x_{k+1}$. Then $\bigvee_{u,v\in U}d(u,v)\leq b$
in view of inequality \ref{eq:temp-244}. At the same time, applying
inequality \ref{eq:temp-243-2-1} to $k+2$ in place of $k$, we obtain
\begin{equation}
2^{-m(k+2)}\leq2^{-m(k+1)-2}e^{-b}\delta_{\log@1}(2^{-k-2}).\label{eq:temp-325-3}
\end{equation}
Hence we can apply Lemma \ref{Lem. d_D(x,y) well defined for certain simple cadlags}
to the quintuple $(m_{k+1},m_{k+2},k+2,x_{k},x_{k+1})$ in place of
the quintuple $(m,m',k,x,y)$ there, to construct some $\lambda_{k+1}\in\Lambda_{m(k+1),m(k+2)}$
such that 
\begin{equation}
(d_{D}(x_{k},x_{k+1})+2^{-k-1},\lambda_{k+1})\in A_{x(k),x(k+1)}.\label{eq:temp-299-1}
\end{equation}
Since, from inequality \ref{eq:temp-322}, we have 
\[
d_{D}(x_{k},x_{k+1})+2^{-k-1}<(2^{-k}+2^{-k}+2^{-k-1})+2^{-k-1}<2^{-k+2},
\]
relation \ref{eq:temp-299-1} trivially implies the desired relation
\ref{eq:temp-394}. Consequently, 
\begin{equation}
|\log\frac{\lambda_{k+1}t-\lambda_{k+1}s}{t-s}|\leq2^{-k+2}\label{eq:temp-319}
\end{equation}
for each $s,t\in[0,1]$ with $s<t$, and 
\begin{equation}
d(x_{k+1}(\lambda_{k+1}t),x_{k}(t))\leq2^{-k+2}\label{eq:temp-391}
\end{equation}
for each $t\in domain(x_{k})\cap\lambda_{k+1}^{-1}domain(x_{k+1})$. 

7. For each $k\geq0$, define the composite admissible function
\[
\mu_{k}\equiv\lambda_{k}\lambda_{k-1}\cdots\lambda_{0}\in\Lambda.
\]
We will prove that $\mu_{k}\rightarrow\mu$ uniformly on $[0,1]$
for some $\mu\in\Lambda$. To that end, let $h>k\geq0$ be arbitrary.
By Lemma \ref{Lem. B_xy basics}, relation \ref{eq:temp-394} implies
that
\[
(2^{-k+2}+2^{-k+1}+\cdots+2^{-h+3},\lambda_{h}\lambda_{h-1}\cdots\lambda_{k+1})\in A_{x(k),x(h)}.
\]
Hence, since $2^{-k+2}+2^{-k+1}+\cdots+2^{-h+3}<2^{-k+3}$, we also
have 
\begin{equation}
(2^{-k+3},\lambda_{h}\lambda_{h-1}\cdots\lambda_{k+1})\in A_{x(k),x(h)}.\label{eq:temp-304}
\end{equation}
Let $t,s\in[0,1]$ be arbitrary with $s<t$. Write $t'\equiv\mu_{k}t$
and $s'\equiv\mu_{k}s$. Then 
\[
|\log\frac{\mu_{h}t-\mu_{h}s}{\mu_{k}t-\mu_{k}s}|
\]
\begin{equation}
=|\log\frac{\lambda_{h}\lambda_{h-1}\cdots\lambda_{k+1}t'-\lambda_{h}\lambda_{h-1}\cdots\lambda_{k+1}s'}{t'-s'}|<2^{-k+3}.\label{eq:temp-398}
\end{equation}
Equivalently,
\[
\exp(-2^{-k+3})(t'-s')<\lambda_{h}\lambda_{h-1}\cdots\lambda_{k+1}t'-\lambda_{h}\lambda_{h-1}\cdots\lambda_{k+1}s'
\]
\begin{equation}
<\exp(2^{-k+3})(t'-s')\label{eq:temp-393}
\end{equation}
for each $t',s'\in[0,1]$ with $s'<t'$. In the special case where
$s=0$, inequality \ref{eq:temp-398} reduces to 
\begin{equation}
|\log\mu_{h}t-\log\mu_{k}t|=|\log\frac{\mu_{h}t}{\mu_{k}t}|<2^{-k+3}.\label{eq:temp-392}
\end{equation}
Hence the limit 
\[
\mu t\equiv\lim_{h\rightarrow\infty}\mu_{h}t
\]
exists, where $t\in(0,1]$ is arbitrary. Moreover, letting $k=0$
and $h\rightarrow\infty$ in inequality \ref{eq:temp-398}, we obtain
\begin{equation}
|\log\frac{\mu t-\mu s}{t-s}|\leq2^{3}=8.\label{eq:temp-321}
\end{equation}
Therefore $\mu$ is an increasing function which is uniformly continuous
on $(0,1]$. Furthermore,
\[
te^{-8}\leq\mu t\leq te^{8}
\]
where $t\in(0,1]$ is arbitrary. Hence $\mu$ can be extended to a
continuous increasing function on $[0,1]$, with $\mu0=0$. Since
$\mu_{k}1=1$ for each $k\geq0$, we have $\mu1=1$. In view of inequality
\ref{eq:temp-321}, we conclude that $\mu\in\Lambda$.

8. By letting $h\rightarrow\infty$ in inequality \ref{eq:temp-398},
we obtain
\begin{equation}
|\log\frac{\mu t-\mu s}{\mu_{k}t-\mu_{k}s}|\leq2^{-k+3}\label{eq:temp-395}
\end{equation}
where $t,s\in[0,1]$ are arbitrary with $s<t$. Replacing $t,s$ by
$\mu^{-1}t,\mu^{-1}s$ respectively, we obtain 
\begin{equation}
|\log\frac{\mu_{k}\mu^{-1}t-\mu_{k}\mu^{-1}s}{t-s}|=|\log\frac{t-s}{\mu_{k}\mu^{-1}t-\mu_{k}\mu^{-1}s}|\leq2^{-k+3},\label{eq:temp-395-2}
\end{equation}
where $t,s\in[0,1]$ are arbitrary with $s<t$, and where $k\geq0$
is arbitrary.

9. Let $k\geq0$ be arbitrary. Recall from Step 5 that $n_{k}\equiv p_{m(k)}\equiv2^{m(k)},$
and that 
\[
\tau_{k}\equiv(\tau_{k,0},\tau_{k,1},\cdots,\tau_{k,n(k)-1})\equiv q_{m(k)}\equiv(q_{m(k),0},q_{m(k),1},\cdots,q_{m(k),n(k)-1})
\]
 is a sequence of division points of the simple càdlàg function $x_{k}$.
Define the set 
\begin{equation}
A\equiv\bigcap_{h=0}^{\infty}\mu\mu_{h}^{-1}([\tau_{h,0},\tau_{h,1})\cup[\tau_{h,1},\tau_{h,2})\cup\cdots\cup[\tau_{h,n(h)-1},\tau_{h,n(h)}]).\label{eq:temp-387}
\end{equation}
Then $A$ contains all but countably many points in $[0,1]$, and
is therefore dense in $[0,1]$. Define
\begin{equation}
u_{k}\equiv x_{k}\circ\mu_{k}\mu^{-1}.\label{eq:temp-402}
\end{equation}
Then $u_{k}\in D[0,1]$ by Lemma \ref{Lem.  If x=00003Dy on dense A then x=00003Dy. If f(x,y)<=00003Dc on A then...-1}.
Moreover,
\begin{equation}
u_{k+1}=x_{k+1}\circ\lambda_{k+1}\mu_{k}\mu^{-1}.\label{eq:temp-402-1}
\end{equation}
Now let $r\in A$ and $h\geq k$ be arbitrary. Then, by the defining
equality \ref{eq:temp-387} of the set $A$, there exists $i=0,\cdots,n_{h}-1$
such that 
\begin{equation}
\mu_{h}\mu^{-1}r\in[\tau_{h,i},\tau_{h,i+1})\subset domain(x_{h}).\label{eq:temp-411}
\end{equation}
Hence, 
\[
\lambda_{h+1}\mu_{h}\mu^{-1}r\in\lambda_{h+1}[\tau_{h,i},\tau_{h,i+1})
\]
Moreover, 
\begin{equation}
r\in domain(x_{h}\circ\mu_{h}\mu^{-1})\equiv domain(u_{h}).\label{eq:temp-300-1}
\end{equation}
From equalities \ref{eq:temp-402-1} and \ref{eq:temp-402}, and inequality
\ref{eq:temp-391}, we obtain
\begin{equation}
d(u_{k+1}(r),u_{k}(r))\equiv d(x_{k+1}(\lambda_{k+1}\mu_{k}\mu^{-1}r),x_{k}(\mu_{k}\mu^{-1}r))\leq2^{-k+2}.\label{eq:temp-390}
\end{equation}
Hence, since $r\in A$ and $k\geq0$ are arbitrary, we conclude that
$\lim_{k\rightarrow\infty}u_{k}$ exists on $A$. Define the function
$u:[0,1]\rightarrow S$ by $domain(u)\equiv A$ and by $u(t)\equiv\lim_{k\rightarrow\infty}u_{k}(t)$
for each $t\in domain(u)$. Inequality \ref{eq:temp-390} then implies
that 
\begin{equation}
d(u(r),u_{k}(r))\leq2^{-k+3}\label{eq:temp-412-1}
\end{equation}
where $r\in A$ and $k\geq0$ are arbitrary. We proceed to verify
the conditions in Proposition \ref{Prop. Cadlag completion} for the
function $u$ to have a càdlàg completion.

10. For that purpose, let $r\in A$ and $h\geq k$ be arbitrary. Then,
as observed in Step 9, we have $r\in domain(u_{h})$. Hence, since
$u_{h}$ is càdlàg, it is right continuous at $r$. Therefore there
exists $c_{h}>0$ such that 
\[
d(u_{h}(t),u_{h}(r))<2^{-h+3}
\]
for each $t\in[r,r+c_{h})\cap A$. In view of inequality \ref{eq:temp-412-1},
it follows that 
\[
d(u(t),u(r))\leq d(u(t),u_{h}(t))+d(u_{h}(t),u_{h}(r))+d(u_{h}(r),u(r))<3\cdot2^{-h+3}
\]
for each $t\in[r,r+c_{h})\cap A$. Thus $u$ is right continuous at
each point $r\in A\equiv domain(u)$. Condition 1 in Definition \ref{Def. Cadlag functions, D=00005B0,1=00005D}
has been verified for the function $u$.

11. We will now verify Conditions 3 in Definition \ref{Def. Cadlag functions, D=00005B0,1=00005D}.
To that end, let $\varepsilon>0$ be arbitrary. Take $k\geq0$ so
large that $2^{-k+4}<\varepsilon$, and define $\delta(\varepsilon)\equiv\exp(-2^{-k+3})2^{-m(k)}$.
Let $j=0,\cdots,n_{k}$ be arbitrary. For brevity, define $\tau'_{j}\equiv\mu\mu_{k}^{-1}\tau_{k,j}$.
By inequality \ref{eq:temp-319}, we have 
\[
\lambda_{h+1}t-\lambda_{h+1}s\geq(t-s)\exp(-2^{-h+2})
\]
for each $s,t\in[0,1]$ with $s<t$, for each $h\geq0$. Hence, for
each $j=0,\cdots,n_{k}-1$, we have 
\[
\tau'_{j+1}-\tau'_{j}\equiv\mu\mu_{k}^{-1}\tau_{k,j+1}-\mu\mu_{k}^{-1}\tau_{k,j}
\]
\[
=\lim_{h\rightarrow\infty}(\mu_{h}\mu_{k}^{-1}\tau_{k,j+1}-\mu_{h}\mu_{k}^{-1}\tau_{k,j})
\]
\[
=\lim_{h\rightarrow\infty}(\lambda_{h}\cdots\lambda_{k+1}\tau_{k,j+1}-\lambda_{h}\cdots\lambda_{k+1}\tau_{k,j})
\]
\begin{equation}
\geq\exp(-2^{-k+3})(\tau_{k,j+1}-\tau_{k,j})=\exp(-2^{-k+3})2^{-m(k)}\equiv\delta(\varepsilon)\label{eq:temp-396}
\end{equation}
where the inequality is by inequality \ref{eq:temp-393}. Now consider
each $j=0,\cdots,n_{k}-1$ and each 
\[
t'\in domain(u)\cap[\tau'_{j},\tau'_{j+1})\equiv A\cap[\mu\mu_{k}^{-1}\tau_{k,j},\mu\mu_{k}^{-1}\tau_{k,j+1}).
\]
We will show that 
\begin{equation}
d(u(t'),u(\tau'_{j})))\leq\varepsilon.\label{eq:temp-412-1-1}
\end{equation}
To that end, write $t\equiv\mu^{-1}t'$, and write $s\equiv\mu_{k}t\equiv\mu_{k}\mu^{-1}t'\in[\tau_{k,j},\tau_{k,j+1})$.
Then 
\begin{equation}
x_{k}(s)=x_{k}(\tau_{k,j})\label{eq:temp-301-1}
\end{equation}
since $x_{k}$ is a simple càdlàg function with $(\tau_{k,0},\tau_{k,1},\cdots,\tau_{k,n(k)-1})$
as a sequence of division points. Let $h>k$ be arbitrary, and define
\[
r\equiv\mu_{h}\mu^{-1}t'=\mu_{h}t=\lambda_{h}\lambda_{h-1}\cdots\lambda_{k+1}s.
\]
Then 
\[
u_{h}(t')\equiv x_{h}(\mu_{h}\mu^{-1}t')\equiv x_{h}(r)\equiv x_{h}(\lambda_{h}\lambda_{h-1}\cdots\lambda_{k+1}s).
\]
Combining with equality \ref{eq:temp-301-1}, we obtain 
\[
d(u_{h}(t'),x_{k}(\tau_{k,j}))=d(x_{h}(\lambda_{h}\lambda_{h-1}\cdots\lambda_{k+1}s),x_{k}(s))\leq2^{-k+3},
\]
where the last inequality is a consequence of relation \ref{eq:temp-304},
and Consequently
\begin{equation}
d(u_{h}(t'),u_{h}(\tau'_{j}))\leq d(u_{h}(t'),x_{k}(\tau_{k,j}))+d(u_{h}(\tau'_{j}),x_{k}(\tau_{k,j}))\leq2^{-k+4}<\varepsilon.\label{eq:temp-397}
\end{equation}
where $t'\in domain(u)\cap[\tau'_{j},\tau'_{j+1})$ is arbitrary.
Letting $h\rightarrow\infty$, we obtain the desired inequality \ref{eq:temp-412-1-1}.
Inequalities \ref{eq:temp-412-1-1} and \ref{eq:temp-396} together
say that $(\tau'_{j})_{j=0,\cdots,n(k)+1}$ is a sequence of $\varepsilon$-division
points with separation at least $\delta(\varepsilon)$. 

Thus Conditions 1 and 3 in Definition \ref{Def. Cadlag functions, D=00005B0,1=00005D}
have been verified for the objects $u,\tau'$, and $\delta_{cdlg}\equiv\delta$.
Therefore Proposition \ref{Prop. Cadlag completion} implies that
(i') the completion $y\in D[0,1]$ of $u$ is well-defined, (ii')
$y|domain(u)=u$, and (iii') $\delta_{cdlg}$ is a modulus of càdlàg
of $y$.

12. Finally, we will prove that $d_{D}(y_{h},y)\rightarrow0$ as $h\rightarrow\infty$.
To that end, let $h\geq0$ and
\[
r\in A\subset domain(u)\subset domain(y)
\]
be arbitrary. By the above Condition (ii') we have $u(r)=y(r)$. Hence
\[
d(y(r),x_{h}\circ\mu_{h}\mu^{-1}(r))\equiv d(y(r),u_{h}(r))=d(u(r),u_{h}(r))\leq2^{-h+3}
\]
by inequality \ref{eq:temp-412-1}. Consequently, since $A$ is a
dense subset of $[0,1]$, Lemma \ref{Lem.  If x=00003Dy on dense A then x=00003Dy. If f(x,y)<=00003Dc on A then...-1},
applied to $h$ in the place of $k$, yields 
\begin{equation}
d(y(r),x_{h}\circ\mu_{h}\mu^{-1}(r))\leq2^{-h+3}\label{eq:temp-397-1}
\end{equation}
for each $r\in domain(y)\cap\mu\mu_{h}^{-1}domain(x_{h})$. Inequalities
\ref{eq:temp-397-1} and \ref{eq:temp-395-2} together imply that
$(2^{-h+3},\mu_{h}\mu^{-1})\in A_{y,x(h)}$, whence $2^{-h+3}\in B_{y,x(h)}$.
Accordingly,
\[
d_{D}(y,x_{h})\equiv\inf B_{y,x(h)}\leq2^{-h+3}.
\]
Thus $d_{D}(x_{h},y)\rightarrow0$ as $h\rightarrow\infty$. In view
of inequality \ref{eq:temp-323}, we conclude that $d_{D}(y_{h},y)\rightarrow0$
as $h\rightarrow\infty$. 

13. Summing up, for each Cauchy sequence $(y_{h})_{h=0,1,\cdots}$,
there exists $y\in D[0,1]$ such that $d_{D}(y_{h},y)\rightarrow0$
as $h\rightarrow\infty$. In other words, $(D[0,1],d_{D})$ is complete,
as alleged. 
\end{proof}

\section{a.u. Càdlàg Processes}

Let $(\Omega,L,E)$ be a probability space, and let $(S,d)$ be a
locally compact metric space. Let $(D[0,1],d_{D})$ be the Skorokhod
space of càdlàg functions on the unit interval $[0,1]$ with values
in $(S,d)$, as introduced in the previous section. 

Recall Definition \ref{Def. Notations for dyadic rationals} for notations
related to the enumerated set of dyadic rationals $Q_{\infty}$ in
$[0,1]$.
\begin{defn}
\label{Def. a.u. Random cadlag function} \textbf{(a.u. random càdlàg
function).} A r.v. $Y:\Omega\rightarrow(D[0,1],d_{D})$ with values
in the Skorokhod space is called an \emph{almost uniform} (a.u.) \emph{random
càdlàg function}\index{random càdlàg function}, if, for each $\varepsilon>0$,
there exists a $\mathrm{measurable}$ set $A$ with $P(A)<\varepsilon$
such that members of the family $\{Y(\omega):\omega\in A^{c}\}$ are
càdlàg and share a common modulus of càdlàg\emph{.}
\end{defn}
Of special interest is the subclass of the a.u. random càdlàg functions
corresponding to a.u. càdlàg processes on $[0,1]$, defined next. 
\begin{defn}
\label{Def. a.u. cadlag process} \textbf{(a.u. càdlàg process).}
Let $X:[0,1]\times\Omega\rightarrow S$ be a stochastic process which
is continuous in probability on $[0,1].$ with a modulus of continuity
in probability $\delta_{cp}$. Suppose there exists a full set $B\subset\bigcap_{t\in Q(\infty)}domain(X_{t})$
with the following properties.

1 (Right continuity). For each $\omega\in B$, the function $X(\cdot,\omega)$
is right continuous at each $t\in domain(X(\cdot,\omega))$.

2 (Right completeness).  Let $\omega\in B$ and $t\in[0,1]$ be arbitrary.
If $\lim_{r\rightarrow t;r>t}X(r,\omega)$ exists, then $t\in domain(X(\cdot,\omega))$.

3 (Approximation by step functions). Let $\varepsilon>0$ be arbitrary.
Then there exist (i) $\delta_{aucl}(\varepsilon)>0$, (ii) a $\mathrm{measurable}$
set $A\subset B$ with $P(A^{c})<\varepsilon$, (iii) an integer $h\geq1$,
and (iv) a sequence of r.r.v.'s 
\begin{equation}
0=\tau_{0}<\tau_{1}<\cdots<\tau_{h-1}<\tau_{h}=1,\label{eq:temp-404}
\end{equation}
such that, for each $i=0,\cdots,h-1$, the function $X_{\tau(i)}$
is a r.v., and such that, (v) for each $\omega\in A$, we have
\begin{equation}
\bigwedge_{i=0}^{h-1}(\tau_{i+1}(\omega)-\tau_{i}(\omega))\geq\delta_{aucl}(\varepsilon),\label{eq:temp-342}
\end{equation}
with
\begin{equation}
d(X(\tau_{i}(\omega),\omega),X(\cdot,\omega))\leq\varepsilon,\label{eq:temp-307-4-1}
\end{equation}
on the interval $\theta_{i}(\omega)\equiv[\tau_{i}(\omega),\tau_{i+1}(\omega))$
or $\theta_{i}(\omega)\equiv[\tau_{i}(\omega),\tau_{i+1}(\omega)]$
according as $0\leq i\leq h-2$ or $i=h-1$. 

Then the process $X:[0,1]\times\Omega\rightarrow S$ is called an
\index{a.u. càdlàg process}\emph{a.u. càdlàg process}, with $\delta_{cp}$
as a modulus of continuity in probability and with $\delta_{aucl}$
as a \index{modulus of a.u. càdlàg}\emph{modulus of a.u. càdlàg}.
We will let $\widehat{D}_{\delta(aucl),\delta(cp)}[0,1]$ denote the
set of all such processes. 

We will let $\widehat{D}[0,1]$ denote the set of all a.u. càdlàg
processes. Two members $X,Y$ of $\widehat{D}[0,1]$ are considered
equal if there exists a full set $B'$ such that for each $\omega\in B'$
we have $X(\cdot,\omega)=Y(\cdot,\omega)$ as functions on $[0,1]$.
\end{defn}
$\square$
\begin{defn}
\label{Def. Random cadlag functions} \textbf{(a.u. Random càdlàg
function by extension of an a.u. càdlàg process).} Let $X\in\widehat{D}[0,1]$
be arbitrary. Define a function
\[
X^{*}:\Omega\rightarrow D[0,1]
\]
by 
\[
domain(X^{*})\equiv\{\omega\in\Omega:X(\cdot,\omega)\in D[0,1]\}
\]
and
\[
X^{*}(\omega)\equiv X(\cdot,\omega)
\]
for each $\omega\in domain(X^{*})$. We call $X^{*}$ the \emph{a.u.
random càdlàg function}\index{random càdlàg function} associated
with the a.u. càdlàg process $X$. $\square$
\end{defn}
\begin{prop}
\label{Prop. a.u. cadlag process is r.v. with values in D=00005B0,1=00005D}
\textbf{\emph{(Each a.u. càdlàg process extends to a random}} \textbf{\emph{càdlàg
function). }}\emph{Let }$X\in\widehat{D}[0,1]$ be an arbitrary a.u.
càdlàg process. Let $\varepsilon>0$ be arbitrary. Then there exists
a $\mathrm{measurable}$ set $G$ with $P(G^{c})<\varepsilon$, such
that members of the set $\{X(\cdot,\omega):\omega\in G\}$ are càdlàg
functions which share a common modulus of càdlàg. 

Moreover, $X^{*}$ is a r.v. with values in the complete metric space
$(D[0,1],d_{D}).$
\end{prop}
\begin{proof}
As in Definition \ref{Def. a.u. Random cadlag function}, define the
full subset $B\equiv\bigcap_{t\in Q(\infty)}domain(X_{t})$ of $\Omega$.
By Conditions (1) and (2) in Definition \ref{Def. a.u. Random cadlag function},
members of the set $\{X(\cdot,\omega):\omega\in B\}$ of functions
satisfy the corresponding Conditions (1) and (2) in Definition \ref{Def. Cadlag functions, D=00005B0,1=00005D}. 

Let $n\geq1$ be arbitrary. By Condition (3) in Definition \ref{Def. a.u. Random cadlag function},
there exist (i) $\delta_{aucl}(2^{-n})>0$, (ii) a $\mathrm{measurable}$
set $A_{n}\subset B$ with $P(A_{n}^{c})<2^{-n}$, (iii) an integer
$h_{n}\geq0$, and (iv) a sequence of r.r.v.'s
\begin{equation}
0=\tau_{n,0}<\tau_{n,1}<\cdots<\tau_{n,h(n)}=1,\label{eq:temp-405}
\end{equation}
such that, for each $i=0,\cdots,h_{n}-1$, the function $X_{\tau(n,i)}$
is a r.v., and such that, for each $\omega\in A_{n}$, we have
\begin{equation}
\bigwedge_{i=0}^{h(n)-1}(\tau_{n,i+1}(\omega)-\tau_{n,i}(\omega))\geq\delta_{aucl}(2^{-n}),\label{eq:temp-408}
\end{equation}
with
\begin{equation}
d(X(\tau_{n,i}(\omega),\omega),X(\cdot,\omega))\leq2^{-n}\label{eq:temp-409}
\end{equation}
on the interval $\theta_{n,i}\equiv[\tau_{n,i}(\omega),\tau_{n,i+1}(\omega))$
or $\theta_{n,i}\equiv[\tau_{n,i}(\omega),\tau_{n,i+1}(\omega)]$
according as $i\leq h_{n}-2$ or $i=h_{n}-1$. 

Now let $\varepsilon>0$ be arbitrary. Take $j\geq1$ so large that
$2^{-j}<\varepsilon$. Define $G\equiv\bigcap_{n=j+1}^{\infty}A_{n}\subset B$.
Then $P(G^{c})<\sum_{n=j+1}^{\infty}2^{-n}=2^{-j}<\varepsilon$. Consider
each $\omega\in G$. Let $\varepsilon'>0$ be arbitrary. Consider
each $n\geq j+1$ so large that $2^{-n}<\varepsilon'$. Define $\delta_{cdlg}(\varepsilon')\equiv\delta_{aucl}(2^{-n})$.
Then $\omega\in A_{n}$. Hence inequalities \ref{eq:temp-408} and
\ref{eq:temp-409} hold, and imply that the sequence
\begin{equation}
0=\tau_{n,0}(\omega)<\tau_{n,1}(\omega)<\cdots<\tau_{n,h(n)}(\omega)=1\label{eq:temp-405-1}
\end{equation}
is a sequence of $\varepsilon'$-division points of the function $X(\cdot,\omega)$,
with separation at least $\delta_{cdlg}(\varepsilon')$. Summing up,
all the conditions in Definition \ref{Def. Cadlag functions, D=00005B0,1=00005D}
have been verified for the function $X(\cdot,\omega)$ to be càdlàg,
with the modulus of càdlàg $\delta_{cdlg}$, where $\omega\in G$
is arbitrary. This proves the first part of the proposition.

Let $H_{h(n),\delta(0)}$ be the subset of the product metric space
$([0,1]^{h(n)},d_{ecld}^{h(n)})\otimes(S^{h(n)},d^{h(n)})$ consisting
of elements $\alpha\equiv((a_{0},\cdots,a_{h(n)-1}),(x_{0},\cdots,x_{h(n)-1}))$
such that 
\begin{equation}
0\equiv a_{0}<a_{1}<\cdots<a_{h(n)}\equiv1\label{eq:temp-405-2-1-1}
\end{equation}
with 
\begin{equation}
\bigwedge_{i=0}^{h(n)-1}(a_{i+1}-a_{i})\geq\delta_{0}.\label{eq:temp-408-1-1-1}
\end{equation}
Then the function 
\[
\Phi_{smpl}:(H_{h,(n)\delta(0)},d_{ecld}^{h(n)}\otimes d^{h(n)})\rightarrow(D[0,1],d_{D}),
\]
which assigns to each $((a_{0},\cdots,a_{h(n)-1}),(x_{0},\cdots,x_{h(n)-1}))\in H_{h(n),\delta(0)}$
the simple càdlàg function 
\[
\Phi_{smpl}((a_{0},\cdots,a_{h(n)-1}),(x_{0},\cdots,x_{h(n)-1})),
\]
is uniformly continuous. Hence, the function
\[
X^{(n)}\equiv\Phi_{smpl}((\tau_{n,0},\cdots,\tau_{n,h(n)-1}),(X_{\tau(n,0)},\cdots,X_{\tau(n,h(n)-1)}))
\]
is a r.v. with values in $(D[0,1],d_{D})$, according to Proposition
\ref{Prop. X meas,  f unif continuous and bd on bd subsets =00003D> f(X) meas}.
At the same time, inequality \ref{eq:temp-409} implies that 
\[
d_{D}(X^{(n)},X^{*})\leq2^{-n},
\]
on $A_{n}$, where $P(A_{n}^{c})<2^{-n}$. Hence $X^{(n)}\rightarrow X^{*}$
a.u. Consequently $X^{*}$ is a r.v. with values in $(D[0,1],d_{D})$. 
\end{proof}
\begin{defn}
\label{Def. metric rho_D_hat } \textbf{(Metric for the space of a.u.
càdlàg processes).} Define the metric $\rho_{\widehat{D}[0,1]}$ on
$\widehat{D}[0,1]$ by
\begin{equation}
\rho_{\widehat{D}[0,1]}(X,Y)\equiv\int E(d\omega)1\wedge d_{D}(X(\cdot,\omega),Y(\cdot,\omega))\equiv E1\wedge d_{D}(X^{*},Y^{*})\label{eq:temp-202-2-1}
\end{equation}
for each $X,Y\in\widehat{D}[0,1]$. The next lemma justifies the definition.
$\square$
\end{defn}
\begin{lem}
\label{Lem. metric rho_D_hat is well defined.} \textbf{\emph{($(\widehat{D}[0,1],$$\rho_{\widehat{D}[0,1]})$
is a metric space).}} The function $\rho_{\widehat{D}[0,1]}$ is well-defined
and is a metric. 
\end{lem}
\begin{proof}
Let $X,Y\in\widehat{D}[0,1]$ be arbitrary. Then, according to Proposition
\ref{Prop. a.u. cadlag process is r.v. with values in D=00005B0,1=00005D},
the random càdlàg functions $X^{*},Y^{*}$ associated with $X,Y$
respectively are r.v.'s with values in $(D[0,1],d_{D})$. Therefore
the function $1\wedge d_{D}(X^{*},Y^{*})$ is an integrable r.r.v.,
and the defining equality \ref{eq:temp-202-2-1} makes sense.

Symmetry and the triangle inequality can be trivially verified for
$\rho_{\widehat{D}[0,1]}$. Now suppose $X=Y$ in $\widehat{D}[0,1]$.
Then, by the definition of the equality relation for the set $\widehat{D}[0,1]$,
we have $X(\cdot,\omega)=Y(\cdot,\omega)$ in $D[0,1]$, for a.e.
$\omega\in\Omega$. Hence
\[
1\wedge d_{D}(X^{*}(\omega),Y^{*}(\omega))\equiv1\wedge d_{D}(X(\cdot,\omega),Y(\cdot,\omega))=0.
\]
Thus $1\wedge d_{D}(X^{*},Y^{*})=0$ a.s. Consequently, $\rho_{\widehat{D}[0,1]}(X,Y)=0$
according to the defining equality \ref{eq:temp-202-2-1}. The converse
is proved similarly. Combining, we conclude that $\rho_{\widehat{D}[0,1]}$
is a metric.
\end{proof}

\section{$D$-regular Families of f.j.d.'s and $D$-regular Processes}

In this and the following two sections, we give a construction of
an a.u. process from a consistent family F of f.j.d.'s with parameter
set $[0,1]$ which satisfies a certain $D$-regularity condition to
be defined presently. The construction is by (i) taking any process
$Z:Q_{\infty}\times\Omega\rightarrow S$ with marginal distributions
given by $F|Q_{\infty}$, (ii) extending the process $Z$ to an a.u.
càdlàg, process $X:[0,1]\times\Omega\rightarrow S$ by taking right
limits of sample paths. Step (i) can be done by the Daniell-Kolmogorov
Extension or Daniell-Kolmogorv-Skorokhod Extension, for example. As
a matter of fact, we can define $D$-regularity for $F$ as $D$-regularity
of any process $Z$ with marginal distributions given by $F$. The
key Step (ii) is then by proving that that a process $X:[0,1]\times\Omega\rightarrow S$
is a.u. càdlàg iff it is the right-limit extension of some $D$-regular
process $Z:Q_{\infty}\times\Omega\rightarrow S$. In this section,
we prove the ``only if'' part. In the next section, we will prove
the ``if'' part, which is the useful part for the purpose of the
first sentence in this paragraph.
\begin{defn}
\label{Def.  D-regular process and D-regular family of f.j.d.s on Q_inf}
\textbf{(}$D$\textbf{-regular processes and }$D$\textbf{-regular
families of f.j.d.'s with parameter set $Q_{\infty}$).} Let $Z:Q_{\infty}\times\Omega\rightarrow S$
be a stochastic process, with marginal distributions given by the
family $F$ of f.j.d.'s. Let $\overline{m}\equiv(m_{n})_{n=0.1.\cdots}$
be an increasing sequence of nonnegative integers. Suppose the following
conditions are satisfied.

1. Let $n\geq0$ be arbitrary. Let $\beta>2^{-n}$ be arbitrary such
that the set 
\begin{equation}
A_{t,s}^{\beta}\equiv(d(Z_{t},Z_{s})>\beta)\label{eq:temp-372}
\end{equation}
is $\mathrm{measurable}$ for each $s,t\in Q_{\infty}$. Then 
\begin{equation}
P(D_{n})<2^{-n},\label{eq:temp-204-3}
\end{equation}
where we define the exceptional set 
\begin{equation}
D_{n}\equiv\bigcup_{t\in Q(m(n))}\bigcup_{r,s\in(t,t')Q(m(n+1));r\leq s}(A_{t,r}^{\beta}\cup A_{t,s}^{\beta})(A_{t,r}^{\beta}\cup A_{t',s}^{\beta})(A_{t',r}^{\beta}\cup A_{t',s}^{\beta}),\label{eq:temp-201-4-3}
\end{equation}
where for each $t\in Q{}_{m(n)}$ we abuse notations and write $t'\equiv1\wedge(t+2^{-m(n)}).$ 

2. The process $Z$ is continuous in probability on $Q_{\infty},$
with a modulus of continuity in probability $\delta_{Cp}$.

Then the process $Z:Q_{\infty}\times\Omega\rightarrow S$ and the
family $F$ of f.j.d.'s are said to be $D$\emph{-regular} \index{D-regular process on Q_{infty}@D-regular process on $Q_{\infty}$}
\index{D-regular family of f.j.d.'s on Q_{infty}@D-regular family of f.j.d.'s on $Q_{\infty}$},
with the sequence $\overline{m}$ as a \emph{modulus of }$D$\emph{-regularity
}\index{modulus of D-regularity} and with the operation $\delta_{Cp}$
as a modulus of continuity in probability. 

Let 
\[
\widehat{R}_{Dreg,\overline{m},\delta(Cp)}(Q_{\infty}\times\Omega,S)
\]
denote the set of all such processes. Let $\widehat{R}_{Dreg}(Q_{\infty}\times\Omega,S)$
denote the set of all $D$-regular processes. Thus $\widehat{R}_{Dreg,\overline{m},\delta(Cp)}(Q_{\infty}\times\Omega,S)$
and $\widehat{R}_{Dreg}(Q_{\infty}\times\Omega,S)$ are subsets of
the metric space $(\widehat{R}(Q_{\infty}\times\Omega,S),\widehat{\rho}_{Prob,Q(\infty)})$
introduced in Definition \ref{Def. Metric on random fields w/ countable parameters},
and, as such, inherit the metric $\widehat{\rho}_{Prob,Q(\infty)}$.
Thus we have the metric space 
\[
(\widehat{R}_{Dreg}(Q_{\infty}\times\Omega,S),\widehat{\rho}_{Prob,Q(\infty)}).
\]

Let 
\[
\widehat{F}_{Dreg,\overline{m},\delta(Cp)}(Q_{\infty},S)
\]
denote the set of all such families $F$ of f.j.d.'s. Let $\widehat{F}_{Dreg}(Q_{\infty},S)$
denote the set of all $D$-regular families of f.j.d.'s. Then $\widehat{F}_{Dreg,\overline{m},\delta(Cp)}(Q_{\infty},S)$
and $\widehat{F}_{Dreg}(Q_{\infty},S)$ are a subsets of the metric
space $(\widehat{F}(Q_{\infty},S),\widehat{\rho}_{Marg,\xi,Q(\infty)})$
of consistent families of f.j.d.'s introduced in Definition \ref{Def. Marginal metric},
where the metric $\widehat{\rho}_{Marg,\xi,Q(\infty)}$ is defined
relative to an arbitrarily given, but fixed, binary approximation
$\xi\equiv(A_{q})_{q=1,2,\cdots}$ of $(S,d)$.
\end{defn}
$\square$

Condition 1 in Definition \ref{Def.  D-regular process and D-regular family of f.j.d.s on Q_inf}
is, in essence, equivalent to Condition (13.10) in the key Theorem
13.3 of \cite{Billingsley 99}. The crucial difference of our construction
and the last cited theorem is that the latter proves only that (i)
a sequence of distributions on the path space which satisfies said
Condition (13.10) is tight, and (ii) the weak-convergence limit of
any such subsequence sequence of distributions, if such limit exists,
and (iii) existence of weak-convergence limit is then guaranteed by
Prokhorov's Theorem, Theorem 5.1 in \cite{Billingsley 99}, which
says that each subsequence of a tight sequence of distributions contains
a subsequence which has a weak-convergence limit. As we observed earlier,
Prokhorov's Theorem implies the principle of infinite search. This
is in contrast to our simple and direct construction in developed
in this and the next section.

First we extend the definition of $D$-regularity to families of f.j.d.'s
with parameter interval $[0,1]$ which are continuous in probability
on the interval. 
\begin{defn}
\textbf{\label{Def. metric space of D-regular families of f.j.d.'s on =00005B0,1=00005D}(}$D$\textbf{-regular
families of f.j.d.'s with parameter interval $[0,1]$). }Recall from
\ref{Def. Metric on  of continuous in prob families of finite joint distributions.}the
metric space $(\widehat{F}_{Cp}([0,1],S),\widehat{\rho}_{Cp,\xi,[0,1]|Q(\infty)})$
of families of f.j.d.'s which are continuous in probability on $[0,1]$,
where the metric is defined relative to the enumerated, countable,
dense subset $Q_{\infty}$ of $[0,1]$. Define two subsets of $\widehat{F}_{Cp}([0,1],S)$
by\textbf{
\[
\widehat{F}_{Dreg}([0,1],S)\equiv\{F\in\widehat{F}_{Cp}([0,1],S):\;F|Q_{\infty}\in\widehat{F}_{Dreg}(Q_{\infty},S)\},
\]
}and\textbf{
\[
\widehat{F}_{Dreg,\overline{m},\delta(Cp)}([0,1],S)\equiv\{F\in\widehat{F}_{Cp}([0,1],S):\;F|Q_{\infty}\in\widehat{F}_{Dreg,\overline{m},\delta(Cp)}(Q_{\infty},S)\}.
\]
}These subsets inherit the metric $\widehat{\rho}_{Cp,\xi,[0,1]|Q(\infty)}$. 

$\square$
\end{defn}
We will prove that a process $X:[0,1]\times\Omega\rightarrow S$ is
a.u càdlàg iff it is the extension by right limit of a $D$-regular
process $Z:Q_{\infty}\times\Omega\rightarrow S$. The next theorem
proves the ``only if'' part.
\begin{thm}
\label{Thm. Restriction of a.u. cadlag process to dyadics is D-regular}
\textbf{\emph{(Restriction of each a.u. Càdlàg process to $Q_{\infty}$
is }}$D$\textbf{\emph{-regular).}} Let $X:[0,1]\times\Omega\rightarrow S$
be an a.u càdlàg process, with a modulus of a.u. càdlàg $\delta_{aucl}$
and a modulus of continuity in probability $\delta_{Cp}$. Let $\overline{m}\equiv(m_{n})_{n=0,1,2,\cdots}$
be an arbitrary increasing sequence of integers such that 
\begin{equation}
2^{-m(n)}<\delta_{aucl}(2^{-n-1}).\label{eq:temp-369}
\end{equation}
for each $n\geq0$.

Then the process $Z\equiv X|Q_{\infty}$ is $D$-regular, with a modulus
of $D$-regularity $\overline{m}$, and with the same modulus of continuity
in probability $\delta_{Cp}$. 
\end{thm}
\begin{proof}
First note that, by Definition \ref{Def. a.u. Random cadlag function},
$X$ is continuous in probability on $[0,1],$ with some modulus of
continuity in probability $\delta_{Cp}$. Hence so is $Z$ on $Q_{\infty}$,
with the same modulus of continuity in probability $\delta_{Cp}$.
Consequently, Condition 2 in Definition \ref{Def.  D-regular process and D-regular family of f.j.d.s on Q_inf}
is satisfied.

Now let $n\geq0$ be arbitrary. Write $\varepsilon_{n}\equiv2^{-n}$.
By Condition (3) in Definition \ref{Def. a.u. cadlag process}, there
exist (i) $\delta_{aucl}(\frac{1}{2}\varepsilon_{n})>0$, (ii) a $\mathrm{measurable}$
set $A_{n}\subset B\equiv\bigcap_{t\in Q(\infty)}domain(X_{t})$ with
$P(A_{n}^{c})<\frac{1}{2}\varepsilon_{n}$, (iii) an integer $h_{n}\geq0$,
and (iv) a sequence of r.r.v.'s,
\begin{equation}
0=\tau_{n,0}<\tau_{n,1}<\cdots<\tau_{n,h(n)-1}<\tau_{n,h(n)}=1,\label{eq:temp-410}
\end{equation}
such that, for each $i=0,\cdots,h_{n}-1$, the function $X_{\tau(n,i)}$
is a r.v., and such that, for each $\omega\in A_{n}$, we have
\begin{equation}
\bigwedge_{i=0}^{h(n)-1}(\tau_{n,i+1}(\omega)-\tau_{n,i}(\omega))\geq\delta_{aucl}(\frac{1}{2}\varepsilon_{n})\label{eq:temp-400}
\end{equation}
with
\begin{equation}
d(X(\tau_{n,i}(\omega),\omega),X(\cdot,\omega))\leq\frac{1}{2}\varepsilon_{n},\label{eq:temp-401}
\end{equation}
on the interval $\theta_{n,i}\equiv[\tau_{n,i}(\omega),\tau_{n,i+1}(\omega))$
or $\theta_{n,i}\equiv[\tau_{n,i}(\omega),1]$ according as $i\leq h_{n}-2$
or $i=h_{n}-1$.

Take any $\beta>2^{-n}$ such that the set 
\begin{equation}
A_{t,s}^{\beta}\equiv(d(X_{t},X_{s})>\beta)\label{eq:temp-372-2}
\end{equation}
is $\mathrm{measurable}$ for each $s,t\in Q_{\infty}$. Let $D_{n}$
be the exceptional set as in defining equality \ref{eq:temp-201-4-3}.

Suppose, for the sake of a contradiction, that $P(D_{n})>\varepsilon_{n}\equiv2^{-n}$.
Then $P(D_{n})>\varepsilon_{n}>P(A_{n}^{c})$ by Condition (ii) above.
Hence $P(D_{n}A_{n})>0$. Consequently, there exists $\omega\in D_{n}A_{n}$.
Since $\omega\in A_{n}$, inequalities \ref{eq:temp-400} and \ref{eq:temp-401}
hold at $\omega$. At the same time, since $\omega\in D_{n},$ there
exist $t\in Q_{m(n)}$ and $r,s\in Q_{m(n+1)}$, with $t<r\leq s<t'\equiv t+2^{-m(n)}$,
such that 
\begin{equation}
\omega\in(A_{t,r}^{\beta}\cup A_{t,s}^{\beta})(A_{t,r}^{\beta}\cup A_{t',s}^{\beta})(A_{t',r}^{\beta}\cup A_{t',s}^{\beta}).\label{eq:temp-416}
\end{equation}
At the same time, the interval $[0,1)$ is contained in the union
$\bigcup_{i=1}^{h(n)-1}[\tau_{n,i-1}(\omega),\tau_{n,i+1}(\omega))$.
Hence $t\in[\tau_{n,i-1}(\omega),\tau_{n,i+1}(\omega)$) for some
$i=1,\cdots,h_{n}-1$. Now consider each 
\[
u,v\in G\equiv domain(X(\cdot,\omega))\cap(t,t')\cap\{\tau_{n,i}(\omega)\}^{c}
\]
with $u\leq v$. Either (i') $t<\tau_{n,i}(\omega)$, or (ii') $\tau_{n,i}(\omega)<t'$.
Consider Case (i'). Then
\[
t'\equiv t+2^{-m(n)}<\tau_{n,i}(\omega)+2^{-m(n)}
\]
\[
<\tau_{n,i}(\omega)+\delta_{aucl}(2^{-n-1})\equiv\tau_{n,i}(\omega)+\delta_{aucl}(\frac{1}{2}\varepsilon_{n})\leq\tau_{n,i+1}(\omega),
\]
where the second inequality is from inequality \ref{eq:temp-369},
and where the last inequality is due to inequality \ref{eq:temp-400}.
Combining, $\tau_{n,i-1}(\omega)\leq t<t'<\tau_{n,i+1}(\omega)$.
Hence there are three subcases regarding the order of $u,v$ in relation
to the points $\tau_{n,i-1}(\omega),\tau_{n,i}(\omega),$ and $\tau_{n,i+1}(\omega)$
in the interval $[0,1]$: (i'a) 
\[
\tau_{n,i-1}(\omega)\leq t<u\leq v<\tau_{n,i}(\omega),
\]
(i'b)
\[
\tau_{n,i-1}(\omega)\leq t<u<\tau_{n,i}(\omega)<v<t'<\tau_{n,i+1}(\omega),
\]
and (i'c)
\[
\tau_{n,i}(\omega)<u\leq v<t'<\tau_{n,i+1}(\omega).
\]
In Subcase (i'a), we have, in view of inequality \ref{eq:temp-401},
\[
g_{1}(u,v)\equiv d(X_{t}(\omega),X_{u}(\omega))\vee d(X_{t}(\omega),X_{v}(\omega))\leq\varepsilon_{n}<\beta,
\]
where $u,v\in G$ with $u\leq v$ are arbitrary. Since $G$ is dense
in $domain(X(\cdot,\omega))\cap(t,t')$, we therefore obtain, by right
continuity at $r$ and at $s$, 
\[
g_{1}(r,s)\equiv d(X_{t}(\omega),X_{r}(\omega))\vee d(X_{t}(\omega),X_{s}(\omega))\leq\varepsilon_{n}<\beta,
\]
whence $\omega\in(A_{t,r}^{\beta}\cup A_{t,s}^{\beta})^{c}$. Similarly,
in Subcase (i'b), we have 
\[
g_{2}(r,s)\equiv d(X_{t}(\omega),X_{r}(\omega))\vee d(X_{s}(\omega),X_{t'}(\omega))<\beta,
\]
whence $\omega\in(A_{t,r}^{\beta}\cup A_{t',s}^{\beta})^{c}$. Likewise,
in Subcase (i'c), we have 
\[
g_{3}(r,s)\equiv d(X_{t'}(\omega),X_{r}(\omega))\vee d(X_{t'}(\omega),X_{s}(\omega))<\beta,
\]
whence $\omega\in(A_{t',r}^{\beta}\cup A_{t',s}^{\beta})^{c}$ . Combining,
\[
\omega\in(A_{t,r}^{\beta}\cup A_{t,s}^{\beta})^{c}\cup(A_{t,r}^{\beta}\cup A_{t',s}^{\beta})^{c}\cup(A_{t',r}^{\beta}\cup A_{t',s}^{\beta})^{c},
\]
contradicting relation \ref{eq:temp-416}. Thus the assumption that
$P(D_{n})>2^{-n}$ leads to a contradiction. We conclude that $P(D_{n})\leq2^{-n}$,
where $n\geq0$ and $\beta>2^{-n}$ are arbitrary. Thus the process
$Z\equiv X|Q_{\infty}$ satisfies the conditions in Definition \ref{Def.  D-regular process and D-regular family of f.j.d.s on Q_inf}
to be $D$-regular, with the sequence $(m_{n})_{n=0,1,2,\cdots}$
as a modulus of $D$-regularity.
\end{proof}
The converse to Theorem \ref{Thm. Restriction of a.u. cadlag process to dyadics is D-regular}
will be proved in the next section. From a $D$-regular family $F$
of f.j.d.'s with parameter set $[0,1]$ we will construct an a.u.
càdlàg process with marginal distributions given by the family $F$. 

\section{The Right-Limit Extension of $D$-regular Processes area.u. Càdlàg }

Refer to Definition \ref{Def. Notations for dyadic rationals} for
the notations related to the enumerated set $Q_{\infty}$of dyadic
rationals in the interval $[0,1]$. We proved in Theorem \ref{Thm. Restriction of a.u. cadlag process to dyadics is D-regular}
that the restriction to $Q_{\infty}$ of each a.u. càdlàg process
on $[0,1]$ is $D$-regular. In this section, we will prove the more
useful converse theorem, which is the key theorem in this chapter,
that the extension by right limit of each $D$-regular process on
$Q_{\infty}$ is a.u. càdlàg. Then we will prove the easy corollary
that, given an $D$-regular family $F$ of f.j.d.'s with parameter
set $[0,1]$, we can construct an a.u. càdlàg process $X$ with marginal
distributions given by $F$, and with a modulus of a.u. càdlàg in
terms of the modulus of $D$-regularity of $F$. 

In the remainder of the section, we will use the following assumption
and notations.
\begin{defn}
\label{Def. Assumption of D-regulat process} \textbf{(Assumption
of a $D$-regular process).} Let $Z\in\widehat{R}_{Dreg,\overline{m},\delta(Cp)}(Q_{\infty}\times\Omega,S)$
be arbitrary, but fixed for the remainder of this section. In other
words, $Z:Q_{\infty}\times\Omega\rightarrow S$ is a fixed $D$-regular
process. Let $\overline{m}\equiv(m_{k})_{k=0,1,\cdots}$ be a fixed
modulus of $D$-regularity, and let $\delta_{Cp}$ be a fixed modulus
of continuity in probability. 
\end{defn}
$\square$
\begin{defn}
\textbf{\label{De. Notation for range of a sample function} (Notation
for the range of a sample function). }Let $Y:Q\times\Omega\rightarrow S$
be an arbitrary process, let $A\subset Q$ and $\omega\in\Omega$
be arbitrary. Then we write
\[
Y(A,\omega)\equiv\{x\in S:x=Y(t,\omega)\textrm{ for some }t\in A\}.
\]
\end{defn}
$\square$
\begin{defn}
\textbf{\label{Def. Accordian function} (Accordion function).} Let
$(\beta_{h})_{h=0,1,\cdots}$ be an arbitrary sequence of real numbers
such that, for each $k,h\geq0$ with $k\leq h$, and for each $r,s\in Q_{\infty},$
we have
\begin{equation}
\beta_{h}\in(2^{-h+1},2^{-h+2}),\label{eq:temp-420-1}
\end{equation}
 and the set
\begin{equation}
(d(Z_{r},Z_{s})>\beta_{k}+\cdots+\beta_{h})\label{eq:temp-304-1-1}
\end{equation}
is $\mathrm{measurable}$, and, in particular, the set
\begin{equation}
A_{r,s}^{\beta(h)}\equiv(d(Z_{r},Z_{s})>\beta_{h})\label{eq:temp-246}
\end{equation}
is $\mathrm{measurable}$. 

Let $h,n\geq0$ be arbitrary. Define
\begin{equation}
\beta_{n,h}\equiv\sum_{i=n}^{h}\beta_{i},\label{eq:temp-375}
\end{equation}
where, by convention, $\sum_{i=n}^{h}\beta_{i}\equiv0$ if $h<n$.
Define
\begin{equation}
\beta_{n,\infty}\equiv\sum_{i=n}^{\infty}\beta_{i}.\label{eq:temp-376}
\end{equation}
Note that $\beta_{n,\infty}<\sum_{i=n}^{\infty}2^{-i+2}=2^{-n+3}\rightarrow0$
as $n\rightarrow\infty$. For each subset $A$ of $Q_{\infty}$, write
\[
A^{-}\equiv\{t\in Q_{\infty}:t\neq s\mbox{ for each }s\in A\},
\]
the metric complement of $A$ in $Q_{\infty}$. 

Let $s\in Q_{\infty}=\bigcup_{h=0}^{\infty}Q_{m(h)}$ be arbitrary.
Define $\widehat{h}(s)\equiv h\geq0$ to be the smallest integer such
that $s\in Q_{m(h)}$. Let $n\geq0$ be arbitrary. Define 
\begin{equation}
\widehat{\beta}_{n}(s)\equiv\beta_{n,\widehat{h}(s)}\equiv\sum_{i=n}^{\widehat{h}(s)}\beta_{i}.\label{eq:temp-422}
\end{equation}
Thus we have the functions
\[
\widehat{h}:Q_{\infty}\rightarrow\{0,1,2,\cdots\}
\]
and
\[
\widehat{\beta}_{n}:Q_{\infty}\rightarrow(0,\beta_{n,\infty})
\]
for each $n\geq0$. These functions are defined relative to the sequences
$(\beta_{n})_{n=0,1,\cdots}$ and $(m_{n})_{n=0,1,\cdots}$.

For want of a better name, we might call the function $\widehat{\beta}_{n}$
an \index{accordion function}\emph{accordion function}, because its
graph resembles a fractal-like accordion. They will furnish time-varying
boundaries for some simple first exit times in the proof of the main
theorem. Note that, for arbitrary $s\in Q_{m(k)}$ for some $k\geq0$,
we have $\widehat{h}(s)\leq k$ and so $\widehat{\beta}_{n}(s)\leq\beta_{n,k}$. 
\end{defn}
$\square$
\begin{defn}
\label{Def. Exceptional sets D_n} \textbf{(Some small exceptional
sets).} Let $k\geq0$ be arbitrary. Then $\beta_{k+1}>2^{-k}$. Hence,
by the conditions in Definition \ref{Def.  D-regular process and D-regular family of f.j.d.s on Q_inf}
for $\overline{m}\equiv(m_{h})_{h=0,1,\cdots}$ to be a modulus of
$D$-regularity of the process $Z$, we have 
\begin{equation}
P(D_{k})\leq2^{-k},\label{eq:temp-204-3-2}
\end{equation}
where
\[
D_{k}\equiv\bigcup_{u\in Q(m(k))}\bigcup_{r,s\in(u,u')Q(m(k+1));r\leq s}\{
\]
\begin{equation}
(A_{u,r}^{\beta(k+1)}\cup A_{u,s}^{\beta(k+1)})(A_{u,r}^{\beta(k+1)}\cup A_{u',s}^{\beta(k+1)})(A_{u',r}^{\beta(k+1)}\cup A_{u',s}^{\beta(k+1)})\},\label{eq:temp-201-4-3-3}
\end{equation}
where, for each $u\in Q{}_{m(k)}$, we abuse notations and write $u'\equiv u+\Delta_{m(k)}.$ 

For each $h\geq0$, define the small exceptional set
\begin{equation}
D_{h+}\equiv\bigcup_{k=h}^{\infty}D_{k},\label{eq:temp-358-1}
\end{equation}
with
\begin{equation}
P(D_{h+})\leq\sum_{k=h}^{\infty}2^{-k}=2^{-h+1}.\label{eq:temp-424}
\end{equation}
\end{defn}
$\square$
\begin{lem}
\label{Lem. Existence of supremum distance Z} \textbf{\emph{(Existence
of certain supremums as r.r.v.'s).}} Let $Z:Q_{\infty}\times\Omega\rightarrow S$
be a $D$-regular process, with a modulus of $D$-regularity $\overline{m}\equiv(m_{k})_{k=0,1,\cdots}$
. Let $h\geq0$ and $v,\overline{v},v'\in Q_{m(h)}$ be arbitrary.
with $v\leq\overline{v}\leq v'$. Then the following holds.

1. For each $r\in[v,v']Q_{\infty}$, we have 
\begin{equation}
d(Z_{\overline{v}},Z_{r})\leq\bigvee_{u\in[v,v']Q(m(h))}d(Z_{\overline{v}},Z_{u})+\widehat{\beta}_{h+1}(r)\label{eq:temp-366}
\end{equation}
on $D_{h+}^{c}$.

2. The supremum
\[
Y_{v,v'}\equiv\sup_{u\in[v,v']Q(\infty)}d(Z_{v},Z_{u})
\]
exists as a r.r.v. Moreover
\[
0\leq Y_{v,v'}-\bigvee_{u\in[v,v']Q(m(h))}d(Z_{v},Z_{u})\leq\beta_{h+1,\infty}\leq2^{-h+4}
\]
on $D_{h+}^{c}$, where $P(D_{h+})\leq2^{-h+1}.$

3. Write $\widehat{d}\equiv1\wedge d$. Then
\[
0\leq E\sup_{u\in[v,v']Q(\infty)}\widehat{d}(Z_{v},Z_{u})-E\bigvee_{u\in[v,v']Q(m(h))}\widehat{d}(Z_{v},Z_{u})\leq2^{-h+5}.
\]
where $h\geq0$ and $v,v'\in Q_{m(h)}$ are arbitrary with $v\leq v'$.
\end{lem}
\begin{proof}
1. First let $k\geq0$ be arbitrary. Consider each $\omega\in D_{k}^{c}$.
We will prove that
\begin{equation}
0\leq\bigvee_{u\in[v,v']Q(m(k+1))}d(Z_{\overline{v}}(\omega),Z_{u}(\omega))-\bigvee_{u\in[v,v']Q(m(k))}d(Z_{\overline{v}}(\omega),Z_{u}(\omega))\leq\beta_{k+1}.\label{eq:temp-312}
\end{equation}
Consider each $r\in[v,v']Q_{m(k+1)}$. Then $r\in[s,s+\Delta_{m(k)}]Q_{m(k+1)}$
for some $s\in Q_{m(k)}$ such that $[s,s+\Delta_{m(k)}]\subset[v,v']$.
Write $s'\equiv s+\Delta_{m(k)}$. We need to prove that 
\begin{equation}
d(Z_{\overline{v}}(\omega),Z_{r}(\omega))\leq\bigvee_{u\in[v,v']Q(m(k))}d(Z_{\overline{v}}(\omega),Z_{u}(\omega))+\beta_{k+1}.\label{eq:temp-249}
\end{equation}
If $r=s$ or $r=s'$, then $r\in[v,v']Q_{m(k)}$ and inequality \ref{eq:temp-249}
holds trivially. Hence we may assume that $r\in(s,s')Q_{m(k+1)}$.
Since $\omega\in D_{k}^{c}$ by assumption, the defining equality
\ref{eq:temp-201-4-3-3} implies that 
\[
\omega\in(A_{s,r}^{\beta(k+1)})^{c}(A_{s,r}^{\beta(k+1)})^{c}\cup(A_{s,r}^{\beta(k+1)})^{c}(A_{s',r}^{\beta(k+1)})^{c}\cup(A_{s',r}^{\beta(k+1)})^{c}(A_{s',r}^{\beta(k+1)})^{c}.
\]
Consequently, by the defining equality \ref{eq:temp-246} for the
sets in the last displayed expression, we have
\[
d(Z_{s}(\omega),Z_{r}(\omega))\wedge d(Z_{s'}(\omega),Z_{r}(\omega))\leq\beta_{k+1}.
\]
Hence the triangle inequality implies that
\[
d(Z_{\overline{v}}(\omega),Z_{r}(\omega))
\]
\[
\leq(d(Z_{s}(\omega),Z_{r}(\omega))+d(Z_{\overline{v}}(\omega),Z_{s}(\omega)))\wedge(d(Z_{s'}(\omega),Z_{r}(\omega))+d(Z_{\overline{v}}(\omega),Z_{s'}(\omega)))
\]
\[
\leq(d(Z_{s}(\omega),Z_{r}(\omega))+\bigvee_{u\in[v,v']Q(m(k))}d(Z_{\overline{v}}(\omega),Z_{u}(\omega)))\wedge(d(Z_{s'}(\omega),Z_{r}(\omega))
\]
\[
+\bigvee_{u\in[v,v']Q(m(k))}d(Z_{\overline{v}}(\omega),Z_{u}(\omega)))
\]
\[
\leq(\beta_{k+1}+\bigvee_{u\in[v,v']Q(m(k))}d(Z_{\overline{v}}(\omega),Z_{u}(\omega)))\wedge(\beta_{k+1}+\bigvee_{u\in[v,v']Q(m(k))}d(Z_{\overline{v}}(\omega),Z_{u}(\omega)))
\]
\[
\leq\beta_{k+1}+\bigvee_{u\in[v,v']Q(m(k))}d(Z_{\overline{v}}(\omega),Z_{u}(\omega)),
\]
establishing inequality \ref{eq:temp-249} for arbitrary $r\in[v,v']Q_{m(k+1)}$.
The desired inequality \ref{eq:temp-312} follows.

2. Let $h\geq0$ be arbitrary, and consider each $\omega\in D_{h+}^{c}$.
Then $\omega\in D_{k}^{c}$ for each $k\geq h$. Hence, inequality
\ref{eq:temp-312} can be applied to $h,h+1,\cdots,k+1$ to yield
\[
0\leq\bigvee_{u\in[v,v']Q(m(k+1))}d(Z_{\overline{v}}(\omega),Z_{u}(\omega))-\bigvee_{u\in[v,v']Q(m(h))}d(Z_{\overline{v}}(\omega),Z_{u}(\omega))
\]
\begin{equation}
\leq\beta_{k+1}+\cdots+\beta_{h+1}=\beta_{h+1,k+1}<\beta_{h+1,\infty}<2^{-h+2}.\label{eq:temp-326}
\end{equation}

3. Consider each $r\in[v,v']Q_{\infty}$. We will prove that 
\begin{equation}
d(Z_{\overline{v}}(\omega),Z_{r}(\omega))\leq\bigvee_{u\in[v,v']Q(m(h))}d(Z_{\overline{v}}(\omega),Z_{u}(\omega))+\widehat{\beta}_{h+1}(r).\label{eq:temp-367}
\end{equation}
The desired inequality is trivial if $r\in Q_{m(h)}$. Hence we may
assume that $r\in Q_{m(k+1)}Q_{m(k)}^{-}$ for some $k\geq h$. Then
$\widehat{\beta}_{h}(r)=\beta_{h+1,k+1}$. Therefore the first half
of inequality \ref{eq:temp-326} implies that

\[
d(Z_{\overline{v}}(\omega),Z_{r}(\omega))\leq\bigvee_{u\in[v,v']Q(m(h))}d(Z_{\overline{v}}(\omega),Z_{u}(\omega))+\beta_{h+1,k+1}
\]
\[
=\bigvee_{u\in[v,v']Q(m(h))}d(Z_{\overline{v}}(\omega),Z_{u}(\omega))+\widehat{\beta}_{h+1}(r).
\]
Inequality \ref{eq:temp-367} is proved, where $\omega\in D_{h+}^{c}$
is arbitrary. Inequality \ref{eq:temp-366} follows. Assertion 1 is
proved.

4. Next consider the special case where $\overline{v}=v$. Since $P(D_{h+})\leq2^{-h+1}$,
it follows from inequality \ref{eq:temp-326} that the a.u. limit
\[
\overline{Y}_{v,v'}\equiv\lim_{k\rightarrow\infty}\bigvee_{u\in[v,v']Q(m(k+1))}d(Z_{v},Z_{u})
\]
exists and is a r.r.v. Moreover, for each $\omega\in domain(\overline{Y}_{v,v'})$,
it is easy to verify that $\overline{Y}_{v,v'}(\omega)$ gives the
supremum $\sup_{u\in[v,v']Q(\infty)}d(Z_{v}(\omega),Z_{u}(\omega))$.
Thus the supremum $Y_{v,v'}\equiv\sup_{u\in[v,v']Q(\infty)}d(Z_{v},Z_{u})$
is defined and equal to the r.r.v. $\overline{Y}_{v,v'}$ on a full
set, and is therefore itself a r.r.v. Letting $k\rightarrow\infty$
in inequality \ref{eq:temp-326}, we obtain 
\[
0\leq Y_{v,v'}-\bigvee_{u\in[v,v']Q(m(h))}d(Z_{v},Z_{u})\leq\beta_{h+1,\infty}\leq2^{-h+4}
\]
on $D_{h+}^{c}\cap domain(\overline{Y}_{v,v'})$. This proves Assertion
2.

3. Write $\widehat{d}\equiv1\wedge d$. Then
\[
0\leq E\sup_{u\in[v,v']Q(\infty)}\widehat{d}(Z_{v},Z_{u})-E\bigvee_{u\in[v,v']Q(m(h))}\widehat{d}(Z_{v},Z_{u})
\]
\[
=E(1\wedge Y_{v,v'}-1\wedge\bigvee_{u\in[v,v']Q(m(h))}d(Z_{v},Z_{u}))
\]
\[
\leq E1_{D(h+)^{c}}(1\wedge Y_{v,v'}-1\wedge\bigvee_{u\in[v,v']Q(m(h))}d(Z_{v},Z_{u}))+E1_{D(h+)}.
\]
\[
\leq2^{-h+5}+P(D_{h+})\leq2^{-h+4}+2^{-h+1}<2^{-h+5}.
\]
Assertion 3 and the lemma is proved.
\end{proof}
\begin{defn}
\label{Def. Extension of process by right limit} \textbf{(Right-limit
extension of a process with dyadic rational parameters).} Recall the
convention that if $f$ is an arbitrary function, we write $f(x)$
only with the implicit or explicit condition that $x\in domain(f)$. 

1. Let $Q_{\infty}$ stand for the set of dyadic rationals in $[0,1]$.
Let $Y:Q_{\infty}\times\Omega\rightarrow S$ be an arbitrary process.
Define a function $X:[0,1]\times\Omega\rightarrow S$ by 
\[
domain(X)\equiv\{(r,\omega)\in[0,1]\times\Omega:\lim_{u\rightarrow r;u\geq r}Y(u,\omega)\:\mathrm{exists}\},
\]
and by 
\begin{equation}
X(r,\omega)\equiv\lim_{u\rightarrow r;u\geq r}Y(u,\omega)\label{eq:temp-332-1-1}
\end{equation}
for each $(r,\omega)\in domain(X)$. We will call 
\[
\Phi_{rLim}(Y)\equiv X
\]
the right-limit extension of the process $Y$ \index{right-limit extension of process}
to the parameter set $[0,1]$. 

2. Let $\overline{Q}_{\infty}$ stand for the set of dyadic rationals
in $[0,\infty)$. Let $Y:\overline{Q}_{\infty}\times\Omega\rightarrow S$
be an arbitrary process. Define a function $X:[0,\infty)\times\Omega\rightarrow S$
by 
\[
domain(X)\equiv\{(r,\omega)\in[0,\infty)\times\Omega:\lim_{u\rightarrow r;u\geq r}Y(u,\omega)\:\mathrm{exists}\},
\]
and by 
\begin{equation}
X(r,\omega)\equiv\lim_{u\rightarrow r;u\geq r}Y(u,\omega)\label{eq:temp-332-1-1-2}
\end{equation}
for each $(r,\omega)\in domain(X)$. We will call 
\[
\overline{\Phi}_{rLim}(Y)\equiv X
\]
the right-limit extension of the process $Y$ \index{right-limit extension of process}
to the parameter set $[0,\infty)$. 

In general, the right-limit extension $X$ need not be a well-defined
process. $\square$
\end{defn}
In the following proposition, recall that, as in Definition \ref{Def. continuity in prob, continuity a.u., and a.u. continuity},
continuity a.u. is a weaker condition than a.u. continuity.
\begin{prop}
\textbf{\emph{\label{Prop. Right-limit extenaion of D-regular process is continuous a.u.}
(The right-limit extension of a }}$D$\textbf{\emph{-regular process
is a well-defined stochastic process and is continuous a.u.).}} The
right-limit extension $X\equiv\Phi_{rLim}(Z):[0,1]\times\Omega\rightarrow S$
of the $D$-regular process $Z$ is a stochastic process which is
continuous a.u. 

Specifically, the following holds.

1. Let $\varepsilon>0$ be arbitrary. Take $\nu\geq0$ so large that
$2^{-\nu+5}<\varepsilon$. Take an arbitrary $J\geq0$ so large that
\begin{equation}
\Delta_{m(J)}\equiv2^{-m(J)}<2^{-2}\delta_{Cp}(2^{-2\nu+2}).\label{eq:temp-332}
\end{equation}
Define 
\begin{equation}
\delta_{cau}(\varepsilon)\equiv\delta_{cau}(\varepsilon,\overline{m},\delta_{Cp})\equiv\Delta_{m(J)}\in Q_{\infty}.\label{eq:temp-334}
\end{equation}
Then, for each $t\in[0,1]$, there exists an exceptional set $G_{t}$
with $P(G_{t})<\varepsilon$ such that
\[
d(X(t,\omega),X(t',\omega))<\varepsilon
\]
for each $t'\in[t-\delta_{cau}(\varepsilon),t+\delta_{cau}(\varepsilon)]\cap domain(X(\cdot,\omega))$,
for each $\omega\in G_{t}^{c}$.

2. Moreover, $X(\cdot,\omega)|Q_{\infty}=Z(\cdot,\omega)$ for a.e.
$\omega\in\Omega$. 

3. Furthermore, the process $X$ has the same modulus of continuity
in probability $\delta_{Cp}$ as the process $Z$. 

4. (Right continuity). For a.e. $\omega\in\Omega$, the function $X(\cdot,\omega)$
is right continuous at each $t\in domain(X(\cdot,\omega))$.

5. (Right completeness). For a.e. $\omega\in\Omega$, we have $t\in domain(X(\cdot,\omega))$
for each $t\in[0,1]$ such that $\lim_{r\rightarrow t;r>t}X(r,\omega)$
exists.
\end{prop}
\begin{proof}
We will first verify that $X_{t}\equiv X(t,\cdot)$ is a r.v. with
values in $S$, for each $t\in[0,1]$, then prove the continuity a.u.
To that end, let $t\in[0,1]$ and $\varepsilon>0$ be arbitrary. Write
$\Delta\equiv\Delta_{m(J)}$. For ease of reference to previous theorems,
write $n\equiv\nu$.

1. Recall that $\delta_{Cp}$ is the given modulus of continuity in
probability of the process $Z$, and recall that $\beta_{n}\in(2^{-n+1},2^{-n+2})$
is as selected in Definition \ref{Def. Accordian function}. When
there is no risk of confusion, suppress the subscript $m_{J}$, write
$p\equiv p_{m(J)}\equiv2^{m(J)}$, and write $q_{i}\equiv q_{m(J),i}\equiv i\Delta\equiv i2^{-m(J)}$
for each $i=0,1,\cdots,p$. Then 
\[
t\in[0,1]=[q_{0},q_{2}]\cup(q_{1},q_{3}]\cup\cdots\cup[q_{p-3},q_{p-1}]\cup[q_{p-2},q_{p}].
\]
Hence there exists $i=0,\cdots,p-2$ be such that $t\in[q_{i},q_{i+2}]$.
The neighborhood $\theta_{t,\varepsilon}\equiv[t-\Delta,t+\Delta]\cap[0,1]$
of $t$ in $[0,1]$ is a subset of $[q_{(i-1)\vee0},q_{(i+3)\wedge p}]$.
Write $v\equiv q_{(i-1)\vee0}$ and $v'\equiv q_{(i+3)\wedge p}$.
Then (i) $v,v'\in Q_{m(J)}$, (ii) $v<v'$, and (iii) the set 
\[
[v,v']Q_{m(J)}=\{q_{(i-1)\vee0},q_{i},q_{i+1},q_{i+2},q_{(i+3)\wedge p}\}
\]
contains 4 or 5 distinct and consecutive elements of $Q_{m(J)}$.
Therefore, for each $u\in[v,v']Q_{m(J)}$, we have $|v-u|\leq4\Delta<\delta_{Cp}(2^{-2n+2})$,
whence
\[
E1\wedge d(Z_{v},Z_{u}))\leq2^{-2n+2}<\beta_{n}^{2},
\]
and, by Chebychev's inequality,
\[
P(d(Z_{v},Z_{u})>\beta_{n})\leq\beta_{n}.
\]
Hence the $\mathrm{measurable}$ set
\[
A_{n}\equiv\bigcup_{u\in[v,v']Q(m(J))}(d(Z_{v},Z_{u})>\beta_{n})
\]
has probability bounded by $P(A_{n})\leq4\beta_{n}<2^{-n+4.}$. Define
$G_{t,\varepsilon}\equiv D_{n+}\cup A_{n}$. It follows that 
\begin{equation}
P(G_{t,\varepsilon})\leq P(D_{n+})+P(A_{n})<2^{-n+4}+2^{-n+4}=2^{-n+5}<\varepsilon.\label{eq:temp-333}
\end{equation}

2. Next consider each $\omega\in G_{t,\varepsilon}^{c}=D_{n+}^{c}A_{n}^{c}.$
Then, by the definition of the set $A_{n}$, we have
\[
\bigvee_{s\in[v,v']Q(m(J))}d(Z(v,\omega),Z(s,\omega))\leq\beta_{n}.
\]
At the same time, since $\omega\in D_{n+}^{c}\subset D_{J+}^{c}$,
inequality \ref{eq:temp-366} of Lemma \ref{Lem. Existence of supremum distance Z},
where $h,\overline{v}$ are replaced by $J,v$ respectively, implies
that, for each $r\in\theta_{n}Q_{\infty}\subset[v,v']Q_{\infty}$,
we have
\[
d(Z_{r}(\omega),Z_{v}(\omega))\leq\widehat{\beta}_{J+1}(r)+\bigvee_{s\in[v,v']Q(m(J))}d(Z_{v}(\omega),Z_{s}(\omega))
\]
\[
=\widehat{\beta}_{J+1}(r)+\beta_{n}\leq\beta_{J+1,\infty}+\beta_{n}\leq\beta_{n+1,\infty}+\beta_{n}=\beta_{n,\infty}.
\]
Then the triangle inequality yields 
\begin{equation}
d(Z(r,\omega),Z(r',\omega))\leq2\beta_{n,\infty}<2^{-n+4}<\varepsilon\label{eq:temp-309}
\end{equation}
for each $r,r'\in\theta_{t,\varepsilon}Q_{\infty}\equiv[t-\Delta_{m(J)},t+\Delta_{m(J)}]Q_{\infty}$,
where $\omega\in G_{t,\varepsilon}^{c}$ is arbitrary, where $\varepsilon>0$
is arbitrarily.

3. Now write $\varepsilon_{k}\equiv2^{-k}$ for each $k\geq1$. Then
$G_{\kappa}\equiv\bigcup_{k=\kappa}^{\infty}G_{t,\varepsilon(k)}$
has probability bounded by $P(G_{\kappa})\leq2^{-\kappa+1}$. Hence
$H_{t}\equiv\bigcap_{\kappa=1}^{\infty}G_{\kappa}$ is a null set.
Moreover, for each $\omega\in H_{t}^{c}$, we have $\omega\in G_{\kappa}^{c}$
for some $\kappa\geq1$. Therefore the $\lim_{u\rightarrow t}Z(u,\omega)$
exists, whence the right limit $X(t,\omega)$ is well defined and
\begin{equation}
X(t,\omega)\equiv\lim_{u\rightarrow t;u\geq t}Z(u,\omega)=\lim_{u\rightarrow t}Z(u,\omega).\label{eq:temp-332-1-1-1}
\end{equation}
Applied to the special case of an arbitrary $u\in Q_{\infty}$ in
the place of $t$, equality \ref{eq:temp-332-1-1-1} implies that
$X(u,\omega)=Z(u,\omega)$ for each $\omega\in H_{u}^{c}$. Hence
the same equality \ref{eq:temp-332-1-1-1} can be rewritten as
\[
X(t,\omega)\equiv\lim_{u\rightarrow t;u\geq t}X(u,\omega)=\lim_{u\rightarrow t}X(u,\omega)
\]
for each $\omega\in H_{t}^{c}$. Define the null set $H\equiv\bigcup_{u\in Q(\infty)}H_{u}$.
Then $X(u,\omega)=Z(u,\omega)$ for each $u\in Q_{\infty},$ for each
$\omega\in H^{c}$. In short, $X(\cdot,\omega)|Q_{\infty}=Z(\cdot,\omega)$
for each $\omega\in H^{c}$. This verifies Assertion 2. 

4. For each $k\geq1$, fix an arbitrary $r_{k}\in\theta_{t,\varepsilon(k)}Q_{\infty}$.
Consider each $\kappa\geq1$, and each $\omega\in G_{\kappa}^{c}$.
Then, for each $k\geq\kappa$, we have $r_{k},r_{\kappa}\in\theta_{t,\varepsilon(\kappa)}Q_{\infty},$
and so, by inequality \ref{eq:temp-309}, where $\varepsilon$ is
replaced by $\varepsilon_{\kappa}$, we obtain
\begin{equation}
d(Z(r_{k},\omega),Z(r_{\kappa},\omega))\leq\varepsilon_{\kappa}.\label{eq:temp-309-1}
\end{equation}
 Letting $k\rightarrow\infty$, this yields
\begin{equation}
d(X(t,\omega),Z(r_{\kappa},\omega))\leq\varepsilon_{\kappa},\label{eq:temp-309-1-1}
\end{equation}
where $\omega\in G_{\kappa}^{c}$ is arbitrary. Since $P(G_{\kappa})\leq2^{-\kappa+1}$
is arbitrarily small for sufficiently large $\kappa$, we conclude
that $Z_{r(\kappa)}\rightarrow X_{t}$ a.u. Consequently, the function
$X_{t}$ is a r.v. Since $t\in[0,1]$ is arbitrary, the function $X\equiv\Phi_{rLim}(Z):[0,1]\times\Omega\rightarrow S$
is a stochastic process.

5. Let $t'\in[t-\Delta_{m(J)},t+\Delta_{m(J)}]\cap domain(X(\cdot,\omega))$
be arbitrary. Letting $r\downarrow t$ and $r'\downarrow t'$ in inequality
\ref{eq:temp-309} while $r,r'\in\theta_{t,\varepsilon}Q_{\infty}$,
we obtain 
\begin{equation}
d(X(t,\omega),X(t',\omega))<\varepsilon\label{eq:temp-309-2}
\end{equation}
where $\omega\in G_{t,\varepsilon}^{c}$ is arbitrary. Since $t\in[0,1]$
is arbitrary, and since $P(G_{t,\varepsilon})<\varepsilon$ is arbitrarily
small, we see that the process $X$ is continuous a.u. according to
Definition \ref{Def. continuity in prob, continuity a.u., and a.u. continuity}.
Assertions 1 has been proved.

6. Next, we will verify that the process $X$ has the same modulus
of continuity in probability $\delta_{Cp}$ as the process $Z$. To
that end, let $\varepsilon>0$ be arbitrary, and let $t,s\in[0,1]$
be such that $|t-s|<\delta_{Cp}(\varepsilon)$. In Step 5 we saw that
there exist sequences $(r_{k})_{k=1,2,\cdots}$ and $(v_{k})_{k=1,2,\cdots}$in
$Q_{\infty}$ such that $r_{k}\rightarrow t$, $v_{k}\rightarrow s$,
$Z_{r(k)}\rightarrow X_{t}$ a.u. and $Z_{v(k)}\rightarrow X_{s}$
a.u. Then, for sufficiently large $k\geq0$ we have $|r_{k}-v_{k}|<\delta_{Cp}(\varepsilon)$,
whence $E1\wedge d(Z_{r(k)},Z_{v(k)})\leq\varepsilon$. The last cited
a.u. convergence therefore implies that 
\[
E1\wedge d(X_{t},X_{s})\leq\varepsilon.
\]
Summing up, $\delta_{Cp}$ is a modulus of continuity of probability
of the process $X$. Assertion 3 is proved.

7. Consider each $\omega$ in the full set $\bigcap_{u\in Q(\infty)}domain(Z_{u})$.
Then the function $Z(\cdot,\omega):Q_{\infty}\rightarrow S$ has domain
$Q_{\infty}$ which is dense in$[0,1]$. Hence its right-limit extension
$X(\cdot,\omega):[0,1]\rightarrow S$ satisfies the conditions in
Assertions 4 and 5, according to Proposition \ref{Prop. Cadlag completion}.
The present proposition is proved.
\end{proof}
We now prove the main theorem of this chapter. In the proof, refer
to Proposition \ref{Prop. Basics of simple Exit times} for basic
properties of simple first exit times.
\begin{thm}
\label{Thm. Extension of D-regular process by right limit is a.u.cadlag}
\textbf{\emph{(The right-limit extension of a }}$D$\textbf{\emph{-regular
process is a.u. càdlàg). }}The right-limit extension 
\[
X\equiv\Phi_{rLim}(Z):[0,1]\times\Omega\rightarrow S
\]
is a.u. càdlàg. Specifically, \emph{(i) }it has the same  modulus
of continuity in probability $\delta_{Cp}$ as the given $D$-regular
process $Z$, and \emph{(ii)} it has a modulus of a.u. càdlàg $\delta_{aucl}(\cdot,\overline{m},\delta_{Cp})$
defined as follows. Let $\varepsilon>0$ be arbitrary. Let $n\geq0$
be so large that $2^{-n+6}<\varepsilon.$ Let $J\geq n+1$ be so large
that 
\begin{equation}
\Delta_{m(J)}\equiv2^{-m(J)}<2^{-2}\delta_{Cp}(2^{-2m(n)-2n-10}).\label{eq:temp-332-1-2}
\end{equation}
Define
\[
\delta_{aucl}(\varepsilon,\overline{m},\delta_{Cp})\equiv\Delta_{m(J)}.
\]
Note that the operation $\delta_{aucl}(\cdot,\overline{m},\delta_{Cp})$
depends only on $\overline{m}$ and $\delta_{Cp}$.
\end{thm}
\begin{proof}
We will prove that the operation $\delta_{aucl}(\cdot,\overline{m},\delta_{Cp})$
thus defined is a modulus of a.u. càdlàg of the right-limit extension
$X$. 

1. To that end, let $i=1,\cdots,p_{m(n)}$ be arbitrary, but fixed
until further notice. When there is little risk of confusion, we suppress
references to $n$ and $i$, and write 
\[
p\equiv p_{m(n)}\equiv2^{m(n)},
\]
\[
\Delta\equiv\Delta_{m(n)}\equiv2^{-m(n)},
\]
\[
t\equiv q_{i-1}\equiv q_{m(n),i-1}\equiv(i-1)2^{-m(n)},
\]
and
\[
t'\equiv q_{i}\equiv q_{m(n),i}\equiv i2^{-m(n)}.
\]
Thus $t,t'\in Q_{m(n)}$ and $0\leq t<t'=t+\Delta\leq1.$ 

2. Write $\varepsilon_{n}=2^{-m(n)-n-1}$ and $\nu=m_{n}+n+6$. Then
$2^{-\nu+5}<\varepsilon_{n}$, and 
\begin{equation}
\Delta_{m(J)}\equiv2^{-m(J)}<2^{-2}\delta_{Cp}(2^{-2\nu+2}).\label{eq:temp-332-1-2-3}
\end{equation}
Hence
\[
\Delta_{m(J)}=\delta_{cau}(\varepsilon_{n},\overline{m},\delta_{Cp}),
\]
where $\delta_{cau}(\cdot,\overline{m},\delta_{Cp})$ is the modulus
of continuity a.u. of the process $X$ defined in Proposition \ref{Prop. Right-limit extenaion of D-regular process is continuous a.u.}.
Note also that $\varepsilon_{n}\leq\beta_{n}$.

In the next several steps, we will prove that, for each $\omega\in D_{n+}^{c}$,
the set $Z([t,t']Q_{\infty},\omega)$ can be divided into two sections
$Z([t,\tau)Q_{\infty},\omega)$ and $Z([\tau,t']Q_{\infty},\omega)$
each of which is contained in a ball in $(S,d)$ with radius $2^{-n+5}$. 

2. First introduce some simple first exit times of the process $Z$.
Let $k\geq n$ be arbitrary. As in Definition \ref{Def. Simple First Exit time},
define the simple first exit time
\[
\eta_{k}\equiv\eta_{k,i}\equiv\eta_{t,\widehat{\beta}(n),[t,t']Q(m(k))}
\]
for the process $Z|[t,t']Q_{m(k)}$ to exit the time-varying $\widehat{\beta}_{n}$-neighborhood
of $Z_{t}$. In particular, the r.r.v. $\eta_{k}$ has values in $[t+\Delta_{m(k)},t']Q_{m(k)}$.
Thus 
\begin{equation}
t+\Delta_{m(k)}\leq\eta_{k}\leq t'.\label{eq:temp-310}
\end{equation}
In the case where $k=n,$ this yields
\begin{equation}
\eta_{n}\equiv\eta_{n,i}=t'=q_{m(n),i}.\label{eq:temp-336}
\end{equation}
Since $Q_{m(k)}\subset Q_{m(k+1)}$, the more frequently sampled simple
first exit time $\eta_{k+1}$ comes no later than $\eta_{k}$, according
to Assertion 5 of Proposition \ref{Prop. Basics of simple Exit times}.
Thus
\begin{equation}
\eta_{k+1}\leq\eta_{k}.\label{eq:temp-314}
\end{equation}

3. Now let $\kappa\geq k\geq n$ and $\omega\in D_{k+}^{c}$ be arbitrary.
We will prove that 
\begin{equation}
t\leq\eta_{k}(\omega)-2^{-m(k)}\leq\eta_{\kappa}(\omega)-2^{-m(\kappa)}\leq\eta_{\kappa}(\omega)\leq\eta_{k}(\omega).\label{eq:temp-406}
\end{equation}
The first of these inequalities is from the first part of inequality
\ref{eq:temp-310}. The third is trivial, and the last is by a repeated
application of inequality \ref{eq:temp-314}. It remains only to prove
the second.

To that end, write $v\equiv t,$ $s\equiv\eta_{k}(\omega)$ and $v'\equiv s-\Delta_{m(k)}$.
Then $v,t,s,v'\in[t,t']Q_{m(k)}$. Since $v'<s\equiv\eta_{k}(\omega)$,
the sample path $Z(\cdot,\omega)|[t,t']Q_{m(k)}$ has not exited the
time varying $\widehat{\beta}_{n}$-neighborhood of $Z_{t}(\omega)$
at time $v'$. In other words, 
\begin{equation}
d(Z_{t}(\omega),Z_{u}(\omega))\leq\widehat{\beta}_{n}(u)\leq\beta_{n,k}\label{eq:temp-299-1-1}
\end{equation}
for each $u\in[v,v']Q_{m(k)}$. At the same time, $\omega\in D_{n+}^{c}\subset D_{k+}^{c}$.
Hence inequality \ref{eq:temp-366} of Lemma \ref{Lem. Existence of supremum distance Z},
where $h,v,\overline{v}$ are replaced by $k,t,t$ respectively, implies
that 
\[
d(Z_{t}(\omega),Z_{r}(\omega))\leq\widehat{\beta}_{k+1}(r)+\bigvee_{u\in[t,v']Q(m(k))}d(Z_{t}(\omega),Z_{u}(\omega))
\]
\begin{equation}
\leq\widehat{\beta}_{k+1}(r)+\beta_{n,k}=\widehat{\beta}_{n}(r).\label{eq:temp-329-1}
\end{equation}
for each $r\in[t,v']Q_{\infty}$. In particular, for each $r\in[t,v']Q_{m(\kappa+1)}\subset[v,v']Q_{\infty}$,
inequality \ref{eq:temp-329-1} holds. Thus the sample path $Z(\cdot,\omega)|[t,t']Q_{m(\kappa+1)}$
has not exited the time-varying $\widehat{\beta}_{n}$-neighborhood
of $Z(t,\omega)$ at time $v'\equiv\eta_{k}(\omega)-\Delta_{m(k)}$.
In other words,
\[
\eta_{k}(\omega)-\Delta_{m(k)}<\eta_{\kappa+1}(\omega).
\]
Since both sides of this strict inequality are members of $Q_{m(k+1)}$,
it follows that 
\[
\eta_{k}(\omega)-\Delta_{m(k)}\leq\eta_{\kappa+1}(\omega)-\Delta_{m(\kappa+1)}.
\]
Thus inequality \ref{eq:temp-406} has been verified.

4. Inequality \ref{eq:temp-406} implies that the limit 
\[
\tau\equiv\tau_{i}\equiv\lim_{\kappa\rightarrow\infty}\eta_{\kappa}
\]
exists uniformly on $D_{k+}^{c}$, with
\begin{equation}
t\leq\eta_{k}-2^{-m(k)}\leq\tau\leq\eta_{k}\leq t'\label{eq:temp-406-1}
\end{equation}
on  $D_{k+}^{c}$. Since $P(D_{k+})\leq2^{-k+3}$ is arbitrarily small,
we conclude that 
\[
\eta_{\kappa}\downarrow\tau\quad\mbox{ a.u.},
\]
with $t<\tau\leq t'$, or, in fuller notations, with
\begin{equation}
q_{i-1}<\tau_{i}\leq q_{i}\label{eq:temp-406-1-1}
\end{equation}

5. Now let $h\geq n$ be arbitrary, and let $\omega\in D_{h+}^{c}$
be arbitrary. Consider each $u\in[t,\tau(\omega))Q_{\infty}$. Then
$u\in[t,\eta_{k}(\omega))Q_{m(k)}$ for some $k\geq h$. Hence, by
the basic properties of the simple first exit time $\eta_{k}$, we
have
\begin{equation}
d(Z(t,\omega),Z(u,\omega))\leq\widehat{\beta}_{n}(u)<\beta_{n,\infty}.\label{eq:temp-362-2}
\end{equation}
Thus we obtain the bounding relation
\begin{equation}
Z([q_{i-1},\tau_{i}(\omega))Q_{\infty},\omega)\subset(d(\cdot,Z(q_{i-1},\omega))<\beta_{n,\infty})\subset(d(\cdot,Z(q_{i-1},\omega))<2^{-n+3}).\label{eq:temp-356}
\end{equation}

6. To obtain a similar bounding relation for the set $Z([\tau_{i}(\omega),q_{i}],\omega)$,
we will first prove that
\begin{equation}
d(Z(w,\omega),Z(\eta_{h}(\omega),\omega))\leq\beta_{h+1,\infty}\label{eq:temp-247-1-1}
\end{equation}
for each $w\in[\eta_{h+1}(\omega),\eta_{h}(\omega)]Q_{\infty}$. 

To that end, write, for abbreviation, write $u\equiv\eta_{h}(\omega)-\Delta_{m(h)}$,
$r\equiv\eta_{h+1}(\omega)$, and $u'\equiv\eta_{h}(\omega)$. From
inequality \ref{eq:temp-406}, where $k,\kappa$ are replaced by $h,h+1$
respectively, we obtain $u<r\leq u'$. The desired inequality \ref{eq:temp-247-1-1}
holds trivially if $r=u'$. Hence we may assume that $u<r<u'$. Consequently,
since $u,u'$ are consecutive points in the set $Q_{m(h)}$ of dyadic
rationals, we have $r\in Q_{m(h+1)}Q_{m(h)}^{-}$, whence $\widehat{\beta}_{n}(r)=\beta_{n,h+1}$.
Moreover, since $u'\leq t'$ according to inequality \ref{eq:temp-310},
we have 
\[
\eta_{h+1}(\omega)\equiv r<u'\leq t'.
\]
In words, the sample path $Z(\cdot,\omega)|[t,t']Q_{m(h+1)}$ successfully
exits the time-varying $\widehat{\beta}_{n}$-neighborhood of $Z(t,\omega)$,
for the first time at $r$. Therefore 
\begin{equation}
d(Z(t,\omega),Z(r,\omega))>\widehat{\beta}_{n}(r)=\beta_{n,h+1}.\label{eq:temp-315}
\end{equation}
On the other hand, since $u\in Q_{m(h)}\subset Q_{m(h+1)}$ and $u<r$,
exit has not occurred at time $u$. Hence 
\begin{equation}
d(Z(t,\omega),Z(u,\omega))\leq\widehat{\beta}_{n}(u)\leq\beta_{n,h}.\label{eq:temp-316}
\end{equation}
Inequalities \ref{eq:temp-315} and \ref{eq:temp-316} together yield,
by the triangle inequality,
\[
d(Z(u,\omega),Z(r,\omega))>\beta_{n,h+1}-\beta_{n,h}=\beta_{h+1}.
\]
It follows that $\omega\in A_{u,r}^{\beta(h+1)}$ by the definition
of the sets $A_{u,r}^{\beta(h+1)}$. 

Now consider an arbitrary $s\in[r,u')Q_{m(h+1)}$. Then, trivially,
\begin{equation}
\omega\in A_{u,r}^{\beta(h+1)}\subset(A_{u,r}^{\beta(h+1)}\cup A_{u,s}^{\beta(h+1)})(A_{u,r}^{\beta(h+1)}\cup A_{u',s}^{\beta(h+1)}).\label{eq:temp-364}
\end{equation}
At the same time, since $\omega\in D_{\overline{h}+}^{c}\subset D_{h}^{c}$
and since $u,u'\in Q_{m(h)}$ with $u'\equiv u+\Delta_{m(h)}$, we
can apply the defining formula \ref{eq:temp-201-4-3-3} of the exceptional
set $D_{h}$, to obtain
\[
\omega\in D_{h}^{c}
\]
\begin{equation}
\subset(A_{u,r}^{\beta(h+1)}\cup A_{u,s}^{\beta(h+1)})^{c}\cup(A_{u,r}^{\beta(h+1)}\cup A_{u',s}^{\beta(h+1)})^{c}\cup(A_{u',r}^{\beta(h+1)}\cup A_{u',s}^{\beta(h+1)})^{c}.\label{eq:temp-365}
\end{equation}
Relations \ref{eq:temp-364} and \ref{eq:temp-365} together then
imply that
\[
\omega\in(A_{u',r}^{\beta(h+1)}\cup A_{u',s}^{\beta(h+1)})^{c}\subset(A_{u',s}^{\beta(h+1)})^{c}.
\]
Consequently, by the definition of the set $A_{u',s}^{\beta(h+1)}$,
we have
\begin{equation}
d(Z(s,\omega),Z(u',\omega))\leq\beta_{h+1},\label{eq:temp-363}
\end{equation}
where $s\in[r,u')Q_{m(h+1)}$ is arbitrary. The same inequality trivially
holds for $s=u'$. Summing up,
\[
\bigvee_{s\in[r,u']Q(m(h+1))}d(Z_{u'}(\omega),Z_{s}(\omega))\leq\beta_{h+1}.
\]
Then, for each $w\in[r,u']Q_{\infty}$, we can apply inequality \ref{eq:temp-366}
of Lemma \ref{Lem. Existence of supremum distance Z}, where $h,v,v',\overline{v},r,u$
are replaced by $h+1,$$r,u',u',w,s$ respectively, to obtain
\[
d(Z_{w}(\omega),Z_{u'}(\omega))\leq\widehat{\beta}_{h+2}(w)+\bigvee_{s\in[r,u']Q(m(h+1))}d(Z_{u'}(\omega),Z_{s}(\omega))
\]
\begin{equation}
\leq\beta_{h+2,\infty}+\beta_{h+1}=\beta_{h+1,\infty}.\label{eq:temp-329-2}
\end{equation}
In other words, inequality \ref{eq:temp-247-1-1} is verified.

7. Now let $r'\in(\tau(\omega),\eta_{h}(\omega)]Q_{\infty}$be arbitrary.
Then $r'\in[\eta_{k+1}(\omega),\eta_{h}(\omega)]Q_{\infty}$ for some
$k\geq h$. 
\[
d(Z(r',\omega),Z(\eta_{h}(\omega),\omega))
\]
\[
\leq d(Z(r',\omega),Z(\eta_{k}(\omega),\omega))+d(Z(\eta_{k}(\omega),\omega),Z(\eta_{k-1}(\omega),\omega)+
\]
\[
\cdots+d(Z(\eta_{h+1}(\omega),\omega),Z(\eta_{h}(\omega),\omega)
\]
\begin{equation}
\leq\beta_{k+1,\infty}+\beta_{k,\infty}+\cdots+\beta_{h+1,\infty}<\beta_{h+1,\infty}<2^{-h+5},\label{eq:temp-412}
\end{equation}
where the second inequality is by repeated applications of inequality
\ref{eq:temp-247-1-1}. Thus
\begin{equation}
d(Z(r',\omega),Z(\eta_{h}(\omega),\omega))<2^{-h+5}\label{eq:temp-412-3}
\end{equation}
for each $r'\in(\tau(\omega),\eta_{h}(\omega)]Q_{\infty}$, where
$h\geq n$ and $\omega\in D_{h+}^{c}$ are arbitrary. 

8. We will prove that $X_{\tau}$ is a well defined r.r.v. To that
end, let $h\geq n$ be arbitrary. By Proposition \ref{Prop. Right-limit extenaion of D-regular process is continuous a.u.},
there exists $\delta_{h}\equiv\delta_{cau}(2^{-m(h)-1-h})>0$ with
$\delta_{h}\in Q_{\infty}$ such that, for each $r\in[0,1]$, there
exists an exceptional set $H_{r}$ with $P(H_{r})<2^{-m(h)-1-h}$
such that, for each $\omega\in H_{r}^{c}$, we have
\begin{equation}
d(X(r,\omega),X(r',\omega))<2^{-m(h)-1-h}\leq2^{-h+5}\label{eq:temp-358-3}
\end{equation}
for each $r'\in[r-\delta_{h},r+\delta_{h}]\cap domain(X(\cdot,\omega))$.
Define the exceptional set
\begin{equation}
B_{h}\equiv\bigcup_{r\in Q(m(h))}H_{r}.\label{eq:temp-335-1}
\end{equation}
Then
\begin{equation}
P(B_{h})\leq\sum_{r\in Q(m(h))}2^{-m(h)-1-h}<2^{-h}.\label{eq:temp-423}
\end{equation}
Now let $\omega\in B_{h}^{c}D_{h+}^{c}$ be arbitrary, and write $r\equiv\eta_{h}(\omega)\geq\tau(\omega)$.
Let $u\in(\tau(\omega),\tau(\omega)+\delta_{h})Q_{\infty}$ be arbitrary.
Then either (i) $u\in(\tau(\omega),r]Q_{\infty}$ or (ii) $u\in[r,r+\delta_{h})Q_{\infty}.$
Consider Case (i). Then, since $\omega\in B_{h}^{c}\subset D_{h+}^{c}$,
inequality \ref{eq:temp-412-3} applies and yields
\begin{equation}
d(Z_{u}(\omega),Z_{\eta(h)}(\omega))\equiv d(Z(u,\omega),Z(r,\omega))<2^{-h+5}.\label{eq:temp-412-3-1}
\end{equation}
Consider Case (ii). Then, since $\omega\in B_{h}^{c}\subset H_{r}^{c}$,
and since $r,u\in Q_{\infty}\subset domain(X(\cdot,\omega))$, inequality
\ref{eq:temp-358-3} holds with $X$ replaced by $Z$, to yield
\begin{equation}
d(Z_{\eta(h)}(\omega),Z_{u}(\omega))\equiv d(Z(r,\omega),Z(u,\omega))<2^{-h+5}.\label{eq:temp-358-3-1}
\end{equation}
Combining, we see that 
\begin{equation}
d(Z_{\eta(h)}(\omega),Z_{u}(\omega))<2^{-h+5}\label{eq:temp-358-3-1-2}
\end{equation}
for each $u\in(\tau(\omega),\tau(\omega)+\delta_{h})Q_{\infty}$. 

Now consider each $u\in[\tau(\omega),\tau(\omega)+\delta_{h})Q_{\infty}$.
Pick an arbitrary sequence $(u_{k})_{k=1,2,\cdots}$ in $(u,\tau(\omega)+\delta_{h})Q_{\infty}\subset(\tau(\omega),\tau(\omega)+\delta_{h})Q_{\infty}$
with $u_{k}\rightarrow u$. Then, for each $k\geq1$, inequality \ref{eq:temp-358-3-1-2}
holds for $u_{k}$ in the place of $u$. At the same time, since the
process $Z$ is continuous a.u., we have $Z_{u(k)}(\omega)\rightarrow Z_{u}(\omega)$.
Consequently, 
\begin{equation}
d(Z_{\eta(h)}(\omega),Z_{u}(\omega))\leq2^{-h+5}\label{eq:temp-358-3-1-2-1}
\end{equation}
where $u\in[\tau(\omega),\tau(\omega)+\delta_{h})Q_{\infty}$ is arbitrary.
It follows that $\lim_{u\rightarrow\tau(\omega);u\geq\tau(\omega)}Z(u,\omega)$
exists. In other words, $(\tau(\omega),\omega)\in domain(X)$, or,
equivalently, $\omega\in domain(X_{\tau})$. Moreover, letting $u\downarrow\tau(\omega)$
in inequality \ref{eq:temp-358-3-1-2-1}, we obtain
\begin{equation}
d(X_{\tau(i)}(\omega),Z_{\eta(h)}(\omega))\equiv d(X_{\tau}(\omega),Z_{\eta(h)}(\omega))\leq2^{-h+5}\label{eq:temp-337}
\end{equation}
where $\omega\in B_{h}^{c}D_{h+}^{c}$ is arbitrary. Since 
\[
P(B_{h}\cup D_{h+})\leq P(B_{h})+P(D_{h+})<2^{-h}+2^{-h+4}<2^{-h+5},
\]
where $h\geq n$ is arbitrary, we see that $Z_{\eta(h)}\downarrow X_{\tau}$
a.u. as $h\rightarrow\infty$. Since $Z_{\eta(h)}$ is a r.v. for
each $h\geq n$, we conclude that the function $X_{\tau}\equiv X_{\tau(i)}$
is a r.v., where $\tau\equiv\tau_{i}$ and $i=1,\cdots,p$ are arbitrary. 

9. Consider the special case where $h=n$. Consider each $\omega\in B_{n}^{c}D_{n+}^{c}$.
Recall that $\eta_{n}=t'\equiv q_{i}$. Therefore
\[
[\tau(\omega),q_{i}]Q_{\infty}=[\tau(\omega),\eta_{n}(\omega)]Q_{\infty}\subset[\tau(\omega),\tau(\omega)+\delta_{n})Q_{\infty}\cup(\tau(\omega),\eta_{n}(\omega)]Q_{\infty}.
\]
Hence
\[
Z([\tau_{i}(\omega),q_{i}]Q_{\infty},\omega)\subset Z([\tau(\omega),\tau(\omega)+\delta_{n})Q_{\infty},\omega)\cup Z((\tau(\omega),\eta_{n}(\omega)]Q_{\infty},\omega)
\]
\begin{equation}
\subset(d(\cdot,Z_{\eta(n)}(\omega))\leq2^{-n+5})=(d(\cdot,Z(q_{i},\omega))\leq2^{-n+5}),\label{eq:temp-356-1}
\end{equation}
where $\omega\in B_{n}^{c}D_{n+}^{c}$ is arbitrary, and where the
second containment relation is thanks to inequalities \ref{eq:temp-358-3-1-2-1}
and \ref{eq:temp-412-3}.

10. Now consider an arbitrary $\omega\in B_{n}^{c}D_{n+}^{c}$. For
convenience, define the constant r.r.v.'s $\tau_{0}\equiv0$ and $\tau_{p+1}\equiv1$.
It is easy to combine relations \ref{eq:temp-356} and \ref{eq:temp-356-1}
for each of the intervals in the disjoint union
\begin{equation}
\theta_{1}\cup\theta_{2}\cup\cdots\cup\theta_{p`}\cup\theta_{p+1},\label{eq:temp-371-1}
\end{equation}
where
\[
(\theta_{1},\theta_{2},\cdots,\theta_{p},\theta_{p+1})
\]
\begin{equation}
\equiv([\tau_{0}(\omega),\tau_{1}(\omega)),[\tau_{1}(\omega),\tau_{2}(\omega)),\cdots,[\tau_{p-1}(\omega),\tau_{p}(\omega)),[\tau_{p}(\omega),\tau_{p+1}(\omega)]).\label{eq:temp-371-1-1}
\end{equation}
More precisely, 
\[
Z(\theta_{i}Q_{\infty},\omega)\equiv Z([\tau_{i-1}(\omega),\tau_{i}(\omega))Q_{\infty},\omega)
\]
\[
=Z([\tau_{i-1}(\omega),q_{i-1}]Q_{\infty},\omega)\cup Z([q_{i-1},\tau_{i}(\omega))Q_{\infty},\omega)
\]
\[
\subset(d(\cdot,Z(q_{i-1},\omega))\leq2^{-n+5})\cup(d(\cdot,Z(q_{i-1},\omega))\leq2^{-n+3})
\]
\[
=(d(\cdot,Z(q_{i-1},\omega))\leq2^{-n+5}),
\]
where the containment relation $\subset$ holds by application of
condition \ref{eq:temp-356-1} to $i-1$, and condition \ref{eq:temp-356}
to $i$, where $i=1,\cdots,p$ is arbitrary. Similarly
\[
Z(\theta_{p+1}Q_{\infty},\omega)=Z([\tau_{p}(\omega),q_{p}]Q_{\infty},\omega)\subset(d(\cdot,Z(q_{p},\omega))<2^{-n+5})
\]
by condition \ref{eq:temp-356-1}. Summing up, we have
\begin{equation}
Z(\theta_{\kappa}Q_{\infty},\omega)\subset(d(\cdot,Z(q_{\kappa-1},\omega))\leq2^{-n+5}),\label{eq:temp-362}
\end{equation}
for each $\kappa=1,\cdots,p+1$.

11. We wish to prove that each of the intervals $\theta_{1},\cdots,\theta_{p},\theta_{p+1}$
has positive length, provided that we exclude also a third small exceptional
set of $\omega$. To be precise, recall from the beginning of this
proof that 
\begin{equation}
\delta_{cau}(\varepsilon_{n})\equiv\delta_{cau}(\varepsilon_{n},\overline{m},\delta_{Cp})\equiv\Delta_{m(J)}\equiv2^{-m(J)}<4^{-1}\delta_{Cp}(2^{-2\nu(\varepsilon(n))+2}).\label{eq:temp-332-1-2-1}
\end{equation}
For abbreviation write
\[
\overline{\Delta}\equiv\Delta_{m(J)}<4^{-1}\delta_{Cp}(2^{-2\nu(\varepsilon(n))+2}).
\]
By Proposition \ref{Prop. Right-limit extenaion of D-regular process is continuous a.u.},
for a.e. $\omega\in\Omega$, the right-limit extension $X(\cdot,\omega)$
is defined at each point in $Q_{\infty}$, with $X(\cdot,\omega)|Q_{\infty}=Z(\cdot,\omega)$.
Consequently, for a.e. $\omega\in\Omega$, $domain(X(\cdot,\omega))$
is dense in $[0,1]$. Moreover, by the same Proposition \ref{Prop. Right-limit extenaion of D-regular process is continuous a.u.},
for each $r\in[0,1]$, there exists an exceptional set $G_{r}$ with
$P(G_{r})<\varepsilon_{n}$ such that, for each $\omega\in G_{r}^{c}$,
we have
\begin{equation}
d(X(r,\omega),X(r',\omega))<\varepsilon_{n}\leq\beta_{n}\label{eq:temp-358}
\end{equation}
for each $r'\in[r-\overline{\Delta},r+\overline{\Delta}]\cap domain(X(\cdot,\omega))$. 

Define the exceptional set
\begin{equation}
C_{n}\equiv\bigcup_{r\in Q(m(n))}G_{r}.\label{eq:temp-335}
\end{equation}
Then
\[
P(C_{n})\leq\sum_{r\in Q(m(n))}\varepsilon_{n}<2^{m(n)+1}\varepsilon_{n}=2^{-n}.
\]

12. Let $\omega\in B_{n}^{c}D_{n+}^{c}C_{n}^{c}$ be arbitrary, and
let $s\in[t,t+\overline{\Delta})Q_{\infty}$ be arbitrary. Then $\omega\in G_{r}^{c}$
for each $r\in Q_{m(n)}$. Hence inequality \ref{eq:temp-358} holds
for each $r'\in[r-\overline{\Delta},r+\overline{\Delta}]\cap domain(X(\cdot,\omega))$,
for each $r\in Q_{m(n)}$. In particular, $\omega\in G_{t}^{c}$,
and so inequality \ref{eq:temp-358}, with $t$ in the place of $r$,
holds on
\[
[t,s]Q_{\infty}\subset[t,t+\overline{\Delta})Q_{\infty}\subset[t-\overline{\Delta},t+\overline{\Delta}]\cap domain(X(\cdot,\omega)),
\]
whence
\begin{equation}
d(X(t,\omega),X(r',\omega))\leq\beta_{n}\leq\widehat{\beta}_{n}(r')\label{eq:temp-360}
\end{equation}
for each $r'\in[t,s]Q_{\infty}.$ Consequently, for each $k\geq n$,
the sample path $Z(\cdot,\omega)|[t,t']Q_{m(k)}$ stays within the
the time varying $\widehat{\beta}_{n}$-neighborhood of $Z(t,\omega)$
up to and including time $s$. Hence, according to Assertion 4 of
Proposition \ref{Prop. Basics of simple Exit times}, the simple first
exit time $\eta_{k}(\omega)\equiv\eta_{t,\widehat{\beta}(n),[t,t']Q(m(k))}(\omega)$
can come only after time $s$. In other words $s<\eta_{k}(\omega).$
Letting $k\rightarrow\infty,$ we therefore obtain $s\leq\tau(\omega)$.
Since $s\in[t,t+\overline{\Delta})Q_{\infty}$ is arbitrary, it follows
that $t+\overline{\Delta}\leq\tau(\omega)$. Therefore, 
\begin{equation}
|\theta_{i}|\geq\tau_{i}(\omega)-q_{i-1}\equiv\tau(\omega)-t\geq\overline{\Delta},\label{eq:temp-313}
\end{equation}
where $i=1,\cdots,p$ is arbitrary.

13. The interval $\theta_{p+1}\equiv(\tau_{p}(\omega),1]$ however
remains, with length possibly less than $\overline{\Delta}$. To deal
with this nuisance, we will blatantly replace $\tau_{p}$ with the
r.r.v. $\overline{\tau}_{p}\equiv\tau_{p}\wedge(1-\overline{\Delta})$,
while keeping $\overline{\tau}_{i}\equiv\tau_{i}$ if $i\neq p$,
and prove that the other desirable properties of the intervals in
the sequence \ref{eq:temp-371-1-1} are preserved. More precisely,
define the sequence
\[
(\overline{\theta}_{1},\overline{\theta}_{2},\cdots,\overline{\theta}_{p},\overline{\theta}_{p+1})
\]
\begin{equation}
\equiv([\tau_{0}(\omega),\tau_{1}(\omega)),[\tau_{1}(\omega),\tau_{2}(\omega)),\cdots,[\tau_{p-1}(\omega),\overline{\tau}_{p}(\omega)),[\overline{\tau}_{p}(\omega),\tau_{p+1}(\omega)]).\label{eq:temp-371-1-1-1}
\end{equation}
Note that only the last two members of this sequence are affected
by the change. Hence $|\overline{\theta}_{i}|\equiv|\theta_{i}|\geq\overline{\Delta}$,
if $1\leq i\leq p-1$. Moreover, trivially,
\begin{equation}
|\overline{\theta}_{p+1}|=1-\tau_{p}(\omega)\wedge(1-\overline{\Delta})\geq\overline{\Delta}.\label{eq:temp-359}
\end{equation}
Furthermore, 
\[
\overline{\tau}_{p}(\omega)-q_{p-1}\equiv\tau_{p}(\omega)\wedge(1-\overline{\Delta})-q_{p-1}=(\tau_{p}(\omega)-q_{p-1})\wedge(1-\overline{\Delta}-q_{p-1})
\]
\begin{equation}
=(\tau_{p}(\omega)-q_{p-1})\wedge(\Delta-\overline{\Delta})\geq\overline{\Delta}\wedge\overline{\Delta}=\overline{\Delta},\label{eq:temp-346}
\end{equation}
where the last inequality follows from inequality \ref{eq:temp-313}
and from the inequality
\[
\Delta-\overline{\Delta}\equiv2^{-m(n)}-2^{-m(J)}\geq2^{-m(J)+1}-2^{-m(J)}=2^{-m(J)}=\overline{\Delta}.
\]
Hence
\begin{equation}
|\overline{\theta}_{p}|=\overline{\tau}_{p}(\omega)-\tau_{p-1}(\omega)\geq\overline{\tau}_{p}(\omega)-q_{p-1}\geq\overline{\Delta}.\label{eq:temp-361}
\end{equation}
Combining inequalities \ref{eq:temp-313}, \ref{eq:temp-359}, and
\ref{eq:temp-361}, we see that 
\begin{equation}
|\overline{\theta}_{\kappa}|\equiv\overline{\tau}_{\kappa}(\omega)-\overline{\tau}_{\kappa-1}(\omega)\geq\overline{\Delta}\equiv\delta_{aucl}(\varepsilon,\overline{m},\delta_{Cp})\label{eq:temp-341}
\end{equation}
for each $\kappa=1,\cdots,p+1$, where $\omega\in B_{n}^{c}D_{n+}^{c}C_{n}^{c}$
is arbitrary.

14. We will now verify that relation \ref{eq:temp-362} still holds
when $\theta_{\kappa}$ is replaced by $\overline{\theta}_{\kappa}$,
for each $\kappa=1,\cdots,p+1$. For $\kappa\leq1,\cdots,p$, we have
$\overline{\theta}_{\kappa}\subset\theta_{\kappa}$, whence relation
\ref{eq:temp-362} is trivially preserved. It remains to verify that
\begin{equation}
Z([\overline{\tau}_{p}(\omega),1]Q_{\infty},\omega)\subset(d(\cdot,Z(1,\omega))\leq2^{-n+5}).\label{eq:temp-362-1}
\end{equation}
To that end, let $r'\in[\overline{\tau}_{p}(\omega),1]Q_{\infty}\equiv[\tau_{p}(\omega)\wedge(1-\overline{\Delta}),1]Q_{\infty}$
be arbitrary. Either (i') $r'<1-\overline{\Delta}$ or (ii') $r'\geq1-\overline{\Delta}$.
In Case (i'), the assumption $r'<\tau_{p}(\omega)$ would imply $r'<\tau_{p}(\omega)\wedge(1-\overline{\Delta})$,
a contradiction. Therefore, in Case (i') we have $r'\in[\tau_{p}(\omega),1]Q_{\infty}\equiv\theta_{p+1}Q_{\infty}$,
whence 
\[
d(Z(r',\omega),Z(1,\omega))\leq2^{-n+5}
\]
by relation \ref{eq:temp-362}. In Case (ii'), because $1\in Q_{m(n)},$
inequality \ref{eq:temp-358} applies to $r\equiv1$, to yield
\begin{equation}
d(Z(r',\omega),Z(1,\omega))<\varepsilon_{n}\leq\beta_{n}<2^{-n+5}.\label{eq:temp-358-2}
\end{equation}
Thus relation \ref{eq:temp-362-1} has been verified in either case.
Summing up,
\begin{equation}
Z(\overline{\theta}_{\kappa}Q_{\infty},\omega)\subset(d(\cdot,Z(q_{\kappa-1},\omega))\leq2^{-n+5}),\label{eq:temp-362-3}
\end{equation}
for each $\kappa=1,\cdots,p+1$, where $\omega\in B_{n}^{c}D_{n+}^{c}C_{n}^{c}$
is arbitrary. 

15. We will now prove that $X_{\overline{\tau}(i)}$ is a r.v. If
$i<p$, then $\overline{\tau}_{i}\equiv\tau_{i}$ and so $X_{\overline{\tau}(i)}=X_{\tau(i)}$
is a r.v.. Hence it suffices to treat the case where $i=p$. Then
$\overline{\tau}_{p}\equiv\tau_{p}\wedge\overline{q}$ where we write
$\overline{q}\equiv1-\overline{\Delta}$ for short. Again, let $h\geq n$,
and $\omega\in B_{h}^{c}D_{h+}^{c}$ be arbitrary. Note that $\omega\in B_{h}^{c}\subset H_{\overline{q}}^{c}$.
Let 
\begin{equation}
u,v\in[\overline{\tau}_{p}(\omega),\overline{\tau}_{p}(\omega)+\delta_{h})Q_{\infty}\label{eq:temp-340}
\end{equation}
be arbitrary. There are three possibilities: (i'') $\tau_{p}(\omega)<\overline{q}$,
(ii'') $\overline{q}<\tau_{p}(\omega),$ or (iii'') $|\overline{q}-\tau_{p}(\omega)|<2^{-1}\delta_{h}$. 

Consider Case (i''). We then have $\overline{\tau}_{p}(\omega)=\tau_{p}(\omega)$,
whence $u,v\in[\tau_{p}(\omega),\tau_{p}(\omega)+\delta_{h})Q_{\infty}$.
Therefore inequality \ref{eq:temp-358-3-1-2-1} holds for both $u$
and for $v$, which implies, by the triangle inequality, that
\begin{equation}
d(Z_{v}(\omega),Z_{u}(\omega))<2^{-h+5}+2^{-h+5}=2^{-h+6}.\label{eq:temp-339}
\end{equation}

Consider Case (ii''). We then have $\overline{\tau}_{p}(\omega)=\overline{q}$,
whence $u,v\in[\overline{q},\overline{q}+\delta_{h})Q_{\infty}$.
Therefore inequality \ref{eq:temp-358-3} holds with $X,r,r'$ replaced
by $Z,\overline{q},u$ respectively, to yield
\begin{equation}
d(Z(\overline{q},\omega),Z(u,\omega))<2^{-h+5}.\label{eq:temp-358-3-1-3}
\end{equation}
Similarly, the last inequality holds for $v$ in the place of $u$.
Hence the triangle inequality establishes inequality \ref{eq:temp-339}
also for Case (ii''). 

Now consider Case (iii''), where $|\overline{q}-\tau_{p}(\omega)|<2^{-1}\delta_{h}$.
Let $w\equiv\overline{q}+2^{-1}\delta_{h}\in Q_{\infty}$. Then, by
the triangle inequality, 
\[
w\in[\tau_{p}(\omega),\tau_{p}(\omega)+\delta_{h})Q_{\infty}\cap[\overline{q},\overline{q}+\delta_{h})Q_{\infty}.
\]
Since $u,\overline{q}\in Q_{\infty}$, we have either $u<\overline{q}$
or $u\geq\overline{q}$. Suppose $u<\overline{q}$. Then relation
\ref{eq:temp-340} implies that $\tau_{p}(\omega)\wedge\overline{q}\equiv\overline{\tau}_{p}(\omega)<\overline{q}$,
whence $\overline{\tau}_{p}(\omega)=\tau_{p}(\omega)$ and relation
\ref{eq:temp-340} reduces to $u\in[\tau_{p}(\omega),\tau_{p}(\omega)+\delta_{h})Q_{\infty}$.
Hence inequality \ref{eq:temp-358-3-1-2} holds for $u$, and similarly
for $w$. The triangle inequality then implies that
\begin{equation}
d(Z_{w}(\omega),Z_{u}(\omega))<2^{-h+5}+2^{-h+5}=2^{-h+6}.\label{eq:temp-339-2}
\end{equation}
Now suppose $u\geq\overline{q}$. Then relation \ref{eq:temp-340}
implies that $u\in[\overline{q},\overline{q}+\delta_{h})Q_{\infty}.$
Therefore inequality \ref{eq:temp-358-3} holds with $X,r,r'$ replaced
by $Z,\overline{q},u$ respectively, to yield

\begin{equation}
d(Z(\overline{q},\omega),Z(u,\omega))<2^{-h+5}.\label{eq:temp-358-3-1-3-1}
\end{equation}
Similarly,
\[
d(Z(\overline{q},\omega),Z(w,\omega))<2^{-h+5}.
\]
Combining the last two displayed inequalities, inequality \ref{eq:temp-339-2}
is proved for $u$. Repeating this argument with $v$ in the place
of $u$, we see that inequality \ref{eq:temp-339-2} holds also for
$v$. Hence the triangle inequality implies that 
\begin{equation}
d(Z_{v}(\omega),Z_{u}(\omega))<2^{-h+6}+2^{-h+6}=2^{-h+7}.\label{eq:temp-339-2-1}
\end{equation}
Summing up, in each of Cases (i''-iii''), we have inequality \ref{eq:temp-339-2-1}
where 
\[
u,v\in[\overline{\tau}_{p}(\omega),\overline{\tau}_{p}(\omega)+\delta_{h})Q_{\infty}
\]
and $h\geq n$ are arbitrary. It follows that $\lim_{u\rightarrow\overline{\tau}(p,\omega);u\geq\overline{\tau}(p,\omega)}Z(u,\omega)$
exists. In other words, $(\overline{\tau}_{p}(\omega),\omega)\in domain(X)$,
or, equivalently, $\omega\in domain(X_{\overline{\tau}(p)}),$ with
$X_{\overline{\tau}(p)}(\omega)$ by definition equal to this limit.
Moreover, letting $u\downarrow\overline{\tau}_{p}(\omega)$ in inequality
\ref{eq:temp-339-2-1}, we obtain
\begin{equation}
d(Z_{v}(\omega),X_{\overline{\tau}(p)}(\omega))\leq2^{-h+7},\label{eq:temp-349}
\end{equation}
where $v\in[\overline{\tau}_{p}(\omega),\overline{\tau}_{p}(\omega)+\delta_{h})Q_{\infty}$
and $h\geq n$ are arbitrary. Now define 
\[
\overline{\eta}_{h,p}\equiv\eta_{h,p}\wedge\overline{q}.
\]
Since $\eta_{h,p}(\omega),\overline{q}\in Q_{m(h)}$, we have either
$\overline{q}\leq\eta_{h,p}(\omega)-2^{-m(h)}$ or $\eta_{h,p}(\omega)\leq\overline{q}$.
In view of inequality \ref{eq:temp-406-1}, it follows that either
(i''') $\tau_{p}(\omega)\leq\eta_{h,p}(\omega)\leq\overline{q}$ or
(ii''') $\overline{q}\leq\tau_{p}(\omega)\leq\eta_{h,p}(\omega)$.
In Case (i'''), we have $\overline{\tau}_{p}(\omega)=\tau_{p}(\omega)$
and $\overline{\eta}_{h,p}(\omega)=\eta_{h,p}(\omega)$, whence inequality
\ref{eq:temp-337} implies that 
\begin{equation}
d(X_{\overline{\tau}(p)}(\omega),Z_{\overline{\eta}(h,.p)}(\omega))\leq2^{-h+5}.\label{eq:temp-337-1}
\end{equation}
In Case (ii'''), we have $\overline{\tau}_{p}(\omega)=\overline{\eta}_{h,p}(\omega)=\overline{q}$,
whence the same inequality \ref{eq:temp-337-1} holds trivially. Since
$\omega\in B_{h}^{c}D_{h+}^{c}$ is arbitrary, and since $P(B_{h}\cup D_{h+})<2^{-h+5}$,
where $h\geq n$ is arbitrary, we conclude that $Z_{\overline{\eta}(h,.p)}\downarrow X_{\overline{\tau}(p)}$
a.u. as $h\rightarrow\infty$. Note that, for each $h\geq n$, the
function $\overline{\eta}_{h,p}$ is a r.r.v. with values in the finite
set $Q_{m(h)}$. Hence $Z_{\overline{\eta}(h,.p)}$ is a r.v. for
each $h\geq n$. We conclude that the function $X_{\overline{\tau}(p)}$
is a r.v. 

16. Now define $H\equiv B_{n}\cup D_{n+}\cup C_{n}$, with 
\begin{equation}
P(H)=P(B_{n}\cup D_{n+}\cup C_{n})<2^{-n}+2^{-n+3}+2^{-n}<2^{-n+4}<\varepsilon.\label{eq:temp-370}
\end{equation}
Let $\omega\in H^{c}=C_{n}^{c}B_{n}^{c}D_{n+}^{c}$ and $\kappa=1,\cdots,p+1$
be arbitrary. Relation \ref{eq:temp-362-3} implies that 
\[
d(Z(u,\omega),Z(q_{\kappa-1},\omega))\leq2^{-n+5}
\]
for each $u\in\overline{\theta}_{\kappa}Q_{\infty}$. Hence, by the
definition of $X(\cdot,\omega)$ as the right limit of $Z$, it follows
that 
\begin{equation}
d(X(\cdot,\omega),Z(q_{\kappa-1},\omega))\leq2^{-n+5}<\varepsilon\label{eq:temp-350}
\end{equation}
for each $v\in\overline{\theta}_{\kappa}\cap domain(X(\cdot,\omega))$.
In particular, since $\overline{\tau}_{\kappa}(\omega)\in\overline{\theta}_{\kappa}\cap domain(X(\cdot,\omega))$,
we have 
\[
d(X(\overline{\tau}_{\kappa}(\omega),\omega),Z(q_{\kappa-1},\omega))\leq2^{-n+5}.
\]
Combining the last two displayed inequality, we obtain 
\begin{equation}
d(X(\overline{\tau}_{\kappa}(\omega),\omega),X(v,\omega))\leq2^{-n+5}+2^{-n+5}=2^{-n+6}<\varepsilon\label{eq:temp-357}
\end{equation}
for each $v\in\overline{\theta}_{\kappa}\cap domain(X(\cdot,\omega)),$
where $\omega\in H^{c}$ is arbitrary.

17. Let $\omega\in H^{c}$ and $i=1,\cdots,p+1$ be arbitrary. By
inequalities \ref{eq:temp-406-1-1} and \ref{eq:temp-313},
\begin{equation}
q_{i-1}+\overline{\Delta}\leq\overline{\tau}_{i}(\omega)\equiv\tau_{i}(\omega)\leq q_{i}\label{eq:temp-406-1-1-1}
\end{equation}
if $i\neq p$. At the same time, inequality \ref{eq:temp-346} says
\[
q_{p-1}+\overline{\Delta}\leq\overline{\tau}_{p}(\omega)\equiv\tau_{p}\wedge(1-\overline{\Delta})\leq1-\overline{\Delta}.
\]
Thus
\begin{equation}
0\equiv q_{0}<\overline{\tau}_{1}(\omega)\leq q_{1}<\overline{\tau}_{2}(\omega)\leq\cdots<\overline{\tau}_{p-1}(\omega)\leq q_{p-1}<\overline{\tau}_{p}(\omega)\leq1-\overline{\Delta}<q_{p}=1,\label{eq:temp-308-2}
\end{equation}
with
\begin{equation}
\bigwedge_{j=1}^{p+1}(\overline{\tau}_{j}(\omega)-\overline{\tau}{}_{j-1}(\omega))\geq\bigwedge_{j=1}^{p}(\overline{\tau}_{j}(\omega)-q_{j-1})\wedge(\overline{\tau}_{p+1}(\omega)-\overline{\tau}_{p}(\omega))\geq\overline{\Delta}.\label{eq:temp-342-1-2}
\end{equation}
Since $\varepsilon>0$ is arbitrary, inequalities \ref{eq:temp-342-1-2},
\ref{eq:temp-357}, and \ref{eq:temp-370} together show that the
sequence 
\[
0\equiv\overline{\tau}_{0}<\overline{\tau}_{1}<\cdots<\overline{\tau}_{p}<\overline{\tau}_{p+1}\equiv1
\]
of r.r.v.'s, along with the operation $\delta_{aucl}$, satisfy Condition
3 in Definition \ref{Def. a.u. cadlag process} for the process $X$.
At the same time Assertions 4 and 5 in Proposition \ref{Prop. Right-limit extenaion of D-regular process is continuous a.u.}
imply Conditions 1 and 2 respectively in Definition \ref{Def. a.u. cadlag process}
for the process $X$. The same Proposition \ref{Prop. Right-limit extenaion of D-regular process is continuous a.u.}
also says that the process $X$ is continuous in probability, with
modulus of continuity in probability $\delta_{Cp}$. All the conditions
in Definition \ref{Def. a.u. cadlag process} having been verified,
the process $X$ is a.u. càdlàg, with \emph{$\delta_{aucl}$ }as modulus
of a.u. càdlàg\emph{.} The theorem is proved.
\end{proof}

\section{Continuity of the Right-Limit-Extension Construction}

We will prove that the construction of a.u. càdlàg processes by the
right-limit extension of $D$-regular processes, given in the previous
section, is a continuous construction.

Let $(S,d)$ be a locally compact metric space. As usual, we write
$\widehat{d}\equiv1\wedge d$. Refer to Definition \ref{Def. Notations for dyadic rationals}
for notations related to the enumerated set $Q_{\infty}\equiv\{t_{0},t_{1},\cdots\}$
of dyadic rationals in the interval $[0,1]$, and its subset $Q_{m}\equiv\{q_{m,0},q_{m,1},\cdots,q_{m,p(m)}\}=\{t_{0},\cdots,t_{p(m)}\}$
for each $m\geq0$. 

Recall from Definition \ref{Def. Metric on random fields w/ countable parameters}
that $\widehat{R}(Q_{\infty}\times\Omega,S)$ denotes the space of
stochastic processes with parameter set $Q_{\infty}$, sample space
$\Omega$, and state space $(S,d)$, and that it is equipped with
the metric $\widehat{\rho}_{Prob,Q(\infty)}$ defined by
\begin{equation}
\widehat{\rho}_{Prob,Q(\infty)}(Z,Z')\equiv\sum_{j=0}^{\infty}2^{-j-1}E1\wedge d(Z_{t(j)},Z'_{t(j)})\label{eq:temp-287-1}
\end{equation}
for each $Z,Z'\in\widehat{R}(Q_{\infty}\times\Omega,S)$. 

Recall from Definitions \ref{Def. a.u. cadlag process} and \ref{Def. metric rho_D_hat }
the metric space $(\widehat{D}[0,1],\rho_{\widehat{D}[0,1]})$ of
a.u. càdlàg processes on the interval $[0,1],$ with
\[
\rho_{\widehat{D}[0,1]}(X,X')\equiv\int E(d\omega)\dddot{\widehat{d_{D}}}(X(\cdot,\omega),X'(\cdot,\omega))
\]
for each $X,X'\in\widehat{D}[0,1]$, where $\dddot{\widehat{d_{D}}}\equiv1\wedge d_{D}$,
where $d_{D}$ is the Skorokhod metric on the space $D[0,1]$ of càdlàg
functions. Recall that $[\cdot]_{1}$ is an operation which assigns
to each $a\in R$ an integer $[a]_{1}\in(a,a+2)$.
\begin{thm}
\label{Thm.Construction of a.u. cadlag process by right limit of D-regular proceses is continuous.}
\textbf{\emph{(Continuity of the construction of a.u. càdlàg processes
by right-limit extension of $D$-regular processes).}} Let $\widehat{R}_{Dreg,\overline{m},\delta(Cp)}(Q_{\infty}\times\Omega,S)$
denote the subspace of \emph{$(\widehat{R}(Q_{\infty}\times\Omega,S),\widehat{\rho}_{Prob,Q(\infty)})$}
whose members share some common modulus of continuity in probability
$\delta_{Cp}$ and some common modulus of $D$-regularity $\overline{m}\equiv(m_{n})_{n=0.1.\cdots}$.
Let
\begin{equation}
\Phi_{rLim}:(\widehat{R}_{Dreg,\overline{m},\delta(Cp)}(Q_{\infty}\times\Omega,S),\widehat{\rho}_{Prob,Q(\infty)})\rightarrow(\widehat{D}_{\delta(aucl),\delta(Cp)}[0,1],\rho_{\widehat{D}[0,1]})\label{eq:temp-369-2}
\end{equation}
be the right-limit extension as constructed in Theorem \ref{Thm. Extension of D-regular process by right limit is a.u.cadlag},
where 
\[
(\widehat{D}_{\delta(aucl),\delta(Cp)}[0,1],\rho_{\widehat{D}[0,1]})
\]
is the metric subspace of $a.u.c\grave{a}dl\grave{a}g$ processes
which share the common modulus of continuity in probability $\delta_{Cp}$,
and which share the common modulus of a.u. càdlàg $\delta_{aucl}\equiv\delta_{aucl}(\cdot,\overline{m},\delta_{Cp})$
that is defined in Theorem \ref{Thm. Extension of D-regular process by right limit is a.u.cadlag}.

Then the function $\Phi_{rLim}$ is uniformly continuous, with a modulus
of continuity $\delta_{rLim}(\cdot,\overline{m},\delta_{Cp})$ which
depends only on $\overline{m}$ and $\delta_{Cp}$. 
\end{thm}
\begin{proof}
1. Let $\varepsilon_{0}>0$ be arbitrary. Let $n\geq0$ be so large
that $2^{-n+6}<\varepsilon.$ Let $J\geq n+1$ be so large that 
\begin{equation}
\Delta_{m(J)}\equiv2^{-m(J)}<2^{-2}\delta_{Cp}(2^{-2m(n)-2n-10}).\label{eq:temp-332-1-2-2}
\end{equation}
Then, according to Theorem \ref{Thm. Extension of D-regular process by right limit is a.u.cadlag},
we have
\[
\overline{\Delta}\equiv\Delta_{m(J)}=\delta_{aucl}(\varepsilon,\overline{m},\delta_{Cp}).
\]
Write $\varepsilon\equiv2^{-4}\varepsilon_{0}$. Take $k\geq n$ so
large that 
\[
2^{-m(k)+2}<(1-e^{-\varepsilon})\overline{\Delta}.
\]
Write $p\equiv p_{m(n)}\equiv2^{m(n)}$, $\widetilde{p}\equiv p_{m(k)}\equiv2^{m(k)}$
and $\widetilde{\Delta}\equiv2^{-m(k)}$. Define
\[
\delta_{rLim}(\varepsilon_{0},\overline{m},\delta_{Cp})\equiv2^{-\widetilde{p}-1}\varepsilon^{2}.
\]
We will prove that the operation $\delta_{rLim}(\cdot,\overline{m},\delta_{Cp})$
is a modulus of continuity for the the function $\Phi_{rLim}$ in
expression \ref{eq:temp-369-2}. 

To that end, let $Z,Z'\in\widehat{R}_{Dreg,\overline{m},\delta(Cp)}(Q_{\infty}\times\Omega,S)$
be arbitrary such that
\begin{equation}
\widehat{\rho}_{Prob,Q(\infty)}(Z,Z')\equiv E\sum_{j=0}^{\infty}2^{-j-1}\widehat{d}(Z_{t(j)},Z'_{t(j)})<\delta\equiv\delta_{rLim}(\varepsilon_{0},\overline{m},\delta_{Cp}).\label{eq:temp-287-2}
\end{equation}
Let $X\equiv\Phi_{rLim}(Z)$ and $X'\equiv\Phi_{rLim}(Z')$. We need
only verify that 
\begin{equation}
\rho_{\widehat{D}[0,1]}(X,X')<\varepsilon_{0}.\label{eq:temp-345}
\end{equation}

2. According to Steps 16 and 17 in the proof of Theorem \ref{Thm. Extension of D-regular process by right limit is a.u.cadlag},
there exists an increasing sequence 
\[
0\equiv\overline{\tau}_{0}<\overline{\tau}_{1}<\cdots<\overline{\tau}_{p}<\overline{\tau}_{p+1}\equiv1
\]
of r.r.v.'s, and a $\mathrm{measurable}$ set $H$ with $P(H)<\varepsilon$,
such that, for each $\omega\in H^{c}$, we have
\[
0\equiv\overline{\tau}_{0}(\omega)\equiv q_{0}<\overline{\tau}_{1}(\omega)\leq q_{1}<\overline{\tau}_{2}(\omega)\leq\cdots
\]
\begin{equation}
<\overline{\tau}_{p-1}(\omega)\leq q_{p-1}<\overline{\tau}_{p}(\omega)\leq1-\overline{\Delta}<q_{p}\equiv\overline{\tau}_{p+1}(\omega)\equiv1,\label{eq:temp-308}
\end{equation}
with
\begin{equation}
\bigwedge_{j=1}^{p+1}(\overline{\tau}_{j}(\omega)-\overline{\tau}{}_{j-1}(\omega))\geq\bigwedge_{i=1}^{p}(\overline{\tau}_{i}(\omega)-q_{i-1})\wedge(\overline{\tau}_{p+1}(\omega)-\overline{\tau}_{p}(\omega))\geq\overline{\Delta}.\label{eq:temp-342-1}
\end{equation}
Moreover, inequality \ref{eq:temp-350} in the proof of Theorem \ref{Thm. Extension of D-regular process by right limit is a.u.cadlag}
says that, for each $\omega\in H^{c}$ and for each $i=1,\cdots,p+1$,
we have 
\begin{equation}
d(X(\cdot,\omega),Z(q_{i-1},\omega))\leq\varepsilon\label{eq:temp-307-4-1-1}
\end{equation}
on the interval $\overline{\theta}_{i}\equiv[\overline{\tau}_{i-1}(\omega),\overline{\tau}_{i}(\omega))$
or $\overline{\theta}_{i}\equiv[\overline{\tau}_{i-1}(\omega),\overline{\tau}_{i}(\omega)]$,
according as $1\leq i\leq p$ or $i=p+1$. 

3. Relative to the process $X'$, we can similarly construct the increasing
sequence 
\[
0\equiv\overline{\tau}'_{0}<\overline{\tau}'_{1}<\cdots<\overline{\tau}'_{p}<\overline{\tau}'_{p+1}\equiv1
\]
of r.r.v.'s, and $\mathrm{measurable}$ set $H'$ with $P(H')<\varepsilon$,
such that, for each $\omega\in H'^{c}$, we have 
\[
0=\overline{\tau}'_{0}(\omega)=q_{0}<\overline{\tau}'_{1}(\omega)\leq q_{1}<\overline{\tau}'_{2}(\omega)\leq\cdots
\]
\begin{equation}
<\overline{\tau}'_{p-1}(\omega)\leq q_{p-1}<\overline{\tau}'_{p}(\omega)\leq1-\overline{\Delta}<q_{p}=\overline{\tau}'_{p+1}(\omega)\equiv1,\label{eq:temp-308-1}
\end{equation}
and
\begin{equation}
\bigwedge_{j=1}^{p+1}(\overline{\tau}_{j}(\omega)-\overline{\tau}{}_{j-1}(\omega))\geq\bigwedge_{i=1}^{p}(\overline{\tau}'_{i}(\omega)-q_{i-1})\wedge(\overline{\tau}'_{p+1}(\omega)-\overline{\tau}'_{p}(\omega))\geq\overline{\Delta}.\label{eq:temp-342-1-1}
\end{equation}
Moreover, for each $\omega\in H^{c}$ and for each $i=1,\cdots,p+1$,
we have 
\begin{equation}
d(X'(\cdot,\omega),Z'(q_{i-1},,\omega))\leq\varepsilon\label{eq:temp-307-4-1-1-1}
\end{equation}
on the interval $\overline{\theta}'_{i}\equiv[\overline{\tau}'_{i-1}(\omega),\overline{\tau}'_{i}(\omega))$
or $\overline{\theta}_{i}\equiv[\overline{\tau}'_{i-1}(\omega),\overline{\tau}'_{i}(\omega)]$
according as $1\leq i\leq p$ or $i=p+1$. 

4. Then, in view of the bound \ref{eq:temp-287-2}, Chebychev's inequality
implies that there exists a $\mathrm{measurable}$ set $G$ with $P(G)<\varepsilon$
such that, for each $\omega\in G^{c}$, we have 
\[
\bigvee_{r\in Q(m(k))}d(Z(r,\omega),Z'(r,\omega))=\bigvee_{j=0}^{p(m(k))}d(Z_{t(j)}(\omega),Z'_{t(j)}(\omega))
\]
\[
\leq2^{p(m(k))+1}\sum_{j=0}^{\infty}2^{-j-1}\widehat{d}(Z_{t(j)}(\omega),Z'_{t(j)}(\omega))
\]
\begin{equation}
\leq2^{p(m(k))+1}\delta\varepsilon^{-1}\equiv2^{p(m(k))+1}2^{-\widetilde{p}-1}\varepsilon^{2}\varepsilon^{-1}=\varepsilon.\label{eq:temp-343}
\end{equation}

5. Now let $\omega\in H^{c}H'^{c}G^{c}$ be arbitrary, but fixed.
Let $i=0,\cdots,p+1$ be arbitrary. To simplify notations, we will
henceforth suppress the reference to $\omega$, and suppress the overline,
to write $\tau_{i},\tau'_{i},x,x',z,z'$ for $\overline{\tau}_{i}(\omega),\overline{\tau}'_{i}(\omega),X(\cdot,\omega),X'(\cdot,\omega),Z(\cdot,\omega),Z'(\cdot,\omega)$
respectively. Then inequality \ref{eq:temp-307-4-1-1} can be rewritten
more compactly as

\begin{equation}
x(\theta_{i})\subset(d(\cdot,z(q_{i-1}))\leq\varepsilon).\label{eq:temp-338}
\end{equation}
Similarly, inequality \ref{eq:temp-307-4-1-1-1} can be rewritten
as
\begin{equation}
x(\theta'_{i})\subset(d(\cdot,z'(q_{i-1}))\leq\varepsilon).\label{eq:temp-347}
\end{equation}

6. Partition the set $\{0,\cdots,p+1\}$ into two disjoint subsets
$A$ and $B$ such that (i) $|\tau_{i}-\tau'_{i}|<2\widetilde{\Delta}$
if $i\in A$, and (ii) $|\tau_{i}-\tau'_{i}|>\widetilde{\Delta}$
if $i\in B$. Consider each $i\in B$. We will verify that, then,
$1\leq i\leq p$ and 
\begin{equation}
d(z(q_{i}),z(q_{i-1}))\vee d(z'(q_{i}),z(q'_{i-1}))\leq6\varepsilon.\label{eq:temp-311}
\end{equation}
In view of Condition (ii), we may assume, without loss of generality,
that $\tau'_{i}-\tau_{i}>\widetilde{\Delta}\equiv2^{-m(k)}$. Then
there exists $u\in[\tau_{i},\tau'_{i})Q_{m(k)}.$ Since $[\tau_{0},\tau'_{0})=[0,0)=\phi$
and $[\tau_{p+1},\tau'_{p+1})=[1,1)=\phi$, it follows that $1\leq i\leq p$.
Moreover, \ref{eq:temp-308} and \ref{eq:temp-308-1} imply that 
\[
u\in[\tau_{i},\tau'_{i})\subset[q_{i-1},q_{i})\subset[\tau'_{i-1},\tau{}_{i+1}).
\]
Hence $u\in[\tau_{i},\tau{}_{i+1})\cap[\tau'_{i-1},\tau'_{i})\subset\theta_{i+1}\theta'_{i}.$
Therefore, using relations \ref{eq:temp-338} and \ref{eq:temp-347},
and inequality\ref{eq:temp-343}, we obtain
\[
d(z(q_{i}),z(q_{i-1}))
\]
\[
\leq d(z(q_{i}),z(u))+d(z(u),z'(u))+d(z'(u),z'(q_{i-1}))+d(z'(q_{i-1}),z(q_{i-1}))
\]
\[
\leq\varepsilon+\varepsilon+\varepsilon+\varepsilon=4\varepsilon,
\]
and so
\[
d(z'(q_{i}),z(q'_{i-1}))\leq d(z'(q_{i}),z(q_{i}))+d(z(q_{i}),z(q_{i-1}))+d(z(q_{i-1}),z'(q_{i-1}))
\]
\[
\leq\varepsilon+4\varepsilon+\varepsilon=6\varepsilon.
\]
Thus inequality \ref{eq:temp-311} is verified.

7. Define an increasing function $\lambda:[0,1]\rightarrow[0,1]$
by $\lambda\tau_{i}\equiv\tau'_{i}$ or $\lambda\tau_{i}\equiv\tau{}_{i}$
according as $i\in A$ or $i\in B$, for each $i=0,\cdots,p+1$, and
by linearity on $[\tau_{i},\tau_{i+1}]$ for each $i=0,\cdots,p$.
Here we write $\lambda t\equiv\lambda(t)$ for each $t\in[0,1]$ for
brevity. Then, in view of the definition of the index sets $A$ and
$B$, we have $|\tau_{i}-\lambda\tau_{i}|<2\widetilde{\Delta}$ for
each $i=0,\cdots,p+1$. Now consider each $i=0,\cdots,p$, and write
\[
u_{i}\equiv\frac{\lambda\tau_{i+1}-\tau_{i+1}+\tau_{i}-\lambda\tau_{i}}{\tau_{i+1}-\tau_{i}}.
\]
Then, since $\tau_{i+1}-\tau_{i}\geq\overline{\Delta}$ according
to inequality \ref{eq:temp-342-1}, we have
\[
|u_{i}|\leq|\lambda\tau_{i+1}-\tau_{i+1}|\overline{\Delta}^{-1}+|\lambda\tau_{i}-\tau_{i}|\overline{\Delta}^{-1}
\]
\[
\leq2\widetilde{\Delta}\overline{\Delta}^{-1}+2\widetilde{\Delta}\overline{\Delta}^{-1}=2^{-m(k)+2}\overline{\Delta}^{-1}<(1-e^{-\varepsilon}).
\]
Note that the function $\log(1+u)$ of $u\in[-1+e^{-\varepsilon},1-e^{-\varepsilon}]$
vanishes at $u=0$, and has positive first derivative and negative
second derivative on the interval $[-1+e^{-\varepsilon},1-e^{-\varepsilon}]$.
Hence the maximum of its absolute value is attained at the left end
point of the interval. Therefore
\begin{equation}
|\log(1+u_{i})|\leq|\log(1-1+e^{-\varepsilon})|=\varepsilon.\label{eq:temp-348-1}
\end{equation}
Lemma \ref{Lem. log((gt-gs)/(t-s))}, therefore implies the bound
\[
\sup_{0\leq s<t\leq1}|\log\frac{\lambda t-\lambda s}{t-s}|=\bigvee_{i=0}^{p}|\log\frac{\lambda\tau_{i+1}-\lambda\tau_{i}}{\tau_{i+1}-\tau_{i}}|
\]
\begin{equation}
=\bigvee_{i=0}^{p}|\log(1+u_{i})|\leq\varepsilon<9\varepsilon.\label{eq:temp-347-1}
\end{equation}

8. We will next prove that 
\begin{equation}
d(x,x'\circ\lambda)\leq9\varepsilon\label{eq:temp-245-3}
\end{equation}
on $domain(x)\cap domain(x'\circ\lambda)$. Now let $i=0,\cdots,p$
and 
\[
v\in\bigcup_{i=0}^{p}[\tau_{i},\tau_{i+1})\cap\bigcup_{i=0}^{p}(\lambda^{-1}\tau'_{i},\lambda^{-1}\tau'_{i+1})\cap domain(x)\cap domain(x'\circ\lambda)
\]
be arbitrary. There are four possible cases: (i') $i,i+1\in A$, (ii')
$i,i+1\in B$, (iii') $i\in A$ and $i+1\in B$, and (iv') $i\in B$
and $i+1\in A$. 

Consider Case (i'), where $i,i+1\in A$. Then $\lambda\tau_{i}\equiv\tau'_{i}$
and $\lambda\tau_{i+1}\equiv\tau'_{i+1}$. Hence 
\[
\lambda v\in[\lambda\tau_{i},\lambda\tau_{i+1})\subset[\tau'_{i},\tau'_{i+1})\subset\theta'_{i+1}.
\]
Consequently,
\[
d(x(v),x'(\lambda v))\leq d(x(v),z(q_{i}))+d(z(q_{i}),z'(q_{i}))+d(z'(q_{i}),x'(\lambda v))<\varepsilon+\varepsilon+\varepsilon=3\varepsilon,
\]
where we used relations \ref{eq:temp-338} and \ref{eq:temp-347},
and inequality \ref{eq:temp-343}. 

Next consider Case (ii'), where $i,i+1\in B$. Then, according to
Step 6, we have $1\leq i<i+1\leq p$ and 
\begin{equation}
d(z(q_{i}),z(q_{i-1}))\vee d(z'(q_{i}),z(q'_{i-1}))\vee d(z(q_{i+1}),z(q_{i}))\vee d(z'(q_{i+1}),z(q'_{i}))\leq6\varepsilon.\label{eq:temp-311-1}
\end{equation}
Moreover, $\lambda\tau_{i}\equiv\tau{}_{i}$ and $\lambda\tau_{i+1}\equiv\tau{}_{i+1}$.
Hence 
\[
\lambda v\in[\tau_{i},\tau{}_{i+1})\subset[q_{i-1},q_{i+1})\subset[\tau'_{i-1},\tau'_{i+2}).
\]
Since $\lambda v\neq\tau'_{i}$ and $\lambda v\neq\tau'_{i+1}$ by
assumption, we obtain
\[
\lambda v\in[\tau'_{i-1},\tau'_{i})\cup[\tau'_{i},\tau'_{i+1})\cup[\tau'_{i},\tau'_{i+2})\subset\theta'_{i}\cup\theta'_{i+1}\cup\theta'_{i+2}.
\]
Consequently,
\[
d(x(v),x'(\lambda v))
\]
\[
\leq(d(x(v),z(q_{i}))+d(z(q_{i}),z(q_{i-1}))+d(z(q_{i-1}),z'(q_{i-1}))+d(z'(q_{i-1}),x'(\lambda v)))
\]
\[
\vee(d(x(v),z(q_{i}))+d(z(q_{i}),z'(q_{i}))+d(z'(q_{i}),x'(\lambda v)))
\]
\[
\vee(d(x(v),z(q_{i}))+d(z(q_{i}),z(q_{i+1}))+d(z(q_{i+1}),z'(q_{i+1}))+d(z'(q_{i+1}),x'(\lambda v)))
\]
\[
\leq(\varepsilon+6\varepsilon+\varepsilon+\varepsilon)\vee(\varepsilon+\varepsilon+\varepsilon)\vee(\varepsilon+6\varepsilon+\varepsilon+\varepsilon)=9\varepsilon,
\]
where we used relations \ref{eq:temp-338} and \ref{eq:temp-347},
and inequalities \ref{eq:temp-343} and \ref{eq:temp-311-1}, . 

Now consider Case (iii'), where $i\in A$ and $i+1\in B$. Then, according
to Step 6, we have $i+1\leq p$ and 
\begin{equation}
d(z(q_{i+1}),z(q_{i}))\vee d(z'(q_{i+1}),z(q'_{i}))\leq6\varepsilon.\label{eq:temp-311-1-1}
\end{equation}
Moreover, $\lambda\tau_{i}\equiv\tau'{}_{i}$ and $\lambda\tau_{i+1}\equiv\tau{}_{i+1}$.
Hence 
\[
\lambda v\in[\tau'_{i},\tau{}_{i+1})\subset[\tau'_{i},q_{i+1})\subset[\tau'_{i},\tau'_{i+2}).
\]
Since $\lambda v\neq\tau'_{i+1}$ by assumption, we obtain
\[
\lambda v\in[\tau'_{i},\tau'_{i+1})\cup[\tau'_{i+1},\tau'_{i+2})\subset\theta'_{i+1}\cup\theta'_{i+2}.
\]
Consequently,
\[
d(x(v),x'(\lambda v))
\]
\[
\leq(d(x(v),z(q_{i}))+d(z(q_{i}),z'(q_{i}))+d(z'(q_{i}),x'(\lambda v)))
\]
\[
\vee(d(x(v),z(q_{i}))+d(z(q_{i}),z(q_{i+1}))+d(z(q_{i+1}),z'(q_{i+1}))+d(z'(q_{i+1}),x'(\lambda v)))
\]
\[
\leq(\varepsilon+6\varepsilon+\varepsilon+\varepsilon)\vee(\varepsilon+\varepsilon+\varepsilon)\vee(\varepsilon+6\varepsilon+\varepsilon+\varepsilon)=9\varepsilon,
\]
where we used relations \ref{eq:temp-338} and \ref{eq:temp-347},
and inequalities \ref{eq:temp-343} and  \ref{eq:temp-311-1-1}, . 

Finally, consider Case (iv'), where $i\in B$ and $i+1\in A$. Then,
according to Step 6, we have $1\leq i$ and 
\begin{equation}
d(z(q_{i}),z(q_{i=1}))\vee d(z'(q_{i}),z(q'_{i-1}))\leq6\varepsilon.\label{eq:temp-311-1-1-1}
\end{equation}
Moreover, $\lambda\tau_{i}\equiv\tau{}_{i}$ and $\lambda\tau_{i+1}\equiv\tau'{}_{i+1}$.
Hence 
\[
\lambda v\in[\tau_{i},\tau'{}_{i+1})\subset[q_{i-1},\tau'{}_{i+1})\subset[\tau'_{i-1},\tau'_{i+1}).
\]
Since $\lambda v\neq\tau'_{i}$ by assumption, we obtain
\[
\lambda v\in[\tau'_{i-1},\tau'_{i})\cup[\tau'_{i},\tau'_{i+1})\subset\theta'_{i}\cup\theta'_{i+1}.
\]
Consequently,
\[
d(x(v),x'(\lambda v))
\]
\[
\leq(d(x(v),z(q_{i}))+d(z(q_{i}),z(q_{i-1}))+d(z(q_{i-1}),z'(q_{i-1}))+d(z'(q_{i-1}),x'(\lambda v)))
\]
\[
\vee(d(x(v),z(q_{i}))+d(z(q_{i}),z'(q_{i}))+d(z'(q_{i}),x'(\lambda v)))
\]
\[
\leq(\varepsilon+6\varepsilon+\varepsilon+\varepsilon)\vee(\varepsilon+\varepsilon+\varepsilon)=9\varepsilon,
\]
where we used relations \ref{eq:temp-338} and \ref{eq:temp-347},
and inequalities \ref{eq:temp-343} and  \ref{eq:temp-311-1-1-1}.

Summing up, we see that 
\begin{equation}
d(x(v),x'\circ\lambda(v))\leq9\varepsilon\label{eq:temp-344}
\end{equation}
for each 
\begin{equation}
v\in\bigcup_{i=0}^{p}[\tau_{i},\tau_{i+1})\cap\bigcup_{i=0}^{p}(\lambda^{-1}\tau'_{i},\lambda^{-1}\tau'_{i+1})\cap domain(x)\cap domain(x'\circ\lambda).\label{eq:temp-348}
\end{equation}
Note that the set $\bigcup_{i=0}^{p}[\tau_{i},\tau_{i+1})\cap\bigcup_{i=0}^{p}(\lambda^{-1}\tau'_{i},\lambda^{-1}\tau'_{i+1})$
contain all but finitely many points in the interval $[0,1]$, while
the set $domain(x)\cap domain(x'\circ\lambda)$ contains all but countably
many points in $[0,1]$, according to Proposition \ref{Prop. Points of continuitys of cadlag functions-1}.
Hence the set on the right-hand side of expression \ref{eq:temp-348}
contains all but countably many points in $[0,1]$, and is therefore
dense in $domain(x)\cap domain(x'\circ\lambda)$. Therefore, by Lemma
\ref{Lem.  If x=00003Dy on dense A then x=00003Dy. If f(x,y)<=00003Dc on A then...-1},
inequality \ref{eq:temp-344} extends to the desired inequality \ref{eq:temp-245-3}.

9. Therefore, according to Definition \ref{Def. Skorokhod Metric}
of the Skorokhod metric $d_{D}$, inequalities \ref{eq:temp-347-1}
and \ref{eq:temp-245-3} together yield the bound
\[
\widehat{d_{D}}(X(\cdot,\omega),X'(\cdot,\omega))\equiv1\wedge d_{D}(X(\cdot,\omega),X'(\cdot,\omega))\equiv1\wedge d_{D}(x,x')\leq9\varepsilon,
\]
where $\omega\in H^{c}H'^{c}G^{c}$ is arbitrary. It follows that
\[
\rho_{\widehat{D}[0,1]}(X,X')\equiv\int E(d\omega)\widehat{d_{D}}(X(\cdot,\omega),X'(\cdot,\omega))
\]
\[
\leq P(H\cup H'\cup G)+\int E(d\omega)\widehat{d_{D}}(X(\cdot,\omega),X'(\cdot,\omega))1_{H^{c}H'^{c}G^{c}}\leq3\varepsilon+9\varepsilon=12\varepsilon<\varepsilon_{0},
\]
which proves inequality \ref{eq:temp-345} and the theorem.
\end{proof}
The next easy corollary of Theorems \ref{Thm. Extension of D-regular process by right limit is a.u.cadlag}
and \ref{Thm.Construction of a.u. cadlag process by right limit of D-regular proceses is continuous.}
gives a construction of an a.u. càdlàg process from a $D$-regular
family of f.j.d.'s, and shows that the construction given is uniformly
continuous on a set of such $D$-regular families which share a common
modulus of $D$-regularity and a common modulus of continuity in probability.

We will prove the corollary using the Daniell-Kolmogorov-Skorokhod
Extension. For that purpose, we fix the sample space
\[
(\Omega,L,E)\equiv(\Theta_{0},L_{0},E_{0})\equiv([0,1],L_{0},\int\cdot dx)
\]
to be the Lebesgue integration space based on the interval $\Theta_{0}\equiv[0,1]$.
Let $\xi\equiv(A_{q})_{q=1,2,\cdots}$ be an arbitrarily, but fixed,
binary approximation of the locally compact metric space $(S,d)$
relative to some fixed reference point $x_{\circ}\in S$.

Recall from Definition \ref{Def. Marginal metric} the metric space
$(\widehat{F}(Q_{\infty},S),\widehat{\rho}_{Marg,\xi,Q(\infty)})$
of consistent families of f.j.d.'s with parameter set $Q_{\infty}$
and state space $(S,d)$, where the marginal metric $\widehat{\rho}_{Marg,\xi,Q(\infty)}$
is defined relative to the binary approximation $\xi$ of $(S,d)$.
Recall from Definition \ref{Def.  D-regular process and D-regular family of f.j.d.s on Q_inf}
its subset $\widehat{F}_{Dreg}(Q_{\infty},S)$ consisting of $D$-regular
families, and the subset $\widehat{F}_{Dreg,\overline{m},\delta(Cp)}(Q_{\infty},S)$
consisting of $D$-regular families with some given modulus of $D$-regularity
$\overline{m}$ and some given modulus of continuity in probability
$\delta_{Cp}$. 

Similarly, recall from Definition \ref{Def. metric space of D-regular families of f.j.d.'s on =00005B0,1=00005D}
the metric space $(\widehat{F}_{Dreg}([0,1],S),\widehat{\rho}_{Cp,\xi,[0,1]|Q(\infty)})$
and its subset $\widehat{F}_{Dreg,\overline{m},\delta(Cp)}([0,1],S)$.

We re-emphasize that $D$-regularity is a condition on individual
f.j.d.'s. 
\begin{cor}
\label{Cor.  Construction of a.u. calafg processes from D-regular f.j.d.'s, and continuity}
\textbf{\emph{(Construction of a.u. càdlàg process on $[0,1]$ from
$D$-regular family of f.j.d.'s on $[0,1]$, and continuity of construction).}}
Let 
\[
(\Theta_{0},L_{0},I_{0})\equiv([0,1],L_{0},\int\cdot dx)
\]
denote the Lebesgue integration space based on the interval $\Theta_{0}\equiv[0,1]$.
Then there exists a uniformly continuous function 
\begin{equation}
\Phi_{aucl,\xi}:(\widehat{F}_{Dreg,\overline{m},\delta(Cp)}([0,1],S),\widehat{\rho}_{Cp,\xi,[0,1]|Q(\infty)})\rightarrow(\widehat{D}_{\delta(aucl),\delta(Cp)}[0,1],\rho_{\widehat{D}[0,1]})\label{eq:temp-351-3}
\end{equation}
such that, for each $F\in\widehat{F}_{Dreg,\overline{m},\delta(Cp)}([0,1],S)$,
the a.u. càdlàg process $X\equiv\Phi_{aucl,\xi}(F)$ has marginal
distributions given by $F$, and has the modulus of a.u. càdlàg $\delta_{aucl}(\cdot,\overline{m},\delta_{Cp})$
defined as in Theorem \ref{Thm. Extension of D-regular process by right limit is a.u.cadlag}. 

Moreover, the function $\Phi_{aucl,\xi}$ has a modulus of continuity
\[
\overline{\delta}_{aucl,\xi}\equiv\overline{\delta}_{aucl,\xi}(\cdot,\overline{m},\delta_{Cp},\left\Vert \xi\right\Vert )
\]
which depends only on $\overline{m},\delta_{Cp},$ and $\left\Vert \xi\right\Vert $. 
\end{cor}
\begin{proof}
1. Recall from Definition \ref{Def. Metric on  of continuous in prob families of finite joint distributions.}
the isometry
\[
\Phi_{[0,1]|Q(\infty)}:(\widehat{F}_{Cp}([0,1],S),\widehat{\rho}_{Cp,\xi,[0,1]|Q(\infty)})\rightarrow(\widehat{F}_{Cp}(Q_{\infty},S),\widehat{\rho}_{Marg,\xi,Q(\infty)})
\]
defined by $\Phi_{[0,1]|Q(\infty)}(F)\equiv F|Q_{\infty}$ for each
$F\in\widehat{F}_{Cp}([0,1]$. Since $\widehat{F}_{Dreg,\overline{m},\delta(Cp)}([0,1],S)\subset\widehat{F}_{Cp}([0,1],S)$,
we have an isometry
\[
\Phi_{[0,1]|Q(\infty)}:(\widehat{F}_{Dreg,\overline{m},\delta(Cp)}([0,1],S),\widehat{\rho}_{Cp,\xi,[0,1]|Q(\infty)})\rightarrow(\widehat{F}_{Cp}(Q_{\infty},S),\widehat{\rho}_{Marg,\xi,Q(\infty)}).
\]
Moreover, for each $F\in\widehat{F}_{Dreg,\overline{m},\delta(Cp)}([0,1],S)$,
we have $F|Q_{\infty}\in\widehat{F}_{Dreg,\overline{m},\delta(Cp)}(Q_{\infty},S)$,
we obtain the isometry
\begin{equation}
\Phi_{[0,1]|Q(\infty)}:(\widehat{F}_{Dreg,\overline{m},\delta(Cp)}([0,1],S),\widehat{\rho}_{Cp,\xi,[0,1]|Q(\infty)})\rightarrow(\widehat{F}_{Dreg,\overline{m},\delta(Cp)}(Q_{\infty},S),\widehat{\rho}_{Marg,\xi,Q(\infty)}).\label{eq:temp-203}
\end{equation}

2. Let 
\[
\Phi_{DKS,\xi}:(\widehat{F}_{Dreg,\overline{m},\delta(Cp)}(Q_{\infty},S),\widehat{\rho}_{Marg,\xi,Q(\infty)})\rightarrow(\widehat{R}(Q_{\infty}\times\Theta_{0},S),\widehat{\rho}_{Prob,Q(\infty)})
\]
be the Daniell-Kolmogorov-Skorokhod extension relative to the binary
approximation $\xi$ of $(S,d)$. Let $F'\in\widehat{F}_{Dreg,\overline{m},\delta(Cp)}(Q_{\infty},S)$
be arbitrary. Then the process $Z\equiv\Phi_{DKS,\xi}(F')$ has marginal
distributions given by $F'$. Hence, by Definition \ref{Def.  D-regular process and D-regular family of f.j.d.s on Q_inf},
the process $Z$, like its family $F'$ of marginal distributions,
is $D$-regular, with a\emph{ }modulus of $D$-regularity $\overline{m}$
and with a modulus of continuity in probability $\delta_{Cp}$ . In
other words, $\Phi_{DKS,\xi}(F')\in\widehat{R}_{Dreg,\overline{m},\delta(Cp)}(Q_{\infty}\times\Theta_{0},S)$.
Thus we obtain the continuous function
\begin{equation}
\Phi_{DKS,\xi}:(\widehat{F}_{Dreg,\overline{m},\delta(Cp)}(Q_{\infty},S),\widehat{\rho}_{Marg,\xi,Q(\infty)})\rightarrow(\widehat{R}_{Dreg,\overline{m},\delta(Cp)}(Q_{\infty}\times\Theta_{0},S),\widehat{\rho}_{Prob,Q(\infty)})\label{eq:temp-209}
\end{equation}

3. Let 
\begin{equation}
\Phi_{rLim}:(\widehat{R}_{Dreg,\overline{m},\delta(Cp)}(Q_{\infty}\times\Theta_{0},S),\widehat{\rho}_{Prob,Q(\infty)})\rightarrow(\widehat{D}_{\delta(aucl),\delta(Cp)}[0,1],\rho_{\widehat{D}[0,1]})\label{eq:temp-369-1-1}
\end{equation}
denote the extension by right limit, as in Theorem \ref{Thm. Extension of D-regular process by right limit is a.u.cadlag}.
By Theorem \ref{Thm.Construction of a.u. cadlag process by right limit of D-regular proceses is continuous.},
the mapping $\Phi_{rLim}$ is uniformly continuous, with a modulus
of continuity $\delta_{rLim}(\cdot,\overline{m},\delta_{Cp})$ which
depends only on $\overline{m}$ and $\delta_{Cp}$. 

4. Now define the composite $\Phi_{aucl,\xi}\equiv\Phi_{rLim}\circ\Phi_{DKS,\xi}\circ\Phi_{[0,1]|Q(\infty)}$
of the three mappings \ref{eq:temp-203}, \ref{eq:temp-209}, and
\ref{eq:temp-369-1-1}. Thus 
\[
\Phi_{aucl,\xi}:(\widehat{F}_{Dreg,\overline{m},\delta(Cp)}([0,1],S),\widehat{\rho}_{Cp,\xi,[0,1]|Q(\infty)})\rightarrow(\widehat{D}_{\delta(aucl),\delta(Cp)}[0,1],\rho_{\widehat{D}[0,1]})
\]
is the composite of three uniformly continuous functions. As such,
it is uniformly continuous, with a composite modulus of continuity
given by
\[
\overline{\delta}_{aucl,\xi}\equiv\overline{\delta}_{aucl,\xi}(\cdot,\overline{m},\delta_{Cp},\left\Vert \xi\right\Vert )\equiv\delta_{DKS}(\delta_{rLim}(\cdot,\overline{m},\delta_{Cp}),\left\Vert \xi\right\Vert ).
\]

5. Now consider each $F\in\widehat{F}_{Dreg,\overline{m},\delta(Cp)}([0,1],S)$.
Write $F'\equiv\Phi_{[0,1]|Q(\infty)}(F)\equiv F|Q_{\infty}$, write
$Z\equiv\Phi_{DKS,\xi}(F')$, and $X\equiv\Phi_{rLim}(Z)$. Then $X=\Phi_{aucl,\xi}(F)$.
Then the process $Z$ has marginal distributions given by $F'$. Let
$\overline{F}$ be the family of marginal distributions of the a.u.
càdlàg process $X$. Then, since $Z=X|Q_{\infty}$ by Assertion 2
of Proposition \ref{Prop. Right-limit extenaion of D-regular process is continuous a.u.},
we have 
\[
F|Q_{\infty}\equiv F'=\overline{F}|Q_{\infty}
\]
 Moreover, since a.u. càdlàg processes and $D$-regular families of
f.j.d.'s are continuous in probability, by definition, so are the
families $F$ and $\overline{F}$. Hence 
\[
\widehat{\rho}_{Cp,\xi,[0,1]|Q(\infty)}(F,\overline{F})\equiv\widehat{\rho}_{Marg,\xi,Q(\infty)}(F|Q_{\infty},\overline{F}|Q_{\infty})=0
\]
by the definition of the metric $\widehat{\rho}_{Cp,\xi,[0,1]|Q(\infty)}$.
Therefore $\overline{F}=F$. In other words, the a.u. càdlàg process
$X=\Phi_{aucl,\xi}(F)$ has marginal distributions given by the family
$F$, where $F\in\widehat{F}_{Dreg,\overline{m},\delta(Cp)}([0,1],S)$
is arbitrary. 

The Corollary is proved.
\end{proof}

\section{\label{Sec:Strong-Right-Continuity}Strong Right Continuity in Probability}

In this section, let $(S,d)$ be a locally compact metric space.

In the previous sections we proved that $D$-regularity is a necessary
and sufficient condition for a family of f.j.d.'s to be extendable
to an a.u. càdlàg process. We remarked that this method is a generalization
of the treatment of a Markov process with a Markov semigroup in \cite{Chan74}.
In the next chapter, we will show that this indeed is the case; we
will make precise the notion of a Markov process with a Markov semigroup,
and will show, by means of $D$-regularity, that they are a.u. càdlàg.
In preparation, we will presently introduce a sufficient condition
for $D$-regularity, which is easily verifiable for such Markov processes.
In addition, this sufficient condition will reduce, in the next section,
to a simple condition for a martingale with parameter set $[0,1]$
to be a.u. càdlàg. 

Said sufficient condition will consist of two subconditions. The first
subcondition will be called, for lack of a better name, \emph{strong
right continuity in probability}\index{strong right continuity in probability}.
The second subcondition will be called \emph{a.u. boundlessness }\index{a.u. boundlessness},.
We emphasize that both subconditions are on f.j.d.'s, or equivalently,
on finite samples of the process. 

We will note that the a.u. boundlessness condition will always be
satisfied if the locally compact metric state space $(S,d)$ is bounded,
e.g., if the metric $d$ is replaced by the equivalent metric $d(1+d)^{-1}$,
or if $(S,d)$ is embedded in the one-point compactification $(S\cup\{\Delta\},\overline{d})$,
where $\overline{d}\leq1$. 

Recall Definition \ref{Def. Notations for dyadic rationals} for the
enumerated set $Q_{\infty}$ of dyadic rationals in $[0,1]$, and
the enumerated subset
\[
Q_{h}\equiv\{0,\triangle_{h},2\triangle_{h},\cdots,1\}=\{t_{0},t_{1},t_{2},\cdots,t_{p(h)}\},
\]
where $p_{h}\equiv2^{h}$, $\triangle_{h}\equiv2^{-h}$, and where
$q_{h,i}\equiv i\triangle_{h}$ for each $i=0,\cdots,p_{h}$, for
each $h\geq0$. Recall also the miscellaneous notations and conventions
in Definition \ref{Def. Misc notations}, In addition, we will use
the following notations.
\begin{defn}
\label{Def. Notations for some natural filtrations and space of simple stopping times}
\textbf{(Natural filtrations and certain first exit times for a process
sampled at regular intervals).} In the remainder of the present section,
let $Z:Q_{\infty}\times(\Omega,L,E)\rightarrow(S,d)$ be an arbitrary
process with marginal distributions given by the family $F$ of f.j.d.'s. 

Let $h\geq0$ be arbitrary. Considered the process $Z|Q_{h}$, which
is the process $Z$ sampled at the regular interval of $\triangle_{h}$.
Let 
\[
\mathscr{\mathcal{L}}^{(h)}\equiv\{L^{(t,h)}:t\in Q_{h}\}
\]
denote the natural filtration of the process $Z|Q_{h}$. In other
words,
\[
L^{(t,h)}\equiv L(Z_{r}:r\in[0,t]Q_{h})
\]
for each $t\in Q_{h}$. Let $\tau$ be an arbitrary simple stopping
time with values in $Q_{h}$, relative to the filtration $\mathscr{\mathcal{L}}^{(h)}$.
Define the probability space
\[
L^{(\tau,h)}\equiv\{Y\in L:Y1_{(\tau\leq s)}\in L^{(s)}\mbox{ for each \ensuremath{s\in Q_{h}}}\}.
\]

For each $t\in Q_{h}$ and $\alpha>0$, recall from Part 2 of Definition
\ref{Def. Simple First Exit time}, the simple first exit time $\eta_{t,\alpha,[t,1]Q(h)}$ 

\[
\eta_{t,\alpha,[t,1]Q(h)}\equiv\sum_{r\in[t,1]Q(h)}r1_{(d(Z(t),Z(r))>\alpha)}\prod_{s\in[t,r)Q(h)}1_{(d(Z(t),Z(s))\leq\alpha)}
\]
\begin{equation}
+\prod_{s\in[t,1]Q(h)}1_{(d(Z(t),Z(s))\leq\alpha)}\label{eq:temp-142-2-3}
\end{equation}
for the process $Z|[t,1]Q_{h}$ to exit the closed $\alpha$-neighborhood
of $Z_{t}$. Here an empty product is, by convention, equal to $1$. 

Similarly, for each $\gamma>0$, define the r.r.v.
\[
\zeta_{h,\gamma}\equiv\sum_{r\in Q(h)}r1_{(d(x(\circ),Z(r))>\gamma)}\prod_{s\in[0,r)Q(h)}1_{d(x(\circ),Z(s))\leq\gamma}
\]
\begin{equation}
+\prod_{s\in Q(h)}1_{d(x(\circ),Z(s))\leq\gamma},\label{eq:temp-450}
\end{equation}
where we recall that $x_{\circ}$ is an arbitrary, but fixed reference
point in the state space $(S,d)$. It can easily be verified that
$\zeta_{h,\gamma}$ is a simple stopping time relative to the filtration
$\mathscr{\mathcal{L}}^{(h)}$ . Intuitively, $\zeta_{h,\gamma}$
is the first time $r\in Q_{h}$ when the process $Z|Q_{h}$ is outside
the bounded set $(d(x_{\circ},\cdot)\leq\gamma)$, with $\zeta_{h,\gamma}$
set to $1$ if no such $s\in Q_{h}$ exists.

Refer to Propositions \ref{Prop. simple stopping times} and \ref{Prop. Basics of simple Exit times}
for basic properties of simple stopping times and simple first exit
times. 
\end{defn}
$\square$
\begin{defn}
\label{Def. a.u. boundedness on Q_inf}\textbf{ (a.u. boundlessness
on $Q_{\infty}$).} Let $(S,d)$ be a locally compact metric space.
Suppose the process $Z:Q_{\infty}\times(\Omega,L,E)\rightarrow(S,d)$
is such that, for each $\varepsilon>0$, there exists $\beta_{auB}(\varepsilon)>0$
so large that 
\begin{equation}
P(\bigvee_{r\in Q(h)}d(x_{\circ},Z_{r})>\gamma)<\varepsilon\label{eq:temp-477-2}
\end{equation}
for each $h\geq0$, for each $\gamma>\beta_{auB}(\varepsilon)$. Then
we will say that the process $Z$ and the family $F$ of its marginal
distributions are \emph{a.u. bounded }\index{a.u. boundlessness},
with the operation $\beta_{auB}$ as a \emph{modulus of a.u. boundlessness
}\index{modulus of a.u. boundlessness}, relative to the reference
point $x_{\circ}\in S$. Note that this condition is trivially satisfied
if $d\leq1$, in which case we can take $\beta_{auB}(\varepsilon)\equiv1$
for each $\varepsilon>1$.
\end{defn}
$\square$
\begin{defn}
\label{Def. strongly  right continuity in probability on Q_inf}\textbf{
(Strong right continuity in probability on $Q_{\infty}$).} Let $(S,d)$
be a locally compact metric space. Let $Z:Q_{\infty}\times(\Omega,L,E)\rightarrow(S,d)$
be an arbitrary process. Suppose that, for each $\varepsilon,\gamma>0$,
there exists $\delta_{SRcp}(\varepsilon,\gamma)>0$ such that, for
arbitrary $h\geq0$ and $s,r\in Q_{h}$ with $s\leq r<s+\delta_{SRcp}(\varepsilon,\gamma)$,
we have 
\begin{equation}
P_{A}(d(Z_{s},Z_{r})>\alpha)\leq\varepsilon,\label{eq:temp-493-2}
\end{equation}
for each $\alpha>\varepsilon$ and for each $A\in L^{(s,h)}$ with
$A\subset(d(x_{\circ},Z_{s})\leq\gamma)$ and $P(A)>0$. Then we will
say that the process $Z$ and the family $F$ of its its marginal
distributions are \emph{strongly right continuous} \emph{in probability}\index{strong right continuity in probability},
with the operation $\delta_{SRcp}$ as a \emph{modulus of strong right
 continuity} \emph{in probability}\index{modulus of strong right continuity}. 

Note that the operation $\delta_{SRcp}$ has two variables. Suppose
it is independent of the second variable. Equivalently, suppose $\delta_{SRcp}(\varepsilon,\gamma)=\delta_{SRcp}(\varepsilon,1)$
for each $\varepsilon,\gamma>0$. Then we say that the process $Z$
and the family $F$ of its its marginal distributions are \emph{uniformly
strongly right continuous} \emph{in probability}\index{uniformly strong right continuity in probability}, 
\end{defn}
$\square$

The above definition will next be restated without the assumption
of $P(A)>0$ and the reference to the probability $P_{A}$.
\begin{lem}
\label{Lem. Equivalent def of strong right continuity in prob} \textbf{\emph{(Equivalent
definition of strong right continuity in probability). }}Let $(S,d)$
be a locally compact metric space. Then a process $Z:Q_{\infty}\times(\Omega,L,E)\rightarrow(S,d)$
is strongly right continuous in probability, with a modulus of strong
right  continuity in probability $\delta_{SRcp}$ iff, for each $\varepsilon,\gamma>0$,
there exists $\delta_{SRcp}(\varepsilon,\gamma)>0$ such that, for
arbitrary $h\geq0$ and $s,r\in Q_{h}$ with $s\leq r<s+\delta_{SRcp}(\varepsilon,\gamma)$,
we have 
\begin{equation}
P(d(Z_{s},Z_{r})>\alpha;A)\leq\varepsilon P(A),\label{eq:temp-493-2-1}
\end{equation}
for each $\alpha>\varepsilon$ and for each $A\in L^{(s,h)}$ with
$A\subset(d(x_{\circ},Z_{s})\leq\gamma).$
\end{lem}
\begin{proof}
Suppose Z is strongly right continuous in probability, with a modulus
of strong right  continuity in probability $\delta_{SRcp}$. Let $\varepsilon,\gamma>0$,
$h\geq0$, and $s,r\in Q_{h}$ be arbitrary with $s\leq r<s+\delta_{SRcp}(\varepsilon,\gamma)$.
Consider each $\alpha>0$ and $A\in L^{(s,h)}$ with $A\subset(d(x_{\circ},Z_{s})\leq\gamma).$
Suppose, for the sake of a contradiction, that 
\[
P(d(Z_{s},Z_{r})>\alpha;A)>\varepsilon P(A).
\]
Then $P(A)\geq P(d(Z_{s},Z_{r})>\alpha;A)>0$. Hence we can divide
both sides of the last displayed inequality by $P(A)$ to obtain
\[
P_{A}(d(Z_{s},Z_{r})>\alpha)\equiv P(A)^{-1}P(d(Z_{s},Z_{r})>\alpha;A)>\varepsilon,
\]
which contradicts inequality \ref{eq:temp-493-2} in Definition \ref{Def. strongly  right continuity in probability on Q_inf}
for the assumed strong right continuity. Hence inequality \ref{eq:temp-493-2-1}
holds. Thus the ``only if'' part of the lemma is proved. The ``if''
part is equally straightforward and omitted.
\end{proof}
Three more lemmas to prepare for the main theorem of this section.
The next two are elementary.
\begin{lem}
\textbf{\emph{\label{Lem. Min of a real and sum of two reals}(Minimum
of a real number and a sum of two real numbers). }}For each $a,b,c\in R$,
we have $a\wedge(b+c)=b+c\wedge(a-b)$, or, equivalently, $a\wedge(b+c)-b=c\wedge(a-b)$, 
\end{lem}
\begin{proof}
Write $c'\equiv a-b$. Then $a\wedge(b+c)=(b+c')\wedge(b+c)=b+c'\wedge c=b+(a-b)\wedge c$.
\end{proof}
\begin{lem}
\textbf{\emph{\label{Lem. Function on two contiguous intervals} (Function
on two contiguous intervals). }}\emph{Let} $Z:Q_{\infty}\times\Omega\rightarrow S$
be an arbitrary process. Let $\alpha>0$ and $\beta>2\alpha$ be arbitrary,
such that the set 
\begin{equation}
A_{r,s}^{\beta}\equiv(d(Z_{r},Z_{s})>\beta)\label{eq:temp-372-3-1-2-1}
\end{equation}
is \textup{\emph{measurable}}\emph{ }for each $r,s\in Q_{\infty}$.
Let $\omega\in\bigcap_{r,s\in Q(\infty)}(A_{r,s}^{\beta}\cup(A_{r,s}^{\beta})^{c})$
and let $h\geq0$ be arbitrary. Let $A_{\omega},B_{\omega}$ be arbitrary
intervals, with end points in $Q_{h}$ which depend on $\omega$,
such that the right end point of $A_{\omega}$ is equal to the left
end point of $B_{\omega}$. Let $t,t'\in(A_{\omega}\cup B_{\omega})Q_{h}$
be arbitrary such that $t<t'$. Suppose there exist $x_{\omega},y_{\omega}\in S$
with 
\begin{equation}
\bigvee_{r\in A(\omega)Q(h)}d(x_{\omega},Z(r,\omega))\vee\bigvee_{s\in B(\omega)Q(h)}d(y_{\omega},Z(s,\omega))\leq\alpha.\label{eq:temp-519-1-1}
\end{equation}
Then
\begin{equation}
\omega\in\bigcap_{r,s\in(t,t')Q(h);r\leq s}((A_{t,r}^{\beta}\cup A_{t,s}^{\beta})(A_{t,r}^{\beta}\cup A_{t',s}^{\beta})(A_{t',r}^{\beta}\cup A_{t',s}^{\beta}))^{c}.\label{eq:temp-425}
\end{equation}
\end{lem}
\begin{proof}
With $\omega$ fixed, write $z_{r}\equiv Z(r,\omega)\in S$ for each
$r\in Q_{\infty}$, and write $A\equiv A_{\omega}$, $B\equiv B_{\omega}$,
$x\equiv x_{\omega},$ and $y\equiv y_{\omega}$. Then inequality
\ref{eq:temp-519-1-1} can be restated as
\begin{equation}
\bigvee_{r\in AQ(h)}d(x,z_{r})\vee\bigvee_{s\in BQ(h)}d(y,z_{s})\leq\alpha.\label{eq:temp-519-1-2}
\end{equation}
Let $r,s\in(t,t')Q_{h}\subset AQ_{h}\cup BQ_{h}$ be arbitrary with
$r\leq s$. Note that, by assumption, the end points of the intervals
$A$ and $B$ are members of $Q_{h}.$ Then, since the right end point
of $A$ is equal to the left end point of $B$, there are only three
possibilities: (i) $t,r,s\in AQ_{h}$, (ii) $t,r\in AQ_{h}$ and $s,t'\in BQ_{h}$,
or (iii) $r,s,t'\in BQ_{h}$. In Case (i), inequality \ref{eq:temp-519-1-2}
implies that 
\[
d(z_{t},z_{r})\vee d(z_{t},z_{s})\leq(d(x,z_{t})+d(x,z_{r}))\vee(d(x,z_{t})+d(x,z_{s}))\leq2\alpha\vee2\alpha=2\alpha.
\]
Similarly, in Case (ii), we have 
\[
d(z_{t},z_{r})\vee d(z_{t'},z_{s})\leq(d(x,z_{t})+d(x,z_{r}))\vee(d(y,z_{t'})+d(y,z_{s}))\leq2\alpha\vee2\alpha=2\alpha.
\]
Similarly, in Case (iii), we have 
\[
d(z_{t'},z_{r})\vee d(z_{t'},z_{s})\leq(d(y,z_{t'})+d(y,z_{r}))\vee(d(y,z_{t'})+d(y,z_{s}))\leq2\alpha\vee2\alpha=2\alpha.
\]
Thus, in all cases, we have 
\[
(d(z_{t},z_{r})\vee d(z_{t},z_{s}))\wedge(d(z_{t},z_{r})\vee d(z_{t'},z_{s}))\wedge(d(z_{t'},z_{r})\vee d(z_{t'},z_{s}))\leq2\alpha\leq\beta,
\]
where $r,s\in(t,t')Q_{h}$ are arbitrary with $r\leq s$. Equivalently,
the desired relation \ref{eq:temp-425} holds.
\end{proof}
\begin{lem}
\label{Lem. Lower bound for mean time before exit} \textbf{\emph{(Lower
bound for mean waiting time before exit after a simple  stopping time).
}}Let $(S,d)$ be a locally compact metric space. Suppose the process
$Z:Q_{\infty}\times\Omega\rightarrow S$ is  strongly  right continuous
in probability, with a modulus of strong right  continuity in probability
$\delta_{SRcp}$.

Let $\varepsilon,\gamma>0$ be arbitrary, but fixed. Take any $m\geq0$
so large that 
\[
\Delta_{m}\equiv2^{-m}<\delta_{SRcp}(\varepsilon,\gamma).
\]
Let 
\[
h\geq m
\]
be arbitrary. Define the simple stopping time $\zeta_{h,\gamma}$
to be the first time when the process $Z|Q_{h}$ is outside the bounded
set $(d(x_{\circ},\cdot)\leq\gamma)$, as in Definition \ref{Def. Notations for some natural filtrations and space of simple stopping times}.
Then the following holds.

\emph{1.} Let the point $t\in Q_{h}$ be arbitrary. Let $\overline{\eta}$
be an arbitrary simple stopping time with values in $[t,t+\Delta_{m}]Q_{h}$,
relative to the natural filtration $\mathscr{\mathcal{L}}^{(h)}$
of the process $Z|Q_{h}$. Let $A\in L^{(\overline{\eta},h)}$ be
an arbitrary measurable set such that $A\subset(d(x_{\circ},Z_{\overline{\eta}})\leq\gamma)$.
Then
\begin{equation}
P(d(Z_{\overline{\eta}},Z_{1\wedge(t+\Delta(m))})>\alpha;A)\leq\varepsilon P(A)\label{eq:temp-519-2-4}
\end{equation}
for each $\alpha>\varepsilon$.

\emph{2.} Suppose $\varepsilon\leq2^{-2}$ , and let $\alpha>2\varepsilon$
be arbitrary. For abbreviation, write 
\[
\eta_{t}\equiv\eta_{t,\alpha,[t,1]Q(h)}
\]
for each $t\in Q_{h}$. Let $\tau$ be an arbitrary simple stopping
time with values in $Q_{h}$, relative to the filtration $\mathscr{\mathcal{L}}^{(h)}$.
Then the r.r.v.
\begin{equation}
\overline{\eta}_{\tau}\equiv\sum_{t\in Q(h)}(\eta_{t}1_{\eta(t)<\zeta(h,\gamma)}+1_{\zeta(h,\gamma)\leq\eta(t)})1_{(\tau=t)}\label{eq:temp-421}
\end{equation}
is a simple stopping time with values in $Q_{h}$, relative to the
filtration $\mathscr{\mathcal{L}}^{(h,\tau)},.$which we will loosely
call a simple first exit time after the given stopping time $\tau$.
Moreover, we have the upper bound
\begin{equation}
P(\overline{\eta}_{\tau}<1\wedge(\tau+\Delta_{m}))\leq2\varepsilon\label{eq:temp-381-1}
\end{equation}
for the probability of a quick first exit after $\tau$, and we have
the lower bound
\begin{equation}
E(\overline{\eta}_{\tau}-\tau)\geq2^{-1}E((1-\tau)\wedge\Delta_{m}).\label{eq:temp-415-1}
\end{equation}
for the mean waiting time before exit after $\tau$. 

We emphasize that the upper bound \ref{eq:temp-381-1} and the lower
bound \ref{eq:temp-415-1} are independent of $h$.
\end{lem}
\begin{proof}
1. Let $\alpha>\varepsilon$ be arbitrary. Let $t\in Q_{h},$ the
simple stopping time $\overline{\eta}$ with values in $[t,t+\Delta_{m}]Q_{h}$,
and the set $A\in L^{(\overline{\eta},h)}$ with $A\subset(d(x_{\circ},Z_{\overline{\eta}})\leq\gamma)$
be as given. Write $r\equiv1\wedge(t+\Delta_{m})$. Then $\overline{\eta}$
has values in $[t,r]Q_{h}$. Let $s\in[t,r]Q_{h}$ be arbitrary. Then
\[
s\leq r\leq t+\Delta_{m}<s+\delta_{SRcp}(\varepsilon,\gamma)
\]
and $(\overline{\eta}=s;A)\in L^{(s,h)}$. Therefore, we can apply
inequality \ref{eq:temp-493-2-1} in Lemma \ref{Lem. Equivalent def of strong right continuity in prob}
to the modulus of strong right continuity $\delta_{SRcp}$, the points
$s,r\in Q_{h}$, the simple stopping time $\overline{\eta},$ and
the $\mathrm{measurable}$ set $(\overline{\eta}=s;A)\in L^{(s,h)},$
to obtain
\[
P(d(Z_{s},Z_{r})>\alpha;\overline{\eta}=s;A)\leq\varepsilon P(\overline{\eta}=s;A),
\]
where $s\in[t,r]Q_{h}$ is arbitrary. Consequently,
\[
P(d(Z_{\overline{\eta}},Z_{r})>\alpha;A)=\sum_{s\in[t,r]Q(h)}P(d(Z_{s},Z_{r})>\alpha;\overline{\eta}=s;A)
\]
\[
\leq\sum_{s\in[t,r]Q(h)}\varepsilon P(\overline{\eta}=s;A)=\varepsilon P(A).
\]
Assertion 1 is proved. 

2. To prove Assertion 2, suppose $\varepsilon\leq2^{-2}$, and let$\alpha>2\varepsilon$
be arbitrary. First consider the special case where $\tau\equiv t$
for some $t\in Q_{h}$. Define the r.r.v.
\[
\overline{\eta}_{t}\equiv\eta_{t}1_{\eta(t)<\zeta(h,\gamma)}+1_{\zeta(h,\gamma)\leq\eta(t)}
\]
\begin{equation}
=\sum_{s\in[t,1]Q(h)}s1_{\eta(t)=s;s<\zeta(h,\gamma)}+\sum_{s\in[t,1]Q(h)}1_{\eta(t)=s;s\geq\zeta(h,\gamma)},\label{eq:temp-260}
\end{equation}
which has values in $[t,1]Q_{h}$. Note that (i') if $s\in[0,t)Q_{h}$
, then, trivially ,
\[
1_{(\overline{\eta}(t)=s)}=0\in L^{(s,h)},
\]
(ii') if $s\in[t,1)Q_{h}$, then 
\[
1_{(\overline{\eta}(t)=s)}=1_{\eta(t)=s;s<\zeta(h,\gamma)}\in L^{(s,h)}
\]
because $\eta_{t}$ and $\zeta_{h,\gamma}$ are simple stopping times
with values in $Q_{h}$, relative to the filtration $\mathscr{\mathcal{L}}^{(h)}$,
and (iii') if $s=1$, then 
\[
1_{(\overline{\eta}(t)=s)}=1_{\eta(t)=s;s\geq\zeta(h,\gamma)}=1_{\eta(t)=s}\in L^{(s,h)}.
\]
Combining, we see that $1_{(\overline{\eta}(t)=s)}\in L^{(s,h)}$
for each $s\in Q_{h}$. Thus $\overline{\eta}_{t}$ is a simple stopping
time with values in $[t,1]Q_{h}$, relative to the filtration $\mathscr{\mathcal{L}}^{(h)}$.
Intuitively, $\overline{\eta}_{t}$ is the first time $s\in[t,1]Q_{h}$
when the process $Z|Q_{h}$ exits the $\alpha$-neighborhood of $Z_{t}$
while staying in the $\gamma$-neighborhood of $x_{\circ}$ over the
entire time interval $[t,s]Q_{h}$; $\overline{\eta}_{t}$ is set
to $1$ if no such time $s$ exists. 

Continuing, observe that
\begin{equation}
(\overline{\eta}_{t}<1)\subset(\overline{\eta}_{t}=\eta_{t}<\zeta_{h,\gamma}\leq1)\label{eq:temp-462}
\end{equation}
by the defining equality \ref{eq:temp-260}.

Define $\overline{\eta}\equiv\overline{\eta}_{t}\wedge r$, where
$r\equiv1\wedge(t+\Delta_{m})$. Then $r-t=(1-t)\wedge\Delta_{m}$by
Lemma \ref{Lem. Min of a real and sum of two reals}. Moreover, the
simple stopping time $\overline{\eta}$ has values in $[t,r]Q_{h}$.
Let $A\in L^{(t,h)}$ be arbitrary. Then $A\in L^{(t,h)}\subset L^{(\overline{\eta},h)}$.
We estimate an upper bound for the probability
\[
P(\overline{\eta}_{t}<r;A)\leq P(\overline{\eta}=\overline{\eta}_{t}<1;A)
\]
\[
\leq P(\overline{\eta}_{t}<1;\overline{\eta}=\overline{\eta}_{t}=\eta_{t}<\zeta_{h,\gamma}\leq1;A)
\]
\[
=P(\eta_{t}<1;\overline{\eta}=\overline{\eta}_{t}=\eta_{t}<\zeta_{h,\gamma}\leq1;A)
\]
\[
\leq P(d(Z_{t},Z_{\eta(t)})>\alpha;\overline{\eta}=\overline{\eta}_{t}=\eta_{t}<\zeta_{h,\gamma}\leq1;A)
\]
\[
=P(d(Z_{t},Z_{\overline{\eta}})>\alpha;\overline{\eta}=\overline{\eta}_{t}=\eta_{t}<\zeta_{h,\gamma}\leq1;A)
\]
\[
\leq P(d(Z_{t},Z_{\overline{\eta}})>\alpha;t<\zeta_{h,\gamma};A)
\]
\[
\leq P(d(Z_{t},Z_{\overline{\eta}})>\alpha;d(x_{\circ},Z_{t})\leq\gamma;A)
\]
\[
\leq P(d(Z_{t},Z_{\overline{\eta}})>\alpha;d(Z_{t},Z_{r})\leq2^{-1}\alpha;d(x_{\circ},Z_{t})\leq\gamma;A)
\]
\[
+P(d(Z_{t},Z_{r})>2^{-1}\alpha;d(x_{\circ},Z_{t})\leq\gamma;A)
\]
\[
\leq P(d(Z_{\overline{\eta}},Z_{r})>2^{-1}\alpha;d(x_{\circ},Z_{t})\leq\gamma;A)
\]
\begin{equation}
+P(d(Z_{t},Z_{r})>2^{-1}\alpha;d(x_{\circ},Z_{t})\leq\gamma;A),\label{eq:temp-457}
\end{equation}
where the the second inequality is thanks relation \ref{eq:temp-462},
and the third is by the definition of the simple stopping time time
$\eta_{t}$, and where the fifth inequality is by the definition of
the simple stopping time $\zeta_{h,\gamma}$. Next note that, since
$2^{-1}\alpha>\varepsilon$, we have
\[
P(d(Z_{\overline{\eta}},Z_{r})>2^{-1}\alpha;d(x_{\circ},Z_{t})\leq\gamma;A)\leq\varepsilon P(d(x_{\circ},Z_{t})\leq\gamma;A)
\]
by applying inequality \ref{eq:temp-519-2-4} where $\alpha,$$A$
are replaced by $2^{-1}\alpha$, $(d(x_{\circ},Z_{t})\leq\gamma;A)$
respectively. Similarly, we have
\[
P(d(Z_{t},Z_{r})>2^{-1}\alpha;d(x_{\circ},Z_{t})\leq\gamma;A)\leq\varepsilon P(d(x_{\circ},Z_{t})\leq\gamma;A)
\]
by applying inequality \ref{eq:temp-519-2-4} where $\overline{\eta},\alpha,$$A$
are replaced by $t,2^{-1}\alpha$, $(d(x_{\circ},Z_{t})\leq\gamma;A)$
respectively. Combining, inequality \ref{eq:temp-457} can be continued
to yield
\begin{equation}
P(\overline{\eta}_{t}<r;A)\leq2\varepsilon P(d(x_{\circ},Z_{t})\leq\gamma;A).\label{eq:temp-150-3}
\end{equation}
Consequently,
\[
E(\overline{\eta}_{t}-t;A)\geq E(r-t;\overline{\eta}_{t}\geq r;A)=(r-t)P(\overline{\eta}_{t}\geq r;A)
\]
\[
=((1-t)\wedge\Delta_{m})(P(A)-P(\overline{\eta}_{t}<r;A))
\]
\[
\geq((1-t)\wedge\Delta_{m})(P(A)-2\varepsilon P(d(x_{\circ},Z_{t})\leq\gamma;A)),
\]
\begin{equation}
=((1-t)\wedge\Delta_{m})(E1_{A}-2\varepsilon E1_{(d(x(\circ),Z(t))\leq\gamma;A)}),\label{eq:temp-374-1}
\end{equation}
where the second inequality is by inequality \ref{eq:temp-150-3}. 

Taking $A\equiv\Omega$, inequalities \ref{eq:temp-150-3} and \ref{eq:temp-374-1}
become, respectively, 
\[
P(\overline{\eta}_{t}<1\wedge(t+\Delta_{m}))\leq2\varepsilon P(d(x_{\circ},Z_{t})\leq\gamma)\leq2\varepsilon,
\]
and 
\[
E(\overline{\eta}_{t}-t)\geq((1-t)\wedge\Delta_{m})(1-2\varepsilon)\geq2^{-1}((1-t)\wedge\Delta_{m}),
\]
where the second inequality is because $\varepsilon\leq2^{-2}$ by
assumption. Thus Assertion 2 is proved for the special case where
$\tau\equiv t$ for some $t\in Q_{h}$. 

3. To complete the proof of Assertion 2 for the general case, let
the simple stopping time $\tau$ be arbitrary, with values in $Q_{h}$,
relative to the filtration $\mathscr{\mathcal{L}}^{(h)}$. Then 
\[
P(\overline{\eta}_{\tau}<1\wedge(\tau+\Delta_{m}))=\sum_{t\in Q(h)}P(\overline{\eta}_{t}<1\wedge(t+\Delta_{m});\tau=t)
\]
\[
\leq\sum_{t\in Q(h)}2\varepsilon P(d(x_{\circ},Z_{t})\leq\gamma;\tau=t)
\]
\[
=\sum_{t\in Q(h)}2\varepsilon P(d(x_{\circ},Z_{\tau})\leq\gamma;\tau=t)
\]
\begin{equation}
=2\varepsilon P(d(x_{\circ},Z_{\tau})\leq\gamma)\leq2\varepsilon,\label{eq:temp-531-1-1-1}
\end{equation}
where the inequality is by applying inequality \ref{eq:temp-150-3}
to the $\mathrm{measurable}$ set $A\equiv(\tau=t)\in L^{(t,h)}$,
for each $t\in Q_{h}$. Similarly, 
\[
E(\overline{\eta}_{\tau}-\tau)=\sum_{t\in Q(h)}E(\overline{\eta}_{t}-t;\tau=t)
\]
\[
\geq\sum_{t\in Q(h)}((1-t)\wedge\Delta_{m})(E1_{(\tau=t)}-2\varepsilon E1_{(d(x(\circ),Z(t))\leq\gamma;\tau=t)})
\]
\[
=\sum_{t\in Q(h)}(E((1-t)\wedge\Delta_{m})1_{(\tau=t)}-2\varepsilon E((1-t)\wedge\Delta_{m})1_{(d(x(\circ),Z(t))\leq\gamma;\tau=t)})
\]
\[
=\sum_{t\in Q(h)}(E((1-\tau)\wedge\Delta_{m})1_{(\tau=t)}-2\varepsilon E((1-\tau)\wedge\Delta_{m})1_{(d(x(\circ),Z(t))\leq\gamma;\tau=t)})
\]
\[
=(E((1-\tau)\wedge\Delta_{m})-2\varepsilon E((1-\tau)\wedge\Delta_{m})1_{(d(x(\circ),Z(t))\leq\gamma)})
\]
\[
\geq(E((1-\tau)\wedge\Delta_{m})-2\varepsilon E((1-\tau)\wedge\Delta_{m}))
\]
\begin{equation}
\geq2^{-1}E((1-\tau)\wedge\Delta_{m}),\label{eq:temp-602}
\end{equation}
where the first inequality is by applying inequality \ref{eq:temp-374-1}
to the $\mathrm{measurable}$ set $A\equiv(\tau=t)\in L^{(t,h)}$,
for each $t\in Q_{h}$, and where the last inequality is because $1-2\varepsilon\geq2^{-1}$
by the assumption that $\varepsilon\leq2^{-2}$. Summing up, inequalities
\ref{eq:temp-531-1-1-1} and \ref{eq:temp-602} yield, respectively,
the desired inequalities \ref{eq:temp-381-1} and \ref{eq:temp-415-1}.
The lemma is proved. 
\end{proof}
\begin{thm}
\textbf{\emph{\label{Thm. a.u.boundeness and Strongly right continuity impliy D-regular}
(Strong right  continuity}}\textbf{ }\textbf{\emph{in probability
and a.u. boundlessness together imply $D$-regularity and extendability
by right limit to a.u. càdlàg process). }}Let $(S,d)$ be a locally
compact metric space. Then a process $Z:Q_{\infty}\times(\Omega,L,E)\rightarrow(S,d)$
Suppose the process $Z:Q_{\infty}\times\Omega\rightarrow S$ is  \emph{(i)}
a.u. bounded, with a modulus of a.u. boundlessness $\beta_{auB}$,
and \emph{(ii)} strongly  right continuous in probability, with a
modulus of strong right  continuity in probability $\delta_{SRcp}$.
Then the following holds.

\emph{1.} The process $Z$ is $D$-regular, with a modulus of continuity
in probability 
\[
\delta{}_{Cp}(\cdot,\beta_{auB},\delta_{SRcp})
\]
and has a modulus of $D$-regularity $\overline{m}\equiv\overline{m}(\beta_{auB},\delta_{SRcp})$. 

\emph{2.} The right-limit extension $X\equiv\Phi_{rLim}(Z):[0,1]\times\Omega\rightarrow S$
is an a.u. càdlàg process, with the same modulus of continuity in
probability $\delta{}_{Cp}(\cdot,\beta_{auB},\delta_{SRcp})$, and
with a modulus of a.u. càdlàg $\delta_{aucl}(\cdot,\beta_{auB},\delta_{SRcp}).$
Recall here the right-limit extension $\Phi_{rLim}$ from Definition
\ref{Def. Extension of process by right limit}.
\end{thm}
\begin{proof}
1. Condition (i), the a.u. boundlessness condition in the hypothesis,
says that for each $\varepsilon>0$, we have
\begin{equation}
P(\bigvee_{u\in Q(h)}d(x_{\circ},Z_{u})>\gamma)<\varepsilon\label{eq:temp-477}
\end{equation}
for each $h\geq0$, for each $\gamma>\beta_{auB}(\varepsilon)$. 

2. Let $\varepsilon>0$ and $\alpha\in(2^{-2}\varepsilon,2^{-1}\varepsilon)$
be arbitrary. Let 
\[
\gamma\equiv[\beta_{auB}(2^{-2}\varepsilon)]_{1},
\]
\[
m\equiv[-\log_{2}(1\wedge\delta_{SRcp}(2^{-2}\varepsilon,\gamma))]_{1},
\]
and

\[
\delta{}_{Cp}(\varepsilon)\equiv\delta{}_{Cp}(\varepsilon,\beta_{auB},\delta_{SRcp})\equiv\Delta_{m}\equiv2^{-m}<\delta_{SRcp}(2^{-2}\varepsilon,\gamma).
\]
We will verify that $\delta{}_{Cp}(\cdot,\beta_{auB},\delta_{SRcp})$
is a modulus of continuity in probability of the process $Z$. For
each $h\geq0$, define the $\mathrm{measurable}$ set 
\[
\overline{G}_{h}\equiv(\bigvee_{u\in Q(h)}d(x_{\circ},Z_{u})>\gamma).
\]
Then, since $\gamma>\beta_{auB}(2^{-2}\varepsilon)$, inequality \ref{eq:temp-477},
with $\varepsilon$ replaced by $2^{-2}\varepsilon,$ yields
\begin{equation}
P(\overline{G}_{h})<2^{-2}\varepsilon,\label{eq:temp-493}
\end{equation}
where $h\geq0$ is arbitrary.

3. Consider each $s,r\in Q_{\infty}$ with $|s-r|<\delta{}_{Cp}(\varepsilon)$.
Define the measurable set 
\[
D_{s,r}\equiv(d(Z_{s},Z_{r})\leq\alpha)\subset(d(Z_{s},Z_{r})\leq2^{-1}\varepsilon).
\]
First assume that $s\leq r$. Then
\begin{equation}
s\leq r<s+\delta{}_{Cp}(\varepsilon)\leq s+\delta_{SRcp}(2^{-2}\varepsilon,\gamma).\label{eq:temp-458}
\end{equation}
Take $h\geq m$ so large that $s,r\in Q_{h}$. Since $\alpha>2^{-2}\varepsilon$,
we can apply inequality \ref{eq:temp-493-2-1} in Lemma \ref{Lem. Equivalent def of strong right continuity in prob},
where $\varepsilon$ and $A$ are replaced by $2^{-2}\varepsilon$
and $(d(x_{\circ},Z_{s})\leq\gamma)$ respectively, to obtain 
\[
P(D_{s,r}^{c}\overline{G}_{h}^{c})\leq P(d(Z_{s},Z_{r})>\alpha;d(x_{\circ},Z_{s})\leq\gamma)
\]
\[
\leq2^{-2}\varepsilon P(d(x_{\circ},Z_{s})\leq\gamma)\leq2^{-2}\varepsilon.
\]
Combining with inequality \ref{eq:temp-493}, we obtain
\begin{equation}
P(D_{s,r}^{c})\leq P(D_{s,r}^{c}\overline{G}_{h}^{c})+P(\overline{G}_{h})\leq2^{-2}\varepsilon+2^{-2}\varepsilon=2^{-1}\varepsilon.\label{eq:temp-470}
\end{equation}
Consequently, 
\[
E1\wedge d(Z_{s},Z_{r})\leq P(D_{s,r}^{c})+E(1\wedge d(Z_{s},Z_{r});D_{s,r})\leq P(D_{s,r}^{c})+2^{-1}\varepsilon\leq\varepsilon
\]
By symmetry, the same inequality holds for each $s,r\in Q_{\infty}$
with $|s-r|<\delta{}_{Cp}(\varepsilon),$ where $\varepsilon>0$ is
arbitrary. Hence the process $Z$ is continuous in probability, with
a modulus of continuity in probability $\delta{}_{Cp}\equiv\delta{}_{Cp}(\cdot,\beta_{auB},\delta_{SRcp})$.
Thus the process $Z$ satisfies Condition 2 of Definition \ref{Def.  D-regular process and D-regular family of f.j.d.s on Q_inf}
to be $D$-regular. 

4. It remains to prove Condition 1 in Definition \ref{Def.  D-regular process and D-regular family of f.j.d.s on Q_inf}.
To that end, define $m_{-1}\equiv\kappa_{-1}\equiv0$. Let $n\geq0$
be arbitrary. Write
\[
\varepsilon_{n}\equiv2^{-n}.
\]
Take any 
\[
\gamma_{n}\in(\beta_{auB}(2^{-3}\varepsilon_{n}),\beta_{auB}(2^{-3}\varepsilon_{n})+1).
\]
Define the integers
\begin{equation}
\kappa_{n}\equiv[\kappa_{n-1}\vee-\log_{2}(1\wedge\delta_{Rcp}(2^{-3}\varepsilon_{n},\gamma_{n}))]_{1},\label{eq:temp-368-1}
\end{equation}
\begin{equation}
m_{n}\equiv\kappa_{n}\vee[m_{n-1}\vee-\log_{2}(1\wedge\delta_{Rcp}(2^{-\kappa(n)-n-6}\varepsilon_{n},\gamma_{n}))]_{1},\label{eq:temp-352-1}
\end{equation}
\[
K_{n}\equiv2^{\kappa(n)+n+3},
\]
\[
h_{n}\equiv m_{n+1},
\]
and the measurable set 
\begin{equation}
G_{n}\equiv(\bigvee_{u\in Q(h(n))}d(x_{\circ},Z_{u})>\gamma_{n}).\label{eq:temp-263}
\end{equation}
Then, since $\gamma_{n}>\beta_{auB}(2^{-3}\varepsilon_{n})$, inequality
\ref{eq:temp-477} in the hypothesis implies that
\begin{equation}
P(G_{n})<2^{-3}\varepsilon_{n}.\label{eq:temp-482}
\end{equation}
Moreover,
\[
h_{n}\equiv m_{n+1}>m_{n}\geq\kappa_{n}\geq n\geq0,
\]
and equality \ref{eq:temp-352-1} can be rewritten as 
\begin{equation}
m_{n}\equiv\kappa_{n}\vee[m_{n-1}\vee-\log_{2}(1\wedge\delta_{Rcp}(K_{n}^{-1}2^{-3}\varepsilon_{n},\gamma_{n}))]_{1}.\label{eq:temp-530-1-1}
\end{equation}
Furthermore, equalities \ref{eq:temp-368-1} and \ref{eq:temp-530-1-1}
imply, respectively,
\begin{equation}
\Delta_{\kappa(n)}\equiv2^{-\kappa(n)}<\delta_{Rcp}(2^{-3}\varepsilon_{n},\gamma_{n}),\label{eq:temp-351-1}
\end{equation}
and
\begin{equation}
\Delta_{m(n)}\equiv2^{-m(n)}<\delta_{Rcp}(K^{-1}2^{-3}\varepsilon_{n},\gamma_{n}).\label{eq:temp-494}
\end{equation}

5. Now let 
\[
\beta>\varepsilon_{n}
\]
be arbitrary such that the set
\begin{equation}
A_{r,s}^{\beta}\equiv(d(Z_{r},Z_{s})>\beta)\label{eq:temp-372-3-1-2}
\end{equation}
is $\mathrm{measurable}$ for each $r,s\in Q_{\infty}$. Define the
exceptional set 
\begin{equation}
D_{n}\equiv\bigcup_{t\in Q(m(n))}\bigcup_{r,s\in(t,t')Q(m(n+1));r\leq s}(A_{t,r}^{\beta}\cup A_{t,s}^{\beta})(A_{t,r}^{\beta}\cup A_{t',s}^{\beta})(A_{t',r}^{\beta}\cup A_{t',s}^{\beta}),\label{eq:temp-201-4-3-1-1}
\end{equation}
where, for each $t\in[0,1)Q{}_{m(n)}$, we write $t'\equiv t+2^{-m(n)}.$
We need only prove that $P(D_{n})<2^{-n}$, as required in Condition
1 of Definition \ref{Def.  D-regular process and D-regular family of f.j.d.s on Q_inf}.

6. For that purpose, let the simple first stopping time $\zeta\equiv\zeta_{h(n),\gamma(n)}$
be as in Definition \ref{Def. Notations for some natural filtrations and space of simple stopping times}.
Thus $\zeta$ is the first time $s\in Q_{h(n)}$ when the process
$Z|Q_{h(n)}$ is outside the $\gamma_{n}$-neighborhood of the reference
point $x_{\circ}$, with $\zeta$ set to $1$ if no such $s\in Q_{h}$
exists. Take any 
\[
\alpha\in(2^{-1}\varepsilon_{n},2^{-1}\beta).
\]
This is possible because $\beta>\varepsilon_{n}$ by assumption. 

7. Now let $t\in Q_{h(n)}$ be arbitrary, but fixed till further notice.
Define the simple first exit time
\[
\eta_{t}\equiv\eta_{t,\alpha,[t,1]Q(h(n))},
\]
and define the simple stopping time
\begin{equation}
\overline{\eta}_{t}\equiv\eta_{t}1_{\eta(t)<\zeta}+1_{\zeta\leq\eta(t)}=\sum_{s\in[t,1]Q(h(n))}s1_{\eta(t)=s;s<\zeta}+\sum_{s\in[t,1]Q(h(n))}s1_{\eta(t)=s;\zeta\leq s}\label{eq:temp-479}
\end{equation}
as a special case of formula \ref{eq:temp-421} in Lemma \ref{Lem. Lower bound for mean time before exit}. 

Define the trivial simple stopping time $\tau_{0}\equiv0$. Let $k=1,\cdots,K_{n}$
be arbitrary. Define the simple stopping time 
\[
\tau_{k}\equiv\overline{\eta}_{\tau(k-1)}\equiv\sum_{u\in Q(h(n))}\overline{\eta}_{u}1_{(\tau(k-1)=u)}
\]
relative to the natural filtration $\mathcal{L}^{(h)}$ of the process
$Z|Q_{h}$. Then $\tau_{k}\geq\tau_{k-1}$. Intuitively, $\tau_{1},\tau_{2},\cdots,\tau_{K(n)}$
are the successive stopping times when the process $Z|Q_{h}$ moves
away from the previous stopping state $Z_{\tau(k-1)}$ by a distance
of more than $\alpha$ while staying within the bounded set $(d(x_{\circ},\cdot)<\gamma_{n})$. 

In view of the inequality \ref{eq:temp-351-1} and the bound 
\[
\alpha>2^{-1}\varepsilon_{n}>2(2^{-3}\varepsilon_{n}),
\]
we can apply Part 2 of Lemma \ref{Lem. Lower bound for mean time before exit},
where $\varepsilon,\gamma,m,h,\tau,\Delta$ are replaced by 
\[
2^{-3}\varepsilon_{n},\gamma_{n},\kappa_{n},h_{n},\tau_{k-1},\Delta_{\kappa(n)}
\]
respectively. Then Lemma \ref{Lem. Min of a real and sum of two reals}
and inequality \ref{eq:temp-381-1} in Lemma \ref{Lem. Lower bound for mean time before exit}
together imply 
\[
P(\tau_{k}-\tau_{k-1}<(1-\tau_{k-1})\wedge\Delta_{\kappa(n)})\equiv P(\overline{\eta}_{\tau(k-1)}<\tau_{k-1}+(1-\tau_{k-1})\wedge\Delta_{\kappa(n)}))
\]
\begin{equation}
=P(\overline{\eta}_{\tau(k-1)}<1\wedge(\tau_{k-1}+\Delta_{\kappa(n)}))\leq2(2^{-3}\varepsilon_{n}),\label{eq:temp-379-1}
\end{equation}
while \ref{eq:temp-415-1} in Lemma \ref{Lem. Lower bound for mean time before exit}
yields
\begin{equation}
E(\tau_{k}-\tau_{k-1})\equiv E(\overline{\eta}_{\tau(k-1)}-\tau_{k-1})\geq2^{-1}E((1-\tau_{k-1})\wedge\Delta_{\kappa(n)}),\label{eq:temp-530}
\end{equation}
where $k=1,\cdots,K_{n}$ is arbitrary. Hence 
\[
1\geq E(\tau_{K(n)})=\sum_{k=1}^{K(n)}E(\tau_{k}-\tau_{k-1})
\]
\[
\geq2^{-1}\sum_{k=1}^{K(n)}E((1-\tau_{k-1})\wedge\Delta_{\kappa(n)})\geq2^{-1}\sum_{k=1}^{K(n)}E((1-\tau_{K(n)-1})\wedge\Delta_{\kappa(n)})
\]
\[
=2^{-1}K_{n}E((1-\tau_{K(n)-1})\wedge\Delta_{\kappa(n)})
\]
\[
\geq2^{-1}K_{n}E((1-\tau_{K(n)-1})\wedge\Delta_{\kappa(n)};1-\tau_{K(n)-1}>\Delta_{\kappa(n)})
\]
\[
=2^{-1}K_{n}E(\Delta_{\kappa(n)};1-\tau_{K(n)-1}>\Delta_{\kappa(n)})=2^{-1}K_{n}\Delta_{\kappa(n)}P(1-\tau_{K(n)-1}>\Delta_{\kappa(n)})
\]
\[
=2^{-1}K_{n}\Delta_{\kappa(n)}P(\tau_{K(n)-1}<1-\Delta_{\kappa(n)}),
\]
where the second inequality is from inequality \ref{eq:temp-530}.
Dividing by $2^{-1}K_{n}\Delta_{\kappa(n)}$, we obtain
\begin{equation}
P(\tau_{K(n)-1}<1-\Delta_{\kappa(n)})<2K_{n}^{-1}\Delta_{\kappa(n)}^{-1}\equiv2\cdot2^{-\kappa(n)-n-3}2^{\kappa(n)}=2^{-2}\varepsilon_{n}.\label{eq:temp-532-1-1}
\end{equation}
At the same time,
\[
P(\tau_{K(n)}<1;\tau_{K(n)-1}\geq1-\Delta_{\kappa(n)})\equiv P(\overline{\eta}_{\tau(K(n)-1)}<1;\tau_{K(n)-1}\geq1-\Delta_{\kappa(n)})
\]
\begin{equation}
\leq P(\overline{\eta}_{\tau(K(n))-1)}<1\wedge(\tau_{K(n)-1}+\Delta_{\kappa(n)}))\leq2(2^{-3}\varepsilon_{n})=2^{-2}\varepsilon_{n},\label{eq:temp-531-1-1}
\end{equation}
where the last inequality is by applying inequality \ref{eq:temp-379-1}
to $k=K_{n}$. 

8. Next define the exceptional set 
\[
H_{n}\equiv(\tau_{K(n)}<1).
\]
Then, combining inequalities \ref{eq:temp-531-1-1} and \ref{eq:temp-532-1-1},
we obtain
\[
P(H_{n})\equiv P(\tau_{K(n)}<1)
\]
\[
\leq P(\tau_{K(n)}<1;\tau_{K(n)-1}\geq1-\Delta_{\kappa(n)})+P(\tau_{K(n)-1}<1-\Delta_{\kappa(n)})
\]
\[
\leq2^{-2}\varepsilon_{n}+2^{-2}\varepsilon_{n}=2^{-1}\varepsilon_{n}.
\]
Summing up, except for the small exceptional set $G_{n}\cup H_{n}$,
there are at most the finite number $K_{n}$ of simple stopping times
$0<\tau_{1}<\cdots<\tau_{K(n)}=1$ each of which is the first time
in $Q_{h(n)}$ when the process $Z$ strays from the previous stopped
state by a distance greater than $\alpha$, while staying in the bounded
set $(d(x_{\circ},\cdot)<\gamma_{n})$. At the same time, inequality
\ref{eq:temp-379-1} says that each waiting time $\tau_{k}-\tau_{k-1}$
exceeds a certain lower bound with some probability close to $1$.
We wish, however, to be able to say that, with some probability close
to $1$, all these $K_{n}$ waiting times simultaneously exceed a
certain lower bound. For that purpose, we will relax the lower bound
and specify two more small exceptional sets, as follows.

9. Define two more exceptional sets, 
\begin{equation}
B_{n}\equiv\bigcup_{k=1}^{K(n)}(\tau_{k}-\tau_{k-1}<(1-\tau_{k-1})\wedge\Delta_{m(n)}),\label{eq:temp-378-1}
\end{equation}
\begin{equation}
C_{n}\equiv\bigcup_{k=1}^{K(n)}(d(Z_{\tau(k-1)\vee(1-\Delta(m(n))},Z_{1})>\alpha),\label{eq:temp-419-1}
\end{equation}
and proceed to estimate $P(B_{n})$ and $P(C_{n})$. 

First, let $k=1,\cdots,K_{n}$ be arbitrary. Note that, trivially,
\[
\alpha>2^{-1}\varepsilon_{n}>2K_{n}^{-1}2^{-3}\varepsilon_{n}.
\]
and that, as in inequality \ref{eq:temp-494},
\[
\Delta_{m(n)}\equiv2^{-m(n)}<\delta_{Rcp}(K_{n}^{-1}2^{-3}\varepsilon_{n},\gamma_{n}).
\]
Now define the numbers, $m\equiv m_{n}$, $t\equiv1-\Delta_{m(n)}$,
$r\equiv1\wedge(t+\Delta_{m(n)})=1$, the simple stopping times 
\[
\overline{\eta}\equiv\tau_{k-1}\vee t\equiv\tau_{k-1}\vee(1-\Delta_{m(n)})
\]
and $\tau\equiv\tau_{k-1}$ with values in $Q_{h(n)}$, and the $\mathrm{measurable}$
set $(d(x_{\circ},Z_{\overline{\eta}})\leq\gamma_{n})\in L^{(\overline{\eta},h(n))}$.
Then $\alpha>2\varepsilon$ and 
\[
\Delta_{m(n)}\equiv2^{-m(n)}<\delta_{Rcp}(\varepsilon,\gamma)
\]
according to inequality \ref{eq:temp-494}. Moreover, the simple stopping
time $\overline{\eta}$ has values in $[t,r]Q_{h(n)}$, relative to
the filtration $\mathscr{\mathcal{L}}^{(h)}$. Hence we can apply
Lemma \ref{Lem. Lower bound for mean time before exit}, where $\varepsilon,\gamma,\alpha,m,h,t,r,\overline{\eta}\tau,A$
are replaced by $K_{n}^{-1}2^{-3}\varepsilon_{n},\gamma_{n},\alpha,m_{n},h_{n},t,r,\overline{\eta},\tau_{k-1},(d(x_{\circ},Z_{\overline{\eta}})\leq\gamma_{n})$
respectively. Then inequality \ref{eq:temp-519-2-4} of Lemma \ref{Lem. Lower bound for mean time before exit}
yields
\[
P(d(Z_{\overline{\eta}},Z_{1})>\alpha;d(x_{\circ},Z_{\overline{\eta}})\leq\gamma_{n})\leq K_{n}^{-1}2^{-3}\varepsilon_{n}P(d(x_{\circ},Z_{\overline{\eta}})\leq\gamma_{n})\leq K_{n}^{-1}2^{-3}\varepsilon_{n}.
\]
Hence, recalling the defining equalities \ref{eq:temp-419-1} and
\ref{eq:temp-263} for the measurable sets $C_{n}$ and $G_{n}$ respectively,
we immediately obtain 
\[
P(C_{n}G_{n}^{c})\equiv P(\bigcup_{k=1}^{K(n)}(d(Z_{\tau(k-1)\vee(1-\Delta(n))},Z_{1})>\alpha;\bigvee_{u\in Q(h(n))}d(x_{\circ},Z_{u})\leq\gamma_{n}))
\]
\[
\leq P(\bigcup_{k=1}^{K(n)}(d(Z_{\overline{\eta}},Z_{1})>\alpha;d(x_{\circ},Z_{\overline{\eta}})\leq\gamma_{n}))\leq\sum_{k=1}^{K(n)}K_{n}^{-1}2^{-3}\varepsilon_{n}=2^{-3}\varepsilon_{n}.
\]
At the same time, \ref{eq:temp-381-1} in Lemma \ref{Lem. Lower bound for mean time before exit}
yields

\begin{equation}
P(\overline{\eta}_{\tau(k-1)}<1\wedge(\tau_{k-1}+\Delta_{m(n)}))\leq2K_{n}^{-1}2^{-3}\varepsilon_{n}.\label{eq:temp-381-1-1}
\end{equation}
Since $\overline{\eta}_{\tau(k-1)}\equiv\tau_{k}$ ,and since 
\[
1\wedge(\tau_{k-1}+\Delta_{m(n)})=\tau_{k-1}+\Delta_{m(n)}\wedge(1-\tau_{k-1})
\]
according to Lemma \ref{Lem. Min of a real and sum of two reals},
inequality \ref{eq:temp-381-1-1} is equivalent to
\begin{equation}
P(\tau_{k}-\tau_{k-1}<\Delta_{m(n)}\wedge(1-\tau_{k-1}))\leq2K_{n}^{-1}2^{-3}\varepsilon_{n}.\label{eq:temp-531-1-1-2}
\end{equation}
Hence, recalling the defining equality \ref{eq:temp-378-1} for the
measurable sets $B_{n}$, we obtain
\[
P(B_{n})\equiv P(\bigcup_{k=1}^{K(n)}(\tau_{k}-\tau_{k-1}<\Delta_{n}\wedge(1-\tau_{k-1})))\leq\sum_{k=1}^{K(n)}2K_{n}^{-1}2^{-3}\varepsilon_{n}=2^{-3}\varepsilon_{n}.
\]
Combining, we see that 
\[
P(G_{n}\cup H_{n}\cup B_{n}\cup C_{n})=P(G_{n}\cup H_{n}\cup B_{n}\cup C_{n}G_{n}^{c})\leq2^{-n-3}+2^{-1}\varepsilon_{n}+2^{-3}\varepsilon_{n}+2^{-3}\varepsilon_{n}<\varepsilon_{n}.
\]

10. Finally, we will prove that $D_{n}\subset G_{n}\cup H_{n}\cup B_{n}\cup C_{n}$.
To that end, consider each $\omega\in G_{n}^{c}H_{n}^{c}B_{n}^{c}C_{n}^{c}$.
Then, since $\omega\in G_{n}^{c},$ we have $\bigvee_{t\in Q(h(n))}d(x_{\circ},Z_{t}(\omega))\leq\gamma_{n}$.
Consequently, $\zeta_{h(n),\gamma(n)}(\omega)=1$ according to the
defining equality \ref{eq:temp-450}. Hence, by the defining equality
\ref{eq:temp-479}, we have 
\begin{equation}
\overline{\eta}_{t}(\omega)\equiv\eta_{t}(\omega)1_{\eta(t,\omega)<1}+1_{1\leq\eta(t,\omega)}=\eta_{t}(\omega)\label{eq:temp-480}
\end{equation}
for each $t\in Q_{h(n)}$. 

Separately, since $\omega\in H_{n}^{c}$, we have $\tau_{K(n)}(\omega)=1$.
Let $k=1,\cdots,K_{n}$ be arbitrary. Write $t\equiv\tau_{k-1}(\omega)$
. Then
\[
\tau_{k}(\omega)\equiv\overline{\eta}_{\tau(k-1)}(\omega)=\overline{\eta}_{t}(\omega)=\eta_{t}(\omega)\equiv\eta_{t,\alpha,[t,1]Q(h(n))}(\omega),
\]
where the second equality is by equality \ref{eq:temp-480}. Hence
basic properties of the simple first exit time $\eta_{t,\alpha,[t,1]Q(h(n))}$
implies that
\begin{equation}
d(Z(t,\omega),Z(u,\omega))\leq\alpha\label{eq:temp-519-1}
\end{equation}
for each $u\in[t,\tau_{k}(\omega))Q_{h}$. In other words,
\begin{equation}
d(Z(\tau_{k-1}(\omega),\omega),Z(u,\omega))\leq\alpha\label{eq:temp-519}
\end{equation}
for each $u\in[\tau_{k-1}(\omega),\tau_{k}(\omega))Q_{h}$, where
$k=1,\cdots,K_{n}$ is arbitrary. 

Next, let $t\in[0,1)Q_{m(n)}$ be arbitrary, and write $t'\equiv t+\Delta_{m(n)}$.
Consider the following two possible cases (i') and (ii') regarding
the number of members in the sequence $\tau_{1}(\omega),\cdots,\tau_{K(n)-1}(\omega)$
that are in the interval $(t,t']$.

(i'). Suppose the interval $(t,t']$ contains two or more members
in the sequence $\tau_{1}(\omega),\cdots,\tau_{K-1}(\omega)$. Then
there exists $k=1,\cdots,K_{n}-1$ such that 
\begin{equation}
\tau_{k-1}(\omega)\leq t<\tau_{k}(\omega)\leq\tau_{k+1}(\omega)\leq t'.\label{eq:temp-381}
\end{equation}
It follows that
\begin{equation}
\Delta_{m(n)}\equiv t'-t>\tau_{k+1}(\omega)-\tau_{k}(\omega)\geq\Delta_{m(n)}\wedge(1-\tau_{k}(\omega)),\label{eq:temp-606}
\end{equation}
where the last inequality is because $\omega\in B_{n}^{c}$. Consequently
$\Delta_{m(n)}>1-\tau_{k}(\omega),$ Hence $\Delta_{m(n)}\wedge(1-\tau_{k}(\omega))=1-\tau_{k}(\omega)$.
Therefore the second half of inequality \ref{eq:temp-606} yields
\[
\tau_{k+1}(\omega)-\tau_{k}(\omega)\geq1-\tau_{k}(\omega).
\]
which is equivalent to $1=\tau_{k+1}(\omega)$. Consequently, inequality
\ref{eq:temp-381} implies that $t'=1$ and that $\tau_{k}(\omega)>t=1-\Delta_{m(n)}$.
Therefore 
\begin{equation}
d(Z_{\tau(k)}(\omega),Z_{1}(\omega))=d(Z_{\tau(k)\vee(1-\Delta(m(n))}(\omega),Z_{1}(\omega))\leq\alpha,\label{eq:temp-465-1}
\end{equation}
where the inequality is because $\omega\in C_{n}^{c}$. At the same
time, 
\[
[\tau_{k}(\omega),1)Q_{h(n)}=[\tau_{k}(\omega),\tau_{k+1}(\omega))Q_{h(n)}
\]
\[
\equiv[\tau_{k}(\omega),\overline{\eta}_{\tau(k)}(\omega))Q_{h(n)}=[\tau_{k}(\omega),\eta_{\tau(k)}(\omega))Q_{h(n)},
\]
where the last equality follows from equality \ref{eq:temp-480}.
Hence, basic properties of the simple fist exit time $\eta_{\tau(k)}$
implies that
\begin{equation}
d(Z(\tau_{k}(\omega),\omega),Z(u,\omega))\leq\alpha\label{eq:temp-460-1}
\end{equation}
for each $u\in[\tau_{k}(\omega),1)Q_{h(n)}$. Combining with inequality
\ref{eq:temp-465-1} for the end point $1$, we see that 
\begin{equation}
d(Z(\tau_{k}(\omega),\omega),Z(u,\omega))\leq\alpha.\label{eq:temp-529}
\end{equation}
for each $u\in[\tau_{k}(\omega),1]Q_{h(n)}$. 

Note that $\beta>2\alpha$, and that $t,t'\in[\tau_{k-1}(\omega),\tau_{k}(\omega))\cup[\tau_{k}(\omega),t']$,
with $t'=1$. Hence, in view of inequalities \ref{eq:temp-519} and
\ref{eq:temp-529}, we can apply Lemma \ref{Lem. Function on two contiguous intervals}
to the contiguous intervals $A_{\omega}\equiv[\tau_{k-1}(\omega),\tau_{k}(\omega))$
and $B_{\omega}\equiv[\tau_{k}(\omega),t']=[\tau_{k}(\omega),1]$,
to obtain
\[
\omega\in\bigcap_{r,s\in(t,t')Q(h(n));r\leq s}((A_{t,r}^{\beta}\cup A_{t,s}^{\beta})(A_{t,r}^{\beta}\cup A_{t',s}^{\beta})(A_{t',r}^{\beta}\cup A_{t',s}^{\beta}))^{c}
\]
\begin{equation}
\equiv\bigcap_{r,s\in(t,t')Q(m(n+1));r\leq s}((A_{t,r}^{\beta}\cup A_{t,s}^{\beta})(A_{t,r}^{\beta}\cup A_{t',s}^{\beta})(A_{t',r}^{\beta}\cup A_{t',s}^{\beta}))^{c}\subset D_{n}^{c}.\label{eq:temp-425-2}
\end{equation}

(ii') Now suppose the interval $(t,t']$ contains zero or one member
in the sequence $\tau_{1}(\omega),\cdots,\tau_{K(n)-1}(\omega)$.
Then there exists $k=1,\cdots,K_{n}-1$ such that $\tau_{k-1}(\omega)\leq t<t'\leq\tau_{k+1}(\omega)$.
Hence
\[
t,t'\in[\tau_{k-1}(\omega),\tau_{k}(\omega))\cup[\tau_{k}(\omega),\tau_{k+1}(\omega)).
\]
Then inequality \ref{eq:temp-519} holds for both $k$ and $k+1$,
Hence we can apply Lemma \ref{Lem. Function on two contiguous intervals}
to the contiguous intervals $A_{\omega}\equiv[\tau_{k-1}(\omega),\tau_{k}(\omega))$
and $B_{\omega}\equiv[\tau_{k}(\omega),\tau_{k+1}(\omega))$, to obtain,
again, relation \ref{eq:temp-425-2}. 

11. Summing up, we have proved that $\omega\in D_{n}^{c}$ for each
$\omega\in G_{n}^{c}H_{n}^{c}B_{n}^{c}C_{n}^{c}$. Consequently $D_{n}\subset G_{n}\cup H_{n}\cup B_{n}\cup C_{n}$
, whence 
\[
P(D_{n})\leq P(G_{n}\cup H_{n}\cup B_{n}\cup C_{n})<\varepsilon_{n}\equiv2^{-n},
\]
where $n\geq0$ is arbitrary, Condition 1 of Definition \ref{Def.  D-regular process and D-regular family of f.j.d.s on Q_inf}
has also been verified. Accordingly, the process $Z$ is $D$-regular,
with the sequence $\overline{m}\equiv(m_{n})_{n=0.1.\cdots}$ as a
modulus of $D$-regularity. Assertion 1 is proved.

12. Therefore, by Theorem \ref{Thm. Extension of D-regular process by right limit is a.u.cadlag},
the right limit extension $X$ of the process $Z$ is a.u. càdlàg,
with a modulus of a.u. càdlàg $\delta_{aucl}(\cdot,\overline{m},\delta_{cp})\equiv\delta_{aucl}(\cdot,\beta_{auB},\delta_{Rcp})$.
Assertion 2 is proved..
\end{proof}

\section{A Sufficient Condition for an a.u. Càdlàg Martingale}

Using Theorem \ref{Thm. a.u.boundeness and Strongly right continuity impliy D-regular}
in the preceding section, we will prove a sufficient condition for
a martingale $X:[0,1]\times\Omega\rightarrow R$ to be equivalent
to one which is a.u. càdlàg. 

To that end, recall, from Definition \ref{Def. Special convex function},
the special convex function $\overline{\lambda}:R\rightarrow R$ given
by
\begin{equation}
\overline{\lambda}(x)\equiv2x+(e^{-|x|}-1+|x|)\label{eq:temp-524-1-1}
\end{equation}
for each $x\in R$. Theorem \ref{Thm. Bishop's maximal inequality for w.s. submartingales.}
says that the function $\overline{\lambda}$ is increasing and strictly
convex, with 
\begin{equation}
|x|\leq|\overline{\lambda}(x)|\leq3|x|\label{eq:temp-454}
\end{equation}
for each $x\in R$. 
\begin{lem}
\label{Lem. Each w.s. submartingale on  Q_inf is a.u. bounded} \textbf{\emph{(Each
wide-sense submartingale on $Q_{\infty}$ is a.u. bounded).}} Let
$Z:Q_{\infty}\times\Omega\rightarrow R$ be a wide-sense submartingale
relative to some filtration $\mathcal{L}\equiv\{L^{(t)}:t\in Q_{\infty}\}$.
Suppose $b>0$ is an upper bound of $E|Z_{0}|\vee E|Z_{1}|\leq b$.
Then the process $Z$ is a.u. bounded in the sense of Definition \ref{Def. a.u. boundedness on Q_inf},
with a modulus of a.u. boundlessness $\beta_{auB}\equiv\beta_{auB}(\cdot,b)$. 
\end{lem}
\begin{proof}
Let $\varepsilon>0$ be arbitrary. Take an arbitrary $\alpha\in(2^{-2}\varepsilon,2^{-1}\varepsilon)$.
Take an arbitrary real number $K>0$ so large that
\begin{equation}
6b<\frac{1}{6}\alpha^{3}K\exp(-3K^{-1}b\alpha^{-1}).\label{eq:temp-495}
\end{equation}
Such a real number $K$ exists because the right-hand side of the
inequality \ref{eq:temp-495} is arbitrarily large for sufficiently
large $K$. Define 
\[
\beta_{auB}(\varepsilon)\equiv\beta_{auB}(\varepsilon,b)\equiv b\alpha^{-1}+K\alpha.
\]
Then, by inequality \ref{eq:temp-454}, we have
\[
E\overline{\lambda}(K^{-1}Z_{1})-E\overline{\lambda}(K^{-1}Z_{0})\leq E|\overline{\lambda}(K^{-1}Z_{1})|+E|\overline{\lambda}(K^{-1}Z_{0})|
\]
\[
\leq3E|K^{-1}Z_{1}|+3E|K^{-1}Z_{0}|\leq3K^{-1}b+3K^{-1}b=6K^{-1}b
\]
\begin{equation}
<\frac{1}{6}\alpha^{3}\exp(-3(E|K^{-1}Z_{0}|\vee E|K^{-1}Z_{1}|)\alpha^{-1}),\label{eq:temp-523-3}
\end{equation}
where the third inequality is a consequence of inequality \ref{eq:temp-495}.
Hence we can apply Theorem \ref{Thm. Bishop's maximal inequality for w.s. submartingales.}
to the wide-sense submartingale $K^{-1}Z|Q_{h}:Q_{h}\times\Omega\rightarrow R$,
to obtain
\begin{equation}
P(\bigvee_{r\in Q(h)}|Z_{r}-Z_{0}|>K\alpha)=P(\bigvee_{r\in Q(h)}|K^{-1}Z_{r}-K^{-1}Z_{0}|>\alpha)<\alpha,\label{eq:temp-509-1}
\end{equation}
for each $h\geq0$. Separately, Chebychev's inequality implies that
\[
P(|Z_{0}|>b\alpha^{-1})\leq b^{-1}\alpha E|Z_{0}|\leq\alpha.
\]
Combining with inequality \ref{eq:temp-509-1}, we obtain
\[
P(\bigvee_{r\in Q(h)}|Z_{r}|>K\alpha+b\alpha^{-1})<2\alpha<\varepsilon.
\]
Consequently, $P(\bigvee_{r\in Q(h)}|Z_{r}|>\gamma)<\varepsilon$
for each $\gamma>K\alpha+b\alpha^{-1}\equiv\beta_{auB}(\varepsilon)$.
In other words, the process $Z$ is a.u. bounded, with the operation
$\beta_{auB}$ as a modulus of a.u. boundlessness.
\end{proof}
\begin{lem}
\label{Lem. Martingale given an event} \textbf{\emph{(Martingale
after an event observed at time $t$).}} Let $Z:Q_{\infty}\times(\Omega,L,E)\rightarrow R$
be a martingale relative to some filtration $\mathcal{L}\equiv\{L^{(t)}:t\in Q_{\infty}\}$.
Let $t\in Q_{\infty}$ and $A\in L^{(t)}$ be arbitrary with $P(A)>0$.
Recall from Definition \ref{Def: Conditional Expectation} the conditional
probability space $(\Omega,L,E_{A})$. Then the process
\[
Z|[t,1]Q_{\infty}:[t,1]Q_{\infty}\times(\Omega,L,E_{A})\rightarrow R
\]
is a martingale relative to the filtration $\mathcal{L}$. 
\end{lem}
\begin{proof}
Consider each $s,r\in[t,1]Q_{\infty}$ with $s\leq r$. Let $U\in L^{(s)}$
be arbitrary, with $U$ bounded. Then $U1_{A}\in L^{(s)}$. Hence,
since $Z:Q_{\infty}\times(\Omega,L,E)\rightarrow R$ is a martingale
relative to the filtration $\mathcal{L}$, we have 
\[
E_{A}(Z_{r}U)\equiv P(A)^{-1}E(Z_{r}U1_{A})=P(A)^{-1}E(Z_{s}U1_{A})\equiv E_{A}(Z_{s}U),
\]
where $Z_{s}\in L^{(s)}$. Hence $E_{A}(Z_{r}|L^{(s)})=Z_{s}$, where
$s,r\in[t,1]Q_{\infty}$ are arbitrary with $s\leq r$. Thus the process
\[
Z|[t,1]Q_{\infty}:[t,1]Q_{\infty}\times(\Omega,L,E_{A})\rightarrow R
\]
is a martingale relative to the filtration $\mathcal{L}$. 
\end{proof}
\begin{thm}
\textbf{\emph{\label{Thm.  Sufficient-Condition-for martignale to be equivalent to a.u. martinale.}
(Sufficient Condition for of martingale on $Q_{\infty}$ to have an
a.u. càdlàg martingale extension to $[0,1]$). }}Let $Z:Q_{\infty}\times\Omega\rightarrow R$
be an arbitrary martingale relative to some filtration $\mathscr{\mathcal{L}}\equiv\{L^{(t)}:t\in Q_{\infty}\}$.
Suppose the following conditions are satisfied.

\emph{(i).} The real number $b>0$ is such that $E|Z_{1}|\leq b$. 

\emph{(ii).} For each $\alpha,\gamma>0$, there exists $\overline{\delta}_{Rcp}(\alpha,\gamma)>0$
such that, for each $h\geq1$ and $t,s\in Q_{h}$ with $t\leq s<t+\overline{\delta}_{Rcp}(\alpha,\gamma)$,
and for each $A\in L^{(t,h)}$ with $A\subset(|Z_{t}|\leq\gamma)$
and $P(A)>0$, we have 
\begin{equation}
E_{A}|Z_{s}|-E_{A}|Z_{t}|\leq\alpha,\label{eq:temp-455}
\end{equation}
and 
\begin{equation}
|E_{A}e^{-|Z(s)|}-E_{A}e^{-|Z(t)|}|\leq\alpha.\label{eq:temp-455-1}
\end{equation}
Then the following holds.

1. The martingale $Z$ is $D$-regular, with a modulus of continuity
in probability $\delta{}_{Cp}\equiv\delta{}_{Cp}(\cdot,b,\overline{\delta}_{Rcp})$
and with a modulus of $D$-regularity $\overline{m}\equiv\overline{m}(b,\overline{\delta}_{Rcp})$. 

2. Let $X\equiv\Phi_{rLim}(Z):[0,1]\times\Omega\rightarrow R$ be
the right limit extension of $Z$. Then $X$ is an a.u. càdlàg martingale
relative to the right-limit extension $\mathcal{L}^{+}$ of the given
filtration $\mathcal{L}$, with $\delta{}_{Cp}(\cdot,b,\overline{\delta}_{Rcp})$
as a modulus of continuity in probability, and with a modulus of a.u.
càdlàg $\delta_{aucl}(\cdot,b,\overline{\delta}_{Rcp})$. 

3. Recall from Definition \ref{Def. Metric on random fields w/ countable parameters}
the metric space $(\widehat{R}(Q_{\infty}\times\Omega,R),\widehat{\rho}_{Prob,Q(\infty)})$
of stochastic processes with parameter set $Q_{\infty},$ where sequential
convergence relative to the metric $\widehat{\rho}_{Prob,Q(\infty)}$
is equivalent to convergence in probability when stochastic processes
are viewed as r.v.'s. Let $\widehat{R}_{Mtgl,b,\overline{\delta}(Rcp)}(Q_{\infty}\times\Omega,R)$
denote the subspace of \emph{$(\widehat{R}(Q_{\infty}\times\Omega,R),\widehat{\rho}_{Prob,Q(\infty)}))$}
consisting of all martingales $Z:Q_{\infty}\times\Omega\rightarrow R$
satisfying the above conditions \emph{(i-ii) }with the same bound
$b>0$ and same given operation $\overline{\delta}_{Rcp}$. Then the
right limit extension function
\begin{equation}
\Phi_{rLim}:(\widehat{R}_{Mtgl,b,\overline{\delta}(Rcp)}(Q_{\infty}\times\Omega,R),\widehat{\rho}_{Prob,Q(\infty)}))\rightarrow(\widehat{D}_{\delta(aucl),\delta(cp)}[0,1],\rho_{\widehat{D}[0,1]})\label{eq:temp-369-2-2-1}
\end{equation}
is uniformly continuous, where 
\[
(\widehat{D}_{\delta(aucl),\delta(cp)}[0,1],\rho_{\widehat{D}[0,1]})
\]
is the metric space of $a.u.c\grave{a}dl\grave{a}g$ processes which
share the common modulus of continuity in probability $\delta{}_{Cp}\equiv\delta{}_{Cp}(\cdot,b,\overline{\delta}_{Rcp})$,
and which share the common modulus of a.u. càdlàg $\delta_{aucl}\equiv\delta_{aucl}(\cdot,\overline{m},\delta_{cp})$
Moreover, the mapping has a modulus of continuity $\delta_{rLim}(\cdot,b,\overline{\delta}_{Rcp})$
depending only on $b$ and $\overline{\delta}_{Rcp}$. 
\end{thm}
\begin{proof}
1. Note that, because $Z$ is a martingale, we have $E|Z_{0}|\leq E|Z_{1}|$
by Assertion 4 of Proposition \ref{Prop. Martingale basics}. Hence
$E|Z_{0}|\vee E|Z_{1}|=E|Z_{1}|\leq b$. Therefore, according to Lemma
\ref{Lem. Each w.s. submartingale on  Q_inf is a.u. bounded}, the
martingale $Z$ is a.u. bounded, with the modulus of a.u. boundlessness
$\beta_{auB}\equiv\beta_{auB}(\cdot,b)$, in the sense of Definition
\ref{Def. a.u. boundedness on Q_inf},.

2. Let $\varepsilon,\gamma>0$ and $h\geq0$ be arbitrary. Consider
each $A\in L^{(t,h)}\subset L^{(t)}$ with $A\subset(|Z_{t}|\leq\gamma)$
and $P(A)>0$. By Lemma \ref{Lem. Martingale given an event}, the
process $Z|[t,1]Q_{\infty}:[t,1]Q_{\infty}\times(\Omega,L_{A},E_{A})\rightarrow R$
is a martingale relative to the filtration $\mathcal{L}$. Define
\[
\alpha\equiv\alpha_{\varepsilon,\gamma}\equiv1\wedge\frac{1}{12}\varepsilon^{3}\exp(-3(\gamma\vee b+1)\varepsilon^{-1})
\]
and
\[
\delta_{SRcp}(\varepsilon,\gamma)\equiv\overline{\delta}_{Rcp}(\alpha,\gamma).
\]

3. Now let $s\in Q_{h}$ be arbitrary with 
\[
t\leq s<t+\delta_{SRcp}(\varepsilon,\gamma)\equiv t+\overline{\delta}_{Rcp}(\alpha,\gamma).
\]
be arbitrary. Trivially, $E_{A}|Z_{t}|\leq\gamma$. Therefore, in
view of inequality \ref{eq:temp-455} in the hypothesis, we have 
\[
E_{A}|Z_{s}|\leq E_{A}|Z_{t}|+\alpha\leq\gamma+\alpha\leq\gamma+1.
\]
Since $Z|[t,1]Q_{h}:[t,1]Q_{h}\times(\Omega,L,E_{A})\rightarrow R$
is a martingale according to Lemma \ref{eq:temp-524-1-1}, we have
$E_{A}Z_{s}=E_{A}Z_{t}$ and $E_{A}|Z_{s}|\geq E_{A}|Z_{t}|$. Hence
equality \ref{eq:temp-524-1-1}, and inequalities \ref{eq:temp-455}
and \ref{eq:temp-455-1}, together imply that
\[
E_{A}\overline{\lambda}(Z_{s})-E_{A}\overline{\lambda}(Z_{t})\leq|E_{A}e^{-|Z(s)|}-E_{A}e^{-|Z(t)|}|+|E_{A}|Z_{s}|-E_{A}|Z_{t}||
\]
\[
=|E_{A}e^{-|Z(s)|}-E_{A}e^{-|Z(t)|}|+E_{A}|Z_{s}|-E_{A}|Z_{t}|
\]
\[
\leq\alpha+\alpha=2\alpha\leq\frac{1}{6}\varepsilon^{3}\exp(-3(b\vee\gamma+1)\varepsilon^{-1})
\]
\[
<\frac{1}{6}\varepsilon^{3}\exp(-3b\varepsilon^{-1})
\]
\[
\leq\frac{1}{6}\varepsilon^{3}\exp(-3(E_{A}|Z_{s}|\vee E_{A}|Z_{t}|)\varepsilon^{-1}).
\]
Therefore we can apply Theorem \ref{Thm. Bishop's maximal inequality for w.s. submartingales.},
to obtain the bound
\begin{equation}
P_{A}(|Z_{t}-Z_{s}|>\varepsilon)<\varepsilon,\label{eq:temp-509-3}
\end{equation}
where $\varepsilon>0$, $\gamma>0$, $h\geq0$, $s\in Q_{h}$ with
$t\leq s<t+\delta_{SRcp}(\varepsilon,\gamma)$, and 
\[
A\in L^{(t,h)}\equiv L(Z_{r}:r\in[0,t]Q_{h})
\]
with $A\subset(|X_{t}|\leq\gamma)$ and $P(A)>0$ are arbitrary. Thus
the process $Z$ is strongly right continuous in probability in the
sense of Definition \ref{Def. strongly  right continuity in probability on Q_inf},
with the operation $\delta_{SRcp}$ as a modulus of strong right  continuity
in probability. Assertion 1 is proved.

4. Thus the process $Z$ satisfies the conditions of Theorem \ref{Thm. a.u.boundeness and Strongly right continuity impliy D-regular}.
Accordingly, the process $Z$ has a modulus of continuity in probability
$\delta{}_{Cp}(\cdot,\beta_{auB},\delta_{SRcp})\equiv\delta{}_{Cp}(\cdot,b,\overline{\delta}_{Rcp})$
and a modulus of $D$-regularity $\overline{m}\equiv\overline{m}(\beta_{auB},\delta_{SRcp})\equiv\overline{m}(b,\overline{\delta}_{Rcp})$.
Moreover, Theorem \ref{Thm. a.u.boundeness and Strongly right continuity impliy D-regular}
says that its right limit extension $X\equiv\Phi_{rLim}(Z)$ is a.u.
càdlàg, with the modulus of a.u. càdlàg $\delta_{aucl}(\cdot,\beta_{auB},\delta_{SRcp})\equiv\delta_{aucl}(\cdot,b,\overline{\delta}_{Rcp})$
as defined in Theorem \ref{Thm. a.u.boundeness and Strongly right continuity impliy D-regular}. 

5. Separately, since $1\in Q_{\infty}$, Assertion 4 of Proposition
\ref{Prop. Martingale basics} implies that the family 
\[
H\equiv\{Z_{u}:u\in Q_{\infty}\}=\{Z_{t}:t\in[0,1]Q_{\infty}\}
\]
is uniformly integrable. Now let $t\in[0,1]$ be arbitrary. Let $r\in[0,1]$
be arbitrary. Suppose $t<r$. Take any sequence $(u_{k})_{k=1,2,\cdots}$
in $[0,r]Q_{\infty}$ such that $u_{k}\rightarrow t$. Since the a.u.
càdlàg process ${\normalcolor X}$ is continuous in probability, we
have $Z_{u(k)}=X_{u(k)}\rightarrow X_{t}$ in probability. Since the
subfamily $\{Z_{u(k)}:k=1,2,\cdots\}$ inherits the uniform integrability
from the family $H$, it follows from Proposition \ref{Prop. Unif integrable+convgnce in pr-> L1 convg}
that the r.r.v. $X_{t}$ is integrable, with $E|X_{u(k)}-X_{t}|\rightarrow0$.
At the same time, $X_{u(k)}\in L^{(r)}$ for each $k\geq1$, it follows
that $X_{t}\in L^{(r)}$, where $r\in[0,1]$ is arbitrary such that
$t<r$. Hence
\[
X_{t}\in\bigcap_{r\in Q(\infty);r>t}L^{(r)}\equiv L^{(t+)}.
\]

6. We will next show that the process $X$ is a martingale relative
to the filtration $\mathcal{L}^{+}$. To that end, let $t,s\in[0,1]$
be arbitrary with $t<s$. Now let $r,u\in Q_{\infty}$ be arbitrary
such that $t<r\leq u$ and $s\leq u$. Let the indicator $Y\in L^{(t+)}$
be arbitrary. Then $Y,X_{t}\in L^{(t+)}\subset L^{(r)}$. Hence, since
$Z$ is a martingale relative to the filtration $\mathcal{L},$ we
have 
\[
EYZ_{r}=EYZ_{u}.
\]
Let $r\downarrow t$ and $u\downarrow s$. Then $E|Z_{r}-X_{t}|=E|X_{r}-X_{t}|\rightarrow0$
and $E|Z_{u}-X_{s}|=E|X_{u}-X_{s}|\rightarrow0$. It then follows
that 
\begin{equation}
EYX_{t}=EYX_{s},\label{eq:temp-607}
\end{equation}
where $t,s\in[0,1]$ are arbitrary with $t<s$. Consider each $t,s\in[0,1]$
with $t\leq s.$ Suppose $EYX_{t}\neq EYX_{s}$. If $t<s$, then equality
\ref{eq:temp-607} would hold, which is a contradiction. Hence $t=s$.
Then trivially $EYX_{t}=EYX_{s}$, again a contradiction. We conclude
that $EYX_{t}=EYX_{s}$ for each $t,s\in[0,1]$ with $t\leq s,$ and
for each indicator $Y\in L^{(t+)}$. Assertion 2 is proved.

6. Assertion 3 remains. Note that, since $Z\in\widehat{R}_{Mtgl,b,\overline{\delta}(Rcp)}(Q_{\infty}\times\Omega,R)$
is arbitrary, we have proved that 
\[
\widehat{R}_{Mtgl,b,\overline{\delta}(Rcp)}(Q_{\infty}\times\Omega,R)\subset\widehat{R}_{Dreg,\overline{m},\delta(cp)}(Q_{\infty}\times\Omega,S).
\]
At the same time, Theorem \ref{Thm.Construction of a.u. cadlag process by right limit of D-regular proceses is continuous.}
says that the right-limit extension function
\begin{equation}
\Phi_{rLim}:(\widehat{R}_{Dreg,\overline{m},\delta(cp)}(Q_{\infty}\times\Omega,S),\widehat{\rho}_{Prob,Q(\infty)}))\rightarrow(\widehat{D}_{\delta(aucl),\delta(cp)}[0,1],\rho_{\widehat{D}[0,1]})\label{eq:temp-369-2-2-2}
\end{equation}
is uniformly continuous, with a modulus of continuity $\delta_{rLim}(\cdot,\overline{m},\delta_{cp})\equiv\delta_{rLim}(\cdot,b,\overline{\delta}_{Rcp})$.
Assertion 3 and the theorem are proved.
\end{proof}
Theorem \ref{Thm.  Sufficient-Condition-for martignale to be equivalent to a.u. martinale.}
can be restated in terms of the continuous construction of a.u. càdlàg
martingales $X:[0,1]\times\Omega\rightarrow R$ from their marginal
distributions. More precisely, let $\xi$ be an arbitrary, but fixed,
binary approximation of $R$ relative to $0\in R$. Recall from Definition
\ref{Def. Metric on  of continuous in prob families of finite joint distributions.}
the metric space $(\widehat{F}_{cp}([0,1],R),\widehat{\rho}_{Cp,\xi,[0,1]|Q(\infty)})$
of consistent families of f.j.d.'s on $[0,1]$ which are continuous
in probability. Let $\widehat{F}_{Mtgl,b,\overline{\delta}(Rcp)}([0,1],R)$
be the subspace of consistent families $F$ such that $F|Q_{\infty}$
gives the marginal distributions of some martingale $Z\in\widehat{R}_{Mtgl,b,\overline{\delta}(Rcp)}(Q_{\infty}\times\Omega,R)$.
Let 
\[
(\Theta_{0},L_{0},I_{0})\equiv([0,1],L_{0},\int\cdot dx)
\]
denote the Lebesgue integration space based on the interval $[0,1]$.
Recall from Theorem \ref{Thm. DKS Extension, construction and continuity}
that the Daniell-Kolmogorov-Skorokhod Extension 
\[
\Phi_{DKS,\xi}:(\widehat{F}(Q_{\infty},R),\widehat{\rho}_{Marg,\xi,Q(\infty)})\rightarrow(\widehat{R}(Q_{\infty}\times\Theta_{0},R),\rho_{Q(\infty)\times\Theta(0),R})
\]
is uniformly continuous, with a modulus of continuity $\delta_{DKS}(\cdot,\left\Vert \xi\right\Vert )$
dependent only on $\left\Vert \xi\right\Vert $. Hence the composite
mapping 
\[
\Phi_{rLim}\circ\Phi_{DKS,\xi}:(\widehat{F}_{Mtgl,b,\overline{\delta}(Rcp)}([0,1],R)|Q_{\infty},\widehat{\rho}_{Marg,\xi,Q(\infty)})\rightarrow(\widehat{D}_{\delta(aucl),\delta(cp)}[0,1],\rho_{\widehat{D}[0,1]})
\]
is uniformly continuous.

\section{A Sufficient Condition for a Right Hoelder Process}

In this section, we give a sufficient condition, in terms of triple
distributions, for a set $\widehat{F}$ of consistent families of
f.j.d.'s to be $D$-regular. Theorem \ref{Thm. Extension of D-regular process by right limit is a.u.cadlag}
is then applied to construct corresponding a.u. càdlàg processes. 

As an application, we will prove that, under an additional condition
on the modulus of continuity in probability, the a.u. càdlàg process
so constructed has sample functions which are right Hoelder, in a
sense made precise below. This result seems to be hitherto unknown.

In the following, recall Definition \ref{Def. Notations for dyadic rationals}
for the notations associated with the enumerated sets $(Q_{k})_{k=0,1,\cdots}$and
$Q_{\infty}$ of dyadic rationals in $[0,1]$. In particular, $p_{k}\equiv2^{k}$
and $\Delta_{k}\equiv2^{-k}$ for each $k\geq0$. Separately, $(\Omega,L,E)\equiv(\Theta_{0},L_{0},E_{0})\equiv([0,1],L_{0},\int\cdot dx)$
denote the Lebesgue integration space based on the interval $\Theta_{0}\equiv[0,1]$.
This will serve as the sample space. Let $\xi\equiv(A_{q})_{q=1,2,\cdots}$
be an arbitrarily, but fixed, binary approximation of the locally
compact metric space $(S,d)$ relative to some fixed reference point
$x_{\circ}\in S$. Recall the operation $[\cdot]_{1}$ which assigns
to each $a\in R$ an integer $[a]_{1}\in(a,a+2)$. We will repeated
use the inequality $a<[a]_{1}<a+2$ for each $a\in R$, without further
comments.

The following theorem is in essence due to Kolmogorov. 
\begin{thm}
\label{Thm. Sufficient Conditions for D-regularity in terms of triple joiny distributions}
\textbf{\emph{(Sufficient Condition for D-regularity in terms of triple
joint distributions). }}Let $F$ be an arbitrary consistent family
of f.j.d.'s  with parameter set $[0,1]$ and state space $(S,d)$.
Suppose $F$ is continuous in probability. Let $Z:Q\times\Omega\rightarrow S$
be an arbitrary process with marginal distributions given by $F|Q$.
Suppose there exist two sequences $(\gamma_{k})_{k=0,1\cdots}$ and
$(\alpha_{k})_{k=0,1,\cdots}$ of positive real numbers such that
(i) $\sum_{k=0}^{\infty}2^{k}\alpha_{k}<\infty$ and $\sum_{k=0}^{\infty}\gamma_{k}<\infty$,
(ii) the set 
\begin{equation}
A_{r,s}'^{(k)}\equiv(d(Z_{r},Z_{s})>\gamma_{k+1})\label{eq:temp-246-1-1}
\end{equation}
is \textup{\emph{measurable}}\emph{ } for each $s,t\in Q_{\infty}$,
for each $k\geq0$, and (iii) 
\[
P(A_{v,v''}'^{(k)}A_{v'',v'}'^{(k)})<\alpha_{k},
\]
where $v'\equiv v+\Delta_{k}$ and $v''\equiv v+\Delta_{k+1}=v'-\Delta_{k+1}$,
for each $v\in[0,1)Q_{k}$, for each $k\geq0$. 

Then the family $F|Q_{\infty}$ of f.j.d.'s is $D$-regular. Specifically,
let $m_{0}\equiv0$. For each $n\geq1$, let $m_{n}\geq m_{n-1}+1$
be so large that $\sum_{k=m(n)}^{\infty}2^{k}\alpha_{k}<2^{-n}$ and
$\sum_{k=m(n)+1}^{\infty}\gamma_{k}<2^{-n-1}$. Then the sequence
$(m_{n})_{n=0,1,\cdots}$ is a modulus of $D$-regularity of the family
$F|Q_{\infty}$.
\[
.
\]
\end{thm}
\begin{proof}
1. Let $k\geq0$ be arbitrary. Define 
\begin{equation}
D'_{k}\equiv\bigcup_{v\in[0,1)Q(k)}A_{v,v''}'^{(k)}A_{v'',v'}'^{(k)}.\label{eq:temp-140}
\end{equation}
Then $P(D'_{k})\leq2^{k}\alpha_{k}$ according to Condition (iii)
in the hypothesis.

2. Inductively, for each $n\geq0$, take any 
\begin{equation}
\beta_{n}\in(2^{-n},2^{-n+1})\label{eq:temp-420-1-1}
\end{equation}
such that, for each $s,t\in Q_{\infty},$ and for each $k=0,\cdots,n$,
the sets
\begin{equation}
(d(Z_{t},Z_{s})>\beta_{k}+\cdots+\beta_{n})\label{eq:temp-304-1-1-1}
\end{equation}
and
\begin{equation}
A_{t,s}^{(n)}\equiv(d(Z_{t},Z_{s})>\beta_{n+1})\label{eq:temp-246-1}
\end{equation}
are $\mathrm{measurable}$. Note that $\beta_{n,\infty}\leq\sum_{k=n}^{\infty}2^{-k+1}=2^{-n+2}$
for each $n\geq0$. 

3. Let $n\geq0$ be arbitrary, but fixed until further notice. For
ease of notations, suppress some symbols signifying dependence on
$n$, to write $q_{i}\equiv q_{m(n),i}\equiv2^{-p(m(n))}$ for each
$i=0,\cdots,p_{m(n)}$. Let $\beta>2^{-n}>\beta_{n+1}$ be arbitrary
such that the set 
\begin{equation}
A_{t,s}^{\beta}\equiv(d(Z_{t},Z_{s})>\beta)\label{eq:temp-372-1}
\end{equation}
is $\mathrm{measurable}$ for each $s,t\in Q_{\infty}$. Define the
exceptional set 
\[
D_{n}\equiv\bigcup_{t\in[0,1)Q(m(n))}\bigcup_{r,s\in(t,t')Q(m(n+1));r\leq s}(A_{t,r}^{\beta}\cup A_{t,s}^{\beta})(A_{t,r}^{\beta}\cup A_{t',s}^{\beta})(A_{t',r}^{\beta}\cup A_{t',s}^{\beta})
\]
\begin{equation}
\subset\bigcup_{t\in[0,1)Q(m(n))}\bigcup_{r,s\in(t,t')Q(m(n+1));r\leq s}(A_{t,r}^{(n)}\cup A_{t,s}^{(n)})(A_{t,r}^{(n)}\cup A_{t',s}^{(n)})(A_{t',r}^{(n)}\cup A_{t',s}^{(n)}),\label{eq:temp-201-4-3-3-1-1}
\end{equation}
where, for each $t\in[0,1)Q{}_{m(n)}$, we write $t'\equiv t+\Delta_{m(n)}\in Q_{m(n)}.$
To verify Condition 1 in Definition \ref{Def.  D-regular process and D-regular family of f.j.d.s on Q_inf}
for the sequence $(m)_{n=0,1,\cdots}$ and the process $Z$, we need
only show that 
\begin{equation}
P(D_{n})\leq2^{-n}.\label{eq:temp-204-3-2-1-1}
\end{equation}

4. To that end, consider each $\omega\in(\bigcup_{k=m(n)}^{m(n+1)}D'_{k})^{c}$.
Let $t\in[0,1)Q_{m(n)}$ be arbitrary, and write $t'\equiv t+\Delta_{m(n)}$.
We will show inductively that, for each $k=m_{n},\cdots,m_{n+1}$,
there exists $r_{k}\in(t,t']Q_{k}$ such that 
\begin{equation}
\bigvee_{u\in[t,r(k)-\Delta(k)]Q(k)}d(Z_{t}(\omega),Z_{u}(\omega))\leq\sum_{j=m(n)+1}^{k}\gamma_{j},\label{eq:temp-129}
\end{equation}
and 
\begin{equation}
\bigvee_{v\in[r(k),t']Q(k)}d(Z_{v}(\omega),Z_{t'}(\omega))\leq\sum_{j=m(n)+1}^{k}\gamma_{j},\label{eq:temp-139}
\end{equation}
where an empty sum is, by convention, equal to $0$. Start with $k=m_{n}$.
Define $r_{k}\equiv t'$, whence $r_{k}-\Delta_{k}=t$. Then inequalities
\ref{eq:temp-129} and \ref{eq:temp-139} are trivially satisfied. 

Suppose, for some $k=m_{n},\cdots,m_{n+1}-1$, we have constructed
$r_{k}\in(t,t']Q_{k}$ which satisfies inequalities \ref{eq:temp-129}
and \ref{eq:temp-139}. Note that, according to the defining equality
\ref{eq:temp-140}, we have 
\[
\omega\in D'_{k}{}^{c}\subset(A_{r(k)-\Delta(k),r(k)-\Delta(k+1)}'^{(k)})^{c}\cup(A_{r(k)-\Delta(k+1),r(k)}'^{(k)})^{c}.
\]
Hence, by the defining equality \ref{eq:temp-246-1-1}, we have (i')
\begin{equation}
d(Z_{r(k)-\Delta(k)}(\omega),Z_{r(k)-\Delta(k+1)}(\omega))\leq\gamma_{k+1},\label{eq:temp-353}
\end{equation}
or (ii')
\begin{equation}
d(Z_{r(k)-\Delta(k+1)},Z_{r(k)}(\omega))\leq\gamma_{k+1}.\label{eq:temp-354}
\end{equation}
In Case (i'), define $r_{k+1}\equiv r_{k}$. In Case (ii'), define
$r_{k+1}\equiv r_{k}-\Delta_{k+1}$. 

We wish to prove inequalities \ref{eq:temp-129} and \ref{eq:temp-139}
for $k+1$. For that purpose, first consider each 
\[
u\in[t,r_{k+1}-\Delta_{k+1}]Q_{k+1}.
\]
\begin{equation}
=[t,r_{k}-\Delta_{k}]Q_{k}\cup[t,r_{k}-\Delta_{k}]Q_{k+1}Q_{k}^{c}\cup(r_{k}-\Delta_{k},r_{k+1}-\Delta_{k+1}]Q_{k+1}.\label{eq:temp-198}
\end{equation}
Suppose $u\in[t,r_{k}-\Delta_{k}]Q_{k}$. Then inequality \ref{eq:temp-129}
in the induction hypothesis trivially implies 
\begin{equation}
d(Z_{t}(\omega),Z_{u}(\omega))\leq\sum_{j=m(n)+1}^{k+1}\gamma_{j}.\label{eq:temp-129-1}
\end{equation}
Suppose next that $u\in[t,r_{k}-\Delta_{k}]Q_{k+1}Q_{k}^{c}$. Then
$u\leq r_{k}-\Delta_{k}-\Delta_{k+1}$. Let $v\equiv u-\Delta_{k+1}$
and $v'\equiv v+\Delta_{k}=u+\Delta_{k+1}$. Then $v\in[0,1)Q_{k}$,
and so the defining inequality \ref{eq:temp-140} implies that
\[
\omega\in D'_{k}{}^{c}\subset(A_{v,u}'^{(k)})^{c}\cup(A_{u,v'}'^{(k)})^{c},
\]
and therefore that
\[
d(Z_{u}(\omega),Z_{v}(\omega))\wedge d(Z_{u}(\omega),Z_{v'}(\omega))\leq\gamma_{k+1}.
\]
It follows that
\[
d(Z_{u}(\omega),Z_{t}(\omega))
\]
\[
\leq(d(Z_{u}(\omega),Z_{v}(\omega))+d(Z_{v}(\omega),Z_{t}(\omega)))\wedge(d(Z_{u}(\omega),Z_{v'}(\omega))+d(Z_{v'}(\omega),Z_{t}(\omega)))
\]
\[
\leq\gamma_{k+1}+d(Z_{v}(\omega),Z_{t}(\omega))\vee d(Z_{v'}(\omega),Z_{t}(\omega))
\]
\[
\leq\gamma_{k+1}+\sum_{j=m(n)+1}^{k}\gamma_{j}=\sum_{j=m(n)+1}^{k+1}\gamma_{j},
\]
where the last inequality is due to inequality \ref{eq:temp-129}
in the induction hypothesis. Thus we have verified inequality \ref{eq:temp-129-1}
also for each $u\in[t,r_{k}-\Delta_{k}]Q_{k+1}Q_{k}^{c}$. Now suppose
$u\in(r_{k}-\Delta_{k},r_{k+1}-\Delta_{k+1}]Q_{k+1}$. Then
\[
r_{k+1}>r_{k}-\Delta_{k}+\Delta_{k+1}=r_{k}-\Delta_{k+1},
\]
which, by the definition of $r_{k+1}$, rules out Case (ii'). Hence
Case (i') must hold, where $r_{k+1}\equiv r_{k}$. Consequently, 
\[
u\in(r_{k}-\Delta_{k},r_{k}-\Delta_{k+1}]Q_{k+1}=\{r_{k}-\Delta_{k+1}\}
\]
and so $u=r_{k}-\Delta_{k+1}$. Let $v\equiv r_{k}-\Delta_{k}$. Then
inequality \ref{eq:temp-353} implies that
\begin{equation}
d(Z_{v}(\omega),Z_{u}(\omega))\equiv d(Z_{r(k)-\Delta(k)}(\omega),Z_{r(k)-\Delta(k+1)}(\omega))\leq\gamma_{k+1}.\label{eq:temp-353-1}
\end{equation}
Hence
\[
d(Z_{u}(\omega),Z_{t}(\omega))\leq d(Z_{u}(\omega),Z_{v}(\omega))+d(Z_{v}(\omega),Z_{t}(\omega))
\]
\[
\leq\gamma_{k+1}+d(Z_{v}(\omega),Z_{t}(\omega))\leq\gamma_{k+1}+\sum_{j=m(n)+1}^{k}\gamma_{j}=\sum_{j=m(n)+1}^{k+1}\gamma_{j},
\]
where the last inequality is due to inequality \ref{eq:temp-129}
in the induction hypothesis. Combining, we see that inequality \ref{eq:temp-129-1}
holds for each $u\in[t,r_{k+1}-\Delta_{k+1}]Q_{k+1}$. 

Similarly we can verify that
\begin{equation}
d(Z_{u}(\omega),Z_{t'}(\omega))\leq\sum_{j=m(n)+1}^{k+1}\gamma_{j}.\label{eq:temp-129-1-1}
\end{equation}
for each $u\in[r_{k+1},t']Q_{k+1}$. Summing up, inequalities \ref{eq:temp-129}
and \ref{eq:temp-139} have been verified for $k+1$. Induction is
completed. Thus inequalities \ref{eq:temp-129} and \ref{eq:temp-139}
hold for each $k=m_{n},\cdots,m_{n+1}$.

5. Continuing, let $r,s\in(t,t')Q_{m(n+1)}$ be arbitrary with $r\leq s$.
Write $k\equiv m_{n+1}$. Then there are three possibilities: (i'')
$r,s\in[t,r_{k}-\Delta_{k}]Q_{k}$, (ii'') $r\in[t,r_{k}-\Delta_{k}]Q_{k}$
and $s\in[r_{k},t']Q_{k}$, or (iii'') $r,s\in[r_{k},t']Q_{k}$. In
Case (i''), inequality \ref{eq:temp-129} applies to $r,s$, yielding
\[
d(Z_{t}(\omega),Z_{r}(\omega))\vee d(Z_{t}(\omega),Z_{s}(\omega))\leq\sum_{j=m(n)+1}^{k}\gamma_{j}<2^{-n-1}<\beta_{n+1}.
\]
Hence $\omega\in(A_{t,r}^{(n)}\cup A_{t,s}^{(n)})^{c}\subset D_{n}^{c}.$
In Case (ii''), inequalities \ref{eq:temp-129} and \ref{eq:temp-139},
apply to $r$ and $s$ receptively, and yield
\[
d(Z_{t}(\omega),Z_{r}(\omega))\vee d(Z_{s}(\omega),Z_{t'}(\omega))\leq\sum_{j=m(n)+1}^{k}\gamma_{j}<2^{-n-1}<\beta_{n+1}.
\]
In other words, $\omega\in(A_{t,r}^{(n)}\cup A_{s,t'}^{(n)})^{c}\subset D_{n}^{c}.$
Similarly, in Case (iii'), we can prove that $\omega\in(A_{r,t'}^{(n)}\cup A_{s,t'}^{(n)})^{c}\subset D_{n}^{c}.$ 

6. Summing up, we have shown that $\omega\in D_{n}^{c}$ where $\omega\in(\bigcup_{k=m(n)}^{m(n+1)}D'_{k})^{c}$
is arbitrary. Thus $D_{n}\subset\bigcup_{k=m(n)}^{m(n+1)}D'_{k}$.
Hence
\[
P(D_{n})\leq\sum_{k=m(n)}^{m(n+1)}P(D'_{k})\leq\sum_{k=m(n)}^{\infty}2^{k}\alpha_{k}<2^{-n},
\]
where $n\geq0$ is arbitrary. This proves inequality \ref{eq:temp-204-3-2-1-1},
and verifies Condition 1 in Definition \ref{Def.  D-regular process and D-regular family of f.j.d.s on Q_inf}
for the sequence $(m)_{n=0,1,\cdots}$ and the process $Z$. At the
same time, since the family $F$ is continuous in probability by hypothesis,
Condition 2 of . Definition \ref{Def.  D-regular process and D-regular family of f.j.d.s on Q_inf}
follows for $F|Q_{\infty}$ and for $Z$. We conclude that the family
$F|Q_{\infty}$ of f.j.d.'s is $D$-regular, with the sequence $(m)_{n=0,1,\cdots}$
as modulus of $D$-regularity.
\end{proof}
\begin{defn}
\label{Def. right Hoelder cadlag functions and processes}\textbf{
(Right Hoelder process)} Let $C_{0}\geq0$ and $\lambda>0$ be arbitrary.
Let $X:[0,1]\times\Omega\rightarrow S$ be an arbitrary a.u. càdlàg
process. Suppose, for each $\varepsilon>0$, there exist (i) $\widetilde{\delta}>0$,
(ii) a $\mathrm{measurable}$ subset $B\subset\Omega$ with $P(B^{c})<\varepsilon$,
(iii) for each $\omega\in B$, there exists a Lebesgue $\mathrm{measurable}$
subset $\widetilde{\theta}(\omega)$ of $[0,1]$ with Lebesgue measure
$\mu\widetilde{\theta}_{k}(\omega)^{c}<\varepsilon$, such that, for
each $t\in\widetilde{\theta}(\omega)\cap domain(X(\cdot,\omega))$
and for each $s\in[t,t+\widetilde{\delta})\cap domain(x),$ we have
\[
d(X(t,\omega),X(s,\omega))\leq C_{0}(s-t)^{\lambda}.
\]
Then the a.u. càdlàg process $X$ is said to be \index{right Hoelder process}\emph{right
Hoelder}, with \index{right Hoelder exponent}\emph{right Hoelder
exponent} $\lambda$, and with \emph{right Hoelder constant}\index{right Hoelder constant}
$C_{0}$.$\square$

Part 2 of the next theorem is in essence due to Chentsov\cite{Chentsov56}.
Part 3, concerning right Hoelder processes, seems hitherto unknown. 
\end{defn}
\begin{thm}
\label{Thm Sufficient condition  right Hoelder process} \textbf{\emph{(Sufficient
condition for a right Hoelder process).}} Let $u\geq0$ and $w,K>0$
be arbitrary. Let $F$ be an arbitrary consistent family of f.j.d.'s
 with parameter set $[0,1]$ and state space $(S,d)$. Suppose $F$
is continuous in probability, with a modulus of continuity in probability
$\delta_{cp}$, and suppose 
\begin{equation}
F_{t,r,s}\{(x,y,z)\in S^{3}:d(x,y)\wedge d(y,z)>b\}\leq b^{-u}(Ks-Kt)^{1+w}\label{eq:temp-472-3-1}
\end{equation}
for each $b>0$ and for each $t\leq r\leq s$ in $[0,1]$. Then the
following holds.

1. The family $F|Q_{\infty}$ is $D$-regular. 

2. There exists an a.u. càdlàg process $X:[0,1]\times\Omega\rightarrow S$
with marginal distributions given by the family $F$, and with a modulus
of a.u. càdlàg dependent only on $u,w$ and the modulus of continuity
in probability $\delta_{cp}$.

3. Suppose, in addition, that there exist $\overline{u},\overline{\lambda}>0$
such that 
\begin{equation}
F_{t,s}(d>b)\leq b^{-\overline{u}}|Ks-Kt|^{\overline{\lambda}}\label{eq:temp-150}
\end{equation}
for each $b>0$ and for each $t,s\in[0,1]$. Then there exist constants
$\lambda(u,w,\overline{u},\overline{\lambda})>0$ and $C(K,u,w,\overline{u},\overline{\lambda})>0$
such that each a.u. càdlàg process $X:[0,1]\times\Omega\rightarrow S$
with marginal distribution given by the family $F$ is right Hoelder,
with right Hoelder exponent $\lambda$, and with right Hoelder constant
$C$. Specifically
\begin{equation}
\lambda(u,w,\overline{u}\overline{\lambda})\equiv2^{-1}w((1+2u)(1+\overline{\lambda}^{-1}(1+\overline{u}))+w(2+3\overline{\lambda}^{-1}(1+\overline{u})))^{-1}.\label{eq:temp-296}
\end{equation}
\end{thm}
\begin{proof}
1. Let $Z:Q_{\infty}\times\Omega\rightarrow S$ be an arbitrary process
with marginal distributions given by $F|Q_{\infty}$. Define $u_{0}\equiv u+2^{-1}$
and $u_{1}\equiv u+1$. Then $\gamma_{0}\equiv2^{-w/u(0)}<\gamma_{1}\equiv2^{-w/u(1)}$.
Take an arbitrary $\gamma\in(\gamma_{0},\gamma_{1})$ such that the
subset
\[
A_{r,s}'^{(k)}\equiv(d(Z_{r},Z_{s})>\gamma^{k+1})=(d(Z_{r},Z_{s})\vee\gamma_{0}^{k+1}>\gamma^{k+1})
\]
\begin{equation}
=((d(Z_{r},Z_{s})\vee\gamma_{0}^{k+1})^{1/(k+1)}>\gamma)\label{eq:temp-246-1-1-1}
\end{equation}
of $\Omega$ is $\mathrm{measurable}$ for each $r,s\in Q_{\infty}$
and for each $k\geq0$. Then, since $2^{-w/v}$ is a strictly increasing
continuous function of $v\in(u_{0},u_{1})$ with range $(\gamma_{0},\gamma_{1})$,
there exists a unique 
\[
v\in(u_{0},u_{1})\equiv(u+2^{-1},u+1)
\]
such that 
\begin{equation}
\gamma=2^{-w/v}.\label{eq:temp-352}
\end{equation}
Note that $0<\gamma<1$, and that
\begin{equation}
0<\gamma^{-u}2^{-w}=2^{uw/v-w}<2^{w-w}=1.\label{eq:temp-351}
\end{equation}

2. Let $k\geq0$ be arbitrary. Define the positive real numbers
\[
\gamma_{k}\equiv\gamma^{k},
\]
\[
\alpha_{k}\equiv K^{1+w}2^{-k}\gamma^{-(k+1)u}2^{-kw}.
\]
Then
\[
\sum_{k=0}^{\infty}\gamma_{k}=\sum_{k=0}^{\infty}\gamma^{k}<\infty,
\]
and, in view of inequality \ref{eq:temp-351}, 
\[
\sum_{k=0}^{\infty}2^{k}\alpha_{k}=K^{1+w}\gamma^{-u}\sum_{k=0}^{\infty}(\gamma^{-u}2^{-w})^{k}<\infty.
\]
Let $t\in[0,1)Q_{k}$ be arbitrary. Write $t'\equiv t+\Delta_{k}$
and $t''\equiv t+\Delta_{k+1}=t'-\Delta_{k+1}$. Then
\[
P(A_{t,t''}'^{(k)}A_{t'',t'}'^{(k)})=P(d(Z_{t},Z_{t''})\wedge d(Z_{t''},Z_{t'})>\gamma^{k+1})
\]
\[
\leq\gamma^{-(k+1)u}(Kt'-Kt){}^{1+w}=\gamma^{-(k+1)u}(K\Delta_{k})^{1+w}=\gamma^{-(k+1)u}(K2^{-k})^{1+w}
\]
\begin{equation}
=K^{1+w}2^{-k}\gamma^{-(k+1)u}2^{-kw}\equiv\alpha_{k},\label{eq:temp-474-1}
\end{equation}
where the inequality follows from inequality \ref{eq:temp-472-3-1}
in the hypothesis. Thus we have verified the conditions in Theorem
\ref{Thm. Sufficient Conditions for D-regularity in terms of triple joiny distributions}
for the consistent family $F|Q_{\infty}$ of f.j.d.'s to be $D$-regular,
with with a modulus of $D$-regularity $(m_{n})_{n=0,1,\cdots}$ dependent
only on the sequences $(\alpha_{k})_{k=0,1,\cdots}$. and $(\gamma_{k})_{k=0,1,\cdots}$,
which, in turn, depends only on the constants $u,w$.

3. Theorem \ref{Cor.  Construction of a.u. calafg processes from D-regular f.j.d.'s, and continuity}
can therefore be applied to construct an a.u. càdlàg process $X:[0,1]\times\Omega\rightarrow S$
with marginal distributions given by the family $F$, and with modulus
of a.u. càdlàg depending only on the modulus of $D$-regularity $(m_{n})_{n=0,1,\cdots}$
and the modulus of continuity in probability $\delta_{cp}$. Thus
the process $X$ has a modulus of a.u. càdlàg depending only on the
constants $u,w$ and the modulus of continuity in probability $\delta_{cp}$.
Assertions 1 and 2 have been proved.

4. We proceed to prove Assertion 3. Suppose the positive constants
$\overline{u}$ and $\overline{\lambda}$ are given and satisfy inequality
\ref{eq:temp-150}, and let $\lambda\equiv\lambda(u,w,\overline{u}\overline{\lambda})$
be as defined in equality \ref{eq:temp-296}. Let $X:[0,1]\times\Omega\rightarrow S$
be an arbitrary a.u. càdlàg process with marginal distribution given
by the family $F$. We need to show that $\lambda$ is a right Hoelder
exponent of $X$. For abbreviation, define the constants
\[
c_{0}\equiv(1+w)\log_{2}K-\log_{2}(1-\gamma^{-u}2^{-w}),
\]
\[
c\equiv w(1-uv^{-1})=w-uwv^{-1},
\]
\[
c_{1}\equiv-\log_{2}(1-\gamma),
\]
\[
c_{2}\equiv-\log_{2}\gamma=wv^{-1},
\]
\[
\kappa_{0}\equiv\kappa_{0}(K,u,w)\equiv[(c^{-1}c_{0}-1)\vee c_{2}^{-1}(c_{1}+1)]_{1}
\]
and
\[
\kappa\equiv\kappa(u,w)\equiv[c^{-1}]_{1}\geq1.
\]

5. Note that, because $v-u<1$, we have $1-uv^{-1}<v^{-1}$, whence
$c<c_{2}$ and 
\[
\kappa c_{2}>\kappa c>1.
\]
Note also that, because $u+2^{-1}<v<u+1$, we have 
\[
(1+2u)^{-1}=1-u(2^{-1}+u)^{-1}<1-uv^{-1}<1-u(1+u)^{-1}=(1+u)^{-1},
\]
whence
\[
w(1+2u)^{-1}<c\equiv w(1-uv^{-1})<w(1+u)^{-1},
\]
or 
\[
w^{-1}(1+2u)>c^{-1}>w^{-1}(1+u).
\]
Hence, the defining equality \ref{eq:temp-296} yields 
\[
\lambda=2^{-1}(w^{-1}(1+2u)(1+\overline{\lambda}^{-1}(1+\overline{u}))+(2+3\overline{\lambda}^{-1}(1+\overline{u})))^{-1}
\]
\[
<2^{-1}(c^{-1}(1+\overline{\lambda}^{-1}(1+\overline{u}))+(2+3\overline{\lambda}^{-1}(1+\overline{u})))^{-1}
\]
\[
=2^{-1}(c^{-1}+2+\overline{\lambda}^{-1}(1+\overline{u})(c^{-1}+3))^{-1}
\]
\begin{equation}
<2^{-1}(\kappa+\overline{\lambda}^{-1}(1+\overline{u})(\kappa+1))^{-1}.\label{eq:temp-371}
\end{equation}

6. Consider the process $Z\equiv X|Q_{\infty}$. Then $X$ is equal
to the right-limit extension of the process $Z$. Define $m_{0}\equiv0$.
Let $n\geq1$ be arbitrary. Define 
\begin{equation}
m_{n}\equiv\kappa_{0}+\kappa n.\label{eq:temp-292}
\end{equation}
Then $m_{n}\geq m_{n-1}+1$. Moreover, 
\[
\log_{2}\sum_{h=m(n)}^{\infty}2^{h}\alpha_{h}=\log_{2}\sum_{h=m(n)}^{\infty}2^{h}K^{1+w}2^{-h}\gamma^{-(h+1)u}2^{-hw}
\]
\[
=\log_{2}\gamma^{-(m(n)+1)u}2^{-(m(n)+1)w}K^{1+w}(1-\gamma^{-u}2^{-w})^{-1}
\]
\[
=\log_{2}2^{-(m(n)+1)(w-uw/v)}K^{1+w}(1-\gamma^{-u}2^{-w})^{-1}
\]
\[
=-(m_{n}+1)(w-uwv^{-1})+(1+w)\log_{2}K-\log_{2}(1-\gamma^{-u}2^{-w})
\]
\[
\equiv-(\kappa_{0}+\kappa n+1)c+c_{0}<-((c^{-1}c_{0}-1)+1)c-\kappa cn+c_{0}
\]
\[
=-\kappa cn<-n.
\]
where the last inequality is thanks to the earlier observation that
$\kappa c_{2}>\kappa c>1$. Hence
\begin{equation}
\sum_{k=m(n)}^{\infty}2^{k}\alpha_{k}<2^{-n}.\label{eq:temp-460}
\end{equation}
Similarly, 
\[
\log_{2}\sum_{k=m(n)+1}^{\infty}\gamma_{k}<\log_{2}\sum_{k=m(n)}^{\infty}\gamma^{k}=\log_{2}\gamma^{m(n)}(1-\gamma)^{-1}
\]
\[
=m_{n}\log_{2}\gamma-\log_{2}(1-\gamma)=-(\kappa_{0}+\kappa n)c_{2}+c_{1}
\]
\[
<-\kappa_{0}c_{2}-\kappa c_{2}n+c_{1}<-c_{2}^{-1}(c_{1}+1)c_{2}-\kappa c_{2}n+c_{1}=-1-\kappa c_{2}n<-1-n,
\]
whence
\begin{equation}
\sum_{k=m(n)+1}^{\infty}\gamma_{k}<2^{-n-1}.\label{eq:temp-465}
\end{equation}
In view of inequalities \ref{eq:temp-474-1}, \ref{eq:temp-460},
and \ref{eq:temp-465}, Theorem \ref{Thm. Sufficient Conditions for D-regularity in terms of triple joiny distributions}
applies, and says that the sequence $\overline{m}\equiv(m_{n})_{n=0,1,\cdots}$
is a modulus of $D$-regularity of the family $F|Q_{\infty}$ and
of the process $Z$. 

7. More specifically, let $\varepsilon>0$ be arbitrary. Define 
\begin{equation}
\delta_{cp}(\varepsilon)\equiv K^{-1}(2^{-1}\varepsilon){}^{(1+\overline{u})/\overline{\lambda}}.\label{eq:temp-471}
\end{equation}
Then, for each $t,s\in[0,1]$ with $|s-t|\leq\delta_{cp}(\varepsilon)$,
and for each $b>2^{-1}\varepsilon$, we have
\[
F_{t,s}(\widehat{d})\leq b+F_{t,s}(\widehat{d}>b)\leq b+F_{t,s}(d>b)\leq b+b^{-\overline{u}}|Ks-Kt|^{\overline{\lambda}}
\]
\[
\leq b+b{}^{-\overline{u}}K^{\overline{\lambda}}\delta_{cp}(\varepsilon)^{\overline{\lambda}}=b+b{}^{-\overline{u}}K^{\overline{\lambda}}K^{-\overline{\lambda}}(2^{-1}\varepsilon){}^{(1+\overline{u})}=b+b^{-\overline{u}}(2^{-1}\varepsilon){}^{(1+\overline{u})},
\]
where the first inequality is because $\widehat{d}\equiv d\wedge1\leq1$,
and where the third inequality is from inequality \ref{eq:temp-150}
in the hypothesis. Letting $b\downarrow2^{-1}\varepsilon$, we obtain
\begin{equation}
F_{t,s}(\widehat{d})\leq2^{-1}\varepsilon+(2^{-1}\varepsilon)^{-\overline{u}}(2^{-1}\varepsilon){}^{(1+\overline{u})}=2^{-1}\varepsilon+2^{-1}\varepsilon=\varepsilon.\label{eq:temp-150-1}
\end{equation}
Thus the operation $\delta_{cp}$ is a modulus of continuity in probability
of the family $F$ of f.j.d.'s, and of the $D$-regular process $Z$.
Hence Theorem \ref{Thm. Extension of D-regular process by right limit is a.u.cadlag}
says that the right-limit extension $X$ of $Z$ has a modulus of
a.u. càdlàg $\delta_{aucl}$ defined as follows. Let $\varepsilon>0$
be arbitrary. Let $n\geq0$ be so large that $2^{-n+6}<\varepsilon.$
Let $J>n$ be so large that 
\begin{equation}
\Delta_{m(J)}\equiv2^{-m(J)}<2^{-2}\delta_{cp}(2^{-2m(n)-2n-10}).\label{eq:temp-332-1-2-4}
\end{equation}
Then
\begin{equation}
\delta_{aucl}(\varepsilon,\overline{m},\delta_{cp})\equiv\Delta_{m(J)}.\label{eq:temp-473}
\end{equation}

8. Let $h\geq0$ be arbitrary. Write $\varepsilon_{h}\equiv2^{-h}$
for short. By Definition \ref{Def. a.u. cadlag process} for a modulus
of a.u. càdlàg, there exist a $\mathrm{measurable}$ set $A_{h}$
with 
\[
P(A_{h}^{c})<\varepsilon_{h}\equiv2^{-h},
\]
an integer $p_{h}\geq1$, and a sequence of r.r.v.'s 
\begin{equation}
0=\tau_{h,0}<\tau_{h,1}<\cdots<\tau_{h,p(h)-1}<\tau_{h,p(h)}=1,\label{eq:temp-404-1}
\end{equation}
such that, for each $i=0,\cdots,p_{h}-1$, the function $X_{\tau(h,i)}$
is a r.v., and such that, for each $\omega\in A_{h}$, we have
\begin{equation}
\bigwedge_{i=0}^{p(h)-1}(\tau_{h,i+1}(\omega)-\tau_{h,i}(\omega))\geq\delta_{aucl}(\varepsilon_{h},\overline{m},\delta_{cp})\equiv\Delta_{m(J(h))},\label{eq:temp-342-2}
\end{equation}
with
\begin{equation}
d(X(\tau_{h,i}(\omega),\omega),X(\cdot,\omega))\leq\varepsilon_{h},\label{eq:temp-307-4-1-2}
\end{equation}
on the interval $\theta_{h,i}(\omega)\equiv[\tau_{h,i}(\omega),\tau_{h,i+1}(\omega))$
or $\theta_{h,i}(\omega)\equiv[\tau_{h,i}(\omega),\tau_{h,i+1}(\omega)]$
according as $i\leq p_{h}-2$ or $i=p_{h}-1$. 

9. Now let $h\geq3$ and $\omega\in A_{h}$ be arbitrary. Inequality
\ref{eq:temp-371} implies 
\[
\lambda^{-1}>2(\kappa+\overline{\lambda}^{-1}(1+\overline{u})(\kappa+1))=(2\kappa+\overline{\lambda}^{-1}(1+\overline{u})(2\kappa+2)),
\]
whence 
\begin{equation}
\eta\equiv\lambda^{-1}-(2\kappa+\overline{\lambda}^{-1}(1+\overline{u})(2\kappa+2))>0.\label{eq:temp-379}
\end{equation}
Define 
\[
\delta_{h}\equiv2^{-(h-2)\eta}.
\]
Then $\delta_{h}<1$. For each $i=0,\cdots,p_{h}-1$, define the subinterval
\[
\overline{\theta}_{h,i}(\omega)\equiv[\tau_{h,i}(\omega),\tau_{h,i+1}(\omega)-\delta_{h}(\tau_{h,i+1}(\omega)-\tau_{h,i}(\omega))],
\]
of $\theta_{h,i}(\omega),$ with Lebesgue measure
\[
\mu\overline{\theta}_{h,i}(\omega)=(1-\delta_{h})(\tau_{h,i+1}(\omega)-\tau_{h,i}(\omega)).
\]
Note that the intervals $\overline{\theta}_{h,0}(\omega),\cdots,\overline{\theta}_{h,p-1}(\omega)$
are mutually exclusive. Therefore
\[
\mu(\bigcup_{i=0}^{p(h)-1}\overline{\theta}_{h,i}(\omega))^{c}=\sum_{i=0}^{p(h)-1}\delta_{h}(\tau_{h,i+1}(\omega)-\tau_{h,i}(\omega))=\delta_{h}.
\]

10. Define the positive integers
\[
b_{0}\equiv b_{0}(u,w,\overline{u},\overline{\lambda})\equiv7\vee[\kappa^{-1}(-\kappa_{0}+2+\log_{2}K+\overline{\lambda}^{-1}(1+\overline{u})(2\kappa_{0}+(2\kappa+2)7+11))]_{1},
\]
\begin{equation}
b_{1}\equiv b_{1}(u,w,\overline{u},\overline{\lambda})\equiv[\kappa^{-1}\overline{\lambda}^{-1}(1+\overline{u})(2\kappa+2)]_{1},\label{eq:temp-374}
\end{equation}
\begin{equation}
J_{h}\equiv b_{0}+b_{1}h,\label{eq:temp-415}
\end{equation}
and 
\[
n_{h}\equiv h+7.
\]
Then $2^{-n(h)+6}=2^{-h-1}<2^{-h}\equiv\varepsilon_{h}$ and $J_{h}\equiv b_{0}+b_{1}h\geq7+h\equiv n_{h}$.
Moreover, according to the defining equalities \ref{eq:temp-292}
and \ref{eq:temp-415}, we have 
\[
m_{J(h)}\equiv\kappa_{0}+\kappa J_{h}\equiv\kappa_{0}+\kappa b_{0}+\kappa b_{1}h
\]
\[
>\kappa_{0}+\kappa b_{0}+\overline{\lambda}^{-1}(1+\overline{u})(2\kappa+2)h
\]
\[
\geq\kappa_{0}+(-\kappa_{0}+2+\log_{2}K+\overline{\lambda}^{-1}(1+\overline{u})(2\kappa_{0}+(2\kappa+2)7+11))+\overline{\lambda}^{-1}(1+\overline{u})(2\kappa+2)h
\]
\[
=2+\log_{2}K+\overline{\lambda}^{-1}(1+\overline{u})(2\kappa_{0}+(2\kappa+2)n_{h}+11)
\]
\begin{equation}
\equiv2+\log_{2}K+\overline{\lambda}^{-1}(1+\overline{u})(2m_{n(h)}+2n_{h}+11),\label{eq:temp-373}
\end{equation}
whence
\[
\Delta_{m(J(h))}\equiv2^{-m(J(h))}
\]
\[
<2^{-2}K^{-1}(2^{-2m(n(h))-2n(h)-11})^{(1+\overline{u})/\overline{\lambda}}=2^{-2}\delta_{cp}(2^{-2m(n(h))-2n(h)-10}),
\]
where the last equality is from the defining equality \ref{eq:temp-471}.
Therefore 
\[
\delta_{aucl}(\varepsilon_{h},\overline{m},\delta_{cp})\equiv\Delta_{m(J(h))}\equiv2^{-m(J(h))}.
\]
according to the defining equality \ref{eq:temp-473}. 

11. Now let $\varepsilon>0$ arbitrary. Let $k\geq3$ be so large
that $2^{-k+1}<\varepsilon$ and $2^{-(k-2)\eta}(1-2^{-\eta})<\varepsilon$.
Define the $\mathrm{measurable}$ subset 
\[
B_{k}\equiv\bigcap_{h=k}^{\infty}A_{h}
\]
of $\Omega$, where the $\mathrm{measurable}$ set $A_{h}$ was introduced
in Step 8, for each $h\geq k$. Then 
\[
P(B_{k}^{c})\leq\sum_{h=k}^{\infty}P(A_{h}^{c})<\sum_{h=k}^{\infty}2^{-h}=2^{-k+1}<\varepsilon.
\]
Consider each $\omega\in B_{k}$. Then $\omega\in A_{h}$ for each
$h\geq k$. Define the Lebesgue $\mathrm{measurable}$ subset 
\[
\widetilde{\theta}_{k}(\omega)\equiv\bigcap_{h=k}^{\infty}\bigcup_{i=0}^{p(h)-1}\overline{\theta}_{h,i}(\omega)
\]
of $[0,1]$. Then the Lebesgue measure $\mu(\widetilde{\theta}_{k}(\omega)^{c})$
is bounded by 
\[
\mu(\widetilde{\theta}_{k}(\omega)^{c})\leq\sum_{h=k}^{\infty}\mu(\bigcup_{i=0}^{p(h)-1}\overline{\theta}_{h,i}(\omega))^{c}=\sum_{h=k}^{\infty}\delta_{h}\equiv\sum_{h=k}^{\infty}2^{-(h-2)\eta}=2^{-(k-2)\eta}(1-2^{-\eta})<\varepsilon.
\]
Now let $t\in\widetilde{\theta}_{k}(\omega)\cap domain(X(\cdot,\omega)$
and 
\[
s\in\bigcup_{h=k}^{\infty}(t+\delta_{h+1}\Delta_{m(J(h+1))},t+\delta_{h}\Delta_{m(J(h))})\cap domain(X(\cdot,\omega)
\]
 be arbitrary. Then 
\begin{equation}
s\in(t+\delta_{h+1}\Delta_{m(J(h+1))},\:t+\delta_{h}\Delta_{m(J(h))})\label{eq:temp-474}
\end{equation}
for some $h\geq k\geq3$. Since $t\in\widetilde{\theta}_{k}(\omega)\subset\bigcup_{i=0}^{p(h)-1}\overline{\theta}_{h,i}(\omega)$,
there exists $i=0,\cdots,p_{h}-1$ such that 
\[
t\in\overline{\theta}_{h,i}(\omega)\equiv[\tau_{h,i}(\omega),\tau_{h,i+1}(\omega)-\delta_{h}(\tau_{h,i+1}(\omega)-\tau_{h,i}(\omega))].
\]
It follows, in view of inequality \ref{eq:temp-342-2}, that 
\[
s\in(t,t+\delta_{h}\Delta_{m(J(h))})\subset[\tau_{h,i}(\omega),\tau_{h,i+1}(\omega)-\delta_{h}(\tau_{h,i+1}(\omega)-\tau_{h,i}(\omega)+\delta_{h}\Delta_{m(J(h))})
\]
\[
\subset[\tau_{h,i}(\omega),\tau_{h,i+1}(\omega)-\delta_{h}(\tau_{h,i+1}(\omega)-\tau_{h,i}(\omega))+\delta_{h}(\tau_{h,i+1}(\omega)-\tau_{h,i}(\omega)))
\]
\[
=[\tau_{h,i}(\omega),\tau_{h,i+1}(\omega))\subset\theta_{h,i}(\omega).
\]
Hence inequality \ref{eq:temp-307-4-1-2} implies that
\begin{equation}
d(X(\tau_{h,i}(\omega),\omega),X(t,\omega))\vee d(X(\tau_{h,i}(\omega),\omega),X(s,\omega))\leq\varepsilon_{h},\label{eq:temp-307-1}
\end{equation}
and therefore that 
\begin{equation}
d(X(t,\omega),X(s,\omega))\leq2\varepsilon_{h}=2^{-h+1}.\label{eq:temp-380}
\end{equation}
At the same time, relation \ref{eq:temp-474} implies 
\begin{equation}
s-t\geq\delta_{h+1}\Delta_{m(J(h+1))}\equiv2^{-\eta(h-1)}2^{-m(J(h+1))}.\label{eq:temp-413}
\end{equation}
Moreover, from the defining equalities \ref{eq:temp-292}, \ref{eq:temp-415},
and \ref{eq:temp-374}, we have

\[
m_{J(h+1)}\equiv\kappa_{0}+\kappa J_{h+1}\equiv\kappa_{0}+\kappa b_{0}+\kappa b_{1}(h+1)
\]
\[
<(\kappa_{0}+\kappa b_{0}+\kappa b_{1})+\kappa(\kappa^{-1}\overline{\lambda}^{-1}(1+\overline{u})(2\kappa+2)+2)h
\]
\[
=(\kappa_{0}+\kappa b_{0}+\kappa b_{1})+(\overline{\lambda}^{-1}(1+\overline{u})(2\kappa+2)+2\kappa)h.
\]
Hence
\[
(h-1)\eta+m_{J(h+1)}<((\kappa_{0}+\kappa b_{0}+\kappa b_{1})-\eta)+(\eta+\overline{\lambda}^{-1}(1+\overline{u})(2\kappa+2)+2\kappa)h
\]
\[
\equiv((\kappa_{0}+\kappa b_{0}+\kappa b_{1})-\eta)
\]
\[
+(\lambda^{-1}-(2\kappa+\overline{\lambda}^{-1}(1+\overline{u})(2\kappa+2))+\overline{\lambda}^{-1}(1+\overline{u})(2\kappa+2))+2\kappa)h
\]
\[
=(\kappa_{0}+\kappa b_{0}+\kappa b_{1}-\eta)+\lambda^{-1}h,
\]
and so
\[
((h-1)\eta+m_{J(h+1)})\lambda<(\kappa_{0}+\kappa b_{0}+\kappa b_{1}-\eta)\lambda+h.
\]
Consequently, inequalities \ref{eq:temp-380} and \ref{eq:temp-413}
together yield 
\[
d(X(t,\omega),X(s,\omega))(s-t)^{-\lambda}\leq2^{-h+1}(2^{\eta(h-1)}2^{m(J(h+1))})^{\lambda}
\]
\[
<2^{-h+1}2^{(\kappa(0)+\kappa b(0)+\kappa b(1)-\eta)\lambda+h}=C_{0}\equiv2^{\kappa(0)+\kappa b(0)+\kappa b(1)-\eta)\lambda+1}.
\]
Thus 
\begin{equation}
d(X(t,\omega),X(s,\omega))\leq C_{0}(s-t)^{\lambda}\label{eq:temp-419}
\end{equation}
where $\omega\in B_{k}$, $t\in\widetilde{\theta}_{k}(\omega)\cap domain(X(\cdot,\omega)$,
and 
\[
s\in G\equiv\bigcup_{h=k}^{\infty}(t+\delta_{h+1}\Delta_{m(J(h+1))},t+\delta_{h}\Delta_{m(J(h))})\cap domain(X(\cdot,\omega))
\]
are arbitrary. Since the set $G$ is a dense subset of $[t,t+\delta_{k}\Delta_{m(J(k))})\cap domain(X(\cdot,\omega))$,
the right continuity of the function $X(\cdot,\omega)$ implies that
inequality \ref{eq:temp-419} holds for each $\omega\in B_{k}$, $t\in\widetilde{\theta}_{k}(\omega)\cap domain(X(\cdot,\omega))$,
and $s\in[t,t+\delta_{k}\Delta_{m(J(k))})\cap domain(X(\cdot,\omega))$.
Since $P(B_{k}^{c})<\varepsilon$ and $\mu(\widetilde{\theta}_{k}(\omega)^{c})<\varepsilon$,
where $\varepsilon>0$ is arbitrarily small, the process $X$ is right
Hoelder by Definition \ref{Def. right Hoelder cadlag functions and processes}. 
\end{proof}
As an easy corollary, the following Theorem \ref{Thm Sufficient condition for a time scaled right Hoelder process}
generalizes the preceding Theorem \ref{Thm Sufficient condition  right Hoelder process}.
The proof of Part 1 by means of a deterministic time scaling is a
new and simple proof of Theorem 13.6 of \cite{Billingsley 99}. Part
2, a condition for right Hoelder property, seems to be new.
\begin{thm}
\label{Thm Sufficient condition for a time scaled right Hoelder process}
\textbf{\emph{(Sufficient condition for a time-scaled right Hoelder
process).}} Let $u\geq0$ and $w>0$ be arbitrary. Let $G:[0,1]\rightarrow[0,1]$
be a nondecreasing continuous function. Let $F$ be an arbitrary consistent
family of f.j.d.'s  with parameter set $[0,1]$ and state space $(S,d)$.
Suppose $F$ is continuous in probability, and suppose 
\begin{equation}
F_{t,r,s}\{(x,y,z)\in S^{3}:d(x,y)\wedge d(y,z)>b\}\leq b^{-u}(G(s)-G(t))^{1+w}\label{eq:temp-472-3-1-1}
\end{equation}
for each $b>0$ and for each $t\leq r\leq s$ in $[0,1]$. Then the
following holds.

1. There exists an a.u. càdlàg process $Y:[0,1]\times\Omega\rightarrow S$
with marginal distributions given by the family $F$. 

2. Suppose, in addition, that there exist $\overline{u},\overline{\lambda}>0$
such that 
\begin{equation}
F_{t,s}(d>b)\leq b^{-\overline{u}}|G(s)-G(t)|^{\overline{\lambda}}\label{eq:temp-150-2}
\end{equation}
for each $b>0$ and for each $t,s\in[0,1]$. Then $Y(r,\omega)=X(h(r),\omega)$
for each $(r,\omega)\in domain(Y)$, for some right Hoelder process
$X$ and for some continuous increasing function $h:[0,1]\rightarrow[0,1]$.
\end{thm}
\begin{proof}
Write $a_{0}\equiv G(0)$ and $a_{1}\equiv G(1)$. Write $K\equiv a_{1}-a_{0}+1>0$.
Define the continuous increasing function $h:[0,1]\rightarrow[0,1]$
by $h(r)\equiv K^{-1}(G(r)-a_{0}+r)$ for each $r\in[0,1]$. Then
its inverse $g\equiv h^{-1}$ is also continuous and increasing. Moreover,
for each $s,t\in[0,1]$ with $t\leq s$, we have
\[
G(g(s))-G(g(t))=Kh(g(s))-g(s)+a_{0}-Kh(g(t))+g(t)-a_{0}
\]
\begin{equation}
=(Ks-Kt)-(g(s)-g(t))<Ks-Kt.\label{eq:temp-472}
\end{equation}

1. Since the family $F$ of f.j.d.'s is, by hypothesis, continuous
in probability, there exists, according to Theorem \ref{Thm. Measurable Extension of consistent family of fjd.s continuous in prob},
a process $V:[0,1]\times\Omega\rightarrow S$ with marginal distributions
given by the family $F$. Define the function $U:[0,1]\times\Omega\rightarrow S$
by $domain(U)\equiv\{(t,\omega):(g(t),\omega)\in domain(V)\}$and
by 
\begin{equation}
U(t,\omega)\equiv V(g(t),\omega)\label{eq:temp-475}
\end{equation}
for each $(t,\omega)\in domain(U)$. Then $U(t,\cdot)\equiv V(g(t),\cdot)$
is a r.v. for each $t\in[0,1]$. Thus $U$ is a stochastic process.
Let $F'$ denote the family of its marginal distributions. Then, for
each $b>0$, we have 
\[
F'_{t,r,s}\{(x,y,z)\in S^{3}:d(x,y)\wedge d(y,z)>b)=P(d(U_{t},U_{r})\wedge d(U_{r},U_{s})>b)
\]
\[
\equiv P(d(V_{g(t)},V_{g(r)})\wedge d(V_{g(r)},V_{g(s)})>b)=F{}_{g(t),g(r),g(s)}\{(x,y,z)\in S^{3}:d(x,y)\wedge d(y,z)>b)
\]
\[
\leq b^{-u}(G(g(s))-G(g(t)))^{1+w}\leq b^{-u}(Ks-Kt)^{1+w},
\]
where the next to last inequality follows from inequality \ref{eq:temp-472-3-1-1}
in the hypothesis, and the last inequality is from inequality \ref{eq:temp-472}.
Thus the family $F'$ and the constants $K,u,w$ satisfy the hypothesis
of Part 2 of Theorem \ref{Thm Sufficient condition  right Hoelder process}.
Accordingly, there exists an a.u. càdlàg process $X:[0,1]\times\Omega\rightarrow S$
with marginal distributions given by the family $F'$. Now define
a process $Y:[0,1]\times\Omega\rightarrow S$ by $domain(Y)\equiv\{(r,\omega):(h(r),\omega)\in domain(X)\}$and
by 
\[
Y(r,\omega)\equiv X(h(r),\omega)
\]
for each $(r\omega)\in domain(Y)$. Because the function $h$ is continuous,
it can be easily verified that the process $Y$ is a.u. càdlàg. Moreover,
in view of the defining equality \ref{eq:temp-475}, we have
\[
V(r,\omega)=U(h(r),\omega)
\]
for each $(r\omega)\in domain(V)$. Since the processes $X$ and $U$
share the same marginal distributions given by the family $F'$, the
last two displayed equalities imply that the processes $Y$ and $V$
share the same marginal distributions. Since the process $V$ has
marginal distributions given by the family $F$, so does the process
$Y$. Assertion 1 of the theorem is proved.

2. Suppose, in addition, that there exist $\overline{u},\overline{\lambda}>0$
inequality \ref{eq:temp-150-2} holds for each $b>0$ and for each
$t,s\in[0,1]$. Then, for each $b>0$, we have 
\[
F'_{t,s}(d>b)=P(d(U_{t},U_{s})>b)
\]
\[
\equiv P(d(V_{g(t)},V_{g(s)})>b)=F_{g(t),g(s)}(d>b)
\]
\[
\leq b^{-\overline{u}}|G(g(s))-G(g(t))|^{\overline{\lambda}}\leq b^{-\overline{u}}|Ks-Kt|^{\overline{\lambda}},
\]
where the first inequality is from inequality \ref{eq:temp-150-2},
and the last is from inequality \ref{eq:temp-472}. Thus the family
$F'$ and the constants $K,\overline{u},\overline{\lambda}$ satisfy
the hypothesis of Part 3 of Theorem \ref{Thm Sufficient condition  right Hoelder process}.
Accordingly, the a.u. càdlàg process $X:[0,1]\times\Omega\rightarrow S$
is right Hoelder. The theorem is proved.
\end{proof}

\section{a.u. Càdlàg Processes and a.u. Continuous Processes with Parameter
Set $[0,\infty)$}

In the preceding sections, a.u. continuous processes and a.u. càdlàg
processes are studied with the unit interval $[0,1]$ as the parameter
set. Generalization to parameter intervals $[0,\infty)$ is straightforward.
Specifically, we can study the process piecewise on unit subintervals
$[M,M+1]$, for $M=0,1,\cdots$, using the results in the preceding
sections, and then stitching the results back together. 

For a.u. continuous processes, we need only state the next definition.
\begin{defn}
\label{Def. a.u. continuout process on =00005B0,inf)}\textbf{(a.u.
Continuous process on $[0,\infty)$. }A process
\[
X:[0,\infty)\times(\Omega,L,E)\rightarrow S
\]
is said to be a.u. continuous if, for each $M\geq0$, the shifted
process $X^{M}:[0,1]\times\Omega\rightarrow S$, defined by $X^{M}(t)\equiv X(M+t)$
for each $t\in[0,1]$, is a.u. continuous, in the sense of Definition
\ref{Def. continuity in prob, continuity a.u., and a.u. continuity},
with some modulus of a.u. càdlàg $\delta_{aucl}^{M}$ and with some
modulus of continuity in probability $\delta_{Cp}^{M}$. 
\end{defn}
For a.u. càdlàg processes, we will state several related definitions,
and give a straightforward proof of the theorem which extends an arbitrary
$D$-regular process on $\overline{Q}_{\infty}$ to an a.u. càdlàg
process on $[0,\infty)$. Recall here that $Q_{\infty}$ and $\overline{Q}_{\infty}$
are the enumerated sets of dyadic rationals in $[0,1]$ and $[0,\infty)$
respectively. 

Let $(\Omega,L,E)$ be a sample space. Let $(S,d)$ be a locally compact
metric space, which will serve as the state space. Let $\xi\equiv(A_{q})_{q=1,2,\cdots}$
be a binary approximation of $(S,d)$. Recall that $D[0,1]$ stands
for the space of càdlàg functions on $[0,1]$, and that $\widehat{D}[0,1]$
stands for the space of a.u. càdlàg processes on $[0,1]$.
\begin{defn}
\textbf{\label{Def. Metric space of cadlag functions on =00005B0,inf)}
(Skorokhod metric space of càdlàg functions on $[0,\infty)).$} Let
$x:[0,\infty)\rightarrow S$ be an arbitrary function whose domain
contains the enumerated set $\overline{Q}_{\infty}$ of dyadic rationals.
Let $M\geq0$ be arbitrary. Then the function $x$ is said to be \emph{càdlàg}
on the interval $[M,M+1]$ if the shifted function $x^{M}:[0,1]\rightarrow S$,
defined by $x^{M}(t)\equiv x(M+t)$ for each $t$ with $M+t\in domain(x)$,
is a member of $D[0,1]$. The function $x:[0,\infty)\rightarrow S$
is said to be \emph{càdlàg} \index{càdlàg function on [0,infty)@càdlàg function on {[}0,$\infty$)}
if it is càdlàg on the interval $[M,M+1]$ for each $M\geq0$. We
will write $D[0,\infty)$ for the set of  càdlàg functions on $[0,\infty)$. 

Recall from Definition \ref{Def. Skorokhod Metric} the Skorokhod
metric $d_{D[0,1]}$ on $D[0,1]$. Define the \emph{Skorokhod metric}
on $D[0,\infty)$ by 
\[
d_{D[0,\infty)}(x,y)\equiv\sum_{M=0}^{\infty}2^{-M-1}(1\wedge d_{D[0,1]}(x^{M},y^{M}))
\]
for each $x,y\in D[0,\infty)$. We will call $(D[0,\infty),d_{D[0,\infty)})$
the \emph{Skorokhod space} \index{Skorokhod space on [0,infty)@Skorokhod space on $[0,\infty)$}
on $[0,\infty)$.
\end{defn}
$\square$
\begin{defn}
\label{Def. D-regular processes wirh dyadic rationals in =00005B0,oo)as parameters-1}
\textbf{($D$-regular processes with parameter set $\overline{Q}_{\infty}$).}
Let $Z:\overline{Q}_{\infty}\times\Omega\rightarrow S$ be a stochastic
process. Recall from Definition \ref{Def.  D-regular process and D-regular family of f.j.d.s on Q_inf}
the metric space $(\widehat{R}_{Dreg}(Q_{\infty}\times\Omega,S),\widehat{\rho}_{Prob,Q(\infty)})$
of $D$-regular processes with parameter set $Q_{\infty}$. 

Suppose, for each $M\geq0$, (i) the process $Z|\overline{Q}_{\infty}[0,M+1]$
is continuous in probability, with a modulus of continuity in probability
$\delta{}_{Cp,M}$, and (ii) the shifted process $Z^{M}:Q_{\infty}\times\Omega\rightarrow S$,
defined by $Z^{M}(t)\equiv Z(M+t)$ for each $t\in Q_{\infty}$, is
a member of the space $\widehat{R}_{Dreg}(Q_{\infty}\times\Omega,S)$,
with a modulus of $D$-regularity $\overline{m}_{M}$. 

Then the process $Z:\overline{Q}_{\infty}\times\Omega\rightarrow S$
is said to be $D$-regular\index{D-regular process on overline{Q}_{infty}@$D$-regular process on $\overline{Q}_{\infty}$},
with a modulus of continuity in probability $\widetilde{\delta}{}_{Cp}$$\equiv(\delta{}_{Cp,M})_{M=0,1,\cdots},$
and with a modulus of $D$-regularity $\widetilde{m}\equiv(\overline{m}_{M})_{M=0,1,\cdots}$.
Let $\widehat{R}_{Dreg}(\overline{Q}_{\infty}\times\Omega,S)$ denote
the set of all  $D$-regular processes with parameter set $\overline{Q}_{\infty}$.
Let $\widehat{R}_{Dreg,\widetilde{\delta}(Cp),\widetilde{m}.}(\overline{Q}_{\infty}\times\Omega,S)$
denote the subset whose members share the common modulus of continuity
in probability $\widetilde{\delta}{}_{Cp}$ and the common modulus
of $D$-regularity $\widetilde{m}\equiv(\overline{m}_{M})_{M=0,1,\cdots}$.
If, in addition, $\delta{}_{Cp,M}=\delta{}_{Cp,0}$ and $\overline{m}_{M}=\overline{m}_{0}$
for each $M\geq0$, then we say that the process $Z$ \emph{time-uniformly}
\emph{$D$-regular} \index{time-uniformly D-regular process on overline{Q}_{infty}@time-uniformly $D$-regular process on $\overline{Q}_{\infty}$}
on $\overline{Q}_{\infty}$. 
\end{defn}
$\square$
\begin{defn}
\textbf{\label{Def. Metric space of a.u. cadlag process on =00005B0,inf)}
(Metric space of a.u. càdlàg process on $[0,\infty)).$} Let 
\[
X:[0,\infty)\times(\Omega,L,E)\rightarrow S
\]
be an arbitrary process. Suppose, for each $M\geq0$, (i) the process
$X|[0,M+1]$ is continuous in probability, with a modulus of continuity
in probability $\delta{}_{Cp,M}$, and (ii) the shifted process $X^{M}:[0,1]\times\Omega\rightarrow S$,
defined by $X(t)\equiv X(M+t)$ for each $t\in[0,\infty)$, is a member
of the space $\widehat{D}[0,1]$, with some modulus of a.u. càdlàg
$\delta_{aucl}^{M}$. 

Then the process $X$ is said to be \emph{a.u. càdlàg} \index{a.u. càdlàg process on [0,infty)@a.u. càdlàg process on {[}0,$\infty$)}
on the interval $[0,\infty)$, with a modulus of a.u. càdlàg $\widetilde{\delta}_{aucl}\equiv(\delta_{aucl}^{M})_{M=0,1,\cdots}$
and with a modulus of continuity in probability $\widetilde{\delta}_{Cp}\equiv(\delta_{Cp}^{M})_{M=0,1,\cdots}$.
If, in addition, $\delta_{aucl}^{M}=\delta_{aucl}^{0}$ and $\delta_{Cp}^{M}=\delta_{Cp}^{0}$
for each $M\geq0$ , then we say that the process is \emph{time-uniformly}
\emph{a.u. càdlàg} \index{time-uniformly a.u.càdlàg process on [0,infty)@time-uniformly a.u.càdlàg process on {[}0,$\infty$)}
on the interval $[0,\infty)$. 

We will write $\widehat{D}[0,\infty)$ for the set of  a.u. càdlàg
processes on $[0,\infty)$, and equip it ewith the metric $\widetilde{\rho}_{\widehat{D}[0,\infty)}$
defined by
\[
\widetilde{\rho}_{\widehat{D}[0,\infty)}(X,X')\equiv\rho_{Prob,\overline{Q}(\infty)}(X|\overline{Q}_{\infty},X'|\overline{Q}_{\infty})
\]
for each $X,X'\in\widehat{D}[0,\infty)$, where,according to Definition
\ref{Def. Metric on random fields w/ countable parameters}, we have
\emph{
\begin{equation}
\rho_{Prob,\overline{Q}(\infty)}(X|\overline{Q}_{\infty},X'|\overline{Q}_{\infty})\equiv E\sum_{n=0}^{\infty}2^{-n-1}(1\wedge d(X_{u(n)},X'_{u(n)})).\label{eq:temp-445-2-1-1}
\end{equation}
}Thus
\[
\widetilde{\rho}_{\widehat{D}[0,\infty)}(X,X')\equiv E\sum_{n=0}^{\infty}2^{-n-1}(1\wedge d(X_{u(n)},X'_{u(n)})
\]
for each $X,X'\in\widehat{D}[0,\infty)$.

Let $\widehat{D}_{\widetilde{\delta}(aucl),\widetilde{\delta}(Cp)}[0,\infty)$
denote the subspace of the metric space $(\widehat{D}[0,\infty),\widetilde{\rho}_{\widehat{D}[0,\infty)})$
whose members share a common modulus of continuity in probability
$\widetilde{\delta}_{Cp}\equiv(\delta{}_{Cp,M})_{M=0,1,\cdots}$ and
share a common modulus of a.u. càdlàg $\widetilde{\delta}_{aucl}\equiv(\delta_{aucl}(\cdot,\overline{m}_{M},\delta_{Cp,M}))_{M=0,1,\cdots}$.
\end{defn}
$\square$

In the following, recall the right-limit extension functions $\Phi_{rLim}$
and $\overline{\Phi}_{rLim}$ from Definition \ref{Def. Extension of process by right limit}.
\begin{lem}
\textbf{\emph{\label{Lem.  Right limit extension of D-regular processs on QBarInf is contionuou in probabiity}}}
\textbf{\emph{(The right-limit extension of a $D$-regular process
on $\overline{Q}_{\infty}$ is continuous in probability on $[0,\infty)$).}}
Suppose 
\[
Z:\overline{Q}_{\infty}\times(\Omega,L,E)\rightarrow(S,d)
\]
is an arbitrary $D$-regular process on $\overline{Q}_{\infty}$,
with a modulus of continuity in probability $\widetilde{\delta}_{Cp}$$\equiv(\delta{}_{Cp,M})_{M=0,1,\cdots}.$
Then the right-limit extension
\[
X\equiv\overline{\Phi}_{rLim}(Z):[0,\infty)\times(\Omega,L,E)\rightarrow(S,d)
\]
of $Z$ is a well-defined process which is continuous in probability,
with the same modulus of continuity in probability $\widetilde{\delta}_{Cp}\equiv(\delta{}_{Cp,M})_{M=0,1,\cdots}$
as $Z$.
\end{lem}
\begin{proof}
Let $N\geq0$ be arbitrary. Consider each $t,t'\in[0,N+1]$. Let $(r_{k})_{k=1,2,\cdots}$
be a sequence in $[0,N+1]\overline{Q}_{\infty}$such that $r_{k}\downarrow t$.
Since the process $Z$ is continuous in probability on $[0,N+1]\overline{Q}_{\infty}$,
with a modulus of continuity in probability $\delta{}_{Cp,N}$, it
follows that 
\[
E(1\wedge d(Z_{rk)},Z_{r(h)}))\rightarrow0
\]
as $k,h\rightarrow0$. By passing to a subsequenc if necessary, we
can assume 
\[
E(1\wedge d(Z_{r(k)},Z_{r(k+1)}))\leq2^{-k}
\]
for each $k\geq1$. Then
\begin{equation}
E(1\wedge d(Z_{r(k)},Z_{r(h)}))\leq2^{-k+1}\label{eq:temp-535-2}
\end{equation}
for each $h\geq k\geq1$. It follows that $Z_{r(k)}\rightarrow U_{t}$
a.u. for some r.v. $U_{t}$. Consequently $Z_{rk)}(\omega)\rightarrow U_{t}(\omega)$
for each $\omega$ in some full set $B$. Consider each $\omega\in B$.
Since $X\equiv\overline{\Phi}_{rLim}(Z)$ and since $r_{k}\downarrow t$
with $Z_{r(k)}(\omega)\rightarrow U_{t}(\omega)$, we see that $\omega\in domain(X_{t})$
and that $X_{t}(\omega)=U_{t}(\omega)$. In short, $X_{t}=U_{t}$
a.s. Consequently $X_{t}$ is a r.v., where $t\in[0,N+1]$ and $N\geq0$
are arbitrary. Since $[0,\infty)=\bigcup_{M=0}^{\infty}[0,M+1]$,
it follows that $X_{u}$ is a r.v. for each $u\in[0,\infty)$. Thus
$X:[0,\infty)\times(\Omega,L,E)\rightarrow(S,d)$ is a well-define
process.

Letting $h\rightarrow\infty$ in equality \ref{eq:temp-535-2}, we
obtain
\[
E(1\wedge d(Z_{r(k)},X_{t}))=E(1\wedge d(Z_{r(k)},U_{t}))\leq2^{-k+1}
\]
for each $k\geq1$. Similarly we can construct a sequence $(r'_{k})_{k=1,2,\cdots}$
in $[0,N+1]\overline{Q}_{\infty}$such that
\[
E(1\wedge d(Z_{r'(k)},X_{t'}))\leq2^{-k+1}
\]
Now let $\varepsilon>0$ be arbitrary, and suppose $t,t'\in[0,N+1]$
are such that $|t-t'|<\delta{}_{Cp,N}(\varepsilon)$. Then $|r_{k}-r'_{k}|<\delta{}_{Cp,N}(\varepsilon)$,
whence 
\[
E(1\wedge d(Z_{r(k)},Z_{r'(k)}))\leq\varepsilon
\]
for sufficiently large $k\geq1$. Combinning, we obtain
\[
E(1\wedge d(X_{t},X_{t'}))\leq2^{-k+1}+\varepsilon+2^{-k+1}
\]
for sufficiently large $k\geq1$. Hence
\[
E(1\wedge d(X_{t},X_{t'}))\leq\varepsilon
\]
where $t,t'\in[0,N+1]$ are arbitrary with $|t-t'|<\delta{}_{Cp,N}(\varepsilon)$.
Thus we have verified that the process $X$ is continuous in probability,
with the modulus of continuity in probability $\widetilde{\delta}_{Cp}\equiv(\delta_{Cp.N})_{N=0,1,\cdots}$
which is the same as the modulus of continuity in probability of the
given $D$-regular process. 
\end{proof}
\begin{thm}
\textbf{\emph{\label{Thm.  Right limit extension of D-regular processs on QBarInf is a.u. cadlag}}}
\textbf{\emph{(The right-limit extension of a $D$-regular process
on $\overline{Q}_{\infty}$ is an a.u. càdlàg process on $[0,\infty)$).}}
Suppose 
\[
Z:\overline{Q}_{\infty}\times(\Omega,L,E)\rightarrow(S,d)
\]
is an arbitrary $D$-regular process on $\overline{Q}_{\infty}$,
with a modulus of continuity in probability $\widetilde{\delta}_{Cp}$$\equiv(\delta{}_{Cp,M})_{M=0,1,\cdots},$
and with a modulus of $D$-regularity $\widetilde{m}\equiv(\overline{m}_{M})_{M=0,1,\cdots}$.
In symbols, suppose 
\[
Z\in\widehat{R}_{Dreg,\widetilde{\delta}(Cp),\widetilde{m}.}(\overline{Q}_{\infty}\times\Omega,S).
\]
Then the right-limit extension
\[
X\equiv\overline{\Phi}_{rLim}(Z):[0,\infty)\times(\Omega,L,E)\rightarrow(S,d)
\]
of $Z$ is an a.u. càdlàg process, with the same modulus of continuity
in probability $\widetilde{\delta}_{Cp}\equiv(\delta{}_{Cp,M})_{M=0,1,\cdots}$
as $Z$, and with a modulus of a.u. càdlàg $\widetilde{\delta}_{aucl}\equiv\widetilde{\delta}_{aucl}(\widetilde{m},\widetilde{\delta}_{Cp})$.
In other words,

\[
X\equiv\overline{\Phi}_{rLim}(Z)\in\widehat{D}_{\widetilde{\delta}(aucl,\widetilde{m},\widetilde{\delta}(Cp)),\widetilde{\delta}(Cp)}[0,\infty)\subset\widehat{D}[0,\infty).
\]
\end{thm}
\begin{proof}
1. Lemma \ref{Lem.  Right limit extension of D-regular processs on QBarInf is contionuou in probabiity}
says that the process $X$ is continuous in probability, with the
same modulus of continuity in probability $\widetilde{\delta}_{Cp}$$\equiv(\delta{}_{Cp,M})_{M=0,1,\cdots}$
as $Z$. Let $N\geq0$ be arbitrary. Then $Z^{N}:Q_{\infty}\times(\Omega,L,E)\rightarrow(S,d)$
is a $D$-regular process with a modulus of continuity in probability
$\delta{}_{Cp,N},$ and with a modulus of $D$-regularity $\overline{m}_{N}$.
Theorem \ref{Thm. Extension of D-regular process by right limit is a.u.cadlag}
therefore implies that the right-limit extension process 
\[
Y_{N}\equiv\Phi_{rLim}(Z^{N}):[0,1]\times(\Omega,L,E)\rightarrow(S,d)
\]
is a.u. càdlàg, with the same modulus of continuity in probability
$\delta{}_{Cp,N},$ and with a modulus of a.u. càdlàg $\delta_{aucl}(\cdot,\overline{m}_{N},\delta_{Cp,N})$.
Separately, Proposition \ref{Prop. Right-limit extenaion of D-regular process is continuous a.u.}
implies that the process $Y_{N}$ is continuous a.u., with a modulus
of continuity a.u. $\delta_{cau}(\cdot,\overline{m}_{N},\delta_{Cp,N})$.
Here the reader is reminded that continuity a.u. is not to be confused
with the much stronger condition of a.u. continuity. Note that $Y_{N}=Z^{N}$
on $Q_{\infty}$, and that $X=Z$ on $\overline{Q}_{\infty}$. 

2. Let $k\geq1$ be arbitrary. Define
\[
\delta_{k}\equiv\delta_{cau}(2^{-k},\overline{m}_{N},\delta_{Cp,N})\wedge\delta_{cau}(2^{-k},\overline{m}_{N+1},\delta_{Cp,N+1})\wedge2^{-k}.
\]
Then, since $\delta_{cau}(\cdot,\overline{m}_{N},\delta_{Cp,N})$
is a modulus of continuity a.u. of the process $Y_{N}$, there exists,
according to Definition \ref{Def. continuity in prob, continuity a.u., and a.u. continuity},
a $\mathrm{measurable}$ set $D_{1,k}\subset domain(Y_{N,1})$ with
$P(D_{1,k}^{c})<2^{-k}$ such that for each $\omega\in D_{1,k}$ and
for each $r\in domain(Y_{N}(\cdot,\omega))$ with $|r-1|<\delta_{k}$,
we have
\begin{equation}
d(Y_{N}(r,\omega),Z(N+1,\omega))=d(Y_{N}(r,\omega),Y_{N}(1,\omega))\leq2^{-k}.\label{eq:temp-539}
\end{equation}
Likewise, since $\delta_{cau}(\cdot,\overline{m}_{N+1},\delta_{Cp,N+1})$
is a modulus of continuity a.u. of the process $Y_{N+1}$, there exists,
according to Definition \ref{Def. continuity in prob, continuity a.u., and a.u. continuity},
a $\mathrm{measurable}$ set $D_{0,k}\subset domain(Y_{N+1,0})$ with
$P(D_{0,k}^{c})<2^{-k}$ such that for each $\omega\in D_{0,k}$ and
for each $r\in domain(Y_{N+1}(0,\omega))$ with $|r-0|<\delta_{k}$,
we have
\begin{equation}
d(Y_{N+1}(r,\omega),Z(N+1,\omega))=d(Y_{N+1}(r,\omega),Y_{N+1}(0,\omega))\leq2^{-k}.\label{eq:temp-539-4}
\end{equation}
Define $D_{k+}\equiv\bigcap_{h=k}^{\infty}D_{1,h}D_{0,h}$ and $B\equiv\bigcup_{\kappa=1}^{\infty}D_{\kappa+.}$
Then $P(D_{k+}^{c})<2^{-k+2}.$ Hence $P(B)=1$. In words, $B$ is
a full set.

3. Consider each $t\in[N,N+1)$. Since $Y_{N}\equiv\Phi_{rLim}(Z^{N})$
and $X\equiv\overline{\Phi}_{rLim}(Z)$, we have 
\[
domain(Y_{N,t-N})\equiv\{\omega\in\Omega:\lim_{s\rightarrow t-N;s\in[t-N,\infty)Q(\infty)}Z_{s}^{N}(\omega)\;\mathrm{exists}\}
\]
\[
=\{\omega\in\Omega:\lim_{s\rightarrow t-N;s\in[t-N,1]Q(\infty)}Z_{s}^{N}(\omega)\;\mathrm{exists}\}
\]
\[
=\{\omega\in\Omega:\lim_{s\rightarrow t-N;s\in[t-N,1]Q(\infty)}Z(N+s,\omega)\;\mathrm{exists}\}
\]
\[
=\{\omega\in\Omega:\lim_{r\rightarrow t;r\in[t,N+1]\overline{Q}(\infty)}Z(r,\omega)\;\mathrm{exists}\}
\]
\[
=\{\omega\in\Omega:\lim_{r\rightarrow t;r\in[t,\infty)\overline{Q}(\infty)}Z(r,\omega)\;\mathrm{exists}\}
\]
\[
\equiv domain(X_{t}),
\]
because each limit which appears in the previous equality exists iff
all others exist, in which case they are equal. Hence
\[
X_{t}(\omega)=\lim_{r\rightarrow t;r\in[t,\infty)\overline{Q}(\infty)}Z(r,\omega)=\lim_{s\rightarrow t-N;s\in[t-N,\infty)Q(\infty)}Z_{s}^{N}(\omega)=Y_{N,t-N}(\omega)
\]
for each $\omega\in domain(Y_{N,t-N})$. Thus the two functions $X_{t}$
and $Y_{N,t-N}$ have the same domain, and have equal values on the
common domain. In short, 
\begin{equation}
X_{t}=Y_{N,t-N},\label{eq:temp-536}
\end{equation}
where $t\in[N,N+1)$ is arbitrary. As for the end point $t=N+1$,
we have, trivially,
\[
X_{N+1}=Z_{N+1}=Z_{1}^{N}=Y_{N,1}=Y_{N,(N+1)-1}.
\]
Hence 
\begin{equation}
X_{t}=Y_{N,t-N},\label{eq:temp-536-2}
\end{equation}
for each $t\in[N,N+1)\cup\{N+1\}$. 

4. We wish to extend equality \ref{eq:temp-536-2} to each $t\in[N,N+1]$.
To that end, consider each 
\[
t\in[N,N+1].
\]
We will prove that 
\begin{equation}
X_{t}=Y_{N,t-N}\label{eq:temp-552}
\end{equation}
 on the full set $B$.

5. To that end, let $\omega\in B$ be arbitrary. Suppose $\omega\in domain(X_{t})$.
Then, since $X\equiv\overline{\Phi}_{rLim}(Z)$, the limit $\lim_{r\rightarrow t;r\in[t,\infty)\overline{Q}(\infty)}Z(r,\omega)$
exists and is equal to $X_{t}(\omega)$. Consequently, each of the
following limits
\[
\lim_{r\rightarrow t;r\in[t,N+1]\overline{Q}(\infty)}Z(r,\omega)=\lim_{N+s\rightarrow t;s\in[t-N,1]\overline{Q}(\infty)}Z(N+s,\omega)
\]
\[
=\lim_{s\rightarrow t-N;s\in[t-N,1]\overline{Q}(\infty)}Z(N+s,\omega)=\lim_{s\rightarrow t-N;s\in[t-N,1]Q(\infty)}Z_{s}^{N}(\omega)
\]
exists and is equal to $X_{t}(\omega)$. Therefore, since $Y_{N}\equiv\Phi_{rLim}(Z^{N}),$
the existence of the last limit implies that $\omega\in domain(Y_{N,t-N})$,
and that $Y_{N,t-N}(\omega)=X_{t}(\omega)$. Thus 
\begin{equation}
domain(X_{t})\subset domain(Y_{N,t-N})\label{eq:temp-556}
\end{equation}
and 
\begin{equation}
X_{t}=Y_{N,t-N}\label{eq:temp-535}
\end{equation}
on $domain(X_{t})$. We have proved half of the desired equality \ref{eq:temp-552}.

6. Conversely, suppose $\omega\in domain(Y_{N,t-N})$. Then $y\equiv Y_{N,t-N}(\omega)\in S$
is defined. Hence, since $Y_{N}\equiv\Phi_{rLim}(Z^{N})$, the limit
\[
\lim_{s\rightarrow t-N;s\in[t-N,\infty)Q(\infty)}Z^{N}(s,\omega)
\]
exists, and is equal to $y$. Let $\varepsilon>0$ be arbitrary. Then
there exists $\delta'>0$ such that 
\[
d(Z^{N}(s,\omega),y)<\varepsilon
\]
for each 
\[
s\in[t-N,\infty)Q_{\infty}=[t-N,\infty)[0,1]\overline{Q}_{\infty}
\]
such that $s-(t-N)<\delta'$. In other words,
\begin{equation}
d(Z(u,\omega),y)<\varepsilon\label{eq:temp-555}
\end{equation}
for each $u\in[t,t+\delta')[N,N+1]\overline{Q}_{\infty}$. 

7. Separately, recall the assumption that $\omega\in B\equiv\bigcup_{\kappa=1}^{\infty}D_{\kappa+.}$.
Therefore there exists some $\kappa\geq1$ such that $\omega\in D_{\kappa+}\equiv\bigcap_{k=\kappa}^{\infty}D_{1,k}D_{0,k}$.
Take $k\geq\kappa$ so large that $2^{-k}<\varepsilon$. Then $\omega\in D_{1,k}D_{0,k}$.
Therefore, for each $r\in domain(Y_{N+1}(0,\omega))$ with $|r-0|<\delta_{k}$,
we have, according to inequality \ref{eq:temp-539-4}
\begin{equation}
d(Y_{N+1}(r,\omega),Z(N+1,\omega))\leq2^{-k}<\varepsilon.\label{eq:temp-539-3}
\end{equation}
Similarly, for each $r\in domain(Y_{N}(\cdot,\omega))$ with $|r-1|<\delta_{k}$,
we have, according to inequality \ref{eq:temp-539-4},
\begin{equation}
d(Y_{N}(r,\omega),Z(N+1,\omega))\leq2^{-k}<\varepsilon.\label{eq:temp-539-3-2}
\end{equation}

8. Now let $u,v\in[t,t+\delta_{k}\wedge\delta')\overline{Q}_{\infty}$be
arbitrary with $u<v$. Then $u,v\in[t,t+\delta_{k}\wedge\delta')[N,\infty)\overline{Q}_{\infty}.$
Since $u,v$ are dyadic rationals, there are three possibilities:
(i) $u<v\leq N+1$, (ii) $u\leq N+1<v$, or (iii) $N+1<u<v$. Consider
Case (i). Then $u,v\in[t,t+\delta')[N,N+1]\overline{Q}_{\infty}$.
Hence inequality \ref{eq:temp-555} applies to $u$ and $v$, to yield
\[
d(Z(u,\omega),y)\vee d(Z(v,\omega),y)<\varepsilon,
\]
whence
\[
d(Z(u,\omega),Z(v,\omega))<2\varepsilon.
\]

Next consider Case (ii) . Then $|(u-N)-1|<v-t<\delta_{k}$. Hence,
by inequality \ref{eq:temp-539-3-2}, we obtain
\[
d(Z(u,\omega),Z(N+1,\omega))\equiv d(Z^{N}(u-N,\omega),Z(N+1,\omega))
\]
\begin{equation}
=d(Y_{N}(u-N,\omega),Z(N+1,\omega))<\varepsilon.\label{eq:temp-539-5}
\end{equation}
Similarly, $|(v-(N+1))-0|<v-t<\delta_{k}$. Hence, by inequality \ref{eq:temp-539-3},
we obtain 
\[
d(Z(v,\omega),Z(N+1,\omega))\equiv d(Z^{N+1}(v-(N+1),\omega),Z(N+1,\omega))
\]
\begin{equation}
=d(Y_{N+1}(v-(N+1),\omega),Z(N+1,\omega))<\varepsilon.\label{eq:temp-539-4-1}
\end{equation}
Combining \ref{eq:temp-539-5} and \ref{eq:temp-539-4-1}, we obtain
\[
d(Z(u,\omega),Z(v,\omega))<2\varepsilon
\]
also in Case (ii).

Now consider Case (iii). Then $|(u-(N+1))-0|<v-t<\delta_{k}$ and
$|(v-(N+1))-0|<v-t<\delta_{k}$. Hence inequality \ref{eq:temp-539-3}
implies 
\[
d(Z(u,\omega),Z(N+1,\omega))=d(Z^{N+1}(u-(N+1),\omega),Z(N+1,\omega))
\]
\begin{equation}
=d(Y_{N+1}(u-(N+1),\omega),Z(N+1,\omega))<\varepsilon.\label{eq:temp-539-4-1-1}
\end{equation}
and, similarly,
\begin{equation}
d(Z(v,\omega),Z(N+1,\omega))<\varepsilon.\label{eq:temp-539-3-1}
\end{equation}
Hence 
\[
d(Z(u,\omega),Z(v,\omega))<2\varepsilon
\]
also in Case (iii).

9. Summing up, we see that $d(Z(u,\omega),Z(v,\omega))<2\varepsilon$
for each $u,v\in[t,t+\delta_{k}\wedge\delta')\overline{Q}_{\infty}$with
$u<v$. Since $\varepsilon>0$ is arbitrary, we conclude that $\lim_{u\rightarrow t;u\in[t,\infty)\overline{Q}(\infty)}Z(u,\omega)\;\mathrm{exists}$.
Thus 
\[
(t,\omega)\in domain(\overline{\Phi}_{rLim}(Z))\equiv domain(X).
\]
In other words, we have $\omega\in domain(X_{t})$, where $\omega\in B\cap domain(Y_{N,t-N})$
is arbitrary. Hence 
\[
B\cap domain(Y_{N,t-N})\subset domain(X_{t})\subset domain(Y_{N,t-N}),
\]
where the second inclusion is by relation \ref{eq:temp-556} in Step
5 above. Consequently
\[
B\cap domain(Y_{N,t-N})=B\cap domain(X_{t}),
\]
while, according to equality \ref{eq:temp-535},
\begin{equation}
X_{t-N}^{N}=X_{t}=Y_{N,t-N}\label{eq:temp-535-1}
\end{equation}
on $B\cap domain(X_{t})$, In other words, on the full subset $B,$
we have $X_{t-N}^{N}=Y_{N,t-N}$ for each $t\in[N,N+1].$ Equivalently,
$X^{N}=Y_{N}$ on the full subset $B$. Since the process 
\[
Y_{N}\equiv\Phi_{rLim}(Z^{N}):[0,1]\times(\Omega,L,E)\rightarrow(S,d)
\]
is a.u. càdlàg, with a modulus of a.u. càdlàg $\delta_{aucl}(\cdot,\overline{m}_{N},\delta_{Cp,N})$,
so is the process $X^{N}$. Thus, by Definition \ref{Def. Metric space of a.u. cadlag process on =00005B0,inf)},
the process $X$ is a.u. càdlàg, with modulus of continuity in probability
$\widetilde{\delta}_{Cp}\equiv(\delta{}_{Cp,M})_{M=0,1,\cdots}$,
and with a modulus of a.u. càdlàg $\widetilde{\delta}_{aucl}\equiv(\delta_{aucl}(\cdot,\overline{m}_{M},\delta_{Cp,M}))_{M=0,1,\cdots}$.
In other words,
\[
X\in\widehat{D}_{\widetilde{\delta}(aucl),\widetilde{\delta}(Cp)}[0,\infty).
\]
The theorem is proved.
\end{proof}
The next theorem is straightforward, and is verified here for ease
of future reference.
\begin{thm}
\textbf{\emph{\label{Thm. Phi^bar_rLim is isometry on properly restricted domain}($\overline{\Phi}_{rLim}$
is an isometry on a properly restricted domain).}} Recall from Definition
\ref{Def. D-regular processes wirh dyadic rationals in =00005B0,oo)as parameters-1}
the metric space $(\widehat{R}_{Dreg,\widetilde{\delta}(Cp),\widetilde{m}.}(\overline{Q}_{\infty}\times\Omega,S),\widehat{\rho}_{Prob,\overline{Q}(\infty)})$
of $D$-regular processes whose members $Z$ share a given modulus
of continuity in probability $\widetilde{\delta}_{Cp}$$\equiv(\delta{}_{Cp,M})_{M=0,1,\cdots},$
and share a given modulus of $D$-regularity $\widetilde{m}\equiv(\overline{m}_{M})_{M=0,1,\cdots}$. 

Recall from Definition \ref{Def. Metric space of a.u. cadlag process on =00005B0,inf)}
the metric space $(\widehat{D}[0,\infty),\widetilde{\rho}_{\widehat{D}[0,\infty)})$
of a.u. càdlàg processes on $[0,\infty)$, where
\[
\widetilde{\rho}_{\widehat{D}[0,\infty)}(X,X')\equiv\rho_{Prob,\overline{Q}(\infty)}(X|\overline{Q}_{\infty},X'|\overline{Q}_{\infty})
\]
for each $X,X'\in\widehat{D}[0,\infty)$. 

Then the function 
\[
\overline{\Phi}_{rLim}:(\widehat{R}_{Dreg,\widetilde{\delta}(Cp),\widetilde{m}.}(\overline{Q}_{\infty}\times\Omega,S),\widehat{\rho}_{Prob,\overline{Q}(\infty)})\rightarrow\widehat{D}_{\widetilde{\delta}(aucl,\widetilde{m},\widetilde{\delta}(Cp)),\widetilde{\delta}(Cp)}[0,\infty)\subset(\widehat{D}[0,\infty),\widetilde{\rho}_{\widehat{D}[0,\infty)})
\]
is a well-defined isometry on its domain, where the modulus of a.u.
càdlàg $\widetilde{\delta}_{aucl}\equiv\widetilde{\delta}_{aucl}(\widetilde{m},\widetilde{\delta}_{Cp})$
is defined in the proof below.
\end{thm}
\begin{proof}
1. Let $Z\in\widehat{R}_{Dreg,\widetilde{\delta}(Cp),\widetilde{m}.}(\overline{Q}_{\infty}\times\Omega,S)$
be arbitrary. In other words,
\[
Z,:\overline{Q}_{\infty}\times(\Omega,L,E)\rightarrow(S,d)
\]
is a $D$-regular process, with a modulus of continuity in probability
$\widetilde{\delta}_{Cp}$$\equiv(\delta{}_{Cp,M})_{M=0,1,\cdots},$
and with a modulus of $D$-regularity $\widetilde{m}\equiv(\overline{m}_{M})_{M=0,1,\cdots}$.
Consider each $N\geq0$. Then the shifted processes $Z^{N}:Q_{\infty}\times(\Omega,L,E)\rightarrow(S,d)$
is $D$-regular, with modulus of continuity in probability $\delta{}_{Cp,N}$
and modulus of $D$-regularity $\overline{m}_{N}.$ In other words,
\[
Z^{N}\in(\widehat{R}_{Dreg,\overline{m},\delta(Cp)}(Q_{\infty}\times\Omega,S),\widehat{\rho}_{Prob,Q(\infty)})
\]
Hence, by Theorem \ref{Thm.Construction of a.u. cadlag process by right limit of D-regular proceses is continuous.},
the processes $Y^{N}\equiv\Phi_{rLim}(Z^{N})$ is a.u. càdlàg, with
a modulus of continuity in probability $\delta{}_{Cp,N}$, and with
a modulus of a.u. càdlàg $\delta_{aucl}(\cdot,\overline{m}_{N},\delta{}_{Cp,N})$.
It is therefore easily verified that $X\equiv\overline{\Phi}_{rLim}(Z)$
is a well defined process on $[0,\infty)$, wiith 
\[
X^{N}\equiv\overline{\Phi}_{rLim}(Z)^{N}=\Phi_{rLim}(Z^{N})
\]
for each $N\geq0$. In other words, $\overline{\Phi}_{rLim}(Z)\equiv X\in\widehat{D}[0,\infty)$.
Hence the function $\overline{\Phi}_{rLim}$ is well-defined on $\widehat{R}_{Dreg,\widetilde{\delta}(Cp),\widetilde{m}.}(\overline{Q}_{\infty}\times\Omega,S)$.
In other words, $X\equiv\overline{\Phi}_{rLim}(Z)\in\widehat{D}_{\widetilde{\delta}(aucl,),\widetilde{\delta}(Cp)}[0,\infty)$,
where $\widetilde{\delta}_{Cp}\equiv(\delta{}_{Cp,M})_{M=0,1,\cdots}$
and where $\widetilde{\delta}_{aucl}\equiv(\delta_{aucl}(\cdot,\overline{m}_{M},\delta_{Cp,M}))_{M=0,1,\cdots}$.

2. It rremains to prove that the function $\overline{\Phi}_{rLim}$
is uniformly continuous on its domain. To that end, let $Z,Z'\in\widehat{R}_{Dreg,\widetilde{\delta}(Cp),\widetilde{m}.}(\overline{Q}_{\infty}\times\Omega,S)$
be arbitrary. Define $X\equiv\overline{\Phi}_{rLim}(Z)$ and $X'\equiv\overline{\Phi}_{rLim}(Z')$
as in the previous step. Then
\[
\widetilde{\rho}_{\widehat{D}[0,\infty)}(X,X')\equiv E\sum_{n=0}^{\infty}2^{-n-1}(1\wedge d(X_{u(n)},X'_{u(n)})
\]
\[
=E\sum_{n=0}^{\infty}2^{-n-1}(1\wedge d(Z_{u(n)},Z'_{u(n)})\equiv\widehat{\rho}_{Prob,\overline{Q}(\infty)}(Z,Z').
\]
Hence the function $\overline{\Phi}_{rLim}$ is an isometry on its
domain.
\end{proof}

\chapter{Markov Process}

In this chapter, we will construct an a.u. càdlàg Markov process from
a given Mrekov semigroup of transition distributions, and show that
the construction is a continuous mapping. 
\begin{defn}
\label{Def.  Specification of state space} \textbf{(Specificaton
of state space). }In this chapter, let $(S,d)$ be a locally compact
metric space, with an arbitrary, but fixed reference point. Let $\xi\equiv(A_{k})_{k=1,2,\cdots}$
be a given binary approximation of $(S,d)$ relative to $x_{\circ}$.
Let 
\[
\pi\equiv(\{g_{k,x}:x\in A_{k}\})_{k=1,2,\cdots}
\]
be the partition of unity of $(S,d)$ determined by $\xi$, as in
Definition \ref{Def. Partition of unity for locally compact (S,d)}.
\end{defn}
$\square$
\begin{defn}
\textbf{\label{Def. Notations for dyadic rationals-1} (Notations
for dyadic rationals).} Recall from Definition \ref{Def. Notations for dyadic rationals}
the notations related to the enumerated set $Q_{\infty}\equiv\{t_{0},t_{1},\cdots\}$
of dyadic rationals in the interval $[0,1]$, and the enumerated set
$\overline{Q}_{\infty}\equiv\{u_{0},u_{1},\cdots\}$ of dyadic rationals
in the interval $[0,\infty)$. 

In particular, for each $m\geq0$, recall the notations $p_{m}\equiv2^{m}$
and $\Delta_{m}\equiv2^{-m},$ and
\[
Q_{m}\equiv\{t_{0},t_{1},\cdots,t_{p(m)}\}=\{q_{m,0},\cdots,q_{m,p(m)}\}=\{0,\Delta_{m},2\Delta_{m},3\Delta_{m},\cdots,1\},
\]
where the second equality is equality for sets, and

\[
\overline{Q}_{m}\equiv\{u_{0},u_{1},\cdots,u_{p(2m)}\}\equiv\{0,2^{-m},2\cdot2^{-m},\cdots,2^{m}\}\subset[0,2^{m}],
\]
and 
\[
\overline{Q}_{\infty}\equiv\bigcup_{m=0}^{\infty}\overline{Q}_{m}\equiv\{u_{0},u_{1},\cdots\},
\]
where the second equality is equality for sets.

Unless otherwise specified, we will let $Q$ denote one of the three
parameter sets $\{0,1,\cdots\}$, $\overline{Q}_{\infty}$, or $[0,\infty)$.
\end{defn}
$\square$

To ease the burden on notations, for real-valued expressions $a,b,c$,
we will write the expressions $a=b\pm c$ interchangeably with $|a-b|\leq c$. 

\section{Filtrations, Stopping Times, and Markov Processes}

Filtrations, stopping times, and related objects were introduced in
Sections \ref{Section :Filtration} and \ref{Sec. Stopping time}.
Hitherto we have used only simple stopping times, which have values
in a finite and discrete parameter set. In this section, we will prove
the basic properties of stopping times with values in $[0,\infty)$
relative to some right continuous filtration. Then we will define
Markov processes and strong Markov processes.
\begin{prop}
\label{Prop. Stopping times rlative to right continuous filtration}\textbf{\emph{
(Basic properties of stopping times relative to a right continuous
filtration).}} Let $\mathcal{L}\equiv\{L^{(t)}:t\in[0,\infty)\}$
be an arbitrary right continuous filtration on some probability space
$(\Omega,L,E)$. All stopping times in the following will be relative
to this filtration and will have values in $[0,\infty)$. 

Let $\tau',\tau''$ be arbitrary stopping times . Then the following
holds.

\emph{1.} \emph{(Approximating stopping time by stopping times which
have regularly spaced dyadic values).} Let $\tau$ be an arbitrary
stopping time. Then for each regular point $t$ of the r.r.v. $\tau$,
we have $(\tau<t),(\tau=t)\in L^{(t)}$. Moreover, there exists a
nonincreasing sequence $(\eta_{h})_{h=0,1,\cdots}$ of stopping times,
such that, for each $h\geq0$, the stopping time $\eta_{h}$ has values
in 
\[
A_{h}\equiv\{s_{0},s_{1},s_{2},\cdots\}\equiv\{0,\Delta_{h},2\Delta_{h},\cdots\}
\]
 where $\Delta_{h}\equiv2^{-h}$ and $s_{j}\equiv j\Delta_{h}$ for
each $j\geq0$, and such that
\begin{equation}
\tau\leq\eta_{h}<\tau+2^{-h+2},\label{eq:temp-421-1-2-1}
\end{equation}
for each $h\geq0$. Consequently, $\eta_{h}\rightarrow\tau$ \emph{a.u.}
and $\eta_{h}\rightarrow\tau$ in probability.

\emph{2.} \emph{(Construction of stopping time as right-limit of given
stopping times).} Conversely, suppose $(\eta_{h})_{h=0,1,\cdots}$
is a sequence of stopping times such that, for some r.r.v. $\tau$,
we have \emph{(i)} $\tau\leq\eta_{h}$ for each $h\geq0$, and \emph{(ii)}
$\eta_{h}\rightarrow\tau$ in probability. Then $\tau$ is a stopping
time,

\emph{3.} $\tau'\wedge\tau''$, $\tau'\vee\tau'',$ and $\tau'+\tau''$
are stopping times.

\emph{4}. If $\tau'\leq\tau''$ then $L^{(\tau')}\subset L^{(\tau''')}$.

\emph{5}. Suppose $\tau$ is a stopping time. Define $L^{(\tau+)}\equiv\bigcap_{s>0}L^{(\tau+s)}$.
Then $L^{(\tau+)}=L^{(\tau)}$.
\end{prop}
\begin{proof}
1. Suppose $\tau$ is a stopping time. Let $t$ be an arbitrary regular
point of the r.r.v. $\tau$. Then, according to Definition \ref{Def. Regular=000026 Continuity Pts r.r.v.}
there exists an increasing sequence $(r_{k})_{k=1,2,\cdots}$of regular
points of $\tau$ such that $r_{k}\uparrow t$ and such that $P(\tau\leq r_{k})\uparrow P(\tau<t)$.
Consequently, $E|1_{\tau\leq r(k)}-1_{\tau<t}|\rightarrow0$. Since
$1_{\tau\leq r(k)}\in L^{(r(k))}\subset L^{(t)}$ for each $k\geq1$,
we conclude that $1_{\tau<t}\in L^{(t)}.$ In other words, $(\tau<t)\in L^{(t)}$.
Therefore $(\tau=t)=(\tau\leq t)(\tau<t)^{c}\in L^{(t)}$.

Separately, let $h\geq0$ be arbitrary. For convenience, define $r_{0}\equiv-\Delta_{h}$.
For each $j\geq1$, take a regular point 
\[
r_{j}\in(s_{j-1},s_{j})\equiv((j-1)\Delta_{h},j\Delta_{h})
\]
of the r.r.v. $\tau$. Then the r.r.v.
\[
\zeta_{h}\equiv\sum_{j=1}^{\infty}s_{j}1_{(r(j-1)<\tau\leq r(j))}
\]
has values in $A_{h}$ Moreover, for each $j\geq1$, we have, on the
$\mathrm{measurable}$ set $(r_{j-1}<\tau\leq r_{j})$, the inequality
\[
\tau\leq r_{j}<s_{j}=\zeta_{h}<r_{j+1}<r_{j+1}+(\tau-r_{j-1})
\]
\begin{equation}
<s_{j+1}+(\tau-s_{j-2})=\tau+3\Delta_{h}<\tau+2^{-h+2}.\label{eq:temp-210}
\end{equation}
Since $\bigcup_{j=1}^{\infty}(r_{j-1}<\tau\leq r_{j})$ is a full
set because $\tau$ is a nonnegative r.r.v. and because $r_{j}\uparrow\infty$
as $j\rightarrow\infty$, inequality \ref{eq:temp-210} implies that
\begin{equation}
\tau<\zeta_{h}<\tau+2^{-h+2}\label{eq:temp-421-2}
\end{equation}
as r.r.v.'s. Consider each possible value $s_{j}$ of the r.r.v. $\zeta_{h}$,
for some $j\geq1$. Then 
\[
(\zeta_{h}=s_{j})=(r_{j-1}<\tau\leq r_{j})=(\tau\leq r_{j})(\tau\leq r_{j-1})^{c}\in L^{(r(j))}\subset L^{(s(j))},
\]
because $\tau$ is a stopping time relative to $\mathcal{L}$ and
because $r_{j-1}$ and $r_{j}$ are regular points of $\tau$. Hence
\[
(\zeta_{h}\leq s_{j})=\bigcup_{k=1}^{j}(\zeta_{h}=s_{k})\in L^{(s(j))}.
\]
Thus $\zeta_{h}$ is a stopping time, with values in $A_{h}$. Therefore
the sequence $(\zeta_{h})_{h=0,1,\cdots}$ would be the desired result,
except for the lack of monotonicity. 

To fix that, let $h\geq0$ be arbitrary, and define the r.r.v. 
\begin{equation}
\eta_{h}\equiv\zeta_{0}\wedge\zeta_{1}\wedge\cdots\wedge\zeta_{h}\label{eq:temp-450-1}
\end{equation}
with values in $A_{0}\cup A_{1}\cup\cdots\cup A_{h}=A_{h}.$ Then,
for each possible value $t\in A_{h}$ of the r.r.v. $\eta_{h}$, we
have
\[
(\eta_{h}\leq t)=\bigcup_{i=0}^{h}(\zeta_{i}\leq t)\in L^{(t)}.
\]
Hence $\eta_{h}$ is a stopping time with values in $A_{h}$. Moreover,
equality \ref{eq:temp-450-1} implies that $(\eta_{h})_{h=0,1,\cdots}$
is a nonincreasing sequence, while inequality \ref{eq:temp-421-2}
implies that 
\begin{equation}
\tau<\eta_{h}<\tau+2^{-h+2}.\label{eq:temp-421-1-2}
\end{equation}
Assertion 1 is proved.

2. Conversely, suppose $(\eta_{h})_{h=0,1,\cdots}$ is a sequence
of stopping times such that (i) $\tau\leq\eta_{h}$ for each $h\geq0$,
and (ii) $\eta_{h}\rightarrow\tau$ in probability. Let $t\in R$
be a regular point of the r.r.v. $\tau$. Consider each $r>t$. Let
$\varepsilon>0$ be arbitrary. Let $s\in(t,r)$ be an arbitrary regular
point of all the r.r.v.'s in the countable set $\{\tau,\eta_{0},\eta_{1},\cdots\}$.
Then, by Conditions (i) and (ii), there exists $k\geq0$ such that
\[
P(s<\eta_{h})(\tau\leq t)\leq P(s-t<\eta_{h}-\tau)\leq\varepsilon
\]
for each $h\geq k$. Hence, for each $h\geq k$, we have
\[
E|1_{\tau\leq t}-1_{\eta(h)\leq s}|\leq E1_{\tau\leq t<s<\eta(h)}\leq\varepsilon.
\]
Thus $E|1_{\tau\leq t}-1_{\eta(h)\leq s}|\rightarrow0$ as $h\rightarrow\infty$.
Because $\eta_{h}$ is a stopping time we have $1_{\eta(h)\leq s}\in L^{(s)}\subset L^{(r)}.$
Hence $1_{\tau\leq t}\in L^{(r)},$ where $r>t$ is arbitrary. Thus
\[
1_{\tau\leq t}\in\bigcap_{r>t}L^{(r)}\equiv L^{(t+)}=L^{(t)},
\]
where the last equality is thanks to the assumption that the filtration
$\mathcal{L}$ is right continuous. We have verified that $\tau$
is a stopping time relative to the filtration $\mathcal{L}$ . Assertion
2 is proved.

3. By hypothesis, $\tau',\tau'$ are stopping times. Consider each
regular point $t\in R$ of the r.r.v. $\tau'\wedge\tau''.$ Let $r>t$
be arbitrary. Take any common regular point $s>t$ of the three r.r.v.'s
$\tau',\tau'$,$\tau'\wedge\tau''.$ Then $(\tau'\leq s)\cup(\tau'>s)$
and $(\tau''\leq s)\cup(\tau''>s)$ are full sets. Hence
\[
(\tau'\wedge\tau''\leq s)
\]
\[
=(\tau'\wedge\tau''\leq s)\cap((\tau'\leq s)\cup(\tau'>s))\cap((\tau''\leq s)\cup(\tau''>s))
\]
\[
=(\tau'\wedge\tau''\leq s)\cap(\tau'\leq s)\cap(\tau''\leq s)
\]
\begin{equation}
=(\tau'\leq s)\cap(\tau''\leq s)\in L^{(s)}\subset L^{(r)}.\label{eq:temp-220}
\end{equation}
Now let $s\downarrow t$. Then, since $t$ is a regular point of the
r.r.v. $\tau'\wedge\tau''$, we have $E1_{(\tau'\wedge\tau''\leq s)}\downarrow E1_{(\tau'\wedge\tau''\leq t)}$.
The Monotone Convergence Theorem therefore implies that $E|1_{(\tau'\wedge\tau''\leq s)}-1_{(\tau'\wedge\tau''\leq t)}|\rightarrow0$.
Consequently, since $1_{(\tau'\wedge\tau''\leq s)}\in L^{(r)}$ for
each common regular point $s>t$ of the three r.r.v.'s $\tau',\tau'$,$\tau'\wedge\tau''$,
according to relation \ref{eq:temp-220}, it follows that $1_{(\tau'\wedge\tau''\leq t)}\in L^{(r)}$,
where $r>t$ is arbitrary. Hence
\[
1_{(\tau'\wedge\tau''\leq t)}\in\bigcap_{r\in(t,\infty)}L^{(r)}\equiv L^{(t+)}=L^{(t)},
\]
where the last equality is thanks to the right continuity of the filtration
$\mathcal{L}$. Thus $\tau'\wedge\tau''$ is a stopping time relative
to $\mathcal{L}$. Similarly we can prove that $\tau'\vee\tau'',$
and $\tau'+\tau''$ are stopping times relative to $\mathcal{L}$.
Assertion 3 is verified.

4. Suppose the stopping times $\tau',\tau''$ are such that $\tau'\leq\tau''$.
Let $Y\in L^{(\tau')}$ be arbitrary. Consider each regular point
$t''$ of the stopping time $\tau''$ . Take an arbitrary regular
point $t'$ of $\tau'$ such that $t'\neq t''$. Then 
\[
Y1_{(\tau''\leq t'')}=Y1_{(t'<\tau'\leq\tau''\leq t'')}+Y1_{(\tau'\leq t'<\tau''\leq t'')}+Y1_{(\tau'\leq\tau''\leq t'<t'')}+Y1_{(\tau'\leq\tau''\leq t''<t')}
\]
\[
=Y1_{(t'<\tau')}1_{(\tau''\leq t'')}+Y1_{(\tau'\leq t')}1_{(t'<\tau''\leq t'')}+Y1_{(\tau''\leq t')}1_{(t'<t'')}+Y1_{(\tau''\leq t'')}1_{(t''<t')}.
\]
Consider the first summand in the last sum. Since $Y\in L^{(\tau')}$
by assumption, we have $Y1_{(t'<\tau')}\in L^{(t')}\subset L^{(t'')}$.
At the same time $1_{(\tau''\leq t'')}\in L^{(t'')}$. Hence $Y1_{(t'<\tau')}1_{(\tau''\leq t'')}\in L^{(t'')}$,
Similarly all the other summands in the last sum are members of $L^{(t'')}$.
Consequently the sum $Y1_{(\tau''\leq t'')}$ is a member of $L^{(t'')}$,
where $t''$ is an arbitrary regular point of the stopping time $\tau''$.
In other words, $Y\in L^{(\tau'')}$. Since $Y\in L^{(\tau')}$ is
arbitrary, we conclude that $L^{(\tau')}\subset L^{(\tau'')}$, as
alleged in Assertion 4.

5. Suppose $\tau$ is a stopping time relative to the right continuous
filtration $\mathcal{L}$. Let $Y\in L^{(\tau+)}\equiv\bigcap_{s>0}L^{(\tau+s)}$
be arbitrary. Let $t$ be an arbitrary regular point of the stopping
time $\tau$. Consider each $s>t$. Then $Y\in L^{(\tau+s)}$. Hence
$Y1_{(\tau\leq t)}\in L^{(t+s)}$. Consequently,
\[
Y1_{(\tau\leq t)}\in\bigcap_{s>0}L^{(t+s)}\equiv L^{(t+)}=L^{(t)},
\]
where the last equality is because the filtration $\mathcal{L}$ is
right continuous. Thus $Y\in L^{(\tau)}$ for each $Y\in L^{(\tau+)}$.
Equivalently, $L^{(\tau+)}\subset L^{(\tau)}$. In the other direction,
we have, trivially, $L^{(\tau)}\subset L^{(\tau+)}$. Summing up $L^{(\tau)}=L^{(\tau+)}$
as alleged in Assertion 5.
\end{proof}
\begin{defn}
\label{Def. Markov process} \textbf{(Markov process). }Let $(S,d)$
be an abritrary locally compact metric space. Let $Q$ denote one
of the three parameter sets $\{0,1,\cdots\}$, $\overline{Q}_{\infty}$,
or $[0,\infty)$. Let $\mathcal{L}\equiv\{L^{(t)}:t\in Q\}$ denote
an arbitrary right continuous filtration of the sample space $(\Omega,L,E)$.
Let 
\[
X:Q\times(\Omega,L,E)\rightarrow(S,d)
\]
be an arbitrary process which is adapted to the filtration $\mathcal{L}$.

Suppose, for each $t\in Q$, for each nondecreasing sequence $t_{0}\equiv0\leq t_{1}\leq\cdots\leq t_{m}$
in $Q$, and for each function $f\in C(S^{m+1},d^{m+1})$, we have
(i) the conditional expectation of $E(f(X_{t+t(0)},X_{t+t(1)},\cdots,X_{t+t(m)})|L^{(t)})$
exists, (ii) the conditional expectation of $E(f(X_{t+t(0)},X_{t+t(1)},\cdots,X_{t+t(m)})|X_{t})$
exists, and (iii)

\begin{equation}
E(f(X_{t+t(0)},X_{t+t(1)},\cdots,X_{t+t(m)})|L^{(t)})=E(f(X_{t+t(0)},X_{t+t(1)},\cdots,X_{t+t(m)})|X_{t}).\label{eq:temp-459}
\end{equation}
Then the process $X$ is called a \emph{Markov process}\index{Markov process}
relative to the filtration $\mathcal{L}$. We will refer to Conditions
(i-iii) as the \index{Markov property}\emph{Markov property}. In
the special case where $\mathcal{L}$ is the natural filtration of
the process $X$, we will omit the reference to $\mathcal{L}$ and
simply say that $X$ is a Markov process. 
\end{defn}
$\square$
\begin{defn}
\textbf{\label{Def. Strong Markov process}(Strong Markov process).
}Let $(S,d)$ be an abritrary locally compact metric space. Let $Q$
denote one of the three parameter sets $\{0,1,\cdots\}$, $\overline{Q}_{\infty}$,
or $[0,\infty)$. Let $\mathcal{L}\equiv\{L^{(t)}:t\in Q\}$ denote
an arbitrary right continuous filtration of the sample space $(\Omega,L,E)$.
Let 
\[
X:Q\times(\Omega,L,E)\rightarrow(S,d)
\]
be an arbitrary process which is adapted to the filtration $\mathcal{L}$.
Suppose, for each stopping time $\tau$ with values in $Q$ relative
to the filtration\emph{ $\mathcal{L}$, }the following two conditions
hold.

1. The function $X_{\tau}$ is a well defined r.v. relative to $L^{(\tau)}$,
with values in $(S,d)$. 

2. For each nondecreasing sequence $t_{0}\equiv0\leq t_{1}\leq\cdots\leq t_{m}$
in $Q$, and for each function $f\in C(S^{m+1},d^{m+1})$, we have
(i) the conditional expectation of $E(f(X_{\tau+t(0)},X_{\tau+t(1)},\cdots,X_{\tau+t(m)})|L^{(\tau)})$
exists, (ii) the conditional expectation of $(f(X_{\tau+t(0)},X_{\tau+t(1)},\cdots,X_{\tau+t(m)})|X_{\tau})$
exists, and (iii)

\begin{equation}
E(f(X_{\tau+t(0)},X_{\tau+t(1)},\cdots,X_{\tau+t(m)})|L^{(\tau)})=(f(X_{\tau+t(0)},X_{\tau+t(1)},\cdots,X_{\tau+t(m)})|X_{\tau}).\label{eq:temp-459-1}
\end{equation}

Then the process $X$ is called a \emph{strong Markov process}\index{Strong Markov process}
relative to the filtration $\mathcal{L}$. Since each constant time
$t\in Q$ is a stopping time, each strong Markov process is clearly
a Markov process. 
\end{defn}
$\square$

\section{Transition Distributions}
\begin{defn}
\label{Def.  Transition Distribution} \textbf{(Transition distribution).}
Let $(S_{0},d_{0})$ and $(S_{1},d_{1})$ be compact metric spaces,
with $d_{0}\leq1$ and $d_{1}\leq1$, and with fixed reference points
$x_{0,\circ}$ and $x_{1,\circ}$ respectively. Let
\[
T:C(S_{1},d_{1})\rightarrow C(S_{0},d_{0})
\]
be an arbitrary nonnegative linear function. Write 
\[
T^{x}\equiv T(\cdot)(x):C(S_{1},d_{1})\rightarrow R
\]
for each $x\in S_{0}$. Suppose, (i) for each $x\in S_{0}$, the function
$T^{x}$ is a distribution on $(S_{1},d_{1})$, in the sense of Definition
\ref{Def. distributions on complete metric space}, and (ii) for each
$f\in C(S_{1},d_{1})$ with a modulus of continuity $\delta_{f}$
, the function $Tf\in C(S_{0},d_{0})$ has a modulus of continuity
$\alpha(\delta_{f}):(0,\infty)\rightarrow(0,\infty)$ which depends
only on $\delta_{f}$, and otherwise not on the function $f$.

Then the function $T$ is called a \index{transition distribution}\emph{
transition distribution} from $(S_{0},d_{0})$ to $(S_{1},d_{1})$.
The operation $\alpha$ is then called a \emph{modulus of smoothness
of the transition distribution}\index{modulus of smoothness oa transition distribution}
$T$.
\end{defn}
$\square$
\begin{lem}
\label{Lem. Composite transition distributions} \textbf{\emph{(Composite
transition distributions).}} For each $j=0,1,2$, let $(S_{j},d_{j})$
be a compact metric space with $d_{i}\leq1$. For each $j=0,1$, let
$T_{j,j+1}$ be a transition distribution from $(S_{j},d_{j})$ to
$(S_{j+1},d_{j+1})$, with modulus of smoothness $\alpha_{j,j+1}$.
Then the composite function 
\[
T_{0,2}\equiv T_{0,1}T_{1,2}\equiv T_{0,1}\circ T_{1,2}:C(S_{2},d_{2})\rightarrow C(S_{0},d_{0})
\]
is a transition distribution from $(S_{0},d_{0})$ to $(S_{2},d_{2})$,
with a modulus of smoothness $\alpha_{0,2}$ defined by 
\[
\alpha_{0,2}(\delta_{f})\equiv\alpha_{0,1}(\alpha_{1,2}(\delta_{f})):(0,\infty)\rightarrow(0,\infty)
\]
 for each $modulus$ of continuity $\delta_{f}$.. 

We will call $T_{0,2}\equiv T_{0,1}T_{1,2}$ the \emph{composite transition
distribution\index{composite transition distribution}} of $T_{0,1}$
and $T_{1,2}$, and call $\alpha_{0,2}\equiv\alpha_{0,1}\circ\alpha_{1,2}$
the \emph{\index{composite modulus of smoothness}composite modulus
of smoothness} of $T_{0,1}T_{1,2}$.
\end{lem}
\begin{proof}
1. Being the composite of two linear functions, the function $T_{0,2}$
is linear. Let $f\in C(S_{2},d_{2})$ be arbitrary, with a modulus
of continuity $\delta_{f}$ and with $|f|\leq1$. Then $T_{1,2}f\in C(S_{1},d_{1})$
has a modulus of continuity $\alpha_{1,2}(\delta_{f})$. Hence the
function $T_{0,1}(T_{1,2}f)$ has a modulus of continuity $\alpha_{0,1}(\alpha_{1,2}((\delta_{f}))$. 

2. Consider each $x\in S_{0}$. Suppose $f\in C(S_{2},d_{2})$ is
such that $T_{0,2}^{x}f>0$. Then $T_{0,1}^{x}T_{1,2}f=T_{0,2}^{x}f>0$.
Therefore, since $T_{0,1}^{x}$ is a distribution, there exists $y\in S_{1}$
such that $T_{1,2}^{y}f\equiv T_{1,2}f(y)>0.$ Since $T_{1,2}^{y}$
is a distribution, there exists, in turn, some $z\in S_{2}$ such
that $f(z)>0$. Thus $T_{0,2}^{x}$ is an integration on $(S_{2},d_{2})$
in the sense of Definition \ref{Def. integration on loc compact space}.
Since $1\in C(S_{2},d_{2})$ and $d_{2}\leq1$, it follows that $T_{0,2}^{x}$
is a distribution on $(S_{2},d_{2})$ in the sense of Definition \ref{Def. distributions on complete metric space},
where $x\in S_{0}$ is arbitrary. The conditions in Definition \ref{Def.  Transition Distribution}
have been verified for $T_{0,2}$ to be a transition distribution. 
\end{proof}
By Definition \ref{Def.  Transition Distribution}, the domain and
range of a transition distribution are spaces of continuous functions.
We next extend both to spaces of integrable functions. Recall that,
for each $x\in S$, we let $\delta_{x}$ denote the distribution concentrated
at $x.$ 
\begin{prop}
\label{Prop. Completion of a transition distribution relativve to initial distribution}
\textbf{\emph{(Complete extension of a transition distribution relative
to an initial distribution).}} Let $(S_{0},d_{0})$ and $(S_{1},d_{1})$
be compact metric spaces, with $d_{0}\leq1$ and $d_{1}\leq1$. Let
$T$ be a transition distribution from $(S_{0},d_{0})$ to $(S_{1},d_{1})$.
Let $E_{0}$ be a distribution on $(S_{0},d_{0})$. Define the composite
function 
\begin{equation}
E_{1}\equiv E_{0}T:C(S_{1})\rightarrow R\label{eq:temp-428}
\end{equation}
Then the the following holds.

\emph{1. }$E_{1}\equiv E_{0}T$ is a distribution on $(S_{1},d_{1})$. 

\emph{2.} For each $i=0,1,$ let $(S_{i},L_{E(i)},E_{i})$ be the
complete extension of $(S_{i},C(S_{i}),E_{i})$. Let $f\in L_{E(1)}$
be arbitrary. Define the function $Tf$ on $S_{0}$ by 
\begin{equation}
domain(Tf)\equiv\{x\in S_{0}:f\in L_{\delta(x)T}\}\label{eq:temp-533}
\end{equation}
and by 
\begin{equation}
(Tf)(x)\equiv(\delta_{x}T)f\label{eq:temp-534}
\end{equation}
 for each $x\in domain(Tf)$. Then \emph{(i)} $Tf\in L_{E(0)}$, and
\emph{(ii)} $E_{0}(Tf)=E_{1}f\equiv(E_{0}T)f$. 

\emph{3. }Moreover, the extended function 
\[
T:L_{E(1)}\rightarrow L_{E(0)},
\]
thus defined, is a contraction mapping, hence continuous, relative
to the norm $E_{1}|\cdot|$ on $L_{E(1)}$ and the norm $E_{0}|\cdot|$
on $L_{E(0)}$. 
\end{prop}
\begin{proof}
1. By the defining equality \ref{Lem. Composite transition distributions},
the function $E_{1}$ is clearly linear and nonegative. Suppose $E_{1}f\equiv E_{0}Tf>0$
for some $C(S_{1},d_{1})$. Then, since $E_{0}$ is a disribution,
there exists $x\in S_{0}$ such that $T^{x}f>0$. In turn, since $T^{x}$
is a distribution, there exists $y\in S_{1}$ such that $f(y)>0$.
Thus $E_{1}$is an integration. Since $1\in C(S_{1},d_{1})$ and $d_{1}\leq1$,
it follows that $E_{1}$ is a distribution.

2. For each $i=0,1,$ let $(S_{i},L_{E(i)},E_{i})$ be the complete
extension of $(S_{i},C(S_{i}),E_{i})$. Let $f\in L_{E(1)}$ be arbitrary.
Define the function $Tf$ on $S_{0}$ by equalities \ref{eq:temp-533}
and \ref{eq:temp-534}. Then, by Definition \ref{Def. Integrablee functions and Completion of integration space}
of complete extensions, there exists a sequence ($f_{n})_{n=1,2,\cdots}$
in $C(S_{1})$ such that (i') $\sum_{n=1}^{\infty}E_{0}T|f_{n}|=\sum_{n=1}^{\infty}E_{1}|f_{n}|<\infty$,
(ii') 
\[
\{x\in S_{1}:\sum_{n=1}^{\infty}|f_{n}(x)|<\infty\}\subset domain(f),
\]
and (iii') $f(x)=\sum_{n=1}^{\infty}f_{n}(x)$ for each $x\in S_{1}$
with $\sum_{n=1}^{\infty}|f_{n}(x)|<\infty$. Condition (i') implies
that the subset
\[
D_{0}\equiv\{x\in S_{0}:\sum_{n=1}^{\infty}T|f_{n}|(x)<\infty\}\equiv\{x\in S_{0}:\sum_{n=1}^{\infty}T^{x}|f_{n}|<\infty\}
\]
of the probability space $(S_{0},L_{E(0)},E_{0})$ is a full subset.
It implies also that the function $g\equiv\sum_{n=1}^{\infty}Tf_{n}$,
with $domain(g)\equiv D_{0}$, is a member of $L_{E(0)}$. Now consider
an arbitrary $x\in D_{0}$. Then 
\begin{equation}
\sum_{n=1}^{\infty}|T^{x}f_{n}|\leq\sum_{n=1}^{\infty}T^{x}|f_{n}|<\infty.\label{eq:temp-528}
\end{equation}
Together with Condition (iii'), this implies that $f\in L_{\delta(x)T}$,
with $(\delta_{x}T)f=\sum_{n=1}^{\infty}(\delta_{x}T)f_{n}$. Hence,
according to the defining equalities \ref{eq:temp-533} and \ref{eq:temp-534},
we have $x\in domain(Tf)$, with 
\[
(Tf)(x)\equiv(\delta_{x}T)f=(\sum_{n=1}^{\infty}Tf_{n})(x)=g(x).
\]
Thus $Tf=g$ on the full subset $D_{0}$ of $(S_{0},L_{E(0)},E_{0})$.
Since $g\in L_{E(0)}$, it follows that $Tf\in L_{E(0)}$. The desired
Condition (i) is verified. Moreover, 
\[
E_{0}(Tf)=E_{0}g\equiv E_{0}\sum_{n=1}^{\infty}Tf_{n}=\sum_{n=1}^{\infty}E_{0}Tf_{n}=E_{0}T\sum_{n=1}^{\infty}f_{n}=E_{1}f,
\]
where the third and fourth equality are both justified by Condition
(i'). The desired Condtion (ii) is also verified. Thus Assertion 2
is proved.

3. Let $f\in L_{E(1)}$ be arbitrary. Then, in the notations of the
previous step, we have
\[
E_{0}|Tf|=E_{0}|g|=\lim_{N\rightarrow\infty}E_{0}|\sum_{n=1}^{N}Tf_{n}|=\lim_{N\rightarrow\infty}E_{0}|T\sum_{n=1}^{N}f_{n}|
\]
\[
\leq\lim_{N\rightarrow\infty}E_{0}T|\sum_{n=1}^{N}f_{n}|=\lim_{N\rightarrow\infty}E_{1}|\sum_{n=1}^{N}f_{n}|=E_{1}|f|.
\]
In short, $E_{0}|Tf|\leq E_{1}|f|$. Thus the mapping $T:L_{E(1)}\rightarrow L_{E(0)}$
is a contraction, as alleged in Assertion 3.
\end{proof}
\begin{defn}
\label{Def. Assumptions, conventions, and notations re traansition fjd's}\textbf{
(Convention regarding automatic completion of transition distributions).}
We hereby make the convention that, given each transition distribution
$T$ from a compact metric space $(S_{0},d_{0})$ to a compact metric
space $(S_{1},d_{1})$ with $d_{0}\leq1$ and $d_{1}\leq1$, and given
each initial distribution $E_{0}$ on $(S_{0},d_{0})$, the transition
distribution 
\[
T:C(S_{1},d_{1})\rightarrow C(S_{0},d_{0})
\]
 is automaticaly completely extended to the nonnegative linear function
\[
T:L_{E(1)}\rightarrow L_{E(0)}
\]
in the manner of Proposition \ref{Prop. Completion of a transition distribution relativve to initial distribution},
where $E_{1}\equiv E_{0}T$, and where $(S_{i},L_{E(i)},E_{i})$ is
the complete extension of $(S_{i},C(S_{i}),E_{i}),$ for each $i=0,1$, 

Thus $Tf$ is integrable relative to $E_{0}$ for each integrable
function $f$ relative to $E_{0}T$, with $(E_{0}T)f=E_{0}(Tf)$.
In the special case where $E_{0}$ is the point mass distribution
$\delta_{x}$ concentrated at some $x\in S_{0}$, we have $x\in domain(Tf)$
and $T^{x}f=(Tf)^{x}$, for each integrable function $f$ relative
to $T^{x}$. 
\end{defn}
$\square$
\begin{lem}
\label{Lem.One-step transition distributions} \textbf{\emph{(One-step
transition distributions). }}Let $(S,d)$ be a compact metric space
with $d\leq1$. Let $T$ be a transition distribution from $(S,d)$
to $(S,d),$ with a modulus of smoothness $\alpha_{T}$. Define $^{1}T\equiv$T
. Let $m\geq2$ and $f\in C(S^{m},d^{m})$ be arbitrary. Define a
function $(^{m-1}T)f$ on $S^{m-1}$ by 
\begin{equation}
((^{m-1}T)f)(x_{1},\cdots,x_{m-1})\equiv\int T^{x(m-1)}(dx_{m})f(x_{1},\cdots,x_{m-1},x_{m})\label{eq:temp-426}
\end{equation}
for each $x\equiv(x_{1},\cdots,x_{m-1})\in S^{m-1}$. Then the following
holds for each $m\geq1$.

\emph{1. }If $f\in C(S^{m},d^{m})$ has values in $[0,1]$ and has
modulus of continuity $\delta_{f}$, then the function $(^{m-1}T)f$
is a member of $C(S^{m-1},d^{m-1})$, with values in $[0,1]$, and
has the modulus of continuity 
\[
\widetilde{\alpha}_{\alpha(T),\xi}(\delta_{f}):(0,\infty)\rightarrow(0,\infty)
\]
defined by 
\begin{equation}
\widetilde{\alpha}_{\alpha(T),\xi}(\delta_{f})(\varepsilon)\equiv\alpha_{T}(\delta_{f})(2^{-1}\varepsilon)\wedge\delta_{f}(2^{-1}\varepsilon)\label{eq:temp-562}
\end{equation}
for each $\varepsilon>0$.

\noun{2. }For each $x\equiv(x_{1},\cdots,x_{m-1})\in S^{m-1}$, the
function 
\[
^{m-1}T:C(S^{m},d^{m})\rightarrow C(S^{m-1},d^{m-1})
\]
is a transition distribution from $(S^{m-1},d^{m-1})$ to $(S^{m},d^{m}),$
with modulus of smoothness $\widetilde{\alpha}_{\alpha(T),\xi}$.

We will call $^{m-1}T$ the \emph{\index{one-step transition distribution}
one-step transition distribution at step $m$ according to} $T$. 
\end{lem}
\begin{proof}
1. Let $f\in C(S^{m+1},d^{m+1})$ be arbitrary, with values in $[0,1]$
and with a modulus of continuity $\delta_{f}$. Since $T$ is a transition
distribution, $T^{x(m-1)}$ is a distibution on $(S,d)$ for each
$x_{m-1}\in S$. Hence the integration on the right-hand side of equality
\ref{eq:temp-426} makes sense and has values in $[0,1]$. Therefore
the left-hand side is well defined and has values in $[0,1]$. We
need to prove that the function $(^{m-1}T)f$ is a continuous function.

2. To that end, let $\varepsilon>0$ be arbitrary. Define

Let $(x_{1},\cdots,x_{m-1}),(x'_{1},\cdots,x'_{m-1})\in S^{m-1}$
be arbitrary such that
\[
d^{m-1}((x_{1},\cdots,x_{m-1}),(x'_{1},\cdots,x'_{m-1}))<\widetilde{\alpha}_{\alpha(T),\xi}(\delta_{f})(\varepsilon)
\]
For abbreviation, write $x\equiv(x_{1},\cdots,x_{m-2})\in S^{m-2}$,
where the sequence is by convention empty if $m=2$. Write $y\equiv x_{m-1}$.
Define the abbreviations $x',y'$ similarly. 

With $x,y$ fixed, the function $f(x,y,\cdot)$ also has a modulus
of continuity $\delta_{f}$. Hence the function $Tf(x,y,\cdot)$ has
a moudlus of continuuity $\alpha_{T}(\delta_{f})$, by Definition
\ref{Def.  Transition Distribution}. Therefore, since 
\[
d(y,y')\leq d^{m-1}((x_{1},\cdots,x_{m-1}),(x'_{1},\cdots,x'_{m-1}))<\widetilde{\alpha}_{\alpha(T),\xi}(\delta_{f})(\varepsilon)\leq\alpha_{T}(\delta_{f})(2^{-1}\varepsilon),
\]
it follows that 
\[
|(Tf(x,y,\cdot))(y)-(Tf(x,y,\cdot))(y')|<2^{-1}\varepsilon.
\]
In other words
\begin{equation}
\int T^{y}(dz)f(x,y,z)=\int T^{y'}(dz)f(x,y,z)\pm2^{-1}\varepsilon\label{eq:temp-561}
\end{equation}
,At the same time, for each $z\in S$, since 
\[
d^{m}((x,y,z),(x',y',z))=d^{m-1}((x,y),(x',y'))<\delta_{f}(2^{-1}\varepsilon)
\]
we have $|f(x,y,z)-f(x',y',z)|<2^{-1}\varepsilon$. Hence
\[
\int T^{y'}(dz)f(x,y,z)=\int T^{y'}(dz)(f(x',y',z)\pm2^{-1}\varepsilon)=\int T^{y'}(dz)f(x',y',z)\pm2^{-1}\varepsilon.
\]
Combining with equality \ref{eq:temp-561}, we obtain
\[
\int T^{y}(dz)f(x,y,z)=\int T^{y'}(dz)(f(x',y',z)\pm\varepsilon.
\]
In other words,
\[
((^{m-1}T)f)(x,y)=((^{m-1}T)f)(x',y')\pm\varepsilon.
\]
Thus $(^{m-1}T)f$ is continuous, with modulus of continuity $\widetilde{\alpha}_{\alpha(T),\xi}(\delta_{f})$.
Assertion 1 has been proved. 

2. By linearity, we see that $(^{m}T)f\in C(S^{m},d^{m})$ for each
$f\in C(S^{m+1},d^{m+1})$. Therefore the function 
\[
^{m-1}T:C(S^{m},d^{m})\rightarrow C(S^{m-1},d^{m-1})
\]
is a well-defined. It is clearly linear and nonnegative from the defining
formula \ref{eq:temp-426}. Consider each $x\equiv(x_{1},\cdots,x_{m-1})\in S^{m-1}$.
Suppose $(^{m-1}T)^{x}f\equiv\int T^{x(m-1)}(dy)f(x,y)>0$. Then,
since $T^{x(m-1)}$ is a distribution, there exists $y\in S$ such
that $f(x,y)>0.$ Hence $(^{m-1}T)^{x}$ is an integration on $(S^{m-1},d^{m-1})$
in the sense of Definition \ref{Def. integration on loc compact space}.
Since $d^{m-1}\leq1$ and $1\in C(S^{m-1},d^{m-1})$, the function
$(^{m-1}T)^{x}$ is a distribution on $(S^{m-1},d^{m-1})$ in the
sense of Definition \ref{Def. distributions on complete metric space}.
We have verified the conditions in Definition \ref{Def.  Transition Distribution}
for $^{m-1}T$ to be a transition distribution. Assertion 2 is proved.
\end{proof}

\section{Markov Semigroup}

Recall that $Q$ denotes one of the three parameter sets $\{0,1,\cdots\}$,
$\overline{Q}_{\infty}$, or $[0,\infty)$.
\begin{defn}
\label{Def. Markov semigroup  semigroup} \textbf{(Markov semigroup).
}Let $(S,d)$ be a compact metric space with $d\leq1$. Unless otherwise
spacified, The ssymbol $\left\Vert \cdot\right\Vert $ will stand
for the supremum norm for the space $C(S,d)$. Let $\mathbf{T}\equiv\{T_{t}:t\in Q\}$
be a family of transition distributions from $(S,d)$ to $(S,d)$,
such that $T_{0}$ is the identity mapping. Suppose the following
three conditions are satisfied.

1. (Smoothness). For each $N\geq1$, for each $t\in[0,N]Q$, the transition
distrbution $T_{t}$ has some modulus of smoothness $\alpha_{\mathbf{T,\mathit{N}}}$,
in the sense of Definition \ref{Def.  Transition Distribution}. Note
that the modulus of smoohness $\alpha_{\mathbf{T,\mathit{N}}}$ is
dependent on the finite interval $[0,N]$, but is otherwise independent
of $t$.

2. (Semigroup property). For each $s,t\in Q$, we have $T_{t+s}=T_{t}T_{s}$. 

3. (Strong continuity). For each $f\in C(S,d)$ with a modulus of
continuity $\delta_{f}$ and with $\left\Vert f\right\Vert \leq1$,
and for each $\varepsilon>0$, there exists $\delta_{\mathbf{T}}(\varepsilon,\delta_{f})>0$
so small that, for each $t\in[0,\delta_{\mathbf{T}}(\varepsilon,\delta_{f}))Q$,
we have 
\begin{equation}
\left\Vert f-T_{t}f\right\Vert \leq\varepsilon.\label{eq:temp-382-1-2-2}
\end{equation}
Note that this strong continuity condition is trivially satisfied
if $Q=\{0,1,\cdots\}$.

Then we call the family\emph{ }$\mathbf{T}$\emph{ }a \emph{Markov
semigroup}\index{Markov semigroup} of transition distributions with
state space $(S,d)$ and parameter space $Q$. For abbreviaation,
we will simply call $\mathbf{T}$ a \index{semigroup}\emph{semigroup}.
The operation $\delta_{\mathbf{T}}$ is called a \index{modulus of strong continuity of a Markov semigroup}\emph{
modulus of strong continuity} of\emph{ }$\mathbf{T}$. The sequence
$\alpha_{\mathbf{T}}\equiv(\alpha_{\mathbf{T,N}})_{N=1,2,\cdots}$
is called the \emph{modulus of smoothness  of the semgroup}\index{modulus of smoothness of Markov semgroup}
$\mathbf{T}$.$\square$
\end{defn}
\begin{rem}
\textbf{\label{Rem. No loss of genrality}(No loss of generally in
restricting state space to comact space).} Eventhough we have defined
transition distributions and Markov semigroups only for compact metric
spaces $(S,d)$ with $d\leq1$, there is no loss of generality because
a locally compact metric space $(S_{00},d_{00})$ can be embedded
into its one-point compactification $(S,d)\equiv(S_{00}\cup\{\Delta_{00}\},\overline{d}_{00})$
where $\Delta_{00}$ is the point at infinity, and where $\overline{d}_{00}\leq1$.
The assumption of compactness simplifies proofs.
\end{rem}
The next lemma strengthens the continuity of $T_{t\cdot}$ at $t=0$
to uniform continuity over $t\in Q$. Its proof is somewhat longer
than its one-lined classical counterpart.
\begin{lem}
\label{Lem. Uniform strong continuity on parameteer set}\textbf{\emph{(Uniform
strong continuity on the parameter set).}} Suppose $Q$ is one of
the three parameter sets $\{0,1,\cdots\}$, $\overline{Q}_{\infty}$,
or $[0,\infty)$. Let $\mathbf{T}$ be an arbitrary semigroup with
parameter set $Q$, and with a modulus of strong continuity $\delta_{\mathbf{T}}$.
Let $f\in C(S,d)$ be arbitrary, with a modulus of continuity $\delta_{f}$,
and with $|f|\leq1$. Let $\varepsilon>0$ and $r,s\in Q$ be arbitrary
with $|r-s|<\delta_{\mathbf{T}}(\varepsilon,\delta_{f})$. Then $\left\Vert T_{r}f-T_{s}f\right\Vert \leq\varepsilon.$
\end{lem}
\begin{proof}
1. The case where $Q=\{0,1,\cdots\}$ is trivial.

2. Suppose $Q=\overline{Q}_{\infty}$ or $Q=[0,\infty)$. Let $\varepsilon>0$
and $r,s\in Q$ be arbitrary with $0\leq s-r<\delta_{\mathbf{T}}(\varepsilon,\delta_{f})$.
Then, for each $x\in S$, we have
\[
|T_{s}^{x}f-T_{r}^{x}f|=|T_{r}^{x}(T_{s-r}f-f)|\leq\left\Vert T_{s-r}f-f\right\Vert \leq\varepsilon,
\]
where the equality is by the semigroup property, where the first inequality
is because $T_{r}^{x}$ is a distribution on $(S,d)$, and where the
last inequality is by the definition of $\delta_{\mathbf{T}}$ as
a modulus of strong continuity. Thus 
\begin{equation}
\left\Vert T_{r}f-T_{s}f\right\Vert \leq\varepsilon.\label{eq:temp-378}
\end{equation}

3. Let $\varepsilon>0$ and $r,s\in\overline{Q}_{\infty}$ be arbitrary
with $|s-r|<\delta_{\mathbf{T}}(\varepsilon,\delta_{f})$. Either
$0\leq s-r<\delta_{\mathbf{T}}(\varepsilon,\delta_{f})$, in which
case inequality\ref{eq:temp-378} holds according to Step 2, or $0\leq r-s<\delta_{\mathbf{T}}(\varepsilon,\delta_{f})$,
in which case inequality\ref{eq:temp-378} holds similarly. Thus the
lemma is proved if $Q=\overline{Q}_{\infty}$.

4. Now suppose $Q=[0,\infty)$. Let $\varepsilon>0$ and $r,s\in[0,\infty)$
be arbitrary with $|r-s|<\delta_{\mathbf{T}}(\varepsilon,\delta_{f})$.
Let $\varepsilon'>0$ be arbitrary. Let $t,v\in\overline{Q}_{\infty}$be
arbitrary such that (i)\textbf{ $r\leq t<r+\delta_{\mathbf{T}}(\varepsilon',\delta_{f})$},
(ii) $s\leq v<s+\delta_{\mathbf{T}}(\varepsilon',\delta_{f})$ and
(iii) $|t-v|<\delta_{\mathbf{T}}(\varepsilon,\delta_{f})$. Then,
according to inequality \ref{eq:temp-378} in Step 2, we have $\left\Vert T_{t}f-T_{r}f\right\Vert \leq\varepsilon'$,
$\left\Vert T_{v}f-T_{s}f\right\Vert \leq\varepsilon'$. According
to Step 3, we have $\left\Vert T_{t}f-T_{v}f\right\Vert \leq\varepsilon$.
Combining, we obtain
\[
\left\Vert T_{r}f-T_{s}f\right\Vert \leq\left\Vert T_{t}f-T_{r}f\right\Vert +\left\Vert T_{s}f-T_{v}f\right\Vert +\left\Vert T_{t}f-T_{v}f\right\Vert <\varepsilon'+\varepsilon'+\varepsilon.
\]
Letting $\varepsilon'\rightarrow0,$ we obtain $\left\Vert T_{r}f-T_{s}f\right\Vert \leq\varepsilon$.
Thus the lemma is also proved for the case where $Q=[0,\infty)$.
\end{proof}
$\square$

\section{Markov Transition f.j.d.'s }

In this section, we will define consistent family of f.j.d.'s generated
by an initial distribution and a semigroup. The parameter set $Q$
is assumed to be one of the three sets $\{0,1,\cdots\}$, $\overline{Q}_{\infty}$,
or $[0,\infty)$. We will refer informally to the first two as the
discrete parameter sets. The state space $(S,d)$ is assumed to be
compact with $d\leq1$. Let $\xi\equiv(A_{k})_{k=1,2,\cdots}$ be
a binary approximation of $(S,d)$ relative to $x_{\circ}$. Let 
\[
\pi\equiv(\{g_{k,x}:x\in A_{k}\})_{k=1,2,\cdots}
\]
be the partition of unity of $(S,d)$ determined by $\xi$, as in
Definition \ref{Def. Partition of unity for locally compact (S,d)}.
\begin{defn}
\label{Def. Family of f.j.d.s from E0 and Markov semi-group} \textbf{(Family
of transition f.j.d.'s generated by an initial distribution and a
Markov semigroup).} Let $Q$ be one of the three sets $\{0,1,\cdots\}$,
$\overline{Q}_{\infty}$, or $[0,\infty)$. Let $\mathbf{T}$ be an
arbitrary Markov semigroup, with the compact state space $(S,d)$
where $d\leq1$, and with parameter set $Q.$  Let $E_{0}$ be an
arbitrary distribution on $(S,d)$. For arbitrary $m\geq1$, $f\in C(S^{m},d^{m})$,
and nondecreasing sequence $r_{1}\leq\cdots\leq r_{m}$ in $Q$, define
\[
F_{r(1),\cdots,r(m)}^{E(0),\mathbf{T}}f
\]
\begin{equation}
=\int E_{0}(dx_{0})\int T_{r(1)}^{x(0)}(dx_{1})\int T_{r(2)-r(1)}^{x(1)}(dx_{2})\cdots\int T_{r(m)-r(m-1)}^{x(m-1)}(dx_{m})f(x_{1},\cdots,x_{m}).\label{eq:temp-384-1-1}
\end{equation}
In the special case where $E_{0}\equiv\delta_{x}$ is the distribution
which assigns probability $1$ to some point $x\in S$, we will simply
write 
\[
F_{r(1),\cdots,r(m)}^{*,\mathbf{T}}(x)\equiv F_{r(1),\cdots,r(m)}^{x,\mathbf{T}}\equiv F_{r(1),\cdots,r(m)}^{\delta(x),\mathbf{T}}.
\]
The next theorem willl prove that $F_{r(1),\cdots,r(m)}^{*,\mathbf{T}}:C(S^{m},d^{m})\rightarrow C(S,d)$
is then a well-defined transition distribution. 

An arbitrary consistent family 
\[
\{F_{r(1),\cdots,r(m)}f:m\geq0;r_{1},\cdots,r_{m}\;\mathrm{in}\;Q\}
\]
 of f.j.d.'s satisfying Condition \ref{eq:temp-384-1-1} is said to
be \emph{generated by the initial distribution $E_{0}$ and the semigroup\index{consistent family of f.j.d.'s from initial distribution and semigroup}}
\emph{$\mathbf{T}$.} 

A process 
\[
X:Q\times(\Omega,L,E)\rightarrow(S,d),
\]
whose marginal distributions are given by a consistent family generated
by the initial distribution $E_{0}$ and the semigroup \emph{$\mathbf{T}$,}
will be called \emph{a process }\index{process with Markov semigroup}
\emph{generated by the initial distribution $E_{0}$ and semigroup
$\mathbf{T}$. }We will see later that such processes are Markov processes.

$\square$
\end{defn}
\begin{thm}
\label{Thm. Construction of transition  f.j.d.s from E0 and Markov semi-group, discrete parameters-1}
\textbf{\emph{(Construction of family of transition f.j.d.'s from
initial distribution and semigroup, for discrete parameters ). }}Let
$(S,d)$ be a compact metric space with $d\leq1$. Let $Q$ be one
of the three sets $\{0,1,\cdots\}$, $\overline{Q}_{\infty}$, or
$[0,\infty)$. Let $\mathbf{T}$ be an arbitrary semigroup with state
space $(S,d)$ and with parameter set  $Q$, and with a modulus of
strong continuity $\delta_{\mathbf{T}}$ and with a modulus of smoothness
\emph{$\alpha_{\mathbf{T}}\equiv(\alpha_{\mathbf{T,N}})_{N=1,2,\cdots}$}
in the sense of Definition \ref{Def. Markov semigroup  semigroup}.
Let $E_{0}$ be an arbitrary distribution on $(S,d)$. Then the following
holds.

\emph{1. }For each $m\geq1$ and for each nondecreasing sequence $r_{1}\leq\cdots\leq r_{m}$
in {[}0,N{]}$Q$, the function
\[
F_{r(1),\cdots,r(m)}^{*,\mathbf{T}}:C(S^{m},d^{m})\rightarrow C(S,d),
\]
in Definition \ref{Def. Family of f.j.d.s from E0 and Markov semi-group},
is a well-defined transition distribution. More precisely, if $f\in C(S^{m},d^{m})$
has values in $[0,1]$ and has modulus of continuity $\delta_{f}$,
then the function $F_{r(1),\cdots,r(m)}^{*,\mathbf{T}}f$ is a member
of $C(S,d)$, with values in $[0,1]$ and with a modulus of continuity
$\widetilde{\alpha}_{\alpha(\mathbf{T},N),\xi}^{(m)}(\delta_{f})$,
where the operation $\widetilde{\alpha}_{\alpha(\mathbf{T},N),\xi}^{(m)}$
is the $m$-fold composite product operation
\[
\widetilde{\alpha}_{\alpha(\mathbf{T},N),\xi}^{(m)}\equiv\widetilde{\alpha}_{\alpha(\mathbf{T},N),\xi}\circ\cdots\circ\widetilde{\alpha}_{\alpha(\mathbf{T},N),\xi},
\]
where each factor on the right-hand side is as defined in Lemma \ref{Lem.One-step transition distributions}.

\emph{2.} For each $m\geq1$ and $\varepsilon>0$ there exists $\delta_{m}(\varepsilon,\delta_{f},\delta_{\mathbf{T}},\alpha_{\mathbf{T}})>0$
such that, for arbitrary $f\in C(S^{m},d^{m})$ with values in $[0,1]$
and with modulus of continuity $\delta_{f}$, and for arbitrary nondecreasing
sequences $r_{1}\leq\cdots\leq r_{m}$ and $s_{1}\leq\cdots\leq s_{m}$
in $Q$ with $\bigvee_{i=1}^{m}|r_{i}-s_{i}|<\delta_{m}(\varepsilon,\delta_{f},\delta_{\mathbf{T}},\alpha_{\mathbf{T}})$,
we have
\[
\left\Vert F_{r(1),\cdots,r(m)}^{*,\mathbf{T}}f-F_{s(1),\cdots s(m)}^{*,\mathbf{T}}f\right\Vert \leq\varepsilon.
\]

\emph{3. }In the special case where $Q=\{0,1,\cdots\}$ or $Q=\overline{Q}_{\infty}$,
the family
\[
F^{E(0),\mathbf{T}}\equiv\{F_{r(1),\cdots,r(m)}^{E(0),\mathbf{T}}f:m\geq1;r_{1}\leq\cdots\leq r_{m}\;\mathrm{in}\;Q\}
\]
can be uniquely extended to a consistent family 
\[
\Phi_{Sg,fjd}(E_{0},\mathbf{T)}\equiv F^{E(0),\mathbf{T}}\equiv\{F_{r(1),\cdots,r(m)}^{E(0),\mathbf{T}}f:m\geq1;r_{1},\cdots,r_{m}\;\mathrm{in}\;Q\}
\]
of f.j.d.'s with parameterr set $Q$ and state space $(S,d).$ Moreover,
the consistent family $\Phi_{Sg,fjd}(E_{0},\mathbf{T)}$ is generated
by the initial distribution $E_{0}$ and the semigroup \emph{$\mathbf{T}$,
}in the sense of Definition \emph{\ref{Def. Family of f.j.d.s from E0 and Markov semi-group}.
}Thus, in the special case where $Q=\{0,1,\cdots\}$ or $Q=\overline{Q}_{\infty}$,
we have a mapping
\[
\Phi_{Sg,fjd}:\widehat{J}(S,d)\times\mathscr{T}\rightarrow\widehat{F}(Q,S),
\]
where $\widehat{J}(S,d)$ is the space of distributions $E_{0}$ on
$(S,d)$, where $\mathscr{T}$ is the space of semigroups with state
space $(S,d)$ and with parameter set  $Q$, and where $\widehat{F}(Q,S)$
is the set of consistent families of \emph{f.j.d.'s }with parameter
set\emph{ $Q$} and state space $S$.

4. In the special case where\emph{ }$E_{0}\equiv\delta_{x}$ for some
$x\in S$, we simply write 
\[
\Phi_{Sg,fjd}(x,\mathbf{T)\equiv}\Phi_{Sg,fjd}(\delta_{x},\mathbf{T)\equiv\mathit{F}^{\mathit{x},\mathbf{T}}}.
\]
Then we have a function
\[
\Phi_{Sg,fjd}:S\times\mathscr{T}\rightarrow\widehat{F}(Q,S).
\]

5. Let $m\geq1$ and let the sequence $0\equiv r_{0}\leq r_{1}\leq\cdots\leq r_{m}$
in $Q$ be arbitrary. Then we have

\begin{equation}
F_{r(1),\cdots,r(m)}^{*,\mathbf{T}}=(^{1}T_{r(1)-r(0)})(^{2}T_{r(2)-r(1)})\cdots(^{m}T_{r(m)-r(m-1)}).\label{eq:temp-377-2}
\end{equation}
where the factors on the right hand side are one-step transition distributions
defined in Lemma \ref{Lem.One-step transition distributions}. In
other words,
\[
F_{r(1),\cdots,r(m)}^{x,\mathbf{T}}f
\]
\begin{equation}
=\int T_{r(1)}^{x}(dx_{1})\int T_{r(2)-r(1)}^{x(1)}(dx_{2})\cdots\int T_{r(m)-r(m-1)}^{x(m-1)}(dx_{m})f(x_{1},\cdots,x_{m}),\label{eq:temp-384-1-1-1-1}
\end{equation}
for each $x\in S$.
\end{thm}
\begin{proof}
1. Let $N\geq0$, $m\geq1$ and let the sequence $0\equiv r_{0}\leq r_{1}\leq\cdots\leq r_{m}$
in $Q[0,N]$ be arbitrary. By the defining equality \ref{eq:temp-384-1-1},
we have
\[
F_{r(1),\cdots,r(m)}^{E(0),\mathbf{T}}f
\]
\begin{equation}
=\int E_{0}(dx_{0})\int T_{r(1)}^{x(0)}(dx_{1})\int T_{r(2)-r(1)}^{x(1)}(dx_{2})\cdots\int T_{r(m)-r(m-1)}^{x(m-1)}(dx_{m})f(x_{1},\cdots,x_{m}).\label{eq:temp-384-1-1-1}
\end{equation}
By the defining equality \ref{eq:temp-426} of Lemma \ref{Lem.One-step transition distributions},
the right-most integral is equal to 
\begin{equation}
\int T_{r(n)-r(n-1)}^{x(m-1)}(dx_{m})f(x_{1},\cdots,x_{m-1},x_{m})\equiv((^{m-1}T_{r(m)-r(m-1)})f)(x_{1},\cdots,x_{m-1}).\label{eq:temp-426-6}
\end{equation}
Rcursively backward, we obtain
\[
F_{r(1),\cdots,r(m)}^{E(0),\mathbf{T}}f=\int E_{0}(dx_{0})(^{0}T_{r(1)-r(0)})\cdots(^{m-1}T_{r(m)-r(m-1)})f.
\]
In particular,
\[
F_{r(1),\cdots,r(m)}^{x,\mathbf{T}}f=(^{0}T_{r(1)-r(0)}^{x})\cdots(^{m-1}T_{r(m)-r(m-1)})f.
\]
for each $x\in S$. In other words,

\begin{equation}
F_{r(1),\cdots,r(m)}^{*,\mathbf{T}}=(^{1}T_{r(1)-r(0)})(^{2}T_{r(2)-r(1)})\cdots(^{m}T_{r(m)-r(m-1)}).\label{eq:temp-377}
\end{equation}
This proves Assertion 5 of the present Theorem. 

2. Now Lemma \ref{Lem.One-step transition distributions} says the
factors $(^{1}T_{r(1)-r(0)}),(^{2}T_{r(2)-r(1)}),\cdots,(^{m}T_{r(m)-r(m-1)})$
on the right-hand side are transition distibution with the modulus
of smoothness $\widetilde{\alpha}_{\alpha(T,N),\xi}$. Hence the composite
$F_{r(1),\cdots,r(m)}^{*,\mathbf{T}}$ is a composite transition distribution,
according to repeated applications of Lemma \ref{Lem. Composite transition distributions},
with the modulus of smoothness which is the $m$-fold composite operation
\[
\widetilde{\alpha}_{\alpha(T,N),\xi}^{(m)}\equiv\widetilde{\alpha}_{\alpha(\mathbf{T},N),\xi}\circ\cdots\circ\widetilde{\alpha}_{\alpha(\mathbf{T},N),\xi}.
\]
 This proves Assertion 1 of the present theorem. 

3. Proceed to prove Assertion 2 by induction on $m$. In the case
where $m=1$, define $\delta_{1}\equiv\delta_{1}(\varepsilon,\delta_{f},\delta_{\mathbf{T}},\alpha_{\mathbf{T}})\equiv\delta_{\mathbf{T}}(\varepsilon,\delta_{f})$.
Suppose $r_{1},s_{1}$ in $Q$ are such that 
\[
|r_{1}-s_{1}|<\delta_{1}(\varepsilon,\delta_{f},\delta_{\mathbf{T}},\alpha_{\mathbf{T}})\equiv\delta_{\mathbf{T}}(\varepsilon,\delta_{f}).
\]
Then 
\[
\left\Vert F_{r(1)}^{*,\mathbf{T}}f-F_{s(1)}^{*,\mathbf{T}}f\right\Vert =\left\Vert ^{1}T_{r(1)-r(0)}f-{}^{1}T_{s(1)-s(0)}f\right\Vert =\left\Vert T_{r(1)}f-T_{s(1)}f\right\Vert \leq\varepsilon
\]
where the inequality is by Lemma \ref{Lem. Uniform strong continuity on parameteer set}.
Assertion 2 is thus proved for the starting case $m=1$. 

4. Suppose it has been proved for $m-1$ for some $m\geq2,$ and the
operation
\[
\delta_{m-1}(\cdot,\delta_{f},\delta_{\mathbf{T}},\alpha_{\mathbf{T}})
\]
has been constructed with the desired propertiess. Let $f\in C(S^{m},d^{m})$
be arbitrary with values in $[0,1]$ and with modulus of continuity
$\delta_{f}$. Define
\begin{equation}
\delta_{m}(\varepsilon,\delta_{f},\delta_{\mathbf{T}},\alpha_{\mathbf{T}})\equiv2^{-1}\delta_{m-1}(2^{-1}\varepsilon,\delta_{f},\delta_{\mathbf{T}},\alpha_{\mathbf{T}})\wedge\delta_{m-1}(\varepsilon,\delta_{f},\delta_{\mathbf{T}},\alpha_{\mathbf{T}})\label{eq:temp-427}
\end{equation}
Suppose 
\begin{equation}
\bigvee_{i=1}^{m}|r_{i}-s_{i}|<\delta_{m}(\varepsilon,\delta_{f},\delta_{\mathbf{T}},\alpha_{\mathbf{T}}).\label{eq:temp-448}
\end{equation}
Deffine the function 
\begin{equation}
h\equiv(^{2}T_{r(2)-r(1)})\cdots(^{m}T_{r(m)-r(m-1)})f\in C(S,d).\label{eq:temp-563}
\end{equation}
Then, by the induction hypothesis for $(m-1)$-fold composite, the
function $h$ has modulus of continuity $\delta_{1}(\cdot,\delta_{f},\delta_{\mathbf{T}},\alpha_{\mathbf{T}})$.
We emphasize here that, as $m\geq2$, the modulus of smothness of
the one-step transition distribution $^{2}T_{r(2)-r(1)},$ on the
right-hand side of equality \ref{eq:temp-563} actually depends on
the modulus $\alpha_{\mathbf{T}}$, according to as Lemma \ref{Lem.One-step transition distributions}.
Hence the modulus of continuity of function $h$ indeed depends on
$\alpha_{\mathbf{T}}$, which justifies the notation.

At the same time, inequality \ref{eq:temp-448} and the defining equality
\ref{eq:temp-427} together imply that
\[
|r_{1}-s_{1}|<\delta_{m}(\varepsilon,\delta_{f},\delta_{\mathbf{T}},\alpha_{\mathbf{T}})\leq\cdots\leq\delta_{1}(\varepsilon,\delta_{f},\delta_{\mathbf{T}},\alpha_{\mathbf{T}}).
\]
Hence 
\begin{equation}
\left\Vert (^{1}T_{r(1)-r(0)})h-(^{1}T_{s(1)-s(0)})h\right\Vert =\left\Vert F_{r(1)}^{*,\mathbf{T}}h-F_{s(1)}^{*,\mathbf{T}}h\right\Vert \leq2^{-1}\varepsilon,\label{eq:temp-492}
\end{equation}
where the inequality is by the induction hypothesis for the starting
case where $m=1$. 

5. Similarly, inequality \ref{eq:temp-448} and the defining equality
\ref{eq:temp-427} together imply that
\[
\bigvee_{i=2}^{m}|(r_{i}-r_{i-1})-(s_{i}-s_{i-1})|\leq2\bigvee_{i=1}^{m}(r_{i}-s_{i}|<\delta_{m-1}(2^{-1}\varepsilon,\delta_{f},\delta_{\mathbf{T}},\alpha_{\mathbf{T}}).
\]
Hence 
\[
\left\Vert (^{2}T_{r(2)-r(1)})\cdots(^{m}T_{r(m)-r(m-1)})f-(^{2}T_{s(2)-s(1)})\cdots(^{m}T_{s(m)-s(m-1)})f\right\Vert 
\]
\begin{equation}
\equiv\left\Vert F_{r(2),\cdots,r(m)}^{*,\mathbf{T}}f-F_{r(s),\cdots,r(s)}^{*,\mathbf{T}}f\right\Vert <2^{-1}\varepsilon,\label{eq:temp-497}
\end{equation}
where the indequality isby the induction hypothesis for the case where
$m-1$. 

5. Combining, we estimate, for each $x\in S$, the bound
\[
|F_{r(1),\cdots,r(m)}^{x,\mathbf{T}}f-F_{s(1),\cdots,s(m)}^{x,\mathbf{T}}f|
\]
\[
=|(^{1}T_{r(1)-r(0)}^{x})(^{2}T_{r(2)-r(1)})\cdots(^{m}T_{r(m)-r(m-1)})f-(^{1}T_{s(1)-s(0)}^{x})(^{2}T_{s(2)-s(1)})\cdots(^{m}T_{s(m)-s(m-1)})f|
\]
\[
\equiv|(^{1}T_{r(1)-r(0)}^{x})h-(^{1}T_{s(1)-s(0)}^{x})(^{2}T_{s(2)-s(1)})\cdots(^{m}T_{s(m)-s(m-1)})f|
\]
\[
\leq|(^{1}T_{r(1)-r(0)}^{x})h-(^{1}T_{s(1)-s(0)}^{x})h|
\]
\[
+|(^{1}T_{s(1)-s(0)}^{x})h-(^{1}T_{s(1)-s(0)}^{x})(^{2}T_{s(2)-s(1)})\cdots(^{m}T_{s(m)-s(m-1)})f|
\]
\[
\leq2^{-1}\varepsilon+|(^{1}T_{s(1)-s(0)}^{x})(^{2}T_{r(2)-r(1)})\cdots(^{m}T_{r(m)-r(m-1)})f-(^{1}T_{s(1)-s(0)}^{x})(^{2}T_{s(2)-s(1)})\cdots(^{m}T_{s(m)-s(m-1)})f|
\]
\[
\leq2^{-1}\varepsilon+(^{1}T_{s(1)-s(0)}^{x})|(^{2}T_{r(2)-r(1)})\cdots(^{m}T_{r(m)-r(m-1)})f-(^{2}T_{s(2)-s(1)})\cdots(^{m}T_{s(m)-s(m-1)})f|
\]
\[
<2^{-1}\varepsilon+2^{-1}\varepsilon=\varepsilon,
\]
where the second inequality is by inequality \ref{eq:temp-492}, and
where the last inequality is by inquality \ref{eq:temp-497}. Since
$x\in S$ is arbitrary, it follows that 
\[
\left\Vert F_{r(1),\cdots,r(m)}^{*,\mathbf{T}}f-F_{s(1),\cdots s(m)}^{*,\mathbf{T}}f\right\Vert \leq\varepsilon.
\]
Induction is completed, and Assertion 2 is proved. 

6. We still need to prove Assertion 3. In other words, assuming that
$Q=\{0,1,\cdots\}$ or $Q=\overline{Q}_{\infty}$, we need to prove
that the family
\[
\{F_{r(1),\cdots,r(m)}^{E(0),\mathbf{T}}:m\geq1;r_{1},\cdots,r_{m}\in Q;r_{1}\leq\cdots\leq r_{m}\}
\]
can be uniquely extended to a consistent family 
\[
\{F_{s(1),\cdots,s(m)}^{E(0),\mathbf{T}}:m\geq1;s_{1},\cdots,s_{m}\in Q\}
\]
of f.j.d.'s with parameter set $Q$. We will give the proof only for
case where $Q=\overline{Q}_{\infty}$, the case of $\{0,1,\cdots\}$
being similar. So assume in the following that $Q=\overline{Q}_{\infty}$.

7. Because $F_{r(1),\cdots,r(m)}^{*,\mathbf{T}}$ is a transition
distribution, Proposition \ref{Prop. Completion of a transition distribution relativve to initial distribution}
says that the composite function
\begin{equation}
F_{r(1),\cdots,r(m)}^{E(0),\mathbf{T}}=E_{0}F_{r(1),\cdots,r(m)}^{*,\mathbf{T}}=E_{0}(^{1}T_{r(1)-r(0)})(^{2}T_{r(2)-r(1)})\cdots(^{m}T_{r(m)-r(m-1)}).\label{eq:temp-384-2-1}
\end{equation}
is a distribution on $(S^{m},d^{m})$, for each $m\geq1$, and for
each sequence $0\equiv r_{0}\leq r_{1}\leq\cdots\leq r_{m}$ in $\overline{Q}_{\infty}$
is arbitrary. 

8. To proceed, let $m\geq2$ and $r_{1},\cdots,r_{m}\in\overline{Q}_{\infty}$
be arbitrary, with $r_{1}\leq\cdots\leq r_{m}.$ Let $n=1,\cdots,m$
be arbitrary. Define the sequence
\[
\kappa\equiv\kappa_{n,m}\equiv(\kappa_{1},\cdots,\kappa_{m-1})\equiv(1,\cdots,\widehat{n},\cdots,m),
\]
where the caret on the top of an element in a sequence signifies the
omission of that element in the sequence. Let $\kappa^{*}\equiv\kappa_{n,m}^{*}:S^{m}\rightarrow S^{m-1}$
denote the dual function of the sequence $\kappa$, defined by 
\[
\kappa^{*}(x_{1},\cdots,x_{m})\equiv\kappa^{*}(x)\equiv x\circ\kappa=(x_{\kappa(1)},\cdots,x_{\kappa(m-1)})=(x_{1},\cdots,\widehat{x_{n}},\cdots,x_{m})
\]
for each $x\equiv(x_{1},\cdots,x_{m})\in S^{m}$. Let $f\in C(S^{m-1})$
be arbitrary. We will prove that 
\begin{equation}
F_{r(1),\cdots,\widehat{r(n)}\cdots,r(m)}^{E(0),\mathbf{T}}f=F_{r(1),\cdots,r(m)}^{E(0),\mathbf{T}}f\circ\kappa_{n,m}^{*}.\label{eq:temp-429-2}
\end{equation}
To that end, note that equality \ref{eq:temp-384-1-1-1} yields 
\[
F_{r(1),\cdots,\widehat{r(n)}\cdots,r(m)}^{E(0),\mathbf{T}}f
\]
\[
\equiv\int E_{0}(dx_{0})\int T_{r(1)}^{x(0)}(dx_{1})\cdots\int T_{r(n+1)-r(n-1)}^{x(n-1)}(dy_{n})\{
\]
\[
\int T_{r(n+2)-r(n+1)}^{y(n)}(dy_{n+1})\cdots\int T_{r(m)-r(m-1)}^{y(m-2)}(dy_{m-1})f(x_{1},\cdots,x_{n-1},y_{n},\cdots,y_{m-1})\}
\]
For each fixed $(x_{1},\cdots,x_{n-1})$, the expression in braces
is a continuous function of the one variable $y_{n}$. Call this function
$g_{x(1),\cdots,x(n-1)}\in C(S,d)$. Then the last displayed equality
can be continued as 
\[
\equiv\int E_{0}(dx_{0})\int T_{r(1)}^{x(0)}(dx_{1})\cdots(\int T_{r(n+1)-r(n-1)}^{x(n-1)}(dy_{n})g_{x(1),\cdots,x(n-1)}(y_{n}))
\]
\[
=\int E_{0}(dx_{0})\int T_{r(1)}^{x(0)}(dx_{1})\cdots(\int T_{r(n)-r(n-1)}^{x(n-1)}(dx_{n})\int T_{r(n+1)-r(n)}^{x(n)}(dy_{n})g_{x(1),\cdots,x(n-1)}(y_{n}))
\]
where the last equality is thanks to the semigroup property of $\mathbf{T}$.
Combining, we obtain
\[
F_{r(1),\cdots,\widehat{r(n)}\cdots,r(m)}^{E(0),\mathbf{T}}f
\]
\[
=\int E_{0}(dx_{0})\int T_{r(1)}^{x(0)}(dx_{1})\cdots\int T_{r(n)-r(n-1)}^{x(n-1)}(dx_{n})\int T_{r(n+1)-r(n)}^{x(n)}(dy_{n})g_{x(1),\cdots,x(n-1)}(y_{n})
\]
\[
\equiv\int E_{0}(dx_{0})\int T_{r(1)}^{x(0)}(dx_{1})\cdots\int T_{r(n)-r(n-1)}^{x(n-1)}(dx_{n})\int T_{r(n+1)-r(n)}^{x(n)}(dy_{n})\{
\]
\[
\int T_{r(n+2)-r(n+1)}^{y(n)}(dy_{n+1})\cdots\int T_{r(m)-r(m-1)}^{y(m-2)}(dy_{m-1})f(x_{1},\cdots,x_{n-1},y_{n},\cdots,y_{m-1})\}
\]
\[
=\int E_{0}(dx_{0})\int T_{r(1)}^{x(0)}(dx_{1})\cdots\int T_{r(n)-r(n-1)}^{x(n-1)}(dx_{n})\int T_{r(n+1)-r(n)}^{x(n)}(dx_{n+1})\{
\]
\[
\int T_{r(n+2)-r(n+1)}^{x(n+1)}(dx_{n+2})\cdots\int T_{r(m)-r(m-1)}^{x(m-1)}(dx_{m})f(x_{1},\cdots,x_{n-1},x_{n+1},\cdots,x_{m})\}
\]
\[
=F_{r(1),\cdots,r(m)}^{E(0),\mathbf{T}}(f\circ\kappa_{n,m}^{*}),
\]
where the third equality is by a trivial change of the dummy integration
variables $y_{n},\cdots,y_{m-1}$ to $x_{n+1},\cdots,x_{m}$ respectively.
Here the function 

Thus equality \ref{eq:temp-429-2} has been proved for the family
\[
\{F_{r(1),\cdots,r(m)}^{E(0),\mathbf{T}}:m\geq1;r_{1},\cdots,r_{m}\in Q;r_{1}\leq\cdots\leq r_{m}\}
\]
of f.j.d.'s. Consequently, the conditions in Lemma \ref{Lem. Consistency when parameter set is discrete subset of R}
are satisfied, to yield a unique extension of this family to a consistent
family 
\[
\{F_{s(1),\cdots,s(m)}^{E(0),\mathbf{T}}:m\geq0;s_{0},\cdots,s_{m}\in Q\}
\]
of f.j.d.'s with parameter set $Q$. Equality \ref{eq:temp-384-1-1-1}
says that $F^{E(0),\mathbf{T}}$ \emph{is} generated by the initial
distribution $E_{0}$ and the semigroup \emph{$\mathbf{T}$.} Assertion
3 is proved.
\end{proof}
$\square$

\section{Construction of Markov Process from Semigroup}

In this section, we construct Markov processs from a Markov semigroup
and an initial distribution. First the discrete parameters.
\begin{thm}
\textbf{\emph{\label{Thm. Markov process from initial distribution and Markov semigroup}
(Construction of Markov Process with initial distribution and a semigroup,
assuming discrete parameters).}} Let $(S,d)$ be a compact metric
space with $d\leq1$. Suppose $Q$ is one of the two parameter sets
$\{0,1,\cdots\}$ or $\overline{Q}_{\infty}$. Let 
\[
(\Theta_{0},L_{0},I_{0})\equiv([0,1],L_{0},\int\cdot dx)
\]
denote the Lebesgue integration space based on the unit interval $\Theta_{0}$.

Let $E_{0}$ be an arbitrary initial distribution on $(S,d)$, and
let $\mathbf{T}\equiv\{T_{t}:t\in Q\}$ be an arbitrary Markov semigroup
with state space $(S,d)$, with a modulus of strong continuity $\delta_{\mathbf{T}}$.
Let 
\[
F^{E(0),\mathbf{T}}\equiv\Phi_{Sg,fjd}(E_{0},\mathbf{T)}
\]
be the corresponding consistent family of f.j.d.'s constructed in
Theorem \ref{Thm. Construction of transition  f.j.d.s from E0 and Markov semi-group, discrete parameters-1}.
Let 
\[
Z^{E(0),\mathbf{T}}\equiv\overline{\Phi}_{DKS,\xi}(F^{E(0),\mathbf{T}}):Q\times\Theta_{0}\rightarrow S
\]
be the Compact Daniell-Kolmogorov-Skorokhod Extension of the consistent
family $F^{E(0),\mathbf{T}}$, as constructed in Theorem \ref{Thm. Compact Daniell-Kolmogorov-Skorohod}.

Then the following holds.

\emph{1.} The process \textbf{$Z^{E(0),\mathbf{T}}$ }is generated
by the initial distribution $E_{0}$ and semigroup $\mathbf{T}$,
in the sense of Definition \ref{Def. Family of f.j.d.s from E0 and Markov semi-group}.

\noun{2.} The process \textbf{$Z\equiv Z^{E(0),\mathbf{T}}$ }is a
Markov\emph{ }process relative to its natural filtration $\mathcal{L}\equiv\{L^{(t)}:t\in Q\}$.
More precisely, let $t\in Q$ be arbitrary, and define
\[
L^{(t)}\equiv L(Z_{s}^{E(0),\mathbf{T}}:s\in[0,t]Q\}\subset L_{0}.
\]
Let the nondecreasing sequence $0\equiv s_{0}\leq s_{1}\leq\cdots\leq s_{m}$
in $Q$, the function $f\in C(S^{m+1},d^{m+1})$, and $t\in Q$ be
arbitrary. Then
\[
I_{0}(f(Z_{t+s(0)},Z_{t+s(1)},\cdots,Z_{t+s(m)})|L^{(t)})
\]
\begin{equation}
=I_{0}(f(Z_{t+s(0)},Z_{t+s(1)},\cdots,Z_{t+s(m)})|Z_{t})=F_{s(0),\cdots,s(m)}^{Z(t),\mathbf{T}}(f)\label{eq:temp-575}
\end{equation}
as r.r.v.'s, where $F_{s(0),\cdots,s(m)}^{*,\mathbf{T}}$ is the transition
distribution in Definition \ref{Def. Family of f.j.d.s from E0 and Markov semi-group}.

\emph{3.} The process \textbf{$Z^{E(0),\mathbf{T}}$} is continuous
in probability, with a modulus of continuity in probability $\delta_{Cp,\delta(\mathbf{T})}$
which is completely determined by $\delta_{\mathbf{T}}$. 
\end{thm}
\begin{proof}
1. Let $t\in Q$ be arbitrary. Let $0\equiv r_{0}\leq r_{1}\leq\cdots\leq r_{n}\equiv t$
and $0\equiv s_{0}\leq s_{1}\leq\cdots\leq s_{m}$ be arbitrary sequences
in $Q$. Write $r_{n+j}\equiv t+s_{j}$ for each $j=0,\cdots,m$.
Thus $s_{j}=r_{n+j}-r_{n}$ for each $j=0,\cdots,m$. Consider each
$f\in C(S^{m+1},d^{m+1})$. Let $h\in C(S^{n+1},S^{n+1})$ be arbitrary.
Then, since the process $Z\equiv Z^{E(0),\mathbf{T}}\equiv\overline{\Phi}_{DKS,\xi}(F^{E(0),\mathbf{T}})$has
marginal distributions given by the family $F^{E(0),\mathbf{T}}$,
we have
\[
I_{0}h(Z_{r(0)},\cdots,Z_{r(n)})f(Z_{r(n)},\cdots,Z_{r(n+m)})
\]
\[
=\int F_{r(0),\cdots,r(n+m)}^{E(0),\mathbf{T}}d(x_{0},\cdots,x_{n+m})h(x_{0},\cdots,x_{n})f(x_{n},\cdots,x_{n+m})
\]
\[
\equiv\int E_{0}(dx_{0})\int T_{r(1)}^{x(0)}(dx_{1})\int T_{r(2)-r(1)}^{x(1)}(dx_{2})\cdots\int T^{x(n+m-1)}(dx_{n+m})h(x_{1},\cdots,x_{n})f(x_{n},\cdots,x_{n+m})
\]
\[
=\int E_{0}(dx_{0})\int T_{r(1)}^{x(0)}(dx_{1})\int T_{r(2)-r(1)}^{x(1)}(dx_{2})\cdots\int T_{r(n)-r(n-1)}^{x(n-1)}(dx_{n})h(x_{1},\cdots,x_{n})\{
\]
\begin{equation}
\int T_{r(n+1)-r(n)}^{x(n)}(dx_{n+1})\cdots\int T_{r(n+m)-r(n+m-1)}^{x(n+m-1)}(dx_{n+m})f(x_{n},x_{n+1}\cdots,x_{n+m})\}.\label{eq:temp-516-1}
\end{equation}
The term inside the braces in the last expression is, by changing
the names of the dummy integration variables, equal to
\[
\int T_{r(n+1)-r(n)}^{x(n)}(dy_{1})\int T_{r(n+2)-r(n+1)}^{y(1)}(dy_{2})\cdots\int T_{r(n+m)-r(n+m-1)}^{y(m-1)}(dy_{m})f(x_{n},y_{1},\cdots,y_{m})
\]
\[
=\int T_{x(1)}^{x(n)}(dy_{0})\int T_{s(2)}^{y(1)}(dy_{2})\cdots\int T_{s(m)}^{y(m-1)}(dy_{m})f(x_{n},y_{1},\cdots,y_{m})
\]
\[
=F_{s(1),\cdots,s(m)}^{x(n),\mathbf{T}}f.
\]
Substituting back into equality \ref{eq:temp-516-1}, we obtain
\[
I_{0}h(Z_{r(0)},\cdots,Z_{r(n)})f(Z_{r(n)},\cdots,Z_{r(n+m)})
\]
\[
=\int E_{0}(dx_{0})\int T_{r(1)}^{x(0)}(dx_{1})\int T_{r(2)-r(1)}^{x(1)}(dx_{2})\cdots\int T_{r(n)-r(n-1)}^{x(n-1)}(dx_{n})F_{s(1),\cdots,s(m)}^{x(n),\mathbf{T}}f
\]
\begin{equation}
=I_{0}h(Z_{r(0)},\cdots,Z_{r(n)})F_{s(1),\cdots,s(m)}^{Z(r(n)),\mathbf{T}}f,\label{eq:temp-508}
\end{equation}
where $F_{s(1),\cdots,s(m)}^{*,\mathbf{T}}f\in C(S,d)$ because $F_{s(1),\cdots,s(m)}^{*,\mathbf{T}}$
is a transition distribution according to Theorem \ref{Thm. Construction of transition  f.j.d.s from E0 and Markov semi-group, discrete parameters-1}.

2. Next note that the set of r.r.v.'s $h(Z_{r(0)},\cdots,Z_{r(n)})$,
with arbitrary $0\equiv r_{0}\leq r_{1}\leq\cdots\leq r_{n}\equiv t$
and arbitrary $h\in C(S^{n+1},d^{n+1})$, is dense in $L^{(t)}$relative
to the norm $I_{0}|\cdot|$. Hence equality \ref{eq:temp-508} implies,
by continuity relative to the norm $I_{0}|\cdot|$, that 
\begin{equation}
I_{0}Yf(Z_{r(n)},\cdots,Z_{r(n+m)})=I_{0}YF_{s(0),\cdots,s(m)}^{Z(r(n)),\mathbf{T}}(f)\label{eq:temp-541-2}
\end{equation}
for each $Y\in L^{(t)}$. It follows that 
\begin{equation}
I_{0}(f(Z_{r(n)},\cdots,Z_{r(n+m)})|L^{(t)})=F_{s(0),\cdots,s(m)}^{Z(r(n)),\mathbf{T}}(f),\label{eq:temp-459-4}
\end{equation}
or, equivalently,
\begin{equation}
I_{0}(f(Z_{t},Z_{t+s(1)},\cdots,Z_{t+s(n)})|L^{(t)})=F_{s(0),\cdots,s(m)}^{Z(t),\mathbf{T}}(f).\label{eq:temp-253}
\end{equation}
In the special case where $Y$ is arbitrary in $L(Z_{t})\subset L^{(t)},$
inequality \ref{eq:temp-541-2} holds, whence
\begin{equation}
I_{0}(f(Z_{t},Z_{t+s(1)},\cdots,Z_{t+s(n)})|Z_{t})=F_{s(0),\cdots,s(m)}^{Z(t),\mathbf{T}}(f).\label{eq:temp-253-1}
\end{equation}
Equalities \ref{eq:temp-253} and \ref{eq:temp-253-1} together prove
Assertions 1 and 2.

2. It remains to prove that the process $Z$ is continuous in probability.
We need only give the proof in the case where $Q=\overline{Q}_{\infty}$,
the case where $Q=\{0,1,\cdots\}$ being trivial. To that end, recall
that $d\leq1$. Consider each $x\in S$. Then the function $f\equiv1-d(\cdot,x)\in C(S,d)$
has a modulus of continuity $\delta_{f}$ defined by $\delta_{f}(\varepsilon)\equiv\iota(\varepsilon)\equiv\varepsilon$
for each $\varepsilon>0$. Now let $\varepsilon>0$ be arbitrary.
Let $t\in\overline{Q}_{\infty}$ be arbitrary with 
\[
t<\delta_{Cp,\delta(\mathbf{T})}(\varepsilon)\equiv\delta_{\mathbf{T}}(\varepsilon,\iota)\equiv\delta_{\mathbf{T}}(\varepsilon,\delta_{f}).
\]
Then, by Lemma \ref{Lem. Uniform strong continuity on parameteer set},
we have
\begin{equation}
|d(\cdot,x)-T_{t}d(\cdot,x)|=|f-T_{t}f|\leq\varepsilon,\label{eq:temp-382-1-2-2-1-1}
\end{equation}
where $x\in S$ is arbitrary. Now let $r_{1},r_{2}\in Q$ be arbitrary
with $|r_{2}-r_{1}|<\delta_{Cp}(\varepsilon)$. Then
\[
I_{0}d(Z_{r(1)},Z_{r(2)})=I_{0}d(Z_{r(1)\wedge r(2)},Z_{r(1)\vee r(2)})=F_{0,r(1)\wedge r(2),r(1)\vee r(2)}^{E(0),\mathbf{T}}d
\]
\[
=\int E_{0}(dx_{0})\int T_{r(1)\wedge r(2)}^{x(0)}(dx_{1})\int T_{r(1)\vee r(2)-r(1)\wedge r(2)}^{x(1)}(dx_{2})d(x_{1},x_{2})
\]
\[
=\int E_{0}(dx_{0})\int T_{r(1)\wedge r(2)}^{x(0)}(dx_{1})\int T_{|r(2)-r(1)|}^{x(1)}(dx_{2})d(x_{1},x_{2})
\]
\[
=\int E_{0}(dx_{0})\int T_{r(1)\wedge r(2)}^{x(0)}(dx_{1})T_{|r(2)-r(1)|}^{x(1)}d(x_{1},\cdot)
\]
\[
\leq\int E_{0}(dx_{0})\int T_{r(1)\wedge r(2)}^{x(0)}(dx_{1})(d(x_{1},x_{1})+\varepsilon)
\]
\[
=\int E_{0}(dx_{0})\int T_{r(1)\wedge r(2)}^{x(0)}(dx_{1})\varepsilon=\varepsilon,
\]
where the inequality is by applying inequality \ref{eq:temp-382-1-2-2-1-1}
to $t\equiv|r_{2}-r_{1}|$ and $x=x_{1}$. Thus we have shown that
$\delta_{Cp,\delta(\boldsymbol{\mathbf{T)}}}\equiv\delta_{\mathbf{T}}(\cdot,\iota)$
is a modulus of continuity in probability for the process $Z$, according
to Definition \ref{Def. continuity in prob, continuity a.u., and a.u. continuity}.
Note here that $\delta_{Cp,\delta(\boldsymbol{\mathbf{T)}}}$ is an
opration dependent only on $\delta_{\mathbf{T}}$. Assertion 3 is
proved.
\end{proof}
The next proposition says that each Markov process $Z:\overline{Q}_{\infty}\times(\Omega,L,E)\rightarrow(S,d)$
with a semigroup can be extended by right limit to an a.u. càdlàg
process $X:[0,\infty)\times(\Omega,L,E)\rightarrow(S,d)$. 
\begin{prop}
\textbf{\emph{\label{Prop. Process with Markov semigroup on dyadic rationals in =00005B0,oo) is D-regular}
(Process with semigroup on dyadic rationals is $D$-regular, extendable
to time-uniformly a.u. càdlàg process). }}Let $(S,d)$ be a compact
metric space with $d\leq1$. Suppose $Q=\overline{Q}_{\infty}$. Let
$E_{0}$ be an arbitrary initial distribution on $(S,d)$, and let
$\mathbf{T}\equiv\{T_{t}:t\in\overline{Q}_{\infty}\}$ be an arbitrary
Markov semigroup with state space $(S,d)$ and with a modulus of strong
continuity $\delta_{\mathbf{T}}$. Let 
\[
Z:\overline{Q}_{\infty}\times(\Omega,L,E)\rightarrow(S,d)
\]
be an arbitrary process generated by the  initial distribution $E_{0}$
and semigroup $\mathbf{T}$. 

Let $N\geq0$ be arbitrary. Define the shifted process $Z^{N}:Q_{\infty}\times\Omega\rightarrow S$
by $Z^{N}(t,\cdot)\equiv Z(N+t,\cdot)$ for each $t\in Q_{\infty}$.
Then the following holds.

\emph{1}. The process $Z^{N}$ is continuous in probability, with
a modulus of continuity in probability $\delta_{Cp,\delta(\mathbf{T})}\equiv\delta_{Cp}(\cdot,\delta_{\mathbf{T}})$
whih is completely determined by $\delta_{\mathbf{T}}$. 

\emph{2.} The process $Z^{N}$ is strongly right continuous in probability,
in the sense of Definition \ref{Def. strongly  right continuity in probability on Q_inf},
with a modulus of strong right continuity in probability given by
the operation $\delta_{SRcp,\delta(\mathbf{T})}$ defined by
\[
\delta_{\gamma,\delta(\mathbf{T})}(\varepsilon,\gamma)\equiv\delta_{Cp}(\varepsilon^{2},\delta_{\mathbf{T}})
\]
for each $\varepsilon\gamma>0$. Note that $\delta_{SRcp,\delta(\mathbf{T})}(\varepsilon,\gamma)$
is actually independent of $\gamma$ or $N$.

\emph{3. }The process $Z^{N}:Q_{\infty}\times\Omega\rightarrow S$
is $D$-regular, with a modulus of $D$-regularity $\overline{m}_{\delta(\mathbf{T})}\equiv\overline{m}_{\delta(\mathbf{T})}$. 

\emph{4.} The process $Z$ is time-uniformly $D$-regular in the sense
of Definition \ref{Def. D-regular processes wirh dyadic rationals in =00005B0,oo)as parameters-1},
with a modulus of continuity in probability $\widetilde{\delta}{}_{Cp,\delta(\mathbf{T)}}\equiv(\delta_{Cp,\delta(\mathbf{T})},\delta_{Cp,\delta(\mathbf{T})},\cdots)$
and with a modulus of $D$-regularity $\widetilde{m}_{\delta(\mathbf{T})}\equiv(\overline{m}_{\delta(\mathbf{T})},\overline{m}_{\delta(\mathbf{T})},\cdots)$. 

\emph{5.} The right-limit extension
\[
X\equiv\overline{\Phi}_{rLim}(Z):[0,\infty)\times(\Omega,L,E)\rightarrow(S,d)
\]
is a time-uniformly a.u. càdlàg process in the sense of Definition
\ref{Def. Metric space of a.u. cadlag process on =00005B0,inf)},
with a modulus of a.u. càdlàg $\widetilde{\delta}_{aucl,\delta(\mathbf{T})}\equiv(\delta_{aucl,\delta(\mathbf{T})},\delta_{aucl,\delta(\mathbf{T})},\cdots)$
and with a modulus of continuity in probability $\widetilde{\delta}{}_{Cp,\delta(\mathbf{T)}}$.
Here we recall the function $\overline{\Phi}_{rLim}$ from Definition
\ref{Def. Extension of process by right limit}. In other words, 
\[
X\equiv\overline{\Phi}_{rLim}(Z)\in\widehat{D}_{\widetilde{\delta}(aucl,\delta(\mathbf{T})),\widetilde{\delta}(Cp,\delta(\mathbf{T)})}[0,\infty).
\]
\end{prop}
\begin{proof}
1. By hypothesis, the process $Z$ has initial distribution $E_{0}$
and Markov semigroup $\mathbf{\mathbf{T}}$. In other words, the process
$Z$ has marginal distributions given by the consistent family $F^{E(0),\mathbf{T}}$
of f.j.d.'s. as in Definition \ref{Def. Family of f.j.d.s from E0 and Markov semi-group}.
Hence the process $Z$ is equivalent to the process $Z^{E(0),\mathbf{T}}$
constructed in Theorem \ref{Thm. Markov process from initial distribution and Markov semigroup}.
Therefore, by Assertion 3 of the Theorem \ref{Thm. Markov process from initial distribution and Markov semigroup},
the process $Z$ is continuous in probability, with a modulus of continuity
in probability $\delta_{Cp,\delta(\mathbf{T})}\equiv\delta_{Cp}(\cdot,\delta_{\mathbf{T}})$
which is is completely determined by $\delta_{\mathbf{T}}$. 

2. Hence, trivially, the shifted process
\[
Y\equiv Z^{N}:Q_{\infty}\times\Omega\rightarrow S
\]
is continuous in probability, with a modulus of continuity in probability
$\delta_{Cp,\delta(\mathbf{T})}\equiv\delta_{Cp}(\cdot,\delta_{\mathbf{T}})$
which is completely determined by $\delta_{\mathbf{T}}$. Assertion
1 is proved. 

3. To prove Assertion 2, let $\varepsilon,\gamma>0$ be arbitrary.
Define 
\begin{equation}
\delta_{SRcp,\delta(\mathbf{T})}(\varepsilon,\gamma)\equiv\delta_{\mathbf{T}}(\varepsilon^{2},\iota)\label{eq:temp-252}
\end{equation}
where the operation $\iota$ is defined by $\iota(\varepsilon')\equiv\varepsilon'$
for each \textbf{$\varepsilon'>0.$}

Let $h\geq0$ and $s,r\in Q_{h}$ be arbitrary with $s\leq r<s+\delta_{SRcp,\delta(\mathbf{T})}(\varepsilon,\gamma)$.
Then 
\[
Q_{h}=\{0,\Delta_{h},2\Delta_{h},\cdots,1\}\equiv\{q_{h,0},\cdots,q_{h,p(h)}\}
\]
 where $\Delta\equiv2^{-h},p_{h}\equiv2^{h},.$and $q_{h,i}\equiv i\Delta,$
for each $i=0,\cdots,p_{h}.$ Moreover, $s=q_{h,i}$ and $r=q_{h,j}$
for some $i,j=0,\cdots,p_{h}$ with $i\leq j$. Now let $g\in C(S^{i+1},d^{i+1})$
be arbitrary. Then
\[
Eg(Y_{q(h,0)},\cdots,Y_{q(h,i)})d(Y_{s},Y_{r})=Eg(Y_{q(h,0)},\cdots,Y_{q(h,i)})d(Y_{q(h,i)},Y_{q(h,j)})
\]
\[
\equiv Eg(Z_{N+q(h,0)},\cdots,Z_{N+q(h,i)})d(Z_{N+s},Z_{N+r})
\]
\[
=Eg(Z_{N+0\Delta},\cdots,Z_{N+i\Delta})d(Z_{N+i\Delta},Z_{N+j\Delta})
\]
\[
=\int E_{0}(dx_{0})\int F_{N+0\Delta,\cdots,N+i\Delta,N+j\Delta}^{x(0),\mathbf{T}}(d(x_{0},\cdots,x_{i},x_{j}))g(x_{0},\cdots,x_{i})d(x_{i},x_{j})
\]
\[
=\int E_{0}(dx_{0})\int T_{N+0\Delta}^{x(0)}(dx_{1})\int T_{\Delta}^{x(1)}(dx_{1})\cdots\int T_{\Delta}^{x(i-1)}(dx_{i})\int T_{(j-i)\Delta}^{x(i)}(dx_{j})g(x_{0},\cdots,x_{i})d(x_{i},x_{j})
\]
\begin{equation}
=\int E_{0}(dx_{0})\int T_{\Delta}^{x(0)}(dx_{1})\cdots\int T_{\Delta}^{x(i-1)}(dx_{i})g(x_{0},\cdots,x_{i})T_{(j-i)\Delta}^{x(i)}d(x_{i},\cdot),\label{eq:temp-147}
\end{equation}
where the next-to-last equality is because the family $F^{x(0),\mathbf{T}}$
of f.j.d.'s is generated by the initial state $x_{0}$ and Markov
semigroup $\mathbf{T}$, for each $x_{0}\in S$. Now consider each
$x_{i}\in S$. Since the function $d(x_{i},\cdot)$ has modulus of
continuity $\iota$, and since 
\[
(j-i)\Delta=r-s<\delta_{SRcp,\delta(\mathbf{T})}(\varepsilon,1)\equiv\delta_{\mathbf{T}}(\varepsilon^{2},\iota),
\]
we have, according to the definition of $\delta_{\mathbf{T}}$ as
the modulus of strong continuity of $\mathbf{T}$, the bound
\begin{equation}
|d(x_{i},\cdot)-T_{(j-i)\Delta}d(x_{i},\cdot)|\leq\varepsilon^{2}\label{eq:temp-382-1-2-2-1}
\end{equation}
as functions on $S$. In particular, 
\begin{equation}
|d(x_{i},x_{i})-T_{(j-i)\Delta}^{x(i)}d(x_{i},\cdot)|\leq\varepsilon^{2}.\label{eq:temp-382-1-2-2-1-2}
\end{equation}
Consequently $T_{(j-i)\Delta}^{x(i)}d(x_{i},\cdot)\leq\varepsilon^{2}$
where $x_{i}\in S$ is arbitrary. Hence equality \ref{eq:temp-147}
can be continued to yield
\[
Eg(Y_{q(h,0)},\cdots,Y_{q(h,i)})d(Y_{s},Y_{r})\leq\int E_{0}(dx_{0})\int T_{\Delta}^{x(0)}(dx_{1})\cdots\int T_{\Delta}^{x(i-1)}(dx_{i})g(x_{0},\cdots,x_{i})\varepsilon^{2}
\]
\[
\equiv\varepsilon^{2}Eg(Y_{q(h,0)},\cdots,Y_{q(h,i)}),
\]
where $g\in C(S^{i+1},d^{i+1})$ is arbitrary. It follows that 
\begin{equation}
EUd(Y_{s},Y_{r})\leq\varepsilon^{2}EU\label{eq:temp-417}
\end{equation}
for each $U\in L(Y_{q(h,0)},Y_{q(h,1)},\cdots,Y_{q(h,i)}).$ Now let
$\gamma>0$ be arbitrary, and take an arbitrary measurable set
\begin{equation}
A\in L^{(s,h)}\equiv L(Y_{r}:r\in[0,s]Q_{h})=L(Y_{q(h,0)},Y_{q(h,1)},\cdots,Y_{q(h,i)})\label{eq:temp-429}
\end{equation}
with $A\subset(d(x_{\circ},Y_{s})\leq\gamma)$, and with $P(A)>0$.
Let $U\equiv1_{A}$ denote the indicator of $A$. Then, in view of
the relation \ref{eq:temp-429}, we have $U\equiv1_{A}\in L(Y_{q(h,0)},Y_{q(h,1)},\cdots,Y_{q(h,i)}).$
Hence equality \ref{eq:temp-417} holds, to yield
\begin{equation}
E1_{A}d(Y_{s},Y_{r})\leq\varepsilon^{2}E1_{A}\label{eq:temp-417-2}
\end{equation}
Equivalently
\[
E_{A}d(Y_{s},Y_{r})\leq\varepsilon^{2}
\]
where $E_{A}$ is the conditional expectation given the event $A$.
Chebychev's inequality therefore implies
\[
P_{A}(d(Y_{s},Y_{r})>\alpha)\leq\varepsilon
\]
for each $\alpha>0$. Here $h\geq0,$ $\varepsilon,\gamma>0$, and
$s,r\in Q_{h}$ are arbitrary with $s\leq r<s+\delta_{SRcp,\delta(\mathbf{T})}(\varepsilon,\gamma)$.
Summing up, the process $Y$ is strongly right continuous in the sense
of Definition \ref{Def. strongly  right continuity in probability on Q_inf},
with a modulus of strong right continuity $\delta_{SRcp,\delta(\mathbf{T})}$.
Assertion 2 is proved. 

4. Proceed to prove Assertions 3-5. To that end, recall that $d\leq1$
by hypothesis. Hence the process $Y\equiv Z^{N}:Q_{\infty}\times\Omega\rightarrow(S,d)$
is trivially a.u. bounded, with a trivial modulus of a.u. boundlessness
\[
\beta_{auB}(\cdot,\delta_{\mathbf{T}})\equiv1.
\]
Combining with Assertion 2, we see that the conditions for Theorem
\ref{Thm. a.u.boundeness and Strongly right continuity impliy D-regular}
are satisfied for the process 
\[
Y\equiv Z^{N}:Q_{\infty}\times\Omega\rightarrow(S,d)
\]
to be $D$-regular, with a modulus of $D$-regularity 
\[
\overline{m}_{\delta(\mathbf{T})}\equiv\overline{m}(\beta_{auB},\delta_{SRcp,\delta(\mathbf{T})})\equiv\overline{m}(1,\delta_{SRcp,\delta(\mathbf{T})})
\]
and a modulus of continuity in probability $\delta{}_{Cp,\mathbf{\delta(\mathbf{T})}}$,
and for the process $Y\equiv Z^{N}$ to be extendable by right-limit
to an a.u. càdlàg process
\[
\Phi_{rLim}(Z^{N}):[0,1]\times\Omega\rightarrow S
\]
which has the same modulus of continuity in probability $\delta{}_{Cp}(\cdot,\beta_{auB},\delta_{SRcp})$
as $Z^{N}$, and which has a modulus of a.u. càdlàg 
\[
\delta_{aucl,\mathbf{\delta(\mathbf{T})}}\equiv\delta_{aucl}(\cdot,\beta_{auB},\delta_{SRcp,\delta(\mathbf{T})})\equiv\delta_{aucl}(\cdot,1,\delta_{SRcp,\delta(\mathbf{T})}).
\]
Recall here the right-limit extension $\Phi_{rLim}$ from Definition
\ref{Def. Extension of process by right limit}. Assertion 3 is proved.

5. Moreover, since the moduli $\delta{}_{Cp,\mathbf{\delta(\mathbf{T})}}$
and $\overline{m}_{\delta(\mathbf{T})}$ are independent of the interval
$\{N,N+1]$, we see that the process 
\[
Z:\overline{Q}_{\infty}\times\Omega\rightarrow(S,d)
\]
is time-uniformly $D$-regular in the sense of Definition \ref{Def. D-regular processes wirh dyadic rationals in =00005B0,oo)as parameters-1},
with modulus of continuity in probability $\widetilde{\delta}{}_{Cp,\mathbf{\delta(\mathbf{T})}}$$\equiv(\delta{}_{Cp,\mathbf{\delta(\mathbf{T})}},\delta{}_{Cp,\mathbf{\delta(\mathbf{T})}},\cdots),$
and with a modulus of $D$-regularity $\widetilde{m}_{\delta(\mathbf{T})}\equiv(\overline{m}_{\delta(\mathbf{T})},\overline{m}_{\delta(\mathbf{T})},\cdots)$.
This proves Assertion 4 of the present theorem is proved. 

6. To prove the remaining Asssertion 5. Consider the right-limit extension
process 
\[
X\equiv\overline{\Phi}_{rLim}(Z):[0,\infty)\times(\Omega,L,E)\rightarrow(S,d).
\]
Consider each $N\geq0$. Then clearly $X^{N}=\Phi_{rLim}(Z^{N})$
on the interval $[0,1)$. Near the end point $1$, things are a bit
more subtle. We need to recall that, since $\Phi_{rLim}(Z^{N})$ is
a.u. càdlàg, it is continuous a.u. on $[0,1]$. Hence, for a.s.$\omega\in\Omega$,
the funtion $Z(\cdot,\omega)$ is continuous at $1\in[0,1]$. It therefore
follows that $X^{N}=\Phi_{rLim}(Z^{N})$ on the interval $[0,1]$.
We saw in Step 4 that the process $\Phi_{rLim}(Z^{N})$ is a.u. càdlàg,
with a modulus of a.u. càdlàg $\delta_{aucl,\mathbf{\delta(\mathbf{T})}}$. 

7. Now note that, as an immediate consequence of Step 5, the process
$Z|]0,N+1]$ is continuous in probability with modulus of continuity
continuity in probability $\delta{}_{Cp,\mathbf{\delta(\mathbf{T})}}.$
It follows that the processs $X|[N+1]$ is continuous in probability
with modulus of continuity continuity in probability $\delta{}_{Cp,\mathbf{\delta(\mathbf{T})}}.$ 

8.. Summing up the results in Steps 7 and 6, we see that the process
$X$ time-uniformly a.u. càdlàg in the sense of Definition \ref{Def. Metric space of a.u. cadlag process on =00005B0,inf)},
with modulus of continuity in probability $\widetilde{\delta}{}_{Cp,\mathbf{\delta(\mathbf{T})}}$$\equiv(\delta{}_{Cp,\mathbf{\delta(\mathbf{T})}},\delta{}_{Cp,\mathbf{\delta(\mathbf{T})}},\cdots),$
and with a modulus of a.u. càdlàg $\widetilde{\delta}_{aucl,\delta(\mathbf{T})}\equiv(\delta_{aucl,\mathbf{\delta(\mathbf{T})}},\delta_{aucl,\mathbf{\delta(\mathbf{T})}},\cdots)$
. Assertions 5 and the theorem are proved. 
\end{proof}
\begin{thm}
\textbf{\emph{\label{Thm. Markov process from initial distribution and Markov semigroup-1}
(Construction of a.u. càdlàg}} \textbf{\emph{Markov Process from an
initial distribution and a semigroup, with parameter set $[0,\infty)$).}}
Let $(S,d)$ be a compact metric space with $d\leq1$. Let $E_{0}$
be an arbitrary initial distribution on $(S,d)$, and let $\mathbf{T}\equiv\{T_{t}:t\in[0,\infty)\}$
be an arbitrary Markov semigroup with state space $(S,d)$, with a
modulus of strong continuity $\delta_{\mathbf{T}}$, and with a modulus
of \emph{\noun{smoothness}}\emph{ $\alpha_{\mathbf{T}}\equiv(\alpha_{\mathbf{T,N}})_{N=1,2,\cdots}$}
in the sense of Definition \ref{Def. Markov semigroup  semigroup}..
Let 
\[
(\Theta_{0},L_{0},I_{0})\equiv([0,1],L_{0},\int\cdot dx)
\]
denote the Lebesgue integration space based on the interval $\Theta_{0}$.
Let 
\[
Z^{E(0),\mathbf{T|\mathrm{\overline{Q}(\infty)}}}\equiv\overline{\Phi}_{DKS,\xi}(F^{E(0),\mathbf{T|\mathrm{\overline{Q}(\infty)}}}):\overline{Q}_{\infty}\times\Theta_{0}\rightarrow S
\]
be the Markov process with initial distribution $E_{0}$ and Markov
semigroup $\mathbf{T}|Q_{\infty}$, as constructed in Assertion 2
of Theorem \ref{Thm. Markov process from initial distribution and Markov semigroup}. 

Then the following holds.

1. The right-limit extension
\[
X^{E(0),\mathbf{T}}\equiv\overline{\Phi}_{rLim}(Z^{E(0),\mathbf{T|\mathrm{Q(\infty)}}}):[0,\infty)\times\Theta_{0}\rightarrow(S,d)
\]
is a time-uniformly a.u. càdlàg\emph{ }process. 

2. The the process $X^{E(0),\mathbf{T}}$ is Markov relative to its
natural filtration $\overline{\mathcal{L}}\equiv\{\overline{L}^{(t)}:t\in[0,\infty)\}$.
Specifically, let $v\geq0$ and let $t_{0}\equiv0\leq t_{1}\leq\cdots\leq t_{m}$
be an arbitrary sequence in $[0,\infty),$with $m\geq1$. Then
\begin{equation}
I_{0}(f(X_{v+t(0)},\cdots,X_{v+t(m)})|\overline{L}^{(v)})=I_{0}(f(X_{v+t(0)},\cdots,X_{v+t(m)})|X_{v})=F_{0,t(1)\cdots,t(m)}^{X(v),\mathbf{T}}f,\label{eq:temp-559}
\end{equation}
where $F_{0,t(1)\cdots,t(m)}^{*,\mathbf{T}}$ is the transition distribution
in Definition \ref{Def. Family of f.j.d.s from E0 and Markov semi-group}.

3. Let $F^{E(0),\mathbf{T}}$ denote the family of marginal f.j.d.'s
of the process $X$. \textup{\emph{In the special case where $E_{0}\equiv\delta_{x}$
assigns probability $1$ to some point $x\in S$, write }}\emph{$F^{x,\mathbf{T}}\equiv F^{\delta(x),\mathbf{T}}$.}
Then the family\textup{ }$F^{E(0),\mathbf{T}}$ \textup{\emph{is}}\emph{
}generated by the initial distribution $E_{0}$ and the semigroup\emph{$\mathbf{T}$,
in the sense of }\textup{\emph{Definition \ref{Def. Family of f.j.d.s from E0 and Markov semi-group}.
In particular, the family }}\emph{$F^{x,\mathbf{T}}$ }is generated
by the initial state $x$ and the semigroup\emph{$\mathbf{T}$.}
\end{thm}
\begin{proof}
To minimize notational clutter, we will write $a=b\pm c$ to mean
$|a-b|\leq c,$ for arbitrary real-valued expressions $a,b,c.$ Let
$v\geq0$ and the sequence $t_{0}\equiv0\leq t_{1}\leq\cdots\leq t_{m}$
in $[0,\infty)$ be arbitrary, but fixed. Let $n\geq1$ and $v_{0}\equiv0\leq v_{1}\leq\cdots\leq v_{n-1}$
in $[0,v]$ be arbitrary. Define $v_{n+i}\equiv v+t_{i}$ for each
$i=0,\cdots,m$. Note that 
\[
v_{0}\equiv0\leq v_{1}\leq\cdots\leq v_{n+m}.
\]
Fix any integer $N\geq0$ so large that $v_{n+m}\in[0,N-1]$.

1. Assertion 1 of the present theorem is merely a restatement of Assertion
5 of Proposition \ref{Prop. Process with Markov semigroup on dyadic rationals in =00005B0,oo) is D-regular},
which says that the right-limit extension $X\equiv X^{E(0),\mathbf{T}}$
is a time-uniformly a.u. càdlàg process, with a modulus of a.u. càdlàg
$\widetilde{\delta}_{aucl,\delta(\mathbf{T})}\equiv(\delta_{aucl,\delta(\mathbf{T})},\delta_{aucl,\delta(\mathbf{T})},\cdots)$
and with a modulus of continuity in probability $\widetilde{\delta}{}_{Cp,\delta(\mathbf{T)}}$$\equiv(\delta_{Cp,\delta(\mathbf{T)}},\delta_{Cp,\delta(\mathbf{T)}},\cdots)$.
In particular, the process $X|[0,N+1]$ is continuous in probability,
with a modulus of continuity in probability $\delta_{Cp,\delta(\mathbf{T)}}$.

2. Proceed to prove Assertion 2, by extending the Markov proeperty
of the process $Z$ to the process $X$. First,let $\mathcal{L}\equiv\{L^{(t)}:t\in\overline{Q}_{\infty}\}$
be the natural filtration of the process $Z\equiv Z^{E(0),\mathbf{T|\mathrm{\overline{Q}(\infty)}}}$.
Assertion 2 of Theorem \ref{Thm. Markov process from initial distribution and Markov semigroup}
says that the process \textbf{$Z\equiv Z^{E(0),\mathbf{T}}$ }is a
Markov\emph{ }process relative to the filtration $\mathcal{L}$. Specifically,
let the nondecreasing sequence $0\equiv s_{0}\leq s_{1}\leq\cdots\leq s_{m}$
in $\overline{Q}_{\infty}$, the function $f\in C(S^{m+1},d^{m+1})$,
and the pointy $t\in\overline{Q}_{\infty}$ be arbitrary. Then
\[
I_{0}(f(Z_{t+s(0)},Z_{t+s(1)},\cdots,Z_{t+s(m)})|L^{(t)}))
\]
\begin{equation}
=I_{0}(f(Z_{t+s(0)},Z_{t+s(1)},\cdots,Z_{t+s(m)})|Z_{t}))=F_{s(0),\cdots,s(m)}^{Z(t),\mathbf{T|\mathrm{\overline{Q}(\infty)}}}(f)\label{eq:temp-575-1}
\end{equation}
as r.r.v.'s, where $F_{s(0),\cdots,s(m)}^{*,T|\mathrm{\overline{Q}(\infty)}}$
is the transition distribution as in Definition \ref{Def. Family of f.j.d.s from E0 and Markov semi-group}. 

3. Next, let $\overline{\mathcal{L}}\equiv\{\overline{L}^{(t)}:t\in[0,\infty)\}$
be the natural filtration of the process $X$. Let the sequence $0\equiv r_{0}\leq r_{1}\leq\cdots\leq r_{n}\leq\cdots\leq r_{n+m}$
in $[0,N]\overline{Q}_{\infty}$ be arbitrary. Let the function $g\in C(S^{n+1},d^{n+1})$
be arbitrary, with values in $[0,1]$ and with modulus of continuity
$\delta_{g}$ and $\delta_{f}$ respectively. Then equality \ref{eq:temp-575-1}
implies that 
\[
I_{0}g(Z_{r(0)},\cdots,Z_{r(n)})f(Z_{r(n)},\cdots,Z_{r(n+m)})
\]
\[
I_{0}(I_{0}g(Z_{r(0)},\cdots,Z_{r(n)})f(Z_{r(n)},\cdots,Z_{r(n+m)})|L^{(r(n))}))=I_{0}g(Z_{r(0)},\cdots,Z_{r(n)})F_{0,r(n+1)-r(n),\cdots,r(n+m)-r(n)}^{Z(r(n)),\mathbf{T|\mathrm{\overline{Q}(\infty)}}}(f).
\]
Since $X_{r}=Z_{r}$ for each $r\in\overline{Q}_{\infty}$, this is
equivalent to
\begin{equation}
I_{0}g(X_{r(0)},\cdots,X_{r(n)})f(X_{r(n)},\cdots,X_{r(n+m)})=I_{0}g(X_{r(0)},\cdots,X_{r(n)})F_{0,r(n+1)-r(n),\cdots,r(n+m)-r(n)}^{X(r(n)),\mathbf{T}}f.\label{eq:temp-254}
\end{equation}
where $F_{0,r(n+1)-r(n),\cdots,r(n+m)-r(n)}^{*,\mathbf{T}}$ is the
transition distribution as in Definition \ref{Def. Family of f.j.d.s from E0 and Markov semi-group}. 

4. Next let $r_{i}\downarrow v_{i}$ for each $i=0,\cdots,n+m$. Then
the left-hand side of equality \ref{eq:temp-254} converges to the
limit
\[
I_{0}g(X_{v(0)},\cdots,X_{v(n)})f(X_{v(n)},\cdots,X_{v(n+m)}),
\]
thanks to the continuity in probability of the process $X|[0,N+1]$. 

5. Consider the right-hand side of equality \ref{eq:temp-254}. Let
$\varepsilon>0$ be arbitrary. As observed in Step 1 above, the process
$X|[0,N+1]$ has a modulus of continuity in probability $\delta_{Cp,\delta(\mathbf{T)}}$.
Consequently, there exists $\delta_{0}>0$ so small that 
\begin{equation}
I_{0}|g(X_{r(0)},\cdots,X_{r(n)})-g(X_{v(0)},\cdots,X_{v(n)})|<\varepsilon\label{eq:temp-564}
\end{equation}
prvided that $\bigvee_{i=0}^{n}(r_{i}-v_{i})<\delta_{0}$.

6. Separately, Assertion 2 of Theorem \ref{Thm. Construction of transition  f.j.d.s from E0 and Markov semi-group, discrete parameters-1}
implies that there exists $\delta_{m+1}(\varepsilon,\delta_{f},\delta_{\mathbf{T}},\alpha_{\mathbf{T}})>0$
such that, for arbitrary nondecreasing sequences $0\equiv\overline{r}_{0}\leq\overline{r}_{1}\leq\cdots\leq\overline{r}_{m}$
and $0\equiv\overline{s}_{0}\leq\overline{s}_{1}\leq\cdots\leq\overline{s}_{m}$
in $[0,\infty)$ with 
\[
\bigvee_{i=0}^{m}|\overline{r}_{i}-\overline{s}_{i}|<\delta_{m+1}(\varepsilon,\delta_{f},\delta_{\mathbf{T}},\alpha_{\mathbf{T}}),
\]
we have
\begin{equation}
\left\Vert F_{\overline{r}(0),\overline{r}(1),\cdots,\overline{r}(m)}^{*,\mathbf{T}}f-F_{\overline{s}(0),\overline{s}(1),\cdots\overline{s(}m)}^{*,\mathbf{T}}f\right\Vert <\varepsilon.\label{eq:temp-548}
\end{equation}

7. Suppose 
\begin{equation}
\bigvee_{i=0}^{n+m}(r_{i}-v_{i})<2^{-1}\delta_{m+1}(\varepsilon,\delta_{f},\delta_{\mathbf{T}},\alpha_{\mathbf{T}})\wedge\delta{}_{Cp,\delta(\mathbf{T)}}(\varepsilon\delta_{f}(\varepsilon))\wedge\delta_{0}.\label{eq:temp-545}
\end{equation}
Let $\overline{r}_{i}\equiv r_{n+i}-r_{n}$ and $\overline{s}_{i}\equiv t_{i}=v_{n+i}-v_{n}$
for each $i=0,\cdots,m$. Then $\overline{r}_{i},\overline{s}_{i}\in[0,N]$,
with
\[
\bigvee_{i=0}^{m}|\overline{r}_{i}-\overline{s}_{i}|\leq2\bigvee_{i=0}^{n+m}(r_{i}-v_{i})<\delta_{m+1}(\varepsilon,\delta_{f},\delta_{\mathbf{T}},\alpha_{\mathbf{T}}).
\]
Hence inequality \ref{eq:temp-548} holds. For abbreviation, define

\[
h\equiv F_{\overline{r}(0),\overline{r}(1),\cdots,\overline{r}(m)}^{*,\mathbf{T}}f
\]
and 
\[
\overline{h}\equiv F_{\overline{s}(0),\overline{s}(1),\cdots\overline{s}(m)}^{*,\mathbf{T}}f.
\]
Then inequality \ref{eq:temp-548} can be rewritten as 
\begin{equation}
\left\Vert h-\overline{h}\right\Vert <\varepsilon.\label{eq:temp-566}
\end{equation}

8. Since $0\equiv\overline{r}_{0}\leq\overline{r}_{1}\leq\cdots\leq\overline{r}_{m}$
is a sequence in $[0,N]\cap[0,\infty)$, Assertion 1 of Theorem \ref{Thm. Construction of transition  f.j.d.s from E0 and Markov semi-group, discrete parameters-1}implies
that the function $h\equiv F_{\overline{r}(0),\overline{r}(1),\cdots,\overline{r}(m)}^{*,\mathbf{T}}f$
is a member of $C(S,d)$, with values in $[0,1]$ and with a modulus
of continuity $\delta_{h}\equiv\widetilde{\alpha}_{\alpha(\mathbf{T},N),\xi}^{(m)}(\delta_{f})$,
where the operation $\widetilde{\alpha}_{\alpha(\mathbf{T},N),\xi}^{(m)}$
is the $m$-fold composite product operation
\[
\widetilde{\alpha}_{\alpha(\mathbf{T},N),\xi}^{(m)}\equiv\widetilde{\alpha}_{\alpha(\mathbf{T},N),\xi}\circ\cdots\circ\widetilde{\alpha}_{\alpha(\mathbf{T},N),\xi},
\]
where each factor on the right-hand side is as defined in Lemma \ref{Lem.One-step transition distributions}. 

9. Now suppose, in addition to the bound \ref{eq:temp-545}, we also
have
\begin{equation}
\bigvee_{i=0}^{n+m}(r_{i}-v_{i})<\delta{}_{Cp,\delta(\mathbf{T)}}(\varepsilon\delta_{h}(\varepsilon))\equiv\delta{}_{Cp,\delta(\mathbf{T)}}(\varepsilon\widetilde{\alpha}_{\alpha(\mathbf{T},N),\xi}^{(m)}(\delta_{f})(\varepsilon)).\label{eq:temp-545-1}
\end{equation}
Then, in particular, we have 
\[
r_{n}-v_{n}<\delta{}_{Cp,\delta(\mathbf{T)}}(\varepsilon\delta_{h}(\varepsilon)).
\]
Hence 
\[
I_{0}d(X_{r(n)},X_{v(n)})<\varepsilon\delta_{h}(\varepsilon)
\]
by the definition of $\delta{}_{Cp,\delta(\mathbf{T)}}$ as a modulus
of continuity in probability of the process $X|[0,N+1]$. Therefore
Chebychev's inequality yields a measurable set $A\subset\Theta_{0}$
with $I_{0}A^{c}<\varepsilon$ such that 
\begin{equation}
A\subset(d(X_{r(n)},X_{v(n)})<\delta_{h}(\varepsilon)).\label{eq:temp-510}
\end{equation}

10. Since the function $h$ has modulus of continuity $\delta_{h}$
, relation \ref{eq:temp-510} immediately extends to
\[
A\subset(d(X_{r(n)},X_{v(n)})<\delta_{f}(\varepsilon))\subset(|h(X_{r(n)})-h(X_{v(n)})|<\varepsilon).
\]
As a reseult, we obtain the estimate
\[
I_{0}g(X_{r(0)},\cdots,X_{r(n)})F_{0,r(n+1)-r(n),\cdots,r(n+m)-r(n)}^{X(r(n)),\mathbf{T}}f
\]
\[
\equiv(I_{0}g(X_{r(0)},\cdots,X_{r(n)})h(X_{r(n)})1_{A})+(I_{0}g(X_{r(0)},\cdots,X_{r(n)})h(X_{r(n)})1_{A^{c}})
\]
\[
=(I_{0}g(X_{r(0)},\cdots,X_{r(n)})h(X_{r(n)})1_{A})\pm\varepsilon
\]
\begin{equation}
=(I_{0}g(X_{r(0)},\cdots,X_{r(n)})h(X_{v(n)})1_{A})\pm\varepsilon.\label{eq:temp-565}
\end{equation}
At this point, note that the bound \ref{eq:temp-545} implies that
$\bigvee_{i=0}^{n}(r_{i}-v_{i})<\delta_{0}$. Hence inequality \ref{eq:temp-564}
holds, and leads to

\[
I_{0}(g(X_{r(0)},\cdots,X_{r(n)})h(X_{v(n)})1_{A})=(I_{0}g(X_{v(0)},\cdots,X_{v(n)})h(X_{v(n)})1_{A})\pm\varepsilon.
\]
Therefore equality \ref{eq:temp-565} can be continued, to yield
\[
I_{0}g(X_{r(0)},\cdots,X_{r(n)})F_{0,r(n+1)-r(n),\cdots,r(n+m)-r(n)}^{X(r(n)),\mathbf{T}}
\]
\[
=(I_{0}g(X_{v(0)},\cdots,X_{v(n)})h(X_{v(n)})1_{A})\pm2\varepsilon
\]
\[
=(I_{0}g(X_{v(0)},\cdots,X_{v(n)})h(X_{v(n)})1_{A})+(I_{0}g(X_{v(0)},\cdots,X_{v(n)})h(X_{v(n)})1_{A^{c}})\pm3\varepsilon
\]
\[
=I_{0}g(X_{v(0)},\cdots,X_{v(n)})h(X_{v(n)})\pm3\varepsilon
\]
\[
=I_{0}g(X_{v(0)},\cdots,X_{v(n)})\overline{h}(X_{v(n)})\pm\varepsilon\pm3\varepsilon
\]
\[
=I_{0}g(X_{v(0)},\cdots,X_{v(n)})\overline{h}(X_{v(n)})\pm4\varepsilon,
\]
where we used the condition that the functions $f,g,h$ have values
in $[0,1]$, and where the fourth equality is thanks to equality \ref{eq:temp-566}.
Summing up, we have proved that 
\[
|I_{0}g(X_{r(0)},\cdots,X_{r(n)})F_{0,r(n+1)-r(n),\cdots,r(n+m)-r(n)}^{X(r(n)),\mathbf{T}}f-I_{0}g(X_{v(0)},\cdots,X_{v(n)})\overline{h}(X_{v(n)})|\leq4\varepsilon
\]
provided that the bounds \ref{eq:temp-545} and \ref{eq:temp-545-1}
are satisfied. Since $\varepsilon>0$ is arbitraily small, we have
proved the convergence of the right-hand side of equality \ref{eq:temp-254},
with, specifically,
\[
I_{0}g(X_{r(0)},\cdots,X_{r(n)})F_{0,r(n+1)-r(n),\cdots,r(n+m)-r(n)}^{X(r(n)),\mathbf{T}}f\rightarrow I_{0}g(X_{v(0)},\cdots,X_{v(n)})\overline{h}(X_{v(n)}),
\]
as $r_{i}\downarrow v_{i}$ for each $i=0,\cdots,n+m$. In view of
the convegence of the left-hand side of the same equality \ref{eq:temp-254},
as observed in Step 4, the two limits are equal. Thus 
\begin{equation}
I_{0}g(X_{v(0)},\cdots,X_{v(n)})f(X_{v(n)},\cdots,X_{v(n+m)})=I_{0}g(X_{v(0)},\cdots,X_{v(n)})\overline{h}(X_{v(n)}),\label{eq:temp-550}
\end{equation}
where the function $g\in C(S^{n+1},d^{n+1})$ is arbitrary with values
in $[0,1],$ and where the integer $n\geq1$ and the sequence $0\equiv v_{0}\leq v_{1}\leq\cdots\leq v_{n-1}$
are arbitrary. Hence, by linearity, equality \ref{eq:temp-550} implies
that
\[
I_{0}Yf(X_{v(n)},\cdots,X_{v(n+m)})=I_{0}Y\overline{h}(X_{v(n)})
\]
for each r.r.v. $Y$ in the family
\[
G_{v(n)}\equiv\{g(X_{v(0)},\cdots,X_{v(n-1)},X_{v(n)}):g\in C(S^{n+1},d^{n+1});(v_{0},\cdots,v_{n-1})\in S^{n};v_{0}\leq v_{1}\leq\cdots\leq v_{n-1}\leq v_{n}\}
\]
Since the family $G_{v(n)}$ is dense in the space $\overline{L}^{(v(n))}\equiv\overline{L}^{(v)}$
relative to the norm $I_{0}|\cdot|$, it follows that 
\begin{equation}
I_{0}Yf(X_{v(n)},\cdots,X_{v(n+m)})=I_{0}Y\overline{h}(X_{v(n)})\label{eq:temp-544}
\end{equation}
for each $Y\in\overline{L}^{(v)}$. Siince $\overline{h}(X_{v(n)})\in G_{v(n)}\subset\overline{L}^{(v)}$,
we obtain the conditional expectation 
\begin{equation}
I_{0}(f(X_{v(n)},\cdots,X_{v(n+m)})|\overline{L}^{(v)})=\overline{h}(X_{v(n)})\equiv F_{\overline{s}(0),\overline{s}(1),\cdots\overline{s(}m))}^{X(v(n)),\mathbf{T}}f.\label{eq:temp-554}
\end{equation}
In the special case of an arbitrary $Y\in L(X_{v})\subset\overline{L}^{(v)}$,
equality \ref{eq:temp-544} holds. Hence 
\begin{equation}
I_{0}(f(X_{v(n)},\cdots,X_{v(n+m)})|X_{v})=\overline{h}(X_{v(n)})\equiv F_{\overline{s}(0),\overline{s}(1),\cdots\overline{s(}m))}^{X(v(n)),\mathbf{T}}f.\label{eq:temp-558}
\end{equation}
Recall that $v_{n+i}\equiv v+t_{i}$ and $\overline{s}_{i}\equiv t_{i}$,
for each $i=0,\cdots,m,$ and recall that $t_{0}\equiv0$. Equalities
\ref{eq:temp-554} and \ref{eq:temp-558} can then be rewritten as
\begin{equation}
I_{0}(f(X_{v+t(0)},\cdots,X_{v+t(m)})|\overline{L}^{(v)})=F_{0,t(1),\cdots,t(m)}^{X(v),\mathbf{T}}f.\label{eq:temp-554-1}
\end{equation}
and 
\begin{equation}
I_{0}(f(X_{v+t(0)},\cdots,X_{v+t(m)})|X_{v})=F_{0,t(1),\cdots,t(m)}^{X(v),\mathbf{T}}.\label{eq:temp-558-1}
\end{equation}
The lasst two equality together yield the desired equality \ref{eq:temp-559}.
Assertion 2 is proved. {[}

11.\emph{ }Finally, note that, by Assertion 2 of Theorem \ref{Thm. Construction of transition  f.j.d.s from E0 and Markov semi-group, discrete parameters-1},the
consistent family $F^{E(0),\mathbf{T|\overline{\mathrm{\mathit{Q}}}(\infty)}}$
is generated by the initial distribution $E_{0}$ and the semigroup
\emph{$\mathbf{T}|\overline{\mathrm{\mathit{Q}}}_{\infty}$, }in the
sense of Definition \emph{\ref{Def. Family of f.j.d.s from E0 and Markov semi-group}.
}Hence, for each sequence $0\equiv r_{0}\leq r_{1}\leq\cdots\leq r_{n}\leq\cdots\leq r_{n+m}$
in $\overline{Q}_{\infty}$ and for each $f\in C(S^{m},d^{m})$, we
have
\[
F_{r(1),\cdots,r(m)}^{E(0),\mathbf{T}}f\equiv I_{0}f(X_{r(1)},\cdots,X_{r(m)})=I_{0}f(Z_{r(1)},\cdots,Z_{r(m)})
\]
\begin{equation}
=\int E_{0}(dx_{0})\int T_{r(1)}^{x(0)}(dx_{1})\int T_{r(2)-r(1)}^{x(1)}(dx_{2})\cdots\int T_{r(m)-r(m-1)}^{x(m-1)}(dx_{m})f(x_{1},\cdots,x_{m}).\label{eq:temp-384-1-1-3}
\end{equation}
Because the process is continuous in probability, and because the
semigroup is strongly continuous, this equality extends to 
\[
F_{r(1),\cdots,r(m)}^{E(0),\mathbf{T}}f
\]
\begin{equation}
=\int E_{0}(dx_{0})\int T_{r(1)}^{x(0)}(dx_{1})\int T_{r(2)-r(1)}^{x(1)}(dx_{2})\cdots\int T_{r(m)-r(m-1)}^{x(m-1)}(dx_{m})f(x_{1},\cdots,x_{m})\label{eq:temp-384-1-1-3-1}
\end{equation}
for each sequence $0\equiv r_{0}\leq r_{1}\leq\cdots\leq r_{n}\leq\cdots\leq r_{n+m}$
in $[0,\infty)$. Thus the family is $F^{E(0),\mathbf{T}}$ \emph{is
}generated by the initial distribution $E_{0}$ and the semigroup\emph{$\mathbf{T}$,
}in the sense of Definition\emph{ \ref{Def. Family of f.j.d.s from E0 and Markov semi-group}.
}Assertion 3 and the theorem are proved.
\end{proof}

\section{Continuity of Construction}

In this section, we will prove that the construction in Theorem \ref{Thm. Markov process from initial distribution and Markov semigroup-1}
of a Markov process from an initial state $x$ and a semigroup $\mathbf{T}$
is uniformly metrically coontinuous over each subspace of semigroups
$\mathbf{T}$ whose members share a common modulus of strong continuity
and share a common modulus of smoothness. We will limit the discussion
to the parameter set $Q=\overline{Q}_{\infty}$or $Q=[0,\infty)$.
The case of $Q=\{0,1,\cdots\}$ being similar.

First we specify a compact state space, and define a metric on the
space of Markov semigroups.
\begin{defn}
\textbf{\label{Def. Specification of compact metric space and binary approx, and partition of unity}
(Specification of state space, its binary approximation, and partition
of unity). }In this section, unless otherwise spacified, $(S,d)$
will denote a given compact metric space, with $d\leq1$, and with
a fixed reference point $x_{\circ}\in S$. Recall that $\xi\equiv(A_{k})_{k=1,2,\cdots}$
is a binary approximation of $(S,d)$ relative to $x_{\circ}$, and
that 
\[
\pi\equiv(\{g_{k,x}:x\in A_{k}\})_{k=1,2,\cdots}
\]
is the partition of unity of $(S,d)$ determined by $\xi$, as in
Definition \ref{Def. Partition of unity for locally compact (S,d)}.
Recall that $|A_{k}|$ denotes the number of elements in the discrete
finite subset $A_{k}\subset S$, for each $k\geq1$.

For each $n\geq1$ let  $\xi^{n}$ denote the $n$-th power of $\xi$,
and let $\pi^{(n)}$denote the corresponding partition of unity for
$(S^{n},d^{n})$, Thus, for each $n\geq1,$ the sequence $\xi^{n}\equiv(A_{k}^{(n)})_{k=1,2,\cdots}$
is the product binary approximation for $(S^{n},d^{n})$ relative
to the reference point $x_{\circ}^{(n)}\equiv(x_{\circ},\cdots x_{\circ})\in S^{n}$,
and 
\[
\pi^{(n)}\equiv(\{g_{k,x}^{(n)}:x\in A_{k}^{(n)}\})_{k=1,2,\cdots}
\]
is the partition of unity of $(S^{n},d^{n})$ determined by the binary
approximation $\xi\equiv(A_{k})_{k=1,2,\cdots}\equiv(A_{k}^{(1)})_{k=1,2,\cdots}$
of $(S,d)$. For each $k\geq1$, the set $A_{k}^{(n)}$ is an $2^{k}$-approximation
of the bounded subset 
\[
(d^{n}(\cdot,(x_{\circ},\cdots,x_{o})\leq2^{k})\subset S^{n}.
\]
To lessen the burden of subscripts, we write $A_{k}^{(n)}$ and $A_{n,k}$
interchangeably, for each $n\geq1$. 

We will use the notations in Definitions \ref{Def. Notations for dyadic rationals-1}
related to the enumerated set $\overline{Q}_{\infty}\equiv\{u_{0},u_{1},\cdots\}$.
\end{defn}
$\square$
\begin{defn}
\textbf{\label{Def. Metric space of Markov semigroups} (Metric space
of Markov semigroups).} Let $(S,d)$ be the specified compact metric
space, with $d\leq1$. Suppose $Q=\overline{Q}_{\infty}$or $Q=[0,\infty)$.
Let $\mathscr{T}$ be  family of all Markov semigroups on the parameter
set $Q$ and with the compact metric state space $(S,d)$. For each
$n\leq0$ write $\Delta_{n}\equiv2^{-n}$ . Define the metric $\rho_{\mathscr{T}}\equiv\rho_{\mathscr{T,\xi}}$
on the family $\mathscr{T}$ by
\[
\rho_{\mathscr{T}}(\mathbf{T},\mathbf{\overline{T}})\equiv\rho_{\mathscr{T,\xi}}(\mathbf{T},\mathbf{\overline{T}})\equiv\sum_{n=0}^{\infty}2^{-n-1}\sum_{k=1}^{\infty}2^{-k}|A_{k}|^{-1}\sum_{z\in A(k)}\left\Vert T_{\Delta(n)}g_{k,z}-T{}_{\Delta(n)}g_{k,z}\right\Vert 
\]
for arbitrary members $\mathbf{T}\equiv\{T_{t}:t\in Q\}$ and $\mathbf{\overline{T}}\equiv\{\overline{T}{}_{t}:t\in Q\}$
of the family $\mathscr{T}$ . Here $\left\Vert \cdot\right\Vert $
stands for the supremum norm on $C(S,d)$.
\end{defn}
It follows easily from the strong continuity of the semigorups $\mathbf{T}$
that $\rho_{\mathscr{T}}$ is in fact a metric. Note that $\rho_{\mathscr{T}}\leq1$.
Let $(S\times\mathscr{T},d\otimes\rho_{\mathscr{T}})$ denote the
product metric metric space of $(S,d)$ and $(\mathscr{T},\rho_{\mathscr{T}})$.

For each $\mathbf{T}\equiv\{T_{t}:t\in Q\}\in\mathscr{T}$, define
the semigroup $\mathbf{T}|\overline{Q}_{\infty}\equiv\{T_{t}:t\in\overline{Q}_{\infty}\}$.
Then
\[
\rho_{\mathscr{T,\xi}}(\mathbf{T}|\overline{Q}_{\infty},\mathbf{\overline{T}}|\overline{Q}_{\infty})=\rho_{\mathscr{T,\xi}}(\mathbf{T},\mathbf{\overline{T}}).
\]
In other words, the mapping 
\[
\Psi:\mathscr{(\mathscr{\mathscr{T}}},\rho_{\mathscr{T,\xi}})\rightarrow\mathscr{(\mathscr{\mathscr{T}}}|\overline{Q}_{\infty},\rho_{\mathscr{T,\xi}}),
\]
defined by $\mathit{\mathrm{\Psi(}}\mathbf{T})\equiv\mathbf{T|}\overline{Q}_{\infty}$
for each $\mathbf{T}\in\mathscr{\mathscr{\mathscr{T}}}$, is an isometry. 

$\square$

The next theorem proves the promised continuity of construction in
the case where the parameter set is $\overline{Q}_{\infty}$. Then
the sunsequent Theorem \ref{Thm. Construction of Markov process from  initial distrinbution and semigrp-1-1}
will prove a similar continuity in the case where the parameter set
$[0,\infty)$.
\begin{thm}
\textbf{\emph{\label{Thm.  Continuity of construction of f.j.d.s from E0 and Markov semigroups-1-1}
(Continuity of construction of family of transition f.j.d.'s from
initial state and semigroup, for discrete parameters). }}Let $(S,d)$
be the specified compact metric space, with $d\leq1$.\textbf{\emph{
}}Let $\mathscr{T}(\overline{\delta},\overline{\alpha})$ be an arbitrary
family of  Markov semigroups with parameter set $\overline{Q}_{\infty}$
and state space $(S,d)$, such that\emph{ }all its members $\mathbf{T}\in\mathscr{T}(\overline{\delta},\overline{\alpha})$
share a common modulus of strong continuity $\delta_{\mathbf{T}}=\overline{\delta}$,
and share a modulus of smoothness $\alpha_{\mathbf{T}}=\overline{\alpha},$
in the sense of Definition \ref{Def. Markov semigroup  semigroup}.
Thus $\mathscr{T}(\overline{\delta},\overline{\alpha})$ is a subset
of $(\mathscr{T},\rho_{\mathscr{T}}),$ and, as such, inherits the
metric $\rho_{\mathscr{T,\xi}}$ introduced in Definition \ref{Def. Metric space of Markov semigroups}.
Recall from Definition \ref{Def. Marginal metric} the metric space
$(\widehat{F}(\overline{Q}_{\infty},S),\widehat{\rho}_{Marg,\xi,\overline{Q}(\infty)})$
of consistent families of f.j.d.'s with paramter set $\overline{Q}_{\infty}$
and state space $(S,d)$. 

Then the mapping
\[
\Phi_{Sg,fjd}:(S\times\mathscr{\mathscr{T}(\overline{\delta}\mathit{,}\overline{\alpha})},d\otimes\rho_{\mathscr{T,\xi}})\rightarrow(\widehat{F}(\overline{Q}_{\infty},S),\widehat{\rho}_{Marg,\xi,\overline{Q}(\infty)}),
\]
constructed in Assertion 4 of Theorem  \ref{Thm. Construction of transition  f.j.d.s from E0 and Markov semi-group, discrete parameters-1},
is uniformly continuous, with a modulus of continuity $\delta_{Sg,fjd}(\cdot,\overline{\delta},\overline{\alpha},\left\Vert \xi\right\Vert )$
determined by the moduli $\overline{\delta},\overline{\alpha}$ and
by the modulus of local compactness $\left\Vert \xi\right\Vert \equiv(|A_{n}|)_{n=1,2,\cdots}$of
$(S,d)$. 
\end{thm}
\begin{proof}
1. Let $(x,\mathbf{T}),(\overline{x},\overline{\mathbf{T}})\in S\times\mathscr{\mathscr{T}(\overline{\delta},\overline{\alpha})}$
be arbitrary. For abbreviation, write $F\equiv F^{x,\mathbf{T}}\equiv\Phi_{Sg,fjd}(x,\mathbf{T})$
and $\overline{F}\equiv F^{\overline{x},\mathbf{\mathbb{\overline{\mathbf{T}}}}}\equiv\Phi_{Sg,fjd}(\overline{x},\overline{\mathbf{T}})$.
Define the distance
\begin{equation}
\overline{\rho}_{0}\equiv(d\otimes\rho_{\mathscr{T}})((x,\mathbf{T}),(\overline{x},\overline{\mathbf{T}}))\equiv d(x,\overline{x})\vee\rho_{\mathscr{T}}(\mathbf{T},\overline{\mathbf{T}})\label{eq:temp-567}
\end{equation}

2. Let $\varepsilon_{0}>0$ be arbitrary, but fixed till further notice.
Take $M\geq0$ so large that $2^{-N-1}<\frac{1}{3}\varepsilon_{0}$,
where $N\equiv p_{2M}\equiv2^{2M}$. Write $\Delta\equiv\Delta_{M}\equiv2^{-M}$
for abbreviation. Then, in the notations of Definitions \ref{Def. Notations for dyadic rationals-1},
we have 

\begin{equation}
\overline{Q}_{M}\equiv\{u_{0},u_{1},\cdots,u_{p(2M)}\}\equiv\{u_{0},u_{1},\cdots,u_{N}\}=\{0,\Delta,2\Delta,\cdots,N\Delta\}\label{eq:temp-221-2}
\end{equation}
as enumerated sets. Thus $u_{n}=n\Delta$ for each $n=0,\cdots,N$.
For abbreviation, also write $K\equiv N+1,$ and
\[
\gamma_{N}\equiv(N+1)^{-1}2^{-N-2}.
\]
Then $2^{-K}=2^{-N-1}<\frac{1}{3}\varepsilon_{0}$. 

3. By Definition \ref{Def. Marginal metric}, we have
\[
\widehat{\rho}_{Marg,\xi,\overline{Q}(\infty)}(F,\overline{F})\equiv\sum_{n=0}^{\infty}2^{-n-1}\rho_{Dist,\xi^{n+1}}(F_{u(0),\cdots,u(n)},\overline{F}{}_{u(0),\cdots,u(n)})
\]
\[
\leq\sum_{n=0}^{N}2^{-n-1}\rho_{Dist,\xi^{n+1}}(F_{u(0),\cdots,u(n)},\overline{F}{}_{u(0),\cdots,u(n)})+2^{-N-1}
\]
\begin{equation}
\leq\sum_{n=0}^{N}2^{-n-1}\rho_{Dist,\xi^{n+1}}(F_{0,\Delta,\cdots,n\Delta},\overline{F}{}_{0,\Delta,\cdots,n\Delta})+\frac{1}{3}\varepsilon_{0}.\label{eq:temp-208-2-2-1}
\end{equation}
For each $n\geq0$, the metric $\rho_{Dist,\xi^{n+1}}$ was introduced
in Definition \ref{Def. Distribution metric} for the space of distributions
on $(S^{n+1},d^{n+1})$, where it is observed that $\rho_{Dist,\xi^{n+1}}\leq1$
and that sequential convergence relative to $\rho_{Dist,\xi^{n+1}}$
is equivalent to weak convergence. 

4. We will prove that the $n$-th summand in the last sum in equality
is bounded by $2^{-n-1}\frac{2}{3}\varepsilon_{0}$, for each $n=0,\cdots,N$,
provided that the the distance $\overline{\rho}_{0}$ is sufficiently
small.

To that end, let $n=0,\cdots,N$ be arbitrary, but fixed till further
notice. Then, in the $n$-th summand on the right-hand side of inequality
\ref{eq:temp-208-2-2-1}, we have 
\[
\rho_{Dist,\xi^{n+1}}(F_{0,\Delta,\cdots,n\Delta},\overline{F}{}_{0,\Delta,\cdots,n\Delta})
\]
\[
\equiv\sum_{k=1}^{\infty}2^{-k}|A_{n+1,k}|^{-1}\sum_{y\in A(n+1,k)}|F_{0,\Delta,\cdots,n\Delta}g_{k,y}^{(n+1)}-\overline{F}{}_{0,\Delta,\cdots,n\Delta}g_{k,y}^{(n+1)}|
\]
\[
\leq\sum_{k=1}^{K}2^{-k}|A_{n+1,k}|^{-1}\sum_{y\in A(n+1,k)}|F_{0,\Delta,\cdots,n\Delta}g_{k,y}^{(n+1)}-\overline{F}{}_{0,\Delta,\cdots,n\Delta})g_{k,y}^{(n+1)}|+2^{-K},
\]
\begin{equation}
<\sum_{k=1}^{K}2^{-k}|A_{n+1,k}|^{-1}\sum_{y\in A(n+1,k)}|F_{0,\Delta,\cdots,n\Delta}g_{k,y}^{(n+1)}-\overline{F}{}_{0,\Delta,\cdots,n\Delta})g_{k,y}^{(n+1)}|+\frac{1}{3}\varepsilon_{0}\label{eq:temp-219-3-1}
\end{equation}
where the first equality is from Definition \ref{Def. Distribution metric}
for the distribution metric $\rho_{Dist,\xi^{n+1}}$, and where, for
each $k\geq1$ and each $y\in A_{n+1,k}$, the basis function $g_{k,y}^{(n+1)}\in C(S^{n+1},d^{n+1})$
is from the partition of unity 
\[
\pi^{(n+1)}\equiv(\{g_{k,y}^{(n+1)}:y\in A_{n+1,k}\})_{k=1,2,\cdots}
\]
of $(S^{n+1},d^{n+1})$ relative to the $2^{-k}$-apprximation $A_{n+1,k}$
of the metric space $(S^{n+1},d^{n+1})$, as specified in Definition
\ref{Def. Specification of compact metric space and binary approx, and partition of unity}.

5. Next, with $n=0,\cdots,N$ arbitrary but fixed, we will show that
\[
|F_{0,\Delta,\cdots,n\Delta}g_{k,y}^{(n+1)}-\overline{F}_{0,\Delta,\cdots,n\Delta}g_{k,y}^{(n+1)}|\leq3^{-1}\varepsilon_{0}
\]
for each $k=0,\cdots,K$ and for each $y\in A_{n+1,k}$. Note that,
according to Assertion 5 of Theorem \ref{Thm. Construction of transition  f.j.d.s from E0 and Markov semi-group, discrete parameters-1},
we have
\[
F_{0,\Delta,\cdots,n\Delta}\equiv F_{0,\Delta,\cdots,n\Delta}^{x,\mathbf{T}}=(^{1}T_{\Delta}^{x})(^{2}T_{\Delta)})\cdots(^{n}T_{\Delta})
\]
where the factors on the right hand side are one-step transition distributions
according to $T_{\Delta}$, as defined in Lemma \ref{Lem.One-step transition distributions}.
A similar equality holds for $\overline{F}_{0,\Delta,\cdots,n\Delta}$
and $\overline{\mathbf{T}}$. Therefore it suffices to prove that
\begin{equation}
|(T_{\Delta}^{x})(^{2}T_{\Delta})\cdots(^{n}T_{\Delta})g_{k,y}^{(n+1)}-(^{1}\overline{T}_{\Delta}^{\overline{x}})(^{2}\overline{T}_{\Delta)})\cdots(^{n}\overline{T}_{\Delta})g_{k,y}^{(n+1)}|\leq3^{-1}\varepsilon_{0}\label{eq:temp-570}
\end{equation}
for each $k=0,\cdots,K$ and for each $y\in A_{n+1,k}$. 

6. To that end, let $p=0,\cdots,n$ be arbitrary. Define the finite
family
\[
G_{p}\equiv\{(^{p+1}\overline{T}_{\Delta)})\cdots(^{n}\overline{T}_{\Delta})g_{k,y}^{(n+1)}:k=0,\cdots,K\;\mathrm{and}\;y\in A_{n+1,k}\}\subset C(S^{p+1},d^{p+1}),
\]
with the understanding that the composite mapping $(^{p+1}\overline{T}_{\Delta)})\cdots(^{n}\overline{T}_{\Delta})$
stands for the identity mapping if $p=n$. Thus 
\[
G_{n}=\{g_{k,y}^{(n+1)}:k=0,\cdots,K\;\mathrm{and}\;y\in A_{n+1,k}\}\subset C(S^{n+1},d^{n+1}).
\]

7. We will prove by backward induction, that, for each $p=n,n-1,\cdots,0$,
there exists an operation $\rho_{p}:(0,\infty)\rightarrow(0,\infty)$
and an operation $\widetilde{\delta}_{p}:(0,\infty)\rightarrow(0,\infty)$
such that (i) $\widetilde{\delta}_{p}$ is a modulus of continuity
of each of the members of $G_{p}$, and (ii) if $p<n$, then for each
$\varepsilon'>0$, we have 
\[
\bigvee_{g\in G(p)}\left\Vert (^{p}T_{\Delta})g-(^{p}\overline{T}_{\Delta})g\right\Vert <\varepsilon'
\]
provided that the distance $\overline{\rho}_{0}$ is bounded by
\[
\overline{\rho}_{0}<\rho_{p}(\varepsilon').
\]

8. Start with $p=n$. Arbitrarily define $\rho_{n}\equiv1$. Then
the operation $\rho_{n}$ trivially satisfies Condition (ii) in Step
7. Consider each $g\in G_{n}.$ Then $g=g_{k,y}^{(n+1)}$ for some
$k=0,\cdots,K$ and $y\in A_{n+1,k}$. According to Proposition \ref{Prop. Properties of  epsilon partition of unity },
the basis function $g_{k,y}^{(n+1)}\in C(S^{n+1},d^{n+1})$ has Lipschitz
constant $2\cdot2^{k}\leq2^{K+1}=2^{N+2}$, and has values in $[0,1]$.
Thus the function $g$ has modulus of continuity $\widetilde{\delta}_{n}$
defined by $\widetilde{\delta}_{n}(\varepsilon')\equiv2^{-N-2}\varepsilon'$
for each $\varepsilon'>0$. Since $g\in G_{n}$ is arbitrary, the
modulus $\widetilde{\delta}_{n}$ satisfies Condition (i) in Step
7. The pair $\widetilde{\delta}_{n},\rho_{n}$ has been constructed
to satisfy Conditions (i) and (ii), in the case where $p=n$.

9. Suppose, for some $p=n,n-1,\cdots,1$, the pair of operations $\widetilde{\delta}_{p},\rho_{p}$
has been constructed to satisfiy the Conditions (i) and (ii) in Step
7. We proceed to construct the pari $\widetilde{\delta}_{p-1},\rho_{p-1}$
. Let $\varepsilon'>0$ be arbitrary. Consider each 
\[
g\in G_{p-1}\equiv\{(^{p}\overline{T}_{\Delta)})\cdots(^{n}\overline{T}_{\Delta})g_{k,y}^{(n+1)}:k=0,\cdots,K\;\mathrm{and}\;y\in A_{n+1,k}\}
\]
\[
=(^{p}\overline{T}_{\Delta})G_{p}\subset C(S^{p},d^{p})
\]
Then $g=(^{p}\overline{T}_{\Delta})\overline{g}$ for some $\overline{g}\in G_{p}$.
By Condition (ii) in the backward induction hypothesis, the function
$\overline{g}$ has modulus of continuity $\widetilde{\delta}_{p}$.
Lemma \ref{Lem.One-step transition distributions} therefore says
that the function $g=(^{p}\overline{T}_{\Delta})\overline{g}$ has
modulus of continuity given by
\[
\widetilde{\delta}_{p-1}\equiv\widetilde{\alpha}_{\alpha(\overline{T}(\Delta)),\xi}(\widetilde{\delta}_{p}),
\]
where the modulus of smoothness $\widetilde{\alpha}_{\alpha(\overline{T}(\Delta)),\xi}$
is as defined in Lemma \ref{Lem.One-step transition distributions}.
Since $g\in G_{p-1}$ is arbitrary, we see that $\widetilde{\delta}_{p-1}$
satisfies Condition (i) in Step 7. It remains to construct $\rho_{p-1}$to
satisfy Condition (ii) in Step 7.

10. To that end, let $\varepsilon'>0$ be arbitrary. Take $h\geq1$
so large that 
\begin{equation}
2^{-h}<\frac{1}{2}\widetilde{\delta}_{p}(\frac{1}{3}\varepsilon').\label{eq:temp-222}
\end{equation}
Define 
\begin{equation}
\rho_{p-1}(\varepsilon')\equiv\rho_{p}(\varepsilon')\wedge2^{-M-1}2^{-h}|A_{h}|^{-1}\varepsilon'.\label{eq:temp-418}
\end{equation}
Consider each $g\in G_{p-1}$. By Step 9 above, the function $g$
has modulus of continuity $\widetilde{\delta}_{p-1}$. Let $(w_{1},\cdots,w_{p-1})\in S^{p-1}$
be arbitrary, and consider the function $\overline{f}\equiv g(w_{1},\cdots,w_{p-1},\cdot)\in C(S,d)$.
Then clearly $\overline{f}$ has the same modulus of continuity $\widetilde{\delta}_{p-1}$. 

Recall that $\{g_{h,z}:z\in A_{h}\}$ is a $2^{-h}$-partition of
unity corresponding to $A_{h}$, where $A_{h}$ is a $2^{-h}$-approximation
of the compact metric space $(S,d)$. Hence, by Definition , we have\ref{Def. Binary approximationt and Modulus of local compactness}
\[
S\subset(d(\cdot,x_{\circ})\leq1)\subset(d(\cdot,x_{\circ})\leq2^{h})\subset\bigcup_{z\in A(h)}(d(\cdot,z)\leq2^{-h})
\]
Therefore, trivially, the function $\overline{f}\in C(S,d)$ has the
set $\bigcup_{z\in A(h)}(d(\cdot,z)\leq2^{-h})$ as support. Hence,
in view of inequality \ref{eq:temp-222}, the conditions in the hypothesis
of Proposition \ref{Prop. Approx  by  interpolation} are satsified.
Accordingly, 

\emph{
\begin{equation}
\left\Vert \overline{f}-\sum_{z\in A(h)}\overline{f}(z)g_{h,z}\right\Vert \leq\varepsilon'\label{eq:temp-151-1}
\end{equation}
}on $S$. Consequently,\emph{
\begin{equation}
\left\Vert T_{\Delta}\overline{f}-\sum_{z\in A(h)}\overline{f}(z)T_{\Delta}g_{h,z}\right\Vert \leq\varepsilon'\label{eq:temp-151-1-1}
\end{equation}
}and, similarly,\emph{
\begin{equation}
\left\Vert \overline{T}_{\Delta}\overline{f}-\sum_{z\in A(h)}\overline{f}(z)\overline{T}_{\Delta}g_{h,z}\right\Vert \leq\varepsilon'.\label{eq:temp-151-1-1-1}
\end{equation}
}

11. Now suppose $\overline{\rho}_{0}<\rho_{p-1}(\varepsilon').$ Then
\[
\sum_{j=0}^{\infty}2^{-j-1}\sum_{h'=1}^{\infty}2^{-h'}|A_{h'}|^{-1}\sum_{z\in A(h')}\left\Vert T_{\Delta(j)}g_{h',z}-\overline{T}{}_{\Delta(j)}g_{h',z}\right\Vert \equiv\rho_{\mathscr{T,\xi}}(\mathbf{T},\mathbf{\overline{T}})\equiv\overline{\rho}_{0}<\rho_{p-1}(\varepsilon').
\]
Consequently,
\[
2^{-M-1}2^{-h}|A_{h}|^{-1}\sum_{z\in A(h)}\left\Vert T_{\Delta(M)}g_{h,z}-\overline{T}{}_{\Delta(M)}g_{h,z}\right\Vert <\rho_{p-1}(\varepsilon').
\]
Recall that $\Delta_{M}\equiv2^{-M}\equiv\Delta.$ The last inequality
can be rewritten as
\[
2^{-M-1}2^{-h}|A_{h}|^{-1}\sum_{z\in A(h)}\left\Vert T_{\Delta}g_{h,z}-\overline{T}{}_{\Delta}g_{h,z}\right\Vert <\rho_{p-1}(\varepsilon').
\]
Therefore, since the function $\overline{f}$ has values in $[0,1]$,
we have 
\[
\left\Vert \sum_{z\in A(h)}\overline{f}(z)T_{\Delta}g_{h,z}-\sum_{z\in A(h)}\overline{f}(z)\overline{T}{}_{\Delta}g_{h,z}\right\Vert 
\]
\[
\leq\sum_{z\in A(h)}\left\Vert \overline{f}(z)T_{\Delta}g_{h,z}-\overline{f}(z)\overline{T}{}_{\Delta}g_{h,z}\right\Vert \leq\sum_{z\in A(h)}\left\Vert T_{\Delta}g_{h,z}-\overline{T}{}_{\Delta}g_{h,z}\right\Vert 
\]
\begin{equation}
<2^{M+1}2^{h}|A_{h}|\rho_{p-1}(\varepsilon')\leq\varepsilon',\label{eq:temp-543}
\end{equation}
where the last inequality follows from the defining formula \ref{eq:temp-418}. 

12. Combining inequalities \ref{eq:temp-543}, \ref{eq:temp-151-1-1-1},
and \ref{eq:temp-151-1-1}, we obtain, by the triangle inequality,
\[
\left\Vert T_{\Delta}\overline{f}-\overline{T}_{\Delta}\overline{f}\right\Vert \leq3\varepsilon'.
\]
In other words,
\[
\left\Vert T_{\Delta}g(w_{1},\cdots,w_{p-1},\cdot)-\overline{T}_{\Delta}g(w_{1},\cdots,w_{p-1},\cdot)\right\Vert \leq3\varepsilon',
\]
where $(w_{1},\cdots,w_{p-1})\in S^{p}$is arbitrary. Consequently
\[
\left\Vert ^{p-1}T_{\Delta}g-{}^{p-1}\overline{T}_{\Delta}g\right\Vert \leq3\varepsilon'.
\]
where $\varepsilon'>0$ and $g\in G_{p-1}$ are arbitrary, provided
that $\overline{\rho}_{0}<\rho_{p-1}(\varepsilon').$ 

13. In short, the operation $\rho_{p-1}$ has been constructed to
satisfy Condtion (ii) in Step 7. The backward induction is completed,
and we have obtained the pair ${\normalcolor (\widetilde{\delta}_{p},\rho_{p})}$
for each $p=n,n-1,\cdots,0$ to satisfy Conditions (i) and (ii) in
Step 7. In particular, we obtained the pair ${\normalcolor (\widetilde{\delta}_{0},\rho_{0})}.$
of operations.

14. Now let $\varepsilon'\equiv3^{-1}(n+2)^{-1}\varepsilon_{0}$.
Suppose 
\[
\overline{\rho}_{0}<\delta_{Sg,fjd}(\varepsilon_{0},\overline{\delta},\overline{\alpha},\left\Vert \xi\right\Vert )\equiv\rho_{0}(\varepsilon')\wedge\widetilde{\delta}_{0}(\varepsilon').
\]
Let $p=1,\cdots n$ and $g\in G_{p,}$ be arbitrary. Then
\begin{equation}
(^{p}\overline{T}_{\Delta})g\in G_{p-1},\label{eq:temp-414}
\end{equation}
Moreover, as a result of the built-in monotonicity of the operations
$\rho_{p}$ in the defining formula\ref{eq:temp-418}, we have $\overline{\rho}_{0}<\rho_{p}(\varepsilon').$
Hence, by Condition (ii) in Step 7, we have 
\[
(^{p}T_{\Delta})g=(^{p}\overline{T}_{\Delta})g\pm\varepsilon'.
\]
Therefore, since the transition distributions $(^{1}T_{\Delta}),\cdots,(^{p-1}T_{\Delta})$
are contraction mappings, it follows that
\[
(^{1}T_{\Delta})(^{2}T_{\Delta})\cdots(^{p-2}T_{\Delta})(^{p-1}T_{\Delta})(^{p}T_{\Delta})g
\]
\begin{equation}
=(^{1}T_{\Delta})(^{2}T_{\Delta})\cdots(^{p-2}T_{\Delta})(^{p-1}T_{\Delta})(^{p}\overline{T}_{\Delta})g\pm\varepsilon'\label{eq:temp-329}
\end{equation}

15. Finally, let $k=0,\cdots,K$ and $y\in A_{n+1,k}$. be arbitrary.
Let $g\equiv g_{k,y}^{(n+1)}$ $\in G_{n}$. Hence, using equality
\ref{eq:temp-329} repeatedly, we obtain

\[
F_{0,\Delta,\cdots,n\Delta}^{*,\mathbf{T}}g=(^{1}T_{\Delta})(^{2}T_{\Delta})\cdots(^{n-1}T_{\Delta})(^{n}T_{\Delta})g
\]
\[
=(^{1}T_{\Delta})(^{2}T_{\Delta})\cdots(^{n-1}T_{\Delta})(^{n}T_{\Delta})(^{n+1}\overline{T}_{\Delta})g\pm\varepsilon'
\]
\[
=(^{1}T_{\Delta})(^{2}T_{\Delta})\cdots(^{n-1}T_{\Delta})(^{n}\overline{T}_{\Delta})(^{n+1}\overline{T}_{\Delta})g\pm2\varepsilon'
\]
\[
=\cdots
\]
\[
=(^{1}\overline{T}_{\Delta})(^{2}\overline{T}_{\Delta})\cdots(^{n-1}T_{\Delta})(^{n}\overline{T}_{\Delta})(^{n+1}\overline{T}_{\Delta})g\pm(n+1)\varepsilon'
\]
\begin{equation}
=F_{0,\Delta,\cdots,n\Delta}^{*,\mathbf{\overline{T}}}g\pm(n+1)\varepsilon'.\label{eq:temp-573}
\end{equation}

16. Moreover, by Condition (i) in Step 7, the function 
\[
F_{0,\Delta,\cdots,n\Delta}^{*,\mathbf{T}}g=(^{1}T_{\Delta})(^{2}T_{\Delta})\cdots(^{n-1}T_{\Delta})(^{n}T_{\Delta})g\in G_{0}
\]
has modulus of continuity $\widetilde{\delta}_{0}$. Therefore
\[
F_{0,\Delta,\cdots,n\Delta}^{x,\mathbf{T}}g=F_{0,\Delta,\cdots,n\Delta}^{\overline{x},\mathbf{T}}g\pm\varepsilon'
\]
because 
\[
d(x,\overline{x})\leq\overline{\rho}_{0}<\delta_{Sg,fjd}(\varepsilon_{0},\overline{\delta},\overline{\alpha},\left\Vert \xi\right\Vert )\equiv\rho_{0}(\varepsilon')\wedge\widetilde{\delta}_{0}(\varepsilon')\leq\widetilde{\delta}_{0}(\varepsilon').
\]
By symmetry,
\begin{equation}
F_{0,\Delta,\cdots,n\Delta}^{x,\mathbf{\overline{T}}}g=F_{0,\Delta,\cdots,n\Delta}^{\overline{x},\mathbf{\overline{T}}}g\pm\varepsilon'.\label{eq:temp-574}
\end{equation}
Equalities \ref{eq:temp-573} and \ref{eq:temp-574} together imply
\[
F_{0,\Delta,\cdots,n\Delta}^{x,\mathbf{T}}g=F_{0,\Delta,\cdots,n\Delta}^{x,\mathbf{\overline{T}}}g\pm(n+1)\varepsilon'=F_{0,\Delta,\cdots,n\Delta}^{\overline{x},\mathbf{\overline{T}}}g\pm(n+1)\varepsilon'\pm\varepsilon'
\]
\[
=F_{0,\Delta,\cdots,n\Delta}^{\overline{x},\mathbf{\overline{T}}}g\pm(n+2)\varepsilon'\equiv F_{0,\Delta,\cdots,n\Delta}^{\overline{x},\mathbf{\overline{T}}}g\pm3^{-1}\varepsilon_{0},
\]
Equivalently,
\begin{equation}
(T_{\Delta}^{x})(^{2}T_{\Delta})\cdots(^{n}T_{\Delta})g_{k,y}^{(n+1)}=(^{1}\overline{T}_{\Delta}^{\overline{x}})(^{2}\overline{T}_{\Delta)})\cdots(^{n}\overline{T}_{\Delta})g_{k,y}^{(n+1)}\pm3^{-1}\varepsilon_{0}\label{eq:temp-570-2}
\end{equation}
where $k=0,\cdots,K$ and $y\in A_{n+1,k}$. are arbitrary. The desired
equality \ref{eq:temp-570} followsfor each $n=0,\cdots,N$. In turn,
inequalities \ref{eq:temp-219-3-1} and \ref{eq:temp-208-2-2-1} then
imply that
\[
\widehat{\rho}_{Marg,\xi,\overline{Q}(\infty)}(\Phi_{Sg,fjd}(x,\mathbf{\mathbf{T}}\mathit{),}\Phi_{Sg,fjd}(\overline{x},\mathbf{\overbrace{\mathbf{T}}}\mathit{)})
\]
\[
=\widehat{\rho}_{Marg,\xi,\overline{Q}(\infty)}(F,\overline{F})\leq3^{-1}\varepsilon_{0}+3^{-1}\varepsilon_{0}+3^{-1}\varepsilon_{0}=\varepsilon_{0},
\]
provide that the distance $\overline{\rho}_{0}\equiv(d\otimes\rho_{\mathscr{T}})((x,\mathbf{T}),(\overline{x},\overline{\mathbf{T}}))$
is bounded by
\[
\overline{\rho}_{0}<\delta_{Sg,fjd}(\varepsilon_{0},\overline{\delta},\overline{\alpha},\left\Vert \xi\right\Vert ).
\]
Summing up, $\delta_{Sg,fjd}(\cdot,\overline{\alpha},\left\Vert \xi\right\Vert )$
is a modulus of continuity of the mapping $\Phi_{Sg,fjd}$. The theorem
is proved.
\end{proof}
Following is the main theorem of the section. It proves the continuity
of construction in the case where the parameter set is $[0,\infty)$.
\begin{thm}
\textbf{\emph{\label{Thm. Construction of Markov process from  initial distrinbution and semigrp-1-1}}}
\textbf{\emph{(Construction of time-uniformly a.u. càdlàg Markov process
from an initial state and a Markov semigroup on $[0,\infty)$, and
continuity of said construction).}} Let $(S,d)$ be the specified
compact metric space, with $d\leq1$.\textbf{\emph{ }}Let $\mathscr{T}(\overline{\delta},\overline{\alpha})$
be an arbitrary family of  Markov semigroups with parameter set $]0,\infty)$
and state space $(S,d)$, such that\emph{ }all its members $\mathbf{T}\in\mathscr{T}(\overline{\delta},\overline{\alpha})$
share a common modulus of strong continuity $\delta_{\mathbf{T}}=\overline{\delta}$,
and share a common modulus of smoothness $\alpha_{\mathbf{T}}=\overline{\alpha}$,
in the sense of Definition \ref{Def. Markov semigroup  semigroup}.
Thus $\mathscr{T}(\overline{\delta},\overline{\alpha})$ is a subset
of $(\mathscr{T},\rho_{\mathscr{T,\xi}}),$ and, as such, inherits
the metric $\rho_{\mathscr{T,\xi}}$ introduced in Definition \ref{Def. Metric space of Markov semigroups}.
Recall from \ref{Def. Continuous in probabillity of Finite Joint Distributions}.
the space $\widehat{F}_{Cp}([0,\infty),S)$ of consistent families
of f.j.d.'s which are continuout in probability, equipped with the
metric $\widehat{\rho}_{Cp,\xi,[0,\infty)|\overline{Q}(\infty)}$
introduced in \textup{\emph{Definition \ref{Def. Metric on  of continuous in prob families of finite joint distributions.}}}.
Then the following holds.

\emph{1.} There exists a uniformly continuous mapping 
\[
\overline{\Phi}_{Sg,fjd}:(S\times\mathscr{T}(\overline{\delta},\overline{\alpha}),d\otimes\rho_{\mathscr{T,\xi}})\rightarrow(\widehat{F}_{Cp}([0,\infty),S),\widehat{\rho}_{Cp,\xi,[0,\infty)|\overline{Q}(\infty)})
\]
such that, for each $(x,\mathbf{T})\in S\times\mathscr{T}(\overline{\delta},\overline{\alpha})$,
the family $\overline{F}\equiv\overline{\Phi}_{Sg,fjd}(x,\mathbf{T})$
of f.j.d.'s is generated by the initial state $x$ and the semigroup
$\mathbf{T}$, in the sense of Definition \textup{\ref{Def. Family of f.j.d.s from E0 and Markov semi-group}}
Roughly speaking, we can generate, from an initial state and a semigroup
the corresponding f.j.d.'s, and the generation is continuous. We will
write $F^{x,\mathbf{T}}\equiv\overline{\Phi}_{Sg,fjd}(x,\mathbf{T})$
and call it \emph{the family of f.j.d.'s generated by the initial
state and the semigroup\index{f.j.d.'s generated by semigroup}} $\mathbf{\mathbf{T}}$.

In particular, for each $\mathbf{T}\in\mathscr{\mathscr{T}(\overline{\delta},\overline{\alpha})}$,
the function 
\[
F^{*,\mathbf{T}}\equiv\overline{\Phi}_{Sg,fjd}(\cdot,\mathbf{T}):(S,d)\rightarrow(\widehat{F}_{Cp}([0,\infty),S),\widehat{\rho}_{Cp,\xi,[0,\infty)|\overline{Q}(\infty)})
\]
is uniformly continuous. Consequently, for each $m\geq1$, $f\in C(S^{m},d^{m})$,
and nondecreasing sequence $r_{1}\leq\cdots\leq r_{m}$ in $[0,\infty)$,
the function 
\[
F_{r(1),\cdots,r(m)}^{*,\mathbf{T}}f:(S,d)\rightarrow R
\]
is uniformly continuous on $(S,d)$.

\emph{2. }There exists a uniformly continuous mapping\emph{
\[
\Phi_{Sg,clMk}:(S\times\mathscr{T}(\overline{\delta},\overline{\alpha}),d\otimes\rho_{\mathscr{T,\xi}})\rightarrow(\widehat{D}[0,\infty),\widetilde{\rho}_{\widehat{D}[0,\infty)}),
\]
}where $(\widehat{D}[0,\infty),\widetilde{\rho}_{\widehat{D}[0,\infty)})$
is the metric space of a.u. càdlàg process with some sample space
$(\Omega,L,E)$ and with parameter set $[0,\infty),$ as defined in
Definition \ref{Def. Metric space of a.u. cadlag process on =00005B0,inf)}.
Specifically
\[
\Phi_{Sg,clMk}\equiv\overline{\Phi}_{rLim}\circ\overline{\Phi}_{DKS,\xi}\circ\Phi_{Sg,fjd}\circ\Psi,
\]
where the component mappings on the right-hand side have been previously
defined and will be recalled in detail in the proof. Roughly speaking,
we can construct, from an initial state and a semigroup, a corresponding
a.u. càdlàg process $X^{x,\mathbf{T}}\equiv\Phi_{Sg,clMk}(x,\mathbf{T})$,
and prove that the construction is continuous.

Moreover, $X^{x,\mathbf{T}}$ has  a modulus of a.u. càdlàg $\widetilde{\delta}_{aucl,\overline{\delta}}\equiv(\delta_{aucl,\overline{\delta}},\delta_{aucl,\overline{\delta}},\cdots)$
and has a modulus of continuity in probability $\widetilde{\delta}{}_{Cp,\overline{\delta}}\equiv(\delta_{Cp,\overline{\delta}},\delta_{Cp,\overline{\delta}},\cdots)$.
In other words, $X^{x,\mathbf{T}}\in\widehat{D}_{\widetilde{\delta}(aucl,\overline{\delta}),\widetilde{\delta}(Cp,\overline{\delta})}[0,\infty).$
Hence the function $\Phi_{Sg,clMk}$ has range in $\widehat{D}_{\widetilde{\delta}(aucl,\overline{\delta}),\widetilde{\delta}(Cp,\overline{\delta})}[0,\infty)$,
and we can regard it as a uniformly continuous mapping\emph{
\[
\Phi_{Sg,clMk}:(S\times\mathscr{T}(\overline{\delta},\overline{\alpha}),d\otimes\rho_{\mathscr{T,\xi}})\rightarrow\widehat{D}_{\widetilde{\delta}(aucl,\overline{\delta}),\widetilde{\delta}(Cp,\overline{\delta})}[0,\infty).
\]
}

\textup{\emph{3. }}Let\emph{ $(x,\mathbf{T})\in S\times\mathscr{T}(\overline{\delta},\overline{\alpha})$
}be arbitrary. Then \textup{\emph{the process
\[
X^{x,\mathbf{T}}\equiv\Phi_{Sg,clMk}(x,\mathbf{T}):[0,\infty)\times(\Omega,L,E)\rightarrow(S,d),
\]
 constructed in Assertion }}\textup{2}\textup{\emph{ above, is Markov
relative to its natural filtration}}\emph{ $\mathcal{\overline{L}}$.}
More precisely, for each $t\geq0$, define 
\[
\overline{L}^{(t)}\equiv L(X_{s}^{x,\mathbf{T}}:s\in[0,t]\}
\]
and define $\mathcal{\overline{L}}\equiv\{\overline{L}^{(t)}:t\geq0$\}.
Let the nondecreasing sequence $0\equiv s_{0}\leq s_{1}\leq\cdots\leq s_{m}$
in $[0,\infty)$, the function $f\in C(S^{m+1},d^{m+1})$, and $t\geq0$
be arbitrary. Then 
\[
E(f(X_{t+s(0)},X_{t+s(1)},\cdots,X_{t+s(m)})|L^{(t)})
\]
\begin{equation}
=E(f(X_{t+s(0)},X_{t+s(1)},\cdots,X_{t+s(m)})|X_{t})=F_{s(0),\cdots,s(m)}^{X(t),\mathbf{T}}(f)\label{eq:temp-575-2}
\end{equation}
as r.r.v.'s. 
\end{thm}
\begin{proof}
We will prove Assertion 2 first, then Assertions 3 and 1. To that
end, let 
\[
(\Omega,L,E)\equiv(\Theta_{0},L_{0},I_{0})\equiv([0,1],L_{0},\int\cdot dx)
\]
denote the Lebesgue integration space based on the unit interval $\Theta_{0}$.

1. Let $\mathbf{T}\in\mathscr{\mathscr{\mathscr{T}(\overline{\delta}\mathit{,}\overline{\alpha})}}$
be arbitrary. Then, since $\mathbf{T}$ is a Markov semigroup with
parameter set $[0,\infty)$, it is trivial to verify the conditions
in Definition \ref{Def. Markov semigroup  semigroup} for its restriction
$\mathbf{T|}\overline{Q}_{\infty}$ to be a Markov semigroup with
parameter set $\overline{Q}_{\infty}$, with the same modulus of strong
continuity $\delta_{\mathbf{T}}=\overline{\delta}$ and modulus of
smoothness  $\alpha_{\mathbf{T}}=\overline{\alpha}$ . 

2. Define the set $\mathscr{\mathscr{T}(\overline{\delta},\overline{\alpha})}|\overline{Q}_{\infty}\equiv\{\mathbf{T|}\overline{Q}_{\infty}:\mathbf{T}\in\mathbf{\mathscr{\mathscr{T}(\overline{\delta}\mathit{,}\overline{\alpha})}}\},$
and equip it with the metric $\rho_{\mathscr{T,\xi}}$ introduced
in Definition \ref{Def. Metric space of Markov semigroups}. As observed
in Step 1 above, members of the family $\mathscr{\mathscr{T}(\overline{\delta}\mathit{,}\overline{\alpha})}|\overline{Q}_{\infty}$
share the the common $\overline{\delta}$ and common$\overline{\alpha}$.As
observed in Definition \ref{Def. Metric space of Markov semigroups},
the mapping
\[
\Psi:\mathscr{(\mathscr{\mathscr{T}(\overline{\delta}\mathit{,}\overline{\alpha})}},\rho_{\mathscr{T,\xi}})\rightarrow\mathscr{(\mathscr{\mathscr{T}(\overline{\delta}\mathit{,}\overline{\alpha})}}|\overline{Q}_{\infty},\rho_{\mathscr{T,\xi}}),
\]
defined by $\mathit{\mathrm{\Psi(}}\mathbf{T})\equiv\mathbf{T|}\overline{Q}_{\infty}$
for each $\mathbf{T}\in\mathscr{\mathscr{\mathscr{T}(\overline{\delta}\mathit{,}\overline{\alpha})}}$,
is trivally an isometry. Hence the mapping 
\[
\overline{\Psi}:(S\times\mathscr{\mathscr{T}(\overline{\delta}\mathit{,}\overline{\alpha})},d\otimes\rho_{\mathscr{T,\xi}})\rightarrow(S\times\mathscr{\mathscr{T}(\overline{\delta}\mathit{,}\overline{\alpha})}|\overline{Q}_{\infty},d\otimes\rho_{\mathscr{T,\xi}}),
\]
defined by $\overline{\Psi}(x,\mathbf{T})\equiv(x,\mathit{\mathrm{\Psi(}}\mathbf{T}))$
for each $(x,\mathbf{T})\in S\times\mathscr{\mathscr{\mathscr{T}(\overline{\delta}\mathit{,}\overline{\alpha})}}$,
is, in turn, trivially uniformly continuous.

3. Theorem \ref{Thm.  Continuity of construction of f.j.d.s from E0 and Markov semigroups-1-1}
herefore says that the mapping 
\[
\Phi_{Sg,fjd}:(S\times\mathscr{\mathscr{T}(\overline{\delta}\mathit{,}\overline{\alpha})}|\overline{Q}_{\infty},d\otimes\rho_{\mathscr{T,\xi}})\rightarrow(\widehat{F}(\overline{Q}_{\infty},S),\widehat{\rho}_{Marg,\xi,\overline{Q}(\infty)}),
\]
constructed in Assertion 4 of Theorem \ref{Thm. Construction of transition  f.j.d.s from E0 and Markov semi-group, discrete parameters-1},
is uniformly continuous, with a modulus of continuity $\delta_{Sg,fjd}(\cdot,\overline{\delta},\overline{\alpha},\left\Vert \xi\right\Vert )$
completely determined by the moduli $\overline{\delta},\overline{\alpha}$
and the modulus of local compactness $\left\Vert \xi\right\Vert \equiv(|A_{n}|)_{n=1,2,\cdots}$of
$(S,d)$. 

4. Separately, Theorem \ref{Thm. Continuity of Compact Daniell-Kolmogorov-Skorohod Extension}
says that the Compact Daniell-Kolmogorov-Skorokhod Extension 
\[
\overline{\Phi}_{DKS,\xi}:(\widehat{F}(\overline{Q}_{\infty},S),\widehat{\rho}_{Marg,\xi,\overline{Q}(\infty)})\rightarrow(\widehat{R}(\overline{Q}_{\infty}\times\Theta_{0},S),\widehat{\rho}_{Prob,\overline{Q}(\infty)})
\]
is uniformly continuous with a modulus of continuity $\overline{\delta}_{DKS}(\cdot,\left\Vert \xi\right\Vert )$
dependent only on $\left\Vert \xi\right\Vert $. 

5. Combining, we see that the composite mapping $\overline{\Phi}_{DKS,\xi}\circ\Phi_{Sg,fjd}\circ\overline{\Psi}$
is uniformly continuous. Now consider the range of this composite
mapping. Specifically, take an arbitrary $(x,\mathbf{T})\in S\times\mathscr{T}$
and consider the image 
\[
Z\equiv\overline{\Phi}_{DKS,\xi}(\Phi_{Sg,fjd}(\overline{\Psi}(x,\mathbf{T}))).
\]
where 
\[
F\equiv\Phi_{Sg,fjd}(x,\mathbf{T}|\overline{Q}_{\infty})
\]
is the consistent family constructed in Assertion 3 of Theorem \ref{Thm. Construction of transition  f.j.d.s from E0 and Markov semi-group, discrete parameters-1},
and is generated by the initial state $x$ and the semigroup \emph{$\mathbf{\mathbf{T}|\overline{Q}_{\infty}}$,
}in the sense of Definition\emph{ \ref{Def. Family of f.j.d.s from E0 and Markov semi-group}.
}In the notations of Theorem \ref{Thm. Construction of transition  f.j.d.s from E0 and Markov semi-group, discrete parameters-1},
we have $F\equiv F^{x,\mathbf{T}|\overline{Q}(\infty)}$. Thus $Z\equiv\overline{\Phi}_{DKS,\xi}(F^{x,\mathbf{T}|\overline{Q}(\infty)})$
and, by the definition of the mapping $\overline{\Phi}_{DKS,\xi}$,
the process $Z$ has marginal distributions given by the faamily $F^{x,\mathbf{T}|\overline{Q}(\infty)}$.
In particular $Z_{0}=x$.

6. Since the semigroup $\mathbf{\mathbf{T}|\overline{Q}_{\infty}}$
has modulus of strong continuuity $\overline{\delta}$, Assertion
1 of Theorem \ref{Thm. Markov process from initial distribution and Markov semigroup}
implies that the process $Z\equiv\overline{\Phi}_{DKS,\xi}(F^{x,\mathbf{T}|\overline{Q}(\infty)})$
is generated by the initial state $x$ and semigroup $\mathbf{\mathbf{T}|\overline{Q}_{\infty}}$,
in the sense of Definition \ref{Def. Family of f.j.d.s from E0 and Markov semi-group}.
Hence, by Assertions 4 of Proposition \ref{Prop. Process with Markov semigroup on dyadic rationals in =00005B0,oo) is D-regular},
the process $Z$ is time-uniformly $D$-regular in the sense of Definition
\ref{Def. D-regular processes wirh dyadic rationals in =00005B0,oo)as parameters-1},
with some modulus of continuity in probability $\widetilde{\delta}{}_{Cp,\overline{\delta}}$
$\equiv(\delta_{Cp,\overline{\delta}},\delta_{Cp,\overline{\delta}},\cdots)$
and with some modulus of $D$-regularity $\widetilde{m}_{\overline{\delta}}\equiv(\overline{m}_{\overline{\delta}},\overline{m}_{\overline{\delta}},\cdots)$.
In the notations of we have $Z\in\widehat{R}_{Dreg,\widetilde{\delta}(Cp,\overline{\delta}),\widetilde{m}(\overline{\delta)}.}(\overline{Q}_{\infty}\times\Omega,S)$.
Summing up, we see that the range of the composite mapping $\overline{\Phi}_{DKS,\xi}\circ\Phi_{Sg,fjd}\circ\overline{\Psi}$
is contained in the subset $\widehat{R}_{Dreg,\widetilde{\delta}(Cp,\overline{\delta}),\widetilde{m}(\overline{\delta)}.}(\overline{Q}_{\infty}\times\Omega,S)$
of $(\widehat{R}(\overline{Q}_{\infty}\times\Theta_{0},S),\widehat{\rho}_{Prob,\overline{Q}(\infty)})$.
Thus we have the uniformly continuous mapping 
\[
\overline{\Phi}_{DKS,\xi}\circ\Phi_{Sg,fjd}\circ\overline{\Psi}:(S\times\mathscr{\mathscr{T}(\overline{\delta}\mathit{,}\overline{\alpha})},d\otimes\rho_{\mathscr{T,\xi}})\rightarrow(\widehat{R}_{Dreg,\widetilde{\delta}(Cp,\overline{\delta}),\widetilde{m}(\overline{\delta)}.}(\overline{Q}_{\infty}\times\Omega,S),\widehat{\rho}_{Prob,\overline{Q}(\infty)}).
\]
In the notations of Theorem \ref{Thm. Markov process from initial distribution and Markov semigroup},
we defined $Z^{x,\mathbf{T}|\overline{Q}(\infty)}\equiv\overline{\Phi}_{DKS,\xi}(.F^{x,\mathbf{T}|\overline{Q}(\infty)})$.
Hence $Z=Z^{x,\mathbf{T}|\overline{Q}(\infty)}$. Therefore, by the
definiton of $\overline{\Phi}_{DKS,\xi}$, the process $Z$ has marginal
distributions given by the family $F^{x,\mathbf{T}|\overline{Q}(\infty)}$.

7. Moreover, by Assertions 5 Proposition \ref{Prop. Process with Markov semigroup on dyadic rationals in =00005B0,oo) is D-regular},
the right-limit extension
\[
X\equiv\overline{\Phi}_{rLim}(Z):[0,\infty)\times(\Omega,L,E)\rightarrow(S,d)
\]
is a well-defined time-uniformly a.u. càdlàg process on $[0,\infty)$,
in the sense of Definition \ref{Def. Metric space of a.u. cadlag process on =00005B0,inf)},
with some modulus of a.u. càdlàg $\widetilde{\delta}_{aucl,\overline{\delta}}\equiv(\delta_{aucl,\overline{\delta}},\delta_{aucl,\overline{\delta}},\cdots)$
and with the modulus of continuity in probability $\widetilde{\delta}{}_{Cp,\delta}$.
In short $X\equiv\overline{\Phi}_{rLim}(Z)\in\widehat{D}[0,\infty).$

8. Since $(x,\mathbf{T})\in S\times\mathscr{T}$ is arbitrary, the
composite mappping $\overline{\Phi}_{rLim}\circ\overline{\Phi}_{DKS,\xi}\circ\Phi_{Sg,fjd}\circ\overline{\Psi}$
is well defined. We have already seen that the mapping $\overline{\Phi}_{DKS,\xi}\circ\Phi_{Sg,fjd}\circ\overline{\Psi}$
is uniformly continuous. In addition, Theorem \ref{Thm. Phi^bar_rLim is isometry on properly restricted domain}
says that 
\[
\overline{\Phi}_{rLim}:(\widehat{R}_{Dreg,\widetilde{\delta}(Cp,\overline{\delta}),\widetilde{m}(\overline{\delta)}.}(\overline{Q}_{\infty}\times\Omega,S),\widehat{\rho}_{Prob,\overline{Q}(\infty)})\rightarrow(\widehat{D}[0,\infty),\widetilde{\rho}_{\widehat{D}[0,\infty)})
\]
is a well-defined isometry. Summing up, we see that the composite
construction mapping 
\[
\Phi_{Sg,clMk}\equiv\overline{\Phi}_{rLim}\circ\overline{\Phi}_{DKS,\xi}\circ\Phi_{Sg,fjd}\circ\overline{\Psi}:(S\times\mathscr{T}(\overline{\delta},\overline{\alpha}),d\otimes\rho_{\mathscr{T,\xi}})\rightarrow(\widehat{D}[0,\infty),\widetilde{\rho}_{\widehat{D}[0,\infty)})
\]
 is well-defined and is uniformly continuous. Assertion 2 of the present
theorem has been proved.

9. Next, note that Assertion 3 of Theorem \ref{Thm. Markov process from initial distribution and Markov semigroup-1}
says that the process $X^{x,\mathbf{T}}\equiv\overline{\Phi}_{rLim}(Z^{x,\mathbf{T}})=X$
is Markov relative to its natural filtration $\overline{\mathcal{L}}\equiv\{\overline{L}^{(t)}:t\in[0,\infty)\}$,
and that, specifically, for each $v\geq0$ and for each sequence $t_{0}\equiv0\leq t_{1}\leq\cdots\leq t_{m}$
in $[0,\infty)$, we have
\begin{equation}
I_{0}(f(X_{v+t(0)},\cdots,X_{v+t(m)})|\overline{L}^{(v)})=I_{0}(f(X_{v+t(0)},\cdots,X_{v+t(m)})|X_{v})=F_{0,t(1)\cdots,t(m)}^{X(v),\mathbf{T}}f.\label{eq:temp-559-1-1}
\end{equation}
Assertion 3 of the present theorem is proved.. 

10. To prove Assertion 1, define $\overline{\Phi}_{Sg,fjd}(x,\mathbf{T})$
to be the family $\overline{F}$ of marginal distributions of the
process $X=X^{x,\mathbf{T}}\equiv\overline{\Phi}_{rLim}(Z^{x,\mathbf{T}})$.
In the special case where $v=0$, equality \ref{eq:temp-559-1-1}
then yields 
\begin{equation}
\overline{F}_{0,t(1)\cdots,t(m)}f\equiv Ef(X_{t(0)},\cdots,X_{t(m)})=E(Ef(X_{t(0)},\cdots,X_{t(m)})|\overline{L}^{(0)}))=EF_{0,t(1)\cdots,t(m)}^{X(0),\mathbf{T}}f=F_{0,t(1)\cdots,t(m)}^{x,\mathbf{T}}f,\label{eq:temp-291}
\end{equation}
 where the last equality holds because $X_{0}=Z_{0}=x$ a.s. If members
of the sequence $t_{0}\equiv0\leq t_{1}\leq\cdots\leq t_{m}$ are
in $\overline{Q}_{\infty}$, then Assertion 5 of Theorem \ref{Thm. Construction of transition  f.j.d.s from E0 and Markov semi-group, discrete parameters-1}
implies
\[
F_{0,t(1)\cdots,t(m)}^{x,\mathbf{T}}f
\]
\begin{equation}
=\int T_{t(1)}^{x}(dx_{1})\int dT_{t(2)-t(1)}^{x(1)}(dx_{2})\cdots\int dT_{t(m)-t(m-1)}^{x(m-1)}(dx_{m})f(x_{0},\cdots,x_{m}).\label{eq:temp-384-1-1-1-1-1}
\end{equation}
Hence, by continuity in the time parameters, this same equality holds
for an arbitrary sequence $t_{0}\equiv0\leq t_{1}\leq\cdots\leq t_{m}$
in $[0,\infty)$. Combining with equality \ref{eq:temp-291}, we obtain
\[
\overline{F}_{0,t(1)\cdots,t(m)}f=\int T_{t(1)}^{x}(dx_{1})\int dT_{t(2)-t(1)}^{x(1)}(dx_{2})\cdots\int dT_{t(m)-t(m-1)}^{x(m-1)}(dx_{m})f(x_{0},\cdots,x_{m}).
\]
Thus the family $\overline{\Phi}_{Sg,fjd}(x,\mathbf{T})\equiv\overline{F}$
of f.j.d.'s is generated by the initial state $x$ and the semigroup
$\mathbf{T}$, in the sense of Definition \ref{Def. Family of f.j.d.s from E0 and Markov semi-group}. 

Continuity of the function$\overline{\Phi}_{Sg,fjd}$ then follows
easily from the continuity of $\overline{\Phi}_{Sg,clMk}$ established
in Step 8. Assertion 1 and the theorem have been proved.
\end{proof}

\section{a.u. Càdlàg Markov Processes are Strongly Markov}
\begin{defn}
\textbf{\label{Def. Specification of compaact state space, semigroup, parameter sets etc}(Specification
of state space, parameter sets, Markov semigroup, and related objects).}
In this section, unless otherwise specified, let $(S,d)$ be a given
compact metric space, with $d\leq1$. Let $\mathbf{T}$ be denote
an arbitrary Markov semigroup with parameter set $[0,\infty)$ and
state space $(S,d)$, with a modulus of strong continuity $\delta_{\mathbf{T}}$
and a modulus of smoothness $\alpha_{\mathbf{T}}\equiv(\alpha_{\mathbf{T,\mathit{N}}})_{N=1,2,\cdots}$
For each $x\in S$, let $F^{x,\mathbf{\mathrm{\mathbf{T}}}}$ denote
the family of transition f.j.d.'s , with parameter set $[0,\infty)$,
generated by the initial state $x$ and semigroup $\mathbf{T},$ as
defined and constructed in Theorem \ref{Thm. Construction of Markov process from  initial distrinbution and semigrp-1-1}.
We will use results in the previous sections reegarding these objects,
without further comments.

Refer also to Definition \ref{Def.  Specification of state space}
for the a specified binary approximation $\xi$ of the compact metric
space $(S,d)$. Refer to Definition \ref{Def. Notations for dyadic rationals-1}
for notations related to the sets $\overline{Q}_{1}$,$\overline{Q}_{2}\cdots,\overline{Q}_{\infty}$
of dyadic rationals in $[0,\infty)$, and to the sets $Q_{1}$,$Q_{2}\cdots,Q_{\infty}$
of dyadic rationals in $[0,1]$.
\end{defn}
\begin{lem}
\textbf{\emph{\label{Lem. Short-term quiecience of a.u. cadlag Feller process after each stopping time with uniformly spaced dyadic values}
(Each a.u.}} \textbf{\emph{càdlàg Markov process}} \textbf{\emph{is
a.u. right continuous at each stopping time with regularly spaced
dyadic values). }}Let $x\in S$ be arbitrary.
\[
X:[0,\infty)\times(\Omega,L,E)\rightarrow(S,d)
\]
be an arbitrary a.u. càdlàg Markov process generated by the initial
state $x$ and with semigroup $\mathbf{T}$ which is adapted to some
filtration $\mathcal{L}\equiv\{L^{(t)}:t\in[0,\infty)\}$.\emph{ }Let
$k\geq0$ be arbitrary. Then there exists $m_{k}\equiv m_{k}(\delta_{\mathbf{T}})\geq0$
with the following properties. 

Let $\eta$ be an arbitrary stopping time with values in $\overline{Q}_{h}$
for some $h\geq0$, relative to the filtration $\mathcal{L}$. Then
the function 
\[
V_{\eta,k}\equiv\sup_{s\in[0,\Delta(m(k))]}d(X_{\eta},X_{\eta+s})
\]
is a well defined r.r.v., with 
\[
EV_{\eta,k}\leq2^{-k+1}.
\]
Recall here that $\Delta_{m(k)}\equiv2^{-m(k)}.$ We emphasize that
$m_{k}$ depends only on $\delta_{\mathbf{T}}$, and is independent
of $\eta,h$ or $the$ initial state $x$. Furthermore, by inductively
replacing $m_{k}$ with $m_{k}\vee(m_{k-1}+1)$ for each $k\geq1$,
we may assume that the sequence $(m_{k})_{k=0,1,\cdots}$ of integers
is increasing.
\end{lem}
\begin{proof}
1. Define $Z\equiv X|\overline{Q}_{\infty}$. Then $X=\overline{\Phi}_{rLim}(Z)$.
Let $M\geq0$ be arbitrary. Define the shifted process 
\[
X^{M}:[0,1]\times(\Omega,L,E)\rightarrow(S,d)
\]
by $X^{M}(t)\equiv X(M+t)$ for each $t\in[0,1]$. Similarly, define
the shifted process 
\[
Z^{M}:Q_{\infty}\times(\Omega,L,E)\rightarrow(S,d)
\]
by $Z^{M}(t)\equiv X(M+t)$ for each $t\in Q_{\infty}$. Then $X^{M}\equiv\Phi_{rLim}(Z^{M})$.

2. Thanks to Theorem \ref{Thm. Construction of Markov process from  initial distrinbution and semigrp-1-1}
the a.u. càdlàg Markov process $X$ is time-uniformly a.u. càdlàg
in the sense of Definition \ref{Def. Metric space of a.u. cadlag process on =00005B0,inf)},
with a modulus of continuity in probability $\widetilde{\delta}{}_{Cp,\delta(\mathscr{\mathrm{\mathbf{T}}})}\equiv(\delta_{Cp,\delta(\mathscr{\mathrm{\mathbf{T}}})},\delta_{Cp,\delta(\mathscr{\mathrm{\mathbf{T}}})},\cdots)$
and a modulus of a.u. càdlàg $\widetilde{\delta}_{aucl,\delta(\mathscr{\mathrm{\mathbf{T}}})}\equiv(\delta_{aucl,\delta(\mathscr{\mathrm{\mathbf{T}}})},\delta_{aucl,\delta(\mathscr{\mathrm{\mathbf{T}}})},\cdots)$.
In other words, the shifted process $X^{M}$ has a modulus of a.u.
càdlàg $\delta_{aucl,\delta(\mathscr{\mathrm{\mathbf{T}}})}$ and
a modulus of continuity in probability $\delta_{Cp,\delta(\mathbf{\mathbf{\mathscr{\mathscr{\mathrm{\mathbf{T}}}}}})}$
which are independent of $M$ and of the initial state $x$. In particular,
the process $X|[0,1]=X^{0}$ has a modulus of a.u. càdlàg $\delta_{aucl,\delta(\mathscr{\mathrm{\mathbf{T}}})}$.
Take an arbitrary increasing sequence $\overline{m}_{\delta(\mathscr{\mathrm{\mathbf{T}}})}\equiv(m_{k})_{k=0,1,2,\cdots}$
of integers such that 
\begin{equation}
\Delta_{m(k)}\equiv2^{-m(k)}<\delta_{aucl,\delta(\mathscr{\mathrm{\mathbf{T}}})}(2^{-k})\label{eq:temp-369-3}
\end{equation}
for each $k\geq0$. Then, by Theorem \ref{Thm. Restriction of a.u. cadlag process to dyadics is D-regular},
the process $Z|Q_{\infty}=Z^{0}$ has a modulus of $D$-regularity
given by the sequence $\overline{m}_{\delta(\mathscr{\mathrm{\mathbf{T}}})}$.
For convenience, write $m_{-1}\equiv m_{0}$. 

3. Let $k\geq0$ be arbitrary, but fixed till further notice. By Condition
3 in Definition \ref{Def. a.u. cadlag process}, there exist a $\mathrm{measurable}$
set $A$ with $P(A^{c})<2^{-k}$ and a r.r.v. $\tau_{1}$ with values
in $[0,1]$ such that, for each $\omega\in A$, we have 
\begin{equation}
\tau_{1}(\omega)\geq\delta_{aucl,\delta(\mathscr{\mathrm{\mathbf{T}}})}(2^{-k})>2^{-m(k)}\equiv\Delta_{m(k)},\label{eq:temp-264}
\end{equation}
and
\begin{equation}
d(X(0,\omega),X(\cdot,\omega))\leq2^{-k},\label{eq:temp-307-4-1-5-1}
\end{equation}
on the interval $\theta_{0}(\omega)\equiv[0,\tau_{1}(\omega))$. Then
inequalities \ref{eq:temp-307-4-1-5-1} and \ref{eq:temp-264} together
imply that, for each $\omega\in A$, we have
\[
d(X_{0}(\omega),X_{u}(\omega))\leq2^{-k}
\]
for each $u\in[0,\Delta_{m(k)}]\cap domain(X(\cdot,\omega))$. Consider
each $\kappa\geq k$. Then 
\[
\bigvee_{j=0}^{J(\kappa)}d(X_{0},X_{j\Delta(\kappa)})1_{A}\leq2^{-k},
\]
where $J_{\kappa}\equiv2^{m(\kappa)-m(k)}$. Therefore, with the function
$f_{\kappa}\in C(S^{J(\kappa)+1},d^{J(\kappa)+1})$ defined by 
\[
f_{\kappa}(x_{0},x_{1},\cdots,x_{J(\kappa)})\equiv\bigvee_{j=0}^{J(\kappa)}d(x_{0},x_{j})
\]
for each $(x_{0},x_{1},\cdots,x_{J(\kappa)})\in S^{J(\kappa)+1}$,
we have
\[
Ef_{\kappa}(X_{0},X_{\Delta(\kappa)},\cdots,X_{J(\kappa)\Delta(\kappa)})\equiv E\bigvee_{j=0}^{J(\kappa)}d(X_{0},X_{j\Delta(\kappa)})
\]
\begin{equation}
\leq E\bigvee_{j=0}^{J(\kappa)}d(X_{0},X_{j\Delta(\kappa)})1_{A}+P(A^{c})\leq2^{-k}+2^{-k}=2^{-k+1},\label{eq:temp-491}
\end{equation}
In terms of marginal distributions, inequality \ref{eq:temp-491}
can be rewritten as 
\begin{equation}
F_{0,\Delta(\kappa),\cdots,J(\kappa)\Delta(\kappa)}^{x,\mathbf{\mathrm{\mathbf{T}}}}f_{\kappa}\leq2^{-k+1},\label{eq:temp-600}
\end{equation}
where the family $F^{x,\mathbf{\mathrm{\mathbf{T}}}}$ of transition
f.j.d.'s are as constructed in \ref{Thm. Construction of Markov process from  initial distrinbution and semigrp-1-1}.
Recall here that $F_{0,\Delta(\kappa),\cdots,J(\kappa)\Delta(\kappa)}^{x,\mathbf{\mathbf{\mathrm{\mathbf{T}}}}}f_{\kappa}$
is continuous in $x\in(S,d)$. 

3. Separately, Assertion 2 of Lemma \ref{Lem. Existence of supremum distance Z}
applies to $v\equiv0$ and $v'\equiv\Delta_{m(k)}$, to show that
the supremum
\[
V_{0,k}\equiv\sup_{u\in[0,\Delta(m(k))]Q(\infty)}d(Z_{0},Z_{u})=\sup_{u\in[0,\Delta(m(k))]Q(\infty)}d(X_{0},X_{u})=\sup_{u\in[0,\Delta(m(k))]}d(X_{0},X_{u})
\]
is a well defined r.r.v, where the last equality is by the right continuity
of the a.u. càdlàg process $X$ . Consider each $\kappa\geq k$, Assertion
3 of Lemma \ref{Lem. Existence of supremum distance Z}then applies
to $v\equiv0$, $v'\equiv\Delta_{m(k)}$, and $h\equiv\kappa$, to
yield

\begin{equation}
0\leq E\sup_{u\in[0,\Delta(m(k))]Q(\infty)}d(Z_{0},Z_{u})-E\bigvee_{u\in[0,\Delta(m(k))]Q(m(\kappa))}d(Z_{0},Z_{u})\leq2^{-\kappa+5}.\label{eq:temp-446}
\end{equation}
Define the r.r.v.
\[
V_{0,k,\kappa}\equiv\bigvee_{u\in[0,\Delta(m(k))]Q(m(\kappa))}d(Z_{0},Z_{u})\equiv\bigvee_{j=0}^{J(\kappa)}d(X_{0},X_{j\Delta(\kappa)})\equiv f_{\kappa}(X_{0},X_{\Delta(\kappa)},\cdots,X_{J(\kappa)\Delta(\kappa)}).
\]
Then inequality \ref{eq:temp-446} can be rewritten compactly as
\begin{equation}
0\leq EV_{0,k}-EV_{0,k,\kappa}\leq2^{-\kappa+5}.\label{eq:temp-599}
\end{equation}
Hence for each $\kappa'\geq\kappa\geq k$, we have $V_{0,k,\kappa+1}\geq V_{0,k,\kappa}$
and 
\begin{equation}
0\leq EV_{0,k,\kappa'}-EV_{0,k,\kappa}\leq2^{-\kappa+5}.\label{eq:temp-599-1}
\end{equation}
Equivalently, 
\begin{equation}
0\leq F_{0,\Delta(\kappa'),\cdots,J(\kappa')\Delta(\kappa')}^{x,\mathbf{T}}f_{\kappa'}-F_{0,\Delta(\kappa),\cdots,J(\kappa)\Delta(\kappa)}^{x,\mathbf{T}}f_{\kappa}\leq2^{-\kappa+5},\label{eq:temp-169}
\end{equation}
where $x\in S$ is the arbitrary initial state. 

4. Now let $\eta$ be an arbitrary stopping time with values in $\overline{Q}_{h}$
for some $h\geq0$, relative to the filtration $\mathcal{L}$. Define
\[
V_{\eta,k,\kappa}\equiv\bigvee_{u\in[0,\Delta(m(k))]Q(m(\kappa))}d(X_{\eta},X_{\eta+u})
\]
\[
\equiv\bigvee_{j=0}^{J(\kappa)}d(X_{\eta},X_{\eta+j\Delta(\kappa)})\equiv f_{\kappa}(X_{\eta},X_{\eta+\Delta(\kappa)},\cdots,X_{\eta+J(\kappa)\Delta(\kappa)}).
\]
Then
\[
0\leq E(V_{\eta,k,\kappa'}-V_{\eta,k,\kappa})
\]
\[
=E(f_{\kappa'}(X_{\eta},X_{\eta+\Delta(\kappa')},\cdots,X_{\eta+J(\kappa')\Delta(\kappa')})-f_{\kappa}(X_{\eta},X_{\eta+\Delta(\kappa)},\cdots,X_{\eta+J(\kappa)\Delta(\kappa)}))
\]
\[
=E(F_{0,\Delta(\kappa'),\cdots,J(\kappa')\Delta(\kappa')}^{X(\eta),\mathbf{T}}f_{\kappa'}-EF_{0,\Delta(\kappa),\cdots,J(\kappa)\Delta(\kappa)}^{X(\eta),\mathbf{T}}f_{\kappa})\leq2^{-\kappa+5},
\]
where the last inequality is thanks to inequality \ref{eq:temp-169}.
Hence

\[
V_{\eta,k,\kappa}\uparrow U\quad\mathrm{a.u}.
\]
for some r.r.v.$U\in L$, as $\kappa\rightarrow\infty$. Let $\omega\in domain(U)$
be arbitrary. Then the limit
\[
\lim_{\kappa\rightarrow\infty}V_{\eta,k,\kappa}(\omega)\equiv\lim_{\kappa\rightarrow\infty}\bigvee_{u\in[0,\Delta(m(k))]Q(m(\kappa))}d(X_{\eta}(\omega),X_{\eta+u}(\omega))
\]
exists and is equal to $U(\omega)$. Hence 
\[
\sup_{u\in[0,\Delta(m(k))]Q(\infty)}d(X_{\eta}(\omega),X_{\eta+u}(\omega))=U(\omega).
\]
Consequently, by the right continuity of $X(\cdot,\omega)$ at $\eta(\omega)$,
we obtain
\[
V_{\eta,k}(\omega)\equiv\sup_{u\in[0,\Delta(m(k))]}d(X_{\eta}(\omega),X_{\eta+u}(\omega))=U(\omega)
\]
We conclude that the function $V_{\eta,k}=U$ on a full set, and therefore
that $V_{\eta,k}$ is an integrable r.r.v. Moreover
\[
EV_{\eta,k}=EU=\lim_{\kappa\rightarrow\infty}EV_{\eta,k,\kappa}==\lim_{\kappa\rightarrow\infty}F_{0,\Delta(\kappa),\cdots,J(\kappa)\Delta(\kappa)}^{X(\eta),\mathbf{T}}f_{\kappa}\leq2^{-k+1},
\]
where the inequality follows from inequality \ref{eq:temp-600}. The
lemma is proved. 
\end{proof}
\begin{lem}
\textbf{\emph{\label{Lem. Observability at stopping time}(Observability
of a.u. càdlàg Markov process at stopping time). }}Let $x\in S$ be
arbitary. Let 
\[
X:[0,\infty)\times(\Omega,L,E)\rightarrow(S,d)
\]
be an arbitrary a.u. càdlàg Markov process generrated by the initial
state $x$ Markov semigroup $\mathbf{T}$, and which is adapted to
some right continuous filtration $\mathcal{L}\equiv\{L^{(t)}:t\in[0,\infty)\}$.
Let $\tau$ be an arbitrary stopping time with values in $[0,\infty)$,
relative to $\mathcal{L}$. Then the following holds.

\emph{1.} The function $X_{\tau}$ is a well defined r.v. which is
measurable relative to $L^{(\tau)}$.

\emph{2.} Specifically, let $(\eta_{h})_{h=0,1,\cdots}$ be an arbitrary
non increasing sequence of stopping times such that, for each $h\geq0$,
the r.r.v. $\eta_{h}$ has values in $\overline{Q}_{h}$ , relative
to the filtration $\mathcal{L}$, and such that 
\begin{equation}
\tau\leq\eta_{h}<\tau+2^{-h+2}.\label{eq:temp-421-1-2-1-2}
\end{equation}
Note that such a sequence exists according to Assertion \emph{1} of
Proposition \ref{Prop. Stopping times rlative to right continuous filtration}.
Then $X_{\eta(h)}\rightarrow X_{\tau}$ a.u. as $h\rightarrow\infty$. 
\end{lem}
\begin{proof}
1. Let $(\eta_{h})_{h=0,1,\cdots}$ be an arbitrary sequence of stopping
times as in the hypothesis of Assertion 2. Let $k\geq0$ be arbitrary.
Let $m_{k}\equiv m_{k}(\delta_{\mathbf{T}})\geq0$ be the integer
constructed in Lemma \ref{Lem. Short-term quiecience of a.u. cadlag Feller process after each stopping time with uniformly spaced dyadic values}
relative to the Markov semigroup $\mathbf{T}$. For abbreviation,
write $h_{k}\equiv(m_{k}+3)$ and $\zeta_{k}\equiv\eta_{h(k)}$. Then
$\zeta_{k}$ is a stopping time with values in $\overline{Q}_{h(k)}$,
relative to the filtration $\mathcal{L}$. Hence Lemma \ref{Lem. Short-term quiecience of a.u. cadlag Feller process after each stopping time with uniformly spaced dyadic values},
applied to the stopping time $\zeta_{k}$, says that the function
\[
V_{k}\equiv V_{\zeta(k),k}\equiv\sup_{s\in[0,\Delta(m(k))]}d(X_{\zeta(k)},X_{\zeta(k)+s})
\]
with
\begin{equation}
EV_{k}\leq2^{-k+1}.\label{eq:temp-447-1}
\end{equation}

2. At the same time, inequality \ref{eq:temp-421-1-2-1-2} implies
that 
\begin{equation}
\tau\leq\eta_{h(k)}\equiv\zeta_{k}<\tau+2^{-h(k)+3}\equiv\tau+2^{-m(k)}\equiv\tau+\Delta_{m(k)},\label{eq:temp-446-2}
\end{equation}
whence $\zeta_{k}\downarrow\tau$ uniformly as $k\rightarrow\infty$.
Note that $\zeta_{k+1}\leq\zeta_{k}$ because $\eta_{h(k+1)}\leq\eta_{h(k)}$
. Hence inequality \ref{eq:temp-446-2} leads to
\begin{equation}
\tau\leq\zeta_{k+1}\leq\zeta_{k}<\tau+\Delta_{m(k)}\leq\zeta_{k+1}+\Delta_{m(k)},\label{eq:temp-449-1}
\end{equation}
An immediate consequence is that, for each $\kappa\geq k$, we have
\[
(\tau,\tau+\Delta_{m(\kappa)}]\overline{Q}_{\infty}\subset(\tau,\zeta_{\kappa+1}+\Delta_{m(\kappa)}]\overline{Q}_{\infty}
\]
\[
=(\bigcup_{j=\kappa}^{\infty}[\zeta_{j+2},\zeta_{j+1}]\overline{Q}_{\infty})\cup[\zeta_{\kappa+1},\zeta_{\kappa+1}+\Delta_{m(\kappa)}]\overline{Q}_{\infty}
\]
\[
\subset(\bigcup_{j=\kappa}^{\infty}[\zeta_{j+2},\zeta_{j+2}+\Delta_{m(j+1)}]\overline{Q}_{\infty})\cup[\zeta_{\kappa+1},\zeta_{\kappa+1}+\Delta_{m(\kappa)}]\overline{Q}_{\infty}
\]
\begin{equation}
=\bigcup_{j=\kappa}^{\infty}[\zeta_{j+1},\zeta_{j+1}+\Delta_{m(j)}]\overline{Q}_{\infty}\label{eq:temp-448-1}
\end{equation}

3. For each $j\geq0$, take any $\varepsilon_{j}\in(2^{-j/2-1},2^{-j/2}]$
such that the set $A_{j}\equiv(V_{j}>\varepsilon_{j})$ is measurable.
Then 
\[
P(A_{j})\leq\varepsilon_{j}^{-1}EV_{j}\leq\varepsilon_{j}^{-1}2^{-j+1}<2^{j/2+1}2^{-j+1}=2^{-j/2+2},
\]
where the first inequality is by Chebychev's inequality, and where
the second inequality is by inequality \ref{eq:temp-447-1}. Therefore,
we can define the measurable set $A_{k+}\equiv\bigcup_{j=k}^{\infty}A_{j}$
, with
\[
P(A_{k+})<\sum_{j=k}^{\infty}P(A_{j})<\sum_{j=k}^{\infty}2^{-j/2+2}=2^{-k/2+2}(1-2^{-1/2})^{-1}<2^{-k/2+4}
\]

4. Now let $\omega\in A_{k+}^{c}$be arbitrary. Consider each $j\geq k$
and each 
\[
t\in[\zeta_{j+1}(\omega),\zeta_{j+1}(\omega)+\Delta_{m(j)}]\overline{Q}_{\infty}.
\]
Then $u\equiv t-\zeta_{j+1}(\omega)\in[0,\Delta_{m(j)}]\overline{Q}_{\infty}$.
Hence
\[
d(X_{\zeta(j+1)}(\omega),X_{t}(\omega))\equiv d(X_{\zeta(j+1)}(\omega),X_{\zeta(j+1)+u}(\omega))
\]
\begin{equation}
\leq\sup_{s\in[0,\Delta(m(j))]}d(X_{\zeta(j+1)}(\omega),X_{\zeta(j+1)+s}(\omega))\equiv V_{j+1}(\omega)\leq\varepsilon_{j+1},\label{eq:temp-604-1}
\end{equation}
where the last inequality is because $\omega\in A_{k+}^{c}\subset A_{,j+1}^{c}$.
Note that 
\[
\zeta_{j}(\omega)\in[\zeta_{j+1}(\omega),\zeta_{j+1}(\omega)+\Delta_{m(j)}]\overline{Q}_{\infty},
\]
thanks to inequality \ref{eq:temp-449-1}. Hence inequality \ref{eq:temp-604-1}
applies to $\zeta_{j}(\omega)$ in the place of $t$, to yield
\begin{equation}
d(X_{\zeta(j+1)}(\omega),X_{\zeta(j)}(\omega))\leq\varepsilon_{j+1},\label{eq:temp-576-3}
\end{equation}
where $j\geq k$ is arbitrary. Since $\sum_{j=0}^{\infty}\varepsilon_{j+1}<\infty$,
it follows that $X_{\zeta(j)}(\omega)\rightarrow x_{\omega}$ for
some $x_{\omega}\in S,$ as $j\rightarrow\infty$, with 
\begin{equation}
d(x_{\omega},X_{\zeta(j)}(\omega))\leq\sum_{i=j}^{\infty}\varepsilon_{i+1}\leq\sum_{i=j}^{\infty}2^{-(i+1)/2}=2^{-(j+1)/2}(1-2^{-1/2})^{-1}<2^{-(j+1)/2+4}.\label{eq:temp-605-1}
\end{equation}

5. Now consider each $t\in(\tau(\omega),\tau(\omega)+\Delta_{m(k)})\overline{Q}_{\infty}$.
Then, according to relation \ref{eq:temp-448-1}, there exists $j\geq k$
such that
\[
t\in[\zeta_{j+1}(\omega),\zeta_{j+1}(\omega)+\Delta_{m(j)}]\overline{Q}_{\infty}.
\]
Therefore
\[
d(x_{\omega},X_{t}(\omega))\leq d(x_{\omega},X_{\zeta(j)}(\omega))+d(X_{\zeta(j)}(\omega),X_{\zeta(j+1)}(\omega))+d(X_{\zeta(j+1)}(\omega),X_{t}(\omega))
\]
\[
\leq\sum_{i=j}^{\infty}\varepsilon_{i+1}+\varepsilon_{j+1}+\varepsilon_{j+1}<3\sum_{i=j}^{\infty}\varepsilon_{i+1}<2^{-(J+1)/2+6}\leq2^{-(k+1)/2+6}
\]
where $t\in(\tau(\omega),\tau(\omega)+\Delta_{m(k)})\overline{Q}_{\infty}$
is arbitrary. Since the set $(\tau(\omega),\tau(\omega)+\Delta_{m(k)})\overline{Q}_{\infty}$
is a dense subset of $(\tau(\omega),\tau(\omega)+\Delta_{m(k)})\cap domain(X(\cdot,\omega))$,
and since the function $X(\cdot,\omega)$ is right continuous on $domain(X(\cdot,\omega))$
according to Condition 1 in Definition \ref{Def. Metric space of a.u. cadlag process on =00005B0,inf)},
it follows that 
\begin{equation}
d(x_{\omega},X_{t}(\omega))\leq2^{-(k+1)/2+6}\label{eq:temp-219}
\end{equation}
for each $t\in(\tau(\omega),\tau(\omega)+\Delta_{m(k)})\cap domain(X(\cdot,\omega))$,
where $2^{-(k+1)/2+6}$ is arbitrarily small. Thus $\lim_{t\rightarrow\tau(\omega);t>\tau(\omega)}X(t,\omega)$
exists and is equal $x_{\omega}\in S$. Therefore, by Condition 2
in Definition \ref{Def. Metric space of a.u. cadlag process on =00005B0,inf)},
the right-completeness condition, we have $\tau(\omega)\in domain(X(\cdot,\omega))$.
Then, by Condition 1 in Definition \ref{Def. Metric space of a.u. cadlag process on =00005B0,inf)},
the right-continuity condition, we have $X_{\tau}(\omega)\equiv X(\tau(\omega),\omega)=x_{\omega}.$ 

Inequality \ref{eq:temp-605-1} can be rewritten as 
\[
d(X_{\tau}(\omega),X_{\zeta(j)}(\omega))\leq2^{-(j+1)/2+4}
\]
where $j\geq k$ and $\omega\in A_{k+}^{c}$are arbitrary. Thus $X_{\zeta(j)}\rightarrow X_{\tau}$
uniformly on $A_{k+}^{c},$ as $j\rightarrow\infty$, where $P(A_{k+}^{c})<2^{-k/2+4}$
is arbitrarily small for sufficiently large $k\geq0$. We conclude
that $X_{\zeta(j)}\rightarrow X_{\tau}$ a.u. as $j\rightarrow\infty$.
It follows that $X_{\tau}$ is a r.v.

6. We will next verify that $X_{\tau}$ is measurable relative to
the probability subspace $L^{(\tau)}.$ To that end, let $s>0$ be
arbitrary. Take $k\geq0$ so large that $\Delta_{m(k)}<s$. Then,
for each $j\geq k$, we have 
\[
\zeta_{j}\leq\tau+\Delta_{m(j)}\leq\tau+\Delta_{m(k)}<\tau+s,
\]
whence the r.v. $X_{\zeta(j)}$ is measurable relative to the probability
subspace $L^{(\tau+s)}.$ Hence, as $j\rightarrow\infty$, the limiting
r.v. $X_{\tau}$ is measurable relative to $L^{(\tau+s)}$, where
$s>0$ is arbitrary. Therefor $X_{\tau}$ is measurable relative to
the probability subspace
\[
\bigcap_{s>0}L^{(\tau+s)}\equiv L^{(\tau+)}=L^{(\tau)},
\]
where the last equality follows from \ref{Prop. Stopping times rlative to right continuous filtration},
in view of the assumed right continuity of the filtration $\mathcal{L}$.
Assertion 1 is proved.

7. Inequality \ref{eq:temp-219} can be rewritten as
\begin{equation}
d(X_{\tau}(\omega),X_{t}(\omega))\leq2^{-(k+1)/2+6}\label{eq:temp-219-1}
\end{equation}
for each $t\in(\tau(\omega),\tau(\omega)+\Delta_{m(k)})\cap domain(X(\cdot,\omega))$.
By the right continuity of the function $X(\cdot,\omega)$, we therefore
obtain
\begin{equation}
d(X_{\tau}(\omega),X_{t}(\omega))\leq2^{-(k+1)/2+6}\label{eq:temp-219-1-2}
\end{equation}
 for each $t\in[\tau(\omega),\tau(\omega)+\Delta_{m(k)})\cap domain(X(\cdot,\omega)).$
Now let $h\geq h_{k}$ be arbitrary. Then 
\begin{equation}
\tau\leq\eta_{h}\leq\eta_{h(k)}<\tau+\Delta_{m(k)},\label{eq:temp-446-2-1}
\end{equation}
according to inequality \ref{eq:temp-446-2}. Hence inequality \ref{eq:temp-219-1-2}
implies that 
\begin{equation}
d(X_{\tau}(\omega),X_{\eta(h)}(\omega))\leq2^{-(k+1)/2+6},\label{eq:temp-219-1-2-1}
\end{equation}
for each $h\geq h_{k},$ where $k$ and $\omega\in A_{k+}^{c}$are
arbitrary. Since $2^{-(k+1)/2+6}$ and $P(A_{k+})<2^{-k/2+4}$ are
arbitrarily small, we conclude that $X_{\eta(h)}\rightarrow X_{\tau}$
a.u. The lemma is proved.
\end{proof}
\begin{thm}
\textbf{\emph{\label{Thm. Each a.u. cadlag Markov process  is a strong Markov process}
(Each a.u. càdlàg Markov process is a strong Markov process).}} Let
$x\in S$ be arbitrary. Let 
\[
X:[0,\infty)\times(\Omega,L,E)\rightarrow(S,d)
\]
be an arbitrary a.u. càdlàg Markov process generated by the initial
state $x$ and the Markov semigroup $\mathbf{T}$, and adapted to
some right continuous filtration $\mathcal{L}\equiv\{L^{(t)}:t\in[0,\infty)\}$.
All stopping times will be understood to be relative to this filtration
$\mathcal{L}$. 

Then the following holds.

1. Let $\tau$ be an arbitrary stopping time with values in $[0,\infty)$.
Then the function $X_{\tau}$ is a r.v. relative to $L^{(\tau)}$. 

2. The process $X$ is strongly Markov\emph{ }relative to  filtration
$\mathcal{L}$. Specifically, let $\tau$ be an arbitrary stopping
time with values in $[0,\infty)$. Let $0\equiv r_{0}\leq r_{1}\leq\cdots\leq r_{m}$
be an arbitrary nondecreasing sequence in $[0,\infty)$, and let $f\in C(S^{m+1},d^{m+1})$
be arbitrary. Then
\[
E(f(X_{\tau+r(0)},X_{\tau+r(1)},\cdots,X_{\tau+r(m)})|L^{(\tau)})
\]
\begin{equation}
=E(f(X_{\tau+r(0)},X_{\tau+r(1)},\cdots,X_{\tau+r(m)})|X_{\tau})=F_{r(0),\cdots,r(m)}^{X(\tau),\mathbf{T}}(f)\label{eq:temp-576-1-1-2}
\end{equation}
as r.r.v.'s
\end{thm}
\begin{proof}
For ease of notations, we present the proof only for the case where
$n=2$, the general case being similar.

1. Assertion 1 is merely a restatement of Assertion 1 of Lemma \ref{Lem. Observability at stopping time}

2. To prove Assertion 2, let $\tau$ be an arbitrary stopping time
with values in $[0,\infty)$. Let $0\equiv r_{0}\leq r_{1}\leq r_{2}$
be an arbitrary nondecreasing sequence in $[0,\infty)$, and let $f\in C(S^{3},d^{3})$
be arbitrary. As in Assertion 1 of Proposition \ref{Prop. Stopping times rlative to right continuous filtration},
construct a nonincreasing sequence $(\eta_{h})_{h=0,1,\cdots}$ of
stopping times such that, for each $h\geq0$, the r.r.v. $\eta_{h}$
is a stopping time with values in $\overline{Q}_{h}$, and such that
\begin{equation}
\tau\leq\eta_{h}<\tau+2^{-h+1}.\label{eq:temp-234}
\end{equation}
Let $h\geq0$ and $i=1,2$ be arbitrary. Take any $s_{h,i}\in[r_{i},r_{i}+2^{-h+1})\overline{Q}_{h}$.
Take $s_{h,0}\equiv0$. Then
\begin{equation}
\tau+r_{i}\leq\eta_{h}+s_{h,i}<\tau+r_{i}+2^{-h+2}.\label{eq:temp-234-1}
\end{equation}
where $i=0,1,2$ is arbitrary.

3. Consider each indicator $U\in L^{(\tau)}$. Then $U\in L^{(\eta(h))}$
because $\tau\leq\eta_{h}$. Hence $U1_{(\eta(h)=r)}\in L^{(r))}$for
each $r\in\overline{Q}_{h}$. Therefore
\[
Ef(X_{\eta(h)},X_{\eta(h)+s(h,1)},X_{\eta(h)+s(h,2)})U
\]
\[
=\sum_{r\in\overline{Q}(h)}Ef(X_{\eta(h)},X_{\eta(h)+s(h,1)},X_{\eta(h)+s(h,2)})U1_{(\eta(h)=r)}
\]
\[
=\sum_{r\in\overline{Q}(h)}Ef(X_{r},X_{r+s(h,1)},X_{r+s(h,2)})U1_{(\eta(h)=r)}
\]
\[
=\sum_{r\in\overline{Q}(h)}E(F_{0,s(h,1),s(h,2)}^{X(r),\mathbf{T}}f)U1_{(\eta(h)=r)}
\]
\[
=\sum_{r\in\overline{Q}(h)}E(F_{0,s(h,1),s(h,2)}^{X(\eta(h)),\mathbf{T}}f)U1_{(\eta(h)=r)}
\]
\[
=E(F_{0,s(h,1),s(h,2)}^{X(\eta(h)),\mathbf{T}}f)U\sum_{r\in\overline{Q}(h)}1_{(\eta(h)=r)}
\]
\begin{equation}
=E(F_{0,s(h,1),s(h,2)}^{X(\eta(h)),\mathbf{T}}f)U.\label{eq:temp-576-1-1-1}
\end{equation}

4. Consider each $i=0,1,2$. Then we can apply Assertion 2 of Lemma
\ref{Lem. Observability at stopping time}, where $\tau,(\eta_{h})_{h=0,1,\cdots}$
are replaced by $\tau+r_{i},(\eta_{h}+s_{h,i})_{h=0,1,\cdots}$ respectively,
to obtain 
\begin{equation}
X_{\eta(h)+s(h,i)}\rightarrow X_{\tau+r(i)}\quad\mathrm{a.u.}\label{eq:temp-443}
\end{equation}
as $h\rightarrow\infty$. Hence, we have, for the left-hand sid of
equality \ref{eq:temp-576-1-1-1},
\[
Ef(X_{\eta(h)},X_{\eta(h)+s(h,1)},X_{\eta(h)+s(h,2)})U\rightarrow Ef(X_{\tau},X_{\tau+r(1)},X_{\tau+r(2)})U.
\]

5. For convergence of the right-hand side, let $\varepsilon>0$ be
arbitrary. There is no loss of genrality to assume that the function
$f\in C(S^{3},d^{3})$ has values in $[0,1]$ and has a modulus of
continuity $\delta_{f}$. Then Theorem \ref{Thm. Construction of transition  f.j.d.s from E0 and Markov semi-group, discrete parameters-1}
yields some $\delta_{3}(\varepsilon,\delta_{f},\delta_{\mathbf{T}},\alpha_{\mathbf{T}})>0$
such that
\[
\left\Vert F_{s(h,0),s(h,1),s(h,2)}^{*,\mathbf{T}}f-F_{r(0),r(,1),r(2)}^{*,\mathbf{T}}f\right\Vert \leq\varepsilon
\]
provided that 
\begin{equation}
\bigvee_{i=0}^{2}|s_{h,i}-r_{i}|<\delta_{3}(\varepsilon,\delta_{f},\delta_{\mathbf{T}},\alpha_{\mathbf{T}}).\label{eq:temp-551}
\end{equation}
 Hence
\begin{equation}
|E(F_{0,s(h,1),s(h,2)}^{X(\eta(h)),\mathbf{T}}f)U-E(F_{r(0),r(,1),r(2)}^{X(\eta(h)),\mathbf{T}}f)U|\leq|E(F_{0,s(h,1),s(h,2)}^{X(\eta(h)),\mathbf{T}}f)-E(F_{r(0),r(,1),r(2)}^{X(\eta(h)),\mathbf{T}}f)|\leq\varepsilon\label{eq:temp-557}
\end{equation}
provided that inequality \ref{eq:temp-551} holds.

6. At the same time, Theorem \ref{Thm. Construction of transition  f.j.d.s from E0 and Markov semi-group, discrete parameters-1}
implies that the function $F_{r(0),r(,1),r(2)}^{*,\mathbf{T}}f$ is
a member of $C(S,d)$. Hence, since $X_{\eta(h)}\rightarrow X_{\tau}\quad\mathrm{a.u.}$,
we have
\[
E(F_{r(0),r(,1),r(2)}^{X(\eta(h)),\mathbf{T}}f)U\rightarrow E(F_{r(0),r(,1),r(2)}^{X(\tau),\mathbf{T}}f)U
\]
as $h\rightarrow\infty$. Therefore there exists $h_{0}\geq0$ so
large that 
\[
|E(F_{r(0),r(,1),r(2)}^{X(\eta(h)),\mathbf{T}}f)U-E(F_{r(0),r(,1),r(2)}^{X(\tau),\mathbf{T}}f)U|<\varepsilon
\]
for each $h\geq h_{0}$. Combining with inequality \ref{eq:temp-557},
we obtain
\[
|E(F_{0,s(h,1),s(h,2)}^{X(\eta(h)),\mathbf{T}}f)U-E(F_{r(0),r(,1),r(2)}^{X(\tau),\mathbf{T}}f)U|<2\varepsilon
\]
provided that $h\geq h_{0}$ is so large that inequality \ref{eq:temp-551}
holds. Since $\varepsilon>0$ is arbitrary, we see that 
\[
E(F_{0,s(h,1),s(h,2)}^{X(\eta(h)),\mathbf{T}}f)U\rightarrow E(F_{r(0),r(,1),r(2)}^{X(\tau),\mathbf{T}}f)U
\]
as $h\rightarrow\infty$. 

7. Thus we have convergence also of the right-hand sied of equality
\ref{eq:temp-576-1-1-1}. The limits on each side, found in Steps
4 and 6 respectively, must therefore be equal, namely
\begin{equation}
Ef(X_{\tau},X_{\tau+r(1)},X_{\tau+r(2)})U=E(F_{r(0),r(,1),r(2)}^{X(\tau),\mathbf{T}}f)U,\label{eq:temp-516}
\end{equation}
where the indicator $U\in L^{(\tau)}$ is arbitrary. Hence, since
$F_{r(0),r(,1),r(2)}^{X(\tau),\mathbf{T}}f\in L^{(\tau)}$, we have
\[
E(f(X_{\tau},X_{\tau+r(1)},X_{\tau+r(2)})|L^{(\tau)})=F_{r(0),r(,1),r(2)}^{X(\tau),\mathbf{T}}f.
\]
In particular, equality \ref{eq:temp-516} holds for each indicator
$U\in L(X_{\tau})$. Hence
\[
E(f(X_{\tau},X_{\tau+r(1)},X_{\tau+r(2)})|X_{\tau})=F_{r(0),r(,1),r(2)}^{X(\tau),\mathbf{T}}f.
\]
Thus the desired equality \ref{eq:temp-576-1-1-2} is proved for the
case $m=2$, with a simillar proof for the general case. We conclude
that the process $X$ is strongly Markov. 
\end{proof}

\section{Abundance of Stopping Times for a.u. càdlàg Markov Processes}

..
\begin{defn}
\textbf{\emph{\label{Def. First exit time} }}\textbf{(First exit
time).} Let $(S,d)$ be an arbitrary locally compact metric space.
Let $X:[0,\infty)\times(\Omega,L,E)\rightarrow(S,d)$ be an arbitrary
a.u. càdlàg process which is adapted to some right continuous filtration
$\mathcal{L}\equiv\{L^{(t)}:t\in[0,\infty)\}$. 

Let $f:(S,d)\rightarrow R$ be an arbitrary function which is continuous
on compact subsets of $(S,d)$, such that $f(X_{0})\leq a_{0}$ for
some $a_{0}\in R$. Let $a\in(a_{0},\infty)$ and $M\geq1$ be arbitrary.
Suppose $\tau$ is a stopping time relative to $\mathcal{L}$ and
with values in $[0,M]$ such that the function $X_{\tau}$ is a well-defined
r.v relative to $L^{(\tau)}$. Suppose, for each $\omega\in domain(X_{\tau})$,
we have 

(i). $f(X(\cdot,\omega))<a$ on the interval $[0,\tau(\omega))$,
and 

(ii).\emph{ $f(X_{\tau}(\omega))\geq a$ }if\emph{ $\tau(\omega)<M$.}

Then we say that $\tau$ is the\index{first exit time} \emph{first
exit time }in\emph{ }$[0,M]$ of the open subset $(f<a)$ by the process
$X$, and write $\overline{\tau}_{f,a,M}\equiv\overline{\tau}_{f,a,M}(X)\equiv\tau$. 

Note that there is no requirement that the process actually exits
$(f<a)$ ever. It is stopped at time $M$ if it does not exit by then.
$\square$

Before proving the existence of these first exit times, the next proposition
make precise some  intuitions.
\end{defn}
\begin{lem}
\label{Lem. Basics of first exit times}\textbf{\emph{(Basics of first
exit times).}} Let $(S,d)$ be an arbitrary locally compact metric
space. Let $X:[0,\infty)\times(\Omega,L,E)\rightarrow(S,d)$ be an
arbitrary a.u. càdlàg process which is adapted to some right continuous
filtration $\mathcal{L}\equiv\{L^{(t)}:t\in[0,\infty)\}$. Let $f:(S,d)\rightarrow R$
be an arbitrary function which is continuous on compact subsets of
$(S,d)$, such that $f(X_{0})\leq a_{0}$ for some $a_{0}\in R$.
Let $a\in(a_{0},\infty)$ and $M\geq1$ be arbitrary. 

Let $a\in(a_{0},\infty)$ be such that the first exit time $\overline{\tau}_{f,a,M}$
exists for each $M\geq1$. Consider each $M\geq1$. Let $r\in(0,M)$
and $N\geq M$ be arbitrary. Then the following holds.

\emph{1.} $\overline{\tau}_{f,a,M}\leq\overline{\tau}_{f,a,N}$.

\textup{2. }$(\overline{\tau}_{f,a,M}<M)\subset(\overline{\tau}_{f,a,N}=\overline{\tau}_{f,a,M})$. 

\emph{3}. $(\overline{\tau}_{f,a,N}\leq r)=(\overline{\tau}_{f,a,M}\leq r)$.
\[
.
\]
\end{lem}
\begin{proof}
1. Let $\omega\in domain(\overline{\tau}_{f,a,M})\cap domain(\overline{\tau}_{f,a,M})$
be arbitrary. For the sake of a contradiction, suppose $t\equiv\overline{\tau}_{f,a,M}(\omega)>s\equiv\overline{\tau}_{f,a,N}(\omega)$.
Then $s<\overline{\tau}_{f,a,M}(\omega)$. Hence we can apply Condition
(i) to $M$, to obtain $f(X_{s}(\omega))<a$. At the same time, $\overline{\tau}_{f,a,N}(\omega)<t\leq N$.
Hence we can apply Condition (ii) to $N$, to obtain $f(X_{\overline{\tau}(f,a,N)}(\omega))\geq a$.
In other words, $f(X_{s}(\omega))\geq a$, a contradiction. We conclude
that $\overline{\tau}_{f,a,N}(\omega)\geq\overline{\tau}_{f,a,M}(\omega)$
where $\omega\in domain(\overline{\tau}_{f,a,M})\cap domain(\overline{\tau}_{f,a,M})$
is arbitrary/ Assertion 1 is proved.

2. Next, suppose $t\equiv\overline{\tau}_{f,a,M}(\omega)<M$. Then
Condition (ii) implies that $f(X_{t}(\omega))\geq a$. For the sake
of a contradiction, suppose $t<\overline{\tau}_{f,a,N}(\omega).$
Then Condition (i) implies that $f(X_{t}(\omega))<a$, a contradiction.
we conclude that $\overline{\tau}_{f,a,M}(\omega)\equiv t\geq\overline{\tau}_{f,a,N}(\omega)$.
Combining, with Assertion 1, we obtain $\overline{\tau}_{f,a,M}(\omega)=\overline{\tau}_{f,a,N}(\omega)$.
Assertion 2 is proved.

3. Note that 
\[
(\overline{\tau}_{f,a,N}\leq r)\subset(\overline{\tau}_{f,a,M}\leq r)=(\overline{\tau}_{f,a,M}\leq r)(\overline{\tau}_{f,a,M}<M)
\]
\begin{equation}
\subset(\overline{\tau}_{f,a,M}\leq r)(\overline{\tau}_{f,a,M}=\overline{\tau}_{f,a,N})\subset(\overline{\tau}_{f,a,M}\leq r)(\overline{\tau}_{f,a,N}\leq r)=(\overline{\tau}_{f,a,N}\leq r),\label{eq:temp-597}
\end{equation}
where we used the just established Assertions 1 and 2 repeatedly.
Since the left-most set and the right-most set in set relation \ref{eq:temp-597}
are the same, the inclusion relations can be replaced by equality.
Assertion 3 then follows.
\end{proof}
\begin{defn}
\textbf{\label{Def. Specification of a filtered a.u. cadlag Feller Process}
(Specification of a filtered a.u.càdlàg Markov process, and notations
for related objects).} In the remainder of this section, let $(S,d)$
denote a compact metric space with $d\leq1$. Let $\mathbf{T}$ be
a Markov semigroup with parameter set $[0,\infty)$ and state space
$(S,d)$, and with a modulus of strong continuity $\delta_{\mathbf{T}}$.
Let $\mathcal{L}\equiv\{L^{(t)}:t\in[0,\infty)\}$ be a right continuous
filtration on a given sample space $(\Omega,L,E)$. 

Let $x\in S$ be arbitrary. Let 
\[
X:[0,\infty)\times(\Omega,L,E)\rightarrow(S,d)
\]
be an a.u. càdlàg Markov process generated by the initial state $x$
and the Markov semigroup $\mathbf{T}$ which is adapted. to the filtration
$\mathcal{L}$.

Note that,then, the process $X$ is equivalent to the process $X^{x,\mathbf{T}}$
constructed in has a modulus of continuity in Theorem \ref{Thm. Construction of Markov process from  initial distrinbution and semigrp-1-1},
and is therefore continuous in probability with a modulus of continuity
in probability $\delta_{Cp,\delta(\mathbf{\mathbf{\mathscr{\mathbf{T}}}})}$.
Moreover, for each $N\geq0$, the shifted process $X^{N}:[0,1]\times(\Omega,L,E)\rightarrow(S,d)$,
defined by $X_{r}^{N}\equiv X_{N+r}$ for each $r\in[0,1]$, has a
modulus of a.u. càdlàg $\delta_{aucl,\delta(\mathbf{\mathbf{\mathscr{\mathbf{T}}}})}$.
Define the operation $\delta_{cp,\delta(\mathscr{\mathbf{T}})}$ by\emph{
}
\[
\delta_{cp,\delta(\mathscr{\mathbf{T}})}(\varepsilon)\equiv\delta_{Cp,\delta(\mathscr{\mathbf{T}})}(2^{-2}(1\wedge\varepsilon)^{2})>0
\]
for each $\varepsilon>0$. Note the lower-case in the subscript of
$\delta_{cp,\delta(\mathscr{\mathbf{T}})}$, in contrast to the subscript
in $\delta_{Cp,\delta(\mathbf{\mathbf{\mathscr{\mathbf{T}}}})}.$

Proposition \ref{Prop. Alternative characterization of r.f. continuity}
then says that, for each $N\geq0$, $\varepsilon>0$, and $t,s\in[0,1]$
with $|s-t|<\delta_{cp,\delta(\mathscr{\mathbf{T}})}(\varepsilon)$,
there exists a $\mathrm{measurable}$ set $D_{t,s}$ with $P(D_{t,s}^{c})\leq\varepsilon$
such that
\[
d(X^{N}(t,\omega),X^{N}(s,\omega))\leq\varepsilon
\]
for each $\omega\in D_{t,s}$.
\end{defn}
\begin{thm}
\textbf{\emph{\label{Thm. Abundance of first exit times for Markov process}
(Abundance of first exit times for Markov processes).}} Let the objects
$(S,d)$, $\mathbf{T}$, $(\Omega,L,E)$, $\mathcal{L}$, x, and $X$
be as specified. In particular $(S,d)$ is a compact metric space
with $d\leq1$. Let $f\in C(S,d)$ be arbitrary such that $f(X_{0})=f(x)\leq a_{0}$
for some $a_{0}\in R$. 

Then there exists a countable subset $G$ of $R$ such that, for each
$M\geq1$ and for each $a\in(a_{0},\infty)G_{c}$, the first exit
time $\overline{\tau}_{f,a,M}(X)$ exists as defined in Definition
\ref{Def. First exit time}. Here $G_{c}$ denotes the metric complement
of $G$ in $R$. Moreover, the set $G$ of exceptional points is completely
determined by the function $f$ and the marginal f.j.d.'s of the process
$X$. 
\end{thm}
\begin{proof}
1. Let $\delta_{f}$ be a modulus of continuity of the function $f$.
Let $N\geq0$ be arbitrary, but fixed till further notice. Consider
the a.u. càdlàg\emph{ }process $X^{N}:[0,1]\times\Omega\rightarrow S,$
with modulus of continuity in probability $\delta_{Cp,\delta(\mathbf{\mathbf{\mathscr{\mathbf{T}}}})}$
and with modulus of a.u. càdlàg $\delta_{aucl,\delta(\mathbf{\mathbf{\mathscr{\mathbf{T}}}})}$.
Definition \ref{Def. a.u. cadlag process}, says that there exists
a full set $B\subset\bigcap_{t\in Q(\infty)}domain(X_{t})$ with the
following properties. Each $\omega\in B$ satisfies the right-continuity
condition and the right-completeness condition in Definition \ref{Def. a.u. cadlag process}.
Moreover, for each $k\geq0$ and $\varepsilon_{k}>0$, there exist
(i) $\delta_{k}\equiv\delta_{aucl,\delta(\mathscr{\mathbf{T}})}(\varepsilon_{k})>0$,
(ii) a $\mathrm{measurable}$ set $A_{k}\subset B$ with $P(A_{k}^{c})<\varepsilon_{k}$,
(iii) an integer $h_{k}\geq1$, and (iv) a sequence of r.r.v.'s 
\begin{equation}
0=\tau_{k,0}<\tau_{k,1}<\cdots<\tau_{k,h-1}<\tau_{k,h}=1,\label{eq:temp-404-2}
\end{equation}
such that, for each $i=0,\cdots,h_{k}-1$, the function $X_{\tau(k,i)}^{N}$
is a r.v., and such that, (v) for each $\omega\in A_{k}$, we have
\begin{equation}
\bigwedge_{i=0}^{h(k)-1}(\tau_{k,i+1}(\omega)-\tau_{k,i}(\omega))\geq\delta_{k},\label{eq:temp-342-4}
\end{equation}
with
\begin{equation}
d(X^{N}(\tau_{k,i}(\omega),\omega),X^{N}(\cdot,\omega))\leq\varepsilon_{k}\label{eq:temp-307-4-1-4}
\end{equation}
on the interval $\theta_{k,i}(\omega)\equiv[\tau_{k,i}(\omega),\tau_{k,i+1}(\omega))$
or $\theta_{k,i}(\omega)\equiv[\tau_{k,i}(\omega),\tau_{k,i+1}(\omega)]$
according as $0\leq i\leq h_{k}-2$ or $i=h_{k}-1$. 

2. For each $k\geq0$, let 
\[
\varepsilon_{k}\equiv2^{-k}\wedge2^{-2}\delta_{f}(2^{-k}).
\]
Then Conditions (i-v) in the previous step hold. Separately, write
$n_{-1}\equiv0$. Inductively, for each $k\geq0$, take an integer
$n_{k}\geq n_{k-1}+1$ so large that 
\[
2^{-n(k)}<\delta_{k}\equiv\delta_{aucl,\delta(\mathscr{\mathbf{T}})}(\varepsilon_{k}).
\]

3. Now let $s\in(0,1]Q_{\infty}$ be arbitrary. Then there exists
$k\geq0$ so large that $s\in(0,1]Q_{n(k)}$ . Consider each $\omega\in A_{k}$.
Let 
\[
r\in G_{k,\omega}\equiv[0,s]\cap\bigcup_{i=0}^{h(k)-1}\theta_{k,i}(\omega)\cap Q_{\infty}
\]
be arbitrary. Then there exists $i=0,\cdots,h_{k}-1$ such that $r\in[0,s]\cap\theta_{k,i}(\omega)$.
Since, according to \ref{eq:temp-342-4}, we have
\[
|\theta_{k,i}(\omega)|\equiv\tau_{k,i+1}(\omega)-\tau_{k,i}(\omega)\geq\delta_{k}>2^{-n(k)},
\]
there exists some $t\in\theta_{k,i}(\omega)Q_{n(k)}$. Either $t\leq s$
or $s<t.$ 

Consider the case where $t\leq s$. Then $t,r\in\theta_{k,i}(\omega).$
Hence 
\[
d(X^{N}(\tau_{k,i}(\omega),\omega),X^{N}(t,\omega))\leq\varepsilon_{k}\leq2^{-2}\delta_{f}(2^{-k})
\]
and
\[
d(X^{N}(\tau_{k,i}(\omega),\omega),X^{N}(r,\omega))\leq\varepsilon_{k}\leq2^{-2}\delta_{f}(2^{-k}),
\]
according to inequality \ref{eq:temp-307-4-1-4}. Consequently 
\[
d(X_{t}^{N}(\omega),X_{r}^{N}(\omega))\equiv d(X^{N}(t,\omega),X^{N}(r,\omega))\leq2^{-1}\delta_{f}(2^{-k})<\delta_{f}(2^{-k}).
\]
Therefore
\[
f(X_{r}^{N}(\omega))<f(X_{t}^{N}(\omega))+2^{-k}\leq\bigvee_{u\in[0,s]Q(n(k)}f(X_{u}^{N}(\omega))+2^{-k},
\]
where the last inequality is because $t\in[0,s]Q_{n(k)}$. 

Now consider the case where $t>s$. Then $r,s\in\theta_{k,i}(\omega).$
Hence
\[
d(X^{N}(\tau_{k,i}(\omega),\omega),X^{N}(s,\omega))\leq\varepsilon_{k}\leq2^{-2}\delta_{f}(2^{-k})
\]
and
\[
d(X^{N}(\tau_{k,i}(\omega),\omega),X^{N}(r,\omega))\leq\varepsilon_{k}\leq2^{-2}\delta_{f}(2^{-k}),
\]
according to inequality \ref{eq:temp-307-4-1-4}. Consequently 
\[
d(X_{s}^{N}(\omega),X_{r}^{N}(\omega))\equiv d(X^{N}(s,\omega),X^{N}(r,\omega))\leq2^{-1}\delta_{f}(2^{-k})<\delta_{f}(2^{-k}).
\]
Therefore
\[
f(X_{r}^{N}(\omega))<f(X_{s}^{N}(\omega))+2^{-k}\leq\bigvee_{u\in[0,s]Q(n(k)}f(X_{u}^{N}(\omega))+2^{-k},
\]
where the last inequality is because $s\in[0,s]Q_{n(k)}$. 

Either way, we obtain
\[
f(X_{r}^{N}(\omega))\leq\bigvee_{u\in[0,s]Q(n(k)}f(X_{u}^{N}(\omega))+2^{-k},
\]
where $r\in G_{k,\omega}$ is arbitrary. Since the set $G_{k,\omega}$
is dense in 
\[
G_{\omega}\equiv[0,s]\cap domain(X^{N}(\cdot,\omega)),
\]
the last displayed inequality holds for each $r\in G_{\omega}$, thanks
to the right continuity of the function $X^{N}(\cdot,\omega)$. In
particular, 
\[
f(X_{r}^{N}(\omega))\leq\bigvee_{u\in[0,s]Q(n(k)}f(X_{u}^{N}(\omega))+2^{-k},
\]
for each $r\in[0,s]Q_{n(k+1)}\subset G_{\omega},$where $\omega\in A_{k}$
is arbitrary. Thus 
\[
\bigvee_{r\in[0,s]Q(n(k+1))}f(X_{r}^{N})\leq\bigvee_{u\in[0,s]Q(n(k)}f(X_{u}^{N})+2^{-k}
\]
on $A_{k}$. Consequently, 
\[
0\leq\bigvee_{r\in[0,s]Q(n(k+1))}f(X_{r}^{N})-\bigvee_{u\in[0,s]Q(n(k)}f(X_{u}^{N})\leq2^{-k+1}
\]
on $A_{k}$, where $P(A_{k}^{c})<\varepsilon_{k}\leq2^{-k}$, and
where $k\geq0$ is arbitrary. Since $\sum_{\kappa=0}^{\infty}2^{-\kappa}<\infty$,
it follows that the a.u.- and $L_{1}$-limit 
\begin{equation}
Y_{N,s}\equiv\lim_{\kappa\rightarrow\infty}\bigvee_{u\in[0,s]Q(n(\kappa))}f(X_{u}^{N})\label{eq:temp-273}
\end{equation}
exists is a r.r.v. Hence, for each $\omega$ in the full set $domain(Y_{N,s})$,
the supremum 
\[
\sup_{u\in[0,s]Q(\infty)}f(X_{u}^{N}(\omega))
\]
exists and is given by $Y_{N,s}(\omega)$. Now the function $X^{N}(\cdot,\omega)$
is right continuous for each $\omega$ in the full set $B$. Hence,
by right continuity, we have
\begin{equation}
\sup_{u\in[0,s]}f(X_{u}^{N})=\sup_{u\in[0,s]Q(\infty)}f(X_{u}^{N})\label{eq:temp-541}
\end{equation}
on the full set $domain(\sup_{u\in[0,s]Q(\infty)}f(X_{u}^{N}(\omega)))$.
Therefore $\sup_{u\in[0,s]}f(X_{u}^{N})$ is a well-defined r.r.v.,
where $N\geq0$ and $s\in(N,N+1]Q_{\infty}$are arbitrary. Moreover,
from equality \ref{eq:temp-273}, we see that $Y_{N,s}$ is the $L_{1}$-limit
of a sequence in $L^{(s)}$. Hence $\sup_{u\in[0,s]}f(X_{u}^{N})\in L^{(s)}.$ 

4. Next let $s\in\overline{Q}_{\infty}$ be arbitrary. Then $s\in[N,N+1]\overline{Q}_{\infty}$
for some $N\geq0$. There are three possibilities: (i) $s=0$, in
which case, trivially
\begin{equation}
\sup_{u\in[N,s]}f(X_{u})\equiv f(X_{0})\label{eq:temp-444}
\end{equation}
is a r.r.v., (ii) $N=0$ and $s\in(N,N+1]$, in which case, by Step
3 above, we have
\begin{equation}
\sup_{u\in[N,s]}f(X_{u})\equiv\sup_{u\in[N,s]}f(X_{u}^{N})\label{eq:temp-447}
\end{equation}
is a r.r.v, or (iii) where $N\geq1$, in which case, again by Step
3 above, we have
\[
\sup_{u\in[0,s])}f(X_{u})=\sup_{u\in[0,1]}f(X_{u})\vee\sup_{u\in[1,2]}f(X_{u})\cdots\vee\sup_{u\in[N-1,N]}f(X_{u})\vee\sup_{u\in[N,s])}f(X_{u})
\]
\begin{equation}
\equiv\sup_{u\in[0,1]}f(X_{u}^{0})\vee\sup_{u\in[0,1]}f(X_{u}^{1})\cdots\vee\sup_{u\in[0,1]}f(X_{u}^{N-1})\vee\sup_{u\in[0,s-N])}f(X_{u}^{N})\label{eq:temp-553}
\end{equation}
 is a r.r.v. In all three case, we see that 
\begin{equation}
V_{s}\equiv\sup_{u\in[0,s]}f(X_{u})\label{eq:temp-582}
\end{equation}
is a r.r.v., with 
\begin{equation}
V_{s}\in L^{(s)},\label{eq:temp-576}
\end{equation}
for each $s\in\overline{Q}_{\infty}$. Summing up, $V:\overline{Q}_{\infty}\times(\Omega,L,E)\rightarrow R$
is a nondecreasing real-valued process adapted to the filtration $\mathcal{L}$.

5. Since the set $\{V_{s}:s\in\overline{Q}_{\infty}\}$ of r.r.v.'s
is countable, there exists a countable subset $G$ of $(0,\infty)$
such that each point $a\in(0,\infty)G_{c}$ is a regular point of
the r.r.v.'s in the set $\{V_{s}:s\in\overline{Q}_{\infty}\}$, where
\[
G_{c}\equiv\{a\in(0,\infty):|a-b|>0\quad\mathrm{for\;each}\;b\in G\}
\]
is the metric complement of $G$ in $(0,\infty)$. Note from the equalities
\ref{eq:temp-541}, \ref{eq:temp-444}, \ref{eq:temp-447}, \ref{eq:temp-553},
and \ref{eq:temp-582}, that the distribution of the r.r.v. $V_{s}$
depends only on the function $f$ and on the family of marginal f.j.d.'s
of the process $X$, for each $s\in\overline{Q}_{\infty}$. Therefore
the set $G$ of exceptional points is completely determined by the
function $f$ and the marginal f.j.d.'s of the process $X$. 

6. Consider each $a\in G_{c}$. Then, by the definition of the countable
set $G$, we see that $a\in(a_{0},\infty)G_{c}$ is a regular point
of the r.r.v. $V_{s}$ for each $s\in\overline{Q}_{\infty}$. Hence
the set $(V_{s}<a)$ is $\mathrm{measurable}$ for each $s\in\overline{Q}_{\infty}$. 

7. Now let $M\geq1$ and $k\geq0$ be arbitrary. Recall that $\Delta_{n(k)}\equiv2^{-n(k)}$
and that $\overline{Q}_{n(k)}\equiv\{0,\Delta_{n(k)},2\Delta_{n(k)},\cdots\}$.
According to Definition \ref{Def. Simple First Exit time}, define
the r.r.v.
\begin{equation}
\eta_{k}\equiv\sum_{u\in(0,M]\overline{Q}(n(k))}u1_{(V(u)\geq a)}\prod_{s\in(0,M]\overline{Q}(n(k));s<u}1_{(V(s)<a)}+M\prod_{s\in(0,M]\overline{Q}(n(k))}1_{(V(s)<a)}\label{eq:temp-449}
\end{equation}
Because the real-valued process $V$ is nondecreasing, the defining
equality \ref{eq:temp-449} simplifies to 
\begin{equation}
\eta_{k}=\sum_{u\in(0,M]\overline{Q}(n(k))}u1_{(V(u)\geq a)}1_{(V(u-\Delta(n(k))<a)}+M1_{(V(M)<a)}.\label{eq:temp-142-2-1-1}
\end{equation}
In words, $\eta_{k}$ is the first time in $[0,M]\overline{Q}_{n(k)}$
for the real-valued nondecreasing process $V$ to exit the interval
$(a_{0},a)$, with $\eta_{k}$ set to $M$ if no such time exists.
Note that $\eta_{k}$ is a r.r.v. with values in the finite set $[0,M]\overline{Q}_{n(k)}$.
Moreover, from the defining equality \ref{eq:temp-142-2-1-1}, we
see that $(\eta_{k}=u)\in L^{(u)}$ for each $u\in[0,M]\overline{Q}_{n(k)}$.
Thus $\eta_{k}$ is a stopping time relative to the filtration $\mathcal{L}$. 

8. Since $\eta_{k}-\Delta_{n(k)}$ has values in $[0,M)\overline{Q}_{n(k)}\subset[0,M)\overline{Q}_{n(k+1)}$,
we have $V_{\eta(k)-\Delta(n(k))}<a$ by equality \ref{eq:temp-142-2-1-1}.
Consequently, equality \ref{eq:temp-142-2-1-1}, applied to the stopping
time $\eta_{k+1}$ in the place of $\eta_{k},$ yields
\begin{equation}
\eta_{k}-\Delta_{n(k)}<\eta_{k+1}.\label{eq:temp-420-2}
\end{equation}
Moreover, $\eta_{k}-\Delta_{n(k)},\eta_{k+1}$ have values in $\overline{Q}_{n(k+1)}$.
Therefore the strict inequality \ref{eq:temp-420-2} implies that
\[
\eta_{k}-\Delta_{n(k)}\leq\eta_{k+1}-\Delta_{n(k+1)}.
\]
Suppose, for the sake of a contradiction, that $u\equiv\eta_{k}(\omega)<\eta_{k+1}(\omega)$
for some $\omega\in\Omega$. Then $\eta_{k+1}<M$, whence $V_{u}(\omega)\geq a$
by the defining formula \ref{eq:temp-142-2-1-1}. Consequently, since
$u\in(0,M)\overline{Q}_{n(k+1)}$, we have $\eta_{k+1}(\omega)\leq u$
by applying the defining formula \ref{eq:temp-142-2-1-1} to $\eta_{k+1}(\omega)$.
This is a contradiction. Thus $\eta_{k+1}(\omega)\leq\eta_{k}(\omega)$.
Roughly speaking, if we sample more frequently, we observe any exit
no later, 

9. Combining, we get
\begin{equation}
\eta_{k}-\Delta_{n(k)}\leq\eta_{k+1}-\Delta_{n(k+1)}<\eta_{k+1}\leq\eta_{k}.\label{eq:temp-271}
\end{equation}
Iterating, we obtain 
\begin{equation}
\eta_{k}-\Delta_{n(k)}\leq\eta_{\kappa}-\Delta_{n(\kappa)}\leq\eta_{\kappa}\leq\eta_{k}\label{eq:temp-293}
\end{equation}
for each $\kappa\geq k+1$. It follows that $\eta_{\kappa}\downarrow\tau$
a.u. for some r.r.v. $\tau$ with 
\begin{equation}
\eta_{k}-\Delta_{n(k)}\leq\tau\leq\eta_{k}\label{eq:temp-293-1}
\end{equation}
where $k\geq0$ is arbitrary,. Consequently, Assertion 2 of Proposition
\ref{Prop. Stopping times rlative to right continuous filtration}
implies that $\tau$ is a stopping time relative to the filtration
$\mathcal{L}$. Moreover, note that $\Delta_{n(k)}\equiv2^{-n(k)}\leq2^{-k}<2^{-k+2}.$
Hence inequality \ref{eq:temp-293-1} implies that 
\begin{equation}
\tau\leq\eta_{k}<\tau+2^{-k+2}.\label{eq:temp-287}
\end{equation}
Therefore Lemma \ref{Lem. Observability at stopping time} implies
$X_{\tau}$ is a well defined r.r.v., and that $X_{\eta(h)}\rightarrow X_{\tau}$
a.u. as $k\rightarrow\infty$. 

10. Now consider each $\omega\in domain(X_{\tau})$. Then $\omega\in domain(\tau)$
and 
\[
\tau(\omega)\in domain(X(\cdot,\omega)).
\]
Let $t\in domain(X(\cdot,\omega))\cap[0,\tau(\omega))$ be arbitrary.
Then, for some sufficiently large $k\geq0$, we have 
\[
t<\tau(\omega)-2^{-n(k)+1}\leq\eta_{k}(\omega)-2\Delta_{n(k)}
\]
Hence there exists $s\in(t,\eta_{k}(\omega))\overline{Q}_{n(k)}$.
Because $s\in\overline{Q}_{n(k)}$ is before the simple first exit
time $\eta_{k}(\omega)$, we have $V_{s}(\omega)<a$, as can more
precisely be verified from to the defining equality \ref{eq:temp-449}.
Consequently,
\[
f(X(t,\omega))\leq\sup_{u\in[0,s]}f(X_{u}(\omega))\equiv V_{s}(\omega)<a.
\]
where $t\in domain(X(\cdot,\omega))\cap[0,\tau(\omega))$ is arbitrary.
Condition (i) of Definition \ref{Def. First exit time} is proved
for the stopping time $\tau$ to be the first exit time\emph{ }in\emph{
}$[0,M]$ of the open subset $(f<a)$. 

12. We will next verify Condition (ii) of Definition \ref{Def. First exit time}.
To that end, suppose \emph{$\tau(\omega)<M$. }Then, in view of inequality
\ref{eq:temp-293-1}, there exists $k\geq0$ so large that $\tau(\omega)\leq r\equiv\eta_{k}(\omega)<M$.
Because successful exit occurs at $r\in[0,M)\overline{Q}_{n(k)}$,
we have $V_{r}(\omega)\geq a$, as can more precisely be verified
from to the defining equality \ref{eq:temp-449}. Consequently,
\[
f(X_{\tau}(\omega))\equiv f(X_{\tau(\omega)}(\omega))\equiv f(X_{r}(\omega))\leq\sup_{u\in[0,r]}f(X_{u}(\omega))\equiv V_{r}(\omega)\equiv V_{\eta(k)}(\omega)\geq a,
\]
Condition (ii) of Definition \ref{Def. First exit time} is also verified
for the stopping time $\tau$ to be the first exit time\emph{ }in\emph{
}$[0,M]$ of the open subset $(f<a)$. Accordingly, $\overline{\tau}_{f,a,M}(X)\equiv\tau$
is the first exit time\emph{ }in\emph{ }$[0,M]$ of the open subset
$(f<a)$. The theorem is proved.
\end{proof}

\section{Feller Semigroup and Feller Process}
\begin{defn}
\label{Specification of locally compact metric space etc)} \textbf{(Specification
of locally compact state space and related objects).} In this section,
let $(S,d)$ be a locally compact metric space, as specified in Definition
\ref{Def.  Specification of state space}, along with related objects
including a reference point $x_{\circ}\in S$ and a binary approximation
$\xi$. In addition, for each $n\geq0$ and for each $y\in S$, define
the function 
\[
h_{y,n}\equiv(1\wedge(1+n-d(\cdot,y))_{+}\in C(S,d).
\]
Clearly, for each fixed $y\in S$ , we have $h_{y,n}\uparrow1$ as
$n\rightarrow\infty$, uniformly on compact subsets of $(S,d)$. Define
\[
\overline{h}_{y,n}\equiv1-h_{y,n}\in C_{u,b}(S,d).
\]
The continuous functions $h_{y,n}$ and $\overline{h}_{y,n}$ are
surrogates for the indicators $1_{(d(y,\cdot)\leq n)}$ and $1_{(d(y,\cdot)>n)}$
respectively. 

Let $(\overline{S},\overline{d})$$\equiv(S\cup\{\Delta\},\overline{d})$
be a one-point compactification of the metric space \emph{(S,d}),
where $\overline{d}\leq1$ and where $\Delta$ is called the point
at infinity. For ease of reference, we list almost verbatim from Definition
\ref{Def. One point compactification}, the conditions for the one-point
compactification  $(\overline{S},\overline{d})$. 

\emph{Condition} $1$. $S\cup\{\Delta\}$ is dense in $(\overline{S},\overline{d})$.
Moreover, $\overline{d}\leq1$.

\emph{Condition} 2. For each compact subset $K$ of $(S,d)$, there
exists $c>0$ such that $\overline{d}(x,\Delta)\geq c$ for each $x\in K$. 

\emph{Condition} 3. Let $K$ be an arbitrary compact subset  of $(S,d)$.
Let $\varepsilon>0$ be arbitrary. Then there exists $\delta_{K}(\varepsilon)>0$
such that for each $y\in K$ and $z\in S$ with $\overline{d}(y,z)<\delta_{K}(\varepsilon)$,
we have $d(y,z)<\varepsilon$. In particular,  the identity mapping
$\bar{\iota}:(S,\overline{d})\rightarrow(S,d)$ is uniformly continuous
on each compact subset of $S$.

\emph{Condition} 4. The identity mapping $\iota:(S,d)\rightarrow(S,\overline{d})$,
defined by $\iota(x)\equiv x$ for each $x\in S$, is uniformly continuous
on $(S,d)$. In other words, for each $\varepsilon>0$, there exists
$\delta_{\overline{d}}(\varepsilon)>0$ such that $\overline{d}(x,y)<\varepsilon$
for each $x,y\in S$ with $d(x,y)<\delta_{\overline{d}}(\varepsilon)$.

\emph{Condition} 5. For each $n\geq1$, we have 
\[
(d(\cdot,x_{\circ})>2^{n+1})\subset(\overline{d}(\cdot,\Delta)\leq2^{-n}).
\]

Refer to Definition \ref{Def. Notations for dyadic rationals-1} for
notations related to the enumerated sets $\overline{Q}_{0}$,$\overline{Q}_{1}\cdots,\overline{Q}_{\infty}$
of dyadic rationals in $[0,\infty)$, and to the enumerated sets $Q_{0}$,$Q_{1}\cdots,Q_{\infty}$
of dyadic rationals in $[0,1]$.
\end{defn}
$\square$
\begin{defn}
\label{Def. Feller-Semigroup}\textbf{(Feller Semigroup).} Let $\mathbf{V}\equiv\{V_{t}:t\in[0,\infty)\}$
be an arbitrary family of nonnegative linear mappings from from $C_{ub}(S,d)$
to $C_{ub}(S,d)$ such that $V_{0}$ is the identity mapping. Suppose,
for each $t\in[0,\infty)$ and for each $y\in S$, that the function
\[
V_{t}^{y}\equiv V_{t}(\cdot)(y):C_{ub}(S,d)\rightarrow R
\]
is a distribution on the locally compact space $(S,d)$. Suppose,
in addition, the following four conditions are satisfied.

1. (Smoothness). For each $N\geq1$, for each $t\in[0,N]$, and for
each $f\in C_{ub}(S,d)$ with a modulus of continuity $\delta_{f}$
and with $|f|\leq1$, the function $V_{t}f\in C_{ub}(S,d)$ has a
modulus of continuity $\alpha_{\mathbf{V,\mathit{N}}}(\delta_{f})$
which depends on the finite interval $[0,N]$ and on $\delta_{f}$,
and otherwise not on the function $f$. 

2. (Semigroup property). For each $s,t\in[0,\infty)$, we have $V_{t+s}=V_{t}V_{s}$. 

3. (Strong continuity). For each $f\in C_{ub}(S,d)$ with a modulus
of continuity $\delta_{f}$ and with $|f|\leq1$, and for each $\varepsilon>0$,
there exists $\delta_{\mathbf{V}}(\varepsilon,\delta_{f})>0$ so small
that, for each $t\in[0,\delta_{\mathbf{V}}(\varepsilon,\delta_{f}))$,
we have 
\begin{equation}
|f-V_{t}f|\leq\varepsilon\label{eq:temp-382-1-2-2-3}
\end{equation}
as functions on $(S,d)$.

4. (Non-explosion). For each $N\geq1$, for each $t\in[0,N]$, and
for each $\varepsilon>0$, there exists an integer $\kappa_{\mathbf{V,\mathit{N}}}(\varepsilon)>0$
so large that, if $n\geq\kappa_{\mathbf{V,\mathit{N}}}(\varepsilon)$
then 
\[
V_{t}^{y}\overline{h}_{y,n}\leq\varepsilon
\]
for each $y\in S$.

Then we call the family\emph{ }$\mathbf{V}$\emph{ }a \emph{Feller
semigroup}\index{Feller semigroup}. The operation $\delta_{\mathbf{V}}$
is called a \index{modulus of strong continuity of Feller semigroup}\emph{
modulus of strong continuity} of\emph{ }$\mathbf{V}$. The sequence
$\alpha_{\mathbf{V}}\equiv(\alpha_{\mathbf{V,\mathrm{N}}})_{N=1,2,\cdots}$
of operations is called a \emph{modulus of smoothness of }\index{modulus of smoothness of Feller semigroup}
$\mathbf{V}$. The sequence $\kappa_{\mathbf{V}}\equiv(\kappa_{\mathbf{V,\mathrm{N}}})_{N=1,2,\cdots}$
of operations is called \emph{a modulus of non-explosion}\index{modulus of non-explosion of Feller semigroup}
of $\mathbf{V}$.$\square$

In order to use results developed in previous sections for Markov
semigroups and Markov processes, where the state space is assumed
to be compact, we embed each given Feller semigroup on the locally
compact state space $(S,d)$ into a Markov semigroup on the one-point
compactification $(\overline{S},\overline{d})$ state space, in the
following sense.
\end{defn}
\begin{thm}
\label{Thmp. Compactification-of-Feller semigoup into Markoc semgroup}
\textbf{\emph{(Compactification of Feller semigroup into a Markov
semigroup with a compact state space). }}Let $\mathbf{V}\equiv\{V_{t}:t\in[0,\infty)\}$
be an arbitrary Feller semigroup on the locally compact metric space
$(S,d),$ with moduli $\delta_{\mathbf{V}}$, $\alpha_{\mathbf{V}}$,
$\kappa_{\mathbf{V}}$ as in Definition \ref{Def. Feller-Semigroup}.
Let $t\in[0,\infty)$ be arbitrary. Define the function
\[
T_{t}:C(\overline{S},\overline{d})\rightarrow C(\overline{S},\overline{d})
\]
by
\[
(T_{t}g)(\Delta)\equiv T_{t}^{\Delta}g\equiv g(\Delta),
\]
and by
\begin{equation}
(T_{t}g)(y)\equiv T_{t}^{y}g\equiv\intop_{z\in\overline{S}}T_{t}^{y}(dz)g(z)\equiv\intop_{z\in S}V_{t}^{y}(dz)g(z)\label{eq:temp-542}
\end{equation}
for each $y\in S$, for each $g\in C(\overline{S},\overline{d})$.
Equality \ref{eq:temp-542} is equivalent to 
\[
T_{t}^{y}g\equiv V_{t}^{y}(g|S)\equiv V_{t}(g|S)(y).
\]
for each $y\in S$, for each $g\in C(\overline{S},\overline{d})$. 

Then the family $\mathbf{T}\equiv\{T_{t}:t\in[0,\infty)\}$ is a Markov
semigroup with state space $(\overline{S},\overline{d})$. This Markov
semigroup $\mathbf{T}$ will be called the \emph{\index{compactification of Feller semigroup}
compactification }of the given Feller semigroup $\mathbf{V}$.
\end{thm}
\begin{proof}
Let $t\in[0,\infty)$ be arbitrary. Let $N\geq1$ be such that $t\in[0,N]$.
Let $g\in C(\overline{S},\overline{d})$ be arbitrary, with modulus
of continuity $\overline{\delta}_{g}$. There is no loss of generality
in assuming that $g$ has values in $[0,1]$. For abbreviation, write
$f\equiv g|S\in C_{ub}(S,d)$.

Let $\varepsilon>0.$ be arbitrary. 

1. Let $\delta_{\overline{d}}$ be the operation listed in Condition
4 of Definition \ref{Specification of locally compact metric space etc)}.
Take arbitrary points ${\normalcolor y,z\in S}$ with $d(y,z)<\delta_{\overline{d}}(\overline{\delta}_{g}(\varepsilon))$.
Then, according to Condition 4 of Definition \ref{Specification of locally compact metric space etc)},
we have $\overline{d}(y,z)<\overline{\delta}_{g}(\varepsilon)$. Hence
\[
|f(y)-f(z)|=|g(y)-g(z)|<\varepsilon.
\]
Thus the function $f\equiv g|S$ has a modulus of continuity $\delta_{\overline{d}}\circ\overline{\delta}_{g}$.
Therefore, according to the definition of the modulus of smoothness
$\alpha_{\mathbf{V}}$ in Definition \ref{Def. Feller-Semigroup},
the function $V_{t}f\in C_{ub}(S,d)$ has modulus of continuity $\alpha_{\mathbf{V,\mathit{N}}}(\delta_{\overline{d}}\circ\overline{\delta}_{g})$.

2. Let $k\geq0$ be so large that 
\begin{equation}
2^{-k+1}<\overline{\delta}_{g}(\varepsilon)\label{eq:temp-540}
\end{equation}
Let 
\[
n\equiv2^{k+1}\vee\kappa_{\mathbf{V,\mathit{N}}}(\varepsilon).
\]
Then 
\begin{equation}
V_{t}^{y}\overline{h}_{y,n}=\pm\varepsilon\label{eq:temp-560}
\end{equation}
for each $y\in S$.

3. By Condition 5 of Definition \ref{Specification of locally compact metric space etc)},
we have, for each $u\in S$, if $d(u,x_{\circ})>n=2^{k+1}$ then $\overline{d}(u,\Delta)\leq2^{-k}<\overline{\delta}_{g}(\varepsilon)$,
whence $f(u)=g(u)=g(\Delta)\pm\varepsilon$. 

4. Take an arbitrary $a\in(2n+1,2n+2)$ such that the set 
\[
K\equiv\{u\in S:d(x_{\circ},u)\leq a\}
\]
is a compact subset of $(S,d)$. Let 
\begin{equation}
K^{c}\equiv\{u\in S:d(x_{\circ},u)>2n+1\}\cup\{\Delta\}.\label{eq:temp-593}
\end{equation}
Then $\overline{S}=K\cup K^{c}$.

5. Define 
\[
\varepsilon'\equiv\alpha_{\mathbf{V,\mathit{N}}}(\delta_{\overline{d}}\circ\overline{\delta}_{g})(\varepsilon).
\]
Then, by Condition 3 of Definition \ref{Specification of locally compact metric space etc)},
there exists $\delta_{K}(\varepsilon')>0$ such that, for each $y\in K$
and $z\in S$ with 
\[
\overline{d}(y,z)<\delta_{K}(\varepsilon'),
\]
we have $d(y,z)<\varepsilon'$, whence, in view of the last statement
in Step 1, we have
\begin{equation}
|V_{t}^{y}f-V_{t}^{y}f|<\varepsilon.\label{eq:temp-592}
\end{equation}

6. Let $y\in K^{c}$ be arbitrary. Suppose $y\in S$. Then $d(x_{\circ},y)>2n+1$
by the defining equality \ref{eq:temp-593}. Therefore, for each point
$u\in S$ with $h_{y,n}(u)>0$, we have $d(y,u)\leq n+1,$ and so
$d(u,x_{\circ})>n$ . In view of Step 3, it follow that, for each
point $u\in S$ with $h_{y,n}(u)>0$, we have $f(u)=g(\Delta)\pm\varepsilon$.
Therefore
\[
(T_{t}g)(y)=V_{t}^{y}f=V_{t}^{y}h_{y,n}f+V_{t}^{y}\overline{h}_{y,n}f=V_{t}^{y}h_{y,n}f\pm V_{t}^{y}\overline{h}_{y,n}
\]
\[
=V_{t}^{y}h_{y,n}f\pm\varepsilon=V_{t}^{y}((g(\Delta)\pm\varepsilon)h_{y,n}\pm\varepsilon=g(\Delta)V_{t}^{y}h_{y,n}\pm2\varepsilon
\]
\[
=g(\Delta)V_{t}^{y}(1-\overline{h}_{y,n})\pm2\varepsilon=g(\Delta)\pm\varepsilon\pm2\varepsilon
\]
\[
=g(\Delta)\pm3\varepsilon,
\]
where we have used equality \ref{eq:temp-560} twice. In short
\begin{equation}
(T_{t}g)(y)=g(\Delta)\pm3\varepsilon\label{eq:temp-590}
\end{equation}
Suppose $y=\Delta.$ Then trivially equality \ref{eq:temp-590} also
holds. Combining equality \ref{eq:temp-590} also holds for each $y\in K^{c}$
.

7. Proceed to examine a pair of points $y,z\in\overline{S}$ with
\begin{equation}
\overline{d}(y,z)<\delta_{K}(\alpha_{\mathbf{V,\mathit{N}}}(\delta_{\overline{d}}\circ\overline{\delta}_{g})(\varepsilon))\equiv\delta_{K}(\varepsilon'),\label{eq:temp-594}
\end{equation}
 Suppose first that $y,z\in K^{c}.$ Then, by inequality \ref{eq:temp-590},
we have
\[
|(T_{t}g)(y)-(T_{t}g)(y)|\leq6\varepsilon.
\]

8. Suppose, on the other hand, that $y\in K$ or $z\in K$ . Then,
by inequality \ref{eq:temp-592}, we obtain.
\begin{equation}
|(T_{t}g)(y)-(T_{t}g)(y)|=|V_{t}^{y}f-V_{t}^{y}f|<\varepsilon.\label{eq:temp-592-1}
\end{equation}

9. Combining, we see that, in either case, we have
\[
(T_{t}g)(y)-(T_{t}g)(z)|=|V_{t}^{y}f-V_{t}^{y}f|\leq6\varepsilon,
\]
provided that the bound \ref{eq:temp-594} holds. 

10. Summing up, the function $T_{t}g$ is continuous on $(\overline{S},\overline{d})$,
with a modulus of continuity $\alpha_{\mathbf{T,N}}(\overline{\delta}_{g})$
defined by
\[
\alpha_{\mathbf{T,\mathit{N}}}(\overline{\delta}_{g})\equiv\delta_{K}\circ\alpha_{\mathbf{V,\mathit{N}}}\circ\delta_{\overline{d}}(\overline{\delta}_{g})
\]
where $\overline{\delta}_{g}$ is the modulus of continuity of the
arbitrary function $g\in C(\overline{S},\overline{d})$.

11. In particular, $T_{t}:C(\overline{S},\overline{d})\rightarrow C(\overline{S},\overline{d})$
is a well-defined function. From the defining equality \ref{eq:temp-542},
it is obvious that it is a nonnegative linear function, with $T_{t}1=1$.
Hence, for each ${\normalcolor y\in\overline{S}}$, the linear, nonnegative,
function $T_{t}^{y}$ is an integration with $T_{t}^{y}1=1$. Moreover,
$T_{t}g$ for each $s,t\in[0,\infty)$ with modulus of continuity
$\overline{\delta}_{g}$, the function $T_{t}g$ has a modulus of
continuity $\alpha_{\mathbf{T,\mathit{N}}}(\overline{\delta}_{g})$
defined in Step 10. We conclude that $T_{t}$ is a transition distribution
from $(\overline{S},\overline{d})$ to $(\overline{S},\overline{d})$,
where $N\geq1$ and $t\in[0,N]$ are arbitrary. It is also clear from
the defining equality \ref{eq:temp-542} that $T_{0}$ is the identity
mapping. 

12. It remains to verify the conditions in Definition \ref{Def. Markov semigroup  semigroup}
for the family $\mathbf{T}\equiv\{T_{t}:t\in[0,\infty)\}$ to be a
Markov semigroup. The smoothness condition follows immediately from
Step 10, where we found that the operation $\alpha_{\mathbf{T,\mathit{N}}}$
is a modulus of smoothness for the transition distribution $T_{t}$,
for each each $N\geq1$, for each $t\in[0,N]$. 

13. For the semigroup property, consider each $s,t\in[0,\infty)$.
Let $y\in S$ be arbitrary. Then, by inequality \ref{eq:temp-560},
we have 
\begin{equation}
T_{s}^{y}h_{y,k}=V_{s}^{y}h_{y,k}\uparrow1\label{eq:temp-560-1}
\end{equation}
and $h_{y,k}\uparrow1_{S}$ as $k\rightarrow\infty$. Consequently,
$S$ is a full subset, and $\{\Delta\}$ is a null subset of $\overline{S}$
relative to the distribution $T_{s}^{y}.$ Hence according to the
defining equality \ref{eq:temp-542}, we have 
\[
T_{t}^{y}g\equiv V_{t}^{y}(g|S)+g(\Delta)1_{\{\Delta\}}=V_{t}^{y}(g|S)+0
\]
as integrable functions on $(\overline{S},\overline{d})$ relative
to the distribution $T_{s}^{y}$. Therefore
\[
T_{s}^{y}(T_{t}g)\equiv T_{s}^{y}(V_{t}(g|S)1_{S}+g(\Delta)1_{\{\Delta\}})=T_{s}^{y}(V_{t}(g|S))
\]
\begin{equation}
=V_{s}^{y}(V_{t}(g|S)|S))=V_{s}^{y}(V_{t}(g|S))=V_{s+t}^{y}(g|S),\label{eq:temp-571}
\end{equation}
where the last equality is by the semigroup property of the Feller
semigroup $V$. Applying equality \ref{eq:temp-571}, with $t,s,$
replaced by $0,t+s$ respectively, we obtain $T_{s+t}^{y}(g)=V_{s+t}^{y}(g|S)$.
Combining with equality \ref{eq:temp-571}, we in turn obtain
\[
T_{s}^{y}(T_{t}g)=T_{s+t}^{y}(g)
\]
where $y\in S$ is arbitrary. At the same time, we have, trivially
\[
T_{s}^{\Delta}(T_{t}g)\equiv(T_{t}g)(\Delta)\equiv g(\Delta)=T_{s+t}^{\Delta}(g).
\]
Thus we proved that $T_{s}(T_{t}g)=T_{s+t}(g)$ on $(\overline{S},\overline{d})$,
where $g\in C(\overline{S},\overline{d})$ is arbitrary. The semigroup
property is proved for the family $\mathbf{T}$. 

14. It remains to verify strong continuity of the family $\mathbf{T}$.
To that end, let $\varepsilon>0$ be arbitrary, and let $g\in C(\overline{S},\overline{d})$
be arbitrary, with a modulus of continuity $\overline{\delta}_{g}$
and with $\left\Vert g\right\Vert \leq1$.We wish to prove that 
\begin{equation}
\left\Vert g-T_{t}g\right\Vert \leq\varepsilon,\label{eq:temp-382-1-2-2-2-1}
\end{equation}
provided that $t\in[0,\delta_{\mathbf{T}}(\varepsilon,\overline{\delta}_{g}))$.
First note that, trivially, 
\begin{equation}
g(\Delta)-(T_{t}g)(\Delta)=0.\label{eq:temp-572}
\end{equation}
Next, recall from Step 1 that the function $g|S$ has a modulus of
continuity $\delta_{\overline{d}}\circ\overline{\delta}_{g}$. Hence,
by the strong continuity of the Feller semigroup , there exists $\delta_{\mathbf{V}}(\varepsilon,\delta_{\overline{d}}\circ\overline{\delta}_{g})>0$
so small that, for each $t\in[0,\delta_{\mathbf{V}}(\varepsilon,\delta_{\overline{d}}\circ\overline{\delta}_{g}))$,
we have 
\begin{equation}
|(g|S)-V_{t}(g|S)|\leq\varepsilon.\label{eq:temp-382-1-2-2-3-1-1}
\end{equation}
Define
\[
\delta_{\mathbf{T}}(\varepsilon,\overline{\delta}_{g})\equiv\delta_{\mathbf{V}}(\varepsilon,\delta_{\overline{d}}\circ\overline{\delta}_{g})
\]
Then, for each $y\in S$, we have
\[
|T_{t}^{y}g-g(y)|\equiv|V_{t}^{y}(g|S)-g(y)|\leq\varepsilon.
\]
Combining with equality \ref{eq:temp-572}, we obtain
\[
\left\Vert T_{t}g-g\right\Vert \leq\varepsilon,
\]
provided that $t\in[0,\delta_{\mathbf{T}}(\varepsilon,\overline{\delta}_{g}))$,
where $g\in C(\overline{S},\overline{d})$ is arbitrary, with a modulus
of continuity $\overline{\delta}_{g}$ and with $\left\Vert g\right\Vert \leq1$.Thus
we have verified also the strong continuity condition in Definition
\ref{Def. Markov semigroup  semigroup} for the family $\mathbf{T}\equiv\{T_{t}:t\in[0,\infty)\}$
to be a Markov semigroup. 
\end{proof}
\begin{cor}
\emph{\label{Cor:(Nonexplosion-of-Feller process at stopping time}}
\textbf{\emph{(Nonexplosion of Feller process in finite time intervals).}}
Let $\mathbf{V}\equiv\{V_{t}:t\in[0,\infty)\}$ be an arbitrary Feller
semigroup on the locally compact metric space $(S,d),$ with moduli
$\delta_{\mathbf{V}}$, $\alpha_{\mathbf{V}}$, $\kappa_{\mathbf{V}}$
as in Definition \ref{Def. Feller-Semigroup}. Let $\mathbf{T}\equiv\{T_{t}:t\in[0,\infty)\}$
be the\emph{ }compactification\emph{ }of $\mathbf{V}$. Let 
\[
X\equiv X^{x,\mathbf{T}}:[0,\infty)\times(\Omega,L,E)\rightarrow(\overline{S},\overline{d})
\]
be an a.u. càdlàg Markov process generated by the initial state x
and semigroup $\mathbf{T}$, as constructed in \ref{Thm. Construction of Markov process from  initial distrinbution and semigrp-1-1}.
Let $\mathcal{L}\equiv\{L^{(t)}:t\in[0,\infty)\}$ an arbitrary right
continuous filtration to which $X^{x,\mathbf{T}}$ is adapted. All
stopping times will be understood to be relative to this filtration
$\mathcal{L}$. Theorem \ref{Thm. Each a.u. cadlag Markov process  is a strong Markov process}
says that $X^{x,\mathbf{T}}$ is a strong Markov process relative
to the filtration $\mathcal{L}$.

Let $\tau$ be an arbitrary stopping time with values in $[0,M]$.
Then the following holds.

1. Let $\varepsilon>0$ be arbitrary. Then there exists $c>0$ so
small that 
\[
P(\overline{d}(\Delta,X_{\tau})>c)>1-\varepsilon.
\]
And there exist a compact subset $K$ of $(S,d)$ which depends on
$\varepsilon,x$ and $\mathbf{V}$, such that
\[
P(X_{\tau}\in K)>1-\varepsilon.
\]
Since $\varepsilon>0$ is arbitrary, we have $P(X_{\tau}\in S)=1$.

2. For each $y\in S$, and for each $N\geq M$, we have
\begin{equation}
T_{M-\tau}^{y}\overline{h}_{y,N}=V_{M-\tau}^{y}\overline{h}_{y,N}\label{eq:temp-578-1-1-1}
\end{equation}
as r.r.v.'s

3. For each $y\in S$, and for each $N\geq M\vee\kappa_{\mathbf{V,\mathit{N}}}(\varepsilon)$,
we have
\begin{equation}
T_{M-\tau}^{y}\overline{h}_{y,N}=V_{M-\tau}^{y}\overline{h}_{y,N}\leq\varepsilon\label{eq:temp-578-1-1-1-1}
\end{equation}
as r.r.v.'s

4. For each $v\in[0,\infty)$, the function $X_{v}:(\Omega,L,E)\rightarrow(S,d)$
is a r.v.
\end{cor}
\begin{proof}
1. Lemma \ref{Lem. Observability at stopping time} says that the
function $X_{\tau}$ is a well defined r.v. with values in $(\overline{S},\overline{d})$
which is measurable relative to $L^{(\tau)}$, and that there exists
a nonincreasing sequence $(\eta_{j})_{j=0,1,\cdots}$ of stopping
times such that, (i) for each $j\geq0$, the r.r.v. $\eta_{j}$ has
values in $\overline{Q}_{j}$, (ii) 
\begin{equation}
\tau\leq\eta_{j}<\tau+2^{-j+2}.\label{eq:temp-421-1-2-1-2-2}
\end{equation}
and (iii) $X_{\eta(j)}\rightarrow X_{\tau}$ a.u. as $j\rightarrow\infty$.
By assumption, we have $\tau\leq M$. Hence, by replacing $\eta_{j}$
with $\eta_{j}\wedge M$ if necessary, we may assume that $\eta_{j}$
has values in the finite set $\overline{Q}_{j}[0,M]$, for each $j\geq0$.

2. Since $X_{\eta(j)}\rightarrow X_{\tau}$ a.u. and since $\overline{d}\leq1$,
we have $E\overline{d}(X_{\eta(j)},X_{\tau})\rightarrow0$. 

3. Let $\varepsilon>0$ be arbitrary. Consider each $n\geq\kappa_{\mathbf{V,\mathit{M}}}(\varepsilon).$
Let $v\in[0,M]$ be arbitrary. Since $x\in S$, we have
\[
E(h_{x,n}(X_{v}))=T_{v}^{x}h_{x,n}=V_{v}^{x}h_{x,n}=1-V_{v}^{x}\overline{h}_{x,n}\geq1-\varepsilon.
\]
Now let $v\in\overline{Q}_{j}[0,M]$ be arbitrary. Take any regular
point $a\in(n-1,n)$ of the r.r.v.'s in the finite family $\{\overline{d}(\Delta,X_{v}):\in\overline{Q}_{j}[0,M]\}$,
such that the subset $K\equiv(d(\cdot,x)\leq a)\subset(S,d)$ is a
compact subset of $(S,d),$ and such that the set $(d(X_{v},x)\leq a)$
is measurable. Then
\[
P(X_{v}\in K)=P(d(X_{v},x)\leq a)\geq E(h_{x,n}(X_{v}))\geq1-\varepsilon.
\]

4. Condition 2 of Definition \ref{Specification of locally compact metric space etc)}
implies that there exists $c>0$ such that $\overline{d}(x,\Delta)>c$
for each $x\in K$. Then
\[
P(\overline{d}(\Delta,X_{v})>c)\geq P(X_{v}\in K)>1-\varepsilon.
\]
Hence, for each $j\geq0$, we have 
\[
P(\overline{d}(\Delta,{\normalcolor X_{\eta(j)}})>c)=\sum_{v\in\overline{Q}(j)[0,M]}P(\eta_{j}=v;\overline{d}(\Delta,{\normalcolor X_{\eta(j)}})>c)
\]
\[
=\sum_{v\in\overline{Q}(j)[0,M]}P(\eta_{j}=v;\overline{d}(\Delta,X_{v})>c)
\]
\[
>\sum_{v\in\overline{Q}(j)[0,M]}P(\eta_{j}=v)(1-\varepsilon)=1-\varepsilon.
\]
Since $\overline{d}(X_{\eta(j)},X_{\tau})\rightarrow0$ a.u. and since
$\overline{d}\leq1$, it follows that
\[
P(\overline{d}(\Delta,X_{\tau})>c)\geq1-\varepsilon.
\]
Similarly, $P(X_{\tau}\in K)>1-\varepsilon,$ where $\varepsilon>0$
is arbitrary. Consequently $P(X_{\tau}\in S)=1$. Assertion 1 of the
corollary is proved.

5. Proceed to prove Assertion 2. Let $y\in S$ and $N\geq M$ be arbitrary.
Then
\[
T_{M-\eta(j)}^{y}\overline{h}_{y,N}=\sum_{v\in\overline{Q}(j)[0,M]}1_{\eta(j)=v}T_{M-v}^{y}\overline{h}_{y,N}
\]
\begin{equation}
=\sum_{v\in\overline{Q}(j)[0,M]}1_{\eta(j)=v}V_{M-v}^{y}\overline{h}_{y,N}=V_{M-\eta(j)}^{y}\overline{h}_{y,N}.\label{eq:temp-578-1-1}
\end{equation}
Letting $j\rightarrow\infty$ and $\eta_{j}\downarrow\tau,$ the uniform
continuity of $T_{w}^{y}$ in the subscript variable $w\in[0,M]$
yields 
\[
T_{M-\tau}^{y}\overline{h}_{y,N}=V_{M-\tau}^{y}\overline{h}_{y,N}.
\]
Assertion 2 is proved.

6. Now suppose, in addition, that $N\geq M\vee\kappa_{\mathbf{V,\mathit{N}}}(\varepsilon).$
Then $V_{M-v}^{y}\overline{h}_{y,N}\leq\varepsilon$ for each $v\in\overline{Q}_{j}[0,M]$.
Hence equality \ref{eq:temp-578-1-1} yields $T_{M-\eta(j)}^{y}\overline{h}_{y,N}=V_{M-\eta(j)}^{y}\overline{h}_{y,N}\leq\varepsilon$.
With $j\rightarrow\infty$, we obtain the desired inequality \ref{eq:temp-578-1-1-1-1}.
Assertion 3 and is proved.

7. Finally, let $v\in[0,\infty)$ be arbitrary. Let $M\geq1$ be so
large that $v\in[0,M]$. Then, in view of Assertion 1, on a full set,
$X_{v}$ is a function with values in $S$. Let $f\in C_{ub}(S,d)$
be arbitrary with $1\geq f\geq0$. Then 
\[
E((\overline{h}_{x(\circ),n}f)(X_{v}))=T_{v}^{x}(\overline{h}_{x(\circ),n}f)=V_{v}^{x}(\overline{h}_{x(\circ),n}f)\leq V_{v}^{x}(\overline{h}_{x(\circ),n})\downarrow0.
\]
Therefore the limit $E((h_{x(\circ),n}f)(X_{v}))$ exists as $n\rightarrow\infty$.
Hence the Monotone Convergence Theorem implies that the function $f(X_{v})$
is integrable. By linearity, the function $g(X_{v})$ is integrable
for each $g\in C_{ub}(S,d)$. In particular, $Eh_{x(\circ),n}(X_{v})\uparrow1.$
Thus $X_{v}$ is a r.v. with values in $(S,d)$. Assertion 4 and the
Corollary are proved.
\end{proof}
\begin{thm}
\textbf{\emph{\label{Thm. Construction of Feller process from semigroup}(Construction
of Feller process from Feller semigroup). }}Let $\mathbf{V}\equiv\{V_{t}:t\in[0,\infty)\}$
be an arbitrary Feller semigroup on the locally compact metric space
$(S,d)$. Then the following holds.

1. Let $x\in(S,d)$ be arbitrary. Then there exists an a.u. càdlàg
process
\[
X\equiv X^{x,\mathbf{V}}:[0,\infty)\times(\Omega,L,E)\rightarrow(S,d),
\]
whose family of marginal distributions $F^{x,\mathbf{V}}$ is generated
by the initial state $x$ and Feller semigroup $\mathbf{V}$ in the
following sense. For arbitrary $m\geq1$, $f\in C_{ub}(S^{m},d^{m})$,
and nondecreasing sequence $r_{1}\leq\cdots\leq r_{m}$ in $[0,\infty)$,
we have
\[
F_{r(1),\cdots,r(m)}^{x,\mathbf{V}}f
\]
\begin{equation}
=\int V_{r(1)}^{x}(dx_{1})\int V_{r(2)-r(1)}^{x(1)}(dx_{2})\cdots\int V_{r(m)-r(m-1)}^{x(m-1)}(dx_{m})f(x_{1},\cdots,x_{m}).\label{eq:temp-384-1-1-2}
\end{equation}

2. For arbitrary $m\geq1$, $f\in C_{ub}(S^{m},d^{m})$, and nondecreasing
sequence $r_{1}\leq\cdots\leq r_{m}$ in $[0,\infty)$, the function
$F_{r(1),\cdots,r(m)}^{*,\mathbf{V}}f$ is continuous on compact subsets
$K\subset(S,d)$ , relative to the metric $d$.

3. Let $\mathcal{L}\equiv\{L^{(t)}:t\in[0,\infty)\}$ be an arbitrary
right continuous filtration. Let $x\in(S,d)$ be arbitrary. Suppose
the a.u. càdlàg Markov process $X^{x,\mathbf{V}}$ is adapted to the
filtration $\mathcal{L}$ . Then the process $X\equiv X^{x,\mathbf{V}}$is
strongly Markov relative to the Feller semigroup $\mathbf{V}$ in
the following sense. 

Let $\tau$ be an arbitrary stopping time with values in $[0,\infty)$.
Let $0\equiv r_{0}\leq r_{1}\leq\cdots\leq r_{m}$ be an arbitrary
nondecreasing sequence in $[0,\infty)$, and let $f\in C(S^{m+1},d^{m+1})$
be arbitrary. Then \emph{(i) } $X_{\tau}$ is a r.v. relative to $L^{(\tau)}$,
and \emph{(ii)}
\[
E(f(X_{\tau+r(0)},X_{\tau+r(1)},\cdots,X_{\tau+r(m)})|L^{(\tau)})
\]
\begin{equation}
=E(f(X_{\tau+r(0)},X_{\tau+r(1)},\cdots,X_{\tau+r(m)})|X_{\tau})=F_{r(0),\cdots,r(m)}^{X(\tau),\mathbf{V}}(f)\label{eq:temp-576-1-1-2-2}
\end{equation}
as r.r.v.'s
\end{thm}
\begin{proof}
Let
\[
(\Omega,L,E)\equiv(\Theta_{0},L_{0},I_{0})\equiv([0,1],L_{0},\int\cdot dx)
\]
denote the Lebesgue integration space based on the unit interval $\Theta_{0}$.

1. Let $\mathbf{T}\equiv\{T_{t}:t\in[0,\infty)\}$ be the compactification
of the Feller semigroup as constructed in Theorem \ref{Thmp. Compactification-of-Feller semigoup into Markoc semgroup}.
Thus $\mathbf{T}$ is a Markov semigroup with the compact metric state
space $(\overline{S},\overline{d})$, where $\overline{d}\leq1$.
Moreover, $\mathbf{T}$ has the modulus of smoothness $\alpha_{\mathbf{T}}\equiv(\alpha_{\mathbf{T,N}})_{N=1,2,\cdots}$
and the modulus of strong continuity $\delta_{\mathbf{T}}$ as constructed
in the proof of Theorem \ref{Thmp. Compactification-of-Feller semigoup into Markoc semgroup}.

2. Let $x\in(S,d)$ be arbitrary. Let 
\begin{equation}
\overline{X}\equiv\overline{X}^{x,\mathbf{T}}:[0,\infty)\times(\Theta_{0},L_{0},I_{0})\rightarrow(\overline{S},\overline{d})\label{eq:temp-585}
\end{equation}
be the a.u. càdlàg\emph{ }Markov\emph{ }process generated by the initial
state $x$ and semigroup $\mathbf{T}$, as constructed in Theorem
\ref{Thm. Markov process from initial distribution and Markov semigroup-1}.

3. Let $t\in[0,\infty)$ be arbitrary. Assertion 4 of Corollary \ref{Cor:(Nonexplosion-of-Feller process at stopping time}
implies that, $(\overline{X}_{t}\in S)$ is a full set. Define the
function $X_{t}:(\Theta_{0},L_{0},I_{0})\rightarrow(S,d)$ by $domain(X_{t})=(\overline{X}_{t}\in S)$
and by $X_{t}=\overline{X}_{t}$ on $domain(X_{t})$. Then $X_{t}$
is a r.v. with values in $(S,\overline{d})$. Since the identity mapping
$\iota:(S\overline{,d})\rightarrow(S,d)$ is continuous, the function
$X_{t}$ is a r.v. with values in $(S,d)$. Thus the function 
\begin{equation}
X\equiv X^{x,\mathbf{V}}:[0,\infty)\times(\Theta_{0},L_{0},I_{0})\rightarrow(S,d)\label{eq:temp-585-2}
\end{equation}
is a well-defined process with state space $(S,d)$. Let $F^{x,\mathbf{V}}$
denote the family of marginal f.j.d.'s of this process $X\equiv X^{x,\mathbf{V}}$.

4. Let the integer $m\geq1$, the function $f\in C(\overline{S}^{m},\overline{d}^{m})$,
and the nondecreasing sequence $r_{1}\leq\cdots\leq r_{m}$ in $[0,\infty)$
be arbitrary. Then Theorem \ref{Thm. Construction of Markov process from  initial distrinbution and semigrp-1-1}
says that the function 
\[
F_{r(1),\cdots,r(m)}^{*,\mathbf{T}}f:(\overline{S},\overline{d})\rightarrow R
\]
is uniformly continuous on $(\overline{S},\overline{d})$. Hence it
is uniformly continuous on compact subsets of $(S,d),$relative to
the metric $d$. At the same time,
\[
F_{r(1),\cdots,r(m)}^{x,\mathbf{V}}\equiv I_{0}f(X_{r(1)},\cdots,X_{r(m)})=I_{0}f(\overline{X}_{r(1)},\cdots,\overline{X}_{r(m)})\equiv F_{r(1),\cdots,r(m)}^{x,\mathbf{T}}f
\]
\[
=\int_{x(1)\in\overline{S}}T_{r(1)}^{x}(dx_{1})\int_{x(2)\in\overline{S}}T_{r(2)-r(1)}^{x(1)}(dx_{2})\cdots\int_{x(m)\in\overline{S}}T_{r(m)-r(m-1)}^{x(m-1)}(dx_{m})f(x_{1},\cdots,x_{m})
\]
\[
=\int_{x(1)\in S}V_{r(1)}^{x}(dx_{1})\int_{x(2)\in S}V_{r(2)-r(1)}^{x(1)}(dx_{2})\cdots\int_{x(m)\in S}V_{r(m)-r(m-1)}^{x(m-1)}(dx_{m})f(x_{1},\cdots,x_{m})
\]
\begin{equation}
=\int V_{r(1)}^{x}(dx_{1})\int V_{r(2)-r(1)}^{x(1)}(dx_{2})\cdots\int V_{r(m)-r(m-1)}^{x(m-1)}(dx_{m})f(x_{1},\cdots,x_{m}).\label{eq:temp-584}
\end{equation}
where the fifth equality is thanks to equality \ref{eq:temp-542}.
We have proved the desired equality \ref{eq:temp-384-1-1-2} in Assertion
1. We have yet to complete the proof of Assertion 1, pending the a.u.
càdlàg property of the process $X$.

5. It was observed in Step 4 that the function $F_{r(1),\cdots,r(m)}^{*,\mathbf{T}}$
is continuous on compact subsets $K$ of $(S,d)$, relative to the
metric $d$. Hence, in view of equality \ref{eq:temp-584}, the function
$F_{r(1),\cdots,r(m)}^{*,\mathbf{V}}f$ is continuous on compact subsets
$K$ of $(S,d)$, relative to the metric $d$. This proves Assertion
2.

5. By Theorem \ref{Thm. Each a.u. cadlag Markov process  is a strong Markov process},
the a.u. càdlàg Markov process $\overline{X}$ is strongly Markov.
To be precise, the following two conditions hold.

(i) Let $\tau$ be an arbitrary stopping time with values in $[0,\infty),$
relative to the filtration $\mathcal{L}$. Then the function $\overline{X}_{\tau}$
is a r.v. relative to $L^{(\tau)}$, and 

(ii). The process $\overline{X}$ is strongly Markov\emph{ }relative
to  filtration $\mathcal{\mathcal{L}}$. Specifically, let $\tau$
be an arbitrary stopping time with values in $[0,\infty),$ relative
to the filtration $\mathcal{L}$. Let $0\equiv r_{0}\leq r_{1}\leq\cdots\leq r_{m}$
be an arbitrary nondecreasing sequence in $[0,\infty),$ and let $f\in C(S^{m+1},d^{m+1})$
be arbitrary. Then
\[
I_{0}(f(\overline{X}_{\tau+r(0)},\overline{X}_{\tau+r(1)},\cdots,\overline{X}_{\tau+r(m)})|L^{(\tau)}
\]
\begin{equation}
=I_{0}(f(\overline{X}_{\tau+r(0)},\overline{X}_{\tau+r(1)},\cdots,\overline{X}_{\tau+r(m)})|\overline{X}_{\tau})=F_{r(0),\cdots,r(m)}^{\overline{X}(\tau),\mathbf{T}}(f)\label{eq:temp-576-1-1-2-3}
\end{equation}
as r.r.v.'s. By Corollary \ref{Cor:(Nonexplosion-of-Feller process at stopping time},
we have $\overline{X}_{\tau+r(1)}=X_{\tau+r(1)}$ on some full set,
for each$i=1,\cdots,m$. Hence equality \ref{eq:temp-576-1-1-2-3}
implies that 
\[
I_{0}(f(X_{\tau+r(0)},X_{\tau+r(1)},\cdots,X_{\tau+r(m)})|L^{(\tau)})
\]
\begin{equation}
=I_{0}(f(X_{\tau+r(0)},X_{\tau+r(1)},\cdots,X_{\tau+r(m)})|X_{\tau})=F_{r(0),\cdots,r(m)}^{\overline{X}(\tau),\mathbf{T}}(f)=F_{r(0),\cdots,r(m)}^{X(\tau),\mathbf{V}}(f)\label{eq:temp-576-1-1-2-3-2}
\end{equation}
Assertion 3 has been proved. This shows that the process $X$ is also
strongly Markov.

6. To complete the proof, we need to show that the process $X$ is
a.u. càdlàg relative to the metric $d$, as in Definition \ref{Def. Metric space of a.u. cadlag process on =00005B0,inf)}.
In other words, we need to prove that it is continuous in probability
on each interval $[0,N]$, and that the shifted process $Y\equiv X^{N}:[0,1]\times(\Omega,L,E)\rightarrow(S,d)$
is a.u. càdlàg, for each $N\geq0$. Definition \ref{Def. Metric space of a.u. cadlag process on =00005B0,inf)}
also defines the shifted processes. We will give the proof only for
the case where $N=0$. The other cases are similar. 

7. To that end, let $Y\equiv X^{0}\equiv X|[0,1]$. Let $Z\equiv Y|Q_{\infty}=X|Q_{\infty}$.
Then, thanks to the strong continuity of the Feller semigroup $\mathbf{V}$
and to the strong Markov property of the process $X$, the process
\textbf{$Z$} is strongly right continuous in probability, in the
sense of Definition \ref{Def. strongly  right continuity in probability on Q_inf}.
We will prove that, in addition, the process $Z$ is also a.u. bounded,
in the sense of Definition \ref{Def. a.u. boundedness on Q_inf} 

To that end, Let $\varepsilon_{0}>0$ be arbitrary. Let $\varepsilon\equiv2^{-1}\varepsilon_{0}$.
Write $n\equiv\kappa_{\mathbf{V,\mathit{1}}}(\varepsilon)$. Take
any $b_{1}\in(n+1,n+2)$. Let 
\[
\beta_{auB}(\varepsilon_{0})\equiv2b_{1}+d(x_{\circ},x).
\]
 Now consider each $k\geq0$ and each $\gamma'>\beta_{auB}(\varepsilon)$.
For abbreviation, write$\gamma\equiv\gamma'-d(x_{\circ},x)$. Then
$\gamma>2b_{1}$.

Let $\eta\equiv\eta_{0,\gamma,Q(k)}$ be the simple first exit time
in the sense of Definition \ref{Def. Simple First Exit time}. Then
\[
P(\bigvee_{r\in Q(h)}d(x,Z_{r})>\gamma;d(x,Z_{1})\leq b_{1})\;\leq P(d(x,Z_{\eta})>\gamma;d(x,Z_{1})\leq b_{1})
\]
\[
\leq P(d(Z_{\eta},Z_{1})>\gamma-b_{1})\;\leq P(d(Z_{\eta},Z_{1})>b_{1})
\]
\[
\leq P(d(X_{\eta},X_{1})>b_{1})\;\leq E(\overline{h}_{X(\eta),n}(X_{1}))
\]
\[
=\int V_{\eta}^{x}(dy)\int V_{1-\eta}^{y}(dz)\overline{h}_{y,n}(z)=\int V_{\eta}^{x}(dy)V_{1-\eta}^{y}d\overline{h}_{y,n}
\]
\[
\leq\int V_{\eta}^{x}(dy)\varepsilon=\varepsilon.
\]
Therefore
\[
P(\bigvee_{r\in Q(h)}d(x,Z_{r})>\gamma)
\]
\[
\leq P(\bigvee_{r\in Q(h)}d(x,Z_{r})>\gamma;d(x,Z_{1})\leq b_{1})+P(d(x,Z_{1})>b_{1})\leq\varepsilon+\varepsilon=2\varepsilon.
\]
Hence 
\[
P(\bigvee_{r\in Q(h)}d(x_{\circ},Z_{r})>\gamma')\equiv P(\bigvee_{r\in Q(h)}d(x_{\circ},Z_{r})>d(x_{\circ},x)+\gamma)
\]
\[
\leq P(\bigvee_{r\in Q(h)}d(x,Z_{r})>\gamma)\leq2\varepsilon.
\]
Summing up, we have proved that the process $Z=X|Q_{\infty}$ is a.u.
bounded, with a modulus of a.u. boundlessness $\beta_{auB}.$ Since
it is also strongly right continuous, Theorem \ref{Thm. a.u.boundeness and Strongly right continuity impliy D-regular}
is applicable to $Z$, and implies that right-limit extension $\widetilde{X}\equiv\Phi_{rLim}(Z):[0,1]\times\Omega\rightarrow S$
relative to the metric $d$, is an a.u. càdlàg process relative to
the metric $d$. At the same time $X$ is a.u. càdlàg relative to
the metric $\overline{d}$. Therefore $X|[0,1]=\Phi_{rLim}(Z)$ relative
to the metric $\overline{d}$. Since sequential convergence relative
to $d$ is equivalent to sequential convergence relative to $\overline{d}$,
we have $X|[0,1]=\widetilde{X}$. Therefore $X|[0,1]$ is a.u. càdlàg
process relative to the metric $d$. 

Similarly we can prove that, for each $N\geq0$, the shifted process
$X^{N}:[0,1]\times(\Omega,L,E)\rightarrow(S,d)$ defined by $X_{t}^{N}\equiv X_{N+t},$
for each $t\in[0,1]$ is a.u. càdlàg. Summing up, the process is a.u.
càdlàg in the sense of Definition \ref{Def. Metric space of a.u. cadlag process on =00005B0,inf)}.

The theorem is proved.
\end{proof}
\begin{cor}
\textbf{\emph{\label{Cor.  a.u.boundedness of Feller process}(Modulus
of a.u. boundlessness on finite intervals).}} The process $X\equiv X^{x,\mathbf{V}}$
constructed in Theorem \ref{Thm. Construction of Feller process from semigroup}
is \emph{a.u. bounded\index{a.u. boundlessness of Feller process}}
in the following sense. Let $\varepsilon_{0}>0$ be arbitrary. Let
$M>N\geq0$ be arbitrary integers. Then there exists $\beta\equiv\widetilde{\beta}_{auB}(\varepsilon_{0},N,M)>0$
so large that, for some measurable set $A$ with $P(A^{c})<\varepsilon_{0}$,
we have 
\[
d(x_{\circ},X(t,\omega))\leq\beta
\]
for each $t\in[N,M]\cap domain(X(\cdot,\omega)$, for each $\omega\in A$.
The operation is then called a modulus of a.u. boundlessness of the
Feller process $X$.
\end{cor}
\begin{proof}
1. We will prove only that the process $X^{0}\equiv X|[0,1]$ is a.u.
bounded, with a modulus $\widetilde{\beta}_{auB}(\cdot,0,1)$. The
prooffor the shifted processes $X^{N}$ with $N\geq1$ is similar
with modulus $\widetilde{\beta}_{auB}(\cdot,N,N+1)$. Then the a.u.
boundlessness of the shifted processes $X^{0},X^{1},\cdots,X^{M-1}$
together imply the a.u. boundlessness of the process $X$ on $[0,M].$
The moduli $\widetilde{\beta}_{auB}(\cdot,N,N+1),\cdots,\widetilde{\beta}_{auB}(\cdot,M-1,M)$
will together yield the modulus $\widetilde{\beta}_{auB}(\cdot,N,M)$

To proceed, let $Y\equiv X|[0,1]$. Let $Z\equiv Y|Q_{\infty}=X|Q_{\infty}.$
In Step 7 of the proof of Theorem \ref{Thm. Construction of Feller process from semigroup},
we saw that the process $Z$ is a.u. bounded in the sense of Definition
\ref{Def. a.u. boundedness on Q_inf}, with some modulus of a.u. boundlessness
$\beta_{auB}$. We also saw that the process $Y\equiv X|[0,1]$ is
a.u. càdlàg, and that $Y=\Phi_{rLim}(Z)$ relative to the metric $d$.

Now let $\varepsilon_{0}>0$ be arbitrary. Write $\varepsilon\equiv2^{-1}\varepsilon_{0}$.
Take any $\gamma>\beta_{auB}(\varepsilon).$ Define $\widetilde{\beta}_{auB,1}(\varepsilon_{0},0,1)\equiv\gamma+2\varepsilon$.

2. By Condition 3 in Definition \ref{Def. a.u. cadlag process} for
a.u. càdlàg processes on $[0,1]$, there exists (i) $\delta_{aucl}(\varepsilon)>0$,
(ii) a $\mathrm{measurable}$ set $A_{1}$ with $P(A_{1}^{c})<\varepsilon$,
(iii) an integer $h\geq1$, and (iv) a sequence of r.r.v.'s 
\begin{equation}
0=\tau_{0}<\tau_{1}<\cdots<\tau_{h-1}<\tau_{h}=1,\label{eq:temp-404-3}
\end{equation}
such that, for each $i=0,\cdots,h-1$, the function $X_{\tau(i)}$
is a r.v., and such that, (v) for each $\omega\in A_{1}$, we have
\begin{equation}
\bigwedge_{i=0}^{h-1}(\tau_{i+1}(\omega)-\tau_{i}(\omega))\geq\delta_{aucl}(\varepsilon),\label{eq:temp-342-3}
\end{equation}
with
\begin{equation}
d(X(\tau_{i}(\omega),\omega),X(\cdot,\omega))\leq\varepsilon,\label{eq:temp-307-4-1-3}
\end{equation}
on the interval $\theta_{i}(\omega)\equiv[\tau_{i}(\omega),\tau_{i+1}(\omega))$
or $\theta_{i}(\omega)\equiv[\tau_{i}(\omega),\tau_{i+1}(\omega)]$
according as $0\leq i\leq h-2$ or $i=h-1$. 

3. Now take $k\geq0$ so large that $2^{-k}<\delta_{aucl}(\varepsilon)$.
Then,because $\gamma>\beta_{auB}(\varepsilon),$ we have 

\begin{equation}
P(A_{0}^{c})<\varepsilon,\label{eq:temp-477-2-1}
\end{equation}
where 
\[
A_{0}\equiv(\bigvee_{r\in Q(k)}d(x_{\circ},Z_{r})\leq\gamma)
\]

4. Consider each $\omega\in A\equiv A_{0}A_{1}$ and each $i=0,\cdots,h-1$.
Then, inequality \ref{eq:temp-342-3} implies that the interval $\theta_{i}(\omega)$
contains some point $r\in Q_{k}\equiv\{0,2^{-k},2\cdot2^{-k},\cdots,1\}$.
Hence inequality \ref{eq:temp-342-3} implies that, for each $t\in\theta_{i}(\omega),$we
have
\[
d(x_{\circ},X(t,\omega))
\]
\[
\leq d(x_{\circ},X(r,\omega))+d(X(r,\omega),X(\tau_{i}(\omega),\omega))+d(X(\tau_{i}(\omega),\omega),X(t,\omega))
\]
\[
\leq d(x_{\circ},Z(r,\omega))+d(X(r,\omega),X(\tau_{i}(\omega),\omega))+d(X(\tau_{i}(\omega),\omega),X(t,\omega))
\]
\begin{equation}
\leq\gamma+\varepsilon+\varepsilon,\label{eq:temp-307-4-1-3-1}
\end{equation}
Combining, we see that $d(x_{\circ},X(t,\omega))\leq\gamma+2\varepsilon$
for each $t\in\theta_{0}(\omega)\cup\cdots\cup\theta_{h-1}(\omega).$
Since this last union of intervals is dense in $[0,1]$, and since
the function $X(\cdot,\omega)$ is right continuous, it follows that
$d(x_{\circ},X(t,\omega))\leq\gamma+2\varepsilon\equiv\widetilde{\beta}_{auB,1}(\varepsilon_{0},0,1)$
for each $t\in[0,1]\cap domain(X(\cdot\omega)),$ for each $\omega\in A.$
Note that $P(A^{c})\equiv P(A_{0}A_{1})^{c}<2\varepsilon\equiv\varepsilon_{0}$.
We have verified that the operation $\widetilde{\beta}_{auB,1}(\cdot,0,1)$
is a modulus of a.u. boundlessness of the process $X|[0,1]$.
\end{proof}
\begin{defn}
\label{Def. Feller process} \textbf{(Feller process). }Let $\mathbf{V}\equiv\{V_{t}:t\in[0,\infty)\}$
be an arbitrary Feller semigroup on the locally compact metric space
$(S,d)$. Let $x\in(S,d)$ be arbitrary. Let $X$ be an arbitrary
a.u. càdlàg process whose marginal distributions are given by the
family $F^{x,\mathbf{V}}$, and which is adapted to some right continuous
filtration $\mathcal{L}$, and which is strongly Markov relative to
the filtration $\mathcal{L}$ and the Feller semigroup $\mathbf{V}$,
in the sense of Conditions 3 in Theorem \ref{Thm. Construction of Feller process from semigroup}.
Then $X$ is called a \emph{Feller process}\index{Feller process}
generated by the initial  state $x$ and the Feller semigroup $\mathbf{V}.$
\end{defn}
\begin{thm}
\textbf{\emph{\label{Thm. Abundance of first exit times for Feller process}
(Abundance of first exit times for Feller process over finite intervals).}}
Let $\mathbf{V}\equiv\{V_{t}:t\in[0,\infty)\}$ be an arbitrary Feller
semi$grou$p on the locally compact metric space $(S,d)$. Let $x\in(S,d)$
be arbitrary. Let 
\[
X:[0,\infty)\times(\Omega,L,E)\rightarrow(S,d),
\]
be the Feller process generated by the initial state $x$ and Feller
semigroup $\mathbf{V}$, in the sense of Theorem \ref{Thm. Construction of Feller process from semigroup}.
Let$\mathcal{L}\equiv\{L^{(t)}:t\in[0,\infty)\}$ be an arbitrary
right continuous filtration relative to which the process $X$ is
adapted.. Let $f:(S,d)\rightarrow R$ be an arbitrary nonnegative
function which is continuous on compact subsets of $(S,d)$, such
that $f(X_{0})\leq a_{0}$ for some $a_{0}\in R$. Let $a\in(a_{0},\infty)$
and be arbitrary. 

Then there exists a countable subset $G$ of $R$ such that, for each
$a\in(a_{0},\infty)G_{c}$, and for each $M\geq1$, the first exit
time $\overline{\tau}_{f,a,M}(X)$ of the open subset $(f<a)$ exists,
in the sense of Definition \textup{\ref{Def. First exit time}.} Here
$G_{c}$ denotes the metric complement of $G$ in $R$. 

Note that, in general, there is no requirement that the process actually
exits $(f<a)$ ever. It is stopped at time $M$ if it does not exit
by then.
\end{thm}
\begin{proof}
Let $x\in(S,d)$ and the continuous function $f$ be as given in the
hypothesis, Thus $f(X_{0})=f(x)\leq a_{0}$ for some $a_{0}\in R$.
Let $b\geq0$ be such that $|f|\leq b$. Let $M\geq1$ be arbitrary.

1. Use the terminology and the notations in the proof of Theorem \ref{Thm. Construction of Feller process from semigroup}.
There is no loss of generality in assuming that $X=X^{x,\mathbf{V}is}$
is the process constructed in the proof of Theorem \ref{Thm. Construction of Feller process from semigroup},
where we saw that, on a full set, we have $X^{x,\mathbf{V}}=\overline{X}^{x,\mathbf{T}},$
where $\mathbf{T}$ is the one-point compactification of the Feller
semigroup $\mathbf{V},$ and where the process $\overline{X}^{x,\mathbf{T}}$
is constructed in Step 2 of said proof. 

2. Let $k\geq0$ be arbitrary. Then Corollary \ref{Cor.  a.u.boundedness of Feller process}
says that there exists $\beta_{k}\equiv\widetilde{\beta}_{auB}(2^{-k},0,M)>0$
so large that, for some measurable set $A_{k}$ with $P(A_{k}^{c})<2^{-k}$,
we have 
\[
d(x_{\circ},X(t,\omega))\leq\beta_{k}
\]
for each $t\in[0,M]\cap domain(X(\cdot,\omega)$, for each $\omega\in A_{k}$.
Define the measurable set
\[
A_{k+}\equiv\bigcap_{i=k}^{\infty}A_{i}.
\]
Then $P(A_{k+}^{c})<2^{-k+1}.$

3. Let $(n_{k})_{k=0,1,\cdots}$ be an increasing sequence such that
$n_{k}>\beta_{k}\equiv\widetilde{\beta}(2^{-k},0,M)$. Let $k\geq0$
be arbitrary. Define $f_{k}\equiv fh_{x(\circ),n(k)}\in C(S,d)$.
Then $|f_{k}|\leq b$, and $f_{k}$ has the same compact support as
$h_{x(\circ),n(k)}$. Hence $f_{k}\in C(S,d)\subset C(\overline{S},\overline{d}),$
with the understanding that $f_{k}(\Delta)\equiv0$. Let $y\in S$
be arbitrary such that $h_{x(\circ),n(k)}(y)=1$.Then
\[
f_{k}(y)=f(y).
\]
Note that. since $f\geq0$, we have $f_{k}\uparrow f$ .

4. First consider the process
\begin{equation}
X=\overline{X}\equiv\overline{X}^{x,\mathbf{T}}:[0,\infty)\times(\Theta_{0},L_{0},I_{0})\rightarrow(\overline{S},\overline{d}).\label{eq:temp-585-1}
\end{equation}
By Theorem \ref{Thm. Abundance of first exit times for Markov process},
there exists a countable subset $G$ of $R$ such that, for each $a\in(a_{0},\infty)G_{c}$,
the first exit time $\tau_{k}\equiv\overline{\tau}_{f(k),a,M}(X)$
exists. By Definition \ref{Def. First exit time}, $\tau_{k}$ is
a stopping time relative to $\mathcal{L}$ , with values in $[0,M]$,
such that the function $X_{\tau(k)}$ is a well-defined r.v relative
to $L^{(\tau(k))}$, with values in $\overline{S}$, and such that,
for each $\omega\in domain(X_{\tau(k)})$, we have 

(i). $f_{k}(X(\cdot,\omega))<a$ on the interval $[0,\tau_{k}(\omega))$,
and 

(ii).\emph{ $f_{k}(X_{\tau(k)}(\omega))\geq a$ }if\emph{ $\tau_{k}(\omega)<M$.}

5. We will verify that the sequence $(\tau_{k})_{k=0,1,\cdots}$ is
nonincreasing. Suppose $\tau_{k}(\omega)<\tau_{k+1}(\omega)$ for
some $\omega\in\Omega$. Then, according to Condition (i), applied
to $f_{k+1}$, we have 
\[
f_{k+1}(X(\tau_{k}(\omega),\omega))<a.
\]
Moreover,\emph{ $\tau_{k}(\omega)<M$ because $\tau_{k+1}(\omega)\leq M$.
}Hence\emph{ $f_{k}(X_{\tau(k)}(\omega))\geq a$ }according to Condition
(ii). This is a contradiction because $f_{k}\leq f_{k+1}$. We conclude
that $\tau_{k}\geq\tau_{k+1}.$

6. Next, let $i>k\geq0$ be arbitrary. Consider each $\omega\in\Omega$
such that $\tau_{i}(\omega)>\tau_{k}(\omega).$ Then, according to
Condition (i), we have $f_{k}(X(\tau_{i}(\omega),\omega))<a$. Moreover,
$\tau_{i}(\omega)<M$, whence,\emph{ $f_{i}(X_{\tau(i)}(\omega))\geq a$,
}according to Conditon (ii) applied to $f_{i}$. Combining,\emph{
\[
f_{k}(X_{\tau(i)}(\omega))<f_{i}(X_{\tau(i)}(\omega)).
\]
}In other words,\emph{
\[
f(X_{\tau(i)}(\omega))h_{x(\circ),n(k)}(X_{\tau(i)}(\omega))<f(X_{\tau(i)}(\omega))h_{x(\circ),n(i)}(X_{\tau(i)}(\omega).
\]
}Cancelling the common factor, we obtain\emph{
\[
h_{x(\circ),n(k)}(X_{\tau(i)}(\omega))<h_{x(\circ),n(i)}(X_{\tau(i)}(\omega))\leq1
\]
}

The strict inequality implies that $d(x_{\circ},X_{\tau(i)}(\omega))\geq n_{k}>\beta_{k}$.
Therefore, according to Step 2, we have $\omega\in A_{k}^{c}\subset A_{k+}^{c}.$
Summing up, $(\tau_{i}>\tau_{k})\subset A_{k+}^{c}$. Hence $\tau_{i}=\tau_{k}$
on $A_{k+}$. Since $P(A_{k+}^{c})\leq2^{-k+1}$ is arbitrarily small
for sufficiently large $k\geq0$, it follows that $\tau_{k}\downarrow\tau$
a.u. for some r.r.v. $\tau$ with values in $[0,M]$. Moreover, 
\begin{equation}
\tau=\tau_{k}\label{eq:temp-578}
\end{equation}
on $A_{k+}$, for each $k\geq0$. 

7. Next we will verify that $\tau$ is a stopping time relative to
the right continuous filtration $\mathcal{L}\equiv\{L^{(t)}:t\in[0,\infty)\}$.
To that end, consider each regular point $t\in[0,M]$ of $\tau$.
Let $s\in(t,\infty)$ be arbitrary. Let $(s_{k})_{k=0,1,\cdots}$
be a sequence such that $s_{k}\in(t,s\wedge(t+2^{-k}))$ is a regular
point of the r.r.v.'s $\tau$ and $\tau_{k}$, for each $k\geq0$.
Then 
\[
1_{(\tau(k)\leq s(k))}\geq1_{(\tau\leq s(k))}\geq1_{(\tau\leq t)}
\]
for each $k\geq0$. Moreover, for each $k\geq0$, we have $\tau=\tau_{k}$
on the set $A_{k+}$. Hence
\[
E1_{(\tau(k)\leq s(k))}=E1_{(\tau(k)\leq s(k))}1_{A(k+)}+E1_{(\tau(k)\leq s(k))}1_{A(k+)^{c}}
\]
\[
\leq E1_{(\tau\leq s(k))}1_{A(k+)}+2^{-k+1}
\]
\[
\leq E1_{(\tau\leq s(k))}+2^{-k+1}\downarrow E1_{(\tau\leq t)}.
\]
At the same time, $1_{(\tau(k)\leq s(k))}\in L^{(s(k))}\subset L^{(s)}$
because $\tau_{k}$ is a stopping time relative to the filtration
$\mathcal{L}$, for each $k\geq0$, according to Step 4. Therefore
the Monotone Convergence Theorem implies that $1_{(\tau\leq t)}\in L^{(s)},$
where $s\in(t,\infty)$ is arbitrary. We conclude that
\[
1_{(\tau\leq t)}\in L^{(t+)}=L^{(t)},
\]
 where the equality is thanks to the right continuity of the filtration
$\mathcal{L}$. Thus we have verified that $\tau$ is a stopping time
relative to the filtration $\mathcal{L}$.

8. Observe that, for each $i\geq k\geq0$, we have $(X_{\tau(k)}=X_{\tau})$
on $A_{k+}$, thanks to equality \ref{eq:temp-578}. Hence,, trivially,
$X_{\tau(k)}\rightarrow X_{\tau}$ a.u. Since $X_{\tau(k)}$ is a
r.v. with values in $(\overline{S},\overline{d})$, so is $X_{\tau}$.
Because there exists a full set on which $\overline{X}_{t}=X_{t}$
for each $t\in[0,1]$, the r.v. $X_{\tau}$ has values in $(S,\overline{d})$
. Since the identity mapping $\iota:(S\overline{,d})\rightarrow(S,d)$
is continuous according to Condition 3 of Definition \ref{Specification of locally compact metric space etc)},
the function $X_{t}$ is actually a r.v. with values in $(S,d)$.
By restricting the r.v. $X_{\tau}$ to the full set $D\equiv\bigcap_{\kappa=0}^{\infty}\bigcup_{k=\kappa}^{\infty}A_{k+}$,
we may assume that $domain(X_{\tau})\subset D$. 

9. Now consider each $\omega\in domain(X_{\tau})$. Let $t\in[0,\tau(\omega))$
be arbitrary. Take $\kappa\geq0$ be so large that
\[
h_{x(\circ),n(\kappa)}(X_{\tau}(\omega))=1=h_{x(\circ),n(\kappa)}(X_{t}(\omega)),
\]
which leads to $f_{k}(X_{\tau}(\omega))=f(X_{\tau}(\omega))$ and
$f_{k}(X_{t}(\omega))=f(X_{t}(\omega))$, for each $k\geq\kappa.$
Separately, because $\omega\in D$, there exists $k\geq\kappa$ such
that $\omega\in A_{k+}.$ Hence $\tau(\omega)=\tau_{k}(\omega)$ and
$X_{\tau(k)}(\omega)=X_{\tau}(\omega)$. Consequently, $t\in[0,\tau_{k}(\omega))$.
Conditions (i) and (ii) in Step 4 therefore implies

(i'). $f(X(t,\omega))<a$ , and 

(ii').\emph{ $f(X_{\tau}(\omega))\geq a$ }if\emph{ $\tau(\omega)<M$.}

We have thus verified all the conditions in Definition \ref{Def. First exit time}
for $\tau$ to be the \emph{first exit time }in\emph{ }$[0,M]$ of
the open subset $(f<a)$ by the process $X$. In symbols, $\overline{\tau}_{f,a,M}\equiv\overline{\tau}_{f,a,M}(X)\equiv\tau$.
The theorem is proved.
\end{proof}

\section{The Brownian motion in $R^{m}$}

In the remainder of this chapter, let $m\geq1$ be arbitrary, but
fixed, and let $(R^{m},d^{m})$ denote the $m$-dimensional Euclidean
metric space. Thus $d^{m}(x,y)=|x-y|$ for each $x,y\in R^{m}$. 
\begin{defn}
\textbf{\label{Def. Notations for matrices}(Notations for matrices
and vectors).} We will use the notations of matrices and vectors in
Section \ref{sec:Normal-Distribution}, with one exception: for the
transpose of an arbitrary matrix $A$, we will write $A'$, rather
than the symbol $A^{T}$ used in Section \ref{sec:Normal-Distribution},
so that we can reserve the symbol $T$ for transition distributions.
We will regard each point $x\equiv(x_{1},\cdots,x_{m})\in R^{m}$
as a column vector, i.e. a $m\times1$ matrix $\left[\begin{array}{c}
x_{1}\\
\vdots\\
x_{m}
\end{array}\right]$. Then $x'x$ is a $1\times1$ matrix which we identify with its sole
entry $\sum_{i=1}^{m}x_{i}^{2}$, the inner product of the vector
$x$. At the risk of some confusion, we will write $0\equiv(0,\cdots,0)\in R^{m},$
and write $I$ for the the $m\times m$ identity matrix.
\end{defn}
$\square$
\begin{defn}
\textbf{\label{Def. Notations-1}(Notations for normal distributions
and normal p.d.f.'s on $R^{m}$).} We will use the notations and results
for normal distributions in Section \ref{sec:Normal-Distribution}.
In particular, for each $t\geq0$, recall, from Definition \ref{Def. n-dim Normal Distribution, nonegative definite},
the normal distribution $\Phi_{x,t}\equiv\Phi_{x,tI}$ on $R^{m}$
with mean vector $x\equiv(x_{1},\cdots,x_{m})\in R^{m}$ and covariance
matrix $tI$. If \textbf{$t>0$}, then $\Phi_{x,t}$ has the p.d.f.
$\varphi_{x,t}\equiv\varphi_{x,tI}$ defined by

\[
\varphi_{x,t}(y)\equiv(2\pi)^{-\frac{m}{2}}t^{-\frac{m}{2}}\exp(-\frac{1}{2}t^{-1}(y-x)'(y-x))
\]
\begin{equation}
=(2\pi)^{-\frac{m}{2}}t^{-\frac{m}{2}}\exp(-\frac{1}{2}t^{-1}\sum_{i=1}^{m}(y_{i}-x_{i})^{2})\label{eq:temp-36-1}
\end{equation}
for each $y\equiv(y_{1},\cdots,y_{m})\in R^{m}$. A r.v. $U$ with
values in $R^{m}$ is said to be \index{standard normal r.v.}\emph{standard
normal} if it has the normal distribution $\Phi_{0,1}$.

In the remainder of this chapter, let $U,W,V,U_{1},U_{2},\cdots$
be an arbitrary independent sequence of standard normal r.v.'s with
values in $R^{m}$, on some probability space $(\widetilde{\Omega},\widetilde{L},\widetilde{E})$.
These r.v.'s are used only for the compact notations for their joint
distributions; the sample space $(\widetilde{\Omega},\widetilde{L},\widetilde{E})$
is immaterial. For example, 
\[
\widetilde{E}f(x+\sqrt{t}U)=\intop_{R^{m}}\Phi_{x,t}(du)f(u)=\intop_{R^{m}}\Phi_{0,t}(du)f(x+u),
\]
where the expectation on the left-hand side gives a compact notation
of two integrals on the right.
\end{defn}
$\square$

We next generalize the definition of Brownian motion introduced in
Definition \ref{Def. Brownian Motion in R}. 
\begin{defn}
\textbf{\label{Def.  Brownian Motion inR^m}(Brownian motion in $R^{m}$).}
An a.u. continuous process $B:[0,\infty)\times(\varOmega,L,E)\rightarrow R^{m}$
is called a \emph{Brownian Motion} \emph{in} $R^{m}$ \index{Brownian Motion in $R^{m}$}
if (i) $B_{0}=0$, (ii) for each sequence $0\equiv t_{0}\leq t_{1}\leq\cdots\leq t_{n-1}\leq t_{n}$
in $[0,\infty)$, the r.v.'s $B_{t(1)}-B_{t(0)},\cdots,B_{t(n)}-B_{t(n-1)}$
are independent, and (iii) for each $s,t\in[0,\infty)$, the r.v.
$B_{t}-B_{s}$ is normal with mean $0$ and covariance matrix $|t-s|I$. 

Let $k=1,\cdots,m$ be arbitrary. We will write $B_{k;t}$ for the
$k$-th coordinate of the r.v. $B_{t}$, for each $t\in[0,\infty)$.
Thus $B_{t}\equiv(B_{1;t},\cdots,B_{m;t})$, for each $t\in[0,\infty)$.
The process $B_{k;}:[0,\infty)\times(\varOmega,L,E)\rightarrow R$
defined by
\[
(B_{k;})_{t}\equiv B_{k;t}
\]
for each $t\in[0,\infty)$, is called the $k$-th coordinate process
of $B$. It is easily verified that if $B$ is a Brownian Motion in
$R^{m},$ then $B_{k;}$ is a Brownian motion in $R.$

Suppose $B:[0,\infty)\times(\varOmega,L,E)\rightarrow R^{m}$ is a
Brownian Motion \emph{in} $R^{m}$. Let 
\[
x\equiv(x_{1},\cdots,x_{m})\in R^{m}
\]
be arbitrary. Then we will write $B^{x}$ for the process $x+B:[0,\infty)\times(\varOmega,L,E)\rightarrow R^{m}$,
and call $B^{x}$ a Brownian motion in $R^{m}$ starting at the point
$x$. Thus, for each $k=1,\cdots,m$, the process $B_{k}^{x(k)}$
is a Brownian motion in $R$ starting at the point $x_{k}$.
\end{defn}
$\square$ 

To construct a Brownian motion, we will define a certain Feller semigroup
$\mathbf{T}$ on $[0,\infty)$ with state space $(R^{m},d^{m})$,
and then prove that any a.u. càdlàg Feller process with initial state
$0$ and with this Feller semigroup $\mathbf{T}$ is a Brownian motion.
This way, we can use the theorems in the preceding sections for a.u.
càdlàg Feller processes, to infer strong Markov property, and abundance
of first exit times.
\begin{lem}
\label{Lem. Brownian semigroup} \textbf{\emph{(Semigroup for the
Brownian motion).}} Let $t\geq0$ be arbitrary. Define the function
$T_{t}:C_{ub}(R^{m},d^{m})\rightarrow C_{ub}(R^{m},d^{m})$ by 

\begin{equation}
T_{t}(f)(x)\equiv T_{t}^{x}(f)\equiv\widetilde{E}f(x+\sqrt{t}U)=\intop_{R^{m}}\Phi_{0,t}(du)f(x+u)\label{eq:temp-517}
\end{equation}
for each $x\in R^{m}$, for each $f\in C_{ub}(R^{m},d^{m})$. Then
the family $\mathbf{T}\equiv\{T_{t}:t\in[0,\infty)\}$ is a Feller
semigroup in the sense of Definition \ref{Def. Feller-Semigroup}.
We will call $\mathbf{T}$ the \emph{$m$-dimensional Brownian semigroup.}
\end{lem}
\begin{proof}
We need to verify the conditions in Definition \ref{Def. Feller-Semigroup}.
for the family $\mathbf{T}$ to be a Feller semigroup. It is trivial
to verify that, for each $t\in[0,\infty)$ and for each $y\in S$,
the function $T_{t}^{x}$ is is a distribution on $(R^{m},d^{m})$.

1. Let $N\geq1$, $t\in[0,N]$, and $f\in C_{ub}(R^{m},d^{m})$ be
arbitrary with a modulus of continuity $\delta_{f}$ and with $|f|\leq1$.
Consider the function $T_{t}f$ defined in equality \ref{eq:temp-517}.
Let $\varepsilon>0$ be arbitrary. Let $x,y\in R^{m}$ be arbitrary
with $|x-y|<\delta_{f}(\varepsilon)$. Then $|(x+u)-(y+u)|=|x-y|<\delta_{f}(\varepsilon)$.
Hence

\[
|T_{t}(f)(x)-T_{t}(f)(y)|\equiv|\intop_{R^{m}}\Phi_{0,t}(du)f(x+u)-\intop_{R^{m}}\Phi_{0,t}(du)f(y+u)|
\]
\[
\leq\intop_{R^{m}}\Phi_{0,t}(du)|f(x+u)-f(y+u)|<\varepsilon.
\]
Hence $T_{t}(f)$ has the same modulus of continuity $\alpha_{\mathbf{V,\mathit{N}}}(\delta_{f})\equiv\delta_{f}$
as the function $f$. Thus the family $\mathbf{T}$ has a modulus
of smoothness $\alpha_{\mathbf{T}}\equiv(\iota,\iota,\cdots)$ where
$\iota$ is the identity operation $\iota$. The smoothness condition,
Condition 1 of Definition \ref{Def. Feller-Semigroup} has been verified
for the family $\mathbf{T}$. In particular, we have a function $T_{t}:C_{ub}(R^{m},d^{m})\rightarrow C_{ub}(R^{m},d^{m})$. 

2. We will next prove that $T_{t}^{x}(f)$ is continuous in $t$.
To that end, let $\varepsilon>0$, $x\in R^{m}$ be arbitrary. Note
that the standard normal r.v. $U$ is a r.v. Hence there exists $a_{\varepsilon}>0$
so large that $\widetilde{E}1_{(|U|>a(\varepsilon))}<2^{-1}\varepsilon$.
Note also that $\sqrt{t}$ is a uniformly continuous function of $t\in[0,\infty)$,
with some modulus of continuity $\delta_{sqrt}$. Let $t,s\geq0$
be such that 
\begin{equation}
|t-s|<\delta_{\mathbf{T}}(\varepsilon,\delta_{f})\equiv\delta_{sqrt}(a_{\varepsilon}^{-1}\delta_{f}(2^{-1}\varepsilon)).\label{eq:temp-547}
\end{equation}
Then $|\sqrt{t}-\sqrt{s}|<a_{\varepsilon}^{-1}\delta_{f}(2^{-1}\varepsilon)$.
Hence 
\[
|(x+\sqrt{t}u)-(x+\sqrt{s}u)|=|\sqrt{t}u-\sqrt{s}u|<\delta_{f}(2^{-1}\varepsilon)
\]
 for each $u\in R^{m}$ with $|u|\leq a_{\varepsilon}$. Therefore
the defining equality \ref{eq:temp-517} leads to 
\[
|T_{t}^{x}(f)-T_{s}^{x}(f)|\leq\widetilde{E}|f(x+\sqrt{t}U)-f(x+\sqrt{s}U)|
\]
\[
\leq\widetilde{E}|f(x+\sqrt{t}U)-f(x+\sqrt{s}U)|1_{(|U|\leq a(\varepsilon))}+\widetilde{E}1_{(|U|>a(\varepsilon))}
\]
\[
\leq\widetilde{E}2^{-1}\varepsilon1_{(|U|\leq a(\varepsilon))}+\widetilde{E}1_{(|U|>a(\varepsilon))}
\]
\begin{equation}
\leq2^{-1}\varepsilon+2^{-1}\varepsilon=\varepsilon.\label{eq:temp-496}
\end{equation}
Thus the function $T_{t\cdot}^{x}f$ is uniformly continuous in $t\in[0,\infty)$,
with a modulus of continuity independent of $x$. In the special case
where $s=0$, we have equality $T_{s}^{x}(f)=\widetilde{E}f(x+\sqrt{0}U)=f(x)$,
whence inequality \ref{eq:temp-496} yields
\[
|T_{t}(f)-f|\leq\varepsilon
\]
for each $t\in[0,\delta_{\mathbf{T}}(\varepsilon,\delta_{f}))$. The
strong-continuity condition, Condition 3 in Definition \ref{Def. Feller-Semigroup},
is also verified for the family $\mathbf{T}$, with $\delta_{\mathbf{T}}$
being the modulus of strong continuity. 

3. Next let $t,s\geq0$ be arbitrary. First assume that $s>0$. Then
$(t+s)^{-\frac{1}{2}}(\sqrt{t}W+\sqrt{s}V)$ is a standard normal
r.v. with values in $R^{m}$. Hence, for each $x\in R^{m}$, we have
\[
T_{t+s}^{x}(f)=\widetilde{E}f(x+\sqrt{t+s}U)=\widetilde{E}f(x+\sqrt{t}W+\sqrt{s}V)
\]
\[
=\int_{w\in R^{m}}\int_{u\in R^{m}}\varphi_{0,t}(w)\varphi_{0,s}(v)f(x+w+v)dwdv
\]
\[
=\int_{v\in R^{m}}(\int_{w\in R^{m}}\varphi_{0,t}(w)f(x+w+v)dw)\varphi_{0,s}(v)dv
\]
\[
=\int_{v\in R^{m}}((T_{t}f)(x+v))\varphi_{0,s}(v)dv
\]
\begin{equation}
=T_{s}^{x}(T_{t}f),\label{eq:temp-420}
\end{equation}
provided that $s\in(0,\infty)$. At the same time, inequality \ref{eq:temp-496}
shows that both ends of equality \ref{eq:temp-420} are continuous
functions of $s\in[0,\infty)$. Hence equality \ref{eq:temp-420}
can be extended, by continuity, to
\[
T_{t+s}^{x}(f)=T_{s}^{x}(T_{t}f)
\]
for each $s\in[0,\infty)$. It follows that $T_{s+t}=T_{s}T_{t}$
for each $t,s\in[0,\infty)$. Thus the semigroup property, Condition
2 in Definition \ref{Def. Feller-Semigroup} is proved for the family
$\mathbf{T}$.

4. The non-explosion condition, Condition 4 in Definition \ref{Def. Feller-Semigroup},
is also straightforward, and is left as an exercise.

5. Summing up, all the conditions in Definition \ref{Def. Feller-Semigroup}
hold for the family $\mathbf{T}$ to be a Feller semigroup with state
space $(R^{m},d^{m})$.
\end{proof}
\begin{thm}
\textbf{\emph{\label{Thm. Construction of -Brownian motion in R^m}(Construction
of Brownian motion in}}\textbf{ $R^{m}$}\textbf{\emph{)}}. A Brownian
motion in $R^{m}$ exists which is an a.u. continuous Feller process
generated by the initial state $x=0$ and the Brownian semigroup $\mathbf{T\equiv}\{T_{t}:t\geq0\}$
constructed in Lemma \ref{Lem. Brownian semigroup}. 
\end{thm}
\begin{proof}
1. By Theorem \ref{Thm. Construction of Feller process from semigroup},
there exists a a time-uniformly a.u. càdlàg, strongly Markov, Feller
process 
\[
B:[0,\infty)\times(\Omega,L,E)\rightarrow(R^{m},d^{m})
\]
with initial state \emph{$0\in R^{m}$} and the given Feller semigroup
$\mathbf{T}$. In particular, the family $F^{0,\mathbf{T}}$ of marginal
distributions of the process $B$ is generated by the initial state
$0\in R^{m}$ and the Feller semigroup \emph{$\mathbf{T}$, }in the
sense of equality \ref{eq:temp-384-1-1-2} in Theorem \ref{Thm. Construction of Feller process from semigroup}.
Furthermore, the process $B$ has a modulus of continuity in probability
$\delta_{Cp,\delta(\mathbf{\mathbf{T}})}$, where $\delta_{\mathbf{T}}$
is the modulus of strong continuity of $\mathbf{T}$ defined in equality
\ref{eq:temp-547} in Lemma \ref{Lem. Brownian semigroup}. We will
prove that the process $B$ is a Brownian motion. 

2. Since the process $B$ has initial state \emph{$0\in R^{m}$, }we
have trivially $B_{0}=0$.

3. Let $\mathcal{L}\equiv\{L^{(t)}:t\in[0,\infty)\}$ be a right continuous
filtration to which the process $B$ is adapted. Let $0\equiv t_{0}\leq t_{1}\leq\cdots\leq t_{n-1}\leq t_{n}$
in $[0,\infty)$. Let $f_{1},\cdots,f_{n}\in C_{ub}(R^{m},d^{m})$
be arbitrary. We will prove, inductively on $k=0,1,\cdots,n$, that
\begin{equation}
E\prod_{i=1}^{k}f_{i}(B_{t(i)}-B_{t(i-1)})=\prod_{i=1}^{k}\intop_{R^{m}}\Phi_{0,t(i)-t(i-1)}(du)f_{i}(u).\label{eq:temp-501}
\end{equation}
To start, with $k=0,$ equality \ref{eq:temp-501} is trivially valid
because, by convention, the empty product on each side of the equality
\ref{eq:temp-501} is equal to $1$. Inductively, suppose equality
\ref{eq:temp-501} has been proved for $k-1$, where $1\leq k\leq n$.
For abbreviation, define the function $g_{k}\in C_{ub}(R^{m}\times R^{m},d^{m}\otimes d^{m})$
by 
\[
g_{k}(y_{1},y_{2})\equiv f_{k}(y_{2}-y_{1})
\]
for each $(y_{1},y_{2})\in R^{m}\times R^{m}$. Then, by the Markov
property of the a.u. càdlàg Feller process $B$, proved in Theorem
\ref{Thm. Construction of Feller process from semigroup}, we have
\begin{equation}
E(g_{k}(B_{t(k-1)},B_{t(k)})|L^{(t(k-1))})=F_{0,t(k)-t(k-1)}^{B(t(k-1)),\mathbf{T}}(g_{k})\label{eq:temp-576-2-1}
\end{equation}
as r.r.v.'s, where $F_{0,t(k)-t(k-1)}^{*,\mathbf{T}}$ is the function
defined in Assertion 2 of Theorem \ref{Thm. Construction of Feller process from semigroup}.
Note that
\[
F_{0,t(k)-t(k-1)}^{x,\mathbf{T}}g_{k}\equiv\int T_{t(k)-t(k-1)}^{x}(dx_{1})g_{k}(x,x_{1})
\]
\begin{equation}
=\intop_{R^{m}}\Phi_{0,t(k)-t(k-1)}(du)g_{k}(x,x+u)\equiv\intop_{R^{m}}\Phi_{0,t(k)-t(k-1)}(du)f_{k}(u),\label{eq:temp-546}
\end{equation}
where the second equality is from the defining equality \ref{eq:temp-517},
and where the last equality is from the definition of the function
$g_{k}$.

At the same time, Since the function $\prod_{i=1}^{k-1}f_{i}(B_{t(i)}-B_{t(i-1)})$
is a member of $L^{(t(k-1))}$ and is bounded, equality \ref{eq:temp-576-2-1}
implies, according to basic properties of conditional expectations
in Assertion 4 of Proposition \ref{Prop. Basics of Conditional expectations},
that

\[
E(\prod_{i=1}^{k-1}f_{i}(B_{t(i)}-B_{t(i-1)}))g_{k}(B_{t(k-1)},B_{t(k)})=E(\prod_{i=1}^{k-1}f_{i}(B_{t(i)}-B_{t(i-1)}))F_{0,t(k)-t(k-1)}^{B(t(k-1)),\mathbf{T}}(g_{k})
\]
\[
=E(\prod_{i=1}^{k-1}f_{i}(B_{t(i)}-B_{t(i-1)}))\intop_{R^{m}}\Phi_{0,t(k)-t(k-1)}(du)f_{k}(u))
\]
\[
=(\prod_{i=1}^{k-1}\intop_{R^{m}}\Phi_{0,t(i)-t(i-1)}(du)f_{i}(u))\intop_{R^{m}}\Phi_{0,t(k)-t(k-1)}(du)f_{k}(u))
\]
\[
=\prod_{i=1}^{k}\intop_{R^{m}}\Phi_{0,t(i)-t(i-1)}(du)f_{i}(u),
\]
where the second equality is from equality \ref{eq:temp-546}, and
where the third equality is from the induction hypothesis that equality
\ref{eq:temp-501} has been proved for $k-1$. Thus equality \ref{eq:temp-501}
is proved for $k$. The induction is completed, resulting in 
\begin{equation}
E\prod_{i=1}^{n}f_{i}(B_{t(i)}-B_{t(i-1)})=\prod_{i=1}^{n}\intop_{R^{m}}\Phi_{0,t(i)-t(i-1)}(du)f_{i}(u).\label{eq:temp-501-1}
\end{equation}

4. Let $j=1,\cdots,n$ be arbitrary. In the special case where $f_{i}\equiv1$
for each $i=1,\cdots,n$ with $i\neq j$, equality \ref{eq:temp-501-1}
reduces to
\begin{equation}
Ef_{j}(B_{t(j)}-B_{t(j-1)})=\intop_{R^{m}}\Phi_{0,t(j)-t(j-1)}(du)f_{j}(u)\label{eq:temp-541-3}
\end{equation}
Hence the r.v. $B_{t(j)}-B_{t(j-1)}$ has the normal distribution
$\Phi_{0,t(j)-t(j-1)}$ where $j\geq1$, and $t_{j-1},t_{j}\in[0,\infty)$
are arbitrary with $t_{j-1}\leq t_{j}$. 

5. In the general case, substituting equality \ref{eq:temp-541-3}
back into equality \ref{eq:temp-501-1} yields
\[
E\prod_{i=1}^{n}f_{i}(B_{t(i)}-B_{t(i-1)})=\prod_{i=1}^{n}Ef_{j}(B_{t(j)}-B_{t(j-1)}),
\]
where the functions $f_{1},\cdots,f_{n}\in C_{ub}(R^{m},d^{m})$ are
arbitrary. It follows that the r.v.'s $B_{t(1)}-B_{t(0)},\cdots,B_{t(n)}-B_{t(n-1)}$
are independent.

6. Let $f\in C_{ub}(R^{m},d^{m})$ be arbitrary. Equality \ref{eq:temp-541-3}
implies that
\begin{equation}
Ef(B_{t}-B_{s})=\intop_{R^{m}}\Phi_{0,|t-s|}(du)f(u)\label{eq:temp-541-3-1}
\end{equation}
for each
\[
(t,s)\in C\equiv\{(r,v)\in[0,\infty)^{2}:r\leq v\}
\]
Similarly, 
\begin{equation}
Ef(B_{t}-B_{s})=Ef(-(B_{s}-B_{t}))=\intop_{R^{m}}\Phi_{0,|t-s|}(du)f(-u)=\intop_{R^{m}}\Phi_{0,|t-s|}(du)f(u).\label{eq:temp-541-3-1-1}
\end{equation}
for each
\[
(t,s)\in C'\equiv\{(r,v)\in[0,\infty)^{2}:r\geq v\}
\]
Since both ends of equality \ref{eq:temp-541-3-1} are continuous
in $(t,s)$, and since $C\cup C'$ is dense in $[0,\infty)^{2}$,
we obtain
\begin{equation}
Ef(B_{t}-B_{s})=\intop_{R^{m}}\Phi_{0,|t-s|}(du)f(u)\label{eq:temp-541-3-1-2}
\end{equation}
for each $(t,s)\in[0,\infty)^{2}$. Summing up, Conditions (i-iii)
in Definition \ref{Def.  Brownian Motion inR^m} have been proved
for the process $B$. It remains to prove a.u. continuity in Definition
for $B$ to be a Brownian motion in $R^{m}$.

7. To that end, let $k=1,\cdots,m$ be arbitrary. As in Definition
\ref{Def.  Brownian Motion inR^m}, write $B_{k;t}$ for the $k$-th
coordinate of the r.v. $B_{t}$, for each $t\in[0,\infty)$. Thus
$B_{t}\equiv(B_{1;t},\cdots,B_{m;t})$, for each $t\in[0,\infty)$.
The process $B_{k;}:[0,\infty)\times(\varOmega,L,E)\rightarrow R$
defined by
\[
(B_{k;})_{t}\equiv B_{k;t}
\]
for each $t\in[0,\infty)$, is called the $k$-th coordinate process
of $B$. 

8. Consider the first-coordinate process $X\equiv B_{1;}$. Note that
the process $X$ inherits continuity in probability from the from
the process $B.$ Since $(B_{1;0},\cdots,B_{m;0})\equiv B_{0}=(0,\cdots,0)\in R^{m}$,
we have $B_{1;0}=0\in R$. Let $t,s\geq0$ and $g\in C_{ub}(R^{1},d^{1})$
be arbitrary. Define $f\in C_{ub}(R^{m},d^{m})$ by $f(x_{1},\cdots,x_{m})\equiv g(x_{1})$
for each $(x_{1},\cdots,x_{m})\in R^{m}$. Then equality \ref{eq:temp-541-3-1-2}
implies that
\[
Eg(X_{t}-X_{s})=Eg(B_{1;t}-B_{1;s})=Ef(B_{t}-B_{s})=\intop_{R^{m}}\Phi_{0,|t-s|}(du)f(u)
\]
\[
=\intop\cdots\int\Phi_{0,|t-s|}^{1}(du_{1})\cdots\Phi_{0,|t-s|}^{1}(du_{n})f(u_{1},\cdots,u_{m})
\]
\[
=\intop\cdots\int\Phi_{0,|t-s|}^{1}(du_{1})\cdots\Phi_{0,|t-s|}^{1}(du_{n})g(u_{1})
\]
\begin{equation}
=\int\Phi_{0,|t-s|}^{1}(du_{1})g(u_{1}),\label{eq:temp-541-3-1-3}
\end{equation}
where $\Phi_{0,|t-s|}^{1}$ stands for the normal distribution on
$(R^{1},d^{1})$ with mean $0$ and variance $|t-s|$, and where $g\in C_{ub}(R^{1},d^{1})$
is arbitrary. Thus we see that $X_{t}-X_{s}$ is normally distributed
r.r.v. with mean $0$ and variance $|t-s|$.

9. Next, let the sequence $0\equiv t_{0}\leq t_{1}\leq\cdots\leq t_{n-1}\leq t_{n}$
in $[0,\infty)$ be arbitrary. For each $i=1,\cdots,m$, let $g_{i}\in C_{ub}(R^{1},d^{1})$
be arbitrary, and define $f_{i}\in C_{ub}(R^{m},d^{m})$ by $f_{i}(x_{1},\cdots,x_{m})\equiv g_{i}(x_{1})$
for each $(x_{1},\cdots,x_{m})\in R^{m}$. Then 
\[
E\prod_{i=1}^{n}g_{i}(X_{t(i)}-X_{t(i-1)})=E\prod_{i=1}^{n}f_{i}(B_{t(i)}-B_{t(i-1)})
\]
\[
=\prod_{i=1}^{n}Ef_{i}(B_{t(i)}-B_{t(i-1)})=\prod_{i=1}^{n}Eg_{i}(B_{1;t(i)}-B_{1;t(i-1)})=\prod_{i=1}^{n}Eg_{i}(X_{t(i)}-X_{t(i-1)})
\]
Consequently, the r.r.v.'s $X_{t(1)}-X_{t(0)},\cdots,X_{t(n)}-X_{t(n-1)}$
are independent, and normally distributed with mean $0$ and variance
$t_{1}-t_{0},\cdots,t_{n}-t_{n-1}$ respectively. 

10. It follows that the restricted process 
\[
Z\equiv X|\overline{Q}_{\infty}:\overline{Q}_{\infty}\times(\Omega,L,E)\rightarrow R
\]
has the property that (i) $Z_{0}=0$, (ii) for each sequence $0\equiv t_{0}\leq t_{1}\leq\cdots\leq t_{n-1}\leq t_{n}$
in $\overline{Q}_{\infty}$, the r.r.v.'s $Z_{t(1)}-Z_{t(0)},\cdots,Z_{t(n)}-Z_{t(n-1)}$
are independent, and (iii) for each $s,t\in\overline{Q}_{\infty}$,
the r.r.v. $Z_{t}-Z_{s}$ is normal with mean $0$ and variance $|t-s|$.
In other words, the process $Z$ satisfies all the conditions in Assertion
1 of \ref{Thm.Construction of  Brownian Motion on =00005B0,inf)}.
Accordingly, the extension-by-limit 
\[
\widetilde{X}_{1}\equiv\Phi_{Lim}(Z):[0,\infty)\times\Omega\rightarrow R
\]
 is a Brownian motion in $R^{1}$. In particular the process $\widetilde{X}$
is a.u. continuous. Hence there exists a full subset $H_{1}$ of $(\Omega,L,E)$
such that
\[
\widetilde{X}_{1}(t,\omega)=\lim_{r\rightarrow t;r\in\overline{Q}(\infty)}Z(r,\omega)=\lim_{r\rightarrow t;r\in\overline{Q}(\infty)}X(r,\omega)\equiv\lim_{r\rightarrow t;r\in\overline{Q}(\infty)}B_{1;}(r,\omega).
\]
for each $t\in[0,\infty),$ for each $\omega\in H_{1}$.

11. For each $k=2,\cdots,m$, repeating the arguments in Steps 8-10
with the process $B_{k;}$in the place of $B_{1;}$, we see that there
exists a full subset $H_{k}$ of $(\Omega,L,E)$ such that
\[
\widetilde{X}_{k}(t,\omega)=\lim_{r\rightarrow t;r\in\overline{Q}(\infty)}B_{k;}(r,\omega)
\]
for each $t\in[0,\infty),$ for each $\omega\in H_{k}$. Combining,
we see that, on the full set $H=H_{1}\cap\cdots\cap H_{m}$, we have
\[
(\widetilde{X}_{1}(t,\omega),\cdots,\widetilde{X}_{m}(t,\omega))=\lim_{r\rightarrow t;r\in\overline{Q}(\infty)}(B_{1;}(r,\omega),\cdots,B_{m;}(r,\omega))\equiv\lim_{r\rightarrow t;r\in\overline{Q}(\infty)}B(r,\omega).
\]
Consequently, the right limit $\lim_{r\rightarrow t;r\in\overline{Q}(\infty)[t,\infty)}B(r,\omega)$
exists and
\[
(\widetilde{X}_{1}(t,\omega),\cdots,\widetilde{X}_{m}(t,\omega))=\lim_{r\rightarrow t;r\in\overline{Q}(\infty)[t,\infty)}B(r,\omega).
\]
for each $t\in[0,\infty),$ for each $\omega\in H$. At the same time,
since $B$ is a.u. càdlàg, there exists a full subset $G\subset\bigcap_{t\in\overline{Q(}\infty)}domain(B_{t})$
of $(\Omega,L,E)$ such that the if the right limit $\lim_{r\rightarrow t;r\in\overline{Q}(\infty)[t,\infty)}B(r,\omega)$
exists, then $t\in domain(B(\cdot,\omega))$ and
\[
B(t,\omega)=\lim_{r\rightarrow t;r\in\overline{Q}(\infty)[t,\infty)}B(r,\omega).
\]
Combining, we see that, for each $\omega\in H\cap G$, we have $domain(B(\cdot,\omega))=[0,\infty)$
and 
\[
B(\cdot,\omega)=(\widetilde{X}_{1}(\cdot,\omega),\cdots,\widetilde{X}_{m}(\cdot,\omega)).
\]
Since the process $(\widetilde{X}_{1},\cdots,\widetilde{X}_{m}):[0,\infty)\times\Omega\rightarrow R^{m}$
is a.u. continuous, it follows that the process $B$ is a.u. continuous.

Summing up, all the conditions in Definition \ref{Def.  Brownian Motion inR^m}
are verified for the process $B$ to be a Brownian motion in $R^{m}$.
$\square$
\end{proof}
\begin{cor}
\label{Cor. Basic-properties-of Brownian Motion in R^m} \textbf{\emph{(Basic
properties of Brownian Motion in $R^{m}$). }}\emph{Let} $B:[0,\infty)\times(\varOmega,L,E)\rightarrow R^{m}$
be an arbitrary Brownian motion in $R^{m}$. Let $\mathcal{L}$ be
an arbitrary right continuous filtration to which the Brownian motion
$B$ is adapted. Then the following holds.

\emph{1.} $B$ is equivalent to the Brownian motion constructed in
Theorem \ref{Thm. Construction of -Brownian motion in R^m}, and is
an a.u. continuous Feller process with Feller semigroup $\mathbf{T}$
defined in Lemma \ref{Lem. Brownian semigroup}.

\emph{2. }The Brownian motion $B$ is strongly Markov relative to
$\mathcal{L}$. Specifically, given any stopping time $\tau$ values
in $[0,\infty)$ relative to $\mathcal{L}$, equality \ref{eq:temp-576-1-1-2-2}
in Theorem \ref{Thm. Construction of Feller process from semigroup}
holds for the process $B$ and the stropping time $\tau$. 

\emph{3.} Let $A$ be an arbitrary orthogonal $k\times m$ matrix.
Thus $AA'=I$ is the $k\times k$ identity matrix. Then the process
$AB:[0,\infty)\times(\varOmega,L,E)\rightarrow R^{k}$ is a Brownian
Motion.

\emph{4.} Let $b$ be an arbitrary unit vector. Then the process $b'B:[0,\infty)\times(\varOmega,L,E)\rightarrow R$
is a Brownian Motion in $R$.

\emph{5.} Let $\gamma>0$ be arbitrary. Define the process $\widetilde{B}:[0,\infty)\times(\varOmega,L,E)\rightarrow R^{m}$
by $\widetilde{B}_{t}\equiv\gamma^{-1/2}B_{\gamma t}$ for each $t\in[0,\infty)$.
Then $\widetilde{B}$ is a Brownian Motion adapted to the right continuous
filtration $\mathcal{L}_{\gamma}\equiv\{L^{(\gamma t)}:t\in[0,\infty)\}$. 
\end{cor}
\begin{proof}
1. Let $\widehat{B}:[0,\infty)\times(\widehat{\varOmega},\widehat{L},\widehat{E})\rightarrow R^{m}$
denote the Brownian motion constructed in Theorem \ref{Thm. Construction of -Brownian motion in R^m}.
Thus $\widehat{B}$ is an a.u. continuous Feller process. Hence it
has marginal f.j.d.'s given by the family $F^{0,\mathsf{\mathbf{T}}}$
generated by the Feller semigroup $\mathbf{T}$. Define the processes
$Z\equiv B|\overline{Q}_{\infty}$and $\widehat{Z}\equiv\widehat{B}|\overline{Q}_{\infty}$.
Let $m\geq0$ and the sequence $0\equiv r_{0}\leq r_{1}\leq\cdots\leq r_{m}$
in $\overline{Q}_{\infty}$, and $f\in C(R^{m+1})$ be arbitrary.
Then
\[
Ef(Z_{0},Z_{r(1)},\cdots,Z_{r(m)})=Ef(B_{0},B_{r(1)},\cdots,B_{r(m)})
\]
\[
=Ef(0,B_{r(1)}-B_{r(0)},(B_{r(1)}-B_{r(0)})
\]
\[
+(B_{r(2)}-B_{r(1)}),\cdots,(B_{r(1)}-B_{r(0)})+\cdots+(B_{r(m)}-B_{r(m-1)}))
\]
\[
=\widehat{E}f(0,\widehat{B}_{r(1)}-\widehat{B}_{r(0)},(\widehat{B}_{r(1)}-\widehat{B}_{r(0)})
\]
\[
+(\widehat{B}_{r(2)}-\widehat{B}_{r(1)}),\cdots,(\widehat{B}_{r(1)}-\widehat{B}_{r(0)})+\cdots+(\widehat{B}_{r(m)}-\widehat{B}_{r(m-1)}))
\]
\[
=\widehat{E}f(\widehat{Z}_{0},\widehat{Z}_{r(1)},\cdots,\widehat{Z}_{r(m)})=F_{r(0),\cdots,r(m)}^{0,\mathbf{T}}f,
\]
where the third equality follows from Conditions (i) and (ii) of Definition
\ref{Def.  Brownian Motion inR^m}. In other words, the marginal f.j.d.'s
of $B|\overline{Q}_{\infty}$and $\widehat{B}|\overline{Q}_{\infty}$are
equal. Since the Brownian motions $B$ and $\widehat{B}$ are a.u.
continuous, it follows that the marginal f.j.d.'s of $B,\widehat{B}$
are equal, and given by the family $F^{0,\mathbf{T}}$. Thus the Brownian
motion $B$ also is generated by the initial state $0$ and Feller
semigroup $\mathbf{T}$. Assertion 1 is proved.

2. Consequently, Assertion 2 then follows from Definition \ref{Def. Feller process}. 

3. Let $A$ be an orthogonal $k\times m$ matrix. Then trivially $AB_{0}=A0=0\in R^{k}$.
Thus Condition (i) of Definition \ref{Def.  Brownian Motion inR^m}
holds for the process $AB$. Next, let the sequence $0\equiv t_{0}\leq t_{1}\leq\cdots\leq t_{n-1}\leq t_{n}$
in $[0,\infty)$ be arbitrary. Then the r.v.'s $B_{t(1)}-B_{t(0)},\cdots,B_{t(n)}-B_{t(n-1)}$
are independent. Hence the r.v.'s $A(B_{t(1)}-B_{t(0)}),\cdots,A(B_{t(n)}-B_{t(n-1)})$
are independent, establishing Condition (ii) for the process $AB$.
Now let $s,t\in[0,\infty)$ be arbitrary. Then the r.v. $B_{t}-B_{s}$
is normal with mean $0$ and covariance matrix $|t-s|I$, where $I$
stands for the $k\times k$ identity matrix. Hence $A(B_{t}-B_{s})$
is normal with mean $0$, with covariance matrix
\[
E(A(B_{t}-B_{s})(B_{t}-B_{s})'A')=A(E(B_{t}-B_{s})(B_{t}-B_{s})')A'
\]
\[
=A(|t-s|IA')=|t-s|AA'=|t-s|I.
\]
This proves Condition (iii) for the process $AB$. Assertion 3 is
proved.

4. Let $b$ be an arbitrary unit vector. Then $u'$ is a $1\times m$
matrix. Hence, according to Assertion 3, the process $u'B:[0,\infty)\times(\varOmega,L,E)\rightarrow R$
is a Brownian Motion in $R$.

5. Define the process $\widetilde{B}:[0,\infty)\times(\varOmega,L,E)\rightarrow R^{m}$
by $\widetilde{B}_{t}\equiv\gamma^{-1/2}B_{\gamma t}$ for each $t\in[0,\infty)$.
Trivially, Conditions (i) and (ii) of Definition \ref{Def.  Brownian Motion inR^m}
holds for the process $\widetilde{B}$. Let $s,t\in[0,\infty)$ be
arbitrary. Then the r.v. $B_{\gamma t}-B_{\gamma s}$ is normal with
mean $0$ and covariance matrix $|\gamma t-\gamma s|I$. Hence the
r.v. $\widetilde{B}_{t}-\widetilde{B}_{s}\equiv\gamma^{-1/2}B_{\gamma t}-\gamma^{-1/2}B_{\gamma s}$
has covariance matrix 
\[
(\gamma^{-1/2})^{2}|\gamma t-\gamma s|I=|t-s|I.
\]
Thus Conditions (iii) of Definition \ref{Def.  Brownian Motion inR^m}
is also verified for the process $\widetilde{B}$ to be a Brownian
motion. Moreover, for each $t\in[0,\infty)$, we have $\widetilde{B}_{t}\equiv\gamma^{-1/2}B_{\gamma t}\in L^{(\gamma t)}$.
Hence the process $\widetilde{B}$ is adapted to the filtration $\mathcal{L}_{\gamma}\equiv\{L^{(\gamma t)}:t\in[0,\infty)\}$.
Assertion 4 and the corollary are proved.
\end{proof}

\section{First Exit Times from Spheres by the Brownian Motion in\emph{ $R^{m}$}}

In this section, let $m\geq1$ be arbitrary, but fixed. We will use
the terminology and notations in the previous section. 
\begin{defn}
\label{Def. Notations for Brownian motions in R^m}\textbf{ (Notations
for a filtered Brownian motion in $R^{m}$).} Let $B:[0,\infty)\times(\Omega,L,E)\rightarrow R^{m}$
be an arbitrary Brownian motion. Let $\mathcal{L}$ be an arbitrary
right continuous filtration to which the process $B$ is adapted.
Unless otherwise stated, all stopping times are relative to this right
continuous filtration $\mathcal{L}$. Let $x\in R^{m}$ be arbitrary.
Recall that the process $B^{x}\equiv x+B$ is then called a Brownian
motion starting at \emph{$x$}. Trivially, the process $B^{x}$ is
an a.u. càdlàg Feller process with the same Brownian semigroup $\mathbf{T}\equiv\{T_{t}:t\in[0,\infty)\}$,
defined in Lemma \ref{Lem. Brownian semigroup}, and is adapted to
the filtration $\mathcal{L}$ like $ $. Corollary \ref{Cor. Basic-properties-of Brownian Motion in R^m}
then says that the process $B^{x}$ is an a.u. continuous Feller process
relative to the filtration $\mathcal{L}$ and the Feller semigroup
$\mathbf{T}$, and therefore that it is strongly Markov relative to
the filtration $\mathcal{L}$ and the Feller semigroup $\mathbf{V}$,
in the sense of Conditions 3 in Theorem \ref{Thm. Construction of Feller process from semigroup}. 
\end{defn}
$\square$
\begin{defn}
\textbf{\label{Def. Notations for spheres and boundaries} (Notations
for spheres and their boundaries).} Let $x\in R^{m}$ and $r>0$ be
arbitrary. Define the open sphere $D_{x,r}\equiv\{z\in R^{m}:|z-x|<r\}$,
its closure 
\[
\overline{D}_{x,r}\equiv\{z\in R^{m}:|z-x|\leq r\},
\]
and its boundary 
\[
\partial D_{x,r}\equiv\{z\in R^{m}:|z-x|=r\}
\]
Define $D_{x,r}^{c}\equiv\{z\in R^{n}:|z-x|\geq r\}$. We will write,
for abbreviation, $D,\partial D,D^{c}$ for $D_{0,1},\partial D_{0,1},D_{0,1}^{c}$
respectively. 
\end{defn}
$\square$
\begin{defn}
\textbf{\emph{\label{Def:first exit time of open spheres by Brownian motions in R^m-1}}}\textbf{(First
exit times from spheres by Brownian Motion} \textbf{in $R^{m}$).}
Let $B:[0,\infty)\times(\Omega,L,E)\rightarrow R^{m}$ be an arbitrary
Brownian motion. Let $\mathcal{L}$ be an arbitrary right continuous
filtration to which the process $B$ is adapted. Let $x,y\in R^{m}$
and $r>0$ be arbitrary such that $|x-y|<r$. Suppose there exists
a stopping time $\tau_{x,y,r;B}$ with values in $(0,\infty)$, relative
to the right continuous filtration $\mathcal{L}$, such that, for
a.e. $\omega\in\Omega$, we have (i) $B^{x}(\cdot,\omega)\in D_{y,r}$
on the interval $[0,\tau_{x,y,r;B}(\omega))$, and (ii)\emph{ $B_{\tau(x,y,r;B)}^{x}(\omega)\in\partial D_{y,r}$. }

Then the stopping time $\tau_{x,y,r;B}$ is called\emph{ }the\emph{
\index{first exit time from open sphere}first exit time from the
open sphere $D_{y,r}$ }by the\emph{ }Brownian motion $B^{x}$ starting
at $x$. Note that, if $\tau_{x,y,r;B}$ exists, then, because the
Brownian motion $B$ is an a.u. càdlàg Feller process, Assertion 3
of Theorem \ref{Thm. Construction of Feller process from semigroup}
implies that $B_{\tau(x,y,r;B)}^{x}$ is a r.v. with values in $R^{m}$,
and is measurable relative to $L^{(\tau(x,y,r;B))}$.

Note that, in contrast to the case of a general Feller process, the
exit time here is over the infinite interval. The following Theorem
\ref{Thm:first exit time of open spheres by Brownian motions in R^m}
will prove tat it exists as a r.r.v.; the Brownian motion exits any
open ball in some finite time.

When the Brownian motion $B$ is understood from context, we will
write simply $\tau_{x,y,r}\equiv\tau_{x,y,r;B}$ . $\square$
\end{defn}
\begin{lem}
\emph{\label{Lem. Translational invariance of certain exit times}
}\textbf{\emph{(Translation-invariance of certain first exit times).
}}Let $B:[0,\infty)\times(\Omega,L,E)\rightarrow R^{m}$ be an arbitrary
Brownian motion. Let $\mathcal{L}$ be an arbitrary right continuous
filtration to which the process $B$ is adapted. Let $x,y,z\in R^{m}$
and $a,r>0$ be arbitrary such that $|x-y|<r$. If the first exit
time $\tau_{x,y,r}$ exists, then \emph{(i)} the first exit time $\tau_{x-z,y-z,r}$
exists, in which case $\tau_{x-z,y-z,r}=\tau_{x,y,r}$, with $u+B_{t}^{x-u,}=B_{t}^{x}$
for each t$\in]0,\tau_{x,y,r}]$, and \emph{(ii)} the first exit time
$\tau_{ax,ay,ar}$ exists, in which case $\tau_{ax,ay,ar}=\tau_{x,y,r}$.
\end{lem}
\begin{proof}
Suppose the first exit time $\tau_{x,y,r}$ exists. Then, by Definition
\ref{Def:first exit time of open spheres by Brownian motions in R^m-1},
for a.e. $\omega\in\Omega$, we have (i) $x+B(\cdot,\omega)\equiv B^{x}(\cdot,\omega)\in D_{y,r}$
on the interval $[0,\tau_{x,y,r}(\omega))$, and (ii)\emph{ $x+B_{\tau(x,y,r)}(\omega)\in\partial D_{y,r}$.
}Equivalently, (i') $B^{x-z}(\cdot,\omega)\in D_{y-z,r}$ on the interval
$[0,\tau_{x,y,r}(\omega))$, and (ii')\emph{ $B_{\tau(x,y,r)}^{x-z}(\omega)\in\partial D_{0,r}$.}
Hence, again by Definition \ref{Def:first exit time of open spheres by Brownian motions in R^m-1},\emph{
}the first exit time $\tau_{x-z,y-z,r}$ exists and is equal to $\tau_{x,y,r}$.
Assertion (i) is proved. The proof of Assertion (ii) is similar.
\end{proof}
\begin{lem}
\textbf{\emph{\label{Lem:first exit time from certain open spheres by Brownian motions in R^m-1}
(Existence of first exit times from certain open spheres by Brownian
motions).}} Let $B:[0,\infty)\times(\Omega,L,E)\rightarrow R^{m}$
be an arbitrary Brownian motion. Let $\mathcal{L}$ be an arbitrary
right continuous filtration to which the process $B$ is adapted.
Then there exists a countable subset $H$ of $R$ such that for each
$x\in R^{m}$ and $a\in(|x|,\infty)H_{c}$, the\emph{ }first exit
time $\tau_{x,0,a}$ exists. Here $H_{c}$ denotes the metric complement
of $H$ in $R$. 
\end{lem}
\begin{proof}
1. Define the function $f:R^{m}\rightarrow R$ by 
\begin{equation}
f(y)\equiv|y|\label{eq:temp-586}
\end{equation}
for each $y\in R^{m}$. Then $f$ is a uniformly continuous function
on each compact subset of $R^{m}$. Moreover, $(f<a)=D_{0,a}$ and
$(f=a)=\partial D_{0,a}$ for each $a>0$. 

2. For each $x\in R^{m}$, the Brownian motion $B^{x}$ is an a.u.
càdlàg Feller process, with $f(B_{0}^{x})\equiv|B_{0}^{x}|=|x|$.
Theorem \ref{Thm. Abundance of first exit times for Feller process}
therefore implies that there exists a countable subset $H$ of $R$
such that, for each $x\in R^{m}$, for each $a\in(|x|,\infty)H_{c}$,
and $M\geq1$for each , the first exit time $\overline{\tau}_{f,a,M}$
exists for the process $B^{x}$. 

3. Consider each $x\in R^{m}$,
\[
a\in(|x|,\infty)H_{c},
\]
and $M\geq1$. Then the first exit time $\overline{\tau}_{f,a,M}$
exists according to Step 2. Recall that, by Definition \ref{Def. First exit time},
$\overline{\tau}_{f,a,M}$ is the\emph{ }first exit time in\emph{
}$[0,M]$ of the open subset $(f<a)=D_{0,a}$ by the process $B^{x}$.
Moreover, by Definition \ref{Def. First exit time}, the function
$B_{\overline{\tau}(f,a,M)}^{x}$ is a well-defined r.v. Furthermore,
Lemma \ref{Lem. Basics of first exit times} says that, for each $r\in(0,M)$
and $N\geq M$ , we have

\begin{equation}
\overline{\tau}_{f,a,M}\leq\overline{\tau}_{f,a,N},\label{eq:temp-587}
\end{equation}
\begin{equation}
(\overline{\tau}_{f,a,M}<M)\subset(\overline{\tau}_{f,a,N}=\overline{\tau}_{f,a,M}),\label{eq:temp-588}
\end{equation}
and
\begin{equation}
(\overline{\tau}_{f,a,N}\leq r)=(\overline{\tau}_{f,a,M}\leq r).\label{eq:temp-589}
\end{equation}
We will prove that the a.u. limit 
\[
\tau_{a}\equiv\lim_{M\rightarrow\infty}\overline{\tau}_{f,a,M}
\]
exists, is a stopping time relative to the filtration $\mathcal{L}$,
and satisfies the conditions in Definition \ref{Def:first exit time of open spheres by Brownian motions in R^m-1}
to be the first exit time $\tau_{x,0,a}$ from the open sphere $D_{0,a}$
by the\emph{ }Brownian motion $B^{x}$.

Recall that $B_{1;}:[0,\infty)\times(\Omega,L,E)\rightarrow R$ denotes
the first component of the Brownian motion $B$. As such, $B_{1;}$
is a Brownian motion in $R$. Suppose $M\geq2$. Then the r.r.v. $B_{1;M-1}$
has a normal distribution on $R$, with mean $0$ and variance $M-1$.
Since said normal distribution has a p.d.f. $\varphi_{0,M-1}$ which
is bounded by $(2\pi(M-1))^{-1/2}$, the set
\[
A_{M}\equiv(|B_{1;M-1}|<a)
\]
is measurable, with
\[
P(A_{M})=\int_{u=-2a}^{2a}\varphi_{0,M-1}(du)\leq4a(2\pi(M-1))^{-1/2}\rightarrow0
\]
as $M\rightarrow\infty$. At the same time,
\[
A_{M}^{c}=(|B_{1;M-1}|\geq a)\subset(|B_{M-1}^{x}|\geq a)\equiv(f(B_{M-1}^{x})\geq a).
\]
Therefore, by Condition (i) of Definition \ref{Def. First exit time}
for the first exit time $\overline{\tau}_{f,a,M}$, we have 
\[
A_{M}^{c}\subset(f(B_{M-1}^{x})\geq a)\subset(M-1\geq\overline{\tau}_{f,a,M})
\]
\begin{equation}
\subset(\overline{\tau}_{f,a,M}<M)\subset\bigcap_{N=M}^{\infty}(\overline{\tau}_{f,a,N}=\overline{\tau}_{f,a,M})\label{eq:temp-522-1}
\end{equation}
where the next-to-last inclusion relation is thanks to relation \ref{eq:temp-588},
and where the last inclusion relation is trivial. Hence, on $A_{M}^{c}$,
we have $\overline{\tau}_{f,a,N}\uparrow\tau_{a}$ uniformly for some
function $\tau_{a}$. Since $P(A_{M}^{c})$ is arbitrarily close to
$1$ if $M$ is sufficiently large, we conclude that $\overline{\tau}_{f,a,N}\uparrow\tau_{a}$
a.u. and therefore that $\tau_{a}$ is a r.r.v. Note that $\overline{\tau}_{f,a,1}>0$.
Hence $\tau_{a}>0$. In other words, $\tau_{a}$ is a r.r.v. with
values in $(0,\infty)$.

3. To prove that the r.r.v. $\tau_{a}$ is a stopping time, consider
each $t\in(0,\infty)$ be an arbitrary regular point of $\tau_{a}$.
Let $s\in(t,\infty)$ be arbitrary. Let $r\in(t,s)$ be an arbitrary
regular point of the r.r.v.'s in the sequence $\tau_{a},\overline{\tau}_{f,a,1},\overline{\tau}_{f,a,2},\cdots$.
Take any $M>r$. Then, according to equality \ref{eq:temp-589}, we
have $(\overline{\tau}_{f,a,M}\leq r)=(\overline{\tau}_{f,a,N}\leq r)$
for each $N\geq M$. Letting $N\rightarrow\infty,$we obtain $(\overline{\tau}_{f,a,M}\leq r)=(\tau_{a}\leq r)$.
It follows that 
\[
(\tau_{a}\leq r)=(\overline{\tau}_{f,a,M}\leq r)\in L^{(r)}\subset L^{(s)}
\]
because $\overline{\tau}_{f,a,M}$ is a stopping time relative to
the right continuous filtration $\mathcal{L}$. Therefore, if we take
a sequence $(r_{k})_{k=1,2,\cdots}$ of regular points in $(t,s)$
of the r.r.v.'s in the sequence $\tau_{a},\tau_{f,a,1},\tau_{f,a,2},\cdots$,
such that $r_{k}\downarrow t$, then we obtain
\[
(\tau_{a}\leq t)=\bigcap_{k=1}^{\infty}(\tau_{a}\leq r_{k})\in L^{(s)}.
\]
Since $s\in(t,\infty)$ is arbitrary, it follows that 
\[
(\tau_{a}\leq t)\in\bigcap_{s\in(t,\infty)}L^{(s)}\equiv L^{(t+)}=L^{(t)},
\]
where the last equality is because the filtration $\mathcal{L}$ is
right continuous by assumption. Thus the r.r.v. $\tau_{a}$ is a stopping
time.

4. It remains to show that $\tau_{a}$ is the desired first exit time.
To that end, first note that, by Assertion 1 of Lemma \ref{Lem. Observability at stopping time},
$B_{\tau(a)}^{x}$ is a well defined r.v. and is measurable relative
to $L^{(\tau(a))}.$ Define the null set
\[
A\equiv\bigcap_{M=2}^{\infty}A_{M}.
\]
Consider each $\omega\in A^{c}$. Then $\omega\in A_{M}^{c}$ for
some $M\geq2$ with 
\[
M>u\equiv\tau_{a}(\omega).
\]
Hence relation \ref{eq:temp-522-1} implies that 
\begin{equation}
\omega\in(\tau_{a}=\overline{\tau}_{f,a,M})\subset(B_{\tau(a)}^{x}=B_{\overline{\tau}(f,a,M)}^{x}).\label{eq:temp-522-1-1}
\end{equation}
Consequently, $\overline{\tau}_{f,a,M}(\omega)=u<M$. By Definition
\ref{Def. First exit time} of the first exit time $\overline{\tau}_{f,a,M}$,
we then have (i) $f(B^{x}(\cdot,\omega))<a$ on the interval $[0,u)$,
and (ii)\emph{ $f(B_{u}^{x}(\omega))\geq a$.} Conditions (i) and
(ii) can be rewritten as conditions (i') $f(B^{x}(\cdot,\omega))<a$
on the interval $[0,\tau_{a}(\omega))$, and (ii')\emph{ $f(B_{\tau(a)}^{x}(\omega))\geq a$.}
Since $B^{x}(\cdot,\omega)$ is a continuous function, we obtain (i'')
$f(B^{x}(\cdot,\omega))<a$ on the interval $[0,\tau_{a}(\omega))$,
and (ii'') $f(B_{\tau(a)}^{x}(\omega))=a$. Therefore, by the observation
in Step 1, we have (i''') $B^{x}(\cdot,\omega))\in D_{0,a}$ on the
interval $[0,\tau_{a}(\omega))$, and (ii''')\emph{ $f(B_{\tau(a)}^{x}(\omega))\in\partial D_{0,a}$,}
where\emph{ }$\omega\in A^{c}$ is arbitrary. By restricting the domain
of the r.r.v. $\tau_{a}$ to the full set $A^{c}$, we obtain conditions
(i''') and (iii''') for each\emph{ }$\omega\in domain(\tau_{a})$. 

We have verified all the conditions in Definition \ref{Def:first exit time of open spheres by Brownian motions in R^m-1}
for the stopping time $\tau_{a}$ to be the first exit time $\tau_{x,0,a}$
from the open sphere $D_{0,a}$ by the\emph{ }Brownian motion $B^{x}$,
where $a\in(|x|,\infty)H_{c}$ are arbitrary. The lemma is proved.
\end{proof}
\begin{lem}
\textbf{\emph{\label{Lem. First level-crossing time for Brownian motion in R}(First
level-crossing time for Brownian motion in $R^{1}$, and the \index{Reflection Principle}Reflection
Principle).}} Suppose $m=1$. Let $x\in R^{1}$ be arbitrary. Then
there exists a countable subset $G$ of $R$ with the following properties.

Let $\overline{B}:[0,\infty)\times(\Omega,L,E)\rightarrow R^{1}$
be an arbitrary Brownian motion adapted to some right continuous filtration
$\mathcal{\overline{L}}$. Let $a\in(x,\infty)G_{c}$ be arbitrary.
Then there exists a stopping time $\overline{\tau}_{x,a}$ with values
in $[0,1]$ relative to $\mathcal{\overline{L}}$, such that for a.e.
$\omega\in\Omega$, we have \emph{(i)} $\overline{B}^{x}(\cdot,\omega)<a$
on the interval $[0,\overline{\tau}_{x,a}(\omega))$, and \emph{(ii)
$\overline{B}_{\overline{\tau}(x,a)}^{x}(\omega)=a$ if $\overline{\tau}_{x,a}(\omega)<1$.
}Moreover, each $t\in[0,1]$ is a regular point of the r.r.v. $\overline{\tau}_{x,a}$,
with 
\begin{equation}
P(\overline{\tau}_{x,a}\leq t)=2P(\overline{B}_{t}>a-x)=2(1-\Phi_{0,t}(a-x)).\label{eq:temp-601}
\end{equation}
We will call $\overline{\tau}_{x,a}$ a \emph{\index{level crossing time}level-crossing
time} by the Brownian motion $\overline{B}$.

This lemma is sufficient for our immediate application in the next
theorem. We note that the lemma can easily be strengthened to have
an empty set $G$ of exceptional points. 
\end{lem}
\begin{proof}
1. Define the function $f:R\rightarrow R$ by 
\begin{equation}
f(y)\equiv y\label{eq:temp-586-1}
\end{equation}
for each $y\in R$. Then, trivially, $f$ is uniformly continuous
and bounded on bounded sets. Let $x\in R^{1}$ be arbitrary. Then
the Brownian motion $\overline{B}^{x}$ is an a.u. càdlàg Feller process,
with $f(\overline{B}_{0}^{x})\equiv\overline{B}_{0}^{x}=x$. Therefore
Assertion 1 of Theorem \ref{Thm. Abundance of first exit times for Feller process}
says that there exists a countable subset $G$ of $R$ such that,
for each point $a\in(x,\infty)G_{c}$, and for $M\equiv1$, the first
exit time $\overline{\tau}_{x,a}\equiv\overline{\tau}_{f,a,1}$ exists
for the process $\overline{B}^{x}$. Consider each 
\[
a\in(x,\infty)G_{c}.
\]
By Definition \ref{Def. First exit time} for first exit times, we
have, for a.e. $\omega\in\Omega$, the conditions (i') $f(\overline{B}^{x}(\cdot,\omega))<a$
on the interval $[0,\overline{\tau}_{x,a}(\omega))$, and (ii')\emph{
$f(\overline{B}_{\overline{\tau}(x,a)}^{x}(\omega))=a$ }if\emph{
$\overline{\tau}_{x,a}(\omega)<1$. }In view of the defining formula
\ref{eq:temp-586-1}, Conditions (i') and (ii') are equivalent to
the desired Conditions (i) and (ii). It remains to verify equality
\ref{eq:temp-601}.

2. To that end, let $t\in(0,1]$ be arbitrary. First assume that $t$
is a regular point of the r.r.v. $\overline{\tau}_{x,a}$, such that
$t\neq u$ for each $u\in\overline{Q}_{\infty}$. Write, for abbreviation,
$\tau\equiv\overline{\tau}_{x,a}$. Let $(\eta_{h})_{h=0,1,\cdots}$
be an arbitrary non increasing sequence of stopping times such that,
for each $h\geq0$, the r.r.v. $\eta_{h}$ has values in $\overline{Q}_{h}$
, and such that 
\begin{equation}
\tau\leq\eta_{h}<\tau+2^{-h+2}.\label{eq:temp-421-1-2-1-2-1}
\end{equation}
Note that such a sequence exists by Assertion 1 of Proposition \ref{Prop. Stopping times rlative to right continuous filtration}.
Then $\eta_{h}\rightarrow\tau$ a.u. Hence $\overline{B}_{\eta(h)}\rightarrow\overline{B}_{\tau}$
a.u., as $h\rightarrow\infty$, thanks to the a.u. continuity of the
Brownian motion $\overline{B}$. Define the indicator function $g:R^{2}\rightarrow R$
by $g(z,w)\equiv1_{(w-z>0)}$ for each $(z,w)\in R^{2}$. For convenience
of notations, let $\widetilde{B}:[0,\infty)\times(\Omega,L,E)\rightarrow R^{1}$
be a second Brownian motion adapted to some right continuous filtration
$\mathcal{\widetilde{L}}$, such that each measurable function relative
to $L(\widetilde{B}_{t}:t\geq0)$ and each measurable function relative
to $L(\overline{B}_{t}:t\geq0)$ are independent. The reader can easily
verify that such an independent copy of $B$ exists, by using replacing
$(\Omega,L,E)$ with the product space $(\Omega,L,E)\otimes(\Omega,L,E)$,
and by identifying $\overline{B}$ with $\overline{B}\otimes1$, and
identifying $\widetilde{B}$ with $1\otimes$. Using basic properties
of first exit times and a.u. convergence, we calculate

\[
P(\tau<t;\overline{B}_{t}^{x}>a)=P(\tau<t;\overline{B}_{\tau}^{x}=a;\overline{B}_{t}^{x}>a)
\]
\[
=P(\tau<t;\overline{B}_{\tau}^{x}=a;\overline{B}_{t}^{x}-\overline{B}_{\tau}^{x}>0)=P(\tau<t;\overline{B}_{t}^{x}-\overline{B}_{\tau}^{x}>0)
\]
\[
=P(\tau<t;\overline{B}_{t}-\overline{B}_{\tau}>0)=\lim_{h\rightarrow\infty}P(\eta_{h}<t;\overline{B}_{t}-\overline{B}_{\eta(h)}>0)
\]
\[
=\lim_{h\rightarrow\infty}\sum_{u\in\overline{Q}(h)[0,t)}P(\eta_{h}=u;\overline{B}_{t}-\overline{B}_{u}>0)=\lim_{h\rightarrow\infty}\sum_{u\in\overline{Q}(h)[0,t)}Eg(\overline{B}_{u},\overline{B}_{t})1_{\eta(h)=u}
\]
\begin{equation}
=\lim_{h\rightarrow\infty}\sum_{u\in\overline{Q}(h)[0,t)}E(F_{0,t-u}^{\overline{B}(u),\mathbf{T}}g)1_{\eta(h)=u},\label{eq:temp-608}
\end{equation}
where $\mathbf{T}$ is the Brownian semigroup defined in Lemma \ref{Lem. Brownian semigroup},
and where, for each $z\in R$, the family $F^{z,\mathbf{T}}$ is generated
by the initial state $z$ and the Feller semigroup \emph{$\mathbf{T}$,
}in the sense of Theorem \ref{Thm. Construction of Feller process from semigroup}.
Thus, for each $z\in R$ and $u<t$, we have
\[
F_{0,t-u}^{z,\mathbf{T}}g\equiv\int T_{t-u}^{z}(dw)g(z,w)=\int_{w\in R}\Phi_{z,t-u}(dw)1_{(w-z>0)}=\int_{w=z}^{\infty}\Phi_{z,t-u}(dw)=\frac{1}{2}.
\]
Hence equality \ref{eq:temp-608} yields
\[
P(\tau<t;\overline{B}_{t}^{x}>a)=\lim_{h\rightarrow\infty}\sum_{u\in\overline{Q}(h)[0,t)}E(\frac{1}{2};\eta_{h}=u)=\frac{1}{2}P(\tau<t),
\]
Note that $(\overline{B}_{t}^{x}>a)\subset(\overline{\tau}_{x,a}<t)$
by the definition of the first level-crossing time $\tau\equiv\overline{\tau}_{x,a}$.
Combining,
\[
P(\overline{\tau}_{x,a}<t)\equiv P(\tau<t)=2P(\tau<t;\overline{B}_{t}^{x}>a)=2P(\overline{B}_{t}^{x}>a)
\]
\begin{equation}
=2P(\overline{B}_{t}\geq a-x)=2(1-\Phi_{0,t}(a-x)),\label{eq:temp-579}
\end{equation}
where $t\in(0,1]$ is an arbitrary regular point of the r.r.v. $\overline{\tau}_{x,a}$,
with $t\neq u$ for each $u\in\overline{Q}_{h}$. Since such points
are dense in the interval $(0,1]$, the continuity of the right-hand
side of equality \ref{eq:temp-579} as a function of $t$ implies
that each $t\in[0,1]$ is a regular point of the r.r.v. $\overline{\tau}_{x,a}$,
with 
\[
P(\overline{\tau}_{x,a}<t)=2(1-\Phi_{0,t}(a-x)),
\]
where $a\in(x,\infty)G_{c}$ is arbitrary. The lemma is proved.
\end{proof}
.The following theorem enables us to dispense with the exceptional
set $H$ in Lemma \ref{Lem:first exit time from certain open spheres by Brownian motions in R^m-1},
and prove that the first exit time exists, from any open sphere by
a Brownian motion started in the interior. It also estimates some
bounds between successive first exit times from two concentric open
spheres with approximately equal radii. 
\begin{thm}
\textbf{\emph{\label{Thm:first exit time of open spheres by Brownian motions in R^m}
(Existence of first exit times from open spheres by Brownian motions).}}
Let $B:[0,\infty)\times(\Omega,L,E)\rightarrow R^{m}$ be an arbitrary
Brownian motion. Let $\mathcal{L}$ be an arbitrary right continuous
filtration to which the process $B$ is adapted. Let $x\in R^{m}$
be arbitrary. Then the following holds.

\emph{1. }Let $b>|x|$ be arbitrary. Then the\emph{ }first exit time
$\tau_{x,0,b}$ from the open sphere $D_{0,b}$ by the Brownian motion
$B^{x}$ exists. 

\emph{2}. Let $r>0$ be arbitrary. Then the first exit time $\tau_{x,x,r}$
from the open sphere $D_{x,r}$ by the Brownian motion $B^{x}$exists,
and is equal to $\tau_{0,0,r}$. In short $\tau_{x,x,r}=\tau_{0,0,r}$.

\emph{3}. Now let $k\geq1$ and $b'>b''$ be arbitrary with $b'-b''<2^{-k}.$
Then 
\[
\tau_{x,0,b'}-\tau_{x,0,b''}<2^{-k}
\]
on some measurable set $A$ with $P(A^{c})<2^{-k/2}$. 
\end{thm}
\begin{proof}
1. By Lemma \ref{Lem:first exit time from certain open spheres by Brownian motions in R^m-1},
there exists a countable subset $H$ of $R$ such that the first exit
time $\tau_{x,0,a}$ exists for each $a\in(|x|,\infty)H_{c}$. By
Lemma \ref{Lem. First level-crossing time for Brownian motion in R},
there exists a countable subset $G$ of $R$, such that for each $a>0$
and $c\in(a,\infty)G_{c}$, we have (i) the first level-crossing time
$\overline{\tau}_{a,c}$ with values in $(0,1]$ exists relative to
$\overline{B}$, and (ii) each $t\in[0,1]$ is a regular point of
$\overline{\tau}_{a,c}$, with
\begin{equation}
P(\overline{\tau}_{a,c}<t)=2(1-\Phi_{0,t}(c-a)).\label{eq:temp-251}
\end{equation}

2. To prove Assertion 1, let $b>|x|$ be arbitrary. Let $n\geq1$
be arbitrary, but fixed till further notice. For abbreviation, write
$\varepsilon_{n}\equiv t_{n}\equiv2^{-n}$. Take $a_{n}\in(|x|,b)H_{c}$
such that $b-a_{n}<\varepsilon_{n}$. Take $c_{n}\in(b,b+\varepsilon_{n})H_{c}G_{c}.$
Then $a_{n}<b<c_{n}$, and the first exit times $\tau_{x,0,a(n)},\tau_{x,0,c(n)}$,
and the first level-crossing time $\overline{\tau}_{a(n),c(n)}$ exist.
In the following, for abbreviation, write $\tau\equiv\tau_{x,0,a(n)}$
and $\widehat{\tau}\equiv\tau_{x,0,c(n)}$.

3. For later reference, we estimate some normal probability. Note
that $c_{n}-a_{n}<2\varepsilon_{n}$. Hence
\[
\frac{c_{n}-a_{n}}{\sqrt{t_{n}}}<\frac{2\varepsilon_{n}}{\sqrt{\varepsilon_{n}}}=2\sqrt{\varepsilon_{n}}.
\]
Therefore
\[
\Phi_{0,t(n)}(c_{n}-a_{n})=\frac{1}{2}+\int_{0}^{\frac{c(n)-a(n)}{\sqrt{t(n)}}}\varphi_{0,1}(v)dv\leq\frac{1}{2}+\int_{0}^{2\sqrt{\varepsilon_{n}}}\varphi_{0,1}(v)dv
\]
\begin{equation}
\leq\frac{1}{2}+\frac{2\sqrt{\varepsilon_{n}}}{\sqrt{2\pi}}\leq\frac{1}{2}+2^{-n/2}.\label{eq:temp-522}
\end{equation}

4. Recall that each point $y\in R^{m}$ is regarded as a column vector,
with transpose denoted by $y'$. Let $u\in R^{m}$ be arbitrary such
that $|u|=1$. In other words, $u$ is a unit column vector. Define
the process 
\[
\overline{B}:[0,\infty)\times(\Omega,L,E)\rightarrow R
\]
by 
\[
\overline{B}_{s}\equiv a_{n}^{-1}(B_{\tau}^{x})'(B_{s+\tau}^{x}-B_{\tau}^{x})
\]
for each $s\in[0,\infty)$. In other words, we reset the clock to
$0$ at the stopping time $\tau$, and starts tracking the projection,
in the direction $u$, of the Brownian motion $B^{x}$ relative to
the new origin $B_{\tau}^{x}$. Note that
\[
\overline{B}_{s}\equiv a_{n}^{-1}(B_{\tau}^{x})'(B_{s+\tau}^{x}-B_{\tau}^{x})
\]
for each $s\in[0,\infty)$.

5. We will next prove that, with $x$ and $a_{n}$ fixed, the process
$\overline{B}$ is a Brownian motion in $R^{1}$. To that end, let
$\overline{\mathbf{T}}\equiv\{\mathbb{T}_{t}:t\geq0$\} denote the
$1$-dimensional Brownian semigroup in the sense of Lemma \ref{Lem. Brownian semigroup}.
Trivially $\overline{B}_{0}=0$, and the process $\overline{B}$ inherits
a.u. continuity from the process $B^{x}$. Write, for abbreviation,
$Y_{s}\equiv B_{\tau(x,0,a(n))+s}^{x}$ for each s$\geq0$, and write
$U\equiv a_{n}^{-1}Y_{0}$. Then
\begin{equation}
U'Y_{0}=a_{n}^{-1}Y_{0}'Y_{0}=a_{n}^{-1}|Y_{0}|^{2}=a_{n}^{-1}|B_{\tau(x,0,a(n))}^{x}|^{2}=a_{n}^{-1}a_{n}^{2}=a_{n},\label{eq:temp-591-1-1}
\end{equation}
where the fourth equality is because$B_{\tau(x,0,a(n))}^{x}\in\partial D_{0,a(n)}$
by the definition of the first exit time $\tau_{x,0,a(n)}$. Next
let $s_{1},\cdots,s_{k}$ be an arbitrary sequence in $[0,\infty),$
and let $f\in C_{ub}(R^{k})$ be arbitrary. Then
\[
Ef(\overline{B}_{s(1)},\cdots,\overline{B}_{s(1)})\equiv E(Ef(\overline{B}_{s(1)},\cdots,\overline{B}_{s(k)})|L^{(\tau(x,0,a(n))}))
\]
\[
=E(Ef(a_{n}^{-1}(B_{\tau}^{x})'(B_{\tau+s(1)}^{x}-B_{\tau}^{x}),\cdots,a_{n}^{-1}(B_{\tau}^{x})'(B_{\tau+s(1)}^{x}-B_{\tau}^{x}))|B_{\tau}^{x}))
\]
\[
=E(E(f(a_{n}^{-1}(B_{0}^{Y(0)})'\widetilde{B}_{s(1)}^{Y(0)}-a_{n},\cdots,a_{n}^{-1}(B_{0}^{Y(0)})'\widetilde{B}_{s(k)}^{Y(0)}-a_{n})|Y_{0}))
\]
\[
=E(E(f(U'\widetilde{B}_{s(1)}^{Y(0)}-a_{n},\cdots,U'\widetilde{B}_{s(k)}^{Y(0)}-a_{n})|U))
\]
\[
=E(E(f(U'\widetilde{B}_{s(1)},\cdots,U'\widetilde{B}_{s(k)})|U)).
\]
\[
=Ef(U'\widetilde{B}_{s(1)},\cdots,U'\widetilde{B}_{s(k)})
\]
\[
=\int_{u\in\partial D(0,1)}E_{U}(du)Ef(u'\widetilde{B}_{s(1)},\cdots,u'\widetilde{B}_{s(k)})
\]
where the last equality is by Fubini's Theorem. At the same time,
for each unit vector $u\in R^{m}$, we have, according to Assertion
4 of Corollary \ref{Cor. Basic-properties-of Brownian Motion in R^m},
\[
Ef(u'\widetilde{B}_{s(1)},\cdots,u'\widetilde{B}_{s(k)})=Ef(\widehat{B}_{s(1)},\cdots,\widehat{B}_{s(k)}),
\]
where $\widehat{B}$ is some arbitrarily fixed Brownian motion in
$R^{1}$. Hence
\[
Ef(\overline{B}_{s(1)},\cdots,\overline{B}_{s(1)})=\int_{u\in\partial D(0,1)}E_{U}(du)Ef(\widehat{B}_{s(1)},\cdots,\widehat{B}_{s(k)})=Ef(\widehat{B}_{s(1)},\cdots,\widehat{B}_{s(k)}).
\]
Thus the process $\overline{B}$ has the same marginal f.j.d.'s as
the Brownian motion $\widehat{B}$. We conclude that $\overline{B}$
is a Brownian motion, as alleged.

6. Since $c_{n}\in(a_{n},\infty)G_{c}$, it follows that the first
level-crossing time $\overline{\tau}_{a(n),c(n)}$ relative to the
Brownian motion $\overline{B}^{a(n)}$ exists, as remarked in Step
2. Moreover, equalities \ref{eq:temp-251} and \ref{eq:temp-522}
together imply that
\[
P(\overline{\tau}_{a(n),c(n)}<t_{n})=2(1-\Phi_{0,t(n)}(c_{n}-a_{n}))\geq2-(1+2^{-n/2})=1-2^{-n/2}.
\]
Recall from above the r.v.'s $Y_{0}\equiv B_{\tau(x,0,a(n))}^{x}$
and $U\equiv a_{n}^{-1}Y$. Then 
\begin{equation}
U'Y=a_{n}^{-1}Y'Y=a_{n}^{-1}|Y|^{2}=a_{n}^{-1}|B_{\tau(x,0,a(n))}^{x}|^{2}=a_{n}^{-1}a_{n}^{2}=a_{n}\label{eq:temp-591}
\end{equation}
Now consider each $\omega\in(\overline{\tau}_{a(n),c(n)}<t_{n})$.
Then 
\[
s\equiv\overline{\tau}_{a(n),c(n)}(\omega)<t_{n}\equiv\varepsilon_{n}\equiv2^{-n}.
\]
Hence $\overline{B}_{s}^{a(n)}(\omega)=c_{n}$. In other words,
\[
a_{n}+U'(\omega)(B_{s+\tau(x,0,a(n))}^{x}(\omega)-B_{\tau(x,0,a(n))}^{x}(\omega))=c_{n},
\]
which, in view of equality \ref{eq:temp-591}, simplifies to
\[
U'(\omega)B_{s+\tau(x,0,a(n))}^{x}(\omega)=c_{n}.
\]
Since $U'(\omega)$ is a unit vector, it follows that 
\[
|B_{s+\tau(x,0,a(n))}^{x}(\omega)|\geq|U'(\omega)B_{s+\tau(x,0,a(n))}^{x}(\omega)|=c_{n}.
\]
Consequently, by the defining properties of the first exit time $\tau_{x,0,c(n)}$,
we obtain
\begin{equation}
\tau_{x,0,a(n)}(\omega)<\tau_{x,0,c(n)}(\omega)<2^{-n}+\tau_{x,0,a(n)}(\omega),\label{eq:temp-583}
\end{equation}
where $\omega\in(\overline{\tau}_{a(n),c(n)}<t_{n})$ is arbitrary,
where $P(\overline{\tau}_{a(n),c(n)}<t_{n})^{c}<2^{-n/2}$. It follows
that that $\tau_{x,0,a(n)}\uparrow\widehat{\tau}$ a.u. and $\tau_{x,0,c(n)}\downarrow\widehat{\tau}$
a.u. for some r.r.v. $\widehat{\tau}$.

6. We will prove that $\widehat{\tau}$ satisfies the conditions in
Definition \ref{Def:first exit time of open spheres by Brownian motions in R^m-1}
for the\emph{ }first exit time $\tau_{x,0,b}$ to exist and be equal
to $\widehat{\tau}$.

To that end, Let $t\in(0,\infty)$ be an arbitrary regular point of
the r.r.v. $\widehat{\tau}$. Let $(t_{j})_{j=0,1,\cdots}$ be a decreasing
sequence in $(t,\infty)$ such that $t_{j}\downarrow t$ and such
that $t_{j}$ is a regular point of the stopping times $\tau_{x,0,a(n)},\tau_{x,0,c(n)}$,
for each $n,j\geq0.$ Then
\[
(\widehat{\tau}<t_{j})=\bigcup_{n=0}^{\infty}\bigcap_{k=n}^{\infty}(\tau_{x,0,c(n)}\leq t_{j})\in L^{(t(j))}
\]
for each $j\geq0$. Consequently
\[
(\widehat{\tau}\leq t)=\bigcap_{j=0}^{\infty}(\widehat{\tau}<t_{j})\in L^{(t+)}=L^{(t)},
\]
where the last equality is due to the right continuity of the filtration
$\mathcal{L}$. Thus we see that the r.r.v. $\widehat{\tau}$ is a
stopping time relative to the filtration $\mathcal{L}$. 

7. Because \emph{$\tau_{x,0,c(n)}\downarrow\widehat{\tau}$ and} \emph{$\tau_{x,0,a(n)}\uparrow\widehat{\tau}$}
a.u., as proved in Step 5, and because the process $B^{x}$ is a.u.
continuous, we see that $B_{\widehat{\tau}}^{x}$ is a well defined
r.v. and that $B_{\tau(x,0,c(k))}^{x}\rightarrow B_{\widehat{\tau}}^{x}$
a.u. as $k\rightarrow\infty$. Since $B_{\tau(x,0,c(k))}^{x}\in D_{0,c(k)}$
and $c_{k}\downarrow b$, the distance $d(B_{\tau(x,0,c(k))}^{x},\partial D_{0,b})\rightarrow0$
a.u. as $k\rightarrow\infty$. Consequently, $B_{\widehat{\tau}}^{x}$
has values in $\partial D_{0,b}$. Now consider each $\omega\in domain(B_{\widehat{\tau}}^{x})$
and consider each $t\in[0,\widehat{\tau}(\omega))$. Then $t\in[0,\tau_{x,0,a(n)}(\omega))$
for some sufficiently large $k\geq0$. Hence $B_{t}^{x}\in D_{0,a(k)}\subset D_{0,b}$.

Thus we have verified the conditions in Definition \ref{Def:first exit time of open spheres by Brownian motions in R^m-1}
for the\emph{ }first exit time $\tau_{x,0,b}$ to exist and be equal
to $\widehat{\tau}$. Since $b\in(|x|,\infty)$ is arbitrary, Assertion
1 of the theorem is proved.

8. Proceed to prove Assertion 2. To that end, let $r>0$ be arbitrary.
By the just-established Assertion 1 of this theorem, applied to the
case where $x=0$ and $b=r>0$, we see that the first exit time $\tau_{0,0,r}$
exists. Therefore Lemma \ref{Lem. Translational invariance of certain exit times}
says that the stopping time $\tau_{x,x,r}$ exists and is equal to
$\tau_{0,0,r}$. Assertion 2 of this theorem is also proved.

9. It remains to prove Assertion 3. To that end, let $k\geq1$ and
$b'>b''$ be arbitrary, with $b'-b''<2^{-k}.$ Then there exist $a''<b''$and
$c'>b'$ such that $a'',c'\in H_{c}G_{c}$ and such that $c'-a''<2^{-k}$.
Repeating the above Steps 3-6, with $a'',c'$ in the place of $a_{n},c_{n}$
respectively, we obtain, in analogy to inequality \ref{eq:temp-583},
some measurable set $A$ with $P(A^{c})<2^{-k/2}$ such that, on $A$,
we have
\[
\tau_{x,0,a''}<\tau_{x,0,c'}<2^{-k}+\tau_{x,0,a''},
\]
whence
\[
0<\tau_{x,0,c'}-\tau_{x,0,a''}<2^{-k}.
\]
Therefore
\[
\tau_{x,0,b'}-\tau_{x,0,b''}<\tau_{x,0,c'}-\tau_{x,0,a''}<2^{-k}
\]
on the measurable set $A$ with $P(A^{c})<2^{-k/2}$. Assertion 3
and the theorem are proved.
\end{proof}
.

.

\part{Appendix}

For completeness and for ease of reference, we include in this Appendix
some theorems in the subject of Several Real Variables, especially
the Theorem for the Change of Integration Variables, which we use
many times in this book. For the proofs of these theorems, we assume
no knowledge of axiomatic Euclidean Geometry. For exampe, in the proof
of the Theorem for the change of Integration Vairables, we do not
assume any prior knowledge that a rotation of the plane preserves
the area of a trianlge, because we consider this a speical case of
the thoerem we are proving. 

The reader who is familiar with the subject as presented in Chapters
9 and 10 of \cite{Rudin13}, can safely skip this Appendix. 

\chapter{The Inverse Function Theorem}

Before proceeding, recall some notations and facts about vectors and
matrices. 

Let $x\in R^{n}$ be arbitrary. Unless otherwise specified, $x_{i}$
will denote the $i$-th component of $x$. Thus $x=(x_{1},\cdots,x_{n})$.
Similarly, if $f:A\rightarrow R^{n}$ is a function, then, unless
otherwise specified, $f_{i}:A\rightarrow R$ will denote the $i$-th
component of $f$. Thus $f_{i}(u)\equiv(f(u))_{i}$ for each $i=1,\cdots,n$
and for each $u\in A$. Let $G$ be an $m\times n$ matrix with entries
in $R$, and let $x=(x_{1},\cdots,x_{n})\in R^{n}$. Then we will
regard $x$ as a column vector, i.e. an $n\times1$ matrix, and write
$Gx$ for the matrix product of $G$ and $x$. Thus $Gx$ is a column
vector in $R^{m}.$ Consider the Euclidean space $R^{n}$. Let $u\in R^{n}$
be arbitrary. We will write $u_{k}$ for the $k-th$ coordinate of
$u$, and define $\left\Vert u\right\Vert \equiv\sqrt{u_{1}^{2}+\cdots+u_{n}^{2}}$.
Then $\frac{1}{n}\sum_{j=1}^{n}|v_{j}|\leq(\frac{1}{n}\sum_{j=1}^{n}v_{j}^{2})^{\frac{1}{2}}=\frac{1}{\sqrt{n}}\left\Vert v\right\Vert $
by Lyapunov's inequality. Suppose $G=[G_{i,j}]$ is an $n\times n$
matrix. The determinant of $G$ is defined as
\[
\det G\equiv\sum_{\pi\in\Pi}\mbox{sign}(\pi)G_{1,\pi(1)}G_{2,\pi(2)}\cdots G_{n,\pi(n)},
\]
where $\Pi$ is the set of all permutations on $\{1,\cdots,n\}$ and
where sign$(\pi)$ is $+1$ or $-1$ according as $\pi$ is an even
or odd permutation. If $n>1$, for each $i,j=1,\cdots,n$ let $G'_{i,j}$
be the $(n-1)\times(n-1)$ matrix obtained by deleting the $i$-th
row and $j$-th column from $G$. If $n=1$, define $G'_{1,1}\equiv[1]$
so that $\mbox{det}G'_{1,1}=1$. For later reference, the next lemma
collects some convenient bounds for matrices, and lists some basic
facts from Matrix Algebra without giving the proofs.
\begin{lem}
\label{Lem. Matrix Basics} \textbf{\emph{(Matrix Basics). }}Suppose
the $m\times n$ matrices $G=[G_{i,j}]$ and $\bar{G}=[\bar{G}_{i,j}]$
are such that $|G_{i,j}|\leq b$ and $|\bar{G}_{i,j}|\leq b$ for
each $i=1,\cdots,m$, and for each $j=1,\cdots,n$, where $b>0$.
Then the following holds.
\end{lem}
\begin{enumerate}
\item \emph{$\left\Vert Gw\right\Vert \leq\sqrt{mn}b\left\Vert w\right\Vert $
for each $w\in R^{n}$.}
\item \emph{Suppose $m=n$. Then $|\det G|\leq n!b^{n}$. If, in addition,
$\mbox{|det}G|\geq c$ where $c>0$, then, }for each $i=1,\cdots,n$,\emph{
there exists $j=1,\cdots,n$ with $|G_{i,j}|\geq(n!b^{n-1})^{-1}c$.}
\item \emph{Suppose $m=n$ and $\mbox{|det}G|\geq c$ where $c>0$. Then
$G$ has an inverse $G^{-1}\equiv F\equiv[F_{i,j}]$ where $F_{i,j}\equiv(\mbox{det}G)^{-1}(-1)^{i+j}\mbox{det}G'_{j,i}$.
Further}mo\emph{re, for each $w\in R^{n}$, we have 
\begin{equation}
\left\Vert G^{-1}w\right\Vert \leq\beta\left\Vert w\right\Vert .\label{eq:temp-32}
\end{equation}
where $\beta\equiv n!c^{-1}b^{n-1}$. }
\item \emph{Suppose $m=n$. Let }$i=1,\cdots,n$\emph{ be arbitrary. Then
\[
\det G=\sum_{j=1}^{n}(-1)^{i+j}G_{i,j}\det G'_{i,j}=\sum_{j=1}^{n}(-1)^{i+j}G_{j,i}\det G'_{j,i}.
\]
In particular, if there exists }$k=1,\cdots,n$ \emph{such that} \emph{$G_{i,j}=0$
for each }$j=1,\cdots,n$ \emph{with $j\neq k$, then $\det G=(-1)^{i+k}G_{i,k}\det G'_{i,k}$.
S}im\emph{ilarly}, \emph{if there exists }$k=1,\cdots,n$ \emph{such
that} \emph{$G'_{j,i}=0$ for each }$j=1,\cdots,n$\emph{ with $j\neq k$,
then }$\det G=(-1)^{i+k}G_{k,i}\det G'_{k,i}$. 
\item Suppose \emph{$m=n$ and $G$ and $F$ are $n\times n$ }ma\emph{trices.
Then $\det(GF)=\det(G)\det(F)$.}
\end{enumerate}
\begin{proof}
$\,$

1. For each $w\in R^{n}$, we have
\[
\left\Vert Gw\right\Vert \equiv\sqrt{\sum_{i=1}^{m}(\sum_{j=1}^{n}G_{i,j}w_{j})^{2}}\leq b\sqrt{\sum_{i=1}^{m}(\sum_{j=1}^{n}|w_{j}|)^{2}}
\]
\begin{equation}
\leq b\sqrt{\sum_{i=1}^{m}(\sqrt{n}\left\Vert w\right\Vert )^{2}}=\sqrt{mn}b\left\Vert w\right\Vert .\label{eq:temp-30}
\end{equation}

2. Let $\Pi$ be the set of all permutations on $\{1,\cdots,n\}$.
Since the set $\Pi$ has $n!$ elements, we obtain, from the definition
of the determinant,

\[
|\det G|\leq\sum_{\pi\in\Pi}|G_{1,\pi(1)}G_{2,\pi(2)}\cdots G_{n,\pi(n)}|\leq n!b^{n}.
\]
Suppose, in addition, that $|\det G|\geq c$. Then there exists $\pi\in\Pi$
such that
\[
|G_{1,\pi(1)}G_{2,\pi(2)}\cdots G_{n,\pi(n)}|\geq(n!)^{-1}c.
\]
Hence $|G_{i,\pi(i)}|\geq(n!b^{n-1})^{-1}c$ for each $i=1,\cdots,n$.

3. Cramer's Rule, whose straightforward proof using the definition
of the determinant is omitted, says that $F$ is the inverse of $G$.
If $n=1$, then $|F_{1,1}|=|(\mbox{det}G)^{-1}|\leq c^{-1}$. If $n>1,$
then $\mbox{|det}G'_{j,i}|\leq(n-1)!b^{n-1}$ for each $i,j=1,\cdots,n$,
according to Assertion 2, and so
\[
|F_{i,j}|=|(\mbox{det}G)^{-1}\mbox{det}G'_{j,i}|\leq c^{-1}(n-1)!b^{n-1}.
\]
Hence, by Assertion 1, $\left\Vert G^{-1}w\right\Vert \equiv\left\Vert Fw\right\Vert \leq n(n-1)!c^{-1}b^{n-1}\left\Vert w\right\Vert $
for each $w\in R^{n}$ if $n>1$. Combining, $\left\Vert G^{-1}w\right\Vert \leq\beta\left\Vert w\right\Vert $
for each $w\in R^{n}$ and for each $n\geq1$, where $\beta\equiv n!c^{-1}b^{n-1}$.

4,5. Proof omitted.
\end{proof}
Let $d$ denote the Euclidean metric defined by $d(u,v)\equiv\left\Vert u-v\right\Vert $
for all $u,v\in R^{n}$. Let $u\in R^{n}$, and let $r\geq0$. For
each compact subset $K$ of $R^{n}$ let $d(u,K)$ denote the distance
from $u$ to $K$. Define $B(u,r)\equiv\{v:\left\Vert v-u\right\Vert \leq r\}$,
$B^{\circ}(u,r)\equiv\{v:\left\Vert v-u\right\Vert <r\}$, and $\partial B(u,r)\equiv\{v:\left\Vert v-u\right\Vert =r\}$.
Suppose $K$ is a compact subset of $R^{n}$. Since the function $d(\cdot,K)$
is continuous on $R^{n}$, the closed $r$-neighborhood $(d(\cdot,K)\leq r)\equiv\{u\in R^{n}:d(u,K)\leq r\}$
is compact for all but countably many $r>0$. 

A compact subset $K$ is said to be \emph{well contained} \index{well contained compact subset}
in a subset $A$ of $R^{n}$ if $K_{r}\subset A$ for some $r>0$.
In that case, we write $K\Subset A$. More generally, a subset $B$
is said to be \emph{well contained}\index{well contained subset}
in $A$ if $B\subset K\Subset A$ for some compact subset $K$. 

Suppose $K\Subset A$. Let $r>0$ be such that $K_{r}\subset A$.
Let $t\in(0,r)$ be arbitrary such that $K_{t}$ is compact. Let $s\in(0,r-t)$
be arbitrary such that $(K_{t})_{s}$ is compact. Then $(K_{t})_{s}\subset K_{r}\subset A$.
Hence $K_{t}\Subset A$. In words, if a compact $K$ is well contained
in the set $A\subset R^{n}$, then some compact neighborhood of $K$
is well contained in $A$. 
\begin{defn}
\label{Def. Differentiable Function} \textbf{(Derivative and Jacobian).
}Let $n,m\geq1$ and let $M$ be the linear space of $m\times n$
matrices with real number components. Let $A$ be any open subset
of $R^{n}$. A function $g:A\rightarrow R^{m}$ is said to be \index{differentiable function}d\emph{ifferentiable}
on $A$ if there exists a function $G:A\rightarrow M$ such that,
for each compact subset $K$ with $K\Subset A$, there exists an operation
$\delta_{K}:(0,\infty)\rightarrow(0,\infty)$ such that $\left\Vert g(v)-g(u)-G(u)(v-u)\right\Vert \leq\varepsilon\left\Vert v-u\right\Vert $
for each $u,v\in K$ with $\left\Vert u-v\right\Vert \leq\delta_{K}(\varepsilon)$.
Here $v-u$ is regarded as a column vector with $n$ rows, and $G(u)(v-u)\in R^{m}$
is the matrix-vector product. The matrix valued function $G$ on $A$
is then called a \index{derivative}\emph{derivative} of $g$ on $A$,
and, for each compact subset $K$ with $K\Subset A$, the operation
$\delta_{K}$ is called the mo\emph{dulus of} \emph{differentiability}\index{modulus of differentiability}
of $g$ on $K$. If $\delta:(0,\infty)\rightarrow(0,\infty)$ is an
operation such that $\delta$ is a modulus of differentiability of
$g$ on each compact subset well contained in $A$, then we say that
$\delta$ is a modulus of differentiability of $g$ on the open set
$A$. To emphasize that $\delta_{K}$ is independent of $u\in K$,
we sometimes say that $g$ is uniform\emph{ly differentiable} on $K$.
For each $u\in A$ and each $i=1,\cdots,m$ and $j=1,\cdots,n$, the
component of $G(u)$ at the $i$-th row and $j$-th column is called
the first order \emph{partial derivative}\index{partial derivative}
of $g_{i}$ relative to the $j$-th component variable, and is denoted
by $G_{i,j}(u)\equiv G(u)_{i,j}\equiv\frac{\partial g_{i}}{\partial v_{j}}(u)$.
In the case where $n=m$, for each $u\in A$, the determinant $\det G(u)\equiv\det[\frac{\partial g_{i}}{\partial v_{j}}(u)]$
is called the Jacobian of $g$ at $u$. $\square$
\end{defn}
Note that we write $G_{i,j}(u)$, $G(u)_{i,j}$, and $\frac{\partial g_{i}}{\partial v_{j}}(u)$
interchangeably. In the last expression, $v_{j}$ is a dummy variable.
For example, the expressions $\frac{\partial g_{i}}{\partial v_{j}}(u)$,
$\frac{\partial g_{i}}{\partial w_{j}}(u)$ or $\frac{\partial g_{i}}{\partial u_{j}}(u)$,
with different dummy variables, all have the same value as $G(u)_{i,j}$.
\begin{prop}
\textbf{\emph{\label{Prop.:Uniqueness of derivative}(Uniqueness of
derivative).}} Suppose $g:A\rightarrow R^{m}$ is differentiable on
the open subset $A$ of $R^{n}$. Then the derivative of g on $A$
is unique.
\end{prop}
\begin{proof}
Suppose $G$ and $H$ are both derivatives of $g$ on $A$. Consider
any $u\in A$. Let $r>0$ be such that $B(u,2r)\subset A$ and let
$K\equiv B(u,r)$. Then $K_{r}\subset A$ and so $K\Subset A$. Let
$\varepsilon>0$ be arbitrary. There exists $\delta_{0}\in(0,r)$
such that $\left\Vert g(v)-g(u)-G(u)(v-u)\right\Vert \leq\frac{\varepsilon}{2}\left\Vert v-u\right\Vert $
and $\left\Vert g(v)-g(u)-H(u)(v-u)\right\Vert \leq\frac{\varepsilon}{2}\left\Vert v-u\right\Vert $
for each $v\in B(u,\delta_{0})$. Write $Q\equiv G(u)-H(u)$. Then
we have 
\begin{equation}
\left\Vert Q(v-u)\right\Vert \leq\varepsilon\left\Vert v-u\right\Vert \label{eq:temp-73}
\end{equation}
 for each $v\in B(u,\delta_{0})$. Consider any $j=1,\cdots,n$. Let
\[
v\equiv(u_{1},\cdots,u_{j-1},u_{j}+\delta_{0},u_{j+1},\cdots,u_{n}).
\]
Then $v-u=(0,\cdots,0,\delta_{0},0,\cdots,0)$ has entries $0$ except
that the $j$-th entry is $\delta_{0}$, whence $Q(v-u)=\delta_{0}(Q_{1,j},\cdots,Q_{m,j})$.
Therefore inequality \ref{eq:temp-73} yields 
\[
\delta_{0}\left\Vert (Q_{1,j},\cdots,Q_{m,j})\right\Vert \leq\varepsilon\delta_{0}.
\]
Canceling, we see that $\left\Vert (Q_{1,j},\cdots,Q_{m,j})\right\Vert \leq\varepsilon$
for arbitrary $\varepsilon>0$. Consequently 
\[
\left\Vert (Q_{1,j},\cdots,Q_{m,j})\right\Vert =0
\]
and so $Q_{i,j}=0$ for each $i=1,\cdots,m$. Equivalently $G(u)=H(u)$
for any $u\in A$.
\end{proof}
\begin{prop}
\label{Prop.g  differenitable =00003D> G unif cont on well contained compacts}
\textbf{\emph{(Differentiability and uniform continuity). }}Suppose
$g:A\rightarrow R^{m}$ is differentiable on the open subset $A$
of $R^{n}$. Then the derivative $G$ of $g$ is uniformly continuous
on every compact subset $K$ with $K\Subset A$. More specifically,
if $K_{t}\Subset A$ for some $t>0$ and if $\delta$ is a modulus
of differentiability of $g$ on $K_{t}$, then the operation  $\delta_{1}:(0,\infty)\rightarrow(0,\infty)$
defined by $\delta_{1}(\varepsilon)\equiv\frac{1}{2}\delta(\frac{\varepsilon}{4})\wedge t$
for each $\varepsilon>0$ is a modulus of continuity of $G_{i,j}$
on $K$ for each $i=1,\cdots,m$ and each $j=1,\cdots,n$.
\end{prop}
\begin{proof}
Let $K$ be a compact subset of $A$ with $K\Subset A$. Let $t>0$
be such that $K_{t}\Subset A$. Let $\delta$ be a modulus of differentiability
of $g$ on $K_{t}$. Let $u\in K$ be arbitrary but fixed. Define
a function $f:A\rightarrow R^{m}$ by $f(v)\equiv g(v)-g(u)-G(u)(v-u)$
for each $v\in A$. Define $F(v)\equiv G(v)-G(u)$ for each $v\in A$.
Then, for arbitrary $\varepsilon>0$ and for each $w,v\in K_{t}$
with $\left\Vert w-v\right\Vert \leq\delta(\varepsilon)$, we have
\[
\left\Vert f(w)-f(v)-F(v)(w-v)\right\Vert \equiv\left\Vert f(w)-f(v)-(G(v)-G(u))(w-v)\right\Vert 
\]
\[
\equiv\left\Vert g(w)-g(u)-G(u)(w-u)-g(v)+g(u)+G(u)(v-u)-(G(v)-G(u))(w-v)\right\Vert 
\]
\begin{equation}
=\left\Vert g(w)-g(v)-G(v)(w-v)\right\Vert \leq\varepsilon\left\Vert w-v\right\Vert .\label{eq:temp-66}
\end{equation}
Clearly $f(u)=0$, and $F(u)=0$, the $m\times n$ matrix whose entries
are zeros. By setting $v=u$ in inequality \ref{eq:temp-66} we have
\begin{equation}
\left\Vert F(w)\right\Vert \leq\varepsilon\left\Vert w-u\right\Vert \label{eq:temp-67}
\end{equation}
for each $w\in K_{t}$ with $\left\Vert w-u\right\Vert \leq\delta(\varepsilon)$
and for each $\varepsilon>0$. 

Now let $\varepsilon>0$ be arbitrary and let $a\equiv\frac{1}{2}\delta(\varepsilon)\wedge t$.
Let $v\in K$ be arbitrary such that $\left\Vert v-u\right\Vert \leq a$.
Fix any $i=1,\cdots,m$ and $j=1,\cdots,n$. Define $w\equiv(v_{1},\cdots,v_{j-1},v_{j}+a,v_{j+1},\cdots,v_{n})$.
Note that $w-v\equiv(0,\cdots,0,a,0,\cdots,0)$ has all components
equal to $0$ except the $j$-th. Hence $\left\Vert w-v\right\Vert =a$
and so $w\in K_{t}$. Moreover$\left\Vert w-u\right\Vert \leq\left\Vert w-v\right\Vert +\left\Vert v-u\right\Vert \leq2a\leq\delta(\varepsilon)$
and so $\left\Vert f(w)\right\Vert \leq\varepsilon\left\Vert w-u\right\Vert $
by inequality \ref{eq:temp-67}. Therefore
\[
|F(v)_{i,j}a|=|\sum_{k=1}^{n}F(v)_{i.k}(w-v)_{k}|=|(F(v)(w-v))_{i}|
\]
 
\[
\leq\left\Vert F(v)(w-v)\right\Vert \leq\left\Vert f(w)-f(v)-F(v)(w-v)\right\Vert +\left\Vert f(w)\right\Vert +\left\Vert f(v)\right\Vert 
\]
\[
\leq\varepsilon\left\Vert v-w\right\Vert +\varepsilon\left\Vert w-u\right\Vert +\varepsilon\left\Vert v-u\right\Vert 
\]
\[
\leq\varepsilon a+2\varepsilon a+\varepsilon a=4\varepsilon a.
\]
Canceling $a$ on both ends of the inequality, we obtain $|F(v)_{i,j}|\leq4\varepsilon$.
Equivalently, $|G(v)_{i,j}-G(u)_{i,j}|\leq4\varepsilon$. The proposition
is proved.
\end{proof}
\begin{prop}
\label{Prop. g differentiable on open set=00003D> g Lipsitz on compact subsets}
\textbf{\emph{(Differentiability and Lipschitz continuity). }}Suppose
$g:A\rightarrow R^{m}$ is a differentiable function on the open subset
$A$ of $R^{n}$. Let $K$ be any compact subset with $K\Subset A$.
Then $g$ is Lipschitz continuous on $K$. 
\end{prop}
\begin{proof}
Let $K$ be a compact subset of $A$ with $K\Subset A$. Let $t>0$
be such that $K_{t}\Subset A$. Let $\delta$ be a modulus of differentiability
of $g$ on $K_{t}$. By Proposition Let $a\equiv\frac{1}{3}\delta(1)$.
Let $u_{1},\cdots,u_{p}\in K$ be an $a$-approximation of $K$. Consider
any $u,v\in K$. Let $k=1,\cdots,p$ be such that $\left\Vert u-u_{k}\right\Vert <a$.
Either (i) $\left\Vert v-u_{k}\right\Vert <3a$ or (ii) $\left\Vert v-u_{k}\right\Vert >2a$.
In case (i), we have $\left\Vert v-u\right\Vert <3a\leq\delta(1)$
whence $\left\Vert g(u)-g(v)-G(u)(u-v)\right\Vert \leq\left\Vert u-v\right\Vert $.
According to Proposition \ref{Prop.g  differenitable =00003D> G unif cont on well contained compacts},
$|G_{i,j}|$ is uniformly continuous, and so bounded by some $b>0$
on $K$, for each $i=1,\cdots,m$ and for each $j=1,\cdots,n$. Therefore,
by Lemma \ref{Lem. Matrix Basics}, there exists $c>0$ such that
$\left\Vert G(u)z\right\Vert \leq c\left\Vert z\right\Vert $ for
each $z\in R^{n}$. Hence $\left\Vert g(u)-g(v)\right\Vert \leq\left\Vert G(u)(u-v)\right\Vert +\left\Vert u-v\right\Vert \leq(1+c)\left\Vert u-v\right\Vert $.
Setting $v=u_{k}$ we conclude that $\left\Vert g(u)-g(u_{k})\right\Vert \leq(1+c)a$.
Thus $\left\Vert g(u)\right\Vert \leq b'\equiv\vee_{k=1}^{p}\left\Vert g(u_{k})\right\Vert +(1+c)a$.
The same consideration yields $\left\Vert g(v)\right\Vert \leq b'$
for arbitrary $v\in K$. Next consider case (ii), where $\left\Vert v-u_{k}\right\Vert >2a$.
Then $\left\Vert v-u\right\Vert >a$. Hence $\left\Vert g(v)-g(u)\right\Vert \leq2b'\leq\frac{2b'}{a}\left\Vert v-u\right\Vert $.
Combining, we see that $\left\Vert g(v)-g(u)\right\Vert \leq b''\left\Vert v-u\right\Vert $
for each $u,v\in K$, where $b''\equiv(1+c)\vee\frac{2b'}{a}$. In
other words, $b''$ is a Lipschitz constant of $g$ on $K$.
\end{proof}
\begin{prop}
\textbf{\emph{\label{Prop. Chain Rule} (Chain Rule)}} Let $f:A\rightarrow R^{m}$
and $g:B\rightarrow R^{n}$ be differentiable functions on the open
subsets $A$ and B of $R^{p}$ and $R^{m}$ respectively, such that
for each compact subset $K$ of $A$ with $K\Subset A$ we have $f(K)\subset K'\Subset B$
for some compact subset $K'$ of $B$. Let $F$ and $G$ denote the
derivatives of $f$ and $g$ respectively. Then the composite function
$g(f):A\rightarrow R^{n}$ is differentiable on $A$, with derivative
$G(f)F$. 
\end{prop}
\begin{proof}
Let $K$ be any compact subset of $A$ with $K\Subset A$. By hypothesis
there exists a compact subset $K'$ of $B$ with $f(K)\subset K'\Subset B$.
Let $\delta$ and $\delta'$ denote the moduli of differentiability
of $f$ and $g$ on $K$ and $K'$ respectively. By Proposition \ref{Prop. g differentiable on open set=00003D> g Lipsitz on compact subsets},
$f$ is Lipschitz continuous on $K$, with a Lipschitz constant $c>0$.
By Proposition \ref{Prop.g  differenitable =00003D> G unif cont on well contained compacts},
the partial derivatives $G_{i,j}$ are uniformly continuous, hence
bounded in absolute value, on $K'$ for each $i=1,\cdots,n$ and each
$j=1,\cdots,m$. Let $b>0$ be such that $|G_{i,j}|\leq b$ on $K'$
for each $i=1,\cdots,n$ and each $j=1,\cdots,m$. 

Let $\varepsilon>0$ be arbitrary. Let $x,y\in K$ be arbitrary with
$\left\Vert y-x\right\Vert <\delta_{1}(\varepsilon)\equiv\delta(\varepsilon)\wedge c^{-1}\delta'(\varepsilon)$.
Let $u\equiv f(x)$ and $v\equiv f(y)$. Then $u,v\in K'$, and $\left\Vert v-u\right\Vert \equiv\left\Vert f(y)-f(x)\right\Vert \leq c\left\Vert y-x\right\Vert \leq\delta'(\varepsilon)$,
whence
\[
\left\Vert g(v)-g(u)-G(u)(v-u)\right\Vert \leq\varepsilon\left\Vert v-u\right\Vert \leq c\varepsilon\left\Vert y-x\right\Vert .
\]
At the same time, 
\[
\left\Vert f(y)-f(x)-F(x)(y-x)\right\Vert \leq\varepsilon\left\Vert y-x\right\Vert .
\]
Hence, according to Lemma \ref{Lem. Matrix Basics},
\[
\left\Vert G(u)(f(y)-f(x)-F(x)(y-x))\right\Vert 
\]
\[
\leq\sqrt{mn}b\left\Vert f(y)-f(x)-F(x)(y-x)\right\Vert \leq\sqrt{mn}b\varepsilon\left\Vert y-x\right\Vert .
\]
Combining, we obtain
\[
\left\Vert g(f(y))-g(f(x))-G(f(x))F(x)(y-x)\right\Vert 
\]
\[
=\left\Vert g(f(y))-g(f(x))-G(u)F(x)(y-x)\right\Vert 
\]
\[
\leq\left\Vert g(f(y))-g(f(x))-G(u)(f(y)-f(x))\right\Vert +\sqrt{mn}b\varepsilon\left\Vert y-x\right\Vert 
\]
\[
\equiv\left\Vert g(v)-g(u)-G(u)(v-u)\right\Vert +\sqrt{mn}b\varepsilon\left\Vert y-x\right\Vert 
\]
\begin{equation}
\leq c\varepsilon\left\Vert y-x\right\Vert +\sqrt{mn}b\varepsilon\left\Vert y-x\right\Vert .\label{eq:temp-74}
\end{equation}
Since $\varepsilon>0$ is arbitrary with the operation $\delta_{2}:(0,\infty)\rightarrow(0,\infty)$
defined by $\delta_{2}(\varepsilon)\equiv\delta_{1}((c+\sqrt{mn}b)^{-1}\varepsilon)$
for each $\varepsilon>0$ as a modulus of differentiability on $K$.
Since $K\Subset A$ is arbitrary, we see from inequality \ref{eq:temp-74}
that $G(f)F$ is the derivative of $g(f)$ on $A$.
\end{proof}
The proof below for the Inverse Function Theorem is by the method
of contraction mapping.
\begin{thm}
\textbf{\emph{(Inverse Function Theorem)}} \label{Thm. Inverse Function Thm}
Let $g:B\rightarrow R^{n}$ be a differentiable function on an open
subset $B$ of $R^{n}$, with derivative $G$. Define the Jacobian
$J(v)\equiv\mbox{det}G(v)$ for each $v\in B$. Suppose $K$ is a
compact subset of $B$ with $K_{r}\Subset B$ for some $r>0$, and
suppose $|J|\geq c$ on $K_{r}$ for some $c>0$. Let $\delta$ be
a modulus of differentiability of $g$ on $K_{r}$, and let $b>0$
be such that $|G_{i,j}|\leq b$ on $K_{r}$, for each $i,j=1,\cdots,n$.
Then there exists $s=s(n,r,b,c,\delta)>0$ such that for each $u\in K$
we have (i) there exists a function $f:C\equiv B^{\circ}(g(u),s)\rightarrow B(u,r)$
such that $g(f(y))=y$ for each $y\in C$, and (ii) the inverse function
$f$ is differentiable on $C$, with derivative $F\equiv G(f)^{-1}$,
and with a modulus of differentiability $\delta'=\delta'(n,r,b,c,\delta)$
on the open set $C$. Note that we write $s=s(n,r,b,c,\delta)$ to
emphasize that $s$ depends only on $n,r,b,c,$ and the operation
$\delta$, and is otherwise independent of $g$,K, or $u$. Similarly
we emphasize that $\delta'$ depends only on $n,r,b,c,$ and $\delta$.
\end{thm}
\begin{proof}
Note that by Proposition \ref{Prop.g  differenitable =00003D> G unif cont on well contained compacts},
$G_{i,j}$ is uniformly continuous on $K_{r}$ for each $i,j=1,\cdots,n$,
with a modulus of continuity $\delta_{1}$ defined by $\delta_{1}(\varepsilon)\equiv\frac{1}{2}\delta(\frac{\varepsilon}{4})\wedge r$
for each $\varepsilon>0$. By the same token, there exists $b>0$
be such that $|G_{i,j}|\leq b$ on $K_{r}$, for each $i,j=1,\cdots,n$.
Define the operation $\delta_{2}:(0,\infty)\rightarrow(0,\infty)$
by $\delta_{2}(\varepsilon)\equiv\delta_{1}(\frac{\varepsilon\wedge\varepsilon_{0}}{n})\wedge\delta(\frac{\varepsilon\wedge\varepsilon_{0}}{4n})$
for each $\varepsilon>0$. Let $\beta\equiv n!c^{-1}b^{n-1}$ and
$\varepsilon_{0}\equiv(4\beta)^{-1}$. Let $a\equiv\frac{1}{2}\delta_{2}(\varepsilon_{0})\wedge r$
and $s\equiv\frac{a}{2\beta}$. We will show that $s$ has the desired
properties. 

To that end, consider any $u\in K$. Let $x\equiv g(u)$. Consider
any $v,w\in B(u,a)$. Then $v,w\in K_{r}$. Moreover $\left\Vert w-v\right\Vert \leq2a\leq\delta_{2}(\varepsilon_{0})$. 

Let $\varepsilon>0$ be arbitrary. Suppose $\left\Vert w-v\right\Vert \leq\delta_{2}(\varepsilon)\equiv\delta_{1}(\frac{\varepsilon\wedge\varepsilon_{0}}{n})\wedge\delta(\frac{\varepsilon\wedge\varepsilon_{0}}{4n})$.
Then, since$\delta$ is a modulus of differentiability of $g$, we
have 
\begin{equation}
\left\Vert g(w)-g(v)-G(v)(w-v)\right\Vert \leq\frac{\varepsilon}{4n}\left\Vert w-v\right\Vert \leq\varepsilon\left\Vert w-v\right\Vert .\label{eq:temp-33}
\end{equation}
At the same time, for each $i,j=1,\cdots,n$, since $\delta_{1}$
is a modulus of continuity of $G_{i,j}$, we have 
\begin{equation}
|G_{i,j}(w)-G_{i,j}(v)|\leq\frac{\varepsilon}{n}.\label{eq:temp-64}
\end{equation}
Hence, by Lemma \ref{Lem. Matrix Basics}, we have 
\begin{equation}
\left\Vert (G(w)-G(v))z\right\Vert \leq\varepsilon\left\Vert z\right\Vert \label{eq:temp-61}
\end{equation}
for each $z\in R^{n}$. Moreover, since $|\mbox{det}G(v)|\geq c>0$
and $|G_{i,j}(v)|\leq b$ for $i,j=1,\cdots,n$. Therefore we see,
by Lemma \ref{Lem. Matrix Basics}, that the inverse matrix $G(v)^{-1}$
exists, with
\begin{equation}
\left\Vert G(v)^{-1}z\right\Vert \leq\beta\left\Vert z\right\Vert \label{eq:temp-69}
\end{equation}
for each $z\in R^{n}$. 

In particular, inequalities \ref{eq:temp-33} through \ref{eq:temp-69}
hold for any $v,w\in B(u,a)$ if $\varepsilon$ is replaced by $\varepsilon_{0}$.

Next, let $y\in B(g(u),s)$ be arbitrary but fixed. For each $v\in B(u,a)$
define $\Phi(v)\equiv v+G^{-1}(u)(y-g(v))$. Then $\left\Vert \Phi(u)-u\right\Vert \equiv\left\Vert G^{-1}(u)(y-g(u))\right\Vert \leq\beta s\equiv\frac{a}{2}$.
Moreover, for any $v,w\in B(u,a)$, we have 
\[
\left\Vert \Phi(w)-\Phi(v)\right\Vert =\left\Vert w-v-G^{-1}(u)(g(w)-g(v))\right\Vert 
\]
\[
=\left\Vert G^{-1}(u)(G(u)(w-v)-g(w)+g(v))\right\Vert 
\]
\[
\leq\beta\left\Vert G(u)(w-v)-g(w)+g(v)\right\Vert 
\]
\[
\leq\beta\left\Vert G(v)(w-v)-g(w)+g(v)\right\Vert +\beta\left\Vert (G(u)-G(v))(w-v)\right\Vert 
\]
\begin{equation}
\leq\beta\varepsilon_{0}\left\Vert w-v\right\Vert +\beta\varepsilon_{0}\left\Vert w-v\right\Vert =2\beta\varepsilon_{0}\left\Vert w-v\right\Vert \equiv\frac{1}{2}\left\Vert w-v\right\Vert \label{eq:temp-39}
\end{equation}
where the first inequality follows from inequality \ref{eq:temp-69}
applied to $u$, and where the third inequality follows from inequalities
\ref{eq:temp-33} and \ref{eq:temp-61} with $\varepsilon$ replaced
by $\varepsilon_{0}$. It follows that, for arbitrary $v\in B(u,a)$,
\[
\left\Vert \Phi(v)-u\right\Vert \leq\left\Vert \Phi(v)-\Phi(u)\right\Vert +\left\Vert \Phi(u)-u\right\Vert \leq\frac{1}{2}\left\Vert v-u\right\Vert +\frac{a}{2}\leq\frac{a}{2}+\frac{a}{2}=a,
\]
and so $\Phi(v)\in B(u,a)$. Thus we see that $\Phi$ is a function
mapping $B(u,a)$ into $B(u,a)$. Now define $u^{(0)}=u$, and inductively
define $u^{(k)}\equiv\Phi(u^{(k-1)})$ for each $k\geq1$. Then, according
to inequality \ref{eq:temp-39},
\[
\left\Vert u^{(k)}-u^{(k-1)}\right\Vert \equiv\left\Vert \Phi(u^{(k-1)})-\Phi(u^{(k-2)})\right\Vert 
\]
\[
\leq\frac{1}{2}\left\Vert u^{(k-1)}-u^{(k-2)}\right\Vert \leq\cdots\leq2^{-k+1}\left\Vert u^{(1)}-u^{(0)}\right\Vert 
\]
for each $k\geq1$. Therefore $u^{(k)}\rightarrow v$ for some $v\in B(u,a)$.
By the Lipschitz continuity of $\Phi$ displayed in inequality \ref{eq:temp-39},
we have $\Phi(u^{(k)})\rightarrow\Phi(v)$. Equivalently $u^{(k+1)}\rightarrow v+G^{-1}(u)(y-g(v))$.
Therefore $v=v+G^{-1}(u)(y-g(v))$ and so $G^{-1}(u)(y-g(v))=0$.
Multiplying the last equality from the left by the matrix $G(u)$,
we obtain $y-g(v)=0$, or $y=g(v)$. We will show that $v$ is unique.
Suppose $y=g(w)$ for some other $w\in B(u,a)$. Then $g(w)=g(v)$,
and so 
\[
\left\Vert G(v)(w-v)\right\Vert =\left\Vert g(w)-g(v)+G(v)(w-v)\right\Vert \leq\varepsilon_{0}\left\Vert w-v\right\Vert 
\]
according to inequality \ref{eq:temp-33} with $\varepsilon$ replaced
by $\varepsilon_{0}$. Hence
\[
\left\Vert w-v\right\Vert =\left\Vert G(v)^{-1}G(v)(w-v)\right\Vert 
\]
\[
\leq\beta\left\Vert G(v)(w-v)\right\Vert \leq\beta\varepsilon_{0}\left\Vert w-v\right\Vert \leq\frac{1}{4}\left\Vert w-v\right\Vert .
\]
Consequently $\left\Vert w-v\right\Vert =0$ and $w=v$. Summing up,
for each $y\in B(g(u),s)$, there exists a unique $v\in B(u,a)$ such
that $y=g(v)$. Hence $B(g(u),s)\subset g(B(u,a))\subset g(B(u,r))$.
We can therefore define a function $f$ on $B(g(u),s)$ by $f(y)=v$
where $v\in B(u,a)$ is such that $y=g(v)$, for each $y\in B(g(u),s)$.
By definition, $g(f(y))=y$ for each $y\in B(g(u),s)$. In other words,
$f$ is the inverse of $g$ on $B(x,s)$, with values in $B(u,a)$.

Next, we will prove the differentiability of $f$ on $C\equiv B^{\circ}(x,s)$.
Consider any $y,z\in B(x,s)$. Let $v\equiv f(y)$ and $w\equiv f(z)$.
Then $y=g(v)$ and $z=g(w)$ by the definition of the inverse function
$f$. Using inequality \ref{eq:temp-33}, we estimate 
\[
\left\Vert v-w-G(w)^{-1}(y-z)\right\Vert =\left\Vert G(w)^{-1}(G(w)(v-w)-g(v)+g(w))\right\Vert 
\]
\begin{equation}
\leq\beta\left\Vert G(w)(v-w)-g(v)+g(w)\right\Vert \leq\beta\varepsilon_{0}\left\Vert v-w\right\Vert .\label{eq:temp-68}
\end{equation}
Hence
\[
\left\Vert v-w\right\Vert \leq\left\Vert G(w)^{-1}(y-z)\right\Vert +\beta\varepsilon_{0}\left\Vert v-w\right\Vert .
\]
It follows that
\[
\left\Vert v-w\right\Vert \leq\frac{1}{1-\beta\varepsilon_{0}}\left\Vert G(w)^{-1}(y-z)\right\Vert 
\]
\begin{equation}
\leq\frac{\beta}{1-\beta\varepsilon_{0}}\left\Vert y-z\right\Vert =\frac{4}{3}\beta\left\Vert y-z\right\Vert .\label{eq:temp-70}
\end{equation}
Thus $\left\Vert f(y)-f(z)\right\Vert \leq\frac{4}{3}\beta\left\Vert y-z\right\Vert $
for each y $y,z\in B(x,s)$. In other words, $f$ is Lipschitz continuous
on $B(x,s)$, with Lipschitz constant $\frac{4}{3}\beta$. 

Now let $\varepsilon>0$ be arbitrary, and let $y,z\in B(x,s)$ be
such that $\left\Vert y-z\right\Vert \leq\delta'(\varepsilon)\equiv\frac{3}{4}\beta^{-1}\delta_{2}(\frac{3}{4}\beta^{-2}\varepsilon)$.
Then inequality \ref{eq:temp-70} implies that $\left\Vert v-w\right\Vert \leq\delta_{2}(\frac{3}{4}\beta^{-2}\varepsilon)$.
Therefore inequalities \ref{eq:temp-33} through \ref{eq:temp-69}
hold for $v,w$. Using inequality \ref{eq:temp-33}, we now obtain
\[
\left\Vert v-w-G(w)^{-1}(y-z)\right\Vert =\left\Vert G(w)^{-1}(G(w)(v-w)-g(v)+g(w))\right\Vert 
\]
\[
\leq\beta\left\Vert G(w)(v-w)-g(v)+g(w)\right\Vert \leq\frac{3}{4}\beta^{-1}\varepsilon\left\Vert v-w\right\Vert \leq\frac{3}{4}\beta^{-1}\varepsilon\frac{4}{3}\beta\left\Vert y-z\right\Vert .
\]
In other words,
\begin{equation}
\left\Vert f(y)-f(z)-G(f(z))^{-1}(y-z)\right\Vert \leq\varepsilon\left\Vert y-z\right\Vert \label{eq:temp-71}
\end{equation}
for arbitrary $y,z\in B(x,s)$ with $\left\Vert y-z\right\Vert \leq\delta'(\varepsilon)$.
In particular, inequality \ref{eq:temp-71} holds for each $y,z\in C\equiv B^{0}(x,s)$
with $\left\Vert y-z\right\Vert \leq\delta'(\varepsilon)$. We have
thus proved that $f$ is differentiable on $C$, with derivative $F=G(h)^{-1}$,
and with modulus of differentiability $\delta'$ for each compact
subset well contained in $C$. The Reader can trace the definitions
of $s$ and $\delta'$ to verify that they depend only on $n,r,b,c,$
and $\delta$. The theorem is proved.
\end{proof}
\begin{cor}
\label{Cor. Condition for Inverse Func to be differentiable} \textbf{\emph{(Condition
for inverse function to be differentiable).}} Let $A,B$ be open subsets
of $R^{n}$. Let $g:B\rightarrow A$ be a differentiable function
on $B$, with derivative $G$. Define the Jacobian $J(v)\equiv\mbox{det}G(v)$
for each $v\in B$. Suppose $|J|\geq c$ for some $c>0$ on each compact
subset $K$ with $K\Subset B$. Suppose the inverse function $f:A\rightarrow B$
of $g$ exists, such that for each compact subset $H$ of $A$ with
$H\Subset A$ we have $f(H)\Subset B$. Then the inverse function
$f$ is differentiable on $A$, with derivative $F\equiv G(f)^{-1}$.
In particular, $f$ is uniformly continuous on compact subsets well
contained in $A$.
\end{cor}
\begin{proof}
Let $H$ be an arbitrary compact subset of $A$ with $H\Subset A$.
Then $H_{a}\Subset A$ for some $a>0$. Define $K\equiv f(H)$. By
hypothesis $K_{r'}\Subset B$ for some $r'>0$. Hence $|J|\geq c$
for some $c>0$ on $K_{r'}.$ Let $\delta$ be a modulus of differentiability
of $g$ on $K_{r'}$, and let $b>0$ be such that $|G_{i,j}|\leq b$
on $K_{r'}$, for each $i,j=1,\cdots,n$. Let $\delta_{g}$ be a modulus
of continuity of $g$ on $K_{r'}$, and let $r>0$ be such that $r<r'\wedge\delta_{g}(a)$.
Then $g(K_{r})\subset H_{a}$. By the Inverse Function Theorem \ref{Thm. Inverse Function Thm},
there exists $s=s(n,r,b,c,\delta)>0$ such that for each $u\in K$
we have (i) there exists a function $\tilde{f}:C\equiv B^{\circ}(g(u),s)\rightarrow B(u,r)$
such that $g(\tilde{f}(y))=y$ for each $y\in C$, and (ii) the inverse
function $\tilde{f}$ is differentiable on $C$, with derivative $F\equiv G(\tilde{f})^{-1}$,
and with a modulus of differentiability $\delta'=\delta'(n,r,b,c,\delta)$
on the open set $C$. Consider any $u\in K$ and let $\tilde{f},C$
be as in conditions (i) and (ii). Then, for each $y\in C$, we have
$\tilde{f}(y)\in B(u,r)\subset K_{r}$. Hence $y=g(\tilde{f}(y))\in g(K_{r})\subset H_{a}$.
Moreover $f(y)=f(g(\tilde{f}(y)))=\tilde{f}(y)$ for each $y\in C$.
Consequently, condition (ii) implies that $f$ is differentiable on
$C$, with derivative $F\equiv G(f)^{-1}$, and with $\delta'$ as
a modulus of differentiability on $C$. Now let $x_{1},\cdots,x_{m}$
be an $\frac{s}{2}$-approximation of $H$. For each $i=1,\cdots,m$
define $u_{i}\equiv f(x_{i})$ and $C_{i}\equiv B^{\circ}(g(u_{i}),s)\equiv B^{\circ}(x_{i},s)$.
Let $\varepsilon>0$ be arbitrary. Let $x,y\in H$ be such that $\left\Vert x-y\right\Vert <\delta'(\varepsilon)\wedge\frac{s}{2}$.
Then $\left\Vert x-x_{i}\right\Vert <\frac{s}{2}$ for some $i=1,\cdots,m$,
whence $x,y\in C_{i}$. Since $\delta'$ is a modulus of differentiability
of $f$ on $C_{i}$, and since $\left\Vert x-y\right\Vert <\delta'(\varepsilon)\wedge\frac{s}{2}$,
we have $\left\Vert f(y)-f(x)-F(x)(y-x)\right\Vert \leq\varepsilon\left\Vert y-x\right\Vert $.
Thus the operation $\delta'(\cdot)\wedge\frac{s}{2}$ is a modulus
of differentiability of $f$ on $H$. Since the compact subset $H\Subset A$
is arbitrary, we see that $f$ is differentiable on $A$.
\end{proof}

\chapter{Change of Integration Variables}

In this section let $n\geq1$ be a fixed integer. Let $\mu_{1}$ and
$\mu$ denote the the measures with respect to the Lebesgue integrations
$\int\cdot dx$ and $\int\cdots\int\cdot dx_{1}\cdots dx_{n}$ respectively.
All measure-theoretic terms will be with respect to the Lebesgue integrations.

First an easy lemma, classically trivial, gives a sufficient condition
for a subset of $R$ to be an open interval. Note that a non-empty
set that is the intersection of two open intervals $(a,b)$ and $(a'b')$
is again an open interval, because $(a,b)\cap(a',b')=(a\vee a',b\wedge b')$.
\begin{lem}
\label{Lem. Condition for  non empty set in R to be an interval}\textbf{\emph{
(Condition of a nonempty subset of $R$ to be an interval).}} Let
$n\geq2$ be arbitrary. Let $z_{k},z'_{k},c_{k},a_{k}\in R$ for $k=1,\cdots,n$.
Suppose the set $\Gamma\equiv\{x\in R:z_{k}<c_{k}+a_{k}x<z'_{k}\mbox{ for }k=1,\cdots,n-1\}\cap(z_{n},z'_{n})$
is non-empty. Then $\Gamma$ is an open interval in $R$. 
\end{lem}
\begin{proof}
Since $\Gamma=\bigcap_{k=1}^{n-1}(\{x\in R:z_{k}<c_{k}+a_{k}x<z'_{k}\}\cap(z_{n},z'_{n}))$,
it suffices, in view of the remark preceding this lemma, to prove
the lemma for the case $n=2$. By hypothesis, there exists $y\in R$
such that $z_{1}<c_{1}+a_{1}y<z'_{1}$ and $z_{2}<y<z'_{2}$. Let
$\varepsilon>0$ be so small that $z_{1}+\varepsilon<c_{1}+a_{1}y<z'_{1}-\varepsilon$.
Then either (i) $|a_{1}y|<\frac{\varepsilon}{2}$ and $|a_{1}|(|z_{2}|\vee|z'_{2}|)<\frac{\varepsilon}{2}$,
or (ii) $|a_{1}y|>0$ or $|a_{1}|(|z_{2}|\vee|z'_{2}|)>0$. Suppose
condition (i) holds. Then, for each $x\in(z_{2},z'_{2})$, we have
\[
|c_{1}+a_{1}x-(c_{1}+a_{1}y)|\leq|a_{1}y|+|a_{1}x|<|a_{1}y|+|a_{1}|(|z_{2}|\vee|z'_{2}|)<\varepsilon
\]
whence $z_{1}<c_{1}+a_{1}x<z'_{1}$. It follows that $(z_{2},z'_{2})\subset\{x\in R:z_{1}<c_{1}+a_{1}x<z'_{1}\}$
and so $\Gamma\equiv\{x\in R:z_{1}<c_{1}+a_{1}x<z'_{1}\}\cap(z_{2},z'_{2})=(z_{2},z'_{2})$.
On the other hand, suppose condition (ii) holds. Then $|a_{1}|>0$.
Hence either $a_{1}>0$ or $a_{1}<0$. In the first case $\Gamma=(a_{1}^{-1}(z_{1}-c_{1}),a_{1}^{-1}(z'_{1}-c_{1}))\cap(z_{2},z'_{2})$.
In the second case $\Gamma=(a_{1}^{-1}(z'_{1}-c_{1}),a_{1}^{-1}(z{}_{1}-c_{1}))\cap(z_{2},z'_{2})$.
Thus we see that $\Gamma$ is an open interval under either of the
conditions (i) and (ii).
\end{proof}
In the following, an interval $\Delta$ in $R$ will mean a non-empty
subset of $R$ that is equal to one of $(a,b)$, $(a,b]$, $[a,b)$,
or $[a,b]$ for some $a,b\in R$ with $a\leq b$, and the length of
$\Delta$ is defined to be $|\Delta|\equiv b-a$. More generally,
the Cartesian product $\Delta\equiv\Delta_{1}\times\cdots\times\Delta_{n}$,
where $\Delta_{i}$ is an interval in $R$ for each $i=1,\cdots,n$,
is called an $n$-\emph{interval}\index{n-interval@\emph{n}-interval},
and the \emph{length\index{length of an n-interval}} of $\Delta$
is defined to be $|\Delta|\equiv\bigvee_{i=1}^{n}|\Delta_{i}|$ while
the \emph{di}am\emph{eter\index{diameter of an n-interval}} of $\Delta$
is defined to be $\left\Vert \Delta\right\Vert \equiv\sqrt{\sum_{i=1}^{n}|\Delta_{i}|^{2}}$.
The intervals $\Delta_{1},\cdots,\Delta_{n}$ in $R$ are then called
the \index{factors of an n-interval@factors of an \emph{n-}interval}\emph{factors}
of the $n$-interval $\Delta$. If each of the factors of $\Delta$
is an open interval, then $\Delta$ is called an \emph{open $n$-interval}\index{open n-interval@open \emph{n}-interval}.
If all the factors of $\Delta$ are of the same length, then $\Delta$
is called an \index{n-cube@\emph{n}-cube}$n$-\emph{cube}. The \emph{center}
$x$ of $\Delta$ is defined to be $x=(x_{1},\cdots,x_{n})$ where
$x_{i}$ is the mid-point of $\Delta_{i}$ for each $i=1,\cdots,n$.
By Fubini's Theorem, we have $\mu(\Delta)=\prod_{i=1}^{n}\mu_{1}(\Delta_{i})$. 

The next lemma is the formula for the change of integration variables
in the special case of a linear transformation.
\begin{lem}
\textbf{\emph{\label{Lem.  Volumn of a Parallelopiped} (Volume of
a parallelepiped)}}. Let $a\in R^{n}$ be arbitrary. Let $F$ be an
$n\times n$ matrix with $|\det F|>0$. Let $f:R^{n}\rightarrow R^{n}$
be the linear function defined by $f(y)=a+Fy$ for each $y\in R^{n}.$
Let $\Delta\equiv\Delta_{1}\times\cdots\times\Delta_{n}$ be any $n$-interval
in $R^{n}$. Then $f(\Delta)$ is an integrable set, and $\mu(f(\Delta))=|\det F|\mu(\Delta)$. 
\end{lem}
\begin{proof}
Since $|\det F|>0$, the matrix $F$ has an inverse $G\equiv F^{-1}$,
with $|\det G|=|\det F|^{-1}$. Define the linear function $g:R^{n}\rightarrow R^{n}$
by $g(v)=b+Gv$ for each $v\in R^{n}$, where $b\equiv-Ga$. Then
$g$ is the inverse of $f$ on $R^{n}$. Note that the desired equality
is equivalent to
\[
\mu(\Delta)=|\det F|^{-1}\mu(g^{-1}(\Delta)),
\]
or
\begin{equation}
\mu(\Delta)=|\det G|\int\cdots\int1_{\bigcap_{k=1}^{n}\{(v_{1},\cdots v_{n}):b_{k}+\sum_{i=1}^{n}G_{k,i}v_{i}\in\Delta_{k}\}}dv_{1}\cdots dv_{n}.\label{eq:temp-77}
\end{equation}

First assume that the lemma holds for each open $n$-interval. Let
$\Delta\equiv\Delta_{1}\times\cdots\times\Delta_{n}$ be an arbitrary
$n$-interval. For each $i=1,\cdots,n$ and $k\geq1$ define the open
interval $\Delta_{i}^{(k)}\equiv(r,s),$ $(r-\frac{1}{k},s),$ $(r,s+\frac{1}{k}),$
or $(r-\frac{1}{k},s+\frac{1}{k})$ according as $\Delta_{i}=(r,s)$,
$[r,s)$, $(r,s]$, or $[r,s]$ for some $r,s\in R$. For each $k\geq1$
define the open $n$-interval $\Delta^{(k)}\equiv\Delta_{1}^{(k)}\times\cdots\times\Delta_{n}^{(k)}$.
 Then (i) $\Delta^{(k)}\supset\Delta^{(k+1)}$ for each $k\geq1$,
(ii) $\bigcap_{k=1}^{\infty}\Delta^{(k)}=\Delta$ and (iii) $\mu(\Delta^{(k)})\downarrow\mu(\Delta)$.
Here condition (iii) follows from Proposition \ref{Prop. Intervals are Lebesgue integrable},
which says that every interval in $R$ is Lebesgue integrable and
has measure equal to its length. By assumption, the lemma holds for
$\Delta^{(k)}$ for each $k\geq1$. Hence $f(\Delta^{(k)})$ is an
integrable set, and $\mu(f(\Delta^{(k)}))=|\det F|\mu(\Delta^{(k)})$
for each $k\geq1$. This implies that $\mu(f(\Delta^{(k)}))\downarrow|\det F|\mu(\Delta)$
while $\bigcap_{k=1}^{\infty}f(\Delta^{(k)})=f(\Delta)$. Hence $f(\Delta)$
is an integrable set, with
\[
\mu(f(\Delta))=\lim_{k\rightarrow\infty}\mu(f(\Delta^{(k)}))=|\det F|\mu(\Delta).
\]
Thus the lemma holds also for the arbitrary $n$-interval $\Delta$.
We see that we need only prove the lemma for open $n$-intervals.

To that end, let $\Delta\equiv(z_{1},z'_{1})\times\cdots\times(z_{n},z'_{n})$
be an arbitrary open $n$-interval. Proceed by induction. 

First assume that $n=1$. The $1\times1$ matrix $F$ can be regarded
as a real number. By hypothesis $|F|=|\det F|>0$. We will assume
$F<0$, the positive case being similar. Then $f(\Delta)=(a+Fz'_{1},a+Fz_{1})$,
an interval. Therefore, by Proposition \ref{Prop. Intervals are Lebesgue integrable},
we have 
\[
\mu(f(\Delta))=a+Fz_{1}-a-Fz'_{1}=F(z_{1}-z'_{1})=|F|\mu(\Delta),
\]
establishing the lemma for $n=1$.

Suppose the lemma has been proved for $n=1,\cdots,m-1$ for some $m>1$.
We will give the proof for $n=m$. Since $|\det G|=|\det F|^{-1}>0$,
there exists, by Lemma \ref{Lem. Matrix Basics}, some $j=1,\cdots,n$
such that $|G_{n,j}|>0$. Without loss of generality, and for ease
of notations, we assume that $j=n$ and that $G_{n,n}>0$. 

Let $(v_{1},\cdots,v_{n-1})\in R^{n-1}$ be arbitrary. Define the
function $\varphi\equiv\varphi_{(v_{1},\cdots,v_{n-1})}:R\rightarrow R$
by
\[
\varphi(v_{n})\equiv\varphi_{(v_{1},\cdots,v_{n-1})}(v_{n})\equiv g_{n}(v_{1},\cdots,v_{n})=b_{n}+\sum_{i=1}^{n-1}G_{n,i}v_{i}+G_{n,n}v_{n}
\]
for each $v_{n}\in R$, Since $G_{n,n}>0$, a direct substitution
shows that $\varphi$ has an inverse $\psi\equiv\psi_{(v_{1},\cdots,v_{n-1})}:R\rightarrow R$
given by
\begin{equation}
\psi(x_{n})\equiv\psi_{(v_{1},\cdots,v_{n-1})}(x_{n})\equiv-G_{n,n}^{-1}b_{n}-\sum_{i=1}^{n-1}G_{n,n}^{-1}G_{n,i}v_{i}+G_{n,n}^{-1}x_{n}\label{eq:temp-84}
\end{equation}
for each $x_{n}\in R$. Define the linear function
\[
\lambda_{k}(v_{1},\cdots,v_{n-1},x_{n})\equiv g_{k}(v_{1},\cdots,v_{n-1},\psi_{(v_{1},\cdots,v_{n-1})}(x_{n}))
\]
\[
=b_{k}+\sum_{i=1}^{n-1}G_{k,i}v_{i}+G_{k,n}(-G_{n,n}^{-1}(b_{n}+\sum_{i=1}^{n-1}G_{n,i}v_{i})+G_{n,n}^{-1}x_{n})
\]
\begin{equation}
\equiv c_{k}+\sum_{i=1}^{n-1}\Lambda_{k,i}v_{i}+\Lambda_{k,n}x_{n}\label{eq:temp-24}
\end{equation}
for each $k=1,\cdots,n$, and for each $(v_{1},\cdots v_{n-1},x_{n})\in R^{n}$,
where $c_{k}\equiv b_{k}-G_{n,n}^{-1}G_{k,n}b_{n}$, $\Lambda_{k,i}\equiv G_{k,i}-G_{n,n}^{-1}G_{k,n}G_{n,i}$,
and $\Lambda_{k,n}\equiv G_{k,n}G_{n,n}^{-1}$ for each $k=1,\cdots,n$
and for each $i=1,\cdots,n-1$. Now define an $n\times n$ matrix
$\Psi$ by (i$'$) $\Psi_{k,i}\equiv1$ or $\Psi_{k,i}\equiv0$ according
as $k=i$ or $k\neq i$, for $k,i=1,\cdots,n-1,$ (ii$'$) $\Psi_{n,i}\equiv-G_{n,n}^{-1}G_{n,i}$
and $\Psi_{i,n}=0$ for $i=1,\cdots,n-1,$ (iii$'$) $\Psi_{n,n}\equiv G_{n,n}^{-1}$.
Then $\psi(v_{1},\cdots,v_{n-1},x_{n})\equiv-G_{n,n}^{-1}b_{n}+\sum_{i=1}^{n-1}\Psi_{n,i}v_{i}+\Psi_{n,n}x_{n}$.
Note that $\Psi$ is a triangular matrix, with each entry above the
diagonal equal to $0$. Hence $\det\Psi=\Psi_{1,1}\cdots\Psi_{n,n}=G_{n,n}^{-1}$.
Moreover, a direct matrix multiplication verifies that $\Lambda\equiv G\Psi$.
Consequently $\det\Lambda=\det G\cdot\det\Psi=G_{n,n}^{-1}\det G$.
Hence $|\det\Lambda|>0$. 

Note from their definition that $\Lambda_{n,i}=0$ for $i=1,\cdots,n-1,$
and that $\Lambda_{n,n}=1$. By Lemma \ref{Lem. Matrix Basics}, it
follows that $\det\Lambda=\det\Lambda'$ where $\Lambda'$ is the
$(n-1)\times(n-1)$ matrix obtained by deleting the $n$-th row and
$n$-th column in $\Lambda$. 

Let $A$ be an arbitrary subset of $R$ which is dense in $R$. Assume
that the lemma holds for each open $n$-interval each of whose factors
has endpoints in $A$. Now let $\Delta$ be an arbitrary open $n$-interval.
Let $(\Delta^{(k)})_{k=1,2,\cdots}$ be a sequence of open $n$-intervals
whose factors have endpoints in $A$, such that (i$''$) $\Delta^{(k)}\subset\Delta^{(k+1)}$
for each $k\geq1$, (ii$''$) $\bigcup_{k=1}^{\infty}\Delta^{(k)}=\Delta$
and (iii$''$) $\mu(\Delta^{(k)})\uparrow\mu(\Delta)$. By assumption,
the lemma holds for $\Delta^{(k)}$ for each $k\geq1$. Hence $f(\Delta^{(k)})$
is an integrable set, and $\mu(f(\Delta^{(k)}))=|\det F|\mu(\Delta^{(k)})$
for each $k\geq1$. This implies that $\mu(f(\Delta^{(k)}))\uparrow\left|\det F\right|\mu(\Delta)$
while $\bigcup_{k=1}^{\infty}f(\Delta^{(k)})=f(\Delta)$. Hence $f(\Delta)$
is integrable with
\[
\mu(f(\Delta))=\lim_{k\rightarrow\infty}\mu(f(\Delta^{(k)}))=|\det F|\mu(\Delta).
\]
Thus the lemma holds also for the arbitrary $n$-interval $\Delta$.
We see that we need only prove the lemma for open $n$-intervals with
endpoints in $A$. In the remainder of this proof, we will let $A$
be the set of continuity points of the $\mathrm{measurable}$ functions
$g_{k}$ and $\lambda_{k}$ for each $k=1,\cdots,n$. By Proposition
\ref{Prop.  Countable Exceptional Pts for Meas X int A}, the set
$A$ is the metric complement of a countable subset of $R$, whence
$A$ is dense in $R$. We will prove the lemma for an arbitrary open
$n$-interval whose factors have endpoints in $A$, thereby completing
the proof. 

To that end, let $\Delta\equiv(z_{1},z'_{1})\times\cdots\times(z_{n},z'_{n})$
be an open $n$-interval with $z_{1},\cdots,z_{n},$ $z'_{1},\cdots,z'_{n}$
$\in A$. Then the set $g^{-1}(\Delta)=\bigcap_{k=1}^{n}(z_{k}<g_{k}<z'_{k})$
is a $\mathrm{measurable}$ subset in $R^{n}$. Since $g^{-1}\equiv f$
is continuous, $g^{-1}(\Delta)$ is also bounded. Hence $g^{-1}(\Delta)$
is an integrable set. Similarly $\lambda^{-1}(\Delta)$ is an integrable
set. By Fubini's Theorem, there exists a full subset $D$ of $R^{n-1}$
such that  $1_{\lambda^{-1}(\Delta)}(v_{1},\cdots,v_{n-1},\cdot)$
is an integrable indicator on $R$ for each $(v_{1},\cdots,v_{n-1})\in D$.
In other words, in terms of equality \ref{eq:temp-24}, the set
\begin{equation}
\Gamma_{v_{1},\cdots,v_{n-1}}\equiv\bigcap_{k=1}^{n}\{x_{n}:c_{k}+\sum_{i=1}^{n-1}\Lambda_{k,i}v_{i}+\Lambda_{k,n}x_{n}\in\Delta_{k}\}\label{eq:temp-82}
\end{equation}
is an integrable subset of $R$ for each $(v_{1},\cdots,v_{n-1})\in D$.
In view of equality \ref{eq:temp-24}, we have 
\begin{equation}
\Gamma_{v_{1},\cdots,v_{n-1}}=\{x_{n}:g(v_{1},\cdots,v_{n-1},\psi(x_{n}))\in\Delta\}.\label{eq:temp-81}
\end{equation}
Moreover, since $\Lambda_{n,i}=0$ for $i=1,\cdots,n-1,$ and $\Lambda_{n,n}=1$
while $c_{n}=0$, equality \ref{eq:temp-82} reduces to 
\begin{equation}
\Gamma_{v_{1},\cdots,v_{n-1}}=\bigcap_{k=1}^{n-1}\{x_{n}:c_{k}+\sum_{i=1}^{n-1}\Lambda_{k,i}v_{i}+\Lambda_{k,n}x_{n}\in\Delta_{k}\}\cap(z_{n},z'_{n}).\label{eq:temp-85}
\end{equation}

Define $\Delta'\equiv\Delta_{1}\times\cdots\times\Delta_{n-1}$. Using
Fubini's Theorem, and applying the induction hypothesis, in the for$m$
of equality \ref{eq:temp-77}, to the $(n-1)$-interval $\Delta'$,
we obtain 
\[
\int\cdots\int_{\Delta}dx_{1}\cdots dx_{n}=\int1_{(z_{n},z'_{n})}(x_{n})\mu'(\Delta')dx_{n}
\]
\[
=\int1_{(z_{n},z'_{n})}(x_{n})(|\det\Lambda'|\qquad\qquad\qquad\qquad\qquad\qquad\qquad\qquad\qquad\qquad\qquad\qquad
\]
\[
\qquad\times\int\cdots\int1_{\bigcap_{k=1}^{n-1}\{(v_{1},\cdots v_{n-1}):c_{k}+\sum_{i=1}^{n-1}\Lambda_{k,i}v_{i}+\Lambda_{k,n}x_{n}\in\Delta_{k}\}}dv_{1}\cdots dv_{n-1})dx_{n}
\]
\[
=|\det\Lambda'|\int\cdots\int\qquad\qquad\qquad\qquad\qquad\qquad\qquad\qquad\qquad\qquad\qquad\qquad\qquad
\]
\[
\qquad1_{\bigcap_{k=1}^{n-1}\{(v_{1},\cdots v_{n-1},x_{n}):c_{k}+\sum_{i=1}^{n-1}\Lambda_{k,i}v_{i}+\Lambda_{k,n}x_{n}\in\Delta_{k};x_{n}\in(z_{n},z'_{n})\}}dv_{1}\cdots dv_{n-1}dx_{n}
\]
\begin{equation}
=|\det\Lambda|\int\cdots\int\mu_{1}(\Gamma_{v_{1},\cdots,v_{n-1}})dv_{1}\cdots dv_{n-1},\label{eq:temp-23}
\end{equation}
where $\Lambda'$ is the $(n-1)\times(n-1)$ submatrix of $\Lambda$
defined previously, and where $\mu'$ is the Lebesgue measure on $R^{n-1}$.

Let $(v_{1},\cdots,v_{n-1})\in D$ be arbitrary. Write $\Gamma\equiv\Gamma_{v_{1},\cdots,v_{n-1}}$,
$\varphi\equiv\varphi_{v_{1},\cdots,v_{n-1}}$, and $\psi\equiv\psi_{v_{1},\cdots,v_{n-1}}$
for short. Define
\[
\Theta\equiv\Theta_{v_{1},\cdots,v_{n-1}}\equiv\{v_{n}:g(v_{1},\cdots,v_{n-1},v_{n})\in\Delta\}.
\]
Since $\psi$ has an inverse function $\varphi$, it follows from
equality \ref{eq:temp-81} that 
\begin{equation}
\psi(\Gamma)=\{v_{n}:g(v_{1},\cdots,v_{n-1},v_{n})\in\Delta\}\equiv\Theta.\label{eq:temp-79}
\end{equation}
Let $\varepsilon>0$ be arbitrary. At lease one of the three following
conditions must hold: (a) $\mu_{1}(\Gamma)<\varepsilon$ and $\mu_{1}(\Theta)<\varepsilon$,
(b) $\mu_{1}(\Gamma)>0$, or (c) $\mu_{1}(\Theta)>0$. If condition
(a) holds, then
\begin{equation}
|\mu_{1}(\Gamma)-|G_{n,n}|\mu_{1}(\Theta)|<\varepsilon(1+|G_{n,n}|).\label{eq:temp-80}
\end{equation}
Suppose condition (b) holds. Then $\Gamma$ is non-empty. Suppose
condition (c) holds. Then $\Theta$ is non-empty, whence $\Gamma$
is also non-empty, thanks to equality \ref{eq:temp-79}. Therefore,
if either of conditions (b) and (c) holds, then $\Gamma$ is non-empty,
and Lemma \ref{Lem. Condition for  non empty set in R to be an interval},
together with equality \ref{eq:temp-85}, implies that $\Gamma$ is
an open interval. The induction hypothesis therefore applies, in either
case, to the interval $\Gamma$ and the linear function $\psi$, yielding
\begin{equation}
\mu_{1}(\Gamma)=|G_{n,n}|\mu_{1}(\psi(\Gamma))\equiv|G_{n,n}|\mu_{1}(\Theta).\label{eq:temp-83}
\end{equation}
Inequality \ref{eq:temp-80} is thus established under each of conditions
(a-c). Since $\varepsilon>0$ is arbitrary, we conclude therefore
that $\mu_{1}(\Gamma)=|G_{n,n}|\mu_{1}(\Theta)$ for each $(v_{1},\cdots,v_{n-1})$
in the full set $D\subset R^{n-1}$. Combining with equality \ref{eq:temp-23},
and recalling that $\det\Lambda=G_{n,n}^{-1}\det G$, we obtain 
\[
\int\cdots\int_{\Delta}dx_{1}\cdots dx_{n}=|\det\Lambda|\int\cdots\int\mu_{1}(\Gamma_{v_{1},\cdots,v_{n-1}})dv_{1}\cdots dv_{n-1}
\]
\[
=|\det\Lambda|\cdot|G_{n,n}|\int\cdots\int\mu_{1}(\Theta)dv_{1}\cdots dv_{n-1}
\]
\[
=|\det G|\int\cdots\int(\int1_{\Theta_{v_{1},\cdots,v_{n-1}}}(v_{n})dv_{n})dv_{1}\cdots dv_{n-1}
\]
\[
\equiv|\det G|\int\cdots\int(\int1_{\{v_{n}:g(v_{1},\cdots,v_{n})\in\Delta\}}dv_{n})dv_{1}\cdots dv_{n-1}
\]
\[
=|\det G|\int\cdots\int1_{\{(v_{1},\cdots,v_{n}):g(v_{1},\cdots,v_{n})\in\Delta\}}dv_{1}\cdots dv_{n}
\]
\[
\equiv|\det G|\int\cdots\int1_{g^{-1}(\Delta)}dv_{1}\cdots dv_{n}=|\det G|\mu(g^{-1}(\Delta)).
\]
Induction is completed.
\end{proof}
In the following, if $\Delta\equiv[a_{1},a_{1}+t]\times\cdots\times[a_{n},a_{n}+t]$
is any closed $n$-cube with length $|\Delta|=t$, then for each $r\in(-1,1)$
we define $\Delta^{r}\equiv[a_{1}-\frac{r}{2}t,a_{1}+t+\frac{r}{2}t]\times\cdots\times[a_{n}-\frac{r}{2}t,a_{n}+t+\frac{r}{2}t]$.
Similar notations are used for half-open $n$-cubes. Thus $\Delta^{r}$
has the same center as $\Delta$ and is similar to $\Delta$ with
a scale of $1+r$. In particular $|\Delta^{r}|=(1+r)|\Delta|$ and
$\Delta^{0}=\Delta$. 

Recall that a subset $H$ of an open set $A$ is said to be well contained
in $A$, or $H\Subset A$ in sy$m$bols, if $H\subset K\subset K_{a}\subset A$
for some compact subset $K$ and some $a>0$. Recall that $K_{a}$
stands for the compact $a$-neighborhood of $K$.
\begin{lem}
\label{Lem. Parallelopiped Sandwich} \textbf{\emph{(Parallelepiped
sandwich).}} Let $A,B$ be open subsets of $R^{n}$. Let $f:A\rightarrow B$
be a differentiable function on $A$, with derivative $F$. Suppose
$f$ has an inverse function $g:B\rightarrow A$ which is differentiable
on $B$, with derivative $G$, and which is Lipschitz continuous on
$B$. Suppose $b,c>0$ are such that $|\det G|\geq c$ and $|G_{i,j}|\leq b$
for $i,j=1,\cdots,n$ on $B$. Let $\delta:(0,\infty)\rightarrow(0,\infty)$
be an operation. Then for each $\varepsilon>0$ there exists $\tau=\tau(\varepsilon,n,b,c,\delta)>0$
with the following properties. Let $t\in(0,\tau)$ be arbitrary. Let
$\Delta\equiv(x_{1}-t,x_{1}+t]\times\cdots\times(x_{n}-t,x_{n}+t]$
be any half open $n$-cube with center $x\equiv(x_{1},\cdots,x_{n})$,
such that (i) $\Delta\Subset A$, (ii) $f(\Delta)$ and $f(\overline{\Delta})$
are integrable sets with $\mu(f(\Delta))=\mu(f(\overline{\Delta}))$,
where $\overline{\Delta}\equiv[x_{1}-t,x_{1}+t]\times\cdots\times[x_{n}-t,x_{n}+t]$,
(iii) $f(\overline{\Delta})\Subset B$, and (iv) $\delta$ is a modulus
of differentiability of $f$ on $\overline{\Delta}$. Then

\begin{equation}
|\mu(f(\Delta))(|\det F(x)|^{-1})-\mu(\Delta)|\leq n(2^{n-1}+1)\varepsilon\mu(\Delta).\label{eq:temp-92}
\end{equation}
\end{lem}
\begin{proof}
By hypothesis $g$ is Lipschitz continuous on $B$ with some Lipschitz
constant $c_{g}>0$. 

Let $x\in A$ be arbitrary, and let $u\equiv f(x)$. Define the linear
function $\bar{f}:R^{n}\rightarrow R^{n}$ by $\bar{f}(y)\equiv f(x)+F(x)(y-x)$
for each $y\in R^{n}$. By the Chain Rule, we have $F(x)=G(u)^{-1}$.
By hypothesis $|\det G(u)|\geq c$ and $|G_{i,j}(x)|\leq b$ for each
$i,j=1,\cdots,n$. By Lemma \ref{Lem. Matrix Basics},
\begin{equation}
\left\Vert F(x)w\right\Vert =\left\Vert G(u)^{-1}w\right\Vert \leq\beta\left\Vert w\right\Vert \label{eq:temp-113}
\end{equation}
and\emph{ $\left\Vert F(x)^{-1}w\right\Vert =\left\Vert G(u)w\right\Vert \leq nb\left\Vert w\right\Vert $
}for each\emph{ $w\in R^{n}$, }where\emph{ $\beta\equiv n!c^{-1}b^{n-1}$.}
Replacing $w$ by $F(x)w$ in the last inequality, we obtain $\left\Vert F(x)w\right\Vert \geq n^{-1}b^{-1}\left\Vert w\right\Vert $
for each \emph{$w\in R^{n}$.} Let $\varepsilon\in(0,1)$ be arbitrary.
Define $\tau\equiv n^{-\frac{1}{2}}\delta(\frac{1}{2}n^{-\frac{3}{2}}b^{-1}\varepsilon)$.
We will show that $\tau$ has the desired properties. 

To that end, let $t\in(0,\tau)$ be arbitrary. Let $\Delta\equiv(x_{1}-t,x_{1}+t]\times\cdots\times(x_{n}-t,x_{n}+t]$
be any half open $n$-cube with center $x\equiv(x_{1},\cdots,x_{n})$,
such that conditions (i-iii) hold. We will prove inequality \ref{eq:temp-87}.
First note that since $\Delta\Subset A$, and $f(\overline{\Delta})\Subset B$,
the functions $f$ and $g$ are uniformly continuous on $\Delta$
and $f(\overline{\Delta})$ respectively. Consequently $f(\overline{\Delta})$
is a compact set.

Consider any $y\in\Delta$. We have $\left\Vert y-x\right\Vert \leq\sqrt{n}t<\sqrt{n}\tau\leq\delta(\frac{1}{2}n^{-\frac{3}{2}}b^{-1}\varepsilon)$,
whence
\[
\left\Vert f(y)-\bar{f}(y)\right\Vert \equiv\left\Vert f(y)-f(x)-F(x)(y-x)\right\Vert \leq\frac{1}{2}n^{-\frac{3}{2}}b^{-1}\varepsilon\left\Vert y-x\right\Vert 
\]
\begin{equation}
\leq\frac{1}{2}n^{-\frac{3}{2}}b^{-1}\varepsilon\sqrt{n}t=\frac{1}{2}n^{-1}b^{-1}\varepsilon t.\label{eq:temp-112}
\end{equation}
Write $z\equiv f(y)-f(x)-F(x)(y-x)$ and $v\equiv F(x)^{-1}z$. Then
$\left\Vert v\right\Vert \leq nb\left\Vert z\right\Vert \leq\frac{1}{2}nbn^{-1}b^{-1}\varepsilon t$
$<\varepsilon t$. Hence $y_{k}+v_{k}\in(x_{k}-t(1+\varepsilon),x_{k}+t(1+\varepsilon)]$
for each $k=1,\cdots,n$. Thus $y+v\in\Delta^{\varepsilon}$, in the
notations introduced before this lemma. Moreover
\[
f(y)=f(x)+F(x)(y-x)+z
\]
\[
=f(x)+F(x)(y-x+v)\equiv\bar{f}(y+v)\in\bar{f}(\Delta^{\varepsilon}).
\]
Since $y\in\Delta$ is arbitrary, we see that $f(\Delta)\subset\bar{f}(\Delta^{\varepsilon})$
and so $\mu(f(\Delta))\leq\mu(\bar{f}(\Delta^{\varepsilon}))$. On
the other hand, by Lemma \ref{Lem.  Volumn of a Parallelopiped},
we have $\mu(\bar{f}(\Delta^{\varepsilon}))=\mu(\Delta^{\varepsilon})|\det F(x)|$.
Combining, we obtain
\begin{equation}
\mu(f(\Delta))\leq\mu(\Delta^{\varepsilon})|\det F(x)|=(1+\varepsilon)^{n}\mu(\Delta)|\det F(x)|.\label{eq:temp-78}
\end{equation}
At the same time, since $\varepsilon\in[0,1]$, we have $|n(1+\varepsilon)^{n-1}-n|\leq n2^{n-1}$.
Taylor's Theorem \ref{Thm. Tayor's Theorem} therefore yields the
bound $|(1+\varepsilon)^{n}-(1+n\varepsilon)|\leq n2^{n-1}\varepsilon$
whence $|(1+\varepsilon)^{n}-1|\leq n(2^{n-1}+1)\varepsilon$. Therefore
inequality \ref{eq:temp-78} implies
\[
\mu(f(\Delta))|\det F(x)|^{-1}-\mu(\Delta)
\]
 
\begin{equation}
\leq((1+\varepsilon)^{n}-1)\mu(\Delta)\leq n(2^{n-1}+1)\varepsilon\mu(\Delta).\label{eq:temp-89}
\end{equation}

Next consider any $y\in\Delta^{-\varepsilon}$. We will prove that
$\bar{f}(y)\in f(\overline{\Delta})$. For the sake of a contradiction,
suppose $a\equiv d(\bar{f}(y),f(\overline{\Delta}))>0$. Write $u\equiv f(x)$.
Since $f(\overline{\Delta})\Subset B$, there exists $r>0$ such that
$f(\overline{\Delta})_{r}\subset B$. Let $\varepsilon'>0$ be so
small that $3\varepsilon'<r\wedge a$ and $3\varepsilon'+3\beta c_{g}\varepsilon'<\frac{1}{2}n^{-1}b^{-1}\varepsilon t$.
Let $p\geq1$ be so large that $\frac{1}{p}\beta\sqrt{n}t\leq\varepsilon'$.
For each $k=0,\cdots,p$ define $y_{k}\equiv(1-\frac{k}{p})x+\frac{k}{p}y$
and define $v_{k}\equiv\bar{f}(y_{k})=u+\frac{k}{p}F(x)(y-x)$. Then
$y_{k}\in\Delta^{-\varepsilon}$ for each $k=0,\cdots,p$ because,
being an $n$-interval, $\Delta^{-\varepsilon}$ is convex. Moreover,
$d(v_{0},f(\overline{\Delta}))=d(f(x),f(\overline{\Delta}))=0$, and
$d(v_{p},f(\overline{\Delta}))=d(\bar{f}(y),f(\overline{\Delta}))=a>0$
by assumption. Hence there exists $j=1,\cdots,p$ such that $d(v_{j-1},f(\overline{\Delta}))<2\varepsilon'<r$
and $d(v_{j},f(\overline{\Delta}))>\varepsilon'$. Since $\left\Vert v_{j}-v_{j-1}\right\Vert =\frac{1}{p}\left\Vert F(x)(y-x)\right\Vert \leq\frac{1}{p}\beta\left\Vert y-x\right\Vert \leq\frac{1}{p}\beta\sqrt{n}t\leq\varepsilon'$,
it follows that $d(v_{j},f(\overline{\Delta}))<3\varepsilon'\leq r$.
Consequently $v_{j-1},v_{j}\in f(\overline{\Delta})_{r}\subset B$.
Define $z'\equiv g(v_{j-1})$ and $z''\equiv g(v_{j})$. Then $v_{j-1}=f(z')$
and $v_{j}=f(z'')$. Since $d(v_{j},f(\overline{\Delta}))>\varepsilon'$,
we have $d(z'',\overline{\Delta})>0$, thanks to the continuity of
$f$. On the other hand, since $d(v_{j-1},f(\overline{\Delta}))<2\varepsilon'$,
there exists $z\in\Delta$ such that $\left\Vert v_{j-1}-f(z)\right\Vert <2\varepsilon'$.
Since the function $g$ is Lipschitz continuous on $B$ with Lipschitz
constant $c_{g}$, the last inequality yields $\left\Vert g(v_{j-1})-g(f(z))\right\Vert <2c_{g}\varepsilon'$,
or equivalently $\left\Vert z'-z\right\Vert <2c_{g}\varepsilon'$.
Similarly, from $\left\Vert v_{j}-v_{j-1}\right\Vert \leq\varepsilon'$
we deduce $\left\Vert z''-z'\right\Vert \leq c_{g}\varepsilon'$.
Combining, we have $\left\Vert z''-z\right\Vert <3c_{g}\varepsilon'$.
Since $y_{j}\in\Delta^{1-\varepsilon}$, and since $d(z'',\overline{\Delta})>0$,
we have $\left\Vert y_{j}-z''\right\Vert >\varepsilon t$. Hence
\[
\left\Vert \bar{f}(y_{j})-\bar{f}(z'')\right\Vert =\left\Vert F(x)(y_{j}-z'')\right\Vert 
\]
\begin{equation}
\geq n^{-1}b^{-1}\left\Vert y_{j}-z''\right\Vert >n^{-1}b^{-1}\varepsilon t.\label{eq:temp-86}
\end{equation}
On the other hand,
\[
\left\Vert \bar{f}(y_{j})-\bar{f}(z'')\right\Vert \equiv\left\Vert v_{j}-\bar{f}(z'')\right\Vert 
\]
\[
\leq\left\Vert v_{j}-v_{j-1}\right\Vert +\left\Vert v_{j-1}-f(z)\right\Vert +\left\Vert f(z)-\bar{f}(z)\right\Vert +\left\Vert \bar{f}(z)-\bar{f}(z'')\right\Vert 
\]
\[
\leq\varepsilon'+2\varepsilon'+\left\Vert f(z)-\bar{f}(z)\right\Vert +\left\Vert F(x)(z-z'')\right\Vert .
\]
In view of inequalities \ref{eq:temp-112} and \ref{eq:temp-113},
the last expression is bounded by
\[
3\varepsilon'+\frac{1}{2}n^{-1}b^{-1}\varepsilon t+\beta\left\Vert z-z''\right\Vert 
\]
\[
\leq3\varepsilon'+\frac{1}{2}n^{-1}b^{-1}\varepsilon t+3\beta c_{g}\varepsilon'
\]
\[
<\frac{1}{2}n^{-1}b^{-1}\varepsilon t+\frac{1}{2}n^{-1}b^{-1}\varepsilon t=n^{-1}b^{-1}\varepsilon t,
\]
which contradicts inequality \ref{eq:temp-86}. Thus we see that $d(\bar{f}(y),f(\overline{\Delta}))=0$.
We conclude that $\bar{f}(y)\in f(\overline{\Delta})$. Since $y\in\Delta^{-\varepsilon}$
is arbitrary, we have $\bar{f}(\Delta^{-\varepsilon})\subset f(\overline{\Delta})$.
On the other hand, by Lemma \ref{Lem.  Volumn of a Parallelopiped},
we have $\mu(\bar{f}(\Delta^{-\varepsilon}))=\mu(\Delta^{-\varepsilon})|\det F(x)|$.
Combining, we obtain
\begin{equation}
(1-\varepsilon)^{n}\mu(\Delta)|\det F(x)|=\mu(\Delta^{-\varepsilon})|\det F(x)|=\mu(\bar{f}(\Delta^{-\varepsilon}))\leq\mu(f(\overline{\Delta}))=\mu(f(\Delta)).\label{eq:temp-90}
\end{equation}
At the same time, since $\varepsilon\in[0,1]$, we have $|n(1-\varepsilon)^{n-1}-n|\leq n\leq n2^{n-1}$.
Taylor's Theorem yields $|(1-\varepsilon)^{n}-(1-n\varepsilon)|\leq n2^{n-1}\varepsilon$
whence $|(1-\varepsilon)^{n}-1|\leq n(2^{n-1}+1)\varepsilon$. Inequality
\ref{eq:temp-90} therefore yields
\[
\mu(f(\Delta))|\det F(x)|^{-1}-\mu(\Delta)\geq((1-\varepsilon)^{n}-1)\mu(\Delta)
\]
\begin{equation}
\geq-n(2^{n-1}+1)\varepsilon\mu(\Delta).\label{eq:temp-91}
\end{equation}
Combining inequalities \ref{eq:temp-89} and \ref{eq:temp-91}, we
obtain the desired inequality \ref{eq:temp-92}. The lemma is proved.
\end{proof}
Let $K$ be an arbitrary compact subset of $R^{n}$. The closed $r$-neighborhood
$(d(\cdot,K)\leq r)\equiv\{u\in R^{n}:d(u,K)\leq r\}$ is compact
for all but countably many $r>0$. At the same time, all but countably
many $r>0$ are regular points of the $\mathrm{measurable}$ function
$d(\cdot,K)$. For the remainder of this section, we will write $K_{r}$
for $(d(\cdot,K)\leq r)$ only with the implicit assumption that $(d(\cdot,K)\leq r)$
is compact and that $r>0$ is a regular point of $d(\cdot,K)$. Thus
$K_{r}$ is compact and $\mathrm{measurable}$, hence integrable.
Furthermore $K_{r}^{c}\equiv(d(\cdot,K)\leq r)^{c}=(d(\cdot,K)>r)$
according to Corollary \ref{Cor. In sigma-finite space, (X<=00003Dt)^c=00003D(t<X) etc}.
\begin{lem}
\label{Lem. Special half open n-intervals} \textbf{\emph{(Special
half open $n$-interval). }}Let $A,B$ be an open subset of $R^{n}$.
Let $f:A\rightarrow B$ be a function which is uniformly continuous
on compact subsets well contained in $A$. Suppose $g:B\rightarrow A$
is an inverse function of $f$ and is uniformly continuous on compact
subsets well contained in $B$. Let $H,K$ be compact subsets of $R^{n}$
with $H\Subset A$ and $K\Subset B$. Then there exists \textup{$\alpha\in R$}\textup{\emph{
with the following properties. Suppose $\overline{\Delta}$ is a closed
$n$-interval with $\overline{\Delta}=[\alpha+q_{1},\alpha+q'_{1}]\times\cdots\times[\alpha+q_{n},\alpha+q'_{n}]$
where $q_{i},q'_{i}$ are arbitrary rational n}}um\textup{\emph{bers
with $q_{i}<q'_{i}$ for each }}\emph{$i=1,\cdots,n$. }Suppose \textup{\emph{$\overline{\Delta}\subset H$}}
and\emph{ }\textup{\emph{$f(\overline{\Delta})\subset K$. Let $\Delta\equiv(\alpha+q_{1},\alpha+q'_{1}]\times\cdots\times(\alpha+q_{n},\alpha+q'_{n}]$.
Then $f(\Delta)$ and $f(\overline{\Delta})$ are integrable sets
with $\mu(f(\Delta))=\mu(f(\overline{\Delta}))$. Moreover $f(\overline{\Delta})^{c}=(d(\cdot,f(\overline{\Delta}))>0)$
a.e. In other words, the }}meas\textup{\emph{ure theoretic c}}om\textup{\emph{pl}}em\textup{\emph{ent
of $f(\overline{\Delta})$ is equal to its }}me\textup{\emph{tric
c}}om\textup{\emph{pl}}em\textup{\emph{ent.}}
\end{lem}
\begin{proof}
Since $f(\Delta)=g^{-1}(\Delta)$, the measurability of $f(\Delta)$
for all rational numbers \emph{$q_{i},q'_{i}$ }and for all but countably
many $\alpha\in R$ would follow if $g$ is a $\mathrm{measurable}$
function on $R^{n}$. However, $g$ is not necessarily defined a.e.
on $R^{n}$. Hence we introduce a function $\bar{g}$ which is $\mathrm{measurable}$
on $R^{n}$ and which is equal to $g$ on a neighborhood of $K$.

Since $K\Subset B$, there exists $\rho>0$ such that $K_{\rho}\Subset B$.
Write $M\equiv K_{\rho}$. Then $M\subset M_{r}\subset B$ for some
$r>0$. Moreover $M$ is a compact and integrable set, and $M_{r}^{c}\equiv(d(\cdot,M)>r)$,
according to the remarks preceding this lemma. Let $i=1,\cdots,n$
be arbitrary. Define a function $\bar{g_{i}}:M_{r}\cup M_{r}^{c}\rightarrow R$
by (i) $\bar{g_{i}}(u)\equiv g_{i}(u)(1-r{}^{-1}d(u,M))_{+}$ if $u\in M_{r}$,
and (ii) $\bar{g_{i}}(u)\equiv0$ if $u\in M_{r}^{c}$. We will show
that $\bar{g_{i}}$ is uniformly continuous on $M_{r}\cup M_{r}^{c}$.
In view of conditions (i) and (ii), the function $\bar{g_{i}}$ is
uniformly continuous on each of $M_{r}$ and $M_{r}^{c}$. Let $\varepsilon>0$
be arbitrary. It suffices to show that $|\bar{g_{i}}(u)-\bar{g_{i}}(v)|<\varepsilon$
if $u\in M_{r}$ and $v\in M_{r}^{c}$ are such that $\left\Vert u-v\right\Vert <\eta$
for sufficiently small $\eta>0$. Let $\eta\in(0,r)$ be so small
that $|g_{i}|\eta\leq\varepsilon r$ on $M_{r}$. Consider $u\in M_{r}$
and $v\in M_{r}^{c}$ with $\left\Vert u-v\right\Vert <\eta$. Since
$v\in M_{r}^{c}$, we have $d(v,M)>r$. Consequently $d(u,M)>r-\eta$.
It follows from the defining condition (i) that $|\bar{g_{i}}(u)|<|g_{i}(u)|(1-\frac{r-\eta}{r})=|g_{i}(u)|\frac{\eta}{r}\leq\varepsilon$.
Since $v\in M_{r}^{c}$, we have $\bar{g_{i}}(v)=0$. Combining, we
obtain $|\bar{g_{i}}(u)-\bar{g_{i}}(v)|<\varepsilon$. Thus we see
that $\bar{g_{i}}$ is uniformly continuous on $M_{r}\cup M_{r}^{c}$.
Since $M_{r}\cup M_{r}^{c}$ is dense in $R^{n}$, the function $\bar{g_{i}}$
can be extended to a continuous function $\bar{g_{i}}:R^{n}\rightarrow R$.
In view of condition (i), we have $\bar{g_{i}}=g_{i}$ on $M$. Being
continuous on $R^{n}$, the function $\bar{g_{i}}$ is $\mathrm{measurable}$.
Therefore $\bar{g_{i}}-q$ is a $\mathrm{measurable}$ function on
$R^{n}$ for each $i=1,\cdots,n$ and each rational number $q$. Because
$R^{n}$ is $\sigma$-finite, there exists, according to Proposition
\ref{Prop. Continuity Pts for Measurable Function}, a real number
$\alpha$ which is a continuity point of $\bar{g_{i}}-q$ for each
$i=1,\cdots,n$ and each rational number $q$. Let $i=1,\cdots,n$
be arbitrary, and let $q$ be an arbitrary rational number. Then $(\bar{g_{i}}-q\leq\alpha)M$,
$(\bar{g_{i}}-q>\alpha)M$, $(\bar{g_{i}}-q\geq\alpha)M$ are integrable
sets, with $\mu((\bar{g_{i}}-q>\alpha)M)=\mu((\bar{g_{i}}-q\geq\alpha)M)$.
Hence $(\bar{g_{i}}-q>\alpha)M=(\bar{g_{i}}-q\geq\alpha)M$ a.e. In
other words $M(\bar{g_{i}}\leq\alpha+q)$, $M(\bar{g_{i}}>\alpha+q)$,
and $M(\bar{g_{i}}\geq\alpha+q)$ are integrable sets, with $M(\bar{g_{i}}>\alpha+q)=M(\bar{g_{i}}\geq\alpha+q)$
a.e.

Now suppose $\Delta$ is a half open $n$-interval with $\Delta=(\alpha+q_{1},\alpha+q'_{1}]\times\cdots\times(\alpha+q_{n},\alpha+q'_{n}]$
where $q_{i},q'_{i}$ are rational numbers for each $i=1,\cdots,n$,
such that \emph{$\overline{\Delta}\subset H$} and\emph{ $g^{-1}(\overline{\Delta})\subset K$.
}Write $x_{i}\equiv\alpha+q_{i}$ and $x'_{i}\equiv\alpha+q'_{i}$
for each $i=1,\cdots,n$. We saw in the last paragraph that $M(\bar{g_{i}}\leq x'_{i})$
and $M(\bar{g_{i}}\geq x'_{i})$ are integrable sets for each $i=1,\cdots,n$.
Since $g^{-1}(\overline{\Delta})\subset K\subset M$ and $g=\bar{g}$
on $M$, we have
\[
g^{-1}(\overline{\Delta})=\{u\in M:\;x_{i}\leq g_{i}(u)\leq x'_{i}\mbox{ for each }i=1,\cdots n\}
\]
\[
=\{u\in M:\;x_{i}\leq\bar{g_{i}}(u)\leq x'_{i}\mbox{ for each }i=1,\cdots n\}
\]
\begin{equation}
=\bigcap_{i=1}^{n}M(x_{i}\leq\bar{g_{i}})M(\bar{g_{i}}\leq x'_{i}),\label{eq:temp-88}
\end{equation}
which is an integrable set. Similarly $g^{-1}(\Delta)=\bigcap_{i=1}^{n}M(x_{i}<\bar{g_{i}})M(\bar{g_{i}}\leq x'_{i})$
is an integrable set. Since $M(x_{i}\leq\bar{g_{i}})=M(x_{i}<\bar{g_{i}})$
a.e. for each $i=1,\cdots,n$, we have $g^{-1}(\Delta)=g^{-1}(\overline{\Delta})$
a.e. It follows that $\mu(g^{-1}(\Delta))=\mu(g^{-1}(\overline{\Delta}))$.
In other words \emph{$\mu(f(\Delta))=\mu(f(\overline{\Delta}))$.}

It remains to prove that \emph{$f(\overline{\Delta})^{c}=(d(\cdot,f(\overline{\Delta}))>0)$.}

Consider an arbitrary point $u\in M$. Suppose $u\in g^{-1}(\overline{\Delta})^{c}$.
Then, according to equality \ref{eq:temp-88}, we have $u\in(M(x_{i}\leq\bar{g_{i}})M(\bar{g_{i}}\leq x'_{i}))^{c}$
for some $i=1,\cdots,n$. In view of Corollary \ref{Cor. In sigma-finite space, (X<=00003Dt)^c=00003D(t<X) etc},
we therefore have $u\in M^{c}\cup(\bar{g}_{i}<x_{i})\cup(x'_{i}<\bar{g}_{i})$.
Since $u\in M$ and so $\bar{g}(u)=g(u)$, it follows that $u\in(g_{i}<x_{i})\cup(x'_{i}<g_{i})$.
Thus $g_{i}(u)<x_{i}-\varepsilon$ or $x'_{i}+\varepsilon<g_{i}(u)$
for some $\varepsilon>0$. Suppose $d(u,g^{-1}(\overline{\Delta}))<\delta_{g}(\varepsilon)$,
where $\delta_{g}$ is a modulus of continuity of $g$ on $M$. Then
$d(g(u),\overline{\Delta})<\varepsilon$. Consequently $g_{i}(u)\in[x_{i}-\varepsilon,x'_{i}+\varepsilon]$,
a contradiction. Hence $d(u,g^{-1}(\overline{\Delta}))\geq\delta_{g}(\varepsilon)>0$.
We conclude that $Mg^{-1}(\overline{\Delta})^{c}\subset(d(\cdot,g^{-1}(\overline{\Delta}))>0)$.
Next consider any $u\in M^{c}=(d(\cdot,M)>0)$. Since $g^{-1}(\overline{\Delta})\subset M$,
we have $d(u,g^{-1}(\overline{\Delta}))\geq d(u,M)>0$. We conclude
that $M^{c}g^{-1}(\overline{\Delta})^{c}\subset(d(\cdot,g^{-1}(\overline{\Delta}))>0)$.
Combining, we obtain
\begin{equation}
(M\cup M^{c})g^{-1}(\overline{\Delta})^{c}\subset(M\cup M^{c})(d(\cdot,g^{-1}(\overline{\Delta}))>0).\label{eq:temp-103}
\end{equation}

Conversely, consider any $u\in M$. Let $y\equiv g(u)$ and so $u=f(y)$.
Suppose $u\in(d(\cdot,g^{-1}(\overline{\Delta}))>0)$. In other words
$d(u,f(\overline{\Delta}))>0$, or
\begin{equation}
d(f(y),f(\overline{\Delta}))\geq\varepsilon>0\label{eq:temp-102}
\end{equation}
 for some $\varepsilon>0$. By hypothesis, \emph{$\overline{\Delta}\subset H\subset H_{a}\Subset A$
}for some $a>0$, and $f$ is uniformly continuous on $H_{a}$ with
some modulus of continuity $\delta_{f}$. Suppose $d(y,\overline{\Delta})<a$.
Then $d(y,H)\leq d(y,\overline{\Delta})<a$. It follows that $y\in H_{a}$
and so $d(y,\overline{\Delta})\geq\delta_{f}(\varepsilon)>0$, in
view of inequality \ref{eq:temp-102}. Combining, we see that in any
case $d(y,\overline{\Delta})\geq a\wedge\delta_{f}(\varepsilon)>0$.
Therefore there exists $i=1,\cdots,n$ such that $y_{i}<x_{i}$ or
$x'_{i}<y_{i}$. In other words, $g_{i}(u)<x_{i}$ or $x'_{i}<g_{i}(u)$.
Since $u\in M$ and $g=\bar{g}$ on $M$, this implies that
\[
u\in M(\bar{g}_{i}<x_{i})\cup M(x'_{i}<\bar{g}_{i})=M(M(x_{i}\leq\bar{g_{i}}\leq x'_{i}))^{c}
\]
\[
\subset M\cap\{u\in M:\;x_{i}\leq\bar{g_{i}}(u)\leq x'_{i}\mbox{ for each }i=1,\cdots n\}^{c}=Mg^{-1}(\overline{\Delta})^{c},
\]
in view of equality \ref{eq:temp-88}. Thus we see that
\[
M(d(\cdot,g^{-1}(\overline{\Delta}))>0)\subset Mg^{-1}(\overline{\Delta})^{c}\subset(M\cup M^{c})g^{-1}(\overline{\Delta})^{c}
\]
At the same time, since $M\supset g^{-1}(\overline{\Delta})$, we
have $M^{c}\subset g^{-1}(\overline{\Delta})^{c}$
\[
M^{c}(d(\cdot,g^{-1}(\overline{\Delta}))>0)\subset M^{c}=M^{c}g^{-1}(\overline{\Delta})^{c}\subset(M\cup M^{c})g^{-1}(\overline{\Delta})^{c}.
\]
Combining, we obtain
\begin{equation}
(M\cup M^{c})(d(\cdot,g^{-1}(\overline{\Delta}))>0)\subset(M\cup M^{c})g^{-1}(\overline{\Delta})^{c}.\label{eq:temp-104}
\end{equation}

Expressions \ref{eq:temp-103} and \ref{eq:temp-104} together show
that $(M\cup M^{c})(d(\cdot,g^{-1}(\overline{\Delta}))>0)=(M\cup M^{c})g^{-1}(\overline{\Delta})^{c}$.
Since $M\cup M^{c}$ is a full set, it follows that $(d(\cdot,g^{-1}(\overline{\Delta}))>0)=g^{-1}(\overline{\Delta})^{c}$
a.e., as desired. 
\end{proof}
\begin{lem}
\label{Lem. Partition of n-cube} \textbf{\emph{(Partition of $n$-cube).
}}Let $\Delta\equiv\prod_{j=1}^{n}[a_{j},a_{j}+t]$ be an arbitrary
$n-cube.$ Let $p\geq1$. For each $j=1,\cdots,n$ and $k=0,\cdots,p$
define $a_{j,k}\equiv a_{j}+\frac{k}{p}t$. For each $\kappa\in\{1,\cdots,p\}^{n}$
define $\Delta_{\kappa}\equiv\prod_{j=1}^{n}[a_{j,\kappa_{j}-1},a_{j,\kappa_{j}}]$.
Then
\[
\Delta^{\eta/p}=\bigcup_{\kappa\in\{1,\cdots,p\}^{n}}\Delta_{\kappa}^{\eta}
\]
for arbitrary $\eta>0$. 
\end{lem}
\begin{proof}
Let $\eta>0$ be arbitrary. Then, for each $j=1,\cdots,n$, we have
\[
[a_{j},a_{j}+t]^{\eta/p}\equiv[a_{j}-\frac{\eta}{2p}t,a_{j}+t+\frac{\eta}{2p}t]
\]
\[
=\bigcup_{k=1}^{p}[a_{j,k-1}-\frac{\eta}{2p}t,a_{j,k}+\frac{\eta}{2p}t]\equiv\bigcup_{k=1}^{p}[a_{j,k-1},a_{j,k}]^{\eta}.
\]
Hence
\[
\Delta^{\eta/p}=\prod_{j=1}^{n}[a_{j},a_{j}+t]^{\eta/p}=\prod_{j=1}^{n}\bigcup_{k=1}^{p}[a_{j,k-1},a_{j,k}]^{\eta}
\]
\[
=\bigcup_{\kappa\in\{1,\cdots,p\}^{n}}\prod_{j=1}^{n}[a_{j,\kappa_{j}-1},a_{j,\kappa_{j}}]^{\eta}=\bigcup_{\kappa\in\{1,\cdots,p\}^{n}}\Delta_{\kappa}^{\eta}.
\]
\end{proof}
\begin{lem}
\label{Lem. Criterion for compactness of  f(H) where H is compact in R^n}
\textbf{\emph{(Sufficient condition for image of compact set to be
compact). }}Let $A$ and $B$ be open subsets of $R^{n}$. Suppose
the following two conditions hold.
\end{lem}
\begin{enumerate}
\item \emph{$g:B\rightarrow A$ is a function which is differentiable on
$B$, such that for each c}om\emph{pact subset $K\Subset B$ we have
$|J|\equiv|\det G|\geq c$ on $K$ for s}om\emph{e $c>0$, where $G$
is the derivative of $g$.}
\item \emph{There exists a function $f:A\rightarrow B$ which is the inverse
of $g$, such that for each c}om\emph{pact subset $H$ of $A$ with
$H\Subset A$ we have $f(H)\Subset B$. }
\end{enumerate}
\emph{Then, for each} \emph{c}om\emph{pact subset $H\Subset A$, the
set $K\equiv f(H)$ is a c}om\emph{pact subset of $B$.}
\begin{proof}
Let $H$ be a compact subset \emph{$H\Subset A$}. Conditions 1 and
2, in view of Corollary \ref{Cor. Condition for Inverse Func to be differentiable},
implies that $f$ is differentiable on $A$. Consequently $f$ is
uniformly continuous on $H$. It follows that $K$ is totally bounded.
It remains to show that $K$ is complete. By condition 2, the set
$K\equiv f(H)$ is such that $K\Subset B$. In other words $K\subset K'\Subset B$
for some compact subset $K'$ of $B$. Hence it suffices to show that
$K$ is closed. Suppose $(u_{k})_{k=1,2,\cdots}$ is a sequence in
$K$ such that $u_{k}\rightarrow u$ for some $u\in K'$. By Proposition
\ref{Prop. g differentiable on open set=00003D> g Lipsitz on compact subsets},
$g$ is Lipschitz continuous on any compact subset well contained
in $B$. In particular $g$ is uniformly continuous on $K'$. Therefore
$g(u_{k})\rightarrow g(u)$ in $A$. Since $g(u_{k})\in H$ for each
$k\geq1$ and since $H$ is closed, we have $g(u)\in H$. Therefore
$u=f(g(u))\in f(H)=K$. Thus $K$ is closed. 
\end{proof}
\begin{lem}
\label{Lem. Special case Change of Integration Variables} \textbf{\emph{(Change
of integration variables, special case).}} Let $A$ and $B$ be open
subsets of $R^{n}$. Suppose the following four conditions hold.
\end{lem}
\begin{enumerate}
\item \emph{$g:B\rightarrow A$ is a function which is differentiable on
$B$, such that for each c}om\emph{pact subset $K\Subset B$ we have
$|J|\equiv|\det G|\geq c$ on $K$ for s}om\emph{e $c>0$, where $G$
is the derivative of $g$. }
\item \emph{There exists a function $f:A\rightarrow B$ which is the inverse
of $g$, such that for each c}om\emph{pact subset $H$ of $A$ with
$H\Subset A$ we have $f(H)\Subset B$.}
\item \emph{$H$ is a given c}om\emph{pact integrable subset of $A$ with
$H\Subset A$, such that $\mu(H_{a})\downarrow\mu(H)$ as $a\downarrow0$.}
\end{enumerate}
\emph{Then $f(H)$ is an integrable set. Moreover, for each $X\in C(R^{n})$,
the function $1_{f(H)}X(g)|J|$ is integrable, with }
\begin{equation}
\int\cdots\int_{H}X(x)dx_{1}\cdots dx_{n}=\int\cdots\int_{f(H)}X(g(u))|J(u)|du_{1}\cdots du_{n}\label{eq:temp-65}
\end{equation}
\emph{where the functions $g$ and $J$ have been extended to $R^{n}$
by $g=0=J$ on $B^{c}$.}
\begin{proof}
For abbreviation, we will write $\int\cdot dx$ for $\int\cdots\int\cdot dx_{1}\cdots dx_{n}$,
and similarly write$\int\cdot du$ for $\int\cdots\int\cdot du_{1}\cdots du_{n}$.

Let $H$ be as given in Condition 3. Conditions 1 and 2, in view of
Corollary \ref{Cor. Condition for Inverse Func to be differentiable},
imply that $f$ is differentiable on $A$. Therefore $f$ is Lipschitz
continuous on compact subsets well contained in $A$, by Proposition
\ref{Prop. g differentiable on open set=00003D> g Lipsitz on compact subsets}.
At the same time, Lemma \ref{Lem. Criterion for compactness of  f(H) where H is compact in R^n}
implies that $K\equiv f(H)$ is compact. Let $X\in C(R^{n})$ be arbitrary.
In the following proof, we may assume that $X\geq0$. The general
case of $X\in C(R)$ follows from $X=2X_{+}-|X|$ and from linearity.

Since $K\Subset B$ by Condition 3, we have $K_{\rho}\Subset B$ for
some $\rho>0$. Let $\delta$ be a modulus of differentiability of
$g$ on $K_{\rho}$. By Condition 1, $|J|\geq c$ on $K_{\rho}$ for
some $c>0$. By Proposition \ref{Prop. g differentiable on open set=00003D> g Lipsitz on compact subsets},
$g$ is Lipschitz continuous on $K_{\rho}$, with some Lipschitz constant
$c_{g}>0$. By Proposition \ref{Prop.g  differenitable =00003D> G unif cont on well contained compacts},
the partial derivative $G_{i,j}$ of $g$ is uniformly continuous
on $K_{\rho}$ for each $i,j=1,\cdots,n$. Hence the Jacobian $J$
is uniformly continuous on $K_{\rho}$, as is the product $X(g)J$.
Let $\delta_{X}$ denote a modulus of continuity of $X$ on $R^{n}$,
and let $\delta_{XgJ}$ denote a modulus of continuity of $X(g)J$
on $K_{\rho}$. Let $\beta>0$ be such that $|X|\leq\beta$ on $R^{n}$.
Let $b>0$ be such that $|G_{i,j}|\leq b$ on $K_{\rho}$ for each
$i,j=1,\cdots,n$. 

Since $H\Subset A$, we have $H_{a}\Subset A$ for some $a>0$. By
Lemma \ref{Lem. Special half open n-intervals}, there exists $\alpha\in R$
with the following properties. Suppose $\overline{\Delta}$ is a closed
$n$-interval with $\overline{\Delta}\equiv[\alpha+q_{1},\alpha+q'_{1}]\times\cdots\times[\alpha+q_{n},\alpha+q'_{n}]$,
where $q_{i},q'_{i}$ are rational numbers for each $i=1,\cdots,n$.
Suppose $\overline{\Delta}\subset H_{a}$ and suppose $f(\overline{\Delta})\subset K_{\rho}$.
Let $\Delta\equiv(\alpha+q_{1},\alpha+q'_{1}]\times\cdots\times(\alpha+q_{n},\alpha+q'_{n}]$.
Then $f(\Delta)$ and $f(\overline{\Delta})$ are integrable sets
with \emph{$\mu(f(\Delta))=\mu(f(\overline{\Delta}))$}, and
\begin{equation}
f(\overline{\Delta})^{c}=(d(\cdot,f(\overline{\Delta}))>0)\quad\mbox{a.e.}\label{eq:temp-93}
\end{equation}

In the remainder of this proof, $\rho$ and $\alpha$ will be fixed.
We will call a number $\alpha'\in R$ \emph{special} if $\alpha'=\alpha+q$
for some rational number $q$. We will call a half open $n$-interval
$\Delta$ \emph{a special} $n$-interval if $\Delta=(x_{1},x'_{1}]\times\cdots\times(x_{n},x'_{n}]$
where $x_{i},x'_{i}$ are special numbers for each $i=1,\cdots,n$.
A special $n$-interval which is an $n$-cube will be called a special
$n$-cube. Suppose $\Delta$ is a special $n$-interval. If, in addition,
$\overline{\Delta}\subset H_{a}$ and $f(\overline{\Delta})\subset K_{\rho}$,
then, according to the preceding paragraph, the sets $f(\Delta)$
and $f(\overline{\Delta})$ are integrable sets with \emph{$\mu(f(\Delta))=\mu(f(\overline{\Delta}))$},
and equality \ref{eq:temp-93} holds. Let $p\geq1$ be any integer.
Then, for each $\kappa\in\{1,\cdots,p\}^{n}$, the subinterval
\[
\Delta_{\kappa}\equiv\prod_{i=1}^{n}(x_{i}+(\kappa_{i}-1)p^{-1}(x'_{i}-x_{i}),\;x_{i}+\kappa_{i}p^{-1}(x'_{i}-x_{i})]
\]
is again special, because $\kappa_{i}p^{-1}(x'_{i}-x_{i})$ is a rational
number for each $i=1,\cdots,n$. Thus we see that each special $n$-interval
$\Delta$ can be subdivided into arbitrarily small special subintervals,
by taking sufficiently large $p\geq1$. 

As remarked at the opening of this proof, the function $f$ is differentiable
on $A$, and is therefore uniformly continuous on $H_{a}$. Hence
$H\subset H_{a_{0}}\Subset A$ and $f(H_{a_{0}})\subset K_{\rho}\Subset B$
for some sufficiently small $a_{0}\in(0,a)$. Recall that, by convention,
$a_{0}$ is so chosen that $H_{a_{0}}$ is an integrable and compact
subset of $R^{n}$. Lemma \ref{Lem. Criterion for compactness of  f(H) where H is compact in R^n}
implies that $f(H_{a_{0}})$ is compact. Let $\delta$ denote a modulus
of differentiability of $f$ on $H_{a_{0}}$, and let $\delta_{f}$
denote a modulus of continuity of $f$ on $H_{a_{0}}$.

Now let $\Delta_{0}$ be a fixed special $n$-cube such that $H_{a_{0}}\subset\Delta_{0}$.
Set $\Phi_{0}\equiv\{\Delta_{0}\}$. Then, trivially, $H_{a_{0}}\subset\bigcup_{\Delta\in\Phi_{0}}\overline{\Delta}^{\eta}$
for arbitrary $\eta>0$.

For abbreviation, write $\varepsilon_{k}\equiv\frac{1}{k}$ for each
$k\geq1$. We will construct, inductively for each $k\geq0$, a real
number $a_{k}>0$ and a set $\Phi_{k}$ of special $n$-cubes of equal
size, with $a_{k}\downarrow0$, and with the following properties:
(I) if $k\geq1$ then each member of $\Phi_{k}$ is a subset of some
member of $\Phi_{k-1}$, (II) $f(H_{a_{k}})$ is compact, (III) if
$k\geq1$ then $\mu(H_{a_{k}})-\mu(H)<\varepsilon_{k}$, and (IV)
$H_{a_{k}}\subset\bigcup_{\Delta\in\Phi_{k}}\overline{\Delta}^{\eta}$
for arbitrary $\eta>0$. For the case $k=0$, the objects $a_{0}$
and $\Phi_{0}$ have been constructed which trivially satisfy the
Conditions (I-IV).

Suppose $k\geq1$ is such that $a_{j}$ and $\Phi_{j}$ have been
constructed for $j=0,\cdots,k-1$ satisfying Conditions (I-IV). Since
$f$ is uniformly continuous on $H_{a_{k-1}}\subset H_{a_{0}}$, and
since $f(H)\equiv K$, there exists $a_{k}\in(0,\frac{1}{2}a_{k-1})$
so small that $H_{a_{k}}\subset H_{a_{k-1}}\Subset A$ and $f(H_{a_{k}})\subset K_{\rho}\Subset B$.
Moreover, by Condition 3 in the hypothesis, we can choose $a_{k}$
so small that $\mu(H_{a_{k}})-\mu(H)<\varepsilon_{k}.$ Recall that,
by convention, $a_{k}$ is so chosen that $H_{a_{k}}$ is an integrable
and compact subset of $R^{n}$. Lemma \ref{Lem. Criterion for compactness of  f(H) where H is compact in R^n}
implies that $f(H_{a_{k}})$ is compact. Thus $H_{a_{k}}$ satisfies
Condition (II) and (III). We proceed to establish also Conditions
(I) and (IV) for $k$.

By Lemma \ref{Lem. Parallelopiped Sandwich} applied to the function
$f$, there exists $\tau_{k}=\tau(\varepsilon_{k},n,b,c,\delta)>0$
with the following properties. Let $t\in(0,\tau_{k})$ be arbitrary.
Let $\Delta\equiv(x_{1}-t,x_{1}+t]\times\cdots\times(x_{n}-t,x_{n}+t]$
be any half open $n$-cube with center $x\equiv(x_{1},\cdots,x_{n})$,
such that (i) , (ii) $f(\Delta)$ and $f(\overline{\Delta})$ are
integrable sets with $\mu(f(\Delta))=\mu(f(\overline{\Delta}))$,
where $\overline{\Delta}\equiv[x_{1}-t,x_{1}+t]\times\cdots\times[x_{n}-t,x_{n}+t]$,
(iii) , and (iv) $\delta$ is a modulus of differentiability of $f$
on $\overline{\Delta}$. Then
\begin{equation}
|\mu(f(\Delta))(|\det F(x)|^{-1})-\mu(\Delta)|\leq n(2^{n-1}+1)\varepsilon_{k}\mu(\Delta).\label{eq:temp-87}
\end{equation}

Continue with the induction process for the construction of $a_{k}$
and $\Phi_{k}$. By the induction hypothesis, the special $n$-cubes
in $\Phi_{k-1}$ are of equal diameter, which we denote by $t$. Thus
$\left\Vert \Delta\right\Vert =t$ for each $\Delta\in\Phi_{k-1}$.
Let $s_{k}\in(0,\frac{1}{2}(a_{k-1}-a_{k})$) be arbitrary. Let $p\equiv p_{k}\geq1$
be so large that
\begin{equation}
p^{-1}t<\tau_{k}\wedge s_{k}\wedge\delta_{X}(\varepsilon_{k})\wedge\delta_{f}(\delta_{XgJ}(\varepsilon_{k})).\label{eq:temp-99}
\end{equation}
Consider an arbitrary special $n$-cube $\Delta\in\Phi_{k-1}$. Write
$\Delta\equiv(\alpha_{1},\alpha_{1}+t]\times\cdots\times(\alpha_{n},\alpha_{n}+t]$,
and, for each $\kappa\in\{1,\cdots,p\}^{n}$, define the special $n$-cube
\[
\Delta_{\kappa}\equiv\prod_{i=1}^{n}(\alpha_{i}+(\kappa_{i}-1)p^{-1}t,\alpha_{i}+\kappa_{i}p^{-1}t].
\]
Define the set $\Psi_{k}\equiv\{\Delta_{\kappa}:\Delta\in\Phi_{k-1};\kappa\in\{1,\cdots,p\}^{n}\}$.
Thus $\Psi_{k}$ is a set of special $n$-cubes each of which is a
subcube of some member of $\Phi_{k-1}$. Moreover, according to Lemma
\ref{Lem. Partition of n-cube}, we have $\overline{\Delta}^{\eta/p}\subset\bigcup_{\kappa\in\{1,\cdots,p\}^{n}}\overline{\Delta}_{\kappa}^{\eta}$
for each $\Delta\in\Phi_{k-1}$ and for arbitrary$\eta>0$. Hence
\begin{equation}
\bigcup_{\Delta\in\Phi_{k-1}}\overline{\Delta}^{\eta/p}\subset\bigcup_{\Delta\in\Phi_{k-1}}\bigcup_{\kappa\in\{1,\cdots,p\}^{n}}\overline{\Delta}_{\kappa}^{\eta}=\bigcup_{\Gamma\in\Psi_{k}}\overline{\gamma}^{\eta}\label{eq:temp-96}
\end{equation}
for arbitrary $\eta>0$. Furthermore, $\left\Vert \Gamma\right\Vert =p^{-1}t$
for each $\Gamma\in\Psi_{k}$. We can partition the set $\Psi_{k}$
into two subsets $\Phi_{k}$ and $\Phi_{k}'$ such that 
\begin{equation}
\mbox{if }\Gamma\in\Phi_{k}\;\mbox{then }d(\overline{\gamma},H_{a_{k}})<s_{k},\label{eq:temp-100}
\end{equation}
and
\begin{equation}
\mbox{if }\Gamma\in\Phi_{k}'\;\mbox{then }d(\overline{\gamma},H_{a_{k}})>\frac{s_{k}}{2}.\label{eq:temp-101}
\end{equation}
The set $\Phi_{k}$, being a subset of $\Psi_{k}$, automatically
satisfies Condition (I). 

Let $\eta>0$ be arbitrary. We will show that $H_{a_{k}}\subset\bigcup_{\Gamma\in\Phi_{k}}\overline{\gamma}^{\eta}$.
Consider an arbitrary $x\in H_{a_{k}}$. Let $\zeta\equiv\eta\wedge\frac{s_{k}}{2\sqrt{n}}$.
By the induction hypothesis, we have $H_{a_{k}}\subset H_{a_{k-1}}\subset\bigcup_{\Delta\in\Phi_{k-1}}\overline{\Delta}^{\zeta/p}$.
Thus $x\in\overline{\Delta}^{\zeta/p}$ for some $\Delta\in\Phi_{k-1}$.
Therefore, according to expression \ref{eq:temp-96}, we have $x\in\bigcup_{\Gamma\in\Psi_{k}}\overline{\gamma}^{\zeta}$.
Consequently, $x\in\overline{\gamma}^{\zeta}$ for some $\Gamma\in\Psi_{k}$.
Hence $d(x,\overline{\gamma})\leq\sqrt{n}\zeta\leq\frac{s_{k}}{2}$.
Inequality \ref{eq:temp-101} therefore implies that $\Gamma\notin\Phi'_{k}$.
Hence $\Gamma\in\Phi_{k}$. Consequently $x\in\bigcup_{\Gamma\in\Phi_{k}}\overline{\gamma}^{\zeta}$.
We conclude that $H_{a_{k}}\subset\bigcup_{\Gamma\in\Phi_{k}}\overline{\gamma}^{\zeta}\subset\bigcup_{\Gamma\in\Phi_{k}}\overline{\gamma}^{\eta}$,
thereby establishing Condition (IV) for $k$. Induction is completed.

Next let $k\geq1$ be arbitrary. Consider any $\Gamma\in\Phi_{k}$.
We have, on the one hand, $\left\Vert \overline{\gamma}\right\Vert =p^{-1}t<s_{k}$,
and, on the other hand, $d(\overline{\gamma},H_{a_{k}})<s_{k}$ in
view of inequality \ref{eq:temp-100}. Since $2s_{k}<a_{k-1}-a_{k}$,
we obtain 
\[
\overline{\gamma}\subset(H_{a_{k}})_{2s_{k}}\subset H_{a_{k-1}}\subset H_{a_{0}}\subset H_{a}.
\]
It follows that
\begin{equation}
f(\overline{\gamma})\subset f(H_{a_{k-1}}).\label{eq:temp-94}
\end{equation}
Consequently $f(\overline{\gamma})\subset K_{\rho}$. Hence, since
$\Gamma$ is a special $n$-cube, the sets $f(\Gamma)$ and $f(\overline{\gamma})$
are integrable, with \emph{$\mu(f(\Gamma))=\mu(f(\overline{\gamma}))$}.
Moreover $f(\overline{\gamma})^{c}=(d(\cdot,f(\overline{\gamma}))>0)$
on some full set $D_{0}$, according to equality \ref{eq:temp-93}.
Define $M_{k}\equiv\bigcup_{\Gamma\in\Phi_{k}}\overline{\gamma}$.
Then $f(M_{k})=\bigcup_{\Gamma\in\Phi_{k}}f(\overline{\gamma})$.
From relation \ref{eq:temp-94}, we see that $f(M_{k})\subset f(H_{a_{k-1}})$. 

We will next show, in the other direction, that $f(H_{a_{k}})\subset f(M_{k})$
a.e. Consider an arbitrary $v\in D\cap f(H_{a_{k}})$ where $D$ is
the full set
\begin{equation}
D\equiv D_{0}\cap\bigcap_{\Gamma\in\Phi_{k}}(f(\overline{\gamma})\cup f(\overline{\gamma})^{c})=D_{0}\cap\bigcap_{\Gamma\in\Phi_{k}}(f(\overline{\gamma})\cup(d(\cdot,f(\overline{\gamma}))>0)).\label{eq:temp-95}
\end{equation}
Then $v=f(y)$ for some $y\in H_{a_{k}}$. In view of equality \ref{eq:temp-95},
we have either $v\in f(\overline{\gamma})$ for some $\Gamma\in\Phi_{k}$,
or $v\in(d(\cdot,f(\overline{\gamma}))>0)$ for each $\Gamma\in\Phi_{k}$.
In other words either $v\in f(M_{k})$, or $v\in(d(\cdot,f(\overline{\gamma}))>0)$
for each $\Gamma\in\Phi_{k}$. We will show that $v\in f(M_{k})$.
Suppose, for the sake of a contradiction, that $d(v,f(\overline{\gamma}))>0$
for each $\Gamma\in\Phi_{k}$. Then $d(f(y),f(\overline{\gamma}))=d(v,f(\overline{\gamma}))>\zeta$
for each $\Gamma\in\Phi_{k}$, for some $\zeta>0$. Define $\eta\equiv\frac{s_{k}}{2\sqrt{n}}\wedge\frac{1}{2\sqrt{n}}\delta_{f}(\zeta)$,
where $\delta_{f}$ is the previously defined modulus of continuity
of $f$ on the compact set $H_{a_{0}}$. Then, for each $\Gamma\in\Phi_{k}$,
the assumption that $d(y,\overline{\gamma})<2\sqrt{n}\eta$ would
imply $d(y,\overline{\gamma})<\delta_{f}(\zeta)$ and so $d(f(y),f(\overline{\gamma}))<\zeta$,
a contradiction. Hence
\begin{equation}
d(y,\overline{\gamma})\geq2\sqrt{n}\eta\quad\mbox{ for each }\Gamma\in\Phi_{k}.\label{eq:temp-97}
\end{equation}
On the other hand, by Condition (IV), we have $y\in H_{a_{k}}\subset\bigcup_{\Gamma\in\Phi_{k}}\overline{\gamma}^{\eta}$.
Consequently there exists $\Gamma\in\Phi_{k}$ such that $y\in\overline{\gamma}^{\eta}$.
Hence $d(y,\overline{\gamma})\leq\sqrt{n}\eta$, contradicting inequality
\ref{eq:temp-97}. We conclude that $v\in f(M_{k})$. Since $v\in D\cap f(H_{a_{k}})$
is arbitrary, we have established that $D\cap f(H_{a_{k}})\subset f(M_{k})$.
Consequently $f(H_{a_{k}})\subset f(M_{k})$ a.e. 

Summing up, we have proved that
\begin{equation}
f(H_{a_{k}})\subset f(M_{k})\subset f(H_{a_{k-1}})\mbox{ a.e.}\label{eq:temp-75}
\end{equation}
Since $H=\bigcap_{k=1}^{\infty}H_{a_{k}}$, and since $f$ is a bijection,
we have $f(H)=\bigcap_{k=1}^{\infty}f(H_{a_{k}})$. It follows that
\[
f(H)=\bigcap_{k=1}^{\infty}f(H_{a_{k}})\subset\bigcap_{k=1}^{\infty}f(M_{k})\subset\bigcap_{k=1}^{\infty}f(H_{a_{k-1}})=f(H)\quad\mbox{a.e. }
\]
Hence $f(H)=\bigcap_{k=1}^{\infty}f(M_{k})$ a.e.

As a special case, when we apply the arguments in the previous paragraphs
with $f$ replaced by the identity function on $R^{n}$, expression
\ref{eq:temp-75} yields $H_{a_{k}}\subset M_{k}\subset H_{a_{k-1}}$
a.e. 

Consider any $\Gamma\in\Phi_{k}$. Then $d(\overline{\gamma},H)<s{}_{k}$
by Condition \ref{eq:temp-100}. Let $x_{\Gamma}$ be the center of
$\Gamma$, and let $u_{\Gamma}\equiv f(x_{\Gamma})$. By expression
\ref{eq:temp-94}, we have $f(\overline{\gamma})\subset f(H_{a_{0}})\subset B$,
and so $\delta$ is a modulus of differentiability of $f$ on $\overline{\gamma}$.
We saw earlier that $f(\Gamma)$ is an integrable set. Since $|\Gamma|=p^{-1}t\leq\tau_{k}$,
all the conditions are satisfied for inequality \ref{eq:temp-87}
to hold for the $n$-cube $\Gamma$, yielding 
\[
|\mu(f(\Gamma))(|\det F(x_{\Gamma})|^{-1})-\mu(\Gamma)|\leq n(2^{n-1}+1)\varepsilon_{k}\mu(\Gamma),
\]
where $F(x_{\Gamma})$ is the derivative of $f$ at $x_{\Gamma}$.
Since , it follows that 
\[
|\mu(\Gamma)-\mu(f(\Gamma))(|\det G(u_{\Gamma})|)|\leq n(2^{n-1}+1)\varepsilon_{k}\mu(\Gamma).
\]
Equivalently
\[
|\int_{\Gamma}dx-\int_{f(\Gamma)}|J(u_{\Gamma})|du|\leq n(2^{n-1}+1)\varepsilon_{k}\mu(\Gamma).
\]
Multiplying by $X(x_{\Gamma})\equiv X(g(u_{\Gamma}))$, and noting
that $|X|\leq\beta$, we obtain
\begin{equation}
|\int_{\Gamma}X(x_{\Gamma})dx-\int_{f(\Gamma)}X(x_{\Gamma})|J(u_{\Gamma})|du|\leq\beta n(2^{n-1}+1)\varepsilon_{k}\mu(\Gamma).\label{eq:temp-76}
\end{equation}
At the same time $\left\Vert \overline{\gamma}\right\Vert =p^{-1}t<\delta_{X}(\varepsilon_{k})\wedge\delta_{f}(\delta_{XgJ}(\varepsilon_{k}))$
by inequality \ref{eq:temp-99}. Hence $|X(x)-X(x_{\Gamma})|<\varepsilon_{k}$
and $\left\Vert f(x)-f(x_{\Gamma})\right\Vert <\delta_{XgJ}(\varepsilon_{k})$
for each $x\in\Gamma$. Consequently, $\left\Vert u_{\Gamma}-u\right\Vert <\delta_{XgJ}(\varepsilon_{k})$
for each $u\in f(\Gamma)$, and so $|X(g(u_{\Gamma}))(|J(u_{\Gamma})|)-X(g(u))(|J(u)|)|\leq\varepsilon_{k}$
for each $u\in f(\Gamma)$. Combining with inequality \ref{eq:temp-76},
we obtain
\[
|\int_{\Gamma}X(x)dx-\int_{f(\Gamma)}X(g(u))|J(u)|du|
\]
\[
\leq\beta n(2^{n-1}+1)\varepsilon_{k}\mu(\Gamma)+\varepsilon_{k}\mu(\Gamma)+\varepsilon_{k}\mu(f(\Gamma)).
\]
Summation over all $\Gamma\in\Phi_{k}$ yields
\[
|\int_{M_{k}}X(x)dx-\int_{f(M_{k})}X(g(u))|J(u)|du|
\]
\[
\leq\beta n(2^{n-1}+1)\varepsilon_{k}\mu(M_{k})+\varepsilon_{k}\mu(M_{k})+\varepsilon_{k}\mu(f(M_{k}))
\]
\begin{equation}
\leq\beta n(2^{n-1}+1)\varepsilon_{k}\mu(H_{a_{0}})+\varepsilon_{k}\mu(H_{a_{0}})+\varepsilon_{k}\mu(K_{\rho}).\label{eq:temp-72}
\end{equation}
On the other hand, recalling that $H\subset M_{k}\subset H_{a_{k-1}}$
a.e. and that $|X|\leq\beta$, we have
\[
|\int_{M_{k}}X(x)dx-\int_{H}X(x)dx|
\]
\[
\leq\beta(\mu(M_{k})-\mu(H))\leq\beta(\mu(H_{a_{k-1}})-\mu(H))\leq\beta\varepsilon_{k-1}.
\]
Hence inequality \ref{eq:temp-72} yields
\[
|\int_{H}X(x)dx-\int_{f(M_{k})}X(g(u))|J(u)|du|
\]
\[
\leq\beta n(2^{n-1}+1)\varepsilon_{k}\mu(H_{a_{0}})+\varepsilon_{k}\mu(H_{a_{0}})+\varepsilon_{k}\mu(K_{\rho})+\beta\varepsilon_{k-1}.
\]
Since $\varepsilon_{k}\rightarrow0$, it follows that 
\begin{equation}
\int1_{f(M_{k})}(u)X(g(u))|J(u)|du\rightarrow\int1_{H}(x)X(x)dx\label{eq:temp-98}
\end{equation}
as . At the same time, we have $M_{k}\subset M_{k-1}$, and so $1_{f(M_{k})}\leq1_{f(M_{k-1})}$
for each $k\geq1$. Recall the assumption that $X\geq0$. The sequence
$(1_{f(M_{k})}X(g)|J|)_{k=1,2,\cdots}$ is nonincreasing. By the Monotone
Convergence Theorem, expression \ref{eq:temp-98} implies that $Z\equiv\lim_{k\rightarrow\infty}1_{f(M_{k})}X(g)|J|$
is an integrable function, with 
\[
\int Z(u)du=\int1_{H}X(x)dx.
\]

Now let $U\in C(R^{n})$ be such that $U=1$ on $H_{a_{0}}$. For
each $u\in f(M_{k})$ we then have $g(u)\in M_{k}\subset H_{a_{0}}$
and so $U(g(u))=1$. Thus $U(g)=1$ on $f(M_{k})$ for each $k\geq1$.
Define a $\mathrm{measurable}$ function $V$ on $R^{n}$ by $V\equiv|J|^{-1}1_{B}+1_{B^{c}}.$
By the arguments in the previous paragraphs, $Y\equiv\lim_{k\rightarrow\infty}1_{f(M_{k})}U(g)|J|$
is an integrable function. Therefore $VY\equiv\lim_{k\rightarrow\infty}V1_{f(M_{k})}U(g)|J|=\lim_{k\rightarrow\infty}1_{f(M_{k})}$
is a $\mathrm{measurable}$ function. We have seen earlier that $f(H)=\bigcap_{k=1}^{\infty}f(M_{k})$
a.e. Combining, we have $1_{f(H)}=VY$ a.e. and so $1_{f(H)}$ is
a $\mathrm{measurable}$ function, with $1_{f(H)}=\lim_{k\rightarrow\infty}1_{f(M_{k})}$
a.e. Since $1_{f(H)}\leq1_{f(M_{0})}$, it follows that $1_{f(H)}$
is an integrable function. In other words, $f(H)$ is an integrable
set. Moreover, $Z=1_{f(H)}X(g)|J|$. Thus

\[
\int_{f(H)}X(g(u))|J(u)|du=\int Z(u)du=\int_{H}X(x)dx.
\]
The lemma is proved.
\end{proof}
\begin{lem}
\label{Lem. Change of Integration Variables, continuous integrands}
\textbf{\emph{(Change of integration variables for continuous integrands).}}
Let $A$ and $B$ be \textup{\emph{measurable}}\emph{ }open subsets
of $R^{n}$. Suppose the following four conditions hold.

1. $g:B\rightarrow A$ is a function which is differentiable on $B$,
such that for each compact subset $K\Subset B$ we have $|J|\equiv|\det G|\geq c$
on $K$ for some $c>0$, where $G$ is the derivative of $g$. 

2. There exists a function $f:A\rightarrow B$ which is the inverse
of $g$.

3. For each compact subset $H$ of $A$ with $H\Subset A$ we have
$f(H)\Subset B$. 

4. There exists a sequence $(H_{k})_{k=1,2,\cdots}$ of compact integrable
subsets of $A$ such that (i) $H_{k}\subset H_{k+1}\Subset A$ for
each $k\geq1$, (ii) $A=\bigcup_{k=1}^{\infty}H_{k}$ and $1_{H_{k}}\uparrow1_{A}$
in measure.

Then, for each $X\in C(R^{n})$, the function $1_{B}X(g)|J|$ is integrable,
and 
\begin{equation}
\int\cdots\int_{A}X(x)dx_{1}\cdots dx_{n}=\int\cdots\int_{B}X(g(u))|J(u)|du_{1}\cdots du_{n},\label{eq:temp-107}
\end{equation}
 where t\emph{he functions $g$ and $J$ have been extended to the
set $B^{c}$ by setting $g=0=J$ on $B^{c}$. }
\end{lem}
\begin{proof}
For abbreviation, we will write $\int\cdot dx$ for $\int\cdots\int\cdot dx_{1}\cdots dx_{n}$,
and similarly $\int\cdot du$ for $\int\cdots\int\cdot du_{1}\cdots du_{n}$.

Since $H_{1}\Subset A$, there exists $a_{1}>0$ such that $H'_{1}\equiv(H_{1})_{a_{1}}\Subset A$.
Inductively, suppose $H'_{1},\cdots,H'_{k}$ have been constructed
with $H'_{1}\subset\cdots\subset H'_{k}\Subset A$ for some $k\geq1$.
Define $M_{k+1}\equiv(H'_{k}\cup H_{k+1})^{-},$ the closure of $H'_{k}\cup H_{k+1}$.
Then $M_{k+1}$ is compact. Since $H'_{k}\Subset A$ and $H_{k+1}\Subset A$,
we also have $M_{k+1}\Subset A$. Hence $H'_{k+1}\equiv(M_{k+1})_{a_{k+1}}\Subset A$
for some $a_{k+1}>0$. We have constructed inductively the sequence
\emph{$(H'_{k})_{k=1,2,\cdots}$ }such that \emph{$H'_{k}\subset H'_{k+1}\Subset A$
}for each $k\geq1$. Moreover, for each $k\geq1$ we have \emph{$H'_{k}\equiv(M_{k})_{a_{k}}$}
for some compact set $M_{k}$ and for some $a_{k}$ which, by convention,
is chosen to be a regular point of $d(\cdot,M_{k})$. It follows that
\emph{$\mu((H'_{k})_{a})\downarrow\mu(H'_{k})$} as $a\downarrow0$,
for each $k\geq1$. At the same time, since $H_{k}\subset H'_{k}\subset A$
for each $k\geq1$, we have $A=\bigcup_{k=1}^{\infty}H_{k}\subset\bigcup_{k=1}^{\infty}H'_{k}\subset A$
and $1_{H'_{k}}\uparrow1_{A}$ in measure. Thus the sequence $(H'_{k})_{k=1,2,\cdots}$
satisfies conditions 4(i) and 4(ii) in the hypothesis. Therefore,
we can replace the given sequence $(H_{k})_{k=1,2,\cdots}$ with $(H'_{k})_{k=1,2,\cdots}$
and assume that, in addition to conditions 4(i) and 4(ii), we have
\[
\mu((H{}_{k})_{a})\downarrow\mu(H{}_{k})\quad\mbox{ as }a\downarrow0
\]
 for each $k\geq1$.

Let $X\in C(R^{n})$ be arbitrary. First suppose $X\geq0$. For each
$k\geq1$, the conditions in the hypothesis of Lemma \ref{Lem. Special case Change of Integration Variables}
are satisfied by $H_{k}$. Accordingly, $f(H_{k})$ is integrable,
and
\begin{equation}
\int1_{H_{k}}(x)X(x)dx=\int1_{f(H_{k})}(u)X(g(u))|J(u)|du\label{eq:temp-105}
\end{equation}
for each $k\geq1$. Since \emph{$1_{H_{k}}\uparrow1_{A}$ }in measure,
we have \emph{$1_{H_{k}}X\uparrow1_{A}X$ }in measure. Since $X$
is integrable, the Dominated Convergence Theorem implies that 
\begin{equation}
\int1_{H_{k}}(x)X(x)dx\uparrow\int1_{A}(x)X(x)dx.\label{eq:temp-106}
\end{equation}
It follows from equalities \ref{eq:temp-105} and \ref{eq:temp-106}
that $\int1_{f(H_{k})}(u)X(g(u))|J(u)|du$ also converges to $\alpha\equiv1_{A}(x)X(x)dx$.
Hence, according to the Monotone Convergence Theorem, the function
$Z\equiv\lim_{k\rightarrow\infty}1_{f(H_{k})}X(g)|J|$ is integrable,
with integral $\alpha$. On the other hand, because $B=f(A)=\bigcup_{k=1}^{\infty}f(H_{k})$,
we have $1_{f(H_{k})}\uparrow1_{B}$ a.e. and so $1_{f(H_{k})}X(g)|J|\uparrow1_{B}X(g)|J|$
a.e. Combining, we obtain $1_{B}X(g)|J|=Z$ a.e. Accordingly $1_{B}X(g)|J|$
is integrable, with integral $\alpha$. The desired equality \ref{eq:temp-107}
is thus proved for the case of a non-negative$X\in C(R)$. The general
case of $X\in C(R)$ follows from $X=2X_{+}-|X|$ and from linearity.
\end{proof}
The following proposition generalizes the change of integration variables
for continuous integrand $X$ in Lemma \ref{Lem. Change of Integration Variables, continuous integrands}
to allow for a general integrable $X$. 
\begin{thm}
\label{Theorem. Change of Variables for integrable integrands} \textbf{\emph{(Change
of integration variables).}} Let $A$ and $B$ be \textup{\emph{measurable}}\emph{
}open subsets of $R^{n}$. Suppose the following four conditions hold.

\emph{1. $g:B\rightarrow A$ }is a function which is differentiable
on $B$, such that for each compact subset $K\Subset B$ we have $|J|\equiv|\det G|\geq c$
on $K$ for some $c>0$, where $G$ is the derivative of $g$. 

\emph{2.} There exists a function $f:A\rightarrow B$ which is the
inverse of $g$.

\emph{3}. For each compact subset $H$ of $A$ with $H\Subset A$
we have $f(H)\Subset B$. 

\emph{4}. There exists a sequence $(H_{k})_{k=1,2,\cdots}$ of compact
integrable subsets of $A$ such that (i) $H_{k}\subset H_{k+1}\Subset A$
for each $k\geq1$, (ii) $A=\bigcup_{k=1}^{\infty}H_{k}$ and $1_{H_{k}}\uparrow1_{A}$
in measure.

Then, for each Lebesgue integrable function $X$ on $R^{n}$ such
that the function $1_{B}X(g)|J|$ is integrable, we have
\begin{equation}
\int\cdots\int_{A}X(x)dx_{1}\cdots dx_{n}=\int\cdots\int_{B}X(g(u))|J(u)|du_{1}\cdots du_{n},\label{eq:temp-107-2}
\end{equation}
 where the functions \emph{$g$ and $J$ have been extended to the
set $B^{c}$ by setting $g=0=J$ on $B^{c}$. }
\end{thm}
\begin{proof}
By hypothesis, the functions $g$ and $J$ are defined on some full
subset $D_{0}$ of $R^{n}$. Let $D\equiv D_{0}\cap(B\cup B^{c})$.
Then $D$ is a full set. Let $X$ be an arbitrary integrable function
on $R^{n}$. By Proposition\ref{Prop. Completion of C(S1xS2) =00003D compC(S1) x CompC(S2)}
there exists a sequence $(X_{j})_{j=1,2,\cdots}$ in $C(R^{n})$ which
is a representation of $X$ relative to the Lebesgue integration on
$R^{n}$. In other words, (i) $\sum_{j=1}^{\infty}\int|X_{j}(x)|dx<\infty$,
and (ii) if $x\in R^{n}$ is such that $\sum_{j=1}^{\infty}|X_{j}(x)|<\infty$
then $x\in domain(X)$ and $X(x)=\sum_{j=1}^{\infty}X_{j}(x)$. By
hypothesis, for each $j\geq1$, the function \emph{$1_{B}X_{j}(g)|J|$
}is integrable, and equality \ref{eq:temp-107} in Lemma \ref{Lem. Change of Integration Variables, continuous integrands}
holds for each $|X_{j}|$. Hence 
\begin{equation}
\sum_{j=1}^{\infty}\int1_{B}(u)|X_{j}(g(u))J(u)|du=\sum_{j=1}^{\infty}\int_{A}|X_{j}(x)|dx<\infty.\label{eq:temp-109}
\end{equation}
Define a function $Y$ on $R^{n}$ by
\[
domain(Y)\equiv\{u\in D:\sum_{j=1}^{\infty}1_{B}(u)|X_{j}(g(u))J(u)|<\infty\}
\]
and $Y(u)\equiv\sum_{j=1}^{\infty}1_{B}(u)X_{j}(g(u))|J(u)|$ for
each $u\in domain(Y)$. Then, in view of inequality \ref{eq:temp-109},
the function $Y$ is integrable, the sequence $(1_{B}X_{j}(g)|J|)_{j=1,2,\cdots}$
being a representation. Moreover 
\[
\int Y(u)du=\sum_{j=1}^{\infty}\int1_{B}(u)X_{j}(g(u))|J(u)|du=\sum_{j=1}^{\infty}\int_{A}X_{j}(x)dx=\int_{A}X(x)dx.
\]

Consider an arbitrary $u\in DB$. Then
\[
\sum_{j=1}^{\infty}|X_{j}(g(u))J(u)|=\sum_{j=1}^{\infty}1_{B}(u)|X_{j}(g(u))J(u)|<\infty.
\]
Since  by hypothesis, we see that $\sum_{j=1}^{\infty}|X_{j}(g(u))|<\infty$.
It follows from condition (ii) above that $g(u)\in domain(X)$ and
$X(g(u))=\sum_{j=1}^{\infty}X_{j}(g(u))$. Consequently 
\[
1_{B}(u)X(g(u))|J(u)|=\sum_{j=1}^{\infty}1_{B}(u)X_{j}(g(u))|J(u)|=Y(u).
\]
We conclude that $DB\subset domain(X(g))$ and that $1_{B}X(g)|J|=Y$
on $DB$. On the other hand, for each $u\in DB^{c}$, we have . Combining,
we see that $1_{B}X(g)|J|=Y$ on the full set $D(B\cup B^{c})$. It
follows that $1_{B}X(g)|J|$ is integrable, which integral equal to
$\int Y(u)du=\int_{A}X(x)dx$. The proposition is thus proved.
\end{proof}
The next corollary is a formula for changing integration variables
from rectangular coordinates in $R^{2}$ to polar coordinates. It
is all we needed in the text. First an elementary lemma.
\begin{lem}
\label{Lem. bijction h:(0,2pi)->unit circle minus a pt} \textbf{\emph{(bijection
of an arc to an interval).}} Let $C\equiv\{(u,v)\in R^{2}:u^{2}+v^{2}=1\mbox{ and }(u,v)\neq(1,0)\}$,
where the inequality signifies a positive Euclidean distance. Then
the function $h:(0,2\pi)\rightarrow C$ defined by $h(\theta)\equiv(\cos\theta,\sin\theta)$
is a bijection. More specifically, for arbitrarily small $\theta_{0}\in(0,\frac{\pi}{2})$,
the function  is a bijection.
\end{lem}
\begin{proof}
Take arbitrary $\theta_{0},\theta_{1}\in(0,\frac{\pi}{2})$ such that
$\theta_{0}<\theta_{1}$. Define $\delta_{0}\equiv\sin\theta_{0}$
and $\delta_{1}\equiv\sin\theta_{1}$. Then $0<\delta_{0}<\delta_{1}$.
Let $(u,v)\in C_{\theta_{0}}$ be arbitrary. Then $u^{2}+v^{2}=1$
and $(u-1)^{2}+v^{2}\geq\delta_{0}^{2}\equiv\sin^{2}\theta_{0}$.
Either $v^{2}\geq\delta_{0}^{2}$ or $v^{2}\leq\delta_{1}^{2}$. In
the first case, we have in turn $v\geq\delta_{0}$ or $v\leq-\delta_{0}$.
Suppose $v\geq\delta_{0}$. Then $u^{2}\leq1-\delta_{0}^{2}=\cos^{2}\theta_{0}$,
and so $u\in[-\cos\theta_{0},\cos\theta_{0}]\subset(-1,1)$. At the
same time, the function $\cos:[\theta_{0},\pi-\theta_{0}]\rightarrow[-\cos\theta_{0},\cos\theta_{0}]$
has a strictly negative derivative. Hence there exists a unique $\theta\in[\theta_{0},\pi-\theta_{0}]\subset(0,\pi)$
such that $u=\cos\theta$ and $v=\sqrt{1-u^{2}}=\sqrt{1-\cos^{2}\theta}=\sin\theta$,
the last equality holding because $\sin\geq0$ on $(0,\pi)$. Next,
suppose $v\leq-\delta_{0}$. Then, again, $u^{2}\leq1-\delta_{0}^{2}$,
and so $u\in[-\cos\theta_{0},\cos\theta_{0}]\subset(-1,1)$. At the
same time, the function $\cos:[\pi+\theta_{0},2\pi-\theta_{0}]\rightarrow[-\cos\theta_{0},\cos\theta_{0}]$
has a strictly positive derivative. Hence there exists a unique $\theta\in[\pi+\theta_{0},2\pi-\theta_{0}]\subset(\pi,2\pi)$
such that $u=\cos\theta$ and $v=-\sqrt{1-u^{2}}=-\sqrt{1-\cos^{2}\theta}=\sin\theta$,
the last equality holding because $\sin\leq0$ on $(\pi,2\pi)$. Now
consider the second case, where $v^{2}\leq\delta_{1}^{2}\equiv\sin^{2}\theta_{1}$.
Thus $v\in[-\sin\theta_{1},\sin\theta_{1}]\subset(-1,1)$. At the
same time, the function $\sin:[\pi-\theta_{1},\pi+\theta_{1}]\rightarrow[-\sin\theta_{1},\sin\theta_{1}]$
has a strictly negative derivative. Hence there exists a unique $\theta\in[\pi-\theta_{1},\pi+\theta_{1}]\subset(\frac{\pi}{2},\frac{3\pi}{2})$
such that $v=\sin\theta$. Then $u=-\sqrt{1-v^{2}}=-\sqrt{1-\sin\theta^{2}}=\cos\theta$,
the last equality holding because $\cos\leq0$ on $(\frac{\pi}{2},\frac{3\pi}{2})$.
Summing up, we see that for each $(u,v)\in C_{\theta_{0}}$ there
exists $\theta\in[\theta_{0},2\pi-\theta_{0}]$ such that $(u,v)=(\cos\theta,\sin\theta)\equiv h(\theta)$.
Since each $(u,v)\in C$ belongs to $C_{\theta_{0}}$ for sufficiently
small $\theta_{0}\in(0,\frac{\pi}{4})$, we have $(u,v)=h(\theta)$
for some $\theta\in(0,2\pi)$.

Suppose $\theta,\theta'\in(0,2\pi)$ are such that $(\cos\theta,\sin\theta)=(\cos\theta',\sin\theta')$.
Suppose $\theta\neq\theta'$. There are three possibilities: (i) $\theta\in(0,\pi)$,
(ii) $\theta\in(\frac{\pi}{2},\frac{3\pi}{2})$, and (iii) $\theta\in(\pi,2\pi)$.
Consider case (i). Then $\theta'\in[\pi,2\pi)$ and so $\sin\theta'\leq0<\sin\theta$,
a contradiction. Similarly case (iii) would lead to a contradiction.
Now consider case (ii). Then $\theta'\in(0,\frac{\pi}{2}]\cup[\frac{3\pi}{2},2\pi)$
and we have $\cos\theta'\geq0>\cos\theta$, again a contradiction.
Hence $\theta=\theta'$. Thus we conclude that $h$ is one-to-one. 
\end{proof}
\begin{cor}
\label{Cor. Integrating in Polar Coordinates} \textbf{\emph{(Change
of integration variables from rectangular coordinates to polar coordinates).}}
For each integrable function $X$ on $R^{2}$, the function 
\[
X(r\cos\theta,r\sin\theta)r
\]
is integrable on $R^{2}$, with 
\[
\int_{-\infty}^{\infty}\int_{-\infty}^{\infty}X(x,y)dxdy=\int_{0}^{2\pi}\int_{0}^{\infty}X(r\cos\theta,r\sin\theta)rdrd\theta.
\]
\end{cor}
\begin{proof}
Let $B\equiv(0,\infty)\times(0,2\pi)$. Define$g:B\rightarrow R^{2}$
by $g(r,\theta)\equiv(r\cos\theta,r\sin\theta)$ for each $(r,\theta)\in B$.
Then $g$ is differentiable on the open set $B$, with derivative
$G\equiv\left[\begin{array}{cc}
\cos\theta & \sin\theta\\
-r\sin\theta & r\cos\theta
\end{array}\right]$ and Jacobian $\det G=r$. Let $R_{+}\equiv[0,\infty)\times\{0\}\subset R^{2}$.
Define $A\equiv\{(x,y)\in R^{2}:d((x,y),R_{+})>0\}$, the metric complement
of $R_{+}$. It is easily verified that $A$ is an open set. We proceed
to verify the conditions in the hypothesis of Lemma  \ref{Lem. Change of Integration Variables, continuous integrands}.

We will first prove that $g(B)\subset A$. Let $\delta\in(0,\frac{\pi}{2})$
and $0<a<b$ be arbitrary. Consider any $(r,\theta)\in[a,b]\times[\delta,2\pi-\delta]$.
Let $(x,y)\equiv g(r,\theta)$. Consider any $(z,0)\in R^{+}$. Write
$\alpha\equiv d((x,y),(z,0))$. Then $\alpha^{2}=\left\Vert (r\cos\theta,r\sin\theta)-(z,0)\right\Vert ^{2}=r^{2}-2rz\cos\theta+z^{2}$.
If $\theta\in[\frac{\pi}{2},\frac{3\pi}{2}]$ then $\cos\theta\leq0$
and so $\alpha^{2}\geq r^{2}\geq a^{2}$. If $\theta\in[\delta,\frac{\pi}{2})\cup(\frac{3\pi}{2},2\pi-\delta]$
then $\alpha\geq|y|=r|\sin\theta|\geq a\sin\delta$. Hence $\alpha\geq a\sin\delta$
if $\theta\in[\delta,\frac{\pi}{2})\cup[\frac{\pi}{2},\frac{3\pi}{2}]\cup(\frac{3\pi}{2},2\pi-\delta]$.
By continuity, we therefore have $\alpha\geq a\sin\delta$ for the
given $(r,\theta)\in[a,b]\times[\delta,2\pi-\delta]$. Since $(z,0)\in R^{+}$
is arbitrary, we have $d((x,y),R_{+})\geq a\sin\delta>0$ and so $g(r,\theta)\in A$.
It follows that
\begin{equation}
g([a,b]\times[\delta,2\pi-\delta])\subset\{(x,y)\in R^{2}:d((x,y),R_{+})\geq a\sin\delta\}\subset A.\label{eq:temp-111}
\end{equation}
Since $a>0$ and $\delta>0$ are arbitrarily small, and $b$ is arbitrarily
large, it follows that $g(B)\equiv g((0,\infty)\times(0,2\pi))\subset A$.

Next, we will show that the function $g:B\rightarrow A$ has an inverse.
Let $C\equiv\{(u,v)\in R^{2}:u^{2}+v^{2}=1\mbox{ and }(u,v)\neq(1,0)\}$.
Let $h:(0,2\pi)\rightarrow C$ be the bijection defined by $h(\theta)\equiv(\cos\theta,\sin\theta)$
as in Lemma \ref{Lem. bijction h:(0,2pi)->unit circle minus a pt},
and let $h^{-1}:C\rightarrow(0,2\pi)$ denote its inverse. Let $(x,y)\in A$
be arbitrary. Then $r\equiv\sqrt{x^{2}+y^{2}}>0$. Define $u\equiv\frac{x}{r}$
and $v\equiv\frac{y}{r}$. Then $d((u,v),(1,0))=\frac{1}{r}d((x,y),(r,0))\geq\frac{1}{r}d((x,y),R_{+})>0$.
Hence $(u,v)\in C$. Define $\theta\equiv h^{-1}(u,v)\in(0,2\pi)$
and define $f(x,y)\equiv(r,\theta)\in(0,\infty)\times(0,2\pi)\equiv B$.
Then $(u,v)=h(\theta)\equiv(\cos\theta,\sin\theta)$, and so $g(f(x,y))\equiv g(r,\theta)=(ru,rv)=(x,y)$.
In other words, $g(f)$ is the identity function on $A$, and $f:A\rightarrow B$
is the inverse of $g$. 

Let $(a_{k})$,$(b_{k})$,$(\theta_{k})$ be strictly monotone sequences
in $R$ such that $0<a_{k}<b_{k}$, $0<\theta_{k}<\frac{\pi}{4}$
for each $k\geq1$, and such that $a_{k}\downarrow0$, $b_{k}\uparrow\infty$,
$\theta_{k}\downarrow0$. For each $k\geq1$ define $K_{k}\equiv[a_{k},b_{k}]\times[\theta_{k},2\pi-\theta_{k}]\Subset B$
and $H_{k}\equiv g(K_{k})$. Being the product of two closed intervals,
$K_{k}$ is compact and integrable for each $k\geq1$. Moreover $B=\bigcup_{k=1}^{\infty}K_{k}$
and so $A\equiv g(B)=\bigcup_{k=1}^{\infty}g(K_{k})\equiv\bigcup_{k=1}^{\infty}H_{k}$. 

Next let $H$ be an arbitrary compact subset with $H\Subset A$. Since
$H\Subset A$, there exists $a>0$ such that $H_{a}\Subset A$. We
will show that \emph{$f(H)\Subset B$} by proving that\emph{ }$f(H)\subset K_{k}$
for sufficiently large $k$. Let $j\geq1$ be so large that $a_{k}<a$
and $H\subset(d(\cdot,(0,0))\leq b_{k})$ for each $k\geq j$. Consider
an arbitrary $(x,y)\in H$. Let $(r,\theta)\equiv f(x,y)$ and $(u,v)\equiv(\frac{x}{r},\frac{y}{r})$.
Then $r\leq b_{j}$. Moreover, the assumption that $r<a$ would imply
that $d((0,0),H)\leq d((0,0),(x,y))<a$, and so $(0,0)\in H_{a}\subset A$,
whence $0=d((0,0),R_{+})>0$, a contradiction. Therefore $r\geq a\geq a_{j}$.
Next, we have $(u,v)=(\cos\theta,\sin\theta)$. Let $k\geq j$ be
so large that $\sin^{2}\theta_{k}\leq b_{j}^{2}a^{2}$. Then
\begin{equation}
(u-1)^{2}+v^{2}=r^{-2}d((x,y),(r,0))^{2}\geq r^{-2}a^{2}\geq b_{j}^{2}a^{2}\geq\sin^{2}\theta_{k}.\label{eq:temp-110}
\end{equation}
Hence $(u,v)\in C_{\theta_{k}}\equiv\{(u,v)\in C:(u-1)^{2}+v^{2}\geq\sin^{2}\theta_{k}\}$.
By Lemma \ref{Lem. bijction h:(0,2pi)->unit circle minus a pt}, the
function $h:[\theta_{k},2\pi-\theta_{k}]\rightarrow C_{\theta_{k}}$
is a bijection. Therefore we have $\theta=h^{-1}(u,v)\in[\theta_{k},2\pi-\theta_{k}]$.
Accordingly
\[
f(x,y)\equiv(r,\theta)\in[a_{j},b_{j}]\times[\theta_{k},2\pi-\theta_{k}]\subset[a_{k},b_{k}]\times[\theta_{k},2\pi-\theta_{k}]\equiv K_{k}.
\]
Since $(x,y)\in H$ is arbitrary, it follows that $f(H)\subset K_{k}$.
On the other hand $(K_{k})_{b}\subset K_{k+1}$ where $b=(a_{k}-a_{k+1})\wedge(b_{k+1}-b_{k})\wedge(\theta_{k}-\theta_{k+1})$.
Hence $K_{k}\Subset K_{k+1}\subset B$. This proves that $f(H)\Subset B$.

We next prove that $H_{k}\Subset A$ for each $k\geq1$. Let $k\geq1$
and $a\in(0,a_{k}\sin\theta_{k})$ be arbitrary. From expression \ref{eq:temp-111}
we see that $d((x,y),R_{+})\geq a_{k}\sin\theta_{k}$ for each $(x,y)\in H_{k}\equiv g(K_{k})$.
Thus $d(H_{k},R_{+})\geq a_{k}\sin\theta_{k}$. Therefore $d((H_{k})_{a},R_{+})\geq a_{k}\sin\theta_{k}-a>0$.
Hence $(H_{k})_{a}\subset A$It follows that $H_{k}\Subset A$.

We will now prove that $1_{H_{k}}\uparrow1_{A}$ in measure. To that
end, let $M$ be any integrable subset of $R^{2}$. Let $\varepsilon>0$
be arbitrary. Then there exists $q\geq1$ so large that $\mu(MB_{q}^{c})<\varepsilon$,
where $B_{q}\equiv[-q,q]^{2}$. Let $k\geq1$ be so large that $4qa_{k}\sin\theta_{k}<\varepsilon$.
From the previous paragraph we see that $H_{k}\subset(d(\cdot,R_{+})\geq a_{k}\sin\theta_{k})$.
Therefore $\mu(B_{q}H_{k}^{c})\leq\mu([-q,q]\times[-a_{k}\sin\theta_{k},a_{k}\sin\theta_{k}])=4qa_{k}\sin\theta_{k}<\varepsilon$.
Hence
\[
\mu((A-H_{k})M)\leq\mu(MH_{k}^{c})\leq\mu(MB_{q}^{c})+\mu(B_{q}H_{k}^{c})\leq\varepsilon+\varepsilon=2\varepsilon.
\]
Consequently $\mu(|1_{A}-1_{H_{k}}|>\varepsilon)M<2\varepsilon$.
Since $\varepsilon>0$ and the integrable set $M$ are arbitrary,
we have $1_{H_{k}}\uparrow1_{A}$ in measure.

All the conditions in the hypothesis of Lemma  \ref{Lem. Change of Integration Variables, continuous integrands}
have thus been verified. Accordingly,\emph{ }for each $X\in C(R^{2})$,
the function $1_{B}X(g)|J|$ is integrable, with 
\begin{equation}
\int\int_{A}X(x,y)dxdy=\int\int_{B}X(r\cos\theta,r\sin\theta)rdrd\theta.\label{eq:temp-107-1}
\end{equation}
Since $A$ and $B$ are full sets of $R^{2}$ and $[0,\infty)\times[0,2\pi]$
respectively, equality \ref{eq:temp-107-1} yields
\begin{equation}
\int_{-\infty}^{\infty}\int_{-\infty}^{\infty}X(x,y)dxdy=\int_{0}^{2\pi}\int_{0}^{\infty}X(r\cos\theta,r\sin\theta)rdrd\theta\label{eq:temp-108}
\end{equation}
for each $X\in C(R^{2})$. By Theorem \ref{Theorem. Change of Variables for integrable integrands},
equality \ref{eq:temp-108} therefore holds also for each integrable
function $X$ on $H^{n}$.
\end{proof}
\begin{cor}
\label{Cor. If X,Y integrable, then  X(u-v)Y(v) integrable} \textbf{\emph{(Integrability
of convolution of two integrable functions). }}Let $X,Y$ be integrable
functions on $R^{n}$. Define the function $Z:R^{2n}\rightarrow R$
by $Z(u,v)\equiv X(u-v)Y(v)$ where 
\[
domain(Z)\equiv\{(u,v)\in R^{2n}:u-v\in domain(X)\mbox{ and }v\in domain(Y)\}.
\]
Then $Z$ is an integrable function on $R^{2n}$. Hence $X\star Y\equiv\int X(\cdot-v)Y(v)dv$
is an integrable function on $R^{n}$ by Fubini's Theorem.
\end{cor}
\begin{proof}
Define the function $W:R^{2n}\rightarrow R$ by $W(x,y)\equiv X(x)Y(y)$
for each $(x,y)$ in the full set $domain(W)\equiv domain(X)\times domain(Y)$.
By Proposition \ref{Prop. Cartesian product of integrable funcs is integrable},
the function $W$ is integrable on $R^{2n}$. Now define the function
$g:R^{2n}\rightarrow R^{2n}$ by $g(u,v)\equiv(u-v,v)$. In particular
$g$ has the Jacobian $J\equiv1$ on $R^{2n}$ and has an inverse
function $f$ defined by $f(x,y)=(x+y,y)$. The conditions in Lemma
 \ref{Lem. Change of Integration Variables, continuous integrands}
are routinely verified. Accordingly $W(g)\equiv W\circ g$ is an integrable
function on $R^{2n}$, where
\[
domain(W(g))\equiv\{(u,v)\in R^{2}:g(u,v)\in domain(W)\}
\]
\[
=\{(u,v)\in R^{2}:(u-v,v)\in domain(X)\times domain(Y)\}=domain(Z)
\]
Moreover, for each $(u,v)\in domain(W(g))$, we have $W(g(u,v))\equiv W(u-v,v)\equiv X(u-v)Y(v)\equiv Z(u,v)$.
Thus $Z=W(g)$ and so $Z$ is integrable. The corollary is proved.
\end{proof}

\chapter{Taylor's Theorem}

We present the proof of Taylor's Theorem from \cite{BishopBridges85},
which is then extended to higher dimension in the next corollary. 
\begin{thm}
\textbf{\emph{(Taylor's Theorem)}} \label{Thm. Tayor's Theorem} Let
$D$ be a non-empty open interval in $R$. Let $f$ be a complex-valued
function on $D$. Let $n\geq0$ be arbitrary. Suppose $f$ has continuous
derivatives up to order $n$ on $D$. For $k=1,\cdots,n$ write $f^{(k)}$
for the $k$-th derivative of $f$. Let $t_{0}\in D$ be arbitrary,
and define 
\[
r_{n}(t)\equiv f(t)-\sum_{k=0}^{n}f^{(k)}(t_{0})(t-t_{0})^{k}/k!
\]
for each $t\in D$. Then the following holds.
\end{thm}
\begin{enumerate}
\item \emph{If $|f^{(n)}(t)-f^{(n)}(t_{0})|\leq M$ on $D$ for some $M>0$,
then $|r_{n}(t)|\leq M|t-t_{0}|^{n}/n!$ }
\item \emph{$r_{n}(t)=o(|t-t_{0}|^{n})$ as $t\rightarrow t_{0}$. More
precisely, suppose $\delta_{f,n}$ is a modulus of continuity of $f^{(n)}$
at the point $t_{0}$. Let $\varepsilon>0$ be arbitrary. Then $|r_{n}(t)|<\varepsilon|t-t_{0}|^{n}$
for each $t\in R$ with $|t-t_{0}|<\delta_{f,n}(n!\varepsilon)$.}
\item \emph{If $f^{(n+1)}$ exists on $D$ and $|f^{(n+1)}|\leq M$ for
some $M>0$, then $|r_{n}(t)|\leq M|t-t_{0}|^{n+1}/(n+1)!$.}
\end{enumerate}
\begin{proof}
1. Let $t\in D$ be arbitrary. Integrating repeatedly, we obtain 
\[
|r_{n}(t)|=|\int_{t_{0}}^{t}\int_{t_{0}}^{t_{1}}\cdots\int_{t_{0}}^{t_{n-1}}(f^{(n)}(t_{n})-f^{(n)}(t_{0}))dt_{n}\cdots dt_{2}dt_{1}|
\]
\[
\leq|\int_{t_{0}}^{t}\int_{t_{0}}^{t_{1}}\cdots\int_{t_{0}}^{t_{n-1}}Mdt_{n}\cdots dt_{2}dt_{1}|=M|t-t_{0}|^{n}/n!,
\]
where we note that all the integration variables $t_{1},\cdots,t_{n}$
are in the interval $D$.

2. Let $\varepsilon>0$ be arbitrary.\emph{ }Consider each \emph{$t\in R$
with $|t-t_{0}|<\delta_{f,n}(n!\varepsilon)$. }Then $|f^{(n)}(t)-f^{(n)}(t_{0})|\leq n!\varepsilon$
for each $t\in D\equiv(t_{0}-\delta_{f,n}(n!\varepsilon),t_{0}+\delta_{f,n}(n!\varepsilon)).$
By Assertion 1, we then have $|r_{n}(t)|\leq\varepsilon|t-t_{0}|^{n}$
for $t\in D$.

3. Integrating repeatedly, we obtain
\[
r_{n}(t)=|\int_{t_{0}}^{t}\int_{t_{0}}^{t_{1}}\cdots\int_{t_{0}}^{t_{n}}f^{(n+1)}(t_{n+1})dt_{n+1}\cdots dt_{2}dt_{1}|
\]
\[
\leq|\int_{t_{0}}^{t}\int_{t_{0}}^{t_{1}}\cdots\int_{t_{0}}^{t_{n}}Mdt_{n+1}\cdots dt_{2}dt_{1}|=M|t-t_{0}|^{n+1}/(n+1)!.
\]
\end{proof}
\begin{cor}
\label{Cor. Tayor's Theorem in dimension m, for twice coninuoutly diff functon}
(\textbf{\emph{Taylor's Theorem for functions with $m$ variables
which have continuous partial derivatives up to second order).}} For
ease of notations, we state and prove the corollary only for degree
$n=2$; in other words, only for a twice continuously differentiable,
real-valued, function $f$. Let $D$ be a non-empty convex open subset
in $R^{m}$. Let $f$ be a real-valued function on $D$. Suppose $f$
has continuous partial derivatives up to order $2$ on $D$. In other
words, the partial derivatives
\[
f_{i}\equiv\frac{\partial f}{\partial x_{i}}
\]
and
\[
f_{i,j}\equiv\frac{\partial^{2}f}{\partial x_{i}\partial x_{j}}
\]
are continuous functions on $D$, for arbitrary $i,j$ in $\{1,\cdots,m\}$.
Let $\delta_{f,2}$ be a common modulus of the these second order
partial derivatives $f_{i,j}$. Let $\overline{x}\equiv(x_{1},\cdots,x_{m})\in D$
and $\overline{y}\equiv(y_{1},\cdots,y_{m})\in D$ be arbitrary, and
define 
\[
r_{f,2}(\overline{y})\equiv f(\overline{y})-\{f(\overline{x})+\sum_{i=1}^{m}f_{i}(\overline{x})(y_{i}-x_{i}))+\frac{1}{2}\sum_{j=1}^{m}\sum_{i=1}^{m}f_{j,i}(\overline{x})(y_{j}-x_{j})(y_{i}-x_{i})\}
\]
Let $\varepsilon>0$. Suppose 
\[
|\overline{y}-\overline{x}|<\delta_{f,2}(m^{-2}\varepsilon).
\]
Then
\[
|r_{f,2}(\overline{y})|\leq2^{-1}\varepsilon|\overline{y}-\overline{x}|^{2}.
\]
\end{cor}
\begin{proof}
Define the function 
\[
g(t)\equiv f(t(\overline{y}-\overline{x})+\overline{x})
\]
for each $t\in[0,1]$. Then the Chain Rule yields
\[
g'(t)=\sum_{i=1}^{m}f_{i}(t(\overline{y}-\overline{x})+\overline{x})(y_{i}-x_{i})
\]
and
\[
g''(t)=\sum_{j=1}^{m}\sum_{i=1}^{m}f_{j,i}(t(\overline{y}-\overline{x})+\overline{x})(y_{j}-x_{j})(y_{i}-x_{i}).
\]
By hypothesis, the partial derivatives $f_{j,i}$ have a common modulus
of continuity $\delta_{f,2}$ at the point $\overline{x}$, for each
$i,j=1,\cdots,m$. Let $\varepsilon>0$ be arbitrary. 
\[
|\overline{y}-\overline{x}|<\delta_{f,2}(m^{-2}\varepsilon).
\]
Hence 
\[
|g''(1)-g''(0)|=|\sum_{j=1}^{m}\sum_{i=1}^{m}(f_{j,i}(\overline{y})-f_{j,i}(\overline{x}))(y_{j}-x_{j})(y_{i}-x_{i})|
\]
\[
\leq|\sum_{j=1}^{m}\sum_{i=1}^{m}m^{-2}\varepsilon(y_{j}-x_{j})(y_{i}-x_{i})|.
\]
\[
\leq m^{-2}\varepsilon m^{2}|\overline{y}-\overline{x}|^{2}=\varepsilon|\overline{y}-\overline{x}|^{2}.
\]
Therefore
\[
|r_{f,2}(\overline{y})|\equiv|f(\overline{y})-f(\overline{x})-\sum_{i=1}^{m}f_{i}(\overline{x})(y_{i}-x_{i})-\sum_{j=1}^{m}\sum_{i=1}^{m}f_{j,i}(\overline{x})(y_{j}-x_{j})(y_{i}-x_{i})|
\]
\[
=|g(1)-g(0)-g'(0)-2^{-1}g''(0)=|\int_{s=0}^{1}(g'(s)-g'(0)-g''(0)s)ds|
\]
\[
=|\int_{s=0}^{1}\int_{u=0}^{s}(g''(u)-g''(0))duds|\leq\varepsilon|\overline{y}-\overline{x}|^{2}\int_{s=0}^{1}\int_{u=0}^{s}duds
\]
\[
=2^{-1}\varepsilon|\overline{y}-\overline{x}|^{2},
\]
as alleged. 
\end{proof}
$\,$

\printindex{}
\end{document}